%
%
%
\newdimen\rotdimen
\def\vspec#1{\special{ps:#1}}
\def\rotstart#1{\vspec{gsave currentpoint currentpoint translate
   #1 neg exch neg exch translate}}
\def\rotfinish{\vspec{currentpoint grestore moveto}}
%
%
\def\rotr#1{\rotdimen=\ht#1\advance\rotdimen by\dp#1%
   \hbox to\rotdimen{\hskip\ht#1\vbox to\wd#1{\rotstart{90 rotate}%
   \box#1\vss}\hss}\rotfinish}
%
%
\def\rotl#1{\rotdimen=\ht#1\advance\rotdimen by\dp#1%
   \hbox to\rotdimen{\vbox to\wd#1{\vskip\wd#1\rotstart{270 rotate}%
   \box#1\vss}\hss}\rotfinish}%
%
%
\def\rotu#1{\rotdimen=\ht#1\advance\rotdimen by\dp#1%
   \hbox to\wd#1{\hskip\wd#1\vbox to\rotdimen{\vskip\rotdimen
   \rotstart{-1 dup scale}\box#1\vss}\hss}\rotfinish}%
%
%
\def\rotf#1{\hbox to\wd#1{\hskip\wd#1\rotstart{-1 1 scale}%
   \box#1\hss}\rotfinish}%
%
%
\def\flip#1{\hbox{\rotstart{270 rotate}#1}\rotfinish}%
\def\flippts{0.9\colpts}
%
%
\def\twist#1{\hbox{\rotstart{315 rotate}#1}\rotfinish}%
\def\twistpts{0.6\colpts}
%
%
\def\tableopen#1{%
\tabletop%
\vbox to \hdpts{%
\vss
\tabcenterline{%
{\bf Table {\tenrm\number\tablenumber}.\enspace}%
{\sl #1}%
}%
\nobreak%
\vskip 2.5 true pt%
\hrule width\tabhsize%
}%
\immediate\write\toctxt{%
{\bf\qquad Table \the\tablenumber.\enspace}{\sl #1}%
}%
\write\tocnum{\the\pageno}%
\tablestyle
\nobreak
\nointerlineskip
}%
%
%
\def\tablecont{%
\tabletop%
\vbox to \hdpts{%
\vss
\tabcenterline{%
{\bf Table {\tenrm\number\tablenumber},\enspace}%
{\sl continued.}%
}%
\nobreak%
\vskip 2.5 true pt%
\hrule width\tabhsize%
}%
\tablestyle%
\nobreak%
\nointerlineskip
}%
%
%
%
\def\tableclose{%
\global\advance\tablenumber by 1
}%
%
%
\def\lotsofdots{%
\leaders
\hbox to 0.7em{\hss.\hss}%
\hfil
}%
%
%
\newcount\tablenumber \tablenumber=1
\magnification 1200
\hoffset -0.25 true in 
\voffset -0.4 true in 
\hsize=7.0 true in
\vsize=10.1 true in
\newwrite\toctxt
\immediate\openout\toctxt toctxt
\newwrite\tocnum
\immediate\openout\tocnum tocnum
\headline={%
\ifnum\pageno>1{%
\vbox to \headlinepts{\line{\hss\tenrm\folio\hss}\vss}%
}%
\else{%
\ifnum\pageno<-1{%
\vbox to \headlinepts{\line{\hss\tenrm\folio\hss}\vss}%
}%
\else{%
\vbox to \headlinepts{\line{\hfil}\vss}%
}%
\fi
}%
\fi
}%
\footline={}%
\pageno=0  
\mathsurround=0.5 true pt
\def\incltabskip{0.16667em}  
\newdimen\rowpts
\newdimen\colpts
\newdimen\hdpts
\hdpts=30 true pt
\newdimen\flippts
\newdimen\twistpts
\newdimen\rowlabpts
\newdimen\collabpts
\def\cg{\hskip 0pt plus 1pt minus 1pt}%
\def\olskip{\hskip 1.25 true pt}%
\def\ol#1{{\olskip\overline{#1}\olskip}}%
\def\cx#1{\hbox to\colpts{#1}}%
\def\rx{\vbox to\rowpts}%
\def\rlx#1{\hbox to\rowlabpts{#1}}%
\def\clx{\vbox to \collabpts}
\newskip\dblskip
\dblskip=24 true pt
\baselineskip=\dblskip
\newskip\tablebaseskip
\tablebaseskip=18 true pt plus 5 true pt
\def\makeheadline{\vbox to 0pt{\vskip -20 true pt\line{\vbox to \headlinepts{}\the\headline}\vss
}\nointerlineskip}%
\newdimen\headlinepts
\headlinepts=8 true pt
\newskip\tocskip
\tocskip=22 true pt
\def\tabletop{}%
\def\sectionstyle{%
\baselineskip=\dblskip
}%
\def\tablestyle{%
\def\inbaseskip{\tablebaseskip}%
\baselineskip=\tablebaseskip%
}%
\def\sectionhead#1{%
\tableftline{\sectionfont \the\sectionnumber.\quad #1}%
\immediate\write\toctxt{%
{\bf \the\sectionnumber.\enspace #1}%
}%
\write\tocnum{\the\pageno}%
\sectionstyle%
\global\advance\sectionnumber by 1
}%
\def\otherhead#1{%
\tableftline{\sectionfont #1}%
\immediate\write\toctxt{{\bf #1}}%
\nobreak
\write\tocnum{\the\pageno}%
\sectionstyle
}%
\newcount\sectionnumber \sectionnumber=1
\newskip\inclskip
\vbadness=10000
\hbadness=10000
\def\rginit{\nobreak\nointerlineskip\vskip 0pt plus 1 true pt}
\def\rg{\rginit}
\def\eol{\break}%
\def\e#1{\cx{\hss{$#1$}}\cg}%
\def\hfull{\hbox to \tabhsize}%
\def\eop{}%
\def\tabline{\hbox to \tabhsize}%
\def\tableftline#1{\tabline{{#1}\hss}}%
\def\tabcenterline#1{\tabline{\hss{#1}\hss}}%
\newdimen\tabhsize
\tabhsize=\hsize
\font\titlefont=cmbx10 scaled\magstep 3
\font\authorfont=cmr10 scaled\magstep 2
\font\sectionfont=cmbx10 scaled\magstep 2
\def\dn#1{\down{2.5 true pt}{#1}}
\def\up#1#2{\raise #1\hbox{$\scriptstyle{#2}$}}
\def\down#1#2{\lower #1\hbox{$\scriptstyle{#2}$}}
\font\bulf=cmsy10 scaled\magstep 1
\def\entry#1{$\vcenter{\smash{\hbox to 0pt{\hss{#1}\hss}}}$}
\def\spot{\lower 3.5 true pt\hbox{\entry{\bulf\char'017}}}
\def\hbar{\hbox to 25
true pt{\hss\leaders\hrule\hskip 45 true pt\hss}}
\def\hb{\vcenter{\smash{\hbar}}}
\def\vb{\vcenter{\vbar}}
%
%
\vbox to 0pt {}
\vfil
\centerline{\titlefont Induce/Restrict Matrices}
\vskip 0.75 true in
\centerline{\titlefont for Exceptional Weyl Groups}
\vfil
\centerline{
\hss
\vbox{
\authorfont
\hbox{\hss Dean Alvis\hss}
\vskip 10 true pt
\hbox{\hss Department of Mathematics and Computer Science\hss}
\vskip 10 true pt
\hbox{\hss Indiana University at South Bend\hss}
\vskip 10 true pt
\hbox{\hss South Bend, Indiana, 46634\hss}
}}
\vfil
\eject
\sectionhead{Introduction}

\bigskip
The ordinary character theory of the Weyl groups
$W(A_{\ell})$ of type $A_{\ell}$,
that is, the finite symmetric groups, is well established.
The irreducible characters
are indexed
by combinatorial objects
such as partitions or Young diagrams.
A combinatorial algorithm, the
Murnaghan-Nakayama rule, can be used
to compute the values of the
irreducible characters.
Another combinatorial algorithm,
the Littlewood-Richardson rule,
can be used to compute
the multiplicities with which irreducible
characters of
$W(A_{\ell})$
occur in the characters induced from
parabolic subgroups.
A similar combinatorial machinery is available
in types $B_{\ell}$ (or
$C_{\ell}$) and $D_{\ell}$.
In contrast,
the character theory for the
exceptional Weyl groups has developed
in a more ad-hoc fashion.
In particular, no analogue
the Littlewood-Richardson rule
is available
in the exceptional types.

\smallskip
Let $\Gamma$ be a Dynkin diagram of
type $E_6$, $E_7$, $E_8$, $F_4$ or $G_2$,
and let
$\widetilde\Gamma$ denote
the extended diagram of $\Gamma$.
Let $\Gamma_0$ be a
Dynkin diagram
obtained from
$\widetilde\Gamma$ by removing one node,
and suppose
$\Gamma_0$ is not isomorphic to
$\Gamma$.
Denote by $W = W(\Gamma)$ and
$W_0 = W(\Gamma_0)$ the Weyl
groups corresponding
to $\Gamma$ and $\Gamma_0$, respectively.
Our goal is to
present tables containing the
multiplicities with which
the various irreducible characters of
$W$
occur as constituents in the characters
induced from irreducible
characters of $W_0$.
These tables will be
called {\sl induce/restrict matrices}.
The data in these tables has been used,
for example,
by Lusztig and Spaltenstein, in
the study of
Springer's correspondence in
the exceptional types.
This correspondence establishes
a bijection between the conjugacy
classes of unipotent elements
in a connected reductive algebraic
group over the field of complex
numbers and representations of the
associated Weyl group induced from
special representations of
subgroups of the form
$W(\Gamma^{\prime})$,
$\Gamma^{\prime}\,\subset\,\widetilde{\Gamma}$.
See [AL] and [S] for more details on
this application.

\smallskip
We also present
below the induce/restrict
matrices for
$W_1\,\subset\,W$, where
$W_1$
is a maximal parabolic subgroup of
$W$.
Tables such as these have been used
in the computation of the
generic degrees for the
exceptional Weyl groups,
which are certain polynomials
associated with
$W$.
The generic degrees are important in
the representation theory of finite groups
of Lie type since they occur as factors in
the degree formulas for the irreducible
characters.
See [A] for a remark on how
induce/restrict matrices can be
applied in this context.

\smallskip
This work is organized as follows.
In Section 2, notation is
established for the conjugacy
classes of the
various groups $W_0$, $W$.
In Section 3, the class inclusions for
the pairs
$W_0\,\subset\,W$ are indicated.
In Section 4, notation is
given for the irreducible
characters of the groups
$W_0$, $W$.
This information is needed to
interpret the induce/restrict matrices,
which appear in Section 5.
Section 6 contains
tables of inner products for
certain
pairs of induced characters.

\smallskip
Although the situation in type
$G_2$ is extremely
straightforward, the data for this
case is included for the sake
of completeness.

\smallskip
Throughout this discussion
$W(\Gamma)$ will denote the Weyl
group corresponding to the root system
with Dynkin diagram (or type)
$\Gamma$.
The unit and sign
characters of a Weyl group
will be denoted
$1$ and $\varepsilon$,
respectively.
If
$\chi$ is a character of a Weyl group
$W$,
then the {\sl dual} of
$\chi$ is the character
$\chi^{*} = \varepsilon\!\cdot\!\chi$.
The character of
$W$ induced from the character
$\varphi$
of the subgroup
$W_0$ of $W$ will be denoted
$Ind^{W}_{W_0}\varphi$.
If
$\chi$ and $\psi$ are characters of
$W$, then
the usual inner product of
$\chi$ and $\psi$ will
be denoted
$(\chi,\psi)_{\down{2 true pt}{W}}$.
The set of all irreducible characters of
$W$ is denoted
$Irr(W)$.
If $W$ has central involution
$z$, then
the cyclic subgroup ${<}{z}{>}$ of
$W$ is denoted
$\{\pm1\}$.

\smallskip
More information on root systems,
Weyl groups, and
extended Dynkin diagrams can be
found in [B].

\goodbreak
\medskip
\line{\bf{Acknowledgements.}\hss}

\smallskip
The preparation of this manuscript took place over several years,
during portions of which time the author was at the University of
Oregon, the Massachusetts Institute of Technology, the University
of Notre Dame, and Indiana University at South Bend.  The author
wishes to thank his colleagues at
each of these institutions for their assistance in the
preparation of various versions of this manuscript.
The author is indebted to
Charles~Curtis and George~Lusztig
for posing the problems
that led to this investigation,
and for many helpful conversations.

\vfil
\eject
\sectionhead{Conjugacy Classes}

\bigskip
In this section we review some facts concerning the
conjugacy classes of finite Weyl groups.
More information about the classical types can be found
Chapters 1 and 4 of [JK].
The notation established here will be needed
in the discussion of class inclusions
given in the next section.

\medskip
\leftline{\bf Type $A$.}
\smallskip
By a {\sl partition} of a positive integer $n$, we mean an
(unordered) list of positive integers whose sum is $n$.
The entries on the list are called the {\sl parts} of the
partition.
Partitions will be represented
as formal expressions
${1^{e_1}}{2^{e_2}}\dots$,
with
$k^{e_k}$ representing a part $k$ repeated $e_k$ times.
(Exponents equal to 1 are suppressed, as are parts
with exponent 0.)
For example, if $n = 4$, then the partitions of $n$ are
${1^4}$, ${1^2}{2}$, ${2^2}$, ${1}{3}$, and $4$.
(This notation is unambiguous when used below,
since in all cases $n \le 8$.)

\smallskip
Let $V = V_{\ell}$ be
$\ell$-dimensional Euclidean space,
with orthonormal basis
${\cal B} = \{e_1,e_2,\dots,e_{\ell}\}$.
The set of vectors
$\{\,e_{i}-e_{j}, e_{j}-e_{i}\,|\,1 \le i < j \le \ell\,\}$
is a root system of type $A_{\ell-1}$
(in a subspace of $V$).
The group $W(A_{\ell-1})$ acts faithfully
$V$, and permutes the basis
${\cal B}$.
Via this action, we
identify $W(A_{\ell-1})$
with the group
$P_{\ell}$ of all $\ell\times\ell$
permutation matrices.
(A permutation matrix is a matrix
whose
entries are in $\{0,1\}$
and such that each row and each
column contains exactly one nonzero entry.)
If $w \in W(A_{\ell-1})$, then the sizes
of the orbits of the cyclic group
${<}{w}{>}$ on ${\cal B}$ form a partition
of $\ell$, called the {\sl cycle structure}
of $w$.
The conjugacy class of any element
of $W(A_{\ell-1})$ is uniquely
determined by its cycle structure,
and thus we may parametrize the classes
of $W(A_{\ell-1})$ by
partitions of $\ell$.
If $\lambda$ is a partition of
$\ell$, then $(\lambda)$ will
denote the corresponding
conjugacy class of $W(A_{\ell-1})$.
Thus the identity element is in
the class
$(1^{\ell})$,
while the class
$(\ell)$ contains the
Coxeter elements
(that is, the products of complete
sets of fundamental reflections).

\medskip
\leftline{\bf Types $B$ and $C$.}
\smallskip
To parametrize the conjugacy classes of $W(B_{\ell})$,
we utilize the notion of a {\sl signed partition}
of $\ell$, by which we mean
a formal expression
${1^{e_1}}{\ol{1}^{\dn{f_1}}}{2^{e_2}}{\ol{2}^{\dn{f_2}}}\dots$
such that the exponents
$e_1$, $f_1$, $e_2$, $f_2$, \dots
are nonnegative integers
and ${e_1}+{f_1}+{e_2}+{f_2}+\cdots = \ell$.
(The parts $\ol{k}$ are considered to be
{\sl signed}.)
For example, if $\ell = 3$, then the signed partitions
of $\ell$ are
$1^3$,
${1^2}\ol{1}$,
${1}{\ol{1}^{\dn{2}}}$,
$\ol{1}^{\dn{3}}$,
${1}{2}$, $\ol{1}{2}$, ${1}\ol{2}$, $\ol{1}\ol{2}$,
$3$ and $\ol{3}$.

\smallskip
The set of vectors
$$\{\,\pm e_{i} \, | \, 1 \le i \le \ell \, \} \cup
\{\,
e_{i}+e_{j},
e_{i}-e_{j},
e_{j}-e_{i}
\,|\,
1 \le i < j \le \ell
\,\}$$
is a root system in $V = V_{\ell}$
of type $B_{\ell}$.
The Weyl group $W(B_{\ell})$
acts faithfully on
$V_{\ell}$,
and permutes the basis
${\cal B} = \{e_1,e_2,\dots,e_{\ell}\}$
up to sign.
Via this action, we
identify $W(B_{\ell})$ with the group
of all $\ell\times\ell$ signed
permutation matrices.
(A signed permutation matrix is any
matrix whose entries are in $\{-1,0,1\}$
and such that each row and each column
contains exactly one nonzero entry.)
Let $\Delta = \Delta_{\ell}$ denote
the diagonal subgroup of
$W(B_{\ell})$.
Note $W(B_{\ell})$ is
the semidirect product
of $\Delta$ and
$P = P_{\ell}$.
Let $w = {\delta}{\pi} \in W(B_{\ell})$,
where $\delta \in \Delta$
and $\pi \in P$.
Let $\Omega$ be an orbit
of ${<}{\pi}{>}$
on ${\cal B}$.
We consider $\Omega$
to be {\sl positive} ({\sl negative})
if the number of entries $-1$
in $\delta$ corresponding
to $\Omega$ is even
(odd, respectively).
The signed cycle type of
$w$ is then defined to
be
${1^{e_1}}{\ol{1}^{\dn{f_1}}}{2^{e_2}}{\ol{2}^{\dn{f_2}}}\dots$,
where $e_k$ ($f_k$) is the
number of positive (negative, respectively)
orbits of ${<}{\pi}{>}$ of size $k$
for all $k \ge 1$.
The signed cycle type of $w$ is a signed
partition of $\ell$ which uniquely
determines the conjugacy
class of $w$ in $W(B_{\ell})$.
If $\mu$ is a signed partition of
$\ell$, we denote by
$(\mu)$ the corresponding
conjugacy class of $W(B_{\ell})$.
Thus the central involution
is contained in the
class $(\ol{1}^{\dn{\ell}})$,
while $(\ol{\ell})$ is the class of
Coxeter elements of $W(B_{\ell})$.

\smallskip
If $\pm{e_i}$ is replaced by $\pm{2{e_i}}$ in the
root system above, a root system of type $C_{\ell}$
is obtained.
With this choice of root system, the
action of $W(C_{\ell})$ on $V$
coincides with that of $W(B_{\ell})$.
We identify
$W(C_{\ell})$ with $W(B_{\ell})$, and
use the same notation for
the conjugacy classes
of $W(C_{\ell})$ as for $W(B_{\ell})$.

\medskip
\leftline{\bf Type $D$.}
\smallskip
Let $\Delta^{+} = \Delta^{+}_{\ell}$ be the subgroup
of $\Delta_{\ell}$ containing the elements
with an even number of entries $-1$.
Taking as root system the set
$\{\,e_{i}+e_{j},
e_{i}-e_{j},
e_{j}-e_{i}
\,|\,
1 \le i < j \le \ell\,\}$ in
$V_{\ell}$,
we see that $W(D_{\ell})$
may be identified with the
subgroup $\Delta^{+}\!\cdot\!P$
of $W(B_{\ell})$.  In particular,
$W(D_{\ell})$ has index 2 in
$W(B_{\ell})$.
If $\mu = {1^{e_1}}{\ol{1}^{\dn{f_1}}}{2^{e_2}}{\ol{2}^{\dn{f_2}}}\dots$
is a signed partition of $\ell$, then
$(\mu)\,\subset\,W(D_{\ell})$ if and only if
${f_1}+{f_2}+\cdots$ is even.
If $(\mu)\,\subset\,W(D_{\ell})$, then
$(\mu)$ remains a single conjugacy class
in $W(D_{\ell})$ unless $\ell$
is even, $e_{i} = 0$ for all odd $i \ge 1$,
and $f_{j} = 0$ for all $j \ge 1$.
If $\ell$ is even and
$\mu = {2^{e_2}}{4^{e_4}}\dots$
is a signed partition of $\ell$ with only
unsigned
even parts, then the class $(\mu)$
of $W(B_{\ell})$
splits into two classes in $W(D_{\ell})$;
one, denoted
$({2^{e_2}}{4^{e_4}}\dots)_{+}$,
contains permutation matrices, while the
other, denoted
$({2^{e_2}}{4^{e_4}}\dots)_{-}$,
does not.
The classes $({2^{e_2}}{4^{e_4}}\dots)_{+}$ and
$({2^{e_2}}{4^{e_4}}\dots)_{-}$ are interchanged
by the nonidentity
graph automorphism of
$W(D_{\ell})$
if $\ell$ is even and
$\ell \ge 6$.
We note for future reference that
if
$\ell \equiv 2\!\!\pmod{4}$,
then $({2^{e_2}}{4^{e_4}}\dots)_{+}$ and
$({2^{e_2}}{4^{e_4}}\dots)_{-}$ are interchanged
under multiplication by the
central involution in $W(D_{\ell})$,
while if
$\ell \equiv 0\!\!\pmod{4}$,
then
$({2^{e_2}}{4^{e_4}}\dots)_{+}$
and
$({2^{e_2}}{4^{e_4}}\dots)_{-}$
are
invariant under multiplication by 
the central involution.

\medskip
\leftline{\bf Types $E_6$, $E_7$ and $E_8$.}
\smallskip
In selecting a notation for the conjugacy classes in
the exceptional cases, it was necessary to
find a compromise between the desire for
a uniform notation
on the one hand and consistency
with the notation of
[F1], [F2] and [K] on the other.
The notation used here is similar to
that of [F1],
and any differences with that of [F2] and [K] are
easily reconciled.

\smallskip
By a {\sl virtual partition} of the positive integer $n$, we
mean a formal expression ${1^{e_1}}{2^{e_2}}{3^{e_3}}\dots$,
where the exponents $e_i$ are (possibly negative) integers
and ${e_1}+2{e_2}+3{e_3}+\dots = n$.
Any parts greater than $9$ will be enclosed within
parentheses to avoid ambiguities; for example,
${2^{-1}}{4}{6^{-1}}{(12)}$
is a virtual partition of 8
with largest part $12$.
Let $\Gamma$ be one of the diagrams $E_6$, $E_7$, $E_8$, and put
$W = W(\Gamma)$.
If $w \in W$, then the eigenvalues of $w$ in the
natural reflection representation of $W$
may be described by a
virtual partition
$\nu = {1^{e_1}}{2^{e_2}}\dots$ of $|\Gamma|$,
as follows.
If $e_k > 0$, then $k^{e_k}$
represents $e_k$ full sets of complex $k$-th roots of unity, while if
$e_k < 0$, then $k^{e_k}$ represents $|e_k|$ full sets of complex $2k$-th
roots of unity with $|e_k|$ full sets of complex $k$-th roots of unity 
removed.
Thus ${2^{-1}}{4}{6^{-1}}{(12)}$ represents
the union of a set of
primitive complex 4-th roots of unity
with a set of 12-th roots of unity,
with a set of 6-th roots of unity removed.
The conjugacy
classes of $W$ are, with few exceptions, uniquely determined by
their eigenvalue structure in the reflection representation, and
will be denoted by the corresponding
virtual partition.
When
ambiguity might arise, a subscript `$u$' or `$v$' is
used to indicate
which class is if type $u$ or $v$ as in [F1] and [F2].
The notation used here
is consistent
with that used in [F1],
with the
following notable exceptions.
We denote by
${1}{2}{4}_{u}$ and ${1}{2}{4}_{v}$ the classes of
$W(E_7)$ with eigenvalue type ${1}{2}{4}$ and sizes
${2^4}\!\cdot\!{3^4}\!\cdot\!{5}\!\cdot\!{7}$ and
${2^3}\!\cdot\!{3^3}\!\cdot\!{5}\!\cdot\!{7}$,
respectively.  Also, ${1^{-1}}{2}{6}_{u}$ and
${1^{-1}}{2}{6}_{v}$ will denote the classes of $W(E_7)$
of type ${1^{-1}}{2}{6}$ and sizes
${2^7}\!\cdot\!{3^3}\!\cdot\!{5}\!\cdot\!{7}$ and
${2^7}\!\cdot\!{3^2}\!\cdot\!{5}\!\cdot\!{7}$, respectively.
The notation used here agrees
with that of [F2] for classes not uniquely determined
by their eigenvalue structure, with the following
exception: we denote by
${2}{6}_{vv}$
the class with eigenvalue
structure ${2}{6}$ and centralizer
of order ${2^3}\!\cdot\!{3^3}$; this
class is denoted
${1}\ol{1}{3}\ol{3}$ in [F2].

\smallskip
No two distinct conjugacy classes of $W$ have
the same size and eigenvalue structure; thus each
conjugacy class of $W$ is invariant under all
automorphisms of $W$.

\medskip
\leftline{\bf Type $F_4$.}
\smallskip
Let $W = W(F_4)$.
We identify the generators
of $W$ denoted
$d$, $a$, $\tau$, and ${\tau}{\sigma}$
in [K]
with fundamental
reflections in such a way that $d$ and $a$ correspond to long
roots.
Let $b = \tau^{-1}{a}\tau$,
$c = \sigma^{-1}b\sigma$,
$e = (abcd)^2$, and $z = (abcd)^3$.  We denote the classes of
$W$ by their eigenvalue structure in the reflection representation.
When necessary, an expression
for a representative element
in terms of $a$, \dots, $z$, $\tau$, $\sigma$
is used
as a subscript to distinguish classes with the same
eigenvalue structure.
For example, ${1^2}{2}_{d}$ and
${1^2}{2}_{\tau}$ are the classes of reflections in $W$.  (This
notation differs from that used by Kondo in [K].)

\smallskip
The nontrivial orbits of the
graph automorphisms on the
set of conjugacy classes
of $W(F_4)$ are\hfil
$\{{1}{3}_{e},{1}{3}_{\sigma}\}$,\hfil
$\{-({1}{3}_{e}),-({1}{3}_{\sigma})\}$,\hfil
$\{{1^2}{2}_{d},{1^2}{2}_{\tau}\}$,\hfil
$\{-({1^2}{2}_{d}),-({1^2}{2}_{\tau})\}$,\hfil
$\{{1^{-1}}{2}{3}_{\sigma{d}},{1^{-1}}{2}{3}_{e{\tau}}\}$,%
\hfilneg\goodbreak\noindent
$\{-({1^{-1}}{2}{3}_{\sigma{d}}),-({1^{-1}}{2}{3}_{e{\tau}})\}$,
and
$\{4_{adb},4_{ca\tau}\}$.
(See [K], 2.3.)

\medskip
\leftline{\bf Type $G_2$.}
\smallskip
The classes of $W(G_2)$ will be denoted by their eigenvalue
structure in the reflection representation.  The class $2_{\ell}$
contains the reflections corresponding
to long roots, while
$-(2_{\ell})$
contains the reflections
corresponding to
short roots.
The classes
$2_{\ell}$
and 
$-(2_{\ell})$
are interchanged
by the graph automorphism
of $W(G_2)$.

\bigskip
\leftline{\bf General Remarks.}
\medskip
If $W = W(\Gamma)$ contains a central involution z (acting as
$-1$ in the natural reflection representation) and if $C$
is a conjugacy class of $W$, then
we denote by
$-C$
the conjugacy class
$\{\,{w}{z}\,|\,w \in C\,\}$
of $W$.
For example, $-({1}{2}{4}_{u})$ is the class in $W(E_7)$ with
eigenvalue structure ${1^{-1}}{2^2}{4}$ and size
${2^4}\!\cdot\!{3^4}\!\cdot\!{5}\!\cdot\!{7}$.

\medskip
If $\Gamma_1$,\dots, $\Gamma_k$
are the connected components of
$\Gamma$
and
$C_j$ is a conjugacy
class of $W(\Gamma_j)$ for $1 \le j \le k$,
then we denote by
$C_1\times\dots\times C_k$
the conjugacy class
$\{\,{w_1}\dots{w_k}\,|\,{w_j} \in {C_i}\,for\,1 \le j \le k\,\}$
of
$W(\Gamma) = W(\Gamma_1)\times\dots\times W(\Gamma_k).$
For example, the central involution of $W(E_7)$ is contained
in the conjugacy class $(\ol{1}^6)\times(2)$
of $W({D_6}{A_1})\,\subset\,W(E_7)$.

\vfil
\eject
\sectionhead{Class Inclusions}

\bigskip
We present tables in this section indicating
which conjugacy class $C$ of
$W = W(\Gamma)$
contains a given conjugacy class $C_0$ of
$W_0 = W(\Gamma_0)$.
In many cases, there is only one
conjugacy class $C$ of $W$ with the same eigenvalue structure
as $C_0$.  If the eigenvalue structure of $C_0$ does not
determine a unique conjugacy class of $W$, then the
orders of the centralizers
of elements of
$W_0$ and $W$
can in some cases be used
to resolve the class inclusion.
The effect of multiplication by the central
involution, if present, and a knowledge
of the
powers of the conjugacy
classes
are also useful in some cases.
Choosing as a representative of the conjugacy class
$C_0$ an appropriate product of fundamental
reflections helps to identify other
subgroups $W(\Gamma_1)$,
$\Gamma_1\,\subset\,\widetilde{\Gamma}$,
that have nonempty intersection with $C_0$.
Using
these general techniques and the information given in
[F1], [F2] and [K], the class inclusions can be determined.
We discuss a sample of the class inclusions;
the remaining cases are similar.

\smallskip
If $\widetilde{\Gamma}$ possesses a nontrivial
graph automorphism,
we consider only one representative of each orbit of the
graph automorphisms on the various
subgraphs $\Gamma_0$.
The class inclusions are also given
for $W_1\,\subset\,W$ if
$W_1$ is a maximal parabolic subgroup
of $W$ not contained in any of the
subgroups $W_0$ for $\Gamma_0$
a maximal subdiagram of $\widetilde{\Gamma}$,
$\Gamma_0 \not\approx \Gamma$.

\smallskip
We remark that in the classical types
$A$, $B$, $C$, and $D$, the order of the
centralizer of an element and the
eigenvalues of an element acting in the
reflection representation can easily
be determined from its cycle type or
signed cycle type, as appropriate.
For example, if
$\mu = {1^{e_1}}{\ol{1}^{\dn{f_1}}}{2^{e_2}}{\ol{2}^{\dn{f_2}}}\dots$
is a signed partition of $\ell$,
then the order of the centralizer
in $W(B_{\ell})$ of an element of
the class $(\mu)$ is given by
$$\prod_{k \ge 1}(2k)^{{e_k}+{f_k}}{{e_k}!}{{f_k}!}$$
([JK], Lemma 4.2.10).
Also, the eigenvalue structure of an
element of $(\mu)$ in the reflection
representation of $W(B_{\ell})$
is the virtual partition
obtained from $\mu$ by
replacing powers of signed parts
$\ol{k}^{f_k}$ by ${k}^{\up{1 true pt}{-{f_k}}}{(2k)}^{f_k}$
and combining like powers.

\smallskip
The label
`${*}$' in
${A_2}{A_2^{*}}\,\subset\,F_4$
and
${A_1}{A_1^{*}}\,\subset\,G_2$
is attached to
the connected component
corresponding to short roots.
In other types, the order of the connected components
of $\Gamma_0$
is not important.

\medskip
\goodbreak
\leftline{\bf Type $E_6$.}

\smallskip
The extended Dynkin diagram of type $E_6$ is
given below.
In this and subsequent extended diagrams,
the node labeled `{$\circ$}' corresponds
to the opposite of the highest root.

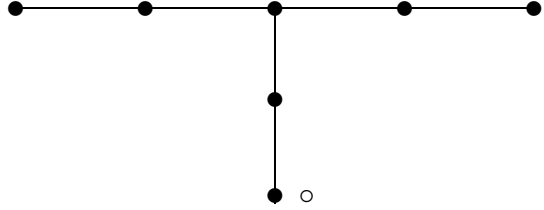
\begingroup
\def\star{\lower 0 true pt\hbox to 0pt{\hss{$\circ$}\hss}}
\def\vbar{\vbox to 25 true pt{\vss\hbox to 0pt{\hss\vrule height 38 true pt\hss}\vss}}
$$
\matrix{
\spot&\hb&\spot&\hb&\spot&\hb&\spot&\hb&\spot\cr
&&&&\vb&&&&\cr
&&&&\spot&&&&\cr
&&&&\vb&&&&\cr
&&&&\spot{\rlap{\hskip 12 true pt{\star}}}&&&&\cr
}
$$
\endgroup

\smallskip
\noindent
The subdiagrams
${A_5}{A_1}$,
${A_2}{A_2}{A_2}\,\subset\,\widetilde{\Gamma}$ and
$D_5\,\subset\Gamma$
are considered.

\smallskip
The conjugacy classes of $W(E_6)$ are
uniquely determined by their eigenvalue structure.
The class inclusions for $W(E_6)$ are given in
Tables~1 through~3.

\def\rg{\rgloc}
\def\bol{\quad}

\def\rg{\noalign{\nobreak\vskip 1 true pt plus 3 true pt}}

\vskip 5 true pt plus 5 true pt
\inclskip=15 true pt
\vbox{%
\tableopen{Class inclusions for $W({A_{5}}{A_{1}})\,\subset\,W(E_{6})$}%
\nobreak\vskip 10 true pt\nobreak%
\baselineskip=\tablebaseskip%
\tabskip=\incltabskip\halign to \hsize{%
\bol\hfil$#$&%
$\,\subset\,#$&&%
$\,\subset\,#$\cr%
{({1^{6}}){\times}({1^{2}})}&
{{1^{6}}}\hskip\inclskip\quad\hfil%
{({3^{2}}){\times}({1^{2}})}&
{{3^{2}}}\hskip\inclskip\quad\hfil%
{({1^{4}}{2}){\times}({2})}&
{{1^{2}}{2^{2}}}\hskip\inclskip\quad\hfil%
{({3^{2}}){\times}({2})}&
{{1^{-2}}{2}{3^{2}}}\hskip\inclskip\hskip-\inclskip\hfil\cr\rg%
{({1^{4}}{2}){\times}({1^{2}})}&
{{1^{4}}{2}}\hskip\inclskip\quad\hfil%
{({1^{2}}{4}){\times}({1^{2}})}&
{{1^{2}}{4}}\hskip\inclskip\quad\hfil%
{({1^{2}}{2^{2}}){\times}({2})}&
{{2^{3}}}\hskip\inclskip\quad\hfil%
{({1^{2}}{4}){\times}({2})}&
{{2}{4}}\hskip\inclskip\hskip-\inclskip\hfil\cr\rg%
{({1^{2}}{2^{2}}){\times}({1^{2}})}&
{{1^{2}}{2^{2}}}\hskip\inclskip\quad\hfil%
{({2}{4}){\times}({1^{2}})}&
{{2}{4}}\hskip\inclskip\quad\hfil%
{({2^{3}}){\times}({2})}&
{{1^{-2}}{2^{4}}}\hskip\inclskip\quad\hfil%
{({2}{4}){\times}({2})}&
{{1^{-2}}{2^{2}}{4}}\hskip\inclskip\hskip-\inclskip\hfil\cr\rg%
{({2^{3}}){\times}({1^{2}})}&
{{2^{3}}}\hskip\inclskip\quad\hfil%
{({1}{5}){\times}({1^{2}})}&
{{1}{5}}\hskip\inclskip\quad\hfil%
{({1^{3}}{3}){\times}({2})}&
{{1}{2}{3}}\hskip\inclskip\quad\hfil%
{({1}{5}){\times}({2})}&
{{1^{-1}}{2}{5}}\hskip\inclskip\hskip-\inclskip\hfil\cr\rg%
{({1^{3}}{3}){\times}({1^{2}})}&
{{1^{3}}{3}}\hskip\inclskip\quad\hfil%
{({6}){\times}({1^{2}})}&
{{6}}\hskip\inclskip\quad\hfil%
{({1}{2}{3}){\times}({2})}&
{{1^{-1}}{2^{2}}{3}}\hskip\inclskip\quad\hfil%
{({6}){\times}({2})}&
{{1^{-2}}{2}{6}}\hskip\inclskip\hskip-\inclskip\hfil\cr\rg%
{({1}{2}{3}){\times}({1^{2}})}&
{{1}{2}{3}}\hskip\inclskip\quad\hfil%
{({1^{6}}){\times}({2})}&
{{1^{4}}{2}}\hskip\inclskip\quad\hfil%
&
\omit\hskip\inclskip\hfil\quad\hfil%
&
\omit\hskip\inclskip\hfil\hskip-\inclskip\hfil\cr\rg%
}}%
\tableclose%

\vskip 5 true pt plus 5 true pt
\inclskip=35 true pt
\vbox{%
\tableopen{Class inclusions for $W({A_{2}}{A_{2}}{A_{2}})\,\subset\,W(E_{6})$}%
\nobreak\vskip 10 true pt\nobreak%
\baselineskip=\tablebaseskip%
\tabskip=\incltabskip\halign to \hsize{%
\bol\hfil$#$&%
$\,\subset\,#$&&%
$\,\subset\,#$\cr%
{({1^{3}}){\times}({1^{3}}){\times}({1^{3}})}&
{{1^{6}}}\hskip\inclskip\quad\hfil%
{({1^{3}}){\times}({1^{3}}){\times}({1}{2})}&
{{1^{4}}{2}}\hskip\inclskip\quad\hfil%
{({1^{3}}){\times}({1^{3}}){\times}({3})}&
{{1^{3}}{3}}\hskip\inclskip\hskip-\inclskip\hfil\cr\rg%
{({1}{2}){\times}({1^{3}}){\times}({1^{3}})}&
{{1^{4}}{2}}\hskip\inclskip\quad\hfil%
{({1}{2}){\times}({1^{3}}){\times}({1}{2})}&
{{1^{2}}{2^{2}}}\hskip\inclskip\quad\hfil%
{({1}{2}){\times}({1^{3}}){\times}({3})}&
{{1}{2}{3}}\hskip\inclskip\hskip-\inclskip\hfil\cr\rg%
{({3}){\times}({1^{3}}){\times}({1^{3}})}&
{{1^{3}}{3}}\hskip\inclskip\quad\hfil%
{({3}){\times}({1^{3}}){\times}({1}{2})}&
{{1}{2}{3}}\hskip\inclskip\quad\hfil%
{({3}){\times}({1^{3}}){\times}({3})}&
{{3^{2}}}\hskip\inclskip\hskip-\inclskip\hfil\cr\rg%
{({1^{3}}){\times}({1}{2}){\times}({1^{3}})}&
{{1^{4}}{2}}\hskip\inclskip\quad\hfil%
{({1^{3}}){\times}({1}{2}){\times}({1}{2})}&
{{1^{2}}{2^{2}}}\hskip\inclskip\quad\hfil%
{({1^{3}}){\times}({1}{2}){\times}({3})}&
{{1}{2}{3}}\hskip\inclskip\hskip-\inclskip\hfil\cr\rg%
{({1}{2}){\times}({1}{2}){\times}({1^{3}})}&
{{1^{2}}{2^{2}}}\hskip\inclskip\quad\hfil%
{({1}{2}){\times}({1}{2}){\times}({1}{2})}&
{{2^{3}}}\hskip\inclskip\quad\hfil%
{({1}{2}){\times}({1}{2}){\times}({3})}&
{{1^{-1}}{2^{2}}{3}}\hskip\inclskip\hskip-\inclskip\hfil\cr\rg%
{({3}){\times}({1}{2}){\times}({1^{3}})}&
{{1}{2}{3}}\hskip\inclskip\quad\hfil%
{({3}){\times}({1}{2}){\times}({1}{2})}&
{{1^{-1}}{2^{2}}{3}}\hskip\inclskip\quad\hfil%
{({3}){\times}({1}{2}){\times}({3})}&
{{1^{-2}}{2}{3^{2}}}\hskip\inclskip\hskip-\inclskip\hfil\cr\rg%
{({1^{3}}){\times}({3}){\times}({1^{3}})}&
{{1^{3}}{3}}\hskip\inclskip\quad\hfil%
{({1^{3}}){\times}({3}){\times}({1}{2})}&
{{1}{2}{3}}\hskip\inclskip\quad\hfil%
{({1^{3}}){\times}({3}){\times}({3})}&
{{3^{2}}}\hskip\inclskip\hskip-\inclskip\hfil\cr\rg%
{({1}{2}){\times}({3}){\times}({1^{3}})}&
{{1}{2}{3}}\hskip\inclskip\quad\hfil%
{({1}{2}){\times}({3}){\times}({1}{2})}&
{{1^{-1}}{2^{2}}{3}}\hskip\inclskip\quad\hfil%
{({1}{2}){\times}({3}){\times}({3})}&
{{1^{-2}}{2}{3^{2}}}\hskip\inclskip\hskip-\inclskip\hfil\cr\rg%
{({3}){\times}({3}){\times}({1^{3}})}&
{{3^{2}}}\hskip\inclskip\quad\hfil%
{({3}){\times}({3}){\times}({1}{2})}&
{{1^{-2}}{2}{3^{2}}}\hskip\inclskip\quad\hfil%
{({3}){\times}({3}){\times}({3})}&
{{1^{-3}}{3^{3}}}\hskip\inclskip\hskip-\inclskip\hfil\cr\rg%
}}%
\tableclose%

\vskip 0pt plus 5 true pt
\eject
\def\rg{\rgloc}

\vskip 0pt plus 50 pt
\eject
\inclskip=32 true pt
\vbox{%
\tableopen{Class inclusions for $W(D_{5})\,\subset\,W(E_{6})$}%
\nobreak\vskip 10 true pt\nobreak%
\baselineskip=\tablebaseskip%
\tabskip=\incltabskip\halign to \hsize{%
\bol\hfil$#$&%
$\,\subset\,#$&&%
$\,\subset\,#$\cr%
{({1^{5}})}&
{{1^{6}}}\hskip\inclskip\quad\hfil%
{({1}{4})}&
{{1^{2}}{4}}\hskip\inclskip\quad\hfil%
{({1^{2}}{\ol{1}}{\ol{2}})}&
{{1^{2}}{4}}\hskip\inclskip\quad\hfil%
{({1}{\ol{1}}{\ol{3}})}&
{{1}{2}{3^{-1}}{6}}\hskip\inclskip\hskip-\inclskip\hfil\cr\rg%
{({1^{3}}{2})}&
{{1^{4}}{2}}\hskip\inclskip\quad\hfil%
{({5})}&
{{1}{5}}\hskip\inclskip\quad\hfil%
{({\ol{1}}{2}{\ol{2}})}&
{{2}{4}}\hskip\inclskip\quad\hfil%
{({\ol{1}^{\dn{3}}}{\ol{2}})}&
{{1^{-2}}{2^{2}}{4}}\hskip\inclskip\hskip-\inclskip\hfil\cr\rg%
{({1}{2^{2}})}&
{{1^{2}}{2^{2}}}\hskip\inclskip\quad\hfil%
{({1^{3}}{\ol{1}^{\dn{2}}})}&
{{1^{2}}{2^{2}}}\hskip\inclskip\quad\hfil%
{({1}{\ol{1}^{\dn{4}}})}&
{{1^{-2}}{2^{4}}}\hskip\inclskip\quad\hfil%
{({\ol{2}}{\ol{3}})}&
{{1}{2^{-1}}{3^{-1}}{4}{6}}\hskip\inclskip\hskip-\inclskip\hfil\cr\rg%
{({1^{2}}{3})}&
{{1^{3}}{3}}\hskip\inclskip\quad\hfil%
{({1}{\ol{1}^{\dn{2}}}{2})}&
{{2^{3}}}\hskip\inclskip\quad\hfil%
{({1}{\ol{2}^{\dn{2}}})}&
{{1^{2}}{2^{-2}}{4^{2}}}\hskip\inclskip\quad\hfil%
{({\ol{1}}{\ol{4}})}&
{{2}{4^{-1}}{8}}\hskip\inclskip\hskip-\inclskip\hfil\cr\rg%
{({2}{3})}&
{{1}{2}{3}}\hskip\inclskip\quad\hfil%
{({\ol{1}^{\dn{2}}}{3})}&
{{1^{-1}}{2^{2}}{3}}\hskip\inclskip\quad\hfil%
&
\omit\hskip\inclskip\hfil\quad\hfil%
&
\omit\hskip\inclskip\hfil\hskip-\inclskip\hfil\cr\rg%
}}%
\tableclose%

\medskip
\leftline{\bf Type $E_7$.}

\smallskip
The extended diagram for $E_7$ is as follows.

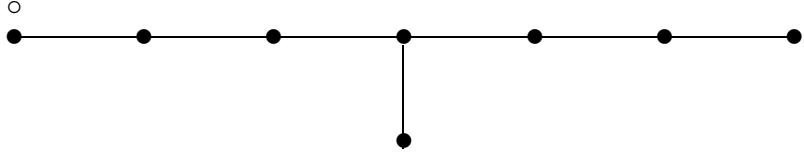
\begingroup
\def\star{\lower 3 true pt\hbox to 0pt{\hss{$\circ$}\hss}}
\def\vbar{\vbox to 30 true pt{\vss\hbox to 0pt{\hss\vrule height 39 true pt\hss}\vss}}
$$
\matrix{
\star&&&&&&&&&&&&&&\cr
\spot&\hb&\spot&\hb&
\spot&\hb&\spot&\hb&\spot&\hb&\spot&\hb&\spot\cr
&&&&&&\hbox{\hss{$\vb$\hskip 1 true pt}\hss}&&&&&&\cr
&&&&&&\hbox{\lower 2 true pt\hbox{$\spot$}}&&&&&&\cr
}
$$
\endgroup

\smallskip
\noindent
The subdiagrams considered for class inclusions are
${D_6}{A_1}$, 
$A_7$,
${A_5}{A_2}$,
${A_3}{A_3}{A_1}\,\subset\,\widetilde{\Gamma}$
and
$E_6\,\subset\,\Gamma$.

\smallskip
The class inclusions for $W(E_6)\,\subset\,W(E_7)$
are given in [F1].
Also, the permutation character of
$W(E_7)/\{\pm{1\}}$ corresponding to the image of $W(A_7)$
is given in [F1], and from this
the class inclusions for $W(A_7)$ in
$W(E_7)$ are easily determined.

\smallskip
We view $W({D_6}{A_1})$ as embedded in $W(E_7)$ in such a way
that the class
$(6)_{+}\times(1^2)$
of $W({D_6}{A_1})$
meets the
class $[6]$ of $W(A_7)$.
The class inclusions for $W({D_6}{A_1})$ can be determined
by using the cases $W(E_6)\,\subset\,W(E_7)$
and $W(A_7)\,\subset\,W(E_7)$ considered earlier.
For
example, the classes
$(6)_{+}\times(1^2)$
and $(6)_{+}\times(2)$
of $W({D_6}{A_1})$
meet the subgroup $W(A_7)$,
and thus are contained in the classes
$-({1^{-1}}{2}{6}_{u})$ and ${1^{-1}}{2}{6}_{v}$
of $W(E_7)$, respectively.
Hence the classes $({6})_{-}\times(1^2) = -(({6})_{+}\times(2))$ and
$({6})_{-}\times(2) = -((6)_{+}\times(1^2))$ are contained
in $-({1^{-1}}{2}{6}_{v})$ and ${1^{-1}}{2}{6}_{u}$, respectively.
As another example, consider the class
$({\ol{1}^{\dn{2}}}{2^2})\times(1^2)$ 
of $W({D_6}{A_1})$.
This class is contained in the same class of $W(E_7)$ as
the class $({2^3})_{-}\times(2)$,
as can be seen by choosing as
representative elements appropriate products
of fundamental reflections.
Further, $(2^3)_{-}\times(2) =
-((2^3)_{+}\times(1^2))$, and
$(2^3)_{+}\times(1^2)$ is contained in the
class
$-({1^{-1}}{2^4}_{\llap{$\scriptstyle{u}$}})$ of $W(E_7)$
(since $(2^3)_{+}\times(1^2)$ meets $W(A_7)$).
Therefore $({\ol{1}^{\dn{2}}}{2^2})\times(1^2)$
and $(2^3)_{-}\times(2)$
are both contained in the class
${1^{-1}}{2^4}_{\llap{$\scriptstyle{u}$}}$.

\smallskip
The remaining cases $W({A_5}{A_2})$
and $W({A_3}{A_3}{A_1})$ follow easily
from the ones already considered.
The class inclusions for type $E_7$ are given in
Tables~4 through~8.
\vskip 0pt plus 30 pt

\medskip
\inclskip=7 true pt
\def\rg{\noalign{\nobreak\vskip 4 true pt plus 5 true pt}}
\vbox{%
\tableopen{Class inclusions for $W({D_{6}}{A_{1}})\,\subset\,W(E_{7})$}%
\nobreak\vskip 10 true pt\nobreak%
\baselineskip=\tablebaseskip%
\tabskip=\incltabskip\halign to \hsize{%
\bol\hfil$#$&%
$\,\subset\,#$&&%
$\,\subset\,#$\cr%
{({1^{6}}){\times}({1^{2}})}&
{{1^{7}}}\hskip\inclskip\quad\hfil%
{({\ol{1}^{\dn{4}}}{2}){\times}({1^{2}})}&
{-({1^{3}}{2^{2}})}\hskip\inclskip\quad\hfil%
{({6})_{-}{\times}({2})}&
{{1^{-1}}{2}{6}_{u}}\hskip\inclskip\hskip-\inclskip\hfil\cr\rg%
{({1^{4}}{2}){\times}({1^{2}})}&
{-({1^{-5}}{2^{6}})}\hskip\inclskip\quad\hfil%
{({2}{\ol{2}^{\dn{2}}}){\times}({1^{2}})}&
{-({1^{-1}}{4^{2}})}\hskip\inclskip\quad\hfil%
{({1^{4}}{\ol{1}^{\dn{2}}}){\times}({2})}&
{-({1^{-1}}{2^{4}}_{\llap{$\scriptstyle{v}$}{}})}\hskip\inclskip\hskip-\inclskip\hfil\cr\rg%
{({1^{2}}{2^{2}}){\times}({1^{2}})}&
{{1^{3}}{2^{2}}}\hskip\inclskip\quad\hfil%
{({\ol{1}}{2}{\ol{3}}){\times}({1^{2}})}&
{-({2^{2}}{3})}\hskip\inclskip\quad\hfil%
{({1^{2}}{\ol{1}^{\dn{2}}}{2}){\times}({2})}&
{{1^{-1}}{2^{4}}_{\llap{$\scriptstyle{u}$}{}}}\hskip\inclskip\hskip-\inclskip\hfil\cr\rg%
{({2^{3}})_{+}{\times}({1^{2}})}&
{-({1^{-1}}{2^{4}}_{\llap{$\scriptstyle{u}$}{}})}\hskip\inclskip\quad\hfil%
{({1}{\ol{1}^{\dn{3}}}{\ol{2}}){\times}({1^{2}})}&
{-({1}{2}{4}_{v})}\hskip\inclskip\quad\hfil%
{({\ol{1}^{\dn{2}}}{2^{2}}){\times}({2})}&
{-({1^{3}}{2^{2}})}\hskip\inclskip\hskip-\inclskip\hfil\cr\rg%
{({2^{3}})_{-}{\times}({1^{2}})}&
{-({1^{-1}}{2^{4}}_{\llap{$\scriptstyle{v}$}{}})}\hskip\inclskip\quad\hfil%
{({1}{\ol{2}}{\ol{3}}){\times}({1^{2}})}&
{-({1^{-2}}{2}{3}{4})}\hskip\inclskip\quad\hfil%
{({1}{\ol{1}^{\dn{2}}}{3}){\times}({2})}&
{-({1^{2}}{2}{3^{-1}}{6})}\hskip\inclskip\hskip-\inclskip\hfil\cr\rg%
{({1^{3}}{3}){\times}({1^{2}})}&
{{1^{4}}{3}}\hskip\inclskip\quad\hfil%
{({1}{\ol{1}}{\ol{4}}){\times}({1^{2}})}&
{-({1^{-1}}{2^{2}}{4^{-1}}{8})}\hskip\inclskip\quad\hfil%
{({\ol{1}^{\dn{2}}}{4}){\times}({2})}&
{{1^{-3}}{2^{3}}{4}}\hskip\inclskip\hskip-\inclskip\hfil\cr\rg%
{({1}{2}{3}){\times}({1^{2}})}&
{-({1^{-2}}{2^{3}}{3^{-1}}{6})}\hskip\inclskip\quad\hfil%
{({\ol{1}^{\dn{6}}}){\times}({1^{2}})}&
{{1^{-5}}{2^{6}}}\hskip\inclskip\quad\hfil%
{({1^{3}}{\ol{1}}{\ol{2}}){\times}({2})}&
{{1}{2}{4}_{v}}\hskip\inclskip\hskip-\inclskip\hfil\cr\rg%
{({3^{2}}){\times}({1^{2}})}&
{{1}{3^{2}}}\hskip\inclskip\quad\hfil%
{({\ol{1}^{\dn{2}}}{\ol{2}^{\dn{2}}}){\times}({1^{2}})}&
{{1^{-1}}{4^{2}}}\hskip\inclskip\quad\hfil%
{({1}{\ol{1}}{2}{\ol{2}}){\times}({2})}&
{-({1}{2}{4}_{u})}\hskip\inclskip\hskip-\inclskip\hfil\cr\rg%
{({1^{2}}{4}){\times}({1^{2}})}&
{-({1^{-3}}{2^{3}}{4})}\hskip\inclskip\quad\hfil%
{({\ol{1}^{\dn{3}}}{\ol{3}}){\times}({1^{2}})}&
{{1^{-2}}{2^{3}}{3^{-1}}{6}}\hskip\inclskip\quad\hfil%
{({\ol{1}}{\ol{2}}{3}){\times}({2})}&
{{1^{-2}}{2}{3}{4}}\hskip\inclskip\hskip-\inclskip\hfil\cr\rg%
{({2}{4})_{+}{\times}({1^{2}})}&
{{1}{2}{4}_{u}}\hskip\inclskip\quad\hfil%
{({\ol{3}^{\dn{2}}}){\times}({1^{2}})}&
{{1}{3^{-2}}{6^{2}}}\hskip\inclskip\quad\hfil%
{({1^{2}}{\ol{1}^{\dn{4}}}){\times}({2})}&
{-({1^{3}}{2^{2}})}\hskip\inclskip\hskip-\inclskip\hfil\cr\rg%
{({2}{4})_{-}{\times}({1^{2}})}&
{{1}{2}{4}_{v}}\hskip\inclskip\quad\hfil%
{({\ol{2}}{\ol{4}}){\times}({1^{2}})}&
{{1}{2^{-1}}{8}}\hskip\inclskip\quad\hfil%
{({1^{2}}{\ol{2}^{\dn{2}}}){\times}({2})}&
{-({1^{-1}}{4^{2}})}\hskip\inclskip\hskip-\inclskip\hfil\cr\rg%
{({1}{5}){\times}({1^{2}})}&
{{1^{2}}{5}}\hskip\inclskip\quad\hfil%
{({\ol{1}}{\ol{5}}){\times}({1^{2}})}&
{{2}{5^{-1}}{(10)}}\hskip\inclskip\quad\hfil%
{({1^{2}}{\ol{1}}{\ol{3}}){\times}({2})}&
{-({2^{2}}{3})}\hskip\inclskip\hskip-\inclskip\hfil\cr\rg%
{({6})_{+}{\times}({1^{2}})}&
{-({1^{-1}}{2}{6}_{u})}\hskip\inclskip\quad\hfil%
{({1^{6}}){\times}({2})}&
{-({1^{-5}}{2^{6}})}\hskip\inclskip\quad\hfil%
{({\ol{1}^{\dn{4}}}{2}){\times}({2})}&
{{1^{-5}}{2^{6}}}\hskip\inclskip\hskip-\inclskip\hfil\cr\rg%
{({6})_{-}{\times}({1^{2}})}&
{-({1^{-1}}{2}{6}_{v})}\hskip\inclskip\quad\hfil%
{({1^{4}}{2}){\times}({2})}&
{{1^{3}}{2^{2}}}\hskip\inclskip\quad\hfil%
{({2}{\ol{2}^{\dn{2}}}){\times}({2})}&
{{1^{-1}}{4^{2}}}\hskip\inclskip\hskip-\inclskip\hfil\cr\rg%
{({1^{4}}{\ol{1}^{\dn{2}}}){\times}({1^{2}})}&
{{1^{3}}{2^{2}}}\hskip\inclskip\quad\hfil%
{({1^{2}}{2^{2}}){\times}({2})}&
{-({1^{-1}}{2^{4}}_{\llap{$\scriptstyle{u}$}{}})}\hskip\inclskip\quad\hfil%
{({\ol{1}}{2}{\ol{3}}){\times}({2})}&
{{1^{-2}}{2^{3}}{3^{-1}}{6}}\hskip\inclskip\hskip-\inclskip\hfil\cr\rg%
{({1^{2}}{\ol{1}^{\dn{2}}}{2}){\times}({1^{2}})}&
{-({1^{-1}}{2^{4}}_{\llap{$\scriptstyle{u}$}{}})}\hskip\inclskip\quad\hfil%
{({2^{3}})_{+}{\times}({2})}&
{{1^{-1}}{2^{4}}_{\llap{$\scriptstyle{v}$}{}}}\hskip\inclskip\quad\hfil%
{({1}{\ol{1}^{\dn{3}}}{\ol{2}}){\times}({2})}&
{{1^{-3}}{2^{3}}{4}}\hskip\inclskip\hskip-\inclskip\hfil\cr\rg%
{({\ol{1}^{\dn{2}}}{2^{2}}){\times}({1^{2}})}&
{{1^{-1}}{2^{4}}_{\llap{$\scriptstyle{u}$}{}}}\hskip\inclskip\quad\hfil%
{({2^{3}})_{-}{\times}({2})}&
{{1^{-1}}{2^{4}}_{\llap{$\scriptstyle{u}$}{}}}\hskip\inclskip\quad\hfil%
{({1}{\ol{2}}{\ol{3}}){\times}({2})}&
{{3^{-1}}{4}{6}}\hskip\inclskip\hskip-\inclskip\hfil\cr\rg%
{({1}{\ol{1}^{\dn{2}}}{3}){\times}({1^{2}})}&
{{2^{2}}{3}}\hskip\inclskip\quad\hfil%
{({1^{3}}{3}){\times}({2})}&
{-({1^{-2}}{2^{3}}{3^{-1}}{6})}\hskip\inclskip\quad\hfil%
{({1}{\ol{1}}{\ol{4}}){\times}({2})}&
{{1^{-1}}{2^{2}}{4^{-1}}{8}}\hskip\inclskip\hskip-\inclskip\hfil\cr\rg%
{({\ol{1}^{\dn{2}}}{4}){\times}({1^{2}})}&
{-({1}{2}{4}_{u})}\hskip\inclskip\quad\hfil%
{({1}{2}{3}){\times}({2})}&
{{2^{2}}{3}}\hskip\inclskip\quad\hfil%
{({\ol{1}^{\dn{6}}}){\times}({2})}&
{-({1^{7}})}\hskip\inclskip\hskip-\inclskip\hfil\cr\rg%
{({1^{3}}{\ol{1}}{\ol{2}}){\times}({1^{2}})}&
{-({1^{-3}}{2^{3}}{4})}\hskip\inclskip\quad\hfil%
{({3^{2}}){\times}({2})}&
{-({1}{3^{-2}}{6^{2}})}\hskip\inclskip\quad\hfil%
{({\ol{1}^{\dn{2}}}{\ol{2}^{\dn{2}}}){\times}({2})}&
{-({1^{3}}{2^{-2}}{4^{2}})}\hskip\inclskip\hskip-\inclskip\hfil\cr\rg%
{({1}{\ol{1}}{2}{\ol{2}}){\times}({1^{2}})}&
{{1}{2}{4}_{u}}\hskip\inclskip\quad\hfil%
{({1^{2}}{4}){\times}({2})}&
{{1}{2}{4}_{u}}\hskip\inclskip\quad\hfil%
{({\ol{1}^{\dn{3}}}{\ol{3}}){\times}({2})}&
{-({1^{4}}{3})}\hskip\inclskip\hskip-\inclskip\hfil\cr\rg%
{({\ol{1}}{\ol{2}}{3}){\times}({1^{2}})}&
{-({3^{-1}}{4}{6})}\hskip\inclskip\quad\hfil%
{({2}{4})_{+}{\times}({2})}&
{-({1}{2}{4}_{v})}\hskip\inclskip\quad\hfil%
{({\ol{3}^{\dn{2}}}){\times}({2})}&
{-({1}{3^{2}})}\hskip\inclskip\hskip-\inclskip\hfil\cr\rg%
{({1^{2}}{\ol{1}^{\dn{4}}}){\times}({1^{2}})}&
{{1^{-1}}{2^{4}}_{\llap{$\scriptstyle{v}$}{}}}\hskip\inclskip\quad\hfil%
{({2}{4})_{-}{\times}({2})}&
{-({1}{2}{4}_{u})}\hskip\inclskip\quad\hfil%
{({\ol{2}}{\ol{4}}){\times}({2})}&
{-({1}{2^{-1}}{8})}\hskip\inclskip\hskip-\inclskip\hfil\cr\rg%
{({1^{2}}{\ol{2}^{\dn{2}}}){\times}({1^{2}})}&
{{1^{3}}{2^{-2}}{4^{2}}}\hskip\inclskip\quad\hfil%
{({1}{5}){\times}({2})}&
{-({2}{5^{-1}}{(10)})}\hskip\inclskip\quad\hfil%
{({\ol{1}}{\ol{5}}){\times}({2})}&
{-({1^{2}}{5})}\hskip\inclskip\hskip-\inclskip\hfil\cr\rg%
{({1^{2}}{\ol{1}}{\ol{3}}){\times}({1^{2}})}&
{{1^{2}}{2}{3^{-1}}{6}}\hskip\inclskip\quad\hfil%
{({6})_{+}{\times}({2})}&
{{1^{-1}}{2}{6}_{v}}\hskip\inclskip\quad\hfil%
&
\omit\hskip\inclskip\hfil\hskip-\inclskip\hfil\cr\rg%
}}%
\tableclose%
\eject
\def\rg{\rgloc}

\inclskip=4.5 true pt
\def\rg{\noalign{\nobreak\vskip 4 true pt plus 2 true pt}}
\vbox{%
\tableopen{Class inclusions for $W(A_{7})\,\subset\,W(E_{7})$}%
\nobreak\vskip 10 true pt\nobreak%
\baselineskip=\tablebaseskip%
\tabskip=\incltabskip\halign to \hsize{%
\bol\hfil$#$&%
$\,\subset\,#$&&%
$\,\subset\,#$\cr%
{({1^{8}})}&
{{1^{7}}}\hskip\inclskip\quad\hfil%
{({1^{3}}{2}{3})}&
{-({1^{-2}}{2^{3}}{3^{-1}}{6})}\hskip\inclskip\quad\hfil%
{({2^{2}}{4})}&
{-({1}{2}{4}_{v})}\hskip\inclskip\quad\hfil%
{({3}{5})}&
{{1^{-1}}{3}{5}}\hskip\inclskip\hskip-\inclskip\hfil\cr\rg%
{({1^{6}}{2})}&
{-({1^{-5}}{2^{6}})}\hskip\inclskip\quad\hfil%
{({1}{2^{2}}{3})}&
{{2^{2}}{3}}\hskip\inclskip\quad\hfil%
{({1}{3}{4})}&
{-({3^{-1}}{4}{6})}\hskip\inclskip\quad\hfil%
{({1^{2}}{6})}&
{-({1^{-1}}{2}{6}_{u})}\hskip\inclskip\hskip-\inclskip\hfil\cr\rg%
{({1^{4}}{2^{2}})}&
{{1^{3}}{2^{2}}}\hskip\inclskip\quad\hfil%
{({1^{2}}{3^{2}})}&
{{1}{3^{2}}}\hskip\inclskip\quad\hfil%
{({4^{2}})}&
{{1^{-1}}{4^{2}}}\hskip\inclskip\quad\hfil%
{({2}{6})}&
{{1^{-1}}{2}{6}_{v}}\hskip\inclskip\hskip-\inclskip\hfil\cr\rg%
{({1^{2}}{2^{3}})}&
{-({1^{-1}}{2^{4}}_{\llap{$\scriptstyle{u}$}{}})}\hskip\inclskip\quad\hfil%
{({2}{3^{2}})}&
{-({1}{3^{-2}}{6^{2}})}\hskip\inclskip\quad\hfil%
{({1^{3}}{5})}&
{{1^{2}}{5}}\hskip\inclskip\quad\hfil%
{({1}{7})}&
{{7}}\hskip\inclskip\hskip-\inclskip\hfil\cr\rg%
{({2^{4}})}&
{{1^{-1}}{2^{4}}_{\llap{$\scriptstyle{v}$}{}}}\hskip\inclskip\quad\hfil%
{({1^{4}}{4})}&
{-({1^{-3}}{2^{3}}{4})}\hskip\inclskip\quad\hfil%
{({1}{2}{5})}&
{-({2}{5^{-1}}{(10)})}\hskip\inclskip\quad\hfil%
{({8})}&
{-({1}{2^{-1}}{8})}\hskip\inclskip\hskip-\inclskip\hfil\cr\rg%
{({1^{5}}{3})}&
{{1^{4}}{3}}\hskip\inclskip\quad\hfil%
{({1^{2}}{2}{4})}&
{{1}{2}{4}_{u}}\hskip\inclskip\quad\hfil%
&
\omit\hskip\inclskip\hfil\quad\hfil%
&
\omit\hskip\inclskip\hfil\hskip-\inclskip\hfil\cr\rg%
}}%
\tableclose%

\vskip 60 true pt

\inclskip=3.5 true pt
\vbox{%
\tableopen{Class inclusions for $W({A_{5}}{A_{2}})\,\subset\,W(E_{7})$}%
\nobreak\vskip 10 true pt\nobreak%
\baselineskip=\tablebaseskip%
\tabskip=\incltabskip\halign to \hsize{%
\bol\hfil$#$&%
$\,\subset\,#$&&%
$\,\subset\,#$\cr%
{({1^{6}}){\times}({1^{3}})}&
{{1^{7}}}\hskip\inclskip\quad\hfil%
{({1^{6}}){\times}({1}{2})}&
{-({1^{-5}}{2^{6}})}\hskip\inclskip\quad\hfil%
{({1^{6}}){\times}({3})}&
{{1^{4}}{3}}\hskip\inclskip\hskip-\inclskip\hfil\cr\rg%
{({1^{4}}{2}){\times}({1^{3}})}&
{-({1^{-5}}{2^{6}})}\hskip\inclskip\quad\hfil%
{({1^{4}}{2}){\times}({1}{2})}&
{{1^{3}}{2^{2}}}\hskip\inclskip\quad\hfil%
{({1^{4}}{2}){\times}({3})}&
{-({1^{-2}}{2^{3}}{3^{-1}}{6})}\hskip\inclskip\hskip-\inclskip\hfil\cr\rg%
{({1^{2}}{2^{2}}){\times}({1^{3}})}&
{{1^{3}}{2^{2}}}\hskip\inclskip\quad\hfil%
{({1^{2}}{2^{2}}){\times}({1}{2})}&
{-({1^{-1}}{2^{4}}_{\llap{$\scriptstyle{u}$}{}})}\hskip\inclskip\quad\hfil%
{({1^{2}}{2^{2}}){\times}({3})}&
{{2^{2}}{3}}\hskip\inclskip\hskip-\inclskip\hfil\cr\rg%
{({2^{3}}){\times}({1^{3}})}&
{-({1^{-1}}{2^{4}}_{\llap{$\scriptstyle{v}$}{}})}\hskip\inclskip\quad\hfil%
{({2^{3}}){\times}({1}{2})}&
{{1^{-1}}{2^{4}}_{\llap{$\scriptstyle{u}$}{}}}\hskip\inclskip\quad\hfil%
{({2^{3}}){\times}({3})}&
{-({1^{2}}{2}{3^{-1}}{6})}\hskip\inclskip\hskip-\inclskip\hfil\cr\rg%
{({1^{3}}{3}){\times}({1^{3}})}&
{{1^{4}}{3}}\hskip\inclskip\quad\hfil%
{({1^{3}}{3}){\times}({1}{2})}&
{-({1^{-2}}{2^{3}}{3^{-1}}{6})}\hskip\inclskip\quad\hfil%
{({1^{3}}{3}){\times}({3})}&
{{1}{3^{2}}}\hskip\inclskip\hskip-\inclskip\hfil\cr\rg%
{({1}{2}{3}){\times}({1^{3}})}&
{-({1^{-2}}{2^{3}}{3^{-1}}{6})}\hskip\inclskip\quad\hfil%
{({1}{2}{3}){\times}({1}{2})}&
{{2^{2}}{3}}\hskip\inclskip\quad\hfil%
{({1}{2}{3}){\times}({3})}&
{-({1}{3^{-2}}{6^{2}})}\hskip\inclskip\hskip-\inclskip\hfil\cr\rg%
{({3^{2}}){\times}({1^{3}})}&
{{1}{3^{2}}}\hskip\inclskip\quad\hfil%
{({3^{2}}){\times}({1}{2})}&
{-({1}{3^{-2}}{6^{2}})}\hskip\inclskip\quad\hfil%
{({3^{2}}){\times}({3})}&
{{1^{-2}}{3^{3}}}\hskip\inclskip\hskip-\inclskip\hfil\cr\rg%
{({1^{2}}{4}){\times}({1^{3}})}&
{-({1^{-3}}{2^{3}}{4})}\hskip\inclskip\quad\hfil%
{({1^{2}}{4}){\times}({1}{2})}&
{{1}{2}{4}_{u}}\hskip\inclskip\quad\hfil%
{({1^{2}}{4}){\times}({3})}&
{-({3^{-1}}{4}{6})}\hskip\inclskip\hskip-\inclskip\hfil\cr\rg%
{({2}{4}){\times}({1^{3}})}&
{{1}{2}{4}_{v}}\hskip\inclskip\quad\hfil%
{({2}{4}){\times}({1}{2})}&
{-({1}{2}{4}_{u})}\hskip\inclskip\quad\hfil%
{({2}{4}){\times}({3})}&
{{1^{-2}}{2}{3}{4}}\hskip\inclskip\hskip-\inclskip\hfil\cr\rg%
{({1}{5}){\times}({1^{3}})}&
{{1^{2}}{5}}\hskip\inclskip\quad\hfil%
{({1}{5}){\times}({1}{2})}&
{-({2}{5^{-1}}{(10)})}\hskip\inclskip\quad\hfil%
{({1}{5}){\times}({3})}&
{{1^{-1}}{3}{5}}\hskip\inclskip\hskip-\inclskip\hfil\cr\rg%
{({6}){\times}({1^{3}})}&
{-({1^{-1}}{2}{6}_{v})}\hskip\inclskip\quad\hfil%
{({6}){\times}({1}{2})}&
{{1^{-1}}{2}{6}_{u}}\hskip\inclskip\quad\hfil%
{({6}){\times}({3})}&
{-({1^{2}}{2^{-2}}{3^{-1}}{6^{2}})}\hskip\inclskip\hskip-\inclskip\hfil\cr\rg%
}}%
\tableclose%

\eject
\def\rg{\rgloc}

\def\bol{\qquad}
\inclskip=60 true pt
\def\rg{\noalign{\nobreak\vskip 0.5 true pt plus 1 true pt}}
\vbox{%
\tableopen{Class inclusions for $W({A_{3}}{A_{3}}{A_{1}})\,\subset\,W(E_{7})$}%
\nobreak\vskip 10 true pt\nobreak%
\baselineskip=\tablebaseskip%
\tabskip=\incltabskip\halign to \hsize{%
\bol\hfil$#$&%
$\,\subset\,#$&&%
$\,\subset\,#$\cr%
{({1^{4}}){\times}({1^{4}}){\times}({1^{2}})}&
{{1^{7}}}\hskip\inclskip\quad\hfil%
{({1^{4}}){\times}({1^{4}}){\times}({2})}&
{-({1^{-5}}{2^{6}})}\hskip\inclskip\hskip-\inclskip\hfil\cr\rg%
{({1^{2}}{2}){\times}({1^{4}}){\times}({1^{2}})}&
{-({1^{-5}}{2^{6}})}\hskip\inclskip\quad\hfil%
{({1^{2}}{2}){\times}({1^{4}}){\times}({2})}&
{{1^{3}}{2^{2}}}\hskip\inclskip\hskip-\inclskip\hfil\cr\rg%
{({2^{2}}){\times}({1^{4}}){\times}({1^{2}})}&
{{1^{3}}{2^{2}}}\hskip\inclskip\quad\hfil%
{({2^{2}}){\times}({1^{4}}){\times}({2})}&
{-({1^{-1}}{2^{4}}_{\llap{$\scriptstyle{v}$}{}})}\hskip\inclskip\hskip-\inclskip\hfil\cr\rg%
{({1}{3}){\times}({1^{4}}){\times}({1^{2}})}&
{{1^{4}}{3}}\hskip\inclskip\quad\hfil%
{({1}{3}){\times}({1^{4}}){\times}({2})}&
{-({1^{-2}}{2^{3}}{3^{-1}}{6})}\hskip\inclskip\hskip-\inclskip\hfil\cr\rg%
{({4}){\times}({1^{4}}){\times}({1^{2}})}&
{-({1^{-3}}{2^{3}}{4})}\hskip\inclskip\quad\hfil%
{({4}){\times}({1^{4}}){\times}({2})}&
{{1}{2}{4}_{v}}\hskip\inclskip\hskip-\inclskip\hfil\cr\rg%
{({1^{4}}){\times}({1^{2}}{2}){\times}({1^{2}})}&
{-({1^{-5}}{2^{6}})}\hskip\inclskip\quad\hfil%
{({1^{4}}){\times}({1^{2}}{2}){\times}({2})}&
{{1^{3}}{2^{2}}}\hskip\inclskip\hskip-\inclskip\hfil\cr\rg%
{({1^{2}}{2}){\times}({1^{2}}{2}){\times}({1^{2}})}&
{{1^{3}}{2^{2}}}\hskip\inclskip\quad\hfil%
{({1^{2}}{2}){\times}({1^{2}}{2}){\times}({2})}&
{-({1^{-1}}{2^{4}}_{\llap{$\scriptstyle{u}$}{}})}\hskip\inclskip\hskip-\inclskip\hfil\cr\rg%
{({2^{2}}){\times}({1^{2}}{2}){\times}({1^{2}})}&
{-({1^{-1}}{2^{4}}_{\llap{$\scriptstyle{u}$}{}})}\hskip\inclskip\quad\hfil%
{({2^{2}}){\times}({1^{2}}{2}){\times}({2})}&
{{1^{-1}}{2^{4}}_{\llap{$\scriptstyle{u}$}{}}}\hskip\inclskip\hskip-\inclskip\hfil\cr\rg%
{({1}{3}){\times}({1^{2}}{2}){\times}({1^{2}})}&
{-({1^{-2}}{2^{3}}{3^{-1}}{6})}\hskip\inclskip\quad\hfil%
{({1}{3}){\times}({1^{2}}{2}){\times}({2})}&
{{2^{2}}{3}}\hskip\inclskip\hskip-\inclskip\hfil\cr\rg%
{({4}){\times}({1^{2}}{2}){\times}({1^{2}})}&
{{1}{2}{4}_{u}}\hskip\inclskip\quad\hfil%
{({4}){\times}({1^{2}}{2}){\times}({2})}&
{-({1}{2}{4}_{u})}\hskip\inclskip\hskip-\inclskip\hfil\cr\rg%
{({1^{4}}){\times}({2^{2}}){\times}({1^{2}})}&
{{1^{3}}{2^{2}}}\hskip\inclskip\quad\hfil%
{({1^{4}}){\times}({2^{2}}){\times}({2})}&
{-({1^{-1}}{2^{4}}_{\llap{$\scriptstyle{v}$}{}})}\hskip\inclskip\hskip-\inclskip\hfil\cr\rg%
{({1^{2}}{2}){\times}({2^{2}}){\times}({1^{2}})}&
{-({1^{-1}}{2^{4}}_{\llap{$\scriptstyle{u}$}{}})}\hskip\inclskip\quad\hfil%
{({1^{2}}{2}){\times}({2^{2}}){\times}({2})}&
{{1^{-1}}{2^{4}}_{\llap{$\scriptstyle{u}$}{}}}\hskip\inclskip\hskip-\inclskip\hfil\cr\rg%
{({2^{2}}){\times}({2^{2}}){\times}({1^{2}})}&
{{1^{-1}}{2^{4}}_{\llap{$\scriptstyle{v}$}{}}}\hskip\inclskip\quad\hfil%
{({2^{2}}){\times}({2^{2}}){\times}({2})}&
{-({1^{3}}{2^{2}})}\hskip\inclskip\hskip-\inclskip\hfil\cr\rg%
{({1}{3}){\times}({2^{2}}){\times}({1^{2}})}&
{{2^{2}}{3}}\hskip\inclskip\quad\hfil%
{({1}{3}){\times}({2^{2}}){\times}({2})}&
{-({1^{2}}{2}{3^{-1}}{6})}\hskip\inclskip\hskip-\inclskip\hfil\cr\rg%
{({4}){\times}({2^{2}}){\times}({1^{2}})}&
{-({1}{2}{4}_{v})}\hskip\inclskip\quad\hfil%
{({4}){\times}({2^{2}}){\times}({2})}&
{{1^{-3}}{2^{3}}{4}}\hskip\inclskip\hskip-\inclskip\hfil\cr\rg%
{({1^{4}}){\times}({1}{3}){\times}({1^{2}})}&
{{1^{4}}{3}}\hskip\inclskip\quad\hfil%
{({1^{4}}){\times}({1}{3}){\times}({2})}&
{-({1^{-2}}{2^{3}}{3^{-1}}{6})}\hskip\inclskip\hskip-\inclskip\hfil\cr\rg%
{({1^{2}}{2}){\times}({1}{3}){\times}({1^{2}})}&
{-({1^{-2}}{2^{3}}{3^{-1}}{6})}\hskip\inclskip\quad\hfil%
{({1^{2}}{2}){\times}({1}{3}){\times}({2})}&
{{2^{2}}{3}}\hskip\inclskip\hskip-\inclskip\hfil\cr\rg%
{({2^{2}}){\times}({1}{3}){\times}({1^{2}})}&
{{2^{2}}{3}}\hskip\inclskip\quad\hfil%
{({2^{2}}){\times}({1}{3}){\times}({2})}&
{-({1^{2}}{2}{3^{-1}}{6})}\hskip\inclskip\hskip-\inclskip\hfil\cr\rg%
{({1}{3}){\times}({1}{3}){\times}({1^{2}})}&
{{1}{3^{2}}}\hskip\inclskip\quad\hfil%
{({1}{3}){\times}({1}{3}){\times}({2})}&
{-({1}{3^{-2}}{6^{2}})}\hskip\inclskip\hskip-\inclskip\hfil\cr\rg%
{({4}){\times}({1}{3}){\times}({1^{2}})}&
{-({3^{-1}}{4}{6})}\hskip\inclskip\quad\hfil%
{({4}){\times}({1}{3}){\times}({2})}&
{{1^{-2}}{2}{3}{4}}\hskip\inclskip\hskip-\inclskip\hfil\cr\rg%
{({1^{4}}){\times}({4}){\times}({1^{2}})}&
{-({1^{-3}}{2^{3}}{4})}\hskip\inclskip\quad\hfil%
{({1^{4}}){\times}({4}){\times}({2})}&
{{1}{2}{4}_{v}}\hskip\inclskip\hskip-\inclskip\hfil\cr\rg%
{({1^{2}}{2}){\times}({4}){\times}({1^{2}})}&
{{1}{2}{4}_{u}}\hskip\inclskip\quad\hfil%
{({1^{2}}{2}){\times}({4}){\times}({2})}&
{-({1}{2}{4}_{u})}\hskip\inclskip\hskip-\inclskip\hfil\cr\rg%
{({2^{2}}){\times}({4}){\times}({1^{2}})}&
{-({1}{2}{4}_{v})}\hskip\inclskip\quad\hfil%
{({2^{2}}){\times}({4}){\times}({2})}&
{{1^{-3}}{2^{3}}{4}}\hskip\inclskip\hskip-\inclskip\hfil\cr\rg%
{({1}{3}){\times}({4}){\times}({1^{2}})}&
{-({3^{-1}}{4}{6})}\hskip\inclskip\quad\hfil%
{({1}{3}){\times}({4}){\times}({2})}&
{{1^{-2}}{2}{3}{4}}\hskip\inclskip\hskip-\inclskip\hfil\cr\rg%
{({4}){\times}({4}){\times}({1^{2}})}&
{{1^{-1}}{4^{2}}}\hskip\inclskip\quad\hfil%
{({4}){\times}({4}){\times}({2})}&
{-({1^{3}}{2^{-2}}{4^{2}})}\hskip\inclskip\hskip-\inclskip\hfil\cr\rg%
}}%
\tableclose%

\def\bol{\quad}
\medskip
\inclskip=15 true pt
\vbox{%
\tableopen{Class inclusions for $W(E_{6})\,\subset\,W(E_{7})$}%
\nobreak\vskip 10 true pt\nobreak%
\baselineskip=\tablebaseskip%
\tabskip=\incltabskip\halign to \hsize{%
\bol\hfil$#$&%
$\,\subset\,#$&&%
$\,\subset\,#$\cr%
{{1^{6}}}&
{{1^{7}}}\hskip\inclskip\quad\hfil%
{{1^{-2}}{2}{6}}&
{{1^{-1}}{2}{6}_{v}}\hskip\inclskip\quad\hfil%
{{1^{2}}{4}}&
{-({1^{-3}}{2^{3}}{4})}\hskip\inclskip\hskip-\inclskip\hfil\cr\rg%
{{1^{2}}{2^{2}}}&
{{1^{3}}{2^{2}}}\hskip\inclskip\quad\hfil%
{{1^{-3}}{3^{3}}}&
{{1^{-2}}{3^{3}}}\hskip\inclskip\quad\hfil%
{{1^{-2}}{2^{2}}{4}}&
{-({1}{2}{4}_{v})}\hskip\inclskip\hskip-\inclskip\hfil\cr\rg%
{{1^{-2}}{2^{4}}}&
{{1^{-1}}{2^{4}}_{\llap{$\scriptstyle{v}$}{}}}\hskip\inclskip\quad\hfil%
{{1}{2^{-2}}{3^{-1}}{6^{2}}}&
{{1^{2}}{2^{-2}}{3^{-1}}{6^{2}}}\hskip\inclskip\quad\hfil%
{{2}{4^{-1}}{8}}&
{-({1^{-1}}{2^{2}}{4^{-1}}{8})}\hskip\inclskip\hskip-\inclskip\hfil\cr\rg%
{{1^{2}}{2^{-2}}{4^{2}}}&
{{1^{3}}{2^{-2}}{4^{2}}}\hskip\inclskip\quad\hfil%
{{1^{-1}}{2}{3}{4^{-1}}{6^{-1}}{(12)}}&
{{2}{3}{4^{-1}}{6^{-1}}{(12)}}\hskip\inclskip\quad\hfil%
{{1}{2}{3}}&
{-({1^{-2}}{2^{3}}{3^{-1}}{6})}\hskip\inclskip\hskip-\inclskip\hfil\cr\rg%
{{2}{4}}&
{{1}{2}{4}_{u}}\hskip\inclskip\quad\hfil%
{{3^{-1}}{9}}&
{{1}{3^{-1}}{9}}\hskip\inclskip\quad\hfil%
{{1^{-2}}{2}{3^{2}}}&
{-({1}{3^{-2}}{6^{2}})}\hskip\inclskip\hskip-\inclskip\hfil\cr\rg%
{{1^{3}}{3}}&
{{1^{4}}{3}}\hskip\inclskip\quad\hfil%
{{1}{5}}&
{{1^{2}}{5}}\hskip\inclskip\quad\hfil%
{{6}}&
{-({1^{-1}}{2}{6}_{u})}\hskip\inclskip\hskip-\inclskip\hfil\cr\rg%
{{1}{2}{3^{-1}}{6}}&
{{1^{2}}{2}{3^{-1}}{6}}\hskip\inclskip\quad\hfil%
{{1^{4}}{2}}&
{-({1^{-5}}{2^{6}})}\hskip\inclskip\quad\hfil%
{{1}{2^{-1}}{3^{-1}}{4}{6}}&
{-({1^{-2}}{2}{3}{4})}\hskip\inclskip\hskip-\inclskip\hfil\cr\rg%
{{1^{-1}}{2^{2}}{3}}&
{{2^{2}}{3}}\hskip\inclskip\quad\hfil%
{{2^{3}}}&
{-({1^{-1}}{2^{4}}_{\llap{$\scriptstyle{u}$}{}})}\hskip\inclskip\quad\hfil%
{{1^{-1}}{2}{5}}&
{-({2}{5^{-1}}{(10)})}\hskip\inclskip\hskip-\inclskip\hfil\cr\rg%
{{3^{2}}}&
{{1}{3^{2}}}\hskip\inclskip\quad\hfil%
&
\omit\hskip\inclskip\hfil\quad\hfil%
&
\omit\hskip\inclskip\hfil\hskip-\inclskip\hfil\cr\rg%
}}%
\tableclose%
\eject
\def\rg{\rgloc}
\leftline{\bf Type $E_8$.}

\smallskip
The extended diagram for $E_8$ appears below.

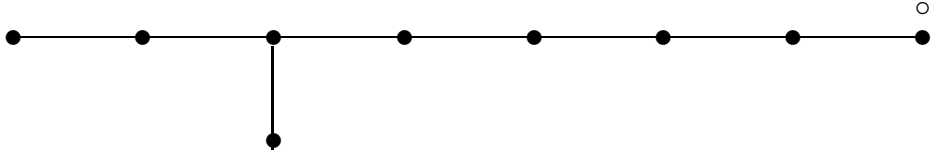
\begingroup
\def\star{\lower 3 true pt\hbox to 0pt{\hss{$\circ$}\hss}}
\def\vbar{\vbox to 30 true pt{\vss\hbox to 0pt{\hss\vrule height 39 true pt\hss}\vss}}
$$
\matrix{
&&&&&&&&&&&&&&\star\cr
\spot&\hb&\spot&\hb&\spot&\hb&
\spot&\hb&\spot&\hb&\spot&\hb&\spot&\hb&\spot\cr
&&&&\hbox{\hss{$\vb$\hskip 1 true pt}\hss}&&&&&&&&&&\cr
&&&&\hbox{\lower 2 true pt\hbox{$\spot$}}&&&&&&&&&&\cr
}
$$
\endgroup

\medskip
\noindent
We consider all maximal subdiagrams $\Gamma_0
\not= \Gamma$
of $\widetilde{\Gamma}$.

\smallskip
The class inclusions for
$W(A_8)$
can be found using
the permutation character
$960_p$ given in [F2];
this character corresponds
to the
image of
$W(A_8)$ in
$W(E_8)/\{\pm1\}$.
In particular, the classes
$(2^4)$,
$({2^2}{4})$,
$(4^2)$,
$({2}{6})$ and
$(8)$ of $W(A_8)$ are contained in
the classes
${2^4}_{\llap{$\scriptstyle{u}$}}$, ${2^2}{4}_{u}$,
${4^2}_{\llap{$\scriptstyle{u}$}}$,
${2}{6}_{u}$,
and ${8}_{u}$ of $W(E_8)$.

\smallskip
We view $W(D_8)$ as embedded in $W(E_8)$ in such a way
that the class $({8})_{-}$ of $W(D_8)$ meets $W(A_8)$.
(The group $W(D_8)$ coincides with the monomial
subgroup denoted $M$ in [F2]; this choice of
embedding ensures
that the class
$(8)_{+}$ is represented in $M$ by a
permutation matrix.)
The classes
$({1^4}{\ol{1}^{\dn{4}}})$,
$(2^4)_{+}$,
$({2^2}{4})_{+}$,
$(4^2)_{+}$,
$({2}{6})_{+}$,
$({1}{\ol{1}}{3}{\ol{3}})$ and
$(8)_{+}$ of $W(D_8)$
are contained in the classes
${2^4}_{\llap{$\scriptstyle{v}$}}$,
${2^4}_{\llap{$\scriptstyle{v}$}}$,
${2^2}{4}_{v}$,
${4^2}_{\llap{$\scriptstyle{v}$}}$,
${2}{6}_{v}$,
${2}{6}_{vv}$
and ${8}_{v}$, respectively,
while the classes
$({2^4})_{-}$,
$({2^2}{4})_{-}$,
$({4^2})_{-}$,
$({2}{6})_{-}$, and
$({8})_{-}$
are contained in
${2^4}_{\llap{$\scriptstyle{u}$}}$,
${2^2}{4}_{u}$,
${4^2_{u}}$,
${2}{6}_{u}$, and
$8_{u}$, respectively
([F2], pp. 113-114).
The remaining class inclusions
are easily determined using
the general techniques described earlier
and the case
$W(A_8)\,\subset\,W(E_8)$.

\smallskip
The permutation character $120_p$
corresponding to the image of $W({E_7}{A_1}$)
in $W(E_8)/\{\pm{1}\}$
is given in [F2],
and from this the class
inclusions for $W({E_7}{A_1})$
are easily resolved.

\smallskip
The class inclusions for
$W({A_7}{A_1})$, $W({A_5}{A_2}{A_1})$, $W({A_4}{A_4})$, $W({D_5}{A_3})$, and
$W({E_6}{A_2})$ in $W(E_8)$
follow
from the cases already considered.
For example, consider $W({A_7}{A_1})$.  A class of the form
$C\times(1^2)$
whose eigenvalues do not determine
a unique class of $W(E_8)$
meets a class of
type `${+}$' in $W(D_8)$,
and hence has type `$v$' in $W(E_8)$.
(In particular,
$({2}{6})\times(1^2)\,\subset\,{2}{6}_{v}$.)
On the other hand,
a class
$C\times(2)$
whose eigenvalues do not
determine a unique class in $W(E_8)$
meets $W(A_8)$, and so
has type `$u$' in $W(E_8)$.
As another example, one sees that the class
${1}{2}{3^{-1}}{6}\times(3)$ of $W({E_6}{A_2})$ meets the class
$({1}{\ol{1}}{3}{\ol{3}})$ of $W(D_8)$ by choosing as a representative the
product of an appropriate subset of the set of fundamental reflections.
Thus ${1}{2}{3^{-1}}{6}\times(3)$ is contained in the class ${2}{6}_{vv}$ of
$W(E_8)$.  The class inclusions for type $E_8$ are given in 
Tables~9 through~16.

\vfil
\eject
\inclskip=20 true pt
\vbox{%
\tableopen{Class inclusions for $W(D_{8})\,\subset\,W(E_{8})$}%
\nobreak\vskip 10 true pt\nobreak%
\baselineskip=\tablebaseskip%
\tabskip=\incltabskip\halign to \hsize{%
\bol\hfil$#$&%
$\,\subset\,#$&&%
$\,\subset\,#$\cr%
{({1^{8}})}&
{{1^{8}}}\hskip\inclskip\quad\hfil%
{({1^{2}}{\ol{1}^{\dn{2}}}{4})}&
{{2^{2}}{4}_{u}}\hskip\inclskip\quad\hfil%
{({\ol{2}}{3}{\ol{3}})}&
{{2^{-1}}{4}{6}}\hskip\inclskip\hskip-\inclskip\hfil\cr\rg%
{({1^{6}}{2})}&
{{1^{6}}{2}}\hskip\inclskip\quad\hfil%
{({\ol{1}^{\dn{2}}}{2}{4})}&
{-({1^{2}}{2}{4})}\hskip\inclskip\quad\hfil%
{({\ol{1}}{3}{\ol{4}})}&
{-({1}{3^{-1}}{4^{-1}}{6}{8})}\hskip\inclskip\hskip-\inclskip\hfil\cr\rg%
{({1^{4}}{2^{2}})}&
{{1^{4}}{2^{2}}}\hskip\inclskip\quad\hfil%
{({1}{\ol{1}^{\dn{2}}}{5})}&
{-({1}{2}{5^{-1}}{(10)})}\hskip\inclskip\quad\hfil%
{({1^{2}}{\ol{1}^{\dn{6}}})}&
{-({1^{4}}{2^{2}})}\hskip\inclskip\hskip-\inclskip\hfil\cr\rg%
{({1^{2}}{2^{3}})}&
{{1^{2}}{2^{3}}}\hskip\inclskip\quad\hfil%
{({\ol{1}^{\dn{2}}}{6})}&
{-({1^{2}}{6})}\hskip\inclskip\quad\hfil%
{({1^{2}}{\ol{1}^{\dn{2}}}{\ol{2}^{\dn{2}}})}&
{{4^{2}}_{\llap{$\scriptstyle{v}$}{}}}\hskip\inclskip\hskip-\inclskip\hfil\cr\rg%
{({2^{4}})_{+}}&
{{2^{4}}_{\llap{$\scriptstyle{v}$}{}}}\hskip\inclskip\quad\hfil%
{({1^{5}}{\ol{1}}{\ol{2}})}&
{{1^{4}}{4}}\hskip\inclskip\quad\hfil%
{({1^{2}}{\ol{1}^{\dn{3}}}{\ol{3}})}&
{-({1}{2^{2}}{3})}\hskip\inclskip\hskip-\inclskip\hfil\cr\rg%
{({2^{4}})_{-}}&
{{2^{4}}_{\llap{$\scriptstyle{u}$}{}}}\hskip\inclskip\quad\hfil%
{({1^{3}}{\ol{1}}{2}{\ol{2}})}&
{{1^{2}}{2}{4}}\hskip\inclskip\quad\hfil%
{({1^{2}}{\ol{3}^{\dn{2}}})}&
{{1^{2}}{3^{-2}}{6^{2}}}\hskip\inclskip\hskip-\inclskip\hfil\cr\rg%
{({1^{5}}{3})}&
{{1^{5}}{3}}\hskip\inclskip\quad\hfil%
{({1}{\ol{1}}{2^{2}}{\ol{2}})}&
{{2^{2}}{4}_{u}}\hskip\inclskip\quad\hfil%
{({1^{2}}{\ol{2}}{\ol{4}})}&
{{1^{2}}{2^{-1}}{8}}\hskip\inclskip\hskip-\inclskip\hfil\cr\rg%
{({1^{3}}{2}{3})}&
{{1^{3}}{2}{3}}\hskip\inclskip\quad\hfil%
{({1^{2}}{\ol{1}}{\ol{2}}{3})}&
{{1}{3}{4}}\hskip\inclskip\quad\hfil%
{({1^{2}}{\ol{1}}{\ol{5}})}&
{{1}{2}{5^{-1}}{(10)}}\hskip\inclskip\hskip-\inclskip\hfil\cr\rg%
{({1}{2^{2}}{3})}&
{{1}{2^{2}}{3}}\hskip\inclskip\quad\hfil%
{({\ol{1}}{2}{\ol{2}}{3})}&
{-({1}{3^{-1}}{4}{6})}\hskip\inclskip\quad\hfil%
{({\ol{1}^{\dn{6}}}{2})}&
{-({1^{6}}{2})}\hskip\inclskip\hskip-\inclskip\hfil\cr\rg%
{({1^{2}}{3^{2}})}&
{{1^{2}}{3^{2}}}\hskip\inclskip\quad\hfil%
{({1}{\ol{1}}{\ol{2}}{4})}&
{{4^{2}}_{\llap{$\scriptstyle{u}$}{}}}\hskip\inclskip\quad\hfil%
{({\ol{1}^{\dn{2}}}{2}{\ol{2}^{\dn{2}}})}&
{-({1^{2}}{2^{-1}}{4^{2}})}\hskip\inclskip\hskip-\inclskip\hfil\cr\rg%
{({2}{3^{2}})}&
{{2}{3^{2}}}\hskip\inclskip\quad\hfil%
{({\ol{1}}{\ol{2}}{5})}&
{-({1}{2^{-1}}{4}{5^{-1}}{(10)})}\hskip\inclskip\quad\hfil%
{({\ol{1}^{\dn{3}}}{2}{\ol{3}})}&
{-({1^{3}}{2}{3})}\hskip\inclskip\hskip-\inclskip\hfil\cr\rg%
{({1^{4}}{4})}&
{{1^{4}}{4}}\hskip\inclskip\quad\hfil%
{({1^{4}}{\ol{1}^{\dn{4}}})}&
{{2^{4}}_{\llap{$\scriptstyle{v}$}{}}}\hskip\inclskip\quad\hfil%
{({2}{\ol{3}^{\dn{2}}})}&
{-({2}{3^{2}})}\hskip\inclskip\hskip-\inclskip\hfil\cr\rg%
{({1^{2}}{2}{4})}&
{{1^{2}}{2}{4}}\hskip\inclskip\quad\hfil%
{({1^{4}}{\ol{2}^{\dn{2}}})}&
{{1^{4}}{2^{-2}}{4^{2}}}\hskip\inclskip\quad\hfil%
{({2}{\ol{2}}{\ol{4}})}&
{{8}_{v}}\hskip\inclskip\hskip-\inclskip\hfil\cr\rg%
{({2^{2}}{4})_{+}}&
{{2^{2}}{4}_{v}}\hskip\inclskip\quad\hfil%
{({1^{4}}{\ol{1}}{\ol{3}})}&
{{1^{3}}{2}{3^{-1}}{6}}\hskip\inclskip\quad\hfil%
{({\ol{1}}{2}{\ol{5}})}&
{-({1}{2}{5})}\hskip\inclskip\hskip-\inclskip\hfil\cr\rg%
{({2^{2}}{4})_{-}}&
{{2^{2}}{4}_{u}}\hskip\inclskip\quad\hfil%
{({1^{2}}{\ol{1}^{\dn{4}}}{2})}&
{-({1^{2}}{2^{3}})}\hskip\inclskip\quad\hfil%
{({1}{\ol{1}^{\dn{5}}}{\ol{2}})}&
{-({1^{4}}{4})}\hskip\inclskip\hskip-\inclskip\hfil\cr\rg%
{({1}{3}{4})}&
{{1}{3}{4}}\hskip\inclskip\quad\hfil%
{({1^{2}}{2}{\ol{2}^{\dn{2}}})}&
{{1^{2}}{2^{-1}}{4^{2}}}\hskip\inclskip\quad\hfil%
{({1}{\ol{1}}{\ol{2}^{\dn{3}}})}&
{{2^{-2}}{4^{3}}}\hskip\inclskip\hskip-\inclskip\hfil\cr\rg%
{({4^{2}})_{+}}&
{{4^{2}}_{\llap{$\scriptstyle{v}$}{}}}\hskip\inclskip\quad\hfil%
{({1^{2}}{\ol{1}}{2}{\ol{3}})}&
{{1}{2^{2}}{3^{-1}}{6}}\hskip\inclskip\quad\hfil%
{({1}{\ol{1}^{\dn{2}}}{\ol{2}}{\ol{3}})}&
{-({1}{3}{4})}\hskip\inclskip\hskip-\inclskip\hfil\cr\rg%
{({4^{2}})_{-}}&
{{4^{2}}_{\llap{$\scriptstyle{u}$}{}}}\hskip\inclskip\quad\hfil%
{({\ol{1}^{\dn{4}}}{2^{2}})}&
{-({1^{4}}{2^{2}})}\hskip\inclskip\quad\hfil%
{({1}{\ol{1}^{\dn{3}}}{\ol{4}})}&
{-({1^{2}}{2}{4^{-1}}{8})}\hskip\inclskip\hskip-\inclskip\hfil\cr\rg%
{({1^{3}}{5})}&
{{1^{3}}{5}}\hskip\inclskip\quad\hfil%
{({2^{2}}{\ol{2}^{\dn{2}}})}&
{{4^{2}}_{\llap{$\scriptstyle{v}$}{}}}\hskip\inclskip\quad\hfil%
{({1}{\ol{3}}{\ol{4}})}&
{{1}{3^{-1}}{4^{-1}}{6}{8}}\hskip\inclskip\hskip-\inclskip\hfil\cr\rg%
{({1}{2}{5})}&
{{1}{2}{5}}\hskip\inclskip\quad\hfil%
{({\ol{1}}{2^{2}}{\ol{3}})}&
{-({1}{2^{2}}{3})}\hskip\inclskip\quad\hfil%
{({1}{\ol{2}}{\ol{5}})}&
{{1}{2^{-1}}{4}{5^{-1}}{(10)}}\hskip\inclskip\hskip-\inclskip\hfil\cr\rg%
{({3}{5})}&
{{3}{5}}\hskip\inclskip\quad\hfil%
{({1}{\ol{1}^{\dn{4}}}{3})}&
{-({1^{3}}{2}{3^{-1}}{6})}\hskip\inclskip\quad\hfil%
{({1}{\ol{1}}{\ol{6}})}&
{{2}{6^{-1}}{(12)}}\hskip\inclskip\hskip-\inclskip\hfil\cr\rg%
{({1^{2}}{6})}&
{{1^{2}}{6}}\hskip\inclskip\quad\hfil%
{({1}{\ol{2}^{\dn{2}}}{3})}&
{{1}{2^{-2}}{3}{4^{2}}}\hskip\inclskip\quad\hfil%
{({\ol{1}^{\dn{8}}})}&
{-({1^{8}})}\hskip\inclskip\hskip-\inclskip\hfil\cr\rg%
{({2}{6})_{+}}&
{{2}{6}_{v}}\hskip\inclskip\quad\hfil%
{({1}{\ol{1}}{3}{\ol{3}})}&
{{2}{6}_{vv}}\hskip\inclskip\quad\hfil%
{({\ol{1}^{\dn{4}}}{\ol{2}^{\dn{2}}})}&
{-({1^{4}}{2^{-2}}{4^{2}})}\hskip\inclskip\hskip-\inclskip\hfil\cr\rg%
{({2}{6})_{-}}&
{{2}{6}_{u}}\hskip\inclskip\quad\hfil%
{({\ol{1}^{\dn{4}}}{4})}&
{-({1^{4}}{4})}\hskip\inclskip\quad\hfil%
{({\ol{2}^{\dn{4}}})}&
{{2^{-4}}{4^{4}}}\hskip\inclskip\hskip-\inclskip\hfil\cr\rg%
{({1}{7})}&
{{1}{7}}\hskip\inclskip\quad\hfil%
{({\ol{2}^{\dn{2}}}{4})}&
{{2^{-2}}{4^{3}}}\hskip\inclskip\quad\hfil%
{({\ol{1}^{\dn{5}}}{\ol{3}})}&
{-({1^{5}}{3})}\hskip\inclskip\hskip-\inclskip\hfil\cr\rg%
{({8})_{+}}&
{{8}_{v}}\hskip\inclskip\quad\hfil%
{({\ol{1}}{\ol{3}}{4})}&
{-({1}{3}{4})}\hskip\inclskip\quad\hfil%
{({\ol{1}}{\ol{2}^{\dn{2}}}{\ol{3}})}&
{-({1}{2^{-2}}{3}{4^{2}})}\hskip\inclskip\hskip-\inclskip\hfil\cr\rg%
{({8})_{-}}&
{{8}_{u}}\hskip\inclskip\quad\hfil%
{({1^{3}}{\ol{1}^{\dn{3}}}{\ol{2}})}&
{{2^{2}}{4}_{v}}\hskip\inclskip\quad\hfil%
{({\ol{1}^{\dn{2}}}{\ol{3}^{\dn{2}}})}&
{-({1^{2}}{3^{2}})}\hskip\inclskip\hskip-\inclskip\hfil\cr\rg%
{({1^{6}}{\ol{1}^{\dn{2}}})}&
{{1^{4}}{2^{2}}}\hskip\inclskip\quad\hfil%
{({1^{3}}{\ol{2}}{\ol{3}})}&
{{1^{3}}{2^{-1}}{3^{-1}}{4}{6}}\hskip\inclskip\quad\hfil%
{({\ol{1}^{\dn{2}}}{\ol{2}}{\ol{4}})}&
{-({1^{2}}{2^{-1}}{8})}\hskip\inclskip\hskip-\inclskip\hfil\cr\rg%
{({1^{4}}{\ol{1}^{\dn{2}}}{2})}&
{{1^{2}}{2^{3}}}\hskip\inclskip\quad\hfil%
{({1^{3}}{\ol{1}}{\ol{4}})}&
{{1^{2}}{2}{4^{-1}}{8}}\hskip\inclskip\quad\hfil%
{({\ol{4}^{\dn{2}}})}&
{{4^{-2}}{8^{2}}}\hskip\inclskip\hskip-\inclskip\hfil\cr\rg%
{({1^{2}}{\ol{1}^{\dn{2}}}{2^{2}})}&
{{2^{4}}_{\llap{$\scriptstyle{u}$}{}}}\hskip\inclskip\quad\hfil%
{({1}{\ol{1}^{\dn{3}}}{2}{\ol{2}})}&
{-({1^{2}}{2}{4})}\hskip\inclskip\quad\hfil%
{({\ol{1}^{\dn{3}}}{\ol{5}})}&
{-({1^{3}}{5})}\hskip\inclskip\hskip-\inclskip\hfil\cr\rg%
{({\ol{1}^{\dn{2}}}{2^{3}})}&
{-({1^{2}}{2^{3}})}\hskip\inclskip\quad\hfil%
{({1}{2}{\ol{2}}{\ol{3}})}&
{{1}{3^{-1}}{4}{6}}\hskip\inclskip\quad\hfil%
{({\ol{3}}{\ol{5}})}&
{-({3}{5})}\hskip\inclskip\hskip-\inclskip\hfil\cr\rg%
{({1^{3}}{\ol{1}^{\dn{2}}}{3})}&
{{1}{2^{2}}{3}}\hskip\inclskip\quad\hfil%
{({1}{\ol{1}}{2}{\ol{4}})}&
{{2^{2}}{4^{-1}}{8}}\hskip\inclskip\quad\hfil%
{({\ol{2}}{\ol{6}})}&
{{2^{-1}}{4}{6^{-1}}{(12)}}\hskip\inclskip\hskip-\inclskip\hfil\cr\rg%
{({1}{\ol{1}^{\dn{2}}}{2}{3})}&
{-({1}{2^{2}}{3^{-1}}{6})}\hskip\inclskip\quad\hfil%
{({\ol{1}^{\dn{3}}}{\ol{2}}{3})}&
{-({1^{3}}{2^{-1}}{3^{-1}}{4}{6})}\hskip\inclskip\quad\hfil%
{({\ol{1}}{\ol{7}})}&
{-({1}{7})}\hskip\inclskip\hskip-\inclskip\hfil\cr\rg%
{({\ol{1}^{\dn{2}}}{3^{2}})}&
{-({1^{2}}{3^{-2}}{6^{2}})}\hskip\inclskip\quad\hfil%
&
\omit\hskip\inclskip\hfil\quad\hfil%
&
\omit\hskip\inclskip\hfil\hskip-\inclskip\hfil\cr\rg%
}}%
\tableclose%

\inclskip=30 true pt
\def\rg{\noalign{\nobreak\vskip 1.5 true pt plus 2 true pt}}
\vbox{%
\tableopen{Class inclusions for $W(A_{8})\,\subset\,W(E_{8})$}%
\nobreak\vskip 10 true pt\nobreak%
\baselineskip=\tablebaseskip%
\tabskip=\incltabskip\halign to \hsize{%
\bol\hfil$#$&%
$\,\subset\,#$&&%
$\,\subset\,#$\cr%
{({1^{9}})}&
{{1^{8}}}\hskip\inclskip\quad\hfil%
{({1}{2}{3^{2}})}&
{{2}{3^{2}}}\hskip\inclskip\quad\hfil%
{({2^{2}}{5})}&
{-({1}{2}{5^{-1}}{(10)})}\hskip\inclskip\hskip-\inclskip\hfil\cr\rg%
{({1^{7}}{2})}&
{{1^{6}}{2}}\hskip\inclskip\quad\hfil%
{({3^{3}})}&
{-({1}{2^{-1}}{3^{-3}}{6^{3}})}\hskip\inclskip\quad\hfil%
{({1}{3}{5})}&
{{3}{5}}\hskip\inclskip\hskip-\inclskip\hfil\cr\rg%
{({1^{5}}{2^{2}})}&
{{1^{4}}{2^{2}}}\hskip\inclskip\quad\hfil%
{({1^{5}}{4})}&
{{1^{4}}{4}}\hskip\inclskip\quad\hfil%
{({4}{5})}&
{-({1}{2^{-1}}{4}{5^{-1}}{(10)})}\hskip\inclskip\hskip-\inclskip\hfil\cr\rg%
{({1^{3}}{2^{3}})}&
{{1^{2}}{2^{3}}}\hskip\inclskip\quad\hfil%
{({1^{3}}{2}{4})}&
{{1^{2}}{2}{4}}\hskip\inclskip\quad\hfil%
{({1^{3}}{6})}&
{{1^{2}}{6}}\hskip\inclskip\hskip-\inclskip\hfil\cr\rg%
{({1}{2^{4}})}&
{{2^{4}}_{\llap{$\scriptstyle{u}$}{}}}\hskip\inclskip\quad\hfil%
{({1}{2^{2}}{4})}&
{{2^{2}}{4}_{u}}\hskip\inclskip\quad\hfil%
{({1}{2}{6})}&
{{2}{6}_{u}}\hskip\inclskip\hskip-\inclskip\hfil\cr\rg%
{({1^{6}}{3})}&
{{1^{5}}{3}}\hskip\inclskip\quad\hfil%
{({1^{2}}{3}{4})}&
{{1}{3}{4}}\hskip\inclskip\quad\hfil%
{({3}{6})}&
{-({1}{2^{-1}}{3^{-1}}{6^{2}})}\hskip\inclskip\hskip-\inclskip\hfil\cr\rg%
{({1^{4}}{2}{3})}&
{{1^{3}}{2}{3}}\hskip\inclskip\quad\hfil%
{({2}{3}{4})}&
{-({1}{3^{-1}}{4}{6})}\hskip\inclskip\quad\hfil%
{({1^{2}}{7})}&
{{1}{7}}\hskip\inclskip\hskip-\inclskip\hfil\cr\rg%
{({1^{2}}{2^{2}}{3})}&
{{1}{2^{2}}{3}}\hskip\inclskip\quad\hfil%
{({1}{4^{2}})}&
{{4^{2}}_{\llap{$\scriptstyle{u}$}{}}}\hskip\inclskip\quad\hfil%
{({2}{7})}&
{-({1}{7^{-1}}{(14)})}\hskip\inclskip\hskip-\inclskip\hfil\cr\rg%
{({2^{3}}{3})}&
{-({1}{2^{2}}{3^{-1}}{6})}\hskip\inclskip\quad\hfil%
{({1^{4}}{5})}&
{{1^{3}}{5}}\hskip\inclskip\quad\hfil%
{({1}{8})}&
{{8}_{u}}\hskip\inclskip\hskip-\inclskip\hfil\cr\rg%
{({1^{3}}{3^{2}})}&
{{1^{2}}{3^{2}}}\hskip\inclskip\quad\hfil%
{({1^{2}}{2}{5})}&
{{1}{2}{5}}\hskip\inclskip\quad\hfil%
{({9})}&
{-({1}{2^{-1}}{9^{-1}}{(18)})}\hskip\inclskip\hskip-\inclskip\hfil\cr\rg%
}}%
\tableclose%

\vskip 40 true pt

\inclskip=30true pt
\vbox{%
\tableopen{Class inclusions for $W({A_{7}}{A_{1}})\,\subset\,W(E_{8})$}%
\nobreak\vskip 10 true pt\nobreak%
\baselineskip=\tablebaseskip%
\tabskip=\incltabskip\halign to \hsize{%
\bol\hfil$#$&%
$\,\subset\,#$&&%
$\,\subset\,#$\cr%
{({1^{8}}){\times}({1^{2}})}&
{{1^{8}}}\hskip\inclskip\quad\hfil%
{({1^{3}}{5}){\times}({1^{2}})}&
{{1^{3}}{5}}\hskip\inclskip\quad\hfil%
{({1^{2}}{3^{2}}){\times}({2})}&
{{2}{3^{2}}}\hskip\inclskip\hskip-\inclskip\hfil\cr\rg%
{({1^{6}}{2}){\times}({1^{2}})}&
{{1^{6}}{2}}\hskip\inclskip\quad\hfil%
{({1}{2}{5}){\times}({1^{2}})}&
{{1}{2}{5}}\hskip\inclskip\quad\hfil%
{({2}{3^{2}}){\times}({2})}&
{-({1^{2}}{3^{-2}}{6^{2}})}\hskip\inclskip\hskip-\inclskip\hfil\cr\rg%
{({1^{4}}{2^{2}}){\times}({1^{2}})}&
{{1^{4}}{2^{2}}}\hskip\inclskip\quad\hfil%
{({3}{5}){\times}({1^{2}})}&
{{3}{5}}\hskip\inclskip\quad\hfil%
{({1^{4}}{4}){\times}({2})}&
{{1^{2}}{2}{4}}\hskip\inclskip\hskip-\inclskip\hfil\cr\rg%
{({1^{2}}{2^{3}}){\times}({1^{2}})}&
{{1^{2}}{2^{3}}}\hskip\inclskip\quad\hfil%
{({1^{2}}{6}){\times}({1^{2}})}&
{{1^{2}}{6}}\hskip\inclskip\quad\hfil%
{({1^{2}}{2}{4}){\times}({2})}&
{{2^{2}}{4}_{u}}\hskip\inclskip\hskip-\inclskip\hfil\cr\rg%
{({2^{4}}){\times}({1^{2}})}&
{{2^{4}}_{\llap{$\scriptstyle{v}$}{}}}\hskip\inclskip\quad\hfil%
{({2}{6}){\times}({1^{2}})}&
{{2}{6}_{v}}\hskip\inclskip\quad\hfil%
{({2^{2}}{4}){\times}({2})}&
{-({1^{2}}{2}{4})}\hskip\inclskip\hskip-\inclskip\hfil\cr\rg%
{({1^{5}}{3}){\times}({1^{2}})}&
{{1^{5}}{3}}\hskip\inclskip\quad\hfil%
{({1}{7}){\times}({1^{2}})}&
{{1}{7}}\hskip\inclskip\quad\hfil%
{({1}{3}{4}){\times}({2})}&
{-({1}{3^{-1}}{4}{6})}\hskip\inclskip\hskip-\inclskip\hfil\cr\rg%
{({1^{3}}{2}{3}){\times}({1^{2}})}&
{{1^{3}}{2}{3}}\hskip\inclskip\quad\hfil%
{({8}){\times}({1^{2}})}&
{{8}_{v}}\hskip\inclskip\quad\hfil%
{({4^{2}}){\times}({2})}&
{-({1^{2}}{2^{-1}}{4^{2}})}\hskip\inclskip\hskip-\inclskip\hfil\cr\rg%
{({1}{2^{2}}{3}){\times}({1^{2}})}&
{{1}{2^{2}}{3}}\hskip\inclskip\quad\hfil%
{({1^{8}}){\times}({2})}&
{{1^{6}}{2}}\hskip\inclskip\quad\hfil%
{({1^{3}}{5}){\times}({2})}&
{{1}{2}{5}}\hskip\inclskip\hskip-\inclskip\hfil\cr\rg%
{({1^{2}}{3^{2}}){\times}({1^{2}})}&
{{1^{2}}{3^{2}}}\hskip\inclskip\quad\hfil%
{({1^{6}}{2}){\times}({2})}&
{{1^{4}}{2^{2}}}\hskip\inclskip\quad\hfil%
{({1}{2}{5}){\times}({2})}&
{-({1}{2}{5^{-1}}{(10)})}\hskip\inclskip\hskip-\inclskip\hfil\cr\rg%
{({2}{3^{2}}){\times}({1^{2}})}&
{{2}{3^{2}}}\hskip\inclskip\quad\hfil%
{({1^{4}}{2^{2}}){\times}({2})}&
{{1^{2}}{2^{3}}}\hskip\inclskip\quad\hfil%
{({3}{5}){\times}({2})}&
{-({1^{2}}{2^{-1}}{3^{-1}}{5^{-1}}{6}{(10)})}\hskip\inclskip\hskip-\inclskip\hfil\cr\rg%
{({1^{4}}{4}){\times}({1^{2}})}&
{{1^{4}}{4}}\hskip\inclskip\quad\hfil%
{({1^{2}}{2^{3}}){\times}({2})}&
{{2^{4}}_{\llap{$\scriptstyle{u}$}{}}}\hskip\inclskip\quad\hfil%
{({1^{2}}{6}){\times}({2})}&
{{2}{6}_{u}}\hskip\inclskip\hskip-\inclskip\hfil\cr\rg%
{({1^{2}}{2}{4}){\times}({1^{2}})}&
{{1^{2}}{2}{4}}\hskip\inclskip\quad\hfil%
{({2^{4}}){\times}({2})}&
{-({1^{2}}{2^{3}})}\hskip\inclskip\quad\hfil%
{({2}{6}){\times}({2})}&
{-({1^{2}}{6})}\hskip\inclskip\hskip-\inclskip\hfil\cr\rg%
{({2^{2}}{4}){\times}({1^{2}})}&
{{2^{2}}{4}_{v}}\hskip\inclskip\quad\hfil%
{({1^{5}}{3}){\times}({2})}&
{{1^{3}}{2}{3}}\hskip\inclskip\quad\hfil%
{({1}{7}){\times}({2})}&
{-({1}{7^{-1}}{(14)})}\hskip\inclskip\hskip-\inclskip\hfil\cr\rg%
{({1}{3}{4}){\times}({1^{2}})}&
{{1}{3}{4}}\hskip\inclskip\quad\hfil%
{({1^{3}}{2}{3}){\times}({2})}&
{{1}{2^{2}}{3}}\hskip\inclskip\quad\hfil%
{({8}){\times}({2})}&
{-({1^{2}}{2^{-1}}{8})}\hskip\inclskip\hskip-\inclskip\hfil\cr\rg%
{({4^{2}}){\times}({1^{2}})}&
{{4^{2}}_{\llap{$\scriptstyle{v}$}{}}}\hskip\inclskip\quad\hfil%
{({1}{2^{2}}{3}){\times}({2})}&
{-({1}{2^{2}}{3^{-1}}{6})}\hskip\inclskip\quad\hfil%
&
\omit\hskip\inclskip\hfil\hskip-\inclskip\hfil\cr\rg%
}}%
\tableclose%
\eject
\def\rg{\rgloc}

\def\bol{\qquad}
\def\rg{\noalign{\nobreak\vskip 0.5 true pt plus 2 true pt}}
\inclskip=70 true pt
\vbox{%
\tableopen{Class inclusions for $W({A_{5}}{A_{2}}{A_{1}})\,\subset\,W(E_{8})$}%
\nobreak\vskip 10 true pt\nobreak%
\baselineskip=\tablebaseskip%
\tabskip=\incltabskip\halign to \hsize{%
\bol\hfil$#$&%
$\,\subset\,#$&&%
$\,\subset\,#$\cr%
{({1^{6}}){\times}({1^{3}}){\times}({1^{2}})}&
{{1^{8}}}\hskip\inclskip\quad\hfil%
{({1^{6}}){\times}({1^{3}}){\times}({2})}&
{{1^{6}}{2}}\hskip\inclskip\hskip-\inclskip\hfil\cr\rg%
{({1^{4}}{2}){\times}({1^{3}}){\times}({1^{2}})}&
{{1^{6}}{2}}\hskip\inclskip\quad\hfil%
{({1^{4}}{2}){\times}({1^{3}}){\times}({2})}&
{{1^{4}}{2^{2}}}\hskip\inclskip\hskip-\inclskip\hfil\cr\rg%
{({1^{2}}{2^{2}}){\times}({1^{3}}){\times}({1^{2}})}&
{{1^{4}}{2^{2}}}\hskip\inclskip\quad\hfil%
{({1^{2}}{2^{2}}){\times}({1^{3}}){\times}({2})}&
{{1^{2}}{2^{3}}}\hskip\inclskip\hskip-\inclskip\hfil\cr\rg%
{({2^{3}}){\times}({1^{3}}){\times}({1^{2}})}&
{{1^{2}}{2^{3}}}\hskip\inclskip\quad\hfil%
{({2^{3}}){\times}({1^{3}}){\times}({2})}&
{{2^{4}}_{\llap{$\scriptstyle{v}$}{}}}\hskip\inclskip\hskip-\inclskip\hfil\cr\rg%
{({1^{3}}{3}){\times}({1^{3}}){\times}({1^{2}})}&
{{1^{5}}{3}}\hskip\inclskip\quad\hfil%
{({1^{3}}{3}){\times}({1^{3}}){\times}({2})}&
{{1^{3}}{2}{3}}\hskip\inclskip\hskip-\inclskip\hfil\cr\rg%
{({1}{2}{3}){\times}({1^{3}}){\times}({1^{2}})}&
{{1^{3}}{2}{3}}\hskip\inclskip\quad\hfil%
{({1}{2}{3}){\times}({1^{3}}){\times}({2})}&
{{1}{2^{2}}{3}}\hskip\inclskip\hskip-\inclskip\hfil\cr\rg%
{({3^{2}}){\times}({1^{3}}){\times}({1^{2}})}&
{{1^{2}}{3^{2}}}\hskip\inclskip\quad\hfil%
{({3^{2}}){\times}({1^{3}}){\times}({2})}&
{{2}{3^{2}}}\hskip\inclskip\hskip-\inclskip\hfil\cr\rg%
{({1^{2}}{4}){\times}({1^{3}}){\times}({1^{2}})}&
{{1^{4}}{4}}\hskip\inclskip\quad\hfil%
{({1^{2}}{4}){\times}({1^{3}}){\times}({2})}&
{{1^{2}}{2}{4}}\hskip\inclskip\hskip-\inclskip\hfil\cr\rg%
{({2}{4}){\times}({1^{3}}){\times}({1^{2}})}&
{{1^{2}}{2}{4}}\hskip\inclskip\quad\hfil%
{({2}{4}){\times}({1^{3}}){\times}({2})}&
{{2^{2}}{4}_{v}}\hskip\inclskip\hskip-\inclskip\hfil\cr\rg%
{({1}{5}){\times}({1^{3}}){\times}({1^{2}})}&
{{1^{3}}{5}}\hskip\inclskip\quad\hfil%
{({1}{5}){\times}({1^{3}}){\times}({2})}&
{{1}{2}{5}}\hskip\inclskip\hskip-\inclskip\hfil\cr\rg%
{({6}){\times}({1^{3}}){\times}({1^{2}})}&
{{1^{2}}{6}}\hskip\inclskip\quad\hfil%
{({6}){\times}({1^{3}}){\times}({2})}&
{{2}{6}_{v}}\hskip\inclskip\hskip-\inclskip\hfil\cr\rg%
{({1^{6}}){\times}({1}{2}){\times}({1^{2}})}&
{{1^{6}}{2}}\hskip\inclskip\quad\hfil%
{({1^{6}}){\times}({1}{2}){\times}({2})}&
{{1^{4}}{2^{2}}}\hskip\inclskip\hskip-\inclskip\hfil\cr\rg%
{({1^{4}}{2}){\times}({1}{2}){\times}({1^{2}})}&
{{1^{4}}{2^{2}}}\hskip\inclskip\quad\hfil%
{({1^{4}}{2}){\times}({1}{2}){\times}({2})}&
{{1^{2}}{2^{3}}}\hskip\inclskip\hskip-\inclskip\hfil\cr\rg%
{({1^{2}}{2^{2}}){\times}({1}{2}){\times}({1^{2}})}&
{{1^{2}}{2^{3}}}\hskip\inclskip\quad\hfil%
{({1^{2}}{2^{2}}){\times}({1}{2}){\times}({2})}&
{{2^{4}}_{\llap{$\scriptstyle{u}$}{}}}\hskip\inclskip\hskip-\inclskip\hfil\cr\rg%
{({2^{3}}){\times}({1}{2}){\times}({1^{2}})}&
{{2^{4}}_{\llap{$\scriptstyle{u}$}{}}}\hskip\inclskip\quad\hfil%
{({2^{3}}){\times}({1}{2}){\times}({2})}&
{-({1^{2}}{2^{3}})}\hskip\inclskip\hskip-\inclskip\hfil\cr\rg%
{({1^{3}}{3}){\times}({1}{2}){\times}({1^{2}})}&
{{1^{3}}{2}{3}}\hskip\inclskip\quad\hfil%
{({1^{3}}{3}){\times}({1}{2}){\times}({2})}&
{{1}{2^{2}}{3}}\hskip\inclskip\hskip-\inclskip\hfil\cr\rg%
{({1}{2}{3}){\times}({1}{2}){\times}({1^{2}})}&
{{1}{2^{2}}{3}}\hskip\inclskip\quad\hfil%
{({1}{2}{3}){\times}({1}{2}){\times}({2})}&
{-({1}{2^{2}}{3^{-1}}{6})}\hskip\inclskip\hskip-\inclskip\hfil\cr\rg%
{({3^{2}}){\times}({1}{2}){\times}({1^{2}})}&
{{2}{3^{2}}}\hskip\inclskip\quad\hfil%
{({3^{2}}){\times}({1}{2}){\times}({2})}&
{-({1^{2}}{3^{-2}}{6^{2}})}\hskip\inclskip\hskip-\inclskip\hfil\cr\rg%
{({1^{2}}{4}){\times}({1}{2}){\times}({1^{2}})}&
{{1^{2}}{2}{4}}\hskip\inclskip\quad\hfil%
{({1^{2}}{4}){\times}({1}{2}){\times}({2})}&
{{2^{2}}{4}_{u}}\hskip\inclskip\hskip-\inclskip\hfil\cr\rg%
{({2}{4}){\times}({1}{2}){\times}({1^{2}})}&
{{2^{2}}{4}_{u}}\hskip\inclskip\quad\hfil%
{({2}{4}){\times}({1}{2}){\times}({2})}&
{-({1^{2}}{2}{4})}\hskip\inclskip\hskip-\inclskip\hfil\cr\rg%
{({1}{5}){\times}({1}{2}){\times}({1^{2}})}&
{{1}{2}{5}}\hskip\inclskip\quad\hfil%
{({1}{5}){\times}({1}{2}){\times}({2})}&
{-({1}{2}{5^{-1}}{(10)})}\hskip\inclskip\hskip-\inclskip\hfil\cr\rg%
{({6}){\times}({1}{2}){\times}({1^{2}})}&
{{2}{6}_{u}}\hskip\inclskip\quad\hfil%
{({6}){\times}({1}{2}){\times}({2})}&
{-({1^{2}}{6})}\hskip\inclskip\hskip-\inclskip\hfil\cr\rg%
{({1^{6}}){\times}({3}){\times}({1^{2}})}&
{{1^{5}}{3}}\hskip\inclskip\quad\hfil%
{({1^{6}}){\times}({3}){\times}({2})}&
{{1^{3}}{2}{3}}\hskip\inclskip\hskip-\inclskip\hfil\cr\rg%
{({1^{4}}{2}){\times}({3}){\times}({1^{2}})}&
{{1^{3}}{2}{3}}\hskip\inclskip\quad\hfil%
{({1^{4}}{2}){\times}({3}){\times}({2})}&
{{1}{2^{2}}{3}}\hskip\inclskip\hskip-\inclskip\hfil\cr\rg%
{({1^{2}}{2^{2}}){\times}({3}){\times}({1^{2}})}&
{{1}{2^{2}}{3}}\hskip\inclskip\quad\hfil%
{({1^{2}}{2^{2}}){\times}({3}){\times}({2})}&
{-({1}{2^{2}}{3^{-1}}{6})}\hskip\inclskip\hskip-\inclskip\hfil\cr\rg%
{({2^{3}}){\times}({3}){\times}({1^{2}})}&
{-({1}{2^{2}}{3^{-1}}{6})}\hskip\inclskip\quad\hfil%
{({2^{3}}){\times}({3}){\times}({2})}&
{-({1^{3}}{2}{3^{-1}}{6})}\hskip\inclskip\hskip-\inclskip\hfil\cr\rg%
{({1^{3}}{3}){\times}({3}){\times}({1^{2}})}&
{{1^{2}}{3^{2}}}\hskip\inclskip\quad\hfil%
{({1^{3}}{3}){\times}({3}){\times}({2})}&
{{2}{3^{2}}}\hskip\inclskip\hskip-\inclskip\hfil\cr\rg%
{({1}{2}{3}){\times}({3}){\times}({1^{2}})}&
{{2}{3^{2}}}\hskip\inclskip\quad\hfil%
{({1}{2}{3}){\times}({3}){\times}({2})}&
{-({1^{2}}{3^{-2}}{6^{2}})}\hskip\inclskip\hskip-\inclskip\hfil\cr\rg%
{({3^{2}}){\times}({3}){\times}({1^{2}})}&
{-({1}{2^{-1}}{3^{-3}}{6^{3}})}\hskip\inclskip\quad\hfil%
{({3^{2}}){\times}({3}){\times}({2})}&
{-({1^{3}}{2^{-2}}{3^{-3}}{6^{3}})}\hskip\inclskip\hskip-\inclskip\hfil\cr\rg%
{({1^{2}}{4}){\times}({3}){\times}({1^{2}})}&
{{1}{3}{4}}\hskip\inclskip\quad\hfil%
{({1^{2}}{4}){\times}({3}){\times}({2})}&
{-({1}{3^{-1}}{4}{6})}\hskip\inclskip\hskip-\inclskip\hfil\cr\rg%
{({2}{4}){\times}({3}){\times}({1^{2}})}&
{-({1}{3^{-1}}{4}{6})}\hskip\inclskip\quad\hfil%
{({2}{4}){\times}({3}){\times}({2})}&
{-({1^{3}}{2^{-1}}{3^{-1}}{4}{6})}\hskip\inclskip\hskip-\inclskip\hfil\cr\rg%
{({1}{5}){\times}({3}){\times}({1^{2}})}&
{{3}{5}}\hskip\inclskip\quad\hfil%
{({1}{5}){\times}({3}){\times}({2})}&
{-({1^{2}}{2^{-1}}{3^{-1}}{5^{-1}}{6}{(10)})}\hskip\inclskip\hskip-\inclskip\hfil\cr\rg%
{({6}){\times}({3}){\times}({1^{2}})}&
{-({1}{2^{-1}}{3^{-1}}{6^{2}})}\hskip\inclskip\quad\hfil%
{({6}){\times}({3}){\times}({2})}&
{-({1^{3}}{2^{-2}}{3^{-1}}{6^{2}})}\hskip\inclskip\hskip-\inclskip\hfil\cr\rg%
}}%
\tableclose%
\eject
\def\rg{\rgloc}

\def\bol{\quad}
\inclskip=10 true pt
\def\rg{\noalign{\nobreak\vskip 5 true pt plus 2 true pt}}
\vbox{%
\tableopen{Class inclusions for $W({A_{4}}{A_{4}})\,\subset\,W(E_{8})$}%
\nobreak\vskip 10 true pt\nobreak%
\baselineskip=\tablebaseskip%
\tabskip=\incltabskip\halign to \hsize{%
\bol\hfil$#$&%
$\,\subset\,#$&&%
$\,\subset\,#$\cr%
{({1^{5}}){\times}({1^{5}})}&
{{1^{8}}}\hskip\inclskip\quad\hfil%
{({1^{2}}{3}){\times}({1}{2^{2}})}&
{{1}{2^{2}}{3}}\hskip\inclskip\quad\hfil%
{({1}{4}){\times}({2}{3})}&
{-({1}{3^{-1}}{4}{6})}\hskip\inclskip\hskip-\inclskip\hfil\cr\rg%
{({1^{3}}{2}){\times}({1^{5}})}&
{{1^{6}}{2}}\hskip\inclskip\quad\hfil%
{({2}{3}){\times}({1}{2^{2}})}&
{-({1}{2^{2}}{3^{-1}}{6})}\hskip\inclskip\quad\hfil%
{({5}){\times}({2}{3})}&
{-({1^{2}}{2^{-1}}{3^{-1}}{5^{-1}}{6}{(10)})}\hskip\inclskip\hskip-\inclskip\hfil\cr\rg%
{({1}{2^{2}}){\times}({1^{5}})}&
{{1^{4}}{2^{2}}}\hskip\inclskip\quad\hfil%
{({1}{4}){\times}({1}{2^{2}})}&
{{2^{2}}{4}_{u}}\hskip\inclskip\quad\hfil%
{({1^{5}}){\times}({1}{4})}&
{{1^{4}}{4}}\hskip\inclskip\hskip-\inclskip\hfil\cr\rg%
{({1^{2}}{3}){\times}({1^{5}})}&
{{1^{5}}{3}}\hskip\inclskip\quad\hfil%
{({5}){\times}({1}{2^{2}})}&
{-({1}{2}{5^{-1}}{(10)})}\hskip\inclskip\quad\hfil%
{({1^{3}}{2}){\times}({1}{4})}&
{{1^{2}}{2}{4}}\hskip\inclskip\hskip-\inclskip\hfil\cr\rg%
{({2}{3}){\times}({1^{5}})}&
{{1^{3}}{2}{3}}\hskip\inclskip\quad\hfil%
{({1^{5}}){\times}({1^{2}}{3})}&
{{1^{5}}{3}}\hskip\inclskip\quad\hfil%
{({1}{2^{2}}){\times}({1}{4})}&
{{2^{2}}{4}_{u}}\hskip\inclskip\hskip-\inclskip\hfil\cr\rg%
{({1}{4}){\times}({1^{5}})}&
{{1^{4}}{4}}\hskip\inclskip\quad\hfil%
{({1^{3}}{2}){\times}({1^{2}}{3})}&
{{1^{3}}{2}{3}}\hskip\inclskip\quad\hfil%
{({1^{2}}{3}){\times}({1}{4})}&
{{1}{3}{4}}\hskip\inclskip\hskip-\inclskip\hfil\cr\rg%
{({5}){\times}({1^{5}})}&
{{1^{3}}{5}}\hskip\inclskip\quad\hfil%
{({1}{2^{2}}){\times}({1^{2}}{3})}&
{{1}{2^{2}}{3}}\hskip\inclskip\quad\hfil%
{({2}{3}){\times}({1}{4})}&
{-({1}{3^{-1}}{4}{6})}\hskip\inclskip\hskip-\inclskip\hfil\cr\rg%
{({1^{5}}){\times}({1^{3}}{2})}&
{{1^{6}}{2}}\hskip\inclskip\quad\hfil%
{({1^{2}}{3}){\times}({1^{2}}{3})}&
{{1^{2}}{3^{2}}}\hskip\inclskip\quad\hfil%
{({1}{4}){\times}({1}{4})}&
{{4^{2}}_{\llap{$\scriptstyle{u}$}{}}}\hskip\inclskip\hskip-\inclskip\hfil\cr\rg%
{({1^{3}}{2}){\times}({1^{3}}{2})}&
{{1^{4}}{2^{2}}}\hskip\inclskip\quad\hfil%
{({2}{3}){\times}({1^{2}}{3})}&
{{2}{3^{2}}}\hskip\inclskip\quad\hfil%
{({5}){\times}({1}{4})}&
{-({1}{2^{-1}}{4}{5^{-1}}{(10)})}\hskip\inclskip\hskip-\inclskip\hfil\cr\rg%
{({1}{2^{2}}){\times}({1^{3}}{2})}&
{{1^{2}}{2^{3}}}\hskip\inclskip\quad\hfil%
{({1}{4}){\times}({1^{2}}{3})}&
{{1}{3}{4}}\hskip\inclskip\quad\hfil%
{({1^{5}}){\times}({5})}&
{{1^{3}}{5}}\hskip\inclskip\hskip-\inclskip\hfil\cr\rg%
{({1^{2}}{3}){\times}({1^{3}}{2})}&
{{1^{3}}{2}{3}}\hskip\inclskip\quad\hfil%
{({5}){\times}({1^{2}}{3})}&
{{3}{5}}\hskip\inclskip\quad\hfil%
{({1^{3}}{2}){\times}({5})}&
{{1}{2}{5}}\hskip\inclskip\hskip-\inclskip\hfil\cr\rg%
{({2}{3}){\times}({1^{3}}{2})}&
{{1}{2^{2}}{3}}\hskip\inclskip\quad\hfil%
{({1^{5}}){\times}({2}{3})}&
{{1^{3}}{2}{3}}\hskip\inclskip\quad\hfil%
{({1}{2^{2}}){\times}({5})}&
{-({1}{2}{5^{-1}}{(10)})}\hskip\inclskip\hskip-\inclskip\hfil\cr\rg%
{({1}{4}){\times}({1^{3}}{2})}&
{{1^{2}}{2}{4}}\hskip\inclskip\quad\hfil%
{({1^{3}}{2}){\times}({2}{3})}&
{{1}{2^{2}}{3}}\hskip\inclskip\quad\hfil%
{({1^{2}}{3}){\times}({5})}&
{{3}{5}}\hskip\inclskip\hskip-\inclskip\hfil\cr\rg%
{({5}){\times}({1^{3}}{2})}&
{{1}{2}{5}}\hskip\inclskip\quad\hfil%
{({1}{2^{2}}){\times}({2}{3})}&
{-({1}{2^{2}}{3^{-1}}{6})}\hskip\inclskip\quad\hfil%
{({2}{3}){\times}({5})}&
{-({1^{2}}{2^{-1}}{3^{-1}}{5^{-1}}{6}{(10)})}\hskip\inclskip\hskip-\inclskip\hfil\cr\rg%
{({1^{5}}){\times}({1}{2^{2}})}&
{{1^{4}}{2^{2}}}\hskip\inclskip\quad\hfil%
{({1^{2}}{3}){\times}({2}{3})}&
{{2}{3^{2}}}\hskip\inclskip\quad\hfil%
{({1}{4}){\times}({5})}&
{-({1}{2^{-1}}{4}{5^{-1}}{(10)})}\hskip\inclskip\hskip-\inclskip\hfil\cr\rg%
{({1^{3}}{2}){\times}({1}{2^{2}})}&
{{1^{2}}{2^{3}}}\hskip\inclskip\quad\hfil%
{({2}{3}){\times}({2}{3})}&
{-({1^{2}}{3^{-2}}{6^{2}})}\hskip\inclskip\quad\hfil%
{({5}){\times}({5})}&
{-({1^{2}}{2^{-2}}{5^{-2}}{(10)^{2}})}\hskip\inclskip\hskip-\inclskip\hfil\cr\rg%
{({1}{2^{2}}){\times}({1}{2^{2}})}&
{{2^{4}}_{\llap{$\scriptstyle{u}$}{}}}\hskip\inclskip\quad\hfil%
&
\omit\hskip\inclskip\hfil\quad\hfil%
&
\omit\hskip\inclskip\hfil\hskip-\inclskip\hfil\cr\rg%
}}%
\tableclose%
\eject
\def\rg{\rgloc}

\def\bol{\hskip 0.2em \relax}
\inclskip=1 true pt
\def\rg{\noalign{\nobreak\vskip 2 true pt plus 1 true pt}}
\vbox{%
\tableopen{Class inclusions for $W({D_{5}}{A_{3}})\,\subset\,W(E_{8})$}%
\nobreak\vskip 10 true pt\nobreak%
\baselineskip=\tablebaseskip%
\tabskip=\incltabskip\halign to \hsize{%
\bol\hfil$#$&%
$\,\subset\,#$&&%
$\,\subset\,#$\cr%
{({1^{5}}){\times}({1^{4}})}&
{{1^{8}}}\hskip\inclskip\quad\hfil%
{({1}{\ol{1}^{\dn{4}}}){\times}({1^{2}}{2})}&
{-({1^{2}}{2^{3}})}\hskip\inclskip\quad\hfil%
{({5}){\times}({1}{3})}&
{{3}{5}}\hskip\inclskip\hskip-\inclskip\hfil\cr\rg%
{({1^{3}}{2}){\times}({1^{4}})}&
{{1^{6}}{2}}\hskip\inclskip\quad\hfil%
{({1}{\ol{2}^{\dn{2}}}){\times}({1^{2}}{2})}&
{{1^{2}}{2^{-1}}{4^{2}}}\hskip\inclskip\quad\hfil%
{({1^{3}}{\ol{1}^{\dn{2}}}){\times}({1}{3})}&
{{1}{2^{2}}{3}}\hskip\inclskip\hskip-\inclskip\hfil\cr\rg%
{({1}{2^{2}}){\times}({1^{4}})}&
{{1^{4}}{2^{2}}}\hskip\inclskip\quad\hfil%
{({1}{\ol{1}}{\ol{3}}){\times}({1^{2}}{2})}&
{{1}{2^{2}}{3^{-1}}{6}}\hskip\inclskip\quad\hfil%
{({1}{\ol{1}^{\dn{2}}}{2}){\times}({1}{3})}&
{-({1}{2^{2}}{3^{-1}}{6})}\hskip\inclskip\hskip-\inclskip\hfil\cr\rg%
{({1^{2}}{3}){\times}({1^{4}})}&
{{1^{5}}{3}}\hskip\inclskip\quad\hfil%
{({\ol{1}^{\dn{3}}}{\ol{2}}){\times}({1^{2}}{2})}&
{-({1^{2}}{2}{4})}\hskip\inclskip\quad\hfil%
{({\ol{1}^{\dn{2}}}{3}){\times}({1}{3})}&
{-({1^{2}}{3^{-2}}{6^{2}})}\hskip\inclskip\hskip-\inclskip\hfil\cr\rg%
{({2}{3}){\times}({1^{4}})}&
{{1^{3}}{2}{3}}\hskip\inclskip\quad\hfil%
{({\ol{2}}{\ol{3}}){\times}({1^{2}}{2})}&
{{1}{3^{-1}}{4}{6}}\hskip\inclskip\quad\hfil%
{({1^{2}}{\ol{1}}{\ol{2}}){\times}({1}{3})}&
{{1}{3}{4}}\hskip\inclskip\hskip-\inclskip\hfil\cr\rg%
{({1}{4}){\times}({1^{4}})}&
{{1^{4}}{4}}\hskip\inclskip\quad\hfil%
{({\ol{1}}{\ol{4}}){\times}({1^{2}}{2})}&
{{2^{2}}{4^{-1}}{8}}\hskip\inclskip\quad\hfil%
{({\ol{1}}{2}{\ol{2}}){\times}({1}{3})}&
{-({1}{3^{-1}}{4}{6})}\hskip\inclskip\hskip-\inclskip\hfil\cr\rg%
{({5}){\times}({1^{4}})}&
{{1^{3}}{5}}\hskip\inclskip\quad\hfil%
{({1^{5}}){\times}({2^{2}})}&
{{1^{4}}{2^{2}}}\hskip\inclskip\quad\hfil%
{({1}{\ol{1}^{\dn{4}}}){\times}({1}{3})}&
{-({1^{3}}{2}{3^{-1}}{6})}\hskip\inclskip\hskip-\inclskip\hfil\cr\rg%
{({1^{3}}{\ol{1}^{\dn{2}}}){\times}({1^{4}})}&
{{1^{4}}{2^{2}}}\hskip\inclskip\quad\hfil%
{({1^{3}}{2}){\times}({2^{2}})}&
{{1^{2}}{2^{3}}}\hskip\inclskip\quad\hfil%
{({1}{\ol{2}^{\dn{2}}}){\times}({1}{3})}&
{{1}{2^{-2}}{3}{4^{2}}}\hskip\inclskip\hskip-\inclskip\hfil\cr\rg%
{({1}{\ol{1}^{\dn{2}}}{2}){\times}({1^{4}})}&
{{1^{2}}{2^{3}}}\hskip\inclskip\quad\hfil%
{({1}{2^{2}}){\times}({2^{2}})}&
{{2^{4}}_{\llap{$\scriptstyle{u}$}{}}}\hskip\inclskip\quad\hfil%
{({1}{\ol{1}}{\ol{3}}){\times}({1}{3})}&
{{2}{6}_{vv}}\hskip\inclskip\hskip-\inclskip\hfil\cr\rg%
{({\ol{1}^{\dn{2}}}{3}){\times}({1^{4}})}&
{{1}{2^{2}}{3}}\hskip\inclskip\quad\hfil%
{({1^{2}}{3}){\times}({2^{2}})}&
{{1}{2^{2}}{3}}\hskip\inclskip\quad\hfil%
{({\ol{1}^{\dn{3}}}{\ol{2}}){\times}({1}{3})}&
{-({1^{3}}{2^{-1}}{3^{-1}}{4}{6})}\hskip\inclskip\hskip-\inclskip\hfil\cr\rg%
{({1^{2}}{\ol{1}}{\ol{2}}){\times}({1^{4}})}&
{{1^{4}}{4}}\hskip\inclskip\quad\hfil%
{({2}{3}){\times}({2^{2}})}&
{-({1}{2^{2}}{3^{-1}}{6})}\hskip\inclskip\quad\hfil%
{({\ol{2}}{\ol{3}}){\times}({1}{3})}&
{{2^{-1}}{4}{6}}\hskip\inclskip\hskip-\inclskip\hfil\cr\rg%
{({\ol{1}}{2}{\ol{2}}){\times}({1^{4}})}&
{{1^{2}}{2}{4}}\hskip\inclskip\quad\hfil%
{({1}{4}){\times}({2^{2}})}&
{{2^{2}}{4}_{u}}\hskip\inclskip\quad\hfil%
{({\ol{1}}{\ol{4}}){\times}({1}{3})}&
{-({1}{3^{-1}}{4^{-1}}{6}{8})}\hskip\inclskip\hskip-\inclskip\hfil\cr\rg%
{({1}{\ol{1}^{\dn{4}}}){\times}({1^{4}})}&
{{2^{4}}_{\llap{$\scriptstyle{v}$}{}}}\hskip\inclskip\quad\hfil%
{({5}){\times}({2^{2}})}&
{-({1}{2}{5^{-1}}{(10)})}\hskip\inclskip\quad\hfil%
{({1^{5}}){\times}({4})}&
{{1^{4}}{4}}\hskip\inclskip\hskip-\inclskip\hfil\cr\rg%
{({1}{\ol{2}^{\dn{2}}}){\times}({1^{4}})}&
{{1^{4}}{2^{-2}}{4^{2}}}\hskip\inclskip\quad\hfil%
{({1^{3}}{\ol{1}^{\dn{2}}}){\times}({2^{2}})}&
{{2^{4}}_{\llap{$\scriptstyle{v}$}{}}}\hskip\inclskip\quad\hfil%
{({1^{3}}{2}){\times}({4})}&
{{1^{2}}{2}{4}}\hskip\inclskip\hskip-\inclskip\hfil\cr\rg%
{({1}{\ol{1}}{\ol{3}}){\times}({1^{4}})}&
{{1^{3}}{2}{3^{-1}}{6}}\hskip\inclskip\quad\hfil%
{({1}{\ol{1}^{\dn{2}}}{2}){\times}({2^{2}})}&
{-({1^{2}}{2^{3}})}\hskip\inclskip\quad\hfil%
{({1}{2^{2}}){\times}({4})}&
{{2^{2}}{4}_{u}}\hskip\inclskip\hskip-\inclskip\hfil\cr\rg%
{({\ol{1}^{\dn{3}}}{\ol{2}}){\times}({1^{4}})}&
{{2^{2}}{4}_{v}}\hskip\inclskip\quad\hfil%
{({\ol{1}^{\dn{2}}}{3}){\times}({2^{2}})}&
{-({1^{3}}{2}{3^{-1}}{6})}\hskip\inclskip\quad\hfil%
{({1^{2}}{3}){\times}({4})}&
{{1}{3}{4}}\hskip\inclskip\hskip-\inclskip\hfil\cr\rg%
{({\ol{2}}{\ol{3}}){\times}({1^{4}})}&
{{1^{3}}{2^{-1}}{3^{-1}}{4}{6}}\hskip\inclskip\quad\hfil%
{({1^{2}}{\ol{1}}{\ol{2}}){\times}({2^{2}})}&
{{2^{2}}{4}_{v}}\hskip\inclskip\quad\hfil%
{({2}{3}){\times}({4})}&
{-({1}{3^{-1}}{4}{6})}\hskip\inclskip\hskip-\inclskip\hfil\cr\rg%
{({\ol{1}}{\ol{4}}){\times}({1^{4}})}&
{{1^{2}}{2}{4^{-1}}{8}}\hskip\inclskip\quad\hfil%
{({\ol{1}}{2}{\ol{2}}){\times}({2^{2}})}&
{-({1^{2}}{2}{4})}\hskip\inclskip\quad\hfil%
{({1}{4}){\times}({4})}&
{{4^{2}}_{\llap{$\scriptstyle{u}$}{}}}\hskip\inclskip\hskip-\inclskip\hfil\cr\rg%
{({1^{5}}){\times}({1^{2}}{2})}&
{{1^{6}}{2}}\hskip\inclskip\quad\hfil%
{({1}{\ol{1}^{\dn{4}}}){\times}({2^{2}})}&
{-({1^{4}}{2^{2}})}\hskip\inclskip\quad\hfil%
{({5}){\times}({4})}&
{-({1}{2^{-1}}{4}{5^{-1}}{(10)})}\hskip\inclskip\hskip-\inclskip\hfil\cr\rg%
{({1^{3}}{2}){\times}({1^{2}}{2})}&
{{1^{4}}{2^{2}}}\hskip\inclskip\quad\hfil%
{({1}{\ol{2}^{\dn{2}}}){\times}({2^{2}})}&
{{4^{2}}_{\llap{$\scriptstyle{v}$}{}}}\hskip\inclskip\quad\hfil%
{({1^{3}}{\ol{1}^{\dn{2}}}){\times}({4})}&
{{2^{2}}{4}_{v}}\hskip\inclskip\hskip-\inclskip\hfil\cr\rg%
{({1}{2^{2}}){\times}({1^{2}}{2})}&
{{1^{2}}{2^{3}}}\hskip\inclskip\quad\hfil%
{({1}{\ol{1}}{\ol{3}}){\times}({2^{2}})}&
{-({1}{2^{2}}{3})}\hskip\inclskip\quad\hfil%
{({1}{\ol{1}^{\dn{2}}}{2}){\times}({4})}&
{-({1^{2}}{2}{4})}\hskip\inclskip\hskip-\inclskip\hfil\cr\rg%
{({1^{2}}{3}){\times}({1^{2}}{2})}&
{{1^{3}}{2}{3}}\hskip\inclskip\quad\hfil%
{({\ol{1}^{\dn{3}}}{\ol{2}}){\times}({2^{2}})}&
{-({1^{4}}{4})}\hskip\inclskip\quad\hfil%
{({\ol{1}^{\dn{2}}}{3}){\times}({4})}&
{-({1^{3}}{2^{-1}}{3^{-1}}{4}{6})}\hskip\inclskip\hskip-\inclskip\hfil\cr\rg%
{({2}{3}){\times}({1^{2}}{2})}&
{{1}{2^{2}}{3}}\hskip\inclskip\quad\hfil%
{({\ol{2}}{\ol{3}}){\times}({2^{2}})}&
{-({1}{3}{4})}\hskip\inclskip\quad\hfil%
{({1^{2}}{\ol{1}}{\ol{2}}){\times}({4})}&
{{4^{2}}_{\llap{$\scriptstyle{v}$}{}}}\hskip\inclskip\hskip-\inclskip\hfil\cr\rg%
{({1}{4}){\times}({1^{2}}{2})}&
{{1^{2}}{2}{4}}\hskip\inclskip\quad\hfil%
{({\ol{1}}{\ol{4}}){\times}({2^{2}})}&
{-({1^{2}}{2}{4^{-1}}{8})}\hskip\inclskip\quad\hfil%
{({\ol{1}}{2}{\ol{2}}){\times}({4})}&
{-({1^{2}}{2^{-1}}{4^{2}})}\hskip\inclskip\hskip-\inclskip\hfil\cr\rg%
{({5}){\times}({1^{2}}{2})}&
{{1}{2}{5}}\hskip\inclskip\quad\hfil%
{({1^{5}}){\times}({1}{3})}&
{{1^{5}}{3}}\hskip\inclskip\quad\hfil%
{({1}{\ol{1}^{\dn{4}}}){\times}({4})}&
{-({1^{4}}{4})}\hskip\inclskip\hskip-\inclskip\hfil\cr\rg%
{({1^{3}}{\ol{1}^{\dn{2}}}){\times}({1^{2}}{2})}&
{{1^{2}}{2^{3}}}\hskip\inclskip\quad\hfil%
{({1^{3}}{2}){\times}({1}{3})}&
{{1^{3}}{2}{3}}\hskip\inclskip\quad\hfil%
{({1}{\ol{2}^{\dn{2}}}){\times}({4})}&
{{2^{-2}}{4^{3}}}\hskip\inclskip\hskip-\inclskip\hfil\cr\rg%
{({1}{\ol{1}^{\dn{2}}}{2}){\times}({1^{2}}{2})}&
{{2^{4}}_{\llap{$\scriptstyle{u}$}{}}}\hskip\inclskip\quad\hfil%
{({1}{2^{2}}){\times}({1}{3})}&
{{1}{2^{2}}{3}}\hskip\inclskip\quad\hfil%
{({1}{\ol{1}}{\ol{3}}){\times}({4})}&
{-({1}{3}{4})}\hskip\inclskip\hskip-\inclskip\hfil\cr\rg%
{({\ol{1}^{\dn{2}}}{3}){\times}({1^{2}}{2})}&
{-({1}{2^{2}}{3^{-1}}{6})}\hskip\inclskip\quad\hfil%
{({1^{2}}{3}){\times}({1}{3})}&
{{1^{2}}{3^{2}}}\hskip\inclskip\quad\hfil%
{({\ol{1}^{\dn{3}}}{\ol{2}}){\times}({4})}&
{-({1^{4}}{2^{-2}}{4^{2}})}\hskip\inclskip\hskip-\inclskip\hfil\cr\rg%
{({1^{2}}{\ol{1}}{\ol{2}}){\times}({1^{2}}{2})}&
{{1^{2}}{2}{4}}\hskip\inclskip\quad\hfil%
{({2}{3}){\times}({1}{3})}&
{{2}{3^{2}}}\hskip\inclskip\quad\hfil%
{({\ol{2}}{\ol{3}}){\times}({4})}&
{-({1}{2^{-2}}{3}{4^{2}})}\hskip\inclskip\hskip-\inclskip\hfil\cr\rg%
{({\ol{1}}{2}{\ol{2}}){\times}({1^{2}}{2})}&
{{2^{2}}{4}_{u}}\hskip\inclskip\quad\hfil%
{({1}{4}){\times}({1}{3})}&
{{1}{3}{4}}\hskip\inclskip\quad\hfil%
{({\ol{1}}{\ol{4}}){\times}({4})}&
{-({1^{2}}{2^{-1}}{8})}\hskip\inclskip\hskip-\inclskip\hfil\cr\rg%
}}%
\tableclose%
\eject
\def\bol{\quad}
\def\rg{\rgloc}

\inclskip=20 true pt
\vbox{%
\tableopen{Class inclusions for $W({E_{6}}{A_{2}})\,\subset\,W(E_{8})$}%
\nobreak\vskip 10 true pt\nobreak%
\baselineskip=\tablebaseskip%
\tabskip=\incltabskip\halign to \hsize{%
\bol\hfil$#$&%
$\,\subset\,#$&&%
$\,\subset\,#$\cr%
{{1^{6}}{\times}({1^{3}})}&
{{1^{8}}}\hskip\inclskip\quad\hfil%
{{3^{-1}}{9}{\times}({1}{2})}&
{{2}{3^{-1}}{9}}\hskip\inclskip\hskip-\inclskip\hfil\cr\rg%
{{1^{2}}{2^{2}}{\times}({1^{3}})}&
{{1^{4}}{2^{2}}}\hskip\inclskip\quad\hfil%
{{1}{5}{\times}({1}{2})}&
{{1}{2}{5}}\hskip\inclskip\hskip-\inclskip\hfil\cr\rg%
{{1^{-2}}{2^{4}}{\times}({1^{3}})}&
{{2^{4}}_{\llap{$\scriptstyle{v}$}{}}}\hskip\inclskip\quad\hfil%
{{1^{4}}{2}{\times}({1}{2})}&
{{1^{4}}{2^{2}}}\hskip\inclskip\hskip-\inclskip\hfil\cr\rg%
{{1^{2}}{2^{-2}}{4^{2}}{\times}({1^{3}})}&
{{1^{4}}{2^{-2}}{4^{2}}}\hskip\inclskip\quad\hfil%
{{2^{3}}{\times}({1}{2})}&
{{2^{4}}_{\llap{$\scriptstyle{u}$}{}}}\hskip\inclskip\hskip-\inclskip\hfil\cr\rg%
{{2}{4}{\times}({1^{3}})}&
{{1^{2}}{2}{4}}\hskip\inclskip\quad\hfil%
{{1^{2}}{4}{\times}({1}{2})}&
{{1^{2}}{2}{4}}\hskip\inclskip\hskip-\inclskip\hfil\cr\rg%
{{1^{3}}{3}{\times}({1^{3}})}&
{{1^{5}}{3}}\hskip\inclskip\quad\hfil%
{{1^{-2}}{2^{2}}{4}{\times}({1}{2})}&
{-({1^{2}}{2}{4})}\hskip\inclskip\hskip-\inclskip\hfil\cr\rg%
{{1}{2}{3^{-1}}{6}{\times}({1^{3}})}&
{{1^{3}}{2}{3^{-1}}{6}}\hskip\inclskip\quad\hfil%
{{2}{4^{-1}}{8}{\times}({1}{2})}&
{{2^{2}}{4^{-1}}{8}}\hskip\inclskip\hskip-\inclskip\hfil\cr\rg%
{{1^{-1}}{2^{2}}{3}{\times}({1^{3}})}&
{{1}{2^{2}}{3}}\hskip\inclskip\quad\hfil%
{{1}{2}{3}{\times}({1}{2})}&
{{1}{2^{2}}{3}}\hskip\inclskip\hskip-\inclskip\hfil\cr\rg%
{{3^{2}}{\times}({1^{3}})}&
{{1^{2}}{3^{2}}}\hskip\inclskip\quad\hfil%
{{1^{-2}}{2}{3^{2}}{\times}({1}{2})}&
{-({1^{2}}{3^{-2}}{6^{2}})}\hskip\inclskip\hskip-\inclskip\hfil\cr\rg%
{{1^{-2}}{2}{6}{\times}({1^{3}})}&
{{2}{6}_{v}}\hskip\inclskip\quad\hfil%
{{6}{\times}({1}{2})}&
{{2}{6}_{u}}\hskip\inclskip\hskip-\inclskip\hfil\cr\rg%
{{1^{-3}}{3^{3}}{\times}({1^{3}})}&
{-({1}{2^{-1}}{3^{-3}}{6^{3}})}\hskip\inclskip\quad\hfil%
{{1}{2^{-1}}{3^{-1}}{4}{6}{\times}({1}{2})}&
{{1}{3^{-1}}{4}{6}}\hskip\inclskip\hskip-\inclskip\hfil\cr\rg%
{{1}{2^{-2}}{3^{-1}}{6^{2}}{\times}({1^{3}})}&
{{1^{3}}{2^{-2}}{3^{-1}}{6^{2}}}\hskip\inclskip\quad\hfil%
{{1^{-1}}{2}{5}{\times}({1}{2})}&
{-({1}{2}{5^{-1}}{(10)})}\hskip\inclskip\hskip-\inclskip\hfil\cr\rg%
{{1^{-1}}{2}{3}{4^{-1}}{6^{-1}}{(12)}{\times}({1^{3}})}&
{{1}{2}{3}{4^{-1}}{6^{-1}}{(12)}}\hskip\inclskip\quad\hfil%
{{1^{6}}{\times}({3})}&
{{1^{5}}{3}}\hskip\inclskip\hskip-\inclskip\hfil\cr\rg%
{{3^{-1}}{9}{\times}({1^{3}})}&
{{1^{2}}{3^{-1}}{9}}\hskip\inclskip\quad\hfil%
{{1^{2}}{2^{2}}{\times}({3})}&
{{1}{2^{2}}{3}}\hskip\inclskip\hskip-\inclskip\hfil\cr\rg%
{{1}{5}{\times}({1^{3}})}&
{{1^{3}}{5}}\hskip\inclskip\quad\hfil%
{{1^{-2}}{2^{4}}{\times}({3})}&
{-({1^{3}}{2}{3^{-1}}{6})}\hskip\inclskip\hskip-\inclskip\hfil\cr\rg%
{{1^{4}}{2}{\times}({1^{3}})}&
{{1^{6}}{2}}\hskip\inclskip\quad\hfil%
{{1^{2}}{2^{-2}}{4^{2}}{\times}({3})}&
{{1}{2^{-2}}{3}{4^{2}}}\hskip\inclskip\hskip-\inclskip\hfil\cr\rg%
{{2^{3}}{\times}({1^{3}})}&
{{1^{2}}{2^{3}}}\hskip\inclskip\quad\hfil%
{{2}{4}{\times}({3})}&
{-({1}{3^{-1}}{4}{6})}\hskip\inclskip\hskip-\inclskip\hfil\cr\rg%
{{1^{2}}{4}{\times}({1^{3}})}&
{{1^{4}}{4}}\hskip\inclskip\quad\hfil%
{{1^{3}}{3}{\times}({3})}&
{{1^{2}}{3^{2}}}\hskip\inclskip\hskip-\inclskip\hfil\cr\rg%
{{1^{-2}}{2^{2}}{4}{\times}({1^{3}})}&
{{2^{2}}{4}_{v}}\hskip\inclskip\quad\hfil%
{{1}{2}{3^{-1}}{6}{\times}({3})}&
{{2}{6}_{vv}}\hskip\inclskip\hskip-\inclskip\hfil\cr\rg%
{{2}{4^{-1}}{8}{\times}({1^{3}})}&
{{1^{2}}{2}{4^{-1}}{8}}\hskip\inclskip\quad\hfil%
{{1^{-1}}{2^{2}}{3}{\times}({3})}&
{-({1^{2}}{3^{-2}}{6^{2}})}\hskip\inclskip\hskip-\inclskip\hfil\cr\rg%
{{1}{2}{3}{\times}({1^{3}})}&
{{1^{3}}{2}{3}}\hskip\inclskip\quad\hfil%
{{3^{2}}{\times}({3})}&
{-({1}{2^{-1}}{3^{-3}}{6^{3}})}\hskip\inclskip\hskip-\inclskip\hfil\cr\rg%
{{1^{-2}}{2}{3^{2}}{\times}({1^{3}})}&
{{2}{3^{2}}}\hskip\inclskip\quad\hfil%
{{1^{-2}}{2}{6}{\times}({3})}&
{-({1^{3}}{2^{-2}}{3^{-1}}{6^{2}})}\hskip\inclskip\hskip-\inclskip\hfil\cr\rg%
{{6}{\times}({1^{3}})}&
{{1^{2}}{6}}\hskip\inclskip\quad\hfil%
{{1^{-3}}{3^{3}}{\times}({3})}&
{-({1^{4}}{2^{-4}}{3^{-4}}{6^{4}})}\hskip\inclskip\hskip-\inclskip\hfil\cr\rg%
{{1}{2^{-1}}{3^{-1}}{4}{6}{\times}({1^{3}})}&
{{1^{3}}{2^{-1}}{3^{-1}}{4}{6}}\hskip\inclskip\quad\hfil%
{{1}{2^{-2}}{3^{-1}}{6^{2}}{\times}({3})}&
{{2^{-2}}{6^{2}}}\hskip\inclskip\hskip-\inclskip\hfil\cr\rg%
{{1^{-1}}{2}{5}{\times}({1^{3}})}&
{{1}{2}{5}}\hskip\inclskip\quad\hfil%
{{1^{-1}}{2}{3}{4^{-1}}{6^{-1}}{(12)}{\times}({3})}&
{-({1^{2}}{2^{-1}}{3^{-2}}{4^{-1}}{6}{(12)})}\hskip\inclskip\hskip-\inclskip\hfil\cr\rg%
{{1^{6}}{\times}({1}{2})}&
{{1^{6}}{2}}\hskip\inclskip\quad\hfil%
{{3^{-1}}{9}{\times}({3})}&
{-({1}{2^{-1}}{9^{-1}}{(18)})}\hskip\inclskip\hskip-\inclskip\hfil\cr\rg%
{{1^{2}}{2^{2}}{\times}({1}{2})}&
{{1^{2}}{2^{3}}}\hskip\inclskip\quad\hfil%
{{1}{5}{\times}({3})}&
{{3}{5}}\hskip\inclskip\hskip-\inclskip\hfil\cr\rg%
{{1^{-2}}{2^{4}}{\times}({1}{2})}&
{-({1^{2}}{2^{3}})}\hskip\inclskip\quad\hfil%
{{1^{4}}{2}{\times}({3})}&
{{1^{3}}{2}{3}}\hskip\inclskip\hskip-\inclskip\hfil\cr\rg%
{{1^{2}}{2^{-2}}{4^{2}}{\times}({1}{2})}&
{{1^{2}}{2^{-1}}{4^{2}}}\hskip\inclskip\quad\hfil%
{{2^{3}}{\times}({3})}&
{-({1}{2^{2}}{3^{-1}}{6})}\hskip\inclskip\hskip-\inclskip\hfil\cr\rg%
{{2}{4}{\times}({1}{2})}&
{{2^{2}}{4}_{u}}\hskip\inclskip\quad\hfil%
{{1^{2}}{4}{\times}({3})}&
{{1}{3}{4}}\hskip\inclskip\hskip-\inclskip\hfil\cr\rg%
{{1^{3}}{3}{\times}({1}{2})}&
{{1^{3}}{2}{3}}\hskip\inclskip\quad\hfil%
{{1^{-2}}{2^{2}}{4}{\times}({3})}&
{-({1^{3}}{2^{-1}}{3^{-1}}{4}{6})}\hskip\inclskip\hskip-\inclskip\hfil\cr\rg%
{{1}{2}{3^{-1}}{6}{\times}({1}{2})}&
{{1}{2^{2}}{3^{-1}}{6}}\hskip\inclskip\quad\hfil%
{{2}{4^{-1}}{8}{\times}({3})}&
{-({1}{3^{-1}}{4^{-1}}{6}{8})}\hskip\inclskip\hskip-\inclskip\hfil\cr\rg%
{{1^{-1}}{2^{2}}{3}{\times}({1}{2})}&
{-({1}{2^{2}}{3^{-1}}{6})}\hskip\inclskip\quad\hfil%
{{1}{2}{3}{\times}({3})}&
{{2}{3^{2}}}\hskip\inclskip\hskip-\inclskip\hfil\cr\rg%
{{3^{2}}{\times}({1}{2})}&
{{2}{3^{2}}}\hskip\inclskip\quad\hfil%
{{1^{-2}}{2}{3^{2}}{\times}({3})}&
{-({1^{3}}{2^{-2}}{3^{-3}}{6^{3}})}\hskip\inclskip\hskip-\inclskip\hfil\cr\rg%
{{1^{-2}}{2}{6}{\times}({1}{2})}&
{-({1^{2}}{6})}\hskip\inclskip\quad\hfil%
{{6}{\times}({3})}&
{-({1}{2^{-1}}{3^{-1}}{6^{2}})}\hskip\inclskip\hskip-\inclskip\hfil\cr\rg%
{{1^{-3}}{3^{3}}{\times}({1}{2})}&
{-({1^{3}}{2^{-2}}{3^{-3}}{6^{3}})}\hskip\inclskip\quad\hfil%
{{1}{2^{-1}}{3^{-1}}{4}{6}{\times}({3})}&
{{2^{-1}}{4}{6}}\hskip\inclskip\hskip-\inclskip\hfil\cr\rg%
{{1}{2^{-2}}{3^{-1}}{6^{2}}{\times}({1}{2})}&
{{1}{2^{-1}}{3^{-1}}{6^{2}}}\hskip\inclskip\quad\hfil%
{{1^{-1}}{2}{5}{\times}({3})}&
{-({1^{2}}{2^{-1}}{3^{-1}}{5^{-1}}{6}{(10)})}\hskip\inclskip\hskip-\inclskip\hfil\cr\rg%
{{1^{-1}}{2}{3}{4^{-1}}{6^{-1}}{(12)}{\times}({1}{2})}&
{-({1}{2}{3^{-1}}{4^{-1}}{(12)})}\hskip\inclskip\quad\hfil%
&
\omit\hskip\inclskip\hfil\hskip-\inclskip\hfil\cr\rg%
}}%
\tableclose%

\inclskip=20 true pt
\def\rg{\noalign{\nobreak\vskip 2 true pt plus 3 true pt}}
\vbox{%
\tableopen{Class inclusions for $W({E_{7}}{A_{1}})\,\subset\,W(E_{8})$}%
\nobreak\vskip 10 true pt\nobreak%
\baselineskip=\tablebaseskip%
\tabskip=\incltabskip\halign to \hsize{%
\bol\hfil$#$&%
$\,\subset\,#$&&%
$\,\subset\,#$\cr%
{{1^{7}}{\times}({1^{2}})}&
{{1^{8}}}\hskip\inclskip\quad\hfil%
{-({1^{7}}){\times}({1^{2}})}&
{-({1^{6}}{2})}\hskip\inclskip\hskip-\inclskip\hfil\cr\rg%
{{1^{3}}{2^{2}}{\times}({1^{2}})}&
{{1^{4}}{2^{2}}}\hskip\inclskip\quad\hfil%
{-({1^{3}}{2^{2}}){\times}({1^{2}})}&
{-({1^{2}}{2^{3}})}\hskip\inclskip\hskip-\inclskip\hfil\cr\rg%
{{1^{-1}}{2^{4}}_{\llap{$\scriptstyle{v}$}{}}{\times}({1^{2}})}&
{{2^{4}}_{\llap{$\scriptstyle{v}$}{}}}\hskip\inclskip\quad\hfil%
{-({1^{-1}}{2^{4}}_{\llap{$\scriptstyle{v}$}{}}){\times}({1^{2}})}&
{{1^{2}}{2^{3}}}\hskip\inclskip\hskip-\inclskip\hfil\cr\rg%
{{1^{3}}{2^{-2}}{4^{2}}{\times}({1^{2}})}&
{{1^{4}}{2^{-2}}{4^{2}}}\hskip\inclskip\quad\hfil%
{-({1^{3}}{2^{-2}}{4^{2}}){\times}({1^{2}})}&
{-({1^{2}}{2^{-1}}{4^{2}})}\hskip\inclskip\hskip-\inclskip\hfil\cr\rg%
{{1}{2}{4}_{u}{\times}({1^{2}})}&
{{1^{2}}{2}{4}}\hskip\inclskip\quad\hfil%
{-({1}{2}{4}_{u}){\times}({1^{2}})}&
{{2^{2}}{4}_{u}}\hskip\inclskip\hskip-\inclskip\hfil\cr\rg%
{{1^{4}}{3}{\times}({1^{2}})}&
{{1^{5}}{3}}\hskip\inclskip\quad\hfil%
{-({1^{4}}{3}){\times}({1^{2}})}&
{-({1^{3}}{2}{3})}\hskip\inclskip\hskip-\inclskip\hfil\cr\rg%
{{1^{2}}{2}{3^{-1}}{6}{\times}({1^{2}})}&
{{1^{3}}{2}{3^{-1}}{6}}\hskip\inclskip\quad\hfil%
{-({1^{2}}{2}{3^{-1}}{6}){\times}({1^{2}})}&
{-({1}{2^{2}}{3^{-1}}{6})}\hskip\inclskip\hskip-\inclskip\hfil\cr\rg%
{{2^{2}}{3}{\times}({1^{2}})}&
{{1}{2^{2}}{3}}\hskip\inclskip\quad\hfil%
{-({2^{2}}{3}){\times}({1^{2}})}&
{{1}{2^{2}}{3^{-1}}{6}}\hskip\inclskip\hskip-\inclskip\hfil\cr\rg%
{{1}{3^{2}}{\times}({1^{2}})}&
{{1^{2}}{3^{2}}}\hskip\inclskip\quad\hfil%
{-({1}{3^{2}}){\times}({1^{2}})}&
{-({2}{3^{2}})}\hskip\inclskip\hskip-\inclskip\hfil\cr\rg%
{{1^{-1}}{2}{6}_{v}{\times}({1^{2}})}&
{{2}{6}_{v}}\hskip\inclskip\quad\hfil%
{-({1^{-1}}{2}{6}_{v}){\times}({1^{2}})}&
{{1^{2}}{6}}\hskip\inclskip\hskip-\inclskip\hfil\cr\rg%
{{1^{-2}}{3^{3}}{\times}({1^{2}})}&
{-({1}{2^{-1}}{3^{-3}}{6^{3}})}\hskip\inclskip\quad\hfil%
{-({1^{-2}}{3^{3}}){\times}({1^{2}})}&
{{1^{3}}{2^{-2}}{3^{-3}}{6^{3}}}\hskip\inclskip\hskip-\inclskip\hfil\cr\rg%
{{1^{2}}{2^{-2}}{3^{-1}}{6^{2}}{\times}({1^{2}})}&
{{1^{3}}{2^{-2}}{3^{-1}}{6^{2}}}\hskip\inclskip\quad\hfil%
{-({1^{2}}{2^{-2}}{3^{-1}}{6^{2}}){\times}({1^{2}})}&
{-({1}{2^{-1}}{3^{-1}}{6^{2}})}\hskip\inclskip\hskip-\inclskip\hfil\cr\rg%
{{2}{3}{4^{-1}}{6^{-1}}{(12)}{\times}({1^{2}})}&
{{1}{2}{3}{4^{-1}}{6^{-1}}{(12)}}\hskip\inclskip\quad\hfil%
{-({2}{3}{4^{-1}}{6^{-1}}{(12)}){\times}({1^{2}})}&
{{1}{2}{3^{-1}}{4^{-1}}{(12)}}\hskip\inclskip\hskip-\inclskip\hfil\cr\rg%
{{1}{3^{-1}}{9}{\times}({1^{2}})}&
{{1^{2}}{3^{-1}}{9}}\hskip\inclskip\quad\hfil%
{-({1}{3^{-1}}{9}){\times}({1^{2}})}&
{-({2}{3^{-1}}{9})}\hskip\inclskip\hskip-\inclskip\hfil\cr\rg%
{{1^{2}}{5}{\times}({1^{2}})}&
{{1^{3}}{5}}\hskip\inclskip\quad\hfil%
{-({1^{2}}{5}){\times}({1^{2}})}&
{-({1}{2}{5})}\hskip\inclskip\hskip-\inclskip\hfil\cr\rg%
{{1^{-5}}{2^{6}}{\times}({1^{2}})}&
{-({1^{4}}{2^{2}})}\hskip\inclskip\quad\hfil%
{-({1^{-5}}{2^{6}}){\times}({1^{2}})}&
{{1^{6}}{2}}\hskip\inclskip\hskip-\inclskip\hfil\cr\rg%
{{1^{-1}}{2^{4}}_{\llap{$\scriptstyle{u}$}{}}{\times}({1^{2}})}&
{{2^{4}}_{\llap{$\scriptstyle{u}$}{}}}\hskip\inclskip\quad\hfil%
{-({1^{-1}}{2^{4}}_{\llap{$\scriptstyle{u}$}{}}){\times}({1^{2}})}&
{{1^{2}}{2^{3}}}\hskip\inclskip\hskip-\inclskip\hfil\cr\rg%
{{1^{-3}}{2^{3}}{4}{\times}({1^{2}})}&
{-({1^{2}}{2}{4})}\hskip\inclskip\quad\hfil%
{-({1^{-3}}{2^{3}}{4}){\times}({1^{2}})}&
{{1^{4}}{4}}\hskip\inclskip\hskip-\inclskip\hfil\cr\rg%
{{1}{2}{4}_{v}{\times}({1^{2}})}&
{{1^{2}}{2}{4}}\hskip\inclskip\quad\hfil%
{-({1}{2}{4}_{v}){\times}({1^{2}})}&
{{2^{2}}{4}_{v}}\hskip\inclskip\hskip-\inclskip\hfil\cr\rg%
{{1^{-1}}{2^{2}}{4^{-1}}{8}{\times}({1^{2}})}&
{{2^{2}}{4^{-1}}{8}}\hskip\inclskip\quad\hfil%
{-({1^{-1}}{2^{2}}{4^{-1}}{8}){\times}({1^{2}})}&
{{1^{2}}{2}{4^{-1}}{8}}\hskip\inclskip\hskip-\inclskip\hfil\cr\rg%
{{1^{-2}}{2^{3}}{3^{-1}}{6}{\times}({1^{2}})}&
{-({1}{2^{2}}{3})}\hskip\inclskip\quad\hfil%
{-({1^{-2}}{2^{3}}{3^{-1}}{6}){\times}({1^{2}})}&
{{1^{3}}{2}{3}}\hskip\inclskip\hskip-\inclskip\hfil\cr\rg%
{{1}{3^{-2}}{6^{2}}{\times}({1^{2}})}&
{{1^{2}}{3^{-2}}{6^{2}}}\hskip\inclskip\quad\hfil%
{-({1}{3^{-2}}{6^{2}}){\times}({1^{2}})}&
{{2}{3^{2}}}\hskip\inclskip\hskip-\inclskip\hfil\cr\rg%
{{1^{-1}}{2}{6}_{u}{\times}({1^{2}})}&
{{2}{6}_{u}}\hskip\inclskip\quad\hfil%
{-({1^{-1}}{2}{6}_{u}){\times}({1^{2}})}&
{{1^{2}}{6}}\hskip\inclskip\hskip-\inclskip\hfil\cr\rg%
{{1^{-2}}{2}{3}{4}{\times}({1^{2}})}&
{-({1}{3^{-1}}{4}{6})}\hskip\inclskip\quad\hfil%
{-({1^{-2}}{2}{3}{4}){\times}({1^{2}})}&
{{1^{3}}{2^{-1}}{3^{-1}}{4}{6}}\hskip\inclskip\hskip-\inclskip\hfil\cr\rg%
{{2}{5^{-1}}{(10)}{\times}({1^{2}})}&
{{1}{2}{5^{-1}}{(10)}}\hskip\inclskip\quad\hfil%
{-({2}{5^{-1}}{(10)}){\times}({1^{2}})}&
{{1}{2}{5}}\hskip\inclskip\hskip-\inclskip\hfil\cr\rg%
{{1^{-1}}{3}{5}{\times}({1^{2}})}&
{{3}{5}}\hskip\inclskip\quad\hfil%
{-({1^{-1}}{3}{5}){\times}({1^{2}})}&
{{1^{2}}{2^{-1}}{3^{-1}}{5^{-1}}{6}{(10)}}\hskip\inclskip\hskip-\inclskip\hfil\cr\rg%
{{3^{-1}}{4}{6}{\times}({1^{2}})}&
{{1}{3^{-1}}{4}{6}}\hskip\inclskip\quad\hfil%
{-({3^{-1}}{4}{6}){\times}({1^{2}})}&
{{1}{3}{4}}\hskip\inclskip\hskip-\inclskip\hfil\cr\rg%
{{1^{-1}}{4^{2}}{\times}({1^{2}})}&
{{4^{2}}_{\llap{$\scriptstyle{v}$}{}}}\hskip\inclskip\quad\hfil%
{-({1^{-1}}{4^{2}}){\times}({1^{2}})}&
{{1^{2}}{2^{-1}}{4^{2}}}\hskip\inclskip\hskip-\inclskip\hfil\cr\rg%
{{1}{2^{-1}}{8}{\times}({1^{2}})}&
{{1^{2}}{2^{-1}}{8}}\hskip\inclskip\quad\hfil%
{-({1}{2^{-1}}{8}){\times}({1^{2}})}&
{{8}_{v}}\hskip\inclskip\hskip-\inclskip\hfil\cr\rg%
{{7}{\times}({1^{2}})}&
{{1}{7}}\hskip\inclskip\quad\hfil%
{-({7}){\times}({1^{2}})}&
{{1}{7^{-1}}{(14)}}\hskip\inclskip\hskip-\inclskip\hfil\cr\rg%
}}%
\vbox{%
\tablecont%
\nobreak\vskip 10 true pt\nobreak%
\baselineskip=\tablebaseskip%
\tabskip=\incltabskip\halign to \hsize{%
\bol\hfil$#$&%
$\,\subset\,#$&&%
$\,\subset\,#$\cr%
{{1^{7}}{\times}({2})}&
{{1^{6}}{2}}\hskip\inclskip\quad\hfil%
{-({1^{7}}){\times}({2})}&
{-({1^{8}})}\hskip\inclskip\hskip-\inclskip\hfil\cr\rg%
{{1^{3}}{2^{2}}{\times}({2})}&
{{1^{2}}{2^{3}}}\hskip\inclskip\quad\hfil%
{-({1^{3}}{2^{2}}){\times}({2})}&
{-({1^{4}}{2^{2}})}\hskip\inclskip\hskip-\inclskip\hfil\cr\rg%
{{1^{-1}}{2^{4}}_{\llap{$\scriptstyle{v}$}{}}{\times}({2})}&
{-({1^{2}}{2^{3}})}\hskip\inclskip\quad\hfil%
{-({1^{-1}}{2^{4}}_{\llap{$\scriptstyle{v}$}{}}){\times}({2})}&
{{2^{4}}_{\llap{$\scriptstyle{v}$}{}}}\hskip\inclskip\hskip-\inclskip\hfil\cr\rg%
{{1^{3}}{2^{-2}}{4^{2}}{\times}({2})}&
{{1^{2}}{2^{-1}}{4^{2}}}\hskip\inclskip\quad\hfil%
{-({1^{3}}{2^{-2}}{4^{2}}){\times}({2})}&
{-({1^{4}}{2^{-2}}{4^{2}})}\hskip\inclskip\hskip-\inclskip\hfil\cr\rg%
{{1}{2}{4}_{u}{\times}({2})}&
{{2^{2}}{4}_{u}}\hskip\inclskip\quad\hfil%
{-({1}{2}{4}_{u}){\times}({2})}&
{-({1^{2}}{2}{4})}\hskip\inclskip\hskip-\inclskip\hfil\cr\rg%
{{1^{4}}{3}{\times}({2})}&
{{1^{3}}{2}{3}}\hskip\inclskip\quad\hfil%
{-({1^{4}}{3}){\times}({2})}&
{-({1^{5}}{3})}\hskip\inclskip\hskip-\inclskip\hfil\cr\rg%
{{1^{2}}{2}{3^{-1}}{6}{\times}({2})}&
{{1}{2^{2}}{3^{-1}}{6}}\hskip\inclskip\quad\hfil%
{-({1^{2}}{2}{3^{-1}}{6}){\times}({2})}&
{-({1^{3}}{2}{3^{-1}}{6})}\hskip\inclskip\hskip-\inclskip\hfil\cr\rg%
{{2^{2}}{3}{\times}({2})}&
{-({1}{2^{2}}{3^{-1}}{6})}\hskip\inclskip\quad\hfil%
{-({2^{2}}{3}){\times}({2})}&
{-({1}{2^{2}}{3})}\hskip\inclskip\hskip-\inclskip\hfil\cr\rg%
{{1}{3^{2}}{\times}({2})}&
{{2}{3^{2}}}\hskip\inclskip\quad\hfil%
{-({1}{3^{2}}){\times}({2})}&
{-({1^{2}}{3^{2}})}\hskip\inclskip\hskip-\inclskip\hfil\cr\rg%
{{1^{-1}}{2}{6}_{v}{\times}({2})}&
{-({1^{2}}{6})}\hskip\inclskip\quad\hfil%
{-({1^{-1}}{2}{6}_{v}){\times}({2})}&
{{2}{6}_{v}}\hskip\inclskip\hskip-\inclskip\hfil\cr\rg%
{{1^{-2}}{3^{3}}{\times}({2})}&
{-({1^{3}}{2^{-2}}{3^{-3}}{6^{3}})}\hskip\inclskip\quad\hfil%
{-({1^{-2}}{3^{3}}){\times}({2})}&
{{1}{2^{-1}}{3^{-3}}{6^{3}}}\hskip\inclskip\hskip-\inclskip\hfil\cr\rg%
{{1^{2}}{2^{-2}}{3^{-1}}{6^{2}}{\times}({2})}&
{{1}{2^{-1}}{3^{-1}}{6^{2}}}\hskip\inclskip\quad\hfil%
{-({1^{2}}{2^{-2}}{3^{-1}}{6^{2}}){\times}({2})}&
{-({1^{3}}{2^{-2}}{3^{-1}}{6^{2}})}\hskip\inclskip\hskip-\inclskip\hfil\cr\rg%
{{2}{3}{4^{-1}}{6^{-1}}{(12)}{\times}({2})}&
{-({1}{2}{3^{-1}}{4^{-1}}{(12)})}\hskip\inclskip\quad\hfil%
{-({2}{3}{4^{-1}}{6^{-1}}{(12)}){\times}({2})}&
{-({1}{2}{3}{4^{-1}}{6^{-1}}{(12)})}\hskip\inclskip\hskip-\inclskip\hfil\cr\rg%
{{1}{3^{-1}}{9}{\times}({2})}&
{{2}{3^{-1}}{9}}\hskip\inclskip\quad\hfil%
{-({1}{3^{-1}}{9}){\times}({2})}&
{-({1^{2}}{3^{-1}}{9})}\hskip\inclskip\hskip-\inclskip\hfil\cr\rg%
{{1^{2}}{5}{\times}({2})}&
{{1}{2}{5}}\hskip\inclskip\quad\hfil%
{-({1^{2}}{5}){\times}({2})}&
{-({1^{3}}{5})}\hskip\inclskip\hskip-\inclskip\hfil\cr\rg%
{{1^{-5}}{2^{6}}{\times}({2})}&
{-({1^{6}}{2})}\hskip\inclskip\quad\hfil%
{-({1^{-5}}{2^{6}}){\times}({2})}&
{{1^{4}}{2^{2}}}\hskip\inclskip\hskip-\inclskip\hfil\cr\rg%
{{1^{-1}}{2^{4}}_{\llap{$\scriptstyle{u}$}{}}{\times}({2})}&
{-({1^{2}}{2^{3}})}\hskip\inclskip\quad\hfil%
{-({1^{-1}}{2^{4}}_{\llap{$\scriptstyle{u}$}{}}){\times}({2})}&
{{2^{4}}_{\llap{$\scriptstyle{u}$}{}}}\hskip\inclskip\hskip-\inclskip\hfil\cr\rg%
{{1^{-3}}{2^{3}}{4}{\times}({2})}&
{-({1^{4}}{4})}\hskip\inclskip\quad\hfil%
{-({1^{-3}}{2^{3}}{4}){\times}({2})}&
{{1^{2}}{2}{4}}\hskip\inclskip\hskip-\inclskip\hfil\cr\rg%
{{1}{2}{4}_{v}{\times}({2})}&
{{2^{2}}{4}_{v}}\hskip\inclskip\quad\hfil%
{-({1}{2}{4}_{v}){\times}({2})}&
{-({1^{2}}{2}{4})}\hskip\inclskip\hskip-\inclskip\hfil\cr\rg%
{{1^{-1}}{2^{2}}{4^{-1}}{8}{\times}({2})}&
{-({1^{2}}{2}{4^{-1}}{8})}\hskip\inclskip\quad\hfil%
{-({1^{-1}}{2^{2}}{4^{-1}}{8}){\times}({2})}&
{{2^{2}}{4^{-1}}{8}}\hskip\inclskip\hskip-\inclskip\hfil\cr\rg%
{{1^{-2}}{2^{3}}{3^{-1}}{6}{\times}({2})}&
{-({1^{3}}{2}{3})}\hskip\inclskip\quad\hfil%
{-({1^{-2}}{2^{3}}{3^{-1}}{6}){\times}({2})}&
{{1}{2^{2}}{3}}\hskip\inclskip\hskip-\inclskip\hfil\cr\rg%
{{1}{3^{-2}}{6^{2}}{\times}({2})}&
{-({2}{3^{2}})}\hskip\inclskip\quad\hfil%
{-({1}{3^{-2}}{6^{2}}){\times}({2})}&
{-({1^{2}}{3^{-2}}{6^{2}})}\hskip\inclskip\hskip-\inclskip\hfil\cr\rg%
{{1^{-1}}{2}{6}_{u}{\times}({2})}&
{-({1^{2}}{6})}\hskip\inclskip\quad\hfil%
{-({1^{-1}}{2}{6}_{u}){\times}({2})}&
{{2}{6}_{u}}\hskip\inclskip\hskip-\inclskip\hfil\cr\rg%
{{1^{-2}}{2}{3}{4}{\times}({2})}&
{-({1^{3}}{2^{-1}}{3^{-1}}{4}{6})}\hskip\inclskip\quad\hfil%
{-({1^{-2}}{2}{3}{4}){\times}({2})}&
{{1}{3^{-1}}{4}{6}}\hskip\inclskip\hskip-\inclskip\hfil\cr\rg%
{{2}{5^{-1}}{(10)}{\times}({2})}&
{-({1}{2}{5})}\hskip\inclskip\quad\hfil%
{-({2}{5^{-1}}{(10)}){\times}({2})}&
{-({1}{2}{5^{-1}}{(10)})}\hskip\inclskip\hskip-\inclskip\hfil\cr\rg%
{{1^{-1}}{3}{5}{\times}({2})}&
{-({1^{2}}{2^{-1}}{3^{-1}}{5^{-1}}{6}{(10)})}\hskip\inclskip\quad\hfil%
{-({1^{-1}}{3}{5}){\times}({2})}&
{-({3}{5})}\hskip\inclskip\hskip-\inclskip\hfil\cr\rg%
{{3^{-1}}{4}{6}{\times}({2})}&
{-({1}{3}{4})}\hskip\inclskip\quad\hfil%
{-({3^{-1}}{4}{6}){\times}({2})}&
{-({1}{3^{-1}}{4}{6})}\hskip\inclskip\hskip-\inclskip\hfil\cr\rg%
{{1^{-1}}{4^{2}}{\times}({2})}&
{-({1^{2}}{2^{-1}}{4^{2}})}\hskip\inclskip\quad\hfil%
{-({1^{-1}}{4^{2}}){\times}({2})}&
{{4^{2}}_{\llap{$\scriptstyle{v}$}{}}}\hskip\inclskip\hskip-\inclskip\hfil\cr\rg%
{{1}{2^{-1}}{8}{\times}({2})}&
{{8}_{v}}\hskip\inclskip\quad\hfil%
{-({1}{2^{-1}}{8}){\times}({2})}&
{-({1^{2}}{2^{-1}}{8})}\hskip\inclskip\hskip-\inclskip\hfil\cr\rg%
{{7}{\times}({2})}&
{-({1}{7^{-1}}{(14)})}\hskip\inclskip\quad\hfil%
{-({7}){\times}({2})}&
{-({1}{7})}\hskip\inclskip\hskip-\inclskip\hfil\cr\rg%
}}%
\tableclose%
\eject
\def\rg{\rgloc}

\medskip
\leftline{\bf Type $F_4$.}

\smallskip
The extended diagram for $F_4$ appears below.
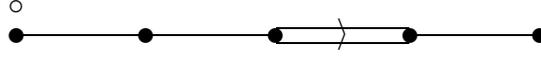
\begingroup
\def\star{\lower 3 true pt\hbox to 0pt{\hss{$\circ$}\hss}}
\def\bgt{\lower 2 true pt\hbox{$\rangle$}}
\def\hbarloc{\hbox to 27 true pt{\hss\leaders\hrule\hskip 50 true
pt\hss}}
\def\vbar{\vbox to 30 true pt{\vss\hbox to 0pt{\hss\vrule height
38 true pt\hss}\vss}}
\def\double{\lower 0.1 true pt\vbox to 6.3 true pt{%
\hbarloc\vss \hbox to 27 true pt{\hss\bgt\hss}\vss\hbarloc}}
$$
\matrix{
\star&&&&&&&&\cr
\spot&\hb&
\spot&\hb&
\spot&\double&
\spot&
\hb&
\spot\cr
}
$$
\endgroup

\medskip
\vskip 5 true pt
\noindent
All maximal subdiagrams
$\Gamma_0 \not= \Gamma$
of
$\widetilde{\Gamma}$ are considered.
Recall that the component
$A_2^{*}$ of ${A_2}{A_2^{*}}\,\subset\,F_4$
corresponds to short roots.

\smallskip
Given our assumptions concerning the generators of $W(F_4)$, the
subgroup $W(B_4)$ of $W(F_4)$ is
generated by $c$, $d$, $a$, and $\tau$
(and hence coincides with the
subgroup H of [K]).
The class inclusions for $W(B_4)$ in
$W(F_4)$ are given in Proposition 5 and Lemma 6 of [K].
The remaining cases $W({C_3}{A_1})$, $W({A_2}{A_2})$ and $W({A_3}{A_1})$
can be analyzed for the most part using the case $W(B_4)$ and the graph
automorphism of $W(F_4)$.
Note that the classes
$({1^2}{2})$ and
$({1^3}{\ol{1}})$ of $W(B_4)$
are contained in ${1^2}{2}_{d}$ and ${1^2}{2}_{\tau}$, respectively,
while the classes
$({1}{2})\times(1^2)$ and
$({1^2}{\ol{1}})\times(1^2)$
of $W({C_3}{A_1})$ are contained in
${1^2}{2}_{\tau}$ and ${1^2}{2}_{d}$, respectively.
Perhaps the most interesting case is that of the
class
$(2^2)\times(2)$ of $W({A_3}{A_1})$.
This class has eigenvalue structure
$-({1^2}{2})$, and so is
contained in the class of ${d}{z}$
or the class of ${\tau}{z}$
in $W(F_4)$.
However, ${d}{z}$ is an element of the subgroup $W(D_4)$ of $W(B_4)$,
$W(D_4)$ is normal in $W(F_4)$, and
$(2^2)\times(2)$ contains
${c}{a}{\tau}{\sigma}$, which is not an element of $W(D_4)$.  Therefore
$(2^2)\times(2)$ is contained in the class of
${\tau}{z}$, that is,
$-({1^2}{2}_{\tau})$.
The class inclusions for type $F_4$ appear in 
Tables~17 through~20.

\def\rg{\noalign{\nobreak\vskip 0pt plus 2 true pt}}

\vskip 10 true pt plus 30 true pt

\inclskip=40 true pt
\vbox{%
\tableopen{Class inclusions for $W({C_{3}}{A_{1}})\,\subset\,W(F_{4})$}%
\nobreak\vskip 10 true pt\nobreak%
\baselineskip=\tablebaseskip%
\tabskip=\incltabskip\halign to \hsize{%
\bol\hfil$#$&%
$\,\subset\,#$&&%
$\,\subset\,#$\cr%
{({1^{3}}){\times}({1^{2}})}&
{{1^{4}}}\hskip\inclskip\quad\hfil%
{({\ol{1}^{\dn{3}}}){\times}({1^{2}})}&
{-({1^{2}}{2}_{d})}\hskip\inclskip\quad\hfil%
{({\ol{1}}{2}){\times}({2})}&
{-({1^{2}}{2}_{\tau})}\hskip\inclskip\hskip-\inclskip\hfil\cr\rg%
{({1}{2}){\times}({1^{2}})}&
{{1^{2}}{2}_{\tau}}\hskip\inclskip\quad\hfil%
{({\ol{1}}{\ol{2}}){\times}({1^{2}})}&
{{4}_{ca{\tau}}}\hskip\inclskip\quad\hfil%
{({1}{\ol{1}^{\dn{2}}}){\times}({2})}&
{-({1^{2}}{2}_{d})}\hskip\inclskip\hskip-\inclskip\hfil\cr\rg%
{({3}){\times}({1^{2}})}&
{{1}{3}_{\sigma}}\hskip\inclskip\quad\hfil%
{({\ol{3}}){\times}({1^{2}})}&
{-({1^{-1}}{2}{3}_{\sigma{d}})}\hskip\inclskip\quad\hfil%
{({1}{\ol{2}}){\times}({2})}&
{{4}_{ca{\tau}}}\hskip\inclskip\hskip-\inclskip\hfil\cr\rg%
{({1^{2}}{\ol{1}}){\times}({1^{2}})}&
{{1^{2}}{2}_{d}}\hskip\inclskip\quad\hfil%
{({1^{3}}){\times}({2})}&
{{1^{2}}{2}_{d}}\hskip\inclskip\quad\hfil%
{({\ol{1}^{\dn{3}}}){\times}({2})}&
{-({1^{4}})}\hskip\inclskip\hskip-\inclskip\hfil\cr\rg%
{({\ol{1}}{2}){\times}({1^{2}})}&
{{2^{2}}_{\llap{$\scriptstyle{\tau}$}{{d}}}}\hskip\inclskip\quad\hfil%
{({1}{2}){\times}({2})}&
{{2^{2}}_{\llap{$\scriptstyle{\tau}$}{{d}}}}\hskip\inclskip\quad\hfil%
{({\ol{1}}{\ol{2}}){\times}({2})}&
{-({1^{2}}{2^{-1}}{4})}\hskip\inclskip\hskip-\inclskip\hfil\cr\rg%
{({1}{\ol{1}^{\dn{2}}}){\times}({1^{2}})}&
{{2^{2}}_{\llap{$\scriptstyle{a}$}{b}}}\hskip\inclskip\quad\hfil%
{({3}){\times}({2})}&
{{1^{-1}}{2}{3}_{\sigma{d}}}\hskip\inclskip\quad\hfil%
{({\ol{3}}){\times}({2})}&
{-({1}{3}_{\sigma})}\hskip\inclskip\hskip-\inclskip\hfil\cr\rg%
{({1}{\ol{2}}){\times}({1^{2}})}&
{{1^{2}}{2^{-1}}{4}}\hskip\inclskip\quad\hfil%
{({1^{2}}{\ol{1}}){\times}({2})}&
{{2^{2}}_{\llap{$\scriptstyle{a}$}{b}}}\hskip\inclskip\quad\hfil%
&
\omit\hskip\inclskip\hfil\hskip-\inclskip\hfil\cr\rg%
}}%
\tableclose%

\vskip 5 true pt plus 15 true pt
\goodbreak
\inclskip=15 true pt
\vbox{%
\tableopen{Class inclusions for $W({A_{2}}{A_{2}^{*}})\,\subset\,W(F_{4})$}%
\nobreak\vskip 10 true pt\nobreak%
\baselineskip=\tablebaseskip%
\tabskip=\incltabskip\halign to \hsize{%
\bol\hfil$#$&%
$\,\subset\,#$&&%
$\,\subset\,#$\cr%
{({1^{3}}){\times}({1^{3}})}&
{{1^{4}}}\hskip\inclskip\quad\hfil%
{({1^{3}}){\times}({1}{2})}&
{{1^{2}}{2}_{\tau}}\hskip\inclskip\quad\hfil%
{({3}){\times}({1}{2})}&
{{1^{-1}}{2}{3}_{{e}\tau}}\hskip\inclskip\quad\hfil%
{({1}{2}){\times}({3})}&
{{1^{-1}}{2}{3}_{\sigma{d}}}\hskip\inclskip\hskip-\inclskip\hfil\cr\rg%
{({1}{2}){\times}({1^{3}})}&
{{1^{2}}{2}_{d}}\hskip\inclskip\quad\hfil%
{({1}{2}){\times}({1}{2})}&
{{2^{2}}_{\llap{$\scriptstyle{\tau}$}{{d}}}}\hskip\inclskip\quad\hfil%
{({1^{3}}){\times}({3})}&
{{1}{3}_{\sigma}}\hskip\inclskip\quad\hfil%
{({3}){\times}({3})}&
{{1^{-2}}{3^{2}}}\hskip\inclskip\hskip-\inclskip\hfil\cr\rg%
{({3}){\times}({1^{3}})}&
{{1}{3}_{e}}\hskip\inclskip\quad\hfil%
&
\omit\hskip\inclskip\hfil\quad\hfil%
&
\omit\hskip\inclskip\hfil\quad\hfil%
&
\omit\hskip\inclskip\hfil\hskip-\inclskip\hfil\cr\rg%
}}%
\tableclose%

\inclskip=14 true pt
\vbox{%
\tableopen{Class inclusions for $W({A_{3}}{A_{1}})\,\subset\,W(F_{4})$}%
\nobreak\vskip 10 true pt\nobreak%
\baselineskip=\tablebaseskip%
\tabskip=\incltabskip\halign to \hsize{%
\bol\hfil$#$&%
$\,\subset\,#$&&%
$\,\subset\,#$\cr%
{({1^{4}}){\times}({1^{2}})}&
{{1^{4}}}\hskip\inclskip\quad\hfil%
{({1}{3}){\times}({1^{2}})}&
{{1}{3}_{e}}\hskip\inclskip\quad\hfil%
{({1^{2}}{2}){\times}({2})}&
{{2^{2}}_{\llap{$\scriptstyle{\tau}$}{{d}}}}\hskip\inclskip\quad\hfil%
{({1}{3}){\times}({2})}&
{{1^{-1}}{2}{3}_{{e}\tau}}\hskip\inclskip\hskip-\inclskip\hfil\cr\rg%
{({1^{2}}{2}){\times}({1^{2}})}&
{{1^{2}}{2}_{d}}\hskip\inclskip\quad\hfil%
{({4}){\times}({1^{2}})}&
{{4}_{adb}}\hskip\inclskip\quad\hfil%
{({2^{2}}){\times}({2})}&
{-({1^{2}}{2}_{\tau})}\hskip\inclskip\quad\hfil%
{({4}){\times}({2})}&
{-({1^{2}}{2^{-1}}{4})}\hskip\inclskip\hskip-\inclskip\hfil\cr\rg%
{({2^{2}}){\times}({1^{2}})}&
{{2^{2}}_{\llap{$\scriptstyle{a}$}{b}}}\hskip\inclskip\quad\hfil%
{({1^{4}}){\times}({2})}&
{{1^{2}}{2}_{\tau}}\hskip\inclskip\quad\hfil%
&
\omit\hskip\inclskip\hfil\quad\hfil%
&
\omit\hskip\inclskip\hfil\hskip-\inclskip\hfil\cr\rg%
}}%
\tableclose%

\def\bol{\qquad}
\inclskip=20 true pt
\vbox{%
\tableopen{Class inclusions for $W(B_{4})\,\subset\,W(F_{4})$}%
\nobreak\vskip 10 true pt\nobreak%
\baselineskip=\tablebaseskip%
\tabskip=\incltabskip\halign to \hsize{%
\bol\hfil$#$&%
$\,\subset\,#$&&%
$\,\subset\,#$\cr%
{({1^{4}})}&
{{1^{4}}}\hskip\inclskip\quad\hfil%
{({1^{3}}{\ol{1}})}&
{{1^{2}}{2}_{\tau}}\hskip\inclskip\quad\hfil%
{({\ol{1}^{\dn{2}}}{2})}&
{-({1^{2}}{2}_{d})}\hskip\inclskip\quad\hfil%
{({\ol{1}^{\dn{4}}})}&
{-({1^{4}})}\hskip\inclskip\hskip-\inclskip\hfil\cr\rg%
{({1^{2}}{2})}&
{{1^{2}}{2}_{d}}\hskip\inclskip\quad\hfil%
{({1}{\ol{1}}{2})}&
{{2^{2}}_{\llap{$\scriptstyle{\tau}$}{{d}}}}\hskip\inclskip\quad\hfil%
{({2}{\ol{2}})}&
{{4}_{ca{\tau}}}\hskip\inclskip\quad\hfil%
{({\ol{1}^{\dn{2}}}{\ol{2}})}&
{-({1^{2}}{2^{-1}}{4})}\hskip\inclskip\hskip-\inclskip\hfil\cr\rg%
{({2^{2}})}&
{{2^{2}}_{\llap{$\scriptstyle{a}$}{b}}}\hskip\inclskip\quad\hfil%
{({\ol{1}}{3})}&
{{1^{-1}}{2}{3}_{{e}\tau}}\hskip\inclskip\quad\hfil%
{({1}{\ol{1}^{\dn{3}}})}&
{-({1^{2}}{2}_{\tau})}\hskip\inclskip\quad\hfil%
{({\ol{2}^{\dn{2}}})}&
{{2^{-2}}{4^{2}}}\hskip\inclskip\hskip-\inclskip\hfil\cr\rg%
{({1}{3})}&
{{1}{3}_{e}}\hskip\inclskip\quad\hfil%
{({1^{2}}{\ol{1}^{\dn{2}}})}&
{{2^{2}}_{\llap{$\scriptstyle{a}$}{b}}}\hskip\inclskip\quad\hfil%
{({1}{\ol{1}}{\ol{2}})}&
{{4}_{adb}}\hskip\inclskip\quad\hfil%
{({\ol{1}}{\ol{3}})}&
{-({1}{3}_{e})}\hskip\inclskip\hskip-\inclskip\hfil\cr\rg%
{({4})}&
{{4}_{adb}}\hskip\inclskip\quad\hfil%
{({1^{2}}{\ol{2}})}&
{{1^{2}}{2^{-1}}{4}}\hskip\inclskip\quad\hfil%
{({1}{\ol{3}})}&
{-({1^{-1}}{2}{3}_{{e}\tau})}\hskip\inclskip\quad\hfil%
{({\ol{4}})}&
{{4^{-1}}{8}}\hskip\inclskip\hskip-\inclskip\hfil\cr\rg%
}}%
\tableclose%

\def\bol{\quad}

\bigskip
\vskip 10 true pt plus 20 true pt
\leftline{\bf Type $G_2$.}

\smallskip
The extended diagram for $G_2$ appears below.
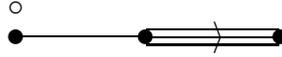
\begingroup
\def\star{\lower 3 true pt\hbox to 0pt{\hss{$\circ$}\hss}}
\def\bgt{\lower 2 true pt\hbox{$\rangle$}}
\def\hbarloc{\hbox to 27 true pt{\hss\leaders\hrule\hskip 50 true
pt\hss}}
\def\vbar{\vbox to 30 true pt{\vss\hbox to 0pt{\hss\vrule height
38 true pt\hss}\vss}}
\def\triple{\lower 0.2 true pt\vbox to 6.3 true pt{%
\hbarloc\vss\hbarloc\vss\hbarloc}}
\def\tri{\hbox{\triple}\llap{$\rangle$\hskip 9 true pt}}
$$
\matrix{
\star&&&&\cr
\spot&\hb&
\spot&\tri&
\spot\cr
}
$$
\endgroup

\vskip 5 true pt
\medskip
\noindent
The class inclusions for type $G_2$ are listed in 
Tables~21 and~22.
Recall that the component

\medskip

\inclskip=30 true pt
\vbox{%
\tableopen{Class inclusions for $W({A_{1}}{A_{1}^{*}})\,\subset\,W(G_{2})$}%
\nobreak\vskip 10 true pt\nobreak%
\baselineskip=\tablebaseskip%
\tabskip=\incltabskip\halign to \hsize{%
\bol\hfil$#$&%
$\,\subset\,#$&&%
$\,\subset\,#$\cr%
{({1^{2}}){\times}({1^{2}})}&
{{1^{2}}}\hskip\inclskip\quad\hfil%
{({2}){\times}({1^{2}})}&
{{2}_{\ell}}\hskip\inclskip\quad\hfil%
{({1^{2}}){\times}({2})}&
{-({2}_{\ell})}\hskip\inclskip\quad\hfil%
{({2}){\times}({2})}&
{-({1^{2}})}\hskip\inclskip\hskip-\inclskip\hfil\cr\rg%
}}%
\tableclose%

\def\bol{\qquad\qquad}
\inclskip=110 true pt

\medskip
\vbox{%
\tableopen{Class inclusions for $W(A_{2})\,\subset\,W(G_{2})$}%
\nobreak\vskip 10 true pt\nobreak%
\baselineskip=\tablebaseskip%
\tabskip=\incltabskip\halign to \hsize{%
\bol\hfil$#$&%
$\,\subset\,#$&&%
$\,\subset\,#$\cr%
{({1^{3}})}&
{{1^{2}}}\hskip\inclskip\quad\hfil%
{({1}{2})}&
{{2}_{\ell}}\hskip\inclskip\quad\hfil%
{({3})}&
{{1^{-1}}{3}}\hskip\inclskip\hskip-\inclskip\hfil\cr\rg%
}}%
\tableclose%
\def\rg{\rginit}

\vfil
\eject
\sectionhead{Irreducible Characters}

\bigskip
This section contains a brief review of the classification
of irreducible characters of finite
Weyl groups.
More information about the classical cases
can be found in Chapters 2 and 4 of [JK], for example.
(See also [Y].)
The notation established here is required in
the next section.

\medskip
\leftline{\bf Type $A$.}
\smallskip
The irreducible characters of $W = W(A_{\ell})$ are parametrized by
partitions of $\ell+1$, as follows.
Let $\alpha$ be a partition of $\ell+1$, and
denote by
$W(\alpha)$
a parabolic subgroup
of $W$ whose Coxeter elements lie in the
class $(\alpha)$ of $W$.
There is a unique irreducible character
of $W$ appearing as a constituent in
both
$Ind^{W}_{W(\alpha)}1$
and
$Ind^{W}_{W(\alpha^{*})}\varepsilon$,
where
$\alpha^{*}$ is the conjugate of $\alpha$;
this irreducible character
will be denoted
$[\alpha]$.
The unit character of
$W$ is
$[\ell+1]$,
and
$[\ell\,1]$ is
the character of the natural
reflection representation of $W$.

\medskip
\leftline{\bf Types $B$ and $C$.}
\smallskip
For $0\, \le \,k\, \le \,\ell$,
define a linear character
$\lambda_{k}$
of $\Delta = \Delta_{\ell}$ by
$$\lambda_{k}(diag(\varepsilon_1,\dots,\varepsilon_{\ell})) =
\varepsilon_{k+1}\cdots\varepsilon_{\ell}{.}$$
The stabilizer of $\lambda_{k}$ in $W(B_{\ell})$ is
the semi-direct product
$\Delta\cdot(P_{k}\times P_{\ell-k})$, where
$P_{k}\times P_{\ell-k}$ 
is a parabolic subgroup
of $P = P_{\ell} = W(A_{\ell-1})$.
Let $\alpha$ and $\beta$ be partitions of $k$ and
$\ell-k$, respectively.
Denote by $\widetilde{\lambda_{k}}$ the extension of
$\lambda_{k}$ to
$\Delta\cdot(P_{k}\times P_{\ell-k})$
such that $P_{k}\times P_{\ell-k}\,\subset\,ker(\widetilde{\lambda_{k}})$.
Let $([\alpha]\times[\beta])\widetilde{\ \ }$ be the extension
of $[\alpha]\times[\beta]$ to
$\Delta\cdot(P_{k}\times P_{\ell-k})$
with $\Delta$ in the kernel.
Define
$[\alpha:\beta]$
to be the character of
$W(B_{\ell})$ induced from the character
$\widetilde{\lambda_{k}}\cdot([\alpha]\times[\beta])\widetilde{\ \ }$
of $\Delta\cdot(P_{k}\times P_{\ell-k})$.
(In the extreme cases $k = 0$ and $k = \ell$, we allow
$\alpha$ or $\beta$ to be the
empty partition of 0,
denoted by `$-$'.
Thus $[\alpha:{-}]$ is the
extension of the irreducible
character $[\alpha]$ of $P$
to $W(B_{\ell})$ with $\Delta$ in the kernel.
Also,
$[{-}:\beta] = \widetilde{\lambda_{0}}\cdot[\beta:{-}] =
[\beta^{*}:{-}]^{*}$.)
By Clifford's Theorem,
$[\alpha:\beta]$ is an
irreducible character of $W(B_{\ell})$;
moreover, each irreducible
character of $W(B_{\ell})$ is obtained exactly once
via this construction.
(See [CR], Proposition 1.18.)
For example, the
character of the natural reflection
representation is $[\ell-1:1]$, and
the sign character of $W(B_{\ell})$
is
$[{-}:1^{\ell}]$.

\smallskip
The same notation will be used
for the irreducible characters
of $W(C_{\ell})$ as is used
for those of $W(B_{\ell})$.

\medskip
\goodbreak
\leftline{\bf Type $D$.}
\smallskip
The irreducible characters
$[\alpha:\beta]$ and
$[\beta:\alpha]$ of $W(B_{\ell})$ have the
same restriction to the subgroup $W(D_{\ell})$.  Since
$W(D_{\ell})$ has index 2 in $W(B_{\ell})$, the
restriction of
$[\alpha:\beta]$ to $W(D_{\ell})$ is
therefore irreducible unless $\alpha = \beta$.
If $\ell$ is even and $\alpha$ is a partition of
$k = \ell/2$, then the restriction of
$[\alpha:\alpha]$
to $W(D_{\ell})$ is the sum of two distinct irreducible
characters of $W(D_{\ell})$,
which we denote by
$[\alpha:\alpha]^{+}$ and
$[\alpha:\alpha]^{-}$.
These characters of $W(D_{\ell})$ agree
on any class of $W(B_{\ell})$ that
does not split in $W(D_{\ell})$.
Clifford's Theorem can be used to compute
the values of
$[\alpha:\alpha]^{\pm}$
on split classes.
Assume $\ell$ is even, and let
$\mu = {2^{e_2}}{4^{e_4}}\dots$
be a
signed
partition of $\ell$ with only even unsigned parts.
Let
$\mu_0 = {1^{e_2}}{2^{e_4}}\dots$ be the partition of $k$
obtained from $\mu$ by dividing each part by 2.
We may choose notation in such a
way that the values of
$[\alpha:\alpha]^{+}$
on $(\mu)_{\pm}$
are given by
$$2[\alpha:\alpha]^{+}\big|_{(\mu)_{+}} =
[\alpha:\alpha]\big|_{(\mu)} +
{2^{e}}[\alpha]\big|_{(\mu_0)}$$
and
$$2[\alpha:\alpha]^{+}\big|_{(\mu)_{-}} =
[\alpha:\alpha]\big|_{(\mu)} -
{2^{e}}[\alpha]\big|_{(\mu_0)}$$
where $e = {e_2}+{e_4}+\dots$
and $[\alpha]$ is the irreducible
character of
$W(A_{k-1})$ corresponding
to $\alpha$.
The values of
$[\alpha:\alpha]^{-}$
are then completely determined,
since
$[\alpha:\alpha]$
agrees with
$[\alpha:\alpha]^{+}+[\alpha:\alpha]^{-}$
on $W(D_{\ell})$.

\medskip
\leftline{\bf Types $E_6$, $E_7$ and $E_8$.}
\smallskip
The irreducible characters of
$W(E_6)$, $W(E_7)$,
and $W(E_8)$ were determined by
Frame in [F1] and [F2].
(See also the minor corrections noted in [BL].)
We adhere to the notation for irreducible
characters that appear in the
Tables II and III of [F1] and Tables 2 through 5 of [F2],
denoting
each character by its degree, with
subscripts
used to distinguish characters of equal degree.
Fifteen of the 25 irreducible characters of
$W(E_6)$ appear in Table II of [F1];
the remaining 10 irreducible characters
are the duals of
$1_p$, $6_p$, \dots, $60_p$, and will
be denoted
$1_p^{*}$, $6_p^{*}$, \dots, $60_p^{*}$.
Similar remarks apply
to $W(E_8)$.
The irreducible characters of the
derived group $W(E_7)^{\prime}$
appear in Table III of [F1].
We use the same notation
for an irreducible character of
$W(E_7)^{\prime}$
and its
extension to $W(E_7)$ with
$\{\pm1\}$ in the kernel,
obtaining 30 irreducible characters
$1_a$, $7_a$, \dots, $70_a$
of $W(E_7)$;
the
remaining 30 irreducible characters
are $1_a^{*}$, $7_{a}^{*},$, \dots, $70_a^{*}$.

\smallskip
In this notation, the unit (sign; reflection)
characters of $W(E_6)$,
$W(E_7)$ and $W(E_8)$
are
$1_p$, $1_a$, and $1_x$
($1_p^{*}$, $1_a^{*}$, and $1_x^{*}$;
$6_p$, $7_a^{*}$, and $8_z$),
respectively.

\medskip
\goodbreak
\leftline{\bf Type $F_4$.}
\smallskip
We denote by $\chi_{d,k}$ the $k$-th irreducible character of
degree $d$ according to the arrangement in the character
table of [K].
Thus $\chi_{1,1}$ is the unit character of $W(F_4)$, while
$\chi_{4,2}$ is the character of the natural reflection
representation.

\medskip
\goodbreak
\leftline{\bf Type $G_2$.}
\smallskip
The irreducible characters of
$W(G_2)$ will be denoted
$\chi_{1,1}$,
$\chi_{1,2}$,
$\chi_{1,3}$,
$\chi_{1,4}$,
$\chi_{2,1}$, and
$\chi_{2,2}$.
The unit and sign characters
are $\chi_{1,1}$ and
$\chi_{1,2}$, respectively.
The linear characters
$\chi_{1,3}$ and
$\chi_{1,4}$ assume the
values $-1$ and $1$ on
$2_{\ell}$, respectively.
The character of the natural
reflection representation
is $\chi_{2,1}$.
Finally,
$\chi_{2,2} =
\chi_{1,3}
\!\cdot\!
\chi_{2,1}$.

\vfil
\eject
\sectionhead{Induce/Restrict Matrices}

\bigskip
This section contains the induce/restrict matrices for
$W_0 = W(\Gamma_0)\,\subset\,W = W(\Gamma)$, where
$\Gamma_0$ is a maximal
subdiagram of either
$\widetilde{\Gamma}$ or
$\Gamma$.
The rows of the induce/restrict matrices
are parametrized by the
irreducible characters $\chi$ of
$W$, while the columns are
parametrized by the irreducible
characters $\varphi$ of $W_0$.
The entry in the row
corresponding to $\chi$ and the
column corresponding to
$\varphi$
is the inner product
$(\chi,Ind^{W}_{W_0}\varphi)_{\down{2 true pt}{W}}$.
These inner products were evaluated using a computer.
A computer was also used to process the data into
the form of the tables below.

\smallskip
In types $E_6$, $E_7$, and $E_8$ the size of each of the tables
was reduced considerably by including only one row for
the pair
$\chi$, $\chi^{*} \in Irr(W)$
whenever
$\chi \not= \chi^{*}$.
Note that
$$(\chi^{*},Ind^{W}_{W_0}\varphi)_{\down{2 true pt}{W}} =
(\chi,Ind^{W}_{W_0}{\varphi^{*}})_{\down{2 true pt}{W}}{,}$$
and hence any entries not given below can be recovered using a
knowledge of the duality mapping $\varphi \mapsto\varphi^{*}$
on $Irr(W_0)$.
We now describe this mapping
in more detail in the necessary cases.

\smallskip
In type $A_{\ell}$, we have
$[\alpha]^{*} = [\alpha^{*}]$, where
$\alpha^{*}$ is the conjugate partition of $\alpha$.
In type $B_{\ell}$ (or $C_{\ell}$),
$[\alpha:\beta]^{*} = [\beta^{*}:\alpha^{*}]$.
The same formula applies in type $D_{\ell}$, where in
addition we have $[\beta^{*}:\alpha^{*}] =
[\alpha^{*}:\beta^{*}]$.
In case $\ell$ is even and $\alpha$ is a partition
of $\ell/2$, we have
\medskip
$$([\alpha:\alpha]^{\pm})^{*} = \cases{
[\alpha^{*}:\alpha^{*}]^{\pm},&if $\ell\equiv0\!\pmod{4}$;\cr
\vbox to 8 true pt{}&{}\cr
[\alpha^{*}:\alpha^{*}]^{\mp},&if $\ell\equiv2\!\pmod{4}$.\cr}$$
\bigskip
\noindent
The duality mapping has been described in the previous section in
types $E_6$ and $E_7$.
Finally, if $\Gamma_1$, \dots, $\Gamma_k$ are the
connected components of $\Gamma$ and
$\chi = \chi_1\times\dots\times\chi_k$ where
$\chi_i \in Irr(W(\Gamma_i))$
for $1 \le i \le k$, then
$\chi^{*} = \chi_1^{*}\times\dots\times\chi_k^{*}$.

\goodbreak

\smallskip
In case $W$ possesses a
nonidentity graph automorphism 
corresponding to a symmetry of
$\widetilde{\Gamma}$,
only one representative of each orbit
of the graph automorphisms
on the set of maximal
subdiagrams of
$\widetilde{\Gamma}$
is considered.
For example, only one subgroup of type
${D_6}{A_1}$ is treated in type $E_7$.
Similarly, if
the Coxeter graph
corresponding to $\Gamma$
possesses a nontrivial symmetry,
then only one representative from
each orbit of the graph
automorphisms on the maximal
subdiagrams of
$\Gamma$
is treated.
For example,
in type $F_4$,
only one parabolic subgroup of type
${A_2}{A_1}$ is considered;
this subgroup is
$W({A_2}{A_1^{*}}) =
W(A_2)\times W(A_1^{*})$,
where
$W(A_1^{*})$
contains a reflection with respect
to a short root.
In type $G_2$ only one
parabolic subgroup of type
$A_1$ is treated; this subgroup
contains a reflection with respect
to a long root.

\smallskip
To compute the induce/restrict matrix for
a case not appearing below,
it suffices to
know the action of the graph automorphisms
on $Irr(W)$.
This action is trivial in
types $E_6$, $E_7$ and $E_8$.
The nontrivial orbits of the graph automorphisms of
$W(F_4)$ on $Irr(W(F_4))$ are
$\{\chi_{1,2},\chi_{1,3}\}$,
$\{\chi_{2,1},\chi_{2,3}\}$,
$\{\chi_{2,2},\chi_{2,4}\}$,
$\{\chi_{9,2},\chi_{9,3}\}$,
$\{\chi_{4,3},\chi_{4,4}\}$,
$\{\chi_{8,1},\chi_{8,3}\}$,
and
$\{\chi_{8,2},\chi_{8,4}\}$.
In type $G_2$, the nonidentity
graph automorphism interchanges
the linear characters
$\chi_{1,3}$ and
$\chi_{1,4}$, leaving
invariant the other
irreducible characters
of $W(G_2)$.

\vfil
\eject
%
%
%
%
%
%
\tableopen{Induce/restrict matrix for $W({A_{5}}{A_{1}})\,\subset\,W(E_{6})$}%
%
%
%
%
%
%
\rowpts=18 true pt%
\colpts=18 true pt%
\rowlabpts=40 true pt%
\collabpts=70 true pt%
\clx{\vss\hfull{%
\rlx{\hss{$ $}}\cg%
\cx{\hskip 16 true pt\flip{$[{6}]{\times}[{2}]$}\hss}\cg%
\cx{\hskip 16 true pt\flip{$[{5}{1}]{\times}[{2}]$}\hss}\cg%
\cx{\hskip 16 true pt\flip{$[{4}{2}]{\times}[{2}]$}\hss}\cg%
\cx{\hskip 16 true pt\flip{$[{4}{1^{2}}]{\times}[{2}]$}\hss}\cg%
\cx{\hskip 16 true pt\flip{$[{3^{2}}]{\times}[{2}]$}\hss}\cg%
\cx{\hskip 16 true pt\flip{$[{3}{2}{1}]{\times}[{2}]$}\hss}\cg%
\cx{\hskip 16 true pt\flip{$[{3}{1^{3}}]{\times}[{2}]$}\hss}\cg%
\cx{\hskip 16 true pt\flip{$[{2^{3}}]{\times}[{2}]$}\hss}\cg%
\cx{\hskip 16 true pt\flip{$[{2^{2}}{1^{2}}]{\times}[{2}]$}\hss}\cg%
\cx{\hskip 16 true pt\flip{$[{2}{1^{4}}]{\times}[{2}]$}\hss}\cg%
\cx{\hskip 16 true pt\flip{$[{1^{6}}]{\times}[{2}]$}\hss}\cg%
\cx{\hskip 16 true pt\flip{$[{6}]{\times}[{1^{2}}]$}\hss}\cg%
\cx{\hskip 16 true pt\flip{$[{5}{1}]{\times}[{1^{2}}]$}\hss}\cg%
\cx{\hskip 16 true pt\flip{$[{4}{2}]{\times}[{1^{2}}]$}\hss}\cg%
\cx{\hskip 16 true pt\flip{$[{4}{1^{2}}]{\times}[{1^{2}}]$}\hss}\cg%
\cx{\hskip 16 true pt\flip{$[{3^{2}}]{\times}[{1^{2}}]$}\hss}\cg%
\cx{\hskip 16 true pt\flip{$[{3}{2}{1}]{\times}[{1^{2}}]$}\hss}\cg%
\cx{\hskip 16 true pt\flip{$[{3}{1^{3}}]{\times}[{1^{2}}]$}\hss}\cg%
\cx{\hskip 16 true pt\flip{$[{2^{3}}]{\times}[{1^{2}}]$}\hss}\cg%
\cx{\hskip 16 true pt\flip{$[{2^{2}}{1^{2}}]{\times}[{1^{2}}]$}\hss}\cg%
\cx{\hskip 16 true pt\flip{$[{2}{1^{4}}]{\times}[{1^{2}}]$}\hss}\cg%
\cx{\hskip 16 true pt\flip{$[{1^{6}}]{\times}[{1^{2}}]$}\hss}\cg%
\eol}}\rg%
%
%
\rx{\vss\hfull{%
\rlx{\hss{$1_p$}}\cg%
\e{1}%
\e{0}%
\e{0}%
\e{0}%
\e{0}%
\e{0}%
\e{0}%
\e{0}%
\e{0}%
\e{0}%
\e{0}%
\e{0}%
\e{0}%
\e{0}%
\e{0}%
\e{0}%
\e{0}%
\e{0}%
\e{0}%
\e{0}%
\e{0}%
\e{0}%
\eol}\vss}\rg%
%
%
\rx{\vss\hfull{%
\rlx{\hss{$6_p$}}\cg%
\e{0}%
\e{1}%
\e{0}%
\e{0}%
\e{0}%
\e{0}%
\e{0}%
\e{0}%
\e{0}%
\e{0}%
\e{0}%
\e{1}%
\e{0}%
\e{0}%
\e{0}%
\e{0}%
\e{0}%
\e{0}%
\e{0}%
\e{0}%
\e{0}%
\e{0}%
\eol}\vss}\rg%
%
%
\rx{\vss\hfull{%
\rlx{\hss{$15_p$}}\cg%
\e{0}%
\e{0}%
\e{0}%
\e{1}%
\e{0}%
\e{0}%
\e{0}%
\e{0}%
\e{0}%
\e{0}%
\e{0}%
\e{0}%
\e{1}%
\e{0}%
\e{0}%
\e{0}%
\e{0}%
\e{0}%
\e{0}%
\e{0}%
\e{0}%
\e{0}%
\eol}\vss}\rg%
%
%
\rx{\vss\hfull{%
\rlx{\hss{$20_p$}}\cg%
\e{1}%
\e{1}%
\e{1}%
\e{0}%
\e{0}%
\e{0}%
\e{0}%
\e{0}%
\e{0}%
\e{0}%
\e{0}%
\e{0}%
\e{1}%
\e{0}%
\e{0}%
\e{0}%
\e{0}%
\e{0}%
\e{0}%
\e{0}%
\e{0}%
\e{0}%
\eol}\vss}\rg%
%
%
\rx{\vss\hfull{%
\rlx{\hss{$30_p$}}\cg%
\e{0}%
\e{1}%
\e{0}%
\e{1}%
\e{1}%
\e{0}%
\e{0}%
\e{0}%
\e{0}%
\e{0}%
\e{0}%
\e{1}%
\e{0}%
\e{1}%
\e{0}%
\e{0}%
\e{0}%
\e{0}%
\e{0}%
\e{0}%
\e{0}%
\e{0}%
\eol}\vss}\rg%
%
%
\rx{\vss\hfull{%
\rlx{\hss{$64_p$}}\cg%
\e{0}%
\e{1}%
\e{1}%
\e{1}%
\e{0}%
\e{1}%
\e{0}%
\e{0}%
\e{0}%
\e{0}%
\e{0}%
\e{0}%
\e{1}%
\e{1}%
\e{1}%
\e{0}%
\e{0}%
\e{0}%
\e{0}%
\e{0}%
\e{0}%
\e{0}%
\eol}\vss}\rg%
%
%
\rx{\vss\hfull{%
\rlx{\hss{$81_p$}}\cg%
\e{0}%
\e{0}%
\e{1}%
\e{1}%
\e{0}%
\e{1}%
\e{1}%
\e{0}%
\e{0}%
\e{0}%
\e{0}%
\e{0}%
\e{1}%
\e{0}%
\e{1}%
\e{1}%
\e{1}%
\e{0}%
\e{0}%
\e{0}%
\e{0}%
\e{0}%
\eol}\vss}\rg%
%
%
\rx{\vss\hfull{%
\rlx{\hss{$15_q$}}\cg%
\e{1}%
\e{0}%
\e{1}%
\e{0}%
\e{0}%
\e{0}%
\e{0}%
\e{0}%
\e{0}%
\e{0}%
\e{0}%
\e{0}%
\e{0}%
\e{0}%
\e{0}%
\e{1}%
\e{0}%
\e{0}%
\e{0}%
\e{0}%
\e{0}%
\e{0}%
\eol}\vss}\rg%
%
%
\rx{\vss\hfull{%
\rlx{\hss{$24_p$}}\cg%
\e{0}%
\e{0}%
\e{1}%
\e{0}%
\e{0}%
\e{0}%
\e{0}%
\e{1}%
\e{0}%
\e{0}%
\e{0}%
\e{0}%
\e{0}%
\e{0}%
\e{1}%
\e{0}%
\e{0}%
\e{0}%
\e{0}%
\e{0}%
\e{0}%
\e{0}%
\eol}\vss}\rg%
%
%
\rx{\vss\hfull{%
\rlx{\hss{$60_p$}}\cg%
\e{0}%
\e{1}%
\e{1}%
\e{0}%
\e{1}%
\e{1}%
\e{0}%
\e{0}%
\e{0}%
\e{0}%
\e{0}%
\e{0}%
\e{0}%
\e{1}%
\e{0}%
\e{0}%
\e{1}%
\e{0}%
\e{0}%
\e{0}%
\e{0}%
\e{0}%
\eol}\vss}\rg%
%
%
\rx{\vss\hfull{%
\rlx{\hss{$20_s$}}\cg%
\e{0}%
\e{0}%
\e{0}%
\e{0}%
\e{0}%
\e{0}%
\e{1}%
\e{0}%
\e{0}%
\e{0}%
\e{0}%
\e{0}%
\e{0}%
\e{0}%
\e{1}%
\e{0}%
\e{0}%
\e{0}%
\e{0}%
\e{0}%
\e{0}%
\e{0}%
\eol}\vss}\rg%
%
%
\rx{\vss\hfull{%
\rlx{\hss{$90_s$}}\cg%
\e{0}%
\e{0}%
\e{0}%
\e{1}%
\e{0}%
\e{1}%
\e{1}%
\e{0}%
\e{1}%
\e{0}%
\e{0}%
\e{0}%
\e{0}%
\e{1}%
\e{1}%
\e{0}%
\e{1}%
\e{1}%
\e{0}%
\e{0}%
\e{0}%
\e{0}%
\eol}\vss}\rg%
%
%
\rx{\vss\hfull{%
\rlx{\hss{$80_s$}}\cg%
\e{0}%
\e{0}%
\e{0}%
\e{1}%
\e{1}%
\e{1}%
\e{0}%
\e{0}%
\e{1}%
\e{0}%
\e{0}%
\e{0}%
\e{0}%
\e{1}%
\e{0}%
\e{0}%
\e{1}%
\e{1}%
\e{1}%
\e{0}%
\e{0}%
\e{0}%
\eol}\vss}\rg%
%
%
\rx{\vss\hfull{%
\rlx{\hss{$60_s$}}\cg%
\e{0}%
\e{0}%
\e{1}%
\e{0}%
\e{0}%
\e{1}%
\e{0}%
\e{1}%
\e{0}%
\e{0}%
\e{0}%
\e{0}%
\e{0}%
\e{0}%
\e{0}%
\e{1}%
\e{1}%
\e{0}%
\e{0}%
\e{1}%
\e{0}%
\e{0}%
\eol}\vss}\rg%
%
%
\rx{\vss\hfull{%
\rlx{\hss{$10_s$}}\cg%
\e{0}%
\e{0}%
\e{0}%
\e{0}%
\e{1}%
\e{0}%
\e{0}%
\e{0}%
\e{0}%
\e{0}%
\e{0}%
\e{0}%
\e{0}%
\e{0}%
\e{0}%
\e{0}%
\e{0}%
\e{0}%
\e{1}%
\e{0}%
\e{0}%
\e{0}%
\eol}\vss}\rg%
\tableclose%
%
%
%
%
%
%
\tableopen{Induce/restrict matrix for $W({A_{2}}{A_{2}}{A_{2}})\,\subset\,W(E_{6})$}%
%
%
%
%
%
%
\rowpts=18 true pt%
\colpts=18 true pt%
\rowlabpts=40 true pt%
\collabpts=90 true pt%
\clx{\vss\hfull{%
\rlx{\hss{$ $}}\cg%
\cx{\hskip 16 true pt\flip{$[{3}]{\times}[{3}]{\times}[{3}]$}\hss}\cg%
\cx{\hskip 16 true pt\flip{$[{2}{1}]{\times}[{3}]{\times}[{3}]$}\hss}\cg%
\cx{\hskip 16 true pt\flip{$[{1^{3}}]{\times}[{3}]{\times}[{3}]$}\hss}\cg%
\cx{\hskip 16 true pt\flip{$[{3}]{\times}[{2}{1}]{\times}[{3}]$}\hss}\cg%
\cx{\hskip 16 true pt\flip{$[{2}{1}]{\times}[{2}{1}]{\times}[{3}]$}\hss}\cg%
\cx{\hskip 16 true pt\flip{$[{1^{3}}]{\times}[{2}{1}]{\times}[{3}]$}\hss}\cg%
\cx{\hskip 16 true pt\flip{$[{3}]{\times}[{1^{3}}]{\times}[{3}]$}\hss}\cg%
\cx{\hskip 16 true pt\flip{$[{2}{1}]{\times}[{1^{3}}]{\times}[{3}]$}\hss}\cg%
\cx{\hskip 16 true pt\flip{$[{1^{3}}]{\times}[{1^{3}}]{\times}[{3}]$}\hss}\cg%
\cx{\hskip 16 true pt\flip{$[{3}]{\times}[{3}]{\times}[{2}{1}]$}\hss}\cg%
\cx{\hskip 16 true pt\flip{$[{2}{1}]{\times}[{3}]{\times}[{2}{1}]$}\hss}\cg%
\cx{\hskip 16 true pt\flip{$[{1^{3}}]{\times}[{3}]{\times}[{2}{1}]$}\hss}\cg%
\cx{\hskip 16 true pt\flip{$[{3}]{\times}[{2}{1}]{\times}[{2}{1}]$}\hss}\cg%
\cx{\hskip 16 true pt\flip{$[{2}{1}]{\times}[{2}{1}]{\times}[{2}{1}]$}\hss}\cg%
\eol}}\rg%
%
%
\rx{\vss\hfull{%
\rlx{\hss{$1_p$}}\cg%
\e{1}%
\e{0}%
\e{0}%
\e{0}%
\e{0}%
\e{0}%
\e{0}%
\e{0}%
\e{0}%
\e{0}%
\e{0}%
\e{0}%
\e{0}%
\e{0}%
\eol}\vss}\rg%
%
%
\rx{\vss\hfull{%
\rlx{\hss{$6_p$}}\cg%
\e{0}%
\e{1}%
\e{0}%
\e{1}%
\e{0}%
\e{0}%
\e{0}%
\e{0}%
\e{0}%
\e{1}%
\e{0}%
\e{0}%
\e{0}%
\e{0}%
\eol}\vss}\rg%
%
%
\rx{\vss\hfull{%
\rlx{\hss{$15_p$}}\cg%
\e{0}%
\e{0}%
\e{1}%
\e{0}%
\e{1}%
\e{0}%
\e{1}%
\e{0}%
\e{0}%
\e{0}%
\e{1}%
\e{0}%
\e{1}%
\e{0}%
\eol}\vss}\rg%
%
%
\rx{\vss\hfull{%
\rlx{\hss{$20_p$}}\cg%
\e{2}%
\e{1}%
\e{0}%
\e{1}%
\e{1}%
\e{0}%
\e{0}%
\e{0}%
\e{0}%
\e{1}%
\e{1}%
\e{0}%
\e{1}%
\e{0}%
\eol}\vss}\rg%
%
%
\rx{\vss\hfull{%
\rlx{\hss{$30_p$}}\cg%
\e{1}%
\e{1}%
\e{1}%
\e{1}%
\e{1}%
\e{0}%
\e{1}%
\e{0}%
\e{0}%
\e{1}%
\e{1}%
\e{0}%
\e{1}%
\e{1}%
\eol}\vss}\rg%
%
%
\rx{\vss\hfull{%
\rlx{\hss{$64_p$}}\cg%
\e{0}%
\e{2}%
\e{0}%
\e{2}%
\e{2}%
\e{1}%
\e{0}%
\e{1}%
\e{0}%
\e{2}%
\e{2}%
\e{1}%
\e{2}%
\e{2}%
\eol}\vss}\rg%
%
%
\rx{\vss\hfull{%
\rlx{\hss{$81_p$}}\cg%
\e{0}%
\e{1}%
\e{1}%
\e{1}%
\e{2}%
\e{1}%
\e{1}%
\e{1}%
\e{0}%
\e{1}%
\e{2}%
\e{1}%
\e{2}%
\e{3}%
\eol}\vss}\rg%
%
%
\rx{\vss\hfull{%
\rlx{\hss{$15_q$}}\cg%
\e{1}%
\e{1}%
\e{0}%
\e{1}%
\e{0}%
\e{0}%
\e{0}%
\e{0}%
\e{0}%
\e{1}%
\e{0}%
\e{0}%
\e{0}%
\e{1}%
\eol}\vss}\rg%
%
%
\rx{\vss\hfull{%
\rlx{\hss{$24_p$}}\cg%
\e{1}%
\e{0}%
\e{0}%
\e{0}%
\e{1}%
\e{0}%
\e{0}%
\e{0}%
\e{1}%
\e{0}%
\e{1}%
\e{0}%
\e{1}%
\e{1}%
\eol}\vss}\rg%
%
%
\rx{\vss\hfull{%
\rlx{\hss{$60_p$}}\cg%
\e{2}%
\e{1}%
\e{0}%
\e{1}%
\e{2}%
\e{0}%
\e{0}%
\e{0}%
\e{0}%
\e{1}%
\e{2}%
\e{0}%
\e{2}%
\e{2}%
\eol}\vss}\rg%
\eop
\eject
\tablecont%
%
%
%
%
%
%
\rowpts=18 true pt%
\colpts=18 true pt%
\rowlabpts=40 true pt%
\collabpts=90 true pt%
\clx{\vss\hfull{%
\rlx{\hss{$ $}}\cg%
\cx{\hskip 16 true pt\flip{$[{3}]{\times}[{3}]{\times}[{3}]$}\hss}\cg%
\cx{\hskip 16 true pt\flip{$[{2}{1}]{\times}[{3}]{\times}[{3}]$}\hss}\cg%
\cx{\hskip 16 true pt\flip{$[{1^{3}}]{\times}[{3}]{\times}[{3}]$}\hss}\cg%
\cx{\hskip 16 true pt\flip{$[{3}]{\times}[{2}{1}]{\times}[{3}]$}\hss}\cg%
\cx{\hskip 16 true pt\flip{$[{2}{1}]{\times}[{2}{1}]{\times}[{3}]$}\hss}\cg%
\cx{\hskip 16 true pt\flip{$[{1^{3}}]{\times}[{2}{1}]{\times}[{3}]$}\hss}\cg%
\cx{\hskip 16 true pt\flip{$[{3}]{\times}[{1^{3}}]{\times}[{3}]$}\hss}\cg%
\cx{\hskip 16 true pt\flip{$[{2}{1}]{\times}[{1^{3}}]{\times}[{3}]$}\hss}\cg%
\cx{\hskip 16 true pt\flip{$[{1^{3}}]{\times}[{1^{3}}]{\times}[{3}]$}\hss}\cg%
\cx{\hskip 16 true pt\flip{$[{3}]{\times}[{3}]{\times}[{2}{1}]$}\hss}\cg%
\cx{\hskip 16 true pt\flip{$[{2}{1}]{\times}[{3}]{\times}[{2}{1}]$}\hss}\cg%
\cx{\hskip 16 true pt\flip{$[{1^{3}}]{\times}[{3}]{\times}[{2}{1}]$}\hss}\cg%
\cx{\hskip 16 true pt\flip{$[{3}]{\times}[{2}{1}]{\times}[{2}{1}]$}\hss}\cg%
\cx{\hskip 16 true pt\flip{$[{2}{1}]{\times}[{2}{1}]{\times}[{2}{1}]$}\hss}\cg%
\eol}}\rg%
%
%
\rx{\vss\hfull{%
\rlx{\hss{$20_s$}}\cg%
\e{0}%
\e{0}%
\e{0}%
\e{0}%
\e{0}%
\e{1}%
\e{0}%
\e{1}%
\e{0}%
\e{0}%
\e{0}%
\e{1}%
\e{0}%
\e{1}%
\eol}\vss}\rg%
%
%
\rx{\vss\hfull{%
\rlx{\hss{$90_s$}}\cg%
\e{0}%
\e{0}%
\e{1}%
\e{0}%
\e{2}%
\e{1}%
\e{1}%
\e{1}%
\e{1}%
\e{0}%
\e{2}%
\e{1}%
\e{2}%
\e{3}%
\eol}\vss}\rg%
%
%
\rx{\vss\hfull{%
\rlx{\hss{$80_s$}}\cg%
\e{0}%
\e{1}%
\e{0}%
\e{1}%
\e{1}%
\e{1}%
\e{0}%
\e{1}%
\e{0}%
\e{1}%
\e{1}%
\e{1}%
\e{1}%
\e{4}%
\eol}\vss}\rg%
%
%
\rx{\vss\hfull{%
\rlx{\hss{$60_s$}}\cg%
\e{0}%
\e{1}%
\e{0}%
\e{1}%
\e{1}%
\e{0}%
\e{0}%
\e{0}%
\e{0}%
\e{1}%
\e{1}%
\e{0}%
\e{1}%
\e{3}%
\eol}\vss}\rg%
%
%
\rx{\vss\hfull{%
\rlx{\hss{$10_s$}}\cg%
\e{1}%
\e{0}%
\e{0}%
\e{0}%
\e{0}%
\e{0}%
\e{0}%
\e{0}%
\e{0}%
\e{0}%
\e{0}%
\e{0}%
\e{0}%
\e{1}%
\eol}\vss}\rg%
%
%
%
%
%
%
\rowpts=18 true pt%
\colpts=18 true pt%
\rowlabpts=40 true pt%
\collabpts=90 true pt%
\clx{\vss\hfull{%
\rlx{\hss{$ $}}\cg%
\cx{\hskip 16 true pt\flip{$[{1^{3}}]{\times}[{2}{1}]{\times}[{2}{1}]$}\hss}\cg%
\cx{\hskip 16 true pt\flip{$[{3}]{\times}[{1^{3}}]{\times}[{2}{1}]$}\hss}\cg%
\cx{\hskip 16 true pt\flip{$[{2}{1}]{\times}[{1^{3}}]{\times}[{2}{1}]$}\hss}\cg%
\cx{\hskip 16 true pt\flip{$[{1^{3}}]{\times}[{1^{3}}]{\times}[{2}{1}]$}\hss}\cg%
\cx{\hskip 16 true pt\flip{$[{3}]{\times}[{3}]{\times}[{1^{3}}]$}\hss}\cg%
\cx{\hskip 16 true pt\flip{$[{2}{1}]{\times}[{3}]{\times}[{1^{3}}]$}\hss}\cg%
\cx{\hskip 16 true pt\flip{$[{1^{3}}]{\times}[{3}]{\times}[{1^{3}}]$}\hss}\cg%
\cx{\hskip 16 true pt\flip{$[{3}]{\times}[{2}{1}]{\times}[{1^{3}}]$}\hss}\cg%
\cx{\hskip 16 true pt\flip{$[{2}{1}]{\times}[{2}{1}]{\times}[{1^{3}}]$}\hss}\cg%
\cx{\hskip 16 true pt\flip{$[{1^{3}}]{\times}[{2}{1}]{\times}[{1^{3}}]$}\hss}\cg%
\cx{\hskip 16 true pt\flip{$[{3}]{\times}[{1^{3}}]{\times}[{1^{3}}]$}\hss}\cg%
\cx{\hskip 16 true pt\flip{$[{2}{1}]{\times}[{1^{3}}]{\times}[{1^{3}}]$}\hss}\cg%
\cx{\hskip 16 true pt\flip{$[{1^{3}}]{\times}[{1^{3}}]{\times}[{1^{3}}]$}\hss}\cg%
\eol}}\rg%
%
%
\rx{\vss\hfull{%
\rlx{\hss{$1_p$}}\cg%
\e{0}%
\e{0}%
\e{0}%
\e{0}%
\e{0}%
\e{0}%
\e{0}%
\e{0}%
\e{0}%
\e{0}%
\e{0}%
\e{0}%
\e{0}%
\eol}\vss}\rg%
%
%
\rx{\vss\hfull{%
\rlx{\hss{$6_p$}}\cg%
\e{0}%
\e{0}%
\e{0}%
\e{0}%
\e{0}%
\e{0}%
\e{0}%
\e{0}%
\e{0}%
\e{0}%
\e{0}%
\e{0}%
\e{0}%
\eol}\vss}\rg%
%
%
\rx{\vss\hfull{%
\rlx{\hss{$15_p$}}\cg%
\e{0}%
\e{0}%
\e{0}%
\e{0}%
\e{1}%
\e{0}%
\e{0}%
\e{0}%
\e{0}%
\e{0}%
\e{0}%
\e{0}%
\e{0}%
\eol}\vss}\rg%
%
%
\rx{\vss\hfull{%
\rlx{\hss{$20_p$}}\cg%
\e{0}%
\e{0}%
\e{0}%
\e{0}%
\e{0}%
\e{0}%
\e{0}%
\e{0}%
\e{0}%
\e{0}%
\e{0}%
\e{0}%
\e{0}%
\eol}\vss}\rg%
%
%
\rx{\vss\hfull{%
\rlx{\hss{$30_p$}}\cg%
\e{0}%
\e{0}%
\e{0}%
\e{0}%
\e{1}%
\e{0}%
\e{0}%
\e{0}%
\e{0}%
\e{0}%
\e{0}%
\e{0}%
\e{0}%
\eol}\vss}\rg%
%
%
\rx{\vss\hfull{%
\rlx{\hss{$64_p$}}\cg%
\e{0}%
\e{1}%
\e{0}%
\e{0}%
\e{0}%
\e{1}%
\e{0}%
\e{1}%
\e{0}%
\e{0}%
\e{0}%
\e{0}%
\e{0}%
\eol}\vss}\rg%
%
%
\rx{\vss\hfull{%
\rlx{\hss{$81_p$}}\cg%
\e{1}%
\e{1}%
\e{1}%
\e{0}%
\e{1}%
\e{1}%
\e{0}%
\e{1}%
\e{1}%
\e{0}%
\e{0}%
\e{0}%
\e{0}%
\eol}\vss}\rg%
%
%
\rx{\vss\hfull{%
\rlx{\hss{$15_q$}}\cg%
\e{0}%
\e{0}%
\e{0}%
\e{0}%
\e{0}%
\e{0}%
\e{0}%
\e{0}%
\e{0}%
\e{0}%
\e{0}%
\e{0}%
\e{0}%
\eol}\vss}\rg%
%
%
\rx{\vss\hfull{%
\rlx{\hss{$24_p$}}\cg%
\e{0}%
\e{0}%
\e{0}%
\e{0}%
\e{0}%
\e{0}%
\e{1}%
\e{0}%
\e{0}%
\e{0}%
\e{1}%
\e{0}%
\e{0}%
\eol}\vss}\rg%
%
%
\rx{\vss\hfull{%
\rlx{\hss{$60_p$}}\cg%
\e{1}%
\e{0}%
\e{1}%
\e{0}%
\e{0}%
\e{0}%
\e{0}%
\e{0}%
\e{1}%
\e{0}%
\e{0}%
\e{0}%
\e{0}%
\eol}\vss}\rg%
%
%
\rx{\vss\hfull{%
\rlx{\hss{$20_s$}}\cg%
\e{0}%
\e{1}%
\e{0}%
\e{0}%
\e{0}%
\e{1}%
\e{0}%
\e{1}%
\e{0}%
\e{0}%
\e{0}%
\e{0}%
\e{0}%
\eol}\vss}\rg%
%
%
\rx{\vss\hfull{%
\rlx{\hss{$90_s$}}\cg%
\e{2}%
\e{1}%
\e{2}%
\e{0}%
\e{1}%
\e{1}%
\e{1}%
\e{1}%
\e{2}%
\e{0}%
\e{1}%
\e{0}%
\e{0}%
\eol}\vss}\rg%
%
%
\rx{\vss\hfull{%
\rlx{\hss{$80_s$}}\cg%
\e{1}%
\e{1}%
\e{1}%
\e{1}%
\e{0}%
\e{1}%
\e{0}%
\e{1}%
\e{1}%
\e{1}%
\e{0}%
\e{1}%
\e{0}%
\eol}\vss}\rg%
%
%
\rx{\vss\hfull{%
\rlx{\hss{$60_s$}}\cg%
\e{1}%
\e{0}%
\e{1}%
\e{1}%
\e{0}%
\e{0}%
\e{0}%
\e{0}%
\e{1}%
\e{1}%
\e{0}%
\e{1}%
\e{0}%
\eol}\vss}\rg%
%
%
\rx{\vss\hfull{%
\rlx{\hss{$10_s$}}\cg%
\e{0}%
\e{0}%
\e{0}%
\e{0}%
\e{0}%
\e{0}%
\e{0}%
\e{0}%
\e{0}%
\e{0}%
\e{0}%
\e{0}%
\e{1}%
\eol}\vss}\rg%
\tableclose%
%
%
%
%
%
%
\eop
\eject
\tableopen{Induce/restrict matrix for $W(D_{5})\,\subset\,W(E_{6})$}%
%
%
%
%
%
%
\rowpts=18 true pt%
\colpts=18 true pt%
\rowlabpts=40 true pt%
\collabpts=55 true pt%
\clx{\vss\hfull{%
\rlx{\hss{$ $}}\cg%
\cx{\hskip 16 true pt\flip{$[{5}:-]$}\hss}\cg%
\cx{\hskip 16 true pt\flip{$[{4}{1}:-]$}\hss}\cg%
\cx{\hskip 16 true pt\flip{$[{3}{2}:-]$}\hss}\cg%
\cx{\hskip 16 true pt\flip{$[{3}{1^{2}}:-]$}\hss}\cg%
\cx{\hskip 16 true pt\flip{$[{2^{2}}{1}:-]$}\hss}\cg%
\cx{\hskip 16 true pt\flip{$[{2}{1^{3}}:-]$}\hss}\cg%
\cx{\hskip 16 true pt\flip{$[{1^{5}}:-]$}\hss}\cg%
\cx{\hskip 16 true pt\flip{$[{4}:{1}]$}\hss}\cg%
\cx{\hskip 16 true pt\flip{$[{3}{1}:{1}]$}\hss}\cg%
\cx{\hskip 16 true pt\flip{$[{2^{2}}:{1}]$}\hss}\cg%
\cx{\hskip 16 true pt\flip{$[{2}{1^{2}}:{1}]$}\hss}\cg%
\cx{\hskip 16 true pt\flip{$[{1^{4}}:{1}]$}\hss}\cg%
\cx{\hskip 16 true pt\flip{$[{3}:{2}]$}\hss}\cg%
\cx{\hskip 16 true pt\flip{$[{3}:{1^{2}}]$}\hss}\cg%
\cx{\hskip 16 true pt\flip{$[{2}{1}:{2}]$}\hss}\cg%
\cx{\hskip 16 true pt\flip{$[{2}{1}:{1^{2}}]$}\hss}\cg%
\cx{\hskip 16 true pt\flip{$[{1^{3}}:{2}]$}\hss}\cg%
\cx{\hskip 16 true pt\flip{$[{1^{3}}:{1^{2}}]$}\hss}\cg%
\eol}}\rg%
%
%
\rx{\vss\hfull{%
\rlx{\hss{$1_p$}}\cg%
\e{1}%
\e{0}%
\e{0}%
\e{0}%
\e{0}%
\e{0}%
\e{0}%
\e{0}%
\e{0}%
\e{0}%
\e{0}%
\e{0}%
\e{0}%
\e{0}%
\e{0}%
\e{0}%
\e{0}%
\e{0}%
\eol}\vss}\rg%
%
%
\rx{\vss\hfull{%
\rlx{\hss{$6_p$}}\cg%
\e{1}%
\e{0}%
\e{0}%
\e{0}%
\e{0}%
\e{0}%
\e{0}%
\e{1}%
\e{0}%
\e{0}%
\e{0}%
\e{0}%
\e{0}%
\e{0}%
\e{0}%
\e{0}%
\e{0}%
\e{0}%
\eol}\vss}\rg%
%
%
\rx{\vss\hfull{%
\rlx{\hss{$15_p$}}\cg%
\e{0}%
\e{0}%
\e{0}%
\e{0}%
\e{0}%
\e{0}%
\e{0}%
\e{1}%
\e{0}%
\e{0}%
\e{0}%
\e{0}%
\e{0}%
\e{1}%
\e{0}%
\e{0}%
\e{0}%
\e{0}%
\eol}\vss}\rg%
%
%
\rx{\vss\hfull{%
\rlx{\hss{$20_p$}}\cg%
\e{1}%
\e{1}%
\e{0}%
\e{0}%
\e{0}%
\e{0}%
\e{0}%
\e{1}%
\e{0}%
\e{0}%
\e{0}%
\e{0}%
\e{1}%
\e{0}%
\e{0}%
\e{0}%
\e{0}%
\e{0}%
\eol}\vss}\rg%
%
%
\rx{\vss\hfull{%
\rlx{\hss{$30_p$}}\cg%
\e{0}%
\e{0}%
\e{0}%
\e{0}%
\e{0}%
\e{0}%
\e{0}%
\e{1}%
\e{1}%
\e{0}%
\e{0}%
\e{0}%
\e{1}%
\e{0}%
\e{0}%
\e{0}%
\e{0}%
\e{0}%
\eol}\vss}\rg%
%
%
\rx{\vss\hfull{%
\rlx{\hss{$64_p$}}\cg%
\e{0}%
\e{1}%
\e{0}%
\e{0}%
\e{0}%
\e{0}%
\e{0}%
\e{1}%
\e{1}%
\e{0}%
\e{0}%
\e{0}%
\e{1}%
\e{1}%
\e{1}%
\e{0}%
\e{0}%
\e{0}%
\eol}\vss}\rg%
%
%
\rx{\vss\hfull{%
\rlx{\hss{$81_p$}}\cg%
\e{0}%
\e{0}%
\e{0}%
\e{1}%
\e{0}%
\e{0}%
\e{0}%
\e{0}%
\e{1}%
\e{0}%
\e{0}%
\e{0}%
\e{1}%
\e{1}%
\e{1}%
\e{1}%
\e{0}%
\e{0}%
\eol}\vss}\rg%
%
%
\rx{\vss\hfull{%
\rlx{\hss{$15_q$}}\cg%
\e{0}%
\e{0}%
\e{1}%
\e{0}%
\e{0}%
\e{0}%
\e{0}%
\e{0}%
\e{0}%
\e{0}%
\e{0}%
\e{0}%
\e{1}%
\e{0}%
\e{0}%
\e{0}%
\e{0}%
\e{0}%
\eol}\vss}\rg%
%
%
\rx{\vss\hfull{%
\rlx{\hss{$24_p$}}\cg%
\e{0}%
\e{1}%
\e{0}%
\e{0}%
\e{0}%
\e{0}%
\e{0}%
\e{0}%
\e{0}%
\e{0}%
\e{0}%
\e{0}%
\e{0}%
\e{0}%
\e{1}%
\e{0}%
\e{0}%
\e{0}%
\eol}\vss}\rg%
%
%
\rx{\vss\hfull{%
\rlx{\hss{$60_p$}}\cg%
\e{0}%
\e{0}%
\e{1}%
\e{0}%
\e{0}%
\e{0}%
\e{0}%
\e{0}%
\e{1}%
\e{1}%
\e{0}%
\e{0}%
\e{1}%
\e{0}%
\e{1}%
\e{0}%
\e{0}%
\e{0}%
\eol}\vss}\rg%
%
%
\rx{\vss\hfull{%
\rlx{\hss{$20_s$}}\cg%
\e{0}%
\e{0}%
\e{0}%
\e{0}%
\e{0}%
\e{0}%
\e{0}%
\e{0}%
\e{0}%
\e{0}%
\e{0}%
\e{0}%
\e{0}%
\e{1}%
\e{0}%
\e{0}%
\e{1}%
\e{0}%
\eol}\vss}\rg%
%
%
\rx{\vss\hfull{%
\rlx{\hss{$90_s$}}\cg%
\e{0}%
\e{0}%
\e{0}%
\e{0}%
\e{0}%
\e{0}%
\e{0}%
\e{0}%
\e{1}%
\e{0}%
\e{1}%
\e{0}%
\e{0}%
\e{1}%
\e{1}%
\e{1}%
\e{1}%
\e{0}%
\eol}\vss}\rg%
%
%
\rx{\vss\hfull{%
\rlx{\hss{$80_s$}}\cg%
\e{0}%
\e{0}%
\e{0}%
\e{0}%
\e{0}%
\e{0}%
\e{0}%
\e{0}%
\e{1}%
\e{1}%
\e{1}%
\e{0}%
\e{0}%
\e{0}%
\e{1}%
\e{1}%
\e{0}%
\e{0}%
\eol}\vss}\rg%
%
%
\rx{\vss\hfull{%
\rlx{\hss{$60_s$}}\cg%
\e{0}%
\e{0}%
\e{1}%
\e{0}%
\e{1}%
\e{0}%
\e{0}%
\e{0}%
\e{0}%
\e{1}%
\e{0}%
\e{0}%
\e{0}%
\e{0}%
\e{1}%
\e{1}%
\e{0}%
\e{0}%
\eol}\vss}\rg%
%
%
\rx{\vss\hfull{%
\rlx{\hss{$10_s$}}\cg%
\e{0}%
\e{0}%
\e{0}%
\e{0}%
\e{0}%
\e{0}%
\e{0}%
\e{0}%
\e{0}%
\e{1}%
\e{0}%
\e{0}%
\e{0}%
\e{0}%
\e{0}%
\e{0}%
\e{0}%
\e{0}%
\eol}\vss}\rg%
\tableclose%
%
%
%
%
%
%
\tableopen{Induce/restrict matrix for $W(A_{5})\,\subset\,W(E_{6})$}%
%
%
%
%
%
%
\rowpts=18 true pt%
\colpts=18 true pt%
\rowlabpts=40 true pt%
\collabpts=45 true pt%
\clx{\vss\hfull{%
\rlx{\hss{$ $}}\cg%
\cx{\hskip 16 true pt\flip{$[{6}]$}\hss}\cg%
\cx{\hskip 16 true pt\flip{$[{5}{1}]$}\hss}\cg%
\cx{\hskip 16 true pt\flip{$[{4}{2}]$}\hss}\cg%
\cx{\hskip 16 true pt\flip{$[{4}{1^{2}}]$}\hss}\cg%
\cx{\hskip 16 true pt\flip{$[{3^{2}}]$}\hss}\cg%
\cx{\hskip 16 true pt\flip{$[{3}{2}{1}]$}\hss}\cg%
\cx{\hskip 16 true pt\flip{$[{3}{1^{3}}]$}\hss}\cg%
\cx{\hskip 16 true pt\flip{$[{2^{3}}]$}\hss}\cg%
\cx{\hskip 16 true pt\flip{$[{2^{2}}{1^{2}}]$}\hss}\cg%
\cx{\hskip 16 true pt\flip{$[{2}{1^{4}}]$}\hss}\cg%
\cx{\hskip 16 true pt\flip{$[{1^{6}}]$}\hss}\cg%
\eol}}\rg%
%
%
\rx{\vss\hfull{%
\rlx{\hss{$1_p$}}\cg%
\e{1}%
\e{0}%
\e{0}%
\e{0}%
\e{0}%
\e{0}%
\e{0}%
\e{0}%
\e{0}%
\e{0}%
\e{0}%
\eol}\vss}\rg%
%
%
\rx{\vss\hfull{%
\rlx{\hss{$6_p$}}\cg%
\e{1}%
\e{1}%
\e{0}%
\e{0}%
\e{0}%
\e{0}%
\e{0}%
\e{0}%
\e{0}%
\e{0}%
\e{0}%
\eol}\vss}\rg%
%
%
\rx{\vss\hfull{%
\rlx{\hss{$15_p$}}\cg%
\e{0}%
\e{1}%
\e{0}%
\e{1}%
\e{0}%
\e{0}%
\e{0}%
\e{0}%
\e{0}%
\e{0}%
\e{0}%
\eol}\vss}\rg%
%
%
\rx{\vss\hfull{%
\rlx{\hss{$20_p$}}\cg%
\e{1}%
\e{2}%
\e{1}%
\e{0}%
\e{0}%
\e{0}%
\e{0}%
\e{0}%
\e{0}%
\e{0}%
\e{0}%
\eol}\vss}\rg%
%
%
\rx{\vss\hfull{%
\rlx{\hss{$30_p$}}\cg%
\e{1}%
\e{1}%
\e{1}%
\e{1}%
\e{1}%
\e{0}%
\e{0}%
\e{0}%
\e{0}%
\e{0}%
\e{0}%
\eol}\vss}\rg%
%
%
\rx{\vss\hfull{%
\rlx{\hss{$64_p$}}\cg%
\e{0}%
\e{2}%
\e{2}%
\e{2}%
\e{0}%
\e{1}%
\e{0}%
\e{0}%
\e{0}%
\e{0}%
\e{0}%
\eol}\vss}\rg%
%
%
\rx{\vss\hfull{%
\rlx{\hss{$81_p$}}\cg%
\e{0}%
\e{1}%
\e{1}%
\e{2}%
\e{1}%
\e{2}%
\e{1}%
\e{0}%
\e{0}%
\e{0}%
\e{0}%
\eol}\vss}\rg%
%
%
\rx{\vss\hfull{%
\rlx{\hss{$15_q$}}\cg%
\e{1}%
\e{0}%
\e{1}%
\e{0}%
\e{1}%
\e{0}%
\e{0}%
\e{0}%
\e{0}%
\e{0}%
\e{0}%
\eol}\vss}\rg%
%
%
\rx{\vss\hfull{%
\rlx{\hss{$24_p$}}\cg%
\e{0}%
\e{0}%
\e{1}%
\e{1}%
\e{0}%
\e{0}%
\e{0}%
\e{1}%
\e{0}%
\e{0}%
\e{0}%
\eol}\vss}\rg%
%
%
\rx{\vss\hfull{%
\rlx{\hss{$60_p$}}\cg%
\e{0}%
\e{1}%
\e{2}%
\e{0}%
\e{1}%
\e{2}%
\e{0}%
\e{0}%
\e{0}%
\e{0}%
\e{0}%
\eol}\vss}\rg%
%
%
\rx{\vss\hfull{%
\rlx{\hss{$20_s$}}\cg%
\e{0}%
\e{0}%
\e{0}%
\e{1}%
\e{0}%
\e{0}%
\e{1}%
\e{0}%
\e{0}%
\e{0}%
\e{0}%
\eol}\vss}\rg%
%
%
\rx{\vss\hfull{%
\rlx{\hss{$90_s$}}\cg%
\e{0}%
\e{0}%
\e{1}%
\e{2}%
\e{0}%
\e{2}%
\e{2}%
\e{0}%
\e{1}%
\e{0}%
\e{0}%
\eol}\vss}\rg%
%
%
\rx{\vss\hfull{%
\rlx{\hss{$80_s$}}\cg%
\e{0}%
\e{0}%
\e{1}%
\e{1}%
\e{1}%
\e{2}%
\e{1}%
\e{1}%
\e{1}%
\e{0}%
\e{0}%
\eol}\vss}\rg%
%
%
\rx{\vss\hfull{%
\rlx{\hss{$60_s$}}\cg%
\e{0}%
\e{0}%
\e{1}%
\e{0}%
\e{1}%
\e{2}%
\e{0}%
\e{1}%
\e{1}%
\e{0}%
\e{0}%
\eol}\vss}\rg%
%
%
\rx{\vss\hfull{%
\rlx{\hss{$10_s$}}\cg%
\e{0}%
\e{0}%
\e{0}%
\e{0}%
\e{1}%
\e{0}%
\e{0}%
\e{1}%
\e{0}%
\e{0}%
\e{0}%
\eol}\vss}\rg%
\tableclose%
%
%
%
%
%
%
\eop
\eject
\tableopen{Induce/restrict matrix for $W({A_{4}}{A_{1}})\,\subset\,W(E_{6})$}%
%
%
%
%
%
%
\rowpts=18 true pt%
\colpts=18 true pt%
\rowlabpts=40 true pt%
\collabpts=70 true pt%
\clx{\vss\hfull{%
\rlx{\hss{$ $}}\cg%
\cx{\hskip 16 true pt\flip{$[{5}]{\times}[{2}]$}\hss}\cg%
\cx{\hskip 16 true pt\flip{$[{4}{1}]{\times}[{2}]$}\hss}\cg%
\cx{\hskip 16 true pt\flip{$[{3}{2}]{\times}[{2}]$}\hss}\cg%
\cx{\hskip 16 true pt\flip{$[{3}{1^{2}}]{\times}[{2}]$}\hss}\cg%
\cx{\hskip 16 true pt\flip{$[{2^{2}}{1}]{\times}[{2}]$}\hss}\cg%
\cx{\hskip 16 true pt\flip{$[{2}{1^{3}}]{\times}[{2}]$}\hss}\cg%
\cx{\hskip 16 true pt\flip{$[{1^{5}}]{\times}[{2}]$}\hss}\cg%
\cx{\hskip 16 true pt\flip{$[{5}]{\times}[{1^{2}}]$}\hss}\cg%
\cx{\hskip 16 true pt\flip{$[{4}{1}]{\times}[{1^{2}}]$}\hss}\cg%
\cx{\hskip 16 true pt\flip{$[{3}{2}]{\times}[{1^{2}}]$}\hss}\cg%
\cx{\hskip 16 true pt\flip{$[{3}{1^{2}}]{\times}[{1^{2}}]$}\hss}\cg%
\cx{\hskip 16 true pt\flip{$[{2^{2}}{1}]{\times}[{1^{2}}]$}\hss}\cg%
\cx{\hskip 16 true pt\flip{$[{2}{1^{3}}]{\times}[{1^{2}}]$}\hss}\cg%
\cx{\hskip 16 true pt\flip{$[{1^{5}}]{\times}[{1^{2}}]$}\hss}\cg%
\eol}}\rg%
%
%
\rx{\vss\hfull{%
\rlx{\hss{$1_p$}}\cg%
\e{1}%
\e{0}%
\e{0}%
\e{0}%
\e{0}%
\e{0}%
\e{0}%
\e{0}%
\e{0}%
\e{0}%
\e{0}%
\e{0}%
\e{0}%
\e{0}%
\eol}\vss}\rg%
%
%
\rx{\vss\hfull{%
\rlx{\hss{$6_p$}}\cg%
\e{1}%
\e{1}%
\e{0}%
\e{0}%
\e{0}%
\e{0}%
\e{0}%
\e{1}%
\e{0}%
\e{0}%
\e{0}%
\e{0}%
\e{0}%
\e{0}%
\eol}\vss}\rg%
%
%
\rx{\vss\hfull{%
\rlx{\hss{$15_p$}}\cg%
\e{0}%
\e{1}%
\e{0}%
\e{1}%
\e{0}%
\e{0}%
\e{0}%
\e{1}%
\e{1}%
\e{0}%
\e{0}%
\e{0}%
\e{0}%
\e{0}%
\eol}\vss}\rg%
%
%
\rx{\vss\hfull{%
\rlx{\hss{$20_p$}}\cg%
\e{2}%
\e{2}%
\e{1}%
\e{0}%
\e{0}%
\e{0}%
\e{0}%
\e{1}%
\e{1}%
\e{0}%
\e{0}%
\e{0}%
\e{0}%
\e{0}%
\eol}\vss}\rg%
%
%
\rx{\vss\hfull{%
\rlx{\hss{$30_p$}}\cg%
\e{1}%
\e{2}%
\e{1}%
\e{1}%
\e{0}%
\e{0}%
\e{0}%
\e{1}%
\e{1}%
\e{1}%
\e{0}%
\e{0}%
\e{0}%
\e{0}%
\eol}\vss}\rg%
%
%
\rx{\vss\hfull{%
\rlx{\hss{$64_p$}}\cg%
\e{1}%
\e{3}%
\e{2}%
\e{2}%
\e{1}%
\e{0}%
\e{0}%
\e{1}%
\e{3}%
\e{1}%
\e{1}%
\e{0}%
\e{0}%
\e{0}%
\eol}\vss}\rg%
%
%
\rx{\vss\hfull{%
\rlx{\hss{$81_p$}}\cg%
\e{0}%
\e{2}%
\e{2}%
\e{3}%
\e{1}%
\e{1}%
\e{0}%
\e{1}%
\e{2}%
\e{2}%
\e{2}%
\e{1}%
\e{0}%
\e{0}%
\eol}\vss}\rg%
%
%
\rx{\vss\hfull{%
\rlx{\hss{$15_q$}}\cg%
\e{1}%
\e{1}%
\e{1}%
\e{0}%
\e{0}%
\e{0}%
\e{0}%
\e{0}%
\e{0}%
\e{1}%
\e{0}%
\e{0}%
\e{0}%
\e{0}%
\eol}\vss}\rg%
%
%
\rx{\vss\hfull{%
\rlx{\hss{$24_p$}}\cg%
\e{0}%
\e{1}%
\e{1}%
\e{0}%
\e{1}%
\e{0}%
\e{0}%
\e{0}%
\e{1}%
\e{0}%
\e{1}%
\e{0}%
\e{0}%
\e{0}%
\eol}\vss}\rg%
%
%
\rx{\vss\hfull{%
\rlx{\hss{$60_p$}}\cg%
\e{1}%
\e{2}%
\e{3}%
\e{1}%
\e{1}%
\e{0}%
\e{0}%
\e{0}%
\e{1}%
\e{2}%
\e{1}%
\e{1}%
\e{0}%
\e{0}%
\eol}\vss}\rg%
%
%
\rx{\vss\hfull{%
\rlx{\hss{$20_s$}}\cg%
\e{0}%
\e{0}%
\e{0}%
\e{1}%
\e{0}%
\e{1}%
\e{0}%
\e{0}%
\e{1}%
\e{0}%
\e{1}%
\e{0}%
\e{0}%
\e{0}%
\eol}\vss}\rg%
%
%
\rx{\vss\hfull{%
\rlx{\hss{$90_s$}}\cg%
\e{0}%
\e{1}%
\e{1}%
\e{3}%
\e{2}%
\e{2}%
\e{0}%
\e{0}%
\e{2}%
\e{2}%
\e{3}%
\e{1}%
\e{1}%
\e{0}%
\eol}\vss}\rg%
%
%
\rx{\vss\hfull{%
\rlx{\hss{$80_s$}}\cg%
\e{0}%
\e{1}%
\e{2}%
\e{2}%
\e{2}%
\e{1}%
\e{0}%
\e{0}%
\e{1}%
\e{2}%
\e{2}%
\e{2}%
\e{1}%
\e{0}%
\eol}\vss}\rg%
%
%
\rx{\vss\hfull{%
\rlx{\hss{$60_s$}}\cg%
\e{0}%
\e{1}%
\e{2}%
\e{1}%
\e{2}%
\e{0}%
\e{0}%
\e{0}%
\e{0}%
\e{2}%
\e{1}%
\e{2}%
\e{1}%
\e{0}%
\eol}\vss}\rg%
%
%
\rx{\vss\hfull{%
\rlx{\hss{$10_s$}}\cg%
\e{0}%
\e{0}%
\e{1}%
\e{0}%
\e{0}%
\e{0}%
\e{0}%
\e{0}%
\e{0}%
\e{0}%
\e{0}%
\e{1}%
\e{0}%
\e{0}%
\eol}\vss}\rg%
\tableclose%
%
%
%
%
%
%
\tableopen{Induce/restrict matrix for $W({A_{2}}{A_{2}}{A_{1}})\,\subset\,W(E_{6})$}%
%
%
%
%
%
%
\rowpts=18 true pt%
\colpts=18 true pt%
\rowlabpts=40 true pt%
\collabpts=90 true pt%
\clx{\vss\hfull{%
\rlx{\hss{$ $}}\cg%
\cx{\hskip 16 true pt\flip{$[{3}]{\times}[{3}]{\times}[{2}]$}\hss}\cg%
\cx{\hskip 16 true pt\flip{$[{2}{1}]{\times}[{3}]{\times}[{2}]$}\hss}\cg%
\cx{\hskip 16 true pt\flip{$[{1^{3}}]{\times}[{3}]{\times}[{2}]$}\hss}\cg%
\cx{\hskip 16 true pt\flip{$[{3}]{\times}[{2}{1}]{\times}[{2}]$}\hss}\cg%
\cx{\hskip 16 true pt\flip{$[{2}{1}]{\times}[{2}{1}]{\times}[{2}]$}\hss}\cg%
\cx{\hskip 16 true pt\flip{$[{1^{3}}]{\times}[{2}{1}]{\times}[{2}]$}\hss}\cg%
\cx{\hskip 16 true pt\flip{$[{3}]{\times}[{1^{3}}]{\times}[{2}]$}\hss}\cg%
\cx{\hskip 16 true pt\flip{$[{2}{1}]{\times}[{1^{3}}]{\times}[{2}]$}\hss}\cg%
\cx{\hskip 16 true pt\flip{$[{1^{3}}]{\times}[{1^{3}}]{\times}[{2}]$}\hss}\cg%
\cx{\hskip 16 true pt\flip{$[{3}]{\times}[{3}]{\times}[{1^{2}}]$}\hss}\cg%
\cx{\hskip 16 true pt\flip{$[{2}{1}]{\times}[{3}]{\times}[{1^{2}}]$}\hss}\cg%
\cx{\hskip 16 true pt\flip{$[{1^{3}}]{\times}[{3}]{\times}[{1^{2}}]$}\hss}\cg%
\cx{\hskip 16 true pt\flip{$[{3}]{\times}[{2}{1}]{\times}[{1^{2}}]$}\hss}\cg%
\cx{\hskip 16 true pt\flip{$[{2}{1}]{\times}[{2}{1}]{\times}[{1^{2}}]$}\hss}\cg%
\cx{\hskip 16 true pt\flip{$[{1^{3}}]{\times}[{2}{1}]{\times}[{1^{2}}]$}\hss}\cg%
\cx{\hskip 16 true pt\flip{$[{3}]{\times}[{1^{3}}]{\times}[{1^{2}}]$}\hss}\cg%
\cx{\hskip 16 true pt\flip{$[{2}{1}]{\times}[{1^{3}}]{\times}[{1^{2}}]$}\hss}\cg%
\cx{\hskip 16 true pt\flip{$[{1^{3}}]{\times}[{1^{3}}]{\times}[{1^{2}}]$}\hss}\cg%
\eol}}\rg%
%
%
\rx{\vss\hfull{%
\rlx{\hss{$1_p$}}\cg%
\e{1}%
\e{0}%
\e{0}%
\e{0}%
\e{0}%
\e{0}%
\e{0}%
\e{0}%
\e{0}%
\e{0}%
\e{0}%
\e{0}%
\e{0}%
\e{0}%
\e{0}%
\e{0}%
\e{0}%
\e{0}%
\eol}\vss}\rg%
%
%
\rx{\vss\hfull{%
\rlx{\hss{$6_p$}}\cg%
\e{1}%
\e{1}%
\e{0}%
\e{1}%
\e{0}%
\e{0}%
\e{0}%
\e{0}%
\e{0}%
\e{1}%
\e{0}%
\e{0}%
\e{0}%
\e{0}%
\e{0}%
\e{0}%
\e{0}%
\e{0}%
\eol}\vss}\rg%
%
%
\rx{\vss\hfull{%
\rlx{\hss{$15_p$}}\cg%
\e{0}%
\e{1}%
\e{1}%
\e{1}%
\e{1}%
\e{0}%
\e{1}%
\e{0}%
\e{0}%
\e{1}%
\e{1}%
\e{0}%
\e{1}%
\e{0}%
\e{0}%
\e{0}%
\e{0}%
\e{0}%
\eol}\vss}\rg%
%
%
\rx{\vss\hfull{%
\rlx{\hss{$20_p$}}\cg%
\e{3}%
\e{2}%
\e{0}%
\e{2}%
\e{1}%
\e{0}%
\e{0}%
\e{0}%
\e{0}%
\e{1}%
\e{1}%
\e{0}%
\e{1}%
\e{0}%
\e{0}%
\e{0}%
\e{0}%
\e{0}%
\eol}\vss}\rg%
%
%
\rx{\vss\hfull{%
\rlx{\hss{$30_p$}}\cg%
\e{2}%
\e{2}%
\e{1}%
\e{2}%
\e{2}%
\e{0}%
\e{1}%
\e{0}%
\e{0}%
\e{2}%
\e{1}%
\e{0}%
\e{1}%
\e{1}%
\e{0}%
\e{0}%
\e{0}%
\e{0}%
\eol}\vss}\rg%
%
%
\rx{\vss\hfull{%
\rlx{\hss{$64_p$}}\cg%
\e{2}%
\e{4}%
\e{1}%
\e{4}%
\e{4}%
\e{1}%
\e{1}%
\e{1}%
\e{0}%
\e{2}%
\e{3}%
\e{1}%
\e{3}%
\e{2}%
\e{0}%
\e{1}%
\e{0}%
\e{0}%
\eol}\vss}\rg%
%
%
\rx{\vss\hfull{%
\rlx{\hss{$81_p$}}\cg%
\e{1}%
\e{3}%
\e{2}%
\e{3}%
\e{5}%
\e{2}%
\e{2}%
\e{2}%
\e{0}%
\e{2}%
\e{3}%
\e{1}%
\e{3}%
\e{4}%
\e{1}%
\e{1}%
\e{1}%
\e{0}%
\eol}\vss}\rg%
%
%
\rx{\vss\hfull{%
\rlx{\hss{$15_q$}}\cg%
\e{2}%
\e{1}%
\e{0}%
\e{1}%
\e{1}%
\e{0}%
\e{0}%
\e{0}%
\e{0}%
\e{1}%
\e{0}%
\e{0}%
\e{0}%
\e{1}%
\e{0}%
\e{0}%
\e{0}%
\e{0}%
\eol}\vss}\rg%
%
%
\rx{\vss\hfull{%
\rlx{\hss{$24_p$}}\cg%
\e{1}%
\e{1}%
\e{0}%
\e{1}%
\e{2}%
\e{0}%
\e{0}%
\e{0}%
\e{1}%
\e{0}%
\e{1}%
\e{1}%
\e{1}%
\e{1}%
\e{0}%
\e{1}%
\e{0}%
\e{0}%
\eol}\vss}\rg%
%
%
\rx{\vss\hfull{%
\rlx{\hss{$60_p$}}\cg%
\e{3}%
\e{3}%
\e{0}%
\e{3}%
\e{4}%
\e{1}%
\e{0}%
\e{1}%
\e{0}%
\e{1}%
\e{2}%
\e{0}%
\e{2}%
\e{3}%
\e{1}%
\e{0}%
\e{1}%
\e{0}%
\eol}\vss}\rg%
\eop
\eject
\tablecont%
%
%
%
%
%
%
\rowpts=18 true pt%
\colpts=18 true pt%
\rowlabpts=40 true pt%
\collabpts=90 true pt%
\clx{\vss\hfull{%
\rlx{\hss{$ $}}\cg%
\cx{\hskip 16 true pt\flip{$[{3}]{\times}[{3}]{\times}[{2}]$}\hss}\cg%
\cx{\hskip 16 true pt\flip{$[{2}{1}]{\times}[{3}]{\times}[{2}]$}\hss}\cg%
\cx{\hskip 16 true pt\flip{$[{1^{3}}]{\times}[{3}]{\times}[{2}]$}\hss}\cg%
\cx{\hskip 16 true pt\flip{$[{3}]{\times}[{2}{1}]{\times}[{2}]$}\hss}\cg%
\cx{\hskip 16 true pt\flip{$[{2}{1}]{\times}[{2}{1}]{\times}[{2}]$}\hss}\cg%
\cx{\hskip 16 true pt\flip{$[{1^{3}}]{\times}[{2}{1}]{\times}[{2}]$}\hss}\cg%
\cx{\hskip 16 true pt\flip{$[{3}]{\times}[{1^{3}}]{\times}[{2}]$}\hss}\cg%
\cx{\hskip 16 true pt\flip{$[{2}{1}]{\times}[{1^{3}}]{\times}[{2}]$}\hss}\cg%
\cx{\hskip 16 true pt\flip{$[{1^{3}}]{\times}[{1^{3}}]{\times}[{2}]$}\hss}\cg%
\cx{\hskip 16 true pt\flip{$[{3}]{\times}[{3}]{\times}[{1^{2}}]$}\hss}\cg%
\cx{\hskip 16 true pt\flip{$[{2}{1}]{\times}[{3}]{\times}[{1^{2}}]$}\hss}\cg%
\cx{\hskip 16 true pt\flip{$[{1^{3}}]{\times}[{3}]{\times}[{1^{2}}]$}\hss}\cg%
\cx{\hskip 16 true pt\flip{$[{3}]{\times}[{2}{1}]{\times}[{1^{2}}]$}\hss}\cg%
\cx{\hskip 16 true pt\flip{$[{2}{1}]{\times}[{2}{1}]{\times}[{1^{2}}]$}\hss}\cg%
\cx{\hskip 16 true pt\flip{$[{1^{3}}]{\times}[{2}{1}]{\times}[{1^{2}}]$}\hss}\cg%
\cx{\hskip 16 true pt\flip{$[{3}]{\times}[{1^{3}}]{\times}[{1^{2}}]$}\hss}\cg%
\cx{\hskip 16 true pt\flip{$[{2}{1}]{\times}[{1^{3}}]{\times}[{1^{2}}]$}\hss}\cg%
\cx{\hskip 16 true pt\flip{$[{1^{3}}]{\times}[{1^{3}}]{\times}[{1^{2}}]$}\hss}\cg%
\eol}}\rg%
%
%
\rx{\vss\hfull{%
\rlx{\hss{$20_s$}}\cg%
\e{0}%
\e{0}%
\e{1}%
\e{0}%
\e{1}%
\e{1}%
\e{1}%
\e{1}%
\e{0}%
\e{0}%
\e{1}%
\e{1}%
\e{1}%
\e{1}%
\e{0}%
\e{1}%
\e{0}%
\e{0}%
\eol}\vss}\rg%
%
%
\rx{\vss\hfull{%
\rlx{\hss{$90_s$}}\cg%
\e{0}%
\e{2}%
\e{2}%
\e{2}%
\e{5}%
\e{3}%
\e{2}%
\e{3}%
\e{1}%
\e{1}%
\e{3}%
\e{2}%
\e{3}%
\e{5}%
\e{2}%
\e{2}%
\e{2}%
\e{0}%
\eol}\vss}\rg%
%
%
\rx{\vss\hfull{%
\rlx{\hss{$80_s$}}\cg%
\e{1}%
\e{2}%
\e{1}%
\e{2}%
\e{5}%
\e{2}%
\e{1}%
\e{2}%
\e{1}%
\e{1}%
\e{2}%
\e{1}%
\e{2}%
\e{5}%
\e{2}%
\e{1}%
\e{2}%
\e{1}%
\eol}\vss}\rg%
%
%
\rx{\vss\hfull{%
\rlx{\hss{$60_s$}}\cg%
\e{1}%
\e{2}%
\e{0}%
\e{2}%
\e{4}%
\e{1}%
\e{0}%
\e{1}%
\e{1}%
\e{1}%
\e{1}%
\e{0}%
\e{1}%
\e{4}%
\e{2}%
\e{0}%
\e{2}%
\e{1}%
\eol}\vss}\rg%
%
%
\rx{\vss\hfull{%
\rlx{\hss{$10_s$}}\cg%
\e{1}%
\e{0}%
\e{0}%
\e{0}%
\e{1}%
\e{0}%
\e{0}%
\e{0}%
\e{0}%
\e{0}%
\e{0}%
\e{0}%
\e{0}%
\e{1}%
\e{0}%
\e{0}%
\e{0}%
\e{1}%
\eol}\vss}\rg%
\tableclose%
%
%
%
%
%
%
\tableopen{Induce/restrict matrix for $W({D_{6}}{A_{1}})\,\subset\,W(E_{7})$}%
%
%
%
%
%
%
\rowpts=18 true pt%
\colpts=18 true pt%
\rowlabpts=40 true pt%
\collabpts=90 true pt%
\clx{\vss\hfull{%
\rlx{\hss{$ $}}\cg%
\cx{\hskip 16 true pt\flip{$[{6}:-]{\times}[{2}]$}\hss}\cg%
\cx{\hskip 16 true pt\flip{$[{5}{1}:-]{\times}[{2}]$}\hss}\cg%
\cx{\hskip 16 true pt\flip{$[{4}{2}:-]{\times}[{2}]$}\hss}\cg%
\cx{\hskip 16 true pt\flip{$[{4}{1^{2}}:-]{\times}[{2}]$}\hss}\cg%
\cx{\hskip 16 true pt\flip{$[{3^{2}}:-]{\times}[{2}]$}\hss}\cg%
\cx{\hskip 16 true pt\flip{$[{3}{2}{1}:-]{\times}[{2}]$}\hss}\cg%
\cx{\hskip 16 true pt\flip{$[{3}{1^{3}}:-]{\times}[{2}]$}\hss}\cg%
\cx{\hskip 16 true pt\flip{$[{2^{3}}:-]{\times}[{2}]$}\hss}\cg%
\cx{\hskip 16 true pt\flip{$[{2^{2}}{1^{2}}:-]{\times}[{2}]$}\hss}\cg%
\cx{\hskip 16 true pt\flip{$[{2}{1^{4}}:-]{\times}[{2}]$}\hss}\cg%
\cx{\hskip 16 true pt\flip{$[{1^{6}}:-]{\times}[{2}]$}\hss}\cg%
\cx{\hskip 16 true pt\flip{$[{5}:{1}]{\times}[{2}]$}\hss}\cg%
\cx{\hskip 16 true pt\flip{$[{4}{1}:{1}]{\times}[{2}]$}\hss}\cg%
\cx{\hskip 16 true pt\flip{$[{3}{2}:{1}]{\times}[{2}]$}\hss}\cg%
\cx{\hskip 16 true pt\flip{$[{3}{1^{2}}:{1}]{\times}[{2}]$}\hss}\cg%
\cx{\hskip 16 true pt\flip{$[{2^{2}}{1}:{1}]{\times}[{2}]$}\hss}\cg%
\cx{\hskip 16 true pt\flip{$[{2}{1^{3}}:{1}]{\times}[{2}]$}\hss}\cg%
\cx{\hskip 16 true pt\flip{$[{1^{5}}:{1}]{\times}[{2}]$}\hss}\cg%
\cx{\hskip 16 true pt\flip{$[{4}:{2}]{\times}[{2}]$}\hss}\cg%
\cx{\hskip 16 true pt\flip{$[{4}:{1^{2}}]{\times}[{2}]$}\hss}\cg%
\cx{\hskip 16 true pt\flip{$[{3}{1}:{2}]{\times}[{2}]$}\hss}\cg%
\cx{\hskip 16 true pt\flip{$[{3}{1}:{1^{2}}]{\times}[{2}]$}\hss}\cg%
\cx{\hskip 16 true pt\flip{$[{2^{2}}:{2}]{\times}[{2}]$}\hss}\cg%
\cx{\hskip 16 true pt\flip{$[{2^{2}}:{1^{2}}]{\times}[{2}]$}\hss}\cg%
\cx{\hskip 16 true pt\flip{$[{2}{1^{2}}:{2}]{\times}[{2}]$}\hss}\cg%
\eol}}\rg%
%
%
\rx{\vss\hfull{%
\rlx{\hss{$1_a$}}\cg%
\e{1}%
\e{0}%
\e{0}%
\e{0}%
\e{0}%
\e{0}%
\e{0}%
\e{0}%
\e{0}%
\e{0}%
\e{0}%
\e{0}%
\e{0}%
\e{0}%
\e{0}%
\e{0}%
\e{0}%
\e{0}%
\e{0}%
\e{0}%
\e{0}%
\e{0}%
\e{0}%
\e{0}%
\e{0}%
\eol}\vss}\rg%
%
%
\rx{\vss\hfull{%
\rlx{\hss{$7_a$}}\cg%
\e{0}%
\e{0}%
\e{0}%
\e{0}%
\e{0}%
\e{0}%
\e{0}%
\e{0}%
\e{0}%
\e{0}%
\e{1}%
\e{0}%
\e{0}%
\e{0}%
\e{0}%
\e{0}%
\e{0}%
\e{0}%
\e{0}%
\e{0}%
\e{0}%
\e{0}%
\e{0}%
\e{0}%
\e{0}%
\eol}\vss}\rg%
%
%
\rx{\vss\hfull{%
\rlx{\hss{$27_a$}}\cg%
\e{1}%
\e{1}%
\e{0}%
\e{0}%
\e{0}%
\e{0}%
\e{0}%
\e{0}%
\e{0}%
\e{0}%
\e{0}%
\e{0}%
\e{0}%
\e{0}%
\e{0}%
\e{0}%
\e{0}%
\e{0}%
\e{1}%
\e{0}%
\e{0}%
\e{0}%
\e{0}%
\e{0}%
\e{0}%
\eol}\vss}\rg%
%
%
\rx{\vss\hfull{%
\rlx{\hss{$21_a$}}\cg%
\e{0}%
\e{0}%
\e{0}%
\e{0}%
\e{0}%
\e{0}%
\e{0}%
\e{0}%
\e{0}%
\e{0}%
\e{0}%
\e{0}%
\e{0}%
\e{0}%
\e{0}%
\e{0}%
\e{0}%
\e{0}%
\e{0}%
\e{1}%
\e{0}%
\e{0}%
\e{0}%
\e{0}%
\e{0}%
\eol}\vss}\rg%
%
%
\rx{\vss\hfull{%
\rlx{\hss{$35_a$}}\cg%
\e{0}%
\e{0}%
\e{0}%
\e{0}%
\e{0}%
\e{0}%
\e{0}%
\e{0}%
\e{0}%
\e{0}%
\e{0}%
\e{0}%
\e{0}%
\e{0}%
\e{0}%
\e{0}%
\e{0}%
\e{0}%
\e{0}%
\e{0}%
\e{0}%
\e{0}%
\e{0}%
\e{0}%
\e{0}%
\eol}\vss}\rg%
%
%
\rx{\vss\hfull{%
\rlx{\hss{$105_a$}}\cg%
\e{0}%
\e{0}%
\e{0}%
\e{0}%
\e{0}%
\e{0}%
\e{0}%
\e{0}%
\e{0}%
\e{1}%
\e{0}%
\e{0}%
\e{0}%
\e{0}%
\e{0}%
\e{0}%
\e{0}%
\e{0}%
\e{0}%
\e{0}%
\e{0}%
\e{0}%
\e{0}%
\e{0}%
\e{0}%
\eol}\vss}\rg%
%
%
\rx{\vss\hfull{%
\rlx{\hss{$189_a$}}\cg%
\e{0}%
\e{0}%
\e{0}%
\e{0}%
\e{0}%
\e{0}%
\e{0}%
\e{0}%
\e{0}%
\e{0}%
\e{0}%
\e{0}%
\e{0}%
\e{0}%
\e{0}%
\e{0}%
\e{0}%
\e{0}%
\e{0}%
\e{1}%
\e{0}%
\e{1}%
\e{0}%
\e{0}%
\e{1}%
\eol}\vss}\rg%
%
%
\rx{\vss\hfull{%
\rlx{\hss{$21_b$}}\cg%
\e{0}%
\e{0}%
\e{0}%
\e{0}%
\e{0}%
\e{0}%
\e{0}%
\e{0}%
\e{0}%
\e{1}%
\e{0}%
\e{0}%
\e{0}%
\e{0}%
\e{0}%
\e{0}%
\e{0}%
\e{0}%
\e{0}%
\e{0}%
\e{0}%
\e{0}%
\e{0}%
\e{0}%
\e{0}%
\eol}\vss}\rg%
%
%
\rx{\vss\hfull{%
\rlx{\hss{$35_b$}}\cg%
\e{1}%
\e{0}%
\e{1}%
\e{0}%
\e{0}%
\e{0}%
\e{0}%
\e{0}%
\e{0}%
\e{0}%
\e{0}%
\e{0}%
\e{0}%
\e{0}%
\e{0}%
\e{0}%
\e{0}%
\e{0}%
\e{1}%
\e{0}%
\e{0}%
\e{0}%
\e{0}%
\e{0}%
\e{0}%
\eol}\vss}\rg%
%
%
\rx{\vss\hfull{%
\rlx{\hss{$189_b$}}\cg%
\e{0}%
\e{0}%
\e{0}%
\e{0}%
\e{0}%
\e{0}%
\e{0}%
\e{0}%
\e{1}%
\e{0}%
\e{0}%
\e{0}%
\e{0}%
\e{0}%
\e{0}%
\e{0}%
\e{0}%
\e{0}%
\e{0}%
\e{0}%
\e{0}%
\e{0}%
\e{0}%
\e{0}%
\e{0}%
\eol}\vss}\rg%
%
%
\rx{\vss\hfull{%
\rlx{\hss{$189_c$}}\cg%
\e{0}%
\e{0}%
\e{0}%
\e{0}%
\e{0}%
\e{0}%
\e{1}%
\e{0}%
\e{0}%
\e{1}%
\e{0}%
\e{0}%
\e{0}%
\e{0}%
\e{0}%
\e{0}%
\e{0}%
\e{0}%
\e{0}%
\e{0}%
\e{0}%
\e{0}%
\e{0}%
\e{0}%
\e{0}%
\eol}\vss}\rg%
%
%
\rx{\vss\hfull{%
\rlx{\hss{$15_a$}}\cg%
\e{0}%
\e{0}%
\e{0}%
\e{0}%
\e{0}%
\e{0}%
\e{0}%
\e{1}%
\e{0}%
\e{0}%
\e{0}%
\e{0}%
\e{0}%
\e{0}%
\e{0}%
\e{0}%
\e{0}%
\e{0}%
\e{0}%
\e{0}%
\e{0}%
\e{0}%
\e{0}%
\e{0}%
\e{0}%
\eol}\vss}\rg%
%
%
\rx{\vss\hfull{%
\rlx{\hss{$105_b$}}\cg%
\e{0}%
\e{0}%
\e{0}%
\e{0}%
\e{1}%
\e{0}%
\e{0}%
\e{0}%
\e{0}%
\e{0}%
\e{0}%
\e{0}%
\e{0}%
\e{0}%
\e{0}%
\e{0}%
\e{0}%
\e{0}%
\e{1}%
\e{0}%
\e{1}%
\e{0}%
\e{0}%
\e{0}%
\e{0}%
\eol}\vss}\rg%
%
%
\rx{\vss\hfull{%
\rlx{\hss{$105_c$}}\cg%
\e{0}%
\e{0}%
\e{0}%
\e{0}%
\e{0}%
\e{0}%
\e{1}%
\e{0}%
\e{0}%
\e{0}%
\e{0}%
\e{0}%
\e{0}%
\e{0}%
\e{0}%
\e{0}%
\e{0}%
\e{0}%
\e{0}%
\e{0}%
\e{0}%
\e{1}%
\e{0}%
\e{0}%
\e{0}%
\eol}\vss}\rg%
%
%
\rx{\vss\hfull{%
\rlx{\hss{$315_a$}}\cg%
\e{0}%
\e{0}%
\e{0}%
\e{0}%
\e{0}%
\e{0}%
\e{0}%
\e{0}%
\e{0}%
\e{0}%
\e{0}%
\e{0}%
\e{0}%
\e{0}%
\e{0}%
\e{0}%
\e{0}%
\e{0}%
\e{0}%
\e{0}%
\e{0}%
\e{0}%
\e{0}%
\e{1}%
\e{1}%
\eol}\vss}\rg%
%
%
\rx{\vss\hfull{%
\rlx{\hss{$405_a$}}\cg%
\e{0}%
\e{0}%
\e{0}%
\e{0}%
\e{0}%
\e{0}%
\e{0}%
\e{0}%
\e{0}%
\e{0}%
\e{0}%
\e{0}%
\e{0}%
\e{0}%
\e{0}%
\e{0}%
\e{0}%
\e{0}%
\e{0}%
\e{1}%
\e{2}%
\e{1}%
\e{0}%
\e{1}%
\e{1}%
\eol}\vss}\rg%
%
%
\rx{\vss\hfull{%
\rlx{\hss{$168_a$}}\cg%
\e{0}%
\e{1}%
\e{1}%
\e{0}%
\e{0}%
\e{0}%
\e{0}%
\e{0}%
\e{0}%
\e{0}%
\e{0}%
\e{0}%
\e{0}%
\e{0}%
\e{0}%
\e{0}%
\e{0}%
\e{0}%
\e{1}%
\e{0}%
\e{1}%
\e{0}%
\e{1}%
\e{0}%
\e{0}%
\eol}\vss}\rg%
%
%
\rx{\vss\hfull{%
\rlx{\hss{$56_a$}}\cg%
\e{0}%
\e{0}%
\e{0}%
\e{0}%
\e{0}%
\e{0}%
\e{0}%
\e{0}%
\e{0}%
\e{0}%
\e{1}%
\e{0}%
\e{0}%
\e{0}%
\e{0}%
\e{0}%
\e{0}%
\e{0}%
\e{0}%
\e{0}%
\e{0}%
\e{0}%
\e{0}%
\e{0}%
\e{0}%
\eol}\vss}\rg%
%
%
\rx{\vss\hfull{%
\rlx{\hss{$120_a$}}\cg%
\e{0}%
\e{1}%
\e{0}%
\e{0}%
\e{0}%
\e{0}%
\e{0}%
\e{0}%
\e{0}%
\e{0}%
\e{0}%
\e{0}%
\e{0}%
\e{0}%
\e{0}%
\e{0}%
\e{0}%
\e{0}%
\e{1}%
\e{1}%
\e{1}%
\e{0}%
\e{0}%
\e{0}%
\e{0}%
\eol}\vss}\rg%
%
%
\rx{\vss\hfull{%
\rlx{\hss{$210_a$}}\cg%
\e{0}%
\e{0}%
\e{0}%
\e{1}%
\e{0}%
\e{0}%
\e{0}%
\e{0}%
\e{0}%
\e{0}%
\e{0}%
\e{0}%
\e{0}%
\e{0}%
\e{0}%
\e{0}%
\e{0}%
\e{0}%
\e{1}%
\e{1}%
\e{1}%
\e{1}%
\e{0}%
\e{0}%
\e{0}%
\eol}\vss}\rg%
\eop
\eject
\tablecont%
%
%
%
%
%
%
\rowpts=18 true pt%
\colpts=18 true pt%
\rowlabpts=40 true pt%
\collabpts=90 true pt%
\clx{\vss\hfull{%
\rlx{\hss{$ $}}\cg%
\cx{\hskip 16 true pt\flip{$[{6}:-]{\times}[{2}]$}\hss}\cg%
\cx{\hskip 16 true pt\flip{$[{5}{1}:-]{\times}[{2}]$}\hss}\cg%
\cx{\hskip 16 true pt\flip{$[{4}{2}:-]{\times}[{2}]$}\hss}\cg%
\cx{\hskip 16 true pt\flip{$[{4}{1^{2}}:-]{\times}[{2}]$}\hss}\cg%
\cx{\hskip 16 true pt\flip{$[{3^{2}}:-]{\times}[{2}]$}\hss}\cg%
\cx{\hskip 16 true pt\flip{$[{3}{2}{1}:-]{\times}[{2}]$}\hss}\cg%
\cx{\hskip 16 true pt\flip{$[{3}{1^{3}}:-]{\times}[{2}]$}\hss}\cg%
\cx{\hskip 16 true pt\flip{$[{2^{3}}:-]{\times}[{2}]$}\hss}\cg%
\cx{\hskip 16 true pt\flip{$[{2^{2}}{1^{2}}:-]{\times}[{2}]$}\hss}\cg%
\cx{\hskip 16 true pt\flip{$[{2}{1^{4}}:-]{\times}[{2}]$}\hss}\cg%
\cx{\hskip 16 true pt\flip{$[{1^{6}}:-]{\times}[{2}]$}\hss}\cg%
\cx{\hskip 16 true pt\flip{$[{5}:{1}]{\times}[{2}]$}\hss}\cg%
\cx{\hskip 16 true pt\flip{$[{4}{1}:{1}]{\times}[{2}]$}\hss}\cg%
\cx{\hskip 16 true pt\flip{$[{3}{2}:{1}]{\times}[{2}]$}\hss}\cg%
\cx{\hskip 16 true pt\flip{$[{3}{1^{2}}:{1}]{\times}[{2}]$}\hss}\cg%
\cx{\hskip 16 true pt\flip{$[{2^{2}}{1}:{1}]{\times}[{2}]$}\hss}\cg%
\cx{\hskip 16 true pt\flip{$[{2}{1^{3}}:{1}]{\times}[{2}]$}\hss}\cg%
\cx{\hskip 16 true pt\flip{$[{1^{5}}:{1}]{\times}[{2}]$}\hss}\cg%
\cx{\hskip 16 true pt\flip{$[{4}:{2}]{\times}[{2}]$}\hss}\cg%
\cx{\hskip 16 true pt\flip{$[{4}:{1^{2}}]{\times}[{2}]$}\hss}\cg%
\cx{\hskip 16 true pt\flip{$[{3}{1}:{2}]{\times}[{2}]$}\hss}\cg%
\cx{\hskip 16 true pt\flip{$[{3}{1}:{1^{2}}]{\times}[{2}]$}\hss}\cg%
\cx{\hskip 16 true pt\flip{$[{2^{2}}:{2}]{\times}[{2}]$}\hss}\cg%
\cx{\hskip 16 true pt\flip{$[{2^{2}}:{1^{2}}]{\times}[{2}]$}\hss}\cg%
\cx{\hskip 16 true pt\flip{$[{2}{1^{2}}:{2}]{\times}[{2}]$}\hss}\cg%
\eol}}\rg%
%
%
\rx{\vss\hfull{%
\rlx{\hss{$280_a$}}\cg%
\e{0}%
\e{0}%
\e{0}%
\e{0}%
\e{0}%
\e{0}%
\e{0}%
\e{0}%
\e{0}%
\e{0}%
\e{0}%
\e{0}%
\e{0}%
\e{0}%
\e{0}%
\e{0}%
\e{0}%
\e{0}%
\e{0}%
\e{0}%
\e{0}%
\e{0}%
\e{0}%
\e{0}%
\e{1}%
\eol}\vss}\rg%
%
%
\rx{\vss\hfull{%
\rlx{\hss{$336_a$}}\cg%
\e{0}%
\e{0}%
\e{0}%
\e{0}%
\e{0}%
\e{0}%
\e{1}%
\e{0}%
\e{0}%
\e{0}%
\e{0}%
\e{0}%
\e{0}%
\e{0}%
\e{0}%
\e{0}%
\e{0}%
\e{0}%
\e{0}%
\e{0}%
\e{0}%
\e{1}%
\e{0}%
\e{0}%
\e{1}%
\eol}\vss}\rg%
%
%
\rx{\vss\hfull{%
\rlx{\hss{$216_a$}}\cg%
\e{0}%
\e{0}%
\e{0}%
\e{0}%
\e{0}%
\e{1}%
\e{0}%
\e{1}%
\e{0}%
\e{0}%
\e{0}%
\e{0}%
\e{0}%
\e{0}%
\e{0}%
\e{0}%
\e{0}%
\e{0}%
\e{0}%
\e{0}%
\e{0}%
\e{0}%
\e{1}%
\e{0}%
\e{0}%
\eol}\vss}\rg%
%
%
\rx{\vss\hfull{%
\rlx{\hss{$512_a$}}\cg%
\e{0}%
\e{0}%
\e{0}%
\e{0}%
\e{0}%
\e{1}%
\e{0}%
\e{0}%
\e{0}%
\e{0}%
\e{0}%
\e{0}%
\e{0}%
\e{0}%
\e{0}%
\e{0}%
\e{0}%
\e{0}%
\e{0}%
\e{0}%
\e{1}%
\e{1}%
\e{1}%
\e{1}%
\e{1}%
\eol}\vss}\rg%
%
%
\rx{\vss\hfull{%
\rlx{\hss{$378_a$}}\cg%
\e{0}%
\e{0}%
\e{0}%
\e{0}%
\e{0}%
\e{0}%
\e{0}%
\e{0}%
\e{1}%
\e{0}%
\e{0}%
\e{0}%
\e{0}%
\e{0}%
\e{0}%
\e{0}%
\e{0}%
\e{0}%
\e{0}%
\e{0}%
\e{0}%
\e{1}%
\e{0}%
\e{1}%
\e{1}%
\eol}\vss}\rg%
%
%
\rx{\vss\hfull{%
\rlx{\hss{$84_a$}}\cg%
\e{0}%
\e{0}%
\e{1}%
\e{0}%
\e{0}%
\e{0}%
\e{0}%
\e{1}%
\e{0}%
\e{0}%
\e{0}%
\e{0}%
\e{0}%
\e{0}%
\e{0}%
\e{0}%
\e{0}%
\e{0}%
\e{0}%
\e{0}%
\e{0}%
\e{0}%
\e{1}%
\e{0}%
\e{0}%
\eol}\vss}\rg%
%
%
\rx{\vss\hfull{%
\rlx{\hss{$420_a$}}\cg%
\e{0}%
\e{0}%
\e{0}%
\e{1}%
\e{0}%
\e{0}%
\e{0}%
\e{0}%
\e{0}%
\e{0}%
\e{0}%
\e{0}%
\e{0}%
\e{0}%
\e{0}%
\e{0}%
\e{0}%
\e{0}%
\e{0}%
\e{0}%
\e{1}%
\e{1}%
\e{1}%
\e{0}%
\e{1}%
\eol}\vss}\rg%
%
%
\rx{\vss\hfull{%
\rlx{\hss{$280_b$}}\cg%
\e{0}%
\e{0}%
\e{1}%
\e{0}%
\e{0}%
\e{1}%
\e{0}%
\e{0}%
\e{0}%
\e{0}%
\e{0}%
\e{0}%
\e{0}%
\e{0}%
\e{0}%
\e{0}%
\e{0}%
\e{0}%
\e{1}%
\e{0}%
\e{1}%
\e{1}%
\e{1}%
\e{0}%
\e{0}%
\eol}\vss}\rg%
%
%
\rx{\vss\hfull{%
\rlx{\hss{$210_b$}}\cg%
\e{0}%
\e{0}%
\e{0}%
\e{0}%
\e{1}%
\e{0}%
\e{0}%
\e{0}%
\e{0}%
\e{0}%
\e{0}%
\e{0}%
\e{0}%
\e{0}%
\e{0}%
\e{0}%
\e{0}%
\e{0}%
\e{0}%
\e{0}%
\e{1}%
\e{0}%
\e{1}%
\e{1}%
\e{0}%
\eol}\vss}\rg%
%
%
\rx{\vss\hfull{%
\rlx{\hss{$70_a$}}\cg%
\e{0}%
\e{0}%
\e{0}%
\e{0}%
\e{0}%
\e{0}%
\e{0}%
\e{0}%
\e{0}%
\e{0}%
\e{0}%
\e{0}%
\e{0}%
\e{0}%
\e{0}%
\e{0}%
\e{0}%
\e{0}%
\e{0}%
\e{0}%
\e{0}%
\e{0}%
\e{0}%
\e{1}%
\e{0}%
\eol}\vss}\rg%
%
%
%
%
%
%
\rowpts=18 true pt%
\colpts=18 true pt%
\rowlabpts=40 true pt%
\collabpts=90 true pt%
\clx{\vss\hfull{%
\rlx{\hss{$ $}}\cg%
\cx{\hskip 16 true pt\flip{$[{2}{1^{2}}:{1^{2}}]{\times}[{2}]$}\hss}\cg%
\cx{\hskip 16 true pt\flip{$[{1^{4}}:{2}]{\times}[{2}]$}\hss}\cg%
\cx{\hskip 16 true pt\flip{$[{1^{4}}:{1^{2}}]{\times}[{2}]$}\hss}\cg%
\cx{\hskip 16 true pt\flip{$[{3}:{3}]^{+}{\times}[{2}]$}\hss}\cg%
\cx{\hskip 16 true pt\flip{$[{3}:{3}]^{-}{\times}[{2}]$}\hss}\cg%
\cx{\hskip 16 true pt\flip{$[{3}:{2}{1}]{\times}[{2}]$}\hss}\cg%
\cx{\hskip 16 true pt\flip{$[{3}:{1^{3}}]{\times}[{2}]$}\hss}\cg%
\cx{\hskip 16 true pt\flip{$[{2}{1}:{2}{1}]^{+}{\times}[{2}]$}\hss}\cg%
\cx{\hskip 16 true pt\flip{$[{2}{1}:{2}{1}]^{-}{\times}[{2}]$}\hss}\cg%
\cx{\hskip 16 true pt\flip{$[{2}{1}:{1^{3}}]{\times}[{2}]$}\hss}\cg%
\cx{\hskip 16 true pt\flip{$[{1^{3}}:{1^{3}}]^{+}{\times}[{2}]$}\hss}\cg%
\cx{\hskip 16 true pt\flip{$[{1^{3}}:{1^{3}}]^{-}{\times}[{2}]$}\hss}\cg%
\cx{\hskip 16 true pt\flip{$[{6}:-]{\times}[{1^{2}}]$}\hss}\cg%
\cx{\hskip 16 true pt\flip{$[{5}{1}:-]{\times}[{1^{2}}]$}\hss}\cg%
\cx{\hskip 16 true pt\flip{$[{4}{2}:-]{\times}[{1^{2}}]$}\hss}\cg%
\cx{\hskip 16 true pt\flip{$[{4}{1^{2}}:-]{\times}[{1^{2}}]$}\hss}\cg%
\cx{\hskip 16 true pt\flip{$[{3^{2}}:-]{\times}[{1^{2}}]$}\hss}\cg%
\cx{\hskip 16 true pt\flip{$[{3}{2}{1}:-]{\times}[{1^{2}}]$}\hss}\cg%
\cx{\hskip 16 true pt\flip{$[{3}{1^{3}}:-]{\times}[{1^{2}}]$}\hss}\cg%
\cx{\hskip 16 true pt\flip{$[{2^{3}}:-]{\times}[{1^{2}}]$}\hss}\cg%
\cx{\hskip 16 true pt\flip{$[{2^{2}}{1^{2}}:-]{\times}[{1^{2}}]$}\hss}\cg%
\cx{\hskip 16 true pt\flip{$[{2}{1^{4}}:-]{\times}[{1^{2}}]$}\hss}\cg%
\cx{\hskip 16 true pt\flip{$[{1^{6}}:-]{\times}[{1^{2}}]$}\hss}\cg%
\cx{\hskip 16 true pt\flip{$[{5}:{1}]{\times}[{1^{2}}]$}\hss}\cg%
\cx{\hskip 16 true pt\flip{$[{4}{1}:{1}]{\times}[{1^{2}}]$}\hss}\cg%
\eol}}\rg%
%
%
\rx{\vss\hfull{%
\rlx{\hss{$1_a$}}\cg%
\e{0}%
\e{0}%
\e{0}%
\e{0}%
\e{0}%
\e{0}%
\e{0}%
\e{0}%
\e{0}%
\e{0}%
\e{0}%
\e{0}%
\e{0}%
\e{0}%
\e{0}%
\e{0}%
\e{0}%
\e{0}%
\e{0}%
\e{0}%
\e{0}%
\e{0}%
\e{0}%
\e{0}%
\e{0}%
\eol}\vss}\rg%
%
%
\rx{\vss\hfull{%
\rlx{\hss{$7_a$}}\cg%
\e{0}%
\e{0}%
\e{0}%
\e{0}%
\e{0}%
\e{0}%
\e{0}%
\e{0}%
\e{0}%
\e{0}%
\e{0}%
\e{0}%
\e{0}%
\e{0}%
\e{0}%
\e{0}%
\e{0}%
\e{0}%
\e{0}%
\e{0}%
\e{0}%
\e{0}%
\e{0}%
\e{0}%
\e{0}%
\eol}\vss}\rg%
%
%
\rx{\vss\hfull{%
\rlx{\hss{$27_a$}}\cg%
\e{0}%
\e{0}%
\e{0}%
\e{0}%
\e{0}%
\e{0}%
\e{0}%
\e{0}%
\e{0}%
\e{0}%
\e{0}%
\e{0}%
\e{0}%
\e{0}%
\e{0}%
\e{0}%
\e{0}%
\e{0}%
\e{0}%
\e{0}%
\e{0}%
\e{0}%
\e{0}%
\e{1}%
\e{0}%
\eol}\vss}\rg%
%
%
\rx{\vss\hfull{%
\rlx{\hss{$21_a$}}\cg%
\e{0}%
\e{0}%
\e{0}%
\e{0}%
\e{0}%
\e{0}%
\e{0}%
\e{0}%
\e{0}%
\e{0}%
\e{0}%
\e{0}%
\e{0}%
\e{0}%
\e{0}%
\e{0}%
\e{0}%
\e{0}%
\e{0}%
\e{0}%
\e{0}%
\e{0}%
\e{0}%
\e{1}%
\e{0}%
\eol}\vss}\rg%
%
%
\rx{\vss\hfull{%
\rlx{\hss{$35_a$}}\cg%
\e{0}%
\e{1}%
\e{0}%
\e{0}%
\e{0}%
\e{0}%
\e{0}%
\e{0}%
\e{0}%
\e{0}%
\e{0}%
\e{0}%
\e{0}%
\e{0}%
\e{0}%
\e{0}%
\e{0}%
\e{0}%
\e{0}%
\e{0}%
\e{0}%
\e{0}%
\e{0}%
\e{0}%
\e{0}%
\eol}\vss}\rg%
%
%
\rx{\vss\hfull{%
\rlx{\hss{$105_a$}}\cg%
\e{0}%
\e{1}%
\e{1}%
\e{0}%
\e{0}%
\e{0}%
\e{0}%
\e{0}%
\e{0}%
\e{0}%
\e{0}%
\e{0}%
\e{0}%
\e{0}%
\e{0}%
\e{0}%
\e{0}%
\e{0}%
\e{0}%
\e{0}%
\e{0}%
\e{0}%
\e{0}%
\e{0}%
\e{0}%
\eol}\vss}\rg%
%
%
\rx{\vss\hfull{%
\rlx{\hss{$189_a$}}\cg%
\e{0}%
\e{0}%
\e{0}%
\e{0}%
\e{0}%
\e{0}%
\e{0}%
\e{0}%
\e{0}%
\e{0}%
\e{0}%
\e{0}%
\e{0}%
\e{0}%
\e{0}%
\e{0}%
\e{0}%
\e{0}%
\e{0}%
\e{0}%
\e{0}%
\e{0}%
\e{0}%
\e{0}%
\e{1}%
\eol}\vss}\rg%
%
%
\rx{\vss\hfull{%
\rlx{\hss{$21_b$}}\cg%
\e{0}%
\e{0}%
\e{0}%
\e{0}%
\e{0}%
\e{0}%
\e{0}%
\e{0}%
\e{0}%
\e{0}%
\e{0}%
\e{0}%
\e{0}%
\e{0}%
\e{0}%
\e{0}%
\e{0}%
\e{0}%
\e{0}%
\e{0}%
\e{0}%
\e{0}%
\e{0}%
\e{0}%
\e{0}%
\eol}\vss}\rg%
%
%
\rx{\vss\hfull{%
\rlx{\hss{$35_b$}}\cg%
\e{0}%
\e{0}%
\e{0}%
\e{0}%
\e{0}%
\e{0}%
\e{0}%
\e{0}%
\e{0}%
\e{0}%
\e{0}%
\e{0}%
\e{0}%
\e{0}%
\e{0}%
\e{0}%
\e{0}%
\e{0}%
\e{0}%
\e{0}%
\e{0}%
\e{0}%
\e{0}%
\e{0}%
\e{0}%
\eol}\vss}\rg%
%
%
\rx{\vss\hfull{%
\rlx{\hss{$189_b$}}\cg%
\e{1}%
\e{0}%
\e{1}%
\e{0}%
\e{0}%
\e{0}%
\e{0}%
\e{0}%
\e{0}%
\e{0}%
\e{0}%
\e{0}%
\e{0}%
\e{0}%
\e{0}%
\e{0}%
\e{0}%
\e{0}%
\e{0}%
\e{0}%
\e{0}%
\e{0}%
\e{0}%
\e{0}%
\e{0}%
\eol}\vss}\rg%
%
%
\rx{\vss\hfull{%
\rlx{\hss{$189_c$}}\cg%
\e{1}%
\e{1}%
\e{0}%
\e{0}%
\e{0}%
\e{0}%
\e{0}%
\e{0}%
\e{0}%
\e{0}%
\e{0}%
\e{0}%
\e{0}%
\e{0}%
\e{0}%
\e{0}%
\e{0}%
\e{0}%
\e{0}%
\e{0}%
\e{0}%
\e{0}%
\e{0}%
\e{0}%
\e{0}%
\eol}\vss}\rg%
%
%
\rx{\vss\hfull{%
\rlx{\hss{$15_a$}}\cg%
\e{0}%
\e{0}%
\e{0}%
\e{0}%
\e{0}%
\e{0}%
\e{0}%
\e{0}%
\e{0}%
\e{0}%
\e{0}%
\e{0}%
\e{0}%
\e{0}%
\e{0}%
\e{0}%
\e{0}%
\e{0}%
\e{0}%
\e{0}%
\e{0}%
\e{0}%
\e{0}%
\e{0}%
\e{0}%
\eol}\vss}\rg%
%
%
\rx{\vss\hfull{%
\rlx{\hss{$105_b$}}\cg%
\e{0}%
\e{0}%
\e{0}%
\e{0}%
\e{0}%
\e{0}%
\e{0}%
\e{0}%
\e{0}%
\e{0}%
\e{0}%
\e{0}%
\e{0}%
\e{0}%
\e{0}%
\e{0}%
\e{0}%
\e{0}%
\e{0}%
\e{0}%
\e{0}%
\e{0}%
\e{0}%
\e{0}%
\e{0}%
\eol}\vss}\rg%
%
%
\rx{\vss\hfull{%
\rlx{\hss{$105_c$}}\cg%
\e{0}%
\e{0}%
\e{0}%
\e{0}%
\e{0}%
\e{0}%
\e{0}%
\e{0}%
\e{0}%
\e{0}%
\e{0}%
\e{0}%
\e{0}%
\e{0}%
\e{0}%
\e{0}%
\e{0}%
\e{0}%
\e{0}%
\e{0}%
\e{0}%
\e{0}%
\e{0}%
\e{0}%
\e{0}%
\eol}\vss}\rg%
%
%
\rx{\vss\hfull{%
\rlx{\hss{$315_a$}}\cg%
\e{1}%
\e{0}%
\e{1}%
\e{0}%
\e{0}%
\e{0}%
\e{0}%
\e{0}%
\e{0}%
\e{0}%
\e{0}%
\e{0}%
\e{0}%
\e{0}%
\e{0}%
\e{0}%
\e{0}%
\e{0}%
\e{0}%
\e{0}%
\e{0}%
\e{0}%
\e{0}%
\e{0}%
\e{0}%
\eol}\vss}\rg%
\eop
\eject
\tablecont%
%
%
%
%
%
%
\rowpts=18 true pt%
\colpts=18 true pt%
\rowlabpts=40 true pt%
\collabpts=90 true pt%
\clx{\vss\hfull{%
\rlx{\hss{$ $}}\cg%
\cx{\hskip 16 true pt\flip{$[{2}{1^{2}}:{1^{2}}]{\times}[{2}]$}\hss}\cg%
\cx{\hskip 16 true pt\flip{$[{1^{4}}:{2}]{\times}[{2}]$}\hss}\cg%
\cx{\hskip 16 true pt\flip{$[{1^{4}}:{1^{2}}]{\times}[{2}]$}\hss}\cg%
\cx{\hskip 16 true pt\flip{$[{3}:{3}]^{+}{\times}[{2}]$}\hss}\cg%
\cx{\hskip 16 true pt\flip{$[{3}:{3}]^{-}{\times}[{2}]$}\hss}\cg%
\cx{\hskip 16 true pt\flip{$[{3}:{2}{1}]{\times}[{2}]$}\hss}\cg%
\cx{\hskip 16 true pt\flip{$[{3}:{1^{3}}]{\times}[{2}]$}\hss}\cg%
\cx{\hskip 16 true pt\flip{$[{2}{1}:{2}{1}]^{+}{\times}[{2}]$}\hss}\cg%
\cx{\hskip 16 true pt\flip{$[{2}{1}:{2}{1}]^{-}{\times}[{2}]$}\hss}\cg%
\cx{\hskip 16 true pt\flip{$[{2}{1}:{1^{3}}]{\times}[{2}]$}\hss}\cg%
\cx{\hskip 16 true pt\flip{$[{1^{3}}:{1^{3}}]^{+}{\times}[{2}]$}\hss}\cg%
\cx{\hskip 16 true pt\flip{$[{1^{3}}:{1^{3}}]^{-}{\times}[{2}]$}\hss}\cg%
\cx{\hskip 16 true pt\flip{$[{6}:-]{\times}[{1^{2}}]$}\hss}\cg%
\cx{\hskip 16 true pt\flip{$[{5}{1}:-]{\times}[{1^{2}}]$}\hss}\cg%
\cx{\hskip 16 true pt\flip{$[{4}{2}:-]{\times}[{1^{2}}]$}\hss}\cg%
\cx{\hskip 16 true pt\flip{$[{4}{1^{2}}:-]{\times}[{1^{2}}]$}\hss}\cg%
\cx{\hskip 16 true pt\flip{$[{3^{2}}:-]{\times}[{1^{2}}]$}\hss}\cg%
\cx{\hskip 16 true pt\flip{$[{3}{2}{1}:-]{\times}[{1^{2}}]$}\hss}\cg%
\cx{\hskip 16 true pt\flip{$[{3}{1^{3}}:-]{\times}[{1^{2}}]$}\hss}\cg%
\cx{\hskip 16 true pt\flip{$[{2^{3}}:-]{\times}[{1^{2}}]$}\hss}\cg%
\cx{\hskip 16 true pt\flip{$[{2^{2}}{1^{2}}:-]{\times}[{1^{2}}]$}\hss}\cg%
\cx{\hskip 16 true pt\flip{$[{2}{1^{4}}:-]{\times}[{1^{2}}]$}\hss}\cg%
\cx{\hskip 16 true pt\flip{$[{1^{6}}:-]{\times}[{1^{2}}]$}\hss}\cg%
\cx{\hskip 16 true pt\flip{$[{5}:{1}]{\times}[{1^{2}}]$}\hss}\cg%
\cx{\hskip 16 true pt\flip{$[{4}{1}:{1}]{\times}[{1^{2}}]$}\hss}\cg%
\eol}}\rg%
%
%
\rx{\vss\hfull{%
\rlx{\hss{$405_a$}}\cg%
\e{0}%
\e{0}%
\e{0}%
\e{0}%
\e{0}%
\e{0}%
\e{0}%
\e{0}%
\e{0}%
\e{0}%
\e{0}%
\e{0}%
\e{0}%
\e{0}%
\e{0}%
\e{0}%
\e{0}%
\e{0}%
\e{0}%
\e{0}%
\e{0}%
\e{0}%
\e{0}%
\e{0}%
\e{1}%
\eol}\vss}\rg%
%
%
\rx{\vss\hfull{%
\rlx{\hss{$168_a$}}\cg%
\e{0}%
\e{0}%
\e{0}%
\e{0}%
\e{0}%
\e{0}%
\e{0}%
\e{0}%
\e{0}%
\e{0}%
\e{0}%
\e{0}%
\e{0}%
\e{0}%
\e{0}%
\e{0}%
\e{0}%
\e{0}%
\e{0}%
\e{0}%
\e{0}%
\e{0}%
\e{0}%
\e{0}%
\e{1}%
\eol}\vss}\rg%
%
%
\rx{\vss\hfull{%
\rlx{\hss{$56_a$}}\cg%
\e{0}%
\e{0}%
\e{1}%
\e{0}%
\e{0}%
\e{0}%
\e{0}%
\e{0}%
\e{0}%
\e{0}%
\e{0}%
\e{0}%
\e{0}%
\e{0}%
\e{0}%
\e{0}%
\e{0}%
\e{0}%
\e{0}%
\e{0}%
\e{0}%
\e{0}%
\e{0}%
\e{0}%
\e{0}%
\eol}\vss}\rg%
%
%
\rx{\vss\hfull{%
\rlx{\hss{$120_a$}}\cg%
\e{0}%
\e{0}%
\e{0}%
\e{0}%
\e{0}%
\e{0}%
\e{0}%
\e{0}%
\e{0}%
\e{0}%
\e{0}%
\e{0}%
\e{0}%
\e{0}%
\e{0}%
\e{0}%
\e{0}%
\e{0}%
\e{0}%
\e{0}%
\e{0}%
\e{0}%
\e{0}%
\e{1}%
\e{1}%
\eol}\vss}\rg%
%
%
\rx{\vss\hfull{%
\rlx{\hss{$210_a$}}\cg%
\e{0}%
\e{0}%
\e{0}%
\e{0}%
\e{0}%
\e{0}%
\e{0}%
\e{0}%
\e{0}%
\e{0}%
\e{0}%
\e{0}%
\e{0}%
\e{0}%
\e{0}%
\e{0}%
\e{0}%
\e{0}%
\e{0}%
\e{0}%
\e{0}%
\e{0}%
\e{0}%
\e{1}%
\e{1}%
\eol}\vss}\rg%
%
%
\rx{\vss\hfull{%
\rlx{\hss{$280_a$}}\cg%
\e{1}%
\e{1}%
\e{1}%
\e{0}%
\e{0}%
\e{0}%
\e{0}%
\e{0}%
\e{0}%
\e{0}%
\e{0}%
\e{0}%
\e{0}%
\e{0}%
\e{0}%
\e{0}%
\e{0}%
\e{0}%
\e{0}%
\e{0}%
\e{0}%
\e{0}%
\e{0}%
\e{0}%
\e{0}%
\eol}\vss}\rg%
%
%
\rx{\vss\hfull{%
\rlx{\hss{$336_a$}}\cg%
\e{1}%
\e{1}%
\e{0}%
\e{0}%
\e{0}%
\e{0}%
\e{0}%
\e{0}%
\e{0}%
\e{0}%
\e{0}%
\e{0}%
\e{0}%
\e{0}%
\e{0}%
\e{0}%
\e{0}%
\e{0}%
\e{0}%
\e{0}%
\e{0}%
\e{0}%
\e{0}%
\e{0}%
\e{0}%
\eol}\vss}\rg%
%
%
\rx{\vss\hfull{%
\rlx{\hss{$216_a$}}\cg%
\e{1}%
\e{0}%
\e{0}%
\e{0}%
\e{0}%
\e{0}%
\e{0}%
\e{0}%
\e{0}%
\e{0}%
\e{0}%
\e{0}%
\e{0}%
\e{0}%
\e{0}%
\e{0}%
\e{0}%
\e{0}%
\e{0}%
\e{0}%
\e{0}%
\e{0}%
\e{0}%
\e{0}%
\e{0}%
\eol}\vss}\rg%
%
%
\rx{\vss\hfull{%
\rlx{\hss{$512_a$}}\cg%
\e{1}%
\e{0}%
\e{0}%
\e{0}%
\e{0}%
\e{0}%
\e{0}%
\e{0}%
\e{0}%
\e{0}%
\e{0}%
\e{0}%
\e{0}%
\e{0}%
\e{0}%
\e{0}%
\e{0}%
\e{0}%
\e{0}%
\e{0}%
\e{0}%
\e{0}%
\e{0}%
\e{0}%
\e{0}%
\eol}\vss}\rg%
%
%
\rx{\vss\hfull{%
\rlx{\hss{$378_a$}}\cg%
\e{1}%
\e{0}%
\e{0}%
\e{0}%
\e{0}%
\e{0}%
\e{0}%
\e{0}%
\e{0}%
\e{0}%
\e{0}%
\e{0}%
\e{0}%
\e{0}%
\e{0}%
\e{0}%
\e{0}%
\e{0}%
\e{0}%
\e{0}%
\e{0}%
\e{0}%
\e{0}%
\e{0}%
\e{0}%
\eol}\vss}\rg%
%
%
\rx{\vss\hfull{%
\rlx{\hss{$84_a$}}\cg%
\e{0}%
\e{0}%
\e{0}%
\e{0}%
\e{0}%
\e{0}%
\e{0}%
\e{0}%
\e{0}%
\e{0}%
\e{0}%
\e{0}%
\e{0}%
\e{0}%
\e{0}%
\e{0}%
\e{0}%
\e{0}%
\e{0}%
\e{0}%
\e{0}%
\e{0}%
\e{0}%
\e{0}%
\e{0}%
\eol}\vss}\rg%
%
%
\rx{\vss\hfull{%
\rlx{\hss{$420_a$}}\cg%
\e{1}%
\e{0}%
\e{0}%
\e{0}%
\e{0}%
\e{0}%
\e{0}%
\e{0}%
\e{0}%
\e{0}%
\e{0}%
\e{0}%
\e{0}%
\e{0}%
\e{0}%
\e{0}%
\e{0}%
\e{0}%
\e{0}%
\e{0}%
\e{0}%
\e{0}%
\e{0}%
\e{0}%
\e{1}%
\eol}\vss}\rg%
%
%
\rx{\vss\hfull{%
\rlx{\hss{$280_b$}}\cg%
\e{0}%
\e{0}%
\e{0}%
\e{0}%
\e{0}%
\e{0}%
\e{0}%
\e{0}%
\e{0}%
\e{0}%
\e{0}%
\e{0}%
\e{0}%
\e{0}%
\e{0}%
\e{0}%
\e{0}%
\e{0}%
\e{0}%
\e{0}%
\e{0}%
\e{0}%
\e{0}%
\e{0}%
\e{0}%
\eol}\vss}\rg%
%
%
\rx{\vss\hfull{%
\rlx{\hss{$210_b$}}\cg%
\e{0}%
\e{0}%
\e{0}%
\e{0}%
\e{0}%
\e{0}%
\e{0}%
\e{0}%
\e{0}%
\e{0}%
\e{0}%
\e{0}%
\e{0}%
\e{0}%
\e{0}%
\e{0}%
\e{0}%
\e{0}%
\e{0}%
\e{0}%
\e{0}%
\e{0}%
\e{0}%
\e{0}%
\e{0}%
\eol}\vss}\rg%
%
%
\rx{\vss\hfull{%
\rlx{\hss{$70_a$}}\cg%
\e{0}%
\e{0}%
\e{0}%
\e{0}%
\e{0}%
\e{0}%
\e{0}%
\e{0}%
\e{0}%
\e{0}%
\e{0}%
\e{0}%
\e{0}%
\e{0}%
\e{0}%
\e{0}%
\e{0}%
\e{0}%
\e{0}%
\e{0}%
\e{0}%
\e{0}%
\e{0}%
\e{0}%
\e{0}%
\eol}\vss}\rg%
%
%
%
%
%
%
\rowpts=18 true pt%
\colpts=18 true pt%
\rowlabpts=40 true pt%
\collabpts=90 true pt%
\clx{\vss\hfull{%
\rlx{\hss{$ $}}\cg%
\cx{\hskip 16 true pt\flip{$[{3}{2}:{1}]{\times}[{1^{2}}]$}\hss}\cg%
\cx{\hskip 16 true pt\flip{$[{3}{1^{2}}:{1}]{\times}[{1^{2}}]$}\hss}\cg%
\cx{\hskip 16 true pt\flip{$[{2^{2}}{1}:{1}]{\times}[{1^{2}}]$}\hss}\cg%
\cx{\hskip 16 true pt\flip{$[{2}{1^{3}}:{1}]{\times}[{1^{2}}]$}\hss}\cg%
\cx{\hskip 16 true pt\flip{$[{1^{5}}:{1}]{\times}[{1^{2}}]$}\hss}\cg%
\cx{\hskip 16 true pt\flip{$[{4}:{2}]{\times}[{1^{2}}]$}\hss}\cg%
\cx{\hskip 16 true pt\flip{$[{4}:{1^{2}}]{\times}[{1^{2}}]$}\hss}\cg%
\cx{\hskip 16 true pt\flip{$[{3}{1}:{2}]{\times}[{1^{2}}]$}\hss}\cg%
\cx{\hskip 16 true pt\flip{$[{3}{1}:{1^{2}}]{\times}[{1^{2}}]$}\hss}\cg%
\cx{\hskip 16 true pt\flip{$[{2^{2}}:{2}]{\times}[{1^{2}}]$}\hss}\cg%
\cx{\hskip 16 true pt\flip{$[{2^{2}}:{1^{2}}]{\times}[{1^{2}}]$}\hss}\cg%
\cx{\hskip 16 true pt\flip{$[{2}{1^{2}}:{2}]{\times}[{1^{2}}]$}\hss}\cg%
\cx{\hskip 16 true pt\flip{$[{2}{1^{2}}:{1^{2}}]{\times}[{1^{2}}]$}\hss}\cg%
\cx{\hskip 16 true pt\flip{$[{1^{4}}:{2}]{\times}[{1^{2}}]$}\hss}\cg%
\cx{\hskip 16 true pt\flip{$[{1^{4}}:{1^{2}}]{\times}[{1^{2}}]$}\hss}\cg%
\cx{\hskip 16 true pt\flip{$[{3}:{3}]^{+}{\times}[{1^{2}}]$}\hss}\cg%
\cx{\hskip 16 true pt\flip{$[{3}:{3}]^{-}{\times}[{1^{2}}]$}\hss}\cg%
\cx{\hskip 16 true pt\flip{$[{3}:{2}{1}]{\times}[{1^{2}}]$}\hss}\cg%
\cx{\hskip 16 true pt\flip{$[{3}:{1^{3}}]{\times}[{1^{2}}]$}\hss}\cg%
\cx{\hskip 16 true pt\flip{$[{2}{1}:{2}{1}]^{+}{\times}[{1^{2}}]$}\hss}\cg%
\cx{\hskip 16 true pt\flip{$[{2}{1}:{2}{1}]^{-}{\times}[{1^{2}}]$}\hss}\cg%
\cx{\hskip 16 true pt\flip{$[{2}{1}:{1^{3}}]{\times}[{1^{2}}]$}\hss}\cg%
\cx{\hskip 16 true pt\flip{$[{1^{3}}:{1^{3}}]^{+}{\times}[{1^{2}}]$}\hss}\cg%
\cx{\hskip 16 true pt\flip{$[{1^{3}}:{1^{3}}]^{-}{\times}[{1^{2}}]$}\hss}\cg%
\eol}}\rg%
%
%
\rx{\vss\hfull{%
\rlx{\hss{$1_a$}}\cg%
\e{0}%
\e{0}%
\e{0}%
\e{0}%
\e{0}%
\e{0}%
\e{0}%
\e{0}%
\e{0}%
\e{0}%
\e{0}%
\e{0}%
\e{0}%
\e{0}%
\e{0}%
\e{0}%
\e{0}%
\e{0}%
\e{0}%
\e{0}%
\e{0}%
\e{0}%
\e{0}%
\e{0}%
\eol}\vss}\rg%
%
%
\rx{\vss\hfull{%
\rlx{\hss{$7_a$}}\cg%
\e{0}%
\e{0}%
\e{0}%
\e{0}%
\e{1}%
\e{0}%
\e{0}%
\e{0}%
\e{0}%
\e{0}%
\e{0}%
\e{0}%
\e{0}%
\e{0}%
\e{0}%
\e{0}%
\e{0}%
\e{0}%
\e{0}%
\e{0}%
\e{0}%
\e{0}%
\e{0}%
\e{0}%
\eol}\vss}\rg%
%
%
\rx{\vss\hfull{%
\rlx{\hss{$27_a$}}\cg%
\e{0}%
\e{0}%
\e{0}%
\e{0}%
\e{0}%
\e{0}%
\e{0}%
\e{0}%
\e{0}%
\e{0}%
\e{0}%
\e{0}%
\e{0}%
\e{0}%
\e{0}%
\e{0}%
\e{0}%
\e{0}%
\e{0}%
\e{0}%
\e{0}%
\e{0}%
\e{0}%
\e{0}%
\eol}\vss}\rg%
%
%
\rx{\vss\hfull{%
\rlx{\hss{$21_a$}}\cg%
\e{0}%
\e{0}%
\e{0}%
\e{0}%
\e{0}%
\e{0}%
\e{0}%
\e{0}%
\e{0}%
\e{0}%
\e{0}%
\e{0}%
\e{0}%
\e{0}%
\e{0}%
\e{0}%
\e{0}%
\e{0}%
\e{0}%
\e{0}%
\e{0}%
\e{0}%
\e{0}%
\e{0}%
\eol}\vss}\rg%
%
%
\rx{\vss\hfull{%
\rlx{\hss{$35_a$}}\cg%
\e{0}%
\e{0}%
\e{0}%
\e{0}%
\e{0}%
\e{0}%
\e{0}%
\e{0}%
\e{0}%
\e{0}%
\e{0}%
\e{0}%
\e{0}%
\e{0}%
\e{0}%
\e{0}%
\e{0}%
\e{0}%
\e{1}%
\e{0}%
\e{0}%
\e{0}%
\e{0}%
\e{0}%
\eol}\vss}\rg%
%
%
\rx{\vss\hfull{%
\rlx{\hss{$105_a$}}\cg%
\e{0}%
\e{0}%
\e{0}%
\e{1}%
\e{1}%
\e{0}%
\e{0}%
\e{0}%
\e{0}%
\e{0}%
\e{0}%
\e{0}%
\e{0}%
\e{0}%
\e{0}%
\e{0}%
\e{0}%
\e{0}%
\e{0}%
\e{0}%
\e{0}%
\e{1}%
\e{0}%
\e{0}%
\eol}\vss}\rg%
%
%
\rx{\vss\hfull{%
\rlx{\hss{$189_a$}}\cg%
\e{0}%
\e{0}%
\e{0}%
\e{0}%
\e{0}%
\e{0}%
\e{0}%
\e{0}%
\e{0}%
\e{0}%
\e{0}%
\e{0}%
\e{0}%
\e{0}%
\e{0}%
\e{0}%
\e{0}%
\e{1}%
\e{1}%
\e{0}%
\e{0}%
\e{0}%
\e{0}%
\e{0}%
\eol}\vss}\rg%
%
%
\rx{\vss\hfull{%
\rlx{\hss{$21_b$}}\cg%
\e{0}%
\e{0}%
\e{0}%
\e{0}%
\e{1}%
\e{0}%
\e{0}%
\e{0}%
\e{0}%
\e{0}%
\e{0}%
\e{0}%
\e{0}%
\e{0}%
\e{0}%
\e{0}%
\e{0}%
\e{0}%
\e{0}%
\e{0}%
\e{0}%
\e{0}%
\e{0}%
\e{1}%
\eol}\vss}\rg%
%
%
\rx{\vss\hfull{%
\rlx{\hss{$35_b$}}\cg%
\e{0}%
\e{0}%
\e{0}%
\e{0}%
\e{0}%
\e{0}%
\e{0}%
\e{0}%
\e{0}%
\e{0}%
\e{0}%
\e{0}%
\e{0}%
\e{0}%
\e{0}%
\e{0}%
\e{1}%
\e{0}%
\e{0}%
\e{0}%
\e{0}%
\e{0}%
\e{0}%
\e{0}%
\eol}\vss}\rg%
%
%
\rx{\vss\hfull{%
\rlx{\hss{$189_b$}}\cg%
\e{0}%
\e{0}%
\e{1}%
\e{1}%
\e{1}%
\e{0}%
\e{0}%
\e{0}%
\e{0}%
\e{0}%
\e{0}%
\e{0}%
\e{0}%
\e{0}%
\e{0}%
\e{0}%
\e{0}%
\e{0}%
\e{0}%
\e{0}%
\e{0}%
\e{1}%
\e{1}%
\e{1}%
\eol}\vss}\rg%
\eop
\eject
\tablecont%
%
%
%
%
%
%
\rowpts=18 true pt%
\colpts=18 true pt%
\rowlabpts=40 true pt%
\collabpts=90 true pt%
\clx{\vss\hfull{%
\rlx{\hss{$ $}}\cg%
\cx{\hskip 16 true pt\flip{$[{3}{2}:{1}]{\times}[{1^{2}}]$}\hss}\cg%
\cx{\hskip 16 true pt\flip{$[{3}{1^{2}}:{1}]{\times}[{1^{2}}]$}\hss}\cg%
\cx{\hskip 16 true pt\flip{$[{2^{2}}{1}:{1}]{\times}[{1^{2}}]$}\hss}\cg%
\cx{\hskip 16 true pt\flip{$[{2}{1^{3}}:{1}]{\times}[{1^{2}}]$}\hss}\cg%
\cx{\hskip 16 true pt\flip{$[{1^{5}}:{1}]{\times}[{1^{2}}]$}\hss}\cg%
\cx{\hskip 16 true pt\flip{$[{4}:{2}]{\times}[{1^{2}}]$}\hss}\cg%
\cx{\hskip 16 true pt\flip{$[{4}:{1^{2}}]{\times}[{1^{2}}]$}\hss}\cg%
\cx{\hskip 16 true pt\flip{$[{3}{1}:{2}]{\times}[{1^{2}}]$}\hss}\cg%
\cx{\hskip 16 true pt\flip{$[{3}{1}:{1^{2}}]{\times}[{1^{2}}]$}\hss}\cg%
\cx{\hskip 16 true pt\flip{$[{2^{2}}:{2}]{\times}[{1^{2}}]$}\hss}\cg%
\cx{\hskip 16 true pt\flip{$[{2^{2}}:{1^{2}}]{\times}[{1^{2}}]$}\hss}\cg%
\cx{\hskip 16 true pt\flip{$[{2}{1^{2}}:{2}]{\times}[{1^{2}}]$}\hss}\cg%
\cx{\hskip 16 true pt\flip{$[{2}{1^{2}}:{1^{2}}]{\times}[{1^{2}}]$}\hss}\cg%
\cx{\hskip 16 true pt\flip{$[{1^{4}}:{2}]{\times}[{1^{2}}]$}\hss}\cg%
\cx{\hskip 16 true pt\flip{$[{1^{4}}:{1^{2}}]{\times}[{1^{2}}]$}\hss}\cg%
\cx{\hskip 16 true pt\flip{$[{3}:{3}]^{+}{\times}[{1^{2}}]$}\hss}\cg%
\cx{\hskip 16 true pt\flip{$[{3}:{3}]^{-}{\times}[{1^{2}}]$}\hss}\cg%
\cx{\hskip 16 true pt\flip{$[{3}:{2}{1}]{\times}[{1^{2}}]$}\hss}\cg%
\cx{\hskip 16 true pt\flip{$[{3}:{1^{3}}]{\times}[{1^{2}}]$}\hss}\cg%
\cx{\hskip 16 true pt\flip{$[{2}{1}:{2}{1}]^{+}{\times}[{1^{2}}]$}\hss}\cg%
\cx{\hskip 16 true pt\flip{$[{2}{1}:{2}{1}]^{-}{\times}[{1^{2}}]$}\hss}\cg%
\cx{\hskip 16 true pt\flip{$[{2}{1}:{1^{3}}]{\times}[{1^{2}}]$}\hss}\cg%
\cx{\hskip 16 true pt\flip{$[{1^{3}}:{1^{3}}]^{+}{\times}[{1^{2}}]$}\hss}\cg%
\cx{\hskip 16 true pt\flip{$[{1^{3}}:{1^{3}}]^{-}{\times}[{1^{2}}]$}\hss}\cg%
\eol}}\rg%
%
%
\rx{\vss\hfull{%
\rlx{\hss{$189_c$}}\cg%
\e{0}%
\e{0}%
\e{0}%
\e{1}%
\e{0}%
\e{0}%
\e{0}%
\e{0}%
\e{0}%
\e{0}%
\e{0}%
\e{0}%
\e{0}%
\e{0}%
\e{0}%
\e{0}%
\e{0}%
\e{0}%
\e{0}%
\e{0}%
\e{1}%
\e{1}%
\e{0}%
\e{1}%
\eol}\vss}\rg%
%
%
\rx{\vss\hfull{%
\rlx{\hss{$15_a$}}\cg%
\e{0}%
\e{0}%
\e{0}%
\e{0}%
\e{0}%
\e{0}%
\e{0}%
\e{0}%
\e{0}%
\e{0}%
\e{0}%
\e{0}%
\e{0}%
\e{0}%
\e{0}%
\e{0}%
\e{0}%
\e{0}%
\e{0}%
\e{0}%
\e{0}%
\e{0}%
\e{0}%
\e{1}%
\eol}\vss}\rg%
%
%
\rx{\vss\hfull{%
\rlx{\hss{$105_b$}}\cg%
\e{1}%
\e{0}%
\e{0}%
\e{0}%
\e{0}%
\e{0}%
\e{0}%
\e{0}%
\e{0}%
\e{0}%
\e{0}%
\e{0}%
\e{0}%
\e{0}%
\e{0}%
\e{1}%
\e{0}%
\e{0}%
\e{0}%
\e{0}%
\e{0}%
\e{0}%
\e{0}%
\e{0}%
\eol}\vss}\rg%
%
%
\rx{\vss\hfull{%
\rlx{\hss{$105_c$}}\cg%
\e{0}%
\e{0}%
\e{0}%
\e{0}%
\e{0}%
\e{0}%
\e{0}%
\e{0}%
\e{0}%
\e{0}%
\e{0}%
\e{0}%
\e{0}%
\e{0}%
\e{0}%
\e{0}%
\e{1}%
\e{0}%
\e{0}%
\e{0}%
\e{1}%
\e{0}%
\e{0}%
\e{0}%
\eol}\vss}\rg%
%
%
\rx{\vss\hfull{%
\rlx{\hss{$315_a$}}\cg%
\e{0}%
\e{1}%
\e{1}%
\e{1}%
\e{0}%
\e{0}%
\e{0}%
\e{0}%
\e{0}%
\e{0}%
\e{0}%
\e{0}%
\e{0}%
\e{0}%
\e{0}%
\e{0}%
\e{0}%
\e{0}%
\e{0}%
\e{1}%
\e{0}%
\e{1}%
\e{1}%
\e{0}%
\eol}\vss}\rg%
%
%
\rx{\vss\hfull{%
\rlx{\hss{$405_a$}}\cg%
\e{1}%
\e{1}%
\e{0}%
\e{0}%
\e{0}%
\e{0}%
\e{0}%
\e{0}%
\e{0}%
\e{0}%
\e{0}%
\e{0}%
\e{0}%
\e{0}%
\e{0}%
\e{1}%
\e{0}%
\e{1}%
\e{0}%
\e{1}%
\e{0}%
\e{0}%
\e{0}%
\e{0}%
\eol}\vss}\rg%
%
%
\rx{\vss\hfull{%
\rlx{\hss{$168_a$}}\cg%
\e{0}%
\e{0}%
\e{0}%
\e{0}%
\e{0}%
\e{0}%
\e{0}%
\e{0}%
\e{0}%
\e{0}%
\e{0}%
\e{0}%
\e{0}%
\e{0}%
\e{0}%
\e{0}%
\e{0}%
\e{1}%
\e{0}%
\e{0}%
\e{0}%
\e{0}%
\e{0}%
\e{0}%
\eol}\vss}\rg%
%
%
\rx{\vss\hfull{%
\rlx{\hss{$56_a$}}\cg%
\e{0}%
\e{0}%
\e{0}%
\e{1}%
\e{1}%
\e{0}%
\e{0}%
\e{0}%
\e{0}%
\e{0}%
\e{0}%
\e{0}%
\e{0}%
\e{0}%
\e{0}%
\e{0}%
\e{0}%
\e{0}%
\e{0}%
\e{0}%
\e{0}%
\e{0}%
\e{1}%
\e{0}%
\eol}\vss}\rg%
%
%
\rx{\vss\hfull{%
\rlx{\hss{$120_a$}}\cg%
\e{0}%
\e{0}%
\e{0}%
\e{0}%
\e{0}%
\e{0}%
\e{0}%
\e{0}%
\e{0}%
\e{0}%
\e{0}%
\e{0}%
\e{0}%
\e{0}%
\e{0}%
\e{1}%
\e{0}%
\e{0}%
\e{0}%
\e{0}%
\e{0}%
\e{0}%
\e{0}%
\e{0}%
\eol}\vss}\rg%
%
%
\rx{\vss\hfull{%
\rlx{\hss{$210_a$}}\cg%
\e{0}%
\e{0}%
\e{0}%
\e{0}%
\e{0}%
\e{0}%
\e{0}%
\e{0}%
\e{0}%
\e{0}%
\e{0}%
\e{0}%
\e{0}%
\e{0}%
\e{0}%
\e{0}%
\e{1}%
\e{1}%
\e{0}%
\e{0}%
\e{0}%
\e{0}%
\e{0}%
\e{0}%
\eol}\vss}\rg%
%
%
\rx{\vss\hfull{%
\rlx{\hss{$280_a$}}\cg%
\e{0}%
\e{1}%
\e{0}%
\e{1}%
\e{0}%
\e{0}%
\e{0}%
\e{0}%
\e{0}%
\e{0}%
\e{0}%
\e{0}%
\e{0}%
\e{0}%
\e{0}%
\e{0}%
\e{0}%
\e{0}%
\e{1}%
\e{1}%
\e{0}%
\e{1}%
\e{0}%
\e{0}%
\eol}\vss}\rg%
%
%
\rx{\vss\hfull{%
\rlx{\hss{$336_a$}}\cg%
\e{0}%
\e{1}%
\e{0}%
\e{0}%
\e{0}%
\e{0}%
\e{0}%
\e{0}%
\e{0}%
\e{0}%
\e{0}%
\e{0}%
\e{0}%
\e{0}%
\e{0}%
\e{0}%
\e{0}%
\e{1}%
\e{1}%
\e{0}%
\e{1}%
\e{1}%
\e{0}%
\e{0}%
\eol}\vss}\rg%
%
%
\rx{\vss\hfull{%
\rlx{\hss{$216_a$}}\cg%
\e{0}%
\e{0}%
\e{1}%
\e{0}%
\e{0}%
\e{0}%
\e{0}%
\e{0}%
\e{0}%
\e{0}%
\e{0}%
\e{0}%
\e{0}%
\e{0}%
\e{0}%
\e{0}%
\e{0}%
\e{0}%
\e{0}%
\e{0}%
\e{1}%
\e{1}%
\e{0}%
\e{1}%
\eol}\vss}\rg%
%
%
\rx{\vss\hfull{%
\rlx{\hss{$512_a$}}\cg%
\e{1}%
\e{1}%
\e{1}%
\e{0}%
\e{0}%
\e{0}%
\e{0}%
\e{0}%
\e{0}%
\e{0}%
\e{0}%
\e{0}%
\e{0}%
\e{0}%
\e{0}%
\e{0}%
\e{0}%
\e{1}%
\e{0}%
\e{1}%
\e{1}%
\e{1}%
\e{0}%
\e{0}%
\eol}\vss}\rg%
%
%
\rx{\vss\hfull{%
\rlx{\hss{$378_a$}}\cg%
\e{1}%
\e{0}%
\e{1}%
\e{1}%
\e{0}%
\e{0}%
\e{0}%
\e{0}%
\e{0}%
\e{0}%
\e{0}%
\e{0}%
\e{0}%
\e{0}%
\e{0}%
\e{0}%
\e{0}%
\e{0}%
\e{0}%
\e{1}%
\e{1}%
\e{1}%
\e{0}%
\e{0}%
\eol}\vss}\rg%
%
%
\rx{\vss\hfull{%
\rlx{\hss{$84_a$}}\cg%
\e{0}%
\e{0}%
\e{0}%
\e{0}%
\e{0}%
\e{0}%
\e{0}%
\e{0}%
\e{0}%
\e{0}%
\e{0}%
\e{0}%
\e{0}%
\e{0}%
\e{0}%
\e{0}%
\e{0}%
\e{0}%
\e{0}%
\e{0}%
\e{1}%
\e{0}%
\e{0}%
\e{0}%
\eol}\vss}\rg%
%
%
\rx{\vss\hfull{%
\rlx{\hss{$420_a$}}\cg%
\e{0}%
\e{1}%
\e{0}%
\e{0}%
\e{0}%
\e{0}%
\e{0}%
\e{0}%
\e{0}%
\e{0}%
\e{0}%
\e{0}%
\e{0}%
\e{0}%
\e{0}%
\e{0}%
\e{0}%
\e{1}%
\e{1}%
\e{1}%
\e{1}%
\e{0}%
\e{0}%
\e{0}%
\eol}\vss}\rg%
%
%
\rx{\vss\hfull{%
\rlx{\hss{$280_b$}}\cg%
\e{1}%
\e{0}%
\e{0}%
\e{0}%
\e{0}%
\e{0}%
\e{0}%
\e{0}%
\e{0}%
\e{0}%
\e{0}%
\e{0}%
\e{0}%
\e{0}%
\e{0}%
\e{0}%
\e{1}%
\e{1}%
\e{0}%
\e{0}%
\e{1}%
\e{0}%
\e{0}%
\e{0}%
\eol}\vss}\rg%
%
%
\rx{\vss\hfull{%
\rlx{\hss{$210_b$}}\cg%
\e{1}%
\e{0}%
\e{1}%
\e{0}%
\e{0}%
\e{0}%
\e{0}%
\e{0}%
\e{0}%
\e{0}%
\e{0}%
\e{0}%
\e{0}%
\e{0}%
\e{0}%
\e{0}%
\e{0}%
\e{0}%
\e{0}%
\e{1}%
\e{0}%
\e{0}%
\e{0}%
\e{0}%
\eol}\vss}\rg%
%
%
\rx{\vss\hfull{%
\rlx{\hss{$70_a$}}\cg%
\e{0}%
\e{0}%
\e{1}%
\e{0}%
\e{0}%
\e{0}%
\e{0}%
\e{0}%
\e{0}%
\e{0}%
\e{0}%
\e{0}%
\e{0}%
\e{0}%
\e{0}%
\e{0}%
\e{0}%
\e{0}%
\e{0}%
\e{0}%
\e{0}%
\e{0}%
\e{1}%
\e{0}%
\eol}\vss}\rg%
\tableclose%
%
%
%
%
%
%
\tableopen{Induce/restrict matrix for $W(A_{7})\,\subset\,W(E_{7})$}%
%
%
%
%
%
%
\rowpts=18 true pt%
\colpts=18 true pt%
\rowlabpts=40 true pt%
\collabpts=50 true pt%
\clx{\vss\hfull{%
\rlx{\hss{$ $}}\cg%
\cx{\hskip 16 true pt\flip{$[{8}]$}\hss}\cg%
\cx{\hskip 16 true pt\flip{$[{7}{1}]$}\hss}\cg%
\cx{\hskip 16 true pt\flip{$[{6}{2}]$}\hss}\cg%
\cx{\hskip 16 true pt\flip{$[{6}{1^{2}}]$}\hss}\cg%
\cx{\hskip 16 true pt\flip{$[{5}{3}]$}\hss}\cg%
\cx{\hskip 16 true pt\flip{$[{5}{2}{1}]$}\hss}\cg%
\cx{\hskip 16 true pt\flip{$[{5}{1^{3}}]$}\hss}\cg%
\cx{\hskip 16 true pt\flip{$[{4^{2}}]$}\hss}\cg%
\cx{\hskip 16 true pt\flip{$[{4}{3}{1}]$}\hss}\cg%
\cx{\hskip 16 true pt\flip{$[{4}{2^{2}}]$}\hss}\cg%
\cx{\hskip 16 true pt\flip{$[{4}{2}{1^{2}}]$}\hss}\cg%
\cx{\hskip 16 true pt\flip{$[{4}{1^{4}}]$}\hss}\cg%
\cx{\hskip 16 true pt\flip{$[{3^{2}}{2}]$}\hss}\cg%
\cx{\hskip 16 true pt\flip{$[{3^{2}}{1^{2}}]$}\hss}\cg%
\cx{\hskip 16 true pt\flip{$[{3}{2^{2}}{1}]$}\hss}\cg%
\cx{\hskip 16 true pt\flip{$[{3}{2}{1^{3}}]$}\hss}\cg%
\cx{\hskip 16 true pt\flip{$[{3}{1^{5}}]$}\hss}\cg%
\cx{\hskip 16 true pt\flip{$[{2^{4}}]$}\hss}\cg%
\cx{\hskip 16 true pt\flip{$[{2^{3}}{1^{2}}]$}\hss}\cg%
\cx{\hskip 16 true pt\flip{$[{2^{2}}{1^{4}}]$}\hss}\cg%
\cx{\hskip 16 true pt\flip{$[{2}{1^{6}}]$}\hss}\cg%
\cx{\hskip 16 true pt\flip{$[{1^{8}}]$}\hss}\cg%
\eol}}\rg%
%
%
\rx{\vss\hfull{%
\rlx{\hss{$1_a$}}\cg%
\e{1}%
\e{0}%
\e{0}%
\e{0}%
\e{0}%
\e{0}%
\e{0}%
\e{0}%
\e{0}%
\e{0}%
\e{0}%
\e{0}%
\e{0}%
\e{0}%
\e{0}%
\e{0}%
\e{0}%
\e{0}%
\e{0}%
\e{0}%
\e{0}%
\e{0}%
\eol}\vss}\rg%
%
%
\rx{\vss\hfull{%
\rlx{\hss{$7_a$}}\cg%
\e{0}%
\e{0}%
\e{0}%
\e{0}%
\e{0}%
\e{0}%
\e{0}%
\e{0}%
\e{0}%
\e{0}%
\e{0}%
\e{0}%
\e{0}%
\e{0}%
\e{0}%
\e{0}%
\e{0}%
\e{0}%
\e{0}%
\e{0}%
\e{1}%
\e{0}%
\eol}\vss}\rg%
%
%
\rx{\vss\hfull{%
\rlx{\hss{$27_a$}}\cg%
\e{0}%
\e{1}%
\e{1}%
\e{0}%
\e{0}%
\e{0}%
\e{0}%
\e{0}%
\e{0}%
\e{0}%
\e{0}%
\e{0}%
\e{0}%
\e{0}%
\e{0}%
\e{0}%
\e{0}%
\e{0}%
\e{0}%
\e{0}%
\e{0}%
\e{0}%
\eol}\vss}\rg%
%
%
\rx{\vss\hfull{%
\rlx{\hss{$21_a$}}\cg%
\e{0}%
\e{0}%
\e{0}%
\e{1}%
\e{0}%
\e{0}%
\e{0}%
\e{0}%
\e{0}%
\e{0}%
\e{0}%
\e{0}%
\e{0}%
\e{0}%
\e{0}%
\e{0}%
\e{0}%
\e{0}%
\e{0}%
\e{0}%
\e{0}%
\e{0}%
\eol}\vss}\rg%
%
%
\rx{\vss\hfull{%
\rlx{\hss{$35_a$}}\cg%
\e{0}%
\e{0}%
\e{0}%
\e{0}%
\e{0}%
\e{0}%
\e{0}%
\e{0}%
\e{0}%
\e{0}%
\e{0}%
\e{1}%
\e{0}%
\e{0}%
\e{0}%
\e{0}%
\e{0}%
\e{0}%
\e{0}%
\e{0}%
\e{0}%
\e{0}%
\eol}\vss}\rg%
\eop
\eject
\tablecont%
%
%
%
%
%
%
\rowpts=18 true pt%
\colpts=18 true pt%
\rowlabpts=40 true pt%
\collabpts=50 true pt%
\clx{\vss\hfull{%
\rlx{\hss{$ $}}\cg%
\cx{\hskip 16 true pt\flip{$[{8}]$}\hss}\cg%
\cx{\hskip 16 true pt\flip{$[{7}{1}]$}\hss}\cg%
\cx{\hskip 16 true pt\flip{$[{6}{2}]$}\hss}\cg%
\cx{\hskip 16 true pt\flip{$[{6}{1^{2}}]$}\hss}\cg%
\cx{\hskip 16 true pt\flip{$[{5}{3}]$}\hss}\cg%
\cx{\hskip 16 true pt\flip{$[{5}{2}{1}]$}\hss}\cg%
\cx{\hskip 16 true pt\flip{$[{5}{1^{3}}]$}\hss}\cg%
\cx{\hskip 16 true pt\flip{$[{4^{2}}]$}\hss}\cg%
\cx{\hskip 16 true pt\flip{$[{4}{3}{1}]$}\hss}\cg%
\cx{\hskip 16 true pt\flip{$[{4}{2^{2}}]$}\hss}\cg%
\cx{\hskip 16 true pt\flip{$[{4}{2}{1^{2}}]$}\hss}\cg%
\cx{\hskip 16 true pt\flip{$[{4}{1^{4}}]$}\hss}\cg%
\cx{\hskip 16 true pt\flip{$[{3^{2}}{2}]$}\hss}\cg%
\cx{\hskip 16 true pt\flip{$[{3^{2}}{1^{2}}]$}\hss}\cg%
\cx{\hskip 16 true pt\flip{$[{3}{2^{2}}{1}]$}\hss}\cg%
\cx{\hskip 16 true pt\flip{$[{3}{2}{1^{3}}]$}\hss}\cg%
\cx{\hskip 16 true pt\flip{$[{3}{1^{5}}]$}\hss}\cg%
\cx{\hskip 16 true pt\flip{$[{2^{4}}]$}\hss}\cg%
\cx{\hskip 16 true pt\flip{$[{2^{3}}{1^{2}}]$}\hss}\cg%
\cx{\hskip 16 true pt\flip{$[{2^{2}}{1^{4}}]$}\hss}\cg%
\cx{\hskip 16 true pt\flip{$[{2}{1^{6}}]$}\hss}\cg%
\cx{\hskip 16 true pt\flip{$[{1^{8}}]$}\hss}\cg%
\eol}}\rg%
%
%
\rx{\vss\hfull{%
\rlx{\hss{$105_a$}}\cg%
\e{0}%
\e{0}%
\e{0}%
\e{0}%
\e{0}%
\e{0}%
\e{0}%
\e{0}%
\e{0}%
\e{0}%
\e{0}%
\e{0}%
\e{0}%
\e{0}%
\e{0}%
\e{1}%
\e{1}%
\e{0}%
\e{0}%
\e{1}%
\e{0}%
\e{0}%
\eol}\vss}\rg%
%
%
\rx{\vss\hfull{%
\rlx{\hss{$189_a$}}\cg%
\e{0}%
\e{0}%
\e{0}%
\e{0}%
\e{0}%
\e{1}%
\e{1}%
\e{0}%
\e{0}%
\e{0}%
\e{1}%
\e{0}%
\e{0}%
\e{0}%
\e{0}%
\e{0}%
\e{0}%
\e{0}%
\e{0}%
\e{0}%
\e{0}%
\e{0}%
\eol}\vss}\rg%
%
%
\rx{\vss\hfull{%
\rlx{\hss{$21_b$}}\cg%
\e{0}%
\e{0}%
\e{0}%
\e{0}%
\e{0}%
\e{0}%
\e{0}%
\e{0}%
\e{0}%
\e{0}%
\e{0}%
\e{0}%
\e{0}%
\e{0}%
\e{0}%
\e{0}%
\e{0}%
\e{0}%
\e{0}%
\e{1}%
\e{0}%
\e{1}%
\eol}\vss}\rg%
%
%
\rx{\vss\hfull{%
\rlx{\hss{$35_b$}}\cg%
\e{1}%
\e{0}%
\e{1}%
\e{0}%
\e{0}%
\e{0}%
\e{0}%
\e{1}%
\e{0}%
\e{0}%
\e{0}%
\e{0}%
\e{0}%
\e{0}%
\e{0}%
\e{0}%
\e{0}%
\e{0}%
\e{0}%
\e{0}%
\e{0}%
\e{0}%
\eol}\vss}\rg%
%
%
\rx{\vss\hfull{%
\rlx{\hss{$189_b$}}\cg%
\e{0}%
\e{0}%
\e{0}%
\e{0}%
\e{0}%
\e{0}%
\e{0}%
\e{0}%
\e{0}%
\e{0}%
\e{0}%
\e{0}%
\e{0}%
\e{0}%
\e{1}%
\e{1}%
\e{0}%
\e{0}%
\e{1}%
\e{1}%
\e{1}%
\e{0}%
\eol}\vss}\rg%
%
%
\rx{\vss\hfull{%
\rlx{\hss{$189_c$}}\cg%
\e{0}%
\e{0}%
\e{0}%
\e{0}%
\e{0}%
\e{0}%
\e{0}%
\e{0}%
\e{0}%
\e{0}%
\e{0}%
\e{1}%
\e{0}%
\e{1}%
\e{0}%
\e{1}%
\e{0}%
\e{1}%
\e{0}%
\e{1}%
\e{0}%
\e{0}%
\eol}\vss}\rg%
%
%
\rx{\vss\hfull{%
\rlx{\hss{$15_a$}}\cg%
\e{0}%
\e{0}%
\e{0}%
\e{0}%
\e{0}%
\e{0}%
\e{0}%
\e{0}%
\e{0}%
\e{0}%
\e{0}%
\e{0}%
\e{0}%
\e{0}%
\e{0}%
\e{0}%
\e{0}%
\e{1}%
\e{0}%
\e{0}%
\e{0}%
\e{1}%
\eol}\vss}\rg%
%
%
\rx{\vss\hfull{%
\rlx{\hss{$105_b$}}\cg%
\e{0}%
\e{1}%
\e{0}%
\e{0}%
\e{1}%
\e{0}%
\e{0}%
\e{0}%
\e{1}%
\e{0}%
\e{0}%
\e{0}%
\e{0}%
\e{0}%
\e{0}%
\e{0}%
\e{0}%
\e{0}%
\e{0}%
\e{0}%
\e{0}%
\e{0}%
\eol}\vss}\rg%
%
%
\rx{\vss\hfull{%
\rlx{\hss{$105_c$}}\cg%
\e{0}%
\e{0}%
\e{0}%
\e{0}%
\e{0}%
\e{0}%
\e{1}%
\e{1}%
\e{0}%
\e{0}%
\e{0}%
\e{0}%
\e{0}%
\e{1}%
\e{0}%
\e{0}%
\e{0}%
\e{0}%
\e{0}%
\e{0}%
\e{0}%
\e{0}%
\eol}\vss}\rg%
%
%
\rx{\vss\hfull{%
\rlx{\hss{$315_a$}}\cg%
\e{0}%
\e{0}%
\e{0}%
\e{0}%
\e{0}%
\e{0}%
\e{0}%
\e{0}%
\e{0}%
\e{0}%
\e{1}%
\e{0}%
\e{1}%
\e{0}%
\e{1}%
\e{1}%
\e{1}%
\e{0}%
\e{1}%
\e{0}%
\e{0}%
\e{0}%
\eol}\vss}\rg%
%
%
\rx{\vss\hfull{%
\rlx{\hss{$405_a$}}\cg%
\e{0}%
\e{0}%
\e{0}%
\e{1}%
\e{1}%
\e{1}%
\e{0}%
\e{0}%
\e{1}%
\e{0}%
\e{2}%
\e{0}%
\e{1}%
\e{0}%
\e{0}%
\e{0}%
\e{0}%
\e{0}%
\e{0}%
\e{0}%
\e{0}%
\e{0}%
\eol}\vss}\rg%
%
%
\rx{\vss\hfull{%
\rlx{\hss{$168_a$}}\cg%
\e{0}%
\e{0}%
\e{1}%
\e{0}%
\e{1}%
\e{1}%
\e{0}%
\e{0}%
\e{0}%
\e{1}%
\e{0}%
\e{0}%
\e{0}%
\e{0}%
\e{0}%
\e{0}%
\e{0}%
\e{0}%
\e{0}%
\e{0}%
\e{0}%
\e{0}%
\eol}\vss}\rg%
%
%
\rx{\vss\hfull{%
\rlx{\hss{$56_a$}}\cg%
\e{0}%
\e{0}%
\e{0}%
\e{0}%
\e{0}%
\e{0}%
\e{0}%
\e{0}%
\e{0}%
\e{0}%
\e{0}%
\e{0}%
\e{0}%
\e{0}%
\e{0}%
\e{0}%
\e{1}%
\e{0}%
\e{1}%
\e{0}%
\e{1}%
\e{0}%
\eol}\vss}\rg%
%
%
\rx{\vss\hfull{%
\rlx{\hss{$120_a$}}\cg%
\e{0}%
\e{1}%
\e{0}%
\e{1}%
\e{1}%
\e{1}%
\e{0}%
\e{0}%
\e{0}%
\e{0}%
\e{0}%
\e{0}%
\e{0}%
\e{0}%
\e{0}%
\e{0}%
\e{0}%
\e{0}%
\e{0}%
\e{0}%
\e{0}%
\e{0}%
\eol}\vss}\rg%
%
%
\rx{\vss\hfull{%
\rlx{\hss{$210_a$}}\cg%
\e{0}%
\e{0}%
\e{1}%
\e{1}%
\e{0}%
\e{1}%
\e{1}%
\e{0}%
\e{1}%
\e{0}%
\e{0}%
\e{0}%
\e{0}%
\e{0}%
\e{0}%
\e{0}%
\e{0}%
\e{0}%
\e{0}%
\e{0}%
\e{0}%
\e{0}%
\eol}\vss}\rg%
%
%
\rx{\vss\hfull{%
\rlx{\hss{$280_a$}}\cg%
\e{0}%
\e{0}%
\e{0}%
\e{0}%
\e{0}%
\e{0}%
\e{0}%
\e{0}%
\e{0}%
\e{0}%
\e{1}%
\e{1}%
\e{0}%
\e{0}%
\e{1}%
\e{1}%
\e{1}%
\e{0}%
\e{0}%
\e{0}%
\e{0}%
\e{0}%
\eol}\vss}\rg%
%
%
\rx{\vss\hfull{%
\rlx{\hss{$336_a$}}\cg%
\e{0}%
\e{0}%
\e{0}%
\e{0}%
\e{0}%
\e{0}%
\e{1}%
\e{0}%
\e{0}%
\e{1}%
\e{1}%
\e{1}%
\e{0}%
\e{1}%
\e{0}%
\e{1}%
\e{0}%
\e{0}%
\e{0}%
\e{0}%
\e{0}%
\e{0}%
\eol}\vss}\rg%
%
%
\rx{\vss\hfull{%
\rlx{\hss{$216_a$}}\cg%
\e{0}%
\e{0}%
\e{0}%
\e{0}%
\e{0}%
\e{0}%
\e{0}%
\e{0}%
\e{0}%
\e{1}%
\e{0}%
\e{0}%
\e{0}%
\e{1}%
\e{1}%
\e{0}%
\e{0}%
\e{1}%
\e{0}%
\e{1}%
\e{0}%
\e{0}%
\eol}\vss}\rg%
%
%
\rx{\vss\hfull{%
\rlx{\hss{$512_a$}}\cg%
\e{0}%
\e{0}%
\e{0}%
\e{0}%
\e{0}%
\e{1}%
\e{0}%
\e{0}%
\e{1}%
\e{1}%
\e{1}%
\e{0}%
\e{1}%
\e{1}%
\e{1}%
\e{1}%
\e{0}%
\e{0}%
\e{0}%
\e{0}%
\e{0}%
\e{0}%
\eol}\vss}\rg%
%
%
\rx{\vss\hfull{%
\rlx{\hss{$378_a$}}\cg%
\e{0}%
\e{0}%
\e{0}%
\e{0}%
\e{0}%
\e{0}%
\e{0}%
\e{0}%
\e{1}%
\e{0}%
\e{1}%
\e{0}%
\e{0}%
\e{1}%
\e{1}%
\e{1}%
\e{0}%
\e{0}%
\e{1}%
\e{0}%
\e{0}%
\e{0}%
\eol}\vss}\rg%
%
%
\rx{\vss\hfull{%
\rlx{\hss{$84_a$}}\cg%
\e{0}%
\e{0}%
\e{0}%
\e{0}%
\e{0}%
\e{0}%
\e{0}%
\e{1}%
\e{0}%
\e{1}%
\e{0}%
\e{0}%
\e{0}%
\e{0}%
\e{0}%
\e{0}%
\e{0}%
\e{1}%
\e{0}%
\e{0}%
\e{0}%
\e{0}%
\eol}\vss}\rg%
%
%
\rx{\vss\hfull{%
\rlx{\hss{$420_a$}}\cg%
\e{0}%
\e{0}%
\e{0}%
\e{0}%
\e{0}%
\e{1}%
\e{1}%
\e{0}%
\e{1}%
\e{1}%
\e{1}%
\e{1}%
\e{0}%
\e{0}%
\e{1}%
\e{0}%
\e{0}%
\e{0}%
\e{0}%
\e{0}%
\e{0}%
\e{0}%
\eol}\vss}\rg%
%
%
\rx{\vss\hfull{%
\rlx{\hss{$280_b$}}\cg%
\e{0}%
\e{0}%
\e{1}%
\e{0}%
\e{0}%
\e{1}%
\e{0}%
\e{1}%
\e{1}%
\e{1}%
\e{0}%
\e{0}%
\e{0}%
\e{1}%
\e{0}%
\e{0}%
\e{0}%
\e{0}%
\e{0}%
\e{0}%
\e{0}%
\e{0}%
\eol}\vss}\rg%
%
%
\rx{\vss\hfull{%
\rlx{\hss{$210_b$}}\cg%
\e{0}%
\e{0}%
\e{0}%
\e{0}%
\e{1}%
\e{0}%
\e{0}%
\e{0}%
\e{1}%
\e{0}%
\e{0}%
\e{0}%
\e{1}%
\e{0}%
\e{1}%
\e{0}%
\e{0}%
\e{0}%
\e{0}%
\e{0}%
\e{0}%
\e{0}%
\eol}\vss}\rg%
%
%
\rx{\vss\hfull{%
\rlx{\hss{$70_a$}}\cg%
\e{0}%
\e{0}%
\e{0}%
\e{0}%
\e{0}%
\e{0}%
\e{0}%
\e{0}%
\e{0}%
\e{0}%
\e{0}%
\e{0}%
\e{1}%
\e{0}%
\e{0}%
\e{0}%
\e{0}%
\e{0}%
\e{1}%
\e{0}%
\e{0}%
\e{0}%
\eol}\vss}\rg%
\tableclose%
%
%
%
%
%
%
\eop
\eject
\tableopen{Induce/restrict matrix for $W({A_{5}}{A_{2}})\,\subset\,W(E_{7})$}%
%
%
%
%
%
%
\rowpts=18 true pt%
\colpts=18 true pt%
\rowlabpts=40 true pt%
\collabpts=70 true pt%
\clx{\vss\hfull{%
\rlx{\hss{$ $}}\cg%
\cx{\hskip 16 true pt\flip{$[{6}]{\times}[{3}]$}\hss}\cg%
\cx{\hskip 16 true pt\flip{$[{5}{1}]{\times}[{3}]$}\hss}\cg%
\cx{\hskip 16 true pt\flip{$[{4}{2}]{\times}[{3}]$}\hss}\cg%
\cx{\hskip 16 true pt\flip{$[{4}{1^{2}}]{\times}[{3}]$}\hss}\cg%
\cx{\hskip 16 true pt\flip{$[{3^{2}}]{\times}[{3}]$}\hss}\cg%
\cx{\hskip 16 true pt\flip{$[{3}{2}{1}]{\times}[{3}]$}\hss}\cg%
\cx{\hskip 16 true pt\flip{$[{3}{1^{3}}]{\times}[{3}]$}\hss}\cg%
\cx{\hskip 16 true pt\flip{$[{2^{3}}]{\times}[{3}]$}\hss}\cg%
\cx{\hskip 16 true pt\flip{$[{2^{2}}{1^{2}}]{\times}[{3}]$}\hss}\cg%
\cx{\hskip 16 true pt\flip{$[{2}{1^{4}}]{\times}[{3}]$}\hss}\cg%
\cx{\hskip 16 true pt\flip{$[{1^{6}}]{\times}[{3}]$}\hss}\cg%
\cx{\hskip 16 true pt\flip{$[{6}]{\times}[{2}{1}]$}\hss}\cg%
\cx{\hskip 16 true pt\flip{$[{5}{1}]{\times}[{2}{1}]$}\hss}\cg%
\cx{\hskip 16 true pt\flip{$[{4}{2}]{\times}[{2}{1}]$}\hss}\cg%
\cx{\hskip 16 true pt\flip{$[{4}{1^{2}}]{\times}[{2}{1}]$}\hss}\cg%
\cx{\hskip 16 true pt\flip{$[{3^{2}}]{\times}[{2}{1}]$}\hss}\cg%
\cx{\hskip 16 true pt\flip{$[{3}{2}{1}]{\times}[{2}{1}]$}\hss}\cg%
\eol}}\rg%
%
%
\rx{\vss\hfull{%
\rlx{\hss{$1_a$}}\cg%
\e{1}%
\e{0}%
\e{0}%
\e{0}%
\e{0}%
\e{0}%
\e{0}%
\e{0}%
\e{0}%
\e{0}%
\e{0}%
\e{0}%
\e{0}%
\e{0}%
\e{0}%
\e{0}%
\e{0}%
\eol}\vss}\rg%
%
%
\rx{\vss\hfull{%
\rlx{\hss{$7_a$}}\cg%
\e{0}%
\e{0}%
\e{0}%
\e{0}%
\e{0}%
\e{0}%
\e{0}%
\e{0}%
\e{0}%
\e{0}%
\e{0}%
\e{0}%
\e{0}%
\e{0}%
\e{0}%
\e{0}%
\e{0}%
\eol}\vss}\rg%
%
%
\rx{\vss\hfull{%
\rlx{\hss{$27_a$}}\cg%
\e{1}%
\e{1}%
\e{1}%
\e{0}%
\e{0}%
\e{0}%
\e{0}%
\e{0}%
\e{0}%
\e{0}%
\e{0}%
\e{1}%
\e{1}%
\e{0}%
\e{0}%
\e{0}%
\e{0}%
\eol}\vss}\rg%
%
%
\rx{\vss\hfull{%
\rlx{\hss{$21_a$}}\cg%
\e{0}%
\e{0}%
\e{0}%
\e{1}%
\e{0}%
\e{0}%
\e{0}%
\e{0}%
\e{0}%
\e{0}%
\e{0}%
\e{0}%
\e{1}%
\e{0}%
\e{0}%
\e{0}%
\e{0}%
\eol}\vss}\rg%
%
%
\rx{\vss\hfull{%
\rlx{\hss{$35_a$}}\cg%
\e{0}%
\e{0}%
\e{0}%
\e{0}%
\e{0}%
\e{0}%
\e{0}%
\e{0}%
\e{0}%
\e{1}%
\e{0}%
\e{0}%
\e{0}%
\e{0}%
\e{0}%
\e{0}%
\e{0}%
\eol}\vss}\rg%
%
%
\rx{\vss\hfull{%
\rlx{\hss{$105_a$}}\cg%
\e{0}%
\e{0}%
\e{0}%
\e{0}%
\e{0}%
\e{0}%
\e{0}%
\e{0}%
\e{0}%
\e{1}%
\e{0}%
\e{0}%
\e{0}%
\e{0}%
\e{0}%
\e{0}%
\e{0}%
\eol}\vss}\rg%
%
%
\rx{\vss\hfull{%
\rlx{\hss{$189_a$}}\cg%
\e{0}%
\e{0}%
\e{0}%
\e{1}%
\e{0}%
\e{1}%
\e{1}%
\e{0}%
\e{1}%
\e{0}%
\e{0}%
\e{0}%
\e{1}%
\e{1}%
\e{2}%
\e{0}%
\e{1}%
\eol}\vss}\rg%
%
%
\rx{\vss\hfull{%
\rlx{\hss{$21_b$}}\cg%
\e{0}%
\e{0}%
\e{0}%
\e{0}%
\e{0}%
\e{0}%
\e{0}%
\e{0}%
\e{0}%
\e{0}%
\e{0}%
\e{0}%
\e{0}%
\e{0}%
\e{0}%
\e{0}%
\e{0}%
\eol}\vss}\rg%
%
%
\rx{\vss\hfull{%
\rlx{\hss{$35_b$}}\cg%
\e{1}%
\e{1}%
\e{1}%
\e{0}%
\e{0}%
\e{0}%
\e{0}%
\e{0}%
\e{0}%
\e{0}%
\e{0}%
\e{1}%
\e{0}%
\e{1}%
\e{0}%
\e{0}%
\e{0}%
\eol}\vss}\rg%
%
%
\rx{\vss\hfull{%
\rlx{\hss{$189_b$}}\cg%
\e{0}%
\e{0}%
\e{0}%
\e{0}%
\e{0}%
\e{0}%
\e{0}%
\e{0}%
\e{1}%
\e{0}%
\e{0}%
\e{0}%
\e{0}%
\e{0}%
\e{0}%
\e{0}%
\e{1}%
\eol}\vss}\rg%
%
%
\rx{\vss\hfull{%
\rlx{\hss{$189_c$}}\cg%
\e{0}%
\e{0}%
\e{0}%
\e{0}%
\e{0}%
\e{0}%
\e{1}%
\e{0}%
\e{0}%
\e{1}%
\e{0}%
\e{0}%
\e{0}%
\e{0}%
\e{0}%
\e{0}%
\e{1}%
\eol}\vss}\rg%
%
%
\rx{\vss\hfull{%
\rlx{\hss{$15_a$}}\cg%
\e{0}%
\e{0}%
\e{0}%
\e{0}%
\e{0}%
\e{0}%
\e{0}%
\e{0}%
\e{0}%
\e{0}%
\e{0}%
\e{0}%
\e{0}%
\e{0}%
\e{0}%
\e{0}%
\e{0}%
\eol}\vss}\rg%
%
%
\rx{\vss\hfull{%
\rlx{\hss{$105_b$}}\cg%
\e{1}%
\e{1}%
\e{1}%
\e{1}%
\e{1}%
\e{0}%
\e{0}%
\e{0}%
\e{0}%
\e{0}%
\e{0}%
\e{0}%
\e{1}%
\e{1}%
\e{0}%
\e{1}%
\e{1}%
\eol}\vss}\rg%
%
%
\rx{\vss\hfull{%
\rlx{\hss{$105_c$}}\cg%
\e{0}%
\e{0}%
\e{0}%
\e{1}%
\e{0}%
\e{0}%
\e{1}%
\e{0}%
\e{0}%
\e{0}%
\e{0}%
\e{0}%
\e{0}%
\e{1}%
\e{0}%
\e{0}%
\e{1}%
\eol}\vss}\rg%
%
%
\rx{\vss\hfull{%
\rlx{\hss{$315_a$}}\cg%
\e{0}%
\e{0}%
\e{0}%
\e{0}%
\e{0}%
\e{1}%
\e{0}%
\e{0}%
\e{1}%
\e{1}%
\e{0}%
\e{0}%
\e{0}%
\e{0}%
\e{1}%
\e{1}%
\e{2}%
\eol}\vss}\rg%
%
%
\rx{\vss\hfull{%
\rlx{\hss{$405_a$}}\cg%
\e{0}%
\e{1}%
\e{1}%
\e{2}%
\e{1}%
\e{2}%
\e{1}%
\e{0}%
\e{1}%
\e{0}%
\e{0}%
\e{0}%
\e{2}%
\e{2}%
\e{3}%
\e{2}%
\e{3}%
\eol}\vss}\rg%
%
%
\rx{\vss\hfull{%
\rlx{\hss{$168_a$}}\cg%
\e{1}%
\e{1}%
\e{2}%
\e{0}%
\e{1}%
\e{1}%
\e{0}%
\e{1}%
\e{0}%
\e{0}%
\e{0}%
\e{0}%
\e{2}%
\e{2}%
\e{1}%
\e{0}%
\e{1}%
\eol}\vss}\rg%
%
%
\rx{\vss\hfull{%
\rlx{\hss{$56_a$}}\cg%
\e{0}%
\e{0}%
\e{0}%
\e{0}%
\e{0}%
\e{0}%
\e{0}%
\e{0}%
\e{0}%
\e{0}%
\e{1}%
\e{0}%
\e{0}%
\e{0}%
\e{0}%
\e{0}%
\e{0}%
\eol}\vss}\rg%
%
%
\rx{\vss\hfull{%
\rlx{\hss{$120_a$}}\cg%
\e{0}%
\e{2}%
\e{1}%
\e{1}%
\e{0}%
\e{1}%
\e{0}%
\e{0}%
\e{0}%
\e{0}%
\e{0}%
\e{1}%
\e{2}%
\e{1}%
\e{1}%
\e{1}%
\e{0}%
\eol}\vss}\rg%
%
%
\rx{\vss\hfull{%
\rlx{\hss{$210_a$}}\cg%
\e{0}%
\e{1}%
\e{1}%
\e{2}%
\e{1}%
\e{1}%
\e{1}%
\e{0}%
\e{0}%
\e{0}%
\e{0}%
\e{1}%
\e{2}%
\e{2}%
\e{2}%
\e{0}%
\e{1}%
\eol}\vss}\rg%
%
%
\rx{\vss\hfull{%
\rlx{\hss{$280_a$}}\cg%
\e{0}%
\e{0}%
\e{0}%
\e{0}%
\e{0}%
\e{0}%
\e{1}%
\e{1}%
\e{1}%
\e{1}%
\e{1}%
\e{0}%
\e{0}%
\e{0}%
\e{1}%
\e{0}%
\e{2}%
\eol}\vss}\rg%
%
%
\rx{\vss\hfull{%
\rlx{\hss{$336_a$}}\cg%
\e{0}%
\e{0}%
\e{0}%
\e{0}%
\e{0}%
\e{1}%
\e{2}%
\e{0}%
\e{1}%
\e{1}%
\e{0}%
\e{0}%
\e{0}%
\e{1}%
\e{2}%
\e{0}%
\e{2}%
\eol}\vss}\rg%
%
%
\rx{\vss\hfull{%
\rlx{\hss{$216_a$}}\cg%
\e{0}%
\e{0}%
\e{0}%
\e{0}%
\e{0}%
\e{1}%
\e{0}%
\e{1}%
\e{0}%
\e{0}%
\e{0}%
\e{0}%
\e{0}%
\e{1}%
\e{0}%
\e{0}%
\e{2}%
\eol}\vss}\rg%
%
%
\rx{\vss\hfull{%
\rlx{\hss{$512_a$}}\cg%
\e{0}%
\e{0}%
\e{1}%
\e{1}%
\e{1}%
\e{2}%
\e{1}%
\e{1}%
\e{1}%
\e{0}%
\e{0}%
\e{0}%
\e{1}%
\e{2}%
\e{2}%
\e{1}%
\e{5}%
\eol}\vss}\rg%
%
%
\rx{\vss\hfull{%
\rlx{\hss{$378_a$}}\cg%
\e{0}%
\e{0}%
\e{0}%
\e{1}%
\e{0}%
\e{1}%
\e{1}%
\e{0}%
\e{1}%
\e{0}%
\e{0}%
\e{0}%
\e{0}%
\e{1}%
\e{1}%
\e{1}%
\e{3}%
\eol}\vss}\rg%
%
%
\rx{\vss\hfull{%
\rlx{\hss{$84_a$}}\cg%
\e{0}%
\e{0}%
\e{1}%
\e{0}%
\e{0}%
\e{0}%
\e{0}%
\e{1}%
\e{0}%
\e{0}%
\e{0}%
\e{0}%
\e{0}%
\e{1}%
\e{0}%
\e{0}%
\e{1}%
\eol}\vss}\rg%
%
%
\rx{\vss\hfull{%
\rlx{\hss{$420_a$}}\cg%
\e{0}%
\e{0}%
\e{1}%
\e{1}%
\e{0}%
\e{2}%
\e{1}%
\e{1}%
\e{1}%
\e{1}%
\e{0}%
\e{0}%
\e{1}%
\e{2}%
\e{3}%
\e{1}%
\e{3}%
\eol}\vss}\rg%
%
%
\rx{\vss\hfull{%
\rlx{\hss{$280_b$}}\cg%
\e{0}%
\e{1}%
\e{2}%
\e{1}%
\e{0}%
\e{2}%
\e{0}%
\e{0}%
\e{0}%
\e{0}%
\e{0}%
\e{1}%
\e{1}%
\e{3}%
\e{1}%
\e{1}%
\e{2}%
\eol}\vss}\rg%
%
%
\rx{\vss\hfull{%
\rlx{\hss{$210_b$}}\cg%
\e{0}%
\e{1}%
\e{1}%
\e{0}%
\e{1}%
\e{1}%
\e{0}%
\e{0}%
\e{0}%
\e{0}%
\e{0}%
\e{0}%
\e{0}%
\e{1}%
\e{1}%
\e{2}%
\e{2}%
\eol}\vss}\rg%
%
%
\rx{\vss\hfull{%
\rlx{\hss{$70_a$}}\cg%
\e{0}%
\e{0}%
\e{0}%
\e{0}%
\e{1}%
\e{0}%
\e{0}%
\e{0}%
\e{0}%
\e{0}%
\e{0}%
\e{0}%
\e{0}%
\e{0}%
\e{0}%
\e{0}%
\e{1}%
\eol}\vss}\rg%
\eop
\eject
\tablecont%
%
%
%
%
%
%
\rowpts=18 true pt%
\colpts=18 true pt%
\rowlabpts=40 true pt%
\collabpts=70 true pt%
\clx{\vss\hfull{%
\rlx{\hss{$ $}}\cg%
\cx{\hskip 16 true pt\flip{$[{3}{1^{3}}]{\times}[{2}{1}]$}\hss}\cg%
\cx{\hskip 16 true pt\flip{$[{2^{3}}]{\times}[{2}{1}]$}\hss}\cg%
\cx{\hskip 16 true pt\flip{$[{2^{2}}{1^{2}}]{\times}[{2}{1}]$}\hss}\cg%
\cx{\hskip 16 true pt\flip{$[{2}{1^{4}}]{\times}[{2}{1}]$}\hss}\cg%
\cx{\hskip 16 true pt\flip{$[{1^{6}}]{\times}[{2}{1}]$}\hss}\cg%
\cx{\hskip 16 true pt\flip{$[{6}]{\times}[{1^{3}}]$}\hss}\cg%
\cx{\hskip 16 true pt\flip{$[{5}{1}]{\times}[{1^{3}}]$}\hss}\cg%
\cx{\hskip 16 true pt\flip{$[{4}{2}]{\times}[{1^{3}}]$}\hss}\cg%
\cx{\hskip 16 true pt\flip{$[{4}{1^{2}}]{\times}[{1^{3}}]$}\hss}\cg%
\cx{\hskip 16 true pt\flip{$[{3^{2}}]{\times}[{1^{3}}]$}\hss}\cg%
\cx{\hskip 16 true pt\flip{$[{3}{2}{1}]{\times}[{1^{3}}]$}\hss}\cg%
\cx{\hskip 16 true pt\flip{$[{3}{1^{3}}]{\times}[{1^{3}}]$}\hss}\cg%
\cx{\hskip 16 true pt\flip{$[{2^{3}}]{\times}[{1^{3}}]$}\hss}\cg%
\cx{\hskip 16 true pt\flip{$[{2^{2}}{1^{2}}]{\times}[{1^{3}}]$}\hss}\cg%
\cx{\hskip 16 true pt\flip{$[{2}{1^{4}}]{\times}[{1^{3}}]$}\hss}\cg%
\cx{\hskip 16 true pt\flip{$[{1^{6}}]{\times}[{1^{3}}]$}\hss}\cg%
\eol}}\rg%
%
%
\rx{\vss\hfull{%
\rlx{\hss{$1_a$}}\cg%
\e{0}%
\e{0}%
\e{0}%
\e{0}%
\e{0}%
\e{0}%
\e{0}%
\e{0}%
\e{0}%
\e{0}%
\e{0}%
\e{0}%
\e{0}%
\e{0}%
\e{0}%
\e{0}%
\eol}\vss}\rg%
%
%
\rx{\vss\hfull{%
\rlx{\hss{$7_a$}}\cg%
\e{0}%
\e{0}%
\e{0}%
\e{0}%
\e{1}%
\e{0}%
\e{0}%
\e{0}%
\e{0}%
\e{0}%
\e{0}%
\e{0}%
\e{0}%
\e{0}%
\e{1}%
\e{0}%
\eol}\vss}\rg%
%
%
\rx{\vss\hfull{%
\rlx{\hss{$27_a$}}\cg%
\e{0}%
\e{0}%
\e{0}%
\e{0}%
\e{0}%
\e{0}%
\e{0}%
\e{0}%
\e{0}%
\e{0}%
\e{0}%
\e{0}%
\e{0}%
\e{0}%
\e{0}%
\e{0}%
\eol}\vss}\rg%
%
%
\rx{\vss\hfull{%
\rlx{\hss{$21_a$}}\cg%
\e{0}%
\e{0}%
\e{0}%
\e{0}%
\e{0}%
\e{1}%
\e{0}%
\e{0}%
\e{0}%
\e{0}%
\e{0}%
\e{0}%
\e{0}%
\e{0}%
\e{0}%
\e{0}%
\eol}\vss}\rg%
%
%
\rx{\vss\hfull{%
\rlx{\hss{$35_a$}}\cg%
\e{1}%
\e{0}%
\e{0}%
\e{0}%
\e{0}%
\e{0}%
\e{0}%
\e{0}%
\e{1}%
\e{0}%
\e{0}%
\e{0}%
\e{0}%
\e{0}%
\e{0}%
\e{0}%
\eol}\vss}\rg%
%
%
\rx{\vss\hfull{%
\rlx{\hss{$105_a$}}\cg%
\e{1}%
\e{0}%
\e{1}%
\e{2}%
\e{1}%
\e{0}%
\e{0}%
\e{0}%
\e{0}%
\e{0}%
\e{1}%
\e{1}%
\e{0}%
\e{1}%
\e{1}%
\e{0}%
\eol}\vss}\rg%
%
%
\rx{\vss\hfull{%
\rlx{\hss{$189_a$}}\cg%
\e{1}%
\e{0}%
\e{0}%
\e{0}%
\e{0}%
\e{0}%
\e{1}%
\e{1}%
\e{1}%
\e{0}%
\e{0}%
\e{0}%
\e{0}%
\e{0}%
\e{0}%
\e{0}%
\eol}\vss}\rg%
%
%
\rx{\vss\hfull{%
\rlx{\hss{$21_b$}}\cg%
\e{0}%
\e{0}%
\e{0}%
\e{1}%
\e{0}%
\e{0}%
\e{0}%
\e{0}%
\e{0}%
\e{0}%
\e{0}%
\e{0}%
\e{1}%
\e{0}%
\e{1}%
\e{1}%
\eol}\vss}\rg%
%
%
\rx{\vss\hfull{%
\rlx{\hss{$35_b$}}\cg%
\e{0}%
\e{0}%
\e{0}%
\e{0}%
\e{0}%
\e{0}%
\e{0}%
\e{0}%
\e{0}%
\e{0}%
\e{0}%
\e{0}%
\e{0}%
\e{0}%
\e{0}%
\e{0}%
\eol}\vss}\rg%
%
%
\rx{\vss\hfull{%
\rlx{\hss{$189_b$}}\cg%
\e{1}%
\e{1}%
\e{2}%
\e{2}%
\e{1}%
\e{0}%
\e{0}%
\e{0}%
\e{0}%
\e{0}%
\e{1}%
\e{1}%
\e{1}%
\e{2}%
\e{2}%
\e{1}%
\eol}\vss}\rg%
%
%
\rx{\vss\hfull{%
\rlx{\hss{$189_c$}}\cg%
\e{2}%
\e{1}%
\e{1}%
\e{2}%
\e{0}%
\e{0}%
\e{0}%
\e{1}%
\e{0}%
\e{0}%
\e{1}%
\e{1}%
\e{1}%
\e{1}%
\e{1}%
\e{0}%
\eol}\vss}\rg%
%
%
\rx{\vss\hfull{%
\rlx{\hss{$15_a$}}\cg%
\e{0}%
\e{1}%
\e{0}%
\e{0}%
\e{0}%
\e{0}%
\e{0}%
\e{0}%
\e{0}%
\e{0}%
\e{0}%
\e{0}%
\e{0}%
\e{0}%
\e{1}%
\e{0}%
\eol}\vss}\rg%
%
%
\rx{\vss\hfull{%
\rlx{\hss{$105_b$}}\cg%
\e{0}%
\e{0}%
\e{0}%
\e{0}%
\e{0}%
\e{0}%
\e{0}%
\e{0}%
\e{0}%
\e{1}%
\e{0}%
\e{0}%
\e{0}%
\e{0}%
\e{0}%
\e{0}%
\eol}\vss}\rg%
%
%
\rx{\vss\hfull{%
\rlx{\hss{$105_c$}}\cg%
\e{1}%
\e{0}%
\e{0}%
\e{0}%
\e{0}%
\e{1}%
\e{0}%
\e{1}%
\e{0}%
\e{0}%
\e{0}%
\e{0}%
\e{1}%
\e{0}%
\e{0}%
\e{0}%
\eol}\vss}\rg%
%
%
\rx{\vss\hfull{%
\rlx{\hss{$315_a$}}\cg%
\e{2}%
\e{1}%
\e{3}%
\e{1}%
\e{1}%
\e{0}%
\e{0}%
\e{0}%
\e{1}%
\e{0}%
\e{2}%
\e{1}%
\e{0}%
\e{2}%
\e{1}%
\e{0}%
\eol}\vss}\rg%
%
%
\rx{\vss\hfull{%
\rlx{\hss{$405_a$}}\cg%
\e{1}%
\e{0}%
\e{1}%
\e{0}%
\e{0}%
\e{0}%
\e{1}%
\e{1}%
\e{1}%
\e{1}%
\e{1}%
\e{0}%
\e{0}%
\e{0}%
\e{0}%
\e{0}%
\eol}\vss}\rg%
%
%
\rx{\vss\hfull{%
\rlx{\hss{$168_a$}}\cg%
\e{0}%
\e{0}%
\e{0}%
\e{0}%
\e{0}%
\e{0}%
\e{0}%
\e{0}%
\e{1}%
\e{0}%
\e{0}%
\e{0}%
\e{0}%
\e{0}%
\e{0}%
\e{0}%
\eol}\vss}\rg%
%
%
\rx{\vss\hfull{%
\rlx{\hss{$56_a$}}\cg%
\e{0}%
\e{0}%
\e{1}%
\e{1}%
\e{1}%
\e{0}%
\e{0}%
\e{0}%
\e{0}%
\e{0}%
\e{0}%
\e{1}%
\e{0}%
\e{1}%
\e{1}%
\e{1}%
\eol}\vss}\rg%
%
%
\rx{\vss\hfull{%
\rlx{\hss{$120_a$}}\cg%
\e{0}%
\e{0}%
\e{0}%
\e{0}%
\e{0}%
\e{0}%
\e{1}%
\e{0}%
\e{0}%
\e{0}%
\e{0}%
\e{0}%
\e{0}%
\e{0}%
\e{0}%
\e{0}%
\eol}\vss}\rg%
%
%
\rx{\vss\hfull{%
\rlx{\hss{$210_a$}}\cg%
\e{0}%
\e{0}%
\e{0}%
\e{0}%
\e{0}%
\e{1}%
\e{1}%
\e{1}%
\e{0}%
\e{0}%
\e{0}%
\e{0}%
\e{0}%
\e{0}%
\e{0}%
\e{0}%
\eol}\vss}\rg%
%
%
\rx{\vss\hfull{%
\rlx{\hss{$280_a$}}\cg%
\e{2}%
\e{0}%
\e{2}%
\e{2}%
\e{0}%
\e{0}%
\e{0}%
\e{0}%
\e{2}%
\e{1}%
\e{1}%
\e{2}%
\e{0}%
\e{1}%
\e{0}%
\e{0}%
\eol}\vss}\rg%
%
%
\rx{\vss\hfull{%
\rlx{\hss{$336_a$}}\cg%
\e{3}%
\e{1}%
\e{1}%
\e{1}%
\e{0}%
\e{0}%
\e{1}%
\e{1}%
\e{1}%
\e{0}%
\e{2}%
\e{1}%
\e{0}%
\e{0}%
\e{0}%
\e{0}%
\eol}\vss}\rg%
%
%
\rx{\vss\hfull{%
\rlx{\hss{$216_a$}}\cg%
\e{1}%
\e{2}%
\e{1}%
\e{1}%
\e{0}%
\e{0}%
\e{0}%
\e{0}%
\e{0}%
\e{0}%
\e{1}%
\e{1}%
\e{1}%
\e{1}%
\e{1}%
\e{0}%
\eol}\vss}\rg%
%
%
\rx{\vss\hfull{%
\rlx{\hss{$512_a$}}\cg%
\e{2}%
\e{1}%
\e{2}%
\e{1}%
\e{0}%
\e{0}%
\e{0}%
\e{1}%
\e{1}%
\e{1}%
\e{2}%
\e{1}%
\e{1}%
\e{1}%
\e{0}%
\e{0}%
\eol}\vss}\rg%
%
%
\rx{\vss\hfull{%
\rlx{\hss{$378_a$}}\cg%
\e{2}%
\e{1}%
\e{3}%
\e{1}%
\e{0}%
\e{0}%
\e{0}%
\e{1}%
\e{0}%
\e{1}%
\e{2}%
\e{1}%
\e{1}%
\e{1}%
\e{1}%
\e{0}%
\eol}\vss}\rg%
%
%
\rx{\vss\hfull{%
\rlx{\hss{$84_a$}}\cg%
\e{0}%
\e{1}%
\e{0}%
\e{0}%
\e{0}%
\e{0}%
\e{0}%
\e{0}%
\e{0}%
\e{0}%
\e{0}%
\e{1}%
\e{0}%
\e{0}%
\e{0}%
\e{0}%
\eol}\vss}\rg%
%
%
\rx{\vss\hfull{%
\rlx{\hss{$420_a$}}\cg%
\e{2}%
\e{1}%
\e{1}%
\e{0}%
\e{0}%
\e{0}%
\e{1}%
\e{1}%
\e{2}%
\e{0}%
\e{1}%
\e{1}%
\e{0}%
\e{0}%
\e{0}%
\e{0}%
\eol}\vss}\rg%
%
%
\rx{\vss\hfull{%
\rlx{\hss{$280_b$}}\cg%
\e{1}%
\e{1}%
\e{0}%
\e{0}%
\e{0}%
\e{0}%
\e{0}%
\e{1}%
\e{0}%
\e{0}%
\e{1}%
\e{0}%
\e{0}%
\e{0}%
\e{0}%
\e{0}%
\eol}\vss}\rg%
%
%
\rx{\vss\hfull{%
\rlx{\hss{$210_b$}}\cg%
\e{0}%
\e{1}%
\e{1}%
\e{0}%
\e{0}%
\e{0}%
\e{0}%
\e{0}%
\e{0}%
\e{0}%
\e{1}%
\e{0}%
\e{0}%
\e{1}%
\e{0}%
\e{0}%
\eol}\vss}\rg%
%
%
\rx{\vss\hfull{%
\rlx{\hss{$70_a$}}\cg%
\e{0}%
\e{0}%
\e{1}%
\e{0}%
\e{0}%
\e{0}%
\e{0}%
\e{0}%
\e{0}%
\e{0}%
\e{0}%
\e{0}%
\e{1}%
\e{1}%
\e{0}%
\e{1}%
\eol}\vss}\rg%
\tableclose%
%
%
%
%
%
%
\eop
\eject
\tableopen{Induce/restrict matrix for $W({A_{3}}{A_{3}}{A_{1}})\,\subset\,W(E_{7})$}%
%
%
%
%
%
%
\rowpts=18 true pt%
\colpts=18 true pt%
\rowlabpts=40 true pt%
\collabpts=90 true pt%
\clx{\vss\hfull{%
\rlx{\hss{$ $}}\cg%
\cx{\hskip 16 true pt\flip{$[{4}]{\times}[{4}]{\times}[{2}]$}\hss}\cg%
\cx{\hskip 16 true pt\flip{$[{3}{1}]{\times}[{4}]{\times}[{2}]$}\hss}\cg%
\cx{\hskip 16 true pt\flip{$[{2^{2}}]{\times}[{4}]{\times}[{2}]$}\hss}\cg%
\cx{\hskip 16 true pt\flip{$[{2}{1^{2}}]{\times}[{4}]{\times}[{2}]$}\hss}\cg%
\cx{\hskip 16 true pt\flip{$[{1^{4}}]{\times}[{4}]{\times}[{2}]$}\hss}\cg%
\cx{\hskip 16 true pt\flip{$[{4}]{\times}[{3}{1}]{\times}[{2}]$}\hss}\cg%
\cx{\hskip 16 true pt\flip{$[{3}{1}]{\times}[{3}{1}]{\times}[{2}]$}\hss}\cg%
\cx{\hskip 16 true pt\flip{$[{2^{2}}]{\times}[{3}{1}]{\times}[{2}]$}\hss}\cg%
\cx{\hskip 16 true pt\flip{$[{2}{1^{2}}]{\times}[{3}{1}]{\times}[{2}]$}\hss}\cg%
\cx{\hskip 16 true pt\flip{$[{1^{4}}]{\times}[{3}{1}]{\times}[{2}]$}\hss}\cg%
\cx{\hskip 16 true pt\flip{$[{4}]{\times}[{2^{2}}]{\times}[{2}]$}\hss}\cg%
\cx{\hskip 16 true pt\flip{$[{3}{1}]{\times}[{2^{2}}]{\times}[{2}]$}\hss}\cg%
\cx{\hskip 16 true pt\flip{$[{2^{2}}]{\times}[{2^{2}}]{\times}[{2}]$}\hss}\cg%
\cx{\hskip 16 true pt\flip{$[{2}{1^{2}}]{\times}[{2^{2}}]{\times}[{2}]$}\hss}\cg%
\cx{\hskip 16 true pt\flip{$[{1^{4}}]{\times}[{2^{2}}]{\times}[{2}]$}\hss}\cg%
\cx{\hskip 16 true pt\flip{$[{4}]{\times}[{2}{1^{2}}]{\times}[{2}]$}\hss}\cg%
\cx{\hskip 16 true pt\flip{$[{3}{1}]{\times}[{2}{1^{2}}]{\times}[{2}]$}\hss}\cg%
\cx{\hskip 16 true pt\flip{$[{2^{2}}]{\times}[{2}{1^{2}}]{\times}[{2}]$}\hss}\cg%
\cx{\hskip 16 true pt\flip{$[{2}{1^{2}}]{\times}[{2}{1^{2}}]{\times}[{2}]$}\hss}\cg%
\cx{\hskip 16 true pt\flip{$[{1^{4}}]{\times}[{2}{1^{2}}]{\times}[{2}]$}\hss}\cg%
\cx{\hskip 16 true pt\flip{$[{4}]{\times}[{1^{4}}]{\times}[{2}]$}\hss}\cg%
\cx{\hskip 16 true pt\flip{$[{3}{1}]{\times}[{1^{4}}]{\times}[{2}]$}\hss}\cg%
\cx{\hskip 16 true pt\flip{$[{2^{2}}]{\times}[{1^{4}}]{\times}[{2}]$}\hss}\cg%
\cx{\hskip 16 true pt\flip{$[{2}{1^{2}}]{\times}[{1^{4}}]{\times}[{2}]$}\hss}\cg%
\cx{\hskip 16 true pt\flip{$[{1^{4}}]{\times}[{1^{4}}]{\times}[{2}]$}\hss}\cg%
\eol}}\rg%
%
%
\rx{\vss\hfull{%
\rlx{\hss{$1_a$}}\cg%
\e{1}%
\e{0}%
\e{0}%
\e{0}%
\e{0}%
\e{0}%
\e{0}%
\e{0}%
\e{0}%
\e{0}%
\e{0}%
\e{0}%
\e{0}%
\e{0}%
\e{0}%
\e{0}%
\e{0}%
\e{0}%
\e{0}%
\e{0}%
\e{0}%
\e{0}%
\e{0}%
\e{0}%
\e{0}%
\eol}\vss}\rg%
%
%
\rx{\vss\hfull{%
\rlx{\hss{$7_a$}}\cg%
\e{0}%
\e{0}%
\e{0}%
\e{0}%
\e{0}%
\e{0}%
\e{0}%
\e{0}%
\e{0}%
\e{0}%
\e{0}%
\e{0}%
\e{0}%
\e{0}%
\e{0}%
\e{0}%
\e{0}%
\e{0}%
\e{0}%
\e{0}%
\e{0}%
\e{0}%
\e{0}%
\e{0}%
\e{1}%
\eol}\vss}\rg%
%
%
\rx{\vss\hfull{%
\rlx{\hss{$27_a$}}\cg%
\e{2}%
\e{1}%
\e{1}%
\e{0}%
\e{0}%
\e{1}%
\e{1}%
\e{0}%
\e{0}%
\e{0}%
\e{1}%
\e{0}%
\e{0}%
\e{0}%
\e{0}%
\e{0}%
\e{0}%
\e{0}%
\e{0}%
\e{0}%
\e{0}%
\e{0}%
\e{0}%
\e{0}%
\e{0}%
\eol}\vss}\rg%
%
%
\rx{\vss\hfull{%
\rlx{\hss{$21_a$}}\cg%
\e{0}%
\e{0}%
\e{0}%
\e{1}%
\e{0}%
\e{0}%
\e{1}%
\e{0}%
\e{0}%
\e{0}%
\e{0}%
\e{0}%
\e{0}%
\e{0}%
\e{0}%
\e{1}%
\e{0}%
\e{0}%
\e{0}%
\e{0}%
\e{0}%
\e{0}%
\e{0}%
\e{0}%
\e{0}%
\eol}\vss}\rg%
%
%
\rx{\vss\hfull{%
\rlx{\hss{$35_a$}}\cg%
\e{0}%
\e{0}%
\e{0}%
\e{0}%
\e{0}%
\e{0}%
\e{0}%
\e{0}%
\e{0}%
\e{1}%
\e{0}%
\e{0}%
\e{0}%
\e{0}%
\e{0}%
\e{0}%
\e{0}%
\e{0}%
\e{1}%
\e{0}%
\e{0}%
\e{1}%
\e{0}%
\e{0}%
\e{0}%
\eol}\vss}\rg%
%
%
\rx{\vss\hfull{%
\rlx{\hss{$105_a$}}\cg%
\e{0}%
\e{0}%
\e{0}%
\e{0}%
\e{0}%
\e{0}%
\e{0}%
\e{0}%
\e{0}%
\e{1}%
\e{0}%
\e{0}%
\e{0}%
\e{0}%
\e{1}%
\e{0}%
\e{0}%
\e{0}%
\e{2}%
\e{1}%
\e{0}%
\e{1}%
\e{1}%
\e{1}%
\e{1}%
\eol}\vss}\rg%
%
%
\rx{\vss\hfull{%
\rlx{\hss{$189_a$}}\cg%
\e{0}%
\e{0}%
\e{0}%
\e{2}%
\e{0}%
\e{0}%
\e{2}%
\e{1}%
\e{2}%
\e{1}%
\e{0}%
\e{1}%
\e{0}%
\e{1}%
\e{0}%
\e{2}%
\e{2}%
\e{1}%
\e{1}%
\e{0}%
\e{0}%
\e{1}%
\e{0}%
\e{0}%
\e{0}%
\eol}\vss}\rg%
%
%
\rx{\vss\hfull{%
\rlx{\hss{$21_b$}}\cg%
\e{0}%
\e{0}%
\e{0}%
\e{0}%
\e{0}%
\e{0}%
\e{0}%
\e{0}%
\e{0}%
\e{0}%
\e{0}%
\e{0}%
\e{0}%
\e{0}%
\e{1}%
\e{0}%
\e{0}%
\e{0}%
\e{0}%
\e{0}%
\e{0}%
\e{0}%
\e{1}%
\e{0}%
\e{1}%
\eol}\vss}\rg%
%
%
\rx{\vss\hfull{%
\rlx{\hss{$35_b$}}\cg%
\e{2}%
\e{1}%
\e{1}%
\e{0}%
\e{0}%
\e{1}%
\e{1}%
\e{0}%
\e{0}%
\e{0}%
\e{1}%
\e{0}%
\e{1}%
\e{0}%
\e{0}%
\e{0}%
\e{0}%
\e{0}%
\e{0}%
\e{0}%
\e{0}%
\e{0}%
\e{0}%
\e{0}%
\e{0}%
\eol}\vss}\rg%
%
%
\rx{\vss\hfull{%
\rlx{\hss{$189_b$}}\cg%
\e{0}%
\e{0}%
\e{0}%
\e{0}%
\e{0}%
\e{0}%
\e{0}%
\e{0}%
\e{1}%
\e{0}%
\e{0}%
\e{0}%
\e{1}%
\e{1}%
\e{1}%
\e{0}%
\e{1}%
\e{1}%
\e{2}%
\e{2}%
\e{0}%
\e{0}%
\e{1}%
\e{2}%
\e{1}%
\eol}\vss}\rg%
%
%
\rx{\vss\hfull{%
\rlx{\hss{$189_c$}}\cg%
\e{0}%
\e{0}%
\e{0}%
\e{0}%
\e{1}%
\e{0}%
\e{0}%
\e{0}%
\e{1}%
\e{1}%
\e{0}%
\e{0}%
\e{1}%
\e{1}%
\e{2}%
\e{0}%
\e{1}%
\e{1}%
\e{2}%
\e{1}%
\e{1}%
\e{1}%
\e{2}%
\e{1}%
\e{1}%
\eol}\vss}\rg%
%
%
\rx{\vss\hfull{%
\rlx{\hss{$15_a$}}\cg%
\e{0}%
\e{0}%
\e{0}%
\e{0}%
\e{0}%
\e{0}%
\e{0}%
\e{0}%
\e{0}%
\e{0}%
\e{0}%
\e{0}%
\e{1}%
\e{0}%
\e{0}%
\e{0}%
\e{0}%
\e{0}%
\e{0}%
\e{0}%
\e{0}%
\e{0}%
\e{0}%
\e{0}%
\e{1}%
\eol}\vss}\rg%
%
%
\rx{\vss\hfull{%
\rlx{\hss{$105_b$}}\cg%
\e{1}%
\e{2}%
\e{0}%
\e{0}%
\e{0}%
\e{2}%
\e{2}%
\e{1}%
\e{1}%
\e{0}%
\e{0}%
\e{1}%
\e{1}%
\e{0}%
\e{0}%
\e{0}%
\e{1}%
\e{0}%
\e{0}%
\e{0}%
\e{0}%
\e{0}%
\e{0}%
\e{0}%
\e{0}%
\eol}\vss}\rg%
%
%
\rx{\vss\hfull{%
\rlx{\hss{$105_c$}}\cg%
\e{0}%
\e{0}%
\e{0}%
\e{1}%
\e{1}%
\e{0}%
\e{1}%
\e{0}%
\e{1}%
\e{0}%
\e{0}%
\e{0}%
\e{1}%
\e{1}%
\e{1}%
\e{1}%
\e{1}%
\e{1}%
\e{0}%
\e{0}%
\e{1}%
\e{0}%
\e{1}%
\e{0}%
\e{0}%
\eol}\vss}\rg%
%
%
\rx{\vss\hfull{%
\rlx{\hss{$315_a$}}\cg%
\e{0}%
\e{0}%
\e{0}%
\e{0}%
\e{0}%
\e{0}%
\e{1}%
\e{1}%
\e{2}%
\e{1}%
\e{0}%
\e{1}%
\e{0}%
\e{2}%
\e{0}%
\e{0}%
\e{2}%
\e{2}%
\e{4}%
\e{2}%
\e{0}%
\e{1}%
\e{0}%
\e{2}%
\e{0}%
\eol}\vss}\rg%
%
%
\rx{\vss\hfull{%
\rlx{\hss{$405_a$}}\cg%
\e{0}%
\e{2}%
\e{0}%
\e{2}%
\e{0}%
\e{2}%
\e{5}%
\e{3}%
\e{4}%
\e{1}%
\e{0}%
\e{3}%
\e{0}%
\e{2}%
\e{0}%
\e{2}%
\e{4}%
\e{2}%
\e{2}%
\e{0}%
\e{0}%
\e{1}%
\e{0}%
\e{0}%
\e{0}%
\eol}\vss}\rg%
%
%
\rx{\vss\hfull{%
\rlx{\hss{$168_a$}}\cg%
\e{2}%
\e{2}%
\e{2}%
\e{0}%
\e{0}%
\e{2}%
\e{3}%
\e{2}%
\e{1}%
\e{0}%
\e{2}%
\e{2}%
\e{1}%
\e{0}%
\e{0}%
\e{0}%
\e{1}%
\e{0}%
\e{1}%
\e{0}%
\e{0}%
\e{0}%
\e{0}%
\e{0}%
\e{0}%
\eol}\vss}\rg%
%
%
\rx{\vss\hfull{%
\rlx{\hss{$56_a$}}\cg%
\e{0}%
\e{0}%
\e{0}%
\e{0}%
\e{0}%
\e{0}%
\e{0}%
\e{0}%
\e{0}%
\e{0}%
\e{0}%
\e{0}%
\e{0}%
\e{0}%
\e{0}%
\e{0}%
\e{0}%
\e{0}%
\e{1}%
\e{1}%
\e{0}%
\e{0}%
\e{0}%
\e{1}%
\e{1}%
\eol}\vss}\rg%
%
%
\rx{\vss\hfull{%
\rlx{\hss{$120_a$}}\cg%
\e{1}%
\e{2}%
\e{1}%
\e{1}%
\e{0}%
\e{2}%
\e{3}%
\e{1}%
\e{1}%
\e{0}%
\e{1}%
\e{1}%
\e{0}%
\e{0}%
\e{0}%
\e{1}%
\e{1}%
\e{0}%
\e{0}%
\e{0}%
\e{0}%
\e{0}%
\e{0}%
\e{0}%
\e{0}%
\eol}\vss}\rg%
%
%
\rx{\vss\hfull{%
\rlx{\hss{$210_a$}}\cg%
\e{0}%
\e{2}%
\e{1}%
\e{2}%
\e{1}%
\e{2}%
\e{4}%
\e{1}%
\e{2}%
\e{0}%
\e{1}%
\e{1}%
\e{1}%
\e{1}%
\e{0}%
\e{2}%
\e{2}%
\e{1}%
\e{0}%
\e{0}%
\e{1}%
\e{0}%
\e{0}%
\e{0}%
\e{0}%
\eol}\vss}\rg%
%
%
\rx{\vss\hfull{%
\rlx{\hss{$280_a$}}\cg%
\e{0}%
\e{0}%
\e{0}%
\e{0}%
\e{0}%
\e{0}%
\e{0}%
\e{1}%
\e{2}%
\e{2}%
\e{0}%
\e{1}%
\e{0}%
\e{1}%
\e{0}%
\e{0}%
\e{2}%
\e{1}%
\e{4}%
\e{2}%
\e{0}%
\e{2}%
\e{0}%
\e{2}%
\e{0}%
\eol}\vss}\rg%
%
%
\rx{\vss\hfull{%
\rlx{\hss{$336_a$}}\cg%
\e{0}%
\e{0}%
\e{0}%
\e{1}%
\e{1}%
\e{0}%
\e{1}%
\e{1}%
\e{3}%
\e{2}%
\e{0}%
\e{1}%
\e{1}%
\e{2}%
\e{1}%
\e{1}%
\e{3}%
\e{2}%
\e{3}%
\e{1}%
\e{1}%
\e{2}%
\e{1}%
\e{1}%
\e{0}%
\eol}\vss}\rg%
%
%
\rx{\vss\hfull{%
\rlx{\hss{$216_a$}}\cg%
\e{0}%
\e{0}%
\e{1}%
\e{0}%
\e{0}%
\e{0}%
\e{1}%
\e{1}%
\e{1}%
\e{0}%
\e{1}%
\e{1}%
\e{3}%
\e{1}%
\e{1}%
\e{0}%
\e{1}%
\e{1}%
\e{2}%
\e{1}%
\e{0}%
\e{0}%
\e{1}%
\e{1}%
\e{1}%
\eol}\vss}\rg%
%
%
\rx{\vss\hfull{%
\rlx{\hss{$512_a$}}\cg%
\e{0}%
\e{1}%
\e{1}%
\e{1}%
\e{0}%
\e{1}%
\e{4}%
\e{3}%
\e{4}%
\e{1}%
\e{1}%
\e{3}%
\e{2}%
\e{3}%
\e{1}%
\e{1}%
\e{4}%
\e{3}%
\e{4}%
\e{1}%
\e{0}%
\e{1}%
\e{1}%
\e{1}%
\e{0}%
\eol}\vss}\rg%
%
%
\rx{\vss\hfull{%
\rlx{\hss{$378_a$}}\cg%
\e{0}%
\e{0}%
\e{0}%
\e{1}%
\e{0}%
\e{0}%
\e{2}%
\e{1}%
\e{3}%
\e{1}%
\e{0}%
\e{1}%
\e{1}%
\e{3}%
\e{1}%
\e{1}%
\e{3}%
\e{3}%
\e{3}%
\e{1}%
\e{0}%
\e{1}%
\e{1}%
\e{1}%
\e{1}%
\eol}\vss}\rg%
%
%
\rx{\vss\hfull{%
\rlx{\hss{$84_a$}}\cg%
\e{1}%
\e{0}%
\e{1}%
\e{0}%
\e{0}%
\e{0}%
\e{1}%
\e{1}%
\e{0}%
\e{0}%
\e{1}%
\e{1}%
\e{2}%
\e{0}%
\e{0}%
\e{0}%
\e{0}%
\e{0}%
\e{1}%
\e{0}%
\e{0}%
\e{0}%
\e{0}%
\e{0}%
\e{1}%
\eol}\vss}\rg%
%
%
\rx{\vss\hfull{%
\rlx{\hss{$420_a$}}\cg%
\e{0}%
\e{1}%
\e{1}%
\e{1}%
\e{1}%
\e{1}%
\e{3}%
\e{3}%
\e{4}%
\e{1}%
\e{1}%
\e{3}%
\e{1}%
\e{2}%
\e{0}%
\e{1}%
\e{4}%
\e{2}%
\e{3}%
\e{1}%
\e{1}%
\e{1}%
\e{0}%
\e{1}%
\e{0}%
\eol}\vss}\rg%
%
%
\rx{\vss\hfull{%
\rlx{\hss{$280_b$}}\cg%
\e{1}%
\e{2}%
\e{2}%
\e{1}%
\e{0}%
\e{2}%
\e{4}%
\e{2}%
\e{2}%
\e{0}%
\e{2}%
\e{2}%
\e{3}%
\e{1}%
\e{1}%
\e{1}%
\e{2}%
\e{1}%
\e{1}%
\e{0}%
\e{0}%
\e{0}%
\e{1}%
\e{0}%
\e{0}%
\eol}\vss}\rg%
%
%
\rx{\vss\hfull{%
\rlx{\hss{$210_b$}}\cg%
\e{1}%
\e{1}%
\e{0}%
\e{0}%
\e{0}%
\e{1}%
\e{3}%
\e{2}%
\e{1}%
\e{0}%
\e{0}%
\e{2}%
\e{1}%
\e{1}%
\e{0}%
\e{0}%
\e{1}%
\e{1}%
\e{2}%
\e{0}%
\e{0}%
\e{0}%
\e{0}%
\e{0}%
\e{0}%
\eol}\vss}\rg%
%
%
\rx{\vss\hfull{%
\rlx{\hss{$70_a$}}\cg%
\e{0}%
\e{0}%
\e{0}%
\e{0}%
\e{0}%
\e{0}%
\e{1}%
\e{0}%
\e{0}%
\e{0}%
\e{0}%
\e{0}%
\e{0}%
\e{1}%
\e{0}%
\e{0}%
\e{0}%
\e{1}%
\e{1}%
\e{0}%
\e{0}%
\e{0}%
\e{0}%
\e{0}%
\e{0}%
\eol}\vss}\rg%
\eop
\eject
\tablecont%
%
%
%
%
%
%
\rowpts=18 true pt%
\colpts=18 true pt%
\rowlabpts=40 true pt%
\collabpts=90 true pt%
\clx{\vss\hfull{%
\rlx{\hss{$ $}}\cg%
\cx{\hskip 16 true pt\flip{$[{4}]{\times}[{4}]{\times}[{1^{2}}]$}\hss}\cg%
\cx{\hskip 16 true pt\flip{$[{3}{1}]{\times}[{4}]{\times}[{1^{2}}]$}\hss}\cg%
\cx{\hskip 16 true pt\flip{$[{2^{2}}]{\times}[{4}]{\times}[{1^{2}}]$}\hss}\cg%
\cx{\hskip 16 true pt\flip{$[{2}{1^{2}}]{\times}[{4}]{\times}[{1^{2}}]$}\hss}\cg%
\cx{\hskip 16 true pt\flip{$[{1^{4}}]{\times}[{4}]{\times}[{1^{2}}]$}\hss}\cg%
\cx{\hskip 16 true pt\flip{$[{4}]{\times}[{3}{1}]{\times}[{1^{2}}]$}\hss}\cg%
\cx{\hskip 16 true pt\flip{$[{3}{1}]{\times}[{3}{1}]{\times}[{1^{2}}]$}\hss}\cg%
\cx{\hskip 16 true pt\flip{$[{2^{2}}]{\times}[{3}{1}]{\times}[{1^{2}}]$}\hss}\cg%
\cx{\hskip 16 true pt\flip{$[{2}{1^{2}}]{\times}[{3}{1}]{\times}[{1^{2}}]$}\hss}\cg%
\cx{\hskip 16 true pt\flip{$[{1^{4}}]{\times}[{3}{1}]{\times}[{1^{2}}]$}\hss}\cg%
\cx{\hskip 16 true pt\flip{$[{4}]{\times}[{2^{2}}]{\times}[{1^{2}}]$}\hss}\cg%
\cx{\hskip 16 true pt\flip{$[{3}{1}]{\times}[{2^{2}}]{\times}[{1^{2}}]$}\hss}\cg%
\cx{\hskip 16 true pt\flip{$[{2^{2}}]{\times}[{2^{2}}]{\times}[{1^{2}}]$}\hss}\cg%
\cx{\hskip 16 true pt\flip{$[{2}{1^{2}}]{\times}[{2^{2}}]{\times}[{1^{2}}]$}\hss}\cg%
\cx{\hskip 16 true pt\flip{$[{1^{4}}]{\times}[{2^{2}}]{\times}[{1^{2}}]$}\hss}\cg%
\cx{\hskip 16 true pt\flip{$[{4}]{\times}[{2}{1^{2}}]{\times}[{1^{2}}]$}\hss}\cg%
\cx{\hskip 16 true pt\flip{$[{3}{1}]{\times}[{2}{1^{2}}]{\times}[{1^{2}}]$}\hss}\cg%
\cx{\hskip 16 true pt\flip{$[{2^{2}}]{\times}[{2}{1^{2}}]{\times}[{1^{2}}]$}\hss}\cg%
\cx{\hskip 16 true pt\flip{$[{2}{1^{2}}]{\times}[{2}{1^{2}}]{\times}[{1^{2}}]$}\hss}\cg%
\cx{\hskip 16 true pt\flip{$[{1^{4}}]{\times}[{2}{1^{2}}]{\times}[{1^{2}}]$}\hss}\cg%
\cx{\hskip 16 true pt\flip{$[{4}]{\times}[{1^{4}}]{\times}[{1^{2}}]$}\hss}\cg%
\cx{\hskip 16 true pt\flip{$[{3}{1}]{\times}[{1^{4}}]{\times}[{1^{2}}]$}\hss}\cg%
\cx{\hskip 16 true pt\flip{$[{2^{2}}]{\times}[{1^{4}}]{\times}[{1^{2}}]$}\hss}\cg%
\cx{\hskip 16 true pt\flip{$[{2}{1^{2}}]{\times}[{1^{4}}]{\times}[{1^{2}}]$}\hss}\cg%
\cx{\hskip 16 true pt\flip{$[{1^{4}}]{\times}[{1^{4}}]{\times}[{1^{2}}]$}\hss}\cg%
\eol}}\rg%
%
%
\rx{\vss\hfull{%
\rlx{\hss{$1_a$}}\cg%
\e{0}%
\e{0}%
\e{0}%
\e{0}%
\e{0}%
\e{0}%
\e{0}%
\e{0}%
\e{0}%
\e{0}%
\e{0}%
\e{0}%
\e{0}%
\e{0}%
\e{0}%
\e{0}%
\e{0}%
\e{0}%
\e{0}%
\e{0}%
\e{0}%
\e{0}%
\e{0}%
\e{0}%
\e{0}%
\eol}\vss}\rg%
%
%
\rx{\vss\hfull{%
\rlx{\hss{$7_a$}}\cg%
\e{0}%
\e{0}%
\e{0}%
\e{0}%
\e{0}%
\e{0}%
\e{0}%
\e{0}%
\e{0}%
\e{0}%
\e{0}%
\e{0}%
\e{0}%
\e{0}%
\e{0}%
\e{0}%
\e{0}%
\e{0}%
\e{0}%
\e{1}%
\e{0}%
\e{0}%
\e{0}%
\e{1}%
\e{0}%
\eol}\vss}\rg%
%
%
\rx{\vss\hfull{%
\rlx{\hss{$27_a$}}\cg%
\e{0}%
\e{1}%
\e{0}%
\e{0}%
\e{0}%
\e{1}%
\e{0}%
\e{0}%
\e{0}%
\e{0}%
\e{0}%
\e{0}%
\e{0}%
\e{0}%
\e{0}%
\e{0}%
\e{0}%
\e{0}%
\e{0}%
\e{0}%
\e{0}%
\e{0}%
\e{0}%
\e{0}%
\e{0}%
\eol}\vss}\rg%
%
%
\rx{\vss\hfull{%
\rlx{\hss{$21_a$}}\cg%
\e{0}%
\e{1}%
\e{0}%
\e{0}%
\e{0}%
\e{1}%
\e{0}%
\e{0}%
\e{0}%
\e{0}%
\e{0}%
\e{0}%
\e{0}%
\e{0}%
\e{0}%
\e{0}%
\e{0}%
\e{0}%
\e{0}%
\e{0}%
\e{0}%
\e{0}%
\e{0}%
\e{0}%
\e{0}%
\eol}\vss}\rg%
%
%
\rx{\vss\hfull{%
\rlx{\hss{$35_a$}}\cg%
\e{0}%
\e{0}%
\e{0}%
\e{0}%
\e{1}%
\e{0}%
\e{0}%
\e{0}%
\e{1}%
\e{0}%
\e{0}%
\e{0}%
\e{0}%
\e{0}%
\e{0}%
\e{0}%
\e{1}%
\e{0}%
\e{0}%
\e{0}%
\e{1}%
\e{0}%
\e{0}%
\e{0}%
\e{0}%
\eol}\vss}\rg%
%
%
\rx{\vss\hfull{%
\rlx{\hss{$105_a$}}\cg%
\e{0}%
\e{0}%
\e{0}%
\e{0}%
\e{0}%
\e{0}%
\e{0}%
\e{0}%
\e{1}%
\e{1}%
\e{0}%
\e{0}%
\e{0}%
\e{1}%
\e{1}%
\e{0}%
\e{1}%
\e{1}%
\e{2}%
\e{2}%
\e{0}%
\e{1}%
\e{1}%
\e{2}%
\e{0}%
\eol}\vss}\rg%
%
%
\rx{\vss\hfull{%
\rlx{\hss{$189_a$}}\cg%
\e{0}%
\e{1}%
\e{1}%
\e{1}%
\e{1}%
\e{1}%
\e{2}%
\e{1}%
\e{2}%
\e{0}%
\e{1}%
\e{1}%
\e{0}%
\e{0}%
\e{0}%
\e{1}%
\e{2}%
\e{0}%
\e{0}%
\e{0}%
\e{1}%
\e{0}%
\e{0}%
\e{0}%
\e{0}%
\eol}\vss}\rg%
%
%
\rx{\vss\hfull{%
\rlx{\hss{$21_b$}}\cg%
\e{0}%
\e{0}%
\e{0}%
\e{0}%
\e{0}%
\e{0}%
\e{0}%
\e{0}%
\e{0}%
\e{0}%
\e{0}%
\e{0}%
\e{0}%
\e{0}%
\e{0}%
\e{0}%
\e{0}%
\e{0}%
\e{1}%
\e{1}%
\e{0}%
\e{0}%
\e{0}%
\e{1}%
\e{1}%
\eol}\vss}\rg%
%
%
\rx{\vss\hfull{%
\rlx{\hss{$35_b$}}\cg%
\e{1}%
\e{0}%
\e{0}%
\e{0}%
\e{0}%
\e{0}%
\e{1}%
\e{0}%
\e{0}%
\e{0}%
\e{0}%
\e{0}%
\e{0}%
\e{0}%
\e{0}%
\e{0}%
\e{0}%
\e{0}%
\e{0}%
\e{0}%
\e{0}%
\e{0}%
\e{0}%
\e{0}%
\e{0}%
\eol}\vss}\rg%
%
%
\rx{\vss\hfull{%
\rlx{\hss{$189_b$}}\cg%
\e{0}%
\e{0}%
\e{0}%
\e{0}%
\e{0}%
\e{0}%
\e{0}%
\e{1}%
\e{1}%
\e{1}%
\e{0}%
\e{1}%
\e{0}%
\e{2}%
\e{1}%
\e{0}%
\e{1}%
\e{2}%
\e{4}%
\e{3}%
\e{0}%
\e{1}%
\e{1}%
\e{3}%
\e{2}%
\eol}\vss}\rg%
%
%
\rx{\vss\hfull{%
\rlx{\hss{$189_c$}}\cg%
\e{0}%
\e{0}%
\e{0}%
\e{0}%
\e{0}%
\e{0}%
\e{1}%
\e{0}%
\e{2}%
\e{1}%
\e{0}%
\e{0}%
\e{1}%
\e{1}%
\e{1}%
\e{0}%
\e{2}%
\e{1}%
\e{4}%
\e{1}%
\e{0}%
\e{1}%
\e{1}%
\e{1}%
\e{1}%
\eol}\vss}\rg%
%
%
\rx{\vss\hfull{%
\rlx{\hss{$15_a$}}\cg%
\e{0}%
\e{0}%
\e{0}%
\e{0}%
\e{0}%
\e{0}%
\e{0}%
\e{0}%
\e{0}%
\e{0}%
\e{0}%
\e{0}%
\e{0}%
\e{0}%
\e{0}%
\e{0}%
\e{0}%
\e{0}%
\e{1}%
\e{0}%
\e{0}%
\e{0}%
\e{0}%
\e{0}%
\e{1}%
\eol}\vss}\rg%
%
%
\rx{\vss\hfull{%
\rlx{\hss{$105_b$}}\cg%
\e{1}%
\e{1}%
\e{0}%
\e{0}%
\e{0}%
\e{1}%
\e{1}%
\e{1}%
\e{0}%
\e{0}%
\e{0}%
\e{1}%
\e{0}%
\e{1}%
\e{0}%
\e{0}%
\e{0}%
\e{1}%
\e{0}%
\e{0}%
\e{0}%
\e{0}%
\e{0}%
\e{0}%
\e{0}%
\eol}\vss}\rg%
%
%
\rx{\vss\hfull{%
\rlx{\hss{$105_c$}}\cg%
\e{1}%
\e{0}%
\e{0}%
\e{0}%
\e{0}%
\e{0}%
\e{2}%
\e{0}%
\e{1}%
\e{0}%
\e{0}%
\e{0}%
\e{1}%
\e{0}%
\e{0}%
\e{0}%
\e{1}%
\e{0}%
\e{1}%
\e{0}%
\e{0}%
\e{0}%
\e{0}%
\e{0}%
\e{0}%
\eol}\vss}\rg%
%
%
\rx{\vss\hfull{%
\rlx{\hss{$315_a$}}\cg%
\e{0}%
\e{0}%
\e{0}%
\e{1}%
\e{0}%
\e{0}%
\e{1}%
\e{2}%
\e{2}%
\e{2}%
\e{0}%
\e{2}%
\e{1}%
\e{3}%
\e{1}%
\e{1}%
\e{2}%
\e{3}%
\e{4}%
\e{2}%
\e{0}%
\e{2}%
\e{1}%
\e{2}%
\e{1}%
\eol}\vss}\rg%
%
%
\rx{\vss\hfull{%
\rlx{\hss{$405_a$}}\cg%
\e{1}%
\e{2}%
\e{1}%
\e{2}%
\e{0}%
\e{2}%
\e{4}%
\e{3}%
\e{2}%
\e{1}%
\e{1}%
\e{3}%
\e{1}%
\e{2}%
\e{0}%
\e{2}%
\e{2}%
\e{2}%
\e{1}%
\e{0}%
\e{0}%
\e{1}%
\e{0}%
\e{0}%
\e{0}%
\eol}\vss}\rg%
%
%
\rx{\vss\hfull{%
\rlx{\hss{$168_a$}}\cg%
\e{0}%
\e{1}%
\e{1}%
\e{1}%
\e{0}%
\e{1}%
\e{2}%
\e{1}%
\e{1}%
\e{0}%
\e{1}%
\e{1}%
\e{0}%
\e{0}%
\e{0}%
\e{1}%
\e{1}%
\e{0}%
\e{0}%
\e{0}%
\e{0}%
\e{0}%
\e{0}%
\e{0}%
\e{0}%
\eol}\vss}\rg%
%
%
\rx{\vss\hfull{%
\rlx{\hss{$56_a$}}\cg%
\e{0}%
\e{0}%
\e{0}%
\e{0}%
\e{0}%
\e{0}%
\e{0}%
\e{0}%
\e{0}%
\e{1}%
\e{0}%
\e{0}%
\e{0}%
\e{1}%
\e{0}%
\e{0}%
\e{0}%
\e{1}%
\e{1}%
\e{2}%
\e{0}%
\e{1}%
\e{0}%
\e{2}%
\e{1}%
\eol}\vss}\rg%
%
%
\rx{\vss\hfull{%
\rlx{\hss{$120_a$}}\cg%
\e{1}%
\e{2}%
\e{0}%
\e{1}%
\e{0}%
\e{2}%
\e{1}%
\e{1}%
\e{0}%
\e{0}%
\e{0}%
\e{1}%
\e{0}%
\e{0}%
\e{0}%
\e{1}%
\e{0}%
\e{0}%
\e{0}%
\e{0}%
\e{0}%
\e{0}%
\e{0}%
\e{0}%
\e{0}%
\eol}\vss}\rg%
%
%
\rx{\vss\hfull{%
\rlx{\hss{$210_a$}}\cg%
\e{1}%
\e{2}%
\e{1}%
\e{1}%
\e{0}%
\e{2}%
\e{3}%
\e{1}%
\e{1}%
\e{0}%
\e{1}%
\e{1}%
\e{0}%
\e{0}%
\e{0}%
\e{1}%
\e{1}%
\e{0}%
\e{0}%
\e{0}%
\e{0}%
\e{0}%
\e{0}%
\e{0}%
\e{0}%
\eol}\vss}\rg%
%
%
\rx{\vss\hfull{%
\rlx{\hss{$280_a$}}\cg%
\e{0}%
\e{0}%
\e{0}%
\e{1}%
\e{1}%
\e{0}%
\e{1}%
\e{1}%
\e{3}%
\e{2}%
\e{0}%
\e{1}%
\e{1}%
\e{2}%
\e{1}%
\e{1}%
\e{3}%
\e{2}%
\e{3}%
\e{1}%
\e{1}%
\e{2}%
\e{1}%
\e{1}%
\e{0}%
\eol}\vss}\rg%
%
%
\rx{\vss\hfull{%
\rlx{\hss{$336_a$}}\cg%
\e{0}%
\e{0}%
\e{1}%
\e{1}%
\e{1}%
\e{0}%
\e{3}%
\e{1}%
\e{4}%
\e{1}%
\e{1}%
\e{1}%
\e{1}%
\e{1}%
\e{1}%
\e{1}%
\e{4}%
\e{1}%
\e{3}%
\e{0}%
\e{1}%
\e{1}%
\e{1}%
\e{0}%
\e{0}%
\eol}\vss}\rg%
%
%
\rx{\vss\hfull{%
\rlx{\hss{$216_a$}}\cg%
\e{0}%
\e{0}%
\e{0}%
\e{0}%
\e{0}%
\e{0}%
\e{1}%
\e{1}%
\e{2}%
\e{0}%
\e{0}%
\e{1}%
\e{1}%
\e{1}%
\e{1}%
\e{0}%
\e{2}%
\e{1}%
\e{4}%
\e{1}%
\e{0}%
\e{0}%
\e{1}%
\e{1}%
\e{1}%
\eol}\vss}\rg%
%
%
\rx{\vss\hfull{%
\rlx{\hss{$512_a$}}\cg%
\e{0}%
\e{1}%
\e{1}%
\e{1}%
\e{0}%
\e{1}%
\e{4}%
\e{3}%
\e{4}%
\e{1}%
\e{1}%
\e{3}%
\e{2}%
\e{3}%
\e{1}%
\e{1}%
\e{4}%
\e{3}%
\e{4}%
\e{1}%
\e{0}%
\e{1}%
\e{1}%
\e{1}%
\e{0}%
\eol}\vss}\rg%
%
%
\rx{\vss\hfull{%
\rlx{\hss{$378_a$}}\cg%
\e{0}%
\e{1}%
\e{0}%
\e{0}%
\e{0}%
\e{1}%
\e{2}%
\e{2}%
\e{3}%
\e{1}%
\e{0}%
\e{2}%
\e{2}%
\e{3}%
\e{1}%
\e{0}%
\e{3}%
\e{3}%
\e{4}%
\e{2}%
\e{0}%
\e{1}%
\e{1}%
\e{2}%
\e{0}%
\eol}\vss}\rg%
%
%
\rx{\vss\hfull{%
\rlx{\hss{$84_a$}}\cg%
\e{0}%
\e{0}%
\e{0}%
\e{0}%
\e{0}%
\e{0}%
\e{1}%
\e{0}%
\e{1}%
\e{0}%
\e{0}%
\e{0}%
\e{1}%
\e{0}%
\e{0}%
\e{0}%
\e{1}%
\e{0}%
\e{1}%
\e{0}%
\e{0}%
\e{0}%
\e{0}%
\e{0}%
\e{0}%
\eol}\vss}\rg%
%
%
\rx{\vss\hfull{%
\rlx{\hss{$420_a$}}\cg%
\e{0}%
\e{1}%
\e{1}%
\e{2}%
\e{1}%
\e{1}%
\e{4}%
\e{2}%
\e{4}%
\e{1}%
\e{1}%
\e{2}%
\e{2}%
\e{1}%
\e{0}%
\e{2}%
\e{4}%
\e{1}%
\e{2}%
\e{0}%
\e{1}%
\e{1}%
\e{0}%
\e{0}%
\e{0}%
\eol}\vss}\rg%
%
%
\rx{\vss\hfull{%
\rlx{\hss{$280_b$}}\cg%
\e{1}%
\e{1}%
\e{1}%
\e{0}%
\e{0}%
\e{1}%
\e{4}%
\e{1}%
\e{2}%
\e{0}%
\e{1}%
\e{1}%
\e{1}%
\e{1}%
\e{0}%
\e{0}%
\e{2}%
\e{1}%
\e{1}%
\e{0}%
\e{0}%
\e{0}%
\e{0}%
\e{0}%
\e{0}%
\eol}\vss}\rg%
%
%
\rx{\vss\hfull{%
\rlx{\hss{$210_b$}}\cg%
\e{0}%
\e{1}%
\e{0}%
\e{0}%
\e{0}%
\e{1}%
\e{1}%
\e{2}%
\e{1}%
\e{0}%
\e{0}%
\e{2}%
\e{1}%
\e{2}%
\e{0}%
\e{0}%
\e{1}%
\e{2}%
\e{1}%
\e{1}%
\e{0}%
\e{0}%
\e{0}%
\e{1}%
\e{0}%
\eol}\vss}\rg%
%
%
\rx{\vss\hfull{%
\rlx{\hss{$70_a$}}\cg%
\e{0}%
\e{0}%
\e{0}%
\e{0}%
\e{0}%
\e{0}%
\e{0}%
\e{1}%
\e{0}%
\e{0}%
\e{0}%
\e{1}%
\e{0}%
\e{1}%
\e{0}%
\e{0}%
\e{0}%
\e{1}%
\e{1}%
\e{1}%
\e{0}%
\e{0}%
\e{0}%
\e{1}%
\e{1}%
\eol}\vss}\rg%
\tableclose%
%
%
%
%
%
%
\eop
\eject
\tableopen{Induce/restrict matrix for $W(D_{6})\,\subset\,W(E_{7})$}%
%
%
%
%
%
%
\rowpts=18 true pt%
\colpts=18 true pt%
\rowlabpts=40 true pt%
\collabpts=60 true pt%
\clx{\vss\hfull{%
\rlx{\hss{$ $}}\cg%
\cx{\hskip 16 true pt\flip{$[{6}:-]$}\hss}\cg%
\cx{\hskip 16 true pt\flip{$[{5}{1}:-]$}\hss}\cg%
\cx{\hskip 16 true pt\flip{$[{4}{2}:-]$}\hss}\cg%
\cx{\hskip 16 true pt\flip{$[{4}{1^{2}}:-]$}\hss}\cg%
\cx{\hskip 16 true pt\flip{$[{3^{2}}:-]$}\hss}\cg%
\cx{\hskip 16 true pt\flip{$[{3}{2}{1}:-]$}\hss}\cg%
\cx{\hskip 16 true pt\flip{$[{3}{1^{3}}:-]$}\hss}\cg%
\cx{\hskip 16 true pt\flip{$[{2^{3}}:-]$}\hss}\cg%
\cx{\hskip 16 true pt\flip{$[{2^{2}}{1^{2}}:-]$}\hss}\cg%
\cx{\hskip 16 true pt\flip{$[{2}{1^{4}}:-]$}\hss}\cg%
\cx{\hskip 16 true pt\flip{$[{1^{6}}:-]$}\hss}\cg%
\cx{\hskip 16 true pt\flip{$[{5}:{1}]$}\hss}\cg%
\cx{\hskip 16 true pt\flip{$[{4}{1}:{1}]$}\hss}\cg%
\cx{\hskip 16 true pt\flip{$[{3}{2}:{1}]$}\hss}\cg%
\cx{\hskip 16 true pt\flip{$[{3}{1^{2}}:{1}]$}\hss}\cg%
\cx{\hskip 16 true pt\flip{$[{2^{2}}{1}:{1}]$}\hss}\cg%
\cx{\hskip 16 true pt\flip{$[{2}{1^{3}}:{1}]$}\hss}\cg%
\cx{\hskip 16 true pt\flip{$[{1^{5}}:{1}]$}\hss}\cg%
\cx{\hskip 16 true pt\flip{$[{4}:{2}]$}\hss}\cg%
\eol}}\rg%
%
%
\rx{\vss\hfull{%
\rlx{\hss{$1_a$}}\cg%
\e{1}%
\e{0}%
\e{0}%
\e{0}%
\e{0}%
\e{0}%
\e{0}%
\e{0}%
\e{0}%
\e{0}%
\e{0}%
\e{0}%
\e{0}%
\e{0}%
\e{0}%
\e{0}%
\e{0}%
\e{0}%
\e{0}%
\eol}\vss}\rg%
%
%
\rx{\vss\hfull{%
\rlx{\hss{$7_a$}}\cg%
\e{0}%
\e{0}%
\e{0}%
\e{0}%
\e{0}%
\e{0}%
\e{0}%
\e{0}%
\e{0}%
\e{0}%
\e{1}%
\e{0}%
\e{0}%
\e{0}%
\e{0}%
\e{0}%
\e{0}%
\e{1}%
\e{0}%
\eol}\vss}\rg%
%
%
\rx{\vss\hfull{%
\rlx{\hss{$27_a$}}\cg%
\e{1}%
\e{1}%
\e{0}%
\e{0}%
\e{0}%
\e{0}%
\e{0}%
\e{0}%
\e{0}%
\e{0}%
\e{0}%
\e{1}%
\e{0}%
\e{0}%
\e{0}%
\e{0}%
\e{0}%
\e{0}%
\e{1}%
\eol}\vss}\rg%
%
%
\rx{\vss\hfull{%
\rlx{\hss{$21_a$}}\cg%
\e{0}%
\e{0}%
\e{0}%
\e{0}%
\e{0}%
\e{0}%
\e{0}%
\e{0}%
\e{0}%
\e{0}%
\e{0}%
\e{1}%
\e{0}%
\e{0}%
\e{0}%
\e{0}%
\e{0}%
\e{0}%
\e{0}%
\eol}\vss}\rg%
%
%
\rx{\vss\hfull{%
\rlx{\hss{$35_a$}}\cg%
\e{0}%
\e{0}%
\e{0}%
\e{0}%
\e{0}%
\e{0}%
\e{0}%
\e{0}%
\e{0}%
\e{0}%
\e{0}%
\e{0}%
\e{0}%
\e{0}%
\e{0}%
\e{0}%
\e{0}%
\e{0}%
\e{0}%
\eol}\vss}\rg%
%
%
\rx{\vss\hfull{%
\rlx{\hss{$105_a$}}\cg%
\e{0}%
\e{0}%
\e{0}%
\e{0}%
\e{0}%
\e{0}%
\e{0}%
\e{0}%
\e{0}%
\e{1}%
\e{0}%
\e{0}%
\e{0}%
\e{0}%
\e{0}%
\e{0}%
\e{1}%
\e{1}%
\e{0}%
\eol}\vss}\rg%
%
%
\rx{\vss\hfull{%
\rlx{\hss{$189_a$}}\cg%
\e{0}%
\e{0}%
\e{0}%
\e{0}%
\e{0}%
\e{0}%
\e{0}%
\e{0}%
\e{0}%
\e{0}%
\e{0}%
\e{0}%
\e{1}%
\e{0}%
\e{0}%
\e{0}%
\e{0}%
\e{0}%
\e{0}%
\eol}\vss}\rg%
%
%
\rx{\vss\hfull{%
\rlx{\hss{$21_b$}}\cg%
\e{0}%
\e{0}%
\e{0}%
\e{0}%
\e{0}%
\e{0}%
\e{0}%
\e{0}%
\e{0}%
\e{1}%
\e{0}%
\e{0}%
\e{0}%
\e{0}%
\e{0}%
\e{0}%
\e{0}%
\e{1}%
\e{0}%
\eol}\vss}\rg%
%
%
\rx{\vss\hfull{%
\rlx{\hss{$35_b$}}\cg%
\e{1}%
\e{0}%
\e{1}%
\e{0}%
\e{0}%
\e{0}%
\e{0}%
\e{0}%
\e{0}%
\e{0}%
\e{0}%
\e{0}%
\e{0}%
\e{0}%
\e{0}%
\e{0}%
\e{0}%
\e{0}%
\e{1}%
\eol}\vss}\rg%
%
%
\rx{\vss\hfull{%
\rlx{\hss{$189_b$}}\cg%
\e{0}%
\e{0}%
\e{0}%
\e{0}%
\e{0}%
\e{0}%
\e{0}%
\e{0}%
\e{1}%
\e{0}%
\e{0}%
\e{0}%
\e{0}%
\e{0}%
\e{0}%
\e{1}%
\e{1}%
\e{1}%
\e{0}%
\eol}\vss}\rg%
%
%
\rx{\vss\hfull{%
\rlx{\hss{$189_c$}}\cg%
\e{0}%
\e{0}%
\e{0}%
\e{0}%
\e{0}%
\e{0}%
\e{1}%
\e{0}%
\e{0}%
\e{1}%
\e{0}%
\e{0}%
\e{0}%
\e{0}%
\e{0}%
\e{0}%
\e{1}%
\e{0}%
\e{0}%
\eol}\vss}\rg%
%
%
\rx{\vss\hfull{%
\rlx{\hss{$15_a$}}\cg%
\e{0}%
\e{0}%
\e{0}%
\e{0}%
\e{0}%
\e{0}%
\e{0}%
\e{1}%
\e{0}%
\e{0}%
\e{0}%
\e{0}%
\e{0}%
\e{0}%
\e{0}%
\e{0}%
\e{0}%
\e{0}%
\e{0}%
\eol}\vss}\rg%
%
%
\rx{\vss\hfull{%
\rlx{\hss{$105_b$}}\cg%
\e{0}%
\e{0}%
\e{0}%
\e{0}%
\e{1}%
\e{0}%
\e{0}%
\e{0}%
\e{0}%
\e{0}%
\e{0}%
\e{0}%
\e{0}%
\e{1}%
\e{0}%
\e{0}%
\e{0}%
\e{0}%
\e{1}%
\eol}\vss}\rg%
%
%
\rx{\vss\hfull{%
\rlx{\hss{$105_c$}}\cg%
\e{0}%
\e{0}%
\e{0}%
\e{0}%
\e{0}%
\e{0}%
\e{1}%
\e{0}%
\e{0}%
\e{0}%
\e{0}%
\e{0}%
\e{0}%
\e{0}%
\e{0}%
\e{0}%
\e{0}%
\e{0}%
\e{0}%
\eol}\vss}\rg%
%
%
\rx{\vss\hfull{%
\rlx{\hss{$315_a$}}\cg%
\e{0}%
\e{0}%
\e{0}%
\e{0}%
\e{0}%
\e{0}%
\e{0}%
\e{0}%
\e{0}%
\e{0}%
\e{0}%
\e{0}%
\e{0}%
\e{0}%
\e{1}%
\e{1}%
\e{1}%
\e{0}%
\e{0}%
\eol}\vss}\rg%
%
%
\rx{\vss\hfull{%
\rlx{\hss{$405_a$}}\cg%
\e{0}%
\e{0}%
\e{0}%
\e{0}%
\e{0}%
\e{0}%
\e{0}%
\e{0}%
\e{0}%
\e{0}%
\e{0}%
\e{0}%
\e{1}%
\e{1}%
\e{1}%
\e{0}%
\e{0}%
\e{0}%
\e{0}%
\eol}\vss}\rg%
%
%
\rx{\vss\hfull{%
\rlx{\hss{$168_a$}}\cg%
\e{0}%
\e{1}%
\e{1}%
\e{0}%
\e{0}%
\e{0}%
\e{0}%
\e{0}%
\e{0}%
\e{0}%
\e{0}%
\e{0}%
\e{1}%
\e{0}%
\e{0}%
\e{0}%
\e{0}%
\e{0}%
\e{1}%
\eol}\vss}\rg%
%
%
\rx{\vss\hfull{%
\rlx{\hss{$56_a$}}\cg%
\e{0}%
\e{0}%
\e{0}%
\e{0}%
\e{0}%
\e{0}%
\e{0}%
\e{0}%
\e{0}%
\e{0}%
\e{1}%
\e{0}%
\e{0}%
\e{0}%
\e{0}%
\e{0}%
\e{1}%
\e{1}%
\e{0}%
\eol}\vss}\rg%
%
%
\rx{\vss\hfull{%
\rlx{\hss{$120_a$}}\cg%
\e{0}%
\e{1}%
\e{0}%
\e{0}%
\e{0}%
\e{0}%
\e{0}%
\e{0}%
\e{0}%
\e{0}%
\e{0}%
\e{1}%
\e{1}%
\e{0}%
\e{0}%
\e{0}%
\e{0}%
\e{0}%
\e{1}%
\eol}\vss}\rg%
%
%
\rx{\vss\hfull{%
\rlx{\hss{$210_a$}}\cg%
\e{0}%
\e{0}%
\e{0}%
\e{1}%
\e{0}%
\e{0}%
\e{0}%
\e{0}%
\e{0}%
\e{0}%
\e{0}%
\e{1}%
\e{1}%
\e{0}%
\e{0}%
\e{0}%
\e{0}%
\e{0}%
\e{1}%
\eol}\vss}\rg%
%
%
\rx{\vss\hfull{%
\rlx{\hss{$280_a$}}\cg%
\e{0}%
\e{0}%
\e{0}%
\e{0}%
\e{0}%
\e{0}%
\e{0}%
\e{0}%
\e{0}%
\e{0}%
\e{0}%
\e{0}%
\e{0}%
\e{0}%
\e{1}%
\e{0}%
\e{1}%
\e{0}%
\e{0}%
\eol}\vss}\rg%
%
%
\rx{\vss\hfull{%
\rlx{\hss{$336_a$}}\cg%
\e{0}%
\e{0}%
\e{0}%
\e{0}%
\e{0}%
\e{0}%
\e{1}%
\e{0}%
\e{0}%
\e{0}%
\e{0}%
\e{0}%
\e{0}%
\e{0}%
\e{1}%
\e{0}%
\e{0}%
\e{0}%
\e{0}%
\eol}\vss}\rg%
%
%
\rx{\vss\hfull{%
\rlx{\hss{$216_a$}}\cg%
\e{0}%
\e{0}%
\e{0}%
\e{0}%
\e{0}%
\e{1}%
\e{0}%
\e{1}%
\e{0}%
\e{0}%
\e{0}%
\e{0}%
\e{0}%
\e{0}%
\e{0}%
\e{1}%
\e{0}%
\e{0}%
\e{0}%
\eol}\vss}\rg%
%
%
\rx{\vss\hfull{%
\rlx{\hss{$512_a$}}\cg%
\e{0}%
\e{0}%
\e{0}%
\e{0}%
\e{0}%
\e{1}%
\e{0}%
\e{0}%
\e{0}%
\e{0}%
\e{0}%
\e{0}%
\e{0}%
\e{1}%
\e{1}%
\e{1}%
\e{0}%
\e{0}%
\e{0}%
\eol}\vss}\rg%
%
%
\rx{\vss\hfull{%
\rlx{\hss{$378_a$}}\cg%
\e{0}%
\e{0}%
\e{0}%
\e{0}%
\e{0}%
\e{0}%
\e{0}%
\e{0}%
\e{1}%
\e{0}%
\e{0}%
\e{0}%
\e{0}%
\e{1}%
\e{0}%
\e{1}%
\e{1}%
\e{0}%
\e{0}%
\eol}\vss}\rg%
%
%
\rx{\vss\hfull{%
\rlx{\hss{$84_a$}}\cg%
\e{0}%
\e{0}%
\e{1}%
\e{0}%
\e{0}%
\e{0}%
\e{0}%
\e{1}%
\e{0}%
\e{0}%
\e{0}%
\e{0}%
\e{0}%
\e{0}%
\e{0}%
\e{0}%
\e{0}%
\e{0}%
\e{0}%
\eol}\vss}\rg%
%
%
\rx{\vss\hfull{%
\rlx{\hss{$420_a$}}\cg%
\e{0}%
\e{0}%
\e{0}%
\e{1}%
\e{0}%
\e{0}%
\e{0}%
\e{0}%
\e{0}%
\e{0}%
\e{0}%
\e{0}%
\e{1}%
\e{0}%
\e{1}%
\e{0}%
\e{0}%
\e{0}%
\e{0}%
\eol}\vss}\rg%
%
%
\rx{\vss\hfull{%
\rlx{\hss{$280_b$}}\cg%
\e{0}%
\e{0}%
\e{1}%
\e{0}%
\e{0}%
\e{1}%
\e{0}%
\e{0}%
\e{0}%
\e{0}%
\e{0}%
\e{0}%
\e{0}%
\e{1}%
\e{0}%
\e{0}%
\e{0}%
\e{0}%
\e{1}%
\eol}\vss}\rg%
%
%
\rx{\vss\hfull{%
\rlx{\hss{$210_b$}}\cg%
\e{0}%
\e{0}%
\e{0}%
\e{0}%
\e{1}%
\e{0}%
\e{0}%
\e{0}%
\e{0}%
\e{0}%
\e{0}%
\e{0}%
\e{0}%
\e{1}%
\e{0}%
\e{1}%
\e{0}%
\e{0}%
\e{0}%
\eol}\vss}\rg%
%
%
\rx{\vss\hfull{%
\rlx{\hss{$70_a$}}\cg%
\e{0}%
\e{0}%
\e{0}%
\e{0}%
\e{0}%
\e{0}%
\e{0}%
\e{0}%
\e{0}%
\e{0}%
\e{0}%
\e{0}%
\e{0}%
\e{0}%
\e{0}%
\e{1}%
\e{0}%
\e{0}%
\e{0}%
\eol}\vss}\rg%
\eop
\eject
\tablecont%
%
%
%
%
%
%
\rowpts=18 true pt%
\colpts=18 true pt%
\rowlabpts=40 true pt%
\collabpts=60 true pt%
\clx{\vss\hfull{%
\rlx{\hss{$ $}}\cg%
\cx{\hskip 16 true pt\flip{$[{4}:{1^{2}}]$}\hss}\cg%
\cx{\hskip 16 true pt\flip{$[{3}{1}:{2}]$}\hss}\cg%
\cx{\hskip 16 true pt\flip{$[{3}{1}:{1^{2}}]$}\hss}\cg%
\cx{\hskip 16 true pt\flip{$[{2^{2}}:{2}]$}\hss}\cg%
\cx{\hskip 16 true pt\flip{$[{2^{2}}:{1^{2}}]$}\hss}\cg%
\cx{\hskip 16 true pt\flip{$[{2}{1^{2}}:{2}]$}\hss}\cg%
\cx{\hskip 16 true pt\flip{$[{2}{1^{2}}:{1^{2}}]$}\hss}\cg%
\cx{\hskip 16 true pt\flip{$[{1^{4}}:{2}]$}\hss}\cg%
\cx{\hskip 16 true pt\flip{$[{1^{4}}:{1^{2}}]$}\hss}\cg%
\cx{\hskip 16 true pt\flip{$[{3}:{3}]^{+}$}\hss}\cg%
\cx{\hskip 16 true pt\flip{$[{3}:{3}]^{-}$}\hss}\cg%
\cx{\hskip 16 true pt\flip{$[{3}:{2}{1}]$}\hss}\cg%
\cx{\hskip 16 true pt\flip{$[{3}:{1^{3}}]$}\hss}\cg%
\cx{\hskip 16 true pt\flip{$[{2}{1}:{2}{1}]^{+}$}\hss}\cg%
\cx{\hskip 16 true pt\flip{$[{2}{1}:{2}{1}]^{-}$}\hss}\cg%
\cx{\hskip 16 true pt\flip{$[{2}{1}:{1^{3}}]$}\hss}\cg%
\cx{\hskip 16 true pt\flip{$[{1^{3}}:{1^{3}}]^{+}$}\hss}\cg%
\cx{\hskip 16 true pt\flip{$[{1^{3}}:{1^{3}}]^{-}$}\hss}\cg%
\eol}}\rg%
%
%
\rx{\vss\hfull{%
\rlx{\hss{$1_a$}}\cg%
\e{0}%
\e{0}%
\e{0}%
\e{0}%
\e{0}%
\e{0}%
\e{0}%
\e{0}%
\e{0}%
\e{0}%
\e{0}%
\e{0}%
\e{0}%
\e{0}%
\e{0}%
\e{0}%
\e{0}%
\e{0}%
\eol}\vss}\rg%
%
%
\rx{\vss\hfull{%
\rlx{\hss{$7_a$}}\cg%
\e{0}%
\e{0}%
\e{0}%
\e{0}%
\e{0}%
\e{0}%
\e{0}%
\e{0}%
\e{0}%
\e{0}%
\e{0}%
\e{0}%
\e{0}%
\e{0}%
\e{0}%
\e{0}%
\e{0}%
\e{0}%
\eol}\vss}\rg%
%
%
\rx{\vss\hfull{%
\rlx{\hss{$27_a$}}\cg%
\e{0}%
\e{0}%
\e{0}%
\e{0}%
\e{0}%
\e{0}%
\e{0}%
\e{0}%
\e{0}%
\e{0}%
\e{0}%
\e{0}%
\e{0}%
\e{0}%
\e{0}%
\e{0}%
\e{0}%
\e{0}%
\eol}\vss}\rg%
%
%
\rx{\vss\hfull{%
\rlx{\hss{$21_a$}}\cg%
\e{1}%
\e{0}%
\e{0}%
\e{0}%
\e{0}%
\e{0}%
\e{0}%
\e{0}%
\e{0}%
\e{0}%
\e{0}%
\e{0}%
\e{0}%
\e{0}%
\e{0}%
\e{0}%
\e{0}%
\e{0}%
\eol}\vss}\rg%
%
%
\rx{\vss\hfull{%
\rlx{\hss{$35_a$}}\cg%
\e{0}%
\e{0}%
\e{0}%
\e{0}%
\e{0}%
\e{0}%
\e{0}%
\e{1}%
\e{0}%
\e{0}%
\e{0}%
\e{0}%
\e{1}%
\e{0}%
\e{0}%
\e{0}%
\e{0}%
\e{0}%
\eol}\vss}\rg%
%
%
\rx{\vss\hfull{%
\rlx{\hss{$105_a$}}\cg%
\e{0}%
\e{0}%
\e{0}%
\e{0}%
\e{0}%
\e{0}%
\e{0}%
\e{1}%
\e{1}%
\e{0}%
\e{0}%
\e{0}%
\e{0}%
\e{0}%
\e{0}%
\e{1}%
\e{0}%
\e{0}%
\eol}\vss}\rg%
%
%
\rx{\vss\hfull{%
\rlx{\hss{$189_a$}}\cg%
\e{1}%
\e{0}%
\e{1}%
\e{0}%
\e{0}%
\e{1}%
\e{0}%
\e{0}%
\e{0}%
\e{0}%
\e{0}%
\e{1}%
\e{1}%
\e{0}%
\e{0}%
\e{0}%
\e{0}%
\e{0}%
\eol}\vss}\rg%
%
%
\rx{\vss\hfull{%
\rlx{\hss{$21_b$}}\cg%
\e{0}%
\e{0}%
\e{0}%
\e{0}%
\e{0}%
\e{0}%
\e{0}%
\e{0}%
\e{0}%
\e{0}%
\e{0}%
\e{0}%
\e{0}%
\e{0}%
\e{0}%
\e{0}%
\e{0}%
\e{1}%
\eol}\vss}\rg%
%
%
\rx{\vss\hfull{%
\rlx{\hss{$35_b$}}\cg%
\e{0}%
\e{0}%
\e{0}%
\e{0}%
\e{0}%
\e{0}%
\e{0}%
\e{0}%
\e{0}%
\e{0}%
\e{1}%
\e{0}%
\e{0}%
\e{0}%
\e{0}%
\e{0}%
\e{0}%
\e{0}%
\eol}\vss}\rg%
%
%
\rx{\vss\hfull{%
\rlx{\hss{$189_b$}}\cg%
\e{0}%
\e{0}%
\e{0}%
\e{0}%
\e{0}%
\e{0}%
\e{1}%
\e{0}%
\e{1}%
\e{0}%
\e{0}%
\e{0}%
\e{0}%
\e{0}%
\e{0}%
\e{1}%
\e{1}%
\e{1}%
\eol}\vss}\rg%
%
%
\rx{\vss\hfull{%
\rlx{\hss{$189_c$}}\cg%
\e{0}%
\e{0}%
\e{0}%
\e{0}%
\e{0}%
\e{0}%
\e{1}%
\e{1}%
\e{0}%
\e{0}%
\e{0}%
\e{0}%
\e{0}%
\e{0}%
\e{1}%
\e{1}%
\e{0}%
\e{1}%
\eol}\vss}\rg%
%
%
\rx{\vss\hfull{%
\rlx{\hss{$15_a$}}\cg%
\e{0}%
\e{0}%
\e{0}%
\e{0}%
\e{0}%
\e{0}%
\e{0}%
\e{0}%
\e{0}%
\e{0}%
\e{0}%
\e{0}%
\e{0}%
\e{0}%
\e{0}%
\e{0}%
\e{0}%
\e{1}%
\eol}\vss}\rg%
%
%
\rx{\vss\hfull{%
\rlx{\hss{$105_b$}}\cg%
\e{0}%
\e{1}%
\e{0}%
\e{0}%
\e{0}%
\e{0}%
\e{0}%
\e{0}%
\e{0}%
\e{1}%
\e{0}%
\e{0}%
\e{0}%
\e{0}%
\e{0}%
\e{0}%
\e{0}%
\e{0}%
\eol}\vss}\rg%
%
%
\rx{\vss\hfull{%
\rlx{\hss{$105_c$}}\cg%
\e{0}%
\e{0}%
\e{1}%
\e{0}%
\e{0}%
\e{0}%
\e{0}%
\e{0}%
\e{0}%
\e{0}%
\e{1}%
\e{0}%
\e{0}%
\e{0}%
\e{1}%
\e{0}%
\e{0}%
\e{0}%
\eol}\vss}\rg%
%
%
\rx{\vss\hfull{%
\rlx{\hss{$315_a$}}\cg%
\e{0}%
\e{0}%
\e{0}%
\e{0}%
\e{1}%
\e{1}%
\e{1}%
\e{0}%
\e{1}%
\e{0}%
\e{0}%
\e{0}%
\e{0}%
\e{1}%
\e{0}%
\e{1}%
\e{1}%
\e{0}%
\eol}\vss}\rg%
%
%
\rx{\vss\hfull{%
\rlx{\hss{$405_a$}}\cg%
\e{1}%
\e{2}%
\e{1}%
\e{0}%
\e{1}%
\e{1}%
\e{0}%
\e{0}%
\e{0}%
\e{1}%
\e{0}%
\e{1}%
\e{0}%
\e{1}%
\e{0}%
\e{0}%
\e{0}%
\e{0}%
\eol}\vss}\rg%
%
%
\rx{\vss\hfull{%
\rlx{\hss{$168_a$}}\cg%
\e{0}%
\e{1}%
\e{0}%
\e{1}%
\e{0}%
\e{0}%
\e{0}%
\e{0}%
\e{0}%
\e{0}%
\e{0}%
\e{1}%
\e{0}%
\e{0}%
\e{0}%
\e{0}%
\e{0}%
\e{0}%
\eol}\vss}\rg%
%
%
\rx{\vss\hfull{%
\rlx{\hss{$56_a$}}\cg%
\e{0}%
\e{0}%
\e{0}%
\e{0}%
\e{0}%
\e{0}%
\e{0}%
\e{0}%
\e{1}%
\e{0}%
\e{0}%
\e{0}%
\e{0}%
\e{0}%
\e{0}%
\e{0}%
\e{1}%
\e{0}%
\eol}\vss}\rg%
%
%
\rx{\vss\hfull{%
\rlx{\hss{$120_a$}}\cg%
\e{1}%
\e{1}%
\e{0}%
\e{0}%
\e{0}%
\e{0}%
\e{0}%
\e{0}%
\e{0}%
\e{1}%
\e{0}%
\e{0}%
\e{0}%
\e{0}%
\e{0}%
\e{0}%
\e{0}%
\e{0}%
\eol}\vss}\rg%
%
%
\rx{\vss\hfull{%
\rlx{\hss{$210_a$}}\cg%
\e{1}%
\e{1}%
\e{1}%
\e{0}%
\e{0}%
\e{0}%
\e{0}%
\e{0}%
\e{0}%
\e{0}%
\e{1}%
\e{1}%
\e{0}%
\e{0}%
\e{0}%
\e{0}%
\e{0}%
\e{0}%
\eol}\vss}\rg%
%
%
\rx{\vss\hfull{%
\rlx{\hss{$280_a$}}\cg%
\e{0}%
\e{0}%
\e{0}%
\e{0}%
\e{0}%
\e{1}%
\e{1}%
\e{1}%
\e{1}%
\e{0}%
\e{0}%
\e{0}%
\e{1}%
\e{1}%
\e{0}%
\e{1}%
\e{0}%
\e{0}%
\eol}\vss}\rg%
%
%
\rx{\vss\hfull{%
\rlx{\hss{$336_a$}}\cg%
\e{0}%
\e{0}%
\e{1}%
\e{0}%
\e{0}%
\e{1}%
\e{1}%
\e{1}%
\e{0}%
\e{0}%
\e{0}%
\e{1}%
\e{1}%
\e{0}%
\e{1}%
\e{1}%
\e{0}%
\e{0}%
\eol}\vss}\rg%
%
%
\rx{\vss\hfull{%
\rlx{\hss{$216_a$}}\cg%
\e{0}%
\e{0}%
\e{0}%
\e{1}%
\e{0}%
\e{0}%
\e{1}%
\e{0}%
\e{0}%
\e{0}%
\e{0}%
\e{0}%
\e{0}%
\e{0}%
\e{1}%
\e{1}%
\e{0}%
\e{1}%
\eol}\vss}\rg%
%
%
\rx{\vss\hfull{%
\rlx{\hss{$512_a$}}\cg%
\e{0}%
\e{1}%
\e{1}%
\e{1}%
\e{1}%
\e{1}%
\e{1}%
\e{0}%
\e{0}%
\e{0}%
\e{0}%
\e{1}%
\e{0}%
\e{1}%
\e{1}%
\e{1}%
\e{0}%
\e{0}%
\eol}\vss}\rg%
%
%
\rx{\vss\hfull{%
\rlx{\hss{$378_a$}}\cg%
\e{0}%
\e{0}%
\e{1}%
\e{0}%
\e{1}%
\e{1}%
\e{1}%
\e{0}%
\e{0}%
\e{0}%
\e{0}%
\e{0}%
\e{0}%
\e{1}%
\e{1}%
\e{1}%
\e{0}%
\e{0}%
\eol}\vss}\rg%
%
%
\rx{\vss\hfull{%
\rlx{\hss{$84_a$}}\cg%
\e{0}%
\e{0}%
\e{0}%
\e{1}%
\e{0}%
\e{0}%
\e{0}%
\e{0}%
\e{0}%
\e{0}%
\e{0}%
\e{0}%
\e{0}%
\e{0}%
\e{1}%
\e{0}%
\e{0}%
\e{0}%
\eol}\vss}\rg%
%
%
\rx{\vss\hfull{%
\rlx{\hss{$420_a$}}\cg%
\e{0}%
\e{1}%
\e{1}%
\e{1}%
\e{0}%
\e{1}%
\e{1}%
\e{0}%
\e{0}%
\e{0}%
\e{0}%
\e{1}%
\e{1}%
\e{1}%
\e{1}%
\e{0}%
\e{0}%
\e{0}%
\eol}\vss}\rg%
%
%
\rx{\vss\hfull{%
\rlx{\hss{$280_b$}}\cg%
\e{0}%
\e{1}%
\e{1}%
\e{1}%
\e{0}%
\e{0}%
\e{0}%
\e{0}%
\e{0}%
\e{0}%
\e{1}%
\e{1}%
\e{0}%
\e{0}%
\e{1}%
\e{0}%
\e{0}%
\e{0}%
\eol}\vss}\rg%
%
%
\rx{\vss\hfull{%
\rlx{\hss{$210_b$}}\cg%
\e{0}%
\e{1}%
\e{0}%
\e{1}%
\e{1}%
\e{0}%
\e{0}%
\e{0}%
\e{0}%
\e{0}%
\e{0}%
\e{0}%
\e{0}%
\e{1}%
\e{0}%
\e{0}%
\e{0}%
\e{0}%
\eol}\vss}\rg%
%
%
\rx{\vss\hfull{%
\rlx{\hss{$70_a$}}\cg%
\e{0}%
\e{0}%
\e{0}%
\e{0}%
\e{1}%
\e{0}%
\e{0}%
\e{0}%
\e{0}%
\e{0}%
\e{0}%
\e{0}%
\e{0}%
\e{0}%
\e{0}%
\e{0}%
\e{1}%
\e{0}%
\eol}\vss}\rg%
\tableclose%
%
%
%
%
%
%
\eop
\eject
\tableopen{Induce/restrict matrix for $W(A_{6})\,\subset\,W(E_{7})$}%
%
%
%
%
%
%
\rowpts=18 true pt%
\colpts=18 true pt%
\rowlabpts=40 true pt%
\collabpts=40 true pt%
\clx{\vss\hfull{%
\rlx{\hss{$ $}}\cg%
\cx{\hskip 16 true pt\flip{$[{7}]$}\hss}\cg%
\cx{\hskip 16 true pt\flip{$[{6}{1}]$}\hss}\cg%
\cx{\hskip 16 true pt\flip{$[{5}{2}]$}\hss}\cg%
\cx{\hskip 16 true pt\flip{$[{5}{1^{2}}]$}\hss}\cg%
\cx{\hskip 16 true pt\flip{$[{4}{3}]$}\hss}\cg%
\cx{\hskip 16 true pt\flip{$[{4}{2}{1}]$}\hss}\cg%
\cx{\hskip 16 true pt\flip{$[{4}{1^{3}}]$}\hss}\cg%
\cx{\hskip 16 true pt\flip{$[{3^{2}}{1}]$}\hss}\cg%
\cx{\hskip 16 true pt\flip{$[{3}{2^{2}}]$}\hss}\cg%
\cx{\hskip 16 true pt\flip{$[{3}{2}{1^{2}}]$}\hss}\cg%
\cx{\hskip 16 true pt\flip{$[{3}{1^{4}}]$}\hss}\cg%
\cx{\hskip 16 true pt\flip{$[{2^{3}}{1}]$}\hss}\cg%
\cx{\hskip 16 true pt\flip{$[{2^{2}}{1^{3}}]$}\hss}\cg%
\cx{\hskip 16 true pt\flip{$[{2}{1^{5}}]$}\hss}\cg%
\cx{\hskip 16 true pt\flip{$[{1^{7}}]$}\hss}\cg%
\eol}}\rg%
%
%
\rx{\vss\hfull{%
\rlx{\hss{$1_a$}}\cg%
\e{1}%
\e{0}%
\e{0}%
\e{0}%
\e{0}%
\e{0}%
\e{0}%
\e{0}%
\e{0}%
\e{0}%
\e{0}%
\e{0}%
\e{0}%
\e{0}%
\e{0}%
\eol}\vss}\rg%
%
%
\rx{\vss\hfull{%
\rlx{\hss{$7_a$}}\cg%
\e{0}%
\e{0}%
\e{0}%
\e{0}%
\e{0}%
\e{0}%
\e{0}%
\e{0}%
\e{0}%
\e{0}%
\e{0}%
\e{0}%
\e{0}%
\e{1}%
\e{1}%
\eol}\vss}\rg%
%
%
\rx{\vss\hfull{%
\rlx{\hss{$27_a$}}\cg%
\e{1}%
\e{2}%
\e{1}%
\e{0}%
\e{0}%
\e{0}%
\e{0}%
\e{0}%
\e{0}%
\e{0}%
\e{0}%
\e{0}%
\e{0}%
\e{0}%
\e{0}%
\eol}\vss}\rg%
%
%
\rx{\vss\hfull{%
\rlx{\hss{$21_a$}}\cg%
\e{0}%
\e{1}%
\e{0}%
\e{1}%
\e{0}%
\e{0}%
\e{0}%
\e{0}%
\e{0}%
\e{0}%
\e{0}%
\e{0}%
\e{0}%
\e{0}%
\e{0}%
\eol}\vss}\rg%
%
%
\rx{\vss\hfull{%
\rlx{\hss{$35_a$}}\cg%
\e{0}%
\e{0}%
\e{0}%
\e{0}%
\e{0}%
\e{0}%
\e{1}%
\e{0}%
\e{0}%
\e{0}%
\e{1}%
\e{0}%
\e{0}%
\e{0}%
\e{0}%
\eol}\vss}\rg%
%
%
\rx{\vss\hfull{%
\rlx{\hss{$105_a$}}\cg%
\e{0}%
\e{0}%
\e{0}%
\e{0}%
\e{0}%
\e{0}%
\e{0}%
\e{0}%
\e{0}%
\e{1}%
\e{2}%
\e{0}%
\e{2}%
\e{2}%
\e{0}%
\eol}\vss}\rg%
%
%
\rx{\vss\hfull{%
\rlx{\hss{$189_a$}}\cg%
\e{0}%
\e{0}%
\e{1}%
\e{2}%
\e{0}%
\e{2}%
\e{2}%
\e{0}%
\e{0}%
\e{1}%
\e{0}%
\e{0}%
\e{0}%
\e{0}%
\e{0}%
\eol}\vss}\rg%
%
%
\rx{\vss\hfull{%
\rlx{\hss{$21_b$}}\cg%
\e{0}%
\e{0}%
\e{0}%
\e{0}%
\e{0}%
\e{0}%
\e{0}%
\e{0}%
\e{0}%
\e{0}%
\e{0}%
\e{0}%
\e{1}%
\e{1}%
\e{1}%
\eol}\vss}\rg%
%
%
\rx{\vss\hfull{%
\rlx{\hss{$35_b$}}\cg%
\e{1}%
\e{1}%
\e{1}%
\e{0}%
\e{1}%
\e{0}%
\e{0}%
\e{0}%
\e{0}%
\e{0}%
\e{0}%
\e{0}%
\e{0}%
\e{0}%
\e{0}%
\eol}\vss}\rg%
%
%
\rx{\vss\hfull{%
\rlx{\hss{$189_b$}}\cg%
\e{0}%
\e{0}%
\e{0}%
\e{0}%
\e{0}%
\e{0}%
\e{0}%
\e{0}%
\e{1}%
\e{2}%
\e{1}%
\e{2}%
\e{3}%
\e{2}%
\e{1}%
\eol}\vss}\rg%
%
%
\rx{\vss\hfull{%
\rlx{\hss{$189_c$}}\cg%
\e{0}%
\e{0}%
\e{0}%
\e{0}%
\e{0}%
\e{0}%
\e{1}%
\e{1}%
\e{0}%
\e{2}%
\e{2}%
\e{1}%
\e{2}%
\e{1}%
\e{0}%
\eol}\vss}\rg%
%
%
\rx{\vss\hfull{%
\rlx{\hss{$15_a$}}\cg%
\e{0}%
\e{0}%
\e{0}%
\e{0}%
\e{0}%
\e{0}%
\e{0}%
\e{0}%
\e{0}%
\e{0}%
\e{0}%
\e{1}%
\e{0}%
\e{0}%
\e{1}%
\eol}\vss}\rg%
%
%
\rx{\vss\hfull{%
\rlx{\hss{$105_b$}}\cg%
\e{1}%
\e{1}%
\e{1}%
\e{0}%
\e{2}%
\e{1}%
\e{0}%
\e{1}%
\e{0}%
\e{0}%
\e{0}%
\e{0}%
\e{0}%
\e{0}%
\e{0}%
\eol}\vss}\rg%
%
%
\rx{\vss\hfull{%
\rlx{\hss{$105_c$}}\cg%
\e{0}%
\e{0}%
\e{0}%
\e{1}%
\e{1}%
\e{0}%
\e{1}%
\e{1}%
\e{0}%
\e{1}%
\e{0}%
\e{0}%
\e{0}%
\e{0}%
\e{0}%
\eol}\vss}\rg%
%
%
\rx{\vss\hfull{%
\rlx{\hss{$315_a$}}\cg%
\e{0}%
\e{0}%
\e{0}%
\e{0}%
\e{0}%
\e{1}%
\e{1}%
\e{1}%
\e{2}%
\e{3}%
\e{2}%
\e{2}%
\e{2}%
\e{1}%
\e{0}%
\eol}\vss}\rg%
%
%
\rx{\vss\hfull{%
\rlx{\hss{$405_a$}}\cg%
\e{0}%
\e{1}%
\e{2}%
\e{2}%
\e{2}%
\e{4}%
\e{2}%
\e{2}%
\e{1}%
\e{2}%
\e{0}%
\e{0}%
\e{0}%
\e{0}%
\e{0}%
\eol}\vss}\rg%
%
%
\rx{\vss\hfull{%
\rlx{\hss{$168_a$}}\cg%
\e{0}%
\e{1}%
\e{3}%
\e{1}%
\e{1}%
\e{2}%
\e{0}%
\e{0}%
\e{1}%
\e{0}%
\e{0}%
\e{0}%
\e{0}%
\e{0}%
\e{0}%
\eol}\vss}\rg%
%
%
\rx{\vss\hfull{%
\rlx{\hss{$56_a$}}\cg%
\e{0}%
\e{0}%
\e{0}%
\e{0}%
\e{0}%
\e{0}%
\e{0}%
\e{0}%
\e{0}%
\e{0}%
\e{1}%
\e{1}%
\e{1}%
\e{2}%
\e{1}%
\eol}\vss}\rg%
%
%
\rx{\vss\hfull{%
\rlx{\hss{$120_a$}}\cg%
\e{1}%
\e{2}%
\e{2}%
\e{2}%
\e{1}%
\e{1}%
\e{0}%
\e{0}%
\e{0}%
\e{0}%
\e{0}%
\e{0}%
\e{0}%
\e{0}%
\e{0}%
\eol}\vss}\rg%
%
%
\rx{\vss\hfull{%
\rlx{\hss{$210_a$}}\cg%
\e{0}%
\e{2}%
\e{2}%
\e{3}%
\e{1}%
\e{2}%
\e{1}%
\e{1}%
\e{0}%
\e{0}%
\e{0}%
\e{0}%
\e{0}%
\e{0}%
\e{0}%
\eol}\vss}\rg%
%
%
\rx{\vss\hfull{%
\rlx{\hss{$280_a$}}\cg%
\e{0}%
\e{0}%
\e{0}%
\e{0}%
\e{0}%
\e{1}%
\e{2}%
\e{0}%
\e{1}%
\e{3}%
\e{3}%
\e{1}%
\e{1}%
\e{1}%
\e{0}%
\eol}\vss}\rg%
%
%
\rx{\vss\hfull{%
\rlx{\hss{$336_a$}}\cg%
\e{0}%
\e{0}%
\e{0}%
\e{1}%
\e{0}%
\e{2}%
\e{3}%
\e{1}%
\e{1}%
\e{3}%
\e{2}%
\e{0}%
\e{1}%
\e{0}%
\e{0}%
\eol}\vss}\rg%
%
%
\rx{\vss\hfull{%
\rlx{\hss{$216_a$}}\cg%
\e{0}%
\e{0}%
\e{0}%
\e{0}%
\e{0}%
\e{1}%
\e{0}%
\e{1}%
\e{2}%
\e{2}%
\e{0}%
\e{2}%
\e{1}%
\e{1}%
\e{0}%
\eol}\vss}\rg%
%
%
\rx{\vss\hfull{%
\rlx{\hss{$512_a$}}\cg%
\e{0}%
\e{0}%
\e{1}%
\e{1}%
\e{1}%
\e{4}%
\e{1}%
\e{3}%
\e{3}%
\e{4}%
\e{1}%
\e{1}%
\e{1}%
\e{0}%
\e{0}%
\eol}\vss}\rg%
%
%
\rx{\vss\hfull{%
\rlx{\hss{$378_a$}}\cg%
\e{0}%
\e{0}%
\e{0}%
\e{0}%
\e{1}%
\e{2}%
\e{1}%
\e{2}%
\e{1}%
\e{4}%
\e{1}%
\e{2}%
\e{2}%
\e{0}%
\e{0}%
\eol}\vss}\rg%
%
%
\rx{\vss\hfull{%
\rlx{\hss{$84_a$}}\cg%
\e{0}%
\e{0}%
\e{0}%
\e{0}%
\e{1}%
\e{1}%
\e{0}%
\e{0}%
\e{1}%
\e{0}%
\e{0}%
\e{1}%
\e{0}%
\e{0}%
\e{0}%
\eol}\vss}\rg%
%
%
\rx{\vss\hfull{%
\rlx{\hss{$420_a$}}\cg%
\e{0}%
\e{0}%
\e{1}%
\e{2}%
\e{1}%
\e{4}%
\e{3}%
\e{1}%
\e{2}%
\e{2}%
\e{1}%
\e{1}%
\e{0}%
\e{0}%
\e{0}%
\eol}\vss}\rg%
%
%
\rx{\vss\hfull{%
\rlx{\hss{$280_b$}}\cg%
\e{0}%
\e{1}%
\e{2}%
\e{1}%
\e{2}%
\e{3}%
\e{0}%
\e{2}%
\e{1}%
\e{1}%
\e{0}%
\e{0}%
\e{0}%
\e{0}%
\e{0}%
\eol}\vss}\rg%
%
%
\rx{\vss\hfull{%
\rlx{\hss{$210_b$}}\cg%
\e{0}%
\e{0}%
\e{1}%
\e{0}%
\e{2}%
\e{1}%
\e{0}%
\e{2}%
\e{2}%
\e{1}%
\e{0}%
\e{1}%
\e{0}%
\e{0}%
\e{0}%
\eol}\vss}\rg%
%
%
\rx{\vss\hfull{%
\rlx{\hss{$70_a$}}\cg%
\e{0}%
\e{0}%
\e{0}%
\e{0}%
\e{0}%
\e{0}%
\e{0}%
\e{1}%
\e{1}%
\e{0}%
\e{0}%
\e{1}%
\e{1}%
\e{0}%
\e{0}%
\eol}\vss}\rg%
\tableclose%
%
%
%
%
%
%
\eop
\eject
\tableopen{Induce/restrict matrix for $W({A_{5}}{A_{1}})\,\subset\,W(E_{7})$}%
%
%
%
%
%
%
\rowpts=18 true pt%
\colpts=18 true pt%
\rowlabpts=40 true pt%
\collabpts=65 true pt%
\clx{\vss\hfull{%
\rlx{\hss{$ $}}\cg%
\cx{\hskip 16 true pt\flip{$[{6}]{\times}[{2}]$}\hss}\cg%
\cx{\hskip 16 true pt\flip{$[{5}{1}]{\times}[{2}]$}\hss}\cg%
\cx{\hskip 16 true pt\flip{$[{4}{2}]{\times}[{2}]$}\hss}\cg%
\cx{\hskip 16 true pt\flip{$[{4}{1^{2}}]{\times}[{2}]$}\hss}\cg%
\cx{\hskip 16 true pt\flip{$[{3^{2}}]{\times}[{2}]$}\hss}\cg%
\cx{\hskip 16 true pt\flip{$[{3}{2}{1}]{\times}[{2}]$}\hss}\cg%
\cx{\hskip 16 true pt\flip{$[{3}{1^{3}}]{\times}[{2}]$}\hss}\cg%
\cx{\hskip 16 true pt\flip{$[{2^{3}}]{\times}[{2}]$}\hss}\cg%
\cx{\hskip 16 true pt\flip{$[{2^{2}}{1^{2}}]{\times}[{2}]$}\hss}\cg%
\cx{\hskip 16 true pt\flip{$[{2}{1^{4}}]{\times}[{2}]$}\hss}\cg%
\cx{\hskip 16 true pt\flip{$[{1^{6}}]{\times}[{2}]$}\hss}\cg%
\cx{\hskip 16 true pt\flip{$[{6}]{\times}[{1^{2}}]$}\hss}\cg%
\cx{\hskip 16 true pt\flip{$[{5}{1}]{\times}[{1^{2}}]$}\hss}\cg%
\cx{\hskip 16 true pt\flip{$[{4}{2}]{\times}[{1^{2}}]$}\hss}\cg%
\cx{\hskip 16 true pt\flip{$[{4}{1^{2}}]{\times}[{1^{2}}]$}\hss}\cg%
\cx{\hskip 16 true pt\flip{$[{3^{2}}]{\times}[{1^{2}}]$}\hss}\cg%
\cx{\hskip 16 true pt\flip{$[{3}{2}{1}]{\times}[{1^{2}}]$}\hss}\cg%
\cx{\hskip 16 true pt\flip{$[{3}{1^{3}}]{\times}[{1^{2}}]$}\hss}\cg%
\cx{\hskip 16 true pt\flip{$[{2^{3}}]{\times}[{1^{2}}]$}\hss}\cg%
\cx{\hskip 16 true pt\flip{$[{2^{2}}{1^{2}}]{\times}[{1^{2}}]$}\hss}\cg%
\cx{\hskip 16 true pt\flip{$[{2}{1^{4}}]{\times}[{1^{2}}]$}\hss}\cg%
\cx{\hskip 16 true pt\flip{$[{1^{6}}]{\times}[{1^{2}}]$}\hss}\cg%
\eol}}\rg%
%
%
\rx{\vss\hfull{%
\rlx{\hss{$1_a$}}\cg%
\e{1}%
\e{0}%
\e{0}%
\e{0}%
\e{0}%
\e{0}%
\e{0}%
\e{0}%
\e{0}%
\e{0}%
\e{0}%
\e{0}%
\e{0}%
\e{0}%
\e{0}%
\e{0}%
\e{0}%
\e{0}%
\e{0}%
\e{0}%
\e{0}%
\e{0}%
\eol}\vss}\rg%
%
%
\rx{\vss\hfull{%
\rlx{\hss{$7_a$}}\cg%
\e{0}%
\e{0}%
\e{0}%
\e{0}%
\e{0}%
\e{0}%
\e{0}%
\e{0}%
\e{0}%
\e{0}%
\e{1}%
\e{0}%
\e{0}%
\e{0}%
\e{0}%
\e{0}%
\e{0}%
\e{0}%
\e{0}%
\e{0}%
\e{1}%
\e{1}%
\eol}\vss}\rg%
%
%
\rx{\vss\hfull{%
\rlx{\hss{$27_a$}}\cg%
\e{2}%
\e{2}%
\e{1}%
\e{0}%
\e{0}%
\e{0}%
\e{0}%
\e{0}%
\e{0}%
\e{0}%
\e{0}%
\e{1}%
\e{1}%
\e{0}%
\e{0}%
\e{0}%
\e{0}%
\e{0}%
\e{0}%
\e{0}%
\e{0}%
\e{0}%
\eol}\vss}\rg%
%
%
\rx{\vss\hfull{%
\rlx{\hss{$21_a$}}\cg%
\e{0}%
\e{1}%
\e{0}%
\e{1}%
\e{0}%
\e{0}%
\e{0}%
\e{0}%
\e{0}%
\e{0}%
\e{0}%
\e{1}%
\e{1}%
\e{0}%
\e{0}%
\e{0}%
\e{0}%
\e{0}%
\e{0}%
\e{0}%
\e{0}%
\e{0}%
\eol}\vss}\rg%
%
%
\rx{\vss\hfull{%
\rlx{\hss{$35_a$}}\cg%
\e{0}%
\e{0}%
\e{0}%
\e{0}%
\e{0}%
\e{0}%
\e{1}%
\e{0}%
\e{0}%
\e{1}%
\e{0}%
\e{0}%
\e{0}%
\e{0}%
\e{1}%
\e{0}%
\e{0}%
\e{1}%
\e{0}%
\e{0}%
\e{0}%
\e{0}%
\eol}\vss}\rg%
%
%
\rx{\vss\hfull{%
\rlx{\hss{$105_a$}}\cg%
\e{0}%
\e{0}%
\e{0}%
\e{0}%
\e{0}%
\e{0}%
\e{1}%
\e{0}%
\e{1}%
\e{3}%
\e{1}%
\e{0}%
\e{0}%
\e{0}%
\e{0}%
\e{0}%
\e{1}%
\e{2}%
\e{0}%
\e{2}%
\e{3}%
\e{1}%
\eol}\vss}\rg%
%
%
\rx{\vss\hfull{%
\rlx{\hss{$189_a$}}\cg%
\e{0}%
\e{1}%
\e{1}%
\e{3}%
\e{0}%
\e{2}%
\e{2}%
\e{0}%
\e{1}%
\e{0}%
\e{0}%
\e{0}%
\e{2}%
\e{2}%
\e{3}%
\e{0}%
\e{1}%
\e{1}%
\e{0}%
\e{0}%
\e{0}%
\e{0}%
\eol}\vss}\rg%
%
%
\rx{\vss\hfull{%
\rlx{\hss{$21_b$}}\cg%
\e{0}%
\e{0}%
\e{0}%
\e{0}%
\e{0}%
\e{0}%
\e{0}%
\e{0}%
\e{0}%
\e{1}%
\e{0}%
\e{0}%
\e{0}%
\e{0}%
\e{0}%
\e{0}%
\e{0}%
\e{0}%
\e{1}%
\e{0}%
\e{2}%
\e{1}%
\eol}\vss}\rg%
%
%
\rx{\vss\hfull{%
\rlx{\hss{$35_b$}}\cg%
\e{2}%
\e{1}%
\e{2}%
\e{0}%
\e{0}%
\e{0}%
\e{0}%
\e{0}%
\e{0}%
\e{0}%
\e{0}%
\e{1}%
\e{0}%
\e{1}%
\e{0}%
\e{0}%
\e{0}%
\e{0}%
\e{0}%
\e{0}%
\e{0}%
\e{0}%
\eol}\vss}\rg%
%
%
\rx{\vss\hfull{%
\rlx{\hss{$189_b$}}\cg%
\e{0}%
\e{0}%
\e{0}%
\e{0}%
\e{0}%
\e{1}%
\e{1}%
\e{1}%
\e{3}%
\e{2}%
\e{1}%
\e{0}%
\e{0}%
\e{0}%
\e{0}%
\e{0}%
\e{2}%
\e{2}%
\e{2}%
\e{4}%
\e{4}%
\e{2}%
\eol}\vss}\rg%
%
%
\rx{\vss\hfull{%
\rlx{\hss{$189_c$}}\cg%
\e{0}%
\e{0}%
\e{0}%
\e{0}%
\e{0}%
\e{1}%
\e{3}%
\e{1}%
\e{1}%
\e{3}%
\e{0}%
\e{0}%
\e{0}%
\e{1}%
\e{0}%
\e{0}%
\e{2}%
\e{3}%
\e{2}%
\e{2}%
\e{3}%
\e{0}%
\eol}\vss}\rg%
%
%
\rx{\vss\hfull{%
\rlx{\hss{$15_a$}}\cg%
\e{0}%
\e{0}%
\e{0}%
\e{0}%
\e{0}%
\e{0}%
\e{0}%
\e{1}%
\e{0}%
\e{0}%
\e{0}%
\e{0}%
\e{0}%
\e{0}%
\e{0}%
\e{0}%
\e{0}%
\e{0}%
\e{1}%
\e{0}%
\e{1}%
\e{0}%
\eol}\vss}\rg%
%
%
\rx{\vss\hfull{%
\rlx{\hss{$105_b$}}\cg%
\e{1}%
\e{2}%
\e{2}%
\e{1}%
\e{2}%
\e{1}%
\e{0}%
\e{0}%
\e{0}%
\e{0}%
\e{0}%
\e{0}%
\e{1}%
\e{1}%
\e{0}%
\e{2}%
\e{1}%
\e{0}%
\e{0}%
\e{0}%
\e{0}%
\e{0}%
\eol}\vss}\rg%
%
%
\rx{\vss\hfull{%
\rlx{\hss{$105_c$}}\cg%
\e{0}%
\e{0}%
\e{1}%
\e{1}%
\e{0}%
\e{1}%
\e{2}%
\e{0}%
\e{0}%
\e{0}%
\e{0}%
\e{1}%
\e{0}%
\e{2}%
\e{0}%
\e{0}%
\e{1}%
\e{1}%
\e{1}%
\e{0}%
\e{0}%
\e{0}%
\eol}\vss}\rg%
%
%
\rx{\vss\hfull{%
\rlx{\hss{$315_a$}}\cg%
\e{0}%
\e{0}%
\e{0}%
\e{1}%
\e{1}%
\e{3}%
\e{2}%
\e{1}%
\e{4}%
\e{2}%
\e{1}%
\e{0}%
\e{0}%
\e{0}%
\e{2}%
\e{1}%
\e{4}%
\e{3}%
\e{1}%
\e{5}%
\e{2}%
\e{1}%
\eol}\vss}\rg%
%
%
\rx{\vss\hfull{%
\rlx{\hss{$405_a$}}\cg%
\e{0}%
\e{3}%
\e{3}%
\e{5}%
\e{3}%
\e{5}%
\e{2}%
\e{0}%
\e{2}%
\e{0}%
\e{0}%
\e{0}%
\e{3}%
\e{3}%
\e{4}%
\e{3}%
\e{4}%
\e{1}%
\e{0}%
\e{1}%
\e{0}%
\e{0}%
\eol}\vss}\rg%
%
%
\rx{\vss\hfull{%
\rlx{\hss{$168_a$}}\cg%
\e{1}%
\e{3}%
\e{4}%
\e{1}%
\e{1}%
\e{2}%
\e{0}%
\e{1}%
\e{0}%
\e{0}%
\e{0}%
\e{0}%
\e{2}%
\e{2}%
\e{2}%
\e{0}%
\e{1}%
\e{0}%
\e{0}%
\e{0}%
\e{0}%
\e{0}%
\eol}\vss}\rg%
%
%
\rx{\vss\hfull{%
\rlx{\hss{$56_a$}}\cg%
\e{0}%
\e{0}%
\e{0}%
\e{0}%
\e{0}%
\e{0}%
\e{0}%
\e{0}%
\e{1}%
\e{1}%
\e{2}%
\e{0}%
\e{0}%
\e{0}%
\e{0}%
\e{0}%
\e{0}%
\e{1}%
\e{0}%
\e{2}%
\e{2}%
\e{2}%
\eol}\vss}\rg%
%
%
\rx{\vss\hfull{%
\rlx{\hss{$120_a$}}\cg%
\e{1}%
\e{4}%
\e{2}%
\e{2}%
\e{1}%
\e{1}%
\e{0}%
\e{0}%
\e{0}%
\e{0}%
\e{0}%
\e{1}%
\e{3}%
\e{1}%
\e{1}%
\e{1}%
\e{0}%
\e{0}%
\e{0}%
\e{0}%
\e{0}%
\e{0}%
\eol}\vss}\rg%
%
%
\rx{\vss\hfull{%
\rlx{\hss{$210_a$}}\cg%
\e{1}%
\e{3}%
\e{3}%
\e{4}%
\e{1}%
\e{2}%
\e{1}%
\e{0}%
\e{0}%
\e{0}%
\e{0}%
\e{2}%
\e{3}%
\e{3}%
\e{2}%
\e{0}%
\e{1}%
\e{0}%
\e{0}%
\e{0}%
\e{0}%
\e{0}%
\eol}\vss}\rg%
%
%
\rx{\vss\hfull{%
\rlx{\hss{$280_a$}}\cg%
\e{0}%
\e{0}%
\e{0}%
\e{1}%
\e{0}%
\e{2}%
\e{3}%
\e{1}%
\e{3}%
\e{3}%
\e{1}%
\e{0}%
\e{0}%
\e{0}%
\e{3}%
\e{1}%
\e{3}%
\e{4}%
\e{0}%
\e{3}%
\e{2}%
\e{0}%
\eol}\vss}\rg%
%
%
\rx{\vss\hfull{%
\rlx{\hss{$336_a$}}\cg%
\e{0}%
\e{0}%
\e{1}%
\e{2}%
\e{0}%
\e{3}%
\e{5}%
\e{1}%
\e{2}%
\e{2}%
\e{0}%
\e{0}%
\e{1}%
\e{2}%
\e{3}%
\e{0}%
\e{4}%
\e{4}%
\e{1}%
\e{1}%
\e{1}%
\e{0}%
\eol}\vss}\rg%
%
%
\rx{\vss\hfull{%
\rlx{\hss{$216_a$}}\cg%
\e{0}%
\e{0}%
\e{1}%
\e{0}%
\e{0}%
\e{3}%
\e{1}%
\e{3}%
\e{1}%
\e{1}%
\e{0}%
\e{0}%
\e{0}%
\e{1}%
\e{0}%
\e{0}%
\e{3}%
\e{2}%
\e{3}%
\e{2}%
\e{2}%
\e{0}%
\eol}\vss}\rg%
%
%
\rx{\vss\hfull{%
\rlx{\hss{$512_a$}}\cg%
\e{0}%
\e{1}%
\e{3}%
\e{3}%
\e{2}%
\e{7}%
\e{3}%
\e{2}%
\e{3}%
\e{1}%
\e{0}%
\e{0}%
\e{1}%
\e{3}%
\e{3}%
\e{2}%
\e{7}%
\e{3}%
\e{2}%
\e{3}%
\e{1}%
\e{0}%
\eol}\vss}\rg%
%
%
\rx{\vss\hfull{%
\rlx{\hss{$378_a$}}\cg%
\e{0}%
\e{0}%
\e{1}%
\e{2}%
\e{1}%
\e{4}%
\e{3}%
\e{1}%
\e{4}%
\e{1}%
\e{0}%
\e{0}%
\e{0}%
\e{2}%
\e{1}%
\e{2}%
\e{5}%
\e{3}%
\e{2}%
\e{4}%
\e{2}%
\e{0}%
\eol}\vss}\rg%
%
%
\rx{\vss\hfull{%
\rlx{\hss{$84_a$}}\cg%
\e{0}%
\e{0}%
\e{2}%
\e{0}%
\e{0}%
\e{1}%
\e{0}%
\e{2}%
\e{0}%
\e{0}%
\e{0}%
\e{0}%
\e{0}%
\e{1}%
\e{0}%
\e{0}%
\e{1}%
\e{1}%
\e{1}%
\e{0}%
\e{0}%
\e{0}%
\eol}\vss}\rg%
%
%
\rx{\vss\hfull{%
\rlx{\hss{$420_a$}}\cg%
\e{0}%
\e{1}%
\e{3}%
\e{4}%
\e{1}%
\e{5}%
\e{3}%
\e{2}%
\e{2}%
\e{1}%
\e{0}%
\e{0}%
\e{2}%
\e{3}%
\e{5}%
\e{1}%
\e{4}%
\e{3}%
\e{1}%
\e{1}%
\e{0}%
\e{0}%
\eol}\vss}\rg%
%
%
\rx{\vss\hfull{%
\rlx{\hss{$280_b$}}\cg%
\e{1}%
\e{2}%
\e{5}%
\e{2}%
\e{1}%
\e{4}%
\e{1}%
\e{1}%
\e{0}%
\e{0}%
\e{0}%
\e{1}%
\e{1}%
\e{4}%
\e{1}%
\e{1}%
\e{3}%
\e{1}%
\e{1}%
\e{0}%
\e{0}%
\e{0}%
\eol}\vss}\rg%
%
%
\rx{\vss\hfull{%
\rlx{\hss{$210_b$}}\cg%
\e{0}%
\e{1}%
\e{2}%
\e{1}%
\e{3}%
\e{3}%
\e{0}%
\e{1}%
\e{1}%
\e{0}%
\e{0}%
\e{0}%
\e{0}%
\e{1}%
\e{1}%
\e{2}%
\e{3}%
\e{0}%
\e{1}%
\e{2}%
\e{0}%
\e{0}%
\eol}\vss}\rg%
%
%
\rx{\vss\hfull{%
\rlx{\hss{$70_a$}}\cg%
\e{0}%
\e{0}%
\e{0}%
\e{0}%
\e{1}%
\e{1}%
\e{0}%
\e{0}%
\e{1}%
\e{0}%
\e{0}%
\e{0}%
\e{0}%
\e{0}%
\e{0}%
\e{0}%
\e{1}%
\e{0}%
\e{1}%
\e{2}%
\e{0}%
\e{1}%
\eol}\vss}\rg%
\tableclose%
%
%
%
%
%
%
\eop
\eject
\tableopen{Induce/restrict matrix for $W({A_{3}}{A_{2}}{A_{1}})\,\subset\,W(E_{7})$}%
%
%
%
%
%
%
\rowpts=18 true pt%
\colpts=18 true pt%
\rowlabpts=40 true pt%
\collabpts=90 true pt%
\clx{\vss\hfull{%
\rlx{\hss{$ $}}\cg%
\cx{\hskip 16 true pt\flip{$[{4}]{\times}[{3}]{\times}[{2}]$}\hss}\cg%
\cx{\hskip 16 true pt\flip{$[{3}{1}]{\times}[{3}]{\times}[{2}]$}\hss}\cg%
\cx{\hskip 16 true pt\flip{$[{2^{2}}]{\times}[{3}]{\times}[{2}]$}\hss}\cg%
\cx{\hskip 16 true pt\flip{$[{2}{1^{2}}]{\times}[{3}]{\times}[{2}]$}\hss}\cg%
\cx{\hskip 16 true pt\flip{$[{1^{4}}]{\times}[{3}]{\times}[{2}]$}\hss}\cg%
\cx{\hskip 16 true pt\flip{$[{4}]{\times}[{2}{1}]{\times}[{2}]$}\hss}\cg%
\cx{\hskip 16 true pt\flip{$[{3}{1}]{\times}[{2}{1}]{\times}[{2}]$}\hss}\cg%
\cx{\hskip 16 true pt\flip{$[{2^{2}}]{\times}[{2}{1}]{\times}[{2}]$}\hss}\cg%
\cx{\hskip 16 true pt\flip{$[{2}{1^{2}}]{\times}[{2}{1}]{\times}[{2}]$}\hss}\cg%
\cx{\hskip 16 true pt\flip{$[{1^{4}}]{\times}[{2}{1}]{\times}[{2}]$}\hss}\cg%
\cx{\hskip 16 true pt\flip{$[{4}]{\times}[{1^{3}}]{\times}[{2}]$}\hss}\cg%
\cx{\hskip 16 true pt\flip{$[{3}{1}]{\times}[{1^{3}}]{\times}[{2}]$}\hss}\cg%
\cx{\hskip 16 true pt\flip{$[{2^{2}}]{\times}[{1^{3}}]{\times}[{2}]$}\hss}\cg%
\cx{\hskip 16 true pt\flip{$[{2}{1^{2}}]{\times}[{1^{3}}]{\times}[{2}]$}\hss}\cg%
\cx{\hskip 16 true pt\flip{$[{1^{4}}]{\times}[{1^{3}}]{\times}[{2}]$}\hss}\cg%
\eol}}\rg%
%
%
\rx{\vss\hfull{%
\rlx{\hss{$1_a$}}\cg%
\e{1}%
\e{0}%
\e{0}%
\e{0}%
\e{0}%
\e{0}%
\e{0}%
\e{0}%
\e{0}%
\e{0}%
\e{0}%
\e{0}%
\e{0}%
\e{0}%
\e{0}%
\eol}\vss}\rg%
%
%
\rx{\vss\hfull{%
\rlx{\hss{$7_a$}}\cg%
\e{0}%
\e{0}%
\e{0}%
\e{0}%
\e{0}%
\e{0}%
\e{0}%
\e{0}%
\e{0}%
\e{0}%
\e{0}%
\e{0}%
\e{0}%
\e{0}%
\e{1}%
\eol}\vss}\rg%
%
%
\rx{\vss\hfull{%
\rlx{\hss{$27_a$}}\cg%
\e{3}%
\e{2}%
\e{1}%
\e{0}%
\e{0}%
\e{2}%
\e{1}%
\e{0}%
\e{0}%
\e{0}%
\e{0}%
\e{0}%
\e{0}%
\e{0}%
\e{0}%
\eol}\vss}\rg%
%
%
\rx{\vss\hfull{%
\rlx{\hss{$21_a$}}\cg%
\e{0}%
\e{1}%
\e{0}%
\e{1}%
\e{0}%
\e{1}%
\e{1}%
\e{0}%
\e{0}%
\e{0}%
\e{1}%
\e{0}%
\e{0}%
\e{0}%
\e{0}%
\eol}\vss}\rg%
%
%
\rx{\vss\hfull{%
\rlx{\hss{$35_a$}}\cg%
\e{0}%
\e{0}%
\e{0}%
\e{0}%
\e{1}%
\e{0}%
\e{0}%
\e{0}%
\e{1}%
\e{1}%
\e{0}%
\e{1}%
\e{0}%
\e{1}%
\e{0}%
\eol}\vss}\rg%
%
%
\rx{\vss\hfull{%
\rlx{\hss{$105_a$}}\cg%
\e{0}%
\e{0}%
\e{0}%
\e{0}%
\e{1}%
\e{0}%
\e{0}%
\e{0}%
\e{2}%
\e{3}%
\e{0}%
\e{1}%
\e{1}%
\e{3}%
\e{2}%
\eol}\vss}\rg%
%
%
\rx{\vss\hfull{%
\rlx{\hss{$189_a$}}\cg%
\e{0}%
\e{2}%
\e{1}%
\e{4}%
\e{1}%
\e{2}%
\e{5}%
\e{2}%
\e{4}%
\e{1}%
\e{2}%
\e{3}%
\e{1}%
\e{1}%
\e{0}%
\eol}\vss}\rg%
%
%
\rx{\vss\hfull{%
\rlx{\hss{$21_b$}}\cg%
\e{0}%
\e{0}%
\e{0}%
\e{0}%
\e{0}%
\e{0}%
\e{0}%
\e{0}%
\e{0}%
\e{1}%
\e{0}%
\e{0}%
\e{1}%
\e{0}%
\e{1}%
\eol}\vss}\rg%
%
%
\rx{\vss\hfull{%
\rlx{\hss{$35_b$}}\cg%
\e{3}%
\e{2}%
\e{1}%
\e{0}%
\e{0}%
\e{2}%
\e{1}%
\e{1}%
\e{0}%
\e{0}%
\e{0}%
\e{0}%
\e{0}%
\e{0}%
\e{0}%
\eol}\vss}\rg%
%
%
\rx{\vss\hfull{%
\rlx{\hss{$189_b$}}\cg%
\e{0}%
\e{0}%
\e{0}%
\e{1}%
\e{0}%
\e{0}%
\e{1}%
\e{2}%
\e{4}%
\e{3}%
\e{0}%
\e{1}%
\e{2}%
\e{4}%
\e{3}%
\eol}\vss}\rg%
%
%
\rx{\vss\hfull{%
\rlx{\hss{$189_c$}}\cg%
\e{0}%
\e{0}%
\e{0}%
\e{1}%
\e{2}%
\e{0}%
\e{1}%
\e{2}%
\e{4}%
\e{4}%
\e{1}%
\e{2}%
\e{3}%
\e{3}%
\e{2}%
\eol}\vss}\rg%
%
%
\rx{\vss\hfull{%
\rlx{\hss{$15_a$}}\cg%
\e{0}%
\e{0}%
\e{0}%
\e{0}%
\e{0}%
\e{0}%
\e{0}%
\e{1}%
\e{0}%
\e{0}%
\e{0}%
\e{0}%
\e{0}%
\e{0}%
\e{1}%
\eol}\vss}\rg%
%
%
\rx{\vss\hfull{%
\rlx{\hss{$105_b$}}\cg%
\e{3}%
\e{4}%
\e{1}%
\e{1}%
\e{0}%
\e{2}%
\e{4}%
\e{2}%
\e{1}%
\e{0}%
\e{0}%
\e{1}%
\e{0}%
\e{0}%
\e{0}%
\eol}\vss}\rg%
%
%
\rx{\vss\hfull{%
\rlx{\hss{$105_c$}}\cg%
\e{0}%
\e{1}%
\e{0}%
\e{2}%
\e{1}%
\e{1}%
\e{2}%
\e{2}%
\e{2}%
\e{1}%
\e{2}%
\e{1}%
\e{2}%
\e{0}%
\e{0}%
\eol}\vss}\rg%
%
%
\rx{\vss\hfull{%
\rlx{\hss{$315_a$}}\cg%
\e{0}%
\e{1}%
\e{1}%
\e{2}%
\e{1}%
\e{0}%
\e{4}%
\e{3}%
\e{8}%
\e{3}%
\e{0}%
\e{3}%
\e{2}%
\e{6}%
\e{2}%
\eol}\vss}\rg%
%
%
\rx{\vss\hfull{%
\rlx{\hss{$405_a$}}\cg%
\e{2}%
\e{7}%
\e{3}%
\e{6}%
\e{1}%
\e{4}%
\e{12}%
\e{5}%
\e{8}%
\e{1}%
\e{2}%
\e{5}%
\e{2}%
\e{2}%
\e{0}%
\eol}\vss}\rg%
%
%
\rx{\vss\hfull{%
\rlx{\hss{$168_a$}}\cg%
\e{4}%
\e{5}%
\e{4}%
\e{1}%
\e{0}%
\e{4}%
\e{6}%
\e{3}%
\e{2}%
\e{0}%
\e{0}%
\e{1}%
\e{0}%
\e{1}%
\e{0}%
\eol}\vss}\rg%
%
%
\rx{\vss\hfull{%
\rlx{\hss{$56_a$}}\cg%
\e{0}%
\e{0}%
\e{0}%
\e{0}%
\e{0}%
\e{0}%
\e{0}%
\e{0}%
\e{1}%
\e{1}%
\e{0}%
\e{0}%
\e{0}%
\e{2}%
\e{2}%
\eol}\vss}\rg%
%
%
\rx{\vss\hfull{%
\rlx{\hss{$120_a$}}\cg%
\e{3}%
\e{5}%
\e{2}%
\e{2}%
\e{0}%
\e{4}%
\e{5}%
\e{1}%
\e{1}%
\e{0}%
\e{1}%
\e{1}%
\e{0}%
\e{0}%
\e{0}%
\eol}\vss}\rg%
%
%
\rx{\vss\hfull{%
\rlx{\hss{$210_a$}}\cg%
\e{2}%
\e{6}%
\e{2}%
\e{4}%
\e{1}%
\e{5}%
\e{7}%
\e{3}%
\e{3}%
\e{0}%
\e{3}%
\e{2}%
\e{1}%
\e{0}%
\e{0}%
\eol}\vss}\rg%
%
%
\rx{\vss\hfull{%
\rlx{\hss{$280_a$}}\cg%
\e{0}%
\e{0}%
\e{1}%
\e{2}%
\e{2}%
\e{0}%
\e{3}%
\e{2}%
\e{7}%
\e{4}%
\e{0}%
\e{4}%
\e{1}%
\e{6}%
\e{2}%
\eol}\vss}\rg%
%
%
\rx{\vss\hfull{%
\rlx{\hss{$336_a$}}\cg%
\e{0}%
\e{1}%
\e{1}%
\e{4}%
\e{3}%
\e{1}%
\e{5}%
\e{4}%
\e{8}%
\e{4}%
\e{2}%
\e{5}%
\e{3}%
\e{4}%
\e{1}%
\eol}\vss}\rg%
%
%
\rx{\vss\hfull{%
\rlx{\hss{$216_a$}}\cg%
\e{0}%
\e{1}%
\e{2}%
\e{1}%
\e{0}%
\e{1}%
\e{3}%
\e{5}%
\e{4}%
\e{2}%
\e{0}%
\e{1}%
\e{2}%
\e{3}%
\e{2}%
\eol}\vss}\rg%
%
%
\rx{\vss\hfull{%
\rlx{\hss{$512_a$}}\cg%
\e{1}%
\e{5}%
\e{4}%
\e{5}%
\e{1}%
\e{3}%
\e{11}%
\e{8}%
\e{11}%
\e{3}%
\e{1}%
\e{5}%
\e{4}%
\e{5}%
\e{1}%
\eol}\vss}\rg%
%
%
\rx{\vss\hfull{%
\rlx{\hss{$378_a$}}\cg%
\e{0}%
\e{2}%
\e{1}%
\e{4}%
\e{1}%
\e{1}%
\e{6}%
\e{5}%
\e{9}%
\e{3}%
\e{1}%
\e{4}%
\e{4}%
\e{4}%
\e{2}%
\eol}\vss}\rg%
%
%
\rx{\vss\hfull{%
\rlx{\hss{$84_a$}}\cg%
\e{1}%
\e{1}%
\e{2}%
\e{0}%
\e{0}%
\e{1}%
\e{2}%
\e{3}%
\e{1}%
\e{0}%
\e{0}%
\e{0}%
\e{0}%
\e{1}%
\e{1}%
\eol}\vss}\rg%
%
%
\rx{\vss\hfull{%
\rlx{\hss{$420_a$}}\cg%
\e{1}%
\e{4}%
\e{4}%
\e{5}%
\e{2}%
\e{3}%
\e{10}%
\e{6}%
\e{9}%
\e{2}%
\e{2}%
\e{5}%
\e{2}%
\e{4}%
\e{1}%
\eol}\vss}\rg%
%
%
\rx{\vss\hfull{%
\rlx{\hss{$280_b$}}\cg%
\e{3}%
\e{6}%
\e{4}%
\e{3}%
\e{0}%
\e{5}%
\e{8}%
\e{6}%
\e{4}%
\e{1}%
\e{1}%
\e{2}%
\e{2}%
\e{1}%
\e{0}%
\eol}\vss}\rg%
%
%
\rx{\vss\hfull{%
\rlx{\hss{$210_b$}}\cg%
\e{2}%
\e{4}%
\e{2}%
\e{1}%
\e{0}%
\e{1}%
\e{6}%
\e{4}%
\e{4}%
\e{0}%
\e{0}%
\e{1}%
\e{1}%
\e{2}%
\e{0}%
\eol}\vss}\rg%
%
%
\rx{\vss\hfull{%
\rlx{\hss{$70_a$}}\cg%
\e{0}%
\e{1}%
\e{0}%
\e{0}%
\e{0}%
\e{0}%
\e{1}%
\e{1}%
\e{2}%
\e{0}%
\e{0}%
\e{0}%
\e{1}%
\e{1}%
\e{0}%
\eol}\vss}\rg%
\eop
\eject
\tablecont%
%
%
%
%
%
%
\rowpts=18 true pt%
\colpts=18 true pt%
\rowlabpts=40 true pt%
\collabpts=90 true pt%
\clx{\vss\hfull{%
\rlx{\hss{$ $}}\cg%
\cx{\hskip 16 true pt\flip{$[{4}]{\times}[{3}]{\times}[{1^{2}}]$}\hss}\cg%
\cx{\hskip 16 true pt\flip{$[{3}{1}]{\times}[{3}]{\times}[{1^{2}}]$}\hss}\cg%
\cx{\hskip 16 true pt\flip{$[{2^{2}}]{\times}[{3}]{\times}[{1^{2}}]$}\hss}\cg%
\cx{\hskip 16 true pt\flip{$[{2}{1^{2}}]{\times}[{3}]{\times}[{1^{2}}]$}\hss}\cg%
\cx{\hskip 16 true pt\flip{$[{1^{4}}]{\times}[{3}]{\times}[{1^{2}}]$}\hss}\cg%
\cx{\hskip 16 true pt\flip{$[{4}]{\times}[{2}{1}]{\times}[{1^{2}}]$}\hss}\cg%
\cx{\hskip 16 true pt\flip{$[{3}{1}]{\times}[{2}{1}]{\times}[{1^{2}}]$}\hss}\cg%
\cx{\hskip 16 true pt\flip{$[{2^{2}}]{\times}[{2}{1}]{\times}[{1^{2}}]$}\hss}\cg%
\cx{\hskip 16 true pt\flip{$[{2}{1^{2}}]{\times}[{2}{1}]{\times}[{1^{2}}]$}\hss}\cg%
\cx{\hskip 16 true pt\flip{$[{1^{4}}]{\times}[{2}{1}]{\times}[{1^{2}}]$}\hss}\cg%
\cx{\hskip 16 true pt\flip{$[{4}]{\times}[{1^{3}}]{\times}[{1^{2}}]$}\hss}\cg%
\cx{\hskip 16 true pt\flip{$[{3}{1}]{\times}[{1^{3}}]{\times}[{1^{2}}]$}\hss}\cg%
\cx{\hskip 16 true pt\flip{$[{2^{2}}]{\times}[{1^{3}}]{\times}[{1^{2}}]$}\hss}\cg%
\cx{\hskip 16 true pt\flip{$[{2}{1^{2}}]{\times}[{1^{3}}]{\times}[{1^{2}}]$}\hss}\cg%
\cx{\hskip 16 true pt\flip{$[{1^{4}}]{\times}[{1^{3}}]{\times}[{1^{2}}]$}\hss}\cg%
\eol}}\rg%
%
%
\rx{\vss\hfull{%
\rlx{\hss{$1_a$}}\cg%
\e{0}%
\e{0}%
\e{0}%
\e{0}%
\e{0}%
\e{0}%
\e{0}%
\e{0}%
\e{0}%
\e{0}%
\e{0}%
\e{0}%
\e{0}%
\e{0}%
\e{0}%
\eol}\vss}\rg%
%
%
\rx{\vss\hfull{%
\rlx{\hss{$7_a$}}\cg%
\e{0}%
\e{0}%
\e{0}%
\e{0}%
\e{0}%
\e{0}%
\e{0}%
\e{0}%
\e{0}%
\e{1}%
\e{0}%
\e{0}%
\e{0}%
\e{1}%
\e{1}%
\eol}\vss}\rg%
%
%
\rx{\vss\hfull{%
\rlx{\hss{$27_a$}}\cg%
\e{1}%
\e{1}%
\e{0}%
\e{0}%
\e{0}%
\e{1}%
\e{0}%
\e{0}%
\e{0}%
\e{0}%
\e{0}%
\e{0}%
\e{0}%
\e{0}%
\e{0}%
\eol}\vss}\rg%
%
%
\rx{\vss\hfull{%
\rlx{\hss{$21_a$}}\cg%
\e{1}%
\e{1}%
\e{0}%
\e{0}%
\e{0}%
\e{1}%
\e{0}%
\e{0}%
\e{0}%
\e{0}%
\e{0}%
\e{0}%
\e{0}%
\e{0}%
\e{0}%
\eol}\vss}\rg%
%
%
\rx{\vss\hfull{%
\rlx{\hss{$35_a$}}\cg%
\e{0}%
\e{0}%
\e{0}%
\e{1}%
\e{1}%
\e{0}%
\e{1}%
\e{0}%
\e{1}%
\e{0}%
\e{1}%
\e{1}%
\e{0}%
\e{0}%
\e{0}%
\eol}\vss}\rg%
%
%
\rx{\vss\hfull{%
\rlx{\hss{$105_a$}}\cg%
\e{0}%
\e{0}%
\e{0}%
\e{1}%
\e{1}%
\e{0}%
\e{1}%
\e{1}%
\e{4}%
\e{4}%
\e{0}%
\e{2}%
\e{2}%
\e{4}%
\e{2}%
\eol}\vss}\rg%
%
%
\rx{\vss\hfull{%
\rlx{\hss{$189_a$}}\cg%
\e{1}%
\e{3}%
\e{2}%
\e{3}%
\e{1}%
\e{3}%
\e{5}%
\e{1}%
\e{2}%
\e{0}%
\e{2}%
\e{2}%
\e{0}%
\e{0}%
\e{0}%
\eol}\vss}\rg%
%
%
\rx{\vss\hfull{%
\rlx{\hss{$21_b$}}\cg%
\e{0}%
\e{0}%
\e{0}%
\e{0}%
\e{0}%
\e{0}%
\e{0}%
\e{0}%
\e{1}%
\e{1}%
\e{0}%
\e{0}%
\e{0}%
\e{2}%
\e{2}%
\eol}\vss}\rg%
%
%
\rx{\vss\hfull{%
\rlx{\hss{$35_b$}}\cg%
\e{1}%
\e{1}%
\e{0}%
\e{0}%
\e{0}%
\e{0}%
\e{1}%
\e{0}%
\e{0}%
\e{0}%
\e{0}%
\e{0}%
\e{0}%
\e{0}%
\e{0}%
\eol}\vss}\rg%
%
%
\rx{\vss\hfull{%
\rlx{\hss{$189_b$}}\cg%
\e{0}%
\e{0}%
\e{1}%
\e{1}%
\e{1}%
\e{0}%
\e{2}%
\e{3}%
\e{7}%
\e{5}%
\e{0}%
\e{2}%
\e{3}%
\e{7}%
\e{5}%
\eol}\vss}\rg%
%
%
\rx{\vss\hfull{%
\rlx{\hss{$189_c$}}\cg%
\e{0}%
\e{1}%
\e{0}%
\e{2}%
\e{1}%
\e{0}%
\e{3}%
\e{2}%
\e{7}%
\e{3}%
\e{0}%
\e{3}%
\e{2}%
\e{5}%
\e{2}%
\eol}\vss}\rg%
%
%
\rx{\vss\hfull{%
\rlx{\hss{$15_a$}}\cg%
\e{0}%
\e{0}%
\e{0}%
\e{0}%
\e{0}%
\e{0}%
\e{0}%
\e{0}%
\e{1}%
\e{0}%
\e{0}%
\e{0}%
\e{0}%
\e{1}%
\e{1}%
\eol}\vss}\rg%
%
%
\rx{\vss\hfull{%
\rlx{\hss{$105_b$}}\cg%
\e{2}%
\e{2}%
\e{1}%
\e{0}%
\e{0}%
\e{1}%
\e{2}%
\e{2}%
\e{1}%
\e{0}%
\e{0}%
\e{0}%
\e{1}%
\e{0}%
\e{0}%
\eol}\vss}\rg%
%
%
\rx{\vss\hfull{%
\rlx{\hss{$105_c$}}\cg%
\e{1}%
\e{2}%
\e{0}%
\e{1}%
\e{0}%
\e{0}%
\e{3}%
\e{1}%
\e{2}%
\e{0}%
\e{0}%
\e{1}%
\e{0}%
\e{1}%
\e{0}%
\eol}\vss}\rg%
%
%
\rx{\vss\hfull{%
\rlx{\hss{$315_a$}}\cg%
\e{0}%
\e{1}%
\e{2}%
\e{3}%
\e{2}%
\e{1}%
\e{5}%
\e{6}%
\e{9}%
\e{5}%
\e{1}%
\e{4}%
\e{4}%
\e{6}%
\e{3}%
\eol}\vss}\rg%
%
%
\rx{\vss\hfull{%
\rlx{\hss{$405_a$}}\cg%
\e{3}%
\e{6}%
\e{4}%
\e{4}%
\e{1}%
\e{5}%
\e{9}%
\e{6}%
\e{5}%
\e{1}%
\e{2}%
\e{3}%
\e{2}%
\e{1}%
\e{0}%
\eol}\vss}\rg%
%
%
\rx{\vss\hfull{%
\rlx{\hss{$168_a$}}\cg%
\e{1}%
\e{3}%
\e{2}%
\e{2}%
\e{0}%
\e{3}%
\e{4}%
\e{1}%
\e{1}%
\e{0}%
\e{1}%
\e{1}%
\e{0}%
\e{0}%
\e{0}%
\eol}\vss}\rg%
%
%
\rx{\vss\hfull{%
\rlx{\hss{$56_a$}}\cg%
\e{0}%
\e{0}%
\e{0}%
\e{0}%
\e{1}%
\e{0}%
\e{0}%
\e{1}%
\e{2}%
\e{3}%
\e{0}%
\e{1}%
\e{1}%
\e{3}%
\e{3}%
\eol}\vss}\rg%
%
%
\rx{\vss\hfull{%
\rlx{\hss{$120_a$}}\cg%
\e{3}%
\e{3}%
\e{1}%
\e{1}%
\e{0}%
\e{3}%
\e{2}%
\e{1}%
\e{0}%
\e{0}%
\e{1}%
\e{0}%
\e{0}%
\e{0}%
\e{0}%
\eol}\vss}\rg%
%
%
\rx{\vss\hfull{%
\rlx{\hss{$210_a$}}\cg%
\e{3}%
\e{5}%
\e{2}%
\e{2}%
\e{0}%
\e{4}%
\e{5}%
\e{1}%
\e{1}%
\e{0}%
\e{1}%
\e{1}%
\e{0}%
\e{0}%
\e{0}%
\eol}\vss}\rg%
%
%
\rx{\vss\hfull{%
\rlx{\hss{$280_a$}}\cg%
\e{0}%
\e{1}%
\e{1}%
\e{4}%
\e{3}%
\e{1}%
\e{5}%
\e{4}%
\e{8}%
\e{4}%
\e{2}%
\e{5}%
\e{3}%
\e{4}%
\e{1}%
\eol}\vss}\rg%
%
%
\rx{\vss\hfull{%
\rlx{\hss{$336_a$}}\cg%
\e{0}%
\e{3}%
\e{2}%
\e{5}%
\e{2}%
\e{2}%
\e{8}%
\e{3}%
\e{8}%
\e{2}%
\e{2}%
\e{5}%
\e{2}%
\e{3}%
\e{0}%
\eol}\vss}\rg%
%
%
\rx{\vss\hfull{%
\rlx{\hss{$216_a$}}\cg%
\e{0}%
\e{1}%
\e{1}%
\e{2}%
\e{0}%
\e{0}%
\e{4}%
\e{3}%
\e{7}%
\e{2}%
\e{0}%
\e{2}%
\e{2}%
\e{5}%
\e{2}%
\eol}\vss}\rg%
%
%
\rx{\vss\hfull{%
\rlx{\hss{$512_a$}}\cg%
\e{1}%
\e{5}%
\e{4}%
\e{5}%
\e{1}%
\e{3}%
\e{11}%
\e{8}%
\e{11}%
\e{3}%
\e{1}%
\e{5}%
\e{4}%
\e{5}%
\e{1}%
\eol}\vss}\rg%
%
%
\rx{\vss\hfull{%
\rlx{\hss{$378_a$}}\cg%
\e{1}%
\e{3}%
\e{2}%
\e{3}%
\e{1}%
\e{1}%
\e{7}%
\e{7}%
\e{10}%
\e{4}%
\e{0}%
\e{4}%
\e{4}%
\e{6}%
\e{2}%
\eol}\vss}\rg%
%
%
\rx{\vss\hfull{%
\rlx{\hss{$84_a$}}\cg%
\e{0}%
\e{1}%
\e{0}%
\e{1}%
\e{0}%
\e{0}%
\e{2}%
\e{1}%
\e{2}%
\e{0}%
\e{0}%
\e{1}%
\e{0}%
\e{1}%
\e{0}%
\eol}\vss}\rg%
%
%
\rx{\vss\hfull{%
\rlx{\hss{$420_a$}}\cg%
\e{1}%
\e{5}%
\e{3}%
\e{6}%
\e{2}%
\e{4}%
\e{10}%
\e{5}%
\e{7}%
\e{1}%
\e{3}%
\e{5}%
\e{1}%
\e{2}%
\e{0}%
\eol}\vss}\rg%
%
%
\rx{\vss\hfull{%
\rlx{\hss{$280_b$}}\cg%
\e{2}%
\e{5}%
\e{2}%
\e{2}%
\e{0}%
\e{2}%
\e{7}%
\e{3}%
\e{4}%
\e{0}%
\e{0}%
\e{2}%
\e{1}%
\e{1}%
\e{0}%
\eol}\vss}\rg%
%
%
\rx{\vss\hfull{%
\rlx{\hss{$210_b$}}\cg%
\e{1}%
\e{2}%
\e{2}%
\e{1}%
\e{0}%
\e{1}%
\e{4}%
\e{5}%
\e{4}%
\e{1}%
\e{0}%
\e{1}%
\e{2}%
\e{2}%
\e{1}%
\eol}\vss}\rg%
%
%
\rx{\vss\hfull{%
\rlx{\hss{$70_a$}}\cg%
\e{0}%
\e{0}%
\e{1}%
\e{0}%
\e{0}%
\e{0}%
\e{1}%
\e{2}%
\e{2}%
\e{1}%
\e{0}%
\e{0}%
\e{1}%
\e{2}%
\e{2}%
\eol}\vss}\rg%
\tableclose%
%
%
%
%
%
%
\eop
\eject
\tableopen{Induce/restrict matrix for $W({A_{4}}{A_{2}})\,\subset\,W(E_{7})$}%
%
%
%
%
%
%
\rowpts=18 true pt%
\colpts=18 true pt%
\rowlabpts=40 true pt%
\collabpts=65 true pt%
\clx{\vss\hfull{%
\rlx{\hss{$ $}}\cg%
\cx{\hskip 16 true pt\flip{$[{5}]{\times}[{3}]$}\hss}\cg%
\cx{\hskip 16 true pt\flip{$[{4}{1}]{\times}[{3}]$}\hss}\cg%
\cx{\hskip 16 true pt\flip{$[{3}{2}]{\times}[{3}]$}\hss}\cg%
\cx{\hskip 16 true pt\flip{$[{3}{1^{2}}]{\times}[{3}]$}\hss}\cg%
\cx{\hskip 16 true pt\flip{$[{2^{2}}{1}]{\times}[{3}]$}\hss}\cg%
\cx{\hskip 16 true pt\flip{$[{2}{1^{3}}]{\times}[{3}]$}\hss}\cg%
\cx{\hskip 16 true pt\flip{$[{1^{5}}]{\times}[{3}]$}\hss}\cg%
\cx{\hskip 16 true pt\flip{$[{5}]{\times}[{2}{1}]$}\hss}\cg%
\cx{\hskip 16 true pt\flip{$[{4}{1}]{\times}[{2}{1}]$}\hss}\cg%
\cx{\hskip 16 true pt\flip{$[{3}{2}]{\times}[{2}{1}]$}\hss}\cg%
\cx{\hskip 16 true pt\flip{$[{3}{1^{2}}]{\times}[{2}{1}]$}\hss}\cg%
\cx{\hskip 16 true pt\flip{$[{2^{2}}{1}]{\times}[{2}{1}]$}\hss}\cg%
\cx{\hskip 16 true pt\flip{$[{2}{1^{3}}]{\times}[{2}{1}]$}\hss}\cg%
\cx{\hskip 16 true pt\flip{$[{1^{5}}]{\times}[{2}{1}]$}\hss}\cg%
\cx{\hskip 16 true pt\flip{$[{5}]{\times}[{1^{3}}]$}\hss}\cg%
\cx{\hskip 16 true pt\flip{$[{4}{1}]{\times}[{1^{3}}]$}\hss}\cg%
\cx{\hskip 16 true pt\flip{$[{3}{2}]{\times}[{1^{3}}]$}\hss}\cg%
\cx{\hskip 16 true pt\flip{$[{3}{1^{2}}]{\times}[{1^{3}}]$}\hss}\cg%
\cx{\hskip 16 true pt\flip{$[{2^{2}}{1}]{\times}[{1^{3}}]$}\hss}\cg%
\cx{\hskip 16 true pt\flip{$[{2}{1^{3}}]{\times}[{1^{3}}]$}\hss}\cg%
\cx{\hskip 16 true pt\flip{$[{1^{5}}]{\times}[{1^{3}}]$}\hss}\cg%
\eol}}\rg%
%
%
\rx{\vss\hfull{%
\rlx{\hss{$1_a$}}\cg%
\e{1}%
\e{0}%
\e{0}%
\e{0}%
\e{0}%
\e{0}%
\e{0}%
\e{0}%
\e{0}%
\e{0}%
\e{0}%
\e{0}%
\e{0}%
\e{0}%
\e{0}%
\e{0}%
\e{0}%
\e{0}%
\e{0}%
\e{0}%
\e{0}%
\eol}\vss}\rg%
%
%
\rx{\vss\hfull{%
\rlx{\hss{$7_a$}}\cg%
\e{0}%
\e{0}%
\e{0}%
\e{0}%
\e{0}%
\e{0}%
\e{0}%
\e{0}%
\e{0}%
\e{0}%
\e{0}%
\e{0}%
\e{0}%
\e{1}%
\e{0}%
\e{0}%
\e{0}%
\e{0}%
\e{0}%
\e{1}%
\e{1}%
\eol}\vss}\rg%
%
%
\rx{\vss\hfull{%
\rlx{\hss{$27_a$}}\cg%
\e{2}%
\e{2}%
\e{1}%
\e{0}%
\e{0}%
\e{0}%
\e{0}%
\e{2}%
\e{1}%
\e{0}%
\e{0}%
\e{0}%
\e{0}%
\e{0}%
\e{0}%
\e{0}%
\e{0}%
\e{0}%
\e{0}%
\e{0}%
\e{0}%
\eol}\vss}\rg%
%
%
\rx{\vss\hfull{%
\rlx{\hss{$21_a$}}\cg%
\e{0}%
\e{1}%
\e{0}%
\e{1}%
\e{0}%
\e{0}%
\e{0}%
\e{1}%
\e{1}%
\e{0}%
\e{0}%
\e{0}%
\e{0}%
\e{0}%
\e{1}%
\e{0}%
\e{0}%
\e{0}%
\e{0}%
\e{0}%
\e{0}%
\eol}\vss}\rg%
%
%
\rx{\vss\hfull{%
\rlx{\hss{$35_a$}}\cg%
\e{0}%
\e{0}%
\e{0}%
\e{0}%
\e{0}%
\e{1}%
\e{1}%
\e{0}%
\e{0}%
\e{0}%
\e{1}%
\e{0}%
\e{1}%
\e{0}%
\e{0}%
\e{1}%
\e{0}%
\e{1}%
\e{0}%
\e{0}%
\e{0}%
\eol}\vss}\rg%
%
%
\rx{\vss\hfull{%
\rlx{\hss{$105_a$}}\cg%
\e{0}%
\e{0}%
\e{0}%
\e{0}%
\e{0}%
\e{1}%
\e{1}%
\e{0}%
\e{0}%
\e{0}%
\e{1}%
\e{1}%
\e{4}%
\e{3}%
\e{0}%
\e{0}%
\e{1}%
\e{2}%
\e{2}%
\e{3}%
\e{1}%
\eol}\vss}\rg%
%
%
\rx{\vss\hfull{%
\rlx{\hss{$189_a$}}\cg%
\e{0}%
\e{1}%
\e{1}%
\e{3}%
\e{2}%
\e{2}%
\e{0}%
\e{1}%
\e{4}%
\e{2}%
\e{4}%
\e{1}%
\e{1}%
\e{0}%
\e{1}%
\e{3}%
\e{1}%
\e{1}%
\e{0}%
\e{0}%
\e{0}%
\eol}\vss}\rg%
%
%
\rx{\vss\hfull{%
\rlx{\hss{$21_b$}}\cg%
\e{0}%
\e{0}%
\e{0}%
\e{0}%
\e{0}%
\e{0}%
\e{0}%
\e{0}%
\e{0}%
\e{0}%
\e{0}%
\e{0}%
\e{1}%
\e{1}%
\e{0}%
\e{0}%
\e{0}%
\e{0}%
\e{1}%
\e{1}%
\e{2}%
\eol}\vss}\rg%
%
%
\rx{\vss\hfull{%
\rlx{\hss{$35_b$}}\cg%
\e{2}%
\e{2}%
\e{1}%
\e{0}%
\e{0}%
\e{0}%
\e{0}%
\e{1}%
\e{1}%
\e{1}%
\e{0}%
\e{0}%
\e{0}%
\e{0}%
\e{0}%
\e{0}%
\e{0}%
\e{0}%
\e{0}%
\e{0}%
\e{0}%
\eol}\vss}\rg%
%
%
\rx{\vss\hfull{%
\rlx{\hss{$189_b$}}\cg%
\e{0}%
\e{0}%
\e{0}%
\e{0}%
\e{1}%
\e{1}%
\e{0}%
\e{0}%
\e{0}%
\e{1}%
\e{2}%
\e{4}%
\e{5}%
\e{3}%
\e{0}%
\e{0}%
\e{1}%
\e{2}%
\e{4}%
\e{5}%
\e{3}%
\eol}\vss}\rg%
%
%
\rx{\vss\hfull{%
\rlx{\hss{$189_c$}}\cg%
\e{0}%
\e{0}%
\e{0}%
\e{1}%
\e{0}%
\e{2}%
\e{1}%
\e{0}%
\e{0}%
\e{1}%
\e{3}%
\e{3}%
\e{5}%
\e{2}%
\e{0}%
\e{1}%
\e{2}%
\e{2}%
\e{3}%
\e{3}%
\e{1}%
\eol}\vss}\rg%
%
%
\rx{\vss\hfull{%
\rlx{\hss{$15_a$}}\cg%
\e{0}%
\e{0}%
\e{0}%
\e{0}%
\e{0}%
\e{0}%
\e{0}%
\e{0}%
\e{0}%
\e{0}%
\e{0}%
\e{1}%
\e{0}%
\e{0}%
\e{0}%
\e{0}%
\e{0}%
\e{0}%
\e{0}%
\e{1}%
\e{1}%
\eol}\vss}\rg%
%
%
\rx{\vss\hfull{%
\rlx{\hss{$105_b$}}\cg%
\e{2}%
\e{3}%
\e{2}%
\e{1}%
\e{0}%
\e{0}%
\e{0}%
\e{1}%
\e{2}%
\e{3}%
\e{1}%
\e{1}%
\e{0}%
\e{0}%
\e{0}%
\e{0}%
\e{1}%
\e{0}%
\e{0}%
\e{0}%
\e{0}%
\eol}\vss}\rg%
%
%
\rx{\vss\hfull{%
\rlx{\hss{$105_c$}}\cg%
\e{0}%
\e{1}%
\e{0}%
\e{2}%
\e{0}%
\e{1}%
\e{0}%
\e{0}%
\e{1}%
\e{2}%
\e{2}%
\e{1}%
\e{1}%
\e{0}%
\e{1}%
\e{1}%
\e{1}%
\e{0}%
\e{1}%
\e{0}%
\e{0}%
\eol}\vss}\rg%
%
%
\rx{\vss\hfull{%
\rlx{\hss{$315_a$}}\cg%
\e{0}%
\e{0}%
\e{1}%
\e{1}%
\e{2}%
\e{2}%
\e{1}%
\e{0}%
\e{1}%
\e{3}%
\e{5}%
\e{6}%
\e{6}%
\e{2}%
\e{0}%
\e{1}%
\e{2}%
\e{4}%
\e{4}%
\e{4}%
\e{1}%
\eol}\vss}\rg%
%
%
\rx{\vss\hfull{%
\rlx{\hss{$405_a$}}\cg%
\e{1}%
\e{4}%
\e{4}%
\e{5}%
\e{3}%
\e{2}%
\e{0}%
\e{2}%
\e{7}%
\e{7}%
\e{7}%
\e{4}%
\e{2}%
\e{0}%
\e{1}%
\e{3}%
\e{3}%
\e{2}%
\e{1}%
\e{0}%
\e{0}%
\eol}\vss}\rg%
%
%
\rx{\vss\hfull{%
\rlx{\hss{$168_a$}}\cg%
\e{2}%
\e{3}%
\e{4}%
\e{1}%
\e{2}%
\e{0}%
\e{0}%
\e{2}%
\e{5}%
\e{3}%
\e{2}%
\e{1}%
\e{0}%
\e{0}%
\e{0}%
\e{1}%
\e{0}%
\e{1}%
\e{0}%
\e{0}%
\e{0}%
\eol}\vss}\rg%
%
%
\rx{\vss\hfull{%
\rlx{\hss{$56_a$}}\cg%
\e{0}%
\e{0}%
\e{0}%
\e{0}%
\e{0}%
\e{0}%
\e{1}%
\e{0}%
\e{0}%
\e{0}%
\e{0}%
\e{1}%
\e{2}%
\e{2}%
\e{0}%
\e{0}%
\e{0}%
\e{1}%
\e{1}%
\e{3}%
\e{2}%
\eol}\vss}\rg%
%
%
\rx{\vss\hfull{%
\rlx{\hss{$120_a$}}\cg%
\e{2}%
\e{4}%
\e{2}%
\e{2}%
\e{1}%
\e{0}%
\e{0}%
\e{3}%
\e{4}%
\e{2}%
\e{1}%
\e{0}%
\e{0}%
\e{0}%
\e{1}%
\e{1}%
\e{0}%
\e{0}%
\e{0}%
\e{0}%
\e{0}%
\eol}\vss}\rg%
%
%
\rx{\vss\hfull{%
\rlx{\hss{$210_a$}}\cg%
\e{1}%
\e{4}%
\e{3}%
\e{4}%
\e{1}%
\e{1}%
\e{0}%
\e{3}%
\e{6}%
\e{3}%
\e{3}%
\e{1}%
\e{0}%
\e{0}%
\e{2}%
\e{2}%
\e{1}%
\e{0}%
\e{0}%
\e{0}%
\e{0}%
\eol}\vss}\rg%
%
%
\rx{\vss\hfull{%
\rlx{\hss{$280_a$}}\cg%
\e{0}%
\e{0}%
\e{0}%
\e{1}%
\e{2}%
\e{3}%
\e{2}%
\e{0}%
\e{1}%
\e{2}%
\e{5}%
\e{4}%
\e{6}%
\e{2}%
\e{0}%
\e{2}%
\e{2}%
\e{5}%
\e{2}%
\e{3}%
\e{0}%
\eol}\vss}\rg%
%
%
\rx{\vss\hfull{%
\rlx{\hss{$336_a$}}\cg%
\e{0}%
\e{0}%
\e{1}%
\e{3}%
\e{2}%
\e{4}%
\e{1}%
\e{0}%
\e{3}%
\e{3}%
\e{7}%
\e{4}%
\e{5}%
\e{1}%
\e{1}%
\e{3}%
\e{3}%
\e{4}%
\e{2}%
\e{1}%
\e{0}%
\eol}\vss}\rg%
%
%
\rx{\vss\hfull{%
\rlx{\hss{$216_a$}}\cg%
\e{0}%
\e{0}%
\e{1}%
\e{1}%
\e{2}%
\e{0}%
\e{0}%
\e{0}%
\e{1}%
\e{3}%
\e{3}%
\e{5}%
\e{3}%
\e{1}%
\e{0}%
\e{0}%
\e{1}%
\e{2}%
\e{3}%
\e{3}%
\e{1}%
\eol}\vss}\rg%
%
%
\rx{\vss\hfull{%
\rlx{\hss{$512_a$}}\cg%
\e{0}%
\e{2}%
\e{4}%
\e{4}%
\e{4}%
\e{2}%
\e{0}%
\e{1}%
\e{5}%
\e{8}%
\e{9}%
\e{8}%
\e{5}%
\e{1}%
\e{0}%
\e{2}%
\e{4}%
\e{4}%
\e{4}%
\e{2}%
\e{0}%
\eol}\vss}\rg%
%
%
\rx{\vss\hfull{%
\rlx{\hss{$378_a$}}\cg%
\e{0}%
\e{1}%
\e{1}%
\e{3}%
\e{2}%
\e{2}%
\e{0}%
\e{0}%
\e{2}%
\e{5}%
\e{6}%
\e{7}%
\e{6}%
\e{1}%
\e{0}%
\e{1}%
\e{4}%
\e{3}%
\e{4}%
\e{3}%
\e{1}%
\eol}\vss}\rg%
%
%
\rx{\vss\hfull{%
\rlx{\hss{$84_a$}}\cg%
\e{0}%
\e{1}%
\e{1}%
\e{0}%
\e{1}%
\e{0}%
\e{0}%
\e{0}%
\e{1}%
\e{2}%
\e{1}%
\e{2}%
\e{0}%
\e{0}%
\e{0}%
\e{0}%
\e{0}%
\e{1}%
\e{0}%
\e{1}%
\e{0}%
\eol}\vss}\rg%
%
%
\rx{\vss\hfull{%
\rlx{\hss{$420_a$}}\cg%
\e{0}%
\e{2}%
\e{3}%
\e{4}%
\e{4}%
\e{3}%
\e{1}%
\e{1}%
\e{6}%
\e{6}%
\e{8}%
\e{5}%
\e{3}%
\e{0}%
\e{1}%
\e{4}%
\e{2}%
\e{4}%
\e{1}%
\e{1}%
\e{0}%
\eol}\vss}\rg%
%
%
\rx{\vss\hfull{%
\rlx{\hss{$280_b$}}\cg%
\e{1}%
\e{4}%
\e{4}%
\e{3}%
\e{2}%
\e{0}%
\e{0}%
\e{2}%
\e{5}%
\e{6}%
\e{4}%
\e{3}%
\e{1}%
\e{0}%
\e{0}%
\e{1}%
\e{2}%
\e{1}%
\e{1}%
\e{0}%
\e{0}%
\eol}\vss}\rg%
%
%
\rx{\vss\hfull{%
\rlx{\hss{$210_b$}}\cg%
\e{1}%
\e{2}%
\e{3}%
\e{1}%
\e{1}%
\e{0}%
\e{0}%
\e{0}%
\e{2}%
\e{5}%
\e{3}%
\e{4}%
\e{1}%
\e{0}%
\e{0}%
\e{0}%
\e{1}%
\e{1}%
\e{2}%
\e{1}%
\e{0}%
\eol}\vss}\rg%
%
%
\rx{\vss\hfull{%
\rlx{\hss{$70_a$}}\cg%
\e{0}%
\e{0}%
\e{1}%
\e{0}%
\e{0}%
\e{0}%
\e{0}%
\e{0}%
\e{0}%
\e{1}%
\e{1}%
\e{2}%
\e{1}%
\e{0}%
\e{0}%
\e{0}%
\e{0}%
\e{0}%
\e{2}%
\e{1}%
\e{1}%
\eol}\vss}\rg%
\tableclose%
%
%
%
%
%
%
\eop
\eject
\tableopen{Induce/restrict matrix for $W({D_{5}}{A_{1}})\,\subset\,W(E_{7})$}%
%
%
%
%
%
%
\rowpts=18 true pt%
\colpts=18 true pt%
\rowlabpts=40 true pt%
\collabpts=85 true pt%
\clx{\vss\hfull{%
\rlx{\hss{$ $}}\cg%
\cx{\hskip 16 true pt\flip{$[{5}:-]{\times}[{2}]$}\hss}\cg%
\cx{\hskip 16 true pt\flip{$[{4}{1}:-]{\times}[{2}]$}\hss}\cg%
\cx{\hskip 16 true pt\flip{$[{3}{2}:-]{\times}[{2}]$}\hss}\cg%
\cx{\hskip 16 true pt\flip{$[{3}{1^{2}}:-]{\times}[{2}]$}\hss}\cg%
\cx{\hskip 16 true pt\flip{$[{2^{2}}{1}:-]{\times}[{2}]$}\hss}\cg%
\cx{\hskip 16 true pt\flip{$[{2}{1^{3}}:-]{\times}[{2}]$}\hss}\cg%
\cx{\hskip 16 true pt\flip{$[{1^{5}}:-]{\times}[{2}]$}\hss}\cg%
\cx{\hskip 16 true pt\flip{$[{4}:{1}]{\times}[{2}]$}\hss}\cg%
\cx{\hskip 16 true pt\flip{$[{3}{1}:{1}]{\times}[{2}]$}\hss}\cg%
\cx{\hskip 16 true pt\flip{$[{2^{2}}:{1}]{\times}[{2}]$}\hss}\cg%
\cx{\hskip 16 true pt\flip{$[{2}{1^{2}}:{1}]{\times}[{2}]$}\hss}\cg%
\cx{\hskip 16 true pt\flip{$[{1^{4}}:{1}]{\times}[{2}]$}\hss}\cg%
\cx{\hskip 16 true pt\flip{$[{3}:{2}]{\times}[{2}]$}\hss}\cg%
\cx{\hskip 16 true pt\flip{$[{3}:{1^{2}}]{\times}[{2}]$}\hss}\cg%
\cx{\hskip 16 true pt\flip{$[{2}{1}:{2}]{\times}[{2}]$}\hss}\cg%
\cx{\hskip 16 true pt\flip{$[{2}{1}:{1^{2}}]{\times}[{2}]$}\hss}\cg%
\cx{\hskip 16 true pt\flip{$[{1^{3}}:{2}]{\times}[{2}]$}\hss}\cg%
\cx{\hskip 16 true pt\flip{$[{1^{3}}:{1^{2}}]{\times}[{2}]$}\hss}\cg%
\eol}}\rg%
%
%
\rx{\vss\hfull{%
\rlx{\hss{$1_a$}}\cg%
\e{1}%
\e{0}%
\e{0}%
\e{0}%
\e{0}%
\e{0}%
\e{0}%
\e{0}%
\e{0}%
\e{0}%
\e{0}%
\e{0}%
\e{0}%
\e{0}%
\e{0}%
\e{0}%
\e{0}%
\e{0}%
\eol}\vss}\rg%
%
%
\rx{\vss\hfull{%
\rlx{\hss{$7_a$}}\cg%
\e{0}%
\e{0}%
\e{0}%
\e{0}%
\e{0}%
\e{0}%
\e{1}%
\e{0}%
\e{0}%
\e{0}%
\e{0}%
\e{0}%
\e{0}%
\e{0}%
\e{0}%
\e{0}%
\e{0}%
\e{0}%
\eol}\vss}\rg%
%
%
\rx{\vss\hfull{%
\rlx{\hss{$27_a$}}\cg%
\e{2}%
\e{1}%
\e{0}%
\e{0}%
\e{0}%
\e{0}%
\e{0}%
\e{1}%
\e{0}%
\e{0}%
\e{0}%
\e{0}%
\e{1}%
\e{0}%
\e{0}%
\e{0}%
\e{0}%
\e{0}%
\eol}\vss}\rg%
%
%
\rx{\vss\hfull{%
\rlx{\hss{$21_a$}}\cg%
\e{0}%
\e{0}%
\e{0}%
\e{0}%
\e{0}%
\e{0}%
\e{0}%
\e{1}%
\e{0}%
\e{0}%
\e{0}%
\e{0}%
\e{0}%
\e{1}%
\e{0}%
\e{0}%
\e{0}%
\e{0}%
\eol}\vss}\rg%
%
%
\rx{\vss\hfull{%
\rlx{\hss{$35_a$}}\cg%
\e{0}%
\e{0}%
\e{0}%
\e{0}%
\e{0}%
\e{0}%
\e{0}%
\e{0}%
\e{0}%
\e{0}%
\e{0}%
\e{1}%
\e{0}%
\e{0}%
\e{0}%
\e{0}%
\e{1}%
\e{0}%
\eol}\vss}\rg%
%
%
\rx{\vss\hfull{%
\rlx{\hss{$105_a$}}\cg%
\e{0}%
\e{0}%
\e{0}%
\e{0}%
\e{0}%
\e{1}%
\e{1}%
\e{0}%
\e{0}%
\e{0}%
\e{0}%
\e{2}%
\e{0}%
\e{0}%
\e{0}%
\e{0}%
\e{1}%
\e{1}%
\eol}\vss}\rg%
%
%
\rx{\vss\hfull{%
\rlx{\hss{$189_a$}}\cg%
\e{0}%
\e{0}%
\e{0}%
\e{0}%
\e{0}%
\e{0}%
\e{0}%
\e{1}%
\e{1}%
\e{0}%
\e{1}%
\e{0}%
\e{0}%
\e{2}%
\e{1}%
\e{1}%
\e{1}%
\e{0}%
\eol}\vss}\rg%
%
%
\rx{\vss\hfull{%
\rlx{\hss{$21_b$}}\cg%
\e{0}%
\e{0}%
\e{0}%
\e{0}%
\e{0}%
\e{1}%
\e{1}%
\e{0}%
\e{0}%
\e{0}%
\e{0}%
\e{0}%
\e{0}%
\e{0}%
\e{0}%
\e{0}%
\e{0}%
\e{0}%
\eol}\vss}\rg%
%
%
\rx{\vss\hfull{%
\rlx{\hss{$35_b$}}\cg%
\e{1}%
\e{1}%
\e{1}%
\e{0}%
\e{0}%
\e{0}%
\e{0}%
\e{1}%
\e{0}%
\e{0}%
\e{0}%
\e{0}%
\e{1}%
\e{0}%
\e{0}%
\e{0}%
\e{0}%
\e{0}%
\eol}\vss}\rg%
%
%
\rx{\vss\hfull{%
\rlx{\hss{$189_b$}}\cg%
\e{0}%
\e{0}%
\e{0}%
\e{0}%
\e{1}%
\e{1}%
\e{0}%
\e{0}%
\e{0}%
\e{0}%
\e{1}%
\e{1}%
\e{0}%
\e{0}%
\e{0}%
\e{1}%
\e{0}%
\e{2}%
\eol}\vss}\rg%
%
%
\rx{\vss\hfull{%
\rlx{\hss{$189_c$}}\cg%
\e{0}%
\e{0}%
\e{0}%
\e{1}%
\e{0}%
\e{2}%
\e{1}%
\e{0}%
\e{0}%
\e{0}%
\e{1}%
\e{1}%
\e{0}%
\e{0}%
\e{0}%
\e{1}%
\e{1}%
\e{1}%
\eol}\vss}\rg%
%
%
\rx{\vss\hfull{%
\rlx{\hss{$15_a$}}\cg%
\e{0}%
\e{0}%
\e{0}%
\e{0}%
\e{1}%
\e{0}%
\e{0}%
\e{0}%
\e{0}%
\e{0}%
\e{0}%
\e{0}%
\e{0}%
\e{0}%
\e{0}%
\e{0}%
\e{0}%
\e{0}%
\eol}\vss}\rg%
%
%
\rx{\vss\hfull{%
\rlx{\hss{$105_b$}}\cg%
\e{0}%
\e{0}%
\e{1}%
\e{0}%
\e{0}%
\e{0}%
\e{0}%
\e{1}%
\e{1}%
\e{0}%
\e{0}%
\e{0}%
\e{2}%
\e{0}%
\e{1}%
\e{0}%
\e{0}%
\e{0}%
\eol}\vss}\rg%
%
%
\rx{\vss\hfull{%
\rlx{\hss{$105_c$}}\cg%
\e{0}%
\e{0}%
\e{0}%
\e{1}%
\e{0}%
\e{1}%
\e{0}%
\e{0}%
\e{1}%
\e{0}%
\e{0}%
\e{0}%
\e{0}%
\e{1}%
\e{0}%
\e{1}%
\e{0}%
\e{0}%
\eol}\vss}\rg%
%
%
\rx{\vss\hfull{%
\rlx{\hss{$315_a$}}\cg%
\e{0}%
\e{0}%
\e{0}%
\e{0}%
\e{0}%
\e{0}%
\e{0}%
\e{0}%
\e{0}%
\e{1}%
\e{2}%
\e{1}%
\e{0}%
\e{0}%
\e{1}%
\e{2}%
\e{1}%
\e{2}%
\eol}\vss}\rg%
%
%
\rx{\vss\hfull{%
\rlx{\hss{$405_a$}}\cg%
\e{0}%
\e{0}%
\e{0}%
\e{0}%
\e{0}%
\e{0}%
\e{0}%
\e{1}%
\e{3}%
\e{1}%
\e{1}%
\e{0}%
\e{2}%
\e{2}%
\e{3}%
\e{2}%
\e{1}%
\e{0}%
\eol}\vss}\rg%
%
%
\rx{\vss\hfull{%
\rlx{\hss{$168_a$}}\cg%
\e{1}%
\e{2}%
\e{1}%
\e{0}%
\e{0}%
\e{0}%
\e{0}%
\e{1}%
\e{1}%
\e{1}%
\e{0}%
\e{0}%
\e{2}%
\e{0}%
\e{2}%
\e{0}%
\e{0}%
\e{0}%
\eol}\vss}\rg%
%
%
\rx{\vss\hfull{%
\rlx{\hss{$56_a$}}\cg%
\e{0}%
\e{0}%
\e{0}%
\e{0}%
\e{0}%
\e{0}%
\e{1}%
\e{0}%
\e{0}%
\e{0}%
\e{0}%
\e{1}%
\e{0}%
\e{0}%
\e{0}%
\e{0}%
\e{0}%
\e{1}%
\eol}\vss}\rg%
%
%
\rx{\vss\hfull{%
\rlx{\hss{$120_a$}}\cg%
\e{1}%
\e{1}%
\e{0}%
\e{0}%
\e{0}%
\e{0}%
\e{0}%
\e{2}%
\e{1}%
\e{0}%
\e{0}%
\e{0}%
\e{2}%
\e{1}%
\e{1}%
\e{0}%
\e{0}%
\e{0}%
\eol}\vss}\rg%
%
%
\rx{\vss\hfull{%
\rlx{\hss{$210_a$}}\cg%
\e{0}%
\e{1}%
\e{0}%
\e{1}%
\e{0}%
\e{0}%
\e{0}%
\e{2}%
\e{2}%
\e{0}%
\e{0}%
\e{0}%
\e{2}%
\e{2}%
\e{1}%
\e{1}%
\e{0}%
\e{0}%
\eol}\vss}\rg%
%
%
\rx{\vss\hfull{%
\rlx{\hss{$280_a$}}\cg%
\e{0}%
\e{0}%
\e{0}%
\e{0}%
\e{0}%
\e{0}%
\e{0}%
\e{0}%
\e{0}%
\e{0}%
\e{2}%
\e{2}%
\e{0}%
\e{0}%
\e{1}%
\e{1}%
\e{2}%
\e{2}%
\eol}\vss}\rg%
%
%
\rx{\vss\hfull{%
\rlx{\hss{$336_a$}}\cg%
\e{0}%
\e{0}%
\e{0}%
\e{1}%
\e{0}%
\e{1}%
\e{0}%
\e{0}%
\e{1}%
\e{0}%
\e{2}%
\e{1}%
\e{0}%
\e{1}%
\e{1}%
\e{2}%
\e{2}%
\e{1}%
\eol}\vss}\rg%
%
%
\rx{\vss\hfull{%
\rlx{\hss{$216_a$}}\cg%
\e{0}%
\e{0}%
\e{1}%
\e{1}%
\e{2}%
\e{0}%
\e{0}%
\e{0}%
\e{0}%
\e{1}%
\e{1}%
\e{0}%
\e{0}%
\e{0}%
\e{1}%
\e{1}%
\e{0}%
\e{1}%
\eol}\vss}\rg%
%
%
\rx{\vss\hfull{%
\rlx{\hss{$512_a$}}\cg%
\e{0}%
\e{0}%
\e{1}%
\e{1}%
\e{1}%
\e{0}%
\e{0}%
\e{0}%
\e{2}%
\e{2}%
\e{2}%
\e{0}%
\e{1}%
\e{1}%
\e{3}%
\e{3}%
\e{1}%
\e{1}%
\eol}\vss}\rg%
%
%
\rx{\vss\hfull{%
\rlx{\hss{$378_a$}}\cg%
\e{0}%
\e{0}%
\e{0}%
\e{0}%
\e{1}%
\e{1}%
\e{0}%
\e{0}%
\e{1}%
\e{1}%
\e{2}%
\e{0}%
\e{0}%
\e{1}%
\e{1}%
\e{3}%
\e{1}%
\e{1}%
\eol}\vss}\rg%
%
%
\rx{\vss\hfull{%
\rlx{\hss{$84_a$}}\cg%
\e{0}%
\e{1}%
\e{1}%
\e{0}%
\e{1}%
\e{0}%
\e{0}%
\e{0}%
\e{0}%
\e{1}%
\e{0}%
\e{0}%
\e{0}%
\e{0}%
\e{1}%
\e{0}%
\e{0}%
\e{0}%
\eol}\vss}\rg%
%
%
\rx{\vss\hfull{%
\rlx{\hss{$420_a$}}\cg%
\e{0}%
\e{1}%
\e{0}%
\e{1}%
\e{0}%
\e{0}%
\e{0}%
\e{0}%
\e{2}%
\e{1}%
\e{2}%
\e{0}%
\e{1}%
\e{1}%
\e{3}%
\e{2}%
\e{1}%
\e{1}%
\eol}\vss}\rg%
%
%
\rx{\vss\hfull{%
\rlx{\hss{$280_b$}}\cg%
\e{0}%
\e{1}%
\e{2}%
\e{1}%
\e{1}%
\e{0}%
\e{0}%
\e{1}%
\e{2}%
\e{1}%
\e{0}%
\e{0}%
\e{2}%
\e{1}%
\e{2}%
\e{1}%
\e{0}%
\e{0}%
\eol}\vss}\rg%
%
%
\rx{\vss\hfull{%
\rlx{\hss{$210_b$}}\cg%
\e{0}%
\e{0}%
\e{1}%
\e{0}%
\e{0}%
\e{0}%
\e{0}%
\e{0}%
\e{1}%
\e{2}%
\e{0}%
\e{0}%
\e{1}%
\e{0}%
\e{2}%
\e{1}%
\e{0}%
\e{0}%
\eol}\vss}\rg%
%
%
\rx{\vss\hfull{%
\rlx{\hss{$70_a$}}\cg%
\e{0}%
\e{0}%
\e{0}%
\e{0}%
\e{0}%
\e{0}%
\e{0}%
\e{0}%
\e{0}%
\e{1}%
\e{0}%
\e{0}%
\e{0}%
\e{0}%
\e{0}%
\e{1}%
\e{0}%
\e{0}%
\eol}\vss}\rg%
\eop
\eject
\tablecont%
%
%
%
%
%
%
\rowpts=18 true pt%
\colpts=18 true pt%
\rowlabpts=40 true pt%
\collabpts=85 true pt%
\clx{\vss\hfull{%
\rlx{\hss{$ $}}\cg%
\cx{\hskip 16 true pt\flip{$[{5}:-]{\times}[{1^{2}}]$}\hss}\cg%
\cx{\hskip 16 true pt\flip{$[{4}{1}:-]{\times}[{1^{2}}]$}\hss}\cg%
\cx{\hskip 16 true pt\flip{$[{3}{2}:-]{\times}[{1^{2}}]$}\hss}\cg%
\cx{\hskip 16 true pt\flip{$[{3}{1^{2}}:-]{\times}[{1^{2}}]$}\hss}\cg%
\cx{\hskip 16 true pt\flip{$[{2^{2}}{1}:-]{\times}[{1^{2}}]$}\hss}\cg%
\cx{\hskip 16 true pt\flip{$[{2}{1^{3}}:-]{\times}[{1^{2}}]$}\hss}\cg%
\cx{\hskip 16 true pt\flip{$[{1^{5}}:-]{\times}[{1^{2}}]$}\hss}\cg%
\cx{\hskip 16 true pt\flip{$[{4}:{1}]{\times}[{1^{2}}]$}\hss}\cg%
\cx{\hskip 16 true pt\flip{$[{3}{1}:{1}]{\times}[{1^{2}}]$}\hss}\cg%
\cx{\hskip 16 true pt\flip{$[{2^{2}}:{1}]{\times}[{1^{2}}]$}\hss}\cg%
\cx{\hskip 16 true pt\flip{$[{2}{1^{2}}:{1}]{\times}[{1^{2}}]$}\hss}\cg%
\cx{\hskip 16 true pt\flip{$[{1^{4}}:{1}]{\times}[{1^{2}}]$}\hss}\cg%
\cx{\hskip 16 true pt\flip{$[{3}:{2}]{\times}[{1^{2}}]$}\hss}\cg%
\cx{\hskip 16 true pt\flip{$[{3}:{1^{2}}]{\times}[{1^{2}}]$}\hss}\cg%
\cx{\hskip 16 true pt\flip{$[{2}{1}:{2}]{\times}[{1^{2}}]$}\hss}\cg%
\cx{\hskip 16 true pt\flip{$[{2}{1}:{1^{2}}]{\times}[{1^{2}}]$}\hss}\cg%
\cx{\hskip 16 true pt\flip{$[{1^{3}}:{2}]{\times}[{1^{2}}]$}\hss}\cg%
\cx{\hskip 16 true pt\flip{$[{1^{3}}:{1^{2}}]{\times}[{1^{2}}]$}\hss}\cg%
\eol}}\rg%
%
%
\rx{\vss\hfull{%
\rlx{\hss{$1_a$}}\cg%
\e{0}%
\e{0}%
\e{0}%
\e{0}%
\e{0}%
\e{0}%
\e{0}%
\e{0}%
\e{0}%
\e{0}%
\e{0}%
\e{0}%
\e{0}%
\e{0}%
\e{0}%
\e{0}%
\e{0}%
\e{0}%
\eol}\vss}\rg%
%
%
\rx{\vss\hfull{%
\rlx{\hss{$7_a$}}\cg%
\e{0}%
\e{0}%
\e{0}%
\e{0}%
\e{0}%
\e{0}%
\e{1}%
\e{0}%
\e{0}%
\e{0}%
\e{0}%
\e{1}%
\e{0}%
\e{0}%
\e{0}%
\e{0}%
\e{0}%
\e{0}%
\eol}\vss}\rg%
%
%
\rx{\vss\hfull{%
\rlx{\hss{$27_a$}}\cg%
\e{1}%
\e{0}%
\e{0}%
\e{0}%
\e{0}%
\e{0}%
\e{0}%
\e{1}%
\e{0}%
\e{0}%
\e{0}%
\e{0}%
\e{0}%
\e{0}%
\e{0}%
\e{0}%
\e{0}%
\e{0}%
\eol}\vss}\rg%
%
%
\rx{\vss\hfull{%
\rlx{\hss{$21_a$}}\cg%
\e{1}%
\e{0}%
\e{0}%
\e{0}%
\e{0}%
\e{0}%
\e{0}%
\e{1}%
\e{0}%
\e{0}%
\e{0}%
\e{0}%
\e{0}%
\e{0}%
\e{0}%
\e{0}%
\e{0}%
\e{0}%
\eol}\vss}\rg%
%
%
\rx{\vss\hfull{%
\rlx{\hss{$35_a$}}\cg%
\e{0}%
\e{0}%
\e{0}%
\e{0}%
\e{0}%
\e{0}%
\e{0}%
\e{0}%
\e{0}%
\e{0}%
\e{0}%
\e{0}%
\e{0}%
\e{1}%
\e{0}%
\e{0}%
\e{1}%
\e{0}%
\eol}\vss}\rg%
%
%
\rx{\vss\hfull{%
\rlx{\hss{$105_a$}}\cg%
\e{0}%
\e{0}%
\e{0}%
\e{0}%
\e{0}%
\e{1}%
\e{1}%
\e{0}%
\e{0}%
\e{0}%
\e{1}%
\e{2}%
\e{0}%
\e{0}%
\e{0}%
\e{1}%
\e{1}%
\e{1}%
\eol}\vss}\rg%
%
%
\rx{\vss\hfull{%
\rlx{\hss{$189_a$}}\cg%
\e{0}%
\e{1}%
\e{0}%
\e{0}%
\e{0}%
\e{0}%
\e{0}%
\e{1}%
\e{1}%
\e{0}%
\e{0}%
\e{0}%
\e{1}%
\e{2}%
\e{1}%
\e{0}%
\e{1}%
\e{0}%
\eol}\vss}\rg%
%
%
\rx{\vss\hfull{%
\rlx{\hss{$21_b$}}\cg%
\e{0}%
\e{0}%
\e{0}%
\e{0}%
\e{0}%
\e{0}%
\e{1}%
\e{0}%
\e{0}%
\e{0}%
\e{0}%
\e{1}%
\e{0}%
\e{0}%
\e{0}%
\e{0}%
\e{0}%
\e{1}%
\eol}\vss}\rg%
%
%
\rx{\vss\hfull{%
\rlx{\hss{$35_b$}}\cg%
\e{0}%
\e{0}%
\e{0}%
\e{0}%
\e{0}%
\e{0}%
\e{0}%
\e{0}%
\e{0}%
\e{0}%
\e{0}%
\e{0}%
\e{1}%
\e{0}%
\e{0}%
\e{0}%
\e{0}%
\e{0}%
\eol}\vss}\rg%
%
%
\rx{\vss\hfull{%
\rlx{\hss{$189_b$}}\cg%
\e{0}%
\e{0}%
\e{0}%
\e{0}%
\e{1}%
\e{1}%
\e{1}%
\e{0}%
\e{0}%
\e{1}%
\e{2}%
\e{2}%
\e{0}%
\e{0}%
\e{0}%
\e{1}%
\e{1}%
\e{3}%
\eol}\vss}\rg%
%
%
\rx{\vss\hfull{%
\rlx{\hss{$189_c$}}\cg%
\e{0}%
\e{0}%
\e{0}%
\e{0}%
\e{0}%
\e{1}%
\e{0}%
\e{0}%
\e{0}%
\e{0}%
\e{1}%
\e{1}%
\e{0}%
\e{0}%
\e{1}%
\e{2}%
\e{1}%
\e{2}%
\eol}\vss}\rg%
%
%
\rx{\vss\hfull{%
\rlx{\hss{$15_a$}}\cg%
\e{0}%
\e{0}%
\e{0}%
\e{0}%
\e{0}%
\e{0}%
\e{0}%
\e{0}%
\e{0}%
\e{0}%
\e{0}%
\e{0}%
\e{0}%
\e{0}%
\e{0}%
\e{0}%
\e{0}%
\e{1}%
\eol}\vss}\rg%
%
%
\rx{\vss\hfull{%
\rlx{\hss{$105_b$}}\cg%
\e{0}%
\e{0}%
\e{1}%
\e{0}%
\e{0}%
\e{0}%
\e{0}%
\e{0}%
\e{1}%
\e{1}%
\e{0}%
\e{0}%
\e{1}%
\e{0}%
\e{0}%
\e{0}%
\e{0}%
\e{0}%
\eol}\vss}\rg%
%
%
\rx{\vss\hfull{%
\rlx{\hss{$105_c$}}\cg%
\e{0}%
\e{0}%
\e{0}%
\e{0}%
\e{0}%
\e{0}%
\e{0}%
\e{0}%
\e{0}%
\e{0}%
\e{0}%
\e{0}%
\e{1}%
\e{0}%
\e{1}%
\e{1}%
\e{0}%
\e{0}%
\eol}\vss}\rg%
%
%
\rx{\vss\hfull{%
\rlx{\hss{$315_a$}}\cg%
\e{0}%
\e{0}%
\e{0}%
\e{1}%
\e{1}%
\e{1}%
\e{0}%
\e{0}%
\e{1}%
\e{1}%
\e{3}%
\e{1}%
\e{0}%
\e{0}%
\e{1}%
\e{2}%
\e{1}%
\e{2}%
\eol}\vss}\rg%
%
%
\rx{\vss\hfull{%
\rlx{\hss{$405_a$}}\cg%
\e{0}%
\e{1}%
\e{1}%
\e{1}%
\e{0}%
\e{0}%
\e{0}%
\e{1}%
\e{3}%
\e{1}%
\e{1}%
\e{0}%
\e{2}%
\e{1}%
\e{2}%
\e{1}%
\e{0}%
\e{0}%
\eol}\vss}\rg%
%
%
\rx{\vss\hfull{%
\rlx{\hss{$168_a$}}\cg%
\e{0}%
\e{1}%
\e{0}%
\e{0}%
\e{0}%
\e{0}%
\e{0}%
\e{1}%
\e{1}%
\e{0}%
\e{0}%
\e{0}%
\e{1}%
\e{1}%
\e{1}%
\e{0}%
\e{0}%
\e{0}%
\eol}\vss}\rg%
%
%
\rx{\vss\hfull{%
\rlx{\hss{$56_a$}}\cg%
\e{0}%
\e{0}%
\e{0}%
\e{0}%
\e{0}%
\e{1}%
\e{1}%
\e{0}%
\e{0}%
\e{0}%
\e{1}%
\e{2}%
\e{0}%
\e{0}%
\e{0}%
\e{0}%
\e{0}%
\e{1}%
\eol}\vss}\rg%
%
%
\rx{\vss\hfull{%
\rlx{\hss{$120_a$}}\cg%
\e{1}%
\e{1}%
\e{0}%
\e{0}%
\e{0}%
\e{0}%
\e{0}%
\e{2}%
\e{1}%
\e{0}%
\e{0}%
\e{0}%
\e{1}%
\e{0}%
\e{0}%
\e{0}%
\e{0}%
\e{0}%
\eol}\vss}\rg%
%
%
\rx{\vss\hfull{%
\rlx{\hss{$210_a$}}\cg%
\e{1}%
\e{1}%
\e{0}%
\e{0}%
\e{0}%
\e{0}%
\e{0}%
\e{2}%
\e{1}%
\e{0}%
\e{0}%
\e{0}%
\e{2}%
\e{1}%
\e{1}%
\e{0}%
\e{0}%
\e{0}%
\eol}\vss}\rg%
%
%
\rx{\vss\hfull{%
\rlx{\hss{$280_a$}}\cg%
\e{0}%
\e{0}%
\e{0}%
\e{1}%
\e{0}%
\e{1}%
\e{0}%
\e{0}%
\e{1}%
\e{0}%
\e{2}%
\e{1}%
\e{0}%
\e{1}%
\e{1}%
\e{2}%
\e{2}%
\e{1}%
\eol}\vss}\rg%
%
%
\rx{\vss\hfull{%
\rlx{\hss{$336_a$}}\cg%
\e{0}%
\e{0}%
\e{0}%
\e{1}%
\e{0}%
\e{0}%
\e{0}%
\e{0}%
\e{1}%
\e{0}%
\e{1}%
\e{0}%
\e{1}%
\e{2}%
\e{2}%
\e{2}%
\e{2}%
\e{1}%
\eol}\vss}\rg%
%
%
\rx{\vss\hfull{%
\rlx{\hss{$216_a$}}\cg%
\e{0}%
\e{0}%
\e{0}%
\e{0}%
\e{1}%
\e{0}%
\e{0}%
\e{0}%
\e{0}%
\e{1}%
\e{1}%
\e{0}%
\e{0}%
\e{0}%
\e{1}%
\e{2}%
\e{1}%
\e{2}%
\eol}\vss}\rg%
%
%
\rx{\vss\hfull{%
\rlx{\hss{$512_a$}}\cg%
\e{0}%
\e{0}%
\e{1}%
\e{1}%
\e{1}%
\e{0}%
\e{0}%
\e{0}%
\e{2}%
\e{2}%
\e{2}%
\e{0}%
\e{1}%
\e{1}%
\e{3}%
\e{3}%
\e{1}%
\e{1}%
\eol}\vss}\rg%
%
%
\rx{\vss\hfull{%
\rlx{\hss{$378_a$}}\cg%
\e{0}%
\e{0}%
\e{1}%
\e{0}%
\e{1}%
\e{1}%
\e{0}%
\e{0}%
\e{1}%
\e{2}%
\e{2}%
\e{1}%
\e{0}%
\e{0}%
\e{2}%
\e{3}%
\e{1}%
\e{1}%
\eol}\vss}\rg%
%
%
\rx{\vss\hfull{%
\rlx{\hss{$84_a$}}\cg%
\e{0}%
\e{0}%
\e{0}%
\e{0}%
\e{0}%
\e{0}%
\e{0}%
\e{0}%
\e{0}%
\e{0}%
\e{0}%
\e{0}%
\e{0}%
\e{0}%
\e{1}%
\e{1}%
\e{0}%
\e{0}%
\eol}\vss}\rg%
%
%
\rx{\vss\hfull{%
\rlx{\hss{$420_a$}}\cg%
\e{0}%
\e{1}%
\e{0}%
\e{1}%
\e{0}%
\e{0}%
\e{0}%
\e{1}%
\e{2}%
\e{0}%
\e{1}%
\e{0}%
\e{1}%
\e{2}%
\e{3}%
\e{2}%
\e{1}%
\e{0}%
\eol}\vss}\rg%
%
%
\rx{\vss\hfull{%
\rlx{\hss{$280_b$}}\cg%
\e{0}%
\e{0}%
\e{1}%
\e{0}%
\e{0}%
\e{0}%
\e{0}%
\e{0}%
\e{1}%
\e{1}%
\e{0}%
\e{0}%
\e{2}%
\e{1}%
\e{2}%
\e{1}%
\e{0}%
\e{0}%
\eol}\vss}\rg%
%
%
\rx{\vss\hfull{%
\rlx{\hss{$210_b$}}\cg%
\e{0}%
\e{0}%
\e{1}%
\e{0}%
\e{1}%
\e{0}%
\e{0}%
\e{0}%
\e{1}%
\e{2}%
\e{1}%
\e{0}%
\e{0}%
\e{0}%
\e{1}%
\e{1}%
\e{0}%
\e{0}%
\eol}\vss}\rg%
%
%
\rx{\vss\hfull{%
\rlx{\hss{$70_a$}}\cg%
\e{0}%
\e{0}%
\e{0}%
\e{0}%
\e{1}%
\e{0}%
\e{0}%
\e{0}%
\e{0}%
\e{1}%
\e{1}%
\e{0}%
\e{0}%
\e{0}%
\e{0}%
\e{0}%
\e{0}%
\e{1}%
\eol}\vss}\rg%
\tableclose%
%
%
%
%
%
%
\eop
\eject
\tableopen{Induce/restrict matrix for $W(E_{6})\,\subset\,W(E_{7})$}%
%
%
%
%
%
%
\rowpts=18 true pt%
\colpts=18 true pt%
\rowlabpts=40 true pt%
\collabpts=40 true pt%
\clx{\vss\hfull{%
\rlx{\hss{$ $}}\cg%
\cx{\hskip 16 true pt\flip{$1_p$}\hss}\cg%
\cx{\hskip 16 true pt\flip{$6_p$}\hss}\cg%
\cx{\hskip 16 true pt\flip{$15_p$}\hss}\cg%
\cx{\hskip 16 true pt\flip{$20_p$}\hss}\cg%
\cx{\hskip 16 true pt\flip{$30_p$}\hss}\cg%
\cx{\hskip 16 true pt\flip{$64_p$}\hss}\cg%
\cx{\hskip 16 true pt\flip{$81_p$}\hss}\cg%
\cx{\hskip 16 true pt\flip{$15_q$}\hss}\cg%
\cx{\hskip 16 true pt\flip{$24_p$}\hss}\cg%
\cx{\hskip 16 true pt\flip{$60_p$}\hss}\cg%
\cx{\hskip 16 true pt\flip{$1_p^{*}$}\hss}\cg%
\cx{\hskip 16 true pt\flip{$6_p^{*}$}\hss}\cg%
\cx{\hskip 16 true pt\flip{$15_p^{*}$}\hss}\cg%
\cx{\hskip 16 true pt\flip{$20_p^{*}$}\hss}\cg%
\cx{\hskip 16 true pt\flip{$30_p^{*}$}\hss}\cg%
\cx{\hskip 16 true pt\flip{$64_p^{*}$}\hss}\cg%
\cx{\hskip 16 true pt\flip{$81_p^{*}$}\hss}\cg%
\cx{\hskip 16 true pt\flip{$15_q^{*}$}\hss}\cg%
\cx{\hskip 16 true pt\flip{$24_p^{*}$}\hss}\cg%
\cx{\hskip 16 true pt\flip{$60_p^{*}$}\hss}\cg%
\cx{\hskip 16 true pt\flip{$20_s$}\hss}\cg%
\cx{\hskip 16 true pt\flip{$90_s$}\hss}\cg%
\cx{\hskip 16 true pt\flip{$80_s$}\hss}\cg%
\cx{\hskip 16 true pt\flip{$60_s$}\hss}\cg%
\cx{\hskip 16 true pt\flip{$10_s$}\hss}\cg%
\eol}}\rg%
%
%
\rx{\vss\hfull{%
\rlx{\hss{$1_a$}}\cg%
\e{1}%
\e{0}%
\e{0}%
\e{0}%
\e{0}%
\e{0}%
\e{0}%
\e{0}%
\e{0}%
\e{0}%
\e{0}%
\e{0}%
\e{0}%
\e{0}%
\e{0}%
\e{0}%
\e{0}%
\e{0}%
\e{0}%
\e{0}%
\e{0}%
\e{0}%
\e{0}%
\e{0}%
\e{0}%
\eol}\vss}\rg%
%
%
\rx{\vss\hfull{%
\rlx{\hss{$7_a$}}\cg%
\e{0}%
\e{0}%
\e{0}%
\e{0}%
\e{0}%
\e{0}%
\e{0}%
\e{0}%
\e{0}%
\e{0}%
\e{1}%
\e{1}%
\e{0}%
\e{0}%
\e{0}%
\e{0}%
\e{0}%
\e{0}%
\e{0}%
\e{0}%
\e{0}%
\e{0}%
\e{0}%
\e{0}%
\e{0}%
\eol}\vss}\rg%
%
%
\rx{\vss\hfull{%
\rlx{\hss{$27_a$}}\cg%
\e{1}%
\e{1}%
\e{0}%
\e{1}%
\e{0}%
\e{0}%
\e{0}%
\e{0}%
\e{0}%
\e{0}%
\e{0}%
\e{0}%
\e{0}%
\e{0}%
\e{0}%
\e{0}%
\e{0}%
\e{0}%
\e{0}%
\e{0}%
\e{0}%
\e{0}%
\e{0}%
\e{0}%
\e{0}%
\eol}\vss}\rg%
%
%
\rx{\vss\hfull{%
\rlx{\hss{$21_a$}}\cg%
\e{0}%
\e{1}%
\e{1}%
\e{0}%
\e{0}%
\e{0}%
\e{0}%
\e{0}%
\e{0}%
\e{0}%
\e{0}%
\e{0}%
\e{0}%
\e{0}%
\e{0}%
\e{0}%
\e{0}%
\e{0}%
\e{0}%
\e{0}%
\e{0}%
\e{0}%
\e{0}%
\e{0}%
\e{0}%
\eol}\vss}\rg%
%
%
\rx{\vss\hfull{%
\rlx{\hss{$35_a$}}\cg%
\e{0}%
\e{0}%
\e{0}%
\e{0}%
\e{0}%
\e{0}%
\e{0}%
\e{0}%
\e{0}%
\e{0}%
\e{0}%
\e{0}%
\e{1}%
\e{0}%
\e{0}%
\e{0}%
\e{0}%
\e{0}%
\e{0}%
\e{0}%
\e{1}%
\e{0}%
\e{0}%
\e{0}%
\e{0}%
\eol}\vss}\rg%
%
%
\rx{\vss\hfull{%
\rlx{\hss{$105_a$}}\cg%
\e{0}%
\e{0}%
\e{0}%
\e{0}%
\e{0}%
\e{0}%
\e{0}%
\e{0}%
\e{0}%
\e{0}%
\e{0}%
\e{1}%
\e{1}%
\e{1}%
\e{0}%
\e{1}%
\e{0}%
\e{0}%
\e{0}%
\e{0}%
\e{0}%
\e{0}%
\e{0}%
\e{0}%
\e{0}%
\eol}\vss}\rg%
%
%
\rx{\vss\hfull{%
\rlx{\hss{$189_a$}}\cg%
\e{0}%
\e{0}%
\e{1}%
\e{0}%
\e{0}%
\e{1}%
\e{0}%
\e{0}%
\e{0}%
\e{0}%
\e{0}%
\e{0}%
\e{0}%
\e{0}%
\e{0}%
\e{0}%
\e{0}%
\e{0}%
\e{0}%
\e{0}%
\e{1}%
\e{1}%
\e{0}%
\e{0}%
\e{0}%
\eol}\vss}\rg%
%
%
\rx{\vss\hfull{%
\rlx{\hss{$21_b$}}\cg%
\e{0}%
\e{0}%
\e{0}%
\e{0}%
\e{0}%
\e{0}%
\e{0}%
\e{0}%
\e{0}%
\e{0}%
\e{1}%
\e{0}%
\e{0}%
\e{1}%
\e{0}%
\e{0}%
\e{0}%
\e{0}%
\e{0}%
\e{0}%
\e{0}%
\e{0}%
\e{0}%
\e{0}%
\e{0}%
\eol}\vss}\rg%
%
%
\rx{\vss\hfull{%
\rlx{\hss{$35_b$}}\cg%
\e{0}%
\e{0}%
\e{0}%
\e{1}%
\e{0}%
\e{0}%
\e{0}%
\e{1}%
\e{0}%
\e{0}%
\e{0}%
\e{0}%
\e{0}%
\e{0}%
\e{0}%
\e{0}%
\e{0}%
\e{0}%
\e{0}%
\e{0}%
\e{0}%
\e{0}%
\e{0}%
\e{0}%
\e{0}%
\eol}\vss}\rg%
%
%
\rx{\vss\hfull{%
\rlx{\hss{$189_b$}}\cg%
\e{0}%
\e{0}%
\e{0}%
\e{0}%
\e{0}%
\e{0}%
\e{0}%
\e{0}%
\e{0}%
\e{0}%
\e{0}%
\e{0}%
\e{0}%
\e{1}%
\e{1}%
\e{1}%
\e{0}%
\e{1}%
\e{0}%
\e{1}%
\e{0}%
\e{0}%
\e{0}%
\e{0}%
\e{0}%
\eol}\vss}\rg%
%
%
\rx{\vss\hfull{%
\rlx{\hss{$189_c$}}\cg%
\e{0}%
\e{0}%
\e{0}%
\e{0}%
\e{0}%
\e{0}%
\e{0}%
\e{0}%
\e{0}%
\e{0}%
\e{0}%
\e{0}%
\e{0}%
\e{1}%
\e{0}%
\e{1}%
\e{1}%
\e{0}%
\e{1}%
\e{0}%
\e{0}%
\e{0}%
\e{0}%
\e{0}%
\e{0}%
\eol}\vss}\rg%
%
%
\rx{\vss\hfull{%
\rlx{\hss{$15_a$}}\cg%
\e{0}%
\e{0}%
\e{0}%
\e{0}%
\e{0}%
\e{0}%
\e{0}%
\e{0}%
\e{0}%
\e{0}%
\e{0}%
\e{0}%
\e{0}%
\e{0}%
\e{0}%
\e{0}%
\e{0}%
\e{1}%
\e{0}%
\e{0}%
\e{0}%
\e{0}%
\e{0}%
\e{0}%
\e{0}%
\eol}\vss}\rg%
%
%
\rx{\vss\hfull{%
\rlx{\hss{$105_b$}}\cg%
\e{0}%
\e{0}%
\e{0}%
\e{0}%
\e{1}%
\e{0}%
\e{0}%
\e{1}%
\e{0}%
\e{1}%
\e{0}%
\e{0}%
\e{0}%
\e{0}%
\e{0}%
\e{0}%
\e{0}%
\e{0}%
\e{0}%
\e{0}%
\e{0}%
\e{0}%
\e{0}%
\e{0}%
\e{0}%
\eol}\vss}\rg%
%
%
\rx{\vss\hfull{%
\rlx{\hss{$105_c$}}\cg%
\e{0}%
\e{0}%
\e{0}%
\e{0}%
\e{0}%
\e{0}%
\e{1}%
\e{0}%
\e{0}%
\e{0}%
\e{0}%
\e{0}%
\e{0}%
\e{0}%
\e{0}%
\e{0}%
\e{0}%
\e{0}%
\e{1}%
\e{0}%
\e{0}%
\e{0}%
\e{0}%
\e{0}%
\e{0}%
\eol}\vss}\rg%
%
%
\rx{\vss\hfull{%
\rlx{\hss{$315_a$}}\cg%
\e{0}%
\e{0}%
\e{0}%
\e{0}%
\e{0}%
\e{0}%
\e{0}%
\e{0}%
\e{0}%
\e{0}%
\e{0}%
\e{0}%
\e{0}%
\e{0}%
\e{1}%
\e{1}%
\e{1}%
\e{0}%
\e{0}%
\e{1}%
\e{0}%
\e{0}%
\e{1}%
\e{0}%
\e{0}%
\eol}\vss}\rg%
%
%
\rx{\vss\hfull{%
\rlx{\hss{$405_a$}}\cg%
\e{0}%
\e{0}%
\e{0}%
\e{0}%
\e{1}%
\e{1}%
\e{1}%
\e{0}%
\e{0}%
\e{1}%
\e{0}%
\e{0}%
\e{0}%
\e{0}%
\e{0}%
\e{0}%
\e{0}%
\e{0}%
\e{0}%
\e{0}%
\e{0}%
\e{1}%
\e{1}%
\e{0}%
\e{0}%
\eol}\vss}\rg%
%
%
\rx{\vss\hfull{%
\rlx{\hss{$168_a$}}\cg%
\e{0}%
\e{0}%
\e{0}%
\e{1}%
\e{0}%
\e{1}%
\e{0}%
\e{0}%
\e{1}%
\e{1}%
\e{0}%
\e{0}%
\e{0}%
\e{0}%
\e{0}%
\e{0}%
\e{0}%
\e{0}%
\e{0}%
\e{0}%
\e{0}%
\e{0}%
\e{0}%
\e{0}%
\e{0}%
\eol}\vss}\rg%
%
%
\rx{\vss\hfull{%
\rlx{\hss{$56_a$}}\cg%
\e{0}%
\e{0}%
\e{0}%
\e{0}%
\e{0}%
\e{0}%
\e{0}%
\e{0}%
\e{0}%
\e{0}%
\e{0}%
\e{1}%
\e{0}%
\e{1}%
\e{1}%
\e{0}%
\e{0}%
\e{0}%
\e{0}%
\e{0}%
\e{0}%
\e{0}%
\e{0}%
\e{0}%
\e{0}%
\eol}\vss}\rg%
%
%
\rx{\vss\hfull{%
\rlx{\hss{$120_a$}}\cg%
\e{0}%
\e{1}%
\e{0}%
\e{1}%
\e{1}%
\e{1}%
\e{0}%
\e{0}%
\e{0}%
\e{0}%
\e{0}%
\e{0}%
\e{0}%
\e{0}%
\e{0}%
\e{0}%
\e{0}%
\e{0}%
\e{0}%
\e{0}%
\e{0}%
\e{0}%
\e{0}%
\e{0}%
\e{0}%
\eol}\vss}\rg%
%
%
\rx{\vss\hfull{%
\rlx{\hss{$210_a$}}\cg%
\e{0}%
\e{0}%
\e{1}%
\e{1}%
\e{1}%
\e{1}%
\e{1}%
\e{0}%
\e{0}%
\e{0}%
\e{0}%
\e{0}%
\e{0}%
\e{0}%
\e{0}%
\e{0}%
\e{0}%
\e{0}%
\e{0}%
\e{0}%
\e{0}%
\e{0}%
\e{0}%
\e{0}%
\e{0}%
\eol}\vss}\rg%
%
%
\rx{\vss\hfull{%
\rlx{\hss{$280_a$}}\cg%
\e{0}%
\e{0}%
\e{0}%
\e{0}%
\e{0}%
\e{0}%
\e{0}%
\e{0}%
\e{0}%
\e{0}%
\e{0}%
\e{0}%
\e{1}%
\e{0}%
\e{1}%
\e{1}%
\e{1}%
\e{0}%
\e{0}%
\e{0}%
\e{0}%
\e{1}%
\e{0}%
\e{0}%
\e{0}%
\eol}\vss}\rg%
%
%
\rx{\vss\hfull{%
\rlx{\hss{$336_a$}}\cg%
\e{0}%
\e{0}%
\e{0}%
\e{0}%
\e{0}%
\e{0}%
\e{1}%
\e{0}%
\e{0}%
\e{0}%
\e{0}%
\e{0}%
\e{0}%
\e{0}%
\e{0}%
\e{1}%
\e{1}%
\e{0}%
\e{0}%
\e{0}%
\e{1}%
\e{1}%
\e{0}%
\e{0}%
\e{0}%
\eol}\vss}\rg%
%
%
\rx{\vss\hfull{%
\rlx{\hss{$216_a$}}\cg%
\e{0}%
\e{0}%
\e{0}%
\e{0}%
\e{0}%
\e{0}%
\e{0}%
\e{0}%
\e{0}%
\e{0}%
\e{0}%
\e{0}%
\e{0}%
\e{0}%
\e{0}%
\e{0}%
\e{1}%
\e{1}%
\e{0}%
\e{1}%
\e{0}%
\e{0}%
\e{0}%
\e{1}%
\e{0}%
\eol}\vss}\rg%
%
%
\rx{\vss\hfull{%
\rlx{\hss{$512_a$}}\cg%
\e{0}%
\e{0}%
\e{0}%
\e{0}%
\e{0}%
\e{0}%
\e{1}%
\e{0}%
\e{0}%
\e{1}%
\e{0}%
\e{0}%
\e{0}%
\e{0}%
\e{0}%
\e{0}%
\e{1}%
\e{0}%
\e{0}%
\e{1}%
\e{0}%
\e{1}%
\e{1}%
\e{1}%
\e{0}%
\eol}\vss}\rg%
%
%
\rx{\vss\hfull{%
\rlx{\hss{$378_a$}}\cg%
\e{0}%
\e{0}%
\e{0}%
\e{0}%
\e{0}%
\e{0}%
\e{0}%
\e{0}%
\e{0}%
\e{0}%
\e{0}%
\e{0}%
\e{0}%
\e{0}%
\e{0}%
\e{1}%
\e{0}%
\e{0}%
\e{1}%
\e{1}%
\e{0}%
\e{1}%
\e{1}%
\e{1}%
\e{0}%
\eol}\vss}\rg%
%
%
\rx{\vss\hfull{%
\rlx{\hss{$84_a$}}\cg%
\e{0}%
\e{0}%
\e{0}%
\e{0}%
\e{0}%
\e{0}%
\e{0}%
\e{0}%
\e{1}%
\e{0}%
\e{0}%
\e{0}%
\e{0}%
\e{0}%
\e{0}%
\e{0}%
\e{0}%
\e{0}%
\e{0}%
\e{0}%
\e{0}%
\e{0}%
\e{0}%
\e{1}%
\e{0}%
\eol}\vss}\rg%
%
%
\rx{\vss\hfull{%
\rlx{\hss{$420_a$}}\cg%
\e{0}%
\e{0}%
\e{0}%
\e{0}%
\e{0}%
\e{1}%
\e{1}%
\e{0}%
\e{1}%
\e{0}%
\e{0}%
\e{0}%
\e{0}%
\e{0}%
\e{0}%
\e{0}%
\e{1}%
\e{0}%
\e{0}%
\e{0}%
\e{0}%
\e{1}%
\e{1}%
\e{0}%
\e{0}%
\eol}\vss}\rg%
%
%
\rx{\vss\hfull{%
\rlx{\hss{$280_b$}}\cg%
\e{0}%
\e{0}%
\e{0}%
\e{0}%
\e{0}%
\e{1}%
\e{1}%
\e{1}%
\e{0}%
\e{1}%
\e{0}%
\e{0}%
\e{0}%
\e{0}%
\e{0}%
\e{0}%
\e{0}%
\e{0}%
\e{0}%
\e{0}%
\e{0}%
\e{0}%
\e{0}%
\e{1}%
\e{0}%
\eol}\vss}\rg%
%
%
\rx{\vss\hfull{%
\rlx{\hss{$210_b$}}\cg%
\e{0}%
\e{0}%
\e{0}%
\e{0}%
\e{0}%
\e{0}%
\e{0}%
\e{0}%
\e{0}%
\e{1}%
\e{0}%
\e{0}%
\e{0}%
\e{0}%
\e{0}%
\e{0}%
\e{0}%
\e{0}%
\e{0}%
\e{0}%
\e{0}%
\e{0}%
\e{1}%
\e{1}%
\e{1}%
\eol}\vss}\rg%
%
%
\rx{\vss\hfull{%
\rlx{\hss{$70_a$}}\cg%
\e{0}%
\e{0}%
\e{0}%
\e{0}%
\e{0}%
\e{0}%
\e{0}%
\e{0}%
\e{0}%
\e{0}%
\e{0}%
\e{0}%
\e{0}%
\e{0}%
\e{0}%
\e{0}%
\e{0}%
\e{0}%
\e{0}%
\e{1}%
\e{0}%
\e{0}%
\e{0}%
\e{0}%
\e{1}%
\eol}\vss}\rg%
\tableclose%
%
%
%
%
%
%
\eop
\eject
\tableopen{Induce/restrict matrix for $W(D_{8})\,\subset\,W(E_{8})$}%
%
%
%
%
%
%
\rowpts=18 true pt%
\colpts=18 true pt%
\rowlabpts=40 true pt%
\collabpts=70 true pt%
\clx{\vss\hfull{%
\rlx{\hss{$ $}}\cg%
\cx{\hskip 16 true pt\flip{$[{8}:-]$}\hss}\cg%
\cx{\hskip 16 true pt\flip{$[{7}{1}:-]$}\hss}\cg%
\cx{\hskip 16 true pt\flip{$[{6}{2}:-]$}\hss}\cg%
\cx{\hskip 16 true pt\flip{$[{6}{1^{2}}:-]$}\hss}\cg%
\cx{\hskip 16 true pt\flip{$[{5}{3}:-]$}\hss}\cg%
\cx{\hskip 16 true pt\flip{$[{5}{2}{1}:-]$}\hss}\cg%
\cx{\hskip 16 true pt\flip{$[{5}{1^{3}}:-]$}\hss}\cg%
\cx{\hskip 16 true pt\flip{$[{4^{2}}:-]$}\hss}\cg%
\cx{\hskip 16 true pt\flip{$[{4}{3}{1}:-]$}\hss}\cg%
\cx{\hskip 16 true pt\flip{$[{4}{2^{2}}:-]$}\hss}\cg%
\cx{\hskip 16 true pt\flip{$[{4}{2}{1^{2}}:-]$}\hss}\cg%
\cx{\hskip 16 true pt\flip{$[{4}{1^{4}}:-]$}\hss}\cg%
\cx{\hskip 16 true pt\flip{$[{3^{2}}{2}:-]$}\hss}\cg%
\cx{\hskip 16 true pt\flip{$[{3^{2}}{1^{2}}:-]$}\hss}\cg%
\cx{\hskip 16 true pt\flip{$[{3}{2^{2}}{1}:-]$}\hss}\cg%
\cx{\hskip 16 true pt\flip{$[{3}{2}{1^{3}}:-]$}\hss}\cg%
\cx{\hskip 16 true pt\flip{$[{3}{1^{5}}:-]$}\hss}\cg%
\cx{\hskip 16 true pt\flip{$[{2^{4}}:-]$}\hss}\cg%
\cx{\hskip 16 true pt\flip{$[{2^{3}}{1^{2}}:-]$}\hss}\cg%
\cx{\hskip 16 true pt\flip{$[{2^{2}}{1^{4}}:-]$}\hss}\cg%
\cx{\hskip 16 true pt\flip{$[{2}{1^{6}}:-]$}\hss}\cg%
\cx{\hskip 16 true pt\flip{$[{1^{8}}:-]$}\hss}\cg%
\cx{\hskip 16 true pt\flip{$[{7}:{1}]$}\hss}\cg%
\cx{\hskip 16 true pt\flip{$[{6}{1}:{1}]$}\hss}\cg%
\cx{\hskip 16 true pt\flip{$[{5}{2}:{1}]$}\hss}\cg%
\eol}}\rg%
%
%
\rx{\vss\hfull{%
\rlx{\hss{$1_{x}$}}\cg%
\e{1}%
\e{0}%
\e{0}%
\e{0}%
\e{0}%
\e{0}%
\e{0}%
\e{0}%
\e{0}%
\e{0}%
\e{0}%
\e{0}%
\e{0}%
\e{0}%
\e{0}%
\e{0}%
\e{0}%
\e{0}%
\e{0}%
\e{0}%
\e{0}%
\e{0}%
\e{0}%
\e{0}%
\e{0}%
\eol}\vss}\rg%
%
%
\rx{\vss\hfull{%
\rlx{\hss{$28_{x}$}}\cg%
\e{0}%
\e{0}%
\e{0}%
\e{0}%
\e{0}%
\e{0}%
\e{0}%
\e{0}%
\e{0}%
\e{0}%
\e{0}%
\e{0}%
\e{0}%
\e{0}%
\e{0}%
\e{0}%
\e{0}%
\e{0}%
\e{0}%
\e{0}%
\e{0}%
\e{0}%
\e{0}%
\e{0}%
\e{0}%
\eol}\vss}\rg%
%
%
\rx{\vss\hfull{%
\rlx{\hss{$35_{x}$}}\cg%
\e{0}%
\e{1}%
\e{0}%
\e{0}%
\e{0}%
\e{0}%
\e{0}%
\e{0}%
\e{0}%
\e{0}%
\e{0}%
\e{0}%
\e{0}%
\e{0}%
\e{0}%
\e{0}%
\e{0}%
\e{0}%
\e{0}%
\e{0}%
\e{0}%
\e{0}%
\e{0}%
\e{0}%
\e{0}%
\eol}\vss}\rg%
%
%
\rx{\vss\hfull{%
\rlx{\hss{$84_{x}$}}\cg%
\e{1}%
\e{0}%
\e{1}%
\e{0}%
\e{0}%
\e{0}%
\e{0}%
\e{0}%
\e{0}%
\e{0}%
\e{0}%
\e{0}%
\e{0}%
\e{0}%
\e{0}%
\e{0}%
\e{0}%
\e{0}%
\e{0}%
\e{0}%
\e{0}%
\e{0}%
\e{0}%
\e{0}%
\e{0}%
\eol}\vss}\rg%
%
%
\rx{\vss\hfull{%
\rlx{\hss{$50_{x}$}}\cg%
\e{1}%
\e{0}%
\e{0}%
\e{0}%
\e{0}%
\e{0}%
\e{0}%
\e{1}%
\e{0}%
\e{0}%
\e{0}%
\e{0}%
\e{0}%
\e{0}%
\e{0}%
\e{0}%
\e{0}%
\e{0}%
\e{0}%
\e{0}%
\e{0}%
\e{0}%
\e{0}%
\e{0}%
\e{0}%
\eol}\vss}\rg%
%
%
\rx{\vss\hfull{%
\rlx{\hss{$350_{x}$}}\cg%
\e{0}%
\e{0}%
\e{0}%
\e{0}%
\e{0}%
\e{0}%
\e{0}%
\e{0}%
\e{0}%
\e{0}%
\e{0}%
\e{0}%
\e{0}%
\e{0}%
\e{0}%
\e{0}%
\e{0}%
\e{0}%
\e{0}%
\e{0}%
\e{0}%
\e{0}%
\e{0}%
\e{0}%
\e{0}%
\eol}\vss}\rg%
%
%
\rx{\vss\hfull{%
\rlx{\hss{$300_{x}$}}\cg%
\e{0}%
\e{0}%
\e{1}%
\e{0}%
\e{0}%
\e{0}%
\e{0}%
\e{0}%
\e{0}%
\e{0}%
\e{0}%
\e{0}%
\e{0}%
\e{0}%
\e{0}%
\e{0}%
\e{0}%
\e{0}%
\e{0}%
\e{0}%
\e{0}%
\e{0}%
\e{0}%
\e{0}%
\e{0}%
\eol}\vss}\rg%
%
%
\rx{\vss\hfull{%
\rlx{\hss{$567_{x}$}}\cg%
\e{0}%
\e{0}%
\e{0}%
\e{1}%
\e{0}%
\e{0}%
\e{0}%
\e{0}%
\e{0}%
\e{0}%
\e{0}%
\e{0}%
\e{0}%
\e{0}%
\e{0}%
\e{0}%
\e{0}%
\e{0}%
\e{0}%
\e{0}%
\e{0}%
\e{0}%
\e{0}%
\e{0}%
\e{0}%
\eol}\vss}\rg%
%
%
\rx{\vss\hfull{%
\rlx{\hss{$210_{x}$}}\cg%
\e{0}%
\e{1}%
\e{0}%
\e{0}%
\e{0}%
\e{0}%
\e{0}%
\e{0}%
\e{0}%
\e{0}%
\e{0}%
\e{0}%
\e{0}%
\e{0}%
\e{0}%
\e{0}%
\e{0}%
\e{0}%
\e{0}%
\e{0}%
\e{0}%
\e{0}%
\e{0}%
\e{0}%
\e{0}%
\eol}\vss}\rg%
%
%
\rx{\vss\hfull{%
\rlx{\hss{$840_{x}$}}\cg%
\e{0}%
\e{0}%
\e{0}%
\e{0}%
\e{1}%
\e{0}%
\e{0}%
\e{0}%
\e{0}%
\e{0}%
\e{0}%
\e{0}%
\e{0}%
\e{0}%
\e{0}%
\e{0}%
\e{0}%
\e{0}%
\e{0}%
\e{0}%
\e{0}%
\e{0}%
\e{0}%
\e{0}%
\e{0}%
\eol}\vss}\rg%
%
%
\rx{\vss\hfull{%
\rlx{\hss{$700_{x}$}}\cg%
\e{0}%
\e{1}%
\e{0}%
\e{0}%
\e{1}%
\e{0}%
\e{0}%
\e{0}%
\e{0}%
\e{0}%
\e{0}%
\e{0}%
\e{0}%
\e{0}%
\e{0}%
\e{0}%
\e{0}%
\e{0}%
\e{0}%
\e{0}%
\e{0}%
\e{0}%
\e{0}%
\e{0}%
\e{0}%
\eol}\vss}\rg%
%
%
\rx{\vss\hfull{%
\rlx{\hss{$175_{x}$}}\cg%
\e{0}%
\e{0}%
\e{0}%
\e{0}%
\e{0}%
\e{0}%
\e{0}%
\e{0}%
\e{0}%
\e{0}%
\e{0}%
\e{0}%
\e{0}%
\e{0}%
\e{0}%
\e{0}%
\e{0}%
\e{0}%
\e{0}%
\e{0}%
\e{0}%
\e{0}%
\e{0}%
\e{0}%
\e{0}%
\eol}\vss}\rg%
%
%
\rx{\vss\hfull{%
\rlx{\hss{$1400_{x}$}}\cg%
\e{0}%
\e{0}%
\e{0}%
\e{0}%
\e{0}%
\e{0}%
\e{0}%
\e{0}%
\e{0}%
\e{0}%
\e{0}%
\e{0}%
\e{0}%
\e{0}%
\e{0}%
\e{0}%
\e{0}%
\e{0}%
\e{0}%
\e{0}%
\e{0}%
\e{0}%
\e{0}%
\e{0}%
\e{0}%
\eol}\vss}\rg%
%
%
\rx{\vss\hfull{%
\rlx{\hss{$1050_{x}$}}\cg%
\e{0}%
\e{0}%
\e{0}%
\e{0}%
\e{0}%
\e{0}%
\e{0}%
\e{0}%
\e{1}%
\e{0}%
\e{0}%
\e{0}%
\e{0}%
\e{0}%
\e{0}%
\e{0}%
\e{0}%
\e{0}%
\e{0}%
\e{0}%
\e{0}%
\e{0}%
\e{0}%
\e{0}%
\e{0}%
\eol}\vss}\rg%
%
%
\rx{\vss\hfull{%
\rlx{\hss{$1575_{x}$}}\cg%
\e{0}%
\e{0}%
\e{0}%
\e{0}%
\e{0}%
\e{0}%
\e{0}%
\e{0}%
\e{0}%
\e{0}%
\e{0}%
\e{0}%
\e{0}%
\e{0}%
\e{0}%
\e{0}%
\e{0}%
\e{0}%
\e{0}%
\e{0}%
\e{0}%
\e{0}%
\e{0}%
\e{0}%
\e{0}%
\eol}\vss}\rg%
%
%
\rx{\vss\hfull{%
\rlx{\hss{$1344_{x}$}}\cg%
\e{0}%
\e{0}%
\e{1}%
\e{0}%
\e{0}%
\e{1}%
\e{0}%
\e{0}%
\e{0}%
\e{0}%
\e{0}%
\e{0}%
\e{0}%
\e{0}%
\e{0}%
\e{0}%
\e{0}%
\e{0}%
\e{0}%
\e{0}%
\e{0}%
\e{0}%
\e{0}%
\e{0}%
\e{0}%
\eol}\vss}\rg%
%
%
\rx{\vss\hfull{%
\rlx{\hss{$2100_{x}$}}\cg%
\e{0}%
\e{0}%
\e{0}%
\e{0}%
\e{0}%
\e{0}%
\e{1}%
\e{0}%
\e{0}%
\e{0}%
\e{0}%
\e{0}%
\e{0}%
\e{0}%
\e{0}%
\e{0}%
\e{0}%
\e{0}%
\e{0}%
\e{0}%
\e{0}%
\e{0}%
\e{0}%
\e{0}%
\e{0}%
\eol}\vss}\rg%
%
%
\rx{\vss\hfull{%
\rlx{\hss{$2268_{x}$}}\cg%
\e{0}%
\e{0}%
\e{0}%
\e{1}%
\e{0}%
\e{0}%
\e{0}%
\e{0}%
\e{0}%
\e{0}%
\e{0}%
\e{0}%
\e{0}%
\e{0}%
\e{0}%
\e{0}%
\e{0}%
\e{0}%
\e{0}%
\e{0}%
\e{0}%
\e{0}%
\e{0}%
\e{0}%
\e{0}%
\eol}\vss}\rg%
%
%
\rx{\vss\hfull{%
\rlx{\hss{$525_{x}$}}\cg%
\e{0}%
\e{0}%
\e{0}%
\e{0}%
\e{0}%
\e{0}%
\e{1}%
\e{0}%
\e{0}%
\e{0}%
\e{0}%
\e{0}%
\e{0}%
\e{0}%
\e{0}%
\e{0}%
\e{0}%
\e{0}%
\e{0}%
\e{0}%
\e{0}%
\e{0}%
\e{0}%
\e{0}%
\e{0}%
\eol}\vss}\rg%
%
%
\rx{\vss\hfull{%
\rlx{\hss{$700_{xx}$}}\cg%
\e{0}%
\e{0}%
\e{0}%
\e{0}%
\e{0}%
\e{0}%
\e{0}%
\e{1}%
\e{0}%
\e{0}%
\e{0}%
\e{0}%
\e{0}%
\e{1}%
\e{0}%
\e{0}%
\e{0}%
\e{0}%
\e{0}%
\e{0}%
\e{0}%
\e{0}%
\e{0}%
\e{0}%
\e{0}%
\eol}\vss}\rg%
%
%
\rx{\vss\hfull{%
\rlx{\hss{$972_{x}$}}\cg%
\e{0}%
\e{0}%
\e{1}%
\e{0}%
\e{0}%
\e{0}%
\e{0}%
\e{1}%
\e{0}%
\e{1}%
\e{0}%
\e{0}%
\e{0}%
\e{0}%
\e{0}%
\e{0}%
\e{0}%
\e{0}%
\e{0}%
\e{0}%
\e{0}%
\e{0}%
\e{0}%
\e{0}%
\e{0}%
\eol}\vss}\rg%
%
%
\rx{\vss\hfull{%
\rlx{\hss{$4096_{x}$}}\cg%
\e{0}%
\e{0}%
\e{0}%
\e{0}%
\e{0}%
\e{1}%
\e{0}%
\e{0}%
\e{0}%
\e{0}%
\e{0}%
\e{0}%
\e{0}%
\e{0}%
\e{0}%
\e{0}%
\e{0}%
\e{0}%
\e{0}%
\e{0}%
\e{0}%
\e{0}%
\e{0}%
\e{0}%
\e{0}%
\eol}\vss}\rg%
%
%
\rx{\vss\hfull{%
\rlx{\hss{$4200_{x}$}}\cg%
\e{0}%
\e{0}%
\e{0}%
\e{0}%
\e{0}%
\e{0}%
\e{0}%
\e{0}%
\e{1}%
\e{0}%
\e{0}%
\e{0}%
\e{0}%
\e{0}%
\e{0}%
\e{0}%
\e{0}%
\e{0}%
\e{0}%
\e{0}%
\e{0}%
\e{0}%
\e{0}%
\e{0}%
\e{0}%
\eol}\vss}\rg%
%
%
\rx{\vss\hfull{%
\rlx{\hss{$2240_{x}$}}\cg%
\e{0}%
\e{0}%
\e{0}%
\e{0}%
\e{1}%
\e{0}%
\e{0}%
\e{0}%
\e{0}%
\e{0}%
\e{0}%
\e{0}%
\e{0}%
\e{0}%
\e{0}%
\e{0}%
\e{0}%
\e{0}%
\e{0}%
\e{0}%
\e{0}%
\e{0}%
\e{0}%
\e{0}%
\e{0}%
\eol}\vss}\rg%
%
%
\rx{\vss\hfull{%
\rlx{\hss{$2835_{x}$}}\cg%
\e{0}%
\e{0}%
\e{0}%
\e{0}%
\e{0}%
\e{0}%
\e{0}%
\e{0}%
\e{0}%
\e{0}%
\e{0}%
\e{0}%
\e{1}%
\e{0}%
\e{0}%
\e{0}%
\e{0}%
\e{0}%
\e{0}%
\e{0}%
\e{0}%
\e{0}%
\e{0}%
\e{0}%
\e{0}%
\eol}\vss}\rg%
%
%
\rx{\vss\hfull{%
\rlx{\hss{$6075_{x}$}}\cg%
\e{0}%
\e{0}%
\e{0}%
\e{0}%
\e{0}%
\e{0}%
\e{0}%
\e{0}%
\e{0}%
\e{0}%
\e{1}%
\e{0}%
\e{0}%
\e{0}%
\e{0}%
\e{0}%
\e{0}%
\e{0}%
\e{0}%
\e{0}%
\e{0}%
\e{0}%
\e{0}%
\e{0}%
\e{0}%
\eol}\vss}\rg%
%
%
\rx{\vss\hfull{%
\rlx{\hss{$3200_{x}$}}\cg%
\e{0}%
\e{0}%
\e{0}%
\e{0}%
\e{0}%
\e{1}%
\e{0}%
\e{0}%
\e{0}%
\e{1}%
\e{0}%
\e{0}%
\e{0}%
\e{0}%
\e{0}%
\e{0}%
\e{0}%
\e{0}%
\e{0}%
\e{0}%
\e{0}%
\e{0}%
\e{0}%
\e{0}%
\e{0}%
\eol}\vss}\rg%
%
%
\rx{\vss\hfull{%
\rlx{\hss{$70_{y}$}}\cg%
\e{0}%
\e{0}%
\e{0}%
\e{0}%
\e{0}%
\e{0}%
\e{0}%
\e{0}%
\e{0}%
\e{0}%
\e{0}%
\e{0}%
\e{0}%
\e{0}%
\e{0}%
\e{0}%
\e{0}%
\e{0}%
\e{0}%
\e{0}%
\e{0}%
\e{0}%
\e{0}%
\e{0}%
\e{0}%
\eol}\vss}\rg%
%
%
\rx{\vss\hfull{%
\rlx{\hss{$1134_{y}$}}\cg%
\e{0}%
\e{0}%
\e{0}%
\e{0}%
\e{0}%
\e{0}%
\e{0}%
\e{0}%
\e{0}%
\e{0}%
\e{0}%
\e{0}%
\e{0}%
\e{0}%
\e{0}%
\e{0}%
\e{0}%
\e{0}%
\e{0}%
\e{0}%
\e{0}%
\e{0}%
\e{0}%
\e{0}%
\e{0}%
\eol}\vss}\rg%
%
%
\rx{\vss\hfull{%
\rlx{\hss{$1680_{y}$}}\cg%
\e{0}%
\e{0}%
\e{0}%
\e{0}%
\e{0}%
\e{0}%
\e{0}%
\e{0}%
\e{0}%
\e{0}%
\e{0}%
\e{0}%
\e{0}%
\e{0}%
\e{0}%
\e{0}%
\e{0}%
\e{0}%
\e{0}%
\e{0}%
\e{0}%
\e{0}%
\e{0}%
\e{0}%
\e{0}%
\eol}\vss}\rg%
%
%
\rx{\vss\hfull{%
\rlx{\hss{$168_{y}$}}\cg%
\e{0}%
\e{0}%
\e{0}%
\e{0}%
\e{0}%
\e{0}%
\e{0}%
\e{1}%
\e{0}%
\e{0}%
\e{0}%
\e{0}%
\e{0}%
\e{0}%
\e{0}%
\e{0}%
\e{0}%
\e{1}%
\e{0}%
\e{0}%
\e{0}%
\e{0}%
\e{0}%
\e{0}%
\e{0}%
\eol}\vss}\rg%
%
%
\rx{\vss\hfull{%
\rlx{\hss{$420_{y}$}}\cg%
\e{0}%
\e{0}%
\e{0}%
\e{0}%
\e{0}%
\e{0}%
\e{0}%
\e{0}%
\e{0}%
\e{0}%
\e{0}%
\e{0}%
\e{0}%
\e{0}%
\e{0}%
\e{0}%
\e{0}%
\e{0}%
\e{0}%
\e{0}%
\e{0}%
\e{0}%
\e{0}%
\e{0}%
\e{0}%
\eol}\vss}\rg%
%
%
\rx{\vss\hfull{%
\rlx{\hss{$3150_{y}$}}\cg%
\e{0}%
\e{0}%
\e{0}%
\e{0}%
\e{0}%
\e{0}%
\e{0}%
\e{0}%
\e{0}%
\e{0}%
\e{0}%
\e{0}%
\e{0}%
\e{0}%
\e{0}%
\e{0}%
\e{0}%
\e{0}%
\e{0}%
\e{0}%
\e{0}%
\e{0}%
\e{0}%
\e{0}%
\e{0}%
\eol}\vss}\rg%
%
%
\rx{\vss\hfull{%
\rlx{\hss{$4200_{y}$}}\cg%
\e{0}%
\e{0}%
\e{0}%
\e{0}%
\e{0}%
\e{0}%
\e{0}%
\e{0}%
\e{1}%
\e{0}%
\e{0}%
\e{0}%
\e{0}%
\e{0}%
\e{1}%
\e{0}%
\e{0}%
\e{0}%
\e{0}%
\e{0}%
\e{0}%
\e{0}%
\e{0}%
\e{0}%
\e{0}%
\eol}\vss}\rg%
\eop
\eject
\tablecont%
%
%
%
%
%
%
\rowpts=18 true pt%
\colpts=18 true pt%
\rowlabpts=40 true pt%
\collabpts=70 true pt%
\clx{\vss\hfull{%
\rlx{\hss{$ $}}\cg%
\cx{\hskip 16 true pt\flip{$[{8}:-]$}\hss}\cg%
\cx{\hskip 16 true pt\flip{$[{7}{1}:-]$}\hss}\cg%
\cx{\hskip 16 true pt\flip{$[{6}{2}:-]$}\hss}\cg%
\cx{\hskip 16 true pt\flip{$[{6}{1^{2}}:-]$}\hss}\cg%
\cx{\hskip 16 true pt\flip{$[{5}{3}:-]$}\hss}\cg%
\cx{\hskip 16 true pt\flip{$[{5}{2}{1}:-]$}\hss}\cg%
\cx{\hskip 16 true pt\flip{$[{5}{1^{3}}:-]$}\hss}\cg%
\cx{\hskip 16 true pt\flip{$[{4^{2}}:-]$}\hss}\cg%
\cx{\hskip 16 true pt\flip{$[{4}{3}{1}:-]$}\hss}\cg%
\cx{\hskip 16 true pt\flip{$[{4}{2^{2}}:-]$}\hss}\cg%
\cx{\hskip 16 true pt\flip{$[{4}{2}{1^{2}}:-]$}\hss}\cg%
\cx{\hskip 16 true pt\flip{$[{4}{1^{4}}:-]$}\hss}\cg%
\cx{\hskip 16 true pt\flip{$[{3^{2}}{2}:-]$}\hss}\cg%
\cx{\hskip 16 true pt\flip{$[{3^{2}}{1^{2}}:-]$}\hss}\cg%
\cx{\hskip 16 true pt\flip{$[{3}{2^{2}}{1}:-]$}\hss}\cg%
\cx{\hskip 16 true pt\flip{$[{3}{2}{1^{3}}:-]$}\hss}\cg%
\cx{\hskip 16 true pt\flip{$[{3}{1^{5}}:-]$}\hss}\cg%
\cx{\hskip 16 true pt\flip{$[{2^{4}}:-]$}\hss}\cg%
\cx{\hskip 16 true pt\flip{$[{2^{3}}{1^{2}}:-]$}\hss}\cg%
\cx{\hskip 16 true pt\flip{$[{2^{2}}{1^{4}}:-]$}\hss}\cg%
\cx{\hskip 16 true pt\flip{$[{2}{1^{6}}:-]$}\hss}\cg%
\cx{\hskip 16 true pt\flip{$[{1^{8}}:-]$}\hss}\cg%
\cx{\hskip 16 true pt\flip{$[{7}:{1}]$}\hss}\cg%
\cx{\hskip 16 true pt\flip{$[{6}{1}:{1}]$}\hss}\cg%
\cx{\hskip 16 true pt\flip{$[{5}{2}:{1}]$}\hss}\cg%
\eol}}\rg%
%
%
\rx{\vss\hfull{%
\rlx{\hss{$2688_{y}$}}\cg%
\e{0}%
\e{0}%
\e{0}%
\e{0}%
\e{0}%
\e{0}%
\e{0}%
\e{0}%
\e{0}%
\e{1}%
\e{0}%
\e{0}%
\e{0}%
\e{1}%
\e{0}%
\e{0}%
\e{0}%
\e{0}%
\e{0}%
\e{0}%
\e{0}%
\e{0}%
\e{0}%
\e{0}%
\e{0}%
\eol}\vss}\rg%
%
%
\rx{\vss\hfull{%
\rlx{\hss{$2100_{y}$}}\cg%
\e{0}%
\e{0}%
\e{0}%
\e{0}%
\e{0}%
\e{0}%
\e{1}%
\e{0}%
\e{0}%
\e{0}%
\e{0}%
\e{1}%
\e{0}%
\e{0}%
\e{0}%
\e{0}%
\e{0}%
\e{0}%
\e{0}%
\e{0}%
\e{0}%
\e{0}%
\e{0}%
\e{0}%
\e{0}%
\eol}\vss}\rg%
%
%
\rx{\vss\hfull{%
\rlx{\hss{$1400_{y}$}}\cg%
\e{0}%
\e{0}%
\e{0}%
\e{0}%
\e{0}%
\e{0}%
\e{0}%
\e{0}%
\e{0}%
\e{0}%
\e{0}%
\e{0}%
\e{0}%
\e{0}%
\e{0}%
\e{0}%
\e{0}%
\e{0}%
\e{0}%
\e{0}%
\e{0}%
\e{0}%
\e{0}%
\e{0}%
\e{0}%
\eol}\vss}\rg%
%
%
\rx{\vss\hfull{%
\rlx{\hss{$4536_{y}$}}\cg%
\e{0}%
\e{0}%
\e{0}%
\e{0}%
\e{0}%
\e{0}%
\e{0}%
\e{0}%
\e{0}%
\e{0}%
\e{0}%
\e{0}%
\e{0}%
\e{0}%
\e{0}%
\e{0}%
\e{0}%
\e{0}%
\e{0}%
\e{0}%
\e{0}%
\e{0}%
\e{0}%
\e{0}%
\e{0}%
\eol}\vss}\rg%
%
%
\rx{\vss\hfull{%
\rlx{\hss{$5670_{y}$}}\cg%
\e{0}%
\e{0}%
\e{0}%
\e{0}%
\e{0}%
\e{0}%
\e{0}%
\e{0}%
\e{0}%
\e{0}%
\e{0}%
\e{0}%
\e{0}%
\e{0}%
\e{0}%
\e{0}%
\e{0}%
\e{0}%
\e{0}%
\e{0}%
\e{0}%
\e{0}%
\e{0}%
\e{0}%
\e{0}%
\eol}\vss}\rg%
%
%
\rx{\vss\hfull{%
\rlx{\hss{$4480_{y}$}}\cg%
\e{0}%
\e{0}%
\e{0}%
\e{0}%
\e{0}%
\e{0}%
\e{0}%
\e{0}%
\e{0}%
\e{0}%
\e{0}%
\e{0}%
\e{0}%
\e{0}%
\e{0}%
\e{0}%
\e{0}%
\e{0}%
\e{0}%
\e{0}%
\e{0}%
\e{0}%
\e{0}%
\e{0}%
\e{0}%
\eol}\vss}\rg%
%
%
\rx{\vss\hfull{%
\rlx{\hss{$8_{z}$}}\cg%
\e{0}%
\e{0}%
\e{0}%
\e{0}%
\e{0}%
\e{0}%
\e{0}%
\e{0}%
\e{0}%
\e{0}%
\e{0}%
\e{0}%
\e{0}%
\e{0}%
\e{0}%
\e{0}%
\e{0}%
\e{0}%
\e{0}%
\e{0}%
\e{0}%
\e{0}%
\e{1}%
\e{0}%
\e{0}%
\eol}\vss}\rg%
%
%
\rx{\vss\hfull{%
\rlx{\hss{$56_{z}$}}\cg%
\e{0}%
\e{0}%
\e{0}%
\e{0}%
\e{0}%
\e{0}%
\e{0}%
\e{0}%
\e{0}%
\e{0}%
\e{0}%
\e{0}%
\e{0}%
\e{0}%
\e{0}%
\e{0}%
\e{0}%
\e{0}%
\e{0}%
\e{0}%
\e{0}%
\e{0}%
\e{0}%
\e{0}%
\e{0}%
\eol}\vss}\rg%
%
%
\rx{\vss\hfull{%
\rlx{\hss{$160_{z}$}}\cg%
\e{0}%
\e{0}%
\e{0}%
\e{0}%
\e{0}%
\e{0}%
\e{0}%
\e{0}%
\e{0}%
\e{0}%
\e{0}%
\e{0}%
\e{0}%
\e{0}%
\e{0}%
\e{0}%
\e{0}%
\e{0}%
\e{0}%
\e{0}%
\e{0}%
\e{0}%
\e{0}%
\e{1}%
\e{0}%
\eol}\vss}\rg%
%
%
\rx{\vss\hfull{%
\rlx{\hss{$112_{z}$}}\cg%
\e{0}%
\e{0}%
\e{0}%
\e{0}%
\e{0}%
\e{0}%
\e{0}%
\e{0}%
\e{0}%
\e{0}%
\e{0}%
\e{0}%
\e{0}%
\e{0}%
\e{0}%
\e{0}%
\e{0}%
\e{0}%
\e{0}%
\e{0}%
\e{0}%
\e{0}%
\e{1}%
\e{1}%
\e{0}%
\eol}\vss}\rg%
%
%
\rx{\vss\hfull{%
\rlx{\hss{$840_{z}$}}\cg%
\e{0}%
\e{0}%
\e{0}%
\e{0}%
\e{0}%
\e{0}%
\e{0}%
\e{0}%
\e{0}%
\e{0}%
\e{0}%
\e{0}%
\e{0}%
\e{0}%
\e{0}%
\e{0}%
\e{0}%
\e{0}%
\e{0}%
\e{0}%
\e{0}%
\e{0}%
\e{0}%
\e{0}%
\e{1}%
\eol}\vss}\rg%
%
%
\rx{\vss\hfull{%
\rlx{\hss{$1296_{z}$}}\cg%
\e{0}%
\e{0}%
\e{0}%
\e{0}%
\e{0}%
\e{0}%
\e{0}%
\e{0}%
\e{0}%
\e{0}%
\e{0}%
\e{0}%
\e{0}%
\e{0}%
\e{0}%
\e{0}%
\e{0}%
\e{0}%
\e{0}%
\e{0}%
\e{0}%
\e{0}%
\e{0}%
\e{0}%
\e{0}%
\eol}\vss}\rg%
%
%
\rx{\vss\hfull{%
\rlx{\hss{$1400_{z}$}}\cg%
\e{0}%
\e{0}%
\e{0}%
\e{0}%
\e{0}%
\e{0}%
\e{0}%
\e{0}%
\e{0}%
\e{0}%
\e{0}%
\e{0}%
\e{0}%
\e{0}%
\e{0}%
\e{0}%
\e{0}%
\e{0}%
\e{0}%
\e{0}%
\e{0}%
\e{0}%
\e{0}%
\e{1}%
\e{1}%
\eol}\vss}\rg%
%
%
\rx{\vss\hfull{%
\rlx{\hss{$1008_{z}$}}\cg%
\e{0}%
\e{0}%
\e{0}%
\e{0}%
\e{0}%
\e{0}%
\e{0}%
\e{0}%
\e{0}%
\e{0}%
\e{0}%
\e{0}%
\e{0}%
\e{0}%
\e{0}%
\e{0}%
\e{0}%
\e{0}%
\e{0}%
\e{0}%
\e{0}%
\e{0}%
\e{0}%
\e{1}%
\e{0}%
\eol}\vss}\rg%
%
%
\rx{\vss\hfull{%
\rlx{\hss{$560_{z}$}}\cg%
\e{0}%
\e{0}%
\e{0}%
\e{0}%
\e{0}%
\e{0}%
\e{0}%
\e{0}%
\e{0}%
\e{0}%
\e{0}%
\e{0}%
\e{0}%
\e{0}%
\e{0}%
\e{0}%
\e{0}%
\e{0}%
\e{0}%
\e{0}%
\e{0}%
\e{0}%
\e{1}%
\e{1}%
\e{1}%
\eol}\vss}\rg%
%
%
\rx{\vss\hfull{%
\rlx{\hss{$1400_{zz}$}}\cg%
\e{0}%
\e{0}%
\e{0}%
\e{0}%
\e{0}%
\e{0}%
\e{0}%
\e{0}%
\e{0}%
\e{0}%
\e{0}%
\e{0}%
\e{0}%
\e{0}%
\e{0}%
\e{0}%
\e{0}%
\e{0}%
\e{0}%
\e{0}%
\e{0}%
\e{0}%
\e{0}%
\e{0}%
\e{0}%
\eol}\vss}\rg%
%
%
\rx{\vss\hfull{%
\rlx{\hss{$4200_{z}$}}\cg%
\e{0}%
\e{0}%
\e{0}%
\e{0}%
\e{0}%
\e{0}%
\e{0}%
\e{0}%
\e{0}%
\e{0}%
\e{0}%
\e{0}%
\e{0}%
\e{0}%
\e{0}%
\e{0}%
\e{0}%
\e{0}%
\e{0}%
\e{0}%
\e{0}%
\e{0}%
\e{0}%
\e{0}%
\e{0}%
\eol}\vss}\rg%
%
%
\rx{\vss\hfull{%
\rlx{\hss{$400_{z}$}}\cg%
\e{0}%
\e{0}%
\e{0}%
\e{0}%
\e{0}%
\e{0}%
\e{0}%
\e{0}%
\e{0}%
\e{0}%
\e{0}%
\e{0}%
\e{0}%
\e{0}%
\e{0}%
\e{0}%
\e{0}%
\e{0}%
\e{0}%
\e{0}%
\e{0}%
\e{0}%
\e{1}%
\e{0}%
\e{0}%
\eol}\vss}\rg%
%
%
\rx{\vss\hfull{%
\rlx{\hss{$3240_{z}$}}\cg%
\e{0}%
\e{0}%
\e{0}%
\e{0}%
\e{0}%
\e{0}%
\e{0}%
\e{0}%
\e{0}%
\e{0}%
\e{0}%
\e{0}%
\e{0}%
\e{0}%
\e{0}%
\e{0}%
\e{0}%
\e{0}%
\e{0}%
\e{0}%
\e{0}%
\e{0}%
\e{0}%
\e{1}%
\e{1}%
\eol}\vss}\rg%
%
%
\rx{\vss\hfull{%
\rlx{\hss{$4536_{z}$}}\cg%
\e{0}%
\e{0}%
\e{0}%
\e{0}%
\e{0}%
\e{0}%
\e{0}%
\e{0}%
\e{0}%
\e{0}%
\e{0}%
\e{0}%
\e{0}%
\e{0}%
\e{0}%
\e{0}%
\e{0}%
\e{0}%
\e{0}%
\e{0}%
\e{0}%
\e{0}%
\e{0}%
\e{0}%
\e{1}%
\eol}\vss}\rg%
%
%
\rx{\vss\hfull{%
\rlx{\hss{$2400_{z}$}}\cg%
\e{0}%
\e{0}%
\e{0}%
\e{0}%
\e{0}%
\e{0}%
\e{0}%
\e{0}%
\e{0}%
\e{0}%
\e{0}%
\e{0}%
\e{0}%
\e{0}%
\e{0}%
\e{0}%
\e{0}%
\e{0}%
\e{0}%
\e{0}%
\e{0}%
\e{0}%
\e{0}%
\e{0}%
\e{0}%
\eol}\vss}\rg%
%
%
\rx{\vss\hfull{%
\rlx{\hss{$3360_{z}$}}\cg%
\e{0}%
\e{0}%
\e{0}%
\e{0}%
\e{0}%
\e{0}%
\e{0}%
\e{0}%
\e{0}%
\e{0}%
\e{0}%
\e{0}%
\e{0}%
\e{0}%
\e{0}%
\e{0}%
\e{0}%
\e{0}%
\e{0}%
\e{0}%
\e{0}%
\e{0}%
\e{0}%
\e{0}%
\e{0}%
\eol}\vss}\rg%
%
%
\rx{\vss\hfull{%
\rlx{\hss{$2800_{z}$}}\cg%
\e{0}%
\e{0}%
\e{0}%
\e{0}%
\e{0}%
\e{0}%
\e{0}%
\e{0}%
\e{0}%
\e{0}%
\e{0}%
\e{0}%
\e{0}%
\e{0}%
\e{0}%
\e{0}%
\e{0}%
\e{0}%
\e{0}%
\e{0}%
\e{0}%
\e{0}%
\e{0}%
\e{0}%
\e{0}%
\eol}\vss}\rg%
%
%
\rx{\vss\hfull{%
\rlx{\hss{$4096_{z}$}}\cg%
\e{0}%
\e{0}%
\e{0}%
\e{0}%
\e{0}%
\e{0}%
\e{0}%
\e{0}%
\e{0}%
\e{0}%
\e{0}%
\e{0}%
\e{0}%
\e{0}%
\e{0}%
\e{0}%
\e{0}%
\e{0}%
\e{0}%
\e{0}%
\e{0}%
\e{0}%
\e{0}%
\e{0}%
\e{1}%
\eol}\vss}\rg%
%
%
\rx{\vss\hfull{%
\rlx{\hss{$5600_{z}$}}\cg%
\e{0}%
\e{0}%
\e{0}%
\e{0}%
\e{0}%
\e{0}%
\e{0}%
\e{0}%
\e{0}%
\e{0}%
\e{0}%
\e{0}%
\e{0}%
\e{0}%
\e{0}%
\e{0}%
\e{0}%
\e{0}%
\e{0}%
\e{0}%
\e{0}%
\e{0}%
\e{0}%
\e{0}%
\e{0}%
\eol}\vss}\rg%
%
%
\rx{\vss\hfull{%
\rlx{\hss{$448_{z}$}}\cg%
\e{0}%
\e{0}%
\e{0}%
\e{0}%
\e{0}%
\e{0}%
\e{0}%
\e{0}%
\e{0}%
\e{0}%
\e{0}%
\e{0}%
\e{0}%
\e{0}%
\e{0}%
\e{0}%
\e{0}%
\e{0}%
\e{0}%
\e{0}%
\e{0}%
\e{0}%
\e{0}%
\e{0}%
\e{1}%
\eol}\vss}\rg%
%
%
\rx{\vss\hfull{%
\rlx{\hss{$448_{w}$}}\cg%
\e{0}%
\e{0}%
\e{0}%
\e{0}%
\e{0}%
\e{0}%
\e{0}%
\e{0}%
\e{0}%
\e{0}%
\e{0}%
\e{0}%
\e{0}%
\e{0}%
\e{0}%
\e{0}%
\e{0}%
\e{0}%
\e{0}%
\e{0}%
\e{0}%
\e{0}%
\e{0}%
\e{0}%
\e{0}%
\eol}\vss}\rg%
%
%
\rx{\vss\hfull{%
\rlx{\hss{$1344_{w}$}}\cg%
\e{0}%
\e{0}%
\e{0}%
\e{0}%
\e{0}%
\e{0}%
\e{0}%
\e{0}%
\e{0}%
\e{0}%
\e{0}%
\e{0}%
\e{0}%
\e{0}%
\e{0}%
\e{0}%
\e{0}%
\e{0}%
\e{0}%
\e{0}%
\e{0}%
\e{0}%
\e{0}%
\e{0}%
\e{0}%
\eol}\vss}\rg%
%
%
\rx{\vss\hfull{%
\rlx{\hss{$5600_{w}$}}\cg%
\e{0}%
\e{0}%
\e{0}%
\e{0}%
\e{0}%
\e{0}%
\e{0}%
\e{0}%
\e{0}%
\e{0}%
\e{0}%
\e{0}%
\e{0}%
\e{0}%
\e{0}%
\e{0}%
\e{0}%
\e{0}%
\e{0}%
\e{0}%
\e{0}%
\e{0}%
\e{0}%
\e{0}%
\e{0}%
\eol}\vss}\rg%
%
%
\rx{\vss\hfull{%
\rlx{\hss{$2016_{w}$}}\cg%
\e{0}%
\e{0}%
\e{0}%
\e{0}%
\e{0}%
\e{0}%
\e{0}%
\e{0}%
\e{0}%
\e{0}%
\e{0}%
\e{0}%
\e{0}%
\e{0}%
\e{0}%
\e{0}%
\e{0}%
\e{0}%
\e{0}%
\e{0}%
\e{0}%
\e{0}%
\e{0}%
\e{0}%
\e{0}%
\eol}\vss}\rg%
%
%
\rx{\vss\hfull{%
\rlx{\hss{$7168_{w}$}}\cg%
\e{0}%
\e{0}%
\e{0}%
\e{0}%
\e{0}%
\e{0}%
\e{0}%
\e{0}%
\e{0}%
\e{0}%
\e{0}%
\e{0}%
\e{0}%
\e{0}%
\e{0}%
\e{0}%
\e{0}%
\e{0}%
\e{0}%
\e{0}%
\e{0}%
\e{0}%
\e{0}%
\e{0}%
\e{0}%
\eol}\vss}\rg%
\eop
\eject
\tablecont%
%
%
%
%
%
%
\rowpts=18 true pt%
\colpts=18 true pt%
\rowlabpts=40 true pt%
\collabpts=70 true pt%
\clx{\vss\hfull{%
\rlx{\hss{$ $}}\cg%
\cx{\hskip 16 true pt\flip{$[{5}{1^{2}}:{1}]$}\hss}\cg%
\cx{\hskip 16 true pt\flip{$[{4}{3}:{1}]$}\hss}\cg%
\cx{\hskip 16 true pt\flip{$[{4}{2}{1}:{1}]$}\hss}\cg%
\cx{\hskip 16 true pt\flip{$[{4}{1^{3}}:{1}]$}\hss}\cg%
\cx{\hskip 16 true pt\flip{$[{3^{2}}{1}:{1}]$}\hss}\cg%
\cx{\hskip 16 true pt\flip{$[{3}{2^{2}}:{1}]$}\hss}\cg%
\cx{\hskip 16 true pt\flip{$[{3}{2}{1^{2}}:{1}]$}\hss}\cg%
\cx{\hskip 16 true pt\flip{$[{3}{1^{4}}:{1}]$}\hss}\cg%
\cx{\hskip 16 true pt\flip{$[{2^{3}}{1}:{1}]$}\hss}\cg%
\cx{\hskip 16 true pt\flip{$[{2^{2}}{1^{3}}:{1}]$}\hss}\cg%
\cx{\hskip 16 true pt\flip{$[{2}{1^{5}}:{1}]$}\hss}\cg%
\cx{\hskip 16 true pt\flip{$[{1^{7}}:{1}]$}\hss}\cg%
\cx{\hskip 16 true pt\flip{$[{6}:{2}]$}\hss}\cg%
\cx{\hskip 16 true pt\flip{$[{6}:{1^{2}}]$}\hss}\cg%
\cx{\hskip 16 true pt\flip{$[{5}{1}:{2}]$}\hss}\cg%
\cx{\hskip 16 true pt\flip{$[{5}{1}:{1^{2}}]$}\hss}\cg%
\cx{\hskip 16 true pt\flip{$[{4}{2}:{2}]$}\hss}\cg%
\cx{\hskip 16 true pt\flip{$[{4}{2}:{1^{2}}]$}\hss}\cg%
\cx{\hskip 16 true pt\flip{$[{4}{1^{2}}:{2}]$}\hss}\cg%
\cx{\hskip 16 true pt\flip{$[{4}{1^{2}}:{1^{2}}]$}\hss}\cg%
\cx{\hskip 16 true pt\flip{$[{3^{2}}:{2}]$}\hss}\cg%
\cx{\hskip 16 true pt\flip{$[{3^{2}}:{1^{2}}]$}\hss}\cg%
\cx{\hskip 16 true pt\flip{$[{3}{2}{1}:{2}]$}\hss}\cg%
\cx{\hskip 16 true pt\flip{$[{3}{2}{1}:{1^{2}}]$}\hss}\cg%
\cx{\hskip 16 true pt\flip{$[{3}{1^{3}}:{2}]$}\hss}\cg%
\eol}}\rg%
%
%
\rx{\vss\hfull{%
\rlx{\hss{$1_{x}$}}\cg%
\e{0}%
\e{0}%
\e{0}%
\e{0}%
\e{0}%
\e{0}%
\e{0}%
\e{0}%
\e{0}%
\e{0}%
\e{0}%
\e{0}%
\e{0}%
\e{0}%
\e{0}%
\e{0}%
\e{0}%
\e{0}%
\e{0}%
\e{0}%
\e{0}%
\e{0}%
\e{0}%
\e{0}%
\e{0}%
\eol}\vss}\rg%
%
%
\rx{\vss\hfull{%
\rlx{\hss{$28_{x}$}}\cg%
\e{0}%
\e{0}%
\e{0}%
\e{0}%
\e{0}%
\e{0}%
\e{0}%
\e{0}%
\e{0}%
\e{0}%
\e{0}%
\e{0}%
\e{0}%
\e{1}%
\e{0}%
\e{0}%
\e{0}%
\e{0}%
\e{0}%
\e{0}%
\e{0}%
\e{0}%
\e{0}%
\e{0}%
\e{0}%
\eol}\vss}\rg%
%
%
\rx{\vss\hfull{%
\rlx{\hss{$35_{x}$}}\cg%
\e{0}%
\e{0}%
\e{0}%
\e{0}%
\e{0}%
\e{0}%
\e{0}%
\e{0}%
\e{0}%
\e{0}%
\e{0}%
\e{0}%
\e{1}%
\e{0}%
\e{0}%
\e{0}%
\e{0}%
\e{0}%
\e{0}%
\e{0}%
\e{0}%
\e{0}%
\e{0}%
\e{0}%
\e{0}%
\eol}\vss}\rg%
%
%
\rx{\vss\hfull{%
\rlx{\hss{$84_{x}$}}\cg%
\e{0}%
\e{0}%
\e{0}%
\e{0}%
\e{0}%
\e{0}%
\e{0}%
\e{0}%
\e{0}%
\e{0}%
\e{0}%
\e{0}%
\e{1}%
\e{0}%
\e{0}%
\e{0}%
\e{0}%
\e{0}%
\e{0}%
\e{0}%
\e{0}%
\e{0}%
\e{0}%
\e{0}%
\e{0}%
\eol}\vss}\rg%
%
%
\rx{\vss\hfull{%
\rlx{\hss{$50_{x}$}}\cg%
\e{0}%
\e{0}%
\e{0}%
\e{0}%
\e{0}%
\e{0}%
\e{0}%
\e{0}%
\e{0}%
\e{0}%
\e{0}%
\e{0}%
\e{0}%
\e{0}%
\e{0}%
\e{0}%
\e{0}%
\e{0}%
\e{0}%
\e{0}%
\e{0}%
\e{0}%
\e{0}%
\e{0}%
\e{0}%
\eol}\vss}\rg%
%
%
\rx{\vss\hfull{%
\rlx{\hss{$350_{x}$}}\cg%
\e{0}%
\e{0}%
\e{0}%
\e{0}%
\e{0}%
\e{0}%
\e{0}%
\e{0}%
\e{0}%
\e{0}%
\e{0}%
\e{0}%
\e{0}%
\e{0}%
\e{0}%
\e{1}%
\e{0}%
\e{0}%
\e{0}%
\e{0}%
\e{0}%
\e{0}%
\e{0}%
\e{0}%
\e{0}%
\eol}\vss}\rg%
%
%
\rx{\vss\hfull{%
\rlx{\hss{$300_{x}$}}\cg%
\e{0}%
\e{0}%
\e{0}%
\e{0}%
\e{0}%
\e{0}%
\e{0}%
\e{0}%
\e{0}%
\e{0}%
\e{0}%
\e{0}%
\e{0}%
\e{0}%
\e{1}%
\e{0}%
\e{0}%
\e{0}%
\e{0}%
\e{0}%
\e{0}%
\e{0}%
\e{0}%
\e{0}%
\e{0}%
\eol}\vss}\rg%
%
%
\rx{\vss\hfull{%
\rlx{\hss{$567_{x}$}}\cg%
\e{0}%
\e{0}%
\e{0}%
\e{0}%
\e{0}%
\e{0}%
\e{0}%
\e{0}%
\e{0}%
\e{0}%
\e{0}%
\e{0}%
\e{1}%
\e{1}%
\e{1}%
\e{1}%
\e{0}%
\e{0}%
\e{0}%
\e{0}%
\e{0}%
\e{0}%
\e{0}%
\e{0}%
\e{0}%
\eol}\vss}\rg%
%
%
\rx{\vss\hfull{%
\rlx{\hss{$210_{x}$}}\cg%
\e{0}%
\e{0}%
\e{0}%
\e{0}%
\e{0}%
\e{0}%
\e{0}%
\e{0}%
\e{0}%
\e{0}%
\e{0}%
\e{0}%
\e{0}%
\e{1}%
\e{1}%
\e{0}%
\e{0}%
\e{0}%
\e{0}%
\e{0}%
\e{0}%
\e{0}%
\e{0}%
\e{0}%
\e{0}%
\eol}\vss}\rg%
%
%
\rx{\vss\hfull{%
\rlx{\hss{$840_{x}$}}\cg%
\e{0}%
\e{0}%
\e{0}%
\e{0}%
\e{0}%
\e{0}%
\e{0}%
\e{0}%
\e{0}%
\e{0}%
\e{0}%
\e{0}%
\e{0}%
\e{0}%
\e{0}%
\e{0}%
\e{1}%
\e{0}%
\e{0}%
\e{0}%
\e{0}%
\e{0}%
\e{0}%
\e{0}%
\e{0}%
\eol}\vss}\rg%
%
%
\rx{\vss\hfull{%
\rlx{\hss{$700_{x}$}}\cg%
\e{0}%
\e{0}%
\e{0}%
\e{0}%
\e{0}%
\e{0}%
\e{0}%
\e{0}%
\e{0}%
\e{0}%
\e{0}%
\e{0}%
\e{1}%
\e{0}%
\e{1}%
\e{0}%
\e{1}%
\e{0}%
\e{0}%
\e{0}%
\e{0}%
\e{0}%
\e{0}%
\e{0}%
\e{0}%
\eol}\vss}\rg%
%
%
\rx{\vss\hfull{%
\rlx{\hss{$175_{x}$}}\cg%
\e{0}%
\e{0}%
\e{0}%
\e{0}%
\e{0}%
\e{0}%
\e{0}%
\e{0}%
\e{0}%
\e{0}%
\e{0}%
\e{0}%
\e{0}%
\e{0}%
\e{0}%
\e{0}%
\e{0}%
\e{0}%
\e{0}%
\e{0}%
\e{1}%
\e{0}%
\e{0}%
\e{0}%
\e{0}%
\eol}\vss}\rg%
%
%
\rx{\vss\hfull{%
\rlx{\hss{$1400_{x}$}}\cg%
\e{0}%
\e{0}%
\e{0}%
\e{0}%
\e{0}%
\e{0}%
\e{0}%
\e{0}%
\e{0}%
\e{0}%
\e{0}%
\e{0}%
\e{0}%
\e{1}%
\e{1}%
\e{0}%
\e{0}%
\e{1}%
\e{1}%
\e{0}%
\e{1}%
\e{0}%
\e{0}%
\e{0}%
\e{0}%
\eol}\vss}\rg%
%
%
\rx{\vss\hfull{%
\rlx{\hss{$1050_{x}$}}\cg%
\e{0}%
\e{0}%
\e{0}%
\e{0}%
\e{0}%
\e{0}%
\e{0}%
\e{0}%
\e{0}%
\e{0}%
\e{0}%
\e{0}%
\e{1}%
\e{0}%
\e{0}%
\e{0}%
\e{1}%
\e{0}%
\e{0}%
\e{0}%
\e{0}%
\e{1}%
\e{0}%
\e{0}%
\e{0}%
\eol}\vss}\rg%
%
%
\rx{\vss\hfull{%
\rlx{\hss{$1575_{x}$}}\cg%
\e{0}%
\e{0}%
\e{0}%
\e{0}%
\e{0}%
\e{0}%
\e{0}%
\e{0}%
\e{0}%
\e{0}%
\e{0}%
\e{0}%
\e{0}%
\e{1}%
\e{1}%
\e{1}%
\e{0}%
\e{1}%
\e{1}%
\e{0}%
\e{0}%
\e{0}%
\e{0}%
\e{0}%
\e{0}%
\eol}\vss}\rg%
%
%
\rx{\vss\hfull{%
\rlx{\hss{$1344_{x}$}}\cg%
\e{0}%
\e{0}%
\e{0}%
\e{0}%
\e{0}%
\e{0}%
\e{0}%
\e{0}%
\e{0}%
\e{0}%
\e{0}%
\e{0}%
\e{1}%
\e{0}%
\e{1}%
\e{1}%
\e{1}%
\e{0}%
\e{0}%
\e{0}%
\e{0}%
\e{0}%
\e{0}%
\e{0}%
\e{0}%
\eol}\vss}\rg%
%
%
\rx{\vss\hfull{%
\rlx{\hss{$2100_{x}$}}\cg%
\e{0}%
\e{0}%
\e{0}%
\e{0}%
\e{0}%
\e{0}%
\e{0}%
\e{0}%
\e{0}%
\e{0}%
\e{0}%
\e{0}%
\e{0}%
\e{0}%
\e{0}%
\e{1}%
\e{0}%
\e{0}%
\e{1}%
\e{1}%
\e{0}%
\e{0}%
\e{0}%
\e{0}%
\e{0}%
\eol}\vss}\rg%
%
%
\rx{\vss\hfull{%
\rlx{\hss{$2268_{x}$}}\cg%
\e{0}%
\e{0}%
\e{0}%
\e{0}%
\e{0}%
\e{0}%
\e{0}%
\e{0}%
\e{0}%
\e{0}%
\e{0}%
\e{0}%
\e{0}%
\e{0}%
\e{1}%
\e{1}%
\e{1}%
\e{0}%
\e{1}%
\e{1}%
\e{0}%
\e{0}%
\e{0}%
\e{0}%
\e{0}%
\eol}\vss}\rg%
%
%
\rx{\vss\hfull{%
\rlx{\hss{$525_{x}$}}\cg%
\e{0}%
\e{0}%
\e{0}%
\e{0}%
\e{0}%
\e{0}%
\e{0}%
\e{0}%
\e{0}%
\e{0}%
\e{0}%
\e{0}%
\e{0}%
\e{0}%
\e{0}%
\e{1}%
\e{0}%
\e{0}%
\e{0}%
\e{0}%
\e{0}%
\e{0}%
\e{0}%
\e{0}%
\e{0}%
\eol}\vss}\rg%
%
%
\rx{\vss\hfull{%
\rlx{\hss{$700_{xx}$}}\cg%
\e{0}%
\e{0}%
\e{0}%
\e{0}%
\e{0}%
\e{0}%
\e{0}%
\e{0}%
\e{0}%
\e{0}%
\e{0}%
\e{0}%
\e{0}%
\e{0}%
\e{0}%
\e{0}%
\e{0}%
\e{0}%
\e{0}%
\e{0}%
\e{0}%
\e{1}%
\e{0}%
\e{0}%
\e{0}%
\eol}\vss}\rg%
%
%
\rx{\vss\hfull{%
\rlx{\hss{$972_{x}$}}\cg%
\e{0}%
\e{0}%
\e{0}%
\e{0}%
\e{0}%
\e{0}%
\e{0}%
\e{0}%
\e{0}%
\e{0}%
\e{0}%
\e{0}%
\e{0}%
\e{0}%
\e{0}%
\e{0}%
\e{1}%
\e{0}%
\e{0}%
\e{0}%
\e{0}%
\e{0}%
\e{0}%
\e{0}%
\e{0}%
\eol}\vss}\rg%
%
%
\rx{\vss\hfull{%
\rlx{\hss{$4096_{x}$}}\cg%
\e{0}%
\e{0}%
\e{0}%
\e{0}%
\e{0}%
\e{0}%
\e{0}%
\e{0}%
\e{0}%
\e{0}%
\e{0}%
\e{0}%
\e{0}%
\e{0}%
\e{1}%
\e{1}%
\e{1}%
\e{1}%
\e{1}%
\e{1}%
\e{0}%
\e{0}%
\e{1}%
\e{0}%
\e{0}%
\eol}\vss}\rg%
%
%
\rx{\vss\hfull{%
\rlx{\hss{$4200_{x}$}}\cg%
\e{0}%
\e{0}%
\e{0}%
\e{0}%
\e{0}%
\e{0}%
\e{0}%
\e{0}%
\e{0}%
\e{0}%
\e{0}%
\e{0}%
\e{0}%
\e{0}%
\e{1}%
\e{0}%
\e{1}%
\e{1}%
\e{1}%
\e{0}%
\e{1}%
\e{0}%
\e{1}%
\e{1}%
\e{0}%
\eol}\vss}\rg%
%
%
\rx{\vss\hfull{%
\rlx{\hss{$2240_{x}$}}\cg%
\e{0}%
\e{0}%
\e{0}%
\e{0}%
\e{0}%
\e{0}%
\e{0}%
\e{0}%
\e{0}%
\e{0}%
\e{0}%
\e{0}%
\e{0}%
\e{0}%
\e{1}%
\e{0}%
\e{1}%
\e{1}%
\e{0}%
\e{0}%
\e{1}%
\e{0}%
\e{1}%
\e{0}%
\e{0}%
\eol}\vss}\rg%
%
%
\rx{\vss\hfull{%
\rlx{\hss{$2835_{x}$}}\cg%
\e{0}%
\e{0}%
\e{0}%
\e{0}%
\e{0}%
\e{0}%
\e{0}%
\e{0}%
\e{0}%
\e{0}%
\e{0}%
\e{0}%
\e{0}%
\e{0}%
\e{0}%
\e{0}%
\e{1}%
\e{0}%
\e{0}%
\e{0}%
\e{1}%
\e{1}%
\e{1}%
\e{1}%
\e{0}%
\eol}\vss}\rg%
%
%
\rx{\vss\hfull{%
\rlx{\hss{$6075_{x}$}}\cg%
\e{0}%
\e{0}%
\e{0}%
\e{0}%
\e{0}%
\e{0}%
\e{0}%
\e{0}%
\e{0}%
\e{0}%
\e{0}%
\e{0}%
\e{0}%
\e{0}%
\e{0}%
\e{1}%
\e{1}%
\e{1}%
\e{1}%
\e{1}%
\e{0}%
\e{1}%
\e{1}%
\e{1}%
\e{1}%
\eol}\vss}\rg%
%
%
\rx{\vss\hfull{%
\rlx{\hss{$3200_{x}$}}\cg%
\e{0}%
\e{0}%
\e{0}%
\e{0}%
\e{0}%
\e{0}%
\e{0}%
\e{0}%
\e{0}%
\e{0}%
\e{0}%
\e{0}%
\e{0}%
\e{0}%
\e{0}%
\e{0}%
\e{1}%
\e{0}%
\e{0}%
\e{1}%
\e{0}%
\e{0}%
\e{1}%
\e{0}%
\e{0}%
\eol}\vss}\rg%
%
%
\rx{\vss\hfull{%
\rlx{\hss{$70_{y}$}}\cg%
\e{0}%
\e{0}%
\e{0}%
\e{0}%
\e{0}%
\e{0}%
\e{0}%
\e{0}%
\e{0}%
\e{0}%
\e{0}%
\e{0}%
\e{0}%
\e{0}%
\e{0}%
\e{0}%
\e{0}%
\e{0}%
\e{0}%
\e{0}%
\e{0}%
\e{0}%
\e{0}%
\e{0}%
\e{0}%
\eol}\vss}\rg%
%
%
\rx{\vss\hfull{%
\rlx{\hss{$1134_{y}$}}\cg%
\e{0}%
\e{0}%
\e{0}%
\e{0}%
\e{0}%
\e{0}%
\e{0}%
\e{0}%
\e{0}%
\e{0}%
\e{0}%
\e{0}%
\e{0}%
\e{0}%
\e{0}%
\e{0}%
\e{0}%
\e{1}%
\e{0}%
\e{0}%
\e{0}%
\e{0}%
\e{0}%
\e{0}%
\e{0}%
\eol}\vss}\rg%
%
%
\rx{\vss\hfull{%
\rlx{\hss{$1680_{y}$}}\cg%
\e{0}%
\e{0}%
\e{0}%
\e{0}%
\e{0}%
\e{0}%
\e{0}%
\e{0}%
\e{0}%
\e{0}%
\e{0}%
\e{0}%
\e{0}%
\e{0}%
\e{0}%
\e{0}%
\e{0}%
\e{0}%
\e{0}%
\e{1}%
\e{0}%
\e{0}%
\e{0}%
\e{0}%
\e{1}%
\eol}\vss}\rg%
%
%
\rx{\vss\hfull{%
\rlx{\hss{$168_{y}$}}\cg%
\e{0}%
\e{0}%
\e{0}%
\e{0}%
\e{0}%
\e{0}%
\e{0}%
\e{0}%
\e{0}%
\e{0}%
\e{0}%
\e{0}%
\e{0}%
\e{0}%
\e{0}%
\e{0}%
\e{0}%
\e{0}%
\e{0}%
\e{0}%
\e{0}%
\e{0}%
\e{0}%
\e{0}%
\e{0}%
\eol}\vss}\rg%
%
%
\rx{\vss\hfull{%
\rlx{\hss{$420_{y}$}}\cg%
\e{0}%
\e{0}%
\e{0}%
\e{0}%
\e{0}%
\e{0}%
\e{0}%
\e{0}%
\e{0}%
\e{0}%
\e{0}%
\e{0}%
\e{0}%
\e{0}%
\e{0}%
\e{0}%
\e{0}%
\e{0}%
\e{0}%
\e{0}%
\e{1}%
\e{0}%
\e{0}%
\e{0}%
\e{0}%
\eol}\vss}\rg%
%
%
\rx{\vss\hfull{%
\rlx{\hss{$3150_{y}$}}\cg%
\e{0}%
\e{0}%
\e{0}%
\e{0}%
\e{0}%
\e{0}%
\e{0}%
\e{0}%
\e{0}%
\e{0}%
\e{0}%
\e{0}%
\e{0}%
\e{0}%
\e{0}%
\e{0}%
\e{0}%
\e{1}%
\e{0}%
\e{0}%
\e{1}%
\e{0}%
\e{1}%
\e{1}%
\e{0}%
\eol}\vss}\rg%
%
%
\rx{\vss\hfull{%
\rlx{\hss{$4200_{y}$}}\cg%
\e{0}%
\e{0}%
\e{0}%
\e{0}%
\e{0}%
\e{0}%
\e{0}%
\e{0}%
\e{0}%
\e{0}%
\e{0}%
\e{0}%
\e{0}%
\e{0}%
\e{0}%
\e{0}%
\e{1}%
\e{0}%
\e{0}%
\e{0}%
\e{0}%
\e{1}%
\e{1}%
\e{1}%
\e{0}%
\eol}\vss}\rg%
\eop
\eject
\tablecont%
%
%
%
%
%
%
\rowpts=18 true pt%
\colpts=18 true pt%
\rowlabpts=40 true pt%
\collabpts=70 true pt%
\clx{\vss\hfull{%
\rlx{\hss{$ $}}\cg%
\cx{\hskip 16 true pt\flip{$[{5}{1^{2}}:{1}]$}\hss}\cg%
\cx{\hskip 16 true pt\flip{$[{4}{3}:{1}]$}\hss}\cg%
\cx{\hskip 16 true pt\flip{$[{4}{2}{1}:{1}]$}\hss}\cg%
\cx{\hskip 16 true pt\flip{$[{4}{1^{3}}:{1}]$}\hss}\cg%
\cx{\hskip 16 true pt\flip{$[{3^{2}}{1}:{1}]$}\hss}\cg%
\cx{\hskip 16 true pt\flip{$[{3}{2^{2}}:{1}]$}\hss}\cg%
\cx{\hskip 16 true pt\flip{$[{3}{2}{1^{2}}:{1}]$}\hss}\cg%
\cx{\hskip 16 true pt\flip{$[{3}{1^{4}}:{1}]$}\hss}\cg%
\cx{\hskip 16 true pt\flip{$[{2^{3}}{1}:{1}]$}\hss}\cg%
\cx{\hskip 16 true pt\flip{$[{2^{2}}{1^{3}}:{1}]$}\hss}\cg%
\cx{\hskip 16 true pt\flip{$[{2}{1^{5}}:{1}]$}\hss}\cg%
\cx{\hskip 16 true pt\flip{$[{1^{7}}:{1}]$}\hss}\cg%
\cx{\hskip 16 true pt\flip{$[{6}:{2}]$}\hss}\cg%
\cx{\hskip 16 true pt\flip{$[{6}:{1^{2}}]$}\hss}\cg%
\cx{\hskip 16 true pt\flip{$[{5}{1}:{2}]$}\hss}\cg%
\cx{\hskip 16 true pt\flip{$[{5}{1}:{1^{2}}]$}\hss}\cg%
\cx{\hskip 16 true pt\flip{$[{4}{2}:{2}]$}\hss}\cg%
\cx{\hskip 16 true pt\flip{$[{4}{2}:{1^{2}}]$}\hss}\cg%
\cx{\hskip 16 true pt\flip{$[{4}{1^{2}}:{2}]$}\hss}\cg%
\cx{\hskip 16 true pt\flip{$[{4}{1^{2}}:{1^{2}}]$}\hss}\cg%
\cx{\hskip 16 true pt\flip{$[{3^{2}}:{2}]$}\hss}\cg%
\cx{\hskip 16 true pt\flip{$[{3^{2}}:{1^{2}}]$}\hss}\cg%
\cx{\hskip 16 true pt\flip{$[{3}{2}{1}:{2}]$}\hss}\cg%
\cx{\hskip 16 true pt\flip{$[{3}{2}{1}:{1^{2}}]$}\hss}\cg%
\cx{\hskip 16 true pt\flip{$[{3}{1^{3}}:{2}]$}\hss}\cg%
\eol}}\rg%
%
%
\rx{\vss\hfull{%
\rlx{\hss{$2688_{y}$}}\cg%
\e{0}%
\e{0}%
\e{0}%
\e{0}%
\e{0}%
\e{0}%
\e{0}%
\e{0}%
\e{0}%
\e{0}%
\e{0}%
\e{0}%
\e{0}%
\e{0}%
\e{0}%
\e{0}%
\e{0}%
\e{0}%
\e{0}%
\e{0}%
\e{0}%
\e{0}%
\e{1}%
\e{1}%
\e{0}%
\eol}\vss}\rg%
%
%
\rx{\vss\hfull{%
\rlx{\hss{$2100_{y}$}}\cg%
\e{0}%
\e{0}%
\e{0}%
\e{0}%
\e{0}%
\e{0}%
\e{0}%
\e{0}%
\e{0}%
\e{0}%
\e{0}%
\e{0}%
\e{0}%
\e{0}%
\e{0}%
\e{0}%
\e{0}%
\e{0}%
\e{0}%
\e{1}%
\e{0}%
\e{0}%
\e{0}%
\e{0}%
\e{1}%
\eol}\vss}\rg%
%
%
\rx{\vss\hfull{%
\rlx{\hss{$1400_{y}$}}\cg%
\e{0}%
\e{0}%
\e{0}%
\e{0}%
\e{0}%
\e{0}%
\e{0}%
\e{0}%
\e{0}%
\e{0}%
\e{0}%
\e{0}%
\e{0}%
\e{0}%
\e{0}%
\e{0}%
\e{0}%
\e{0}%
\e{1}%
\e{0}%
\e{0}%
\e{0}%
\e{0}%
\e{0}%
\e{0}%
\eol}\vss}\rg%
%
%
\rx{\vss\hfull{%
\rlx{\hss{$4536_{y}$}}\cg%
\e{0}%
\e{0}%
\e{0}%
\e{0}%
\e{0}%
\e{0}%
\e{0}%
\e{0}%
\e{0}%
\e{0}%
\e{0}%
\e{0}%
\e{0}%
\e{0}%
\e{0}%
\e{0}%
\e{0}%
\e{0}%
\e{1}%
\e{1}%
\e{0}%
\e{0}%
\e{1}%
\e{1}%
\e{1}%
\eol}\vss}\rg%
%
%
\rx{\vss\hfull{%
\rlx{\hss{$5670_{y}$}}\cg%
\e{0}%
\e{0}%
\e{0}%
\e{0}%
\e{0}%
\e{0}%
\e{0}%
\e{0}%
\e{0}%
\e{0}%
\e{0}%
\e{0}%
\e{0}%
\e{0}%
\e{0}%
\e{0}%
\e{0}%
\e{1}%
\e{1}%
\e{1}%
\e{0}%
\e{0}%
\e{1}%
\e{1}%
\e{1}%
\eol}\vss}\rg%
%
%
\rx{\vss\hfull{%
\rlx{\hss{$4480_{y}$}}\cg%
\e{0}%
\e{0}%
\e{0}%
\e{0}%
\e{0}%
\e{0}%
\e{0}%
\e{0}%
\e{0}%
\e{0}%
\e{0}%
\e{0}%
\e{0}%
\e{0}%
\e{0}%
\e{0}%
\e{0}%
\e{1}%
\e{1}%
\e{0}%
\e{1}%
\e{0}%
\e{1}%
\e{1}%
\e{0}%
\eol}\vss}\rg%
%
%
\rx{\vss\hfull{%
\rlx{\hss{$8_{z}$}}\cg%
\e{0}%
\e{0}%
\e{0}%
\e{0}%
\e{0}%
\e{0}%
\e{0}%
\e{0}%
\e{0}%
\e{0}%
\e{0}%
\e{0}%
\e{0}%
\e{0}%
\e{0}%
\e{0}%
\e{0}%
\e{0}%
\e{0}%
\e{0}%
\e{0}%
\e{0}%
\e{0}%
\e{0}%
\e{0}%
\eol}\vss}\rg%
%
%
\rx{\vss\hfull{%
\rlx{\hss{$56_{z}$}}\cg%
\e{0}%
\e{0}%
\e{0}%
\e{0}%
\e{0}%
\e{0}%
\e{0}%
\e{0}%
\e{0}%
\e{0}%
\e{0}%
\e{0}%
\e{0}%
\e{0}%
\e{0}%
\e{0}%
\e{0}%
\e{0}%
\e{0}%
\e{0}%
\e{0}%
\e{0}%
\e{0}%
\e{0}%
\e{0}%
\eol}\vss}\rg%
%
%
\rx{\vss\hfull{%
\rlx{\hss{$160_{z}$}}\cg%
\e{0}%
\e{0}%
\e{0}%
\e{0}%
\e{0}%
\e{0}%
\e{0}%
\e{0}%
\e{0}%
\e{0}%
\e{0}%
\e{0}%
\e{0}%
\e{0}%
\e{0}%
\e{0}%
\e{0}%
\e{0}%
\e{0}%
\e{0}%
\e{0}%
\e{0}%
\e{0}%
\e{0}%
\e{0}%
\eol}\vss}\rg%
%
%
\rx{\vss\hfull{%
\rlx{\hss{$112_{z}$}}\cg%
\e{0}%
\e{0}%
\e{0}%
\e{0}%
\e{0}%
\e{0}%
\e{0}%
\e{0}%
\e{0}%
\e{0}%
\e{0}%
\e{0}%
\e{0}%
\e{0}%
\e{0}%
\e{0}%
\e{0}%
\e{0}%
\e{0}%
\e{0}%
\e{0}%
\e{0}%
\e{0}%
\e{0}%
\e{0}%
\eol}\vss}\rg%
%
%
\rx{\vss\hfull{%
\rlx{\hss{$840_{z}$}}\cg%
\e{0}%
\e{0}%
\e{0}%
\e{0}%
\e{0}%
\e{0}%
\e{0}%
\e{0}%
\e{0}%
\e{0}%
\e{0}%
\e{0}%
\e{0}%
\e{0}%
\e{0}%
\e{0}%
\e{0}%
\e{0}%
\e{0}%
\e{0}%
\e{0}%
\e{0}%
\e{0}%
\e{0}%
\e{0}%
\eol}\vss}\rg%
%
%
\rx{\vss\hfull{%
\rlx{\hss{$1296_{z}$}}\cg%
\e{1}%
\e{0}%
\e{0}%
\e{0}%
\e{0}%
\e{0}%
\e{0}%
\e{0}%
\e{0}%
\e{0}%
\e{0}%
\e{0}%
\e{0}%
\e{0}%
\e{0}%
\e{0}%
\e{0}%
\e{0}%
\e{0}%
\e{0}%
\e{0}%
\e{0}%
\e{0}%
\e{0}%
\e{0}%
\eol}\vss}\rg%
%
%
\rx{\vss\hfull{%
\rlx{\hss{$1400_{z}$}}\cg%
\e{1}%
\e{0}%
\e{0}%
\e{0}%
\e{0}%
\e{0}%
\e{0}%
\e{0}%
\e{0}%
\e{0}%
\e{0}%
\e{0}%
\e{0}%
\e{0}%
\e{0}%
\e{0}%
\e{0}%
\e{0}%
\e{0}%
\e{0}%
\e{0}%
\e{0}%
\e{0}%
\e{0}%
\e{0}%
\eol}\vss}\rg%
%
%
\rx{\vss\hfull{%
\rlx{\hss{$1008_{z}$}}\cg%
\e{1}%
\e{0}%
\e{0}%
\e{0}%
\e{0}%
\e{0}%
\e{0}%
\e{0}%
\e{0}%
\e{0}%
\e{0}%
\e{0}%
\e{0}%
\e{0}%
\e{0}%
\e{0}%
\e{0}%
\e{0}%
\e{0}%
\e{0}%
\e{0}%
\e{0}%
\e{0}%
\e{0}%
\e{0}%
\eol}\vss}\rg%
%
%
\rx{\vss\hfull{%
\rlx{\hss{$560_{z}$}}\cg%
\e{0}%
\e{0}%
\e{0}%
\e{0}%
\e{0}%
\e{0}%
\e{0}%
\e{0}%
\e{0}%
\e{0}%
\e{0}%
\e{0}%
\e{0}%
\e{0}%
\e{0}%
\e{0}%
\e{0}%
\e{0}%
\e{0}%
\e{0}%
\e{0}%
\e{0}%
\e{0}%
\e{0}%
\e{0}%
\eol}\vss}\rg%
%
%
\rx{\vss\hfull{%
\rlx{\hss{$1400_{zz}$}}\cg%
\e{0}%
\e{1}%
\e{0}%
\e{0}%
\e{1}%
\e{0}%
\e{0}%
\e{0}%
\e{0}%
\e{0}%
\e{0}%
\e{0}%
\e{0}%
\e{0}%
\e{0}%
\e{0}%
\e{0}%
\e{0}%
\e{0}%
\e{0}%
\e{0}%
\e{0}%
\e{0}%
\e{0}%
\e{0}%
\eol}\vss}\rg%
%
%
\rx{\vss\hfull{%
\rlx{\hss{$4200_{z}$}}\cg%
\e{0}%
\e{1}%
\e{0}%
\e{0}%
\e{1}%
\e{0}%
\e{1}%
\e{0}%
\e{0}%
\e{0}%
\e{0}%
\e{0}%
\e{0}%
\e{0}%
\e{0}%
\e{0}%
\e{0}%
\e{0}%
\e{0}%
\e{0}%
\e{0}%
\e{0}%
\e{0}%
\e{0}%
\e{0}%
\eol}\vss}\rg%
%
%
\rx{\vss\hfull{%
\rlx{\hss{$400_{z}$}}\cg%
\e{0}%
\e{1}%
\e{0}%
\e{0}%
\e{0}%
\e{0}%
\e{0}%
\e{0}%
\e{0}%
\e{0}%
\e{0}%
\e{0}%
\e{0}%
\e{0}%
\e{0}%
\e{0}%
\e{0}%
\e{0}%
\e{0}%
\e{0}%
\e{0}%
\e{0}%
\e{0}%
\e{0}%
\e{0}%
\eol}\vss}\rg%
%
%
\rx{\vss\hfull{%
\rlx{\hss{$3240_{z}$}}\cg%
\e{0}%
\e{1}%
\e{1}%
\e{0}%
\e{0}%
\e{0}%
\e{0}%
\e{0}%
\e{0}%
\e{0}%
\e{0}%
\e{0}%
\e{0}%
\e{0}%
\e{0}%
\e{0}%
\e{0}%
\e{0}%
\e{0}%
\e{0}%
\e{0}%
\e{0}%
\e{0}%
\e{0}%
\e{0}%
\eol}\vss}\rg%
%
%
\rx{\vss\hfull{%
\rlx{\hss{$4536_{z}$}}\cg%
\e{0}%
\e{1}%
\e{1}%
\e{0}%
\e{0}%
\e{1}%
\e{0}%
\e{0}%
\e{0}%
\e{0}%
\e{0}%
\e{0}%
\e{0}%
\e{0}%
\e{0}%
\e{0}%
\e{0}%
\e{0}%
\e{0}%
\e{0}%
\e{0}%
\e{0}%
\e{0}%
\e{0}%
\e{0}%
\eol}\vss}\rg%
%
%
\rx{\vss\hfull{%
\rlx{\hss{$2400_{z}$}}\cg%
\e{0}%
\e{0}%
\e{0}%
\e{1}%
\e{0}%
\e{0}%
\e{0}%
\e{0}%
\e{0}%
\e{0}%
\e{0}%
\e{0}%
\e{0}%
\e{0}%
\e{0}%
\e{0}%
\e{0}%
\e{0}%
\e{0}%
\e{0}%
\e{0}%
\e{0}%
\e{0}%
\e{0}%
\e{0}%
\eol}\vss}\rg%
%
%
\rx{\vss\hfull{%
\rlx{\hss{$3360_{z}$}}\cg%
\e{0}%
\e{0}%
\e{1}%
\e{0}%
\e{1}%
\e{0}%
\e{0}%
\e{0}%
\e{0}%
\e{0}%
\e{0}%
\e{0}%
\e{0}%
\e{0}%
\e{0}%
\e{0}%
\e{0}%
\e{0}%
\e{0}%
\e{0}%
\e{0}%
\e{0}%
\e{0}%
\e{0}%
\e{0}%
\eol}\vss}\rg%
%
%
\rx{\vss\hfull{%
\rlx{\hss{$2800_{z}$}}\cg%
\e{1}%
\e{0}%
\e{0}%
\e{1}%
\e{0}%
\e{0}%
\e{0}%
\e{0}%
\e{0}%
\e{0}%
\e{0}%
\e{0}%
\e{0}%
\e{0}%
\e{0}%
\e{0}%
\e{0}%
\e{0}%
\e{0}%
\e{0}%
\e{0}%
\e{0}%
\e{0}%
\e{0}%
\e{0}%
\eol}\vss}\rg%
%
%
\rx{\vss\hfull{%
\rlx{\hss{$4096_{z}$}}\cg%
\e{1}%
\e{0}%
\e{1}%
\e{0}%
\e{0}%
\e{0}%
\e{0}%
\e{0}%
\e{0}%
\e{0}%
\e{0}%
\e{0}%
\e{0}%
\e{0}%
\e{0}%
\e{0}%
\e{0}%
\e{0}%
\e{0}%
\e{0}%
\e{0}%
\e{0}%
\e{0}%
\e{0}%
\e{0}%
\eol}\vss}\rg%
%
%
\rx{\vss\hfull{%
\rlx{\hss{$5600_{z}$}}\cg%
\e{1}%
\e{0}%
\e{1}%
\e{1}%
\e{0}%
\e{0}%
\e{0}%
\e{0}%
\e{0}%
\e{0}%
\e{0}%
\e{0}%
\e{0}%
\e{0}%
\e{0}%
\e{0}%
\e{0}%
\e{0}%
\e{0}%
\e{0}%
\e{0}%
\e{0}%
\e{0}%
\e{0}%
\e{0}%
\eol}\vss}\rg%
%
%
\rx{\vss\hfull{%
\rlx{\hss{$448_{z}$}}\cg%
\e{0}%
\e{0}%
\e{0}%
\e{0}%
\e{0}%
\e{0}%
\e{0}%
\e{0}%
\e{0}%
\e{0}%
\e{0}%
\e{0}%
\e{0}%
\e{0}%
\e{0}%
\e{0}%
\e{0}%
\e{0}%
\e{0}%
\e{0}%
\e{0}%
\e{0}%
\e{0}%
\e{0}%
\e{0}%
\eol}\vss}\rg%
%
%
\rx{\vss\hfull{%
\rlx{\hss{$448_{w}$}}\cg%
\e{0}%
\e{0}%
\e{0}%
\e{0}%
\e{0}%
\e{0}%
\e{0}%
\e{0}%
\e{0}%
\e{0}%
\e{0}%
\e{0}%
\e{0}%
\e{0}%
\e{0}%
\e{0}%
\e{0}%
\e{0}%
\e{0}%
\e{0}%
\e{0}%
\e{0}%
\e{0}%
\e{0}%
\e{0}%
\eol}\vss}\rg%
%
%
\rx{\vss\hfull{%
\rlx{\hss{$1344_{w}$}}\cg%
\e{0}%
\e{1}%
\e{0}%
\e{0}%
\e{0}%
\e{0}%
\e{0}%
\e{0}%
\e{1}%
\e{0}%
\e{0}%
\e{0}%
\e{0}%
\e{0}%
\e{0}%
\e{0}%
\e{0}%
\e{0}%
\e{0}%
\e{0}%
\e{0}%
\e{0}%
\e{0}%
\e{0}%
\e{0}%
\eol}\vss}\rg%
%
%
\rx{\vss\hfull{%
\rlx{\hss{$5600_{w}$}}\cg%
\e{0}%
\e{0}%
\e{1}%
\e{0}%
\e{0}%
\e{0}%
\e{1}%
\e{0}%
\e{0}%
\e{0}%
\e{0}%
\e{0}%
\e{0}%
\e{0}%
\e{0}%
\e{0}%
\e{0}%
\e{0}%
\e{0}%
\e{0}%
\e{0}%
\e{0}%
\e{0}%
\e{0}%
\e{0}%
\eol}\vss}\rg%
%
%
\rx{\vss\hfull{%
\rlx{\hss{$2016_{w}$}}\cg%
\e{0}%
\e{0}%
\e{0}%
\e{0}%
\e{1}%
\e{1}%
\e{0}%
\e{0}%
\e{0}%
\e{0}%
\e{0}%
\e{0}%
\e{0}%
\e{0}%
\e{0}%
\e{0}%
\e{0}%
\e{0}%
\e{0}%
\e{0}%
\e{0}%
\e{0}%
\e{0}%
\e{0}%
\e{0}%
\eol}\vss}\rg%
%
%
\rx{\vss\hfull{%
\rlx{\hss{$7168_{w}$}}\cg%
\e{0}%
\e{0}%
\e{1}%
\e{0}%
\e{1}%
\e{1}%
\e{1}%
\e{0}%
\e{0}%
\e{0}%
\e{0}%
\e{0}%
\e{0}%
\e{0}%
\e{0}%
\e{0}%
\e{0}%
\e{0}%
\e{0}%
\e{0}%
\e{0}%
\e{0}%
\e{0}%
\e{0}%
\e{0}%
\eol}\vss}\rg%
\eop
\eject
\tablecont%
%
%
%
%
%
%
\rowpts=18 true pt%
\colpts=18 true pt%
\rowlabpts=40 true pt%
\collabpts=70 true pt%
\clx{\vss\hfull{%
\rlx{\hss{$ $}}\cg%
\cx{\hskip 16 true pt\flip{$[{3}{1^{3}}:{1^{2}}]$}\hss}\cg%
\cx{\hskip 16 true pt\flip{$[{2^{3}}:{2}]$}\hss}\cg%
\cx{\hskip 16 true pt\flip{$[{2^{3}}:{1^{2}}]$}\hss}\cg%
\cx{\hskip 16 true pt\flip{$[{2^{2}}{1^{2}}:{2}]$}\hss}\cg%
\cx{\hskip 16 true pt\flip{$[{2^{2}}{1^{2}}:{1^{2}}]$}\hss}\cg%
\cx{\hskip 16 true pt\flip{$[{2}{1^{4}}:{2}]$}\hss}\cg%
\cx{\hskip 16 true pt\flip{$[{2}{1^{4}}:{1^{2}}]$}\hss}\cg%
\cx{\hskip 16 true pt\flip{$[{1^{6}}:{2}]$}\hss}\cg%
\cx{\hskip 16 true pt\flip{$[{1^{6}}:{1^{2}}]$}\hss}\cg%
\cx{\hskip 16 true pt\flip{$[{5}:{3}]$}\hss}\cg%
\cx{\hskip 16 true pt\flip{$[{5}:{2}{1}]$}\hss}\cg%
\cx{\hskip 16 true pt\flip{$[{5}:{1^{3}}]$}\hss}\cg%
\cx{\hskip 16 true pt\flip{$[{4}{1}:{3}]$}\hss}\cg%
\cx{\hskip 16 true pt\flip{$[{4}{1}:{2}{1}]$}\hss}\cg%
\cx{\hskip 16 true pt\flip{$[{4}{1}:{1^{3}}]$}\hss}\cg%
\cx{\hskip 16 true pt\flip{$[{3}{2}:{3}]$}\hss}\cg%
\cx{\hskip 16 true pt\flip{$[{3}{2}:{2}{1}]$}\hss}\cg%
\cx{\hskip 16 true pt\flip{$[{3}{2}:{1^{3}}]$}\hss}\cg%
\cx{\hskip 16 true pt\flip{$[{3}{1^{2}}:{3}]$}\hss}\cg%
\cx{\hskip 16 true pt\flip{$[{3}{1^{2}}:{2}{1}]$}\hss}\cg%
\cx{\hskip 16 true pt\flip{$[{3}{1^{2}}:{1^{3}}]$}\hss}\cg%
\cx{\hskip 16 true pt\flip{$[{2^{2}}{1}:{3}]$}\hss}\cg%
\cx{\hskip 16 true pt\flip{$[{2^{2}}{1}:{2}{1}]$}\hss}\cg%
\cx{\hskip 16 true pt\flip{$[{2^{2}}{1}:{1^{3}}]$}\hss}\cg%
\cx{\hskip 16 true pt\flip{$[{2}{1^{3}}:{3}]$}\hss}\cg%
\eol}}\rg%
%
%
\rx{\vss\hfull{%
\rlx{\hss{$1_{x}$}}\cg%
\e{0}%
\e{0}%
\e{0}%
\e{0}%
\e{0}%
\e{0}%
\e{0}%
\e{0}%
\e{0}%
\e{0}%
\e{0}%
\e{0}%
\e{0}%
\e{0}%
\e{0}%
\e{0}%
\e{0}%
\e{0}%
\e{0}%
\e{0}%
\e{0}%
\e{0}%
\e{0}%
\e{0}%
\e{0}%
\eol}\vss}\rg%
%
%
\rx{\vss\hfull{%
\rlx{\hss{$28_{x}$}}\cg%
\e{0}%
\e{0}%
\e{0}%
\e{0}%
\e{0}%
\e{0}%
\e{0}%
\e{0}%
\e{0}%
\e{0}%
\e{0}%
\e{0}%
\e{0}%
\e{0}%
\e{0}%
\e{0}%
\e{0}%
\e{0}%
\e{0}%
\e{0}%
\e{0}%
\e{0}%
\e{0}%
\e{0}%
\e{0}%
\eol}\vss}\rg%
%
%
\rx{\vss\hfull{%
\rlx{\hss{$35_{x}$}}\cg%
\e{0}%
\e{0}%
\e{0}%
\e{0}%
\e{0}%
\e{0}%
\e{0}%
\e{0}%
\e{0}%
\e{0}%
\e{0}%
\e{0}%
\e{0}%
\e{0}%
\e{0}%
\e{0}%
\e{0}%
\e{0}%
\e{0}%
\e{0}%
\e{0}%
\e{0}%
\e{0}%
\e{0}%
\e{0}%
\eol}\vss}\rg%
%
%
\rx{\vss\hfull{%
\rlx{\hss{$84_{x}$}}\cg%
\e{0}%
\e{0}%
\e{0}%
\e{0}%
\e{0}%
\e{0}%
\e{0}%
\e{0}%
\e{0}%
\e{0}%
\e{0}%
\e{0}%
\e{0}%
\e{0}%
\e{0}%
\e{0}%
\e{0}%
\e{0}%
\e{0}%
\e{0}%
\e{0}%
\e{0}%
\e{0}%
\e{0}%
\e{0}%
\eol}\vss}\rg%
%
%
\rx{\vss\hfull{%
\rlx{\hss{$50_{x}$}}\cg%
\e{0}%
\e{0}%
\e{0}%
\e{0}%
\e{0}%
\e{0}%
\e{0}%
\e{0}%
\e{0}%
\e{0}%
\e{0}%
\e{0}%
\e{0}%
\e{0}%
\e{0}%
\e{0}%
\e{0}%
\e{0}%
\e{0}%
\e{0}%
\e{0}%
\e{0}%
\e{0}%
\e{0}%
\e{0}%
\eol}\vss}\rg%
%
%
\rx{\vss\hfull{%
\rlx{\hss{$350_{x}$}}\cg%
\e{0}%
\e{0}%
\e{0}%
\e{0}%
\e{0}%
\e{0}%
\e{0}%
\e{0}%
\e{0}%
\e{0}%
\e{0}%
\e{0}%
\e{0}%
\e{0}%
\e{0}%
\e{0}%
\e{0}%
\e{0}%
\e{0}%
\e{0}%
\e{0}%
\e{0}%
\e{0}%
\e{0}%
\e{0}%
\eol}\vss}\rg%
%
%
\rx{\vss\hfull{%
\rlx{\hss{$300_{x}$}}\cg%
\e{0}%
\e{0}%
\e{0}%
\e{0}%
\e{0}%
\e{0}%
\e{0}%
\e{0}%
\e{0}%
\e{0}%
\e{0}%
\e{0}%
\e{0}%
\e{0}%
\e{0}%
\e{0}%
\e{0}%
\e{0}%
\e{0}%
\e{0}%
\e{0}%
\e{0}%
\e{0}%
\e{0}%
\e{0}%
\eol}\vss}\rg%
%
%
\rx{\vss\hfull{%
\rlx{\hss{$567_{x}$}}\cg%
\e{0}%
\e{0}%
\e{0}%
\e{0}%
\e{0}%
\e{0}%
\e{0}%
\e{0}%
\e{0}%
\e{0}%
\e{0}%
\e{0}%
\e{0}%
\e{0}%
\e{0}%
\e{0}%
\e{0}%
\e{0}%
\e{0}%
\e{0}%
\e{0}%
\e{0}%
\e{0}%
\e{0}%
\e{0}%
\eol}\vss}\rg%
%
%
\rx{\vss\hfull{%
\rlx{\hss{$210_{x}$}}\cg%
\e{0}%
\e{0}%
\e{0}%
\e{0}%
\e{0}%
\e{0}%
\e{0}%
\e{0}%
\e{0}%
\e{0}%
\e{0}%
\e{0}%
\e{0}%
\e{0}%
\e{0}%
\e{0}%
\e{0}%
\e{0}%
\e{0}%
\e{0}%
\e{0}%
\e{0}%
\e{0}%
\e{0}%
\e{0}%
\eol}\vss}\rg%
%
%
\rx{\vss\hfull{%
\rlx{\hss{$840_{x}$}}\cg%
\e{0}%
\e{1}%
\e{0}%
\e{0}%
\e{0}%
\e{0}%
\e{0}%
\e{0}%
\e{0}%
\e{0}%
\e{0}%
\e{0}%
\e{0}%
\e{0}%
\e{0}%
\e{0}%
\e{0}%
\e{0}%
\e{0}%
\e{0}%
\e{0}%
\e{0}%
\e{0}%
\e{0}%
\e{0}%
\eol}\vss}\rg%
%
%
\rx{\vss\hfull{%
\rlx{\hss{$700_{x}$}}\cg%
\e{0}%
\e{0}%
\e{0}%
\e{0}%
\e{0}%
\e{0}%
\e{0}%
\e{0}%
\e{0}%
\e{0}%
\e{0}%
\e{0}%
\e{0}%
\e{0}%
\e{0}%
\e{0}%
\e{0}%
\e{0}%
\e{0}%
\e{0}%
\e{0}%
\e{0}%
\e{0}%
\e{0}%
\e{0}%
\eol}\vss}\rg%
%
%
\rx{\vss\hfull{%
\rlx{\hss{$175_{x}$}}\cg%
\e{0}%
\e{0}%
\e{0}%
\e{0}%
\e{0}%
\e{0}%
\e{0}%
\e{0}%
\e{0}%
\e{0}%
\e{0}%
\e{0}%
\e{0}%
\e{0}%
\e{0}%
\e{0}%
\e{0}%
\e{0}%
\e{0}%
\e{0}%
\e{0}%
\e{0}%
\e{0}%
\e{0}%
\e{0}%
\eol}\vss}\rg%
%
%
\rx{\vss\hfull{%
\rlx{\hss{$1400_{x}$}}\cg%
\e{0}%
\e{0}%
\e{0}%
\e{0}%
\e{0}%
\e{0}%
\e{0}%
\e{0}%
\e{0}%
\e{0}%
\e{0}%
\e{0}%
\e{0}%
\e{0}%
\e{0}%
\e{0}%
\e{0}%
\e{0}%
\e{0}%
\e{0}%
\e{0}%
\e{0}%
\e{0}%
\e{0}%
\e{0}%
\eol}\vss}\rg%
%
%
\rx{\vss\hfull{%
\rlx{\hss{$1050_{x}$}}\cg%
\e{0}%
\e{0}%
\e{0}%
\e{0}%
\e{0}%
\e{0}%
\e{0}%
\e{0}%
\e{0}%
\e{0}%
\e{0}%
\e{0}%
\e{0}%
\e{0}%
\e{0}%
\e{0}%
\e{0}%
\e{0}%
\e{0}%
\e{0}%
\e{0}%
\e{0}%
\e{0}%
\e{0}%
\e{0}%
\eol}\vss}\rg%
%
%
\rx{\vss\hfull{%
\rlx{\hss{$1575_{x}$}}\cg%
\e{0}%
\e{0}%
\e{0}%
\e{0}%
\e{0}%
\e{0}%
\e{0}%
\e{0}%
\e{0}%
\e{0}%
\e{0}%
\e{0}%
\e{0}%
\e{0}%
\e{0}%
\e{0}%
\e{0}%
\e{0}%
\e{0}%
\e{0}%
\e{0}%
\e{0}%
\e{0}%
\e{0}%
\e{0}%
\eol}\vss}\rg%
%
%
\rx{\vss\hfull{%
\rlx{\hss{$1344_{x}$}}\cg%
\e{0}%
\e{0}%
\e{0}%
\e{0}%
\e{0}%
\e{0}%
\e{0}%
\e{0}%
\e{0}%
\e{0}%
\e{0}%
\e{0}%
\e{0}%
\e{0}%
\e{0}%
\e{0}%
\e{0}%
\e{0}%
\e{0}%
\e{0}%
\e{0}%
\e{0}%
\e{0}%
\e{0}%
\e{0}%
\eol}\vss}\rg%
%
%
\rx{\vss\hfull{%
\rlx{\hss{$2100_{x}$}}\cg%
\e{0}%
\e{0}%
\e{0}%
\e{0}%
\e{0}%
\e{0}%
\e{0}%
\e{0}%
\e{0}%
\e{0}%
\e{0}%
\e{0}%
\e{0}%
\e{0}%
\e{0}%
\e{0}%
\e{0}%
\e{0}%
\e{0}%
\e{0}%
\e{0}%
\e{0}%
\e{0}%
\e{0}%
\e{0}%
\eol}\vss}\rg%
%
%
\rx{\vss\hfull{%
\rlx{\hss{$2268_{x}$}}\cg%
\e{0}%
\e{0}%
\e{0}%
\e{0}%
\e{0}%
\e{0}%
\e{0}%
\e{0}%
\e{0}%
\e{0}%
\e{0}%
\e{0}%
\e{0}%
\e{0}%
\e{0}%
\e{0}%
\e{0}%
\e{0}%
\e{0}%
\e{0}%
\e{0}%
\e{0}%
\e{0}%
\e{0}%
\e{0}%
\eol}\vss}\rg%
%
%
\rx{\vss\hfull{%
\rlx{\hss{$525_{x}$}}\cg%
\e{0}%
\e{0}%
\e{0}%
\e{0}%
\e{0}%
\e{0}%
\e{0}%
\e{0}%
\e{0}%
\e{0}%
\e{0}%
\e{0}%
\e{0}%
\e{0}%
\e{0}%
\e{0}%
\e{0}%
\e{0}%
\e{0}%
\e{0}%
\e{0}%
\e{0}%
\e{0}%
\e{0}%
\e{0}%
\eol}\vss}\rg%
%
%
\rx{\vss\hfull{%
\rlx{\hss{$700_{xx}$}}\cg%
\e{0}%
\e{0}%
\e{0}%
\e{0}%
\e{0}%
\e{0}%
\e{0}%
\e{0}%
\e{0}%
\e{0}%
\e{0}%
\e{0}%
\e{0}%
\e{0}%
\e{0}%
\e{0}%
\e{0}%
\e{0}%
\e{0}%
\e{0}%
\e{0}%
\e{0}%
\e{0}%
\e{0}%
\e{0}%
\eol}\vss}\rg%
%
%
\rx{\vss\hfull{%
\rlx{\hss{$972_{x}$}}\cg%
\e{0}%
\e{0}%
\e{0}%
\e{0}%
\e{0}%
\e{0}%
\e{0}%
\e{0}%
\e{0}%
\e{0}%
\e{0}%
\e{0}%
\e{0}%
\e{0}%
\e{0}%
\e{0}%
\e{0}%
\e{0}%
\e{0}%
\e{0}%
\e{0}%
\e{0}%
\e{0}%
\e{0}%
\e{0}%
\eol}\vss}\rg%
%
%
\rx{\vss\hfull{%
\rlx{\hss{$4096_{x}$}}\cg%
\e{0}%
\e{0}%
\e{0}%
\e{0}%
\e{0}%
\e{0}%
\e{0}%
\e{0}%
\e{0}%
\e{0}%
\e{0}%
\e{0}%
\e{0}%
\e{0}%
\e{0}%
\e{0}%
\e{0}%
\e{0}%
\e{0}%
\e{0}%
\e{0}%
\e{0}%
\e{0}%
\e{0}%
\e{0}%
\eol}\vss}\rg%
%
%
\rx{\vss\hfull{%
\rlx{\hss{$4200_{x}$}}\cg%
\e{0}%
\e{0}%
\e{0}%
\e{0}%
\e{0}%
\e{0}%
\e{0}%
\e{0}%
\e{0}%
\e{0}%
\e{0}%
\e{0}%
\e{0}%
\e{0}%
\e{0}%
\e{0}%
\e{0}%
\e{0}%
\e{0}%
\e{0}%
\e{0}%
\e{0}%
\e{0}%
\e{0}%
\e{0}%
\eol}\vss}\rg%
%
%
\rx{\vss\hfull{%
\rlx{\hss{$2240_{x}$}}\cg%
\e{0}%
\e{0}%
\e{0}%
\e{0}%
\e{0}%
\e{0}%
\e{0}%
\e{0}%
\e{0}%
\e{0}%
\e{0}%
\e{0}%
\e{0}%
\e{0}%
\e{0}%
\e{0}%
\e{0}%
\e{0}%
\e{0}%
\e{0}%
\e{0}%
\e{0}%
\e{0}%
\e{0}%
\e{0}%
\eol}\vss}\rg%
%
%
\rx{\vss\hfull{%
\rlx{\hss{$2835_{x}$}}\cg%
\e{0}%
\e{0}%
\e{0}%
\e{0}%
\e{0}%
\e{0}%
\e{0}%
\e{0}%
\e{0}%
\e{0}%
\e{0}%
\e{0}%
\e{0}%
\e{0}%
\e{0}%
\e{0}%
\e{0}%
\e{0}%
\e{0}%
\e{0}%
\e{0}%
\e{0}%
\e{0}%
\e{0}%
\e{0}%
\eol}\vss}\rg%
%
%
\rx{\vss\hfull{%
\rlx{\hss{$6075_{x}$}}\cg%
\e{0}%
\e{0}%
\e{0}%
\e{0}%
\e{0}%
\e{0}%
\e{0}%
\e{0}%
\e{0}%
\e{0}%
\e{0}%
\e{0}%
\e{0}%
\e{0}%
\e{0}%
\e{0}%
\e{0}%
\e{0}%
\e{0}%
\e{0}%
\e{0}%
\e{0}%
\e{0}%
\e{0}%
\e{0}%
\eol}\vss}\rg%
%
%
\rx{\vss\hfull{%
\rlx{\hss{$3200_{x}$}}\cg%
\e{0}%
\e{1}%
\e{0}%
\e{0}%
\e{0}%
\e{0}%
\e{0}%
\e{0}%
\e{0}%
\e{0}%
\e{0}%
\e{0}%
\e{0}%
\e{0}%
\e{0}%
\e{0}%
\e{0}%
\e{0}%
\e{0}%
\e{0}%
\e{0}%
\e{0}%
\e{0}%
\e{0}%
\e{0}%
\eol}\vss}\rg%
%
%
\rx{\vss\hfull{%
\rlx{\hss{$70_{y}$}}\cg%
\e{0}%
\e{0}%
\e{0}%
\e{0}%
\e{0}%
\e{0}%
\e{0}%
\e{0}%
\e{0}%
\e{0}%
\e{0}%
\e{0}%
\e{0}%
\e{0}%
\e{0}%
\e{0}%
\e{0}%
\e{0}%
\e{0}%
\e{0}%
\e{0}%
\e{0}%
\e{0}%
\e{0}%
\e{0}%
\eol}\vss}\rg%
%
%
\rx{\vss\hfull{%
\rlx{\hss{$1134_{y}$}}\cg%
\e{0}%
\e{0}%
\e{0}%
\e{1}%
\e{0}%
\e{0}%
\e{0}%
\e{0}%
\e{0}%
\e{0}%
\e{0}%
\e{0}%
\e{0}%
\e{0}%
\e{0}%
\e{0}%
\e{0}%
\e{0}%
\e{0}%
\e{0}%
\e{0}%
\e{0}%
\e{0}%
\e{0}%
\e{0}%
\eol}\vss}\rg%
%
%
\rx{\vss\hfull{%
\rlx{\hss{$1680_{y}$}}\cg%
\e{0}%
\e{0}%
\e{0}%
\e{0}%
\e{0}%
\e{0}%
\e{0}%
\e{0}%
\e{0}%
\e{0}%
\e{0}%
\e{0}%
\e{0}%
\e{0}%
\e{0}%
\e{0}%
\e{0}%
\e{0}%
\e{0}%
\e{0}%
\e{0}%
\e{0}%
\e{0}%
\e{0}%
\e{0}%
\eol}\vss}\rg%
%
%
\rx{\vss\hfull{%
\rlx{\hss{$168_{y}$}}\cg%
\e{0}%
\e{0}%
\e{0}%
\e{0}%
\e{0}%
\e{0}%
\e{0}%
\e{0}%
\e{0}%
\e{0}%
\e{0}%
\e{0}%
\e{0}%
\e{0}%
\e{0}%
\e{0}%
\e{0}%
\e{0}%
\e{0}%
\e{0}%
\e{0}%
\e{0}%
\e{0}%
\e{0}%
\e{0}%
\eol}\vss}\rg%
%
%
\rx{\vss\hfull{%
\rlx{\hss{$420_{y}$}}\cg%
\e{0}%
\e{0}%
\e{1}%
\e{0}%
\e{0}%
\e{0}%
\e{0}%
\e{0}%
\e{0}%
\e{0}%
\e{0}%
\e{0}%
\e{0}%
\e{0}%
\e{0}%
\e{0}%
\e{0}%
\e{0}%
\e{0}%
\e{0}%
\e{0}%
\e{0}%
\e{0}%
\e{0}%
\e{0}%
\eol}\vss}\rg%
%
%
\rx{\vss\hfull{%
\rlx{\hss{$3150_{y}$}}\cg%
\e{0}%
\e{0}%
\e{1}%
\e{1}%
\e{0}%
\e{0}%
\e{0}%
\e{0}%
\e{0}%
\e{0}%
\e{0}%
\e{0}%
\e{0}%
\e{0}%
\e{0}%
\e{0}%
\e{0}%
\e{0}%
\e{0}%
\e{0}%
\e{0}%
\e{0}%
\e{0}%
\e{0}%
\e{0}%
\eol}\vss}\rg%
%
%
\rx{\vss\hfull{%
\rlx{\hss{$4200_{y}$}}\cg%
\e{0}%
\e{1}%
\e{0}%
\e{0}%
\e{1}%
\e{0}%
\e{0}%
\e{0}%
\e{0}%
\e{0}%
\e{0}%
\e{0}%
\e{0}%
\e{0}%
\e{0}%
\e{0}%
\e{0}%
\e{0}%
\e{0}%
\e{0}%
\e{0}%
\e{0}%
\e{0}%
\e{0}%
\e{0}%
\eol}\vss}\rg%
\eop
\eject
\tablecont%
%
%
%
%
%
%
\rowpts=18 true pt%
\colpts=18 true pt%
\rowlabpts=40 true pt%
\collabpts=70 true pt%
\clx{\vss\hfull{%
\rlx{\hss{$ $}}\cg%
\cx{\hskip 16 true pt\flip{$[{3}{1^{3}}:{1^{2}}]$}\hss}\cg%
\cx{\hskip 16 true pt\flip{$[{2^{3}}:{2}]$}\hss}\cg%
\cx{\hskip 16 true pt\flip{$[{2^{3}}:{1^{2}}]$}\hss}\cg%
\cx{\hskip 16 true pt\flip{$[{2^{2}}{1^{2}}:{2}]$}\hss}\cg%
\cx{\hskip 16 true pt\flip{$[{2^{2}}{1^{2}}:{1^{2}}]$}\hss}\cg%
\cx{\hskip 16 true pt\flip{$[{2}{1^{4}}:{2}]$}\hss}\cg%
\cx{\hskip 16 true pt\flip{$[{2}{1^{4}}:{1^{2}}]$}\hss}\cg%
\cx{\hskip 16 true pt\flip{$[{1^{6}}:{2}]$}\hss}\cg%
\cx{\hskip 16 true pt\flip{$[{1^{6}}:{1^{2}}]$}\hss}\cg%
\cx{\hskip 16 true pt\flip{$[{5}:{3}]$}\hss}\cg%
\cx{\hskip 16 true pt\flip{$[{5}:{2}{1}]$}\hss}\cg%
\cx{\hskip 16 true pt\flip{$[{5}:{1^{3}}]$}\hss}\cg%
\cx{\hskip 16 true pt\flip{$[{4}{1}:{3}]$}\hss}\cg%
\cx{\hskip 16 true pt\flip{$[{4}{1}:{2}{1}]$}\hss}\cg%
\cx{\hskip 16 true pt\flip{$[{4}{1}:{1^{3}}]$}\hss}\cg%
\cx{\hskip 16 true pt\flip{$[{3}{2}:{3}]$}\hss}\cg%
\cx{\hskip 16 true pt\flip{$[{3}{2}:{2}{1}]$}\hss}\cg%
\cx{\hskip 16 true pt\flip{$[{3}{2}:{1^{3}}]$}\hss}\cg%
\cx{\hskip 16 true pt\flip{$[{3}{1^{2}}:{3}]$}\hss}\cg%
\cx{\hskip 16 true pt\flip{$[{3}{1^{2}}:{2}{1}]$}\hss}\cg%
\cx{\hskip 16 true pt\flip{$[{3}{1^{2}}:{1^{3}}]$}\hss}\cg%
\cx{\hskip 16 true pt\flip{$[{2^{2}}{1}:{3}]$}\hss}\cg%
\cx{\hskip 16 true pt\flip{$[{2^{2}}{1}:{2}{1}]$}\hss}\cg%
\cx{\hskip 16 true pt\flip{$[{2^{2}}{1}:{1^{3}}]$}\hss}\cg%
\cx{\hskip 16 true pt\flip{$[{2}{1^{3}}:{3}]$}\hss}\cg%
\eol}}\rg%
%
%
\rx{\vss\hfull{%
\rlx{\hss{$2688_{y}$}}\cg%
\e{0}%
\e{0}%
\e{0}%
\e{0}%
\e{0}%
\e{0}%
\e{0}%
\e{0}%
\e{0}%
\e{0}%
\e{0}%
\e{0}%
\e{0}%
\e{0}%
\e{0}%
\e{0}%
\e{0}%
\e{0}%
\e{0}%
\e{0}%
\e{0}%
\e{0}%
\e{0}%
\e{0}%
\e{0}%
\eol}\vss}\rg%
%
%
\rx{\vss\hfull{%
\rlx{\hss{$2100_{y}$}}\cg%
\e{0}%
\e{0}%
\e{0}%
\e{0}%
\e{0}%
\e{0}%
\e{0}%
\e{0}%
\e{0}%
\e{0}%
\e{0}%
\e{0}%
\e{0}%
\e{0}%
\e{0}%
\e{0}%
\e{0}%
\e{0}%
\e{0}%
\e{0}%
\e{0}%
\e{0}%
\e{0}%
\e{0}%
\e{0}%
\eol}\vss}\rg%
%
%
\rx{\vss\hfull{%
\rlx{\hss{$1400_{y}$}}\cg%
\e{1}%
\e{0}%
\e{0}%
\e{0}%
\e{0}%
\e{0}%
\e{0}%
\e{0}%
\e{0}%
\e{0}%
\e{0}%
\e{0}%
\e{0}%
\e{0}%
\e{0}%
\e{0}%
\e{0}%
\e{0}%
\e{0}%
\e{0}%
\e{0}%
\e{0}%
\e{0}%
\e{0}%
\e{0}%
\eol}\vss}\rg%
%
%
\rx{\vss\hfull{%
\rlx{\hss{$4536_{y}$}}\cg%
\e{1}%
\e{0}%
\e{0}%
\e{0}%
\e{0}%
\e{0}%
\e{0}%
\e{0}%
\e{0}%
\e{0}%
\e{0}%
\e{0}%
\e{0}%
\e{0}%
\e{0}%
\e{0}%
\e{0}%
\e{0}%
\e{0}%
\e{0}%
\e{0}%
\e{0}%
\e{0}%
\e{0}%
\e{0}%
\eol}\vss}\rg%
%
%
\rx{\vss\hfull{%
\rlx{\hss{$5670_{y}$}}\cg%
\e{1}%
\e{0}%
\e{0}%
\e{1}%
\e{0}%
\e{0}%
\e{0}%
\e{0}%
\e{0}%
\e{0}%
\e{0}%
\e{0}%
\e{0}%
\e{0}%
\e{0}%
\e{0}%
\e{0}%
\e{0}%
\e{0}%
\e{0}%
\e{0}%
\e{0}%
\e{0}%
\e{0}%
\e{0}%
\eol}\vss}\rg%
%
%
\rx{\vss\hfull{%
\rlx{\hss{$4480_{y}$}}\cg%
\e{1}%
\e{0}%
\e{1}%
\e{1}%
\e{0}%
\e{0}%
\e{0}%
\e{0}%
\e{0}%
\e{0}%
\e{0}%
\e{0}%
\e{0}%
\e{0}%
\e{0}%
\e{0}%
\e{0}%
\e{0}%
\e{0}%
\e{0}%
\e{0}%
\e{0}%
\e{0}%
\e{0}%
\e{0}%
\eol}\vss}\rg%
%
%
\rx{\vss\hfull{%
\rlx{\hss{$8_{z}$}}\cg%
\e{0}%
\e{0}%
\e{0}%
\e{0}%
\e{0}%
\e{0}%
\e{0}%
\e{0}%
\e{0}%
\e{0}%
\e{0}%
\e{0}%
\e{0}%
\e{0}%
\e{0}%
\e{0}%
\e{0}%
\e{0}%
\e{0}%
\e{0}%
\e{0}%
\e{0}%
\e{0}%
\e{0}%
\e{0}%
\eol}\vss}\rg%
%
%
\rx{\vss\hfull{%
\rlx{\hss{$56_{z}$}}\cg%
\e{0}%
\e{0}%
\e{0}%
\e{0}%
\e{0}%
\e{0}%
\e{0}%
\e{0}%
\e{0}%
\e{0}%
\e{0}%
\e{1}%
\e{0}%
\e{0}%
\e{0}%
\e{0}%
\e{0}%
\e{0}%
\e{0}%
\e{0}%
\e{0}%
\e{0}%
\e{0}%
\e{0}%
\e{0}%
\eol}\vss}\rg%
%
%
\rx{\vss\hfull{%
\rlx{\hss{$160_{z}$}}\cg%
\e{0}%
\e{0}%
\e{0}%
\e{0}%
\e{0}%
\e{0}%
\e{0}%
\e{0}%
\e{0}%
\e{0}%
\e{1}%
\e{0}%
\e{0}%
\e{0}%
\e{0}%
\e{0}%
\e{0}%
\e{0}%
\e{0}%
\e{0}%
\e{0}%
\e{0}%
\e{0}%
\e{0}%
\e{0}%
\eol}\vss}\rg%
%
%
\rx{\vss\hfull{%
\rlx{\hss{$112_{z}$}}\cg%
\e{0}%
\e{0}%
\e{0}%
\e{0}%
\e{0}%
\e{0}%
\e{0}%
\e{0}%
\e{0}%
\e{1}%
\e{0}%
\e{0}%
\e{0}%
\e{0}%
\e{0}%
\e{0}%
\e{0}%
\e{0}%
\e{0}%
\e{0}%
\e{0}%
\e{0}%
\e{0}%
\e{0}%
\e{0}%
\eol}\vss}\rg%
%
%
\rx{\vss\hfull{%
\rlx{\hss{$840_{z}$}}\cg%
\e{0}%
\e{0}%
\e{0}%
\e{0}%
\e{0}%
\e{0}%
\e{0}%
\e{0}%
\e{0}%
\e{0}%
\e{0}%
\e{0}%
\e{0}%
\e{1}%
\e{0}%
\e{0}%
\e{0}%
\e{0}%
\e{0}%
\e{0}%
\e{0}%
\e{1}%
\e{0}%
\e{0}%
\e{0}%
\eol}\vss}\rg%
%
%
\rx{\vss\hfull{%
\rlx{\hss{$1296_{z}$}}\cg%
\e{0}%
\e{0}%
\e{0}%
\e{0}%
\e{0}%
\e{0}%
\e{0}%
\e{0}%
\e{0}%
\e{0}%
\e{1}%
\e{1}%
\e{0}%
\e{1}%
\e{1}%
\e{0}%
\e{0}%
\e{0}%
\e{1}%
\e{0}%
\e{0}%
\e{0}%
\e{0}%
\e{0}%
\e{0}%
\eol}\vss}\rg%
%
%
\rx{\vss\hfull{%
\rlx{\hss{$1400_{z}$}}\cg%
\e{0}%
\e{0}%
\e{0}%
\e{0}%
\e{0}%
\e{0}%
\e{0}%
\e{0}%
\e{0}%
\e{1}%
\e{1}%
\e{0}%
\e{1}%
\e{1}%
\e{0}%
\e{1}%
\e{0}%
\e{0}%
\e{0}%
\e{0}%
\e{0}%
\e{0}%
\e{0}%
\e{0}%
\e{0}%
\eol}\vss}\rg%
%
%
\rx{\vss\hfull{%
\rlx{\hss{$1008_{z}$}}\cg%
\e{0}%
\e{0}%
\e{0}%
\e{0}%
\e{0}%
\e{0}%
\e{0}%
\e{0}%
\e{0}%
\e{0}%
\e{1}%
\e{1}%
\e{1}%
\e{1}%
\e{0}%
\e{0}%
\e{0}%
\e{0}%
\e{0}%
\e{0}%
\e{0}%
\e{0}%
\e{0}%
\e{0}%
\e{0}%
\eol}\vss}\rg%
%
%
\rx{\vss\hfull{%
\rlx{\hss{$560_{z}$}}\cg%
\e{0}%
\e{0}%
\e{0}%
\e{0}%
\e{0}%
\e{0}%
\e{0}%
\e{0}%
\e{0}%
\e{1}%
\e{1}%
\e{0}%
\e{1}%
\e{0}%
\e{0}%
\e{0}%
\e{0}%
\e{0}%
\e{0}%
\e{0}%
\e{0}%
\e{0}%
\e{0}%
\e{0}%
\e{0}%
\eol}\vss}\rg%
%
%
\rx{\vss\hfull{%
\rlx{\hss{$1400_{zz}$}}\cg%
\e{0}%
\e{0}%
\e{0}%
\e{0}%
\e{0}%
\e{0}%
\e{0}%
\e{0}%
\e{0}%
\e{1}%
\e{0}%
\e{0}%
\e{1}%
\e{0}%
\e{0}%
\e{1}%
\e{1}%
\e{0}%
\e{0}%
\e{0}%
\e{0}%
\e{0}%
\e{0}%
\e{0}%
\e{0}%
\eol}\vss}\rg%
%
%
\rx{\vss\hfull{%
\rlx{\hss{$4200_{z}$}}\cg%
\e{0}%
\e{0}%
\e{0}%
\e{0}%
\e{0}%
\e{0}%
\e{0}%
\e{0}%
\e{0}%
\e{0}%
\e{0}%
\e{0}%
\e{1}%
\e{1}%
\e{0}%
\e{0}%
\e{2}%
\e{1}%
\e{1}%
\e{1}%
\e{0}%
\e{0}%
\e{1}%
\e{0}%
\e{0}%
\eol}\vss}\rg%
%
%
\rx{\vss\hfull{%
\rlx{\hss{$400_{z}$}}\cg%
\e{0}%
\e{0}%
\e{0}%
\e{0}%
\e{0}%
\e{0}%
\e{0}%
\e{0}%
\e{0}%
\e{1}%
\e{0}%
\e{0}%
\e{1}%
\e{0}%
\e{0}%
\e{0}%
\e{0}%
\e{0}%
\e{0}%
\e{0}%
\e{0}%
\e{0}%
\e{0}%
\e{0}%
\e{0}%
\eol}\vss}\rg%
%
%
\rx{\vss\hfull{%
\rlx{\hss{$3240_{z}$}}\cg%
\e{0}%
\e{0}%
\e{0}%
\e{0}%
\e{0}%
\e{0}%
\e{0}%
\e{0}%
\e{0}%
\e{1}%
\e{1}%
\e{0}%
\e{2}%
\e{2}%
\e{0}%
\e{1}%
\e{1}%
\e{0}%
\e{1}%
\e{0}%
\e{0}%
\e{0}%
\e{0}%
\e{0}%
\e{0}%
\eol}\vss}\rg%
%
%
\rx{\vss\hfull{%
\rlx{\hss{$4536_{z}$}}\cg%
\e{0}%
\e{0}%
\e{0}%
\e{0}%
\e{0}%
\e{0}%
\e{0}%
\e{0}%
\e{0}%
\e{0}%
\e{0}%
\e{0}%
\e{1}%
\e{1}%
\e{0}%
\e{2}%
\e{2}%
\e{0}%
\e{0}%
\e{1}%
\e{0}%
\e{1}%
\e{1}%
\e{0}%
\e{0}%
\eol}\vss}\rg%
%
%
\rx{\vss\hfull{%
\rlx{\hss{$2400_{z}$}}\cg%
\e{0}%
\e{0}%
\e{0}%
\e{0}%
\e{0}%
\e{0}%
\e{0}%
\e{0}%
\e{0}%
\e{0}%
\e{0}%
\e{1}%
\e{0}%
\e{1}%
\e{1}%
\e{0}%
\e{0}%
\e{1}%
\e{1}%
\e{1}%
\e{0}%
\e{0}%
\e{0}%
\e{0}%
\e{1}%
\eol}\vss}\rg%
%
%
\rx{\vss\hfull{%
\rlx{\hss{$3360_{z}$}}\cg%
\e{0}%
\e{0}%
\e{0}%
\e{0}%
\e{0}%
\e{0}%
\e{0}%
\e{0}%
\e{0}%
\e{0}%
\e{1}%
\e{0}%
\e{1}%
\e{1}%
\e{0}%
\e{1}%
\e{1}%
\e{1}%
\e{1}%
\e{1}%
\e{0}%
\e{0}%
\e{0}%
\e{0}%
\e{0}%
\eol}\vss}\rg%
%
%
\rx{\vss\hfull{%
\rlx{\hss{$2800_{z}$}}\cg%
\e{0}%
\e{0}%
\e{0}%
\e{0}%
\e{0}%
\e{0}%
\e{0}%
\e{0}%
\e{0}%
\e{0}%
\e{0}%
\e{1}%
\e{1}%
\e{1}%
\e{1}%
\e{0}%
\e{1}%
\e{0}%
\e{1}%
\e{1}%
\e{0}%
\e{0}%
\e{0}%
\e{0}%
\e{0}%
\eol}\vss}\rg%
%
%
\rx{\vss\hfull{%
\rlx{\hss{$4096_{z}$}}\cg%
\e{0}%
\e{0}%
\e{0}%
\e{0}%
\e{0}%
\e{0}%
\e{0}%
\e{0}%
\e{0}%
\e{0}%
\e{1}%
\e{0}%
\e{1}%
\e{2}%
\e{1}%
\e{1}%
\e{1}%
\e{0}%
\e{1}%
\e{1}%
\e{0}%
\e{1}%
\e{0}%
\e{0}%
\e{0}%
\eol}\vss}\rg%
%
%
\rx{\vss\hfull{%
\rlx{\hss{$5600_{z}$}}\cg%
\e{0}%
\e{0}%
\e{0}%
\e{0}%
\e{0}%
\e{0}%
\e{0}%
\e{0}%
\e{0}%
\e{0}%
\e{0}%
\e{0}%
\e{0}%
\e{2}%
\e{1}%
\e{1}%
\e{1}%
\e{0}%
\e{1}%
\e{2}%
\e{1}%
\e{1}%
\e{1}%
\e{0}%
\e{1}%
\eol}\vss}\rg%
%
%
\rx{\vss\hfull{%
\rlx{\hss{$448_{z}$}}\cg%
\e{0}%
\e{0}%
\e{0}%
\e{0}%
\e{0}%
\e{0}%
\e{0}%
\e{0}%
\e{0}%
\e{1}%
\e{0}%
\e{0}%
\e{0}%
\e{0}%
\e{0}%
\e{1}%
\e{0}%
\e{0}%
\e{0}%
\e{0}%
\e{0}%
\e{0}%
\e{0}%
\e{0}%
\e{0}%
\eol}\vss}\rg%
%
%
\rx{\vss\hfull{%
\rlx{\hss{$448_{w}$}}\cg%
\e{0}%
\e{0}%
\e{0}%
\e{0}%
\e{0}%
\e{0}%
\e{0}%
\e{0}%
\e{0}%
\e{0}%
\e{0}%
\e{0}%
\e{0}%
\e{0}%
\e{1}%
\e{0}%
\e{0}%
\e{0}%
\e{0}%
\e{0}%
\e{0}%
\e{0}%
\e{0}%
\e{0}%
\e{1}%
\eol}\vss}\rg%
%
%
\rx{\vss\hfull{%
\rlx{\hss{$1344_{w}$}}\cg%
\e{0}%
\e{0}%
\e{0}%
\e{0}%
\e{0}%
\e{0}%
\e{0}%
\e{0}%
\e{0}%
\e{0}%
\e{0}%
\e{0}%
\e{0}%
\e{0}%
\e{0}%
\e{0}%
\e{1}%
\e{0}%
\e{0}%
\e{0}%
\e{0}%
\e{0}%
\e{1}%
\e{0}%
\e{0}%
\eol}\vss}\rg%
%
%
\rx{\vss\hfull{%
\rlx{\hss{$5600_{w}$}}\cg%
\e{0}%
\e{0}%
\e{0}%
\e{0}%
\e{0}%
\e{0}%
\e{0}%
\e{0}%
\e{0}%
\e{0}%
\e{0}%
\e{0}%
\e{0}%
\e{1}%
\e{1}%
\e{0}%
\e{1}%
\e{1}%
\e{1}%
\e{2}%
\e{1}%
\e{1}%
\e{1}%
\e{0}%
\e{1}%
\eol}\vss}\rg%
%
%
\rx{\vss\hfull{%
\rlx{\hss{$2016_{w}$}}\cg%
\e{0}%
\e{0}%
\e{0}%
\e{0}%
\e{0}%
\e{0}%
\e{0}%
\e{0}%
\e{0}%
\e{0}%
\e{0}%
\e{0}%
\e{0}%
\e{0}%
\e{0}%
\e{1}%
\e{1}%
\e{0}%
\e{0}%
\e{0}%
\e{0}%
\e{0}%
\e{1}%
\e{1}%
\e{0}%
\eol}\vss}\rg%
%
%
\rx{\vss\hfull{%
\rlx{\hss{$7168_{w}$}}\cg%
\e{0}%
\e{0}%
\e{0}%
\e{0}%
\e{0}%
\e{0}%
\e{0}%
\e{0}%
\e{0}%
\e{0}%
\e{0}%
\e{0}%
\e{0}%
\e{1}%
\e{0}%
\e{1}%
\e{2}%
\e{1}%
\e{1}%
\e{2}%
\e{1}%
\e{1}%
\e{2}%
\e{1}%
\e{0}%
\eol}\vss}\rg%
\eop
\eject
\tablecont%
%
%
%
%
%
%
\rowpts=18 true pt%
\colpts=18 true pt%
\rowlabpts=40 true pt%
\collabpts=70 true pt%
\clx{\vss\hfull{%
\rlx{\hss{$ $}}\cg%
\cx{\hskip 16 true pt\flip{$[{2}{1^{3}}:{2}{1}]$}\hss}\cg%
\cx{\hskip 16 true pt\flip{$[{2}{1^{3}}:{1^{3}}]$}\hss}\cg%
\cx{\hskip 16 true pt\flip{$[{1^{5}}:{3}]$}\hss}\cg%
\cx{\hskip 16 true pt\flip{$[{1^{5}}:{2}{1}]$}\hss}\cg%
\cx{\hskip 16 true pt\flip{$[{1^{5}}:{1^{3}}]$}\hss}\cg%
\cx{\hskip 16 true pt\flip{$[{4}:{4}]^{+}$}\hss}\cg%
\cx{\hskip 16 true pt\flip{$[{4}:{4}]^{-}$}\hss}\cg%
\cx{\hskip 16 true pt\flip{$[{4}:{3}{1}]$}\hss}\cg%
\cx{\hskip 16 true pt\flip{$[{4}:{2^{2}}]$}\hss}\cg%
\cx{\hskip 16 true pt\flip{$[{4}:{2}{1^{2}}]$}\hss}\cg%
\cx{\hskip 16 true pt\flip{$[{4}:{1^{4}}]$}\hss}\cg%
\cx{\hskip 16 true pt\flip{$[{3}{1}:{3}{1}]^{+}$}\hss}\cg%
\cx{\hskip 16 true pt\flip{$[{3}{1}:{3}{1}]^{-}$}\hss}\cg%
\cx{\hskip 16 true pt\flip{$[{3}{1}:{2^{2}}]$}\hss}\cg%
\cx{\hskip 16 true pt\flip{$[{3}{1}:{2}{1^{2}}]$}\hss}\cg%
\cx{\hskip 16 true pt\flip{$[{3}{1}:{1^{4}}]$}\hss}\cg%
\cx{\hskip 16 true pt\flip{$[{2^{2}}:{2^{2}}]^{+}$}\hss}\cg%
\cx{\hskip 16 true pt\flip{$[{2^{2}}:{2^{2}}]^{-}$}\hss}\cg%
\cx{\hskip 16 true pt\flip{$[{2^{2}}:{2}{1^{2}}]$}\hss}\cg%
\cx{\hskip 16 true pt\flip{$[{2^{2}}:{1^{4}}]$}\hss}\cg%
\cx{\hskip 16 true pt\flip{$[{2}{1^{2}}:{2}{1^{2}}]^{+}$}\hss}\cg%
\cx{\hskip 16 true pt\flip{$[{2}{1^{2}}:{2}{1^{2}}]^{-}$}\hss}\cg%
\cx{\hskip 16 true pt\flip{$[{2}{1^{2}}:{1^{4}}]$}\hss}\cg%
\cx{\hskip 16 true pt\flip{$[{1^{4}}:{1^{4}}]^{+}$}\hss}\cg%
\cx{\hskip 16 true pt\flip{$[{1^{4}}:{1^{4}}]^{-}$}\hss}\cg%
\eol}}\rg%
%
%
\rx{\vss\hfull{%
\rlx{\hss{$1_{x}$}}\cg%
\e{0}%
\e{0}%
\e{0}%
\e{0}%
\e{0}%
\e{0}%
\e{0}%
\e{0}%
\e{0}%
\e{0}%
\e{0}%
\e{0}%
\e{0}%
\e{0}%
\e{0}%
\e{0}%
\e{0}%
\e{0}%
\e{0}%
\e{0}%
\e{0}%
\e{0}%
\e{0}%
\e{0}%
\e{0}%
\eol}\vss}\rg%
%
%
\rx{\vss\hfull{%
\rlx{\hss{$28_{x}$}}\cg%
\e{0}%
\e{0}%
\e{0}%
\e{0}%
\e{0}%
\e{0}%
\e{0}%
\e{0}%
\e{0}%
\e{0}%
\e{0}%
\e{0}%
\e{0}%
\e{0}%
\e{0}%
\e{0}%
\e{0}%
\e{0}%
\e{0}%
\e{0}%
\e{0}%
\e{0}%
\e{0}%
\e{0}%
\e{0}%
\eol}\vss}\rg%
%
%
\rx{\vss\hfull{%
\rlx{\hss{$35_{x}$}}\cg%
\e{0}%
\e{0}%
\e{0}%
\e{0}%
\e{0}%
\e{0}%
\e{0}%
\e{0}%
\e{0}%
\e{0}%
\e{0}%
\e{0}%
\e{0}%
\e{0}%
\e{0}%
\e{0}%
\e{0}%
\e{0}%
\e{0}%
\e{0}%
\e{0}%
\e{0}%
\e{0}%
\e{0}%
\e{0}%
\eol}\vss}\rg%
%
%
\rx{\vss\hfull{%
\rlx{\hss{$84_{x}$}}\cg%
\e{0}%
\e{0}%
\e{0}%
\e{0}%
\e{0}%
\e{1}%
\e{0}%
\e{0}%
\e{0}%
\e{0}%
\e{0}%
\e{0}%
\e{0}%
\e{0}%
\e{0}%
\e{0}%
\e{0}%
\e{0}%
\e{0}%
\e{0}%
\e{0}%
\e{0}%
\e{0}%
\e{0}%
\e{0}%
\eol}\vss}\rg%
%
%
\rx{\vss\hfull{%
\rlx{\hss{$50_{x}$}}\cg%
\e{0}%
\e{0}%
\e{0}%
\e{0}%
\e{0}%
\e{1}%
\e{0}%
\e{0}%
\e{0}%
\e{0}%
\e{0}%
\e{0}%
\e{0}%
\e{0}%
\e{0}%
\e{0}%
\e{0}%
\e{0}%
\e{0}%
\e{0}%
\e{0}%
\e{0}%
\e{0}%
\e{0}%
\e{0}%
\eol}\vss}\rg%
%
%
\rx{\vss\hfull{%
\rlx{\hss{$350_{x}$}}\cg%
\e{0}%
\e{0}%
\e{0}%
\e{0}%
\e{0}%
\e{0}%
\e{0}%
\e{0}%
\e{0}%
\e{1}%
\e{0}%
\e{0}%
\e{0}%
\e{0}%
\e{0}%
\e{0}%
\e{0}%
\e{0}%
\e{0}%
\e{0}%
\e{0}%
\e{0}%
\e{0}%
\e{0}%
\e{0}%
\eol}\vss}\rg%
%
%
\rx{\vss\hfull{%
\rlx{\hss{$300_{x}$}}\cg%
\e{0}%
\e{0}%
\e{0}%
\e{0}%
\e{0}%
\e{0}%
\e{0}%
\e{0}%
\e{1}%
\e{0}%
\e{0}%
\e{0}%
\e{0}%
\e{0}%
\e{0}%
\e{0}%
\e{0}%
\e{0}%
\e{0}%
\e{0}%
\e{0}%
\e{0}%
\e{0}%
\e{0}%
\e{0}%
\eol}\vss}\rg%
%
%
\rx{\vss\hfull{%
\rlx{\hss{$567_{x}$}}\cg%
\e{0}%
\e{0}%
\e{0}%
\e{0}%
\e{0}%
\e{0}%
\e{0}%
\e{1}%
\e{0}%
\e{0}%
\e{0}%
\e{0}%
\e{0}%
\e{0}%
\e{0}%
\e{0}%
\e{0}%
\e{0}%
\e{0}%
\e{0}%
\e{0}%
\e{0}%
\e{0}%
\e{0}%
\e{0}%
\eol}\vss}\rg%
%
%
\rx{\vss\hfull{%
\rlx{\hss{$210_{x}$}}\cg%
\e{0}%
\e{0}%
\e{0}%
\e{0}%
\e{0}%
\e{0}%
\e{1}%
\e{0}%
\e{0}%
\e{0}%
\e{0}%
\e{0}%
\e{0}%
\e{0}%
\e{0}%
\e{0}%
\e{0}%
\e{0}%
\e{0}%
\e{0}%
\e{0}%
\e{0}%
\e{0}%
\e{0}%
\e{0}%
\eol}\vss}\rg%
%
%
\rx{\vss\hfull{%
\rlx{\hss{$840_{x}$}}\cg%
\e{0}%
\e{0}%
\e{0}%
\e{0}%
\e{0}%
\e{0}%
\e{0}%
\e{0}%
\e{0}%
\e{0}%
\e{0}%
\e{0}%
\e{0}%
\e{1}%
\e{0}%
\e{0}%
\e{0}%
\e{0}%
\e{0}%
\e{0}%
\e{0}%
\e{0}%
\e{0}%
\e{0}%
\e{0}%
\eol}\vss}\rg%
%
%
\rx{\vss\hfull{%
\rlx{\hss{$700_{x}$}}\cg%
\e{0}%
\e{0}%
\e{0}%
\e{0}%
\e{0}%
\e{0}%
\e{1}%
\e{1}%
\e{0}%
\e{0}%
\e{0}%
\e{0}%
\e{0}%
\e{0}%
\e{0}%
\e{0}%
\e{0}%
\e{0}%
\e{0}%
\e{0}%
\e{0}%
\e{0}%
\e{0}%
\e{0}%
\e{0}%
\eol}\vss}\rg%
%
%
\rx{\vss\hfull{%
\rlx{\hss{$175_{x}$}}\cg%
\e{0}%
\e{0}%
\e{0}%
\e{0}%
\e{0}%
\e{0}%
\e{1}%
\e{0}%
\e{0}%
\e{0}%
\e{0}%
\e{0}%
\e{0}%
\e{0}%
\e{0}%
\e{0}%
\e{0}%
\e{0}%
\e{0}%
\e{0}%
\e{0}%
\e{0}%
\e{0}%
\e{0}%
\e{0}%
\eol}\vss}\rg%
%
%
\rx{\vss\hfull{%
\rlx{\hss{$1400_{x}$}}\cg%
\e{0}%
\e{0}%
\e{0}%
\e{0}%
\e{0}%
\e{0}%
\e{1}%
\e{1}%
\e{0}%
\e{0}%
\e{0}%
\e{0}%
\e{1}%
\e{0}%
\e{0}%
\e{0}%
\e{0}%
\e{0}%
\e{0}%
\e{0}%
\e{0}%
\e{0}%
\e{0}%
\e{0}%
\e{0}%
\eol}\vss}\rg%
%
%
\rx{\vss\hfull{%
\rlx{\hss{$1050_{x}$}}\cg%
\e{0}%
\e{0}%
\e{0}%
\e{0}%
\e{0}%
\e{1}%
\e{0}%
\e{1}%
\e{0}%
\e{0}%
\e{0}%
\e{1}%
\e{0}%
\e{0}%
\e{0}%
\e{0}%
\e{0}%
\e{0}%
\e{0}%
\e{0}%
\e{0}%
\e{0}%
\e{0}%
\e{0}%
\e{0}%
\eol}\vss}\rg%
%
%
\rx{\vss\hfull{%
\rlx{\hss{$1575_{x}$}}\cg%
\e{0}%
\e{0}%
\e{0}%
\e{0}%
\e{0}%
\e{0}%
\e{0}%
\e{1}%
\e{0}%
\e{1}%
\e{0}%
\e{0}%
\e{1}%
\e{0}%
\e{0}%
\e{0}%
\e{0}%
\e{0}%
\e{0}%
\e{0}%
\e{0}%
\e{0}%
\e{0}%
\e{0}%
\e{0}%
\eol}\vss}\rg%
%
%
\rx{\vss\hfull{%
\rlx{\hss{$1344_{x}$}}\cg%
\e{0}%
\e{0}%
\e{0}%
\e{0}%
\e{0}%
\e{1}%
\e{0}%
\e{1}%
\e{1}%
\e{0}%
\e{0}%
\e{1}%
\e{0}%
\e{0}%
\e{0}%
\e{0}%
\e{0}%
\e{0}%
\e{0}%
\e{0}%
\e{0}%
\e{0}%
\e{0}%
\e{0}%
\e{0}%
\eol}\vss}\rg%
%
%
\rx{\vss\hfull{%
\rlx{\hss{$2100_{x}$}}\cg%
\e{0}%
\e{0}%
\e{0}%
\e{0}%
\e{0}%
\e{0}%
\e{0}%
\e{0}%
\e{1}%
\e{1}%
\e{1}%
\e{1}%
\e{0}%
\e{0}%
\e{1}%
\e{0}%
\e{0}%
\e{0}%
\e{0}%
\e{0}%
\e{0}%
\e{0}%
\e{0}%
\e{0}%
\e{0}%
\eol}\vss}\rg%
%
%
\rx{\vss\hfull{%
\rlx{\hss{$2268_{x}$}}\cg%
\e{0}%
\e{0}%
\e{0}%
\e{0}%
\e{0}%
\e{0}%
\e{0}%
\e{1}%
\e{0}%
\e{1}%
\e{0}%
\e{0}%
\e{1}%
\e{1}%
\e{0}%
\e{0}%
\e{0}%
\e{0}%
\e{0}%
\e{0}%
\e{0}%
\e{0}%
\e{0}%
\e{0}%
\e{0}%
\eol}\vss}\rg%
%
%
\rx{\vss\hfull{%
\rlx{\hss{$525_{x}$}}\cg%
\e{0}%
\e{0}%
\e{0}%
\e{0}%
\e{0}%
\e{1}%
\e{0}%
\e{0}%
\e{0}%
\e{0}%
\e{0}%
\e{1}%
\e{0}%
\e{0}%
\e{0}%
\e{0}%
\e{0}%
\e{0}%
\e{0}%
\e{0}%
\e{0}%
\e{0}%
\e{0}%
\e{0}%
\e{0}%
\eol}\vss}\rg%
%
%
\rx{\vss\hfull{%
\rlx{\hss{$700_{xx}$}}\cg%
\e{0}%
\e{0}%
\e{0}%
\e{0}%
\e{0}%
\e{1}%
\e{0}%
\e{0}%
\e{0}%
\e{0}%
\e{0}%
\e{1}%
\e{0}%
\e{0}%
\e{0}%
\e{0}%
\e{1}%
\e{0}%
\e{0}%
\e{0}%
\e{0}%
\e{0}%
\e{0}%
\e{0}%
\e{0}%
\eol}\vss}\rg%
%
%
\rx{\vss\hfull{%
\rlx{\hss{$972_{x}$}}\cg%
\e{0}%
\e{0}%
\e{0}%
\e{0}%
\e{0}%
\e{1}%
\e{0}%
\e{0}%
\e{1}%
\e{0}%
\e{0}%
\e{1}%
\e{0}%
\e{0}%
\e{0}%
\e{0}%
\e{1}%
\e{0}%
\e{0}%
\e{0}%
\e{0}%
\e{0}%
\e{0}%
\e{0}%
\e{0}%
\eol}\vss}\rg%
%
%
\rx{\vss\hfull{%
\rlx{\hss{$4096_{x}$}}\cg%
\e{0}%
\e{0}%
\e{0}%
\e{0}%
\e{0}%
\e{0}%
\e{0}%
\e{1}%
\e{1}%
\e{1}%
\e{0}%
\e{1}%
\e{1}%
\e{1}%
\e{1}%
\e{0}%
\e{0}%
\e{0}%
\e{0}%
\e{0}%
\e{0}%
\e{0}%
\e{0}%
\e{0}%
\e{0}%
\eol}\vss}\rg%
%
%
\rx{\vss\hfull{%
\rlx{\hss{$4200_{x}$}}\cg%
\e{0}%
\e{0}%
\e{0}%
\e{0}%
\e{0}%
\e{0}%
\e{0}%
\e{1}%
\e{1}%
\e{0}%
\e{0}%
\e{1}%
\e{1}%
\e{1}%
\e{1}%
\e{0}%
\e{0}%
\e{1}%
\e{0}%
\e{0}%
\e{0}%
\e{0}%
\e{0}%
\e{0}%
\e{0}%
\eol}\vss}\rg%
%
%
\rx{\vss\hfull{%
\rlx{\hss{$2240_{x}$}}\cg%
\e{0}%
\e{0}%
\e{0}%
\e{0}%
\e{0}%
\e{0}%
\e{1}%
\e{1}%
\e{0}%
\e{0}%
\e{0}%
\e{0}%
\e{1}%
\e{1}%
\e{0}%
\e{0}%
\e{0}%
\e{0}%
\e{0}%
\e{0}%
\e{0}%
\e{0}%
\e{0}%
\e{0}%
\e{0}%
\eol}\vss}\rg%
%
%
\rx{\vss\hfull{%
\rlx{\hss{$2835_{x}$}}\cg%
\e{0}%
\e{0}%
\e{0}%
\e{0}%
\e{0}%
\e{0}%
\e{0}%
\e{1}%
\e{0}%
\e{0}%
\e{0}%
\e{0}%
\e{1}%
\e{1}%
\e{0}%
\e{0}%
\e{0}%
\e{0}%
\e{1}%
\e{0}%
\e{0}%
\e{0}%
\e{0}%
\e{0}%
\e{0}%
\eol}\vss}\rg%
%
%
\rx{\vss\hfull{%
\rlx{\hss{$6075_{x}$}}\cg%
\e{0}%
\e{0}%
\e{0}%
\e{0}%
\e{0}%
\e{0}%
\e{0}%
\e{1}%
\e{0}%
\e{1}%
\e{0}%
\e{2}%
\e{1}%
\e{1}%
\e{2}%
\e{0}%
\e{0}%
\e{0}%
\e{1}%
\e{0}%
\e{0}%
\e{0}%
\e{0}%
\e{0}%
\e{0}%
\eol}\vss}\rg%
%
%
\rx{\vss\hfull{%
\rlx{\hss{$3200_{x}$}}\cg%
\e{0}%
\e{0}%
\e{0}%
\e{0}%
\e{0}%
\e{0}%
\e{0}%
\e{0}%
\e{1}%
\e{0}%
\e{0}%
\e{1}%
\e{0}%
\e{1}%
\e{1}%
\e{0}%
\e{1}%
\e{0}%
\e{0}%
\e{0}%
\e{1}%
\e{0}%
\e{0}%
\e{0}%
\e{0}%
\eol}\vss}\rg%
%
%
\rx{\vss\hfull{%
\rlx{\hss{$70_{y}$}}\cg%
\e{0}%
\e{0}%
\e{0}%
\e{0}%
\e{0}%
\e{0}%
\e{0}%
\e{0}%
\e{0}%
\e{0}%
\e{1}%
\e{0}%
\e{0}%
\e{0}%
\e{0}%
\e{0}%
\e{0}%
\e{0}%
\e{0}%
\e{0}%
\e{0}%
\e{0}%
\e{0}%
\e{0}%
\e{0}%
\eol}\vss}\rg%
%
%
\rx{\vss\hfull{%
\rlx{\hss{$1134_{y}$}}\cg%
\e{0}%
\e{0}%
\e{0}%
\e{0}%
\e{0}%
\e{0}%
\e{0}%
\e{0}%
\e{0}%
\e{0}%
\e{0}%
\e{0}%
\e{0}%
\e{0}%
\e{1}%
\e{0}%
\e{0}%
\e{0}%
\e{0}%
\e{0}%
\e{0}%
\e{0}%
\e{0}%
\e{0}%
\e{0}%
\eol}\vss}\rg%
%
%
\rx{\vss\hfull{%
\rlx{\hss{$1680_{y}$}}\cg%
\e{0}%
\e{0}%
\e{0}%
\e{0}%
\e{0}%
\e{0}%
\e{0}%
\e{0}%
\e{0}%
\e{1}%
\e{1}%
\e{0}%
\e{0}%
\e{0}%
\e{1}%
\e{1}%
\e{0}%
\e{0}%
\e{0}%
\e{0}%
\e{0}%
\e{0}%
\e{0}%
\e{0}%
\e{0}%
\eol}\vss}\rg%
%
%
\rx{\vss\hfull{%
\rlx{\hss{$168_{y}$}}\cg%
\e{0}%
\e{0}%
\e{0}%
\e{0}%
\e{0}%
\e{0}%
\e{0}%
\e{0}%
\e{0}%
\e{0}%
\e{0}%
\e{0}%
\e{0}%
\e{0}%
\e{0}%
\e{0}%
\e{1}%
\e{0}%
\e{0}%
\e{0}%
\e{0}%
\e{0}%
\e{0}%
\e{0}%
\e{0}%
\eol}\vss}\rg%
%
%
\rx{\vss\hfull{%
\rlx{\hss{$420_{y}$}}\cg%
\e{0}%
\e{0}%
\e{0}%
\e{0}%
\e{0}%
\e{0}%
\e{0}%
\e{0}%
\e{0}%
\e{0}%
\e{0}%
\e{0}%
\e{0}%
\e{0}%
\e{0}%
\e{0}%
\e{0}%
\e{1}%
\e{0}%
\e{0}%
\e{0}%
\e{0}%
\e{0}%
\e{0}%
\e{0}%
\eol}\vss}\rg%
%
%
\rx{\vss\hfull{%
\rlx{\hss{$3150_{y}$}}\cg%
\e{0}%
\e{0}%
\e{0}%
\e{0}%
\e{0}%
\e{0}%
\e{0}%
\e{0}%
\e{0}%
\e{0}%
\e{0}%
\e{0}%
\e{1}%
\e{1}%
\e{0}%
\e{0}%
\e{0}%
\e{0}%
\e{1}%
\e{0}%
\e{0}%
\e{1}%
\e{0}%
\e{0}%
\e{0}%
\eol}\vss}\rg%
%
%
\rx{\vss\hfull{%
\rlx{\hss{$4200_{y}$}}\cg%
\e{0}%
\e{0}%
\e{0}%
\e{0}%
\e{0}%
\e{0}%
\e{0}%
\e{0}%
\e{0}%
\e{0}%
\e{0}%
\e{1}%
\e{0}%
\e{1}%
\e{1}%
\e{0}%
\e{1}%
\e{1}%
\e{1}%
\e{0}%
\e{1}%
\e{0}%
\e{0}%
\e{0}%
\e{0}%
\eol}\vss}\rg%
\eop
\eject
\tablecont%
%
%
%
%
%
%
\rowpts=18 true pt%
\colpts=18 true pt%
\rowlabpts=40 true pt%
\collabpts=70 true pt%
\clx{\vss\hfull{%
\rlx{\hss{$ $}}\cg%
\cx{\hskip 16 true pt\flip{$[{2}{1^{3}}:{2}{1}]$}\hss}\cg%
\cx{\hskip 16 true pt\flip{$[{2}{1^{3}}:{1^{3}}]$}\hss}\cg%
\cx{\hskip 16 true pt\flip{$[{1^{5}}:{3}]$}\hss}\cg%
\cx{\hskip 16 true pt\flip{$[{1^{5}}:{2}{1}]$}\hss}\cg%
\cx{\hskip 16 true pt\flip{$[{1^{5}}:{1^{3}}]$}\hss}\cg%
\cx{\hskip 16 true pt\flip{$[{4}:{4}]^{+}$}\hss}\cg%
\cx{\hskip 16 true pt\flip{$[{4}:{4}]^{-}$}\hss}\cg%
\cx{\hskip 16 true pt\flip{$[{4}:{3}{1}]$}\hss}\cg%
\cx{\hskip 16 true pt\flip{$[{4}:{2^{2}}]$}\hss}\cg%
\cx{\hskip 16 true pt\flip{$[{4}:{2}{1^{2}}]$}\hss}\cg%
\cx{\hskip 16 true pt\flip{$[{4}:{1^{4}}]$}\hss}\cg%
\cx{\hskip 16 true pt\flip{$[{3}{1}:{3}{1}]^{+}$}\hss}\cg%
\cx{\hskip 16 true pt\flip{$[{3}{1}:{3}{1}]^{-}$}\hss}\cg%
\cx{\hskip 16 true pt\flip{$[{3}{1}:{2^{2}}]$}\hss}\cg%
\cx{\hskip 16 true pt\flip{$[{3}{1}:{2}{1^{2}}]$}\hss}\cg%
\cx{\hskip 16 true pt\flip{$[{3}{1}:{1^{4}}]$}\hss}\cg%
\cx{\hskip 16 true pt\flip{$[{2^{2}}:{2^{2}}]^{+}$}\hss}\cg%
\cx{\hskip 16 true pt\flip{$[{2^{2}}:{2^{2}}]^{-}$}\hss}\cg%
\cx{\hskip 16 true pt\flip{$[{2^{2}}:{2}{1^{2}}]$}\hss}\cg%
\cx{\hskip 16 true pt\flip{$[{2^{2}}:{1^{4}}]$}\hss}\cg%
\cx{\hskip 16 true pt\flip{$[{2}{1^{2}}:{2}{1^{2}}]^{+}$}\hss}\cg%
\cx{\hskip 16 true pt\flip{$[{2}{1^{2}}:{2}{1^{2}}]^{-}$}\hss}\cg%
\cx{\hskip 16 true pt\flip{$[{2}{1^{2}}:{1^{4}}]$}\hss}\cg%
\cx{\hskip 16 true pt\flip{$[{1^{4}}:{1^{4}}]^{+}$}\hss}\cg%
\cx{\hskip 16 true pt\flip{$[{1^{4}}:{1^{4}}]^{-}$}\hss}\cg%
\eol}}\rg%
%
%
\rx{\vss\hfull{%
\rlx{\hss{$2688_{y}$}}\cg%
\e{0}%
\e{0}%
\e{0}%
\e{0}%
\e{0}%
\e{0}%
\e{0}%
\e{0}%
\e{1}%
\e{0}%
\e{0}%
\e{1}%
\e{0}%
\e{0}%
\e{1}%
\e{0}%
\e{1}%
\e{0}%
\e{0}%
\e{1}%
\e{1}%
\e{0}%
\e{0}%
\e{0}%
\e{0}%
\eol}\vss}\rg%
%
%
\rx{\vss\hfull{%
\rlx{\hss{$2100_{y}$}}\cg%
\e{0}%
\e{0}%
\e{0}%
\e{0}%
\e{0}%
\e{0}%
\e{0}%
\e{0}%
\e{0}%
\e{0}%
\e{1}%
\e{1}%
\e{0}%
\e{0}%
\e{1}%
\e{0}%
\e{1}%
\e{0}%
\e{0}%
\e{0}%
\e{1}%
\e{0}%
\e{0}%
\e{0}%
\e{0}%
\eol}\vss}\rg%
%
%
\rx{\vss\hfull{%
\rlx{\hss{$1400_{y}$}}\cg%
\e{0}%
\e{0}%
\e{0}%
\e{0}%
\e{0}%
\e{0}%
\e{0}%
\e{0}%
\e{0}%
\e{0}%
\e{1}%
\e{0}%
\e{0}%
\e{0}%
\e{1}%
\e{0}%
\e{0}%
\e{1}%
\e{0}%
\e{0}%
\e{0}%
\e{0}%
\e{0}%
\e{0}%
\e{0}%
\eol}\vss}\rg%
%
%
\rx{\vss\hfull{%
\rlx{\hss{$4536_{y}$}}\cg%
\e{0}%
\e{0}%
\e{0}%
\e{0}%
\e{0}%
\e{0}%
\e{0}%
\e{0}%
\e{0}%
\e{1}%
\e{0}%
\e{0}%
\e{1}%
\e{1}%
\e{1}%
\e{1}%
\e{0}%
\e{0}%
\e{1}%
\e{0}%
\e{0}%
\e{1}%
\e{0}%
\e{0}%
\e{0}%
\eol}\vss}\rg%
%
%
\rx{\vss\hfull{%
\rlx{\hss{$5670_{y}$}}\cg%
\e{0}%
\e{0}%
\e{0}%
\e{0}%
\e{0}%
\e{0}%
\e{0}%
\e{0}%
\e{0}%
\e{1}%
\e{0}%
\e{0}%
\e{1}%
\e{1}%
\e{2}%
\e{1}%
\e{0}%
\e{0}%
\e{1}%
\e{0}%
\e{0}%
\e{1}%
\e{0}%
\e{0}%
\e{0}%
\eol}\vss}\rg%
%
%
\rx{\vss\hfull{%
\rlx{\hss{$4480_{y}$}}\cg%
\e{0}%
\e{0}%
\e{0}%
\e{0}%
\e{0}%
\e{0}%
\e{0}%
\e{0}%
\e{0}%
\e{0}%
\e{0}%
\e{0}%
\e{1}%
\e{1}%
\e{1}%
\e{0}%
\e{0}%
\e{1}%
\e{1}%
\e{0}%
\e{0}%
\e{1}%
\e{0}%
\e{0}%
\e{0}%
\eol}\vss}\rg%
%
%
\rx{\vss\hfull{%
\rlx{\hss{$8_{z}$}}\cg%
\e{0}%
\e{0}%
\e{0}%
\e{0}%
\e{0}%
\e{0}%
\e{0}%
\e{0}%
\e{0}%
\e{0}%
\e{0}%
\e{0}%
\e{0}%
\e{0}%
\e{0}%
\e{0}%
\e{0}%
\e{0}%
\e{0}%
\e{0}%
\e{0}%
\e{0}%
\e{0}%
\e{0}%
\e{0}%
\eol}\vss}\rg%
%
%
\rx{\vss\hfull{%
\rlx{\hss{$56_{z}$}}\cg%
\e{0}%
\e{0}%
\e{0}%
\e{0}%
\e{0}%
\e{0}%
\e{0}%
\e{0}%
\e{0}%
\e{0}%
\e{0}%
\e{0}%
\e{0}%
\e{0}%
\e{0}%
\e{0}%
\e{0}%
\e{0}%
\e{0}%
\e{0}%
\e{0}%
\e{0}%
\e{0}%
\e{0}%
\e{0}%
\eol}\vss}\rg%
%
%
\rx{\vss\hfull{%
\rlx{\hss{$160_{z}$}}\cg%
\e{0}%
\e{0}%
\e{0}%
\e{0}%
\e{0}%
\e{0}%
\e{0}%
\e{0}%
\e{0}%
\e{0}%
\e{0}%
\e{0}%
\e{0}%
\e{0}%
\e{0}%
\e{0}%
\e{0}%
\e{0}%
\e{0}%
\e{0}%
\e{0}%
\e{0}%
\e{0}%
\e{0}%
\e{0}%
\eol}\vss}\rg%
%
%
\rx{\vss\hfull{%
\rlx{\hss{$112_{z}$}}\cg%
\e{0}%
\e{0}%
\e{0}%
\e{0}%
\e{0}%
\e{0}%
\e{0}%
\e{0}%
\e{0}%
\e{0}%
\e{0}%
\e{0}%
\e{0}%
\e{0}%
\e{0}%
\e{0}%
\e{0}%
\e{0}%
\e{0}%
\e{0}%
\e{0}%
\e{0}%
\e{0}%
\e{0}%
\e{0}%
\eol}\vss}\rg%
%
%
\rx{\vss\hfull{%
\rlx{\hss{$840_{z}$}}\cg%
\e{0}%
\e{0}%
\e{0}%
\e{0}%
\e{0}%
\e{0}%
\e{0}%
\e{0}%
\e{0}%
\e{0}%
\e{0}%
\e{0}%
\e{0}%
\e{0}%
\e{0}%
\e{0}%
\e{0}%
\e{0}%
\e{0}%
\e{0}%
\e{0}%
\e{0}%
\e{0}%
\e{0}%
\e{0}%
\eol}\vss}\rg%
%
%
\rx{\vss\hfull{%
\rlx{\hss{$1296_{z}$}}\cg%
\e{0}%
\e{0}%
\e{0}%
\e{0}%
\e{0}%
\e{0}%
\e{0}%
\e{0}%
\e{0}%
\e{0}%
\e{0}%
\e{0}%
\e{0}%
\e{0}%
\e{0}%
\e{0}%
\e{0}%
\e{0}%
\e{0}%
\e{0}%
\e{0}%
\e{0}%
\e{0}%
\e{0}%
\e{0}%
\eol}\vss}\rg%
%
%
\rx{\vss\hfull{%
\rlx{\hss{$1400_{z}$}}\cg%
\e{0}%
\e{0}%
\e{0}%
\e{0}%
\e{0}%
\e{0}%
\e{0}%
\e{0}%
\e{0}%
\e{0}%
\e{0}%
\e{0}%
\e{0}%
\e{0}%
\e{0}%
\e{0}%
\e{0}%
\e{0}%
\e{0}%
\e{0}%
\e{0}%
\e{0}%
\e{0}%
\e{0}%
\e{0}%
\eol}\vss}\rg%
%
%
\rx{\vss\hfull{%
\rlx{\hss{$1008_{z}$}}\cg%
\e{0}%
\e{0}%
\e{0}%
\e{0}%
\e{0}%
\e{0}%
\e{0}%
\e{0}%
\e{0}%
\e{0}%
\e{0}%
\e{0}%
\e{0}%
\e{0}%
\e{0}%
\e{0}%
\e{0}%
\e{0}%
\e{0}%
\e{0}%
\e{0}%
\e{0}%
\e{0}%
\e{0}%
\e{0}%
\eol}\vss}\rg%
%
%
\rx{\vss\hfull{%
\rlx{\hss{$560_{z}$}}\cg%
\e{0}%
\e{0}%
\e{0}%
\e{0}%
\e{0}%
\e{0}%
\e{0}%
\e{0}%
\e{0}%
\e{0}%
\e{0}%
\e{0}%
\e{0}%
\e{0}%
\e{0}%
\e{0}%
\e{0}%
\e{0}%
\e{0}%
\e{0}%
\e{0}%
\e{0}%
\e{0}%
\e{0}%
\e{0}%
\eol}\vss}\rg%
%
%
\rx{\vss\hfull{%
\rlx{\hss{$1400_{zz}$}}\cg%
\e{0}%
\e{0}%
\e{0}%
\e{0}%
\e{0}%
\e{0}%
\e{0}%
\e{0}%
\e{0}%
\e{0}%
\e{0}%
\e{0}%
\e{0}%
\e{0}%
\e{0}%
\e{0}%
\e{0}%
\e{0}%
\e{0}%
\e{0}%
\e{0}%
\e{0}%
\e{0}%
\e{0}%
\e{0}%
\eol}\vss}\rg%
%
%
\rx{\vss\hfull{%
\rlx{\hss{$4200_{z}$}}\cg%
\e{0}%
\e{0}%
\e{0}%
\e{0}%
\e{0}%
\e{0}%
\e{0}%
\e{0}%
\e{0}%
\e{0}%
\e{0}%
\e{0}%
\e{0}%
\e{0}%
\e{0}%
\e{0}%
\e{0}%
\e{0}%
\e{0}%
\e{0}%
\e{0}%
\e{0}%
\e{0}%
\e{0}%
\e{0}%
\eol}\vss}\rg%
%
%
\rx{\vss\hfull{%
\rlx{\hss{$400_{z}$}}\cg%
\e{0}%
\e{0}%
\e{0}%
\e{0}%
\e{0}%
\e{0}%
\e{0}%
\e{0}%
\e{0}%
\e{0}%
\e{0}%
\e{0}%
\e{0}%
\e{0}%
\e{0}%
\e{0}%
\e{0}%
\e{0}%
\e{0}%
\e{0}%
\e{0}%
\e{0}%
\e{0}%
\e{0}%
\e{0}%
\eol}\vss}\rg%
%
%
\rx{\vss\hfull{%
\rlx{\hss{$3240_{z}$}}\cg%
\e{0}%
\e{0}%
\e{0}%
\e{0}%
\e{0}%
\e{0}%
\e{0}%
\e{0}%
\e{0}%
\e{0}%
\e{0}%
\e{0}%
\e{0}%
\e{0}%
\e{0}%
\e{0}%
\e{0}%
\e{0}%
\e{0}%
\e{0}%
\e{0}%
\e{0}%
\e{0}%
\e{0}%
\e{0}%
\eol}\vss}\rg%
%
%
\rx{\vss\hfull{%
\rlx{\hss{$4536_{z}$}}\cg%
\e{0}%
\e{0}%
\e{0}%
\e{0}%
\e{0}%
\e{0}%
\e{0}%
\e{0}%
\e{0}%
\e{0}%
\e{0}%
\e{0}%
\e{0}%
\e{0}%
\e{0}%
\e{0}%
\e{0}%
\e{0}%
\e{0}%
\e{0}%
\e{0}%
\e{0}%
\e{0}%
\e{0}%
\e{0}%
\eol}\vss}\rg%
%
%
\rx{\vss\hfull{%
\rlx{\hss{$2400_{z}$}}\cg%
\e{0}%
\e{0}%
\e{0}%
\e{0}%
\e{0}%
\e{0}%
\e{0}%
\e{0}%
\e{0}%
\e{0}%
\e{0}%
\e{0}%
\e{0}%
\e{0}%
\e{0}%
\e{0}%
\e{0}%
\e{0}%
\e{0}%
\e{0}%
\e{0}%
\e{0}%
\e{0}%
\e{0}%
\e{0}%
\eol}\vss}\rg%
%
%
\rx{\vss\hfull{%
\rlx{\hss{$3360_{z}$}}\cg%
\e{0}%
\e{0}%
\e{0}%
\e{0}%
\e{0}%
\e{0}%
\e{0}%
\e{0}%
\e{0}%
\e{0}%
\e{0}%
\e{0}%
\e{0}%
\e{0}%
\e{0}%
\e{0}%
\e{0}%
\e{0}%
\e{0}%
\e{0}%
\e{0}%
\e{0}%
\e{0}%
\e{0}%
\e{0}%
\eol}\vss}\rg%
%
%
\rx{\vss\hfull{%
\rlx{\hss{$2800_{z}$}}\cg%
\e{0}%
\e{0}%
\e{0}%
\e{0}%
\e{0}%
\e{0}%
\e{0}%
\e{0}%
\e{0}%
\e{0}%
\e{0}%
\e{0}%
\e{0}%
\e{0}%
\e{0}%
\e{0}%
\e{0}%
\e{0}%
\e{0}%
\e{0}%
\e{0}%
\e{0}%
\e{0}%
\e{0}%
\e{0}%
\eol}\vss}\rg%
%
%
\rx{\vss\hfull{%
\rlx{\hss{$4096_{z}$}}\cg%
\e{0}%
\e{0}%
\e{0}%
\e{0}%
\e{0}%
\e{0}%
\e{0}%
\e{0}%
\e{0}%
\e{0}%
\e{0}%
\e{0}%
\e{0}%
\e{0}%
\e{0}%
\e{0}%
\e{0}%
\e{0}%
\e{0}%
\e{0}%
\e{0}%
\e{0}%
\e{0}%
\e{0}%
\e{0}%
\eol}\vss}\rg%
%
%
\rx{\vss\hfull{%
\rlx{\hss{$5600_{z}$}}\cg%
\e{0}%
\e{0}%
\e{0}%
\e{0}%
\e{0}%
\e{0}%
\e{0}%
\e{0}%
\e{0}%
\e{0}%
\e{0}%
\e{0}%
\e{0}%
\e{0}%
\e{0}%
\e{0}%
\e{0}%
\e{0}%
\e{0}%
\e{0}%
\e{0}%
\e{0}%
\e{0}%
\e{0}%
\e{0}%
\eol}\vss}\rg%
%
%
\rx{\vss\hfull{%
\rlx{\hss{$448_{z}$}}\cg%
\e{0}%
\e{0}%
\e{0}%
\e{0}%
\e{0}%
\e{0}%
\e{0}%
\e{0}%
\e{0}%
\e{0}%
\e{0}%
\e{0}%
\e{0}%
\e{0}%
\e{0}%
\e{0}%
\e{0}%
\e{0}%
\e{0}%
\e{0}%
\e{0}%
\e{0}%
\e{0}%
\e{0}%
\e{0}%
\eol}\vss}\rg%
%
%
\rx{\vss\hfull{%
\rlx{\hss{$448_{w}$}}\cg%
\e{0}%
\e{0}%
\e{0}%
\e{0}%
\e{0}%
\e{0}%
\e{0}%
\e{0}%
\e{0}%
\e{0}%
\e{0}%
\e{0}%
\e{0}%
\e{0}%
\e{0}%
\e{0}%
\e{0}%
\e{0}%
\e{0}%
\e{0}%
\e{0}%
\e{0}%
\e{0}%
\e{0}%
\e{0}%
\eol}\vss}\rg%
%
%
\rx{\vss\hfull{%
\rlx{\hss{$1344_{w}$}}\cg%
\e{0}%
\e{0}%
\e{0}%
\e{0}%
\e{0}%
\e{0}%
\e{0}%
\e{0}%
\e{0}%
\e{0}%
\e{0}%
\e{0}%
\e{0}%
\e{0}%
\e{0}%
\e{0}%
\e{0}%
\e{0}%
\e{0}%
\e{0}%
\e{0}%
\e{0}%
\e{0}%
\e{0}%
\e{0}%
\eol}\vss}\rg%
%
%
\rx{\vss\hfull{%
\rlx{\hss{$5600_{w}$}}\cg%
\e{1}%
\e{0}%
\e{0}%
\e{0}%
\e{0}%
\e{0}%
\e{0}%
\e{0}%
\e{0}%
\e{0}%
\e{0}%
\e{0}%
\e{0}%
\e{0}%
\e{0}%
\e{0}%
\e{0}%
\e{0}%
\e{0}%
\e{0}%
\e{0}%
\e{0}%
\e{0}%
\e{0}%
\e{0}%
\eol}\vss}\rg%
%
%
\rx{\vss\hfull{%
\rlx{\hss{$2016_{w}$}}\cg%
\e{0}%
\e{0}%
\e{0}%
\e{0}%
\e{0}%
\e{0}%
\e{0}%
\e{0}%
\e{0}%
\e{0}%
\e{0}%
\e{0}%
\e{0}%
\e{0}%
\e{0}%
\e{0}%
\e{0}%
\e{0}%
\e{0}%
\e{0}%
\e{0}%
\e{0}%
\e{0}%
\e{0}%
\e{0}%
\eol}\vss}\rg%
%
%
\rx{\vss\hfull{%
\rlx{\hss{$7168_{w}$}}\cg%
\e{1}%
\e{0}%
\e{0}%
\e{0}%
\e{0}%
\e{0}%
\e{0}%
\e{0}%
\e{0}%
\e{0}%
\e{0}%
\e{0}%
\e{0}%
\e{0}%
\e{0}%
\e{0}%
\e{0}%
\e{0}%
\e{0}%
\e{0}%
\e{0}%
\e{0}%
\e{0}%
\e{0}%
\e{0}%
\eol}\vss}\rg%
\tableclose%
%
%
%
%
%
%
\eop
\eject
\tableopen{Induce/restrict matrix for $W(A_{8})\,\subset\,W(E_{8})$}%
%
%
%
%
%
%
\rowpts=18 true pt%
\colpts=18 true pt%
\rowlabpts=40 true pt%
\collabpts=45 true pt%
\clx{\vss\hfull{%
\rlx{\hss{$ $}}\cg%
\cx{\hskip 16 true pt\flip{$[{9}]$}\hss}\cg%
\cx{\hskip 16 true pt\flip{$[{8}{1}]$}\hss}\cg%
\cx{\hskip 16 true pt\flip{$[{7}{2}]$}\hss}\cg%
\cx{\hskip 16 true pt\flip{$[{7}{1^{2}}]$}\hss}\cg%
\cx{\hskip 16 true pt\flip{$[{6}{3}]$}\hss}\cg%
\cx{\hskip 16 true pt\flip{$[{6}{2}{1}]$}\hss}\cg%
\cx{\hskip 16 true pt\flip{$[{6}{1^{3}}]$}\hss}\cg%
\cx{\hskip 16 true pt\flip{$[{5}{4}]$}\hss}\cg%
\cx{\hskip 16 true pt\flip{$[{5}{3}{1}]$}\hss}\cg%
\cx{\hskip 16 true pt\flip{$[{5}{2^{2}}]$}\hss}\cg%
\cx{\hskip 16 true pt\flip{$[{5}{2}{1^{2}}]$}\hss}\cg%
\cx{\hskip 16 true pt\flip{$[{5}{1^{4}}]$}\hss}\cg%
\cx{\hskip 16 true pt\flip{$[{4^{2}}{1}]$}\hss}\cg%
\cx{\hskip 16 true pt\flip{$[{4}{3}{2}]$}\hss}\cg%
\cx{\hskip 16 true pt\flip{$[{4}{3}{1^{2}}]$}\hss}\cg%
\eol}}\rg%
%
%
\rx{\vss\hfull{%
\rlx{\hss{$1_{x}$}}\cg%
\e{1}%
\e{0}%
\e{0}%
\e{0}%
\e{0}%
\e{0}%
\e{0}%
\e{0}%
\e{0}%
\e{0}%
\e{0}%
\e{0}%
\e{0}%
\e{0}%
\e{0}%
\eol}\vss}\rg%
%
%
\rx{\vss\hfull{%
\rlx{\hss{$28_{x}$}}\cg%
\e{0}%
\e{0}%
\e{0}%
\e{1}%
\e{0}%
\e{0}%
\e{0}%
\e{0}%
\e{0}%
\e{0}%
\e{0}%
\e{0}%
\e{0}%
\e{0}%
\e{0}%
\eol}\vss}\rg%
%
%
\rx{\vss\hfull{%
\rlx{\hss{$35_{x}$}}\cg%
\e{0}%
\e{1}%
\e{1}%
\e{0}%
\e{0}%
\e{0}%
\e{0}%
\e{0}%
\e{0}%
\e{0}%
\e{0}%
\e{0}%
\e{0}%
\e{0}%
\e{0}%
\eol}\vss}\rg%
%
%
\rx{\vss\hfull{%
\rlx{\hss{$84_{x}$}}\cg%
\e{1}%
\e{1}%
\e{1}%
\e{0}%
\e{1}%
\e{0}%
\e{0}%
\e{0}%
\e{0}%
\e{0}%
\e{0}%
\e{0}%
\e{0}%
\e{0}%
\e{0}%
\eol}\vss}\rg%
%
%
\rx{\vss\hfull{%
\rlx{\hss{$50_{x}$}}\cg%
\e{0}%
\e{1}%
\e{0}%
\e{0}%
\e{0}%
\e{0}%
\e{0}%
\e{1}%
\e{0}%
\e{0}%
\e{0}%
\e{0}%
\e{0}%
\e{0}%
\e{0}%
\eol}\vss}\rg%
%
%
\rx{\vss\hfull{%
\rlx{\hss{$350_{x}$}}\cg%
\e{0}%
\e{0}%
\e{0}%
\e{0}%
\e{0}%
\e{1}%
\e{1}%
\e{0}%
\e{0}%
\e{0}%
\e{1}%
\e{0}%
\e{0}%
\e{0}%
\e{0}%
\eol}\vss}\rg%
%
%
\rx{\vss\hfull{%
\rlx{\hss{$300_{x}$}}\cg%
\e{0}%
\e{0}%
\e{1}%
\e{0}%
\e{1}%
\e{1}%
\e{0}%
\e{0}%
\e{0}%
\e{1}%
\e{0}%
\e{0}%
\e{0}%
\e{0}%
\e{0}%
\eol}\vss}\rg%
%
%
\rx{\vss\hfull{%
\rlx{\hss{$567_{x}$}}\cg%
\e{0}%
\e{1}%
\e{1}%
\e{2}%
\e{1}%
\e{2}%
\e{1}%
\e{0}%
\e{1}%
\e{0}%
\e{0}%
\e{0}%
\e{0}%
\e{0}%
\e{0}%
\eol}\vss}\rg%
%
%
\rx{\vss\hfull{%
\rlx{\hss{$210_{x}$}}\cg%
\e{0}%
\e{1}%
\e{1}%
\e{1}%
\e{0}%
\e{1}%
\e{0}%
\e{1}%
\e{0}%
\e{0}%
\e{0}%
\e{0}%
\e{0}%
\e{0}%
\e{0}%
\eol}\vss}\rg%
%
%
\rx{\vss\hfull{%
\rlx{\hss{$840_{x}$}}\cg%
\e{0}%
\e{0}%
\e{0}%
\e{0}%
\e{1}%
\e{0}%
\e{0}%
\e{1}%
\e{1}%
\e{1}%
\e{0}%
\e{0}%
\e{0}%
\e{1}%
\e{0}%
\eol}\vss}\rg%
%
%
\rx{\vss\hfull{%
\rlx{\hss{$700_{x}$}}\cg%
\e{1}%
\e{1}%
\e{2}%
\e{1}%
\e{2}%
\e{1}%
\e{0}%
\e{1}%
\e{1}%
\e{1}%
\e{0}%
\e{0}%
\e{1}%
\e{0}%
\e{0}%
\eol}\vss}\rg%
%
%
\rx{\vss\hfull{%
\rlx{\hss{$175_{x}$}}\cg%
\e{1}%
\e{0}%
\e{0}%
\e{0}%
\e{1}%
\e{0}%
\e{0}%
\e{0}%
\e{0}%
\e{0}%
\e{0}%
\e{0}%
\e{1}%
\e{0}%
\e{0}%
\eol}\vss}\rg%
%
%
\rx{\vss\hfull{%
\rlx{\hss{$1400_{x}$}}\cg%
\e{0}%
\e{1}%
\e{1}%
\e{1}%
\e{1}%
\e{2}%
\e{1}%
\e{1}%
\e{2}%
\e{0}%
\e{1}%
\e{0}%
\e{1}%
\e{1}%
\e{1}%
\eol}\vss}\rg%
%
%
\rx{\vss\hfull{%
\rlx{\hss{$1050_{x}$}}\cg%
\e{0}%
\e{1}%
\e{1}%
\e{1}%
\e{1}%
\e{1}%
\e{0}%
\e{1}%
\e{2}%
\e{0}%
\e{0}%
\e{0}%
\e{1}%
\e{1}%
\e{1}%
\eol}\vss}\rg%
%
%
\rx{\vss\hfull{%
\rlx{\hss{$1575_{x}$}}\cg%
\e{0}%
\e{0}%
\e{1}%
\e{2}%
\e{1}%
\e{2}%
\e{2}%
\e{0}%
\e{2}%
\e{1}%
\e{2}%
\e{0}%
\e{1}%
\e{0}%
\e{1}%
\eol}\vss}\rg%
%
%
\rx{\vss\hfull{%
\rlx{\hss{$1344_{x}$}}\cg%
\e{0}%
\e{1}%
\e{2}%
\e{1}%
\e{2}%
\e{3}%
\e{0}%
\e{1}%
\e{2}%
\e{1}%
\e{1}%
\e{0}%
\e{0}%
\e{1}%
\e{0}%
\eol}\vss}\rg%
%
%
\rx{\vss\hfull{%
\rlx{\hss{$2100_{x}$}}\cg%
\e{0}%
\e{0}%
\e{0}%
\e{0}%
\e{0}%
\e{2}%
\e{2}%
\e{0}%
\e{1}%
\e{1}%
\e{3}%
\e{2}%
\e{0}%
\e{1}%
\e{1}%
\eol}\vss}\rg%
%
%
\rx{\vss\hfull{%
\rlx{\hss{$2268_{x}$}}\cg%
\e{0}%
\e{0}%
\e{1}%
\e{1}%
\e{1}%
\e{3}%
\e{2}%
\e{1}%
\e{2}%
\e{2}%
\e{2}%
\e{1}%
\e{1}%
\e{1}%
\e{1}%
\eol}\vss}\rg%
%
%
\rx{\vss\hfull{%
\rlx{\hss{$525_{x}$}}\cg%
\e{0}%
\e{0}%
\e{0}%
\e{1}%
\e{1}%
\e{0}%
\e{1}%
\e{0}%
\e{1}%
\e{0}%
\e{1}%
\e{0}%
\e{0}%
\e{0}%
\e{0}%
\eol}\vss}\rg%
%
%
\rx{\vss\hfull{%
\rlx{\hss{$700_{xx}$}}\cg%
\e{0}%
\e{0}%
\e{0}%
\e{1}%
\e{0}%
\e{0}%
\e{0}%
\e{1}%
\e{1}%
\e{0}%
\e{0}%
\e{0}%
\e{1}%
\e{0}%
\e{1}%
\eol}\vss}\rg%
%
%
\rx{\vss\hfull{%
\rlx{\hss{$972_{x}$}}\cg%
\e{0}%
\e{0}%
\e{1}%
\e{0}%
\e{1}%
\e{1}%
\e{0}%
\e{1}%
\e{1}%
\e{1}%
\e{0}%
\e{0}%
\e{1}%
\e{1}%
\e{0}%
\eol}\vss}\rg%
%
%
\rx{\vss\hfull{%
\rlx{\hss{$4096_{x}$}}\cg%
\e{0}%
\e{0}%
\e{1}%
\e{1}%
\e{1}%
\e{4}%
\e{1}%
\e{1}%
\e{4}%
\e{3}%
\e{4}%
\e{1}%
\e{1}%
\e{2}%
\e{2}%
\eol}\vss}\rg%
%
%
\rx{\vss\hfull{%
\rlx{\hss{$4200_{x}$}}\cg%
\e{0}%
\e{0}%
\e{1}%
\e{1}%
\e{2}%
\e{2}%
\e{1}%
\e{1}%
\e{4}%
\e{3}%
\e{2}%
\e{0}%
\e{2}%
\e{3}%
\e{3}%
\eol}\vss}\rg%
%
%
\rx{\vss\hfull{%
\rlx{\hss{$2240_{x}$}}\cg%
\e{0}%
\e{1}%
\e{1}%
\e{0}%
\e{2}%
\e{2}%
\e{0}%
\e{2}%
\e{3}%
\e{1}%
\e{1}%
\e{0}%
\e{1}%
\e{2}%
\e{1}%
\eol}\vss}\rg%
%
%
\rx{\vss\hfull{%
\rlx{\hss{$2835_{x}$}}\cg%
\e{0}%
\e{0}%
\e{1}%
\e{0}%
\e{1}%
\e{1}%
\e{0}%
\e{1}%
\e{2}%
\e{1}%
\e{1}%
\e{0}%
\e{2}%
\e{3}%
\e{2}%
\eol}\vss}\rg%
%
%
\rx{\vss\hfull{%
\rlx{\hss{$6075_{x}$}}\cg%
\e{0}%
\e{0}%
\e{0}%
\e{1}%
\e{1}%
\e{3}%
\e{2}%
\e{1}%
\e{4}%
\e{2}%
\e{5}%
\e{1}%
\e{2}%
\e{3}%
\e{5}%
\eol}\vss}\rg%
%
%
\rx{\vss\hfull{%
\rlx{\hss{$3200_{x}$}}\cg%
\e{0}%
\e{0}%
\e{0}%
\e{0}%
\e{1}%
\e{1}%
\e{0}%
\e{0}%
\e{2}%
\e{3}%
\e{2}%
\e{0}%
\e{1}%
\e{2}%
\e{1}%
\eol}\vss}\rg%
%
%
\rx{\vss\hfull{%
\rlx{\hss{$70_{y}$}}\cg%
\e{0}%
\e{0}%
\e{0}%
\e{0}%
\e{0}%
\e{0}%
\e{0}%
\e{0}%
\e{0}%
\e{0}%
\e{0}%
\e{1}%
\e{0}%
\e{0}%
\e{0}%
\eol}\vss}\rg%
%
%
\rx{\vss\hfull{%
\rlx{\hss{$1134_{y}$}}\cg%
\e{0}%
\e{0}%
\e{0}%
\e{0}%
\e{0}%
\e{0}%
\e{0}%
\e{0}%
\e{1}%
\e{0}%
\e{1}%
\e{0}%
\e{0}%
\e{0}%
\e{1}%
\eol}\vss}\rg%
%
%
\rx{\vss\hfull{%
\rlx{\hss{$1680_{y}$}}\cg%
\e{0}%
\e{0}%
\e{0}%
\e{0}%
\e{0}%
\e{0}%
\e{1}%
\e{0}%
\e{0}%
\e{1}%
\e{2}%
\e{2}%
\e{0}%
\e{0}%
\e{1}%
\eol}\vss}\rg%
%
%
\rx{\vss\hfull{%
\rlx{\hss{$168_{y}$}}\cg%
\e{0}%
\e{0}%
\e{0}%
\e{0}%
\e{0}%
\e{0}%
\e{0}%
\e{0}%
\e{0}%
\e{0}%
\e{0}%
\e{0}%
\e{1}%
\e{0}%
\e{0}%
\eol}\vss}\rg%
%
%
\rx{\vss\hfull{%
\rlx{\hss{$420_{y}$}}\cg%
\e{0}%
\e{0}%
\e{0}%
\e{0}%
\e{0}%
\e{0}%
\e{0}%
\e{1}%
\e{0}%
\e{0}%
\e{0}%
\e{0}%
\e{0}%
\e{1}%
\e{0}%
\eol}\vss}\rg%
%
%
\rx{\vss\hfull{%
\rlx{\hss{$3150_{y}$}}\cg%
\e{0}%
\e{0}%
\e{0}%
\e{0}%
\e{1}%
\e{0}%
\e{0}%
\e{0}%
\e{2}%
\e{1}%
\e{1}%
\e{0}%
\e{1}%
\e{2}%
\e{2}%
\eol}\vss}\rg%
%
%
\rx{\vss\hfull{%
\rlx{\hss{$4200_{y}$}}\cg%
\e{0}%
\e{0}%
\e{0}%
\e{0}%
\e{0}%
\e{1}%
\e{0}%
\e{1}%
\e{2}%
\e{1}%
\e{1}%
\e{0}%
\e{2}%
\e{3}%
\e{3}%
\eol}\vss}\rg%
\eop
\eject
\tablecont%
%
%
%
%
%
%
\rowpts=18 true pt%
\colpts=18 true pt%
\rowlabpts=40 true pt%
\collabpts=45 true pt%
\clx{\vss\hfull{%
\rlx{\hss{$ $}}\cg%
\cx{\hskip 16 true pt\flip{$[{9}]$}\hss}\cg%
\cx{\hskip 16 true pt\flip{$[{8}{1}]$}\hss}\cg%
\cx{\hskip 16 true pt\flip{$[{7}{2}]$}\hss}\cg%
\cx{\hskip 16 true pt\flip{$[{7}{1^{2}}]$}\hss}\cg%
\cx{\hskip 16 true pt\flip{$[{6}{3}]$}\hss}\cg%
\cx{\hskip 16 true pt\flip{$[{6}{2}{1}]$}\hss}\cg%
\cx{\hskip 16 true pt\flip{$[{6}{1^{3}}]$}\hss}\cg%
\cx{\hskip 16 true pt\flip{$[{5}{4}]$}\hss}\cg%
\cx{\hskip 16 true pt\flip{$[{5}{3}{1}]$}\hss}\cg%
\cx{\hskip 16 true pt\flip{$[{5}{2^{2}}]$}\hss}\cg%
\cx{\hskip 16 true pt\flip{$[{5}{2}{1^{2}}]$}\hss}\cg%
\cx{\hskip 16 true pt\flip{$[{5}{1^{4}}]$}\hss}\cg%
\cx{\hskip 16 true pt\flip{$[{4^{2}}{1}]$}\hss}\cg%
\cx{\hskip 16 true pt\flip{$[{4}{3}{2}]$}\hss}\cg%
\cx{\hskip 16 true pt\flip{$[{4}{3}{1^{2}}]$}\hss}\cg%
\eol}}\rg%
%
%
\rx{\vss\hfull{%
\rlx{\hss{$2688_{y}$}}\cg%
\e{0}%
\e{0}%
\e{0}%
\e{0}%
\e{0}%
\e{1}%
\e{0}%
\e{0}%
\e{1}%
\e{1}%
\e{1}%
\e{0}%
\e{0}%
\e{2}%
\e{2}%
\eol}\vss}\rg%
%
%
\rx{\vss\hfull{%
\rlx{\hss{$2100_{y}$}}\cg%
\e{0}%
\e{0}%
\e{0}%
\e{0}%
\e{0}%
\e{0}%
\e{1}%
\e{0}%
\e{1}%
\e{0}%
\e{2}%
\e{2}%
\e{0}%
\e{1}%
\e{1}%
\eol}\vss}\rg%
%
%
\rx{\vss\hfull{%
\rlx{\hss{$1400_{y}$}}\cg%
\e{0}%
\e{0}%
\e{0}%
\e{0}%
\e{0}%
\e{0}%
\e{1}%
\e{0}%
\e{0}%
\e{1}%
\e{1}%
\e{1}%
\e{1}%
\e{0}%
\e{1}%
\eol}\vss}\rg%
%
%
\rx{\vss\hfull{%
\rlx{\hss{$4536_{y}$}}\cg%
\e{0}%
\e{0}%
\e{0}%
\e{0}%
\e{0}%
\e{1}%
\e{1}%
\e{0}%
\e{1}%
\e{2}%
\e{3}%
\e{2}%
\e{1}%
\e{2}%
\e{3}%
\eol}\vss}\rg%
%
%
\rx{\vss\hfull{%
\rlx{\hss{$5670_{y}$}}\cg%
\e{0}%
\e{0}%
\e{0}%
\e{0}%
\e{0}%
\e{1}%
\e{1}%
\e{0}%
\e{2}%
\e{2}%
\e{4}%
\e{2}%
\e{1}%
\e{2}%
\e{4}%
\eol}\vss}\rg%
%
%
\rx{\vss\hfull{%
\rlx{\hss{$4480_{y}$}}\cg%
\e{0}%
\e{0}%
\e{0}%
\e{0}%
\e{0}%
\e{1}%
\e{0}%
\e{1}%
\e{2}%
\e{1}%
\e{2}%
\e{1}%
\e{1}%
\e{3}%
\e{3}%
\eol}\vss}\rg%
%
%
\rx{\vss\hfull{%
\rlx{\hss{$8_{z}$}}\cg%
\e{0}%
\e{1}%
\e{0}%
\e{0}%
\e{0}%
\e{0}%
\e{0}%
\e{0}%
\e{0}%
\e{0}%
\e{0}%
\e{0}%
\e{0}%
\e{0}%
\e{0}%
\eol}\vss}\rg%
%
%
\rx{\vss\hfull{%
\rlx{\hss{$56_{z}$}}\cg%
\e{0}%
\e{0}%
\e{0}%
\e{0}%
\e{0}%
\e{0}%
\e{1}%
\e{0}%
\e{0}%
\e{0}%
\e{0}%
\e{0}%
\e{0}%
\e{0}%
\e{0}%
\eol}\vss}\rg%
%
%
\rx{\vss\hfull{%
\rlx{\hss{$160_{z}$}}\cg%
\e{0}%
\e{0}%
\e{1}%
\e{1}%
\e{0}%
\e{1}%
\e{0}%
\e{0}%
\e{0}%
\e{0}%
\e{0}%
\e{0}%
\e{0}%
\e{0}%
\e{0}%
\eol}\vss}\rg%
%
%
\rx{\vss\hfull{%
\rlx{\hss{$112_{z}$}}\cg%
\e{1}%
\e{1}%
\e{1}%
\e{1}%
\e{1}%
\e{0}%
\e{0}%
\e{0}%
\e{0}%
\e{0}%
\e{0}%
\e{0}%
\e{0}%
\e{0}%
\e{0}%
\eol}\vss}\rg%
%
%
\rx{\vss\hfull{%
\rlx{\hss{$840_{z}$}}\cg%
\e{0}%
\e{0}%
\e{0}%
\e{0}%
\e{1}%
\e{1}%
\e{0}%
\e{0}%
\e{1}%
\e{1}%
\e{1}%
\e{0}%
\e{0}%
\e{0}%
\e{0}%
\eol}\vss}\rg%
%
%
\rx{\vss\hfull{%
\rlx{\hss{$1296_{z}$}}\cg%
\e{0}%
\e{0}%
\e{0}%
\e{1}%
\e{0}%
\e{2}%
\e{2}%
\e{0}%
\e{1}%
\e{1}%
\e{2}%
\e{1}%
\e{0}%
\e{0}%
\e{1}%
\eol}\vss}\rg%
%
%
\rx{\vss\hfull{%
\rlx{\hss{$1400_{z}$}}\cg%
\e{0}%
\e{1}%
\e{2}%
\e{1}%
\e{2}%
\e{3}%
\e{1}%
\e{1}%
\e{2}%
\e{1}%
\e{1}%
\e{0}%
\e{0}%
\e{1}%
\e{0}%
\eol}\vss}\rg%
%
%
\rx{\vss\hfull{%
\rlx{\hss{$1008_{z}$}}\cg%
\e{0}%
\e{0}%
\e{1}%
\e{2}%
\e{1}%
\e{2}%
\e{2}%
\e{0}%
\e{1}%
\e{1}%
\e{1}%
\e{0}%
\e{1}%
\e{0}%
\e{0}%
\eol}\vss}\rg%
%
%
\rx{\vss\hfull{%
\rlx{\hss{$560_{z}$}}\cg%
\e{0}%
\e{2}%
\e{2}%
\e{1}%
\e{1}%
\e{2}%
\e{0}%
\e{1}%
\e{1}%
\e{0}%
\e{0}%
\e{0}%
\e{0}%
\e{0}%
\e{0}%
\eol}\vss}\rg%
%
%
\rx{\vss\hfull{%
\rlx{\hss{$1400_{zz}$}}\cg%
\e{0}%
\e{1}%
\e{1}%
\e{0}%
\e{1}%
\e{1}%
\e{0}%
\e{2}%
\e{2}%
\e{0}%
\e{0}%
\e{0}%
\e{1}%
\e{2}%
\e{1}%
\eol}\vss}\rg%
%
%
\rx{\vss\hfull{%
\rlx{\hss{$4200_{z}$}}\cg%
\e{0}%
\e{0}%
\e{0}%
\e{1}%
\e{1}%
\e{1}%
\e{1}%
\e{1}%
\e{3}%
\e{1}%
\e{2}%
\e{0}%
\e{3}%
\e{2}%
\e{4}%
\eol}\vss}\rg%
%
%
\rx{\vss\hfull{%
\rlx{\hss{$400_{z}$}}\cg%
\e{1}%
\e{1}%
\e{1}%
\e{1}%
\e{1}%
\e{0}%
\e{0}%
\e{1}%
\e{1}%
\e{0}%
\e{0}%
\e{0}%
\e{1}%
\e{0}%
\e{0}%
\eol}\vss}\rg%
%
%
\rx{\vss\hfull{%
\rlx{\hss{$3240_{z}$}}\cg%
\e{0}%
\e{1}%
\e{2}%
\e{2}%
\e{3}%
\e{4}%
\e{1}%
\e{2}%
\e{4}%
\e{2}%
\e{2}%
\e{0}%
\e{2}%
\e{2}%
\e{2}%
\eol}\vss}\rg%
%
%
\rx{\vss\hfull{%
\rlx{\hss{$4536_{z}$}}\cg%
\e{0}%
\e{0}%
\e{1}%
\e{0}%
\e{2}%
\e{2}%
\e{0}%
\e{2}%
\e{4}%
\e{3}%
\e{2}%
\e{0}%
\e{2}%
\e{4}%
\e{2}%
\eol}\vss}\rg%
%
%
\rx{\vss\hfull{%
\rlx{\hss{$2400_{z}$}}\cg%
\e{0}%
\e{0}%
\e{0}%
\e{0}%
\e{0}%
\e{1}%
\e{2}%
\e{0}%
\e{1}%
\e{0}%
\e{3}%
\e{2}%
\e{0}%
\e{1}%
\e{2}%
\eol}\vss}\rg%
%
%
\rx{\vss\hfull{%
\rlx{\hss{$3360_{z}$}}\cg%
\e{0}%
\e{0}%
\e{1}%
\e{1}%
\e{1}%
\e{2}%
\e{1}%
\e{0}%
\e{3}%
\e{2}%
\e{2}%
\e{0}%
\e{2}%
\e{2}%
\e{3}%
\eol}\vss}\rg%
%
%
\rx{\vss\hfull{%
\rlx{\hss{$2800_{z}$}}\cg%
\e{0}%
\e{0}%
\e{0}%
\e{1}%
\e{0}%
\e{2}%
\e{2}%
\e{1}%
\e{2}%
\e{1}%
\e{3}%
\e{2}%
\e{1}%
\e{1}%
\e{2}%
\eol}\vss}\rg%
%
%
\rx{\vss\hfull{%
\rlx{\hss{$4096_{z}$}}\cg%
\e{0}%
\e{0}%
\e{1}%
\e{1}%
\e{1}%
\e{4}%
\e{1}%
\e{1}%
\e{4}%
\e{3}%
\e{4}%
\e{1}%
\e{1}%
\e{2}%
\e{2}%
\eol}\vss}\rg%
%
%
\rx{\vss\hfull{%
\rlx{\hss{$5600_{z}$}}\cg%
\e{0}%
\e{0}%
\e{0}%
\e{0}%
\e{1}%
\e{2}%
\e{2}%
\e{0}%
\e{3}%
\e{3}%
\e{5}%
\e{2}%
\e{1}%
\e{3}%
\e{3}%
\eol}\vss}\rg%
%
%
\rx{\vss\hfull{%
\rlx{\hss{$448_{z}$}}\cg%
\e{1}%
\e{0}%
\e{1}%
\e{0}%
\e{2}%
\e{0}%
\e{0}%
\e{0}%
\e{1}%
\e{1}%
\e{0}%
\e{0}%
\e{0}%
\e{0}%
\e{0}%
\eol}\vss}\rg%
%
%
\rx{\vss\hfull{%
\rlx{\hss{$448_{w}$}}\cg%
\e{0}%
\e{0}%
\e{0}%
\e{0}%
\e{0}%
\e{0}%
\e{0}%
\e{0}%
\e{0}%
\e{0}%
\e{1}%
\e{1}%
\e{0}%
\e{0}%
\e{0}%
\eol}\vss}\rg%
%
%
\rx{\vss\hfull{%
\rlx{\hss{$1344_{w}$}}\cg%
\e{0}%
\e{0}%
\e{0}%
\e{0}%
\e{0}%
\e{0}%
\e{0}%
\e{1}%
\e{1}%
\e{0}%
\e{0}%
\e{0}%
\e{1}%
\e{1}%
\e{1}%
\eol}\vss}\rg%
%
%
\rx{\vss\hfull{%
\rlx{\hss{$5600_{w}$}}\cg%
\e{0}%
\e{0}%
\e{0}%
\e{0}%
\e{0}%
\e{1}%
\e{1}%
\e{0}%
\e{2}%
\e{2}%
\e{4}%
\e{1}%
\e{1}%
\e{2}%
\e{4}%
\eol}\vss}\rg%
%
%
\rx{\vss\hfull{%
\rlx{\hss{$2016_{w}$}}\cg%
\e{0}%
\e{0}%
\e{0}%
\e{0}%
\e{1}%
\e{0}%
\e{0}%
\e{0}%
\e{1}%
\e{1}%
\e{0}%
\e{0}%
\e{1}%
\e{2}%
\e{1}%
\eol}\vss}\rg%
%
%
\rx{\vss\hfull{%
\rlx{\hss{$7168_{w}$}}\cg%
\e{0}%
\e{0}%
\e{0}%
\e{0}%
\e{0}%
\e{2}%
\e{0}%
\e{1}%
\e{3}%
\e{2}%
\e{3}%
\e{1}%
\e{1}%
\e{5}%
\e{5}%
\eol}\vss}\rg%
\eop
\eject
\tablecont%
%
%
%
%
%
%
\rowpts=18 true pt%
\colpts=18 true pt%
\rowlabpts=40 true pt%
\collabpts=45 true pt%
\clx{\vss\hfull{%
\rlx{\hss{$ $}}\cg%
\cx{\hskip 16 true pt\flip{$[{4}{2^{2}}{1}]$}\hss}\cg%
\cx{\hskip 16 true pt\flip{$[{4}{2}{1^{3}}]$}\hss}\cg%
\cx{\hskip 16 true pt\flip{$[{4}{1^{5}}]$}\hss}\cg%
\cx{\hskip 16 true pt\flip{$[{3^{3}}]$}\hss}\cg%
\cx{\hskip 16 true pt\flip{$[{3^{2}}{2}{1}]$}\hss}\cg%
\cx{\hskip 16 true pt\flip{$[{3^{2}}{1^{3}}]$}\hss}\cg%
\cx{\hskip 16 true pt\flip{$[{3}{2^{3}}]$}\hss}\cg%
\cx{\hskip 16 true pt\flip{$[{3}{2^{2}}{1^{2}}]$}\hss}\cg%
\cx{\hskip 16 true pt\flip{$[{3}{2}{1^{4}}]$}\hss}\cg%
\cx{\hskip 16 true pt\flip{$[{3}{1^{6}}]$}\hss}\cg%
\cx{\hskip 16 true pt\flip{$[{2^{4}}{1}]$}\hss}\cg%
\cx{\hskip 16 true pt\flip{$[{2^{3}}{1^{3}}]$}\hss}\cg%
\cx{\hskip 16 true pt\flip{$[{2^{2}}{1^{5}}]$}\hss}\cg%
\cx{\hskip 16 true pt\flip{$[{2}{1^{7}}]$}\hss}\cg%
\cx{\hskip 16 true pt\flip{$[{1^{9}}]$}\hss}\cg%
\eol}}\rg%
%
%
\rx{\vss\hfull{%
\rlx{\hss{$1_{x}$}}\cg%
\e{0}%
\e{0}%
\e{0}%
\e{0}%
\e{0}%
\e{0}%
\e{0}%
\e{0}%
\e{0}%
\e{0}%
\e{0}%
\e{0}%
\e{0}%
\e{0}%
\e{0}%
\eol}\vss}\rg%
%
%
\rx{\vss\hfull{%
\rlx{\hss{$28_{x}$}}\cg%
\e{0}%
\e{0}%
\e{0}%
\e{0}%
\e{0}%
\e{0}%
\e{0}%
\e{0}%
\e{0}%
\e{0}%
\e{0}%
\e{0}%
\e{0}%
\e{0}%
\e{0}%
\eol}\vss}\rg%
%
%
\rx{\vss\hfull{%
\rlx{\hss{$35_{x}$}}\cg%
\e{0}%
\e{0}%
\e{0}%
\e{0}%
\e{0}%
\e{0}%
\e{0}%
\e{0}%
\e{0}%
\e{0}%
\e{0}%
\e{0}%
\e{0}%
\e{0}%
\e{0}%
\eol}\vss}\rg%
%
%
\rx{\vss\hfull{%
\rlx{\hss{$84_{x}$}}\cg%
\e{0}%
\e{0}%
\e{0}%
\e{0}%
\e{0}%
\e{0}%
\e{0}%
\e{0}%
\e{0}%
\e{0}%
\e{0}%
\e{0}%
\e{0}%
\e{0}%
\e{0}%
\eol}\vss}\rg%
%
%
\rx{\vss\hfull{%
\rlx{\hss{$50_{x}$}}\cg%
\e{0}%
\e{0}%
\e{0}%
\e{0}%
\e{0}%
\e{0}%
\e{0}%
\e{0}%
\e{0}%
\e{0}%
\e{0}%
\e{0}%
\e{0}%
\e{0}%
\e{0}%
\eol}\vss}\rg%
%
%
\rx{\vss\hfull{%
\rlx{\hss{$350_{x}$}}\cg%
\e{0}%
\e{0}%
\e{0}%
\e{0}%
\e{0}%
\e{0}%
\e{0}%
\e{0}%
\e{0}%
\e{0}%
\e{0}%
\e{0}%
\e{0}%
\e{0}%
\e{0}%
\eol}\vss}\rg%
%
%
\rx{\vss\hfull{%
\rlx{\hss{$300_{x}$}}\cg%
\e{0}%
\e{0}%
\e{0}%
\e{0}%
\e{0}%
\e{0}%
\e{0}%
\e{0}%
\e{0}%
\e{0}%
\e{0}%
\e{0}%
\e{0}%
\e{0}%
\e{0}%
\eol}\vss}\rg%
%
%
\rx{\vss\hfull{%
\rlx{\hss{$567_{x}$}}\cg%
\e{0}%
\e{0}%
\e{0}%
\e{0}%
\e{0}%
\e{0}%
\e{0}%
\e{0}%
\e{0}%
\e{0}%
\e{0}%
\e{0}%
\e{0}%
\e{0}%
\e{0}%
\eol}\vss}\rg%
%
%
\rx{\vss\hfull{%
\rlx{\hss{$210_{x}$}}\cg%
\e{0}%
\e{0}%
\e{0}%
\e{0}%
\e{0}%
\e{0}%
\e{0}%
\e{0}%
\e{0}%
\e{0}%
\e{0}%
\e{0}%
\e{0}%
\e{0}%
\e{0}%
\eol}\vss}\rg%
%
%
\rx{\vss\hfull{%
\rlx{\hss{$840_{x}$}}\cg%
\e{1}%
\e{0}%
\e{0}%
\e{0}%
\e{0}%
\e{0}%
\e{1}%
\e{0}%
\e{0}%
\e{0}%
\e{0}%
\e{0}%
\e{0}%
\e{0}%
\e{0}%
\eol}\vss}\rg%
%
%
\rx{\vss\hfull{%
\rlx{\hss{$700_{x}$}}\cg%
\e{0}%
\e{0}%
\e{0}%
\e{0}%
\e{0}%
\e{0}%
\e{0}%
\e{0}%
\e{0}%
\e{0}%
\e{0}%
\e{0}%
\e{0}%
\e{0}%
\e{0}%
\eol}\vss}\rg%
%
%
\rx{\vss\hfull{%
\rlx{\hss{$175_{x}$}}\cg%
\e{0}%
\e{0}%
\e{0}%
\e{1}%
\e{0}%
\e{0}%
\e{0}%
\e{0}%
\e{0}%
\e{0}%
\e{0}%
\e{0}%
\e{0}%
\e{0}%
\e{0}%
\eol}\vss}\rg%
%
%
\rx{\vss\hfull{%
\rlx{\hss{$1400_{x}$}}\cg%
\e{0}%
\e{0}%
\e{0}%
\e{0}%
\e{0}%
\e{0}%
\e{0}%
\e{0}%
\e{0}%
\e{0}%
\e{0}%
\e{0}%
\e{0}%
\e{0}%
\e{0}%
\eol}\vss}\rg%
%
%
\rx{\vss\hfull{%
\rlx{\hss{$1050_{x}$}}\cg%
\e{0}%
\e{0}%
\e{0}%
\e{0}%
\e{0}%
\e{0}%
\e{0}%
\e{0}%
\e{0}%
\e{0}%
\e{0}%
\e{0}%
\e{0}%
\e{0}%
\e{0}%
\eol}\vss}\rg%
%
%
\rx{\vss\hfull{%
\rlx{\hss{$1575_{x}$}}\cg%
\e{0}%
\e{0}%
\e{0}%
\e{0}%
\e{0}%
\e{0}%
\e{0}%
\e{0}%
\e{0}%
\e{0}%
\e{0}%
\e{0}%
\e{0}%
\e{0}%
\e{0}%
\eol}\vss}\rg%
%
%
\rx{\vss\hfull{%
\rlx{\hss{$1344_{x}$}}\cg%
\e{0}%
\e{0}%
\e{0}%
\e{0}%
\e{0}%
\e{0}%
\e{0}%
\e{0}%
\e{0}%
\e{0}%
\e{0}%
\e{0}%
\e{0}%
\e{0}%
\e{0}%
\eol}\vss}\rg%
%
%
\rx{\vss\hfull{%
\rlx{\hss{$2100_{x}$}}\cg%
\e{1}%
\e{1}%
\e{0}%
\e{0}%
\e{0}%
\e{0}%
\e{0}%
\e{0}%
\e{0}%
\e{0}%
\e{0}%
\e{0}%
\e{0}%
\e{0}%
\e{0}%
\eol}\vss}\rg%
%
%
\rx{\vss\hfull{%
\rlx{\hss{$2268_{x}$}}\cg%
\e{1}%
\e{0}%
\e{0}%
\e{0}%
\e{0}%
\e{0}%
\e{0}%
\e{0}%
\e{0}%
\e{0}%
\e{0}%
\e{0}%
\e{0}%
\e{0}%
\e{0}%
\eol}\vss}\rg%
%
%
\rx{\vss\hfull{%
\rlx{\hss{$525_{x}$}}\cg%
\e{0}%
\e{0}%
\e{0}%
\e{1}%
\e{0}%
\e{0}%
\e{0}%
\e{0}%
\e{0}%
\e{0}%
\e{0}%
\e{0}%
\e{0}%
\e{0}%
\e{0}%
\eol}\vss}\rg%
%
%
\rx{\vss\hfull{%
\rlx{\hss{$700_{xx}$}}\cg%
\e{0}%
\e{0}%
\e{0}%
\e{0}%
\e{1}%
\e{0}%
\e{0}%
\e{0}%
\e{0}%
\e{0}%
\e{0}%
\e{0}%
\e{0}%
\e{0}%
\e{0}%
\eol}\vss}\rg%
%
%
\rx{\vss\hfull{%
\rlx{\hss{$972_{x}$}}\cg%
\e{1}%
\e{0}%
\e{0}%
\e{0}%
\e{0}%
\e{0}%
\e{0}%
\e{0}%
\e{0}%
\e{0}%
\e{0}%
\e{0}%
\e{0}%
\e{0}%
\e{0}%
\eol}\vss}\rg%
%
%
\rx{\vss\hfull{%
\rlx{\hss{$4096_{x}$}}\cg%
\e{2}%
\e{1}%
\e{0}%
\e{0}%
\e{1}%
\e{0}%
\e{0}%
\e{0}%
\e{0}%
\e{0}%
\e{0}%
\e{0}%
\e{0}%
\e{0}%
\e{0}%
\eol}\vss}\rg%
%
%
\rx{\vss\hfull{%
\rlx{\hss{$4200_{x}$}}\cg%
\e{2}%
\e{1}%
\e{0}%
\e{1}%
\e{1}%
\e{1}%
\e{1}%
\e{0}%
\e{0}%
\e{0}%
\e{0}%
\e{0}%
\e{0}%
\e{0}%
\e{0}%
\eol}\vss}\rg%
%
%
\rx{\vss\hfull{%
\rlx{\hss{$2240_{x}$}}\cg%
\e{1}%
\e{0}%
\e{0}%
\e{0}%
\e{1}%
\e{0}%
\e{0}%
\e{0}%
\e{0}%
\e{0}%
\e{0}%
\e{0}%
\e{0}%
\e{0}%
\e{0}%
\eol}\vss}\rg%
%
%
\rx{\vss\hfull{%
\rlx{\hss{$2835_{x}$}}\cg%
\e{1}%
\e{0}%
\e{0}%
\e{1}%
\e{2}%
\e{1}%
\e{0}%
\e{1}%
\e{0}%
\e{0}%
\e{0}%
\e{0}%
\e{0}%
\e{0}%
\e{0}%
\eol}\vss}\rg%
%
%
\rx{\vss\hfull{%
\rlx{\hss{$6075_{x}$}}\cg%
\e{3}%
\e{3}%
\e{0}%
\e{1}%
\e{2}%
\e{1}%
\e{0}%
\e{1}%
\e{0}%
\e{0}%
\e{0}%
\e{0}%
\e{0}%
\e{0}%
\e{0}%
\eol}\vss}\rg%
%
%
\rx{\vss\hfull{%
\rlx{\hss{$3200_{x}$}}\cg%
\e{3}%
\e{1}%
\e{1}%
\e{1}%
\e{1}%
\e{0}%
\e{1}%
\e{1}%
\e{0}%
\e{0}%
\e{0}%
\e{0}%
\e{0}%
\e{0}%
\e{0}%
\eol}\vss}\rg%
%
%
\rx{\vss\hfull{%
\rlx{\hss{$70_{y}$}}\cg%
\e{0}%
\e{0}%
\e{0}%
\e{0}%
\e{0}%
\e{0}%
\e{0}%
\e{0}%
\e{0}%
\e{0}%
\e{0}%
\e{0}%
\e{0}%
\e{0}%
\e{0}%
\eol}\vss}\rg%
%
%
\rx{\vss\hfull{%
\rlx{\hss{$1134_{y}$}}\cg%
\e{1}%
\e{1}%
\e{0}%
\e{0}%
\e{0}%
\e{0}%
\e{0}%
\e{1}%
\e{0}%
\e{0}%
\e{0}%
\e{0}%
\e{0}%
\e{0}%
\e{0}%
\eol}\vss}\rg%
%
%
\rx{\vss\hfull{%
\rlx{\hss{$1680_{y}$}}\cg%
\e{1}%
\e{2}%
\e{1}%
\e{0}%
\e{0}%
\e{1}%
\e{0}%
\e{0}%
\e{0}%
\e{0}%
\e{0}%
\e{0}%
\e{0}%
\e{0}%
\e{0}%
\eol}\vss}\rg%
%
%
\rx{\vss\hfull{%
\rlx{\hss{$168_{y}$}}\cg%
\e{0}%
\e{0}%
\e{0}%
\e{0}%
\e{0}%
\e{0}%
\e{1}%
\e{0}%
\e{0}%
\e{0}%
\e{0}%
\e{0}%
\e{0}%
\e{0}%
\e{0}%
\eol}\vss}\rg%
%
%
\rx{\vss\hfull{%
\rlx{\hss{$420_{y}$}}\cg%
\e{0}%
\e{0}%
\e{0}%
\e{0}%
\e{1}%
\e{0}%
\e{0}%
\e{0}%
\e{0}%
\e{0}%
\e{1}%
\e{0}%
\e{0}%
\e{0}%
\e{0}%
\eol}\vss}\rg%
%
%
\rx{\vss\hfull{%
\rlx{\hss{$3150_{y}$}}\cg%
\e{2}%
\e{1}%
\e{0}%
\e{2}%
\e{2}%
\e{1}%
\e{1}%
\e{2}%
\e{0}%
\e{0}%
\e{0}%
\e{1}%
\e{0}%
\e{0}%
\e{0}%
\eol}\vss}\rg%
%
%
\rx{\vss\hfull{%
\rlx{\hss{$4200_{y}$}}\cg%
\e{3}%
\e{1}%
\e{0}%
\e{0}%
\e{3}%
\e{1}%
\e{2}%
\e{2}%
\e{1}%
\e{0}%
\e{1}%
\e{0}%
\e{0}%
\e{0}%
\e{0}%
\eol}\vss}\rg%
\eop
\eject
\tablecont%
%
%
%
%
%
%
\rowpts=18 true pt%
\colpts=18 true pt%
\rowlabpts=40 true pt%
\collabpts=45 true pt%
\clx{\vss\hfull{%
\rlx{\hss{$ $}}\cg%
\cx{\hskip 16 true pt\flip{$[{4}{2^{2}}{1}]$}\hss}\cg%
\cx{\hskip 16 true pt\flip{$[{4}{2}{1^{3}}]$}\hss}\cg%
\cx{\hskip 16 true pt\flip{$[{4}{1^{5}}]$}\hss}\cg%
\cx{\hskip 16 true pt\flip{$[{3^{3}}]$}\hss}\cg%
\cx{\hskip 16 true pt\flip{$[{3^{2}}{2}{1}]$}\hss}\cg%
\cx{\hskip 16 true pt\flip{$[{3^{2}}{1^{3}}]$}\hss}\cg%
\cx{\hskip 16 true pt\flip{$[{3}{2^{3}}]$}\hss}\cg%
\cx{\hskip 16 true pt\flip{$[{3}{2^{2}}{1^{2}}]$}\hss}\cg%
\cx{\hskip 16 true pt\flip{$[{3}{2}{1^{4}}]$}\hss}\cg%
\cx{\hskip 16 true pt\flip{$[{3}{1^{6}}]$}\hss}\cg%
\cx{\hskip 16 true pt\flip{$[{2^{4}}{1}]$}\hss}\cg%
\cx{\hskip 16 true pt\flip{$[{2^{3}}{1^{3}}]$}\hss}\cg%
\cx{\hskip 16 true pt\flip{$[{2^{2}}{1^{5}}]$}\hss}\cg%
\cx{\hskip 16 true pt\flip{$[{2}{1^{7}}]$}\hss}\cg%
\cx{\hskip 16 true pt\flip{$[{1^{9}}]$}\hss}\cg%
\eol}}\rg%
%
%
\rx{\vss\hfull{%
\rlx{\hss{$2688_{y}$}}\cg%
\e{2}%
\e{1}%
\e{0}%
\e{0}%
\e{2}%
\e{1}%
\e{0}%
\e{1}%
\e{1}%
\e{0}%
\e{0}%
\e{0}%
\e{0}%
\e{0}%
\e{0}%
\eol}\vss}\rg%
%
%
\rx{\vss\hfull{%
\rlx{\hss{$2100_{y}$}}\cg%
\e{1}%
\e{2}%
\e{1}%
\e{0}%
\e{1}%
\e{0}%
\e{0}%
\e{1}%
\e{0}%
\e{0}%
\e{0}%
\e{0}%
\e{0}%
\e{0}%
\e{0}%
\eol}\vss}\rg%
%
%
\rx{\vss\hfull{%
\rlx{\hss{$1400_{y}$}}\cg%
\e{1}%
\e{1}%
\e{1}%
\e{0}%
\e{0}%
\e{1}%
\e{1}%
\e{0}%
\e{0}%
\e{0}%
\e{0}%
\e{0}%
\e{0}%
\e{0}%
\e{0}%
\eol}\vss}\rg%
%
%
\rx{\vss\hfull{%
\rlx{\hss{$4536_{y}$}}\cg%
\e{3}%
\e{3}%
\e{1}%
\e{0}%
\e{2}%
\e{2}%
\e{1}%
\e{1}%
\e{1}%
\e{0}%
\e{0}%
\e{0}%
\e{0}%
\e{0}%
\e{0}%
\eol}\vss}\rg%
%
%
\rx{\vss\hfull{%
\rlx{\hss{$5670_{y}$}}\cg%
\e{4}%
\e{4}%
\e{1}%
\e{0}%
\e{2}%
\e{2}%
\e{1}%
\e{2}%
\e{1}%
\e{0}%
\e{0}%
\e{0}%
\e{0}%
\e{0}%
\e{0}%
\eol}\vss}\rg%
%
%
\rx{\vss\hfull{%
\rlx{\hss{$4480_{y}$}}\cg%
\e{3}%
\e{2}%
\e{0}%
\e{0}%
\e{3}%
\e{1}%
\e{1}%
\e{2}%
\e{1}%
\e{0}%
\e{1}%
\e{0}%
\e{0}%
\e{0}%
\e{0}%
\eol}\vss}\rg%
%
%
\rx{\vss\hfull{%
\rlx{\hss{$8_{z}$}}\cg%
\e{0}%
\e{0}%
\e{0}%
\e{0}%
\e{0}%
\e{0}%
\e{0}%
\e{0}%
\e{0}%
\e{0}%
\e{0}%
\e{0}%
\e{0}%
\e{0}%
\e{0}%
\eol}\vss}\rg%
%
%
\rx{\vss\hfull{%
\rlx{\hss{$56_{z}$}}\cg%
\e{0}%
\e{0}%
\e{0}%
\e{0}%
\e{0}%
\e{0}%
\e{0}%
\e{0}%
\e{0}%
\e{0}%
\e{0}%
\e{0}%
\e{0}%
\e{0}%
\e{0}%
\eol}\vss}\rg%
%
%
\rx{\vss\hfull{%
\rlx{\hss{$160_{z}$}}\cg%
\e{0}%
\e{0}%
\e{0}%
\e{0}%
\e{0}%
\e{0}%
\e{0}%
\e{0}%
\e{0}%
\e{0}%
\e{0}%
\e{0}%
\e{0}%
\e{0}%
\e{0}%
\eol}\vss}\rg%
%
%
\rx{\vss\hfull{%
\rlx{\hss{$112_{z}$}}\cg%
\e{0}%
\e{0}%
\e{0}%
\e{0}%
\e{0}%
\e{0}%
\e{0}%
\e{0}%
\e{0}%
\e{0}%
\e{0}%
\e{0}%
\e{0}%
\e{0}%
\e{0}%
\eol}\vss}\rg%
%
%
\rx{\vss\hfull{%
\rlx{\hss{$840_{z}$}}\cg%
\e{1}%
\e{0}%
\e{0}%
\e{0}%
\e{0}%
\e{0}%
\e{0}%
\e{0}%
\e{0}%
\e{0}%
\e{0}%
\e{0}%
\e{0}%
\e{0}%
\e{0}%
\eol}\vss}\rg%
%
%
\rx{\vss\hfull{%
\rlx{\hss{$1296_{z}$}}\cg%
\e{0}%
\e{0}%
\e{0}%
\e{0}%
\e{0}%
\e{0}%
\e{0}%
\e{0}%
\e{0}%
\e{0}%
\e{0}%
\e{0}%
\e{0}%
\e{0}%
\e{0}%
\eol}\vss}\rg%
%
%
\rx{\vss\hfull{%
\rlx{\hss{$1400_{z}$}}\cg%
\e{0}%
\e{0}%
\e{0}%
\e{0}%
\e{0}%
\e{0}%
\e{0}%
\e{0}%
\e{0}%
\e{0}%
\e{0}%
\e{0}%
\e{0}%
\e{0}%
\e{0}%
\eol}\vss}\rg%
%
%
\rx{\vss\hfull{%
\rlx{\hss{$1008_{z}$}}\cg%
\e{0}%
\e{0}%
\e{0}%
\e{0}%
\e{0}%
\e{0}%
\e{0}%
\e{0}%
\e{0}%
\e{0}%
\e{0}%
\e{0}%
\e{0}%
\e{0}%
\e{0}%
\eol}\vss}\rg%
%
%
\rx{\vss\hfull{%
\rlx{\hss{$560_{z}$}}\cg%
\e{0}%
\e{0}%
\e{0}%
\e{0}%
\e{0}%
\e{0}%
\e{0}%
\e{0}%
\e{0}%
\e{0}%
\e{0}%
\e{0}%
\e{0}%
\e{0}%
\e{0}%
\eol}\vss}\rg%
%
%
\rx{\vss\hfull{%
\rlx{\hss{$1400_{zz}$}}\cg%
\e{0}%
\e{0}%
\e{0}%
\e{0}%
\e{1}%
\e{0}%
\e{0}%
\e{0}%
\e{0}%
\e{0}%
\e{0}%
\e{0}%
\e{0}%
\e{0}%
\e{0}%
\eol}\vss}\rg%
%
%
\rx{\vss\hfull{%
\rlx{\hss{$4200_{z}$}}\cg%
\e{2}%
\e{1}%
\e{0}%
\e{1}%
\e{2}%
\e{2}%
\e{1}%
\e{1}%
\e{0}%
\e{0}%
\e{0}%
\e{0}%
\e{0}%
\e{0}%
\e{0}%
\eol}\vss}\rg%
%
%
\rx{\vss\hfull{%
\rlx{\hss{$400_{z}$}}\cg%
\e{0}%
\e{0}%
\e{0}%
\e{0}%
\e{0}%
\e{0}%
\e{0}%
\e{0}%
\e{0}%
\e{0}%
\e{0}%
\e{0}%
\e{0}%
\e{0}%
\e{0}%
\eol}\vss}\rg%
%
%
\rx{\vss\hfull{%
\rlx{\hss{$3240_{z}$}}\cg%
\e{1}%
\e{0}%
\e{0}%
\e{0}%
\e{0}%
\e{0}%
\e{0}%
\e{0}%
\e{0}%
\e{0}%
\e{0}%
\e{0}%
\e{0}%
\e{0}%
\e{0}%
\eol}\vss}\rg%
%
%
\rx{\vss\hfull{%
\rlx{\hss{$4536_{z}$}}\cg%
\e{3}%
\e{1}%
\e{0}%
\e{1}%
\e{2}%
\e{0}%
\e{1}%
\e{1}%
\e{0}%
\e{0}%
\e{0}%
\e{0}%
\e{0}%
\e{0}%
\e{0}%
\eol}\vss}\rg%
%
%
\rx{\vss\hfull{%
\rlx{\hss{$2400_{z}$}}\cg%
\e{1}%
\e{2}%
\e{0}%
\e{0}%
\e{0}%
\e{1}%
\e{0}%
\e{0}%
\e{0}%
\e{0}%
\e{0}%
\e{0}%
\e{0}%
\e{0}%
\e{0}%
\eol}\vss}\rg%
%
%
\rx{\vss\hfull{%
\rlx{\hss{$3360_{z}$}}\cg%
\e{1}%
\e{1}%
\e{0}%
\e{1}%
\e{1}%
\e{1}%
\e{0}%
\e{0}%
\e{0}%
\e{0}%
\e{0}%
\e{0}%
\e{0}%
\e{0}%
\e{0}%
\eol}\vss}\rg%
%
%
\rx{\vss\hfull{%
\rlx{\hss{$2800_{z}$}}\cg%
\e{1}%
\e{1}%
\e{0}%
\e{0}%
\e{1}%
\e{0}%
\e{0}%
\e{0}%
\e{0}%
\e{0}%
\e{0}%
\e{0}%
\e{0}%
\e{0}%
\e{0}%
\eol}\vss}\rg%
%
%
\rx{\vss\hfull{%
\rlx{\hss{$4096_{z}$}}\cg%
\e{2}%
\e{1}%
\e{0}%
\e{0}%
\e{1}%
\e{0}%
\e{0}%
\e{0}%
\e{0}%
\e{0}%
\e{0}%
\e{0}%
\e{0}%
\e{0}%
\e{0}%
\eol}\vss}\rg%
%
%
\rx{\vss\hfull{%
\rlx{\hss{$5600_{z}$}}\cg%
\e{4}%
\e{3}%
\e{1}%
\e{1}%
\e{1}%
\e{1}%
\e{1}%
\e{1}%
\e{0}%
\e{0}%
\e{0}%
\e{0}%
\e{0}%
\e{0}%
\e{0}%
\eol}\vss}\rg%
%
%
\rx{\vss\hfull{%
\rlx{\hss{$448_{z}$}}\cg%
\e{0}%
\e{0}%
\e{0}%
\e{1}%
\e{0}%
\e{0}%
\e{0}%
\e{0}%
\e{0}%
\e{0}%
\e{0}%
\e{0}%
\e{0}%
\e{0}%
\e{0}%
\eol}\vss}\rg%
%
%
\rx{\vss\hfull{%
\rlx{\hss{$448_{w}$}}\cg%
\e{0}%
\e{1}%
\e{0}%
\e{0}%
\e{0}%
\e{0}%
\e{0}%
\e{0}%
\e{0}%
\e{0}%
\e{0}%
\e{0}%
\e{0}%
\e{0}%
\e{0}%
\eol}\vss}\rg%
%
%
\rx{\vss\hfull{%
\rlx{\hss{$1344_{w}$}}\cg%
\e{1}%
\e{0}%
\e{0}%
\e{0}%
\e{1}%
\e{0}%
\e{1}%
\e{1}%
\e{0}%
\e{0}%
\e{1}%
\e{0}%
\e{0}%
\e{0}%
\e{0}%
\eol}\vss}\rg%
%
%
\rx{\vss\hfull{%
\rlx{\hss{$5600_{w}$}}\cg%
\e{4}%
\e{4}%
\e{1}%
\e{0}%
\e{2}%
\e{2}%
\e{1}%
\e{2}%
\e{1}%
\e{0}%
\e{0}%
\e{0}%
\e{0}%
\e{0}%
\e{0}%
\eol}\vss}\rg%
%
%
\rx{\vss\hfull{%
\rlx{\hss{$2016_{w}$}}\cg%
\e{1}%
\e{0}%
\e{0}%
\e{2}%
\e{2}%
\e{1}%
\e{1}%
\e{1}%
\e{0}%
\e{0}%
\e{0}%
\e{1}%
\e{0}%
\e{0}%
\e{0}%
\eol}\vss}\rg%
%
%
\rx{\vss\hfull{%
\rlx{\hss{$7168_{w}$}}\cg%
\e{5}%
\e{3}%
\e{0}%
\e{0}%
\e{5}%
\e{2}%
\e{1}%
\e{3}%
\e{2}%
\e{0}%
\e{1}%
\e{0}%
\e{0}%
\e{0}%
\e{0}%
\eol}\vss}\rg%
\tableclose%
%
%
%
%
%
%
\eop
\eject
\tableopen{Induce/restrict matrix for $W({A_{7}}{A_{1}})\,\subset\,W(E_{8})$}%
%
%
%
%
%
%
\rowpts=18 true pt%
\colpts=18 true pt%
\rowlabpts=40 true pt%
\collabpts=65 true pt%
\clx{\vss\hfull{%
\rlx{\hss{$ $}}\cg%
\cx{\hskip 16 true pt\flip{$[{8}]{\times}[{2}]$}\hss}\cg%
\cx{\hskip 16 true pt\flip{$[{7}{1}]{\times}[{2}]$}\hss}\cg%
\cx{\hskip 16 true pt\flip{$[{6}{2}]{\times}[{2}]$}\hss}\cg%
\cx{\hskip 16 true pt\flip{$[{6}{1^{2}}]{\times}[{2}]$}\hss}\cg%
\cx{\hskip 16 true pt\flip{$[{5}{3}]{\times}[{2}]$}\hss}\cg%
\cx{\hskip 16 true pt\flip{$[{5}{2}{1}]{\times}[{2}]$}\hss}\cg%
\cx{\hskip 16 true pt\flip{$[{5}{1^{3}}]{\times}[{2}]$}\hss}\cg%
\cx{\hskip 16 true pt\flip{$[{4^{2}}]{\times}[{2}]$}\hss}\cg%
\cx{\hskip 16 true pt\flip{$[{4}{3}{1}]{\times}[{2}]$}\hss}\cg%
\cx{\hskip 16 true pt\flip{$[{4}{2^{2}}]{\times}[{2}]$}\hss}\cg%
\cx{\hskip 16 true pt\flip{$[{4}{2}{1^{2}}]{\times}[{2}]$}\hss}\cg%
\cx{\hskip 16 true pt\flip{$[{4}{1^{4}}]{\times}[{2}]$}\hss}\cg%
\cx{\hskip 16 true pt\flip{$[{3^{2}}{2}]{\times}[{2}]$}\hss}\cg%
\cx{\hskip 16 true pt\flip{$[{3^{2}}{1^{2}}]{\times}[{2}]$}\hss}\cg%
\cx{\hskip 16 true pt\flip{$[{3}{2^{2}}{1}]{\times}[{2}]$}\hss}\cg%
\cx{\hskip 16 true pt\flip{$[{3}{2}{1^{3}}]{\times}[{2}]$}\hss}\cg%
\cx{\hskip 16 true pt\flip{$[{3}{1^{5}}]{\times}[{2}]$}\hss}\cg%
\cx{\hskip 16 true pt\flip{$[{2^{4}}]{\times}[{2}]$}\hss}\cg%
\cx{\hskip 16 true pt\flip{$[{2^{3}}{1^{2}}]{\times}[{2}]$}\hss}\cg%
\cx{\hskip 16 true pt\flip{$[{2^{2}}{1^{4}}]{\times}[{2}]$}\hss}\cg%
\cx{\hskip 16 true pt\flip{$[{2}{1^{6}}]{\times}[{2}]$}\hss}\cg%
\cx{\hskip 16 true pt\flip{$[{1^{8}}]{\times}[{2}]$}\hss}\cg%
\eol}}\rg%
%
%
\rx{\vss\hfull{%
\rlx{\hss{$1_{x}$}}\cg%
\e{1}%
\e{0}%
\e{0}%
\e{0}%
\e{0}%
\e{0}%
\e{0}%
\e{0}%
\e{0}%
\e{0}%
\e{0}%
\e{0}%
\e{0}%
\e{0}%
\e{0}%
\e{0}%
\e{0}%
\e{0}%
\e{0}%
\e{0}%
\e{0}%
\e{0}%
\eol}\vss}\rg%
%
%
\rx{\vss\hfull{%
\rlx{\hss{$28_{x}$}}\cg%
\e{0}%
\e{0}%
\e{0}%
\e{1}%
\e{0}%
\e{0}%
\e{0}%
\e{0}%
\e{0}%
\e{0}%
\e{0}%
\e{0}%
\e{0}%
\e{0}%
\e{0}%
\e{0}%
\e{0}%
\e{0}%
\e{0}%
\e{0}%
\e{0}%
\e{0}%
\eol}\vss}\rg%
%
%
\rx{\vss\hfull{%
\rlx{\hss{$35_{x}$}}\cg%
\e{1}%
\e{1}%
\e{1}%
\e{0}%
\e{0}%
\e{0}%
\e{0}%
\e{0}%
\e{0}%
\e{0}%
\e{0}%
\e{0}%
\e{0}%
\e{0}%
\e{0}%
\e{0}%
\e{0}%
\e{0}%
\e{0}%
\e{0}%
\e{0}%
\e{0}%
\eol}\vss}\rg%
%
%
\rx{\vss\hfull{%
\rlx{\hss{$84_{x}$}}\cg%
\e{2}%
\e{1}%
\e{2}%
\e{0}%
\e{0}%
\e{0}%
\e{0}%
\e{1}%
\e{0}%
\e{0}%
\e{0}%
\e{0}%
\e{0}%
\e{0}%
\e{0}%
\e{0}%
\e{0}%
\e{0}%
\e{0}%
\e{0}%
\e{0}%
\e{0}%
\eol}\vss}\rg%
%
%
\rx{\vss\hfull{%
\rlx{\hss{$50_{x}$}}\cg%
\e{1}%
\e{0}%
\e{1}%
\e{0}%
\e{0}%
\e{0}%
\e{0}%
\e{1}%
\e{0}%
\e{0}%
\e{0}%
\e{0}%
\e{0}%
\e{0}%
\e{0}%
\e{0}%
\e{0}%
\e{0}%
\e{0}%
\e{0}%
\e{0}%
\e{0}%
\eol}\vss}\rg%
%
%
\rx{\vss\hfull{%
\rlx{\hss{$350_{x}$}}\cg%
\e{0}%
\e{0}%
\e{0}%
\e{1}%
\e{0}%
\e{1}%
\e{1}%
\e{0}%
\e{0}%
\e{0}%
\e{1}%
\e{0}%
\e{0}%
\e{0}%
\e{0}%
\e{0}%
\e{0}%
\e{0}%
\e{0}%
\e{0}%
\e{0}%
\e{0}%
\eol}\vss}\rg%
%
%
\rx{\vss\hfull{%
\rlx{\hss{$300_{x}$}}\cg%
\e{0}%
\e{1}%
\e{2}%
\e{0}%
\e{1}%
\e{1}%
\e{0}%
\e{0}%
\e{0}%
\e{1}%
\e{0}%
\e{0}%
\e{0}%
\e{0}%
\e{0}%
\e{0}%
\e{0}%
\e{0}%
\e{0}%
\e{0}%
\e{0}%
\e{0}%
\eol}\vss}\rg%
%
%
\rx{\vss\hfull{%
\rlx{\hss{$567_{x}$}}\cg%
\e{0}%
\e{2}%
\e{2}%
\e{3}%
\e{1}%
\e{2}%
\e{1}%
\e{0}%
\e{1}%
\e{0}%
\e{0}%
\e{0}%
\e{0}%
\e{0}%
\e{0}%
\e{0}%
\e{0}%
\e{0}%
\e{0}%
\e{0}%
\e{0}%
\e{0}%
\eol}\vss}\rg%
%
%
\rx{\vss\hfull{%
\rlx{\hss{$210_{x}$}}\cg%
\e{0}%
\e{2}%
\e{1}%
\e{1}%
\e{1}%
\e{1}%
\e{0}%
\e{0}%
\e{0}%
\e{0}%
\e{0}%
\e{0}%
\e{0}%
\e{0}%
\e{0}%
\e{0}%
\e{0}%
\e{0}%
\e{0}%
\e{0}%
\e{0}%
\e{0}%
\eol}\vss}\rg%
%
%
\rx{\vss\hfull{%
\rlx{\hss{$840_{x}$}}\cg%
\e{0}%
\e{0}%
\e{1}%
\e{0}%
\e{2}%
\e{1}%
\e{0}%
\e{1}%
\e{1}%
\e{2}%
\e{0}%
\e{0}%
\e{1}%
\e{0}%
\e{1}%
\e{0}%
\e{0}%
\e{1}%
\e{0}%
\e{0}%
\e{0}%
\e{0}%
\eol}\vss}\rg%
%
%
\rx{\vss\hfull{%
\rlx{\hss{$700_{x}$}}\cg%
\e{1}%
\e{3}%
\e{3}%
\e{1}%
\e{3}%
\e{2}%
\e{0}%
\e{1}%
\e{1}%
\e{1}%
\e{0}%
\e{0}%
\e{0}%
\e{0}%
\e{0}%
\e{0}%
\e{0}%
\e{0}%
\e{0}%
\e{0}%
\e{0}%
\e{0}%
\eol}\vss}\rg%
%
%
\rx{\vss\hfull{%
\rlx{\hss{$175_{x}$}}\cg%
\e{0}%
\e{1}%
\e{0}%
\e{0}%
\e{1}%
\e{0}%
\e{0}%
\e{0}%
\e{1}%
\e{0}%
\e{0}%
\e{0}%
\e{0}%
\e{0}%
\e{0}%
\e{0}%
\e{0}%
\e{0}%
\e{0}%
\e{0}%
\e{0}%
\e{0}%
\eol}\vss}\rg%
%
%
\rx{\vss\hfull{%
\rlx{\hss{$1400_{x}$}}\cg%
\e{0}%
\e{2}%
\e{1}%
\e{3}%
\e{3}%
\e{3}%
\e{1}%
\e{0}%
\e{3}%
\e{0}%
\e{2}%
\e{0}%
\e{1}%
\e{0}%
\e{0}%
\e{0}%
\e{0}%
\e{0}%
\e{0}%
\e{0}%
\e{0}%
\e{0}%
\eol}\vss}\rg%
%
%
\rx{\vss\hfull{%
\rlx{\hss{$1050_{x}$}}\cg%
\e{1}%
\e{1}%
\e{3}%
\e{1}%
\e{1}%
\e{2}%
\e{1}%
\e{2}%
\e{3}%
\e{1}%
\e{0}%
\e{0}%
\e{0}%
\e{1}%
\e{0}%
\e{0}%
\e{0}%
\e{0}%
\e{0}%
\e{0}%
\e{0}%
\e{0}%
\eol}\vss}\rg%
%
%
\rx{\vss\hfull{%
\rlx{\hss{$1575_{x}$}}\cg%
\e{0}%
\e{1}%
\e{1}%
\e{4}%
\e{2}%
\e{4}%
\e{2}%
\e{0}%
\e{2}%
\e{0}%
\e{3}%
\e{0}%
\e{1}%
\e{0}%
\e{0}%
\e{0}%
\e{0}%
\e{0}%
\e{0}%
\e{0}%
\e{0}%
\e{0}%
\eol}\vss}\rg%
%
%
\rx{\vss\hfull{%
\rlx{\hss{$1344_{x}$}}\cg%
\e{1}%
\e{2}%
\e{5}%
\e{2}%
\e{2}%
\e{4}%
\e{1}%
\e{2}%
\e{2}%
\e{2}%
\e{0}%
\e{0}%
\e{0}%
\e{1}%
\e{0}%
\e{0}%
\e{0}%
\e{0}%
\e{0}%
\e{0}%
\e{0}%
\e{0}%
\eol}\vss}\rg%
%
%
\rx{\vss\hfull{%
\rlx{\hss{$2100_{x}$}}\cg%
\e{0}%
\e{0}%
\e{1}%
\e{1}%
\e{0}%
\e{3}%
\e{4}%
\e{0}%
\e{2}%
\e{2}%
\e{3}%
\e{2}%
\e{0}%
\e{1}%
\e{1}%
\e{1}%
\e{0}%
\e{0}%
\e{0}%
\e{0}%
\e{0}%
\e{0}%
\eol}\vss}\rg%
%
%
\rx{\vss\hfull{%
\rlx{\hss{$2268_{x}$}}\cg%
\e{0}%
\e{1}%
\e{2}%
\e{3}%
\e{3}%
\e{5}%
\e{2}%
\e{0}%
\e{3}%
\e{2}%
\e{3}%
\e{1}%
\e{1}%
\e{0}%
\e{1}%
\e{0}%
\e{0}%
\e{0}%
\e{0}%
\e{0}%
\e{0}%
\e{0}%
\eol}\vss}\rg%
%
%
\rx{\vss\hfull{%
\rlx{\hss{$525_{x}$}}\cg%
\e{0}%
\e{0}%
\e{1}%
\e{1}%
\e{0}%
\e{1}%
\e{2}%
\e{1}%
\e{1}%
\e{0}%
\e{0}%
\e{0}%
\e{0}%
\e{1}%
\e{0}%
\e{0}%
\e{0}%
\e{0}%
\e{0}%
\e{0}%
\e{0}%
\e{0}%
\eol}\vss}\rg%
%
%
\rx{\vss\hfull{%
\rlx{\hss{$700_{xx}$}}\cg%
\e{0}%
\e{0}%
\e{1}%
\e{0}%
\e{0}%
\e{1}%
\e{1}%
\e{2}%
\e{1}%
\e{1}%
\e{0}%
\e{0}%
\e{0}%
\e{2}%
\e{0}%
\e{0}%
\e{0}%
\e{0}%
\e{0}%
\e{0}%
\e{0}%
\e{0}%
\eol}\vss}\rg%
%
%
\rx{\vss\hfull{%
\rlx{\hss{$972_{x}$}}\cg%
\e{1}%
\e{0}%
\e{3}%
\e{0}%
\e{1}%
\e{2}%
\e{0}%
\e{3}%
\e{1}%
\e{3}%
\e{0}%
\e{0}%
\e{0}%
\e{1}%
\e{0}%
\e{0}%
\e{0}%
\e{1}%
\e{0}%
\e{0}%
\e{0}%
\e{0}%
\eol}\vss}\rg%
%
%
\rx{\vss\hfull{%
\rlx{\hss{$4096_{x}$}}\cg%
\e{0}%
\e{1}%
\e{3}%
\e{3}%
\e{3}%
\e{8}%
\e{3}%
\e{1}%
\e{5}%
\e{4}%
\e{5}%
\e{1}%
\e{2}%
\e{2}%
\e{2}%
\e{1}%
\e{0}%
\e{0}%
\e{0}%
\e{0}%
\e{0}%
\e{0}%
\eol}\vss}\rg%
%
%
\rx{\vss\hfull{%
\rlx{\hss{$4200_{x}$}}\cg%
\e{0}%
\e{1}%
\e{3}%
\e{2}%
\e{4}%
\e{6}%
\e{2}%
\e{1}%
\e{7}%
\e{4}%
\e{4}%
\e{1}%
\e{3}%
\e{2}%
\e{3}%
\e{1}%
\e{0}%
\e{0}%
\e{0}%
\e{0}%
\e{0}%
\e{0}%
\eol}\vss}\rg%
%
%
\rx{\vss\hfull{%
\rlx{\hss{$2240_{x}$}}\cg%
\e{0}%
\e{2}%
\e{2}%
\e{2}%
\e{5}%
\e{4}%
\e{0}%
\e{1}%
\e{4}%
\e{2}%
\e{2}%
\e{0}%
\e{2}%
\e{1}%
\e{1}%
\e{0}%
\e{0}%
\e{0}%
\e{0}%
\e{0}%
\e{0}%
\e{0}%
\eol}\vss}\rg%
%
%
\rx{\vss\hfull{%
\rlx{\hss{$2835_{x}$}}\cg%
\e{0}%
\e{1}%
\e{1}%
\e{1}%
\e{3}%
\e{3}%
\e{0}%
\e{1}%
\e{5}%
\e{2}%
\e{3}%
\e{0}%
\e{3}%
\e{2}%
\e{2}%
\e{1}%
\e{0}%
\e{0}%
\e{0}%
\e{0}%
\e{0}%
\e{0}%
\eol}\vss}\rg%
%
%
\rx{\vss\hfull{%
\rlx{\hss{$6075_{x}$}}\cg%
\e{0}%
\e{0}%
\e{2}%
\e{3}%
\e{2}%
\e{7}%
\e{5}%
\e{2}%
\e{7}%
\e{4}%
\e{9}%
\e{2}%
\e{3}%
\e{5}%
\e{3}%
\e{3}%
\e{0}%
\e{0}%
\e{1}%
\e{0}%
\e{0}%
\e{0}%
\eol}\vss}\rg%
%
%
\rx{\vss\hfull{%
\rlx{\hss{$3200_{x}$}}\cg%
\e{0}%
\e{0}%
\e{2}%
\e{0}%
\e{1}%
\e{4}%
\e{1}%
\e{2}%
\e{3}%
\e{6}%
\e{2}%
\e{1}%
\e{1}%
\e{3}%
\e{3}%
\e{1}%
\e{0}%
\e{2}%
\e{0}%
\e{1}%
\e{0}%
\e{0}%
\eol}\vss}\rg%
%
%
\rx{\vss\hfull{%
\rlx{\hss{$70_{y}$}}\cg%
\e{0}%
\e{0}%
\e{0}%
\e{0}%
\e{0}%
\e{0}%
\e{0}%
\e{0}%
\e{0}%
\e{0}%
\e{0}%
\e{1}%
\e{0}%
\e{0}%
\e{0}%
\e{0}%
\e{0}%
\e{0}%
\e{0}%
\e{0}%
\e{0}%
\e{0}%
\eol}\vss}\rg%
%
%
\rx{\vss\hfull{%
\rlx{\hss{$1134_{y}$}}\cg%
\e{0}%
\e{0}%
\e{0}%
\e{0}%
\e{0}%
\e{1}%
\e{1}%
\e{0}%
\e{1}%
\e{0}%
\e{2}%
\e{0}%
\e{0}%
\e{1}%
\e{1}%
\e{1}%
\e{0}%
\e{0}%
\e{1}%
\e{0}%
\e{0}%
\e{0}%
\eol}\vss}\rg%
%
%
\rx{\vss\hfull{%
\rlx{\hss{$1680_{y}$}}\cg%
\e{0}%
\e{0}%
\e{0}%
\e{0}%
\e{0}%
\e{1}%
\e{2}%
\e{0}%
\e{0}%
\e{1}%
\e{3}%
\e{3}%
\e{0}%
\e{1}%
\e{1}%
\e{2}%
\e{1}%
\e{0}%
\e{0}%
\e{0}%
\e{0}%
\e{0}%
\eol}\vss}\rg%
%
%
\rx{\vss\hfull{%
\rlx{\hss{$168_{y}$}}\cg%
\e{0}%
\e{0}%
\e{0}%
\e{0}%
\e{0}%
\e{0}%
\e{0}%
\e{1}%
\e{0}%
\e{1}%
\e{0}%
\e{0}%
\e{0}%
\e{0}%
\e{0}%
\e{0}%
\e{0}%
\e{1}%
\e{0}%
\e{0}%
\e{0}%
\e{0}%
\eol}\vss}\rg%
%
%
\rx{\vss\hfull{%
\rlx{\hss{$420_{y}$}}\cg%
\e{0}%
\e{0}%
\e{0}%
\e{0}%
\e{1}%
\e{0}%
\e{0}%
\e{0}%
\e{1}%
\e{0}%
\e{0}%
\e{0}%
\e{1}%
\e{0}%
\e{1}%
\e{0}%
\e{0}%
\e{0}%
\e{0}%
\e{0}%
\e{0}%
\e{0}%
\eol}\vss}\rg%
%
%
\rx{\vss\hfull{%
\rlx{\hss{$3150_{y}$}}\cg%
\e{0}%
\e{0}%
\e{0}%
\e{1}%
\e{2}%
\e{2}%
\e{0}%
\e{0}%
\e{4}%
\e{1}%
\e{4}%
\e{0}%
\e{4}%
\e{2}%
\e{3}%
\e{2}%
\e{0}%
\e{0}%
\e{2}%
\e{0}%
\e{0}%
\e{0}%
\eol}\vss}\rg%
%
%
\rx{\vss\hfull{%
\rlx{\hss{$4200_{y}$}}\cg%
\e{0}%
\e{0}%
\e{1}%
\e{0}%
\e{1}%
\e{3}%
\e{1}%
\e{2}%
\e{5}%
\e{5}%
\e{3}%
\e{1}%
\e{2}%
\e{4}%
\e{5}%
\e{2}%
\e{0}%
\e{2}%
\e{1}%
\e{1}%
\e{0}%
\e{0}%
\eol}\vss}\rg%
\eop
\eject
\tablecont%
%
%
%
%
%
%
\rowpts=18 true pt%
\colpts=18 true pt%
\rowlabpts=40 true pt%
\collabpts=65 true pt%
\clx{\vss\hfull{%
\rlx{\hss{$ $}}\cg%
\cx{\hskip 16 true pt\flip{$[{8}]{\times}[{2}]$}\hss}\cg%
\cx{\hskip 16 true pt\flip{$[{7}{1}]{\times}[{2}]$}\hss}\cg%
\cx{\hskip 16 true pt\flip{$[{6}{2}]{\times}[{2}]$}\hss}\cg%
\cx{\hskip 16 true pt\flip{$[{6}{1^{2}}]{\times}[{2}]$}\hss}\cg%
\cx{\hskip 16 true pt\flip{$[{5}{3}]{\times}[{2}]$}\hss}\cg%
\cx{\hskip 16 true pt\flip{$[{5}{2}{1}]{\times}[{2}]$}\hss}\cg%
\cx{\hskip 16 true pt\flip{$[{5}{1^{3}}]{\times}[{2}]$}\hss}\cg%
\cx{\hskip 16 true pt\flip{$[{4^{2}}]{\times}[{2}]$}\hss}\cg%
\cx{\hskip 16 true pt\flip{$[{4}{3}{1}]{\times}[{2}]$}\hss}\cg%
\cx{\hskip 16 true pt\flip{$[{4}{2^{2}}]{\times}[{2}]$}\hss}\cg%
\cx{\hskip 16 true pt\flip{$[{4}{2}{1^{2}}]{\times}[{2}]$}\hss}\cg%
\cx{\hskip 16 true pt\flip{$[{4}{1^{4}}]{\times}[{2}]$}\hss}\cg%
\cx{\hskip 16 true pt\flip{$[{3^{2}}{2}]{\times}[{2}]$}\hss}\cg%
\cx{\hskip 16 true pt\flip{$[{3^{2}}{1^{2}}]{\times}[{2}]$}\hss}\cg%
\cx{\hskip 16 true pt\flip{$[{3}{2^{2}}{1}]{\times}[{2}]$}\hss}\cg%
\cx{\hskip 16 true pt\flip{$[{3}{2}{1^{3}}]{\times}[{2}]$}\hss}\cg%
\cx{\hskip 16 true pt\flip{$[{3}{1^{5}}]{\times}[{2}]$}\hss}\cg%
\cx{\hskip 16 true pt\flip{$[{2^{4}}]{\times}[{2}]$}\hss}\cg%
\cx{\hskip 16 true pt\flip{$[{2^{3}}{1^{2}}]{\times}[{2}]$}\hss}\cg%
\cx{\hskip 16 true pt\flip{$[{2^{2}}{1^{4}}]{\times}[{2}]$}\hss}\cg%
\cx{\hskip 16 true pt\flip{$[{2}{1^{6}}]{\times}[{2}]$}\hss}\cg%
\cx{\hskip 16 true pt\flip{$[{1^{8}}]{\times}[{2}]$}\hss}\cg%
\eol}}\rg%
%
%
\rx{\vss\hfull{%
\rlx{\hss{$2688_{y}$}}\cg%
\e{0}%
\e{0}%
\e{1}%
\e{0}%
\e{0}%
\e{2}%
\e{1}%
\e{1}%
\e{2}%
\e{4}%
\e{2}%
\e{1}%
\e{1}%
\e{4}%
\e{2}%
\e{2}%
\e{0}%
\e{1}%
\e{0}%
\e{1}%
\e{0}%
\e{0}%
\eol}\vss}\rg%
%
%
\rx{\vss\hfull{%
\rlx{\hss{$2100_{y}$}}\cg%
\e{0}%
\e{0}%
\e{0}%
\e{0}%
\e{0}%
\e{1}%
\e{3}%
\e{1}%
\e{1}%
\e{2}%
\e{2}%
\e{3}%
\e{0}%
\e{3}%
\e{1}%
\e{2}%
\e{0}%
\e{1}%
\e{0}%
\e{1}%
\e{0}%
\e{0}%
\eol}\vss}\rg%
%
%
\rx{\vss\hfull{%
\rlx{\hss{$1400_{y}$}}\cg%
\e{0}%
\e{0}%
\e{0}%
\e{0}%
\e{0}%
\e{1}%
\e{1}%
\e{0}%
\e{1}%
\e{1}%
\e{2}%
\e{2}%
\e{0}%
\e{0}%
\e{2}%
\e{1}%
\e{1}%
\e{0}%
\e{0}%
\e{0}%
\e{0}%
\e{0}%
\eol}\vss}\rg%
%
%
\rx{\vss\hfull{%
\rlx{\hss{$4536_{y}$}}\cg%
\e{0}%
\e{0}%
\e{0}%
\e{1}%
\e{1}%
\e{3}%
\e{2}%
\e{0}%
\e{3}%
\e{3}%
\e{7}%
\e{3}%
\e{3}%
\e{2}%
\e{4}%
\e{4}%
\e{2}%
\e{0}%
\e{1}%
\e{0}%
\e{0}%
\e{0}%
\eol}\vss}\rg%
%
%
\rx{\vss\hfull{%
\rlx{\hss{$5670_{y}$}}\cg%
\e{0}%
\e{0}%
\e{0}%
\e{1}%
\e{1}%
\e{4}%
\e{3}%
\e{0}%
\e{4}%
\e{3}%
\e{9}%
\e{3}%
\e{3}%
\e{3}%
\e{5}%
\e{5}%
\e{2}%
\e{0}%
\e{2}%
\e{0}%
\e{0}%
\e{0}%
\eol}\vss}\rg%
%
%
\rx{\vss\hfull{%
\rlx{\hss{$4480_{y}$}}\cg%
\e{0}%
\e{0}%
\e{0}%
\e{1}%
\e{2}%
\e{3}%
\e{1}%
\e{0}%
\e{5}%
\e{2}%
\e{6}%
\e{1}%
\e{4}%
\e{2}%
\e{5}%
\e{3}%
\e{1}%
\e{0}%
\e{2}%
\e{0}%
\e{0}%
\e{0}%
\eol}\vss}\rg%
%
%
\rx{\vss\hfull{%
\rlx{\hss{$8_{z}$}}\cg%
\e{0}%
\e{1}%
\e{0}%
\e{0}%
\e{0}%
\e{0}%
\e{0}%
\e{0}%
\e{0}%
\e{0}%
\e{0}%
\e{0}%
\e{0}%
\e{0}%
\e{0}%
\e{0}%
\e{0}%
\e{0}%
\e{0}%
\e{0}%
\e{0}%
\e{0}%
\eol}\vss}\rg%
%
%
\rx{\vss\hfull{%
\rlx{\hss{$56_{z}$}}\cg%
\e{0}%
\e{0}%
\e{0}%
\e{0}%
\e{0}%
\e{0}%
\e{1}%
\e{0}%
\e{0}%
\e{0}%
\e{0}%
\e{0}%
\e{0}%
\e{0}%
\e{0}%
\e{0}%
\e{0}%
\e{0}%
\e{0}%
\e{0}%
\e{0}%
\e{0}%
\eol}\vss}\rg%
%
%
\rx{\vss\hfull{%
\rlx{\hss{$160_{z}$}}\cg%
\e{0}%
\e{1}%
\e{1}%
\e{1}%
\e{0}%
\e{1}%
\e{0}%
\e{0}%
\e{0}%
\e{0}%
\e{0}%
\e{0}%
\e{0}%
\e{0}%
\e{0}%
\e{0}%
\e{0}%
\e{0}%
\e{0}%
\e{0}%
\e{0}%
\e{0}%
\eol}\vss}\rg%
%
%
\rx{\vss\hfull{%
\rlx{\hss{$112_{z}$}}\cg%
\e{1}%
\e{2}%
\e{1}%
\e{1}%
\e{1}%
\e{0}%
\e{0}%
\e{0}%
\e{0}%
\e{0}%
\e{0}%
\e{0}%
\e{0}%
\e{0}%
\e{0}%
\e{0}%
\e{0}%
\e{0}%
\e{0}%
\e{0}%
\e{0}%
\e{0}%
\eol}\vss}\rg%
%
%
\rx{\vss\hfull{%
\rlx{\hss{$840_{z}$}}\cg%
\e{0}%
\e{0}%
\e{1}%
\e{1}%
\e{1}%
\e{2}%
\e{0}%
\e{0}%
\e{1}%
\e{1}%
\e{1}%
\e{0}%
\e{0}%
\e{0}%
\e{1}%
\e{0}%
\e{0}%
\e{0}%
\e{0}%
\e{0}%
\e{0}%
\e{0}%
\eol}\vss}\rg%
%
%
\rx{\vss\hfull{%
\rlx{\hss{$1296_{z}$}}\cg%
\e{0}%
\e{0}%
\e{1}%
\e{2}%
\e{0}%
\e{3}%
\e{3}%
\e{0}%
\e{1}%
\e{1}%
\e{2}%
\e{1}%
\e{0}%
\e{1}%
\e{0}%
\e{0}%
\e{0}%
\e{0}%
\e{0}%
\e{0}%
\e{0}%
\e{0}%
\eol}\vss}\rg%
%
%
\rx{\vss\hfull{%
\rlx{\hss{$1400_{z}$}}\cg%
\e{1}%
\e{2}%
\e{4}%
\e{3}%
\e{3}%
\e{4}%
\e{1}%
\e{1}%
\e{2}%
\e{1}%
\e{1}%
\e{0}%
\e{1}%
\e{0}%
\e{0}%
\e{0}%
\e{0}%
\e{0}%
\e{0}%
\e{0}%
\e{0}%
\e{0}%
\eol}\vss}\rg%
%
%
\rx{\vss\hfull{%
\rlx{\hss{$1008_{z}$}}\cg%
\e{0}%
\e{1}%
\e{2}%
\e{3}%
\e{1}%
\e{3}%
\e{2}%
\e{1}%
\e{1}%
\e{1}%
\e{1}%
\e{0}%
\e{0}%
\e{0}%
\e{0}%
\e{0}%
\e{0}%
\e{0}%
\e{0}%
\e{0}%
\e{0}%
\e{0}%
\eol}\vss}\rg%
%
%
\rx{\vss\hfull{%
\rlx{\hss{$560_{z}$}}\cg%
\e{1}%
\e{3}%
\e{3}%
\e{2}%
\e{2}%
\e{2}%
\e{0}%
\e{0}%
\e{1}%
\e{0}%
\e{0}%
\e{0}%
\e{0}%
\e{0}%
\e{0}%
\e{0}%
\e{0}%
\e{0}%
\e{0}%
\e{0}%
\e{0}%
\e{0}%
\eol}\vss}\rg%
%
%
\rx{\vss\hfull{%
\rlx{\hss{$1400_{zz}$}}\cg%
\e{1}%
\e{1}%
\e{2}%
\e{1}%
\e{3}%
\e{2}%
\e{0}%
\e{2}%
\e{3}%
\e{1}%
\e{1}%
\e{0}%
\e{2}%
\e{1}%
\e{0}%
\e{0}%
\e{0}%
\e{0}%
\e{0}%
\e{0}%
\e{0}%
\e{0}%
\eol}\vss}\rg%
%
%
\rx{\vss\hfull{%
\rlx{\hss{$4200_{z}$}}\cg%
\e{0}%
\e{0}%
\e{1}%
\e{2}%
\e{2}%
\e{4}%
\e{2}%
\e{2}%
\e{6}%
\e{3}%
\e{5}%
\e{1}%
\e{2}%
\e{4}%
\e{3}%
\e{2}%
\e{0}%
\e{1}%
\e{0}%
\e{0}%
\e{0}%
\e{0}%
\eol}\vss}\rg%
%
%
\rx{\vss\hfull{%
\rlx{\hss{$400_{z}$}}\cg%
\e{1}%
\e{2}%
\e{1}%
\e{1}%
\e{2}%
\e{1}%
\e{0}%
\e{1}%
\e{1}%
\e{0}%
\e{0}%
\e{0}%
\e{0}%
\e{0}%
\e{0}%
\e{0}%
\e{0}%
\e{0}%
\e{0}%
\e{0}%
\e{0}%
\e{0}%
\eol}\vss}\rg%
%
%
\rx{\vss\hfull{%
\rlx{\hss{$3240_{z}$}}\cg%
\e{0}%
\e{3}%
\e{5}%
\e{4}%
\e{5}%
\e{7}%
\e{2}%
\e{2}%
\e{6}%
\e{3}%
\e{3}%
\e{0}%
\e{1}%
\e{1}%
\e{1}%
\e{0}%
\e{0}%
\e{0}%
\e{0}%
\e{0}%
\e{0}%
\e{0}%
\eol}\vss}\rg%
%
%
\rx{\vss\hfull{%
\rlx{\hss{$4536_{z}$}}\cg%
\e{0}%
\e{1}%
\e{3}%
\e{1}%
\e{5}%
\e{6}%
\e{1}%
\e{2}%
\e{7}%
\e{5}%
\e{4}%
\e{0}%
\e{4}%
\e{2}%
\e{4}%
\e{1}%
\e{0}%
\e{0}%
\e{1}%
\e{0}%
\e{0}%
\e{0}%
\eol}\vss}\rg%
%
%
\rx{\vss\hfull{%
\rlx{\hss{$2400_{z}$}}\cg%
\e{0}%
\e{0}%
\e{0}%
\e{1}%
\e{0}%
\e{2}%
\e{4}%
\e{0}%
\e{2}%
\e{1}%
\e{4}%
\e{3}%
\e{0}%
\e{2}%
\e{1}%
\e{2}%
\e{0}%
\e{0}%
\e{0}%
\e{0}%
\e{0}%
\e{0}%
\eol}\vss}\rg%
%
%
\rx{\vss\hfull{%
\rlx{\hss{$3360_{z}$}}\cg%
\e{0}%
\e{1}%
\e{2}%
\e{2}%
\e{2}%
\e{5}%
\e{2}%
\e{1}%
\e{5}%
\e{3}%
\e{4}%
\e{1}%
\e{2}%
\e{3}%
\e{1}%
\e{1}%
\e{0}%
\e{0}%
\e{0}%
\e{0}%
\e{0}%
\e{0}%
\eol}\vss}\rg%
%
%
\rx{\vss\hfull{%
\rlx{\hss{$2800_{z}$}}\cg%
\e{0}%
\e{0}%
\e{1}%
\e{2}%
\e{1}%
\e{4}%
\e{4}%
\e{1}%
\e{3}%
\e{2}%
\e{4}%
\e{2}%
\e{1}%
\e{2}%
\e{1}%
\e{1}%
\e{0}%
\e{0}%
\e{0}%
\e{0}%
\e{0}%
\e{0}%
\eol}\vss}\rg%
%
%
\rx{\vss\hfull{%
\rlx{\hss{$4096_{z}$}}\cg%
\e{0}%
\e{1}%
\e{3}%
\e{3}%
\e{3}%
\e{8}%
\e{3}%
\e{1}%
\e{5}%
\e{4}%
\e{5}%
\e{1}%
\e{2}%
\e{2}%
\e{2}%
\e{1}%
\e{0}%
\e{0}%
\e{0}%
\e{0}%
\e{0}%
\e{0}%
\eol}\vss}\rg%
%
%
\rx{\vss\hfull{%
\rlx{\hss{$5600_{z}$}}\cg%
\e{0}%
\e{0}%
\e{1}%
\e{2}%
\e{2}%
\e{6}%
\e{4}%
\e{1}%
\e{5}%
\e{5}%
\e{8}%
\e{3}%
\e{3}%
\e{3}%
\e{4}%
\e{3}%
\e{1}%
\e{1}%
\e{1}%
\e{0}%
\e{0}%
\e{0}%
\eol}\vss}\rg%
%
%
\rx{\vss\hfull{%
\rlx{\hss{$448_{z}$}}\cg%
\e{1}%
\e{1}%
\e{2}%
\e{0}%
\e{2}%
\e{1}%
\e{0}%
\e{0}%
\e{1}%
\e{0}%
\e{0}%
\e{0}%
\e{1}%
\e{0}%
\e{0}%
\e{0}%
\e{0}%
\e{0}%
\e{0}%
\e{0}%
\e{0}%
\e{0}%
\eol}\vss}\rg%
%
%
\rx{\vss\hfull{%
\rlx{\hss{$448_{w}$}}\cg%
\e{0}%
\e{0}%
\e{0}%
\e{0}%
\e{0}%
\e{0}%
\e{1}%
\e{0}%
\e{0}%
\e{0}%
\e{1}%
\e{1}%
\e{0}%
\e{0}%
\e{0}%
\e{1}%
\e{0}%
\e{0}%
\e{0}%
\e{0}%
\e{0}%
\e{0}%
\eol}\vss}\rg%
%
%
\rx{\vss\hfull{%
\rlx{\hss{$1344_{w}$}}\cg%
\e{0}%
\e{0}%
\e{0}%
\e{0}%
\e{1}%
\e{1}%
\e{0}%
\e{1}%
\e{2}%
\e{1}%
\e{1}%
\e{0}%
\e{1}%
\e{1}%
\e{2}%
\e{0}%
\e{0}%
\e{1}%
\e{1}%
\e{0}%
\e{0}%
\e{0}%
\eol}\vss}\rg%
%
%
\rx{\vss\hfull{%
\rlx{\hss{$5600_{w}$}}\cg%
\e{0}%
\e{0}%
\e{0}%
\e{1}%
\e{1}%
\e{4}%
\e{3}%
\e{0}%
\e{4}%
\e{4}%
\e{8}%
\e{3}%
\e{2}%
\e{4}%
\e{5}%
\e{5}%
\e{1}%
\e{1}%
\e{1}%
\e{1}%
\e{0}%
\e{0}%
\eol}\vss}\rg%
%
%
\rx{\vss\hfull{%
\rlx{\hss{$2016_{w}$}}\cg%
\e{0}%
\e{0}%
\e{1}%
\e{0}%
\e{1}%
\e{1}%
\e{0}%
\e{1}%
\e{3}%
\e{2}%
\e{1}%
\e{0}%
\e{3}%
\e{2}%
\e{2}%
\e{1}%
\e{0}%
\e{0}%
\e{1}%
\e{0}%
\e{0}%
\e{0}%
\eol}\vss}\rg%
%
%
\rx{\vss\hfull{%
\rlx{\hss{$7168_{w}$}}\cg%
\e{0}%
\e{0}%
\e{1}%
\e{1}%
\e{2}%
\e{5}%
\e{2}%
\e{1}%
\e{7}%
\e{6}%
\e{8}%
\e{2}%
\e{5}%
\e{6}%
\e{7}%
\e{5}%
\e{1}%
\e{1}%
\e{2}%
\e{1}%
\e{0}%
\e{0}%
\eol}\vss}\rg%
\eop
\eject
\tablecont%
%
%
%
%
%
%
\rowpts=18 true pt%
\colpts=18 true pt%
\rowlabpts=40 true pt%
\collabpts=65 true pt%
\clx{\vss\hfull{%
\rlx{\hss{$ $}}\cg%
\cx{\hskip 16 true pt\flip{$[{8}]{\times}[{1^{2}}]$}\hss}\cg%
\cx{\hskip 16 true pt\flip{$[{7}{1}]{\times}[{1^{2}}]$}\hss}\cg%
\cx{\hskip 16 true pt\flip{$[{6}{2}]{\times}[{1^{2}}]$}\hss}\cg%
\cx{\hskip 16 true pt\flip{$[{6}{1^{2}}]{\times}[{1^{2}}]$}\hss}\cg%
\cx{\hskip 16 true pt\flip{$[{5}{3}]{\times}[{1^{2}}]$}\hss}\cg%
\cx{\hskip 16 true pt\flip{$[{5}{2}{1}]{\times}[{1^{2}}]$}\hss}\cg%
\cx{\hskip 16 true pt\flip{$[{5}{1^{3}}]{\times}[{1^{2}}]$}\hss}\cg%
\cx{\hskip 16 true pt\flip{$[{4^{2}}]{\times}[{1^{2}}]$}\hss}\cg%
\cx{\hskip 16 true pt\flip{$[{4}{3}{1}]{\times}[{1^{2}}]$}\hss}\cg%
\cx{\hskip 16 true pt\flip{$[{4}{2^{2}}]{\times}[{1^{2}}]$}\hss}\cg%
\cx{\hskip 16 true pt\flip{$[{4}{2}{1^{2}}]{\times}[{1^{2}}]$}\hss}\cg%
\cx{\hskip 16 true pt\flip{$[{4}{1^{4}}]{\times}[{1^{2}}]$}\hss}\cg%
\cx{\hskip 16 true pt\flip{$[{3^{2}}{2}]{\times}[{1^{2}}]$}\hss}\cg%
\cx{\hskip 16 true pt\flip{$[{3^{2}}{1^{2}}]{\times}[{1^{2}}]$}\hss}\cg%
\cx{\hskip 16 true pt\flip{$[{3}{2^{2}}{1}]{\times}[{1^{2}}]$}\hss}\cg%
\cx{\hskip 16 true pt\flip{$[{3}{2}{1^{3}}]{\times}[{1^{2}}]$}\hss}\cg%
\cx{\hskip 16 true pt\flip{$[{3}{1^{5}}]{\times}[{1^{2}}]$}\hss}\cg%
\cx{\hskip 16 true pt\flip{$[{2^{4}}]{\times}[{1^{2}}]$}\hss}\cg%
\cx{\hskip 16 true pt\flip{$[{2^{3}}{1^{2}}]{\times}[{1^{2}}]$}\hss}\cg%
\cx{\hskip 16 true pt\flip{$[{2^{2}}{1^{4}}]{\times}[{1^{2}}]$}\hss}\cg%
\cx{\hskip 16 true pt\flip{$[{2}{1^{6}}]{\times}[{1^{2}}]$}\hss}\cg%
\cx{\hskip 16 true pt\flip{$[{1^{8}}]{\times}[{1^{2}}]$}\hss}\cg%
\eol}}\rg%
%
%
\rx{\vss\hfull{%
\rlx{\hss{$1_{x}$}}\cg%
\e{0}%
\e{0}%
\e{0}%
\e{0}%
\e{0}%
\e{0}%
\e{0}%
\e{0}%
\e{0}%
\e{0}%
\e{0}%
\e{0}%
\e{0}%
\e{0}%
\e{0}%
\e{0}%
\e{0}%
\e{0}%
\e{0}%
\e{0}%
\e{0}%
\e{0}%
\eol}\vss}\rg%
%
%
\rx{\vss\hfull{%
\rlx{\hss{$28_{x}$}}\cg%
\e{0}%
\e{1}%
\e{0}%
\e{0}%
\e{0}%
\e{0}%
\e{0}%
\e{0}%
\e{0}%
\e{0}%
\e{0}%
\e{0}%
\e{0}%
\e{0}%
\e{0}%
\e{0}%
\e{0}%
\e{0}%
\e{0}%
\e{0}%
\e{0}%
\e{0}%
\eol}\vss}\rg%
%
%
\rx{\vss\hfull{%
\rlx{\hss{$35_{x}$}}\cg%
\e{0}%
\e{1}%
\e{0}%
\e{0}%
\e{0}%
\e{0}%
\e{0}%
\e{0}%
\e{0}%
\e{0}%
\e{0}%
\e{0}%
\e{0}%
\e{0}%
\e{0}%
\e{0}%
\e{0}%
\e{0}%
\e{0}%
\e{0}%
\e{0}%
\e{0}%
\eol}\vss}\rg%
%
%
\rx{\vss\hfull{%
\rlx{\hss{$84_{x}$}}\cg%
\e{1}%
\e{0}%
\e{1}%
\e{0}%
\e{0}%
\e{0}%
\e{0}%
\e{0}%
\e{0}%
\e{0}%
\e{0}%
\e{0}%
\e{0}%
\e{0}%
\e{0}%
\e{0}%
\e{0}%
\e{0}%
\e{0}%
\e{0}%
\e{0}%
\e{0}%
\eol}\vss}\rg%
%
%
\rx{\vss\hfull{%
\rlx{\hss{$50_{x}$}}\cg%
\e{1}%
\e{0}%
\e{0}%
\e{0}%
\e{0}%
\e{0}%
\e{0}%
\e{1}%
\e{0}%
\e{0}%
\e{0}%
\e{0}%
\e{0}%
\e{0}%
\e{0}%
\e{0}%
\e{0}%
\e{0}%
\e{0}%
\e{0}%
\e{0}%
\e{0}%
\eol}\vss}\rg%
%
%
\rx{\vss\hfull{%
\rlx{\hss{$350_{x}$}}\cg%
\e{0}%
\e{0}%
\e{1}%
\e{1}%
\e{0}%
\e{1}%
\e{1}%
\e{0}%
\e{0}%
\e{0}%
\e{0}%
\e{0}%
\e{0}%
\e{0}%
\e{0}%
\e{0}%
\e{0}%
\e{0}%
\e{0}%
\e{0}%
\e{0}%
\e{0}%
\eol}\vss}\rg%
%
%
\rx{\vss\hfull{%
\rlx{\hss{$300_{x}$}}\cg%
\e{0}%
\e{0}%
\e{1}%
\e{1}%
\e{0}%
\e{1}%
\e{0}%
\e{0}%
\e{0}%
\e{0}%
\e{0}%
\e{0}%
\e{0}%
\e{0}%
\e{0}%
\e{0}%
\e{0}%
\e{0}%
\e{0}%
\e{0}%
\e{0}%
\e{0}%
\eol}\vss}\rg%
%
%
\rx{\vss\hfull{%
\rlx{\hss{$567_{x}$}}\cg%
\e{1}%
\e{2}%
\e{2}%
\e{2}%
\e{1}%
\e{1}%
\e{0}%
\e{0}%
\e{0}%
\e{0}%
\e{0}%
\e{0}%
\e{0}%
\e{0}%
\e{0}%
\e{0}%
\e{0}%
\e{0}%
\e{0}%
\e{0}%
\e{0}%
\e{0}%
\eol}\vss}\rg%
%
%
\rx{\vss\hfull{%
\rlx{\hss{$210_{x}$}}\cg%
\e{0}%
\e{2}%
\e{0}%
\e{1}%
\e{1}%
\e{0}%
\e{0}%
\e{0}%
\e{0}%
\e{0}%
\e{0}%
\e{0}%
\e{0}%
\e{0}%
\e{0}%
\e{0}%
\e{0}%
\e{0}%
\e{0}%
\e{0}%
\e{0}%
\e{0}%
\eol}\vss}\rg%
%
%
\rx{\vss\hfull{%
\rlx{\hss{$840_{x}$}}\cg%
\e{0}%
\e{0}%
\e{0}%
\e{0}%
\e{1}%
\e{1}%
\e{0}%
\e{0}%
\e{1}%
\e{1}%
\e{1}%
\e{0}%
\e{0}%
\e{0}%
\e{1}%
\e{0}%
\e{0}%
\e{0}%
\e{0}%
\e{0}%
\e{0}%
\e{0}%
\eol}\vss}\rg%
%
%
\rx{\vss\hfull{%
\rlx{\hss{$700_{x}$}}\cg%
\e{0}%
\e{2}%
\e{1}%
\e{1}%
\e{2}%
\e{1}%
\e{0}%
\e{0}%
\e{1}%
\e{0}%
\e{0}%
\e{0}%
\e{0}%
\e{0}%
\e{0}%
\e{0}%
\e{0}%
\e{0}%
\e{0}%
\e{0}%
\e{0}%
\e{0}%
\eol}\vss}\rg%
%
%
\rx{\vss\hfull{%
\rlx{\hss{$175_{x}$}}\cg%
\e{0}%
\e{0}%
\e{0}%
\e{0}%
\e{1}%
\e{0}%
\e{0}%
\e{0}%
\e{0}%
\e{0}%
\e{0}%
\e{0}%
\e{1}%
\e{0}%
\e{0}%
\e{0}%
\e{0}%
\e{0}%
\e{0}%
\e{0}%
\e{0}%
\e{0}%
\eol}\vss}\rg%
%
%
\rx{\vss\hfull{%
\rlx{\hss{$1400_{x}$}}\cg%
\e{0}%
\e{2}%
\e{1}%
\e{2}%
\e{3}%
\e{2}%
\e{0}%
\e{0}%
\e{2}%
\e{0}%
\e{1}%
\e{0}%
\e{1}%
\e{0}%
\e{0}%
\e{0}%
\e{0}%
\e{0}%
\e{0}%
\e{0}%
\e{0}%
\e{0}%
\eol}\vss}\rg%
%
%
\rx{\vss\hfull{%
\rlx{\hss{$1050_{x}$}}\cg%
\e{1}%
\e{1}%
\e{2}%
\e{0}%
\e{1}%
\e{1}%
\e{0}%
\e{2}%
\e{2}%
\e{1}%
\e{0}%
\e{0}%
\e{0}%
\e{1}%
\e{0}%
\e{0}%
\e{0}%
\e{0}%
\e{0}%
\e{0}%
\e{0}%
\e{0}%
\eol}\vss}\rg%
%
%
\rx{\vss\hfull{%
\rlx{\hss{$1575_{x}$}}\cg%
\e{0}%
\e{2}%
\e{2}%
\e{3}%
\e{2}%
\e{3}%
\e{1}%
\e{0}%
\e{2}%
\e{0}%
\e{1}%
\e{0}%
\e{0}%
\e{0}%
\e{0}%
\e{0}%
\e{0}%
\e{0}%
\e{0}%
\e{0}%
\e{0}%
\e{0}%
\eol}\vss}\rg%
%
%
\rx{\vss\hfull{%
\rlx{\hss{$1344_{x}$}}\cg%
\e{1}%
\e{1}%
\e{4}%
\e{1}%
\e{1}%
\e{3}%
\e{1}%
\e{1}%
\e{1}%
\e{1}%
\e{0}%
\e{0}%
\e{0}%
\e{0}%
\e{0}%
\e{0}%
\e{0}%
\e{0}%
\e{0}%
\e{0}%
\e{0}%
\e{0}%
\eol}\vss}\rg%
%
%
\rx{\vss\hfull{%
\rlx{\hss{$2100_{x}$}}\cg%
\e{0}%
\e{0}%
\e{2}%
\e{2}%
\e{0}%
\e{4}%
\e{4}%
\e{1}%
\e{1}%
\e{2}%
\e{2}%
\e{1}%
\e{0}%
\e{1}%
\e{0}%
\e{0}%
\e{0}%
\e{0}%
\e{0}%
\e{0}%
\e{0}%
\e{0}%
\eol}\vss}\rg%
%
%
\rx{\vss\hfull{%
\rlx{\hss{$2268_{x}$}}\cg%
\e{0}%
\e{1}%
\e{2}%
\e{4}%
\e{2}%
\e{4}%
\e{2}%
\e{1}%
\e{2}%
\e{1}%
\e{2}%
\e{0}%
\e{1}%
\e{0}%
\e{0}%
\e{0}%
\e{0}%
\e{0}%
\e{0}%
\e{0}%
\e{0}%
\e{0}%
\eol}\vss}\rg%
%
%
\rx{\vss\hfull{%
\rlx{\hss{$525_{x}$}}\cg%
\e{1}%
\e{0}%
\e{2}%
\e{0}%
\e{0}%
\e{1}%
\e{1}%
\e{1}%
\e{0}%
\e{1}%
\e{0}%
\e{0}%
\e{0}%
\e{0}%
\e{0}%
\e{0}%
\e{0}%
\e{0}%
\e{0}%
\e{0}%
\e{0}%
\e{0}%
\eol}\vss}\rg%
%
%
\rx{\vss\hfull{%
\rlx{\hss{$700_{xx}$}}\cg%
\e{1}%
\e{0}%
\e{1}%
\e{0}%
\e{0}%
\e{0}%
\e{0}%
\e{3}%
\e{1}%
\e{1}%
\e{0}%
\e{0}%
\e{0}%
\e{2}%
\e{0}%
\e{0}%
\e{0}%
\e{1}%
\e{0}%
\e{0}%
\e{0}%
\e{0}%
\eol}\vss}\rg%
%
%
\rx{\vss\hfull{%
\rlx{\hss{$972_{x}$}}\cg%
\e{0}%
\e{0}%
\e{2}%
\e{0}%
\e{0}%
\e{1}%
\e{1}%
\e{2}%
\e{1}%
\e{2}%
\e{0}%
\e{0}%
\e{0}%
\e{1}%
\e{0}%
\e{0}%
\e{0}%
\e{0}%
\e{0}%
\e{0}%
\e{0}%
\e{0}%
\eol}\vss}\rg%
%
%
\rx{\vss\hfull{%
\rlx{\hss{$4096_{x}$}}\cg%
\e{0}%
\e{1}%
\e{3}%
\e{3}%
\e{3}%
\e{7}%
\e{3}%
\e{1}%
\e{4}%
\e{3}%
\e{4}%
\e{1}%
\e{1}%
\e{1}%
\e{1}%
\e{0}%
\e{0}%
\e{0}%
\e{0}%
\e{0}%
\e{0}%
\e{0}%
\eol}\vss}\rg%
%
%
\rx{\vss\hfull{%
\rlx{\hss{$4200_{x}$}}\cg%
\e{0}%
\e{1}%
\e{2}%
\e{2}%
\e{3}%
\e{5}%
\e{1}%
\e{1}%
\e{6}%
\e{3}%
\e{4}%
\e{0}%
\e{2}%
\e{2}%
\e{2}%
\e{1}%
\e{0}%
\e{0}%
\e{0}%
\e{0}%
\e{0}%
\e{0}%
\eol}\vss}\rg%
%
%
\rx{\vss\hfull{%
\rlx{\hss{$2240_{x}$}}\cg%
\e{0}%
\e{1}%
\e{1}%
\e{1}%
\e{4}%
\e{3}%
\e{0}%
\e{0}%
\e{3}%
\e{1}%
\e{2}%
\e{0}%
\e{2}%
\e{0}%
\e{1}%
\e{0}%
\e{0}%
\e{0}%
\e{0}%
\e{0}%
\e{0}%
\e{0}%
\eol}\vss}\rg%
%
%
\rx{\vss\hfull{%
\rlx{\hss{$2835_{x}$}}\cg%
\e{0}%
\e{0}%
\e{1}%
\e{1}%
\e{2}%
\e{2}%
\e{0}%
\e{1}%
\e{4}%
\e{2}%
\e{2}%
\e{0}%
\e{4}%
\e{2}%
\e{2}%
\e{1}%
\e{0}%
\e{0}%
\e{1}%
\e{0}%
\e{0}%
\e{0}%
\eol}\vss}\rg%
%
%
\rx{\vss\hfull{%
\rlx{\hss{$6075_{x}$}}\cg%
\e{0}%
\e{1}%
\e{3}%
\e{2}%
\e{3}%
\e{7}%
\e{4}%
\e{2}%
\e{7}%
\e{5}%
\e{6}%
\e{2}%
\e{2}%
\e{4}%
\e{3}%
\e{2}%
\e{0}%
\e{0}%
\e{0}%
\e{0}%
\e{0}%
\e{0}%
\eol}\vss}\rg%
%
%
\rx{\vss\hfull{%
\rlx{\hss{$3200_{x}$}}\cg%
\e{0}%
\e{0}%
\e{1}%
\e{0}%
\e{1}%
\e{4}%
\e{2}%
\e{1}%
\e{2}%
\e{5}%
\e{3}%
\e{2}%
\e{1}%
\e{2}%
\e{2}%
\e{1}%
\e{0}%
\e{1}%
\e{0}%
\e{0}%
\e{0}%
\e{0}%
\eol}\vss}\rg%
%
%
\rx{\vss\hfull{%
\rlx{\hss{$70_{y}$}}\cg%
\e{0}%
\e{0}%
\e{0}%
\e{0}%
\e{0}%
\e{0}%
\e{1}%
\e{0}%
\e{0}%
\e{0}%
\e{0}%
\e{0}%
\e{0}%
\e{0}%
\e{0}%
\e{0}%
\e{0}%
\e{0}%
\e{0}%
\e{0}%
\e{0}%
\e{0}%
\eol}\vss}\rg%
%
%
\rx{\vss\hfull{%
\rlx{\hss{$1134_{y}$}}\cg%
\e{0}%
\e{0}%
\e{0}%
\e{0}%
\e{1}%
\e{1}%
\e{0}%
\e{0}%
\e{1}%
\e{1}%
\e{2}%
\e{1}%
\e{0}%
\e{0}%
\e{1}%
\e{1}%
\e{0}%
\e{0}%
\e{0}%
\e{0}%
\e{0}%
\e{0}%
\eol}\vss}\rg%
%
%
\rx{\vss\hfull{%
\rlx{\hss{$1680_{y}$}}\cg%
\e{0}%
\e{0}%
\e{0}%
\e{1}%
\e{0}%
\e{2}%
\e{3}%
\e{0}%
\e{1}%
\e{1}%
\e{3}%
\e{2}%
\e{0}%
\e{1}%
\e{0}%
\e{1}%
\e{0}%
\e{0}%
\e{0}%
\e{0}%
\e{0}%
\e{0}%
\eol}\vss}\rg%
%
%
\rx{\vss\hfull{%
\rlx{\hss{$168_{y}$}}\cg%
\e{0}%
\e{0}%
\e{0}%
\e{0}%
\e{0}%
\e{0}%
\e{0}%
\e{1}%
\e{0}%
\e{0}%
\e{0}%
\e{0}%
\e{0}%
\e{1}%
\e{0}%
\e{0}%
\e{0}%
\e{1}%
\e{0}%
\e{0}%
\e{0}%
\e{0}%
\eol}\vss}\rg%
%
%
\rx{\vss\hfull{%
\rlx{\hss{$420_{y}$}}\cg%
\e{0}%
\e{0}%
\e{0}%
\e{0}%
\e{0}%
\e{0}%
\e{0}%
\e{0}%
\e{1}%
\e{0}%
\e{0}%
\e{0}%
\e{1}%
\e{0}%
\e{1}%
\e{0}%
\e{0}%
\e{0}%
\e{1}%
\e{0}%
\e{0}%
\e{0}%
\eol}\vss}\rg%
%
%
\rx{\vss\hfull{%
\rlx{\hss{$3150_{y}$}}\cg%
\e{0}%
\e{0}%
\e{0}%
\e{0}%
\e{2}%
\e{2}%
\e{0}%
\e{0}%
\e{3}%
\e{2}%
\e{4}%
\e{0}%
\e{4}%
\e{1}%
\e{4}%
\e{2}%
\e{1}%
\e{0}%
\e{2}%
\e{0}%
\e{0}%
\e{0}%
\eol}\vss}\rg%
%
%
\rx{\vss\hfull{%
\rlx{\hss{$4200_{y}$}}\cg%
\e{0}%
\e{0}%
\e{1}%
\e{0}%
\e{1}%
\e{2}%
\e{1}%
\e{2}%
\e{5}%
\e{4}%
\e{3}%
\e{1}%
\e{2}%
\e{5}%
\e{5}%
\e{3}%
\e{0}%
\e{2}%
\e{1}%
\e{1}%
\e{0}%
\e{0}%
\eol}\vss}\rg%
\eop
\eject
\tablecont%
%
%
%
%
%
%
\rowpts=18 true pt%
\colpts=18 true pt%
\rowlabpts=40 true pt%
\collabpts=65 true pt%
\clx{\vss\hfull{%
\rlx{\hss{$ $}}\cg%
\cx{\hskip 16 true pt\flip{$[{8}]{\times}[{1^{2}}]$}\hss}\cg%
\cx{\hskip 16 true pt\flip{$[{7}{1}]{\times}[{1^{2}}]$}\hss}\cg%
\cx{\hskip 16 true pt\flip{$[{6}{2}]{\times}[{1^{2}}]$}\hss}\cg%
\cx{\hskip 16 true pt\flip{$[{6}{1^{2}}]{\times}[{1^{2}}]$}\hss}\cg%
\cx{\hskip 16 true pt\flip{$[{5}{3}]{\times}[{1^{2}}]$}\hss}\cg%
\cx{\hskip 16 true pt\flip{$[{5}{2}{1}]{\times}[{1^{2}}]$}\hss}\cg%
\cx{\hskip 16 true pt\flip{$[{5}{1^{3}}]{\times}[{1^{2}}]$}\hss}\cg%
\cx{\hskip 16 true pt\flip{$[{4^{2}}]{\times}[{1^{2}}]$}\hss}\cg%
\cx{\hskip 16 true pt\flip{$[{4}{3}{1}]{\times}[{1^{2}}]$}\hss}\cg%
\cx{\hskip 16 true pt\flip{$[{4}{2^{2}}]{\times}[{1^{2}}]$}\hss}\cg%
\cx{\hskip 16 true pt\flip{$[{4}{2}{1^{2}}]{\times}[{1^{2}}]$}\hss}\cg%
\cx{\hskip 16 true pt\flip{$[{4}{1^{4}}]{\times}[{1^{2}}]$}\hss}\cg%
\cx{\hskip 16 true pt\flip{$[{3^{2}}{2}]{\times}[{1^{2}}]$}\hss}\cg%
\cx{\hskip 16 true pt\flip{$[{3^{2}}{1^{2}}]{\times}[{1^{2}}]$}\hss}\cg%
\cx{\hskip 16 true pt\flip{$[{3}{2^{2}}{1}]{\times}[{1^{2}}]$}\hss}\cg%
\cx{\hskip 16 true pt\flip{$[{3}{2}{1^{3}}]{\times}[{1^{2}}]$}\hss}\cg%
\cx{\hskip 16 true pt\flip{$[{3}{1^{5}}]{\times}[{1^{2}}]$}\hss}\cg%
\cx{\hskip 16 true pt\flip{$[{2^{4}}]{\times}[{1^{2}}]$}\hss}\cg%
\cx{\hskip 16 true pt\flip{$[{2^{3}}{1^{2}}]{\times}[{1^{2}}]$}\hss}\cg%
\cx{\hskip 16 true pt\flip{$[{2^{2}}{1^{4}}]{\times}[{1^{2}}]$}\hss}\cg%
\cx{\hskip 16 true pt\flip{$[{2}{1^{6}}]{\times}[{1^{2}}]$}\hss}\cg%
\cx{\hskip 16 true pt\flip{$[{1^{8}}]{\times}[{1^{2}}]$}\hss}\cg%
\eol}}\rg%
%
%
\rx{\vss\hfull{%
\rlx{\hss{$2688_{y}$}}\cg%
\e{0}%
\e{0}%
\e{1}%
\e{0}%
\e{0}%
\e{2}%
\e{1}%
\e{1}%
\e{2}%
\e{4}%
\e{2}%
\e{1}%
\e{1}%
\e{4}%
\e{2}%
\e{2}%
\e{0}%
\e{1}%
\e{0}%
\e{1}%
\e{0}%
\e{0}%
\eol}\vss}\rg%
%
%
\rx{\vss\hfull{%
\rlx{\hss{$2100_{y}$}}\cg%
\e{0}%
\e{0}%
\e{1}%
\e{0}%
\e{0}%
\e{2}%
\e{3}%
\e{1}%
\e{1}%
\e{3}%
\e{2}%
\e{3}%
\e{0}%
\e{2}%
\e{1}%
\e{1}%
\e{0}%
\e{1}%
\e{0}%
\e{0}%
\e{0}%
\e{0}%
\eol}\vss}\rg%
%
%
\rx{\vss\hfull{%
\rlx{\hss{$1400_{y}$}}\cg%
\e{0}%
\e{0}%
\e{0}%
\e{1}%
\e{0}%
\e{1}%
\e{2}%
\e{0}%
\e{2}%
\e{0}%
\e{2}%
\e{1}%
\e{0}%
\e{1}%
\e{1}%
\e{1}%
\e{0}%
\e{0}%
\e{0}%
\e{0}%
\e{0}%
\e{0}%
\eol}\vss}\rg%
%
%
\rx{\vss\hfull{%
\rlx{\hss{$4536_{y}$}}\cg%
\e{0}%
\e{0}%
\e{0}%
\e{2}%
\e{1}%
\e{4}%
\e{3}%
\e{0}%
\e{4}%
\e{2}%
\e{7}%
\e{2}%
\e{3}%
\e{3}%
\e{3}%
\e{3}%
\e{1}%
\e{0}%
\e{1}%
\e{0}%
\e{0}%
\e{0}%
\eol}\vss}\rg%
%
%
\rx{\vss\hfull{%
\rlx{\hss{$5670_{y}$}}\cg%
\e{0}%
\e{0}%
\e{0}%
\e{2}%
\e{2}%
\e{5}%
\e{3}%
\e{0}%
\e{5}%
\e{3}%
\e{9}%
\e{3}%
\e{3}%
\e{3}%
\e{4}%
\e{4}%
\e{1}%
\e{0}%
\e{1}%
\e{0}%
\e{0}%
\e{0}%
\eol}\vss}\rg%
%
%
\rx{\vss\hfull{%
\rlx{\hss{$4480_{y}$}}\cg%
\e{0}%
\e{0}%
\e{0}%
\e{1}%
\e{2}%
\e{3}%
\e{1}%
\e{0}%
\e{5}%
\e{2}%
\e{6}%
\e{1}%
\e{4}%
\e{2}%
\e{5}%
\e{3}%
\e{1}%
\e{0}%
\e{2}%
\e{0}%
\e{0}%
\e{0}%
\eol}\vss}\rg%
%
%
\rx{\vss\hfull{%
\rlx{\hss{$8_{z}$}}\cg%
\e{1}%
\e{0}%
\e{0}%
\e{0}%
\e{0}%
\e{0}%
\e{0}%
\e{0}%
\e{0}%
\e{0}%
\e{0}%
\e{0}%
\e{0}%
\e{0}%
\e{0}%
\e{0}%
\e{0}%
\e{0}%
\e{0}%
\e{0}%
\e{0}%
\e{0}%
\eol}\vss}\rg%
%
%
\rx{\vss\hfull{%
\rlx{\hss{$56_{z}$}}\cg%
\e{0}%
\e{0}%
\e{0}%
\e{1}%
\e{0}%
\e{0}%
\e{0}%
\e{0}%
\e{0}%
\e{0}%
\e{0}%
\e{0}%
\e{0}%
\e{0}%
\e{0}%
\e{0}%
\e{0}%
\e{0}%
\e{0}%
\e{0}%
\e{0}%
\e{0}%
\eol}\vss}\rg%
%
%
\rx{\vss\hfull{%
\rlx{\hss{$160_{z}$}}\cg%
\e{0}%
\e{1}%
\e{1}%
\e{1}%
\e{0}%
\e{0}%
\e{0}%
\e{0}%
\e{0}%
\e{0}%
\e{0}%
\e{0}%
\e{0}%
\e{0}%
\e{0}%
\e{0}%
\e{0}%
\e{0}%
\e{0}%
\e{0}%
\e{0}%
\e{0}%
\eol}\vss}\rg%
%
%
\rx{\vss\hfull{%
\rlx{\hss{$112_{z}$}}\cg%
\e{1}%
\e{1}%
\e{1}%
\e{0}%
\e{0}%
\e{0}%
\e{0}%
\e{0}%
\e{0}%
\e{0}%
\e{0}%
\e{0}%
\e{0}%
\e{0}%
\e{0}%
\e{0}%
\e{0}%
\e{0}%
\e{0}%
\e{0}%
\e{0}%
\e{0}%
\eol}\vss}\rg%
%
%
\rx{\vss\hfull{%
\rlx{\hss{$840_{z}$}}\cg%
\e{0}%
\e{0}%
\e{1}%
\e{0}%
\e{1}%
\e{2}%
\e{1}%
\e{0}%
\e{0}%
\e{1}%
\e{1}%
\e{0}%
\e{0}%
\e{0}%
\e{0}%
\e{0}%
\e{0}%
\e{0}%
\e{0}%
\e{0}%
\e{0}%
\e{0}%
\eol}\vss}\rg%
%
%
\rx{\vss\hfull{%
\rlx{\hss{$1296_{z}$}}\cg%
\e{0}%
\e{1}%
\e{1}%
\e{3}%
\e{1}%
\e{3}%
\e{2}%
\e{0}%
\e{1}%
\e{0}%
\e{1}%
\e{0}%
\e{0}%
\e{0}%
\e{0}%
\e{0}%
\e{0}%
\e{0}%
\e{0}%
\e{0}%
\e{0}%
\e{0}%
\eol}\vss}\rg%
%
%
\rx{\vss\hfull{%
\rlx{\hss{$1400_{z}$}}\cg%
\e{0}%
\e{2}%
\e{3}%
\e{2}%
\e{2}%
\e{3}%
\e{1}%
\e{0}%
\e{1}%
\e{1}%
\e{0}%
\e{0}%
\e{0}%
\e{0}%
\e{0}%
\e{0}%
\e{0}%
\e{0}%
\e{0}%
\e{0}%
\e{0}%
\e{0}%
\eol}\vss}\rg%
%
%
\rx{\vss\hfull{%
\rlx{\hss{$1008_{z}$}}\cg%
\e{0}%
\e{2}%
\e{2}%
\e{3}%
\e{1}%
\e{2}%
\e{1}%
\e{0}%
\e{1}%
\e{0}%
\e{0}%
\e{0}%
\e{0}%
\e{0}%
\e{0}%
\e{0}%
\e{0}%
\e{0}%
\e{0}%
\e{0}%
\e{0}%
\e{0}%
\eol}\vss}\rg%
%
%
\rx{\vss\hfull{%
\rlx{\hss{$560_{z}$}}\cg%
\e{1}%
\e{2}%
\e{2}%
\e{1}%
\e{1}%
\e{1}%
\e{0}%
\e{1}%
\e{0}%
\e{0}%
\e{0}%
\e{0}%
\e{0}%
\e{0}%
\e{0}%
\e{0}%
\e{0}%
\e{0}%
\e{0}%
\e{0}%
\e{0}%
\e{0}%
\eol}\vss}\rg%
%
%
\rx{\vss\hfull{%
\rlx{\hss{$1400_{zz}$}}\cg%
\e{0}%
\e{1}%
\e{1}%
\e{0}%
\e{2}%
\e{1}%
\e{0}%
\e{1}%
\e{3}%
\e{1}%
\e{0}%
\e{0}%
\e{1}%
\e{1}%
\e{1}%
\e{0}%
\e{0}%
\e{0}%
\e{0}%
\e{0}%
\e{0}%
\e{0}%
\eol}\vss}\rg%
%
%
\rx{\vss\hfull{%
\rlx{\hss{$4200_{z}$}}\cg%
\e{0}%
\e{1}%
\e{1}%
\e{1}%
\e{3}%
\e{3}%
\e{1}%
\e{2}%
\e{6}%
\e{2}%
\e{4}%
\e{0}%
\e{3}%
\e{4}%
\e{3}%
\e{2}%
\e{0}%
\e{0}%
\e{1}%
\e{0}%
\e{0}%
\e{0}%
\eol}\vss}\rg%
%
%
\rx{\vss\hfull{%
\rlx{\hss{$400_{z}$}}\cg%
\e{1}%
\e{1}%
\e{1}%
\e{0}%
\e{1}%
\e{0}%
\e{0}%
\e{1}%
\e{1}%
\e{0}%
\e{0}%
\e{0}%
\e{0}%
\e{0}%
\e{0}%
\e{0}%
\e{0}%
\e{0}%
\e{0}%
\e{0}%
\e{0}%
\e{0}%
\eol}\vss}\rg%
%
%
\rx{\vss\hfull{%
\rlx{\hss{$3240_{z}$}}\cg%
\e{1}%
\e{2}%
\e{4}%
\e{3}%
\e{4}%
\e{5}%
\e{1}%
\e{2}%
\e{4}%
\e{2}%
\e{2}%
\e{0}%
\e{1}%
\e{1}%
\e{0}%
\e{0}%
\e{0}%
\e{0}%
\e{0}%
\e{0}%
\e{0}%
\e{0}%
\eol}\vss}\rg%
%
%
\rx{\vss\hfull{%
\rlx{\hss{$4536_{z}$}}\cg%
\e{0}%
\e{0}%
\e{2}%
\e{1}%
\e{3}%
\e{5}%
\e{1}%
\e{2}%
\e{5}%
\e{5}%
\e{4}%
\e{1}%
\e{3}%
\e{2}%
\e{3}%
\e{1}%
\e{0}%
\e{1}%
\e{0}%
\e{0}%
\e{0}%
\e{0}%
\eol}\vss}\rg%
%
%
\rx{\vss\hfull{%
\rlx{\hss{$2400_{z}$}}\cg%
\e{0}%
\e{0}%
\e{1}%
\e{2}%
\e{1}%
\e{3}%
\e{3}%
\e{0}%
\e{2}%
\e{1}%
\e{4}%
\e{1}%
\e{1}%
\e{1}%
\e{0}%
\e{1}%
\e{0}%
\e{0}%
\e{0}%
\e{0}%
\e{0}%
\e{0}%
\eol}\vss}\rg%
%
%
\rx{\vss\hfull{%
\rlx{\hss{$3360_{z}$}}\cg%
\e{0}%
\e{1}%
\e{2}%
\e{2}%
\e{2}%
\e{4}%
\e{1}%
\e{1}%
\e{5}%
\e{2}%
\e{3}%
\e{0}%
\e{2}%
\e{2}%
\e{1}%
\e{1}%
\e{0}%
\e{0}%
\e{0}%
\e{0}%
\e{0}%
\e{0}%
\eol}\vss}\rg%
%
%
\rx{\vss\hfull{%
\rlx{\hss{$2800_{z}$}}\cg%
\e{0}%
\e{1}%
\e{1}%
\e{3}%
\e{2}%
\e{4}%
\e{3}%
\e{1}%
\e{3}%
\e{1}%
\e{3}%
\e{1}%
\e{1}%
\e{1}%
\e{1}%
\e{0}%
\e{0}%
\e{0}%
\e{0}%
\e{0}%
\e{0}%
\e{0}%
\eol}\vss}\rg%
%
%
\rx{\vss\hfull{%
\rlx{\hss{$4096_{z}$}}\cg%
\e{0}%
\e{1}%
\e{3}%
\e{3}%
\e{3}%
\e{7}%
\e{3}%
\e{1}%
\e{4}%
\e{3}%
\e{4}%
\e{1}%
\e{1}%
\e{1}%
\e{1}%
\e{0}%
\e{0}%
\e{0}%
\e{0}%
\e{0}%
\e{0}%
\e{0}%
\eol}\vss}\rg%
%
%
\rx{\vss\hfull{%
\rlx{\hss{$5600_{z}$}}\cg%
\e{0}%
\e{0}%
\e{2}%
\e{2}%
\e{2}%
\e{7}%
\e{5}%
\e{0}%
\e{5}%
\e{5}%
\e{7}%
\e{3}%
\e{2}%
\e{2}%
\e{3}%
\e{2}%
\e{0}%
\e{0}%
\e{0}%
\e{0}%
\e{0}%
\e{0}%
\eol}\vss}\rg%
%
%
\rx{\vss\hfull{%
\rlx{\hss{$448_{z}$}}\cg%
\e{0}%
\e{0}%
\e{1}%
\e{0}%
\e{1}%
\e{1}%
\e{0}%
\e{0}%
\e{0}%
\e{1}%
\e{0}%
\e{0}%
\e{0}%
\e{0}%
\e{0}%
\e{0}%
\e{0}%
\e{0}%
\e{0}%
\e{0}%
\e{0}%
\e{0}%
\eol}\vss}\rg%
%
%
\rx{\vss\hfull{%
\rlx{\hss{$448_{w}$}}\cg%
\e{0}%
\e{0}%
\e{0}%
\e{0}%
\e{0}%
\e{1}%
\e{1}%
\e{0}%
\e{0}%
\e{0}%
\e{1}%
\e{1}%
\e{0}%
\e{0}%
\e{0}%
\e{0}%
\e{0}%
\e{0}%
\e{0}%
\e{0}%
\e{0}%
\e{0}%
\eol}\vss}\rg%
%
%
\rx{\vss\hfull{%
\rlx{\hss{$1344_{w}$}}\cg%
\e{0}%
\e{0}%
\e{0}%
\e{0}%
\e{1}%
\e{0}%
\e{0}%
\e{1}%
\e{2}%
\e{1}%
\e{1}%
\e{0}%
\e{1}%
\e{1}%
\e{2}%
\e{1}%
\e{0}%
\e{1}%
\e{1}%
\e{0}%
\e{0}%
\e{0}%
\eol}\vss}\rg%
%
%
\rx{\vss\hfull{%
\rlx{\hss{$5600_{w}$}}\cg%
\e{0}%
\e{0}%
\e{1}%
\e{1}%
\e{1}%
\e{5}%
\e{3}%
\e{1}%
\e{5}%
\e{4}%
\e{8}%
\e{3}%
\e{2}%
\e{4}%
\e{4}%
\e{4}%
\e{1}%
\e{0}%
\e{1}%
\e{0}%
\e{0}%
\e{0}%
\eol}\vss}\rg%
%
%
\rx{\vss\hfull{%
\rlx{\hss{$2016_{w}$}}\cg%
\e{0}%
\e{0}%
\e{0}%
\e{0}%
\e{1}%
\e{1}%
\e{0}%
\e{0}%
\e{2}%
\e{2}%
\e{1}%
\e{0}%
\e{3}%
\e{2}%
\e{3}%
\e{1}%
\e{0}%
\e{1}%
\e{1}%
\e{1}%
\e{0}%
\e{0}%
\eol}\vss}\rg%
%
%
\rx{\vss\hfull{%
\rlx{\hss{$7168_{w}$}}\cg%
\e{0}%
\e{0}%
\e{1}%
\e{1}%
\e{2}%
\e{5}%
\e{2}%
\e{1}%
\e{7}%
\e{6}%
\e{8}%
\e{2}%
\e{5}%
\e{6}%
\e{7}%
\e{5}%
\e{1}%
\e{1}%
\e{2}%
\e{1}%
\e{0}%
\e{0}%
\eol}\vss}\rg%
\tableclose%
%
%
%
%
%
%
\eop
\eject
\tableopen{Induce/restrict matrix for $W({A_{5}}{A_{2}}{A_{1}})\,\subset\,W(E_{8})$}%
%
%
%
%
%
%
\rowpts=18 true pt%
\colpts=18 true pt%
\rowlabpts=40 true pt%
\collabpts=90 true pt%
\clx{\vss\hfull{%
\rlx{\hss{$ $}}\cg%
\cx{\hskip 16 true pt\flip{$[{6}]{\times}[{3}]{\times}[{2}]$}\hss}\cg%
\cx{\hskip 16 true pt\flip{$[{5}{1}]{\times}[{3}]{\times}[{2}]$}\hss}\cg%
\cx{\hskip 16 true pt\flip{$[{4}{2}]{\times}[{3}]{\times}[{2}]$}\hss}\cg%
\cx{\hskip 16 true pt\flip{$[{4}{1^{2}}]{\times}[{3}]{\times}[{2}]$}\hss}\cg%
\cx{\hskip 16 true pt\flip{$[{3^{2}}]{\times}[{3}]{\times}[{2}]$}\hss}\cg%
\cx{\hskip 16 true pt\flip{$[{3}{2}{1}]{\times}[{3}]{\times}[{2}]$}\hss}\cg%
\cx{\hskip 16 true pt\flip{$[{3}{1^{3}}]{\times}[{3}]{\times}[{2}]$}\hss}\cg%
\cx{\hskip 16 true pt\flip{$[{2^{3}}]{\times}[{3}]{\times}[{2}]$}\hss}\cg%
\cx{\hskip 16 true pt\flip{$[{2^{2}}{1^{2}}]{\times}[{3}]{\times}[{2}]$}\hss}\cg%
\cx{\hskip 16 true pt\flip{$[{2}{1^{4}}]{\times}[{3}]{\times}[{2}]$}\hss}\cg%
\cx{\hskip 16 true pt\flip{$[{1^{6}}]{\times}[{3}]{\times}[{2}]$}\hss}\cg%
\cx{\hskip 16 true pt\flip{$[{6}]{\times}[{2}{1}]{\times}[{2}]$}\hss}\cg%
\cx{\hskip 16 true pt\flip{$[{5}{1}]{\times}[{2}{1}]{\times}[{2}]$}\hss}\cg%
\cx{\hskip 16 true pt\flip{$[{4}{2}]{\times}[{2}{1}]{\times}[{2}]$}\hss}\cg%
\cx{\hskip 16 true pt\flip{$[{4}{1^{2}}]{\times}[{2}{1}]{\times}[{2}]$}\hss}\cg%
\cx{\hskip 16 true pt\flip{$[{3^{2}}]{\times}[{2}{1}]{\times}[{2}]$}\hss}\cg%
\cx{\hskip 16 true pt\flip{$[{3}{2}{1}]{\times}[{2}{1}]{\times}[{2}]$}\hss}\cg%
\cx{\hskip 16 true pt\flip{$[{3}{1^{3}}]{\times}[{2}{1}]{\times}[{2}]$}\hss}\cg%
\cx{\hskip 16 true pt\flip{$[{2^{3}}]{\times}[{2}{1}]{\times}[{2}]$}\hss}\cg%
\cx{\hskip 16 true pt\flip{$[{2^{2}}{1^{2}}]{\times}[{2}{1}]{\times}[{2}]$}\hss}\cg%
\cx{\hskip 16 true pt\flip{$[{2}{1^{4}}]{\times}[{2}{1}]{\times}[{2}]$}\hss}\cg%
\cx{\hskip 16 true pt\flip{$[{1^{6}}]{\times}[{2}{1}]{\times}[{2}]$}\hss}\cg%
\eol}}\rg%
%
%
\rx{\vss\hfull{%
\rlx{\hss{$1_{x}$}}\cg%
\e{1}%
\e{0}%
\e{0}%
\e{0}%
\e{0}%
\e{0}%
\e{0}%
\e{0}%
\e{0}%
\e{0}%
\e{0}%
\e{0}%
\e{0}%
\e{0}%
\e{0}%
\e{0}%
\e{0}%
\e{0}%
\e{0}%
\e{0}%
\e{0}%
\e{0}%
\eol}\vss}\rg%
%
%
\rx{\vss\hfull{%
\rlx{\hss{$28_{x}$}}\cg%
\e{0}%
\e{0}%
\e{0}%
\e{1}%
\e{0}%
\e{0}%
\e{0}%
\e{0}%
\e{0}%
\e{0}%
\e{0}%
\e{0}%
\e{1}%
\e{0}%
\e{0}%
\e{0}%
\e{0}%
\e{0}%
\e{0}%
\e{0}%
\e{0}%
\e{0}%
\eol}\vss}\rg%
%
%
\rx{\vss\hfull{%
\rlx{\hss{$35_{x}$}}\cg%
\e{2}%
\e{1}%
\e{1}%
\e{0}%
\e{0}%
\e{0}%
\e{0}%
\e{0}%
\e{0}%
\e{0}%
\e{0}%
\e{1}%
\e{1}%
\e{0}%
\e{0}%
\e{0}%
\e{0}%
\e{0}%
\e{0}%
\e{0}%
\e{0}%
\e{0}%
\eol}\vss}\rg%
%
%
\rx{\vss\hfull{%
\rlx{\hss{$84_{x}$}}\cg%
\e{3}%
\e{2}%
\e{2}%
\e{0}%
\e{0}%
\e{0}%
\e{0}%
\e{0}%
\e{0}%
\e{0}%
\e{0}%
\e{2}%
\e{1}%
\e{1}%
\e{0}%
\e{0}%
\e{0}%
\e{0}%
\e{0}%
\e{0}%
\e{0}%
\e{0}%
\eol}\vss}\rg%
%
%
\rx{\vss\hfull{%
\rlx{\hss{$50_{x}$}}\cg%
\e{1}%
\e{1}%
\e{1}%
\e{0}%
\e{0}%
\e{0}%
\e{0}%
\e{0}%
\e{0}%
\e{0}%
\e{0}%
\e{1}%
\e{0}%
\e{1}%
\e{0}%
\e{0}%
\e{0}%
\e{0}%
\e{0}%
\e{0}%
\e{0}%
\e{0}%
\eol}\vss}\rg%
%
%
\rx{\vss\hfull{%
\rlx{\hss{$350_{x}$}}\cg%
\e{0}%
\e{0}%
\e{0}%
\e{2}%
\e{0}%
\e{1}%
\e{1}%
\e{0}%
\e{1}%
\e{0}%
\e{0}%
\e{0}%
\e{2}%
\e{1}%
\e{2}%
\e{0}%
\e{1}%
\e{1}%
\e{0}%
\e{0}%
\e{0}%
\e{0}%
\eol}\vss}\rg%
%
%
\rx{\vss\hfull{%
\rlx{\hss{$300_{x}$}}\cg%
\e{2}%
\e{2}%
\e{3}%
\e{0}%
\e{1}%
\e{1}%
\e{0}%
\e{1}%
\e{0}%
\e{0}%
\e{0}%
\e{1}%
\e{3}%
\e{2}%
\e{1}%
\e{0}%
\e{1}%
\e{0}%
\e{0}%
\e{0}%
\e{0}%
\e{0}%
\eol}\vss}\rg%
%
%
\rx{\vss\hfull{%
\rlx{\hss{$567_{x}$}}\cg%
\e{1}%
\e{4}%
\e{3}%
\e{4}%
\e{1}%
\e{2}%
\e{1}%
\e{0}%
\e{0}%
\e{0}%
\e{0}%
\e{3}%
\e{6}%
\e{3}%
\e{3}%
\e{1}%
\e{1}%
\e{0}%
\e{0}%
\e{0}%
\e{0}%
\e{0}%
\eol}\vss}\rg%
%
%
\rx{\vss\hfull{%
\rlx{\hss{$210_{x}$}}\cg%
\e{1}%
\e{3}%
\e{2}%
\e{1}%
\e{0}%
\e{1}%
\e{0}%
\e{0}%
\e{0}%
\e{0}%
\e{0}%
\e{2}%
\e{3}%
\e{1}%
\e{1}%
\e{1}%
\e{0}%
\e{0}%
\e{0}%
\e{0}%
\e{0}%
\e{0}%
\eol}\vss}\rg%
%
%
\rx{\vss\hfull{%
\rlx{\hss{$840_{x}$}}\cg%
\e{1}%
\e{2}%
\e{4}%
\e{0}%
\e{2}%
\e{2}%
\e{0}%
\e{2}%
\e{0}%
\e{0}%
\e{0}%
\e{0}%
\e{2}%
\e{4}%
\e{2}%
\e{2}%
\e{4}%
\e{0}%
\e{2}%
\e{1}%
\e{0}%
\e{0}%
\eol}\vss}\rg%
%
%
\rx{\vss\hfull{%
\rlx{\hss{$700_{x}$}}\cg%
\e{4}%
\e{6}%
\e{6}%
\e{2}%
\e{2}%
\e{2}%
\e{0}%
\e{1}%
\e{0}%
\e{0}%
\e{0}%
\e{3}%
\e{6}%
\e{5}%
\e{2}%
\e{2}%
\e{2}%
\e{0}%
\e{0}%
\e{0}%
\e{0}%
\e{0}%
\eol}\vss}\rg%
%
%
\rx{\vss\hfull{%
\rlx{\hss{$175_{x}$}}\cg%
\e{1}%
\e{1}%
\e{1}%
\e{1}%
\e{1}%
\e{0}%
\e{0}%
\e{0}%
\e{0}%
\e{0}%
\e{0}%
\e{0}%
\e{1}%
\e{1}%
\e{0}%
\e{1}%
\e{1}%
\e{0}%
\e{0}%
\e{0}%
\e{0}%
\e{0}%
\eol}\vss}\rg%
%
%
\rx{\vss\hfull{%
\rlx{\hss{$1400_{x}$}}\cg%
\e{1}%
\e{5}%
\e{4}%
\e{6}%
\e{3}%
\e{4}%
\e{2}%
\e{0}%
\e{1}%
\e{0}%
\e{0}%
\e{2}%
\e{7}%
\e{6}%
\e{6}%
\e{4}%
\e{5}%
\e{1}%
\e{0}%
\e{1}%
\e{0}%
\e{0}%
\eol}\vss}\rg%
%
%
\rx{\vss\hfull{%
\rlx{\hss{$1050_{x}$}}\cg%
\e{2}%
\e{4}%
\e{5}%
\e{4}%
\e{2}%
\e{3}%
\e{1}%
\e{0}%
\e{0}%
\e{0}%
\e{0}%
\e{3}%
\e{4}%
\e{7}%
\e{3}%
\e{2}%
\e{4}%
\e{1}%
\e{1}%
\e{0}%
\e{0}%
\e{0}%
\eol}\vss}\rg%
%
%
\rx{\vss\hfull{%
\rlx{\hss{$1575_{x}$}}\cg%
\e{0}%
\e{4}%
\e{3}%
\e{7}%
\e{2}%
\e{5}%
\e{3}%
\e{0}%
\e{2}%
\e{0}%
\e{0}%
\e{2}%
\e{8}%
\e{6}%
\e{8}%
\e{3}%
\e{5}%
\e{2}%
\e{0}%
\e{1}%
\e{0}%
\e{0}%
\eol}\vss}\rg%
%
%
\rx{\vss\hfull{%
\rlx{\hss{$1344_{x}$}}\cg%
\e{3}%
\e{7}%
\e{8}%
\e{4}%
\e{2}%
\e{5}%
\e{1}%
\e{1}%
\e{0}%
\e{0}%
\e{0}%
\e{5}%
\e{8}%
\e{9}%
\e{5}%
\e{2}%
\e{4}%
\e{1}%
\e{1}%
\e{0}%
\e{0}%
\e{0}%
\eol}\vss}\rg%
%
%
\rx{\vss\hfull{%
\rlx{\hss{$2100_{x}$}}\cg%
\e{0}%
\e{1}%
\e{2}%
\e{4}%
\e{1}%
\e{5}%
\e{5}%
\e{1}%
\e{3}%
\e{2}%
\e{0}%
\e{1}%
\e{4}%
\e{6}%
\e{9}%
\e{1}%
\e{7}%
\e{6}%
\e{2}%
\e{2}%
\e{1}%
\e{0}%
\eol}\vss}\rg%
%
%
\rx{\vss\hfull{%
\rlx{\hss{$2268_{x}$}}\cg%
\e{1}%
\e{5}%
\e{6}%
\e{6}%
\e{3}%
\e{7}%
\e{3}%
\e{2}%
\e{2}%
\e{1}%
\e{0}%
\e{2}%
\e{9}%
\e{9}%
\e{10}%
\e{4}%
\e{8}%
\e{3}%
\e{1}%
\e{2}%
\e{0}%
\e{0}%
\eol}\vss}\rg%
%
%
\rx{\vss\hfull{%
\rlx{\hss{$525_{x}$}}\cg%
\e{0}%
\e{1}%
\e{1}%
\e{3}%
\e{1}%
\e{1}%
\e{2}%
\e{0}%
\e{0}%
\e{0}%
\e{0}%
\e{1}%
\e{2}%
\e{3}%
\e{2}%
\e{0}%
\e{2}%
\e{1}%
\e{0}%
\e{0}%
\e{0}%
\e{0}%
\eol}\vss}\rg%
%
%
\rx{\vss\hfull{%
\rlx{\hss{$700_{xx}$}}\cg%
\e{0}%
\e{1}%
\e{2}%
\e{2}%
\e{0}%
\e{2}%
\e{1}%
\e{0}%
\e{0}%
\e{0}%
\e{0}%
\e{1}%
\e{1}%
\e{4}%
\e{1}%
\e{1}%
\e{3}%
\e{2}%
\e{1}%
\e{0}%
\e{0}%
\e{0}%
\eol}\vss}\rg%
%
%
\rx{\vss\hfull{%
\rlx{\hss{$972_{x}$}}\cg%
\e{2}%
\e{3}%
\e{6}%
\e{1}%
\e{1}%
\e{3}%
\e{0}%
\e{2}%
\e{0}%
\e{0}%
\e{0}%
\e{2}%
\e{3}%
\e{7}%
\e{2}%
\e{1}%
\e{4}%
\e{1}%
\e{2}%
\e{0}%
\e{0}%
\e{0}%
\eol}\vss}\rg%
%
%
\rx{\vss\hfull{%
\rlx{\hss{$4096_{x}$}}\cg%
\e{1}%
\e{6}%
\e{9}%
\e{9}%
\e{4}%
\e{12}%
\e{5}%
\e{3}%
\e{4}%
\e{1}%
\e{0}%
\e{3}%
\e{12}%
\e{15}%
\e{15}%
\e{6}%
\e{16}%
\e{7}%
\e{3}%
\e{4}%
\e{1}%
\e{0}%
\eol}\vss}\rg%
%
%
\rx{\vss\hfull{%
\rlx{\hss{$4200_{x}$}}\cg%
\e{2}%
\e{6}%
\e{10}%
\e{8}%
\e{6}%
\e{11}%
\e{4}%
\e{3}%
\e{3}%
\e{1}%
\e{0}%
\e{2}%
\e{10}%
\e{15}%
\e{13}%
\e{8}%
\e{18}%
\e{6}%
\e{4}%
\e{5}%
\e{1}%
\e{0}%
\eol}\vss}\rg%
%
%
\rx{\vss\hfull{%
\rlx{\hss{$2240_{x}$}}\cg%
\e{2}%
\e{7}%
\e{8}%
\e{5}%
\e{4}%
\e{7}%
\e{1}%
\e{1}%
\e{1}%
\e{0}%
\e{0}%
\e{2}%
\e{8}%
\e{10}%
\e{7}%
\e{7}%
\e{9}%
\e{2}%
\e{2}%
\e{2}%
\e{0}%
\e{0}%
\eol}\vss}\rg%
%
%
\rx{\vss\hfull{%
\rlx{\hss{$2835_{x}$}}\cg%
\e{1}%
\e{4}%
\e{6}%
\e{5}%
\e{4}%
\e{7}%
\e{2}%
\e{1}%
\e{2}%
\e{0}%
\e{0}%
\e{1}%
\e{5}%
\e{9}%
\e{7}%
\e{7}%
\e{13}%
\e{4}%
\e{3}%
\e{4}%
\e{1}%
\e{0}%
\eol}\vss}\rg%
%
%
\rx{\vss\hfull{%
\rlx{\hss{$6075_{x}$}}\cg%
\e{0}%
\e{4}%
\e{7}%
\e{12}%
\e{4}%
\e{14}%
\e{10}%
\e{2}%
\e{7}%
\e{2}%
\e{0}%
\e{2}%
\e{10}%
\e{17}%
\e{19}%
\e{8}%
\e{24}%
\e{14}%
\e{5}%
\e{9}%
\e{3}%
\e{0}%
\eol}\vss}\rg%
%
%
\rx{\vss\hfull{%
\rlx{\hss{$3200_{x}$}}\cg%
\e{1}%
\e{2}%
\e{7}%
\e{3}%
\e{2}%
\e{8}%
\e{2}%
\e{5}%
\e{2}%
\e{1}%
\e{0}%
\e{1}%
\e{5}%
\e{11}%
\e{7}%
\e{3}%
\e{14}%
\e{6}%
\e{6}%
\e{4}%
\e{2}%
\e{0}%
\eol}\vss}\rg%
\eop
\eject
\tablecont%
%
%
%
%
%
%
\rowpts=18 true pt%
\colpts=18 true pt%
\rowlabpts=40 true pt%
\collabpts=90 true pt%
\clx{\vss\hfull{%
\rlx{\hss{$ $}}\cg%
\cx{\hskip 16 true pt\flip{$[{6}]{\times}[{3}]{\times}[{2}]$}\hss}\cg%
\cx{\hskip 16 true pt\flip{$[{5}{1}]{\times}[{3}]{\times}[{2}]$}\hss}\cg%
\cx{\hskip 16 true pt\flip{$[{4}{2}]{\times}[{3}]{\times}[{2}]$}\hss}\cg%
\cx{\hskip 16 true pt\flip{$[{4}{1^{2}}]{\times}[{3}]{\times}[{2}]$}\hss}\cg%
\cx{\hskip 16 true pt\flip{$[{3^{2}}]{\times}[{3}]{\times}[{2}]$}\hss}\cg%
\cx{\hskip 16 true pt\flip{$[{3}{2}{1}]{\times}[{3}]{\times}[{2}]$}\hss}\cg%
\cx{\hskip 16 true pt\flip{$[{3}{1^{3}}]{\times}[{3}]{\times}[{2}]$}\hss}\cg%
\cx{\hskip 16 true pt\flip{$[{2^{3}}]{\times}[{3}]{\times}[{2}]$}\hss}\cg%
\cx{\hskip 16 true pt\flip{$[{2^{2}}{1^{2}}]{\times}[{3}]{\times}[{2}]$}\hss}\cg%
\cx{\hskip 16 true pt\flip{$[{2}{1^{4}}]{\times}[{3}]{\times}[{2}]$}\hss}\cg%
\cx{\hskip 16 true pt\flip{$[{1^{6}}]{\times}[{3}]{\times}[{2}]$}\hss}\cg%
\cx{\hskip 16 true pt\flip{$[{6}]{\times}[{2}{1}]{\times}[{2}]$}\hss}\cg%
\cx{\hskip 16 true pt\flip{$[{5}{1}]{\times}[{2}{1}]{\times}[{2}]$}\hss}\cg%
\cx{\hskip 16 true pt\flip{$[{4}{2}]{\times}[{2}{1}]{\times}[{2}]$}\hss}\cg%
\cx{\hskip 16 true pt\flip{$[{4}{1^{2}}]{\times}[{2}{1}]{\times}[{2}]$}\hss}\cg%
\cx{\hskip 16 true pt\flip{$[{3^{2}}]{\times}[{2}{1}]{\times}[{2}]$}\hss}\cg%
\cx{\hskip 16 true pt\flip{$[{3}{2}{1}]{\times}[{2}{1}]{\times}[{2}]$}\hss}\cg%
\cx{\hskip 16 true pt\flip{$[{3}{1^{3}}]{\times}[{2}{1}]{\times}[{2}]$}\hss}\cg%
\cx{\hskip 16 true pt\flip{$[{2^{3}}]{\times}[{2}{1}]{\times}[{2}]$}\hss}\cg%
\cx{\hskip 16 true pt\flip{$[{2^{2}}{1^{2}}]{\times}[{2}{1}]{\times}[{2}]$}\hss}\cg%
\cx{\hskip 16 true pt\flip{$[{2}{1^{4}}]{\times}[{2}{1}]{\times}[{2}]$}\hss}\cg%
\cx{\hskip 16 true pt\flip{$[{1^{6}}]{\times}[{2}{1}]{\times}[{2}]$}\hss}\cg%
\eol}}\rg%
%
%
\rx{\vss\hfull{%
\rlx{\hss{$70_{y}$}}\cg%
\e{0}%
\e{0}%
\e{0}%
\e{0}%
\e{0}%
\e{0}%
\e{0}%
\e{0}%
\e{0}%
\e{1}%
\e{0}%
\e{0}%
\e{0}%
\e{0}%
\e{0}%
\e{0}%
\e{0}%
\e{1}%
\e{0}%
\e{0}%
\e{0}%
\e{0}%
\eol}\vss}\rg%
%
%
\rx{\vss\hfull{%
\rlx{\hss{$1134_{y}$}}\cg%
\e{0}%
\e{0}%
\e{0}%
\e{2}%
\e{0}%
\e{2}%
\e{2}%
\e{0}%
\e{2}%
\e{0}%
\e{0}%
\e{0}%
\e{1}%
\e{2}%
\e{3}%
\e{1}%
\e{4}%
\e{3}%
\e{1}%
\e{3}%
\e{1}%
\e{0}%
\eol}\vss}\rg%
%
%
\rx{\vss\hfull{%
\rlx{\hss{$1680_{y}$}}\cg%
\e{0}%
\e{0}%
\e{0}%
\e{1}%
\e{0}%
\e{2}%
\e{4}%
\e{1}%
\e{3}%
\e{3}%
\e{1}%
\e{0}%
\e{1}%
\e{2}%
\e{5}%
\e{0}%
\e{5}%
\e{7}%
\e{1}%
\e{3}%
\e{3}%
\e{0}%
\eol}\vss}\rg%
%
%
\rx{\vss\hfull{%
\rlx{\hss{$168_{y}$}}\cg%
\e{0}%
\e{0}%
\e{1}%
\e{0}%
\e{0}%
\e{0}%
\e{0}%
\e{1}%
\e{0}%
\e{0}%
\e{0}%
\e{0}%
\e{0}%
\e{1}%
\e{0}%
\e{0}%
\e{1}%
\e{0}%
\e{1}%
\e{0}%
\e{0}%
\e{0}%
\eol}\vss}\rg%
%
%
\rx{\vss\hfull{%
\rlx{\hss{$420_{y}$}}\cg%
\e{0}%
\e{1}%
\e{1}%
\e{0}%
\e{1}%
\e{1}%
\e{0}%
\e{0}%
\e{0}%
\e{0}%
\e{0}%
\e{0}%
\e{0}%
\e{1}%
\e{1}%
\e{2}%
\e{2}%
\e{0}%
\e{1}%
\e{1}%
\e{0}%
\e{0}%
\eol}\vss}\rg%
%
%
\rx{\vss\hfull{%
\rlx{\hss{$3150_{y}$}}\cg%
\e{0}%
\e{2}%
\e{3}%
\e{4}%
\e{4}%
\e{6}%
\e{3}%
\e{1}%
\e{3}%
\e{0}%
\e{0}%
\e{0}%
\e{3}%
\e{6}%
\e{7}%
\e{6}%
\e{14}%
\e{5}%
\e{3}%
\e{8}%
\e{2}%
\e{0}%
\eol}\vss}\rg%
%
%
\rx{\vss\hfull{%
\rlx{\hss{$4200_{y}$}}\cg%
\e{0}%
\e{2}%
\e{6}%
\e{4}%
\e{2}%
\e{9}%
\e{3}%
\e{4}%
\e{3}%
\e{1}%
\e{0}%
\e{1}%
\e{3}%
\e{11}%
\e{8}%
\e{6}%
\e{18}%
\e{8}%
\e{8}%
\e{8}%
\e{3}%
\e{0}%
\eol}\vss}\rg%
%
%
\rx{\vss\hfull{%
\rlx{\hss{$2688_{y}$}}\cg%
\e{0}%
\e{1}%
\e{3}%
\e{2}%
\e{1}%
\e{6}%
\e{3}%
\e{2}%
\e{2}%
\e{1}%
\e{0}%
\e{1}%
\e{2}%
\e{7}%
\e{5}%
\e{2}%
\e{11}%
\e{7}%
\e{5}%
\e{4}%
\e{3}%
\e{0}%
\eol}\vss}\rg%
%
%
\rx{\vss\hfull{%
\rlx{\hss{$2100_{y}$}}\cg%
\e{0}%
\e{0}%
\e{1}%
\e{2}%
\e{0}%
\e{3}%
\e{5}%
\e{1}%
\e{2}%
\e{3}%
\e{0}%
\e{0}%
\e{1}%
\e{4}%
\e{5}%
\e{1}%
\e{7}%
\e{8}%
\e{3}%
\e{3}%
\e{3}%
\e{0}%
\eol}\vss}\rg%
%
%
\rx{\vss\hfull{%
\rlx{\hss{$1400_{y}$}}\cg%
\e{0}%
\e{0}%
\e{1}%
\e{1}%
\e{0}%
\e{2}%
\e{2}%
\e{2}%
\e{2}%
\e{2}%
\e{1}%
\e{0}%
\e{1}%
\e{2}%
\e{4}%
\e{1}%
\e{5}%
\e{4}%
\e{1}%
\e{3}%
\e{2}%
\e{0}%
\eol}\vss}\rg%
%
%
\rx{\vss\hfull{%
\rlx{\hss{$4536_{y}$}}\cg%
\e{0}%
\e{1}%
\e{3}%
\e{4}%
\e{2}%
\e{8}%
\e{6}%
\e{3}%
\e{6}%
\e{4}%
\e{1}%
\e{0}%
\e{4}%
\e{7}%
\e{12}%
\e{5}%
\e{17}%
\e{12}%
\e{4}%
\e{10}%
\e{5}%
\e{1}%
\eol}\vss}\rg%
%
%
\rx{\vss\hfull{%
\rlx{\hss{$5670_{y}$}}\cg%
\e{0}%
\e{1}%
\e{3}%
\e{6}%
\e{2}%
\e{10}%
\e{8}%
\e{3}%
\e{8}%
\e{4}%
\e{1}%
\e{0}%
\e{5}%
\e{9}%
\e{15}%
\e{6}%
\e{21}%
\e{15}%
\e{5}%
\e{13}%
\e{6}%
\e{1}%
\eol}\vss}\rg%
%
%
\rx{\vss\hfull{%
\rlx{\hss{$4480_{y}$}}\cg%
\e{0}%
\e{2}%
\e{4}%
\e{5}%
\e{3}%
\e{9}%
\e{4}%
\e{2}%
\e{5}%
\e{2}%
\e{0}%
\e{0}%
\e{4}%
\e{8}%
\e{11}%
\e{8}%
\e{18}%
\e{9}%
\e{5}%
\e{11}%
\e{3}%
\e{1}%
\eol}\vss}\rg%
%
%
\rx{\vss\hfull{%
\rlx{\hss{$8_{z}$}}\cg%
\e{0}%
\e{1}%
\e{0}%
\e{0}%
\e{0}%
\e{0}%
\e{0}%
\e{0}%
\e{0}%
\e{0}%
\e{0}%
\e{1}%
\e{0}%
\e{0}%
\e{0}%
\e{0}%
\e{0}%
\e{0}%
\e{0}%
\e{0}%
\e{0}%
\e{0}%
\eol}\vss}\rg%
%
%
\rx{\vss\hfull{%
\rlx{\hss{$56_{z}$}}\cg%
\e{0}%
\e{0}%
\e{0}%
\e{0}%
\e{0}%
\e{0}%
\e{1}%
\e{0}%
\e{0}%
\e{0}%
\e{0}%
\e{0}%
\e{0}%
\e{0}%
\e{1}%
\e{0}%
\e{0}%
\e{0}%
\e{0}%
\e{0}%
\e{0}%
\e{0}%
\eol}\vss}\rg%
%
%
\rx{\vss\hfull{%
\rlx{\hss{$160_{z}$}}\cg%
\e{0}%
\e{2}%
\e{1}%
\e{1}%
\e{0}%
\e{1}%
\e{0}%
\e{0}%
\e{0}%
\e{0}%
\e{0}%
\e{2}%
\e{2}%
\e{1}%
\e{1}%
\e{0}%
\e{0}%
\e{0}%
\e{0}%
\e{0}%
\e{0}%
\e{0}%
\eol}\vss}\rg%
%
%
\rx{\vss\hfull{%
\rlx{\hss{$112_{z}$}}\cg%
\e{2}%
\e{3}%
\e{1}%
\e{1}%
\e{1}%
\e{0}%
\e{0}%
\e{0}%
\e{0}%
\e{0}%
\e{0}%
\e{2}%
\e{2}%
\e{1}%
\e{0}%
\e{0}%
\e{0}%
\e{0}%
\e{0}%
\e{0}%
\e{0}%
\e{0}%
\eol}\vss}\rg%
%
%
\rx{\vss\hfull{%
\rlx{\hss{$840_{z}$}}\cg%
\e{0}%
\e{2}%
\e{2}%
\e{2}%
\e{1}%
\e{3}%
\e{0}%
\e{1}%
\e{1}%
\e{0}%
\e{0}%
\e{1}%
\e{3}%
\e{4}%
\e{3}%
\e{1}%
\e{3}%
\e{1}%
\e{1}%
\e{1}%
\e{0}%
\e{0}%
\eol}\vss}\rg%
%
%
\rx{\vss\hfull{%
\rlx{\hss{$1296_{z}$}}\cg%
\e{0}%
\e{1}%
\e{2}%
\e{4}%
\e{0}%
\e{4}%
\e{4}%
\e{1}%
\e{1}%
\e{1}%
\e{0}%
\e{1}%
\e{5}%
\e{4}%
\e{7}%
\e{1}%
\e{4}%
\e{3}%
\e{0}%
\e{1}%
\e{0}%
\e{0}%
\eol}\vss}\rg%
%
%
\rx{\vss\hfull{%
\rlx{\hss{$1400_{z}$}}\cg%
\e{3}%
\e{7}%
\e{7}%
\e{5}%
\e{3}%
\e{5}%
\e{1}%
\e{0}%
\e{1}%
\e{0}%
\e{0}%
\e{4}%
\e{9}%
\e{8}%
\e{6}%
\e{3}%
\e{4}%
\e{1}%
\e{1}%
\e{0}%
\e{0}%
\e{0}%
\eol}\vss}\rg%
%
%
\rx{\vss\hfull{%
\rlx{\hss{$1008_{z}$}}\cg%
\e{1}%
\e{3}%
\e{4}%
\e{5}%
\e{1}%
\e{3}%
\e{2}%
\e{1}%
\e{1}%
\e{0}%
\e{0}%
\e{2}%
\e{7}%
\e{5}%
\e{5}%
\e{1}%
\e{3}%
\e{1}%
\e{0}%
\e{0}%
\e{0}%
\e{0}%
\eol}\vss}\rg%
%
%
\rx{\vss\hfull{%
\rlx{\hss{$560_{z}$}}\cg%
\e{3}%
\e{6}%
\e{4}%
\e{3}%
\e{2}%
\e{2}%
\e{0}%
\e{0}%
\e{0}%
\e{0}%
\e{0}%
\e{4}%
\e{6}%
\e{4}%
\e{2}%
\e{1}%
\e{1}%
\e{0}%
\e{0}%
\e{0}%
\e{0}%
\e{0}%
\eol}\vss}\rg%
%
%
\rx{\vss\hfull{%
\rlx{\hss{$1400_{zz}$}}\cg%
\e{2}%
\e{5}%
\e{6}%
\e{3}%
\e{3}%
\e{4}%
\e{1}%
\e{0}%
\e{0}%
\e{0}%
\e{0}%
\e{2}%
\e{4}%
\e{7}%
\e{4}%
\e{5}%
\e{6}%
\e{1}%
\e{1}%
\e{1}%
\e{0}%
\e{0}%
\eol}\vss}\rg%
%
%
\rx{\vss\hfull{%
\rlx{\hss{$4200_{z}$}}\cg%
\e{0}%
\e{3}%
\e{6}%
\e{8}%
\e{3}%
\e{9}%
\e{6}%
\e{3}%
\e{3}%
\e{1}%
\e{0}%
\e{1}%
\e{6}%
\e{12}%
\e{11}%
\e{7}%
\e{18}%
\e{8}%
\e{4}%
\e{7}%
\e{2}%
\e{0}%
\eol}\vss}\rg%
%
%
\rx{\vss\hfull{%
\rlx{\hss{$400_{z}$}}\cg%
\e{2}%
\e{4}%
\e{3}%
\e{2}%
\e{1}%
\e{1}%
\e{0}%
\e{0}%
\e{0}%
\e{0}%
\e{0}%
\e{2}%
\e{3}%
\e{3}%
\e{1}%
\e{2}%
\e{1}%
\e{0}%
\e{0}%
\e{0}%
\e{0}%
\e{0}%
\eol}\vss}\rg%
%
%
\rx{\vss\hfull{%
\rlx{\hss{$3240_{z}$}}\cg%
\e{3}%
\e{10}%
\e{12}%
\e{10}%
\e{5}%
\e{10}%
\e{3}%
\e{2}%
\e{2}%
\e{0}%
\e{0}%
\e{5}%
\e{14}%
\e{16}%
\e{12}%
\e{7}%
\e{12}%
\e{3}%
\e{2}%
\e{2}%
\e{0}%
\e{0}%
\eol}\vss}\rg%
%
%
\rx{\vss\hfull{%
\rlx{\hss{$4536_{z}$}}\cg%
\e{2}%
\e{7}%
\e{11}%
\e{7}%
\e{7}%
\e{12}%
\e{2}%
\e{3}%
\e{4}%
\e{0}%
\e{0}%
\e{2}%
\e{9}%
\e{17}%
\e{12}%
\e{9}%
\e{20}%
\e{6}%
\e{6}%
\e{6}%
\e{1}%
\e{0}%
\eol}\vss}\rg%
%
%
\rx{\vss\hfull{%
\rlx{\hss{$2400_{z}$}}\cg%
\e{0}%
\e{0}%
\e{1}%
\e{4}%
\e{0}%
\e{4}%
\e{7}%
\e{1}%
\e{3}%
\e{3}%
\e{0}%
\e{0}%
\e{3}%
\e{4}%
\e{9}%
\e{2}%
\e{8}%
\e{8}%
\e{1}%
\e{4}%
\e{2}%
\e{0}%
\eol}\vss}\rg%
\eop
\eject
\tablecont%
%
%
%
%
%
%
\rowpts=18 true pt%
\colpts=18 true pt%
\rowlabpts=40 true pt%
\collabpts=90 true pt%
\clx{\vss\hfull{%
\rlx{\hss{$ $}}\cg%
\cx{\hskip 16 true pt\flip{$[{6}]{\times}[{3}]{\times}[{2}]$}\hss}\cg%
\cx{\hskip 16 true pt\flip{$[{5}{1}]{\times}[{3}]{\times}[{2}]$}\hss}\cg%
\cx{\hskip 16 true pt\flip{$[{4}{2}]{\times}[{3}]{\times}[{2}]$}\hss}\cg%
\cx{\hskip 16 true pt\flip{$[{4}{1^{2}}]{\times}[{3}]{\times}[{2}]$}\hss}\cg%
\cx{\hskip 16 true pt\flip{$[{3^{2}}]{\times}[{3}]{\times}[{2}]$}\hss}\cg%
\cx{\hskip 16 true pt\flip{$[{3}{2}{1}]{\times}[{3}]{\times}[{2}]$}\hss}\cg%
\cx{\hskip 16 true pt\flip{$[{3}{1^{3}}]{\times}[{3}]{\times}[{2}]$}\hss}\cg%
\cx{\hskip 16 true pt\flip{$[{2^{3}}]{\times}[{3}]{\times}[{2}]$}\hss}\cg%
\cx{\hskip 16 true pt\flip{$[{2^{2}}{1^{2}}]{\times}[{3}]{\times}[{2}]$}\hss}\cg%
\cx{\hskip 16 true pt\flip{$[{2}{1^{4}}]{\times}[{3}]{\times}[{2}]$}\hss}\cg%
\cx{\hskip 16 true pt\flip{$[{1^{6}}]{\times}[{3}]{\times}[{2}]$}\hss}\cg%
\cx{\hskip 16 true pt\flip{$[{6}]{\times}[{2}{1}]{\times}[{2}]$}\hss}\cg%
\cx{\hskip 16 true pt\flip{$[{5}{1}]{\times}[{2}{1}]{\times}[{2}]$}\hss}\cg%
\cx{\hskip 16 true pt\flip{$[{4}{2}]{\times}[{2}{1}]{\times}[{2}]$}\hss}\cg%
\cx{\hskip 16 true pt\flip{$[{4}{1^{2}}]{\times}[{2}{1}]{\times}[{2}]$}\hss}\cg%
\cx{\hskip 16 true pt\flip{$[{3^{2}}]{\times}[{2}{1}]{\times}[{2}]$}\hss}\cg%
\cx{\hskip 16 true pt\flip{$[{3}{2}{1}]{\times}[{2}{1}]{\times}[{2}]$}\hss}\cg%
\cx{\hskip 16 true pt\flip{$[{3}{1^{3}}]{\times}[{2}{1}]{\times}[{2}]$}\hss}\cg%
\cx{\hskip 16 true pt\flip{$[{2^{3}}]{\times}[{2}{1}]{\times}[{2}]$}\hss}\cg%
\cx{\hskip 16 true pt\flip{$[{2^{2}}{1^{2}}]{\times}[{2}{1}]{\times}[{2}]$}\hss}\cg%
\cx{\hskip 16 true pt\flip{$[{2}{1^{4}}]{\times}[{2}{1}]{\times}[{2}]$}\hss}\cg%
\cx{\hskip 16 true pt\flip{$[{1^{6}}]{\times}[{2}{1}]{\times}[{2}]$}\hss}\cg%
\eol}}\rg%
%
%
\rx{\vss\hfull{%
\rlx{\hss{$3360_{z}$}}\cg%
\e{1}%
\e{4}%
\e{7}%
\e{7}%
\e{3}%
\e{9}%
\e{5}%
\e{2}%
\e{2}%
\e{1}%
\e{0}%
\e{2}%
\e{8}%
\e{11}%
\e{11}%
\e{6}%
\e{14}%
\e{6}%
\e{2}%
\e{4}%
\e{1}%
\e{0}%
\eol}\vss}\rg%
%
%
\rx{\vss\hfull{%
\rlx{\hss{$2800_{z}$}}\cg%
\e{0}%
\e{2}%
\e{4}%
\e{6}%
\e{1}%
\e{7}%
\e{6}%
\e{1}%
\e{3}%
\e{2}%
\e{0}%
\e{1}%
\e{6}%
\e{8}%
\e{11}%
\e{4}%
\e{10}%
\e{7}%
\e{2}%
\e{3}%
\e{1}%
\e{0}%
\eol}\vss}\rg%
%
%
\rx{\vss\hfull{%
\rlx{\hss{$4096_{z}$}}\cg%
\e{1}%
\e{6}%
\e{9}%
\e{9}%
\e{4}%
\e{12}%
\e{5}%
\e{3}%
\e{4}%
\e{1}%
\e{0}%
\e{3}%
\e{12}%
\e{15}%
\e{15}%
\e{6}%
\e{16}%
\e{7}%
\e{3}%
\e{4}%
\e{1}%
\e{0}%
\eol}\vss}\rg%
%
%
\rx{\vss\hfull{%
\rlx{\hss{$5600_{z}$}}\cg%
\e{0}%
\e{3}%
\e{6}%
\e{8}%
\e{4}%
\e{12}%
\e{8}%
\e{4}%
\e{7}%
\e{3}%
\e{1}%
\e{1}%
\e{8}%
\e{14}%
\e{17}%
\e{6}%
\e{22}%
\e{13}%
\e{6}%
\e{9}%
\e{4}%
\e{0}%
\eol}\vss}\rg%
%
%
\rx{\vss\hfull{%
\rlx{\hss{$448_{z}$}}\cg%
\e{3}%
\e{3}%
\e{3}%
\e{1}%
\e{3}%
\e{1}%
\e{0}%
\e{0}%
\e{0}%
\e{0}%
\e{0}%
\e{1}%
\e{3}%
\e{3}%
\e{1}%
\e{1}%
\e{2}%
\e{0}%
\e{0}%
\e{0}%
\e{0}%
\e{0}%
\eol}\vss}\rg%
%
%
\rx{\vss\hfull{%
\rlx{\hss{$448_{w}$}}\cg%
\e{0}%
\e{0}%
\e{0}%
\e{0}%
\e{0}%
\e{0}%
\e{2}%
\e{0}%
\e{1}%
\e{1}%
\e{0}%
\e{0}%
\e{0}%
\e{0}%
\e{2}%
\e{0}%
\e{1}%
\e{2}%
\e{0}%
\e{1}%
\e{1}%
\e{0}%
\eol}\vss}\rg%
%
%
\rx{\vss\hfull{%
\rlx{\hss{$1344_{w}$}}\cg%
\e{0}%
\e{1}%
\e{2}%
\e{2}%
\e{1}%
\e{3}%
\e{0}%
\e{1}%
\e{1}%
\e{0}%
\e{0}%
\e{0}%
\e{1}%
\e{4}%
\e{2}%
\e{3}%
\e{6}%
\e{2}%
\e{3}%
\e{3}%
\e{0}%
\e{0}%
\eol}\vss}\rg%
%
%
\rx{\vss\hfull{%
\rlx{\hss{$5600_{w}$}}\cg%
\e{0}%
\e{1}%
\e{3}%
\e{6}%
\e{2}%
\e{10}%
\e{8}%
\e{4}%
\e{7}%
\e{4}%
\e{0}%
\e{0}%
\e{5}%
\e{10}%
\e{14}%
\e{5}%
\e{21}%
\e{15}%
\e{6}%
\e{12}%
\e{6}%
\e{1}%
\eol}\vss}\rg%
%
%
\rx{\vss\hfull{%
\rlx{\hss{$2016_{w}$}}\cg%
\e{1}%
\e{1}%
\e{4}%
\e{2}%
\e{3}%
\e{4}%
\e{1}%
\e{1}%
\e{1}%
\e{0}%
\e{0}%
\e{0}%
\e{2}%
\e{5}%
\e{3}%
\e{4}%
\e{10}%
\e{3}%
\e{3}%
\e{4}%
\e{1}%
\e{0}%
\eol}\vss}\rg%
%
%
\rx{\vss\hfull{%
\rlx{\hss{$7168_{w}$}}\cg%
\e{0}%
\e{3}%
\e{7}%
\e{7}%
\e{4}%
\e{15}%
\e{7}%
\e{4}%
\e{7}%
\e{3}%
\e{0}%
\e{1}%
\e{6}%
\e{15}%
\e{16}%
\e{10}%
\e{29}%
\e{16}%
\e{10}%
\e{15}%
\e{6}%
\e{1}%
\eol}\vss}\rg%
%
%
%
%
%
%
\rowpts=18 true pt%
\colpts=18 true pt%
\rowlabpts=40 true pt%
\collabpts=90 true pt%
\clx{\vss\hfull{%
\rlx{\hss{$ $}}\cg%
\cx{\hskip 16 true pt\flip{$[{6}]{\times}[{1^{3}}]{\times}[{2}]$}\hss}\cg%
\cx{\hskip 16 true pt\flip{$[{5}{1}]{\times}[{1^{3}}]{\times}[{2}]$}\hss}\cg%
\cx{\hskip 16 true pt\flip{$[{4}{2}]{\times}[{1^{3}}]{\times}[{2}]$}\hss}\cg%
\cx{\hskip 16 true pt\flip{$[{4}{1^{2}}]{\times}[{1^{3}}]{\times}[{2}]$}\hss}\cg%
\cx{\hskip 16 true pt\flip{$[{3^{2}}]{\times}[{1^{3}}]{\times}[{2}]$}\hss}\cg%
\cx{\hskip 16 true pt\flip{$[{3}{2}{1}]{\times}[{1^{3}}]{\times}[{2}]$}\hss}\cg%
\cx{\hskip 16 true pt\flip{$[{3}{1^{3}}]{\times}[{1^{3}}]{\times}[{2}]$}\hss}\cg%
\cx{\hskip 16 true pt\flip{$[{2^{3}}]{\times}[{1^{3}}]{\times}[{2}]$}\hss}\cg%
\cx{\hskip 16 true pt\flip{$[{2^{2}}{1^{2}}]{\times}[{1^{3}}]{\times}[{2}]$}\hss}\cg%
\cx{\hskip 16 true pt\flip{$[{2}{1^{4}}]{\times}[{1^{3}}]{\times}[{2}]$}\hss}\cg%
\cx{\hskip 16 true pt\flip{$[{1^{6}}]{\times}[{1^{3}}]{\times}[{2}]$}\hss}\cg%
\cx{\hskip 16 true pt\flip{$[{6}]{\times}[{3}]{\times}[{1^{2}}]$}\hss}\cg%
\cx{\hskip 16 true pt\flip{$[{5}{1}]{\times}[{3}]{\times}[{1^{2}}]$}\hss}\cg%
\cx{\hskip 16 true pt\flip{$[{4}{2}]{\times}[{3}]{\times}[{1^{2}}]$}\hss}\cg%
\cx{\hskip 16 true pt\flip{$[{4}{1^{2}}]{\times}[{3}]{\times}[{1^{2}}]$}\hss}\cg%
\cx{\hskip 16 true pt\flip{$[{3^{2}}]{\times}[{3}]{\times}[{1^{2}}]$}\hss}\cg%
\cx{\hskip 16 true pt\flip{$[{3}{2}{1}]{\times}[{3}]{\times}[{1^{2}}]$}\hss}\cg%
\cx{\hskip 16 true pt\flip{$[{3}{1^{3}}]{\times}[{3}]{\times}[{1^{2}}]$}\hss}\cg%
\cx{\hskip 16 true pt\flip{$[{2^{3}}]{\times}[{3}]{\times}[{1^{2}}]$}\hss}\cg%
\cx{\hskip 16 true pt\flip{$[{2^{2}}{1^{2}}]{\times}[{3}]{\times}[{1^{2}}]$}\hss}\cg%
\cx{\hskip 16 true pt\flip{$[{2}{1^{4}}]{\times}[{3}]{\times}[{1^{2}}]$}\hss}\cg%
\cx{\hskip 16 true pt\flip{$[{1^{6}}]{\times}[{3}]{\times}[{1^{2}}]$}\hss}\cg%
\eol}}\rg%
%
%
\rx{\vss\hfull{%
\rlx{\hss{$1_{x}$}}\cg%
\e{0}%
\e{0}%
\e{0}%
\e{0}%
\e{0}%
\e{0}%
\e{0}%
\e{0}%
\e{0}%
\e{0}%
\e{0}%
\e{0}%
\e{0}%
\e{0}%
\e{0}%
\e{0}%
\e{0}%
\e{0}%
\e{0}%
\e{0}%
\e{0}%
\e{0}%
\eol}\vss}\rg%
%
%
\rx{\vss\hfull{%
\rlx{\hss{$28_{x}$}}\cg%
\e{1}%
\e{0}%
\e{0}%
\e{0}%
\e{0}%
\e{0}%
\e{0}%
\e{0}%
\e{0}%
\e{0}%
\e{0}%
\e{0}%
\e{1}%
\e{0}%
\e{0}%
\e{0}%
\e{0}%
\e{0}%
\e{0}%
\e{0}%
\e{0}%
\e{0}%
\eol}\vss}\rg%
%
%
\rx{\vss\hfull{%
\rlx{\hss{$35_{x}$}}\cg%
\e{0}%
\e{0}%
\e{0}%
\e{0}%
\e{0}%
\e{0}%
\e{0}%
\e{0}%
\e{0}%
\e{0}%
\e{0}%
\e{0}%
\e{1}%
\e{0}%
\e{0}%
\e{0}%
\e{0}%
\e{0}%
\e{0}%
\e{0}%
\e{0}%
\e{0}%
\eol}\vss}\rg%
%
%
\rx{\vss\hfull{%
\rlx{\hss{$84_{x}$}}\cg%
\e{0}%
\e{0}%
\e{0}%
\e{0}%
\e{0}%
\e{0}%
\e{0}%
\e{0}%
\e{0}%
\e{0}%
\e{0}%
\e{1}%
\e{1}%
\e{0}%
\e{0}%
\e{1}%
\e{0}%
\e{0}%
\e{0}%
\e{0}%
\e{0}%
\e{0}%
\eol}\vss}\rg%
%
%
\rx{\vss\hfull{%
\rlx{\hss{$50_{x}$}}\cg%
\e{0}%
\e{0}%
\e{0}%
\e{0}%
\e{0}%
\e{0}%
\e{0}%
\e{0}%
\e{0}%
\e{0}%
\e{0}%
\e{0}%
\e{1}%
\e{0}%
\e{0}%
\e{0}%
\e{0}%
\e{0}%
\e{0}%
\e{0}%
\e{0}%
\e{0}%
\eol}\vss}\rg%
%
%
\rx{\vss\hfull{%
\rlx{\hss{$350_{x}$}}\cg%
\e{1}%
\e{1}%
\e{1}%
\e{1}%
\e{0}%
\e{0}%
\e{0}%
\e{0}%
\e{0}%
\e{0}%
\e{0}%
\e{0}%
\e{1}%
\e{1}%
\e{1}%
\e{0}%
\e{1}%
\e{1}%
\e{0}%
\e{0}%
\e{0}%
\e{0}%
\eol}\vss}\rg%
%
%
\rx{\vss\hfull{%
\rlx{\hss{$300_{x}$}}\cg%
\e{0}%
\e{0}%
\e{0}%
\e{1}%
\e{0}%
\e{0}%
\e{0}%
\e{0}%
\e{0}%
\e{0}%
\e{0}%
\e{0}%
\e{1}%
\e{1}%
\e{1}%
\e{0}%
\e{1}%
\e{0}%
\e{0}%
\e{0}%
\e{0}%
\e{0}%
\eol}\vss}\rg%
%
%
\rx{\vss\hfull{%
\rlx{\hss{$567_{x}$}}\cg%
\e{2}%
\e{2}%
\e{1}%
\e{0}%
\e{0}%
\e{0}%
\e{0}%
\e{0}%
\e{0}%
\e{0}%
\e{0}%
\e{2}%
\e{4}%
\e{2}%
\e{2}%
\e{1}%
\e{1}%
\e{0}%
\e{0}%
\e{0}%
\e{0}%
\e{0}%
\eol}\vss}\rg%
%
%
\rx{\vss\hfull{%
\rlx{\hss{$210_{x}$}}\cg%
\e{0}%
\e{1}%
\e{0}%
\e{0}%
\e{0}%
\e{0}%
\e{0}%
\e{0}%
\e{0}%
\e{0}%
\e{0}%
\e{1}%
\e{2}%
\e{1}%
\e{1}%
\e{0}%
\e{0}%
\e{0}%
\e{0}%
\e{0}%
\e{0}%
\e{0}%
\eol}\vss}\rg%
%
%
\rx{\vss\hfull{%
\rlx{\hss{$840_{x}$}}\cg%
\e{0}%
\e{0}%
\e{0}%
\e{1}%
\e{0}%
\e{1}%
\e{1}%
\e{0}%
\e{1}%
\e{0}%
\e{0}%
\e{0}%
\e{1}%
\e{1}%
\e{1}%
\e{1}%
\e{2}%
\e{0}%
\e{1}%
\e{1}%
\e{0}%
\e{0}%
\eol}\vss}\rg%
%
%
\rx{\vss\hfull{%
\rlx{\hss{$700_{x}$}}\cg%
\e{0}%
\e{1}%
\e{0}%
\e{1}%
\e{1}%
\e{0}%
\e{0}%
\e{0}%
\e{0}%
\e{0}%
\e{0}%
\e{2}%
\e{3}%
\e{3}%
\e{2}%
\e{1}%
\e{1}%
\e{0}%
\e{0}%
\e{0}%
\e{0}%
\e{0}%
\eol}\vss}\rg%
%
%
\rx{\vss\hfull{%
\rlx{\hss{$175_{x}$}}\cg%
\e{0}%
\e{0}%
\e{0}%
\e{0}%
\e{1}%
\e{0}%
\e{0}%
\e{0}%
\e{0}%
\e{0}%
\e{0}%
\e{1}%
\e{0}%
\e{1}%
\e{0}%
\e{1}%
\e{0}%
\e{0}%
\e{0}%
\e{0}%
\e{0}%
\e{0}%
\eol}\vss}\rg%
%
%
\rx{\vss\hfull{%
\rlx{\hss{$1400_{x}$}}\cg%
\e{1}%
\e{3}%
\e{2}%
\e{1}%
\e{2}%
\e{1}%
\e{0}%
\e{0}%
\e{0}%
\e{0}%
\e{0}%
\e{2}%
\e{4}%
\e{5}%
\e{3}%
\e{1}%
\e{3}%
\e{1}%
\e{0}%
\e{0}%
\e{0}%
\e{0}%
\eol}\vss}\rg%
%
%
\rx{\vss\hfull{%
\rlx{\hss{$1050_{x}$}}\cg%
\e{1}%
\e{1}%
\e{2}%
\e{0}%
\e{1}%
\e{1}%
\e{0}%
\e{0}%
\e{0}%
\e{0}%
\e{0}%
\e{1}%
\e{4}%
\e{3}%
\e{2}%
\e{2}%
\e{2}%
\e{0}%
\e{0}%
\e{0}%
\e{0}%
\e{0}%
\eol}\vss}\rg%
\eop
\eject
\tablecont%
%
%
%
%
%
%
\rowpts=18 true pt%
\colpts=18 true pt%
\rowlabpts=40 true pt%
\collabpts=90 true pt%
\clx{\vss\hfull{%
\rlx{\hss{$ $}}\cg%
\cx{\hskip 16 true pt\flip{$[{6}]{\times}[{1^{3}}]{\times}[{2}]$}\hss}\cg%
\cx{\hskip 16 true pt\flip{$[{5}{1}]{\times}[{1^{3}}]{\times}[{2}]$}\hss}\cg%
\cx{\hskip 16 true pt\flip{$[{4}{2}]{\times}[{1^{3}}]{\times}[{2}]$}\hss}\cg%
\cx{\hskip 16 true pt\flip{$[{4}{1^{2}}]{\times}[{1^{3}}]{\times}[{2}]$}\hss}\cg%
\cx{\hskip 16 true pt\flip{$[{3^{2}}]{\times}[{1^{3}}]{\times}[{2}]$}\hss}\cg%
\cx{\hskip 16 true pt\flip{$[{3}{2}{1}]{\times}[{1^{3}}]{\times}[{2}]$}\hss}\cg%
\cx{\hskip 16 true pt\flip{$[{3}{1^{3}}]{\times}[{1^{3}}]{\times}[{2}]$}\hss}\cg%
\cx{\hskip 16 true pt\flip{$[{2^{3}}]{\times}[{1^{3}}]{\times}[{2}]$}\hss}\cg%
\cx{\hskip 16 true pt\flip{$[{2^{2}}{1^{2}}]{\times}[{1^{3}}]{\times}[{2}]$}\hss}\cg%
\cx{\hskip 16 true pt\flip{$[{2}{1^{4}}]{\times}[{1^{3}}]{\times}[{2}]$}\hss}\cg%
\cx{\hskip 16 true pt\flip{$[{1^{6}}]{\times}[{1^{3}}]{\times}[{2}]$}\hss}\cg%
\cx{\hskip 16 true pt\flip{$[{6}]{\times}[{3}]{\times}[{1^{2}}]$}\hss}\cg%
\cx{\hskip 16 true pt\flip{$[{5}{1}]{\times}[{3}]{\times}[{1^{2}}]$}\hss}\cg%
\cx{\hskip 16 true pt\flip{$[{4}{2}]{\times}[{3}]{\times}[{1^{2}}]$}\hss}\cg%
\cx{\hskip 16 true pt\flip{$[{4}{1^{2}}]{\times}[{3}]{\times}[{1^{2}}]$}\hss}\cg%
\cx{\hskip 16 true pt\flip{$[{3^{2}}]{\times}[{3}]{\times}[{1^{2}}]$}\hss}\cg%
\cx{\hskip 16 true pt\flip{$[{3}{2}{1}]{\times}[{3}]{\times}[{1^{2}}]$}\hss}\cg%
\cx{\hskip 16 true pt\flip{$[{3}{1^{3}}]{\times}[{3}]{\times}[{1^{2}}]$}\hss}\cg%
\cx{\hskip 16 true pt\flip{$[{2^{3}}]{\times}[{3}]{\times}[{1^{2}}]$}\hss}\cg%
\cx{\hskip 16 true pt\flip{$[{2^{2}}{1^{2}}]{\times}[{3}]{\times}[{1^{2}}]$}\hss}\cg%
\cx{\hskip 16 true pt\flip{$[{2}{1^{4}}]{\times}[{3}]{\times}[{1^{2}}]$}\hss}\cg%
\cx{\hskip 16 true pt\flip{$[{1^{6}}]{\times}[{3}]{\times}[{1^{2}}]$}\hss}\cg%
\eol}}\rg%
%
%
\rx{\vss\hfull{%
\rlx{\hss{$1575_{x}$}}\cg%
\e{2}%
\e{4}%
\e{3}%
\e{2}%
\e{1}%
\e{1}%
\e{0}%
\e{0}%
\e{0}%
\e{0}%
\e{0}%
\e{2}%
\e{4}%
\e{5}%
\e{5}%
\e{1}%
\e{3}%
\e{2}%
\e{1}%
\e{0}%
\e{0}%
\e{0}%
\eol}\vss}\rg%
%
%
\rx{\vss\hfull{%
\rlx{\hss{$1344_{x}$}}\cg%
\e{1}%
\e{2}%
\e{2}%
\e{1}%
\e{0}%
\e{1}%
\e{0}%
\e{0}%
\e{0}%
\e{0}%
\e{0}%
\e{2}%
\e{5}%
\e{4}%
\e{3}%
\e{3}%
\e{3}%
\e{0}%
\e{0}%
\e{1}%
\e{0}%
\e{0}%
\eol}\vss}\rg%
%
%
\rx{\vss\hfull{%
\rlx{\hss{$2100_{x}$}}\cg%
\e{1}%
\e{4}%
\e{4}%
\e{4}%
\e{0}%
\e{3}%
\e{2}%
\e{0}%
\e{0}%
\e{0}%
\e{0}%
\e{0}%
\e{2}%
\e{3}%
\e{5}%
\e{1}%
\e{5}%
\e{4}%
\e{1}%
\e{2}%
\e{1}%
\e{0}%
\eol}\vss}\rg%
%
%
\rx{\vss\hfull{%
\rlx{\hss{$2268_{x}$}}\cg%
\e{1}%
\e{4}%
\e{3}%
\e{4}%
\e{1}%
\e{2}%
\e{1}%
\e{0}%
\e{0}%
\e{0}%
\e{0}%
\e{1}%
\e{4}%
\e{6}%
\e{6}%
\e{1}%
\e{5}%
\e{3}%
\e{1}%
\e{1}%
\e{0}%
\e{0}%
\eol}\vss}\rg%
%
%
\rx{\vss\hfull{%
\rlx{\hss{$525_{x}$}}\cg%
\e{2}%
\e{1}%
\e{2}%
\e{0}%
\e{0}%
\e{0}%
\e{0}%
\e{1}%
\e{0}%
\e{0}%
\e{0}%
\e{1}%
\e{2}%
\e{1}%
\e{1}%
\e{2}%
\e{1}%
\e{0}%
\e{0}%
\e{1}%
\e{0}%
\e{0}%
\eol}\vss}\rg%
%
%
\rx{\vss\hfull{%
\rlx{\hss{$700_{xx}$}}\cg%
\e{1}%
\e{0}%
\e{2}%
\e{0}%
\e{0}%
\e{1}%
\e{0}%
\e{1}%
\e{0}%
\e{0}%
\e{0}%
\e{0}%
\e{2}%
\e{1}%
\e{2}%
\e{1}%
\e{1}%
\e{0}%
\e{0}%
\e{0}%
\e{0}%
\e{0}%
\eol}\vss}\rg%
%
%
\rx{\vss\hfull{%
\rlx{\hss{$972_{x}$}}\cg%
\e{0}%
\e{0}%
\e{1}%
\e{1}%
\e{0}%
\e{1}%
\e{1}%
\e{0}%
\e{0}%
\e{0}%
\e{0}%
\e{0}%
\e{2}%
\e{2}%
\e{2}%
\e{2}%
\e{2}%
\e{0}%
\e{0}%
\e{1}%
\e{0}%
\e{0}%
\eol}\vss}\rg%
%
%
\rx{\vss\hfull{%
\rlx{\hss{$4096_{x}$}}\cg%
\e{1}%
\e{5}%
\e{6}%
\e{6}%
\e{2}%
\e{5}%
\e{2}%
\e{1}%
\e{1}%
\e{0}%
\e{0}%
\e{1}%
\e{6}%
\e{8}%
\e{8}%
\e{3}%
\e{10}%
\e{4}%
\e{2}%
\e{3}%
\e{1}%
\e{0}%
\eol}\vss}\rg%
%
%
\rx{\vss\hfull{%
\rlx{\hss{$4200_{x}$}}\cg%
\e{1}%
\e{3}%
\e{5}%
\e{5}%
\e{3}%
\e{6}%
\e{2}%
\e{1}%
\e{2}%
\e{0}%
\e{0}%
\e{1}%
\e{5}%
\e{8}%
\e{7}%
\e{4}%
\e{9}%
\e{4}%
\e{3}%
\e{2}%
\e{0}%
\e{0}%
\eol}\vss}\rg%
%
%
\rx{\vss\hfull{%
\rlx{\hss{$2240_{x}$}}\cg%
\e{0}%
\e{2}%
\e{2}%
\e{2}%
\e{2}%
\e{3}%
\e{0}%
\e{0}%
\e{1}%
\e{0}%
\e{0}%
\e{2}%
\e{4}%
\e{6}%
\e{3}%
\e{3}%
\e{5}%
\e{1}%
\e{1}%
\e{1}%
\e{0}%
\e{0}%
\eol}\vss}\rg%
%
%
\rx{\vss\hfull{%
\rlx{\hss{$2835_{x}$}}\cg%
\e{0}%
\e{1}%
\e{3}%
\e{2}%
\e{3}%
\e{5}%
\e{1}%
\e{1}%
\e{2}%
\e{0}%
\e{0}%
\e{1}%
\e{2}%
\e{6}%
\e{3}%
\e{3}%
\e{6}%
\e{2}%
\e{1}%
\e{1}%
\e{0}%
\e{0}%
\eol}\vss}\rg%
%
%
\rx{\vss\hfull{%
\rlx{\hss{$6075_{x}$}}\cg%
\e{2}%
\e{6}%
\e{10}%
\e{7}%
\e{4}%
\e{10}%
\e{4}%
\e{3}%
\e{2}%
\e{1}%
\e{0}%
\e{1}%
\e{6}%
\e{9}%
\e{10}%
\e{5}%
\e{13}%
\e{7}%
\e{3}%
\e{5}%
\e{2}%
\e{0}%
\eol}\vss}\rg%
%
%
\rx{\vss\hfull{%
\rlx{\hss{$3200_{x}$}}\cg%
\e{0}%
\e{1}%
\e{3}%
\e{4}%
\e{1}%
\e{5}%
\e{4}%
\e{2}%
\e{2}%
\e{1}%
\e{0}%
\e{0}%
\e{2}%
\e{3}%
\e{4}%
\e{4}%
\e{7}%
\e{2}%
\e{2}%
\e{5}%
\e{1}%
\e{1}%
\eol}\vss}\rg%
%
%
\rx{\vss\hfull{%
\rlx{\hss{$70_{y}$}}\cg%
\e{0}%
\e{0}%
\e{0}%
\e{1}%
\e{0}%
\e{0}%
\e{0}%
\e{0}%
\e{0}%
\e{0}%
\e{0}%
\e{0}%
\e{0}%
\e{0}%
\e{0}%
\e{0}%
\e{0}%
\e{1}%
\e{0}%
\e{0}%
\e{0}%
\e{0}%
\eol}\vss}\rg%
%
%
\rx{\vss\hfull{%
\rlx{\hss{$1134_{y}$}}\cg%
\e{0}%
\e{1}%
\e{2}%
\e{1}%
\e{1}%
\e{2}%
\e{1}%
\e{1}%
\e{1}%
\e{1}%
\e{0}%
\e{0}%
\e{1}%
\e{1}%
\e{1}%
\e{1}%
\e{2}%
\e{1}%
\e{1}%
\e{2}%
\e{1}%
\e{0}%
\eol}\vss}\rg%
%
%
\rx{\vss\hfull{%
\rlx{\hss{$1680_{y}$}}\cg%
\e{0}%
\e{2}%
\e{2}%
\e{5}%
\e{1}%
\e{3}%
\e{3}%
\e{0}%
\e{1}%
\e{0}%
\e{0}%
\e{0}%
\e{0}%
\e{1}%
\e{3}%
\e{0}%
\e{3}%
\e{5}%
\e{1}%
\e{2}%
\e{2}%
\e{0}%
\eol}\vss}\rg%
%
%
\rx{\vss\hfull{%
\rlx{\hss{$168_{y}$}}\cg%
\e{0}%
\e{0}%
\e{0}%
\e{0}%
\e{0}%
\e{0}%
\e{1}%
\e{0}%
\e{0}%
\e{0}%
\e{0}%
\e{0}%
\e{0}%
\e{0}%
\e{1}%
\e{0}%
\e{0}%
\e{0}%
\e{0}%
\e{0}%
\e{0}%
\e{0}%
\eol}\vss}\rg%
%
%
\rx{\vss\hfull{%
\rlx{\hss{$420_{y}$}}\cg%
\e{0}%
\e{0}%
\e{0}%
\e{0}%
\e{0}%
\e{1}%
\e{0}%
\e{0}%
\e{1}%
\e{0}%
\e{0}%
\e{0}%
\e{0}%
\e{1}%
\e{0}%
\e{0}%
\e{1}%
\e{0}%
\e{0}%
\e{0}%
\e{0}%
\e{0}%
\eol}\vss}\rg%
%
%
\rx{\vss\hfull{%
\rlx{\hss{$3150_{y}$}}\cg%
\e{0}%
\e{1}%
\e{3}%
\e{2}%
\e{3}%
\e{6}%
\e{2}%
\e{3}%
\e{4}%
\e{1}%
\e{1}%
\e{1}%
\e{1}%
\e{4}%
\e{2}%
\e{3}%
\e{6}%
\e{2}%
\e{3}%
\e{3}%
\e{1}%
\e{0}%
\eol}\vss}\rg%
%
%
\rx{\vss\hfull{%
\rlx{\hss{$4200_{y}$}}\cg%
\e{0}%
\e{1}%
\e{4}%
\e{3}%
\e{2}%
\e{8}%
\e{5}%
\e{3}%
\e{4}%
\e{2}%
\e{0}%
\e{0}%
\e{2}%
\e{4}%
\e{5}%
\e{3}%
\e{8}%
\e{3}%
\e{2}%
\e{4}%
\e{1}%
\e{0}%
\eol}\vss}\rg%
%
%
\rx{\vss\hfull{%
\rlx{\hss{$2688_{y}$}}\cg%
\e{0}%
\e{1}%
\e{3}%
\e{2}%
\e{1}%
\e{6}%
\e{3}%
\e{2}%
\e{2}%
\e{1}%
\e{0}%
\e{0}%
\e{1}%
\e{2}%
\e{3}%
\e{2}%
\e{6}%
\e{2}%
\e{1}%
\e{3}%
\e{1}%
\e{0}%
\eol}\vss}\rg%
%
%
\rx{\vss\hfull{%
\rlx{\hss{$2100_{y}$}}\cg%
\e{1}%
\e{2}%
\e{4}%
\e{3}%
\e{0}%
\e{4}%
\e{3}%
\e{2}%
\e{1}%
\e{1}%
\e{0}%
\e{0}%
\e{1}%
\e{1}%
\e{3}%
\e{2}%
\e{4}%
\e{3}%
\e{0}%
\e{4}%
\e{2}%
\e{1}%
\eol}\vss}\rg%
%
%
\rx{\vss\hfull{%
\rlx{\hss{$1400_{y}$}}\cg%
\e{0}%
\e{1}%
\e{1}%
\e{4}%
\e{1}%
\e{2}%
\e{3}%
\e{0}%
\e{1}%
\e{0}%
\e{0}%
\e{0}%
\e{0}%
\e{1}%
\e{3}%
\e{0}%
\e{2}%
\e{4}%
\e{1}%
\e{1}%
\e{1}%
\e{0}%
\eol}\vss}\rg%
%
%
\rx{\vss\hfull{%
\rlx{\hss{$4536_{y}$}}\cg%
\e{0}%
\e{3}%
\e{4}%
\e{8}%
\e{3}%
\e{9}%
\e{6}%
\e{1}%
\e{4}%
\e{1}%
\e{0}%
\e{0}%
\e{1}%
\e{4}%
\e{6}%
\e{1}%
\e{9}%
\e{8}%
\e{3}%
\e{4}%
\e{3}%
\e{0}%
\eol}\vss}\rg%
%
%
\rx{\vss\hfull{%
\rlx{\hss{$5670_{y}$}}\cg%
\e{0}%
\e{4}%
\e{6}%
\e{9}%
\e{4}%
\e{11}%
\e{7}%
\e{2}%
\e{5}%
\e{2}%
\e{0}%
\e{0}%
\e{2}%
\e{5}%
\e{7}%
\e{2}%
\e{11}%
\e{9}%
\e{4}%
\e{6}%
\e{4}%
\e{0}%
\eol}\vss}\rg%
%
%
\rx{\vss\hfull{%
\rlx{\hss{$4480_{y}$}}\cg%
\e{0}%
\e{2}%
\e{4}%
\e{5}%
\e{3}%
\e{9}%
\e{4}%
\e{2}%
\e{5}%
\e{2}%
\e{0}%
\e{0}%
\e{2}%
\e{5}%
\e{4}%
\e{2}%
\e{9}%
\e{5}%
\e{3}%
\e{4}%
\e{2}%
\e{0}%
\eol}\vss}\rg%
%
%
\rx{\vss\hfull{%
\rlx{\hss{$8_{z}$}}\cg%
\e{0}%
\e{0}%
\e{0}%
\e{0}%
\e{0}%
\e{0}%
\e{0}%
\e{0}%
\e{0}%
\e{0}%
\e{0}%
\e{1}%
\e{0}%
\e{0}%
\e{0}%
\e{0}%
\e{0}%
\e{0}%
\e{0}%
\e{0}%
\e{0}%
\e{0}%
\eol}\vss}\rg%
%
%
\rx{\vss\hfull{%
\rlx{\hss{$56_{z}$}}\cg%
\e{0}%
\e{1}%
\e{0}%
\e{0}%
\e{0}%
\e{0}%
\e{0}%
\e{0}%
\e{0}%
\e{0}%
\e{0}%
\e{0}%
\e{0}%
\e{0}%
\e{1}%
\e{0}%
\e{0}%
\e{0}%
\e{0}%
\e{0}%
\e{0}%
\e{0}%
\eol}\vss}\rg%
%
%
\rx{\vss\hfull{%
\rlx{\hss{$160_{z}$}}\cg%
\e{0}%
\e{1}%
\e{0}%
\e{0}%
\e{0}%
\e{0}%
\e{0}%
\e{0}%
\e{0}%
\e{0}%
\e{0}%
\e{1}%
\e{1}%
\e{1}%
\e{1}%
\e{0}%
\e{0}%
\e{0}%
\e{0}%
\e{0}%
\e{0}%
\e{0}%
\eol}\vss}\rg%
%
%
\rx{\vss\hfull{%
\rlx{\hss{$112_{z}$}}\cg%
\e{1}%
\e{0}%
\e{0}%
\e{0}%
\e{0}%
\e{0}%
\e{0}%
\e{0}%
\e{0}%
\e{0}%
\e{0}%
\e{2}%
\e{1}%
\e{1}%
\e{0}%
\e{0}%
\e{0}%
\e{0}%
\e{0}%
\e{0}%
\e{0}%
\e{0}%
\eol}\vss}\rg%
%
%
\rx{\vss\hfull{%
\rlx{\hss{$840_{z}$}}\cg%
\e{0}%
\e{1}%
\e{1}%
\e{1}%
\e{0}%
\e{1}%
\e{1}%
\e{0}%
\e{0}%
\e{0}%
\e{0}%
\e{1}%
\e{1}%
\e{2}%
\e{1}%
\e{1}%
\e{2}%
\e{1}%
\e{1}%
\e{1}%
\e{0}%
\e{0}%
\eol}\vss}\rg%
\eop
\eject
\tablecont%
%
%
%
%
%
%
\rowpts=18 true pt%
\colpts=18 true pt%
\rowlabpts=40 true pt%
\collabpts=90 true pt%
\clx{\vss\hfull{%
\rlx{\hss{$ $}}\cg%
\cx{\hskip 16 true pt\flip{$[{6}]{\times}[{1^{3}}]{\times}[{2}]$}\hss}\cg%
\cx{\hskip 16 true pt\flip{$[{5}{1}]{\times}[{1^{3}}]{\times}[{2}]$}\hss}\cg%
\cx{\hskip 16 true pt\flip{$[{4}{2}]{\times}[{1^{3}}]{\times}[{2}]$}\hss}\cg%
\cx{\hskip 16 true pt\flip{$[{4}{1^{2}}]{\times}[{1^{3}}]{\times}[{2}]$}\hss}\cg%
\cx{\hskip 16 true pt\flip{$[{3^{2}}]{\times}[{1^{3}}]{\times}[{2}]$}\hss}\cg%
\cx{\hskip 16 true pt\flip{$[{3}{2}{1}]{\times}[{1^{3}}]{\times}[{2}]$}\hss}\cg%
\cx{\hskip 16 true pt\flip{$[{3}{1^{3}}]{\times}[{1^{3}}]{\times}[{2}]$}\hss}\cg%
\cx{\hskip 16 true pt\flip{$[{2^{3}}]{\times}[{1^{3}}]{\times}[{2}]$}\hss}\cg%
\cx{\hskip 16 true pt\flip{$[{2^{2}}{1^{2}}]{\times}[{1^{3}}]{\times}[{2}]$}\hss}\cg%
\cx{\hskip 16 true pt\flip{$[{2}{1^{4}}]{\times}[{1^{3}}]{\times}[{2}]$}\hss}\cg%
\cx{\hskip 16 true pt\flip{$[{1^{6}}]{\times}[{1^{3}}]{\times}[{2}]$}\hss}\cg%
\cx{\hskip 16 true pt\flip{$[{6}]{\times}[{3}]{\times}[{1^{2}}]$}\hss}\cg%
\cx{\hskip 16 true pt\flip{$[{5}{1}]{\times}[{3}]{\times}[{1^{2}}]$}\hss}\cg%
\cx{\hskip 16 true pt\flip{$[{4}{2}]{\times}[{3}]{\times}[{1^{2}}]$}\hss}\cg%
\cx{\hskip 16 true pt\flip{$[{4}{1^{2}}]{\times}[{3}]{\times}[{1^{2}}]$}\hss}\cg%
\cx{\hskip 16 true pt\flip{$[{3^{2}}]{\times}[{3}]{\times}[{1^{2}}]$}\hss}\cg%
\cx{\hskip 16 true pt\flip{$[{3}{2}{1}]{\times}[{3}]{\times}[{1^{2}}]$}\hss}\cg%
\cx{\hskip 16 true pt\flip{$[{3}{1^{3}}]{\times}[{3}]{\times}[{1^{2}}]$}\hss}\cg%
\cx{\hskip 16 true pt\flip{$[{2^{3}}]{\times}[{3}]{\times}[{1^{2}}]$}\hss}\cg%
\cx{\hskip 16 true pt\flip{$[{2^{2}}{1^{2}}]{\times}[{3}]{\times}[{1^{2}}]$}\hss}\cg%
\cx{\hskip 16 true pt\flip{$[{2}{1^{4}}]{\times}[{3}]{\times}[{1^{2}}]$}\hss}\cg%
\cx{\hskip 16 true pt\flip{$[{1^{6}}]{\times}[{3}]{\times}[{1^{2}}]$}\hss}\cg%
\eol}}\rg%
%
%
\rx{\vss\hfull{%
\rlx{\hss{$1296_{z}$}}\cg%
\e{1}%
\e{4}%
\e{2}%
\e{3}%
\e{1}%
\e{1}%
\e{0}%
\e{0}%
\e{0}%
\e{0}%
\e{0}%
\e{0}%
\e{3}%
\e{2}%
\e{5}%
\e{1}%
\e{3}%
\e{2}%
\e{0}%
\e{1}%
\e{0}%
\e{0}%
\eol}\vss}\rg%
%
%
\rx{\vss\hfull{%
\rlx{\hss{$1400_{z}$}}\cg%
\e{1}%
\e{3}%
\e{2}%
\e{1}%
\e{0}%
\e{1}%
\e{0}%
\e{0}%
\e{0}%
\e{0}%
\e{0}%
\e{2}%
\e{5}%
\e{5}%
\e{3}%
\e{2}%
\e{3}%
\e{1}%
\e{1}%
\e{0}%
\e{0}%
\e{0}%
\eol}\vss}\rg%
%
%
\rx{\vss\hfull{%
\rlx{\hss{$1008_{z}$}}\cg%
\e{2}%
\e{3}%
\e{1}%
\e{2}%
\e{1}%
\e{0}%
\e{0}%
\e{0}%
\e{0}%
\e{0}%
\e{0}%
\e{1}%
\e{4}%
\e{3}%
\e{4}%
\e{1}%
\e{2}%
\e{1}%
\e{0}%
\e{0}%
\e{0}%
\e{0}%
\eol}\vss}\rg%
%
%
\rx{\vss\hfull{%
\rlx{\hss{$560_{z}$}}\cg%
\e{1}%
\e{1}%
\e{1}%
\e{0}%
\e{0}%
\e{0}%
\e{0}%
\e{0}%
\e{0}%
\e{0}%
\e{0}%
\e{2}%
\e{4}%
\e{3}%
\e{1}%
\e{0}%
\e{1}%
\e{0}%
\e{0}%
\e{0}%
\e{0}%
\e{0}%
\eol}\vss}\rg%
%
%
\rx{\vss\hfull{%
\rlx{\hss{$1400_{zz}$}}\cg%
\e{0}%
\e{1}%
\e{2}%
\e{0}%
\e{1}%
\e{2}%
\e{0}%
\e{1}%
\e{0}%
\e{0}%
\e{0}%
\e{1}%
\e{3}%
\e{4}%
\e{2}%
\e{2}%
\e{3}%
\e{0}%
\e{0}%
\e{0}%
\e{0}%
\e{0}%
\eol}\vss}\rg%
%
%
\rx{\vss\hfull{%
\rlx{\hss{$4200_{z}$}}\cg%
\e{1}%
\e{3}%
\e{5}%
\e{4}%
\e{5}%
\e{7}%
\e{3}%
\e{1}%
\e{3}%
\e{0}%
\e{0}%
\e{1}%
\e{4}%
\e{6}%
\e{7}%
\e{4}%
\e{8}%
\e{4}%
\e{1}%
\e{3}%
\e{0}%
\e{0}%
\eol}\vss}\rg%
%
%
\rx{\vss\hfull{%
\rlx{\hss{$400_{z}$}}\cg%
\e{1}%
\e{0}%
\e{1}%
\e{0}%
\e{0}%
\e{0}%
\e{0}%
\e{0}%
\e{0}%
\e{0}%
\e{0}%
\e{2}%
\e{2}%
\e{2}%
\e{1}%
\e{1}%
\e{0}%
\e{0}%
\e{0}%
\e{0}%
\e{0}%
\e{0}%
\eol}\vss}\rg%
%
%
\rx{\vss\hfull{%
\rlx{\hss{$3240_{z}$}}\cg%
\e{2}%
\e{4}%
\e{5}%
\e{3}%
\e{2}%
\e{3}%
\e{1}%
\e{0}%
\e{0}%
\e{0}%
\e{0}%
\e{3}%
\e{8}%
\e{9}%
\e{7}%
\e{4}%
\e{7}%
\e{2}%
\e{1}%
\e{1}%
\e{0}%
\e{0}%
\eol}\vss}\rg%
%
%
\rx{\vss\hfull{%
\rlx{\hss{$4536_{z}$}}\cg%
\e{0}%
\e{2}%
\e{5}%
\e{4}%
\e{2}%
\e{7}%
\e{3}%
\e{3}%
\e{2}%
\e{1}%
\e{0}%
\e{1}%
\e{4}%
\e{9}%
\e{5}%
\e{4}%
\e{10}%
\e{3}%
\e{4}%
\e{3}%
\e{1}%
\e{0}%
\eol}\vss}\rg%
%
%
\rx{\vss\hfull{%
\rlx{\hss{$2400_{z}$}}\cg%
\e{1}%
\e{4}%
\e{4}%
\e{5}%
\e{2}%
\e{4}%
\e{2}%
\e{0}%
\e{1}%
\e{0}%
\e{0}%
\e{0}%
\e{2}%
\e{2}%
\e{5}%
\e{2}%
\e{5}%
\e{5}%
\e{0}%
\e{3}%
\e{1}%
\e{0}%
\eol}\vss}\rg%
%
%
\rx{\vss\hfull{%
\rlx{\hss{$3360_{z}$}}\cg%
\e{1}%
\e{3}%
\e{5}%
\e{4}%
\e{3}%
\e{5}%
\e{1}%
\e{1}%
\e{1}%
\e{0}%
\e{0}%
\e{1}%
\e{4}%
\e{6}%
\e{7}%
\e{4}%
\e{7}%
\e{3}%
\e{1}%
\e{2}%
\e{0}%
\e{0}%
\eol}\vss}\rg%
%
%
\rx{\vss\hfull{%
\rlx{\hss{$2800_{z}$}}\cg%
\e{1}%
\e{5}%
\e{4}%
\e{5}%
\e{2}%
\e{4}%
\e{1}%
\e{0}%
\e{1}%
\e{0}%
\e{0}%
\e{0}%
\e{4}%
\e{4}%
\e{7}%
\e{2}%
\e{6}%
\e{4}%
\e{1}%
\e{2}%
\e{1}%
\e{0}%
\eol}\vss}\rg%
%
%
\rx{\vss\hfull{%
\rlx{\hss{$4096_{z}$}}\cg%
\e{1}%
\e{5}%
\e{6}%
\e{6}%
\e{2}%
\e{5}%
\e{2}%
\e{1}%
\e{1}%
\e{0}%
\e{0}%
\e{1}%
\e{6}%
\e{8}%
\e{8}%
\e{3}%
\e{10}%
\e{4}%
\e{2}%
\e{3}%
\e{1}%
\e{0}%
\eol}\vss}\rg%
%
%
\rx{\vss\hfull{%
\rlx{\hss{$5600_{z}$}}\cg%
\e{1}%
\e{5}%
\e{7}%
\e{9}%
\e{3}%
\e{9}%
\e{6}%
\e{2}%
\e{3}%
\e{1}%
\e{0}%
\e{1}%
\e{3}%
\e{7}%
\e{8}%
\e{4}%
\e{12}%
\e{8}%
\e{4}%
\e{6}%
\e{3}%
\e{0}%
\eol}\vss}\rg%
%
%
\rx{\vss\hfull{%
\rlx{\hss{$448_{z}$}}\cg%
\e{0}%
\e{0}%
\e{1}%
\e{0}%
\e{0}%
\e{0}%
\e{0}%
\e{1}%
\e{0}%
\e{0}%
\e{0}%
\e{1}%
\e{1}%
\e{2}%
\e{0}%
\e{1}%
\e{1}%
\e{0}%
\e{1}%
\e{0}%
\e{0}%
\e{0}%
\eol}\vss}\rg%
%
%
\rx{\vss\hfull{%
\rlx{\hss{$448_{w}$}}\cg%
\e{0}%
\e{1}%
\e{1}%
\e{1}%
\e{0}%
\e{1}%
\e{1}%
\e{0}%
\e{0}%
\e{0}%
\e{0}%
\e{0}%
\e{0}%
\e{0}%
\e{1}%
\e{0}%
\e{1}%
\e{1}%
\e{0}%
\e{1}%
\e{1}%
\e{0}%
\eol}\vss}\rg%
%
%
\rx{\vss\hfull{%
\rlx{\hss{$1344_{w}$}}\cg%
\e{0}%
\e{0}%
\e{1}%
\e{1}%
\e{1}%
\e{2}%
\e{1}%
\e{1}%
\e{2}%
\e{1}%
\e{0}%
\e{0}%
\e{1}%
\e{2}%
\e{1}%
\e{1}%
\e{2}%
\e{1}%
\e{1}%
\e{1}%
\e{0}%
\e{0}%
\eol}\vss}\rg%
%
%
\rx{\vss\hfull{%
\rlx{\hss{$5600_{w}$}}\cg%
\e{1}%
\e{3}%
\e{7}%
\e{8}%
\e{3}%
\e{11}%
\e{7}%
\e{2}%
\e{5}%
\e{2}%
\e{0}%
\e{0}%
\e{2}%
\e{5}%
\e{7}%
\e{2}%
\e{11}%
\e{8}%
\e{3}%
\e{7}%
\e{3}%
\e{1}%
\eol}\vss}\rg%
%
%
\rx{\vss\hfull{%
\rlx{\hss{$2016_{w}$}}\cg%
\e{0}%
\e{0}%
\e{1}%
\e{1}%
\e{2}%
\e{4}%
\e{1}%
\e{3}%
\e{2}%
\e{1}%
\e{0}%
\e{0}%
\e{1}%
\e{2}%
\e{1}%
\e{3}%
\e{4}%
\e{1}%
\e{2}%
\e{1}%
\e{0}%
\e{0}%
\eol}\vss}\rg%
%
%
\rx{\vss\hfull{%
\rlx{\hss{$7168_{w}$}}\cg%
\e{0}%
\e{3}%
\e{7}%
\e{7}%
\e{4}%
\e{15}%
\e{7}%
\e{4}%
\e{7}%
\e{3}%
\e{0}%
\e{0}%
\e{3}%
\e{7}%
\e{7}%
\e{4}%
\e{15}%
\e{7}%
\e{4}%
\e{7}%
\e{3}%
\e{0}%
\eol}\vss}\rg%
%
%
%
%
%
%
\rowpts=18 true pt%
\colpts=18 true pt%
\rowlabpts=40 true pt%
\collabpts=90 true pt%
\clx{\vss\hfull{%
\rlx{\hss{$ $}}\cg%
\cx{\hskip 16 true pt\flip{$[{6}]{\times}[{2}{1}]{\times}[{1^{2}}]$}\hss}\cg%
\cx{\hskip 16 true pt\flip{$[{5}{1}]{\times}[{2}{1}]{\times}[{1^{2}}]$}\hss}\cg%
\cx{\hskip 16 true pt\flip{$[{4}{2}]{\times}[{2}{1}]{\times}[{1^{2}}]$}\hss}\cg%
\cx{\hskip 16 true pt\flip{$[{4}{1^{2}}]{\times}[{2}{1}]{\times}[{1^{2}}]$}\hss}\cg%
\cx{\hskip 16 true pt\flip{$[{3^{2}}]{\times}[{2}{1}]{\times}[{1^{2}}]$}\hss}\cg%
\cx{\hskip 16 true pt\flip{$[{3}{2}{1}]{\times}[{2}{1}]{\times}[{1^{2}}]$}\hss}\cg%
\cx{\hskip 16 true pt\flip{$[{3}{1^{3}}]{\times}[{2}{1}]{\times}[{1^{2}}]$}\hss}\cg%
\cx{\hskip 16 true pt\flip{$[{2^{3}}]{\times}[{2}{1}]{\times}[{1^{2}}]$}\hss}\cg%
\cx{\hskip 16 true pt\flip{$[{2^{2}}{1^{2}}]{\times}[{2}{1}]{\times}[{1^{2}}]$}\hss}\cg%
\cx{\hskip 16 true pt\flip{$[{2}{1^{4}}]{\times}[{2}{1}]{\times}[{1^{2}}]$}\hss}\cg%
\cx{\hskip 16 true pt\flip{$[{1^{6}}]{\times}[{2}{1}]{\times}[{1^{2}}]$}\hss}\cg%
\cx{\hskip 16 true pt\flip{$[{6}]{\times}[{1^{3}}]{\times}[{1^{2}}]$}\hss}\cg%
\cx{\hskip 16 true pt\flip{$[{5}{1}]{\times}[{1^{3}}]{\times}[{1^{2}}]$}\hss}\cg%
\cx{\hskip 16 true pt\flip{$[{4}{2}]{\times}[{1^{3}}]{\times}[{1^{2}}]$}\hss}\cg%
\cx{\hskip 16 true pt\flip{$[{4}{1^{2}}]{\times}[{1^{3}}]{\times}[{1^{2}}]$}\hss}\cg%
\cx{\hskip 16 true pt\flip{$[{3^{2}}]{\times}[{1^{3}}]{\times}[{1^{2}}]$}\hss}\cg%
\cx{\hskip 16 true pt\flip{$[{3}{2}{1}]{\times}[{1^{3}}]{\times}[{1^{2}}]$}\hss}\cg%
\cx{\hskip 16 true pt\flip{$[{3}{1^{3}}]{\times}[{1^{3}}]{\times}[{1^{2}}]$}\hss}\cg%
\cx{\hskip 16 true pt\flip{$[{2^{3}}]{\times}[{1^{3}}]{\times}[{1^{2}}]$}\hss}\cg%
\cx{\hskip 16 true pt\flip{$[{2^{2}}{1^{2}}]{\times}[{1^{3}}]{\times}[{1^{2}}]$}\hss}\cg%
\cx{\hskip 16 true pt\flip{$[{2}{1^{4}}]{\times}[{1^{3}}]{\times}[{1^{2}}]$}\hss}\cg%
\cx{\hskip 16 true pt\flip{$[{1^{6}}]{\times}[{1^{3}}]{\times}[{1^{2}}]$}\hss}\cg%
\eol}}\rg%
%
%
\rx{\vss\hfull{%
\rlx{\hss{$1_{x}$}}\cg%
\e{0}%
\e{0}%
\e{0}%
\e{0}%
\e{0}%
\e{0}%
\e{0}%
\e{0}%
\e{0}%
\e{0}%
\e{0}%
\e{0}%
\e{0}%
\e{0}%
\e{0}%
\e{0}%
\e{0}%
\e{0}%
\e{0}%
\e{0}%
\e{0}%
\e{0}%
\eol}\vss}\rg%
%
%
\rx{\vss\hfull{%
\rlx{\hss{$28_{x}$}}\cg%
\e{1}%
\e{0}%
\e{0}%
\e{0}%
\e{0}%
\e{0}%
\e{0}%
\e{0}%
\e{0}%
\e{0}%
\e{0}%
\e{0}%
\e{0}%
\e{0}%
\e{0}%
\e{0}%
\e{0}%
\e{0}%
\e{0}%
\e{0}%
\e{0}%
\e{0}%
\eol}\vss}\rg%
%
%
\rx{\vss\hfull{%
\rlx{\hss{$35_{x}$}}\cg%
\e{1}%
\e{0}%
\e{0}%
\e{0}%
\e{0}%
\e{0}%
\e{0}%
\e{0}%
\e{0}%
\e{0}%
\e{0}%
\e{0}%
\e{0}%
\e{0}%
\e{0}%
\e{0}%
\e{0}%
\e{0}%
\e{0}%
\e{0}%
\e{0}%
\e{0}%
\eol}\vss}\rg%
%
%
\rx{\vss\hfull{%
\rlx{\hss{$84_{x}$}}\cg%
\e{0}%
\e{1}%
\e{0}%
\e{0}%
\e{0}%
\e{0}%
\e{0}%
\e{0}%
\e{0}%
\e{0}%
\e{0}%
\e{0}%
\e{0}%
\e{0}%
\e{0}%
\e{0}%
\e{0}%
\e{0}%
\e{0}%
\e{0}%
\e{0}%
\e{0}%
\eol}\vss}\rg%
%
%
\rx{\vss\hfull{%
\rlx{\hss{$50_{x}$}}\cg%
\e{0}%
\e{0}%
\e{0}%
\e{0}%
\e{1}%
\e{0}%
\e{0}%
\e{0}%
\e{0}%
\e{0}%
\e{0}%
\e{0}%
\e{0}%
\e{0}%
\e{0}%
\e{0}%
\e{0}%
\e{0}%
\e{0}%
\e{0}%
\e{0}%
\e{0}%
\eol}\vss}\rg%
%
%
\rx{\vss\hfull{%
\rlx{\hss{$350_{x}$}}\cg%
\e{1}%
\e{2}%
\e{1}%
\e{2}%
\e{0}%
\e{0}%
\e{0}%
\e{0}%
\e{0}%
\e{0}%
\e{0}%
\e{0}%
\e{2}%
\e{0}%
\e{0}%
\e{0}%
\e{0}%
\e{0}%
\e{0}%
\e{0}%
\e{0}%
\e{0}%
\eol}\vss}\rg%
%
%
\rx{\vss\hfull{%
\rlx{\hss{$300_{x}$}}\cg%
\e{1}%
\e{2}%
\e{1}%
\e{1}%
\e{0}%
\e{0}%
\e{0}%
\e{0}%
\e{0}%
\e{0}%
\e{0}%
\e{0}%
\e{1}%
\e{0}%
\e{0}%
\e{0}%
\e{0}%
\e{0}%
\e{0}%
\e{0}%
\e{0}%
\e{0}%
\eol}\vss}\rg%
\eop
\eject
\tablecont%
%
%
%
%
%
%
\rowpts=18 true pt%
\colpts=18 true pt%
\rowlabpts=40 true pt%
\collabpts=90 true pt%
\clx{\vss\hfull{%
\rlx{\hss{$ $}}\cg%
\cx{\hskip 16 true pt\flip{$[{6}]{\times}[{2}{1}]{\times}[{1^{2}}]$}\hss}\cg%
\cx{\hskip 16 true pt\flip{$[{5}{1}]{\times}[{2}{1}]{\times}[{1^{2}}]$}\hss}\cg%
\cx{\hskip 16 true pt\flip{$[{4}{2}]{\times}[{2}{1}]{\times}[{1^{2}}]$}\hss}\cg%
\cx{\hskip 16 true pt\flip{$[{4}{1^{2}}]{\times}[{2}{1}]{\times}[{1^{2}}]$}\hss}\cg%
\cx{\hskip 16 true pt\flip{$[{3^{2}}]{\times}[{2}{1}]{\times}[{1^{2}}]$}\hss}\cg%
\cx{\hskip 16 true pt\flip{$[{3}{2}{1}]{\times}[{2}{1}]{\times}[{1^{2}}]$}\hss}\cg%
\cx{\hskip 16 true pt\flip{$[{3}{1^{3}}]{\times}[{2}{1}]{\times}[{1^{2}}]$}\hss}\cg%
\cx{\hskip 16 true pt\flip{$[{2^{3}}]{\times}[{2}{1}]{\times}[{1^{2}}]$}\hss}\cg%
\cx{\hskip 16 true pt\flip{$[{2^{2}}{1^{2}}]{\times}[{2}{1}]{\times}[{1^{2}}]$}\hss}\cg%
\cx{\hskip 16 true pt\flip{$[{2}{1^{4}}]{\times}[{2}{1}]{\times}[{1^{2}}]$}\hss}\cg%
\cx{\hskip 16 true pt\flip{$[{1^{6}}]{\times}[{2}{1}]{\times}[{1^{2}}]$}\hss}\cg%
\cx{\hskip 16 true pt\flip{$[{6}]{\times}[{1^{3}}]{\times}[{1^{2}}]$}\hss}\cg%
\cx{\hskip 16 true pt\flip{$[{5}{1}]{\times}[{1^{3}}]{\times}[{1^{2}}]$}\hss}\cg%
\cx{\hskip 16 true pt\flip{$[{4}{2}]{\times}[{1^{3}}]{\times}[{1^{2}}]$}\hss}\cg%
\cx{\hskip 16 true pt\flip{$[{4}{1^{2}}]{\times}[{1^{3}}]{\times}[{1^{2}}]$}\hss}\cg%
\cx{\hskip 16 true pt\flip{$[{3^{2}}]{\times}[{1^{3}}]{\times}[{1^{2}}]$}\hss}\cg%
\cx{\hskip 16 true pt\flip{$[{3}{2}{1}]{\times}[{1^{3}}]{\times}[{1^{2}}]$}\hss}\cg%
\cx{\hskip 16 true pt\flip{$[{3}{1^{3}}]{\times}[{1^{3}}]{\times}[{1^{2}}]$}\hss}\cg%
\cx{\hskip 16 true pt\flip{$[{2^{3}}]{\times}[{1^{3}}]{\times}[{1^{2}}]$}\hss}\cg%
\cx{\hskip 16 true pt\flip{$[{2^{2}}{1^{2}}]{\times}[{1^{3}}]{\times}[{1^{2}}]$}\hss}\cg%
\cx{\hskip 16 true pt\flip{$[{2}{1^{4}}]{\times}[{1^{3}}]{\times}[{1^{2}}]$}\hss}\cg%
\cx{\hskip 16 true pt\flip{$[{1^{6}}]{\times}[{1^{3}}]{\times}[{1^{2}}]$}\hss}\cg%
\eol}}\rg%
%
%
\rx{\vss\hfull{%
\rlx{\hss{$567_{x}$}}\cg%
\e{3}%
\e{4}%
\e{2}%
\e{1}%
\e{0}%
\e{0}%
\e{0}%
\e{0}%
\e{0}%
\e{0}%
\e{0}%
\e{1}%
\e{1}%
\e{0}%
\e{0}%
\e{0}%
\e{0}%
\e{0}%
\e{0}%
\e{0}%
\e{0}%
\e{0}%
\eol}\vss}\rg%
%
%
\rx{\vss\hfull{%
\rlx{\hss{$210_{x}$}}\cg%
\e{2}%
\e{1}%
\e{1}%
\e{0}%
\e{0}%
\e{0}%
\e{0}%
\e{0}%
\e{0}%
\e{0}%
\e{0}%
\e{1}%
\e{0}%
\e{0}%
\e{0}%
\e{0}%
\e{0}%
\e{0}%
\e{0}%
\e{0}%
\e{0}%
\e{0}%
\eol}\vss}\rg%
%
%
\rx{\vss\hfull{%
\rlx{\hss{$840_{x}$}}\cg%
\e{0}%
\e{1}%
\e{3}%
\e{2}%
\e{1}%
\e{3}%
\e{1}%
\e{1}%
\e{1}%
\e{0}%
\e{0}%
\e{0}%
\e{0}%
\e{1}%
\e{1}%
\e{0}%
\e{1}%
\e{1}%
\e{0}%
\e{0}%
\e{0}%
\e{0}%
\eol}\vss}\rg%
%
%
\rx{\vss\hfull{%
\rlx{\hss{$700_{x}$}}\cg%
\e{2}%
\e{3}%
\e{3}%
\e{1}%
\e{1}%
\e{1}%
\e{0}%
\e{0}%
\e{0}%
\e{0}%
\e{0}%
\e{1}%
\e{0}%
\e{1}%
\e{0}%
\e{0}%
\e{0}%
\e{0}%
\e{0}%
\e{0}%
\e{0}%
\e{0}%
\eol}\vss}\rg%
%
%
\rx{\vss\hfull{%
\rlx{\hss{$175_{x}$}}\cg%
\e{0}%
\e{0}%
\e{1}%
\e{0}%
\e{0}%
\e{1}%
\e{0}%
\e{0}%
\e{0}%
\e{0}%
\e{0}%
\e{0}%
\e{0}%
\e{0}%
\e{0}%
\e{0}%
\e{0}%
\e{0}%
\e{1}%
\e{0}%
\e{0}%
\e{0}%
\eol}\vss}\rg%
%
%
\rx{\vss\hfull{%
\rlx{\hss{$1400_{x}$}}\cg%
\e{3}%
\e{4}%
\e{6}%
\e{3}%
\e{2}%
\e{3}%
\e{1}%
\e{1}%
\e{0}%
\e{0}%
\e{0}%
\e{1}%
\e{1}%
\e{2}%
\e{0}%
\e{0}%
\e{1}%
\e{0}%
\e{0}%
\e{0}%
\e{0}%
\e{0}%
\eol}\vss}\rg%
%
%
\rx{\vss\hfull{%
\rlx{\hss{$1050_{x}$}}\cg%
\e{1}%
\e{3}%
\e{3}%
\e{2}%
\e{4}%
\e{3}%
\e{0}%
\e{0}%
\e{1}%
\e{0}%
\e{0}%
\e{0}%
\e{0}%
\e{1}%
\e{0}%
\e{1}%
\e{1}%
\e{0}%
\e{0}%
\e{0}%
\e{0}%
\e{0}%
\eol}\vss}\rg%
%
%
\rx{\vss\hfull{%
\rlx{\hss{$1575_{x}$}}\cg%
\e{3}%
\e{7}%
\e{6}%
\e{4}%
\e{1}%
\e{3}%
\e{1}%
\e{0}%
\e{0}%
\e{0}%
\e{0}%
\e{2}%
\e{2}%
\e{2}%
\e{1}%
\e{1}%
\e{0}%
\e{0}%
\e{0}%
\e{0}%
\e{0}%
\e{0}%
\eol}\vss}\rg%
%
%
\rx{\vss\hfull{%
\rlx{\hss{$1344_{x}$}}\cg%
\e{2}%
\e{7}%
\e{4}%
\e{4}%
\e{2}%
\e{2}%
\e{0}%
\e{0}%
\e{0}%
\e{0}%
\e{0}%
\e{0}%
\e{2}%
\e{1}%
\e{1}%
\e{0}%
\e{0}%
\e{0}%
\e{0}%
\e{0}%
\e{0}%
\e{0}%
\eol}\vss}\rg%
%
%
\rx{\vss\hfull{%
\rlx{\hss{$2100_{x}$}}\cg%
\e{1}%
\e{7}%
\e{5}%
\e{9}%
\e{2}%
\e{5}%
\e{3}%
\e{0}%
\e{1}%
\e{0}%
\e{0}%
\e{1}%
\e{5}%
\e{2}%
\e{4}%
\e{1}%
\e{1}%
\e{0}%
\e{0}%
\e{0}%
\e{0}%
\e{0}%
\eol}\vss}\rg%
%
%
\rx{\vss\hfull{%
\rlx{\hss{$2268_{x}$}}\cg%
\e{3}%
\e{8}%
\e{8}%
\e{7}%
\e{2}%
\e{5}%
\e{2}%
\e{1}%
\e{0}%
\e{0}%
\e{0}%
\e{2}%
\e{4}%
\e{2}%
\e{2}%
\e{1}%
\e{1}%
\e{0}%
\e{0}%
\e{0}%
\e{0}%
\e{0}%
\eol}\vss}\rg%
%
%
\rx{\vss\hfull{%
\rlx{\hss{$525_{x}$}}\cg%
\e{0}%
\e{3}%
\e{1}%
\e{2}%
\e{1}%
\e{1}%
\e{0}%
\e{0}%
\e{0}%
\e{0}%
\e{0}%
\e{0}%
\e{1}%
\e{0}%
\e{1}%
\e{0}%
\e{0}%
\e{0}%
\e{0}%
\e{0}%
\e{0}%
\e{0}%
\eol}\vss}\rg%
%
%
\rx{\vss\hfull{%
\rlx{\hss{$700_{xx}$}}\cg%
\e{0}%
\e{1}%
\e{1}%
\e{1}%
\e{4}%
\e{3}%
\e{0}%
\e{0}%
\e{2}%
\e{0}%
\e{0}%
\e{0}%
\e{0}%
\e{0}%
\e{0}%
\e{2}%
\e{1}%
\e{0}%
\e{0}%
\e{1}%
\e{0}%
\e{0}%
\eol}\vss}\rg%
%
%
\rx{\vss\hfull{%
\rlx{\hss{$972_{x}$}}\cg%
\e{0}%
\e{3}%
\e{2}%
\e{3}%
\e{3}%
\e{3}%
\e{0}%
\e{0}%
\e{1}%
\e{0}%
\e{0}%
\e{0}%
\e{1}%
\e{0}%
\e{1}%
\e{1}%
\e{1}%
\e{0}%
\e{0}%
\e{0}%
\e{0}%
\e{0}%
\eol}\vss}\rg%
%
%
\rx{\vss\hfull{%
\rlx{\hss{$4096_{x}$}}\cg%
\e{3}%
\e{11}%
\e{13}%
\e{13}%
\e{5}%
\e{11}%
\e{5}%
\e{2}%
\e{2}%
\e{0}%
\e{0}%
\e{1}%
\e{5}%
\e{5}%
\e{5}%
\e{1}%
\e{3}%
\e{1}%
\e{0}%
\e{0}%
\e{0}%
\e{0}%
\eol}\vss}\rg%
%
%
\rx{\vss\hfull{%
\rlx{\hss{$4200_{x}$}}\cg%
\e{2}%
\e{8}%
\e{13}%
\e{10}%
\e{6}%
\e{15}%
\e{5}%
\e{3}%
\e{4}%
\e{1}%
\e{0}%
\e{1}%
\e{2}%
\e{5}%
\e{3}%
\e{3}%
\e{5}%
\e{2}%
\e{1}%
\e{1}%
\e{0}%
\e{0}%
\eol}\vss}\rg%
%
%
\rx{\vss\hfull{%
\rlx{\hss{$2240_{x}$}}\cg%
\e{2}%
\e{4}%
\e{9}%
\e{5}%
\e{3}%
\e{7}%
\e{2}%
\e{2}%
\e{1}%
\e{0}%
\e{0}%
\e{0}%
\e{1}%
\e{3}%
\e{1}%
\e{0}%
\e{2}%
\e{1}%
\e{1}%
\e{0}%
\e{0}%
\e{0}%
\eol}\vss}\rg%
%
%
\rx{\vss\hfull{%
\rlx{\hss{$2835_{x}$}}\cg%
\e{1}%
\e{3}%
\e{8}%
\e{5}%
\e{5}%
\e{12}%
\e{4}%
\e{4}%
\e{4}%
\e{1}%
\e{0}%
\e{0}%
\e{1}%
\e{2}%
\e{1}%
\e{2}%
\e{5}%
\e{1}%
\e{3}%
\e{2}%
\e{1}%
\e{0}%
\eol}\vss}\rg%
%
%
\rx{\vss\hfull{%
\rlx{\hss{$6075_{x}$}}\cg%
\e{2}%
\e{11}%
\e{16}%
\e{17}%
\e{9}%
\e{21}%
\e{10}%
\e{4}%
\e{7}%
\e{2}%
\e{0}%
\e{1}%
\e{5}%
\e{7}%
\e{7}%
\e{4}%
\e{8}%
\e{3}%
\e{1}%
\e{2}%
\e{0}%
\e{0}%
\eol}\vss}\rg%
%
%
\rx{\vss\hfull{%
\rlx{\hss{$3200_{x}$}}\cg%
\e{0}%
\e{5}%
\e{7}%
\e{10}%
\e{4}%
\e{12}%
\e{5}%
\e{2}%
\e{5}%
\e{1}%
\e{0}%
\e{0}%
\e{2}%
\e{3}%
\e{6}%
\e{1}%
\e{4}%
\e{3}%
\e{1}%
\e{1}%
\e{0}%
\e{0}%
\eol}\vss}\rg%
%
%
\rx{\vss\hfull{%
\rlx{\hss{$70_{y}$}}\cg%
\e{0}%
\e{0}%
\e{0}%
\e{1}%
\e{0}%
\e{0}%
\e{0}%
\e{0}%
\e{0}%
\e{0}%
\e{0}%
\e{0}%
\e{1}%
\e{0}%
\e{0}%
\e{0}%
\e{0}%
\e{0}%
\e{0}%
\e{0}%
\e{0}%
\e{0}%
\eol}\vss}\rg%
%
%
\rx{\vss\hfull{%
\rlx{\hss{$1134_{y}$}}\cg%
\e{0}%
\e{1}%
\e{3}%
\e{3}%
\e{1}%
\e{4}%
\e{3}%
\e{1}%
\e{2}%
\e{1}%
\e{0}%
\e{0}%
\e{0}%
\e{2}%
\e{2}%
\e{0}%
\e{2}%
\e{2}%
\e{0}%
\e{0}%
\e{0}%
\e{0}%
\eol}\vss}\rg%
%
%
\rx{\vss\hfull{%
\rlx{\hss{$1680_{y}$}}\cg%
\e{0}%
\e{3}%
\e{3}%
\e{7}%
\e{1}%
\e{5}%
\e{5}%
\e{0}%
\e{2}%
\e{1}%
\e{0}%
\e{1}%
\e{3}%
\e{3}%
\e{4}%
\e{1}%
\e{2}%
\e{1}%
\e{0}%
\e{0}%
\e{0}%
\e{0}%
\eol}\vss}\rg%
%
%
\rx{\vss\hfull{%
\rlx{\hss{$168_{y}$}}\cg%
\e{0}%
\e{0}%
\e{0}%
\e{0}%
\e{1}%
\e{1}%
\e{0}%
\e{0}%
\e{1}%
\e{0}%
\e{0}%
\e{0}%
\e{0}%
\e{0}%
\e{0}%
\e{1}%
\e{0}%
\e{0}%
\e{0}%
\e{1}%
\e{0}%
\e{0}%
\eol}\vss}\rg%
%
%
\rx{\vss\hfull{%
\rlx{\hss{$420_{y}$}}\cg%
\e{0}%
\e{0}%
\e{1}%
\e{0}%
\e{1}%
\e{2}%
\e{1}%
\e{2}%
\e{1}%
\e{0}%
\e{0}%
\e{0}%
\e{0}%
\e{0}%
\e{0}%
\e{0}%
\e{1}%
\e{0}%
\e{1}%
\e{1}%
\e{1}%
\e{0}%
\eol}\vss}\rg%
%
%
\rx{\vss\hfull{%
\rlx{\hss{$3150_{y}$}}\cg%
\e{0}%
\e{2}%
\e{8}%
\e{5}%
\e{3}%
\e{14}%
\e{7}%
\e{6}%
\e{6}%
\e{3}%
\e{0}%
\e{0}%
\e{0}%
\e{3}%
\e{3}%
\e{1}%
\e{6}%
\e{4}%
\e{4}%
\e{3}%
\e{2}%
\e{0}%
\eol}\vss}\rg%
%
%
\rx{\vss\hfull{%
\rlx{\hss{$4200_{y}$}}\cg%
\e{0}%
\e{3}%
\e{8}%
\e{8}%
\e{8}%
\e{18}%
\e{8}%
\e{6}%
\e{11}%
\e{3}%
\e{1}%
\e{0}%
\e{1}%
\e{3}%
\e{3}%
\e{4}%
\e{9}%
\e{4}%
\e{2}%
\e{6}%
\e{2}%
\e{0}%
\eol}\vss}\rg%
\eop
\eject
\tablecont%
%
%
%
%
%
%
\rowpts=18 true pt%
\colpts=18 true pt%
\rowlabpts=40 true pt%
\collabpts=90 true pt%
\clx{\vss\hfull{%
\rlx{\hss{$ $}}\cg%
\cx{\hskip 16 true pt\flip{$[{6}]{\times}[{2}{1}]{\times}[{1^{2}}]$}\hss}\cg%
\cx{\hskip 16 true pt\flip{$[{5}{1}]{\times}[{2}{1}]{\times}[{1^{2}}]$}\hss}\cg%
\cx{\hskip 16 true pt\flip{$[{4}{2}]{\times}[{2}{1}]{\times}[{1^{2}}]$}\hss}\cg%
\cx{\hskip 16 true pt\flip{$[{4}{1^{2}}]{\times}[{2}{1}]{\times}[{1^{2}}]$}\hss}\cg%
\cx{\hskip 16 true pt\flip{$[{3^{2}}]{\times}[{2}{1}]{\times}[{1^{2}}]$}\hss}\cg%
\cx{\hskip 16 true pt\flip{$[{3}{2}{1}]{\times}[{2}{1}]{\times}[{1^{2}}]$}\hss}\cg%
\cx{\hskip 16 true pt\flip{$[{3}{1^{3}}]{\times}[{2}{1}]{\times}[{1^{2}}]$}\hss}\cg%
\cx{\hskip 16 true pt\flip{$[{2^{3}}]{\times}[{2}{1}]{\times}[{1^{2}}]$}\hss}\cg%
\cx{\hskip 16 true pt\flip{$[{2^{2}}{1^{2}}]{\times}[{2}{1}]{\times}[{1^{2}}]$}\hss}\cg%
\cx{\hskip 16 true pt\flip{$[{2}{1^{4}}]{\times}[{2}{1}]{\times}[{1^{2}}]$}\hss}\cg%
\cx{\hskip 16 true pt\flip{$[{1^{6}}]{\times}[{2}{1}]{\times}[{1^{2}}]$}\hss}\cg%
\cx{\hskip 16 true pt\flip{$[{6}]{\times}[{1^{3}}]{\times}[{1^{2}}]$}\hss}\cg%
\cx{\hskip 16 true pt\flip{$[{5}{1}]{\times}[{1^{3}}]{\times}[{1^{2}}]$}\hss}\cg%
\cx{\hskip 16 true pt\flip{$[{4}{2}]{\times}[{1^{3}}]{\times}[{1^{2}}]$}\hss}\cg%
\cx{\hskip 16 true pt\flip{$[{4}{1^{2}}]{\times}[{1^{3}}]{\times}[{1^{2}}]$}\hss}\cg%
\cx{\hskip 16 true pt\flip{$[{3^{2}}]{\times}[{1^{3}}]{\times}[{1^{2}}]$}\hss}\cg%
\cx{\hskip 16 true pt\flip{$[{3}{2}{1}]{\times}[{1^{3}}]{\times}[{1^{2}}]$}\hss}\cg%
\cx{\hskip 16 true pt\flip{$[{3}{1^{3}}]{\times}[{1^{3}}]{\times}[{1^{2}}]$}\hss}\cg%
\cx{\hskip 16 true pt\flip{$[{2^{3}}]{\times}[{1^{3}}]{\times}[{1^{2}}]$}\hss}\cg%
\cx{\hskip 16 true pt\flip{$[{2^{2}}{1^{2}}]{\times}[{1^{3}}]{\times}[{1^{2}}]$}\hss}\cg%
\cx{\hskip 16 true pt\flip{$[{2}{1^{4}}]{\times}[{1^{3}}]{\times}[{1^{2}}]$}\hss}\cg%
\cx{\hskip 16 true pt\flip{$[{1^{6}}]{\times}[{1^{3}}]{\times}[{1^{2}}]$}\hss}\cg%
\eol}}\rg%
%
%
\rx{\vss\hfull{%
\rlx{\hss{$2688_{y}$}}\cg%
\e{0}%
\e{3}%
\e{4}%
\e{7}%
\e{5}%
\e{11}%
\e{5}%
\e{2}%
\e{7}%
\e{2}%
\e{1}%
\e{0}%
\e{1}%
\e{2}%
\e{3}%
\e{2}%
\e{6}%
\e{2}%
\e{1}%
\e{3}%
\e{1}%
\e{0}%
\eol}\vss}\rg%
%
%
\rx{\vss\hfull{%
\rlx{\hss{$2100_{y}$}}\cg%
\e{0}%
\e{3}%
\e{3}%
\e{8}%
\e{3}%
\e{7}%
\e{5}%
\e{1}%
\e{4}%
\e{1}%
\e{0}%
\e{0}%
\e{3}%
\e{2}%
\e{5}%
\e{1}%
\e{3}%
\e{2}%
\e{0}%
\e{1}%
\e{0}%
\e{0}%
\eol}\vss}\rg%
%
%
\rx{\vss\hfull{%
\rlx{\hss{$1400_{y}$}}\cg%
\e{0}%
\e{2}%
\e{3}%
\e{4}%
\e{1}%
\e{5}%
\e{4}%
\e{1}%
\e{2}%
\e{1}%
\e{0}%
\e{1}%
\e{2}%
\e{2}%
\e{2}%
\e{2}%
\e{2}%
\e{1}%
\e{0}%
\e{1}%
\e{0}%
\e{0}%
\eol}\vss}\rg%
%
%
\rx{\vss\hfull{%
\rlx{\hss{$4536_{y}$}}\cg%
\e{1}%
\e{5}%
\e{10}%
\e{12}%
\e{4}%
\e{17}%
\e{12}%
\e{5}%
\e{7}%
\e{4}%
\e{0}%
\e{1}%
\e{4}%
\e{6}%
\e{6}%
\e{3}%
\e{8}%
\e{4}%
\e{2}%
\e{3}%
\e{1}%
\e{0}%
\eol}\vss}\rg%
%
%
\rx{\vss\hfull{%
\rlx{\hss{$5670_{y}$}}\cg%
\e{1}%
\e{6}%
\e{13}%
\e{15}%
\e{5}%
\e{21}%
\e{15}%
\e{6}%
\e{9}%
\e{5}%
\e{0}%
\e{1}%
\e{4}%
\e{8}%
\e{8}%
\e{3}%
\e{10}%
\e{6}%
\e{2}%
\e{3}%
\e{1}%
\e{0}%
\eol}\vss}\rg%
%
%
\rx{\vss\hfull{%
\rlx{\hss{$4480_{y}$}}\cg%
\e{1}%
\e{3}%
\e{11}%
\e{9}%
\e{5}%
\e{18}%
\e{11}%
\e{8}%
\e{8}%
\e{4}%
\e{0}%
\e{0}%
\e{2}%
\e{5}%
\e{4}%
\e{2}%
\e{9}%
\e{5}%
\e{3}%
\e{4}%
\e{2}%
\e{0}%
\eol}\vss}\rg%
%
%
\rx{\vss\hfull{%
\rlx{\hss{$8_{z}$}}\cg%
\e{0}%
\e{0}%
\e{0}%
\e{0}%
\e{0}%
\e{0}%
\e{0}%
\e{0}%
\e{0}%
\e{0}%
\e{0}%
\e{0}%
\e{0}%
\e{0}%
\e{0}%
\e{0}%
\e{0}%
\e{0}%
\e{0}%
\e{0}%
\e{0}%
\e{0}%
\eol}\vss}\rg%
%
%
\rx{\vss\hfull{%
\rlx{\hss{$56_{z}$}}\cg%
\e{0}%
\e{1}%
\e{0}%
\e{0}%
\e{0}%
\e{0}%
\e{0}%
\e{0}%
\e{0}%
\e{0}%
\e{0}%
\e{1}%
\e{0}%
\e{0}%
\e{0}%
\e{0}%
\e{0}%
\e{0}%
\e{0}%
\e{0}%
\e{0}%
\e{0}%
\eol}\vss}\rg%
%
%
\rx{\vss\hfull{%
\rlx{\hss{$160_{z}$}}\cg%
\e{1}%
\e{2}%
\e{0}%
\e{0}%
\e{0}%
\e{0}%
\e{0}%
\e{0}%
\e{0}%
\e{0}%
\e{0}%
\e{1}%
\e{0}%
\e{0}%
\e{0}%
\e{0}%
\e{0}%
\e{0}%
\e{0}%
\e{0}%
\e{0}%
\e{0}%
\eol}\vss}\rg%
%
%
\rx{\vss\hfull{%
\rlx{\hss{$112_{z}$}}\cg%
\e{1}%
\e{1}%
\e{0}%
\e{0}%
\e{0}%
\e{0}%
\e{0}%
\e{0}%
\e{0}%
\e{0}%
\e{0}%
\e{0}%
\e{0}%
\e{0}%
\e{0}%
\e{0}%
\e{0}%
\e{0}%
\e{0}%
\e{0}%
\e{0}%
\e{0}%
\eol}\vss}\rg%
%
%
\rx{\vss\hfull{%
\rlx{\hss{$840_{z}$}}\cg%
\e{0}%
\e{3}%
\e{3}%
\e{3}%
\e{0}%
\e{2}%
\e{1}%
\e{0}%
\e{0}%
\e{0}%
\e{0}%
\e{0}%
\e{1}%
\e{1}%
\e{2}%
\e{0}%
\e{0}%
\e{0}%
\e{0}%
\e{0}%
\e{0}%
\e{0}%
\eol}\vss}\rg%
%
%
\rx{\vss\hfull{%
\rlx{\hss{$1296_{z}$}}\cg%
\e{2}%
\e{6}%
\e{4}%
\e{5}%
\e{1}%
\e{2}%
\e{1}%
\e{0}%
\e{0}%
\e{0}%
\e{0}%
\e{2}%
\e{3}%
\e{2}%
\e{1}%
\e{0}%
\e{0}%
\e{0}%
\e{0}%
\e{0}%
\e{0}%
\e{0}%
\eol}\vss}\rg%
%
%
\rx{\vss\hfull{%
\rlx{\hss{$1400_{z}$}}\cg%
\e{3}%
\e{7}%
\e{5}%
\e{4}%
\e{1}%
\e{2}%
\e{0}%
\e{0}%
\e{0}%
\e{0}%
\e{0}%
\e{1}%
\e{2}%
\e{1}%
\e{1}%
\e{0}%
\e{0}%
\e{0}%
\e{0}%
\e{0}%
\e{0}%
\e{0}%
\eol}\vss}\rg%
%
%
\rx{\vss\hfull{%
\rlx{\hss{$1008_{z}$}}\cg%
\e{3}%
\e{6}%
\e{3}%
\e{3}%
\e{1}%
\e{1}%
\e{0}%
\e{0}%
\e{0}%
\e{0}%
\e{0}%
\e{2}%
\e{2}%
\e{1}%
\e{0}%
\e{0}%
\e{0}%
\e{0}%
\e{0}%
\e{0}%
\e{0}%
\e{0}%
\eol}\vss}\rg%
%
%
\rx{\vss\hfull{%
\rlx{\hss{$560_{z}$}}\cg%
\e{3}%
\e{3}%
\e{2}%
\e{1}%
\e{1}%
\e{0}%
\e{0}%
\e{0}%
\e{0}%
\e{0}%
\e{0}%
\e{0}%
\e{1}%
\e{0}%
\e{0}%
\e{0}%
\e{0}%
\e{0}%
\e{0}%
\e{0}%
\e{0}%
\e{0}%
\eol}\vss}\rg%
%
%
\rx{\vss\hfull{%
\rlx{\hss{$1400_{zz}$}}\cg%
\e{1}%
\e{2}%
\e{5}%
\e{2}%
\e{4}%
\e{5}%
\e{1}%
\e{2}%
\e{1}%
\e{0}%
\e{0}%
\e{0}%
\e{0}%
\e{1}%
\e{0}%
\e{1}%
\e{2}%
\e{0}%
\e{0}%
\e{1}%
\e{0}%
\e{0}%
\eol}\vss}\rg%
%
%
\rx{\vss\hfull{%
\rlx{\hss{$4200_{z}$}}\cg%
\e{1}%
\e{5}%
\e{11}%
\e{8}%
\e{8}%
\e{17}%
\e{7}%
\e{4}%
\e{7}%
\e{2}%
\e{0}%
\e{1}%
\e{1}%
\e{5}%
\e{2}%
\e{4}%
\e{7}%
\e{2}%
\e{3}%
\e{3}%
\e{1}%
\e{0}%
\eol}\vss}\rg%
%
%
\rx{\vss\hfull{%
\rlx{\hss{$400_{z}$}}\cg%
\e{1}%
\e{1}%
\e{2}%
\e{0}%
\e{1}%
\e{1}%
\e{0}%
\e{0}%
\e{0}%
\e{0}%
\e{0}%
\e{0}%
\e{0}%
\e{0}%
\e{0}%
\e{1}%
\e{0}%
\e{0}%
\e{0}%
\e{0}%
\e{0}%
\e{0}%
\eol}\vss}\rg%
%
%
\rx{\vss\hfull{%
\rlx{\hss{$3240_{z}$}}\cg%
\e{4}%
\e{10}%
\e{12}%
\e{8}%
\e{5}%
\e{8}%
\e{2}%
\e{1}%
\e{1}%
\e{0}%
\e{0}%
\e{1}%
\e{3}%
\e{3}%
\e{2}%
\e{2}%
\e{2}%
\e{0}%
\e{0}%
\e{0}%
\e{0}%
\e{0}%
\eol}\vss}\rg%
%
%
\rx{\vss\hfull{%
\rlx{\hss{$4536_{z}$}}\cg%
\e{1}%
\e{7}%
\e{13}%
\e{11}%
\e{7}%
\e{17}%
\e{6}%
\e{5}%
\e{5}%
\e{1}%
\e{0}%
\e{0}%
\e{2}%
\e{4}%
\e{5}%
\e{2}%
\e{6}%
\e{3}%
\e{1}%
\e{2}%
\e{0}%
\e{0}%
\eol}\vss}\rg%
%
%
\rx{\vss\hfull{%
\rlx{\hss{$2400_{z}$}}\cg%
\e{1}%
\e{5}%
\e{6}%
\e{9}%
\e{2}%
\e{7}%
\e{5}%
\e{1}%
\e{2}%
\e{1}%
\e{0}%
\e{1}%
\e{4}%
\e{4}%
\e{3}%
\e{1}%
\e{3}%
\e{1}%
\e{0}%
\e{0}%
\e{0}%
\e{0}%
\eol}\vss}\rg%
%
%
\rx{\vss\hfull{%
\rlx{\hss{$3360_{z}$}}\cg%
\e{2}%
\e{7}%
\e{10}%
\e{8}%
\e{5}%
\e{12}%
\e{4}%
\e{2}%
\e{3}%
\e{1}%
\e{0}%
\e{1}%
\e{2}%
\e{4}%
\e{2}%
\e{3}%
\e{4}%
\e{1}%
\e{1}%
\e{1}%
\e{0}%
\e{0}%
\eol}\vss}\rg%
%
%
\rx{\vss\hfull{%
\rlx{\hss{$2800_{z}$}}\cg%
\e{2}%
\e{7}%
\e{8}%
\e{9}%
\e{4}%
\e{8}%
\e{4}%
\e{1}%
\e{2}%
\e{0}%
\e{0}%
\e{2}%
\e{4}%
\e{4}%
\e{3}%
\e{1}%
\e{2}%
\e{1}%
\e{1}%
\e{0}%
\e{0}%
\e{0}%
\eol}\vss}\rg%
%
%
\rx{\vss\hfull{%
\rlx{\hss{$4096_{z}$}}\cg%
\e{3}%
\e{11}%
\e{13}%
\e{13}%
\e{5}%
\e{11}%
\e{5}%
\e{2}%
\e{2}%
\e{0}%
\e{0}%
\e{1}%
\e{5}%
\e{5}%
\e{5}%
\e{1}%
\e{3}%
\e{1}%
\e{0}%
\e{0}%
\e{0}%
\e{0}%
\eol}\vss}\rg%
%
%
\rx{\vss\hfull{%
\rlx{\hss{$5600_{z}$}}\cg%
\e{1}%
\e{10}%
\e{14}%
\e{18}%
\e{5}%
\e{19}%
\e{11}%
\e{4}%
\e{6}%
\e{2}%
\e{0}%
\e{1}%
\e{6}%
\e{7}%
\e{9}%
\e{2}%
\e{7}%
\e{4}%
\e{1}%
\e{1}%
\e{0}%
\e{0}%
\eol}\vss}\rg%
%
%
\rx{\vss\hfull{%
\rlx{\hss{$448_{z}$}}\cg%
\e{0}%
\e{2}%
\e{2}%
\e{1}%
\e{0}%
\e{1}%
\e{0}%
\e{0}%
\e{0}%
\e{0}%
\e{0}%
\e{0}%
\e{0}%
\e{0}%
\e{1}%
\e{0}%
\e{0}%
\e{0}%
\e{0}%
\e{0}%
\e{0}%
\e{0}%
\eol}\vss}\rg%
%
%
\rx{\vss\hfull{%
\rlx{\hss{$448_{w}$}}\cg%
\e{0}%
\e{1}%
\e{1}%
\e{2}%
\e{0}%
\e{1}%
\e{2}%
\e{0}%
\e{0}%
\e{0}%
\e{0}%
\e{0}%
\e{1}%
\e{1}%
\e{2}%
\e{0}%
\e{0}%
\e{0}%
\e{0}%
\e{0}%
\e{0}%
\e{0}%
\eol}\vss}\rg%
%
%
\rx{\vss\hfull{%
\rlx{\hss{$1344_{w}$}}\cg%
\e{0}%
\e{0}%
\e{3}%
\e{2}%
\e{3}%
\e{6}%
\e{2}%
\e{3}%
\e{4}%
\e{1}%
\e{0}%
\e{0}%
\e{0}%
\e{1}%
\e{0}%
\e{1}%
\e{3}%
\e{2}%
\e{1}%
\e{2}%
\e{1}%
\e{0}%
\eol}\vss}\rg%
%
%
\rx{\vss\hfull{%
\rlx{\hss{$5600_{w}$}}\cg%
\e{1}%
\e{6}%
\e{12}%
\e{15}%
\e{6}%
\e{21}%
\e{14}%
\e{5}%
\e{10}%
\e{5}%
\e{0}%
\e{0}%
\e{4}%
\e{7}%
\e{8}%
\e{4}%
\e{10}%
\e{6}%
\e{2}%
\e{3}%
\e{1}%
\e{0}%
\eol}\vss}\rg%
%
%
\rx{\vss\hfull{%
\rlx{\hss{$2016_{w}$}}\cg%
\e{0}%
\e{1}%
\e{4}%
\e{3}%
\e{3}%
\e{10}%
\e{3}%
\e{4}%
\e{5}%
\e{2}%
\e{0}%
\e{0}%
\e{0}%
\e{1}%
\e{1}%
\e{1}%
\e{4}%
\e{2}%
\e{3}%
\e{4}%
\e{1}%
\e{1}%
\eol}\vss}\rg%
%
%
\rx{\vss\hfull{%
\rlx{\hss{$7168_{w}$}}\cg%
\e{1}%
\e{6}%
\e{15}%
\e{16}%
\e{10}%
\e{29}%
\e{16}%
\e{10}%
\e{15}%
\e{6}%
\e{1}%
\e{0}%
\e{3}%
\e{7}%
\e{7}%
\e{4}%
\e{15}%
\e{7}%
\e{4}%
\e{7}%
\e{3}%
\e{0}%
\eol}\vss}\rg%
\tableclose%
%
%
%
%
%
%
\eop
\eject
\tableopen{Induce/restrict matrix for $W({A_{4}}{A_{4}})\,\subset\,W(E_{8})$}%
%
%
%
%
%
%
\rowpts=18 true pt%
\colpts=18 true pt%
\rowlabpts=40 true pt%
\collabpts=65 true pt%
\clx{\vss\hfull{%
\rlx{\hss{$ $}}\cg%
\cx{\hskip 16 true pt\flip{$[{5}]{\times}[{5}]$}\hss}\cg%
\cx{\hskip 16 true pt\flip{$[{4}{1}]{\times}[{5}]$}\hss}\cg%
\cx{\hskip 16 true pt\flip{$[{3}{2}]{\times}[{5}]$}\hss}\cg%
\cx{\hskip 16 true pt\flip{$[{3}{1^{2}}]{\times}[{5}]$}\hss}\cg%
\cx{\hskip 16 true pt\flip{$[{2^{2}}{1}]{\times}[{5}]$}\hss}\cg%
\cx{\hskip 16 true pt\flip{$[{2}{1^{3}}]{\times}[{5}]$}\hss}\cg%
\cx{\hskip 16 true pt\flip{$[{1^{5}}]{\times}[{5}]$}\hss}\cg%
\cx{\hskip 16 true pt\flip{$[{5}]{\times}[{4}{1}]$}\hss}\cg%
\cx{\hskip 16 true pt\flip{$[{4}{1}]{\times}[{4}{1}]$}\hss}\cg%
\cx{\hskip 16 true pt\flip{$[{3}{2}]{\times}[{4}{1}]$}\hss}\cg%
\cx{\hskip 16 true pt\flip{$[{3}{1^{2}}]{\times}[{4}{1}]$}\hss}\cg%
\cx{\hskip 16 true pt\flip{$[{2^{2}}{1}]{\times}[{4}{1}]$}\hss}\cg%
\cx{\hskip 16 true pt\flip{$[{2}{1^{3}}]{\times}[{4}{1}]$}\hss}\cg%
\cx{\hskip 16 true pt\flip{$[{1^{5}}]{\times}[{4}{1}]$}\hss}\cg%
\cx{\hskip 16 true pt\flip{$[{5}]{\times}[{3}{2}]$}\hss}\cg%
\cx{\hskip 16 true pt\flip{$[{4}{1}]{\times}[{3}{2}]$}\hss}\cg%
\cx{\hskip 16 true pt\flip{$[{3}{2}]{\times}[{3}{2}]$}\hss}\cg%
\cx{\hskip 16 true pt\flip{$[{3}{1^{2}}]{\times}[{3}{2}]$}\hss}\cg%
\cx{\hskip 16 true pt\flip{$[{2^{2}}{1}]{\times}[{3}{2}]$}\hss}\cg%
\cx{\hskip 16 true pt\flip{$[{2}{1^{3}}]{\times}[{3}{2}]$}\hss}\cg%
\cx{\hskip 16 true pt\flip{$[{1^{5}}]{\times}[{3}{2}]$}\hss}\cg%
\cx{\hskip 16 true pt\flip{$[{5}]{\times}[{3}{1^{2}}]$}\hss}\cg%
\cx{\hskip 16 true pt\flip{$[{4}{1}]{\times}[{3}{1^{2}}]$}\hss}\cg%
\cx{\hskip 16 true pt\flip{$[{3}{2}]{\times}[{3}{1^{2}}]$}\hss}\cg%
\cx{\hskip 16 true pt\flip{$[{3}{1^{2}}]{\times}[{3}{1^{2}}]$}\hss}\cg%
\eol}}\rg%
%
%
\rx{\vss\hfull{%
\rlx{\hss{$1_{x}$}}\cg%
\e{1}%
\e{0}%
\e{0}%
\e{0}%
\e{0}%
\e{0}%
\e{0}%
\e{0}%
\e{0}%
\e{0}%
\e{0}%
\e{0}%
\e{0}%
\e{0}%
\e{0}%
\e{0}%
\e{0}%
\e{0}%
\e{0}%
\e{0}%
\e{0}%
\e{0}%
\e{0}%
\e{0}%
\e{0}%
\eol}\vss}\rg%
%
%
\rx{\vss\hfull{%
\rlx{\hss{$28_{x}$}}\cg%
\e{0}%
\e{0}%
\e{0}%
\e{1}%
\e{0}%
\e{0}%
\e{0}%
\e{0}%
\e{1}%
\e{0}%
\e{0}%
\e{0}%
\e{0}%
\e{0}%
\e{0}%
\e{0}%
\e{0}%
\e{0}%
\e{0}%
\e{0}%
\e{0}%
\e{1}%
\e{0}%
\e{0}%
\e{0}%
\eol}\vss}\rg%
%
%
\rx{\vss\hfull{%
\rlx{\hss{$35_{x}$}}\cg%
\e{1}%
\e{1}%
\e{1}%
\e{0}%
\e{0}%
\e{0}%
\e{0}%
\e{1}%
\e{1}%
\e{0}%
\e{0}%
\e{0}%
\e{0}%
\e{0}%
\e{1}%
\e{0}%
\e{0}%
\e{0}%
\e{0}%
\e{0}%
\e{0}%
\e{0}%
\e{0}%
\e{0}%
\e{0}%
\eol}\vss}\rg%
%
%
\rx{\vss\hfull{%
\rlx{\hss{$84_{x}$}}\cg%
\e{2}%
\e{2}%
\e{1}%
\e{0}%
\e{0}%
\e{0}%
\e{0}%
\e{2}%
\e{1}%
\e{1}%
\e{0}%
\e{0}%
\e{0}%
\e{0}%
\e{1}%
\e{1}%
\e{0}%
\e{0}%
\e{0}%
\e{0}%
\e{0}%
\e{0}%
\e{0}%
\e{0}%
\e{0}%
\eol}\vss}\rg%
%
%
\rx{\vss\hfull{%
\rlx{\hss{$50_{x}$}}\cg%
\e{1}%
\e{1}%
\e{0}%
\e{0}%
\e{0}%
\e{0}%
\e{0}%
\e{1}%
\e{1}%
\e{0}%
\e{0}%
\e{0}%
\e{0}%
\e{0}%
\e{0}%
\e{0}%
\e{1}%
\e{0}%
\e{0}%
\e{0}%
\e{0}%
\e{0}%
\e{0}%
\e{0}%
\e{0}%
\eol}\vss}\rg%
%
%
\rx{\vss\hfull{%
\rlx{\hss{$350_{x}$}}\cg%
\e{0}%
\e{0}%
\e{0}%
\e{1}%
\e{1}%
\e{1}%
\e{0}%
\e{0}%
\e{1}%
\e{1}%
\e{2}%
\e{1}%
\e{1}%
\e{0}%
\e{0}%
\e{1}%
\e{0}%
\e{1}%
\e{0}%
\e{0}%
\e{0}%
\e{1}%
\e{2}%
\e{1}%
\e{1}%
\eol}\vss}\rg%
%
%
\rx{\vss\hfull{%
\rlx{\hss{$300_{x}$}}\cg%
\e{1}%
\e{1}%
\e{2}%
\e{0}%
\e{1}%
\e{0}%
\e{0}%
\e{1}%
\e{2}%
\e{2}%
\e{1}%
\e{1}%
\e{0}%
\e{0}%
\e{2}%
\e{2}%
\e{1}%
\e{0}%
\e{0}%
\e{0}%
\e{0}%
\e{0}%
\e{1}%
\e{0}%
\e{1}%
\eol}\vss}\rg%
%
%
\rx{\vss\hfull{%
\rlx{\hss{$567_{x}$}}\cg%
\e{0}%
\e{3}%
\e{2}%
\e{3}%
\e{1}%
\e{1}%
\e{0}%
\e{3}%
\e{5}%
\e{3}%
\e{3}%
\e{1}%
\e{0}%
\e{0}%
\e{2}%
\e{3}%
\e{1}%
\e{1}%
\e{0}%
\e{0}%
\e{0}%
\e{3}%
\e{3}%
\e{1}%
\e{0}%
\eol}\vss}\rg%
%
%
\rx{\vss\hfull{%
\rlx{\hss{$210_{x}$}}\cg%
\e{1}%
\e{2}%
\e{1}%
\e{1}%
\e{1}%
\e{0}%
\e{0}%
\e{2}%
\e{3}%
\e{1}%
\e{1}%
\e{0}%
\e{0}%
\e{0}%
\e{1}%
\e{1}%
\e{1}%
\e{0}%
\e{0}%
\e{0}%
\e{0}%
\e{1}%
\e{1}%
\e{0}%
\e{0}%
\eol}\vss}\rg%
%
%
\rx{\vss\hfull{%
\rlx{\hss{$840_{x}$}}\cg%
\e{1}%
\e{1}%
\e{1}%
\e{0}%
\e{0}%
\e{0}%
\e{0}%
\e{1}%
\e{2}%
\e{3}%
\e{1}%
\e{2}%
\e{0}%
\e{0}%
\e{1}%
\e{3}%
\e{4}%
\e{1}%
\e{2}%
\e{0}%
\e{0}%
\e{0}%
\e{1}%
\e{1}%
\e{2}%
\eol}\vss}\rg%
%
%
\rx{\vss\hfull{%
\rlx{\hss{$700_{x}$}}\cg%
\e{3}%
\e{4}%
\e{3}%
\e{1}%
\e{1}%
\e{0}%
\e{0}%
\e{4}%
\e{6}%
\e{4}%
\e{2}%
\e{1}%
\e{0}%
\e{0}%
\e{3}%
\e{4}%
\e{3}%
\e{1}%
\e{1}%
\e{0}%
\e{0}%
\e{1}%
\e{2}%
\e{1}%
\e{1}%
\eol}\vss}\rg%
%
%
\rx{\vss\hfull{%
\rlx{\hss{$175_{x}$}}\cg%
\e{1}%
\e{1}%
\e{0}%
\e{0}%
\e{0}%
\e{0}%
\e{0}%
\e{1}%
\e{1}%
\e{1}%
\e{0}%
\e{0}%
\e{0}%
\e{0}%
\e{0}%
\e{1}%
\e{1}%
\e{1}%
\e{0}%
\e{0}%
\e{0}%
\e{0}%
\e{0}%
\e{1}%
\e{0}%
\eol}\vss}\rg%
%
%
\rx{\vss\hfull{%
\rlx{\hss{$1400_{x}$}}\cg%
\e{1}%
\e{3}%
\e{2}%
\e{3}%
\e{1}%
\e{1}%
\e{0}%
\e{3}%
\e{7}%
\e{5}%
\e{5}%
\e{2}%
\e{1}%
\e{0}%
\e{2}%
\e{5}%
\e{4}%
\e{4}%
\e{2}%
\e{1}%
\e{0}%
\e{3}%
\e{5}%
\e{4}%
\e{2}%
\eol}\vss}\rg%
%
%
\rx{\vss\hfull{%
\rlx{\hss{$1050_{x}$}}\cg%
\e{1}%
\e{3}%
\e{2}%
\e{2}%
\e{0}%
\e{0}%
\e{0}%
\e{3}%
\e{6}%
\e{4}%
\e{3}%
\e{1}%
\e{0}%
\e{0}%
\e{2}%
\e{4}%
\e{4}%
\e{3}%
\e{2}%
\e{1}%
\e{0}%
\e{2}%
\e{3}%
\e{3}%
\e{1}%
\eol}\vss}\rg%
%
%
\rx{\vss\hfull{%
\rlx{\hss{$1575_{x}$}}\cg%
\e{0}%
\e{2}%
\e{2}%
\e{4}%
\e{2}%
\e{2}%
\e{0}%
\e{2}%
\e{7}%
\e{5}%
\e{7}%
\e{3}%
\e{2}%
\e{0}%
\e{2}%
\e{5}%
\e{3}%
\e{4}%
\e{2}%
\e{1}%
\e{0}%
\e{4}%
\e{7}%
\e{4}%
\e{3}%
\eol}\vss}\rg%
%
%
\rx{\vss\hfull{%
\rlx{\hss{$1344_{x}$}}\cg%
\e{2}%
\e{4}%
\e{4}%
\e{3}%
\e{2}%
\e{0}%
\e{0}%
\e{4}%
\e{8}%
\e{7}%
\e{4}%
\e{3}%
\e{1}%
\e{0}%
\e{4}%
\e{7}%
\e{4}%
\e{3}%
\e{1}%
\e{0}%
\e{0}%
\e{3}%
\e{4}%
\e{3}%
\e{2}%
\eol}\vss}\rg%
%
%
\rx{\vss\hfull{%
\rlx{\hss{$2100_{x}$}}\cg%
\e{0}%
\e{0}%
\e{1}%
\e{2}%
\e{2}%
\e{3}%
\e{1}%
\e{0}%
\e{3}%
\e{4}%
\e{7}%
\e{5}%
\e{5}%
\e{1}%
\e{1}%
\e{4}%
\e{3}%
\e{5}%
\e{2}%
\e{2}%
\e{0}%
\e{2}%
\e{7}%
\e{5}%
\e{7}%
\eol}\vss}\rg%
%
%
\rx{\vss\hfull{%
\rlx{\hss{$2268_{x}$}}\cg%
\e{1}%
\e{2}%
\e{3}%
\e{3}%
\e{3}%
\e{2}%
\e{1}%
\e{2}%
\e{8}%
\e{7}%
\e{8}%
\e{5}%
\e{3}%
\e{0}%
\e{3}%
\e{7}%
\e{6}%
\e{5}%
\e{3}%
\e{1}%
\e{0}%
\e{3}%
\e{8}%
\e{5}%
\e{6}%
\eol}\vss}\rg%
%
%
\rx{\vss\hfull{%
\rlx{\hss{$525_{x}$}}\cg%
\e{0}%
\e{1}%
\e{0}%
\e{2}%
\e{0}%
\e{1}%
\e{0}%
\e{1}%
\e{2}%
\e{2}%
\e{2}%
\e{1}%
\e{1}%
\e{0}%
\e{0}%
\e{2}%
\e{1}%
\e{2}%
\e{0}%
\e{0}%
\e{0}%
\e{2}%
\e{2}%
\e{2}%
\e{0}%
\eol}\vss}\rg%
%
%
\rx{\vss\hfull{%
\rlx{\hss{$700_{xx}$}}\cg%
\e{0}%
\e{1}%
\e{0}%
\e{1}%
\e{0}%
\e{0}%
\e{0}%
\e{1}%
\e{3}%
\e{1}%
\e{2}%
\e{0}%
\e{0}%
\e{0}%
\e{0}%
\e{1}%
\e{3}%
\e{2}%
\e{2}%
\e{1}%
\e{0}%
\e{1}%
\e{2}%
\e{2}%
\e{1}%
\eol}\vss}\rg%
%
%
\rx{\vss\hfull{%
\rlx{\hss{$972_{x}$}}\cg%
\e{1}%
\e{2}%
\e{2}%
\e{0}%
\e{1}%
\e{0}%
\e{0}%
\e{2}%
\e{4}%
\e{4}%
\e{2}%
\e{2}%
\e{0}%
\e{0}%
\e{2}%
\e{4}%
\e{4}%
\e{2}%
\e{2}%
\e{0}%
\e{0}%
\e{0}%
\e{2}%
\e{2}%
\e{2}%
\eol}\vss}\rg%
%
%
\rx{\vss\hfull{%
\rlx{\hss{$4096_{x}$}}\cg%
\e{0}%
\e{3}%
\e{4}%
\e{4}%
\e{4}%
\e{2}%
\e{0}%
\e{3}%
\e{10}%
\e{11}%
\e{12}%
\e{9}%
\e{5}%
\e{1}%
\e{4}%
\e{11}%
\e{10}%
\e{10}%
\e{6}%
\e{3}%
\e{0}%
\e{4}%
\e{12}%
\e{10}%
\e{10}%
\eol}\vss}\rg%
%
%
\rx{\vss\hfull{%
\rlx{\hss{$4200_{x}$}}\cg%
\e{1}%
\e{3}%
\e{4}%
\e{3}%
\e{2}%
\e{1}%
\e{0}%
\e{3}%
\e{10}%
\e{11}%
\e{10}%
\e{7}%
\e{3}%
\e{0}%
\e{4}%
\e{11}%
\e{12}%
\e{10}%
\e{8}%
\e{4}%
\e{1}%
\e{3}%
\e{10}%
\e{10}%
\e{10}%
\eol}\vss}\rg%
%
%
\rx{\vss\hfull{%
\rlx{\hss{$2240_{x}$}}\cg%
\e{2}%
\e{4}%
\e{3}%
\e{2}%
\e{1}%
\e{0}%
\e{0}%
\e{4}%
\e{8}%
\e{8}%
\e{5}%
\e{4}%
\e{1}%
\e{0}%
\e{3}%
\e{8}%
\e{8}%
\e{6}%
\e{4}%
\e{1}%
\e{0}%
\e{2}%
\e{5}%
\e{6}%
\e{4}%
\eol}\vss}\rg%
%
%
\rx{\vss\hfull{%
\rlx{\hss{$2835_{x}$}}\cg%
\e{1}%
\e{2}%
\e{2}%
\e{1}%
\e{1}%
\e{0}%
\e{0}%
\e{2}%
\e{6}%
\e{7}%
\e{5}%
\e{3}%
\e{1}%
\e{0}%
\e{2}%
\e{7}%
\e{8}%
\e{8}%
\e{6}%
\e{3}%
\e{0}%
\e{1}%
\e{5}%
\e{8}%
\e{6}%
\eol}\vss}\rg%
%
%
\rx{\vss\hfull{%
\rlx{\hss{$6075_{x}$}}\cg%
\e{0}%
\e{2}%
\e{2}%
\e{5}%
\e{3}%
\e{3}%
\e{0}%
\e{2}%
\e{10}%
\e{11}%
\e{15}%
\e{10}%
\e{8}%
\e{1}%
\e{2}%
\e{11}%
\e{12}%
\e{16}%
\e{10}%
\e{8}%
\e{1}%
\e{5}%
\e{15}%
\e{16}%
\e{15}%
\eol}\vss}\rg%
%
%
\rx{\vss\hfull{%
\rlx{\hss{$3200_{x}$}}\cg%
\e{0}%
\e{1}%
\e{2}%
\e{1}%
\e{2}%
\e{0}%
\e{0}%
\e{1}%
\e{4}%
\e{7}%
\e{5}%
\e{7}%
\e{3}%
\e{1}%
\e{2}%
\e{7}%
\e{8}%
\e{7}%
\e{6}%
\e{2}%
\e{0}%
\e{1}%
\e{5}%
\e{7}%
\e{8}%
\eol}\vss}\rg%
%
%
\rx{\vss\hfull{%
\rlx{\hss{$70_{y}$}}\cg%
\e{0}%
\e{0}%
\e{0}%
\e{0}%
\e{0}%
\e{0}%
\e{1}%
\e{0}%
\e{0}%
\e{0}%
\e{0}%
\e{0}%
\e{1}%
\e{0}%
\e{0}%
\e{0}%
\e{0}%
\e{0}%
\e{0}%
\e{0}%
\e{0}%
\e{0}%
\e{0}%
\e{0}%
\e{1}%
\eol}\vss}\rg%
%
%
\rx{\vss\hfull{%
\rlx{\hss{$1134_{y}$}}\cg%
\e{0}%
\e{0}%
\e{0}%
\e{1}%
\e{0}%
\e{0}%
\e{0}%
\e{0}%
\e{1}%
\e{1}%
\e{2}%
\e{2}%
\e{2}%
\e{0}%
\e{0}%
\e{1}%
\e{1}%
\e{3}%
\e{2}%
\e{2}%
\e{0}%
\e{1}%
\e{2}%
\e{3}%
\e{2}%
\eol}\vss}\rg%
%
%
\rx{\vss\hfull{%
\rlx{\hss{$1680_{y}$}}\cg%
\e{0}%
\e{0}%
\e{0}%
\e{0}%
\e{1}%
\e{2}%
\e{1}%
\e{0}%
\e{0}%
\e{1}%
\e{4}%
\e{3}%
\e{5}%
\e{2}%
\e{0}%
\e{1}%
\e{1}%
\e{3}%
\e{2}%
\e{3}%
\e{1}%
\e{0}%
\e{4}%
\e{3}%
\e{7}%
\eol}\vss}\rg%
%
%
\rx{\vss\hfull{%
\rlx{\hss{$168_{y}$}}\cg%
\e{0}%
\e{0}%
\e{0}%
\e{0}%
\e{0}%
\e{0}%
\e{0}%
\e{0}%
\e{1}%
\e{0}%
\e{0}%
\e{0}%
\e{0}%
\e{0}%
\e{0}%
\e{0}%
\e{1}%
\e{0}%
\e{1}%
\e{0}%
\e{0}%
\e{0}%
\e{0}%
\e{0}%
\e{1}%
\eol}\vss}\rg%
%
%
\rx{\vss\hfull{%
\rlx{\hss{$420_{y}$}}\cg%
\e{1}%
\e{0}%
\e{0}%
\e{0}%
\e{0}%
\e{0}%
\e{0}%
\e{0}%
\e{1}%
\e{1}%
\e{0}%
\e{0}%
\e{0}%
\e{0}%
\e{0}%
\e{1}%
\e{2}%
\e{1}%
\e{1}%
\e{0}%
\e{0}%
\e{0}%
\e{0}%
\e{1}%
\e{1}%
\eol}\vss}\rg%
%
%
\rx{\vss\hfull{%
\rlx{\hss{$3150_{y}$}}\cg%
\e{0}%
\e{1}%
\e{1}%
\e{1}%
\e{0}%
\e{0}%
\e{0}%
\e{1}%
\e{3}%
\e{5}%
\e{4}%
\e{4}%
\e{2}%
\e{0}%
\e{1}%
\e{5}%
\e{7}%
\e{8}%
\e{6}%
\e{4}%
\e{0}%
\e{1}%
\e{4}%
\e{8}%
\e{6}%
\eol}\vss}\rg%
%
%
\rx{\vss\hfull{%
\rlx{\hss{$4200_{y}$}}\cg%
\e{0}%
\e{1}%
\e{1}%
\e{1}%
\e{1}%
\e{0}%
\e{0}%
\e{1}%
\e{5}%
\e{6}%
\e{6}%
\e{5}%
\e{2}%
\e{0}%
\e{1}%
\e{6}%
\e{10}%
\e{9}%
\e{10}%
\e{5}%
\e{1}%
\e{1}%
\e{6}%
\e{9}%
\e{10}%
\eol}\vss}\rg%
\eop
\eject
\tablecont%
%
%
%
%
%
%
\rowpts=18 true pt%
\colpts=18 true pt%
\rowlabpts=40 true pt%
\collabpts=65 true pt%
\clx{\vss\hfull{%
\rlx{\hss{$ $}}\cg%
\cx{\hskip 16 true pt\flip{$[{5}]{\times}[{5}]$}\hss}\cg%
\cx{\hskip 16 true pt\flip{$[{4}{1}]{\times}[{5}]$}\hss}\cg%
\cx{\hskip 16 true pt\flip{$[{3}{2}]{\times}[{5}]$}\hss}\cg%
\cx{\hskip 16 true pt\flip{$[{3}{1^{2}}]{\times}[{5}]$}\hss}\cg%
\cx{\hskip 16 true pt\flip{$[{2^{2}}{1}]{\times}[{5}]$}\hss}\cg%
\cx{\hskip 16 true pt\flip{$[{2}{1^{3}}]{\times}[{5}]$}\hss}\cg%
\cx{\hskip 16 true pt\flip{$[{1^{5}}]{\times}[{5}]$}\hss}\cg%
\cx{\hskip 16 true pt\flip{$[{5}]{\times}[{4}{1}]$}\hss}\cg%
\cx{\hskip 16 true pt\flip{$[{4}{1}]{\times}[{4}{1}]$}\hss}\cg%
\cx{\hskip 16 true pt\flip{$[{3}{2}]{\times}[{4}{1}]$}\hss}\cg%
\cx{\hskip 16 true pt\flip{$[{3}{1^{2}}]{\times}[{4}{1}]$}\hss}\cg%
\cx{\hskip 16 true pt\flip{$[{2^{2}}{1}]{\times}[{4}{1}]$}\hss}\cg%
\cx{\hskip 16 true pt\flip{$[{2}{1^{3}}]{\times}[{4}{1}]$}\hss}\cg%
\cx{\hskip 16 true pt\flip{$[{1^{5}}]{\times}[{4}{1}]$}\hss}\cg%
\cx{\hskip 16 true pt\flip{$[{5}]{\times}[{3}{2}]$}\hss}\cg%
\cx{\hskip 16 true pt\flip{$[{4}{1}]{\times}[{3}{2}]$}\hss}\cg%
\cx{\hskip 16 true pt\flip{$[{3}{2}]{\times}[{3}{2}]$}\hss}\cg%
\cx{\hskip 16 true pt\flip{$[{3}{1^{2}}]{\times}[{3}{2}]$}\hss}\cg%
\cx{\hskip 16 true pt\flip{$[{2^{2}}{1}]{\times}[{3}{2}]$}\hss}\cg%
\cx{\hskip 16 true pt\flip{$[{2}{1^{3}}]{\times}[{3}{2}]$}\hss}\cg%
\cx{\hskip 16 true pt\flip{$[{1^{5}}]{\times}[{3}{2}]$}\hss}\cg%
\cx{\hskip 16 true pt\flip{$[{5}]{\times}[{3}{1^{2}}]$}\hss}\cg%
\cx{\hskip 16 true pt\flip{$[{4}{1}]{\times}[{3}{1^{2}}]$}\hss}\cg%
\cx{\hskip 16 true pt\flip{$[{3}{2}]{\times}[{3}{1^{2}}]$}\hss}\cg%
\cx{\hskip 16 true pt\flip{$[{3}{1^{2}}]{\times}[{3}{1^{2}}]$}\hss}\cg%
\eol}}\rg%
%
%
\rx{\vss\hfull{%
\rlx{\hss{$2688_{y}$}}\cg%
\e{0}%
\e{0}%
\e{1}%
\e{1}%
\e{1}%
\e{0}%
\e{0}%
\e{0}%
\e{2}%
\e{4}%
\e{4}%
\e{4}%
\e{2}%
\e{0}%
\e{1}%
\e{4}%
\e{5}%
\e{6}%
\e{5}%
\e{4}%
\e{1}%
\e{1}%
\e{4}%
\e{6}%
\e{7}%
\eol}\vss}\rg%
%
%
\rx{\vss\hfull{%
\rlx{\hss{$2100_{y}$}}\cg%
\e{0}%
\e{0}%
\e{0}%
\e{1}%
\e{0}%
\e{2}%
\e{1}%
\e{0}%
\e{1}%
\e{2}%
\e{4}%
\e{3}%
\e{5}%
\e{2}%
\e{0}%
\e{2}%
\e{2}%
\e{5}%
\e{3}%
\e{3}%
\e{0}%
\e{1}%
\e{4}%
\e{5}%
\e{6}%
\eol}\vss}\rg%
%
%
\rx{\vss\hfull{%
\rlx{\hss{$1400_{y}$}}\cg%
\e{0}%
\e{0}%
\e{0}%
\e{0}%
\e{1}%
\e{1}%
\e{1}%
\e{0}%
\e{1}%
\e{1}%
\e{3}%
\e{2}%
\e{3}%
\e{1}%
\e{0}%
\e{1}%
\e{2}%
\e{2}%
\e{3}%
\e{2}%
\e{1}%
\e{0}%
\e{3}%
\e{2}%
\e{6}%
\eol}\vss}\rg%
%
%
\rx{\vss\hfull{%
\rlx{\hss{$4536_{y}$}}\cg%
\e{0}%
\e{0}%
\e{1}%
\e{1}%
\e{2}%
\e{2}%
\e{1}%
\e{0}%
\e{3}%
\e{5}%
\e{8}%
\e{7}%
\e{7}%
\e{2}%
\e{1}%
\e{5}%
\e{7}%
\e{9}%
\e{8}%
\e{7}%
\e{2}%
\e{1}%
\e{8}%
\e{9}%
\e{15}%
\eol}\vss}\rg%
%
%
\rx{\vss\hfull{%
\rlx{\hss{$5670_{y}$}}\cg%
\e{0}%
\e{0}%
\e{1}%
\e{2}%
\e{2}%
\e{2}%
\e{1}%
\e{0}%
\e{4}%
\e{6}%
\e{10}%
\e{9}%
\e{9}%
\e{2}%
\e{1}%
\e{6}%
\e{8}%
\e{12}%
\e{10}%
\e{9}%
\e{2}%
\e{2}%
\e{10}%
\e{12}%
\e{17}%
\eol}\vss}\rg%
%
%
\rx{\vss\hfull{%
\rlx{\hss{$4480_{y}$}}\cg%
\e{0}%
\e{1}%
\e{1}%
\e{1}%
\e{1}%
\e{1}%
\e{0}%
\e{1}%
\e{4}%
\e{6}%
\e{7}%
\e{6}%
\e{4}%
\e{1}%
\e{1}%
\e{6}%
\e{9}%
\e{10}%
\e{9}%
\e{6}%
\e{1}%
\e{1}%
\e{7}%
\e{10}%
\e{11}%
\eol}\vss}\rg%
%
%
\rx{\vss\hfull{%
\rlx{\hss{$8_{z}$}}\cg%
\e{0}%
\e{1}%
\e{0}%
\e{0}%
\e{0}%
\e{0}%
\e{0}%
\e{1}%
\e{0}%
\e{0}%
\e{0}%
\e{0}%
\e{0}%
\e{0}%
\e{0}%
\e{0}%
\e{0}%
\e{0}%
\e{0}%
\e{0}%
\e{0}%
\e{0}%
\e{0}%
\e{0}%
\e{0}%
\eol}\vss}\rg%
%
%
\rx{\vss\hfull{%
\rlx{\hss{$56_{z}$}}\cg%
\e{0}%
\e{0}%
\e{0}%
\e{0}%
\e{0}%
\e{1}%
\e{0}%
\e{0}%
\e{0}%
\e{0}%
\e{1}%
\e{0}%
\e{0}%
\e{0}%
\e{0}%
\e{0}%
\e{0}%
\e{0}%
\e{0}%
\e{0}%
\e{0}%
\e{0}%
\e{1}%
\e{0}%
\e{0}%
\eol}\vss}\rg%
%
%
\rx{\vss\hfull{%
\rlx{\hss{$160_{z}$}}\cg%
\e{0}%
\e{1}%
\e{1}%
\e{1}%
\e{1}%
\e{0}%
\e{0}%
\e{1}%
\e{2}%
\e{1}%
\e{1}%
\e{0}%
\e{0}%
\e{0}%
\e{1}%
\e{1}%
\e{0}%
\e{0}%
\e{0}%
\e{0}%
\e{0}%
\e{1}%
\e{1}%
\e{0}%
\e{0}%
\eol}\vss}\rg%
%
%
\rx{\vss\hfull{%
\rlx{\hss{$112_{z}$}}\cg%
\e{2}%
\e{2}%
\e{1}%
\e{1}%
\e{0}%
\e{0}%
\e{0}%
\e{2}%
\e{2}%
\e{1}%
\e{0}%
\e{0}%
\e{0}%
\e{0}%
\e{1}%
\e{1}%
\e{0}%
\e{0}%
\e{0}%
\e{0}%
\e{0}%
\e{1}%
\e{0}%
\e{0}%
\e{0}%
\eol}\vss}\rg%
%
%
\rx{\vss\hfull{%
\rlx{\hss{$840_{z}$}}\cg%
\e{0}%
\e{1}%
\e{1}%
\e{1}%
\e{1}%
\e{0}%
\e{0}%
\e{1}%
\e{2}%
\e{3}%
\e{2}%
\e{3}%
\e{1}%
\e{0}%
\e{1}%
\e{3}%
\e{2}%
\e{2}%
\e{1}%
\e{0}%
\e{0}%
\e{1}%
\e{2}%
\e{2}%
\e{2}%
\eol}\vss}\rg%
%
%
\rx{\vss\hfull{%
\rlx{\hss{$1296_{z}$}}\cg%
\e{0}%
\e{0}%
\e{1}%
\e{3}%
\e{2}%
\e{2}%
\e{1}%
\e{0}%
\e{4}%
\e{3}%
\e{6}%
\e{3}%
\e{3}%
\e{0}%
\e{1}%
\e{3}%
\e{2}%
\e{3}%
\e{1}%
\e{1}%
\e{0}%
\e{3}%
\e{6}%
\e{3}%
\e{4}%
\eol}\vss}\rg%
%
%
\rx{\vss\hfull{%
\rlx{\hss{$1400_{z}$}}\cg%
\e{2}%
\e{4}%
\e{4}%
\e{3}%
\e{2}%
\e{1}%
\e{0}%
\e{4}%
\e{8}%
\e{7}%
\e{5}%
\e{3}%
\e{1}%
\e{0}%
\e{4}%
\e{7}%
\e{4}%
\e{3}%
\e{1}%
\e{0}%
\e{0}%
\e{3}%
\e{5}%
\e{3}%
\e{2}%
\eol}\vss}\rg%
%
%
\rx{\vss\hfull{%
\rlx{\hss{$1008_{z}$}}\cg%
\e{0}%
\e{2}%
\e{2}%
\e{3}%
\e{2}%
\e{2}%
\e{0}%
\e{2}%
\e{6}%
\e{4}%
\e{5}%
\e{2}%
\e{1}%
\e{0}%
\e{2}%
\e{4}%
\e{2}%
\e{2}%
\e{1}%
\e{0}%
\e{0}%
\e{3}%
\e{5}%
\e{2}%
\e{2}%
\eol}\vss}\rg%
%
%
\rx{\vss\hfull{%
\rlx{\hss{$560_{z}$}}\cg%
\e{2}%
\e{4}%
\e{3}%
\e{2}%
\e{1}%
\e{0}%
\e{0}%
\e{4}%
\e{6}%
\e{3}%
\e{2}%
\e{1}%
\e{0}%
\e{0}%
\e{3}%
\e{3}%
\e{2}%
\e{1}%
\e{0}%
\e{0}%
\e{0}%
\e{2}%
\e{2}%
\e{1}%
\e{0}%
\eol}\vss}\rg%
%
%
\rx{\vss\hfull{%
\rlx{\hss{$1400_{zz}$}}\cg%
\e{2}%
\e{3}%
\e{2}%
\e{1}%
\e{0}%
\e{0}%
\e{0}%
\e{3}%
\e{6}%
\e{5}%
\e{3}%
\e{1}%
\e{0}%
\e{0}%
\e{2}%
\e{5}%
\e{6}%
\e{4}%
\e{3}%
\e{1}%
\e{0}%
\e{1}%
\e{3}%
\e{4}%
\e{2}%
\eol}\vss}\rg%
%
%
\rx{\vss\hfull{%
\rlx{\hss{$4200_{z}$}}\cg%
\e{0}%
\e{2}%
\e{1}%
\e{3}%
\e{1}%
\e{1}%
\e{0}%
\e{2}%
\e{8}%
\e{7}%
\e{9}%
\e{5}%
\e{3}%
\e{0}%
\e{1}%
\e{7}%
\e{10}%
\e{11}%
\e{9}%
\e{6}%
\e{1}%
\e{3}%
\e{9}%
\e{11}%
\e{10}%
\eol}\vss}\rg%
%
%
\rx{\vss\hfull{%
\rlx{\hss{$400_{z}$}}\cg%
\e{2}%
\e{3}%
\e{1}%
\e{1}%
\e{0}%
\e{0}%
\e{0}%
\e{3}%
\e{4}%
\e{2}%
\e{1}%
\e{0}%
\e{0}%
\e{0}%
\e{1}%
\e{2}%
\e{2}%
\e{1}%
\e{1}%
\e{0}%
\e{0}%
\e{1}%
\e{1}%
\e{1}%
\e{0}%
\eol}\vss}\rg%
%
%
\rx{\vss\hfull{%
\rlx{\hss{$3240_{z}$}}\cg%
\e{2}%
\e{6}%
\e{5}%
\e{5}%
\e{3}%
\e{1}%
\e{0}%
\e{6}%
\e{14}%
\e{12}%
\e{10}%
\e{6}%
\e{2}%
\e{0}%
\e{5}%
\e{12}%
\e{10}%
\e{8}%
\e{5}%
\e{2}%
\e{0}%
\e{5}%
\e{10}%
\e{8}%
\e{6}%
\eol}\vss}\rg%
%
%
\rx{\vss\hfull{%
\rlx{\hss{$4536_{z}$}}\cg%
\e{2}%
\e{3}%
\e{4}%
\e{2}%
\e{2}%
\e{0}%
\e{0}%
\e{3}%
\e{10}%
\e{12}%
\e{8}%
\e{8}%
\e{3}%
\e{0}%
\e{4}%
\e{12}%
\e{14}%
\e{11}%
\e{9}%
\e{3}%
\e{0}%
\e{2}%
\e{8}%
\e{11}%
\e{10}%
\eol}\vss}\rg%
%
%
\rx{\vss\hfull{%
\rlx{\hss{$2400_{z}$}}\cg%
\e{0}%
\e{0}%
\e{0}%
\e{2}%
\e{1}%
\e{3}%
\e{1}%
\e{0}%
\e{2}%
\e{3}%
\e{7}%
\e{4}%
\e{6}%
\e{1}%
\e{0}%
\e{3}%
\e{2}%
\e{6}%
\e{3}%
\e{4}%
\e{1}%
\e{2}%
\e{7}%
\e{6}%
\e{8}%
\eol}\vss}\rg%
%
%
\rx{\vss\hfull{%
\rlx{\hss{$3360_{z}$}}\cg%
\e{0}%
\e{2}%
\e{3}%
\e{3}%
\e{2}%
\e{1}%
\e{0}%
\e{2}%
\e{8}%
\e{8}%
\e{9}%
\e{5}%
\e{3}%
\e{0}%
\e{3}%
\e{8}%
\e{8}%
\e{9}%
\e{6}%
\e{4}%
\e{1}%
\e{3}%
\e{9}%
\e{9}%
\e{8}%
\eol}\vss}\rg%
%
%
\rx{\vss\hfull{%
\rlx{\hss{$2800_{z}$}}\cg%
\e{0}%
\e{1}%
\e{1}%
\e{3}%
\e{2}%
\e{3}%
\e{1}%
\e{1}%
\e{6}%
\e{5}%
\e{9}%
\e{5}%
\e{5}%
\e{1}%
\e{1}%
\e{5}%
\e{6}%
\e{7}%
\e{4}%
\e{3}%
\e{0}%
\e{3}%
\e{9}%
\e{7}%
\e{8}%
\eol}\vss}\rg%
%
%
\rx{\vss\hfull{%
\rlx{\hss{$4096_{z}$}}\cg%
\e{0}%
\e{3}%
\e{4}%
\e{4}%
\e{4}%
\e{2}%
\e{0}%
\e{3}%
\e{10}%
\e{11}%
\e{12}%
\e{9}%
\e{5}%
\e{1}%
\e{4}%
\e{11}%
\e{10}%
\e{10}%
\e{6}%
\e{3}%
\e{0}%
\e{4}%
\e{12}%
\e{10}%
\e{10}%
\eol}\vss}\rg%
%
%
\rx{\vss\hfull{%
\rlx{\hss{$5600_{z}$}}\cg%
\e{0}%
\e{1}%
\e{2}%
\e{3}%
\e{3}%
\e{3}%
\e{1}%
\e{1}%
\e{6}%
\e{10}%
\e{12}%
\e{11}%
\e{9}%
\e{2}%
\e{2}%
\e{10}%
\e{10}%
\e{13}%
\e{9}%
\e{6}%
\e{1}%
\e{3}%
\e{12}%
\e{13}%
\e{16}%
\eol}\vss}\rg%
%
%
\rx{\vss\hfull{%
\rlx{\hss{$448_{z}$}}\cg%
\e{2}%
\e{2}%
\e{2}%
\e{0}%
\e{0}%
\e{0}%
\e{0}%
\e{2}%
\e{2}%
\e{3}%
\e{1}%
\e{1}%
\e{0}%
\e{0}%
\e{2}%
\e{3}%
\e{2}%
\e{1}%
\e{0}%
\e{0}%
\e{0}%
\e{0}%
\e{1}%
\e{1}%
\e{0}%
\eol}\vss}\rg%
%
%
\rx{\vss\hfull{%
\rlx{\hss{$448_{w}$}}\cg%
\e{0}%
\e{0}%
\e{0}%
\e{0}%
\e{0}%
\e{1}%
\e{0}%
\e{0}%
\e{0}%
\e{0}%
\e{1}%
\e{1}%
\e{2}%
\e{1}%
\e{0}%
\e{0}%
\e{0}%
\e{1}%
\e{0}%
\e{1}%
\e{0}%
\e{0}%
\e{1}%
\e{1}%
\e{2}%
\eol}\vss}\rg%
%
%
\rx{\vss\hfull{%
\rlx{\hss{$1344_{w}$}}\cg%
\e{0}%
\e{1}%
\e{0}%
\e{0}%
\e{0}%
\e{0}%
\e{0}%
\e{1}%
\e{2}%
\e{2}%
\e{2}%
\e{1}%
\e{0}%
\e{0}%
\e{0}%
\e{2}%
\e{4}%
\e{3}%
\e{4}%
\e{1}%
\e{0}%
\e{0}%
\e{2}%
\e{3}%
\e{2}%
\eol}\vss}\rg%
%
%
\rx{\vss\hfull{%
\rlx{\hss{$5600_{w}$}}\cg%
\e{0}%
\e{0}%
\e{1}%
\e{2}%
\e{2}%
\e{2}%
\e{0}%
\e{0}%
\e{4}%
\e{6}%
\e{10}%
\e{9}%
\e{8}%
\e{2}%
\e{1}%
\e{6}%
\e{8}%
\e{12}%
\e{10}%
\e{9}%
\e{2}%
\e{2}%
\e{10}%
\e{12}%
\e{16}%
\eol}\vss}\rg%
%
%
\rx{\vss\hfull{%
\rlx{\hss{$2016_{w}$}}\cg%
\e{0}%
\e{1}%
\e{1}%
\e{0}%
\e{0}%
\e{0}%
\e{0}%
\e{1}%
\e{2}%
\e{4}%
\e{2}%
\e{2}%
\e{0}%
\e{0}%
\e{1}%
\e{4}%
\e{6}%
\e{5}%
\e{4}%
\e{2}%
\e{0}%
\e{0}%
\e{2}%
\e{5}%
\e{4}%
\eol}\vss}\rg%
%
%
\rx{\vss\hfull{%
\rlx{\hss{$7168_{w}$}}\cg%
\e{0}%
\e{1}%
\e{2}%
\e{2}%
\e{2}%
\e{1}%
\e{0}%
\e{1}%
\e{6}%
\e{10}%
\e{11}%
\e{10}%
\e{6}%
\e{1}%
\e{2}%
\e{10}%
\e{14}%
\e{16}%
\e{14}%
\e{10}%
\e{2}%
\e{2}%
\e{11}%
\e{16}%
\e{18}%
\eol}\vss}\rg%
\eop
\eject
\tablecont%
%
%
%
%
%
%
\rowpts=18 true pt%
\colpts=18 true pt%
\rowlabpts=40 true pt%
\collabpts=65 true pt%
\clx{\vss\hfull{%
\rlx{\hss{$ $}}\cg%
\cx{\hskip 16 true pt\flip{$[{2^{2}}{1}]{\times}[{3}{1^{2}}]$}\hss}\cg%
\cx{\hskip 16 true pt\flip{$[{2}{1^{3}}]{\times}[{3}{1^{2}}]$}\hss}\cg%
\cx{\hskip 16 true pt\flip{$[{1^{5}}]{\times}[{3}{1^{2}}]$}\hss}\cg%
\cx{\hskip 16 true pt\flip{$[{5}]{\times}[{2^{2}}{1}]$}\hss}\cg%
\cx{\hskip 16 true pt\flip{$[{4}{1}]{\times}[{2^{2}}{1}]$}\hss}\cg%
\cx{\hskip 16 true pt\flip{$[{3}{2}]{\times}[{2^{2}}{1}]$}\hss}\cg%
\cx{\hskip 16 true pt\flip{$[{3}{1^{2}}]{\times}[{2^{2}}{1}]$}\hss}\cg%
\cx{\hskip 16 true pt\flip{$[{2^{2}}{1}]{\times}[{2^{2}}{1}]$}\hss}\cg%
\cx{\hskip 16 true pt\flip{$[{2}{1^{3}}]{\times}[{2^{2}}{1}]$}\hss}\cg%
\cx{\hskip 16 true pt\flip{$[{1^{5}}]{\times}[{2^{2}}{1}]$}\hss}\cg%
\cx{\hskip 16 true pt\flip{$[{5}]{\times}[{2}{1^{3}}]$}\hss}\cg%
\cx{\hskip 16 true pt\flip{$[{4}{1}]{\times}[{2}{1^{3}}]$}\hss}\cg%
\cx{\hskip 16 true pt\flip{$[{3}{2}]{\times}[{2}{1^{3}}]$}\hss}\cg%
\cx{\hskip 16 true pt\flip{$[{3}{1^{2}}]{\times}[{2}{1^{3}}]$}\hss}\cg%
\cx{\hskip 16 true pt\flip{$[{2^{2}}{1}]{\times}[{2}{1^{3}}]$}\hss}\cg%
\cx{\hskip 16 true pt\flip{$[{2}{1^{3}}]{\times}[{2}{1^{3}}]$}\hss}\cg%
\cx{\hskip 16 true pt\flip{$[{1^{5}}]{\times}[{2}{1^{3}}]$}\hss}\cg%
\cx{\hskip 16 true pt\flip{$[{5}]{\times}[{1^{5}}]$}\hss}\cg%
\cx{\hskip 16 true pt\flip{$[{4}{1}]{\times}[{1^{5}}]$}\hss}\cg%
\cx{\hskip 16 true pt\flip{$[{3}{2}]{\times}[{1^{5}}]$}\hss}\cg%
\cx{\hskip 16 true pt\flip{$[{3}{1^{2}}]{\times}[{1^{5}}]$}\hss}\cg%
\cx{\hskip 16 true pt\flip{$[{2^{2}}{1}]{\times}[{1^{5}}]$}\hss}\cg%
\cx{\hskip 16 true pt\flip{$[{2}{1^{3}}]{\times}[{1^{5}}]$}\hss}\cg%
\cx{\hskip 16 true pt\flip{$[{1^{5}}]{\times}[{1^{5}}]$}\hss}\cg%
\eol}}\rg%
%
%
\rx{\vss\hfull{%
\rlx{\hss{$1_{x}$}}\cg%
\e{0}%
\e{0}%
\e{0}%
\e{0}%
\e{0}%
\e{0}%
\e{0}%
\e{0}%
\e{0}%
\e{0}%
\e{0}%
\e{0}%
\e{0}%
\e{0}%
\e{0}%
\e{0}%
\e{0}%
\e{0}%
\e{0}%
\e{0}%
\e{0}%
\e{0}%
\e{0}%
\e{0}%
\eol}\vss}\rg%
%
%
\rx{\vss\hfull{%
\rlx{\hss{$28_{x}$}}\cg%
\e{0}%
\e{0}%
\e{0}%
\e{0}%
\e{0}%
\e{0}%
\e{0}%
\e{0}%
\e{0}%
\e{0}%
\e{0}%
\e{0}%
\e{0}%
\e{0}%
\e{0}%
\e{0}%
\e{0}%
\e{0}%
\e{0}%
\e{0}%
\e{0}%
\e{0}%
\e{0}%
\e{0}%
\eol}\vss}\rg%
%
%
\rx{\vss\hfull{%
\rlx{\hss{$35_{x}$}}\cg%
\e{0}%
\e{0}%
\e{0}%
\e{0}%
\e{0}%
\e{0}%
\e{0}%
\e{0}%
\e{0}%
\e{0}%
\e{0}%
\e{0}%
\e{0}%
\e{0}%
\e{0}%
\e{0}%
\e{0}%
\e{0}%
\e{0}%
\e{0}%
\e{0}%
\e{0}%
\e{0}%
\e{0}%
\eol}\vss}\rg%
%
%
\rx{\vss\hfull{%
\rlx{\hss{$84_{x}$}}\cg%
\e{0}%
\e{0}%
\e{0}%
\e{0}%
\e{0}%
\e{0}%
\e{0}%
\e{0}%
\e{0}%
\e{0}%
\e{0}%
\e{0}%
\e{0}%
\e{0}%
\e{0}%
\e{0}%
\e{0}%
\e{0}%
\e{0}%
\e{0}%
\e{0}%
\e{0}%
\e{0}%
\e{0}%
\eol}\vss}\rg%
%
%
\rx{\vss\hfull{%
\rlx{\hss{$50_{x}$}}\cg%
\e{0}%
\e{0}%
\e{0}%
\e{0}%
\e{0}%
\e{0}%
\e{0}%
\e{0}%
\e{0}%
\e{0}%
\e{0}%
\e{0}%
\e{0}%
\e{0}%
\e{0}%
\e{0}%
\e{0}%
\e{0}%
\e{0}%
\e{0}%
\e{0}%
\e{0}%
\e{0}%
\e{0}%
\eol}\vss}\rg%
%
%
\rx{\vss\hfull{%
\rlx{\hss{$350_{x}$}}\cg%
\e{0}%
\e{0}%
\e{0}%
\e{1}%
\e{1}%
\e{0}%
\e{0}%
\e{0}%
\e{0}%
\e{0}%
\e{1}%
\e{1}%
\e{0}%
\e{0}%
\e{0}%
\e{0}%
\e{0}%
\e{0}%
\e{0}%
\e{0}%
\e{0}%
\e{0}%
\e{0}%
\e{0}%
\eol}\vss}\rg%
%
%
\rx{\vss\hfull{%
\rlx{\hss{$300_{x}$}}\cg%
\e{0}%
\e{0}%
\e{0}%
\e{1}%
\e{1}%
\e{0}%
\e{0}%
\e{0}%
\e{0}%
\e{0}%
\e{0}%
\e{0}%
\e{0}%
\e{0}%
\e{0}%
\e{0}%
\e{0}%
\e{0}%
\e{0}%
\e{0}%
\e{0}%
\e{0}%
\e{0}%
\e{0}%
\eol}\vss}\rg%
%
%
\rx{\vss\hfull{%
\rlx{\hss{$567_{x}$}}\cg%
\e{0}%
\e{0}%
\e{0}%
\e{1}%
\e{1}%
\e{0}%
\e{0}%
\e{0}%
\e{0}%
\e{0}%
\e{1}%
\e{0}%
\e{0}%
\e{0}%
\e{0}%
\e{0}%
\e{0}%
\e{0}%
\e{0}%
\e{0}%
\e{0}%
\e{0}%
\e{0}%
\e{0}%
\eol}\vss}\rg%
%
%
\rx{\vss\hfull{%
\rlx{\hss{$210_{x}$}}\cg%
\e{0}%
\e{0}%
\e{0}%
\e{1}%
\e{0}%
\e{0}%
\e{0}%
\e{0}%
\e{0}%
\e{0}%
\e{0}%
\e{0}%
\e{0}%
\e{0}%
\e{0}%
\e{0}%
\e{0}%
\e{0}%
\e{0}%
\e{0}%
\e{0}%
\e{0}%
\e{0}%
\e{0}%
\eol}\vss}\rg%
%
%
\rx{\vss\hfull{%
\rlx{\hss{$840_{x}$}}\cg%
\e{2}%
\e{1}%
\e{0}%
\e{0}%
\e{2}%
\e{2}%
\e{2}%
\e{1}%
\e{0}%
\e{0}%
\e{0}%
\e{0}%
\e{0}%
\e{1}%
\e{0}%
\e{1}%
\e{0}%
\e{0}%
\e{0}%
\e{0}%
\e{0}%
\e{0}%
\e{0}%
\e{0}%
\eol}\vss}\rg%
%
%
\rx{\vss\hfull{%
\rlx{\hss{$700_{x}$}}\cg%
\e{0}%
\e{0}%
\e{0}%
\e{1}%
\e{1}%
\e{1}%
\e{0}%
\e{0}%
\e{0}%
\e{0}%
\e{0}%
\e{0}%
\e{0}%
\e{0}%
\e{0}%
\e{0}%
\e{0}%
\e{0}%
\e{0}%
\e{0}%
\e{0}%
\e{0}%
\e{0}%
\e{0}%
\eol}\vss}\rg%
%
%
\rx{\vss\hfull{%
\rlx{\hss{$175_{x}$}}\cg%
\e{0}%
\e{0}%
\e{0}%
\e{0}%
\e{0}%
\e{0}%
\e{0}%
\e{1}%
\e{0}%
\e{0}%
\e{0}%
\e{0}%
\e{0}%
\e{0}%
\e{0}%
\e{0}%
\e{0}%
\e{0}%
\e{0}%
\e{0}%
\e{0}%
\e{0}%
\e{0}%
\e{0}%
\eol}\vss}\rg%
%
%
\rx{\vss\hfull{%
\rlx{\hss{$1400_{x}$}}\cg%
\e{1}%
\e{0}%
\e{0}%
\e{1}%
\e{2}%
\e{2}%
\e{1}%
\e{1}%
\e{0}%
\e{0}%
\e{1}%
\e{1}%
\e{1}%
\e{0}%
\e{0}%
\e{0}%
\e{0}%
\e{0}%
\e{0}%
\e{0}%
\e{0}%
\e{0}%
\e{0}%
\e{0}%
\eol}\vss}\rg%
%
%
\rx{\vss\hfull{%
\rlx{\hss{$1050_{x}$}}\cg%
\e{1}%
\e{0}%
\e{0}%
\e{0}%
\e{1}%
\e{2}%
\e{1}%
\e{1}%
\e{0}%
\e{0}%
\e{0}%
\e{0}%
\e{1}%
\e{0}%
\e{0}%
\e{0}%
\e{0}%
\e{0}%
\e{0}%
\e{0}%
\e{0}%
\e{0}%
\e{0}%
\e{0}%
\eol}\vss}\rg%
%
%
\rx{\vss\hfull{%
\rlx{\hss{$1575_{x}$}}\cg%
\e{1}%
\e{0}%
\e{0}%
\e{2}%
\e{3}%
\e{2}%
\e{1}%
\e{0}%
\e{0}%
\e{0}%
\e{2}%
\e{2}%
\e{1}%
\e{0}%
\e{0}%
\e{0}%
\e{0}%
\e{0}%
\e{0}%
\e{0}%
\e{0}%
\e{0}%
\e{0}%
\e{0}%
\eol}\vss}\rg%
%
%
\rx{\vss\hfull{%
\rlx{\hss{$1344_{x}$}}\cg%
\e{1}%
\e{0}%
\e{0}%
\e{2}%
\e{3}%
\e{1}%
\e{1}%
\e{0}%
\e{0}%
\e{0}%
\e{0}%
\e{1}%
\e{0}%
\e{0}%
\e{0}%
\e{0}%
\e{0}%
\e{0}%
\e{0}%
\e{0}%
\e{0}%
\e{0}%
\e{0}%
\e{0}%
\eol}\vss}\rg%
%
%
\rx{\vss\hfull{%
\rlx{\hss{$2100_{x}$}}\cg%
\e{3}%
\e{2}%
\e{0}%
\e{2}%
\e{5}%
\e{2}%
\e{3}%
\e{1}%
\e{0}%
\e{0}%
\e{3}%
\e{5}%
\e{2}%
\e{2}%
\e{0}%
\e{0}%
\e{0}%
\e{1}%
\e{1}%
\e{0}%
\e{0}%
\e{0}%
\e{0}%
\e{0}%
\eol}\vss}\rg%
%
%
\rx{\vss\hfull{%
\rlx{\hss{$2268_{x}$}}\cg%
\e{2}%
\e{1}%
\e{0}%
\e{3}%
\e{5}%
\e{3}%
\e{2}%
\e{1}%
\e{0}%
\e{0}%
\e{2}%
\e{3}%
\e{1}%
\e{1}%
\e{0}%
\e{0}%
\e{0}%
\e{1}%
\e{0}%
\e{0}%
\e{0}%
\e{0}%
\e{0}%
\e{0}%
\eol}\vss}\rg%
%
%
\rx{\vss\hfull{%
\rlx{\hss{$525_{x}$}}\cg%
\e{1}%
\e{0}%
\e{0}%
\e{0}%
\e{1}%
\e{0}%
\e{1}%
\e{0}%
\e{0}%
\e{0}%
\e{1}%
\e{1}%
\e{0}%
\e{0}%
\e{0}%
\e{0}%
\e{0}%
\e{0}%
\e{0}%
\e{0}%
\e{0}%
\e{0}%
\e{0}%
\e{0}%
\eol}\vss}\rg%
%
%
\rx{\vss\hfull{%
\rlx{\hss{$700_{xx}$}}\cg%
\e{1}%
\e{0}%
\e{0}%
\e{0}%
\e{0}%
\e{2}%
\e{1}%
\e{1}%
\e{1}%
\e{0}%
\e{0}%
\e{0}%
\e{1}%
\e{0}%
\e{1}%
\e{0}%
\e{0}%
\e{0}%
\e{0}%
\e{0}%
\e{0}%
\e{0}%
\e{0}%
\e{0}%
\eol}\vss}\rg%
%
%
\rx{\vss\hfull{%
\rlx{\hss{$972_{x}$}}\cg%
\e{1}%
\e{1}%
\e{0}%
\e{1}%
\e{2}%
\e{2}%
\e{1}%
\e{1}%
\e{0}%
\e{0}%
\e{0}%
\e{0}%
\e{0}%
\e{1}%
\e{0}%
\e{0}%
\e{0}%
\e{0}%
\e{0}%
\e{0}%
\e{0}%
\e{0}%
\e{0}%
\e{0}%
\eol}\vss}\rg%
%
%
\rx{\vss\hfull{%
\rlx{\hss{$4096_{x}$}}\cg%
\e{6}%
\e{3}%
\e{0}%
\e{4}%
\e{9}%
\e{6}%
\e{6}%
\e{2}%
\e{1}%
\e{0}%
\e{2}%
\e{5}%
\e{3}%
\e{3}%
\e{1}%
\e{0}%
\e{0}%
\e{0}%
\e{1}%
\e{0}%
\e{0}%
\e{0}%
\e{0}%
\e{0}%
\eol}\vss}\rg%
%
%
\rx{\vss\hfull{%
\rlx{\hss{$4200_{x}$}}\cg%
\e{7}%
\e{3}%
\e{0}%
\e{2}%
\e{7}%
\e{8}%
\e{7}%
\e{5}%
\e{2}%
\e{0}%
\e{1}%
\e{3}%
\e{4}%
\e{3}%
\e{2}%
\e{1}%
\e{0}%
\e{0}%
\e{0}%
\e{1}%
\e{0}%
\e{0}%
\e{0}%
\e{0}%
\eol}\vss}\rg%
%
%
\rx{\vss\hfull{%
\rlx{\hss{$2240_{x}$}}\cg%
\e{3}%
\e{1}%
\e{0}%
\e{1}%
\e{4}%
\e{4}%
\e{3}%
\e{2}%
\e{1}%
\e{0}%
\e{0}%
\e{1}%
\e{1}%
\e{1}%
\e{1}%
\e{0}%
\e{0}%
\e{0}%
\e{0}%
\e{0}%
\e{0}%
\e{0}%
\e{0}%
\e{0}%
\eol}\vss}\rg%
%
%
\rx{\vss\hfull{%
\rlx{\hss{$2835_{x}$}}\cg%
\e{5}%
\e{2}%
\e{0}%
\e{1}%
\e{3}%
\e{6}%
\e{5}%
\e{6}%
\e{3}%
\e{1}%
\e{0}%
\e{1}%
\e{3}%
\e{2}%
\e{3}%
\e{1}%
\e{0}%
\e{0}%
\e{0}%
\e{0}%
\e{0}%
\e{1}%
\e{0}%
\e{0}%
\eol}\vss}\rg%
%
%
\rx{\vss\hfull{%
\rlx{\hss{$6075_{x}$}}\cg%
\e{11}%
\e{6}%
\e{1}%
\e{3}%
\e{10}%
\e{10}%
\e{11}%
\e{7}%
\e{4}%
\e{0}%
\e{3}%
\e{8}%
\e{8}%
\e{6}%
\e{4}%
\e{1}%
\e{0}%
\e{0}%
\e{1}%
\e{1}%
\e{1}%
\e{0}%
\e{0}%
\e{0}%
\eol}\vss}\rg%
%
%
\rx{\vss\hfull{%
\rlx{\hss{$3200_{x}$}}\cg%
\e{7}%
\e{5}%
\e{1}%
\e{2}%
\e{7}%
\e{6}%
\e{7}%
\e{4}%
\e{2}%
\e{0}%
\e{0}%
\e{3}%
\e{2}%
\e{5}%
\e{2}%
\e{2}%
\e{0}%
\e{0}%
\e{1}%
\e{0}%
\e{1}%
\e{0}%
\e{0}%
\e{0}%
\eol}\vss}\rg%
%
%
\rx{\vss\hfull{%
\rlx{\hss{$70_{y}$}}\cg%
\e{0}%
\e{0}%
\e{0}%
\e{0}%
\e{0}%
\e{0}%
\e{0}%
\e{0}%
\e{0}%
\e{0}%
\e{0}%
\e{1}%
\e{0}%
\e{0}%
\e{0}%
\e{0}%
\e{0}%
\e{1}%
\e{0}%
\e{0}%
\e{0}%
\e{0}%
\e{0}%
\e{0}%
\eol}\vss}\rg%
%
%
\rx{\vss\hfull{%
\rlx{\hss{$1134_{y}$}}\cg%
\e{3}%
\e{2}%
\e{1}%
\e{0}%
\e{2}%
\e{2}%
\e{3}%
\e{1}%
\e{1}%
\e{0}%
\e{0}%
\e{2}%
\e{2}%
\e{2}%
\e{1}%
\e{1}%
\e{0}%
\e{0}%
\e{0}%
\e{0}%
\e{1}%
\e{0}%
\e{0}%
\e{0}%
\eol}\vss}\rg%
%
%
\rx{\vss\hfull{%
\rlx{\hss{$1680_{y}$}}\cg%
\e{3}%
\e{4}%
\e{0}%
\e{1}%
\e{3}%
\e{2}%
\e{3}%
\e{1}%
\e{1}%
\e{0}%
\e{2}%
\e{5}%
\e{3}%
\e{4}%
\e{1}%
\e{0}%
\e{0}%
\e{1}%
\e{2}%
\e{1}%
\e{0}%
\e{0}%
\e{0}%
\e{0}%
\eol}\vss}\rg%
%
%
\rx{\vss\hfull{%
\rlx{\hss{$168_{y}$}}\cg%
\e{0}%
\e{0}%
\e{0}%
\e{0}%
\e{0}%
\e{1}%
\e{0}%
\e{1}%
\e{0}%
\e{0}%
\e{0}%
\e{0}%
\e{0}%
\e{0}%
\e{0}%
\e{1}%
\e{0}%
\e{0}%
\e{0}%
\e{0}%
\e{0}%
\e{0}%
\e{0}%
\e{0}%
\eol}\vss}\rg%
%
%
\rx{\vss\hfull{%
\rlx{\hss{$420_{y}$}}\cg%
\e{1}%
\e{0}%
\e{0}%
\e{0}%
\e{0}%
\e{1}%
\e{1}%
\e{2}%
\e{1}%
\e{0}%
\e{0}%
\e{0}%
\e{0}%
\e{0}%
\e{1}%
\e{1}%
\e{0}%
\e{0}%
\e{0}%
\e{0}%
\e{0}%
\e{0}%
\e{0}%
\e{1}%
\eol}\vss}\rg%
%
%
\rx{\vss\hfull{%
\rlx{\hss{$3150_{y}$}}\cg%
\e{8}%
\e{4}%
\e{1}%
\e{0}%
\e{4}%
\e{6}%
\e{8}%
\e{7}%
\e{5}%
\e{1}%
\e{0}%
\e{2}%
\e{4}%
\e{4}%
\e{5}%
\e{3}%
\e{1}%
\e{0}%
\e{0}%
\e{0}%
\e{1}%
\e{1}%
\e{1}%
\e{0}%
\eol}\vss}\rg%
%
%
\rx{\vss\hfull{%
\rlx{\hss{$4200_{y}$}}\cg%
\e{9}%
\e{6}%
\e{1}%
\e{1}%
\e{5}%
\e{10}%
\e{9}%
\e{10}%
\e{6}%
\e{1}%
\e{0}%
\e{2}%
\e{5}%
\e{6}%
\e{6}%
\e{5}%
\e{1}%
\e{0}%
\e{0}%
\e{1}%
\e{1}%
\e{1}%
\e{1}%
\e{0}%
\eol}\vss}\rg%
\eop
\eject
\tablecont%
%
%
%
%
%
%
\rowpts=18 true pt%
\colpts=18 true pt%
\rowlabpts=40 true pt%
\collabpts=65 true pt%
\clx{\vss\hfull{%
\rlx{\hss{$ $}}\cg%
\cx{\hskip 16 true pt\flip{$[{2^{2}}{1}]{\times}[{3}{1^{2}}]$}\hss}\cg%
\cx{\hskip 16 true pt\flip{$[{2}{1^{3}}]{\times}[{3}{1^{2}}]$}\hss}\cg%
\cx{\hskip 16 true pt\flip{$[{1^{5}}]{\times}[{3}{1^{2}}]$}\hss}\cg%
\cx{\hskip 16 true pt\flip{$[{5}]{\times}[{2^{2}}{1}]$}\hss}\cg%
\cx{\hskip 16 true pt\flip{$[{4}{1}]{\times}[{2^{2}}{1}]$}\hss}\cg%
\cx{\hskip 16 true pt\flip{$[{3}{2}]{\times}[{2^{2}}{1}]$}\hss}\cg%
\cx{\hskip 16 true pt\flip{$[{3}{1^{2}}]{\times}[{2^{2}}{1}]$}\hss}\cg%
\cx{\hskip 16 true pt\flip{$[{2^{2}}{1}]{\times}[{2^{2}}{1}]$}\hss}\cg%
\cx{\hskip 16 true pt\flip{$[{2}{1^{3}}]{\times}[{2^{2}}{1}]$}\hss}\cg%
\cx{\hskip 16 true pt\flip{$[{1^{5}}]{\times}[{2^{2}}{1}]$}\hss}\cg%
\cx{\hskip 16 true pt\flip{$[{5}]{\times}[{2}{1^{3}}]$}\hss}\cg%
\cx{\hskip 16 true pt\flip{$[{4}{1}]{\times}[{2}{1^{3}}]$}\hss}\cg%
\cx{\hskip 16 true pt\flip{$[{3}{2}]{\times}[{2}{1^{3}}]$}\hss}\cg%
\cx{\hskip 16 true pt\flip{$[{3}{1^{2}}]{\times}[{2}{1^{3}}]$}\hss}\cg%
\cx{\hskip 16 true pt\flip{$[{2^{2}}{1}]{\times}[{2}{1^{3}}]$}\hss}\cg%
\cx{\hskip 16 true pt\flip{$[{2}{1^{3}}]{\times}[{2}{1^{3}}]$}\hss}\cg%
\cx{\hskip 16 true pt\flip{$[{1^{5}}]{\times}[{2}{1^{3}}]$}\hss}\cg%
\cx{\hskip 16 true pt\flip{$[{5}]{\times}[{1^{5}}]$}\hss}\cg%
\cx{\hskip 16 true pt\flip{$[{4}{1}]{\times}[{1^{5}}]$}\hss}\cg%
\cx{\hskip 16 true pt\flip{$[{3}{2}]{\times}[{1^{5}}]$}\hss}\cg%
\cx{\hskip 16 true pt\flip{$[{3}{1^{2}}]{\times}[{1^{5}}]$}\hss}\cg%
\cx{\hskip 16 true pt\flip{$[{2^{2}}{1}]{\times}[{1^{5}}]$}\hss}\cg%
\cx{\hskip 16 true pt\flip{$[{2}{1^{3}}]{\times}[{1^{5}}]$}\hss}\cg%
\cx{\hskip 16 true pt\flip{$[{1^{5}}]{\times}[{1^{5}}]$}\hss}\cg%
\eol}}\rg%
%
%
\rx{\vss\hfull{%
\rlx{\hss{$2688_{y}$}}\cg%
\e{6}%
\e{4}%
\e{1}%
\e{1}%
\e{4}%
\e{5}%
\e{6}%
\e{5}%
\e{4}%
\e{1}%
\e{0}%
\e{2}%
\e{4}%
\e{4}%
\e{4}%
\e{2}%
\e{0}%
\e{0}%
\e{0}%
\e{1}%
\e{1}%
\e{1}%
\e{0}%
\e{0}%
\eol}\vss}\rg%
%
%
\rx{\vss\hfull{%
\rlx{\hss{$2100_{y}$}}\cg%
\e{5}%
\e{4}%
\e{1}%
\e{0}%
\e{3}%
\e{3}%
\e{5}%
\e{2}%
\e{2}%
\e{0}%
\e{2}%
\e{5}%
\e{3}%
\e{4}%
\e{2}%
\e{1}%
\e{0}%
\e{1}%
\e{2}%
\e{0}%
\e{1}%
\e{0}%
\e{0}%
\e{0}%
\eol}\vss}\rg%
%
%
\rx{\vss\hfull{%
\rlx{\hss{$1400_{y}$}}\cg%
\e{2}%
\e{3}%
\e{0}%
\e{1}%
\e{2}%
\e{3}%
\e{2}%
\e{2}%
\e{1}%
\e{0}%
\e{1}%
\e{3}%
\e{2}%
\e{3}%
\e{1}%
\e{1}%
\e{0}%
\e{1}%
\e{1}%
\e{1}%
\e{0}%
\e{0}%
\e{0}%
\e{0}%
\eol}\vss}\rg%
%
%
\rx{\vss\hfull{%
\rlx{\hss{$4536_{y}$}}\cg%
\e{9}%
\e{8}%
\e{1}%
\e{2}%
\e{7}%
\e{8}%
\e{9}%
\e{7}%
\e{5}%
\e{1}%
\e{2}%
\e{7}%
\e{7}%
\e{8}%
\e{5}%
\e{3}%
\e{0}%
\e{1}%
\e{2}%
\e{2}%
\e{1}%
\e{1}%
\e{0}%
\e{0}%
\eol}\vss}\rg%
%
%
\rx{\vss\hfull{%
\rlx{\hss{$5670_{y}$}}\cg%
\e{12}%
\e{10}%
\e{2}%
\e{2}%
\e{9}%
\e{10}%
\e{12}%
\e{8}%
\e{6}%
\e{1}%
\e{2}%
\e{9}%
\e{9}%
\e{10}%
\e{6}%
\e{4}%
\e{0}%
\e{1}%
\e{2}%
\e{2}%
\e{2}%
\e{1}%
\e{0}%
\e{0}%
\eol}\vss}\rg%
%
%
\rx{\vss\hfull{%
\rlx{\hss{$4480_{y}$}}\cg%
\e{10}%
\e{7}%
\e{1}%
\e{1}%
\e{6}%
\e{9}%
\e{10}%
\e{9}%
\e{6}%
\e{1}%
\e{1}%
\e{4}%
\e{6}%
\e{7}%
\e{6}%
\e{4}%
\e{1}%
\e{0}%
\e{1}%
\e{1}%
\e{1}%
\e{1}%
\e{1}%
\e{0}%
\eol}\vss}\rg%
%
%
\rx{\vss\hfull{%
\rlx{\hss{$8_{z}$}}\cg%
\e{0}%
\e{0}%
\e{0}%
\e{0}%
\e{0}%
\e{0}%
\e{0}%
\e{0}%
\e{0}%
\e{0}%
\e{0}%
\e{0}%
\e{0}%
\e{0}%
\e{0}%
\e{0}%
\e{0}%
\e{0}%
\e{0}%
\e{0}%
\e{0}%
\e{0}%
\e{0}%
\e{0}%
\eol}\vss}\rg%
%
%
\rx{\vss\hfull{%
\rlx{\hss{$56_{z}$}}\cg%
\e{0}%
\e{0}%
\e{0}%
\e{0}%
\e{0}%
\e{0}%
\e{0}%
\e{0}%
\e{0}%
\e{0}%
\e{1}%
\e{0}%
\e{0}%
\e{0}%
\e{0}%
\e{0}%
\e{0}%
\e{0}%
\e{0}%
\e{0}%
\e{0}%
\e{0}%
\e{0}%
\e{0}%
\eol}\vss}\rg%
%
%
\rx{\vss\hfull{%
\rlx{\hss{$160_{z}$}}\cg%
\e{0}%
\e{0}%
\e{0}%
\e{1}%
\e{0}%
\e{0}%
\e{0}%
\e{0}%
\e{0}%
\e{0}%
\e{0}%
\e{0}%
\e{0}%
\e{0}%
\e{0}%
\e{0}%
\e{0}%
\e{0}%
\e{0}%
\e{0}%
\e{0}%
\e{0}%
\e{0}%
\e{0}%
\eol}\vss}\rg%
%
%
\rx{\vss\hfull{%
\rlx{\hss{$112_{z}$}}\cg%
\e{0}%
\e{0}%
\e{0}%
\e{0}%
\e{0}%
\e{0}%
\e{0}%
\e{0}%
\e{0}%
\e{0}%
\e{0}%
\e{0}%
\e{0}%
\e{0}%
\e{0}%
\e{0}%
\e{0}%
\e{0}%
\e{0}%
\e{0}%
\e{0}%
\e{0}%
\e{0}%
\e{0}%
\eol}\vss}\rg%
%
%
\rx{\vss\hfull{%
\rlx{\hss{$840_{z}$}}\cg%
\e{1}%
\e{1}%
\e{0}%
\e{1}%
\e{3}%
\e{1}%
\e{1}%
\e{0}%
\e{0}%
\e{0}%
\e{0}%
\e{1}%
\e{0}%
\e{1}%
\e{0}%
\e{0}%
\e{0}%
\e{0}%
\e{0}%
\e{0}%
\e{0}%
\e{0}%
\e{0}%
\e{0}%
\eol}\vss}\rg%
%
%
\rx{\vss\hfull{%
\rlx{\hss{$1296_{z}$}}\cg%
\e{1}%
\e{0}%
\e{0}%
\e{2}%
\e{3}%
\e{1}%
\e{1}%
\e{0}%
\e{0}%
\e{0}%
\e{2}%
\e{3}%
\e{1}%
\e{0}%
\e{0}%
\e{0}%
\e{0}%
\e{1}%
\e{0}%
\e{0}%
\e{0}%
\e{0}%
\e{0}%
\e{0}%
\eol}\vss}\rg%
%
%
\rx{\vss\hfull{%
\rlx{\hss{$1400_{z}$}}\cg%
\e{1}%
\e{0}%
\e{0}%
\e{2}%
\e{3}%
\e{1}%
\e{1}%
\e{0}%
\e{0}%
\e{0}%
\e{1}%
\e{1}%
\e{0}%
\e{0}%
\e{0}%
\e{0}%
\e{0}%
\e{0}%
\e{0}%
\e{0}%
\e{0}%
\e{0}%
\e{0}%
\e{0}%
\eol}\vss}\rg%
%
%
\rx{\vss\hfull{%
\rlx{\hss{$1008_{z}$}}\cg%
\e{0}%
\e{0}%
\e{0}%
\e{2}%
\e{2}%
\e{1}%
\e{0}%
\e{0}%
\e{0}%
\e{0}%
\e{2}%
\e{1}%
\e{0}%
\e{0}%
\e{0}%
\e{0}%
\e{0}%
\e{0}%
\e{0}%
\e{0}%
\e{0}%
\e{0}%
\e{0}%
\e{0}%
\eol}\vss}\rg%
%
%
\rx{\vss\hfull{%
\rlx{\hss{$560_{z}$}}\cg%
\e{0}%
\e{0}%
\e{0}%
\e{1}%
\e{1}%
\e{0}%
\e{0}%
\e{0}%
\e{0}%
\e{0}%
\e{0}%
\e{0}%
\e{0}%
\e{0}%
\e{0}%
\e{0}%
\e{0}%
\e{0}%
\e{0}%
\e{0}%
\e{0}%
\e{0}%
\e{0}%
\e{0}%
\eol}\vss}\rg%
%
%
\rx{\vss\hfull{%
\rlx{\hss{$1400_{zz}$}}\cg%
\e{2}%
\e{0}%
\e{0}%
\e{0}%
\e{1}%
\e{3}%
\e{2}%
\e{2}%
\e{1}%
\e{0}%
\e{0}%
\e{0}%
\e{1}%
\e{0}%
\e{1}%
\e{0}%
\e{0}%
\e{0}%
\e{0}%
\e{0}%
\e{0}%
\e{0}%
\e{0}%
\e{0}%
\eol}\vss}\rg%
%
%
\rx{\vss\hfull{%
\rlx{\hss{$4200_{z}$}}\cg%
\e{7}%
\e{4}%
\e{0}%
\e{1}%
\e{5}%
\e{9}%
\e{7}%
\e{8}%
\e{4}%
\e{1}%
\e{1}%
\e{3}%
\e{6}%
\e{4}%
\e{4}%
\e{2}%
\e{0}%
\e{0}%
\e{0}%
\e{1}%
\e{0}%
\e{1}%
\e{0}%
\e{0}%
\eol}\vss}\rg%
%
%
\rx{\vss\hfull{%
\rlx{\hss{$400_{z}$}}\cg%
\e{0}%
\e{0}%
\e{0}%
\e{0}%
\e{0}%
\e{1}%
\e{0}%
\e{0}%
\e{0}%
\e{0}%
\e{0}%
\e{0}%
\e{0}%
\e{0}%
\e{0}%
\e{0}%
\e{0}%
\e{0}%
\e{0}%
\e{0}%
\e{0}%
\e{0}%
\e{0}%
\e{0}%
\eol}\vss}\rg%
%
%
\rx{\vss\hfull{%
\rlx{\hss{$3240_{z}$}}\cg%
\e{3}%
\e{1}%
\e{0}%
\e{3}%
\e{6}%
\e{5}%
\e{3}%
\e{2}%
\e{0}%
\e{0}%
\e{1}%
\e{2}%
\e{2}%
\e{1}%
\e{0}%
\e{0}%
\e{0}%
\e{0}%
\e{0}%
\e{0}%
\e{0}%
\e{0}%
\e{0}%
\e{0}%
\eol}\vss}\rg%
%
%
\rx{\vss\hfull{%
\rlx{\hss{$4536_{z}$}}\cg%
\e{9}%
\e{4}%
\e{1}%
\e{2}%
\e{8}%
\e{9}%
\e{9}%
\e{6}%
\e{3}%
\e{0}%
\e{0}%
\e{3}%
\e{3}%
\e{4}%
\e{3}%
\e{2}%
\e{0}%
\e{0}%
\e{0}%
\e{0}%
\e{1}%
\e{0}%
\e{0}%
\e{0}%
\eol}\vss}\rg%
%
%
\rx{\vss\hfull{%
\rlx{\hss{$2400_{z}$}}\cg%
\e{4}%
\e{3}%
\e{0}%
\e{1}%
\e{4}%
\e{3}%
\e{4}%
\e{2}%
\e{1}%
\e{0}%
\e{3}%
\e{6}%
\e{4}%
\e{3}%
\e{1}%
\e{0}%
\e{0}%
\e{1}%
\e{1}%
\e{1}%
\e{0}%
\e{0}%
\e{0}%
\e{0}%
\eol}\vss}\rg%
%
%
\rx{\vss\hfull{%
\rlx{\hss{$3360_{z}$}}\cg%
\e{5}%
\e{2}%
\e{0}%
\e{2}%
\e{5}%
\e{6}%
\e{5}%
\e{4}%
\e{2}%
\e{0}%
\e{1}%
\e{3}%
\e{4}%
\e{2}%
\e{2}%
\e{0}%
\e{0}%
\e{0}%
\e{0}%
\e{1}%
\e{0}%
\e{0}%
\e{0}%
\e{0}%
\eol}\vss}\rg%
%
%
\rx{\vss\hfull{%
\rlx{\hss{$2800_{z}$}}\cg%
\e{4}%
\e{2}%
\e{0}%
\e{2}%
\e{5}%
\e{4}%
\e{4}%
\e{2}%
\e{1}%
\e{0}%
\e{3}%
\e{5}%
\e{3}%
\e{2}%
\e{1}%
\e{0}%
\e{0}%
\e{1}%
\e{1}%
\e{0}%
\e{0}%
\e{0}%
\e{0}%
\e{0}%
\eol}\vss}\rg%
%
%
\rx{\vss\hfull{%
\rlx{\hss{$4096_{z}$}}\cg%
\e{6}%
\e{3}%
\e{0}%
\e{4}%
\e{9}%
\e{6}%
\e{6}%
\e{2}%
\e{1}%
\e{0}%
\e{2}%
\e{5}%
\e{3}%
\e{3}%
\e{1}%
\e{0}%
\e{0}%
\e{0}%
\e{1}%
\e{0}%
\e{0}%
\e{0}%
\e{0}%
\e{0}%
\eol}\vss}\rg%
%
%
\rx{\vss\hfull{%
\rlx{\hss{$5600_{z}$}}\cg%
\e{11}%
\e{8}%
\e{1}%
\e{3}%
\e{11}%
\e{9}%
\e{11}%
\e{6}%
\e{3}%
\e{0}%
\e{3}%
\e{9}%
\e{6}%
\e{8}%
\e{3}%
\e{2}%
\e{0}%
\e{1}%
\e{2}%
\e{1}%
\e{1}%
\e{0}%
\e{0}%
\e{0}%
\eol}\vss}\rg%
%
%
\rx{\vss\hfull{%
\rlx{\hss{$448_{z}$}}\cg%
\e{1}%
\e{0}%
\e{0}%
\e{0}%
\e{1}%
\e{0}%
\e{1}%
\e{0}%
\e{0}%
\e{0}%
\e{0}%
\e{0}%
\e{0}%
\e{0}%
\e{0}%
\e{0}%
\e{0}%
\e{0}%
\e{0}%
\e{0}%
\e{0}%
\e{0}%
\e{0}%
\e{0}%
\eol}\vss}\rg%
%
%
\rx{\vss\hfull{%
\rlx{\hss{$448_{w}$}}\cg%
\e{1}%
\e{1}%
\e{0}%
\e{0}%
\e{1}%
\e{0}%
\e{1}%
\e{0}%
\e{0}%
\e{0}%
\e{1}%
\e{2}%
\e{1}%
\e{1}%
\e{0}%
\e{0}%
\e{0}%
\e{0}%
\e{1}%
\e{0}%
\e{0}%
\e{0}%
\e{0}%
\e{0}%
\eol}\vss}\rg%
%
%
\rx{\vss\hfull{%
\rlx{\hss{$1344_{w}$}}\cg%
\e{3}%
\e{2}%
\e{0}%
\e{0}%
\e{1}%
\e{4}%
\e{3}%
\e{4}%
\e{2}%
\e{0}%
\e{0}%
\e{0}%
\e{1}%
\e{2}%
\e{2}%
\e{2}%
\e{1}%
\e{0}%
\e{0}%
\e{0}%
\e{0}%
\e{0}%
\e{1}%
\e{0}%
\eol}\vss}\rg%
%
%
\rx{\vss\hfull{%
\rlx{\hss{$5600_{w}$}}\cg%
\e{12}%
\e{10}%
\e{2}%
\e{2}%
\e{9}%
\e{10}%
\e{12}%
\e{8}%
\e{6}%
\e{1}%
\e{2}%
\e{8}%
\e{9}%
\e{10}%
\e{6}%
\e{4}%
\e{0}%
\e{0}%
\e{2}%
\e{2}%
\e{2}%
\e{1}%
\e{0}%
\e{0}%
\eol}\vss}\rg%
%
%
\rx{\vss\hfull{%
\rlx{\hss{$2016_{w}$}}\cg%
\e{5}%
\e{2}%
\e{0}%
\e{0}%
\e{2}%
\e{4}%
\e{5}%
\e{6}%
\e{4}%
\e{1}%
\e{0}%
\e{0}%
\e{2}%
\e{2}%
\e{4}%
\e{2}%
\e{1}%
\e{0}%
\e{0}%
\e{0}%
\e{0}%
\e{1}%
\e{1}%
\e{0}%
\eol}\vss}\rg%
%
%
\rx{\vss\hfull{%
\rlx{\hss{$7168_{w}$}}\cg%
\e{16}%
\e{11}%
\e{2}%
\e{2}%
\e{10}%
\e{14}%
\e{16}%
\e{14}%
\e{10}%
\e{2}%
\e{1}%
\e{6}%
\e{10}%
\e{11}%
\e{10}%
\e{6}%
\e{1}%
\e{0}%
\e{1}%
\e{2}%
\e{2}%
\e{2}%
\e{1}%
\e{0}%
\eol}\vss}\rg%
\tableclose%
%
%
%
%
%
%
\eop
\eject
\tableopen{Induce/restrict matrix for $W({D_{5}}{A_{3}})\,\subset\,W(E_{8})$}%
%
%
%
%
%
%
\rowpts=18 true pt%
\colpts=18 true pt%
\rowlabpts=40 true pt%
\collabpts=85 true pt%
\clx{\vss\hfull{%
\rlx{\hss{$ $}}\cg%
\cx{\hskip 16 true pt\flip{$[{5}:-]{\times}[{4}]$}\hss}\cg%
\cx{\hskip 16 true pt\flip{$[{4}{1}:-]{\times}[{4}]$}\hss}\cg%
\cx{\hskip 16 true pt\flip{$[{3}{2}:-]{\times}[{4}]$}\hss}\cg%
\cx{\hskip 16 true pt\flip{$[{3}{1^{2}}:-]{\times}[{4}]$}\hss}\cg%
\cx{\hskip 16 true pt\flip{$[{2^{2}}{1}:-]{\times}[{4}]$}\hss}\cg%
\cx{\hskip 16 true pt\flip{$[{2}{1^{3}}:-]{\times}[{4}]$}\hss}\cg%
\cx{\hskip 16 true pt\flip{$[{1^{5}}:-]{\times}[{4}]$}\hss}\cg%
\cx{\hskip 16 true pt\flip{$[{4}:{1}]{\times}[{4}]$}\hss}\cg%
\cx{\hskip 16 true pt\flip{$[{3}{1}:{1}]{\times}[{4}]$}\hss}\cg%
\cx{\hskip 16 true pt\flip{$[{2^{2}}:{1}]{\times}[{4}]$}\hss}\cg%
\cx{\hskip 16 true pt\flip{$[{2}{1^{2}}:{1}]{\times}[{4}]$}\hss}\cg%
\cx{\hskip 16 true pt\flip{$[{1^{4}}:{1}]{\times}[{4}]$}\hss}\cg%
\cx{\hskip 16 true pt\flip{$[{3}:{2}]{\times}[{4}]$}\hss}\cg%
\cx{\hskip 16 true pt\flip{$[{3}:{1^{2}}]{\times}[{4}]$}\hss}\cg%
\cx{\hskip 16 true pt\flip{$[{2}{1}:{2}]{\times}[{4}]$}\hss}\cg%
\cx{\hskip 16 true pt\flip{$[{2}{1}:{1^{2}}]{\times}[{4}]$}\hss}\cg%
\cx{\hskip 16 true pt\flip{$[{1^{3}}:{2}]{\times}[{4}]$}\hss}\cg%
\cx{\hskip 16 true pt\flip{$[{1^{3}}:{1^{2}}]{\times}[{4}]$}\hss}\cg%
\cx{\hskip 16 true pt\flip{$[{5}:-]{\times}[{3}{1}]$}\hss}\cg%
\cx{\hskip 16 true pt\flip{$[{4}{1}:-]{\times}[{3}{1}]$}\hss}\cg%
\cx{\hskip 16 true pt\flip{$[{3}{2}:-]{\times}[{3}{1}]$}\hss}\cg%
\cx{\hskip 16 true pt\flip{$[{3}{1^{2}}:-]{\times}[{3}{1}]$}\hss}\cg%
\cx{\hskip 16 true pt\flip{$[{2^{2}}{1}:-]{\times}[{3}{1}]$}\hss}\cg%
\eol}}\rg%
%
%
\rx{\vss\hfull{%
\rlx{\hss{$1_{x}$}}\cg%
\e{1}%
\e{0}%
\e{0}%
\e{0}%
\e{0}%
\e{0}%
\e{0}%
\e{0}%
\e{0}%
\e{0}%
\e{0}%
\e{0}%
\e{0}%
\e{0}%
\e{0}%
\e{0}%
\e{0}%
\e{0}%
\e{0}%
\e{0}%
\e{0}%
\e{0}%
\e{0}%
\eol}\vss}\rg%
%
%
\rx{\vss\hfull{%
\rlx{\hss{$28_{x}$}}\cg%
\e{0}%
\e{0}%
\e{0}%
\e{0}%
\e{0}%
\e{0}%
\e{0}%
\e{0}%
\e{0}%
\e{0}%
\e{0}%
\e{0}%
\e{0}%
\e{1}%
\e{0}%
\e{0}%
\e{0}%
\e{0}%
\e{0}%
\e{0}%
\e{0}%
\e{0}%
\e{0}%
\eol}\vss}\rg%
%
%
\rx{\vss\hfull{%
\rlx{\hss{$35_{x}$}}\cg%
\e{1}%
\e{1}%
\e{0}%
\e{0}%
\e{0}%
\e{0}%
\e{0}%
\e{0}%
\e{0}%
\e{0}%
\e{0}%
\e{0}%
\e{1}%
\e{0}%
\e{0}%
\e{0}%
\e{0}%
\e{0}%
\e{1}%
\e{0}%
\e{0}%
\e{0}%
\e{0}%
\eol}\vss}\rg%
%
%
\rx{\vss\hfull{%
\rlx{\hss{$84_{x}$}}\cg%
\e{2}%
\e{1}%
\e{1}%
\e{0}%
\e{0}%
\e{0}%
\e{0}%
\e{1}%
\e{0}%
\e{0}%
\e{0}%
\e{0}%
\e{1}%
\e{0}%
\e{0}%
\e{0}%
\e{0}%
\e{0}%
\e{1}%
\e{0}%
\e{0}%
\e{0}%
\e{0}%
\eol}\vss}\rg%
%
%
\rx{\vss\hfull{%
\rlx{\hss{$50_{x}$}}\cg%
\e{1}%
\e{1}%
\e{0}%
\e{0}%
\e{0}%
\e{0}%
\e{0}%
\e{1}%
\e{0}%
\e{0}%
\e{0}%
\e{0}%
\e{0}%
\e{0}%
\e{0}%
\e{0}%
\e{0}%
\e{0}%
\e{0}%
\e{0}%
\e{0}%
\e{0}%
\e{0}%
\eol}\vss}\rg%
%
%
\rx{\vss\hfull{%
\rlx{\hss{$350_{x}$}}\cg%
\e{0}%
\e{0}%
\e{0}%
\e{0}%
\e{0}%
\e{0}%
\e{0}%
\e{0}%
\e{0}%
\e{0}%
\e{1}%
\e{0}%
\e{0}%
\e{1}%
\e{0}%
\e{1}%
\e{0}%
\e{0}%
\e{0}%
\e{0}%
\e{0}%
\e{0}%
\e{0}%
\eol}\vss}\rg%
%
%
\rx{\vss\hfull{%
\rlx{\hss{$300_{x}$}}\cg%
\e{1}%
\e{1}%
\e{1}%
\e{0}%
\e{0}%
\e{0}%
\e{0}%
\e{0}%
\e{0}%
\e{1}%
\e{0}%
\e{0}%
\e{1}%
\e{0}%
\e{1}%
\e{0}%
\e{0}%
\e{0}%
\e{1}%
\e{1}%
\e{0}%
\e{0}%
\e{0}%
\eol}\vss}\rg%
%
%
\rx{\vss\hfull{%
\rlx{\hss{$567_{x}$}}\cg%
\e{0}%
\e{1}%
\e{0}%
\e{1}%
\e{0}%
\e{0}%
\e{0}%
\e{1}%
\e{1}%
\e{0}%
\e{0}%
\e{0}%
\e{2}%
\e{2}%
\e{1}%
\e{1}%
\e{0}%
\e{0}%
\e{2}%
\e{1}%
\e{0}%
\e{0}%
\e{0}%
\eol}\vss}\rg%
%
%
\rx{\vss\hfull{%
\rlx{\hss{$210_{x}$}}\cg%
\e{1}%
\e{1}%
\e{0}%
\e{0}%
\e{0}%
\e{0}%
\e{0}%
\e{1}%
\e{0}%
\e{0}%
\e{0}%
\e{0}%
\e{1}%
\e{1}%
\e{1}%
\e{0}%
\e{0}%
\e{0}%
\e{1}%
\e{1}%
\e{0}%
\e{0}%
\e{0}%
\eol}\vss}\rg%
%
%
\rx{\vss\hfull{%
\rlx{\hss{$840_{x}$}}\cg%
\e{1}%
\e{1}%
\e{1}%
\e{0}%
\e{0}%
\e{0}%
\e{0}%
\e{0}%
\e{0}%
\e{1}%
\e{0}%
\e{0}%
\e{1}%
\e{0}%
\e{1}%
\e{0}%
\e{0}%
\e{0}%
\e{0}%
\e{1}%
\e{1}%
\e{0}%
\e{1}%
\eol}\vss}\rg%
%
%
\rx{\vss\hfull{%
\rlx{\hss{$700_{x}$}}\cg%
\e{2}%
\e{2}%
\e{1}%
\e{0}%
\e{0}%
\e{0}%
\e{0}%
\e{2}%
\e{1}%
\e{0}%
\e{0}%
\e{0}%
\e{3}%
\e{0}%
\e{2}%
\e{0}%
\e{0}%
\e{0}%
\e{2}%
\e{2}%
\e{1}%
\e{0}%
\e{0}%
\eol}\vss}\rg%
%
%
\rx{\vss\hfull{%
\rlx{\hss{$175_{x}$}}\cg%
\e{0}%
\e{0}%
\e{0}%
\e{0}%
\e{0}%
\e{0}%
\e{0}%
\e{1}%
\e{0}%
\e{0}%
\e{0}%
\e{0}%
\e{1}%
\e{0}%
\e{0}%
\e{0}%
\e{0}%
\e{0}%
\e{0}%
\e{0}%
\e{1}%
\e{0}%
\e{0}%
\eol}\vss}\rg%
%
%
\rx{\vss\hfull{%
\rlx{\hss{$1400_{x}$}}\cg%
\e{0}%
\e{0}%
\e{0}%
\e{0}%
\e{0}%
\e{0}%
\e{0}%
\e{2}%
\e{2}%
\e{0}%
\e{0}%
\e{0}%
\e{2}%
\e{2}%
\e{2}%
\e{1}%
\e{1}%
\e{0}%
\e{1}%
\e{2}%
\e{1}%
\e{1}%
\e{0}%
\eol}\vss}\rg%
%
%
\rx{\vss\hfull{%
\rlx{\hss{$1050_{x}$}}\cg%
\e{0}%
\e{1}%
\e{1}%
\e{1}%
\e{0}%
\e{0}%
\e{0}%
\e{2}%
\e{2}%
\e{0}%
\e{0}%
\e{0}%
\e{2}%
\e{1}%
\e{1}%
\e{0}%
\e{0}%
\e{0}%
\e{1}%
\e{1}%
\e{1}%
\e{0}%
\e{0}%
\eol}\vss}\rg%
%
%
\rx{\vss\hfull{%
\rlx{\hss{$1575_{x}$}}\cg%
\e{0}%
\e{0}%
\e{0}%
\e{0}%
\e{0}%
\e{0}%
\e{0}%
\e{1}%
\e{2}%
\e{0}%
\e{1}%
\e{0}%
\e{1}%
\e{3}%
\e{2}%
\e{2}%
\e{1}%
\e{0}%
\e{1}%
\e{2}%
\e{0}%
\e{1}%
\e{0}%
\eol}\vss}\rg%
%
%
\rx{\vss\hfull{%
\rlx{\hss{$1344_{x}$}}\cg%
\e{1}%
\e{2}%
\e{2}%
\e{1}%
\e{1}%
\e{0}%
\e{0}%
\e{2}%
\e{2}%
\e{1}%
\e{0}%
\e{0}%
\e{3}%
\e{1}%
\e{2}%
\e{1}%
\e{0}%
\e{0}%
\e{2}%
\e{2}%
\e{1}%
\e{0}%
\e{0}%
\eol}\vss}\rg%
%
%
\rx{\vss\hfull{%
\rlx{\hss{$2100_{x}$}}\cg%
\e{0}%
\e{0}%
\e{0}%
\e{1}%
\e{0}%
\e{1}%
\e{0}%
\e{0}%
\e{1}%
\e{1}%
\e{2}%
\e{1}%
\e{0}%
\e{1}%
\e{1}%
\e{2}%
\e{1}%
\e{1}%
\e{0}%
\e{1}%
\e{0}%
\e{1}%
\e{0}%
\eol}\vss}\rg%
%
%
\rx{\vss\hfull{%
\rlx{\hss{$2268_{x}$}}\cg%
\e{0}%
\e{1}%
\e{0}%
\e{1}%
\e{0}%
\e{0}%
\e{0}%
\e{1}%
\e{2}%
\e{1}%
\e{1}%
\e{0}%
\e{2}%
\e{1}%
\e{3}%
\e{2}%
\e{1}%
\e{1}%
\e{1}%
\e{3}%
\e{1}%
\e{1}%
\e{0}%
\eol}\vss}\rg%
%
%
\rx{\vss\hfull{%
\rlx{\hss{$525_{x}$}}\cg%
\e{0}%
\e{0}%
\e{0}%
\e{1}%
\e{0}%
\e{1}%
\e{0}%
\e{1}%
\e{1}%
\e{0}%
\e{0}%
\e{0}%
\e{0}%
\e{1}%
\e{0}%
\e{1}%
\e{0}%
\e{0}%
\e{0}%
\e{0}%
\e{0}%
\e{0}%
\e{0}%
\eol}\vss}\rg%
%
%
\rx{\vss\hfull{%
\rlx{\hss{$700_{xx}$}}\cg%
\e{0}%
\e{1}%
\e{0}%
\e{1}%
\e{0}%
\e{0}%
\e{0}%
\e{1}%
\e{1}%
\e{0}%
\e{0}%
\e{0}%
\e{0}%
\e{1}%
\e{0}%
\e{0}%
\e{0}%
\e{0}%
\e{0}%
\e{0}%
\e{0}%
\e{0}%
\e{0}%
\eol}\vss}\rg%
%
%
\rx{\vss\hfull{%
\rlx{\hss{$972_{x}$}}\cg%
\e{1}%
\e{2}%
\e{2}%
\e{0}%
\e{1}%
\e{0}%
\e{0}%
\e{1}%
\e{1}%
\e{1}%
\e{0}%
\e{0}%
\e{1}%
\e{0}%
\e{1}%
\e{0}%
\e{0}%
\e{0}%
\e{0}%
\e{1}%
\e{1}%
\e{0}%
\e{0}%
\eol}\vss}\rg%
%
%
\rx{\vss\hfull{%
\rlx{\hss{$4096_{x}$}}\cg%
\e{0}%
\e{1}%
\e{1}%
\e{1}%
\e{1}%
\e{0}%
\e{0}%
\e{1}%
\e{3}%
\e{2}%
\e{2}%
\e{0}%
\e{2}%
\e{2}%
\e{4}%
\e{3}%
\e{1}%
\e{1}%
\e{1}%
\e{3}%
\e{2}%
\e{2}%
\e{1}%
\eol}\vss}\rg%
%
%
\rx{\vss\hfull{%
\rlx{\hss{$4200_{x}$}}\cg%
\e{0}%
\e{1}%
\e{1}%
\e{1}%
\e{0}%
\e{0}%
\e{0}%
\e{1}%
\e{3}%
\e{2}%
\e{1}%
\e{0}%
\e{3}%
\e{1}%
\e{4}%
\e{2}%
\e{1}%
\e{0}%
\e{1}%
\e{3}%
\e{3}%
\e{2}%
\e{1}%
\eol}\vss}\rg%
%
%
\rx{\vss\hfull{%
\rlx{\hss{$2240_{x}$}}\cg%
\e{1}%
\e{1}%
\e{1}%
\e{0}%
\e{0}%
\e{0}%
\e{0}%
\e{2}%
\e{2}%
\e{1}%
\e{0}%
\e{0}%
\e{3}%
\e{1}%
\e{3}%
\e{1}%
\e{0}%
\e{0}%
\e{1}%
\e{2}%
\e{3}%
\e{1}%
\e{1}%
\eol}\vss}\rg%
%
%
\rx{\vss\hfull{%
\rlx{\hss{$2835_{x}$}}\cg%
\e{0}%
\e{0}%
\e{1}%
\e{0}%
\e{0}%
\e{0}%
\e{0}%
\e{1}%
\e{2}%
\e{1}%
\e{0}%
\e{0}%
\e{2}%
\e{1}%
\e{2}%
\e{1}%
\e{0}%
\e{0}%
\e{0}%
\e{1}%
\e{3}%
\e{1}%
\e{1}%
\eol}\vss}\rg%
%
%
\rx{\vss\hfull{%
\rlx{\hss{$6075_{x}$}}\cg%
\e{0}%
\e{0}%
\e{0}%
\e{1}%
\e{1}%
\e{1}%
\e{0}%
\e{1}%
\e{4}%
\e{1}%
\e{3}%
\e{0}%
\e{1}%
\e{3}%
\e{3}%
\e{4}%
\e{2}%
\e{1}%
\e{0}%
\e{2}%
\e{2}%
\e{3}%
\e{1}%
\eol}\vss}\rg%
%
%
\rx{\vss\hfull{%
\rlx{\hss{$3200_{x}$}}\cg%
\e{0}%
\e{1}%
\e{2}%
\e{1}%
\e{2}%
\e{0}%
\e{0}%
\e{0}%
\e{1}%
\e{2}%
\e{1}%
\e{0}%
\e{1}%
\e{0}%
\e{2}%
\e{1}%
\e{0}%
\e{1}%
\e{0}%
\e{1}%
\e{2}%
\e{1}%
\e{2}%
\eol}\vss}\rg%
\eop
\eject
\tablecont%
%
%
%
%
%
%
\rowpts=18 true pt%
\colpts=18 true pt%
\rowlabpts=40 true pt%
\collabpts=85 true pt%
\clx{\vss\hfull{%
\rlx{\hss{$ $}}\cg%
\cx{\hskip 16 true pt\flip{$[{5}:-]{\times}[{4}]$}\hss}\cg%
\cx{\hskip 16 true pt\flip{$[{4}{1}:-]{\times}[{4}]$}\hss}\cg%
\cx{\hskip 16 true pt\flip{$[{3}{2}:-]{\times}[{4}]$}\hss}\cg%
\cx{\hskip 16 true pt\flip{$[{3}{1^{2}}:-]{\times}[{4}]$}\hss}\cg%
\cx{\hskip 16 true pt\flip{$[{2^{2}}{1}:-]{\times}[{4}]$}\hss}\cg%
\cx{\hskip 16 true pt\flip{$[{2}{1^{3}}:-]{\times}[{4}]$}\hss}\cg%
\cx{\hskip 16 true pt\flip{$[{1^{5}}:-]{\times}[{4}]$}\hss}\cg%
\cx{\hskip 16 true pt\flip{$[{4}:{1}]{\times}[{4}]$}\hss}\cg%
\cx{\hskip 16 true pt\flip{$[{3}{1}:{1}]{\times}[{4}]$}\hss}\cg%
\cx{\hskip 16 true pt\flip{$[{2^{2}}:{1}]{\times}[{4}]$}\hss}\cg%
\cx{\hskip 16 true pt\flip{$[{2}{1^{2}}:{1}]{\times}[{4}]$}\hss}\cg%
\cx{\hskip 16 true pt\flip{$[{1^{4}}:{1}]{\times}[{4}]$}\hss}\cg%
\cx{\hskip 16 true pt\flip{$[{3}:{2}]{\times}[{4}]$}\hss}\cg%
\cx{\hskip 16 true pt\flip{$[{3}:{1^{2}}]{\times}[{4}]$}\hss}\cg%
\cx{\hskip 16 true pt\flip{$[{2}{1}:{2}]{\times}[{4}]$}\hss}\cg%
\cx{\hskip 16 true pt\flip{$[{2}{1}:{1^{2}}]{\times}[{4}]$}\hss}\cg%
\cx{\hskip 16 true pt\flip{$[{1^{3}}:{2}]{\times}[{4}]$}\hss}\cg%
\cx{\hskip 16 true pt\flip{$[{1^{3}}:{1^{2}}]{\times}[{4}]$}\hss}\cg%
\cx{\hskip 16 true pt\flip{$[{5}:-]{\times}[{3}{1}]$}\hss}\cg%
\cx{\hskip 16 true pt\flip{$[{4}{1}:-]{\times}[{3}{1}]$}\hss}\cg%
\cx{\hskip 16 true pt\flip{$[{3}{2}:-]{\times}[{3}{1}]$}\hss}\cg%
\cx{\hskip 16 true pt\flip{$[{3}{1^{2}}:-]{\times}[{3}{1}]$}\hss}\cg%
\cx{\hskip 16 true pt\flip{$[{2^{2}}{1}:-]{\times}[{3}{1}]$}\hss}\cg%
\eol}}\rg%
%
%
\rx{\vss\hfull{%
\rlx{\hss{$70_{y}$}}\cg%
\e{0}%
\e{0}%
\e{0}%
\e{0}%
\e{0}%
\e{0}%
\e{0}%
\e{0}%
\e{0}%
\e{0}%
\e{0}%
\e{1}%
\e{0}%
\e{0}%
\e{0}%
\e{0}%
\e{0}%
\e{0}%
\e{0}%
\e{0}%
\e{0}%
\e{0}%
\e{0}%
\eol}\vss}\rg%
%
%
\rx{\vss\hfull{%
\rlx{\hss{$1134_{y}$}}\cg%
\e{0}%
\e{0}%
\e{0}%
\e{0}%
\e{0}%
\e{0}%
\e{0}%
\e{0}%
\e{0}%
\e{0}%
\e{1}%
\e{0}%
\e{0}%
\e{1}%
\e{0}%
\e{1}%
\e{0}%
\e{0}%
\e{0}%
\e{0}%
\e{0}%
\e{0}%
\e{1}%
\eol}\vss}\rg%
%
%
\rx{\vss\hfull{%
\rlx{\hss{$1680_{y}$}}\cg%
\e{0}%
\e{0}%
\e{0}%
\e{0}%
\e{0}%
\e{0}%
\e{0}%
\e{0}%
\e{0}%
\e{0}%
\e{2}%
\e{2}%
\e{0}%
\e{0}%
\e{0}%
\e{1}%
\e{1}%
\e{1}%
\e{0}%
\e{0}%
\e{0}%
\e{1}%
\e{0}%
\eol}\vss}\rg%
%
%
\rx{\vss\hfull{%
\rlx{\hss{$168_{y}$}}\cg%
\e{0}%
\e{1}%
\e{0}%
\e{0}%
\e{0}%
\e{0}%
\e{0}%
\e{0}%
\e{0}%
\e{0}%
\e{0}%
\e{0}%
\e{0}%
\e{0}%
\e{0}%
\e{0}%
\e{0}%
\e{0}%
\e{0}%
\e{0}%
\e{0}%
\e{0}%
\e{0}%
\eol}\vss}\rg%
%
%
\rx{\vss\hfull{%
\rlx{\hss{$420_{y}$}}\cg%
\e{0}%
\e{0}%
\e{0}%
\e{0}%
\e{0}%
\e{0}%
\e{0}%
\e{0}%
\e{0}%
\e{0}%
\e{0}%
\e{0}%
\e{1}%
\e{0}%
\e{0}%
\e{0}%
\e{0}%
\e{0}%
\e{0}%
\e{0}%
\e{1}%
\e{0}%
\e{0}%
\eol}\vss}\rg%
%
%
\rx{\vss\hfull{%
\rlx{\hss{$3150_{y}$}}\cg%
\e{0}%
\e{0}%
\e{0}%
\e{0}%
\e{0}%
\e{0}%
\e{0}%
\e{0}%
\e{1}%
\e{1}%
\e{0}%
\e{0}%
\e{1}%
\e{1}%
\e{1}%
\e{2}%
\e{0}%
\e{0}%
\e{0}%
\e{0}%
\e{2}%
\e{1}%
\e{2}%
\eol}\vss}\rg%
%
%
\rx{\vss\hfull{%
\rlx{\hss{$4200_{y}$}}\cg%
\e{0}%
\e{1}%
\e{1}%
\e{1}%
\e{1}%
\e{0}%
\e{0}%
\e{0}%
\e{1}%
\e{1}%
\e{1}%
\e{0}%
\e{1}%
\e{1}%
\e{2}%
\e{1}%
\e{0}%
\e{0}%
\e{0}%
\e{1}%
\e{2}%
\e{1}%
\e{2}%
\eol}\vss}\rg%
%
%
\rx{\vss\hfull{%
\rlx{\hss{$2688_{y}$}}\cg%
\e{0}%
\e{0}%
\e{1}%
\e{1}%
\e{1}%
\e{0}%
\e{0}%
\e{0}%
\e{1}%
\e{1}%
\e{1}%
\e{0}%
\e{0}%
\e{0}%
\e{1}%
\e{1}%
\e{0}%
\e{0}%
\e{0}%
\e{0}%
\e{1}%
\e{1}%
\e{1}%
\eol}\vss}\rg%
%
%
\rx{\vss\hfull{%
\rlx{\hss{$2100_{y}$}}\cg%
\e{0}%
\e{0}%
\e{0}%
\e{1}%
\e{0}%
\e{2}%
\e{1}%
\e{0}%
\e{1}%
\e{0}%
\e{1}%
\e{1}%
\e{0}%
\e{0}%
\e{0}%
\e{1}%
\e{1}%
\e{1}%
\e{0}%
\e{0}%
\e{0}%
\e{1}%
\e{0}%
\eol}\vss}\rg%
%
%
\rx{\vss\hfull{%
\rlx{\hss{$1400_{y}$}}\cg%
\e{0}%
\e{0}%
\e{0}%
\e{0}%
\e{0}%
\e{0}%
\e{0}%
\e{0}%
\e{0}%
\e{0}%
\e{1}%
\e{1}%
\e{0}%
\e{0}%
\e{1}%
\e{0}%
\e{1}%
\e{1}%
\e{0}%
\e{1}%
\e{0}%
\e{1}%
\e{0}%
\eol}\vss}\rg%
%
%
\rx{\vss\hfull{%
\rlx{\hss{$4536_{y}$}}\cg%
\e{0}%
\e{0}%
\e{0}%
\e{0}%
\e{0}%
\e{0}%
\e{0}%
\e{0}%
\e{1}%
\e{1}%
\e{2}%
\e{1}%
\e{0}%
\e{0}%
\e{2}%
\e{2}%
\e{2}%
\e{2}%
\e{0}%
\e{1}%
\e{1}%
\e{3}%
\e{1}%
\eol}\vss}\rg%
%
%
\rx{\vss\hfull{%
\rlx{\hss{$5670_{y}$}}\cg%
\e{0}%
\e{0}%
\e{0}%
\e{0}%
\e{0}%
\e{0}%
\e{0}%
\e{0}%
\e{1}%
\e{1}%
\e{3}%
\e{1}%
\e{0}%
\e{1}%
\e{2}%
\e{3}%
\e{2}%
\e{2}%
\e{0}%
\e{1}%
\e{1}%
\e{3}%
\e{2}%
\eol}\vss}\rg%
%
%
\rx{\vss\hfull{%
\rlx{\hss{$4480_{y}$}}\cg%
\e{0}%
\e{0}%
\e{0}%
\e{0}%
\e{0}%
\e{0}%
\e{0}%
\e{0}%
\e{1}%
\e{1}%
\e{1}%
\e{0}%
\e{1}%
\e{1}%
\e{2}%
\e{2}%
\e{1}%
\e{1}%
\e{0}%
\e{1}%
\e{2}%
\e{2}%
\e{2}%
\eol}\vss}\rg%
%
%
\rx{\vss\hfull{%
\rlx{\hss{$8_{z}$}}\cg%
\e{0}%
\e{0}%
\e{0}%
\e{0}%
\e{0}%
\e{0}%
\e{0}%
\e{1}%
\e{0}%
\e{0}%
\e{0}%
\e{0}%
\e{0}%
\e{0}%
\e{0}%
\e{0}%
\e{0}%
\e{0}%
\e{1}%
\e{0}%
\e{0}%
\e{0}%
\e{0}%
\eol}\vss}\rg%
%
%
\rx{\vss\hfull{%
\rlx{\hss{$56_{z}$}}\cg%
\e{0}%
\e{0}%
\e{0}%
\e{0}%
\e{0}%
\e{0}%
\e{0}%
\e{0}%
\e{0}%
\e{0}%
\e{0}%
\e{0}%
\e{0}%
\e{0}%
\e{0}%
\e{0}%
\e{1}%
\e{0}%
\e{0}%
\e{0}%
\e{0}%
\e{0}%
\e{0}%
\eol}\vss}\rg%
%
%
\rx{\vss\hfull{%
\rlx{\hss{$160_{z}$}}\cg%
\e{0}%
\e{0}%
\e{0}%
\e{0}%
\e{0}%
\e{0}%
\e{0}%
\e{1}%
\e{1}%
\e{0}%
\e{0}%
\e{0}%
\e{0}%
\e{0}%
\e{1}%
\e{0}%
\e{0}%
\e{0}%
\e{1}%
\e{1}%
\e{0}%
\e{0}%
\e{0}%
\eol}\vss}\rg%
%
%
\rx{\vss\hfull{%
\rlx{\hss{$112_{z}$}}\cg%
\e{1}%
\e{0}%
\e{0}%
\e{0}%
\e{0}%
\e{0}%
\e{0}%
\e{2}%
\e{1}%
\e{0}%
\e{0}%
\e{0}%
\e{1}%
\e{0}%
\e{0}%
\e{0}%
\e{0}%
\e{0}%
\e{2}%
\e{1}%
\e{0}%
\e{0}%
\e{0}%
\eol}\vss}\rg%
%
%
\rx{\vss\hfull{%
\rlx{\hss{$840_{z}$}}\cg%
\e{0}%
\e{0}%
\e{0}%
\e{0}%
\e{1}%
\e{0}%
\e{0}%
\e{1}%
\e{1}%
\e{1}%
\e{0}%
\e{0}%
\e{0}%
\e{0}%
\e{1}%
\e{1}%
\e{0}%
\e{0}%
\e{1}%
\e{1}%
\e{1}%
\e{0}%
\e{0}%
\eol}\vss}\rg%
%
%
\rx{\vss\hfull{%
\rlx{\hss{$1296_{z}$}}\cg%
\e{0}%
\e{0}%
\e{0}%
\e{1}%
\e{0}%
\e{0}%
\e{0}%
\e{0}%
\e{1}%
\e{0}%
\e{1}%
\e{0}%
\e{0}%
\e{1}%
\e{2}%
\e{1}%
\e{2}%
\e{1}%
\e{0}%
\e{1}%
\e{0}%
\e{1}%
\e{0}%
\eol}\vss}\rg%
%
%
\rx{\vss\hfull{%
\rlx{\hss{$1400_{z}$}}\cg%
\e{1}%
\e{1}%
\e{1}%
\e{0}%
\e{0}%
\e{0}%
\e{0}%
\e{2}%
\e{3}%
\e{1}%
\e{1}%
\e{0}%
\e{3}%
\e{1}%
\e{2}%
\e{1}%
\e{0}%
\e{0}%
\e{2}%
\e{3}%
\e{1}%
\e{1}%
\e{0}%
\eol}\vss}\rg%
%
%
\rx{\vss\hfull{%
\rlx{\hss{$1008_{z}$}}\cg%
\e{0}%
\e{1}%
\e{0}%
\e{0}%
\e{0}%
\e{0}%
\e{0}%
\e{1}%
\e{2}%
\e{0}%
\e{1}%
\e{0}%
\e{1}%
\e{1}%
\e{2}%
\e{1}%
\e{1}%
\e{0}%
\e{1}%
\e{2}%
\e{0}%
\e{1}%
\e{0}%
\eol}\vss}\rg%
%
%
\rx{\vss\hfull{%
\rlx{\hss{$560_{z}$}}\cg%
\e{1}%
\e{1}%
\e{0}%
\e{0}%
\e{0}%
\e{0}%
\e{0}%
\e{3}%
\e{2}%
\e{1}%
\e{0}%
\e{0}%
\e{2}%
\e{1}%
\e{1}%
\e{0}%
\e{0}%
\e{0}%
\e{3}%
\e{2}%
\e{1}%
\e{0}%
\e{0}%
\eol}\vss}\rg%
%
%
\rx{\vss\hfull{%
\rlx{\hss{$1400_{zz}$}}\cg%
\e{1}%
\e{1}%
\e{1}%
\e{0}%
\e{0}%
\e{0}%
\e{0}%
\e{1}%
\e{2}%
\e{0}%
\e{0}%
\e{0}%
\e{3}%
\e{1}%
\e{1}%
\e{0}%
\e{0}%
\e{0}%
\e{0}%
\e{1}%
\e{2}%
\e{1}%
\e{0}%
\eol}\vss}\rg%
%
%
\rx{\vss\hfull{%
\rlx{\hss{$4200_{z}$}}\cg%
\e{0}%
\e{1}%
\e{0}%
\e{1}%
\e{0}%
\e{0}%
\e{0}%
\e{1}%
\e{2}%
\e{0}%
\e{1}%
\e{0}%
\e{1}%
\e{2}%
\e{3}%
\e{2}%
\e{1}%
\e{0}%
\e{0}%
\e{1}%
\e{2}%
\e{2}%
\e{1}%
\eol}\vss}\rg%
%
%
\rx{\vss\hfull{%
\rlx{\hss{$400_{z}$}}\cg%
\e{1}%
\e{1}%
\e{0}%
\e{0}%
\e{0}%
\e{0}%
\e{0}%
\e{2}%
\e{1}%
\e{0}%
\e{0}%
\e{0}%
\e{2}%
\e{1}%
\e{0}%
\e{0}%
\e{0}%
\e{0}%
\e{1}%
\e{1}%
\e{1}%
\e{0}%
\e{0}%
\eol}\vss}\rg%
%
%
\rx{\vss\hfull{%
\rlx{\hss{$3240_{z}$}}\cg%
\e{1}%
\e{2}%
\e{1}%
\e{1}%
\e{0}%
\e{0}%
\e{0}%
\e{3}%
\e{4}%
\e{2}%
\e{1}%
\e{0}%
\e{4}%
\e{3}%
\e{4}%
\e{2}%
\e{0}%
\e{0}%
\e{2}%
\e{4}%
\e{3}%
\e{1}%
\e{1}%
\eol}\vss}\rg%
%
%
\rx{\vss\hfull{%
\rlx{\hss{$4536_{z}$}}\cg%
\e{0}%
\e{1}%
\e{2}%
\e{0}%
\e{1}%
\e{0}%
\e{0}%
\e{2}%
\e{3}%
\e{3}%
\e{1}%
\e{0}%
\e{3}%
\e{1}%
\e{3}%
\e{2}%
\e{0}%
\e{0}%
\e{1}%
\e{3}%
\e{4}%
\e{1}%
\e{2}%
\eol}\vss}\rg%
%
%
\rx{\vss\hfull{%
\rlx{\hss{$2400_{z}$}}\cg%
\e{0}%
\e{0}%
\e{0}%
\e{1}%
\e{0}%
\e{1}%
\e{0}%
\e{0}%
\e{0}%
\e{0}%
\e{1}%
\e{1}%
\e{0}%
\e{1}%
\e{1}%
\e{2}%
\e{3}%
\e{1}%
\e{0}%
\e{0}%
\e{0}%
\e{1}%
\e{0}%
\eol}\vss}\rg%
%
%
\rx{\vss\hfull{%
\rlx{\hss{$3360_{z}$}}\cg%
\e{0}%
\e{1}%
\e{1}%
\e{1}%
\e{0}%
\e{0}%
\e{0}%
\e{0}%
\e{2}%
\e{1}%
\e{1}%
\e{0}%
\e{2}%
\e{2}%
\e{3}%
\e{2}%
\e{1}%
\e{0}%
\e{0}%
\e{1}%
\e{2}%
\e{2}%
\e{1}%
\eol}\vss}\rg%
%
%
\rx{\vss\hfull{%
\rlx{\hss{$2800_{z}$}}\cg%
\e{0}%
\e{1}%
\e{0}%
\e{1}%
\e{0}%
\e{0}%
\e{0}%
\e{0}%
\e{1}%
\e{0}%
\e{2}%
\e{1}%
\e{1}%
\e{2}%
\e{2}%
\e{2}%
\e{2}%
\e{1}%
\e{0}%
\e{1}%
\e{0}%
\e{2}%
\e{0}%
\eol}\vss}\rg%
%
%
\rx{\vss\hfull{%
\rlx{\hss{$4096_{z}$}}\cg%
\e{0}%
\e{1}%
\e{1}%
\e{1}%
\e{1}%
\e{0}%
\e{0}%
\e{1}%
\e{3}%
\e{2}%
\e{2}%
\e{0}%
\e{2}%
\e{2}%
\e{4}%
\e{3}%
\e{1}%
\e{1}%
\e{1}%
\e{3}%
\e{2}%
\e{2}%
\e{1}%
\eol}\vss}\rg%
%
%
\rx{\vss\hfull{%
\rlx{\hss{$5600_{z}$}}\cg%
\e{0}%
\e{0}%
\e{1}%
\e{1}%
\e{1}%
\e{1}%
\e{0}%
\e{0}%
\e{2}%
\e{1}%
\e{3}%
\e{1}%
\e{1}%
\e{1}%
\e{3}%
\e{4}%
\e{1}%
\e{2}%
\e{0}%
\e{2}%
\e{1}%
\e{3}%
\e{1}%
\eol}\vss}\rg%
%
%
\rx{\vss\hfull{%
\rlx{\hss{$448_{z}$}}\cg%
\e{1}%
\e{0}%
\e{1}%
\e{0}%
\e{0}%
\e{0}%
\e{0}%
\e{1}%
\e{1}%
\e{1}%
\e{0}%
\e{0}%
\e{2}%
\e{0}%
\e{0}%
\e{0}%
\e{0}%
\e{0}%
\e{1}%
\e{1}%
\e{1}%
\e{0}%
\e{0}%
\eol}\vss}\rg%
\eop
\eject
\tablecont%
%
%
%
%
%
%
\rowpts=18 true pt%
\colpts=18 true pt%
\rowlabpts=40 true pt%
\collabpts=85 true pt%
\clx{\vss\hfull{%
\rlx{\hss{$ $}}\cg%
\cx{\hskip 16 true pt\flip{$[{5}:-]{\times}[{4}]$}\hss}\cg%
\cx{\hskip 16 true pt\flip{$[{4}{1}:-]{\times}[{4}]$}\hss}\cg%
\cx{\hskip 16 true pt\flip{$[{3}{2}:-]{\times}[{4}]$}\hss}\cg%
\cx{\hskip 16 true pt\flip{$[{3}{1^{2}}:-]{\times}[{4}]$}\hss}\cg%
\cx{\hskip 16 true pt\flip{$[{2^{2}}{1}:-]{\times}[{4}]$}\hss}\cg%
\cx{\hskip 16 true pt\flip{$[{2}{1^{3}}:-]{\times}[{4}]$}\hss}\cg%
\cx{\hskip 16 true pt\flip{$[{1^{5}}:-]{\times}[{4}]$}\hss}\cg%
\cx{\hskip 16 true pt\flip{$[{4}:{1}]{\times}[{4}]$}\hss}\cg%
\cx{\hskip 16 true pt\flip{$[{3}{1}:{1}]{\times}[{4}]$}\hss}\cg%
\cx{\hskip 16 true pt\flip{$[{2^{2}}:{1}]{\times}[{4}]$}\hss}\cg%
\cx{\hskip 16 true pt\flip{$[{2}{1^{2}}:{1}]{\times}[{4}]$}\hss}\cg%
\cx{\hskip 16 true pt\flip{$[{1^{4}}:{1}]{\times}[{4}]$}\hss}\cg%
\cx{\hskip 16 true pt\flip{$[{3}:{2}]{\times}[{4}]$}\hss}\cg%
\cx{\hskip 16 true pt\flip{$[{3}:{1^{2}}]{\times}[{4}]$}\hss}\cg%
\cx{\hskip 16 true pt\flip{$[{2}{1}:{2}]{\times}[{4}]$}\hss}\cg%
\cx{\hskip 16 true pt\flip{$[{2}{1}:{1^{2}}]{\times}[{4}]$}\hss}\cg%
\cx{\hskip 16 true pt\flip{$[{1^{3}}:{2}]{\times}[{4}]$}\hss}\cg%
\cx{\hskip 16 true pt\flip{$[{1^{3}}:{1^{2}}]{\times}[{4}]$}\hss}\cg%
\cx{\hskip 16 true pt\flip{$[{5}:-]{\times}[{3}{1}]$}\hss}\cg%
\cx{\hskip 16 true pt\flip{$[{4}{1}:-]{\times}[{3}{1}]$}\hss}\cg%
\cx{\hskip 16 true pt\flip{$[{3}{2}:-]{\times}[{3}{1}]$}\hss}\cg%
\cx{\hskip 16 true pt\flip{$[{3}{1^{2}}:-]{\times}[{3}{1}]$}\hss}\cg%
\cx{\hskip 16 true pt\flip{$[{2^{2}}{1}:-]{\times}[{3}{1}]$}\hss}\cg%
\eol}}\rg%
%
%
\rx{\vss\hfull{%
\rlx{\hss{$448_{w}$}}\cg%
\e{0}%
\e{0}%
\e{0}%
\e{0}%
\e{0}%
\e{1}%
\e{0}%
\e{0}%
\e{0}%
\e{0}%
\e{0}%
\e{0}%
\e{0}%
\e{0}%
\e{0}%
\e{0}%
\e{1}%
\e{1}%
\e{0}%
\e{0}%
\e{0}%
\e{0}%
\e{0}%
\eol}\vss}\rg%
%
%
\rx{\vss\hfull{%
\rlx{\hss{$1344_{w}$}}\cg%
\e{0}%
\e{0}%
\e{0}%
\e{0}%
\e{0}%
\e{0}%
\e{0}%
\e{1}%
\e{1}%
\e{0}%
\e{0}%
\e{0}%
\e{0}%
\e{0}%
\e{1}%
\e{0}%
\e{0}%
\e{0}%
\e{0}%
\e{1}%
\e{1}%
\e{0}%
\e{1}%
\eol}\vss}\rg%
%
%
\rx{\vss\hfull{%
\rlx{\hss{$5600_{w}$}}\cg%
\e{0}%
\e{0}%
\e{0}%
\e{1}%
\e{1}%
\e{1}%
\e{0}%
\e{0}%
\e{1}%
\e{1}%
\e{2}%
\e{0}%
\e{0}%
\e{1}%
\e{2}%
\e{3}%
\e{2}%
\e{2}%
\e{0}%
\e{1}%
\e{1}%
\e{2}%
\e{2}%
\eol}\vss}\rg%
%
%
\rx{\vss\hfull{%
\rlx{\hss{$2016_{w}$}}\cg%
\e{0}%
\e{0}%
\e{1}%
\e{0}%
\e{0}%
\e{0}%
\e{0}%
\e{0}%
\e{1}%
\e{1}%
\e{0}%
\e{0}%
\e{1}%
\e{0}%
\e{1}%
\e{0}%
\e{0}%
\e{0}%
\e{0}%
\e{0}%
\e{2}%
\e{1}%
\e{1}%
\eol}\vss}\rg%
%
%
\rx{\vss\hfull{%
\rlx{\hss{$7168_{w}$}}\cg%
\e{0}%
\e{0}%
\e{1}%
\e{1}%
\e{1}%
\e{0}%
\e{0}%
\e{0}%
\e{2}%
\e{2}%
\e{2}%
\e{0}%
\e{1}%
\e{1}%
\e{3}%
\e{3}%
\e{1}%
\e{1}%
\e{0}%
\e{1}%
\e{3}%
\e{3}%
\e{3}%
\eol}\vss}\rg%
%
%
%
%
%
%
\rowpts=18 true pt%
\colpts=18 true pt%
\rowlabpts=40 true pt%
\collabpts=85 true pt%
\clx{\vss\hfull{%
\rlx{\hss{$ $}}\cg%
\cx{\hskip 16 true pt\flip{$[{2}{1^{3}}:-]{\times}[{3}{1}]$}\hss}\cg%
\cx{\hskip 16 true pt\flip{$[{1^{5}}:-]{\times}[{3}{1}]$}\hss}\cg%
\cx{\hskip 16 true pt\flip{$[{4}:{1}]{\times}[{3}{1}]$}\hss}\cg%
\cx{\hskip 16 true pt\flip{$[{3}{1}:{1}]{\times}[{3}{1}]$}\hss}\cg%
\cx{\hskip 16 true pt\flip{$[{2^{2}}:{1}]{\times}[{3}{1}]$}\hss}\cg%
\cx{\hskip 16 true pt\flip{$[{2}{1^{2}}:{1}]{\times}[{3}{1}]$}\hss}\cg%
\cx{\hskip 16 true pt\flip{$[{1^{4}}:{1}]{\times}[{3}{1}]$}\hss}\cg%
\cx{\hskip 16 true pt\flip{$[{3}:{2}]{\times}[{3}{1}]$}\hss}\cg%
\cx{\hskip 16 true pt\flip{$[{3}:{1^{2}}]{\times}[{3}{1}]$}\hss}\cg%
\cx{\hskip 16 true pt\flip{$[{2}{1}:{2}]{\times}[{3}{1}]$}\hss}\cg%
\cx{\hskip 16 true pt\flip{$[{2}{1}:{1^{2}}]{\times}[{3}{1}]$}\hss}\cg%
\cx{\hskip 16 true pt\flip{$[{1^{3}}:{2}]{\times}[{3}{1}]$}\hss}\cg%
\cx{\hskip 16 true pt\flip{$[{1^{3}}:{1^{2}}]{\times}[{3}{1}]$}\hss}\cg%
\cx{\hskip 16 true pt\flip{$[{5}:-]{\times}[{2^{2}}]$}\hss}\cg%
\cx{\hskip 16 true pt\flip{$[{4}{1}:-]{\times}[{2^{2}}]$}\hss}\cg%
\cx{\hskip 16 true pt\flip{$[{3}{2}:-]{\times}[{2^{2}}]$}\hss}\cg%
\cx{\hskip 16 true pt\flip{$[{3}{1^{2}}:-]{\times}[{2^{2}}]$}\hss}\cg%
\cx{\hskip 16 true pt\flip{$[{2^{2}}{1}:-]{\times}[{2^{2}}]$}\hss}\cg%
\cx{\hskip 16 true pt\flip{$[{2}{1^{3}}:-]{\times}[{2^{2}}]$}\hss}\cg%
\cx{\hskip 16 true pt\flip{$[{1^{5}}:-]{\times}[{2^{2}}]$}\hss}\cg%
\cx{\hskip 16 true pt\flip{$[{4}:{1}]{\times}[{2^{2}}]$}\hss}\cg%
\cx{\hskip 16 true pt\flip{$[{3}{1}:{1}]{\times}[{2^{2}}]$}\hss}\cg%
\cx{\hskip 16 true pt\flip{$[{2^{2}}:{1}]{\times}[{2^{2}}]$}\hss}\cg%
\eol}}\rg%
%
%
\rx{\vss\hfull{%
\rlx{\hss{$1_{x}$}}\cg%
\e{0}%
\e{0}%
\e{0}%
\e{0}%
\e{0}%
\e{0}%
\e{0}%
\e{0}%
\e{0}%
\e{0}%
\e{0}%
\e{0}%
\e{0}%
\e{0}%
\e{0}%
\e{0}%
\e{0}%
\e{0}%
\e{0}%
\e{0}%
\e{0}%
\e{0}%
\e{0}%
\eol}\vss}\rg%
%
%
\rx{\vss\hfull{%
\rlx{\hss{$28_{x}$}}\cg%
\e{0}%
\e{0}%
\e{1}%
\e{0}%
\e{0}%
\e{0}%
\e{0}%
\e{0}%
\e{0}%
\e{0}%
\e{0}%
\e{0}%
\e{0}%
\e{0}%
\e{0}%
\e{0}%
\e{0}%
\e{0}%
\e{0}%
\e{0}%
\e{0}%
\e{0}%
\e{0}%
\eol}\vss}\rg%
%
%
\rx{\vss\hfull{%
\rlx{\hss{$35_{x}$}}\cg%
\e{0}%
\e{0}%
\e{1}%
\e{0}%
\e{0}%
\e{0}%
\e{0}%
\e{0}%
\e{0}%
\e{0}%
\e{0}%
\e{0}%
\e{0}%
\e{1}%
\e{0}%
\e{0}%
\e{0}%
\e{0}%
\e{0}%
\e{0}%
\e{0}%
\e{0}%
\e{0}%
\eol}\vss}\rg%
%
%
\rx{\vss\hfull{%
\rlx{\hss{$84_{x}$}}\cg%
\e{0}%
\e{0}%
\e{1}%
\e{0}%
\e{0}%
\e{0}%
\e{0}%
\e{1}%
\e{0}%
\e{0}%
\e{0}%
\e{0}%
\e{0}%
\e{1}%
\e{1}%
\e{0}%
\e{0}%
\e{0}%
\e{0}%
\e{0}%
\e{0}%
\e{0}%
\e{0}%
\eol}\vss}\rg%
%
%
\rx{\vss\hfull{%
\rlx{\hss{$50_{x}$}}\cg%
\e{0}%
\e{0}%
\e{0}%
\e{0}%
\e{0}%
\e{0}%
\e{0}%
\e{1}%
\e{0}%
\e{0}%
\e{0}%
\e{0}%
\e{0}%
\e{0}%
\e{0}%
\e{1}%
\e{0}%
\e{0}%
\e{0}%
\e{0}%
\e{0}%
\e{0}%
\e{0}%
\eol}\vss}\rg%
%
%
\rx{\vss\hfull{%
\rlx{\hss{$350_{x}$}}\cg%
\e{0}%
\e{0}%
\e{1}%
\e{1}%
\e{0}%
\e{0}%
\e{0}%
\e{0}%
\e{1}%
\e{1}%
\e{0}%
\e{1}%
\e{0}%
\e{0}%
\e{0}%
\e{0}%
\e{0}%
\e{0}%
\e{0}%
\e{0}%
\e{1}%
\e{0}%
\e{0}%
\eol}\vss}\rg%
%
%
\rx{\vss\hfull{%
\rlx{\hss{$300_{x}$}}\cg%
\e{0}%
\e{0}%
\e{1}%
\e{1}%
\e{0}%
\e{0}%
\e{0}%
\e{1}%
\e{0}%
\e{1}%
\e{0}%
\e{0}%
\e{0}%
\e{1}%
\e{1}%
\e{0}%
\e{0}%
\e{0}%
\e{0}%
\e{0}%
\e{1}%
\e{0}%
\e{0}%
\eol}\vss}\rg%
%
%
\rx{\vss\hfull{%
\rlx{\hss{$567_{x}$}}\cg%
\e{0}%
\e{0}%
\e{4}%
\e{2}%
\e{0}%
\e{0}%
\e{0}%
\e{2}%
\e{1}%
\e{1}%
\e{0}%
\e{0}%
\e{0}%
\e{1}%
\e{1}%
\e{0}%
\e{0}%
\e{0}%
\e{0}%
\e{0}%
\e{1}%
\e{0}%
\e{0}%
\eol}\vss}\rg%
%
%
\rx{\vss\hfull{%
\rlx{\hss{$210_{x}$}}\cg%
\e{0}%
\e{0}%
\e{2}%
\e{1}%
\e{0}%
\e{0}%
\e{0}%
\e{1}%
\e{0}%
\e{0}%
\e{0}%
\e{0}%
\e{0}%
\e{1}%
\e{0}%
\e{0}%
\e{0}%
\e{0}%
\e{0}%
\e{0}%
\e{0}%
\e{0}%
\e{0}%
\eol}\vss}\rg%
%
%
\rx{\vss\hfull{%
\rlx{\hss{$840_{x}$}}\cg%
\e{0}%
\e{0}%
\e{1}%
\e{1}%
\e{2}%
\e{0}%
\e{0}%
\e{1}%
\e{0}%
\e{2}%
\e{1}%
\e{0}%
\e{0}%
\e{0}%
\e{1}%
\e{1}%
\e{0}%
\e{0}%
\e{0}%
\e{0}%
\e{0}%
\e{1}%
\e{1}%
\eol}\vss}\rg%
%
%
\rx{\vss\hfull{%
\rlx{\hss{$700_{x}$}}\cg%
\e{0}%
\e{0}%
\e{3}%
\e{2}%
\e{1}%
\e{0}%
\e{0}%
\e{3}%
\e{1}%
\e{1}%
\e{0}%
\e{0}%
\e{0}%
\e{1}%
\e{1}%
\e{1}%
\e{0}%
\e{0}%
\e{0}%
\e{0}%
\e{1}%
\e{0}%
\e{0}%
\eol}\vss}\rg%
%
%
\rx{\vss\hfull{%
\rlx{\hss{$175_{x}$}}\cg%
\e{0}%
\e{0}%
\e{0}%
\e{1}%
\e{0}%
\e{0}%
\e{0}%
\e{1}%
\e{0}%
\e{0}%
\e{0}%
\e{0}%
\e{0}%
\e{0}%
\e{0}%
\e{0}%
\e{0}%
\e{0}%
\e{0}%
\e{0}%
\e{0}%
\e{0}%
\e{0}%
\eol}\vss}\rg%
%
%
\rx{\vss\hfull{%
\rlx{\hss{$1400_{x}$}}\cg%
\e{0}%
\e{0}%
\e{3}%
\e{4}%
\e{1}%
\e{1}%
\e{0}%
\e{4}%
\e{2}%
\e{2}%
\e{1}%
\e{0}%
\e{0}%
\e{0}%
\e{0}%
\e{0}%
\e{0}%
\e{0}%
\e{0}%
\e{0}%
\e{1}%
\e{1}%
\e{0}%
\eol}\vss}\rg%
%
%
\rx{\vss\hfull{%
\rlx{\hss{$1050_{x}$}}\cg%
\e{0}%
\e{0}%
\e{2}%
\e{2}%
\e{1}%
\e{0}%
\e{0}%
\e{4}%
\e{2}%
\e{2}%
\e{1}%
\e{0}%
\e{0}%
\e{0}%
\e{1}%
\e{2}%
\e{1}%
\e{1}%
\e{0}%
\e{0}%
\e{1}%
\e{1}%
\e{0}%
\eol}\vss}\rg%
%
%
\rx{\vss\hfull{%
\rlx{\hss{$1575_{x}$}}\cg%
\e{0}%
\e{0}%
\e{4}%
\e{4}%
\e{1}%
\e{1}%
\e{0}%
\e{3}%
\e{3}%
\e{3}%
\e{1}%
\e{1}%
\e{0}%
\e{0}%
\e{0}%
\e{0}%
\e{0}%
\e{0}%
\e{0}%
\e{0}%
\e{2}%
\e{1}%
\e{0}%
\eol}\vss}\rg%
%
%
\rx{\vss\hfull{%
\rlx{\hss{$1344_{x}$}}\cg%
\e{0}%
\e{0}%
\e{4}%
\e{3}%
\e{1}%
\e{0}%
\e{0}%
\e{5}%
\e{2}%
\e{3}%
\e{1}%
\e{0}%
\e{0}%
\e{2}%
\e{3}%
\e{1}%
\e{1}%
\e{0}%
\e{0}%
\e{0}%
\e{2}%
\e{1}%
\e{0}%
\eol}\vss}\rg%
%
%
\rx{\vss\hfull{%
\rlx{\hss{$2100_{x}$}}\cg%
\e{0}%
\e{0}%
\e{1}%
\e{3}%
\e{0}%
\e{2}%
\e{0}%
\e{2}%
\e{3}%
\e{4}%
\e{3}%
\e{3}%
\e{1}%
\e{0}%
\e{1}%
\e{0}%
\e{1}%
\e{0}%
\e{0}%
\e{0}%
\e{2}%
\e{2}%
\e{0}%
\eol}\vss}\rg%
%
%
\rx{\vss\hfull{%
\rlx{\hss{$2268_{x}$}}\cg%
\e{0}%
\e{0}%
\e{3}%
\e{5}%
\e{1}%
\e{2}%
\e{0}%
\e{4}%
\e{3}%
\e{5}%
\e{2}%
\e{1}%
\e{0}%
\e{1}%
\e{1}%
\e{0}%
\e{0}%
\e{0}%
\e{0}%
\e{0}%
\e{2}%
\e{2}%
\e{1}%
\eol}\vss}\rg%
%
%
\rx{\vss\hfull{%
\rlx{\hss{$525_{x}$}}\cg%
\e{0}%
\e{0}%
\e{1}%
\e{1}%
\e{0}%
\e{0}%
\e{0}%
\e{2}%
\e{1}%
\e{1}%
\e{1}%
\e{0}%
\e{0}%
\e{0}%
\e{1}%
\e{0}%
\e{1}%
\e{0}%
\e{0}%
\e{0}%
\e{0}%
\e{1}%
\e{0}%
\eol}\vss}\rg%
%
%
\rx{\vss\hfull{%
\rlx{\hss{$700_{xx}$}}\cg%
\e{0}%
\e{0}%
\e{0}%
\e{1}%
\e{0}%
\e{0}%
\e{0}%
\e{2}%
\e{1}%
\e{2}%
\e{1}%
\e{0}%
\e{0}%
\e{0}%
\e{0}%
\e{2}%
\e{1}%
\e{1}%
\e{1}%
\e{0}%
\e{0}%
\e{1}%
\e{1}%
\eol}\vss}\rg%
%
%
\rx{\vss\hfull{%
\rlx{\hss{$972_{x}$}}\cg%
\e{0}%
\e{0}%
\e{1}%
\e{1}%
\e{1}%
\e{0}%
\e{0}%
\e{3}%
\e{1}%
\e{3}%
\e{1}%
\e{0}%
\e{0}%
\e{1}%
\e{2}%
\e{2}%
\e{1}%
\e{1}%
\e{0}%
\e{0}%
\e{1}%
\e{1}%
\e{1}%
\eol}\vss}\rg%
\eop
\eject
\tablecont%
%
%
%
%
%
%
\rowpts=18 true pt%
\colpts=18 true pt%
\rowlabpts=40 true pt%
\collabpts=85 true pt%
\clx{\vss\hfull{%
\rlx{\hss{$ $}}\cg%
\cx{\hskip 16 true pt\flip{$[{2}{1^{3}}:-]{\times}[{3}{1}]$}\hss}\cg%
\cx{\hskip 16 true pt\flip{$[{1^{5}}:-]{\times}[{3}{1}]$}\hss}\cg%
\cx{\hskip 16 true pt\flip{$[{4}:{1}]{\times}[{3}{1}]$}\hss}\cg%
\cx{\hskip 16 true pt\flip{$[{3}{1}:{1}]{\times}[{3}{1}]$}\hss}\cg%
\cx{\hskip 16 true pt\flip{$[{2^{2}}:{1}]{\times}[{3}{1}]$}\hss}\cg%
\cx{\hskip 16 true pt\flip{$[{2}{1^{2}}:{1}]{\times}[{3}{1}]$}\hss}\cg%
\cx{\hskip 16 true pt\flip{$[{1^{4}}:{1}]{\times}[{3}{1}]$}\hss}\cg%
\cx{\hskip 16 true pt\flip{$[{3}:{2}]{\times}[{3}{1}]$}\hss}\cg%
\cx{\hskip 16 true pt\flip{$[{3}:{1^{2}}]{\times}[{3}{1}]$}\hss}\cg%
\cx{\hskip 16 true pt\flip{$[{2}{1}:{2}]{\times}[{3}{1}]$}\hss}\cg%
\cx{\hskip 16 true pt\flip{$[{2}{1}:{1^{2}}]{\times}[{3}{1}]$}\hss}\cg%
\cx{\hskip 16 true pt\flip{$[{1^{3}}:{2}]{\times}[{3}{1}]$}\hss}\cg%
\cx{\hskip 16 true pt\flip{$[{1^{3}}:{1^{2}}]{\times}[{3}{1}]$}\hss}\cg%
\cx{\hskip 16 true pt\flip{$[{5}:-]{\times}[{2^{2}}]$}\hss}\cg%
\cx{\hskip 16 true pt\flip{$[{4}{1}:-]{\times}[{2^{2}}]$}\hss}\cg%
\cx{\hskip 16 true pt\flip{$[{3}{2}:-]{\times}[{2^{2}}]$}\hss}\cg%
\cx{\hskip 16 true pt\flip{$[{3}{1^{2}}:-]{\times}[{2^{2}}]$}\hss}\cg%
\cx{\hskip 16 true pt\flip{$[{2^{2}}{1}:-]{\times}[{2^{2}}]$}\hss}\cg%
\cx{\hskip 16 true pt\flip{$[{2}{1^{3}}:-]{\times}[{2^{2}}]$}\hss}\cg%
\cx{\hskip 16 true pt\flip{$[{1^{5}}:-]{\times}[{2^{2}}]$}\hss}\cg%
\cx{\hskip 16 true pt\flip{$[{4}:{1}]{\times}[{2^{2}}]$}\hss}\cg%
\cx{\hskip 16 true pt\flip{$[{3}{1}:{1}]{\times}[{2^{2}}]$}\hss}\cg%
\cx{\hskip 16 true pt\flip{$[{2^{2}}:{1}]{\times}[{2^{2}}]$}\hss}\cg%
\eol}}\rg%
%
%
\rx{\vss\hfull{%
\rlx{\hss{$4096_{x}$}}\cg%
\e{0}%
\e{0}%
\e{4}%
\e{7}%
\e{3}%
\e{3}%
\e{0}%
\e{6}%
\e{5}%
\e{8}%
\e{5}%
\e{2}%
\e{1}%
\e{1}%
\e{2}%
\e{1}%
\e{1}%
\e{0}%
\e{0}%
\e{0}%
\e{3}%
\e{4}%
\e{1}%
\eol}\vss}\rg%
%
%
\rx{\vss\hfull{%
\rlx{\hss{$4200_{x}$}}\cg%
\e{0}%
\e{0}%
\e{3}%
\e{7}%
\e{4}%
\e{3}%
\e{0}%
\e{6}%
\e{4}%
\e{8}%
\e{5}%
\e{1}%
\e{1}%
\e{0}%
\e{1}%
\e{2}%
\e{1}%
\e{1}%
\e{0}%
\e{0}%
\e{2}%
\e{4}%
\e{2}%
\eol}\vss}\rg%
%
%
\rx{\vss\hfull{%
\rlx{\hss{$2240_{x}$}}\cg%
\e{0}%
\e{0}%
\e{3}%
\e{5}%
\e{3}%
\e{1}%
\e{0}%
\e{5}%
\e{2}%
\e{4}%
\e{2}%
\e{0}%
\e{0}%
\e{0}%
\e{1}%
\e{1}%
\e{0}%
\e{0}%
\e{0}%
\e{0}%
\e{1}%
\e{2}%
\e{1}%
\eol}\vss}\rg%
%
%
\rx{\vss\hfull{%
\rlx{\hss{$2835_{x}$}}\cg%
\e{0}%
\e{0}%
\e{1}%
\e{5}%
\e{3}%
\e{2}%
\e{0}%
\e{4}%
\e{2}%
\e{5}%
\e{3}%
\e{1}%
\e{0}%
\e{0}%
\e{0}%
\e{1}%
\e{1}%
\e{1}%
\e{0}%
\e{0}%
\e{1}%
\e{2}%
\e{2}%
\eol}\vss}\rg%
%
%
\rx{\vss\hfull{%
\rlx{\hss{$6075_{x}$}}\cg%
\e{1}%
\e{0}%
\e{3}%
\e{8}%
\e{4}%
\e{5}%
\e{1}%
\e{6}%
\e{7}%
\e{10}%
\e{9}%
\e{4}%
\e{2}%
\e{0}%
\e{1}%
\e{1}%
\e{2}%
\e{1}%
\e{1}%
\e{0}%
\e{2}%
\e{6}%
\e{2}%
\eol}\vss}\rg%
%
%
\rx{\vss\hfull{%
\rlx{\hss{$3200_{x}$}}\cg%
\e{0}%
\e{0}%
\e{1}%
\e{3}%
\e{3}%
\e{2}%
\e{0}%
\e{3}%
\e{2}%
\e{6}%
\e{5}%
\e{1}%
\e{2}%
\e{1}%
\e{3}%
\e{2}%
\e{2}%
\e{1}%
\e{0}%
\e{0}%
\e{1}%
\e{3}%
\e{2}%
\eol}\vss}\rg%
%
%
\rx{\vss\hfull{%
\rlx{\hss{$70_{y}$}}\cg%
\e{0}%
\e{0}%
\e{0}%
\e{0}%
\e{0}%
\e{0}%
\e{0}%
\e{0}%
\e{0}%
\e{0}%
\e{0}%
\e{1}%
\e{0}%
\e{0}%
\e{0}%
\e{0}%
\e{0}%
\e{0}%
\e{0}%
\e{0}%
\e{0}%
\e{0}%
\e{0}%
\eol}\vss}\rg%
%
%
\rx{\vss\hfull{%
\rlx{\hss{$1134_{y}$}}\cg%
\e{1}%
\e{0}%
\e{1}%
\e{1}%
\e{1}%
\e{1}%
\e{0}%
\e{0}%
\e{1}%
\e{1}%
\e{2}%
\e{1}%
\e{1}%
\e{0}%
\e{0}%
\e{0}%
\e{0}%
\e{0}%
\e{0}%
\e{0}%
\e{0}%
\e{1}%
\e{0}%
\eol}\vss}\rg%
%
%
\rx{\vss\hfull{%
\rlx{\hss{$1680_{y}$}}\cg%
\e{1}%
\e{0}%
\e{0}%
\e{1}%
\e{0}%
\e{2}%
\e{1}%
\e{0}%
\e{2}%
\e{2}%
\e{2}%
\e{4}%
\e{2}%
\e{0}%
\e{0}%
\e{0}%
\e{0}%
\e{0}%
\e{0}%
\e{0}%
\e{1}%
\e{1}%
\e{0}%
\eol}\vss}\rg%
%
%
\rx{\vss\hfull{%
\rlx{\hss{$168_{y}$}}\cg%
\e{0}%
\e{0}%
\e{0}%
\e{0}%
\e{0}%
\e{0}%
\e{0}%
\e{0}%
\e{0}%
\e{1}%
\e{0}%
\e{0}%
\e{0}%
\e{0}%
\e{0}%
\e{1}%
\e{0}%
\e{1}%
\e{0}%
\e{0}%
\e{0}%
\e{0}%
\e{1}%
\eol}\vss}\rg%
%
%
\rx{\vss\hfull{%
\rlx{\hss{$420_{y}$}}\cg%
\e{0}%
\e{0}%
\e{0}%
\e{1}%
\e{1}%
\e{0}%
\e{0}%
\e{0}%
\e{0}%
\e{1}%
\e{0}%
\e{0}%
\e{0}%
\e{0}%
\e{0}%
\e{0}%
\e{0}%
\e{0}%
\e{0}%
\e{0}%
\e{0}%
\e{0}%
\e{1}%
\eol}\vss}\rg%
%
%
\rx{\vss\hfull{%
\rlx{\hss{$3150_{y}$}}\cg%
\e{1}%
\e{0}%
\e{1}%
\e{4}%
\e{4}%
\e{3}%
\e{0}%
\e{2}%
\e{1}%
\e{4}%
\e{4}%
\e{1}%
\e{1}%
\e{0}%
\e{0}%
\e{0}%
\e{0}%
\e{0}%
\e{0}%
\e{0}%
\e{0}%
\e{2}%
\e{2}%
\eol}\vss}\rg%
%
%
\rx{\vss\hfull{%
\rlx{\hss{$4200_{y}$}}\cg%
\e{0}%
\e{0}%
\e{1}%
\e{4}%
\e{4}%
\e{3}%
\e{0}%
\e{2}%
\e{2}%
\e{7}%
\e{6}%
\e{2}%
\e{2}%
\e{0}%
\e{1}%
\e{3}%
\e{2}%
\e{3}%
\e{1}%
\e{0}%
\e{0}%
\e{3}%
\e{4}%
\eol}\vss}\rg%
%
%
\rx{\vss\hfull{%
\rlx{\hss{$2688_{y}$}}\cg%
\e{0}%
\e{0}%
\e{0}%
\e{2}%
\e{2}%
\e{2}%
\e{0}%
\e{2}%
\e{2}%
\e{4}%
\e{4}%
\e{2}%
\e{2}%
\e{0}%
\e{1}%
\e{2}%
\e{2}%
\e{2}%
\e{1}%
\e{0}%
\e{1}%
\e{2}%
\e{1}%
\eol}\vss}\rg%
%
%
\rx{\vss\hfull{%
\rlx{\hss{$2100_{y}$}}\cg%
\e{1}%
\e{0}%
\e{0}%
\e{1}%
\e{0}%
\e{2}%
\e{1}%
\e{1}%
\e{2}%
\e{3}%
\e{4}%
\e{2}%
\e{2}%
\e{0}%
\e{1}%
\e{0}%
\e{2}%
\e{0}%
\e{1}%
\e{0}%
\e{0}%
\e{2}%
\e{1}%
\eol}\vss}\rg%
%
%
\rx{\vss\hfull{%
\rlx{\hss{$1400_{y}$}}\cg%
\e{0}%
\e{0}%
\e{0}%
\e{1}%
\e{0}%
\e{2}%
\e{1}%
\e{0}%
\e{1}%
\e{2}%
\e{2}%
\e{2}%
\e{1}%
\e{0}%
\e{0}%
\e{0}%
\e{0}%
\e{0}%
\e{0}%
\e{0}%
\e{0}%
\e{1}%
\e{1}%
\eol}\vss}\rg%
%
%
\rx{\vss\hfull{%
\rlx{\hss{$4536_{y}$}}\cg%
\e{1}%
\e{0}%
\e{0}%
\e{4}%
\e{2}%
\e{6}%
\e{2}%
\e{2}%
\e{3}%
\e{6}%
\e{6}%
\e{4}%
\e{3}%
\e{0}%
\e{0}%
\e{0}%
\e{0}%
\e{0}%
\e{0}%
\e{0}%
\e{1}%
\e{3}%
\e{2}%
\eol}\vss}\rg%
%
%
\rx{\vss\hfull{%
\rlx{\hss{$5670_{y}$}}\cg%
\e{2}%
\e{0}%
\e{1}%
\e{5}%
\e{3}%
\e{7}%
\e{2}%
\e{2}%
\e{4}%
\e{7}%
\e{8}%
\e{5}%
\e{4}%
\e{0}%
\e{0}%
\e{0}%
\e{0}%
\e{0}%
\e{0}%
\e{0}%
\e{1}%
\e{4}%
\e{2}%
\eol}\vss}\rg%
%
%
\rx{\vss\hfull{%
\rlx{\hss{$4480_{y}$}}\cg%
\e{1}%
\e{0}%
\e{1}%
\e{5}%
\e{4}%
\e{5}%
\e{1}%
\e{2}%
\e{2}%
\e{6}%
\e{6}%
\e{2}%
\e{2}%
\e{0}%
\e{0}%
\e{0}%
\e{0}%
\e{0}%
\e{0}%
\e{0}%
\e{0}%
\e{3}%
\e{3}%
\eol}\vss}\rg%
%
%
\rx{\vss\hfull{%
\rlx{\hss{$8_{z}$}}\cg%
\e{0}%
\e{0}%
\e{0}%
\e{0}%
\e{0}%
\e{0}%
\e{0}%
\e{0}%
\e{0}%
\e{0}%
\e{0}%
\e{0}%
\e{0}%
\e{0}%
\e{0}%
\e{0}%
\e{0}%
\e{0}%
\e{0}%
\e{0}%
\e{0}%
\e{0}%
\e{0}%
\eol}\vss}\rg%
%
%
\rx{\vss\hfull{%
\rlx{\hss{$56_{z}$}}\cg%
\e{0}%
\e{0}%
\e{0}%
\e{0}%
\e{0}%
\e{0}%
\e{0}%
\e{0}%
\e{1}%
\e{0}%
\e{0}%
\e{0}%
\e{0}%
\e{0}%
\e{0}%
\e{0}%
\e{0}%
\e{0}%
\e{0}%
\e{0}%
\e{0}%
\e{0}%
\e{0}%
\eol}\vss}\rg%
%
%
\rx{\vss\hfull{%
\rlx{\hss{$160_{z}$}}\cg%
\e{0}%
\e{0}%
\e{1}%
\e{0}%
\e{0}%
\e{0}%
\e{0}%
\e{1}%
\e{1}%
\e{0}%
\e{0}%
\e{0}%
\e{0}%
\e{1}%
\e{0}%
\e{0}%
\e{0}%
\e{0}%
\e{0}%
\e{0}%
\e{1}%
\e{0}%
\e{0}%
\eol}\vss}\rg%
%
%
\rx{\vss\hfull{%
\rlx{\hss{$112_{z}$}}\cg%
\e{0}%
\e{0}%
\e{1}%
\e{0}%
\e{0}%
\e{0}%
\e{0}%
\e{1}%
\e{0}%
\e{0}%
\e{0}%
\e{0}%
\e{0}%
\e{0}%
\e{0}%
\e{0}%
\e{0}%
\e{0}%
\e{0}%
\e{0}%
\e{1}%
\e{0}%
\e{0}%
\eol}\vss}\rg%
%
%
\rx{\vss\hfull{%
\rlx{\hss{$840_{z}$}}\cg%
\e{0}%
\e{0}%
\e{1}%
\e{1}%
\e{1}%
\e{1}%
\e{0}%
\e{1}%
\e{1}%
\e{2}%
\e{1}%
\e{0}%
\e{0}%
\e{0}%
\e{1}%
\e{0}%
\e{0}%
\e{0}%
\e{0}%
\e{0}%
\e{1}%
\e{1}%
\e{0}%
\eol}\vss}\rg%
%
%
\rx{\vss\hfull{%
\rlx{\hss{$1296_{z}$}}\cg%
\e{0}%
\e{0}%
\e{2}%
\e{2}%
\e{0}%
\e{1}%
\e{0}%
\e{2}%
\e{4}%
\e{2}%
\e{2}%
\e{1}%
\e{0}%
\e{1}%
\e{1}%
\e{0}%
\e{0}%
\e{0}%
\e{0}%
\e{0}%
\e{1}%
\e{1}%
\e{0}%
\eol}\vss}\rg%
%
%
\rx{\vss\hfull{%
\rlx{\hss{$1400_{z}$}}\cg%
\e{0}%
\e{0}%
\e{4}%
\e{3}%
\e{1}%
\e{0}%
\e{0}%
\e{5}%
\e{2}%
\e{3}%
\e{1}%
\e{0}%
\e{0}%
\e{1}%
\e{1}%
\e{0}%
\e{0}%
\e{0}%
\e{0}%
\e{0}%
\e{3}%
\e{2}%
\e{0}%
\eol}\vss}\rg%
%
%
\rx{\vss\hfull{%
\rlx{\hss{$1008_{z}$}}\cg%
\e{0}%
\e{0}%
\e{3}%
\e{2}%
\e{0}%
\e{0}%
\e{0}%
\e{3}%
\e{3}%
\e{2}%
\e{1}%
\e{0}%
\e{0}%
\e{1}%
\e{1}%
\e{0}%
\e{0}%
\e{0}%
\e{0}%
\e{0}%
\e{2}%
\e{1}%
\e{0}%
\eol}\vss}\rg%
%
%
\rx{\vss\hfull{%
\rlx{\hss{$560_{z}$}}\cg%
\e{0}%
\e{0}%
\e{3}%
\e{1}%
\e{0}%
\e{0}%
\e{0}%
\e{3}%
\e{1}%
\e{1}%
\e{0}%
\e{0}%
\e{0}%
\e{1}%
\e{0}%
\e{0}%
\e{0}%
\e{0}%
\e{0}%
\e{0}%
\e{2}%
\e{1}%
\e{0}%
\eol}\vss}\rg%
%
%
\rx{\vss\hfull{%
\rlx{\hss{$1400_{zz}$}}\cg%
\e{0}%
\e{0}%
\e{2}%
\e{3}%
\e{2}%
\e{0}%
\e{0}%
\e{4}%
\e{1}%
\e{3}%
\e{1}%
\e{0}%
\e{0}%
\e{0}%
\e{0}%
\e{1}%
\e{0}%
\e{0}%
\e{0}%
\e{0}%
\e{0}%
\e{2}%
\e{2}%
\eol}\vss}\rg%
\eop
\eject
\tablecont%
%
%
%
%
%
%
\rowpts=18 true pt%
\colpts=18 true pt%
\rowlabpts=40 true pt%
\collabpts=85 true pt%
\clx{\vss\hfull{%
\rlx{\hss{$ $}}\cg%
\cx{\hskip 16 true pt\flip{$[{2}{1^{3}}:-]{\times}[{3}{1}]$}\hss}\cg%
\cx{\hskip 16 true pt\flip{$[{1^{5}}:-]{\times}[{3}{1}]$}\hss}\cg%
\cx{\hskip 16 true pt\flip{$[{4}:{1}]{\times}[{3}{1}]$}\hss}\cg%
\cx{\hskip 16 true pt\flip{$[{3}{1}:{1}]{\times}[{3}{1}]$}\hss}\cg%
\cx{\hskip 16 true pt\flip{$[{2^{2}}:{1}]{\times}[{3}{1}]$}\hss}\cg%
\cx{\hskip 16 true pt\flip{$[{2}{1^{2}}:{1}]{\times}[{3}{1}]$}\hss}\cg%
\cx{\hskip 16 true pt\flip{$[{1^{4}}:{1}]{\times}[{3}{1}]$}\hss}\cg%
\cx{\hskip 16 true pt\flip{$[{3}:{2}]{\times}[{3}{1}]$}\hss}\cg%
\cx{\hskip 16 true pt\flip{$[{3}:{1^{2}}]{\times}[{3}{1}]$}\hss}\cg%
\cx{\hskip 16 true pt\flip{$[{2}{1}:{2}]{\times}[{3}{1}]$}\hss}\cg%
\cx{\hskip 16 true pt\flip{$[{2}{1}:{1^{2}}]{\times}[{3}{1}]$}\hss}\cg%
\cx{\hskip 16 true pt\flip{$[{1^{3}}:{2}]{\times}[{3}{1}]$}\hss}\cg%
\cx{\hskip 16 true pt\flip{$[{1^{3}}:{1^{2}}]{\times}[{3}{1}]$}\hss}\cg%
\cx{\hskip 16 true pt\flip{$[{5}:-]{\times}[{2^{2}}]$}\hss}\cg%
\cx{\hskip 16 true pt\flip{$[{4}{1}:-]{\times}[{2^{2}}]$}\hss}\cg%
\cx{\hskip 16 true pt\flip{$[{3}{2}:-]{\times}[{2^{2}}]$}\hss}\cg%
\cx{\hskip 16 true pt\flip{$[{3}{1^{2}}:-]{\times}[{2^{2}}]$}\hss}\cg%
\cx{\hskip 16 true pt\flip{$[{2^{2}}{1}:-]{\times}[{2^{2}}]$}\hss}\cg%
\cx{\hskip 16 true pt\flip{$[{2}{1^{3}}:-]{\times}[{2^{2}}]$}\hss}\cg%
\cx{\hskip 16 true pt\flip{$[{1^{5}}:-]{\times}[{2^{2}}]$}\hss}\cg%
\cx{\hskip 16 true pt\flip{$[{4}:{1}]{\times}[{2^{2}}]$}\hss}\cg%
\cx{\hskip 16 true pt\flip{$[{3}{1}:{1}]{\times}[{2^{2}}]$}\hss}\cg%
\cx{\hskip 16 true pt\flip{$[{2^{2}}:{1}]{\times}[{2^{2}}]$}\hss}\cg%
\eol}}\rg%
%
%
\rx{\vss\hfull{%
\rlx{\hss{$4200_{z}$}}\cg%
\e{1}%
\e{0}%
\e{2}%
\e{6}%
\e{3}%
\e{3}%
\e{0}%
\e{4}%
\e{4}%
\e{7}%
\e{6}%
\e{2}%
\e{1}%
\e{0}%
\e{1}%
\e{2}%
\e{1}%
\e{1}%
\e{0}%
\e{0}%
\e{0}%
\e{3}%
\e{3}%
\eol}\vss}\rg%
%
%
\rx{\vss\hfull{%
\rlx{\hss{$400_{z}$}}\cg%
\e{0}%
\e{0}%
\e{2}%
\e{1}%
\e{0}%
\e{0}%
\e{0}%
\e{2}%
\e{0}%
\e{1}%
\e{0}%
\e{0}%
\e{0}%
\e{0}%
\e{0}%
\e{0}%
\e{0}%
\e{0}%
\e{0}%
\e{0}%
\e{0}%
\e{1}%
\e{1}%
\eol}\vss}\rg%
%
%
\rx{\vss\hfull{%
\rlx{\hss{$3240_{z}$}}\cg%
\e{0}%
\e{0}%
\e{6}%
\e{7}%
\e{2}%
\e{1}%
\e{0}%
\e{8}%
\e{4}%
\e{7}%
\e{3}%
\e{1}%
\e{0}%
\e{1}%
\e{2}%
\e{1}%
\e{0}%
\e{0}%
\e{0}%
\e{0}%
\e{3}%
\e{4}%
\e{2}%
\eol}\vss}\rg%
%
%
\rx{\vss\hfull{%
\rlx{\hss{$4536_{z}$}}\cg%
\e{0}%
\e{0}%
\e{2}%
\e{7}%
\e{6}%
\e{3}%
\e{0}%
\e{6}%
\e{3}%
\e{9}%
\e{5}%
\e{1}%
\e{1}%
\e{0}%
\e{1}%
\e{2}%
\e{1}%
\e{1}%
\e{0}%
\e{0}%
\e{2}%
\e{5}%
\e{3}%
\eol}\vss}\rg%
%
%
\rx{\vss\hfull{%
\rlx{\hss{$2400_{z}$}}\cg%
\e{1}%
\e{0}%
\e{1}%
\e{3}%
\e{0}%
\e{3}%
\e{1}%
\e{1}%
\e{4}%
\e{3}%
\e{4}%
\e{3}%
\e{1}%
\e{0}%
\e{1}%
\e{0}%
\e{1}%
\e{0}%
\e{0}%
\e{0}%
\e{0}%
\e{1}%
\e{0}%
\eol}\vss}\rg%
%
%
\rx{\vss\hfull{%
\rlx{\hss{$3360_{z}$}}\cg%
\e{0}%
\e{0}%
\e{3}%
\e{6}%
\e{2}%
\e{2}%
\e{0}%
\e{5}%
\e{4}%
\e{6}%
\e{4}%
\e{2}%
\e{1}%
\e{1}%
\e{1}%
\e{1}%
\e{1}%
\e{0}%
\e{0}%
\e{0}%
\e{1}%
\e{3}%
\e{2}%
\eol}\vss}\rg%
%
%
\rx{\vss\hfull{%
\rlx{\hss{$2800_{z}$}}\cg%
\e{1}%
\e{0}%
\e{2}%
\e{5}%
\e{1}%
\e{2}%
\e{0}%
\e{3}%
\e{4}%
\e{5}%
\e{4}%
\e{2}%
\e{1}%
\e{0}%
\e{1}%
\e{1}%
\e{1}%
\e{0}%
\e{0}%
\e{0}%
\e{1}%
\e{2}%
\e{0}%
\eol}\vss}\rg%
%
%
\rx{\vss\hfull{%
\rlx{\hss{$4096_{z}$}}\cg%
\e{0}%
\e{0}%
\e{4}%
\e{7}%
\e{3}%
\e{3}%
\e{0}%
\e{6}%
\e{5}%
\e{8}%
\e{5}%
\e{2}%
\e{1}%
\e{1}%
\e{2}%
\e{1}%
\e{1}%
\e{0}%
\e{0}%
\e{0}%
\e{3}%
\e{4}%
\e{1}%
\eol}\vss}\rg%
%
%
\rx{\vss\hfull{%
\rlx{\hss{$5600_{z}$}}\cg%
\e{1}%
\e{0}%
\e{2}%
\e{7}%
\e{4}%
\e{6}%
\e{1}%
\e{4}%
\e{4}%
\e{9}%
\e{8}%
\e{4}%
\e{3}%
\e{0}%
\e{2}%
\e{1}%
\e{2}%
\e{1}%
\e{0}%
\e{0}%
\e{2}%
\e{4}%
\e{1}%
\eol}\vss}\rg%
%
%
\rx{\vss\hfull{%
\rlx{\hss{$448_{z}$}}\cg%
\e{0}%
\e{0}%
\e{1}%
\e{1}%
\e{1}%
\e{0}%
\e{0}%
\e{2}%
\e{0}%
\e{1}%
\e{0}%
\e{0}%
\e{0}%
\e{0}%
\e{0}%
\e{0}%
\e{0}%
\e{0}%
\e{0}%
\e{0}%
\e{1}%
\e{1}%
\e{0}%
\eol}\vss}\rg%
%
%
\rx{\vss\hfull{%
\rlx{\hss{$448_{w}$}}\cg%
\e{0}%
\e{0}%
\e{0}%
\e{0}%
\e{0}%
\e{1}%
\e{1}%
\e{0}%
\e{1}%
\e{0}%
\e{1}%
\e{1}%
\e{0}%
\e{0}%
\e{0}%
\e{0}%
\e{0}%
\e{0}%
\e{0}%
\e{0}%
\e{0}%
\e{0}%
\e{0}%
\eol}\vss}\rg%
%
%
\rx{\vss\hfull{%
\rlx{\hss{$1344_{w}$}}\cg%
\e{0}%
\e{0}%
\e{0}%
\e{1}%
\e{2}%
\e{1}%
\e{0}%
\e{1}%
\e{1}%
\e{2}%
\e{2}%
\e{0}%
\e{0}%
\e{0}%
\e{0}%
\e{1}%
\e{0}%
\e{1}%
\e{0}%
\e{0}%
\e{0}%
\e{1}%
\e{2}%
\eol}\vss}\rg%
%
%
\rx{\vss\hfull{%
\rlx{\hss{$5600_{w}$}}\cg%
\e{1}%
\e{0}%
\e{1}%
\e{5}%
\e{3}%
\e{7}%
\e{2}%
\e{2}%
\e{4}%
\e{7}%
\e{8}%
\e{5}%
\e{4}%
\e{0}%
\e{1}%
\e{1}%
\e{2}%
\e{1}%
\e{1}%
\e{0}%
\e{1}%
\e{3}%
\e{2}%
\eol}\vss}\rg%
%
%
\rx{\vss\hfull{%
\rlx{\hss{$2016_{w}$}}\cg%
\e{0}%
\e{0}%
\e{0}%
\e{2}%
\e{3}%
\e{1}%
\e{0}%
\e{2}%
\e{1}%
\e{3}%
\e{3}%
\e{0}%
\e{0}%
\e{0}%
\e{0}%
\e{1}%
\e{0}%
\e{1}%
\e{0}%
\e{0}%
\e{0}%
\e{2}%
\e{2}%
\eol}\vss}\rg%
%
%
\rx{\vss\hfull{%
\rlx{\hss{$7168_{w}$}}\cg%
\e{1}%
\e{0}%
\e{1}%
\e{7}%
\e{6}%
\e{7}%
\e{1}%
\e{4}%
\e{4}%
\e{10}%
\e{10}%
\e{4}%
\e{4}%
\e{0}%
\e{1}%
\e{2}%
\e{2}%
\e{2}%
\e{1}%
\e{0}%
\e{1}%
\e{5}%
\e{4}%
\eol}\vss}\rg%
%
%
%
%
%
%
\rowpts=18 true pt%
\colpts=18 true pt%
\rowlabpts=40 true pt%
\collabpts=85 true pt%
\clx{\vss\hfull{%
\rlx{\hss{$ $}}\cg%
\cx{\hskip 16 true pt\flip{$[{2}{1^{2}}:{1}]{\times}[{2^{2}}]$}\hss}\cg%
\cx{\hskip 16 true pt\flip{$[{1^{4}}:{1}]{\times}[{2^{2}}]$}\hss}\cg%
\cx{\hskip 16 true pt\flip{$[{3}:{2}]{\times}[{2^{2}}]$}\hss}\cg%
\cx{\hskip 16 true pt\flip{$[{3}:{1^{2}}]{\times}[{2^{2}}]$}\hss}\cg%
\cx{\hskip 16 true pt\flip{$[{2}{1}:{2}]{\times}[{2^{2}}]$}\hss}\cg%
\cx{\hskip 16 true pt\flip{$[{2}{1}:{1^{2}}]{\times}[{2^{2}}]$}\hss}\cg%
\cx{\hskip 16 true pt\flip{$[{1^{3}}:{2}]{\times}[{2^{2}}]$}\hss}\cg%
\cx{\hskip 16 true pt\flip{$[{1^{3}}:{1^{2}}]{\times}[{2^{2}}]$}\hss}\cg%
\cx{\hskip 16 true pt\flip{$[{5}:-]{\times}[{2}{1^{2}}]$}\hss}\cg%
\cx{\hskip 16 true pt\flip{$[{4}{1}:-]{\times}[{2}{1^{2}}]$}\hss}\cg%
\cx{\hskip 16 true pt\flip{$[{3}{2}:-]{\times}[{2}{1^{2}}]$}\hss}\cg%
\cx{\hskip 16 true pt\flip{$[{3}{1^{2}}:-]{\times}[{2}{1^{2}}]$}\hss}\cg%
\cx{\hskip 16 true pt\flip{$[{2^{2}}{1}:-]{\times}[{2}{1^{2}}]$}\hss}\cg%
\cx{\hskip 16 true pt\flip{$[{2}{1^{3}}:-]{\times}[{2}{1^{2}}]$}\hss}\cg%
\cx{\hskip 16 true pt\flip{$[{1^{5}}:-]{\times}[{2}{1^{2}}]$}\hss}\cg%
\cx{\hskip 16 true pt\flip{$[{4}:{1}]{\times}[{2}{1^{2}}]$}\hss}\cg%
\cx{\hskip 16 true pt\flip{$[{3}{1}:{1}]{\times}[{2}{1^{2}}]$}\hss}\cg%
\cx{\hskip 16 true pt\flip{$[{2^{2}}:{1}]{\times}[{2}{1^{2}}]$}\hss}\cg%
\cx{\hskip 16 true pt\flip{$[{2}{1^{2}}:{1}]{\times}[{2}{1^{2}}]$}\hss}\cg%
\cx{\hskip 16 true pt\flip{$[{1^{4}}:{1}]{\times}[{2}{1^{2}}]$}\hss}\cg%
\cx{\hskip 16 true pt\flip{$[{3}:{2}]{\times}[{2}{1^{2}}]$}\hss}\cg%
\cx{\hskip 16 true pt\flip{$[{3}:{1^{2}}]{\times}[{2}{1^{2}}]$}\hss}\cg%
\eol}}\rg%
%
%
\rx{\vss\hfull{%
\rlx{\hss{$1_{x}$}}\cg%
\e{0}%
\e{0}%
\e{0}%
\e{0}%
\e{0}%
\e{0}%
\e{0}%
\e{0}%
\e{0}%
\e{0}%
\e{0}%
\e{0}%
\e{0}%
\e{0}%
\e{0}%
\e{0}%
\e{0}%
\e{0}%
\e{0}%
\e{0}%
\e{0}%
\e{0}%
\eol}\vss}\rg%
%
%
\rx{\vss\hfull{%
\rlx{\hss{$28_{x}$}}\cg%
\e{0}%
\e{0}%
\e{0}%
\e{0}%
\e{0}%
\e{0}%
\e{0}%
\e{0}%
\e{1}%
\e{0}%
\e{0}%
\e{0}%
\e{0}%
\e{0}%
\e{0}%
\e{0}%
\e{0}%
\e{0}%
\e{0}%
\e{0}%
\e{0}%
\e{0}%
\eol}\vss}\rg%
%
%
\rx{\vss\hfull{%
\rlx{\hss{$35_{x}$}}\cg%
\e{0}%
\e{0}%
\e{0}%
\e{0}%
\e{0}%
\e{0}%
\e{0}%
\e{0}%
\e{0}%
\e{0}%
\e{0}%
\e{0}%
\e{0}%
\e{0}%
\e{0}%
\e{0}%
\e{0}%
\e{0}%
\e{0}%
\e{0}%
\e{0}%
\e{0}%
\eol}\vss}\rg%
%
%
\rx{\vss\hfull{%
\rlx{\hss{$84_{x}$}}\cg%
\e{0}%
\e{0}%
\e{0}%
\e{0}%
\e{0}%
\e{0}%
\e{0}%
\e{0}%
\e{0}%
\e{0}%
\e{0}%
\e{0}%
\e{0}%
\e{0}%
\e{0}%
\e{0}%
\e{0}%
\e{0}%
\e{0}%
\e{0}%
\e{0}%
\e{0}%
\eol}\vss}\rg%
%
%
\rx{\vss\hfull{%
\rlx{\hss{$50_{x}$}}\cg%
\e{0}%
\e{0}%
\e{0}%
\e{0}%
\e{0}%
\e{0}%
\e{0}%
\e{0}%
\e{0}%
\e{0}%
\e{0}%
\e{0}%
\e{0}%
\e{0}%
\e{0}%
\e{0}%
\e{0}%
\e{0}%
\e{0}%
\e{0}%
\e{0}%
\e{0}%
\eol}\vss}\rg%
%
%
\rx{\vss\hfull{%
\rlx{\hss{$350_{x}$}}\cg%
\e{0}%
\e{0}%
\e{0}%
\e{1}%
\e{0}%
\e{0}%
\e{0}%
\e{0}%
\e{1}%
\e{1}%
\e{0}%
\e{0}%
\e{0}%
\e{0}%
\e{0}%
\e{1}%
\e{0}%
\e{0}%
\e{0}%
\e{0}%
\e{1}%
\e{1}%
\eol}\vss}\rg%
%
%
\rx{\vss\hfull{%
\rlx{\hss{$300_{x}$}}\cg%
\e{0}%
\e{0}%
\e{1}%
\e{0}%
\e{0}%
\e{0}%
\e{0}%
\e{0}%
\e{0}%
\e{0}%
\e{0}%
\e{0}%
\e{0}%
\e{0}%
\e{0}%
\e{1}%
\e{0}%
\e{0}%
\e{0}%
\e{0}%
\e{0}%
\e{1}%
\eol}\vss}\rg%
\eop
\eject
\tablecont%
%
%
%
%
%
%
\rowpts=18 true pt%
\colpts=18 true pt%
\rowlabpts=40 true pt%
\collabpts=85 true pt%
\clx{\vss\hfull{%
\rlx{\hss{$ $}}\cg%
\cx{\hskip 16 true pt\flip{$[{2}{1^{2}}:{1}]{\times}[{2^{2}}]$}\hss}\cg%
\cx{\hskip 16 true pt\flip{$[{1^{4}}:{1}]{\times}[{2^{2}}]$}\hss}\cg%
\cx{\hskip 16 true pt\flip{$[{3}:{2}]{\times}[{2^{2}}]$}\hss}\cg%
\cx{\hskip 16 true pt\flip{$[{3}:{1^{2}}]{\times}[{2^{2}}]$}\hss}\cg%
\cx{\hskip 16 true pt\flip{$[{2}{1}:{2}]{\times}[{2^{2}}]$}\hss}\cg%
\cx{\hskip 16 true pt\flip{$[{2}{1}:{1^{2}}]{\times}[{2^{2}}]$}\hss}\cg%
\cx{\hskip 16 true pt\flip{$[{1^{3}}:{2}]{\times}[{2^{2}}]$}\hss}\cg%
\cx{\hskip 16 true pt\flip{$[{1^{3}}:{1^{2}}]{\times}[{2^{2}}]$}\hss}\cg%
\cx{\hskip 16 true pt\flip{$[{5}:-]{\times}[{2}{1^{2}}]$}\hss}\cg%
\cx{\hskip 16 true pt\flip{$[{4}{1}:-]{\times}[{2}{1^{2}}]$}\hss}\cg%
\cx{\hskip 16 true pt\flip{$[{3}{2}:-]{\times}[{2}{1^{2}}]$}\hss}\cg%
\cx{\hskip 16 true pt\flip{$[{3}{1^{2}}:-]{\times}[{2}{1^{2}}]$}\hss}\cg%
\cx{\hskip 16 true pt\flip{$[{2^{2}}{1}:-]{\times}[{2}{1^{2}}]$}\hss}\cg%
\cx{\hskip 16 true pt\flip{$[{2}{1^{3}}:-]{\times}[{2}{1^{2}}]$}\hss}\cg%
\cx{\hskip 16 true pt\flip{$[{1^{5}}:-]{\times}[{2}{1^{2}}]$}\hss}\cg%
\cx{\hskip 16 true pt\flip{$[{4}:{1}]{\times}[{2}{1^{2}}]$}\hss}\cg%
\cx{\hskip 16 true pt\flip{$[{3}{1}:{1}]{\times}[{2}{1^{2}}]$}\hss}\cg%
\cx{\hskip 16 true pt\flip{$[{2^{2}}:{1}]{\times}[{2}{1^{2}}]$}\hss}\cg%
\cx{\hskip 16 true pt\flip{$[{2}{1^{2}}:{1}]{\times}[{2}{1^{2}}]$}\hss}\cg%
\cx{\hskip 16 true pt\flip{$[{1^{4}}:{1}]{\times}[{2}{1^{2}}]$}\hss}\cg%
\cx{\hskip 16 true pt\flip{$[{3}:{2}]{\times}[{2}{1^{2}}]$}\hss}\cg%
\cx{\hskip 16 true pt\flip{$[{3}:{1^{2}}]{\times}[{2}{1^{2}}]$}\hss}\cg%
\eol}}\rg%
%
%
\rx{\vss\hfull{%
\rlx{\hss{$567_{x}$}}\cg%
\e{0}%
\e{0}%
\e{1}%
\e{1}%
\e{0}%
\e{0}%
\e{0}%
\e{0}%
\e{2}%
\e{1}%
\e{0}%
\e{0}%
\e{0}%
\e{0}%
\e{0}%
\e{2}%
\e{0}%
\e{0}%
\e{0}%
\e{0}%
\e{1}%
\e{0}%
\eol}\vss}\rg%
%
%
\rx{\vss\hfull{%
\rlx{\hss{$210_{x}$}}\cg%
\e{0}%
\e{0}%
\e{1}%
\e{0}%
\e{0}%
\e{0}%
\e{0}%
\e{0}%
\e{1}%
\e{0}%
\e{0}%
\e{0}%
\e{0}%
\e{0}%
\e{0}%
\e{1}%
\e{0}%
\e{0}%
\e{0}%
\e{0}%
\e{0}%
\e{0}%
\eol}\vss}\rg%
%
%
\rx{\vss\hfull{%
\rlx{\hss{$840_{x}$}}\cg%
\e{0}%
\e{0}%
\e{1}%
\e{0}%
\e{2}%
\e{0}%
\e{0}%
\e{0}%
\e{0}%
\e{0}%
\e{0}%
\e{0}%
\e{0}%
\e{0}%
\e{0}%
\e{0}%
\e{1}%
\e{0}%
\e{1}%
\e{0}%
\e{0}%
\e{1}%
\eol}\vss}\rg%
%
%
\rx{\vss\hfull{%
\rlx{\hss{$700_{x}$}}\cg%
\e{0}%
\e{0}%
\e{2}%
\e{0}%
\e{1}%
\e{0}%
\e{0}%
\e{0}%
\e{0}%
\e{0}%
\e{0}%
\e{0}%
\e{0}%
\e{0}%
\e{0}%
\e{1}%
\e{1}%
\e{0}%
\e{0}%
\e{0}%
\e{1}%
\e{0}%
\eol}\vss}\rg%
%
%
\rx{\vss\hfull{%
\rlx{\hss{$175_{x}$}}\cg%
\e{0}%
\e{0}%
\e{0}%
\e{0}%
\e{1}%
\e{0}%
\e{0}%
\e{0}%
\e{0}%
\e{0}%
\e{0}%
\e{0}%
\e{0}%
\e{0}%
\e{0}%
\e{0}%
\e{0}%
\e{1}%
\e{0}%
\e{0}%
\e{0}%
\e{0}%
\eol}\vss}\rg%
%
%
\rx{\vss\hfull{%
\rlx{\hss{$1400_{x}$}}\cg%
\e{0}%
\e{0}%
\e{2}%
\e{1}%
\e{2}%
\e{1}%
\e{0}%
\e{0}%
\e{1}%
\e{1}%
\e{1}%
\e{0}%
\e{0}%
\e{0}%
\e{0}%
\e{2}%
\e{2}%
\e{1}%
\e{0}%
\e{0}%
\e{2}%
\e{0}%
\eol}\vss}\rg%
%
%
\rx{\vss\hfull{%
\rlx{\hss{$1050_{x}$}}\cg%
\e{0}%
\e{0}%
\e{1}%
\e{0}%
\e{1}%
\e{1}%
\e{0}%
\e{0}%
\e{0}%
\e{0}%
\e{1}%
\e{0}%
\e{0}%
\e{0}%
\e{0}%
\e{0}%
\e{1}%
\e{1}%
\e{0}%
\e{0}%
\e{2}%
\e{0}%
\eol}\vss}\rg%
%
%
\rx{\vss\hfull{%
\rlx{\hss{$1575_{x}$}}\cg%
\e{0}%
\e{0}%
\e{2}%
\e{2}%
\e{1}%
\e{1}%
\e{0}%
\e{0}%
\e{2}%
\e{2}%
\e{1}%
\e{0}%
\e{0}%
\e{0}%
\e{0}%
\e{3}%
\e{2}%
\e{0}%
\e{0}%
\e{0}%
\e{3}%
\e{1}%
\eol}\vss}\rg%
%
%
\rx{\vss\hfull{%
\rlx{\hss{$1344_{x}$}}\cg%
\e{0}%
\e{0}%
\e{2}%
\e{1}%
\e{1}%
\e{0}%
\e{0}%
\e{0}%
\e{1}%
\e{1}%
\e{0}%
\e{0}%
\e{0}%
\e{0}%
\e{0}%
\e{2}%
\e{1}%
\e{0}%
\e{0}%
\e{0}%
\e{2}%
\e{1}%
\eol}\vss}\rg%
%
%
\rx{\vss\hfull{%
\rlx{\hss{$2100_{x}$}}\cg%
\e{1}%
\e{0}%
\e{1}%
\e{2}%
\e{1}%
\e{1}%
\e{0}%
\e{0}%
\e{1}%
\e{2}%
\e{0}%
\e{1}%
\e{0}%
\e{0}%
\e{0}%
\e{3}%
\e{2}%
\e{0}%
\e{0}%
\e{0}%
\e{3}%
\e{4}%
\eol}\vss}\rg%
%
%
\rx{\vss\hfull{%
\rlx{\hss{$2268_{x}$}}\cg%
\e{0}%
\e{0}%
\e{3}%
\e{2}%
\e{2}%
\e{1}%
\e{0}%
\e{0}%
\e{1}%
\e{2}%
\e{0}%
\e{1}%
\e{0}%
\e{0}%
\e{0}%
\e{4}%
\e{3}%
\e{0}%
\e{0}%
\e{0}%
\e{3}%
\e{2}%
\eol}\vss}\rg%
%
%
\rx{\vss\hfull{%
\rlx{\hss{$525_{x}$}}\cg%
\e{0}%
\e{0}%
\e{0}%
\e{1}%
\e{0}%
\e{0}%
\e{0}%
\e{0}%
\e{1}%
\e{1}%
\e{0}%
\e{0}%
\e{0}%
\e{0}%
\e{0}%
\e{1}%
\e{0}%
\e{0}%
\e{0}%
\e{0}%
\e{1}%
\e{0}%
\eol}\vss}\rg%
%
%
\rx{\vss\hfull{%
\rlx{\hss{$700_{xx}$}}\cg%
\e{0}%
\e{0}%
\e{0}%
\e{0}%
\e{0}%
\e{1}%
\e{0}%
\e{0}%
\e{0}%
\e{0}%
\e{1}%
\e{0}%
\e{0}%
\e{0}%
\e{0}%
\e{0}%
\e{0}%
\e{1}%
\e{0}%
\e{0}%
\e{1}%
\e{0}%
\eol}\vss}\rg%
%
%
\rx{\vss\hfull{%
\rlx{\hss{$972_{x}$}}\cg%
\e{0}%
\e{0}%
\e{1}%
\e{0}%
\e{1}%
\e{0}%
\e{0}%
\e{0}%
\e{0}%
\e{0}%
\e{0}%
\e{0}%
\e{0}%
\e{0}%
\e{0}%
\e{0}%
\e{1}%
\e{0}%
\e{0}%
\e{0}%
\e{1}%
\e{1}%
\eol}\vss}\rg%
%
%
\rx{\vss\hfull{%
\rlx{\hss{$4096_{x}$}}\cg%
\e{1}%
\e{0}%
\e{4}%
\e{3}%
\e{4}%
\e{2}%
\e{1}%
\e{0}%
\e{1}%
\e{3}%
\e{1}%
\e{1}%
\e{0}%
\e{0}%
\e{0}%
\e{4}%
\e{5}%
\e{1}%
\e{1}%
\e{0}%
\e{5}%
\e{4}%
\eol}\vss}\rg%
%
%
\rx{\vss\hfull{%
\rlx{\hss{$4200_{x}$}}\cg%
\e{1}%
\e{0}%
\e{4}%
\e{2}%
\e{5}%
\e{3}%
\e{1}%
\e{1}%
\e{0}%
\e{1}%
\e{2}%
\e{1}%
\e{1}%
\e{0}%
\e{0}%
\e{2}%
\e{5}%
\e{3}%
\e{2}%
\e{0}%
\e{4}%
\e{3}%
\eol}\vss}\rg%
%
%
\rx{\vss\hfull{%
\rlx{\hss{$2240_{x}$}}\cg%
\e{0}%
\e{0}%
\e{3}%
\e{1}%
\e{4}%
\e{1}%
\e{1}%
\e{0}%
\e{0}%
\e{1}%
\e{1}%
\e{0}%
\e{0}%
\e{0}%
\e{0}%
\e{1}%
\e{3}%
\e{2}%
\e{1}%
\e{0}%
\e{2}%
\e{1}%
\eol}\vss}\rg%
%
%
\rx{\vss\hfull{%
\rlx{\hss{$2835_{x}$}}\cg%
\e{1}%
\e{0}%
\e{2}%
\e{1}%
\e{4}%
\e{3}%
\e{1}%
\e{1}%
\e{0}%
\e{0}%
\e{2}%
\e{1}%
\e{1}%
\e{0}%
\e{0}%
\e{0}%
\e{3}%
\e{4}%
\e{2}%
\e{0}%
\e{2}%
\e{1}%
\eol}\vss}\rg%
%
%
\rx{\vss\hfull{%
\rlx{\hss{$6075_{x}$}}\cg%
\e{3}%
\e{0}%
\e{3}%
\e{4}%
\e{5}%
\e{5}%
\e{2}%
\e{1}%
\e{1}%
\e{3}%
\e{3}%
\e{2}%
\e{1}%
\e{0}%
\e{0}%
\e{3}%
\e{7}%
\e{3}%
\e{3}%
\e{0}%
\e{7}%
\e{4}%
\eol}\vss}\rg%
%
%
\rx{\vss\hfull{%
\rlx{\hss{$3200_{x}$}}\cg%
\e{2}%
\e{0}%
\e{2}%
\e{1}%
\e{4}%
\e{1}%
\e{1}%
\e{0}%
\e{0}%
\e{1}%
\e{0}%
\e{1}%
\e{0}%
\e{0}%
\e{0}%
\e{1}%
\e{3}%
\e{1}%
\e{2}%
\e{0}%
\e{2}%
\e{3}%
\eol}\vss}\rg%
%
%
\rx{\vss\hfull{%
\rlx{\hss{$70_{y}$}}\cg%
\e{0}%
\e{0}%
\e{0}%
\e{0}%
\e{0}%
\e{0}%
\e{0}%
\e{0}%
\e{0}%
\e{0}%
\e{0}%
\e{0}%
\e{0}%
\e{0}%
\e{0}%
\e{0}%
\e{0}%
\e{0}%
\e{0}%
\e{0}%
\e{0}%
\e{1}%
\eol}\vss}\rg%
%
%
\rx{\vss\hfull{%
\rlx{\hss{$1134_{y}$}}\cg%
\e{1}%
\e{0}%
\e{0}%
\e{1}%
\e{1}%
\e{1}%
\e{1}%
\e{0}%
\e{0}%
\e{1}%
\e{1}%
\e{0}%
\e{0}%
\e{0}%
\e{0}%
\e{0}%
\e{1}%
\e{1}%
\e{1}%
\e{1}%
\e{1}%
\e{1}%
\eol}\vss}\rg%
%
%
\rx{\vss\hfull{%
\rlx{\hss{$1680_{y}$}}\cg%
\e{1}%
\e{1}%
\e{0}%
\e{1}%
\e{1}%
\e{1}%
\e{1}%
\e{0}%
\e{0}%
\e{1}%
\e{0}%
\e{1}%
\e{0}%
\e{0}%
\e{0}%
\e{1}%
\e{2}%
\e{0}%
\e{1}%
\e{0}%
\e{2}%
\e{4}%
\eol}\vss}\rg%
%
%
\rx{\vss\hfull{%
\rlx{\hss{$168_{y}$}}\cg%
\e{0}%
\e{0}%
\e{0}%
\e{0}%
\e{0}%
\e{0}%
\e{0}%
\e{0}%
\e{0}%
\e{0}%
\e{0}%
\e{0}%
\e{0}%
\e{0}%
\e{0}%
\e{0}%
\e{0}%
\e{0}%
\e{0}%
\e{0}%
\e{0}%
\e{0}%
\eol}\vss}\rg%
%
%
\rx{\vss\hfull{%
\rlx{\hss{$420_{y}$}}\cg%
\e{0}%
\e{0}%
\e{0}%
\e{0}%
\e{1}%
\e{1}%
\e{0}%
\e{0}%
\e{0}%
\e{0}%
\e{0}%
\e{0}%
\e{1}%
\e{0}%
\e{0}%
\e{0}%
\e{0}%
\e{1}%
\e{1}%
\e{0}%
\e{0}%
\e{0}%
\eol}\vss}\rg%
%
%
\rx{\vss\hfull{%
\rlx{\hss{$3150_{y}$}}\cg%
\e{2}%
\e{0}%
\e{1}%
\e{2}%
\e{4}%
\e{4}%
\e{2}%
\e{1}%
\e{0}%
\e{1}%
\e{2}%
\e{1}%
\e{2}%
\e{0}%
\e{0}%
\e{0}%
\e{3}%
\e{4}%
\e{4}%
\e{1}%
\e{1}%
\e{1}%
\eol}\vss}\rg%
%
%
\rx{\vss\hfull{%
\rlx{\hss{$4200_{y}$}}\cg%
\e{3}%
\e{0}%
\e{2}%
\e{1}%
\e{4}%
\e{4}%
\e{1}%
\e{2}%
\e{0}%
\e{0}%
\e{2}%
\e{1}%
\e{2}%
\e{1}%
\e{0}%
\e{0}%
\e{3}%
\e{4}%
\e{4}%
\e{1}%
\e{2}%
\e{2}%
\eol}\vss}\rg%
%
%
\rx{\vss\hfull{%
\rlx{\hss{$2688_{y}$}}\cg%
\e{2}%
\e{1}%
\e{1}%
\e{1}%
\e{2}%
\e{2}%
\e{1}%
\e{1}%
\e{0}%
\e{0}%
\e{1}%
\e{1}%
\e{1}%
\e{0}%
\e{0}%
\e{0}%
\e{2}%
\e{2}%
\e{2}%
\e{0}%
\e{2}%
\e{2}%
\eol}\vss}\rg%
%
%
\rx{\vss\hfull{%
\rlx{\hss{$2100_{y}$}}\cg%
\e{2}%
\e{0}%
\e{0}%
\e{1}%
\e{1}%
\e{1}%
\e{1}%
\e{0}%
\e{0}%
\e{1}%
\e{0}%
\e{1}%
\e{0}%
\e{0}%
\e{0}%
\e{1}%
\e{2}%
\e{0}%
\e{1}%
\e{0}%
\e{2}%
\e{2}%
\eol}\vss}\rg%
%
%
\rx{\vss\hfull{%
\rlx{\hss{$1400_{y}$}}\cg%
\e{1}%
\e{0}%
\e{1}%
\e{0}%
\e{1}%
\e{1}%
\e{0}%
\e{1}%
\e{0}%
\e{0}%
\e{0}%
\e{1}%
\e{0}%
\e{1}%
\e{0}%
\e{1}%
\e{2}%
\e{0}%
\e{1}%
\e{0}%
\e{1}%
\e{2}%
\eol}\vss}\rg%
%
%
\rx{\vss\hfull{%
\rlx{\hss{$4536_{y}$}}\cg%
\e{3}%
\e{1}%
\e{2}%
\e{2}%
\e{4}%
\e{4}%
\e{2}%
\e{2}%
\e{0}%
\e{1}%
\e{1}%
\e{3}%
\e{1}%
\e{1}%
\e{0}%
\e{2}%
\e{6}%
\e{2}%
\e{4}%
\e{0}%
\e{3}%
\e{4}%
\eol}\vss}\rg%
%
%
\rx{\vss\hfull{%
\rlx{\hss{$5670_{y}$}}\cg%
\e{4}%
\e{1}%
\e{2}%
\e{3}%
\e{5}%
\e{5}%
\e{3}%
\e{2}%
\e{0}%
\e{2}%
\e{2}%
\e{3}%
\e{1}%
\e{1}%
\e{0}%
\e{2}%
\e{7}%
\e{3}%
\e{5}%
\e{1}%
\e{4}%
\e{5}%
\eol}\vss}\rg%
%
%
\rx{\vss\hfull{%
\rlx{\hss{$4480_{y}$}}\cg%
\e{3}%
\e{0}%
\e{2}%
\e{2}%
\e{5}%
\e{5}%
\e{2}%
\e{2}%
\e{0}%
\e{1}%
\e{2}%
\e{2}%
\e{2}%
\e{1}%
\e{0}%
\e{1}%
\e{5}%
\e{4}%
\e{5}%
\e{1}%
\e{2}%
\e{2}%
\eol}\vss}\rg%
\eop
\eject
\tablecont%
%
%
%
%
%
%
\rowpts=18 true pt%
\colpts=18 true pt%
\rowlabpts=40 true pt%
\collabpts=85 true pt%
\clx{\vss\hfull{%
\rlx{\hss{$ $}}\cg%
\cx{\hskip 16 true pt\flip{$[{2}{1^{2}}:{1}]{\times}[{2^{2}}]$}\hss}\cg%
\cx{\hskip 16 true pt\flip{$[{1^{4}}:{1}]{\times}[{2^{2}}]$}\hss}\cg%
\cx{\hskip 16 true pt\flip{$[{3}:{2}]{\times}[{2^{2}}]$}\hss}\cg%
\cx{\hskip 16 true pt\flip{$[{3}:{1^{2}}]{\times}[{2^{2}}]$}\hss}\cg%
\cx{\hskip 16 true pt\flip{$[{2}{1}:{2}]{\times}[{2^{2}}]$}\hss}\cg%
\cx{\hskip 16 true pt\flip{$[{2}{1}:{1^{2}}]{\times}[{2^{2}}]$}\hss}\cg%
\cx{\hskip 16 true pt\flip{$[{1^{3}}:{2}]{\times}[{2^{2}}]$}\hss}\cg%
\cx{\hskip 16 true pt\flip{$[{1^{3}}:{1^{2}}]{\times}[{2^{2}}]$}\hss}\cg%
\cx{\hskip 16 true pt\flip{$[{5}:-]{\times}[{2}{1^{2}}]$}\hss}\cg%
\cx{\hskip 16 true pt\flip{$[{4}{1}:-]{\times}[{2}{1^{2}}]$}\hss}\cg%
\cx{\hskip 16 true pt\flip{$[{3}{2}:-]{\times}[{2}{1^{2}}]$}\hss}\cg%
\cx{\hskip 16 true pt\flip{$[{3}{1^{2}}:-]{\times}[{2}{1^{2}}]$}\hss}\cg%
\cx{\hskip 16 true pt\flip{$[{2^{2}}{1}:-]{\times}[{2}{1^{2}}]$}\hss}\cg%
\cx{\hskip 16 true pt\flip{$[{2}{1^{3}}:-]{\times}[{2}{1^{2}}]$}\hss}\cg%
\cx{\hskip 16 true pt\flip{$[{1^{5}}:-]{\times}[{2}{1^{2}}]$}\hss}\cg%
\cx{\hskip 16 true pt\flip{$[{4}:{1}]{\times}[{2}{1^{2}}]$}\hss}\cg%
\cx{\hskip 16 true pt\flip{$[{3}{1}:{1}]{\times}[{2}{1^{2}}]$}\hss}\cg%
\cx{\hskip 16 true pt\flip{$[{2^{2}}:{1}]{\times}[{2}{1^{2}}]$}\hss}\cg%
\cx{\hskip 16 true pt\flip{$[{2}{1^{2}}:{1}]{\times}[{2}{1^{2}}]$}\hss}\cg%
\cx{\hskip 16 true pt\flip{$[{1^{4}}:{1}]{\times}[{2}{1^{2}}]$}\hss}\cg%
\cx{\hskip 16 true pt\flip{$[{3}:{2}]{\times}[{2}{1^{2}}]$}\hss}\cg%
\cx{\hskip 16 true pt\flip{$[{3}:{1^{2}}]{\times}[{2}{1^{2}}]$}\hss}\cg%
\eol}}\rg%
%
%
\rx{\vss\hfull{%
\rlx{\hss{$8_{z}$}}\cg%
\e{0}%
\e{0}%
\e{0}%
\e{0}%
\e{0}%
\e{0}%
\e{0}%
\e{0}%
\e{0}%
\e{0}%
\e{0}%
\e{0}%
\e{0}%
\e{0}%
\e{0}%
\e{0}%
\e{0}%
\e{0}%
\e{0}%
\e{0}%
\e{0}%
\e{0}%
\eol}\vss}\rg%
%
%
\rx{\vss\hfull{%
\rlx{\hss{$56_{z}$}}\cg%
\e{0}%
\e{0}%
\e{0}%
\e{0}%
\e{0}%
\e{0}%
\e{0}%
\e{0}%
\e{0}%
\e{0}%
\e{0}%
\e{0}%
\e{0}%
\e{0}%
\e{0}%
\e{1}%
\e{0}%
\e{0}%
\e{0}%
\e{0}%
\e{0}%
\e{0}%
\eol}\vss}\rg%
%
%
\rx{\vss\hfull{%
\rlx{\hss{$160_{z}$}}\cg%
\e{0}%
\e{0}%
\e{0}%
\e{0}%
\e{0}%
\e{0}%
\e{0}%
\e{0}%
\e{1}%
\e{0}%
\e{0}%
\e{0}%
\e{0}%
\e{0}%
\e{0}%
\e{1}%
\e{0}%
\e{0}%
\e{0}%
\e{0}%
\e{0}%
\e{0}%
\eol}\vss}\rg%
%
%
\rx{\vss\hfull{%
\rlx{\hss{$112_{z}$}}\cg%
\e{0}%
\e{0}%
\e{0}%
\e{0}%
\e{0}%
\e{0}%
\e{0}%
\e{0}%
\e{1}%
\e{0}%
\e{0}%
\e{0}%
\e{0}%
\e{0}%
\e{0}%
\e{0}%
\e{0}%
\e{0}%
\e{0}%
\e{0}%
\e{0}%
\e{0}%
\eol}\vss}\rg%
%
%
\rx{\vss\hfull{%
\rlx{\hss{$840_{z}$}}\cg%
\e{0}%
\e{0}%
\e{1}%
\e{1}%
\e{1}%
\e{0}%
\e{0}%
\e{0}%
\e{0}%
\e{1}%
\e{0}%
\e{0}%
\e{0}%
\e{0}%
\e{0}%
\e{1}%
\e{1}%
\e{0}%
\e{0}%
\e{0}%
\e{1}%
\e{1}%
\eol}\vss}\rg%
%
%
\rx{\vss\hfull{%
\rlx{\hss{$1296_{z}$}}\cg%
\e{0}%
\e{0}%
\e{1}%
\e{1}%
\e{1}%
\e{0}%
\e{1}%
\e{0}%
\e{1}%
\e{1}%
\e{0}%
\e{0}%
\e{0}%
\e{0}%
\e{0}%
\e{4}%
\e{2}%
\e{0}%
\e{0}%
\e{0}%
\e{2}%
\e{2}%
\eol}\vss}\rg%
%
%
\rx{\vss\hfull{%
\rlx{\hss{$1400_{z}$}}\cg%
\e{0}%
\e{0}%
\e{2}%
\e{1}%
\e{1}%
\e{0}%
\e{0}%
\e{0}%
\e{2}%
\e{2}%
\e{0}%
\e{0}%
\e{0}%
\e{0}%
\e{0}%
\e{2}%
\e{1}%
\e{0}%
\e{0}%
\e{0}%
\e{2}%
\e{1}%
\eol}\vss}\rg%
%
%
\rx{\vss\hfull{%
\rlx{\hss{$1008_{z}$}}\cg%
\e{0}%
\e{0}%
\e{1}%
\e{0}%
\e{1}%
\e{0}%
\e{0}%
\e{0}%
\e{2}%
\e{1}%
\e{0}%
\e{0}%
\e{0}%
\e{0}%
\e{0}%
\e{3}%
\e{1}%
\e{0}%
\e{0}%
\e{0}%
\e{2}%
\e{1}%
\eol}\vss}\rg%
%
%
\rx{\vss\hfull{%
\rlx{\hss{$560_{z}$}}\cg%
\e{0}%
\e{0}%
\e{1}%
\e{0}%
\e{0}%
\e{0}%
\e{0}%
\e{0}%
\e{1}%
\e{1}%
\e{0}%
\e{0}%
\e{0}%
\e{0}%
\e{0}%
\e{1}%
\e{0}%
\e{0}%
\e{0}%
\e{0}%
\e{1}%
\e{0}%
\eol}\vss}\rg%
%
%
\rx{\vss\hfull{%
\rlx{\hss{$1400_{zz}$}}\cg%
\e{1}%
\e{0}%
\e{2}%
\e{1}%
\e{1}%
\e{1}%
\e{0}%
\e{0}%
\e{0}%
\e{0}%
\e{2}%
\e{0}%
\e{1}%
\e{0}%
\e{0}%
\e{0}%
\e{1}%
\e{1}%
\e{0}%
\e{0}%
\e{1}%
\e{0}%
\eol}\vss}\rg%
%
%
\rx{\vss\hfull{%
\rlx{\hss{$4200_{z}$}}\cg%
\e{3}%
\e{1}%
\e{2}%
\e{1}%
\e{5}%
\e{4}%
\e{1}%
\e{1}%
\e{0}%
\e{0}%
\e{3}%
\e{1}%
\e{2}%
\e{1}%
\e{0}%
\e{1}%
\e{5}%
\e{4}%
\e{2}%
\e{0}%
\e{4}%
\e{2}%
\eol}\vss}\rg%
%
%
\rx{\vss\hfull{%
\rlx{\hss{$400_{z}$}}\cg%
\e{0}%
\e{0}%
\e{1}%
\e{0}%
\e{0}%
\e{0}%
\e{0}%
\e{0}%
\e{0}%
\e{0}%
\e{1}%
\e{0}%
\e{0}%
\e{0}%
\e{0}%
\e{0}%
\e{0}%
\e{0}%
\e{0}%
\e{0}%
\e{1}%
\e{0}%
\eol}\vss}\rg%
%
%
\rx{\vss\hfull{%
\rlx{\hss{$3240_{z}$}}\cg%
\e{1}%
\e{0}%
\e{4}%
\e{2}%
\e{3}%
\e{1}%
\e{0}%
\e{0}%
\e{1}%
\e{2}%
\e{2}%
\e{1}%
\e{0}%
\e{0}%
\e{0}%
\e{3}%
\e{3}%
\e{1}%
\e{0}%
\e{0}%
\e{5}%
\e{2}%
\eol}\vss}\rg%
%
%
\rx{\vss\hfull{%
\rlx{\hss{$4536_{z}$}}\cg%
\e{2}%
\e{0}%
\e{4}%
\e{3}%
\e{5}%
\e{4}%
\e{0}%
\e{0}%
\e{0}%
\e{2}%
\e{2}%
\e{2}%
\e{1}%
\e{0}%
\e{0}%
\e{1}%
\e{4}%
\e{3}%
\e{2}%
\e{0}%
\e{3}%
\e{2}%
\eol}\vss}\rg%
%
%
\rx{\vss\hfull{%
\rlx{\hss{$2400_{z}$}}\cg%
\e{1}%
\e{0}%
\e{1}%
\e{2}%
\e{2}%
\e{1}%
\e{2}%
\e{1}%
\e{0}%
\e{1}%
\e{0}%
\e{1}%
\e{0}%
\e{0}%
\e{0}%
\e{3}%
\e{4}%
\e{1}%
\e{1}%
\e{0}%
\e{3}%
\e{3}%
\eol}\vss}\rg%
%
%
\rx{\vss\hfull{%
\rlx{\hss{$3360_{z}$}}\cg%
\e{2}%
\e{0}%
\e{3}%
\e{2}%
\e{3}%
\e{2}%
\e{1}%
\e{1}%
\e{0}%
\e{1}%
\e{2}%
\e{1}%
\e{1}%
\e{0}%
\e{0}%
\e{2}%
\e{4}%
\e{2}%
\e{1}%
\e{0}%
\e{4}%
\e{2}%
\eol}\vss}\rg%
%
%
\rx{\vss\hfull{%
\rlx{\hss{$2800_{z}$}}\cg%
\e{1}%
\e{0}%
\e{2}%
\e{1}%
\e{3}%
\e{2}%
\e{1}%
\e{0}%
\e{1}%
\e{2}%
\e{0}%
\e{1}%
\e{0}%
\e{0}%
\e{0}%
\e{3}%
\e{4}%
\e{1}%
\e{1}%
\e{0}%
\e{4}%
\e{3}%
\eol}\vss}\rg%
%
%
\rx{\vss\hfull{%
\rlx{\hss{$4096_{z}$}}\cg%
\e{1}%
\e{0}%
\e{4}%
\e{3}%
\e{4}%
\e{2}%
\e{1}%
\e{0}%
\e{1}%
\e{3}%
\e{1}%
\e{1}%
\e{0}%
\e{0}%
\e{0}%
\e{4}%
\e{5}%
\e{1}%
\e{1}%
\e{0}%
\e{5}%
\e{4}%
\eol}\vss}\rg%
%
%
\rx{\vss\hfull{%
\rlx{\hss{$5600_{z}$}}\cg%
\e{2}%
\e{0}%
\e{3}%
\e{4}%
\e{6}%
\e{4}%
\e{2}%
\e{1}%
\e{1}%
\e{3}%
\e{1}%
\e{2}%
\e{0}%
\e{0}%
\e{0}%
\e{3}%
\e{7}%
\e{2}%
\e{4}%
\e{0}%
\e{5}%
\e{6}%
\eol}\vss}\rg%
%
%
\rx{\vss\hfull{%
\rlx{\hss{$448_{z}$}}\cg%
\e{0}%
\e{0}%
\e{1}%
\e{1}%
\e{0}%
\e{0}%
\e{0}%
\e{0}%
\e{0}%
\e{1}%
\e{0}%
\e{0}%
\e{0}%
\e{0}%
\e{0}%
\e{0}%
\e{0}%
\e{0}%
\e{0}%
\e{0}%
\e{0}%
\e{0}%
\eol}\vss}\rg%
%
%
\rx{\vss\hfull{%
\rlx{\hss{$448_{w}$}}\cg%
\e{0}%
\e{0}%
\e{0}%
\e{1}%
\e{0}%
\e{0}%
\e{1}%
\e{0}%
\e{0}%
\e{0}%
\e{0}%
\e{0}%
\e{0}%
\e{0}%
\e{0}%
\e{1}%
\e{1}%
\e{0}%
\e{0}%
\e{0}%
\e{0}%
\e{1}%
\eol}\vss}\rg%
%
%
\rx{\vss\hfull{%
\rlx{\hss{$1344_{w}$}}\cg%
\e{1}%
\e{0}%
\e{0}%
\e{0}%
\e{2}%
\e{2}%
\e{0}%
\e{0}%
\e{0}%
\e{0}%
\e{1}%
\e{0}%
\e{1}%
\e{1}%
\e{0}%
\e{0}%
\e{1}%
\e{2}%
\e{1}%
\e{0}%
\e{0}%
\e{0}%
\eol}\vss}\rg%
%
%
\rx{\vss\hfull{%
\rlx{\hss{$5600_{w}$}}\cg%
\e{3}%
\e{1}%
\e{2}%
\e{3}%
\e{5}%
\e{5}%
\e{3}%
\e{2}%
\e{0}%
\e{1}%
\e{2}%
\e{2}%
\e{1}%
\e{1}%
\e{0}%
\e{2}%
\e{7}%
\e{3}%
\e{5}%
\e{1}%
\e{4}%
\e{5}%
\eol}\vss}\rg%
%
%
\rx{\vss\hfull{%
\rlx{\hss{$2016_{w}$}}\cg%
\e{2}%
\e{0}%
\e{1}%
\e{1}%
\e{2}%
\e{2}%
\e{1}%
\e{1}%
\e{0}%
\e{0}%
\e{1}%
\e{1}%
\e{2}%
\e{0}%
\e{0}%
\e{0}%
\e{1}%
\e{3}%
\e{2}%
\e{0}%
\e{0}%
\e{0}%
\eol}\vss}\rg%
%
%
\rx{\vss\hfull{%
\rlx{\hss{$7168_{w}$}}\cg%
\e{5}%
\e{1}%
\e{3}%
\e{3}%
\e{7}%
\e{7}%
\e{3}%
\e{3}%
\e{0}%
\e{1}%
\e{3}%
\e{3}%
\e{3}%
\e{1}%
\e{0}%
\e{1}%
\e{7}%
\e{6}%
\e{7}%
\e{1}%
\e{4}%
\e{4}%
\eol}\vss}\rg%
\eop
\eject
\tablecont%
%
%
%
%
%
%
\rowpts=18 true pt%
\colpts=18 true pt%
\rowlabpts=40 true pt%
\collabpts=85 true pt%
\clx{\vss\hfull{%
\rlx{\hss{$ $}}\cg%
\cx{\hskip 16 true pt\flip{$[{2}{1}:{2}]{\times}[{2}{1^{2}}]$}\hss}\cg%
\cx{\hskip 16 true pt\flip{$[{2}{1}:{1^{2}}]{\times}[{2}{1^{2}}]$}\hss}\cg%
\cx{\hskip 16 true pt\flip{$[{1^{3}}:{2}]{\times}[{2}{1^{2}}]$}\hss}\cg%
\cx{\hskip 16 true pt\flip{$[{1^{3}}:{1^{2}}]{\times}[{2}{1^{2}}]$}\hss}\cg%
\cx{\hskip 16 true pt\flip{$[{5}:-]{\times}[{1^{4}}]$}\hss}\cg%
\cx{\hskip 16 true pt\flip{$[{4}{1}:-]{\times}[{1^{4}}]$}\hss}\cg%
\cx{\hskip 16 true pt\flip{$[{3}{2}:-]{\times}[{1^{4}}]$}\hss}\cg%
\cx{\hskip 16 true pt\flip{$[{3}{1^{2}}:-]{\times}[{1^{4}}]$}\hss}\cg%
\cx{\hskip 16 true pt\flip{$[{2^{2}}{1}:-]{\times}[{1^{4}}]$}\hss}\cg%
\cx{\hskip 16 true pt\flip{$[{2}{1^{3}}:-]{\times}[{1^{4}}]$}\hss}\cg%
\cx{\hskip 16 true pt\flip{$[{1^{5}}:-]{\times}[{1^{4}}]$}\hss}\cg%
\cx{\hskip 16 true pt\flip{$[{4}:{1}]{\times}[{1^{4}}]$}\hss}\cg%
\cx{\hskip 16 true pt\flip{$[{3}{1}:{1}]{\times}[{1^{4}}]$}\hss}\cg%
\cx{\hskip 16 true pt\flip{$[{2^{2}}:{1}]{\times}[{1^{4}}]$}\hss}\cg%
\cx{\hskip 16 true pt\flip{$[{2}{1^{2}}:{1}]{\times}[{1^{4}}]$}\hss}\cg%
\cx{\hskip 16 true pt\flip{$[{1^{4}}:{1}]{\times}[{1^{4}}]$}\hss}\cg%
\cx{\hskip 16 true pt\flip{$[{3}:{2}]{\times}[{1^{4}}]$}\hss}\cg%
\cx{\hskip 16 true pt\flip{$[{3}:{1^{2}}]{\times}[{1^{4}}]$}\hss}\cg%
\cx{\hskip 16 true pt\flip{$[{2}{1}:{2}]{\times}[{1^{4}}]$}\hss}\cg%
\cx{\hskip 16 true pt\flip{$[{2}{1}:{1^{2}}]{\times}[{1^{4}}]$}\hss}\cg%
\cx{\hskip 16 true pt\flip{$[{1^{3}}:{2}]{\times}[{1^{4}}]$}\hss}\cg%
\cx{\hskip 16 true pt\flip{$[{1^{3}}:{1^{2}}]{\times}[{1^{4}}]$}\hss}\cg%
\eol}}\rg%
%
%
\rx{\vss\hfull{%
\rlx{\hss{$1_{x}$}}\cg%
\e{0}%
\e{0}%
\e{0}%
\e{0}%
\e{0}%
\e{0}%
\e{0}%
\e{0}%
\e{0}%
\e{0}%
\e{0}%
\e{0}%
\e{0}%
\e{0}%
\e{0}%
\e{0}%
\e{0}%
\e{0}%
\e{0}%
\e{0}%
\e{0}%
\e{0}%
\eol}\vss}\rg%
%
%
\rx{\vss\hfull{%
\rlx{\hss{$28_{x}$}}\cg%
\e{0}%
\e{0}%
\e{0}%
\e{0}%
\e{0}%
\e{0}%
\e{0}%
\e{0}%
\e{0}%
\e{0}%
\e{0}%
\e{0}%
\e{0}%
\e{0}%
\e{0}%
\e{0}%
\e{0}%
\e{0}%
\e{0}%
\e{0}%
\e{0}%
\e{0}%
\eol}\vss}\rg%
%
%
\rx{\vss\hfull{%
\rlx{\hss{$35_{x}$}}\cg%
\e{0}%
\e{0}%
\e{0}%
\e{0}%
\e{0}%
\e{0}%
\e{0}%
\e{0}%
\e{0}%
\e{0}%
\e{0}%
\e{0}%
\e{0}%
\e{0}%
\e{0}%
\e{0}%
\e{0}%
\e{0}%
\e{0}%
\e{0}%
\e{0}%
\e{0}%
\eol}\vss}\rg%
%
%
\rx{\vss\hfull{%
\rlx{\hss{$84_{x}$}}\cg%
\e{0}%
\e{0}%
\e{0}%
\e{0}%
\e{0}%
\e{0}%
\e{0}%
\e{0}%
\e{0}%
\e{0}%
\e{0}%
\e{0}%
\e{0}%
\e{0}%
\e{0}%
\e{0}%
\e{0}%
\e{0}%
\e{0}%
\e{0}%
\e{0}%
\e{0}%
\eol}\vss}\rg%
%
%
\rx{\vss\hfull{%
\rlx{\hss{$50_{x}$}}\cg%
\e{0}%
\e{0}%
\e{0}%
\e{0}%
\e{0}%
\e{0}%
\e{0}%
\e{0}%
\e{0}%
\e{0}%
\e{0}%
\e{0}%
\e{0}%
\e{0}%
\e{0}%
\e{0}%
\e{0}%
\e{0}%
\e{0}%
\e{0}%
\e{0}%
\e{0}%
\eol}\vss}\rg%
%
%
\rx{\vss\hfull{%
\rlx{\hss{$350_{x}$}}\cg%
\e{0}%
\e{0}%
\e{0}%
\e{0}%
\e{0}%
\e{0}%
\e{0}%
\e{0}%
\e{0}%
\e{0}%
\e{0}%
\e{1}%
\e{0}%
\e{0}%
\e{0}%
\e{0}%
\e{0}%
\e{0}%
\e{0}%
\e{0}%
\e{0}%
\e{0}%
\eol}\vss}\rg%
%
%
\rx{\vss\hfull{%
\rlx{\hss{$300_{x}$}}\cg%
\e{0}%
\e{0}%
\e{0}%
\e{0}%
\e{0}%
\e{0}%
\e{0}%
\e{0}%
\e{0}%
\e{0}%
\e{0}%
\e{0}%
\e{0}%
\e{0}%
\e{0}%
\e{0}%
\e{0}%
\e{0}%
\e{0}%
\e{0}%
\e{0}%
\e{0}%
\eol}\vss}\rg%
%
%
\rx{\vss\hfull{%
\rlx{\hss{$567_{x}$}}\cg%
\e{0}%
\e{0}%
\e{0}%
\e{0}%
\e{1}%
\e{0}%
\e{0}%
\e{0}%
\e{0}%
\e{0}%
\e{0}%
\e{0}%
\e{0}%
\e{0}%
\e{0}%
\e{0}%
\e{0}%
\e{0}%
\e{0}%
\e{0}%
\e{0}%
\e{0}%
\eol}\vss}\rg%
%
%
\rx{\vss\hfull{%
\rlx{\hss{$210_{x}$}}\cg%
\e{0}%
\e{0}%
\e{0}%
\e{0}%
\e{0}%
\e{0}%
\e{0}%
\e{0}%
\e{0}%
\e{0}%
\e{0}%
\e{0}%
\e{0}%
\e{0}%
\e{0}%
\e{0}%
\e{0}%
\e{0}%
\e{0}%
\e{0}%
\e{0}%
\e{0}%
\eol}\vss}\rg%
%
%
\rx{\vss\hfull{%
\rlx{\hss{$840_{x}$}}\cg%
\e{1}%
\e{1}%
\e{0}%
\e{0}%
\e{0}%
\e{0}%
\e{0}%
\e{0}%
\e{0}%
\e{0}%
\e{0}%
\e{0}%
\e{0}%
\e{0}%
\e{0}%
\e{0}%
\e{0}%
\e{0}%
\e{0}%
\e{0}%
\e{1}%
\e{0}%
\eol}\vss}\rg%
%
%
\rx{\vss\hfull{%
\rlx{\hss{$700_{x}$}}\cg%
\e{0}%
\e{0}%
\e{0}%
\e{0}%
\e{0}%
\e{0}%
\e{0}%
\e{0}%
\e{0}%
\e{0}%
\e{0}%
\e{0}%
\e{0}%
\e{0}%
\e{0}%
\e{0}%
\e{0}%
\e{0}%
\e{0}%
\e{0}%
\e{0}%
\e{0}%
\eol}\vss}\rg%
%
%
\rx{\vss\hfull{%
\rlx{\hss{$175_{x}$}}\cg%
\e{0}%
\e{0}%
\e{0}%
\e{0}%
\e{0}%
\e{0}%
\e{0}%
\e{0}%
\e{0}%
\e{0}%
\e{0}%
\e{0}%
\e{0}%
\e{0}%
\e{0}%
\e{0}%
\e{0}%
\e{0}%
\e{0}%
\e{0}%
\e{0}%
\e{0}%
\eol}\vss}\rg%
%
%
\rx{\vss\hfull{%
\rlx{\hss{$1400_{x}$}}\cg%
\e{1}%
\e{0}%
\e{0}%
\e{0}%
\e{0}%
\e{0}%
\e{0}%
\e{0}%
\e{0}%
\e{0}%
\e{0}%
\e{0}%
\e{0}%
\e{0}%
\e{0}%
\e{0}%
\e{1}%
\e{0}%
\e{0}%
\e{0}%
\e{0}%
\e{0}%
\eol}\vss}\rg%
%
%
\rx{\vss\hfull{%
\rlx{\hss{$1050_{x}$}}\cg%
\e{1}%
\e{0}%
\e{0}%
\e{0}%
\e{0}%
\e{0}%
\e{1}%
\e{0}%
\e{0}%
\e{0}%
\e{0}%
\e{0}%
\e{0}%
\e{0}%
\e{0}%
\e{0}%
\e{0}%
\e{0}%
\e{0}%
\e{0}%
\e{0}%
\e{0}%
\eol}\vss}\rg%
%
%
\rx{\vss\hfull{%
\rlx{\hss{$1575_{x}$}}\cg%
\e{1}%
\e{0}%
\e{0}%
\e{0}%
\e{0}%
\e{0}%
\e{0}%
\e{0}%
\e{0}%
\e{0}%
\e{0}%
\e{1}%
\e{0}%
\e{0}%
\e{0}%
\e{0}%
\e{1}%
\e{0}%
\e{0}%
\e{0}%
\e{0}%
\e{0}%
\eol}\vss}\rg%
%
%
\rx{\vss\hfull{%
\rlx{\hss{$1344_{x}$}}\cg%
\e{1}%
\e{0}%
\e{0}%
\e{0}%
\e{0}%
\e{1}%
\e{0}%
\e{0}%
\e{0}%
\e{0}%
\e{0}%
\e{0}%
\e{0}%
\e{0}%
\e{0}%
\e{0}%
\e{0}%
\e{0}%
\e{0}%
\e{0}%
\e{0}%
\e{0}%
\eol}\vss}\rg%
%
%
\rx{\vss\hfull{%
\rlx{\hss{$2100_{x}$}}\cg%
\e{3}%
\e{1}%
\e{1}%
\e{0}%
\e{1}%
\e{1}%
\e{0}%
\e{0}%
\e{0}%
\e{0}%
\e{0}%
\e{2}%
\e{1}%
\e{0}%
\e{0}%
\e{0}%
\e{1}%
\e{1}%
\e{0}%
\e{0}%
\e{0}%
\e{0}%
\eol}\vss}\rg%
%
%
\rx{\vss\hfull{%
\rlx{\hss{$2268_{x}$}}\cg%
\e{2}%
\e{1}%
\e{0}%
\e{0}%
\e{1}%
\e{0}%
\e{0}%
\e{0}%
\e{0}%
\e{0}%
\e{0}%
\e{1}%
\e{0}%
\e{0}%
\e{0}%
\e{0}%
\e{1}%
\e{1}%
\e{0}%
\e{0}%
\e{0}%
\e{0}%
\eol}\vss}\rg%
%
%
\rx{\vss\hfull{%
\rlx{\hss{$525_{x}$}}\cg%
\e{1}%
\e{0}%
\e{0}%
\e{0}%
\e{1}%
\e{1}%
\e{0}%
\e{0}%
\e{0}%
\e{0}%
\e{0}%
\e{0}%
\e{0}%
\e{0}%
\e{0}%
\e{0}%
\e{0}%
\e{0}%
\e{0}%
\e{0}%
\e{0}%
\e{0}%
\eol}\vss}\rg%
%
%
\rx{\vss\hfull{%
\rlx{\hss{$700_{xx}$}}\cg%
\e{1}%
\e{1}%
\e{0}%
\e{0}%
\e{0}%
\e{0}%
\e{1}%
\e{0}%
\e{1}%
\e{0}%
\e{0}%
\e{0}%
\e{0}%
\e{0}%
\e{0}%
\e{0}%
\e{0}%
\e{0}%
\e{0}%
\e{0}%
\e{0}%
\e{0}%
\eol}\vss}\rg%
%
%
\rx{\vss\hfull{%
\rlx{\hss{$972_{x}$}}\cg%
\e{1}%
\e{1}%
\e{0}%
\e{0}%
\e{0}%
\e{0}%
\e{0}%
\e{1}%
\e{0}%
\e{0}%
\e{0}%
\e{0}%
\e{0}%
\e{0}%
\e{0}%
\e{0}%
\e{0}%
\e{0}%
\e{0}%
\e{0}%
\e{0}%
\e{0}%
\eol}\vss}\rg%
%
%
\rx{\vss\hfull{%
\rlx{\hss{$4096_{x}$}}\cg%
\e{5}%
\e{2}%
\e{1}%
\e{0}%
\e{0}%
\e{1}%
\e{0}%
\e{0}%
\e{0}%
\e{0}%
\e{0}%
\e{1}%
\e{1}%
\e{0}%
\e{0}%
\e{0}%
\e{1}%
\e{1}%
\e{1}%
\e{0}%
\e{0}%
\e{0}%
\eol}\vss}\rg%
%
%
\rx{\vss\hfull{%
\rlx{\hss{$4200_{x}$}}\cg%
\e{5}%
\e{3}%
\e{1}%
\e{0}%
\e{0}%
\e{0}%
\e{1}%
\e{0}%
\e{0}%
\e{0}%
\e{0}%
\e{0}%
\e{1}%
\e{0}%
\e{0}%
\e{0}%
\e{1}%
\e{0}%
\e{1}%
\e{1}%
\e{0}%
\e{0}%
\eol}\vss}\rg%
%
%
\rx{\vss\hfull{%
\rlx{\hss{$2240_{x}$}}\cg%
\e{2}%
\e{1}%
\e{0}%
\e{0}%
\e{0}%
\e{0}%
\e{0}%
\e{0}%
\e{0}%
\e{0}%
\e{0}%
\e{0}%
\e{0}%
\e{0}%
\e{0}%
\e{0}%
\e{0}%
\e{0}%
\e{1}%
\e{0}%
\e{0}%
\e{0}%
\eol}\vss}\rg%
%
%
\rx{\vss\hfull{%
\rlx{\hss{$2835_{x}$}}\cg%
\e{3}%
\e{3}%
\e{0}%
\e{1}%
\e{0}%
\e{0}%
\e{0}%
\e{0}%
\e{1}%
\e{0}%
\e{0}%
\e{0}%
\e{0}%
\e{1}%
\e{0}%
\e{0}%
\e{0}%
\e{0}%
\e{1}%
\e{1}%
\e{0}%
\e{0}%
\eol}\vss}\rg%
%
%
\rx{\vss\hfull{%
\rlx{\hss{$6075_{x}$}}\cg%
\e{9}%
\e{5}%
\e{2}%
\e{1}%
\e{0}%
\e{1}%
\e{1}%
\e{1}%
\e{0}%
\e{0}%
\e{0}%
\e{1}%
\e{2}%
\e{1}%
\e{0}%
\e{0}%
\e{2}%
\e{1}%
\e{2}%
\e{1}%
\e{0}%
\e{0}%
\eol}\vss}\rg%
%
%
\rx{\vss\hfull{%
\rlx{\hss{$3200_{x}$}}\cg%
\e{5}%
\e{4}%
\e{2}%
\e{1}%
\e{0}%
\e{1}%
\e{0}%
\e{1}%
\e{0}%
\e{0}%
\e{0}%
\e{0}%
\e{1}%
\e{0}%
\e{1}%
\e{0}%
\e{0}%
\e{1}%
\e{1}%
\e{0}%
\e{1}%
\e{0}%
\eol}\vss}\rg%
\eop
\eject
\tablecont%
%
%
%
%
%
%
\rowpts=18 true pt%
\colpts=18 true pt%
\rowlabpts=40 true pt%
\collabpts=85 true pt%
\clx{\vss\hfull{%
\rlx{\hss{$ $}}\cg%
\cx{\hskip 16 true pt\flip{$[{2}{1}:{2}]{\times}[{2}{1^{2}}]$}\hss}\cg%
\cx{\hskip 16 true pt\flip{$[{2}{1}:{1^{2}}]{\times}[{2}{1^{2}}]$}\hss}\cg%
\cx{\hskip 16 true pt\flip{$[{1^{3}}:{2}]{\times}[{2}{1^{2}}]$}\hss}\cg%
\cx{\hskip 16 true pt\flip{$[{1^{3}}:{1^{2}}]{\times}[{2}{1^{2}}]$}\hss}\cg%
\cx{\hskip 16 true pt\flip{$[{5}:-]{\times}[{1^{4}}]$}\hss}\cg%
\cx{\hskip 16 true pt\flip{$[{4}{1}:-]{\times}[{1^{4}}]$}\hss}\cg%
\cx{\hskip 16 true pt\flip{$[{3}{2}:-]{\times}[{1^{4}}]$}\hss}\cg%
\cx{\hskip 16 true pt\flip{$[{3}{1^{2}}:-]{\times}[{1^{4}}]$}\hss}\cg%
\cx{\hskip 16 true pt\flip{$[{2^{2}}{1}:-]{\times}[{1^{4}}]$}\hss}\cg%
\cx{\hskip 16 true pt\flip{$[{2}{1^{3}}:-]{\times}[{1^{4}}]$}\hss}\cg%
\cx{\hskip 16 true pt\flip{$[{1^{5}}:-]{\times}[{1^{4}}]$}\hss}\cg%
\cx{\hskip 16 true pt\flip{$[{4}:{1}]{\times}[{1^{4}}]$}\hss}\cg%
\cx{\hskip 16 true pt\flip{$[{3}{1}:{1}]{\times}[{1^{4}}]$}\hss}\cg%
\cx{\hskip 16 true pt\flip{$[{2^{2}}:{1}]{\times}[{1^{4}}]$}\hss}\cg%
\cx{\hskip 16 true pt\flip{$[{2}{1^{2}}:{1}]{\times}[{1^{4}}]$}\hss}\cg%
\cx{\hskip 16 true pt\flip{$[{1^{4}}:{1}]{\times}[{1^{4}}]$}\hss}\cg%
\cx{\hskip 16 true pt\flip{$[{3}:{2}]{\times}[{1^{4}}]$}\hss}\cg%
\cx{\hskip 16 true pt\flip{$[{3}:{1^{2}}]{\times}[{1^{4}}]$}\hss}\cg%
\cx{\hskip 16 true pt\flip{$[{2}{1}:{2}]{\times}[{1^{4}}]$}\hss}\cg%
\cx{\hskip 16 true pt\flip{$[{2}{1}:{1^{2}}]{\times}[{1^{4}}]$}\hss}\cg%
\cx{\hskip 16 true pt\flip{$[{1^{3}}:{2}]{\times}[{1^{4}}]$}\hss}\cg%
\cx{\hskip 16 true pt\flip{$[{1^{3}}:{1^{2}}]{\times}[{1^{4}}]$}\hss}\cg%
\eol}}\rg%
%
%
\rx{\vss\hfull{%
\rlx{\hss{$70_{y}$}}\cg%
\e{0}%
\e{0}%
\e{0}%
\e{0}%
\e{0}%
\e{0}%
\e{0}%
\e{0}%
\e{0}%
\e{0}%
\e{0}%
\e{1}%
\e{0}%
\e{0}%
\e{0}%
\e{0}%
\e{0}%
\e{0}%
\e{0}%
\e{0}%
\e{0}%
\e{0}%
\eol}\vss}\rg%
%
%
\rx{\vss\hfull{%
\rlx{\hss{$1134_{y}$}}\cg%
\e{2}%
\e{1}%
\e{1}%
\e{0}%
\e{0}%
\e{0}%
\e{0}%
\e{0}%
\e{0}%
\e{0}%
\e{0}%
\e{0}%
\e{1}%
\e{0}%
\e{0}%
\e{0}%
\e{0}%
\e{0}%
\e{1}%
\e{0}%
\e{1}%
\e{0}%
\eol}\vss}\rg%
%
%
\rx{\vss\hfull{%
\rlx{\hss{$1680_{y}$}}\cg%
\e{2}%
\e{2}%
\e{2}%
\e{0}%
\e{0}%
\e{0}%
\e{0}%
\e{0}%
\e{0}%
\e{0}%
\e{0}%
\e{2}%
\e{2}%
\e{0}%
\e{0}%
\e{0}%
\e{1}%
\e{1}%
\e{1}%
\e{0}%
\e{0}%
\e{0}%
\eol}\vss}\rg%
%
%
\rx{\vss\hfull{%
\rlx{\hss{$168_{y}$}}\cg%
\e{0}%
\e{1}%
\e{0}%
\e{0}%
\e{0}%
\e{0}%
\e{0}%
\e{0}%
\e{0}%
\e{1}%
\e{0}%
\e{0}%
\e{0}%
\e{0}%
\e{0}%
\e{0}%
\e{0}%
\e{0}%
\e{0}%
\e{0}%
\e{0}%
\e{0}%
\eol}\vss}\rg%
%
%
\rx{\vss\hfull{%
\rlx{\hss{$420_{y}$}}\cg%
\e{0}%
\e{1}%
\e{0}%
\e{0}%
\e{0}%
\e{0}%
\e{0}%
\e{0}%
\e{0}%
\e{0}%
\e{0}%
\e{0}%
\e{0}%
\e{0}%
\e{0}%
\e{0}%
\e{0}%
\e{0}%
\e{0}%
\e{0}%
\e{0}%
\e{1}%
\eol}\vss}\rg%
%
%
\rx{\vss\hfull{%
\rlx{\hss{$3150_{y}$}}\cg%
\e{4}%
\e{4}%
\e{1}%
\e{2}%
\e{0}%
\e{0}%
\e{0}%
\e{0}%
\e{0}%
\e{0}%
\e{0}%
\e{0}%
\e{0}%
\e{1}%
\e{1}%
\e{0}%
\e{0}%
\e{0}%
\e{2}%
\e{1}%
\e{1}%
\e{1}%
\eol}\vss}\rg%
%
%
\rx{\vss\hfull{%
\rlx{\hss{$4200_{y}$}}\cg%
\e{6}%
\e{7}%
\e{2}%
\e{2}%
\e{0}%
\e{0}%
\e{1}%
\e{1}%
\e{1}%
\e{1}%
\e{0}%
\e{0}%
\e{1}%
\e{1}%
\e{1}%
\e{0}%
\e{0}%
\e{0}%
\e{1}%
\e{2}%
\e{1}%
\e{1}%
\eol}\vss}\rg%
%
%
\rx{\vss\hfull{%
\rlx{\hss{$2688_{y}$}}\cg%
\e{4}%
\e{4}%
\e{2}%
\e{2}%
\e{0}%
\e{0}%
\e{1}%
\e{1}%
\e{1}%
\e{0}%
\e{0}%
\e{0}%
\e{1}%
\e{1}%
\e{1}%
\e{0}%
\e{0}%
\e{0}%
\e{1}%
\e{1}%
\e{0}%
\e{0}%
\eol}\vss}\rg%
%
%
\rx{\vss\hfull{%
\rlx{\hss{$2100_{y}$}}\cg%
\e{4}%
\e{3}%
\e{2}%
\e{1}%
\e{1}%
\e{2}%
\e{0}%
\e{1}%
\e{0}%
\e{0}%
\e{0}%
\e{1}%
\e{1}%
\e{0}%
\e{1}%
\e{0}%
\e{1}%
\e{1}%
\e{1}%
\e{0}%
\e{0}%
\e{0}%
\eol}\vss}\rg%
%
%
\rx{\vss\hfull{%
\rlx{\hss{$1400_{y}$}}\cg%
\e{2}%
\e{2}%
\e{1}%
\e{0}%
\e{0}%
\e{0}%
\e{0}%
\e{0}%
\e{0}%
\e{0}%
\e{0}%
\e{1}%
\e{1}%
\e{0}%
\e{0}%
\e{0}%
\e{1}%
\e{1}%
\e{0}%
\e{1}%
\e{0}%
\e{0}%
\eol}\vss}\rg%
%
%
\rx{\vss\hfull{%
\rlx{\hss{$4536_{y}$}}\cg%
\e{6}%
\e{6}%
\e{3}%
\e{2}%
\e{0}%
\e{0}%
\e{0}%
\e{0}%
\e{0}%
\e{0}%
\e{0}%
\e{1}%
\e{2}%
\e{1}%
\e{1}%
\e{0}%
\e{2}%
\e{2}%
\e{2}%
\e{2}%
\e{0}%
\e{0}%
\eol}\vss}\rg%
%
%
\rx{\vss\hfull{%
\rlx{\hss{$5670_{y}$}}\cg%
\e{8}%
\e{7}%
\e{4}%
\e{2}%
\e{0}%
\e{0}%
\e{0}%
\e{0}%
\e{0}%
\e{0}%
\e{0}%
\e{1}%
\e{3}%
\e{1}%
\e{1}%
\e{0}%
\e{2}%
\e{2}%
\e{3}%
\e{2}%
\e{1}%
\e{0}%
\eol}\vss}\rg%
%
%
\rx{\vss\hfull{%
\rlx{\hss{$4480_{y}$}}\cg%
\e{6}%
\e{6}%
\e{2}%
\e{2}%
\e{0}%
\e{0}%
\e{0}%
\e{0}%
\e{0}%
\e{0}%
\e{0}%
\e{0}%
\e{1}%
\e{1}%
\e{1}%
\e{0}%
\e{1}%
\e{1}%
\e{2}%
\e{2}%
\e{1}%
\e{1}%
\eol}\vss}\rg%
%
%
\rx{\vss\hfull{%
\rlx{\hss{$8_{z}$}}\cg%
\e{0}%
\e{0}%
\e{0}%
\e{0}%
\e{0}%
\e{0}%
\e{0}%
\e{0}%
\e{0}%
\e{0}%
\e{0}%
\e{0}%
\e{0}%
\e{0}%
\e{0}%
\e{0}%
\e{0}%
\e{0}%
\e{0}%
\e{0}%
\e{0}%
\e{0}%
\eol}\vss}\rg%
%
%
\rx{\vss\hfull{%
\rlx{\hss{$56_{z}$}}\cg%
\e{0}%
\e{0}%
\e{0}%
\e{0}%
\e{1}%
\e{0}%
\e{0}%
\e{0}%
\e{0}%
\e{0}%
\e{0}%
\e{0}%
\e{0}%
\e{0}%
\e{0}%
\e{0}%
\e{0}%
\e{0}%
\e{0}%
\e{0}%
\e{0}%
\e{0}%
\eol}\vss}\rg%
%
%
\rx{\vss\hfull{%
\rlx{\hss{$160_{z}$}}\cg%
\e{0}%
\e{0}%
\e{0}%
\e{0}%
\e{0}%
\e{0}%
\e{0}%
\e{0}%
\e{0}%
\e{0}%
\e{0}%
\e{0}%
\e{0}%
\e{0}%
\e{0}%
\e{0}%
\e{0}%
\e{0}%
\e{0}%
\e{0}%
\e{0}%
\e{0}%
\eol}\vss}\rg%
%
%
\rx{\vss\hfull{%
\rlx{\hss{$112_{z}$}}\cg%
\e{0}%
\e{0}%
\e{0}%
\e{0}%
\e{0}%
\e{0}%
\e{0}%
\e{0}%
\e{0}%
\e{0}%
\e{0}%
\e{0}%
\e{0}%
\e{0}%
\e{0}%
\e{0}%
\e{0}%
\e{0}%
\e{0}%
\e{0}%
\e{0}%
\e{0}%
\eol}\vss}\rg%
%
%
\rx{\vss\hfull{%
\rlx{\hss{$840_{z}$}}\cg%
\e{1}%
\e{0}%
\e{1}%
\e{0}%
\e{0}%
\e{0}%
\e{0}%
\e{0}%
\e{0}%
\e{0}%
\e{0}%
\e{0}%
\e{0}%
\e{0}%
\e{0}%
\e{0}%
\e{0}%
\e{1}%
\e{0}%
\e{0}%
\e{0}%
\e{0}%
\eol}\vss}\rg%
%
%
\rx{\vss\hfull{%
\rlx{\hss{$1296_{z}$}}\cg%
\e{1}%
\e{0}%
\e{0}%
\e{0}%
\e{1}%
\e{1}%
\e{0}%
\e{0}%
\e{0}%
\e{0}%
\e{0}%
\e{1}%
\e{0}%
\e{0}%
\e{0}%
\e{0}%
\e{1}%
\e{0}%
\e{0}%
\e{0}%
\e{0}%
\e{0}%
\eol}\vss}\rg%
%
%
\rx{\vss\hfull{%
\rlx{\hss{$1400_{z}$}}\cg%
\e{1}%
\e{0}%
\e{0}%
\e{0}%
\e{0}%
\e{0}%
\e{0}%
\e{0}%
\e{0}%
\e{0}%
\e{0}%
\e{1}%
\e{0}%
\e{0}%
\e{0}%
\e{0}%
\e{0}%
\e{0}%
\e{0}%
\e{0}%
\e{0}%
\e{0}%
\eol}\vss}\rg%
%
%
\rx{\vss\hfull{%
\rlx{\hss{$1008_{z}$}}\cg%
\e{0}%
\e{0}%
\e{0}%
\e{0}%
\e{1}%
\e{0}%
\e{0}%
\e{0}%
\e{0}%
\e{0}%
\e{0}%
\e{1}%
\e{0}%
\e{0}%
\e{0}%
\e{0}%
\e{0}%
\e{0}%
\e{0}%
\e{0}%
\e{0}%
\e{0}%
\eol}\vss}\rg%
%
%
\rx{\vss\hfull{%
\rlx{\hss{$560_{z}$}}\cg%
\e{0}%
\e{0}%
\e{0}%
\e{0}%
\e{0}%
\e{0}%
\e{0}%
\e{0}%
\e{0}%
\e{0}%
\e{0}%
\e{0}%
\e{0}%
\e{0}%
\e{0}%
\e{0}%
\e{0}%
\e{0}%
\e{0}%
\e{0}%
\e{0}%
\e{0}%
\eol}\vss}\rg%
%
%
\rx{\vss\hfull{%
\rlx{\hss{$1400_{zz}$}}\cg%
\e{2}%
\e{1}%
\e{0}%
\e{0}%
\e{0}%
\e{0}%
\e{0}%
\e{0}%
\e{0}%
\e{0}%
\e{0}%
\e{0}%
\e{0}%
\e{1}%
\e{0}%
\e{0}%
\e{0}%
\e{0}%
\e{0}%
\e{0}%
\e{0}%
\e{0}%
\eol}\vss}\rg%
%
%
\rx{\vss\hfull{%
\rlx{\hss{$4200_{z}$}}\cg%
\e{5}%
\e{5}%
\e{1}%
\e{1}%
\e{0}%
\e{0}%
\e{1}%
\e{0}%
\e{0}%
\e{0}%
\e{0}%
\e{0}%
\e{1}%
\e{2}%
\e{1}%
\e{0}%
\e{1}%
\e{0}%
\e{1}%
\e{1}%
\e{0}%
\e{0}%
\eol}\vss}\rg%
%
%
\rx{\vss\hfull{%
\rlx{\hss{$400_{z}$}}\cg%
\e{0}%
\e{0}%
\e{0}%
\e{0}%
\e{0}%
\e{0}%
\e{0}%
\e{0}%
\e{0}%
\e{0}%
\e{0}%
\e{0}%
\e{0}%
\e{0}%
\e{0}%
\e{0}%
\e{0}%
\e{0}%
\e{0}%
\e{0}%
\e{0}%
\e{0}%
\eol}\vss}\rg%
%
%
\rx{\vss\hfull{%
\rlx{\hss{$3240_{z}$}}\cg%
\e{3}%
\e{1}%
\e{0}%
\e{0}%
\e{0}%
\e{0}%
\e{0}%
\e{0}%
\e{0}%
\e{0}%
\e{0}%
\e{0}%
\e{1}%
\e{0}%
\e{0}%
\e{0}%
\e{1}%
\e{0}%
\e{0}%
\e{0}%
\e{0}%
\e{0}%
\eol}\vss}\rg%
%
%
\rx{\vss\hfull{%
\rlx{\hss{$4536_{z}$}}\cg%
\e{7}%
\e{4}%
\e{2}%
\e{1}%
\e{0}%
\e{0}%
\e{0}%
\e{0}%
\e{0}%
\e{0}%
\e{0}%
\e{0}%
\e{1}%
\e{0}%
\e{1}%
\e{0}%
\e{0}%
\e{1}%
\e{1}%
\e{1}%
\e{0}%
\e{0}%
\eol}\vss}\rg%
%
%
\rx{\vss\hfull{%
\rlx{\hss{$2400_{z}$}}\cg%
\e{3}%
\e{2}%
\e{1}%
\e{0}%
\e{1}%
\e{1}%
\e{1}%
\e{0}%
\e{0}%
\e{0}%
\e{0}%
\e{1}%
\e{1}%
\e{0}%
\e{0}%
\e{0}%
\e{2}%
\e{1}%
\e{1}%
\e{0}%
\e{0}%
\e{0}%
\eol}\vss}\rg%
%
%
\rx{\vss\hfull{%
\rlx{\hss{$3360_{z}$}}\cg%
\e{4}%
\e{3}%
\e{0}%
\e{0}%
\e{0}%
\e{0}%
\e{1}%
\e{0}%
\e{0}%
\e{0}%
\e{0}%
\e{0}%
\e{1}%
\e{1}%
\e{0}%
\e{0}%
\e{1}%
\e{0}%
\e{1}%
\e{0}%
\e{0}%
\e{0}%
\eol}\vss}\rg%
%
%
\rx{\vss\hfull{%
\rlx{\hss{$2800_{z}$}}\cg%
\e{3}%
\e{2}%
\e{0}%
\e{0}%
\e{1}%
\e{1}%
\e{0}%
\e{0}%
\e{0}%
\e{0}%
\e{0}%
\e{2}%
\e{1}%
\e{0}%
\e{0}%
\e{0}%
\e{1}%
\e{0}%
\e{1}%
\e{0}%
\e{0}%
\e{0}%
\eol}\vss}\rg%
%
%
\rx{\vss\hfull{%
\rlx{\hss{$4096_{z}$}}\cg%
\e{5}%
\e{2}%
\e{1}%
\e{0}%
\e{0}%
\e{1}%
\e{0}%
\e{0}%
\e{0}%
\e{0}%
\e{0}%
\e{1}%
\e{1}%
\e{0}%
\e{0}%
\e{0}%
\e{1}%
\e{1}%
\e{1}%
\e{0}%
\e{0}%
\e{0}%
\eol}\vss}\rg%
%
%
\rx{\vss\hfull{%
\rlx{\hss{$5600_{z}$}}\cg%
\e{8}%
\e{5}%
\e{3}%
\e{1}%
\e{0}%
\e{1}%
\e{0}%
\e{1}%
\e{0}%
\e{0}%
\e{0}%
\e{2}%
\e{2}%
\e{0}%
\e{0}%
\e{0}%
\e{2}%
\e{2}%
\e{2}%
\e{1}%
\e{1}%
\e{0}%
\eol}\vss}\rg%
%
%
\rx{\vss\hfull{%
\rlx{\hss{$448_{z}$}}\cg%
\e{1}%
\e{0}%
\e{0}%
\e{0}%
\e{0}%
\e{0}%
\e{0}%
\e{0}%
\e{0}%
\e{0}%
\e{0}%
\e{0}%
\e{0}%
\e{0}%
\e{0}%
\e{0}%
\e{0}%
\e{0}%
\e{0}%
\e{0}%
\e{0}%
\e{0}%
\eol}\vss}\rg%
\eop
\eject
\tablecont%
%
%
%
%
%
%
\rowpts=18 true pt%
\colpts=18 true pt%
\rowlabpts=40 true pt%
\collabpts=85 true pt%
\clx{\vss\hfull{%
\rlx{\hss{$ $}}\cg%
\cx{\hskip 16 true pt\flip{$[{2}{1}:{2}]{\times}[{2}{1^{2}}]$}\hss}\cg%
\cx{\hskip 16 true pt\flip{$[{2}{1}:{1^{2}}]{\times}[{2}{1^{2}}]$}\hss}\cg%
\cx{\hskip 16 true pt\flip{$[{1^{3}}:{2}]{\times}[{2}{1^{2}}]$}\hss}\cg%
\cx{\hskip 16 true pt\flip{$[{1^{3}}:{1^{2}}]{\times}[{2}{1^{2}}]$}\hss}\cg%
\cx{\hskip 16 true pt\flip{$[{5}:-]{\times}[{1^{4}}]$}\hss}\cg%
\cx{\hskip 16 true pt\flip{$[{4}{1}:-]{\times}[{1^{4}}]$}\hss}\cg%
\cx{\hskip 16 true pt\flip{$[{3}{2}:-]{\times}[{1^{4}}]$}\hss}\cg%
\cx{\hskip 16 true pt\flip{$[{3}{1^{2}}:-]{\times}[{1^{4}}]$}\hss}\cg%
\cx{\hskip 16 true pt\flip{$[{2^{2}}{1}:-]{\times}[{1^{4}}]$}\hss}\cg%
\cx{\hskip 16 true pt\flip{$[{2}{1^{3}}:-]{\times}[{1^{4}}]$}\hss}\cg%
\cx{\hskip 16 true pt\flip{$[{1^{5}}:-]{\times}[{1^{4}}]$}\hss}\cg%
\cx{\hskip 16 true pt\flip{$[{4}:{1}]{\times}[{1^{4}}]$}\hss}\cg%
\cx{\hskip 16 true pt\flip{$[{3}{1}:{1}]{\times}[{1^{4}}]$}\hss}\cg%
\cx{\hskip 16 true pt\flip{$[{2^{2}}:{1}]{\times}[{1^{4}}]$}\hss}\cg%
\cx{\hskip 16 true pt\flip{$[{2}{1^{2}}:{1}]{\times}[{1^{4}}]$}\hss}\cg%
\cx{\hskip 16 true pt\flip{$[{1^{4}}:{1}]{\times}[{1^{4}}]$}\hss}\cg%
\cx{\hskip 16 true pt\flip{$[{3}:{2}]{\times}[{1^{4}}]$}\hss}\cg%
\cx{\hskip 16 true pt\flip{$[{3}:{1^{2}}]{\times}[{1^{4}}]$}\hss}\cg%
\cx{\hskip 16 true pt\flip{$[{2}{1}:{2}]{\times}[{1^{4}}]$}\hss}\cg%
\cx{\hskip 16 true pt\flip{$[{2}{1}:{1^{2}}]{\times}[{1^{4}}]$}\hss}\cg%
\cx{\hskip 16 true pt\flip{$[{1^{3}}:{2}]{\times}[{1^{4}}]$}\hss}\cg%
\cx{\hskip 16 true pt\flip{$[{1^{3}}:{1^{2}}]{\times}[{1^{4}}]$}\hss}\cg%
\eol}}\rg%
%
%
\rx{\vss\hfull{%
\rlx{\hss{$448_{w}$}}\cg%
\e{1}%
\e{0}%
\e{1}%
\e{0}%
\e{0}%
\e{1}%
\e{0}%
\e{0}%
\e{0}%
\e{0}%
\e{0}%
\e{0}%
\e{0}%
\e{0}%
\e{0}%
\e{0}%
\e{1}%
\e{1}%
\e{0}%
\e{0}%
\e{0}%
\e{0}%
\eol}\vss}\rg%
%
%
\rx{\vss\hfull{%
\rlx{\hss{$1344_{w}$}}\cg%
\e{2}%
\e{2}%
\e{1}%
\e{1}%
\e{0}%
\e{0}%
\e{0}%
\e{0}%
\e{0}%
\e{0}%
\e{0}%
\e{0}%
\e{0}%
\e{0}%
\e{1}%
\e{1}%
\e{0}%
\e{0}%
\e{0}%
\e{1}%
\e{0}%
\e{0}%
\eol}\vss}\rg%
%
%
\rx{\vss\hfull{%
\rlx{\hss{$5600_{w}$}}\cg%
\e{8}%
\e{7}%
\e{4}%
\e{2}%
\e{0}%
\e{1}%
\e{1}%
\e{1}%
\e{0}%
\e{0}%
\e{0}%
\e{0}%
\e{2}%
\e{1}%
\e{1}%
\e{0}%
\e{2}%
\e{2}%
\e{3}%
\e{2}%
\e{1}%
\e{0}%
\eol}\vss}\rg%
%
%
\rx{\vss\hfull{%
\rlx{\hss{$2016_{w}$}}\cg%
\e{3}%
\e{3}%
\e{1}%
\e{2}%
\e{0}%
\e{0}%
\e{0}%
\e{0}%
\e{1}%
\e{0}%
\e{0}%
\e{0}%
\e{0}%
\e{1}%
\e{1}%
\e{0}%
\e{0}%
\e{0}%
\e{0}%
\e{1}%
\e{0}%
\e{1}%
\eol}\vss}\rg%
%
%
\rx{\vss\hfull{%
\rlx{\hss{$7168_{w}$}}\cg%
\e{10}%
\e{10}%
\e{4}%
\e{4}%
\e{0}%
\e{0}%
\e{1}%
\e{1}%
\e{1}%
\e{0}%
\e{0}%
\e{0}%
\e{2}%
\e{2}%
\e{2}%
\e{0}%
\e{1}%
\e{1}%
\e{3}%
\e{3}%
\e{1}%
\e{1}%
\eol}\vss}\rg%
\tableclose%
%
%
%
%
%
%
\tableopen{Induce/restrict matrix for $W({E_{6}}{A_{2}})\,\subset\,W(E_{8})$}%
%
%
%
%
%
%
\rowpts=18 true pt%
\colpts=18 true pt%
\rowlabpts=40 true pt%
\collabpts=60 true pt%
\clx{\vss\hfull{%
\rlx{\hss{$ $}}\cg%
\cx{\hskip 16 true pt\flip{$1_p{\times}[{3}]$}\hss}\cg%
\cx{\hskip 16 true pt\flip{$6_p{\times}[{3}]$}\hss}\cg%
\cx{\hskip 16 true pt\flip{$15_p{\times}[{3}]$}\hss}\cg%
\cx{\hskip 16 true pt\flip{$20_p{\times}[{3}]$}\hss}\cg%
\cx{\hskip 16 true pt\flip{$30_p{\times}[{3}]$}\hss}\cg%
\cx{\hskip 16 true pt\flip{$64_p{\times}[{3}]$}\hss}\cg%
\cx{\hskip 16 true pt\flip{$81_p{\times}[{3}]$}\hss}\cg%
\cx{\hskip 16 true pt\flip{$15_q{\times}[{3}]$}\hss}\cg%
\cx{\hskip 16 true pt\flip{$24_p{\times}[{3}]$}\hss}\cg%
\cx{\hskip 16 true pt\flip{$60_p{\times}[{3}]$}\hss}\cg%
\cx{\hskip 16 true pt\flip{$1_p^{*}{\times}[{3}]$}\hss}\cg%
\cx{\hskip 16 true pt\flip{$6_p^{*}{\times}[{3}]$}\hss}\cg%
\cx{\hskip 16 true pt\flip{$15_p^{*}{\times}[{3}]$}\hss}\cg%
\cx{\hskip 16 true pt\flip{$20_p^{*}{\times}[{3}]$}\hss}\cg%
\cx{\hskip 16 true pt\flip{$30_p^{*}{\times}[{3}]$}\hss}\cg%
\cx{\hskip 16 true pt\flip{$64_p^{*}{\times}[{3}]$}\hss}\cg%
\cx{\hskip 16 true pt\flip{$81_p^{*}{\times}[{3}]$}\hss}\cg%
\cx{\hskip 16 true pt\flip{$15_q^{*}{\times}[{3}]$}\hss}\cg%
\cx{\hskip 16 true pt\flip{$24_p^{*}{\times}[{3}]$}\hss}\cg%
\cx{\hskip 16 true pt\flip{$60_p^{*}{\times}[{3}]$}\hss}\cg%
\cx{\hskip 16 true pt\flip{$20_s{\times}[{3}]$}\hss}\cg%
\cx{\hskip 16 true pt\flip{$90_s{\times}[{3}]$}\hss}\cg%
\cx{\hskip 16 true pt\flip{$80_s{\times}[{3}]$}\hss}\cg%
\cx{\hskip 16 true pt\flip{$60_s{\times}[{3}]$}\hss}\cg%
\cx{\hskip 16 true pt\flip{$10_s{\times}[{3}]$}\hss}\cg%
\eol}}\rg%
%
%
\rx{\vss\hfull{%
\rlx{\hss{$1_{x}$}}\cg%
\e{1}%
\e{0}%
\e{0}%
\e{0}%
\e{0}%
\e{0}%
\e{0}%
\e{0}%
\e{0}%
\e{0}%
\e{0}%
\e{0}%
\e{0}%
\e{0}%
\e{0}%
\e{0}%
\e{0}%
\e{0}%
\e{0}%
\e{0}%
\e{0}%
\e{0}%
\e{0}%
\e{0}%
\e{0}%
\eol}\vss}\rg%
%
%
\rx{\vss\hfull{%
\rlx{\hss{$28_{x}$}}\cg%
\e{0}%
\e{0}%
\e{1}%
\e{0}%
\e{0}%
\e{0}%
\e{0}%
\e{0}%
\e{0}%
\e{0}%
\e{0}%
\e{0}%
\e{0}%
\e{0}%
\e{0}%
\e{0}%
\e{0}%
\e{0}%
\e{0}%
\e{0}%
\e{0}%
\e{0}%
\e{0}%
\e{0}%
\e{0}%
\eol}\vss}\rg%
%
%
\rx{\vss\hfull{%
\rlx{\hss{$35_{x}$}}\cg%
\e{1}%
\e{0}%
\e{0}%
\e{1}%
\e{0}%
\e{0}%
\e{0}%
\e{0}%
\e{0}%
\e{0}%
\e{0}%
\e{0}%
\e{0}%
\e{0}%
\e{0}%
\e{0}%
\e{0}%
\e{0}%
\e{0}%
\e{0}%
\e{0}%
\e{0}%
\e{0}%
\e{0}%
\e{0}%
\eol}\vss}\rg%
%
%
\rx{\vss\hfull{%
\rlx{\hss{$84_{x}$}}\cg%
\e{1}%
\e{1}%
\e{0}%
\e{1}%
\e{0}%
\e{0}%
\e{0}%
\e{1}%
\e{0}%
\e{0}%
\e{0}%
\e{0}%
\e{0}%
\e{0}%
\e{0}%
\e{0}%
\e{0}%
\e{0}%
\e{0}%
\e{0}%
\e{0}%
\e{0}%
\e{0}%
\e{0}%
\e{0}%
\eol}\vss}\rg%
%
%
\rx{\vss\hfull{%
\rlx{\hss{$50_{x}$}}\cg%
\e{0}%
\e{0}%
\e{0}%
\e{1}%
\e{0}%
\e{0}%
\e{0}%
\e{0}%
\e{0}%
\e{0}%
\e{0}%
\e{0}%
\e{0}%
\e{0}%
\e{0}%
\e{0}%
\e{0}%
\e{0}%
\e{0}%
\e{0}%
\e{0}%
\e{0}%
\e{0}%
\e{0}%
\e{0}%
\eol}\vss}\rg%
%
%
\rx{\vss\hfull{%
\rlx{\hss{$350_{x}$}}\cg%
\e{0}%
\e{0}%
\e{1}%
\e{0}%
\e{0}%
\e{0}%
\e{0}%
\e{0}%
\e{0}%
\e{0}%
\e{0}%
\e{0}%
\e{0}%
\e{0}%
\e{0}%
\e{0}%
\e{0}%
\e{0}%
\e{0}%
\e{0}%
\e{0}%
\e{1}%
\e{0}%
\e{0}%
\e{0}%
\eol}\vss}\rg%
%
%
\rx{\vss\hfull{%
\rlx{\hss{$300_{x}$}}\cg%
\e{1}%
\e{0}%
\e{0}%
\e{1}%
\e{0}%
\e{0}%
\e{0}%
\e{0}%
\e{1}%
\e{1}%
\e{0}%
\e{0}%
\e{0}%
\e{0}%
\e{0}%
\e{0}%
\e{0}%
\e{0}%
\e{0}%
\e{0}%
\e{0}%
\e{0}%
\e{0}%
\e{0}%
\e{0}%
\eol}\vss}\rg%
%
%
\rx{\vss\hfull{%
\rlx{\hss{$567_{x}$}}\cg%
\e{0}%
\e{1}%
\e{1}%
\e{1}%
\e{1}%
\e{1}%
\e{1}%
\e{0}%
\e{0}%
\e{0}%
\e{0}%
\e{0}%
\e{0}%
\e{0}%
\e{0}%
\e{0}%
\e{0}%
\e{0}%
\e{0}%
\e{0}%
\e{0}%
\e{0}%
\e{0}%
\e{0}%
\e{0}%
\eol}\vss}\rg%
%
%
\rx{\vss\hfull{%
\rlx{\hss{$210_{x}$}}\cg%
\e{0}%
\e{1}%
\e{0}%
\e{1}%
\e{0}%
\e{1}%
\e{0}%
\e{0}%
\e{0}%
\e{0}%
\e{0}%
\e{0}%
\e{0}%
\e{0}%
\e{0}%
\e{0}%
\e{0}%
\e{0}%
\e{0}%
\e{0}%
\e{0}%
\e{0}%
\e{0}%
\e{0}%
\e{0}%
\eol}\vss}\rg%
%
%
\rx{\vss\hfull{%
\rlx{\hss{$840_{x}$}}\cg%
\e{0}%
\e{0}%
\e{0}%
\e{1}%
\e{0}%
\e{0}%
\e{0}%
\e{0}%
\e{1}%
\e{1}%
\e{0}%
\e{0}%
\e{0}%
\e{0}%
\e{0}%
\e{0}%
\e{0}%
\e{0}%
\e{0}%
\e{0}%
\e{0}%
\e{0}%
\e{0}%
\e{1}%
\e{1}%
\eol}\vss}\rg%
%
%
\rx{\vss\hfull{%
\rlx{\hss{$700_{x}$}}\cg%
\e{1}%
\e{1}%
\e{0}%
\e{2}%
\e{1}%
\e{1}%
\e{0}%
\e{1}%
\e{1}%
\e{1}%
\e{0}%
\e{0}%
\e{0}%
\e{0}%
\e{0}%
\e{0}%
\e{0}%
\e{0}%
\e{0}%
\e{0}%
\e{0}%
\e{0}%
\e{0}%
\e{0}%
\e{0}%
\eol}\vss}\rg%
%
%
\rx{\vss\hfull{%
\rlx{\hss{$175_{x}$}}\cg%
\e{0}%
\e{0}%
\e{0}%
\e{0}%
\e{1}%
\e{0}%
\e{0}%
\e{1}%
\e{0}%
\e{0}%
\e{0}%
\e{0}%
\e{0}%
\e{0}%
\e{0}%
\e{0}%
\e{0}%
\e{0}%
\e{0}%
\e{0}%
\e{0}%
\e{0}%
\e{0}%
\e{0}%
\e{0}%
\eol}\vss}\rg%
%
%
\rx{\vss\hfull{%
\rlx{\hss{$1400_{x}$}}\cg%
\e{0}%
\e{0}%
\e{1}%
\e{1}%
\e{2}%
\e{1}%
\e{1}%
\e{0}%
\e{0}%
\e{1}%
\e{0}%
\e{0}%
\e{0}%
\e{0}%
\e{0}%
\e{0}%
\e{0}%
\e{0}%
\e{0}%
\e{0}%
\e{0}%
\e{1}%
\e{0}%
\e{0}%
\e{0}%
\eol}\vss}\rg%
%
%
\rx{\vss\hfull{%
\rlx{\hss{$1050_{x}$}}\cg%
\e{0}%
\e{0}%
\e{1}%
\e{1}%
\e{1}%
\e{1}%
\e{1}%
\e{1}%
\e{0}%
\e{1}%
\e{0}%
\e{0}%
\e{0}%
\e{0}%
\e{0}%
\e{0}%
\e{0}%
\e{0}%
\e{0}%
\e{0}%
\e{0}%
\e{0}%
\e{0}%
\e{0}%
\e{0}%
\eol}\vss}\rg%
%
%
\rx{\vss\hfull{%
\rlx{\hss{$1575_{x}$}}\cg%
\e{0}%
\e{1}%
\e{1}%
\e{0}%
\e{1}%
\e{2}%
\e{1}%
\e{0}%
\e{0}%
\e{0}%
\e{0}%
\e{0}%
\e{0}%
\e{0}%
\e{0}%
\e{0}%
\e{0}%
\e{0}%
\e{0}%
\e{0}%
\e{1}%
\e{1}%
\e{1}%
\e{0}%
\e{0}%
\eol}\vss}\rg%
%
%
\rx{\vss\hfull{%
\rlx{\hss{$1344_{x}$}}\cg%
\e{0}%
\e{1}%
\e{0}%
\e{2}%
\e{1}%
\e{2}%
\e{1}%
\e{1}%
\e{0}%
\e{1}%
\e{0}%
\e{0}%
\e{0}%
\e{0}%
\e{0}%
\e{0}%
\e{0}%
\e{0}%
\e{0}%
\e{0}%
\e{0}%
\e{0}%
\e{0}%
\e{1}%
\e{0}%
\eol}\vss}\rg%
%
%
\rx{\vss\hfull{%
\rlx{\hss{$2100_{x}$}}\cg%
\e{0}%
\e{0}%
\e{0}%
\e{0}%
\e{0}%
\e{1}%
\e{1}%
\e{0}%
\e{0}%
\e{0}%
\e{0}%
\e{0}%
\e{0}%
\e{0}%
\e{0}%
\e{1}%
\e{1}%
\e{0}%
\e{0}%
\e{0}%
\e{1}%
\e{1}%
\e{1}%
\e{0}%
\e{0}%
\eol}\vss}\rg%
%
%
\rx{\vss\hfull{%
\rlx{\hss{$2268_{x}$}}\cg%
\e{0}%
\e{0}%
\e{0}%
\e{1}%
\e{1}%
\e{2}%
\e{1}%
\e{0}%
\e{1}%
\e{1}%
\e{0}%
\e{0}%
\e{0}%
\e{0}%
\e{0}%
\e{0}%
\e{1}%
\e{0}%
\e{0}%
\e{0}%
\e{0}%
\e{1}%
\e{1}%
\e{0}%
\e{0}%
\eol}\vss}\rg%
%
%
\rx{\vss\hfull{%
\rlx{\hss{$525_{x}$}}\cg%
\e{0}%
\e{0}%
\e{1}%
\e{0}%
\e{1}%
\e{0}%
\e{1}%
\e{0}%
\e{0}%
\e{0}%
\e{0}%
\e{0}%
\e{0}%
\e{0}%
\e{0}%
\e{0}%
\e{0}%
\e{0}%
\e{1}%
\e{0}%
\e{0}%
\e{0}%
\e{0}%
\e{0}%
\e{0}%
\eol}\vss}\rg%
%
%
\rx{\vss\hfull{%
\rlx{\hss{$700_{xx}$}}\cg%
\e{0}%
\e{0}%
\e{0}%
\e{0}%
\e{0}%
\e{1}%
\e{1}%
\e{0}%
\e{0}%
\e{0}%
\e{0}%
\e{0}%
\e{0}%
\e{0}%
\e{0}%
\e{0}%
\e{0}%
\e{0}%
\e{0}%
\e{0}%
\e{0}%
\e{0}%
\e{0}%
\e{0}%
\e{0}%
\eol}\vss}\rg%
%
%
\rx{\vss\hfull{%
\rlx{\hss{$972_{x}$}}\cg%
\e{0}%
\e{0}%
\e{0}%
\e{1}%
\e{0}%
\e{1}%
\e{0}%
\e{1}%
\e{1}%
\e{1}%
\e{0}%
\e{0}%
\e{0}%
\e{0}%
\e{0}%
\e{0}%
\e{0}%
\e{0}%
\e{0}%
\e{0}%
\e{0}%
\e{0}%
\e{0}%
\e{1}%
\e{0}%
\eol}\vss}\rg%
\eop
\eject
\tablecont%
%
%
%
%
%
%
\rowpts=18 true pt%
\colpts=18 true pt%
\rowlabpts=40 true pt%
\collabpts=60 true pt%
\clx{\vss\hfull{%
\rlx{\hss{$ $}}\cg%
\cx{\hskip 16 true pt\flip{$1_p{\times}[{3}]$}\hss}\cg%
\cx{\hskip 16 true pt\flip{$6_p{\times}[{3}]$}\hss}\cg%
\cx{\hskip 16 true pt\flip{$15_p{\times}[{3}]$}\hss}\cg%
\cx{\hskip 16 true pt\flip{$20_p{\times}[{3}]$}\hss}\cg%
\cx{\hskip 16 true pt\flip{$30_p{\times}[{3}]$}\hss}\cg%
\cx{\hskip 16 true pt\flip{$64_p{\times}[{3}]$}\hss}\cg%
\cx{\hskip 16 true pt\flip{$81_p{\times}[{3}]$}\hss}\cg%
\cx{\hskip 16 true pt\flip{$15_q{\times}[{3}]$}\hss}\cg%
\cx{\hskip 16 true pt\flip{$24_p{\times}[{3}]$}\hss}\cg%
\cx{\hskip 16 true pt\flip{$60_p{\times}[{3}]$}\hss}\cg%
\cx{\hskip 16 true pt\flip{$1_p^{*}{\times}[{3}]$}\hss}\cg%
\cx{\hskip 16 true pt\flip{$6_p^{*}{\times}[{3}]$}\hss}\cg%
\cx{\hskip 16 true pt\flip{$15_p^{*}{\times}[{3}]$}\hss}\cg%
\cx{\hskip 16 true pt\flip{$20_p^{*}{\times}[{3}]$}\hss}\cg%
\cx{\hskip 16 true pt\flip{$30_p^{*}{\times}[{3}]$}\hss}\cg%
\cx{\hskip 16 true pt\flip{$64_p^{*}{\times}[{3}]$}\hss}\cg%
\cx{\hskip 16 true pt\flip{$81_p^{*}{\times}[{3}]$}\hss}\cg%
\cx{\hskip 16 true pt\flip{$15_q^{*}{\times}[{3}]$}\hss}\cg%
\cx{\hskip 16 true pt\flip{$24_p^{*}{\times}[{3}]$}\hss}\cg%
\cx{\hskip 16 true pt\flip{$60_p^{*}{\times}[{3}]$}\hss}\cg%
\cx{\hskip 16 true pt\flip{$20_s{\times}[{3}]$}\hss}\cg%
\cx{\hskip 16 true pt\flip{$90_s{\times}[{3}]$}\hss}\cg%
\cx{\hskip 16 true pt\flip{$80_s{\times}[{3}]$}\hss}\cg%
\cx{\hskip 16 true pt\flip{$60_s{\times}[{3}]$}\hss}\cg%
\cx{\hskip 16 true pt\flip{$10_s{\times}[{3}]$}\hss}\cg%
\eol}}\rg%
%
%
\rx{\vss\hfull{%
\rlx{\hss{$4096_{x}$}}\cg%
\e{0}%
\e{0}%
\e{1}%
\e{1}%
\e{1}%
\e{2}%
\e{2}%
\e{0}%
\e{1}%
\e{2}%
\e{0}%
\e{0}%
\e{0}%
\e{0}%
\e{0}%
\e{0}%
\e{1}%
\e{0}%
\e{0}%
\e{1}%
\e{0}%
\e{2}%
\e{1}%
\e{1}%
\e{0}%
\eol}\vss}\rg%
%
%
\rx{\vss\hfull{%
\rlx{\hss{$4200_{x}$}}\cg%
\e{0}%
\e{0}%
\e{0}%
\e{1}%
\e{1}%
\e{2}%
\e{2}%
\e{1}%
\e{1}%
\e{2}%
\e{0}%
\e{0}%
\e{0}%
\e{0}%
\e{0}%
\e{0}%
\e{1}%
\e{0}%
\e{0}%
\e{0}%
\e{0}%
\e{1}%
\e{2}%
\e{1}%
\e{1}%
\eol}\vss}\rg%
%
%
\rx{\vss\hfull{%
\rlx{\hss{$2240_{x}$}}\cg%
\e{0}%
\e{1}%
\e{0}%
\e{1}%
\e{1}%
\e{2}%
\e{1}%
\e{1}%
\e{0}%
\e{2}%
\e{0}%
\e{0}%
\e{0}%
\e{0}%
\e{0}%
\e{0}%
\e{0}%
\e{0}%
\e{0}%
\e{0}%
\e{0}%
\e{0}%
\e{1}%
\e{1}%
\e{0}%
\eol}\vss}\rg%
%
%
\rx{\vss\hfull{%
\rlx{\hss{$2835_{x}$}}\cg%
\e{0}%
\e{0}%
\e{0}%
\e{0}%
\e{1}%
\e{1}%
\e{1}%
\e{1}%
\e{0}%
\e{2}%
\e{0}%
\e{0}%
\e{0}%
\e{0}%
\e{0}%
\e{0}%
\e{0}%
\e{0}%
\e{0}%
\e{0}%
\e{0}%
\e{1}%
\e{1}%
\e{1}%
\e{0}%
\eol}\vss}\rg%
%
%
\rx{\vss\hfull{%
\rlx{\hss{$6075_{x}$}}\cg%
\e{0}%
\e{0}%
\e{1}%
\e{0}%
\e{1}%
\e{2}%
\e{3}%
\e{0}%
\e{0}%
\e{1}%
\e{0}%
\e{0}%
\e{0}%
\e{0}%
\e{0}%
\e{1}%
\e{1}%
\e{0}%
\e{1}%
\e{1}%
\e{1}%
\e{3}%
\e{2}%
\e{1}%
\e{0}%
\eol}\vss}\rg%
%
%
\rx{\vss\hfull{%
\rlx{\hss{$3200_{x}$}}\cg%
\e{0}%
\e{0}%
\e{0}%
\e{0}%
\e{0}%
\e{1}%
\e{1}%
\e{1}%
\e{1}%
\e{1}%
\e{0}%
\e{0}%
\e{0}%
\e{0}%
\e{0}%
\e{0}%
\e{1}%
\e{1}%
\e{0}%
\e{1}%
\e{0}%
\e{0}%
\e{1}%
\e{2}%
\e{0}%
\eol}\vss}\rg%
%
%
\rx{\vss\hfull{%
\rlx{\hss{$70_{y}$}}\cg%
\e{0}%
\e{0}%
\e{0}%
\e{0}%
\e{0}%
\e{0}%
\e{0}%
\e{0}%
\e{0}%
\e{0}%
\e{0}%
\e{0}%
\e{1}%
\e{0}%
\e{0}%
\e{0}%
\e{0}%
\e{0}%
\e{0}%
\e{0}%
\e{0}%
\e{0}%
\e{0}%
\e{0}%
\e{0}%
\eol}\vss}\rg%
%
%
\rx{\vss\hfull{%
\rlx{\hss{$1134_{y}$}}\cg%
\e{0}%
\e{0}%
\e{1}%
\e{0}%
\e{0}%
\e{0}%
\e{0}%
\e{0}%
\e{0}%
\e{0}%
\e{0}%
\e{0}%
\e{0}%
\e{0}%
\e{0}%
\e{0}%
\e{0}%
\e{0}%
\e{1}%
\e{1}%
\e{0}%
\e{1}%
\e{0}%
\e{0}%
\e{0}%
\eol}\vss}\rg%
%
%
\rx{\vss\hfull{%
\rlx{\hss{$1680_{y}$}}\cg%
\e{0}%
\e{0}%
\e{0}%
\e{0}%
\e{0}%
\e{0}%
\e{0}%
\e{0}%
\e{0}%
\e{0}%
\e{0}%
\e{0}%
\e{1}%
\e{0}%
\e{1}%
\e{1}%
\e{1}%
\e{0}%
\e{0}%
\e{0}%
\e{1}%
\e{1}%
\e{0}%
\e{0}%
\e{0}%
\eol}\vss}\rg%
%
%
\rx{\vss\hfull{%
\rlx{\hss{$168_{y}$}}\cg%
\e{0}%
\e{0}%
\e{0}%
\e{0}%
\e{0}%
\e{0}%
\e{0}%
\e{0}%
\e{1}%
\e{0}%
\e{0}%
\e{0}%
\e{0}%
\e{0}%
\e{0}%
\e{0}%
\e{0}%
\e{0}%
\e{0}%
\e{0}%
\e{0}%
\e{0}%
\e{0}%
\e{0}%
\e{0}%
\eol}\vss}\rg%
%
%
\rx{\vss\hfull{%
\rlx{\hss{$420_{y}$}}\cg%
\e{0}%
\e{0}%
\e{0}%
\e{0}%
\e{0}%
\e{0}%
\e{0}%
\e{0}%
\e{0}%
\e{1}%
\e{0}%
\e{0}%
\e{0}%
\e{0}%
\e{0}%
\e{0}%
\e{0}%
\e{0}%
\e{0}%
\e{0}%
\e{0}%
\e{0}%
\e{0}%
\e{0}%
\e{0}%
\eol}\vss}\rg%
%
%
\rx{\vss\hfull{%
\rlx{\hss{$3150_{y}$}}\cg%
\e{0}%
\e{0}%
\e{0}%
\e{0}%
\e{1}%
\e{0}%
\e{1}%
\e{0}%
\e{0}%
\e{1}%
\e{0}%
\e{0}%
\e{0}%
\e{0}%
\e{0}%
\e{0}%
\e{0}%
\e{0}%
\e{1}%
\e{1}%
\e{0}%
\e{1}%
\e{1}%
\e{1}%
\e{1}%
\eol}\vss}\rg%
%
%
\rx{\vss\hfull{%
\rlx{\hss{$4200_{y}$}}\cg%
\e{0}%
\e{0}%
\e{0}%
\e{0}%
\e{0}%
\e{1}%
\e{1}%
\e{0}%
\e{1}%
\e{1}%
\e{0}%
\e{0}%
\e{0}%
\e{0}%
\e{0}%
\e{0}%
\e{1}%
\e{0}%
\e{0}%
\e{1}%
\e{0}%
\e{1}%
\e{1}%
\e{2}%
\e{0}%
\eol}\vss}\rg%
%
%
\rx{\vss\hfull{%
\rlx{\hss{$2688_{y}$}}\cg%
\e{0}%
\e{0}%
\e{0}%
\e{0}%
\e{0}%
\e{0}%
\e{1}%
\e{0}%
\e{0}%
\e{1}%
\e{0}%
\e{0}%
\e{0}%
\e{0}%
\e{0}%
\e{0}%
\e{1}%
\e{0}%
\e{0}%
\e{1}%
\e{0}%
\e{1}%
\e{0}%
\e{1}%
\e{0}%
\eol}\vss}\rg%
%
%
\rx{\vss\hfull{%
\rlx{\hss{$2100_{y}$}}\cg%
\e{0}%
\e{0}%
\e{0}%
\e{0}%
\e{0}%
\e{0}%
\e{1}%
\e{0}%
\e{0}%
\e{0}%
\e{0}%
\e{0}%
\e{0}%
\e{1}%
\e{0}%
\e{1}%
\e{1}%
\e{0}%
\e{1}%
\e{0}%
\e{0}%
\e{1}%
\e{0}%
\e{0}%
\e{0}%
\eol}\vss}\rg%
%
%
\rx{\vss\hfull{%
\rlx{\hss{$1400_{y}$}}\cg%
\e{0}%
\e{0}%
\e{0}%
\e{0}%
\e{0}%
\e{0}%
\e{0}%
\e{0}%
\e{1}%
\e{0}%
\e{0}%
\e{0}%
\e{1}%
\e{0}%
\e{1}%
\e{0}%
\e{1}%
\e{0}%
\e{0}%
\e{0}%
\e{0}%
\e{1}%
\e{0}%
\e{0}%
\e{0}%
\eol}\vss}\rg%
%
%
\rx{\vss\hfull{%
\rlx{\hss{$4536_{y}$}}\cg%
\e{0}%
\e{0}%
\e{0}%
\e{0}%
\e{0}%
\e{0}%
\e{1}%
\e{0}%
\e{1}%
\e{1}%
\e{0}%
\e{0}%
\e{1}%
\e{0}%
\e{1}%
\e{1}%
\e{2}%
\e{0}%
\e{0}%
\e{1}%
\e{0}%
\e{2}%
\e{1}%
\e{0}%
\e{0}%
\eol}\vss}\rg%
%
%
\rx{\vss\hfull{%
\rlx{\hss{$5670_{y}$}}\cg%
\e{0}%
\e{0}%
\e{0}%
\e{0}%
\e{0}%
\e{1}%
\e{1}%
\e{0}%
\e{0}%
\e{0}%
\e{0}%
\e{0}%
\e{0}%
\e{0}%
\e{1}%
\e{2}%
\e{2}%
\e{0}%
\e{0}%
\e{1}%
\e{1}%
\e{2}%
\e{2}%
\e{1}%
\e{0}%
\eol}\vss}\rg%
%
%
\rx{\vss\hfull{%
\rlx{\hss{$4480_{y}$}}\cg%
\e{0}%
\e{0}%
\e{0}%
\e{0}%
\e{0}%
\e{1}%
\e{1}%
\e{0}%
\e{0}%
\e{1}%
\e{0}%
\e{0}%
\e{0}%
\e{0}%
\e{0}%
\e{1}%
\e{1}%
\e{0}%
\e{0}%
\e{1}%
\e{0}%
\e{1}%
\e{2}%
\e{1}%
\e{0}%
\eol}\vss}\rg%
%
%
\rx{\vss\hfull{%
\rlx{\hss{$8_{z}$}}\cg%
\e{0}%
\e{1}%
\e{0}%
\e{0}%
\e{0}%
\e{0}%
\e{0}%
\e{0}%
\e{0}%
\e{0}%
\e{0}%
\e{0}%
\e{0}%
\e{0}%
\e{0}%
\e{0}%
\e{0}%
\e{0}%
\e{0}%
\e{0}%
\e{0}%
\e{0}%
\e{0}%
\e{0}%
\e{0}%
\eol}\vss}\rg%
%
%
\rx{\vss\hfull{%
\rlx{\hss{$56_{z}$}}\cg%
\e{0}%
\e{0}%
\e{0}%
\e{0}%
\e{0}%
\e{0}%
\e{0}%
\e{0}%
\e{0}%
\e{0}%
\e{0}%
\e{0}%
\e{0}%
\e{0}%
\e{0}%
\e{0}%
\e{0}%
\e{0}%
\e{0}%
\e{0}%
\e{1}%
\e{0}%
\e{0}%
\e{0}%
\e{0}%
\eol}\vss}\rg%
%
%
\rx{\vss\hfull{%
\rlx{\hss{$160_{z}$}}\cg%
\e{0}%
\e{1}%
\e{0}%
\e{0}%
\e{0}%
\e{1}%
\e{0}%
\e{0}%
\e{0}%
\e{0}%
\e{0}%
\e{0}%
\e{0}%
\e{0}%
\e{0}%
\e{0}%
\e{0}%
\e{0}%
\e{0}%
\e{0}%
\e{0}%
\e{0}%
\e{0}%
\e{0}%
\e{0}%
\eol}\vss}\rg%
%
%
\rx{\vss\hfull{%
\rlx{\hss{$112_{z}$}}\cg%
\e{1}%
\e{1}%
\e{0}%
\e{1}%
\e{1}%
\e{0}%
\e{0}%
\e{0}%
\e{0}%
\e{0}%
\e{0}%
\e{0}%
\e{0}%
\e{0}%
\e{0}%
\e{0}%
\e{0}%
\e{0}%
\e{0}%
\e{0}%
\e{0}%
\e{0}%
\e{0}%
\e{0}%
\e{0}%
\eol}\vss}\rg%
%
%
\rx{\vss\hfull{%
\rlx{\hss{$840_{z}$}}\cg%
\e{0}%
\e{1}%
\e{0}%
\e{0}%
\e{0}%
\e{1}%
\e{0}%
\e{0}%
\e{0}%
\e{0}%
\e{0}%
\e{0}%
\e{0}%
\e{0}%
\e{0}%
\e{0}%
\e{0}%
\e{0}%
\e{0}%
\e{0}%
\e{0}%
\e{0}%
\e{1}%
\e{1}%
\e{0}%
\eol}\vss}\rg%
%
%
\rx{\vss\hfull{%
\rlx{\hss{$1296_{z}$}}\cg%
\e{0}%
\e{0}%
\e{1}%
\e{0}%
\e{0}%
\e{1}%
\e{1}%
\e{0}%
\e{0}%
\e{0}%
\e{0}%
\e{0}%
\e{0}%
\e{0}%
\e{0}%
\e{0}%
\e{1}%
\e{0}%
\e{0}%
\e{0}%
\e{1}%
\e{1}%
\e{0}%
\e{0}%
\e{0}%
\eol}\vss}\rg%
%
%
\rx{\vss\hfull{%
\rlx{\hss{$1400_{z}$}}\cg%
\e{0}%
\e{1}%
\e{0}%
\e{2}%
\e{1}%
\e{2}%
\e{1}%
\e{1}%
\e{0}%
\e{1}%
\e{0}%
\e{0}%
\e{0}%
\e{0}%
\e{0}%
\e{0}%
\e{0}%
\e{0}%
\e{0}%
\e{0}%
\e{0}%
\e{0}%
\e{1}%
\e{0}%
\e{0}%
\eol}\vss}\rg%
%
%
\rx{\vss\hfull{%
\rlx{\hss{$1008_{z}$}}\cg%
\e{0}%
\e{0}%
\e{1}%
\e{1}%
\e{1}%
\e{1}%
\e{1}%
\e{0}%
\e{1}%
\e{0}%
\e{0}%
\e{0}%
\e{0}%
\e{0}%
\e{0}%
\e{0}%
\e{0}%
\e{0}%
\e{0}%
\e{0}%
\e{0}%
\e{1}%
\e{0}%
\e{0}%
\e{0}%
\eol}\vss}\rg%
%
%
\rx{\vss\hfull{%
\rlx{\hss{$560_{z}$}}\cg%
\e{1}%
\e{1}%
\e{1}%
\e{2}%
\e{1}%
\e{1}%
\e{0}%
\e{0}%
\e{0}%
\e{1}%
\e{0}%
\e{0}%
\e{0}%
\e{0}%
\e{0}%
\e{0}%
\e{0}%
\e{0}%
\e{0}%
\e{0}%
\e{0}%
\e{0}%
\e{0}%
\e{0}%
\e{0}%
\eol}\vss}\rg%
%
%
\rx{\vss\hfull{%
\rlx{\hss{$1400_{zz}$}}\cg%
\e{0}%
\e{0}%
\e{0}%
\e{1}%
\e{1}%
\e{1}%
\e{1}%
\e{1}%
\e{0}%
\e{2}%
\e{0}%
\e{0}%
\e{0}%
\e{0}%
\e{0}%
\e{0}%
\e{0}%
\e{0}%
\e{0}%
\e{0}%
\e{0}%
\e{0}%
\e{0}%
\e{0}%
\e{0}%
\eol}\vss}\rg%
%
%
\rx{\vss\hfull{%
\rlx{\hss{$4200_{z}$}}\cg%
\e{0}%
\e{0}%
\e{1}%
\e{0}%
\e{1}%
\e{1}%
\e{2}%
\e{0}%
\e{1}%
\e{1}%
\e{0}%
\e{0}%
\e{0}%
\e{0}%
\e{0}%
\e{0}%
\e{1}%
\e{0}%
\e{1}%
\e{0}%
\e{0}%
\e{2}%
\e{1}%
\e{1}%
\e{0}%
\eol}\vss}\rg%
%
%
\rx{\vss\hfull{%
\rlx{\hss{$400_{z}$}}\cg%
\e{0}%
\e{1}%
\e{0}%
\e{1}%
\e{1}%
\e{1}%
\e{0}%
\e{1}%
\e{0}%
\e{0}%
\e{0}%
\e{0}%
\e{0}%
\e{0}%
\e{0}%
\e{0}%
\e{0}%
\e{0}%
\e{0}%
\e{0}%
\e{0}%
\e{0}%
\e{0}%
\e{0}%
\e{0}%
\eol}\vss}\rg%
%
%
\rx{\vss\hfull{%
\rlx{\hss{$3240_{z}$}}\cg%
\e{0}%
\e{1}%
\e{1}%
\e{2}%
\e{2}%
\e{3}%
\e{2}%
\e{1}%
\e{1}%
\e{2}%
\e{0}%
\e{0}%
\e{0}%
\e{0}%
\e{0}%
\e{0}%
\e{0}%
\e{0}%
\e{0}%
\e{0}%
\e{0}%
\e{1}%
\e{1}%
\e{1}%
\e{0}%
\eol}\vss}\rg%
%
%
\rx{\vss\hfull{%
\rlx{\hss{$4536_{z}$}}\cg%
\e{0}%
\e{0}%
\e{0}%
\e{1}%
\e{1}%
\e{2}%
\e{1}%
\e{1}%
\e{1}%
\e{3}%
\e{0}%
\e{0}%
\e{0}%
\e{0}%
\e{0}%
\e{0}%
\e{0}%
\e{0}%
\e{0}%
\e{1}%
\e{0}%
\e{1}%
\e{2}%
\e{2}%
\e{1}%
\eol}\vss}\rg%
%
%
\rx{\vss\hfull{%
\rlx{\hss{$2400_{z}$}}\cg%
\e{0}%
\e{0}%
\e{1}%
\e{0}%
\e{0}%
\e{0}%
\e{1}%
\e{0}%
\e{0}%
\e{0}%
\e{0}%
\e{0}%
\e{1}%
\e{0}%
\e{0}%
\e{1}%
\e{1}%
\e{0}%
\e{1}%
\e{0}%
\e{1}%
\e{2}%
\e{0}%
\e{0}%
\e{0}%
\eol}\vss}\rg%
\eop
\eject
\tablecont%
%
%
%
%
%
%
\rowpts=18 true pt%
\colpts=18 true pt%
\rowlabpts=40 true pt%
\collabpts=60 true pt%
\clx{\vss\hfull{%
\rlx{\hss{$ $}}\cg%
\cx{\hskip 16 true pt\flip{$1_p{\times}[{3}]$}\hss}\cg%
\cx{\hskip 16 true pt\flip{$6_p{\times}[{3}]$}\hss}\cg%
\cx{\hskip 16 true pt\flip{$15_p{\times}[{3}]$}\hss}\cg%
\cx{\hskip 16 true pt\flip{$20_p{\times}[{3}]$}\hss}\cg%
\cx{\hskip 16 true pt\flip{$30_p{\times}[{3}]$}\hss}\cg%
\cx{\hskip 16 true pt\flip{$64_p{\times}[{3}]$}\hss}\cg%
\cx{\hskip 16 true pt\flip{$81_p{\times}[{3}]$}\hss}\cg%
\cx{\hskip 16 true pt\flip{$15_q{\times}[{3}]$}\hss}\cg%
\cx{\hskip 16 true pt\flip{$24_p{\times}[{3}]$}\hss}\cg%
\cx{\hskip 16 true pt\flip{$60_p{\times}[{3}]$}\hss}\cg%
\cx{\hskip 16 true pt\flip{$1_p^{*}{\times}[{3}]$}\hss}\cg%
\cx{\hskip 16 true pt\flip{$6_p^{*}{\times}[{3}]$}\hss}\cg%
\cx{\hskip 16 true pt\flip{$15_p^{*}{\times}[{3}]$}\hss}\cg%
\cx{\hskip 16 true pt\flip{$20_p^{*}{\times}[{3}]$}\hss}\cg%
\cx{\hskip 16 true pt\flip{$30_p^{*}{\times}[{3}]$}\hss}\cg%
\cx{\hskip 16 true pt\flip{$64_p^{*}{\times}[{3}]$}\hss}\cg%
\cx{\hskip 16 true pt\flip{$81_p^{*}{\times}[{3}]$}\hss}\cg%
\cx{\hskip 16 true pt\flip{$15_q^{*}{\times}[{3}]$}\hss}\cg%
\cx{\hskip 16 true pt\flip{$24_p^{*}{\times}[{3}]$}\hss}\cg%
\cx{\hskip 16 true pt\flip{$60_p^{*}{\times}[{3}]$}\hss}\cg%
\cx{\hskip 16 true pt\flip{$20_s{\times}[{3}]$}\hss}\cg%
\cx{\hskip 16 true pt\flip{$90_s{\times}[{3}]$}\hss}\cg%
\cx{\hskip 16 true pt\flip{$80_s{\times}[{3}]$}\hss}\cg%
\cx{\hskip 16 true pt\flip{$60_s{\times}[{3}]$}\hss}\cg%
\cx{\hskip 16 true pt\flip{$10_s{\times}[{3}]$}\hss}\cg%
\eol}}\rg%
%
%
\rx{\vss\hfull{%
\rlx{\hss{$3360_{z}$}}\cg%
\e{0}%
\e{0}%
\e{0}%
\e{0}%
\e{1}%
\e{2}%
\e{2}%
\e{1}%
\e{0}%
\e{1}%
\e{0}%
\e{0}%
\e{0}%
\e{0}%
\e{0}%
\e{0}%
\e{1}%
\e{0}%
\e{0}%
\e{0}%
\e{1}%
\e{1}%
\e{1}%
\e{1}%
\e{0}%
\eol}\vss}\rg%
%
%
\rx{\vss\hfull{%
\rlx{\hss{$2800_{z}$}}\cg%
\e{0}%
\e{0}%
\e{0}%
\e{0}%
\e{0}%
\e{2}%
\e{2}%
\e{0}%
\e{0}%
\e{0}%
\e{0}%
\e{0}%
\e{0}%
\e{0}%
\e{0}%
\e{1}%
\e{1}%
\e{0}%
\e{0}%
\e{0}%
\e{1}%
\e{1}%
\e{1}%
\e{0}%
\e{0}%
\eol}\vss}\rg%
%
%
\rx{\vss\hfull{%
\rlx{\hss{$4096_{z}$}}\cg%
\e{0}%
\e{0}%
\e{1}%
\e{1}%
\e{1}%
\e{2}%
\e{2}%
\e{0}%
\e{1}%
\e{2}%
\e{0}%
\e{0}%
\e{0}%
\e{0}%
\e{0}%
\e{0}%
\e{1}%
\e{0}%
\e{0}%
\e{1}%
\e{0}%
\e{2}%
\e{1}%
\e{1}%
\e{0}%
\eol}\vss}\rg%
%
%
\rx{\vss\hfull{%
\rlx{\hss{$5600_{z}$}}\cg%
\e{0}%
\e{0}%
\e{0}%
\e{0}%
\e{1}%
\e{1}%
\e{2}%
\e{0}%
\e{1}%
\e{1}%
\e{0}%
\e{0}%
\e{0}%
\e{0}%
\e{1}%
\e{1}%
\e{2}%
\e{0}%
\e{1}%
\e{1}%
\e{0}%
\e{2}%
\e{2}%
\e{1}%
\e{0}%
\eol}\vss}\rg%
%
%
\rx{\vss\hfull{%
\rlx{\hss{$448_{z}$}}\cg%
\e{1}%
\e{0}%
\e{0}%
\e{1}%
\e{1}%
\e{0}%
\e{0}%
\e{1}%
\e{0}%
\e{1}%
\e{0}%
\e{0}%
\e{0}%
\e{0}%
\e{0}%
\e{0}%
\e{0}%
\e{0}%
\e{0}%
\e{0}%
\e{0}%
\e{0}%
\e{0}%
\e{0}%
\e{1}%
\eol}\vss}\rg%
%
%
\rx{\vss\hfull{%
\rlx{\hss{$448_{w}$}}\cg%
\e{0}%
\e{0}%
\e{0}%
\e{0}%
\e{0}%
\e{0}%
\e{0}%
\e{0}%
\e{0}%
\e{0}%
\e{0}%
\e{0}%
\e{0}%
\e{0}%
\e{0}%
\e{1}%
\e{0}%
\e{0}%
\e{0}%
\e{0}%
\e{1}%
\e{0}%
\e{0}%
\e{0}%
\e{0}%
\eol}\vss}\rg%
%
%
\rx{\vss\hfull{%
\rlx{\hss{$1344_{w}$}}\cg%
\e{0}%
\e{0}%
\e{0}%
\e{0}%
\e{0}%
\e{1}%
\e{0}%
\e{0}%
\e{0}%
\e{0}%
\e{0}%
\e{0}%
\e{0}%
\e{0}%
\e{0}%
\e{0}%
\e{0}%
\e{0}%
\e{0}%
\e{0}%
\e{0}%
\e{0}%
\e{1}%
\e{1}%
\e{0}%
\eol}\vss}\rg%
%
%
\rx{\vss\hfull{%
\rlx{\hss{$5600_{w}$}}\cg%
\e{0}%
\e{0}%
\e{0}%
\e{0}%
\e{0}%
\e{1}%
\e{1}%
\e{0}%
\e{0}%
\e{0}%
\e{0}%
\e{0}%
\e{0}%
\e{0}%
\e{0}%
\e{2}%
\e{2}%
\e{1}%
\e{0}%
\e{1}%
\e{1}%
\e{2}%
\e{2}%
\e{1}%
\e{0}%
\eol}\vss}\rg%
%
%
\rx{\vss\hfull{%
\rlx{\hss{$2016_{w}$}}\cg%
\e{0}%
\e{0}%
\e{0}%
\e{0}%
\e{0}%
\e{0}%
\e{1}%
\e{1}%
\e{0}%
\e{1}%
\e{0}%
\e{0}%
\e{0}%
\e{0}%
\e{0}%
\e{0}%
\e{0}%
\e{0}%
\e{0}%
\e{0}%
\e{0}%
\e{0}%
\e{1}%
\e{1}%
\e{1}%
\eol}\vss}\rg%
%
%
\rx{\vss\hfull{%
\rlx{\hss{$7168_{w}$}}\cg%
\e{0}%
\e{0}%
\e{0}%
\e{0}%
\e{0}%
\e{1}%
\e{2}%
\e{0}%
\e{0}%
\e{2}%
\e{0}%
\e{0}%
\e{0}%
\e{0}%
\e{0}%
\e{1}%
\e{2}%
\e{0}%
\e{0}%
\e{2}%
\e{0}%
\e{2}%
\e{2}%
\e{2}%
\e{0}%
\eol}\vss}\rg%
%
%
%
%
%
%
\rowpts=18 true pt%
\colpts=18 true pt%
\rowlabpts=40 true pt%
\collabpts=60 true pt%
\clx{\vss\hfull{%
\rlx{\hss{$ $}}\cg%
\cx{\hskip 16 true pt\flip{$1_p{\times}[{2}{1}]$}\hss}\cg%
\cx{\hskip 16 true pt\flip{$6_p{\times}[{2}{1}]$}\hss}\cg%
\cx{\hskip 16 true pt\flip{$15_p{\times}[{2}{1}]$}\hss}\cg%
\cx{\hskip 16 true pt\flip{$20_p{\times}[{2}{1}]$}\hss}\cg%
\cx{\hskip 16 true pt\flip{$30_p{\times}[{2}{1}]$}\hss}\cg%
\cx{\hskip 16 true pt\flip{$64_p{\times}[{2}{1}]$}\hss}\cg%
\cx{\hskip 16 true pt\flip{$81_p{\times}[{2}{1}]$}\hss}\cg%
\cx{\hskip 16 true pt\flip{$15_q{\times}[{2}{1}]$}\hss}\cg%
\cx{\hskip 16 true pt\flip{$24_p{\times}[{2}{1}]$}\hss}\cg%
\cx{\hskip 16 true pt\flip{$60_p{\times}[{2}{1}]$}\hss}\cg%
\cx{\hskip 16 true pt\flip{$1_p^{*}{\times}[{2}{1}]$}\hss}\cg%
\cx{\hskip 16 true pt\flip{$6_p^{*}{\times}[{2}{1}]$}\hss}\cg%
\cx{\hskip 16 true pt\flip{$15_p^{*}{\times}[{2}{1}]$}\hss}\cg%
\cx{\hskip 16 true pt\flip{$20_p^{*}{\times}[{2}{1}]$}\hss}\cg%
\cx{\hskip 16 true pt\flip{$30_p^{*}{\times}[{2}{1}]$}\hss}\cg%
\cx{\hskip 16 true pt\flip{$64_p^{*}{\times}[{2}{1}]$}\hss}\cg%
\cx{\hskip 16 true pt\flip{$81_p^{*}{\times}[{2}{1}]$}\hss}\cg%
\cx{\hskip 16 true pt\flip{$15_q^{*}{\times}[{2}{1}]$}\hss}\cg%
\cx{\hskip 16 true pt\flip{$24_p^{*}{\times}[{2}{1}]$}\hss}\cg%
\cx{\hskip 16 true pt\flip{$60_p^{*}{\times}[{2}{1}]$}\hss}\cg%
\cx{\hskip 16 true pt\flip{$20_s{\times}[{2}{1}]$}\hss}\cg%
\cx{\hskip 16 true pt\flip{$90_s{\times}[{2}{1}]$}\hss}\cg%
\cx{\hskip 16 true pt\flip{$80_s{\times}[{2}{1}]$}\hss}\cg%
\cx{\hskip 16 true pt\flip{$60_s{\times}[{2}{1}]$}\hss}\cg%
\cx{\hskip 16 true pt\flip{$10_s{\times}[{2}{1}]$}\hss}\cg%
\eol}}\rg%
%
%
\rx{\vss\hfull{%
\rlx{\hss{$1_{x}$}}\cg%
\e{0}%
\e{0}%
\e{0}%
\e{0}%
\e{0}%
\e{0}%
\e{0}%
\e{0}%
\e{0}%
\e{0}%
\e{0}%
\e{0}%
\e{0}%
\e{0}%
\e{0}%
\e{0}%
\e{0}%
\e{0}%
\e{0}%
\e{0}%
\e{0}%
\e{0}%
\e{0}%
\e{0}%
\e{0}%
\eol}\vss}\rg%
%
%
\rx{\vss\hfull{%
\rlx{\hss{$28_{x}$}}\cg%
\e{0}%
\e{1}%
\e{0}%
\e{0}%
\e{0}%
\e{0}%
\e{0}%
\e{0}%
\e{0}%
\e{0}%
\e{0}%
\e{0}%
\e{0}%
\e{0}%
\e{0}%
\e{0}%
\e{0}%
\e{0}%
\e{0}%
\e{0}%
\e{0}%
\e{0}%
\e{0}%
\e{0}%
\e{0}%
\eol}\vss}\rg%
%
%
\rx{\vss\hfull{%
\rlx{\hss{$35_{x}$}}\cg%
\e{1}%
\e{1}%
\e{0}%
\e{0}%
\e{0}%
\e{0}%
\e{0}%
\e{0}%
\e{0}%
\e{0}%
\e{0}%
\e{0}%
\e{0}%
\e{0}%
\e{0}%
\e{0}%
\e{0}%
\e{0}%
\e{0}%
\e{0}%
\e{0}%
\e{0}%
\e{0}%
\e{0}%
\e{0}%
\eol}\vss}\rg%
%
%
\rx{\vss\hfull{%
\rlx{\hss{$84_{x}$}}\cg%
\e{1}%
\e{0}%
\e{0}%
\e{1}%
\e{0}%
\e{0}%
\e{0}%
\e{0}%
\e{0}%
\e{0}%
\e{0}%
\e{0}%
\e{0}%
\e{0}%
\e{0}%
\e{0}%
\e{0}%
\e{0}%
\e{0}%
\e{0}%
\e{0}%
\e{0}%
\e{0}%
\e{0}%
\e{0}%
\eol}\vss}\rg%
%
%
\rx{\vss\hfull{%
\rlx{\hss{$50_{x}$}}\cg%
\e{0}%
\e{0}%
\e{0}%
\e{0}%
\e{0}%
\e{0}%
\e{0}%
\e{1}%
\e{0}%
\e{0}%
\e{0}%
\e{0}%
\e{0}%
\e{0}%
\e{0}%
\e{0}%
\e{0}%
\e{0}%
\e{0}%
\e{0}%
\e{0}%
\e{0}%
\e{0}%
\e{0}%
\e{0}%
\eol}\vss}\rg%
%
%
\rx{\vss\hfull{%
\rlx{\hss{$350_{x}$}}\cg%
\e{0}%
\e{1}%
\e{1}%
\e{0}%
\e{0}%
\e{1}%
\e{0}%
\e{0}%
\e{0}%
\e{0}%
\e{0}%
\e{0}%
\e{0}%
\e{0}%
\e{0}%
\e{0}%
\e{0}%
\e{0}%
\e{0}%
\e{0}%
\e{1}%
\e{0}%
\e{0}%
\e{0}%
\e{0}%
\eol}\vss}\rg%
%
%
\rx{\vss\hfull{%
\rlx{\hss{$300_{x}$}}\cg%
\e{0}%
\e{1}%
\e{0}%
\e{1}%
\e{0}%
\e{1}%
\e{0}%
\e{0}%
\e{0}%
\e{0}%
\e{0}%
\e{0}%
\e{0}%
\e{0}%
\e{0}%
\e{0}%
\e{0}%
\e{0}%
\e{0}%
\e{0}%
\e{0}%
\e{0}%
\e{0}%
\e{0}%
\e{0}%
\eol}\vss}\rg%
%
%
\rx{\vss\hfull{%
\rlx{\hss{$567_{x}$}}\cg%
\e{1}%
\e{2}%
\e{1}%
\e{2}%
\e{1}%
\e{1}%
\e{0}%
\e{0}%
\e{0}%
\e{0}%
\e{0}%
\e{0}%
\e{0}%
\e{0}%
\e{0}%
\e{0}%
\e{0}%
\e{0}%
\e{0}%
\e{0}%
\e{0}%
\e{0}%
\e{0}%
\e{0}%
\e{0}%
\eol}\vss}\rg%
%
%
\rx{\vss\hfull{%
\rlx{\hss{$210_{x}$}}\cg%
\e{1}%
\e{1}%
\e{0}%
\e{1}%
\e{1}%
\e{0}%
\e{0}%
\e{0}%
\e{0}%
\e{0}%
\e{0}%
\e{0}%
\e{0}%
\e{0}%
\e{0}%
\e{0}%
\e{0}%
\e{0}%
\e{0}%
\e{0}%
\e{0}%
\e{0}%
\e{0}%
\e{0}%
\e{0}%
\eol}\vss}\rg%
%
%
\rx{\vss\hfull{%
\rlx{\hss{$840_{x}$}}\cg%
\e{0}%
\e{0}%
\e{0}%
\e{0}%
\e{0}%
\e{1}%
\e{0}%
\e{0}%
\e{1}%
\e{1}%
\e{0}%
\e{0}%
\e{0}%
\e{0}%
\e{0}%
\e{0}%
\e{0}%
\e{0}%
\e{0}%
\e{0}%
\e{0}%
\e{0}%
\e{1}%
\e{1}%
\e{0}%
\eol}\vss}\rg%
%
%
\rx{\vss\hfull{%
\rlx{\hss{$700_{x}$}}\cg%
\e{0}%
\e{1}%
\e{0}%
\e{2}%
\e{1}%
\e{1}%
\e{0}%
\e{1}%
\e{0}%
\e{1}%
\e{0}%
\e{0}%
\e{0}%
\e{0}%
\e{0}%
\e{0}%
\e{0}%
\e{0}%
\e{0}%
\e{0}%
\e{0}%
\e{0}%
\e{0}%
\e{0}%
\e{0}%
\eol}\vss}\rg%
%
%
\rx{\vss\hfull{%
\rlx{\hss{$175_{x}$}}\cg%
\e{0}%
\e{0}%
\e{0}%
\e{0}%
\e{0}%
\e{0}%
\e{0}%
\e{0}%
\e{0}%
\e{1}%
\e{0}%
\e{0}%
\e{0}%
\e{0}%
\e{0}%
\e{0}%
\e{0}%
\e{0}%
\e{0}%
\e{0}%
\e{0}%
\e{0}%
\e{0}%
\e{0}%
\e{0}%
\eol}\vss}\rg%
%
%
\rx{\vss\hfull{%
\rlx{\hss{$1400_{x}$}}\cg%
\e{0}%
\e{1}%
\e{0}%
\e{1}%
\e{2}%
\e{2}%
\e{1}%
\e{1}%
\e{0}%
\e{1}%
\e{0}%
\e{0}%
\e{0}%
\e{0}%
\e{0}%
\e{0}%
\e{0}%
\e{0}%
\e{0}%
\e{0}%
\e{0}%
\e{0}%
\e{1}%
\e{0}%
\e{0}%
\eol}\vss}\rg%
%
%
\rx{\vss\hfull{%
\rlx{\hss{$1050_{x}$}}\cg%
\e{0}%
\e{0}%
\e{0}%
\e{1}%
\e{1}%
\e{1}%
\e{1}%
\e{2}%
\e{0}%
\e{1}%
\e{0}%
\e{0}%
\e{0}%
\e{0}%
\e{0}%
\e{0}%
\e{0}%
\e{0}%
\e{0}%
\e{0}%
\e{0}%
\e{0}%
\e{0}%
\e{1}%
\e{0}%
\eol}\vss}\rg%
%
%
\rx{\vss\hfull{%
\rlx{\hss{$1575_{x}$}}\cg%
\e{0}%
\e{1}%
\e{2}%
\e{2}%
\e{2}%
\e{2}%
\e{1}%
\e{0}%
\e{0}%
\e{1}%
\e{0}%
\e{0}%
\e{0}%
\e{0}%
\e{0}%
\e{0}%
\e{0}%
\e{0}%
\e{0}%
\e{0}%
\e{0}%
\e{1}%
\e{0}%
\e{0}%
\e{0}%
\eol}\vss}\rg%
%
%
\rx{\vss\hfull{%
\rlx{\hss{$1344_{x}$}}\cg%
\e{1}%
\e{1}%
\e{1}%
\e{3}%
\e{1}%
\e{2}%
\e{1}%
\e{1}%
\e{1}%
\e{1}%
\e{0}%
\e{0}%
\e{0}%
\e{0}%
\e{0}%
\e{0}%
\e{0}%
\e{0}%
\e{0}%
\e{0}%
\e{0}%
\e{0}%
\e{0}%
\e{0}%
\e{0}%
\eol}\vss}\rg%
%
%
\rx{\vss\hfull{%
\rlx{\hss{$2100_{x}$}}\cg%
\e{0}%
\e{0}%
\e{2}%
\e{1}%
\e{1}%
\e{2}%
\e{2}%
\e{0}%
\e{1}%
\e{0}%
\e{0}%
\e{0}%
\e{0}%
\e{0}%
\e{0}%
\e{0}%
\e{1}%
\e{0}%
\e{0}%
\e{0}%
\e{1}%
\e{2}%
\e{0}%
\e{0}%
\e{0}%
\eol}\vss}\rg%
%
%
\rx{\vss\hfull{%
\rlx{\hss{$2268_{x}$}}\cg%
\e{0}%
\e{1}%
\e{1}%
\e{2}%
\e{2}%
\e{3}%
\e{2}%
\e{0}%
\e{1}%
\e{1}%
\e{0}%
\e{0}%
\e{0}%
\e{0}%
\e{0}%
\e{0}%
\e{0}%
\e{0}%
\e{0}%
\e{0}%
\e{0}%
\e{1}%
\e{1}%
\e{0}%
\e{0}%
\eol}\vss}\rg%
%
%
\rx{\vss\hfull{%
\rlx{\hss{$525_{x}$}}\cg%
\e{0}%
\e{0}%
\e{0}%
\e{1}%
\e{0}%
\e{1}%
\e{1}%
\e{0}%
\e{0}%
\e{0}%
\e{0}%
\e{0}%
\e{0}%
\e{0}%
\e{0}%
\e{0}%
\e{0}%
\e{0}%
\e{0}%
\e{0}%
\e{0}%
\e{0}%
\e{0}%
\e{0}%
\e{0}%
\eol}\vss}\rg%
%
%
\rx{\vss\hfull{%
\rlx{\hss{$700_{xx}$}}\cg%
\e{0}%
\e{0}%
\e{0}%
\e{0}%
\e{0}%
\e{0}%
\e{1}%
\e{1}%
\e{0}%
\e{1}%
\e{0}%
\e{0}%
\e{0}%
\e{0}%
\e{0}%
\e{0}%
\e{0}%
\e{0}%
\e{1}%
\e{0}%
\e{0}%
\e{0}%
\e{0}%
\e{1}%
\e{0}%
\eol}\vss}\rg%
%
%
\rx{\vss\hfull{%
\rlx{\hss{$972_{x}$}}\cg%
\e{0}%
\e{0}%
\e{0}%
\e{1}%
\e{0}%
\e{1}%
\e{1}%
\e{1}%
\e{1}%
\e{1}%
\e{0}%
\e{0}%
\e{0}%
\e{0}%
\e{0}%
\e{0}%
\e{0}%
\e{0}%
\e{0}%
\e{0}%
\e{0}%
\e{0}%
\e{0}%
\e{1}%
\e{0}%
\eol}\vss}\rg%
\eop
\eject
\tablecont%
%
%
%
%
%
%
\rowpts=18 true pt%
\colpts=18 true pt%
\rowlabpts=40 true pt%
\collabpts=60 true pt%
\clx{\vss\hfull{%
\rlx{\hss{$ $}}\cg%
\cx{\hskip 16 true pt\flip{$1_p{\times}[{2}{1}]$}\hss}\cg%
\cx{\hskip 16 true pt\flip{$6_p{\times}[{2}{1}]$}\hss}\cg%
\cx{\hskip 16 true pt\flip{$15_p{\times}[{2}{1}]$}\hss}\cg%
\cx{\hskip 16 true pt\flip{$20_p{\times}[{2}{1}]$}\hss}\cg%
\cx{\hskip 16 true pt\flip{$30_p{\times}[{2}{1}]$}\hss}\cg%
\cx{\hskip 16 true pt\flip{$64_p{\times}[{2}{1}]$}\hss}\cg%
\cx{\hskip 16 true pt\flip{$81_p{\times}[{2}{1}]$}\hss}\cg%
\cx{\hskip 16 true pt\flip{$15_q{\times}[{2}{1}]$}\hss}\cg%
\cx{\hskip 16 true pt\flip{$24_p{\times}[{2}{1}]$}\hss}\cg%
\cx{\hskip 16 true pt\flip{$60_p{\times}[{2}{1}]$}\hss}\cg%
\cx{\hskip 16 true pt\flip{$1_p^{*}{\times}[{2}{1}]$}\hss}\cg%
\cx{\hskip 16 true pt\flip{$6_p^{*}{\times}[{2}{1}]$}\hss}\cg%
\cx{\hskip 16 true pt\flip{$15_p^{*}{\times}[{2}{1}]$}\hss}\cg%
\cx{\hskip 16 true pt\flip{$20_p^{*}{\times}[{2}{1}]$}\hss}\cg%
\cx{\hskip 16 true pt\flip{$30_p^{*}{\times}[{2}{1}]$}\hss}\cg%
\cx{\hskip 16 true pt\flip{$64_p^{*}{\times}[{2}{1}]$}\hss}\cg%
\cx{\hskip 16 true pt\flip{$81_p^{*}{\times}[{2}{1}]$}\hss}\cg%
\cx{\hskip 16 true pt\flip{$15_q^{*}{\times}[{2}{1}]$}\hss}\cg%
\cx{\hskip 16 true pt\flip{$24_p^{*}{\times}[{2}{1}]$}\hss}\cg%
\cx{\hskip 16 true pt\flip{$60_p^{*}{\times}[{2}{1}]$}\hss}\cg%
\cx{\hskip 16 true pt\flip{$20_s{\times}[{2}{1}]$}\hss}\cg%
\cx{\hskip 16 true pt\flip{$90_s{\times}[{2}{1}]$}\hss}\cg%
\cx{\hskip 16 true pt\flip{$80_s{\times}[{2}{1}]$}\hss}\cg%
\cx{\hskip 16 true pt\flip{$60_s{\times}[{2}{1}]$}\hss}\cg%
\cx{\hskip 16 true pt\flip{$10_s{\times}[{2}{1}]$}\hss}\cg%
\eol}}\rg%
%
%
\rx{\vss\hfull{%
\rlx{\hss{$4096_{x}$}}\cg%
\e{0}%
\e{1}%
\e{1}%
\e{2}%
\e{2}%
\e{5}%
\e{3}%
\e{1}%
\e{1}%
\e{2}%
\e{0}%
\e{0}%
\e{0}%
\e{0}%
\e{0}%
\e{0}%
\e{1}%
\e{0}%
\e{0}%
\e{0}%
\e{1}%
\e{2}%
\e{2}%
\e{1}%
\e{0}%
\eol}\vss}\rg%
%
%
\rx{\vss\hfull{%
\rlx{\hss{$4200_{x}$}}\cg%
\e{0}%
\e{0}%
\e{1}%
\e{1}%
\e{2}%
\e{3}%
\e{3}%
\e{1}%
\e{1}%
\e{4}%
\e{0}%
\e{0}%
\e{0}%
\e{0}%
\e{0}%
\e{0}%
\e{1}%
\e{0}%
\e{0}%
\e{1}%
\e{0}%
\e{2}%
\e{2}%
\e{2}%
\e{0}%
\eol}\vss}\rg%
%
%
\rx{\vss\hfull{%
\rlx{\hss{$2240_{x}$}}\cg%
\e{0}%
\e{0}%
\e{0}%
\e{1}%
\e{2}%
\e{2}%
\e{1}%
\e{1}%
\e{1}%
\e{3}%
\e{0}%
\e{0}%
\e{0}%
\e{0}%
\e{0}%
\e{0}%
\e{0}%
\e{0}%
\e{0}%
\e{0}%
\e{0}%
\e{1}%
\e{1}%
\e{1}%
\e{1}%
\eol}\vss}\rg%
%
%
\rx{\vss\hfull{%
\rlx{\hss{$2835_{x}$}}\cg%
\e{0}%
\e{0}%
\e{0}%
\e{0}%
\e{1}%
\e{1}%
\e{2}%
\e{1}%
\e{0}%
\e{3}%
\e{0}%
\e{0}%
\e{0}%
\e{0}%
\e{0}%
\e{0}%
\e{1}%
\e{0}%
\e{0}%
\e{1}%
\e{0}%
\e{1}%
\e{2}%
\e{2}%
\e{1}%
\eol}\vss}\rg%
%
%
\rx{\vss\hfull{%
\rlx{\hss{$6075_{x}$}}\cg%
\e{0}%
\e{0}%
\e{1}%
\e{1}%
\e{2}%
\e{4}%
\e{5}%
\e{1}%
\e{1}%
\e{3}%
\e{0}%
\e{0}%
\e{0}%
\e{0}%
\e{0}%
\e{1}%
\e{2}%
\e{0}%
\e{1}%
\e{1}%
\e{1}%
\e{4}%
\e{3}%
\e{2}%
\e{0}%
\eol}\vss}\rg%
%
%
\rx{\vss\hfull{%
\rlx{\hss{$3200_{x}$}}\cg%
\e{0}%
\e{0}%
\e{0}%
\e{1}%
\e{0}%
\e{2}%
\e{2}%
\e{0}%
\e{2}%
\e{2}%
\e{0}%
\e{0}%
\e{0}%
\e{0}%
\e{0}%
\e{0}%
\e{2}%
\e{0}%
\e{0}%
\e{1}%
\e{0}%
\e{2}%
\e{1}%
\e{2}%
\e{0}%
\eol}\vss}\rg%
%
%
\rx{\vss\hfull{%
\rlx{\hss{$70_{y}$}}\cg%
\e{0}%
\e{0}%
\e{0}%
\e{0}%
\e{0}%
\e{0}%
\e{0}%
\e{0}%
\e{0}%
\e{0}%
\e{0}%
\e{0}%
\e{0}%
\e{0}%
\e{0}%
\e{0}%
\e{0}%
\e{0}%
\e{0}%
\e{0}%
\e{1}%
\e{0}%
\e{0}%
\e{0}%
\e{0}%
\eol}\vss}\rg%
%
%
\rx{\vss\hfull{%
\rlx{\hss{$1134_{y}$}}\cg%
\e{0}%
\e{0}%
\e{0}%
\e{0}%
\e{0}%
\e{1}%
\e{0}%
\e{0}%
\e{0}%
\e{0}%
\e{0}%
\e{0}%
\e{0}%
\e{0}%
\e{0}%
\e{1}%
\e{0}%
\e{0}%
\e{0}%
\e{0}%
\e{1}%
\e{1}%
\e{1}%
\e{1}%
\e{0}%
\eol}\vss}\rg%
%
%
\rx{\vss\hfull{%
\rlx{\hss{$1680_{y}$}}\cg%
\e{0}%
\e{0}%
\e{1}%
\e{0}%
\e{0}%
\e{1}%
\e{1}%
\e{0}%
\e{0}%
\e{0}%
\e{0}%
\e{0}%
\e{1}%
\e{0}%
\e{0}%
\e{1}%
\e{1}%
\e{0}%
\e{0}%
\e{0}%
\e{2}%
\e{2}%
\e{0}%
\e{0}%
\e{0}%
\eol}\vss}\rg%
%
%
\rx{\vss\hfull{%
\rlx{\hss{$168_{y}$}}\cg%
\e{0}%
\e{0}%
\e{0}%
\e{0}%
\e{0}%
\e{0}%
\e{0}%
\e{0}%
\e{0}%
\e{0}%
\e{0}%
\e{0}%
\e{0}%
\e{0}%
\e{0}%
\e{0}%
\e{0}%
\e{0}%
\e{0}%
\e{0}%
\e{0}%
\e{0}%
\e{0}%
\e{1}%
\e{0}%
\eol}\vss}\rg%
%
%
\rx{\vss\hfull{%
\rlx{\hss{$420_{y}$}}\cg%
\e{0}%
\e{0}%
\e{0}%
\e{0}%
\e{0}%
\e{0}%
\e{0}%
\e{0}%
\e{0}%
\e{0}%
\e{0}%
\e{0}%
\e{0}%
\e{0}%
\e{0}%
\e{0}%
\e{0}%
\e{0}%
\e{0}%
\e{0}%
\e{0}%
\e{0}%
\e{1}%
\e{1}%
\e{1}%
\eol}\vss}\rg%
%
%
\rx{\vss\hfull{%
\rlx{\hss{$3150_{y}$}}\cg%
\e{0}%
\e{0}%
\e{0}%
\e{0}%
\e{0}%
\e{1}%
\e{1}%
\e{0}%
\e{0}%
\e{2}%
\e{0}%
\e{0}%
\e{0}%
\e{0}%
\e{0}%
\e{1}%
\e{1}%
\e{0}%
\e{0}%
\e{2}%
\e{0}%
\e{2}%
\e{3}%
\e{2}%
\e{1}%
\eol}\vss}\rg%
%
%
\rx{\vss\hfull{%
\rlx{\hss{$4200_{y}$}}\cg%
\e{0}%
\e{0}%
\e{0}%
\e{0}%
\e{0}%
\e{1}%
\e{2}%
\e{1}%
\e{1}%
\e{2}%
\e{0}%
\e{0}%
\e{0}%
\e{0}%
\e{0}%
\e{1}%
\e{2}%
\e{1}%
\e{1}%
\e{2}%
\e{0}%
\e{2}%
\e{3}%
\e{4}%
\e{1}%
\eol}\vss}\rg%
%
%
\rx{\vss\hfull{%
\rlx{\hss{$2688_{y}$}}\cg%
\e{0}%
\e{0}%
\e{0}%
\e{0}%
\e{0}%
\e{1}%
\e{2}%
\e{1}%
\e{0}%
\e{1}%
\e{0}%
\e{0}%
\e{0}%
\e{0}%
\e{0}%
\e{1}%
\e{2}%
\e{1}%
\e{0}%
\e{1}%
\e{1}%
\e{1}%
\e{1}%
\e{2}%
\e{0}%
\eol}\vss}\rg%
%
%
\rx{\vss\hfull{%
\rlx{\hss{$2100_{y}$}}\cg%
\e{0}%
\e{0}%
\e{0}%
\e{0}%
\e{0}%
\e{1}%
\e{2}%
\e{0}%
\e{1}%
\e{0}%
\e{0}%
\e{0}%
\e{0}%
\e{0}%
\e{0}%
\e{1}%
\e{2}%
\e{0}%
\e{1}%
\e{0}%
\e{1}%
\e{1}%
\e{1}%
\e{0}%
\e{0}%
\eol}\vss}\rg%
%
%
\rx{\vss\hfull{%
\rlx{\hss{$1400_{y}$}}\cg%
\e{0}%
\e{0}%
\e{0}%
\e{0}%
\e{0}%
\e{1}%
\e{1}%
\e{0}%
\e{0}%
\e{0}%
\e{0}%
\e{0}%
\e{0}%
\e{0}%
\e{0}%
\e{1}%
\e{1}%
\e{0}%
\e{0}%
\e{0}%
\e{0}%
\e{1}%
\e{1}%
\e{0}%
\e{0}%
\eol}\vss}\rg%
%
%
\rx{\vss\hfull{%
\rlx{\hss{$4536_{y}$}}\cg%
\e{0}%
\e{0}%
\e{0}%
\e{0}%
\e{1}%
\e{2}%
\e{3}%
\e{0}%
\e{0}%
\e{1}%
\e{0}%
\e{0}%
\e{0}%
\e{0}%
\e{1}%
\e{2}%
\e{3}%
\e{0}%
\e{0}%
\e{1}%
\e{1}%
\e{3}%
\e{3}%
\e{1}%
\e{0}%
\eol}\vss}\rg%
%
%
\rx{\vss\hfull{%
\rlx{\hss{$5670_{y}$}}\cg%
\e{0}%
\e{0}%
\e{1}%
\e{0}%
\e{1}%
\e{2}%
\e{3}%
\e{0}%
\e{1}%
\e{2}%
\e{0}%
\e{0}%
\e{1}%
\e{0}%
\e{1}%
\e{2}%
\e{3}%
\e{0}%
\e{1}%
\e{2}%
\e{1}%
\e{5}%
\e{3}%
\e{1}%
\e{0}%
\eol}\vss}\rg%
%
%
\rx{\vss\hfull{%
\rlx{\hss{$4480_{y}$}}\cg%
\e{0}%
\e{0}%
\e{0}%
\e{0}%
\e{1}%
\e{1}%
\e{2}%
\e{0}%
\e{1}%
\e{2}%
\e{0}%
\e{0}%
\e{0}%
\e{0}%
\e{1}%
\e{1}%
\e{2}%
\e{0}%
\e{1}%
\e{2}%
\e{0}%
\e{3}%
\e{4}%
\e{2}%
\e{1}%
\eol}\vss}\rg%
%
%
\rx{\vss\hfull{%
\rlx{\hss{$8_{z}$}}\cg%
\e{1}%
\e{0}%
\e{0}%
\e{0}%
\e{0}%
\e{0}%
\e{0}%
\e{0}%
\e{0}%
\e{0}%
\e{0}%
\e{0}%
\e{0}%
\e{0}%
\e{0}%
\e{0}%
\e{0}%
\e{0}%
\e{0}%
\e{0}%
\e{0}%
\e{0}%
\e{0}%
\e{0}%
\e{0}%
\eol}\vss}\rg%
%
%
\rx{\vss\hfull{%
\rlx{\hss{$56_{z}$}}\cg%
\e{0}%
\e{0}%
\e{1}%
\e{0}%
\e{0}%
\e{0}%
\e{0}%
\e{0}%
\e{0}%
\e{0}%
\e{0}%
\e{0}%
\e{0}%
\e{0}%
\e{0}%
\e{0}%
\e{0}%
\e{0}%
\e{0}%
\e{0}%
\e{0}%
\e{0}%
\e{0}%
\e{0}%
\e{0}%
\eol}\vss}\rg%
%
%
\rx{\vss\hfull{%
\rlx{\hss{$160_{z}$}}\cg%
\e{1}%
\e{1}%
\e{1}%
\e{1}%
\e{0}%
\e{0}%
\e{0}%
\e{0}%
\e{0}%
\e{0}%
\e{0}%
\e{0}%
\e{0}%
\e{0}%
\e{0}%
\e{0}%
\e{0}%
\e{0}%
\e{0}%
\e{0}%
\e{0}%
\e{0}%
\e{0}%
\e{0}%
\e{0}%
\eol}\vss}\rg%
%
%
\rx{\vss\hfull{%
\rlx{\hss{$112_{z}$}}\cg%
\e{1}%
\e{1}%
\e{0}%
\e{1}%
\e{0}%
\e{0}%
\e{0}%
\e{0}%
\e{0}%
\e{0}%
\e{0}%
\e{0}%
\e{0}%
\e{0}%
\e{0}%
\e{0}%
\e{0}%
\e{0}%
\e{0}%
\e{0}%
\e{0}%
\e{0}%
\e{0}%
\e{0}%
\e{0}%
\eol}\vss}\rg%
%
%
\rx{\vss\hfull{%
\rlx{\hss{$840_{z}$}}\cg%
\e{0}%
\e{0}%
\e{1}%
\e{1}%
\e{0}%
\e{1}%
\e{0}%
\e{0}%
\e{1}%
\e{1}%
\e{0}%
\e{0}%
\e{0}%
\e{0}%
\e{0}%
\e{0}%
\e{0}%
\e{0}%
\e{0}%
\e{0}%
\e{0}%
\e{1}%
\e{0}%
\e{0}%
\e{0}%
\eol}\vss}\rg%
%
%
\rx{\vss\hfull{%
\rlx{\hss{$1296_{z}$}}\cg%
\e{0}%
\e{1}%
\e{2}%
\e{1}%
\e{1}%
\e{2}%
\e{1}%
\e{0}%
\e{0}%
\e{0}%
\e{0}%
\e{0}%
\e{0}%
\e{0}%
\e{0}%
\e{0}%
\e{0}%
\e{0}%
\e{0}%
\e{0}%
\e{1}%
\e{1}%
\e{0}%
\e{0}%
\e{0}%
\eol}\vss}\rg%
%
%
\rx{\vss\hfull{%
\rlx{\hss{$1400_{z}$}}\cg%
\e{1}%
\e{1}%
\e{1}%
\e{3}%
\e{2}%
\e{2}%
\e{1}%
\e{0}%
\e{1}%
\e{1}%
\e{0}%
\e{0}%
\e{0}%
\e{0}%
\e{0}%
\e{0}%
\e{0}%
\e{0}%
\e{0}%
\e{0}%
\e{0}%
\e{0}%
\e{0}%
\e{0}%
\e{0}%
\eol}\vss}\rg%
%
%
\rx{\vss\hfull{%
\rlx{\hss{$1008_{z}$}}\cg%
\e{0}%
\e{2}%
\e{1}%
\e{2}%
\e{1}%
\e{2}%
\e{1}%
\e{0}%
\e{0}%
\e{0}%
\e{0}%
\e{0}%
\e{0}%
\e{0}%
\e{0}%
\e{0}%
\e{0}%
\e{0}%
\e{0}%
\e{0}%
\e{0}%
\e{0}%
\e{0}%
\e{0}%
\e{0}%
\eol}\vss}\rg%
%
%
\rx{\vss\hfull{%
\rlx{\hss{$560_{z}$}}\cg%
\e{1}%
\e{2}%
\e{0}%
\e{2}%
\e{1}%
\e{1}%
\e{0}%
\e{1}%
\e{0}%
\e{0}%
\e{0}%
\e{0}%
\e{0}%
\e{0}%
\e{0}%
\e{0}%
\e{0}%
\e{0}%
\e{0}%
\e{0}%
\e{0}%
\e{0}%
\e{0}%
\e{0}%
\e{0}%
\eol}\vss}\rg%
%
%
\rx{\vss\hfull{%
\rlx{\hss{$1400_{zz}$}}\cg%
\e{0}%
\e{0}%
\e{0}%
\e{0}%
\e{1}%
\e{1}%
\e{1}%
\e{2}%
\e{0}%
\e{2}%
\e{0}%
\e{0}%
\e{0}%
\e{0}%
\e{0}%
\e{0}%
\e{0}%
\e{0}%
\e{0}%
\e{0}%
\e{0}%
\e{0}%
\e{1}%
\e{1}%
\e{1}%
\eol}\vss}\rg%
%
%
\rx{\vss\hfull{%
\rlx{\hss{$4200_{z}$}}\cg%
\e{0}%
\e{0}%
\e{0}%
\e{0}%
\e{1}%
\e{2}%
\e{3}%
\e{1}%
\e{0}%
\e{3}%
\e{0}%
\e{0}%
\e{0}%
\e{0}%
\e{0}%
\e{1}%
\e{1}%
\e{0}%
\e{1}%
\e{1}%
\e{0}%
\e{2}%
\e{3}%
\e{3}%
\e{0}%
\eol}\vss}\rg%
%
%
\rx{\vss\hfull{%
\rlx{\hss{$400_{z}$}}\cg%
\e{0}%
\e{0}%
\e{0}%
\e{1}%
\e{1}%
\e{0}%
\e{0}%
\e{1}%
\e{0}%
\e{1}%
\e{0}%
\e{0}%
\e{0}%
\e{0}%
\e{0}%
\e{0}%
\e{0}%
\e{0}%
\e{0}%
\e{0}%
\e{0}%
\e{0}%
\e{0}%
\e{0}%
\e{0}%
\eol}\vss}\rg%
%
%
\rx{\vss\hfull{%
\rlx{\hss{$3240_{z}$}}\cg%
\e{0}%
\e{1}%
\e{1}%
\e{3}%
\e{3}%
\e{4}%
\e{2}%
\e{2}%
\e{1}%
\e{3}%
\e{0}%
\e{0}%
\e{0}%
\e{0}%
\e{0}%
\e{0}%
\e{0}%
\e{0}%
\e{0}%
\e{0}%
\e{0}%
\e{1}%
\e{1}%
\e{1}%
\e{0}%
\eol}\vss}\rg%
%
%
\rx{\vss\hfull{%
\rlx{\hss{$4536_{z}$}}\cg%
\e{0}%
\e{0}%
\e{0}%
\e{1}%
\e{1}%
\e{3}%
\e{3}%
\e{1}%
\e{2}%
\e{4}%
\e{0}%
\e{0}%
\e{0}%
\e{0}%
\e{0}%
\e{0}%
\e{1}%
\e{0}%
\e{0}%
\e{1}%
\e{0}%
\e{2}%
\e{3}%
\e{3}%
\e{1}%
\eol}\vss}\rg%
%
%
\rx{\vss\hfull{%
\rlx{\hss{$2400_{z}$}}\cg%
\e{0}%
\e{0}%
\e{1}%
\e{0}%
\e{1}%
\e{2}%
\e{2}%
\e{0}%
\e{0}%
\e{0}%
\e{0}%
\e{0}%
\e{0}%
\e{0}%
\e{0}%
\e{1}%
\e{1}%
\e{0}%
\e{0}%
\e{0}%
\e{2}%
\e{2}%
\e{1}%
\e{0}%
\e{0}%
\eol}\vss}\rg%
\eop
\eject
\tablecont%
%
%
%
%
%
%
\rowpts=18 true pt%
\colpts=18 true pt%
\rowlabpts=40 true pt%
\collabpts=60 true pt%
\clx{\vss\hfull{%
\rlx{\hss{$ $}}\cg%
\cx{\hskip 16 true pt\flip{$1_p{\times}[{2}{1}]$}\hss}\cg%
\cx{\hskip 16 true pt\flip{$6_p{\times}[{2}{1}]$}\hss}\cg%
\cx{\hskip 16 true pt\flip{$15_p{\times}[{2}{1}]$}\hss}\cg%
\cx{\hskip 16 true pt\flip{$20_p{\times}[{2}{1}]$}\hss}\cg%
\cx{\hskip 16 true pt\flip{$30_p{\times}[{2}{1}]$}\hss}\cg%
\cx{\hskip 16 true pt\flip{$64_p{\times}[{2}{1}]$}\hss}\cg%
\cx{\hskip 16 true pt\flip{$81_p{\times}[{2}{1}]$}\hss}\cg%
\cx{\hskip 16 true pt\flip{$15_q{\times}[{2}{1}]$}\hss}\cg%
\cx{\hskip 16 true pt\flip{$24_p{\times}[{2}{1}]$}\hss}\cg%
\cx{\hskip 16 true pt\flip{$60_p{\times}[{2}{1}]$}\hss}\cg%
\cx{\hskip 16 true pt\flip{$1_p^{*}{\times}[{2}{1}]$}\hss}\cg%
\cx{\hskip 16 true pt\flip{$6_p^{*}{\times}[{2}{1}]$}\hss}\cg%
\cx{\hskip 16 true pt\flip{$15_p^{*}{\times}[{2}{1}]$}\hss}\cg%
\cx{\hskip 16 true pt\flip{$20_p^{*}{\times}[{2}{1}]$}\hss}\cg%
\cx{\hskip 16 true pt\flip{$30_p^{*}{\times}[{2}{1}]$}\hss}\cg%
\cx{\hskip 16 true pt\flip{$64_p^{*}{\times}[{2}{1}]$}\hss}\cg%
\cx{\hskip 16 true pt\flip{$81_p^{*}{\times}[{2}{1}]$}\hss}\cg%
\cx{\hskip 16 true pt\flip{$15_q^{*}{\times}[{2}{1}]$}\hss}\cg%
\cx{\hskip 16 true pt\flip{$24_p^{*}{\times}[{2}{1}]$}\hss}\cg%
\cx{\hskip 16 true pt\flip{$60_p^{*}{\times}[{2}{1}]$}\hss}\cg%
\cx{\hskip 16 true pt\flip{$20_s{\times}[{2}{1}]$}\hss}\cg%
\cx{\hskip 16 true pt\flip{$90_s{\times}[{2}{1}]$}\hss}\cg%
\cx{\hskip 16 true pt\flip{$80_s{\times}[{2}{1}]$}\hss}\cg%
\cx{\hskip 16 true pt\flip{$60_s{\times}[{2}{1}]$}\hss}\cg%
\cx{\hskip 16 true pt\flip{$10_s{\times}[{2}{1}]$}\hss}\cg%
\eol}}\rg%
%
%
\rx{\vss\hfull{%
\rlx{\hss{$3360_{z}$}}\cg%
\e{0}%
\e{0}%
\e{1}%
\e{1}%
\e{2}%
\e{2}%
\e{3}%
\e{1}%
\e{0}%
\e{3}%
\e{0}%
\e{0}%
\e{0}%
\e{0}%
\e{0}%
\e{0}%
\e{1}%
\e{0}%
\e{0}%
\e{1}%
\e{0}%
\e{2}%
\e{1}%
\e{1}%
\e{0}%
\eol}\vss}\rg%
%
%
\rx{\vss\hfull{%
\rlx{\hss{$2800_{z}$}}\cg%
\e{0}%
\e{0}%
\e{1}%
\e{1}%
\e{2}%
\e{2}%
\e{3}%
\e{0}%
\e{1}%
\e{1}%
\e{0}%
\e{0}%
\e{0}%
\e{0}%
\e{0}%
\e{0}%
\e{1}%
\e{0}%
\e{1}%
\e{0}%
\e{0}%
\e{2}%
\e{1}%
\e{0}%
\e{0}%
\eol}\vss}\rg%
%
%
\rx{\vss\hfull{%
\rlx{\hss{$4096_{z}$}}\cg%
\e{0}%
\e{1}%
\e{1}%
\e{2}%
\e{2}%
\e{5}%
\e{3}%
\e{1}%
\e{1}%
\e{2}%
\e{0}%
\e{0}%
\e{0}%
\e{0}%
\e{0}%
\e{0}%
\e{1}%
\e{0}%
\e{0}%
\e{0}%
\e{1}%
\e{2}%
\e{2}%
\e{1}%
\e{0}%
\eol}\vss}\rg%
%
%
\rx{\vss\hfull{%
\rlx{\hss{$5600_{z}$}}\cg%
\e{0}%
\e{0}%
\e{1}%
\e{1}%
\e{1}%
\e{4}%
\e{4}%
\e{0}%
\e{2}%
\e{2}%
\e{0}%
\e{0}%
\e{0}%
\e{0}%
\e{0}%
\e{1}%
\e{3}%
\e{0}%
\e{0}%
\e{1}%
\e{1}%
\e{4}%
\e{3}%
\e{1}%
\e{0}%
\eol}\vss}\rg%
%
%
\rx{\vss\hfull{%
\rlx{\hss{$448_{z}$}}\cg%
\e{0}%
\e{0}%
\e{0}%
\e{1}%
\e{0}%
\e{1}%
\e{0}%
\e{0}%
\e{0}%
\e{1}%
\e{0}%
\e{0}%
\e{0}%
\e{0}%
\e{0}%
\e{0}%
\e{0}%
\e{0}%
\e{0}%
\e{0}%
\e{0}%
\e{0}%
\e{0}%
\e{0}%
\e{0}%
\eol}\vss}\rg%
%
%
\rx{\vss\hfull{%
\rlx{\hss{$448_{w}$}}\cg%
\e{0}%
\e{0}%
\e{1}%
\e{0}%
\e{0}%
\e{0}%
\e{0}%
\e{0}%
\e{0}%
\e{0}%
\e{0}%
\e{0}%
\e{1}%
\e{0}%
\e{0}%
\e{0}%
\e{0}%
\e{0}%
\e{0}%
\e{0}%
\e{1}%
\e{1}%
\e{0}%
\e{0}%
\e{0}%
\eol}\vss}\rg%
%
%
\rx{\vss\hfull{%
\rlx{\hss{$1344_{w}$}}\cg%
\e{0}%
\e{0}%
\e{0}%
\e{0}%
\e{0}%
\e{0}%
\e{0}%
\e{0}%
\e{1}%
\e{1}%
\e{0}%
\e{0}%
\e{0}%
\e{0}%
\e{0}%
\e{0}%
\e{0}%
\e{0}%
\e{1}%
\e{1}%
\e{0}%
\e{1}%
\e{1}%
\e{2}%
\e{1}%
\eol}\vss}\rg%
%
%
\rx{\vss\hfull{%
\rlx{\hss{$5600_{w}$}}\cg%
\e{0}%
\e{0}%
\e{1}%
\e{0}%
\e{1}%
\e{2}%
\e{3}%
\e{0}%
\e{1}%
\e{2}%
\e{0}%
\e{0}%
\e{1}%
\e{0}%
\e{1}%
\e{2}%
\e{3}%
\e{0}%
\e{1}%
\e{2}%
\e{1}%
\e{5}%
\e{2}%
\e{2}%
\e{0}%
\eol}\vss}\rg%
%
%
\rx{\vss\hfull{%
\rlx{\hss{$2016_{w}$}}\cg%
\e{0}%
\e{0}%
\e{0}%
\e{0}%
\e{0}%
\e{0}%
\e{1}%
\e{0}%
\e{0}%
\e{2}%
\e{0}%
\e{0}%
\e{0}%
\e{0}%
\e{0}%
\e{0}%
\e{1}%
\e{0}%
\e{0}%
\e{2}%
\e{0}%
\e{1}%
\e{1}%
\e{2}%
\e{1}%
\eol}\vss}\rg%
%
%
\rx{\vss\hfull{%
\rlx{\hss{$7168_{w}$}}\cg%
\e{0}%
\e{0}%
\e{0}%
\e{0}%
\e{1}%
\e{2}%
\e{4}%
\e{1}%
\e{1}%
\e{3}%
\e{0}%
\e{0}%
\e{0}%
\e{0}%
\e{1}%
\e{2}%
\e{4}%
\e{1}%
\e{1}%
\e{3}%
\e{1}%
\e{4}%
\e{5}%
\e{4}%
\e{1}%
\eol}\vss}\rg%
%
%
%
%
%
%
\rowpts=18 true pt%
\colpts=18 true pt%
\rowlabpts=40 true pt%
\collabpts=60 true pt%
\clx{\vss\hfull{%
\rlx{\hss{$ $}}\cg%
\cx{\hskip 16 true pt\flip{$1_p{\times}[{1^{3}}]$}\hss}\cg%
\cx{\hskip 16 true pt\flip{$6_p{\times}[{1^{3}}]$}\hss}\cg%
\cx{\hskip 16 true pt\flip{$15_p{\times}[{1^{3}}]$}\hss}\cg%
\cx{\hskip 16 true pt\flip{$20_p{\times}[{1^{3}}]$}\hss}\cg%
\cx{\hskip 16 true pt\flip{$30_p{\times}[{1^{3}}]$}\hss}\cg%
\cx{\hskip 16 true pt\flip{$64_p{\times}[{1^{3}}]$}\hss}\cg%
\cx{\hskip 16 true pt\flip{$81_p{\times}[{1^{3}}]$}\hss}\cg%
\cx{\hskip 16 true pt\flip{$15_q{\times}[{1^{3}}]$}\hss}\cg%
\cx{\hskip 16 true pt\flip{$24_p{\times}[{1^{3}}]$}\hss}\cg%
\cx{\hskip 16 true pt\flip{$60_p{\times}[{1^{3}}]$}\hss}\cg%
\cx{\hskip 16 true pt\flip{$1_p^{*}{\times}[{1^{3}}]$}\hss}\cg%
\cx{\hskip 16 true pt\flip{$6_p^{*}{\times}[{1^{3}}]$}\hss}\cg%
\cx{\hskip 16 true pt\flip{$15_p^{*}{\times}[{1^{3}}]$}\hss}\cg%
\cx{\hskip 16 true pt\flip{$20_p^{*}{\times}[{1^{3}}]$}\hss}\cg%
\cx{\hskip 16 true pt\flip{$30_p^{*}{\times}[{1^{3}}]$}\hss}\cg%
\cx{\hskip 16 true pt\flip{$64_p^{*}{\times}[{1^{3}}]$}\hss}\cg%
\cx{\hskip 16 true pt\flip{$81_p^{*}{\times}[{1^{3}}]$}\hss}\cg%
\cx{\hskip 16 true pt\flip{$15_q^{*}{\times}[{1^{3}}]$}\hss}\cg%
\cx{\hskip 16 true pt\flip{$24_p^{*}{\times}[{1^{3}}]$}\hss}\cg%
\cx{\hskip 16 true pt\flip{$60_p^{*}{\times}[{1^{3}}]$}\hss}\cg%
\cx{\hskip 16 true pt\flip{$20_s{\times}[{1^{3}}]$}\hss}\cg%
\cx{\hskip 16 true pt\flip{$90_s{\times}[{1^{3}}]$}\hss}\cg%
\cx{\hskip 16 true pt\flip{$80_s{\times}[{1^{3}}]$}\hss}\cg%
\cx{\hskip 16 true pt\flip{$60_s{\times}[{1^{3}}]$}\hss}\cg%
\cx{\hskip 16 true pt\flip{$10_s{\times}[{1^{3}}]$}\hss}\cg%
\eol}}\rg%
%
%
\rx{\vss\hfull{%
\rlx{\hss{$1_{x}$}}\cg%
\e{0}%
\e{0}%
\e{0}%
\e{0}%
\e{0}%
\e{0}%
\e{0}%
\e{0}%
\e{0}%
\e{0}%
\e{0}%
\e{0}%
\e{0}%
\e{0}%
\e{0}%
\e{0}%
\e{0}%
\e{0}%
\e{0}%
\e{0}%
\e{0}%
\e{0}%
\e{0}%
\e{0}%
\e{0}%
\eol}\vss}\rg%
%
%
\rx{\vss\hfull{%
\rlx{\hss{$28_{x}$}}\cg%
\e{1}%
\e{0}%
\e{0}%
\e{0}%
\e{0}%
\e{0}%
\e{0}%
\e{0}%
\e{0}%
\e{0}%
\e{0}%
\e{0}%
\e{0}%
\e{0}%
\e{0}%
\e{0}%
\e{0}%
\e{0}%
\e{0}%
\e{0}%
\e{0}%
\e{0}%
\e{0}%
\e{0}%
\e{0}%
\eol}\vss}\rg%
%
%
\rx{\vss\hfull{%
\rlx{\hss{$35_{x}$}}\cg%
\e{0}%
\e{0}%
\e{0}%
\e{0}%
\e{0}%
\e{0}%
\e{0}%
\e{0}%
\e{0}%
\e{0}%
\e{0}%
\e{0}%
\e{0}%
\e{0}%
\e{0}%
\e{0}%
\e{0}%
\e{0}%
\e{0}%
\e{0}%
\e{0}%
\e{0}%
\e{0}%
\e{0}%
\e{0}%
\eol}\vss}\rg%
%
%
\rx{\vss\hfull{%
\rlx{\hss{$84_{x}$}}\cg%
\e{0}%
\e{0}%
\e{0}%
\e{0}%
\e{0}%
\e{0}%
\e{0}%
\e{0}%
\e{0}%
\e{0}%
\e{0}%
\e{0}%
\e{0}%
\e{0}%
\e{0}%
\e{0}%
\e{0}%
\e{0}%
\e{0}%
\e{0}%
\e{0}%
\e{0}%
\e{0}%
\e{0}%
\e{0}%
\eol}\vss}\rg%
%
%
\rx{\vss\hfull{%
\rlx{\hss{$50_{x}$}}\cg%
\e{0}%
\e{0}%
\e{0}%
\e{0}%
\e{0}%
\e{0}%
\e{0}%
\e{0}%
\e{0}%
\e{0}%
\e{0}%
\e{0}%
\e{0}%
\e{0}%
\e{0}%
\e{0}%
\e{0}%
\e{0}%
\e{0}%
\e{0}%
\e{0}%
\e{0}%
\e{0}%
\e{0}%
\e{0}%
\eol}\vss}\rg%
%
%
\rx{\vss\hfull{%
\rlx{\hss{$350_{x}$}}\cg%
\e{0}%
\e{0}%
\e{1}%
\e{1}%
\e{0}%
\e{0}%
\e{0}%
\e{0}%
\e{0}%
\e{0}%
\e{0}%
\e{0}%
\e{0}%
\e{0}%
\e{0}%
\e{0}%
\e{0}%
\e{0}%
\e{0}%
\e{0}%
\e{0}%
\e{0}%
\e{0}%
\e{0}%
\e{0}%
\eol}\vss}\rg%
%
%
\rx{\vss\hfull{%
\rlx{\hss{$300_{x}$}}\cg%
\e{0}%
\e{0}%
\e{1}%
\e{0}%
\e{0}%
\e{0}%
\e{0}%
\e{0}%
\e{0}%
\e{0}%
\e{0}%
\e{0}%
\e{0}%
\e{0}%
\e{0}%
\e{0}%
\e{0}%
\e{0}%
\e{0}%
\e{0}%
\e{0}%
\e{0}%
\e{0}%
\e{0}%
\e{0}%
\eol}\vss}\rg%
%
%
\rx{\vss\hfull{%
\rlx{\hss{$567_{x}$}}\cg%
\e{1}%
\e{1}%
\e{0}%
\e{1}%
\e{0}%
\e{0}%
\e{0}%
\e{0}%
\e{0}%
\e{0}%
\e{0}%
\e{0}%
\e{0}%
\e{0}%
\e{0}%
\e{0}%
\e{0}%
\e{0}%
\e{0}%
\e{0}%
\e{0}%
\e{0}%
\e{0}%
\e{0}%
\e{0}%
\eol}\vss}\rg%
%
%
\rx{\vss\hfull{%
\rlx{\hss{$210_{x}$}}\cg%
\e{0}%
\e{1}%
\e{0}%
\e{0}%
\e{0}%
\e{0}%
\e{0}%
\e{0}%
\e{0}%
\e{0}%
\e{0}%
\e{0}%
\e{0}%
\e{0}%
\e{0}%
\e{0}%
\e{0}%
\e{0}%
\e{0}%
\e{0}%
\e{0}%
\e{0}%
\e{0}%
\e{0}%
\e{0}%
\eol}\vss}\rg%
%
%
\rx{\vss\hfull{%
\rlx{\hss{$840_{x}$}}\cg%
\e{0}%
\e{0}%
\e{0}%
\e{0}%
\e{0}%
\e{0}%
\e{0}%
\e{0}%
\e{0}%
\e{0}%
\e{0}%
\e{0}%
\e{0}%
\e{0}%
\e{0}%
\e{0}%
\e{0}%
\e{0}%
\e{0}%
\e{0}%
\e{0}%
\e{1}%
\e{0}%
\e{0}%
\e{0}%
\eol}\vss}\rg%
%
%
\rx{\vss\hfull{%
\rlx{\hss{$700_{x}$}}\cg%
\e{0}%
\e{0}%
\e{0}%
\e{0}%
\e{1}%
\e{0}%
\e{0}%
\e{0}%
\e{0}%
\e{0}%
\e{0}%
\e{0}%
\e{0}%
\e{0}%
\e{0}%
\e{0}%
\e{0}%
\e{0}%
\e{0}%
\e{0}%
\e{0}%
\e{0}%
\e{0}%
\e{0}%
\e{0}%
\eol}\vss}\rg%
%
%
\rx{\vss\hfull{%
\rlx{\hss{$175_{x}$}}\cg%
\e{0}%
\e{0}%
\e{0}%
\e{0}%
\e{0}%
\e{0}%
\e{0}%
\e{0}%
\e{0}%
\e{0}%
\e{0}%
\e{0}%
\e{0}%
\e{0}%
\e{0}%
\e{0}%
\e{0}%
\e{0}%
\e{0}%
\e{0}%
\e{0}%
\e{0}%
\e{0}%
\e{0}%
\e{1}%
\eol}\vss}\rg%
%
%
\rx{\vss\hfull{%
\rlx{\hss{$1400_{x}$}}\cg%
\e{0}%
\e{0}%
\e{0}%
\e{1}%
\e{1}%
\e{0}%
\e{0}%
\e{0}%
\e{0}%
\e{1}%
\e{0}%
\e{0}%
\e{0}%
\e{0}%
\e{0}%
\e{0}%
\e{0}%
\e{0}%
\e{0}%
\e{0}%
\e{0}%
\e{0}%
\e{0}%
\e{0}%
\e{0}%
\eol}\vss}\rg%
%
%
\rx{\vss\hfull{%
\rlx{\hss{$1050_{x}$}}\cg%
\e{0}%
\e{0}%
\e{0}%
\e{0}%
\e{0}%
\e{0}%
\e{0}%
\e{1}%
\e{0}%
\e{1}%
\e{0}%
\e{0}%
\e{0}%
\e{0}%
\e{0}%
\e{0}%
\e{0}%
\e{0}%
\e{0}%
\e{0}%
\e{0}%
\e{0}%
\e{0}%
\e{0}%
\e{0}%
\eol}\vss}\rg%
%
%
\rx{\vss\hfull{%
\rlx{\hss{$1575_{x}$}}\cg%
\e{0}%
\e{1}%
\e{0}%
\e{1}%
\e{1}%
\e{1}%
\e{0}%
\e{1}%
\e{0}%
\e{0}%
\e{0}%
\e{0}%
\e{0}%
\e{0}%
\e{0}%
\e{0}%
\e{0}%
\e{0}%
\e{0}%
\e{0}%
\e{0}%
\e{0}%
\e{0}%
\e{0}%
\e{0}%
\eol}\vss}\rg%
%
%
\rx{\vss\hfull{%
\rlx{\hss{$1344_{x}$}}\cg%
\e{0}%
\e{0}%
\e{0}%
\e{1}%
\e{0}%
\e{1}%
\e{0}%
\e{0}%
\e{0}%
\e{0}%
\e{0}%
\e{0}%
\e{0}%
\e{0}%
\e{0}%
\e{0}%
\e{0}%
\e{0}%
\e{0}%
\e{0}%
\e{0}%
\e{0}%
\e{0}%
\e{0}%
\e{0}%
\eol}\vss}\rg%
%
%
\rx{\vss\hfull{%
\rlx{\hss{$2100_{x}$}}\cg%
\e{0}%
\e{1}%
\e{1}%
\e{1}%
\e{0}%
\e{2}%
\e{1}%
\e{0}%
\e{0}%
\e{0}%
\e{0}%
\e{0}%
\e{0}%
\e{0}%
\e{0}%
\e{0}%
\e{0}%
\e{0}%
\e{0}%
\e{0}%
\e{1}%
\e{0}%
\e{0}%
\e{0}%
\e{0}%
\eol}\vss}\rg%
%
%
\rx{\vss\hfull{%
\rlx{\hss{$2268_{x}$}}\cg%
\e{0}%
\e{1}%
\e{1}%
\e{1}%
\e{1}%
\e{1}%
\e{1}%
\e{0}%
\e{0}%
\e{0}%
\e{0}%
\e{0}%
\e{0}%
\e{0}%
\e{0}%
\e{0}%
\e{0}%
\e{0}%
\e{0}%
\e{0}%
\e{0}%
\e{0}%
\e{0}%
\e{0}%
\e{0}%
\eol}\vss}\rg%
%
%
\rx{\vss\hfull{%
\rlx{\hss{$525_{x}$}}\cg%
\e{1}%
\e{0}%
\e{0}%
\e{1}%
\e{0}%
\e{0}%
\e{0}%
\e{0}%
\e{1}%
\e{0}%
\e{0}%
\e{0}%
\e{0}%
\e{0}%
\e{0}%
\e{0}%
\e{0}%
\e{0}%
\e{0}%
\e{0}%
\e{0}%
\e{0}%
\e{0}%
\e{0}%
\e{0}%
\eol}\vss}\rg%
%
%
\rx{\vss\hfull{%
\rlx{\hss{$700_{xx}$}}\cg%
\e{0}%
\e{0}%
\e{0}%
\e{0}%
\e{0}%
\e{0}%
\e{0}%
\e{1}%
\e{0}%
\e{0}%
\e{0}%
\e{0}%
\e{0}%
\e{0}%
\e{0}%
\e{0}%
\e{0}%
\e{0}%
\e{0}%
\e{0}%
\e{0}%
\e{0}%
\e{0}%
\e{1}%
\e{0}%
\eol}\vss}\rg%
%
%
\rx{\vss\hfull{%
\rlx{\hss{$972_{x}$}}\cg%
\e{0}%
\e{0}%
\e{0}%
\e{0}%
\e{0}%
\e{0}%
\e{1}%
\e{0}%
\e{0}%
\e{0}%
\e{0}%
\e{0}%
\e{0}%
\e{0}%
\e{0}%
\e{0}%
\e{0}%
\e{0}%
\e{0}%
\e{0}%
\e{0}%
\e{0}%
\e{0}%
\e{0}%
\e{0}%
\eol}\vss}\rg%
\eop
\eject
\tablecont%
%
%
%
%
%
%
\rowpts=18 true pt%
\colpts=18 true pt%
\rowlabpts=40 true pt%
\collabpts=60 true pt%
\clx{\vss\hfull{%
\rlx{\hss{$ $}}\cg%
\cx{\hskip 16 true pt\flip{$1_p{\times}[{1^{3}}]$}\hss}\cg%
\cx{\hskip 16 true pt\flip{$6_p{\times}[{1^{3}}]$}\hss}\cg%
\cx{\hskip 16 true pt\flip{$15_p{\times}[{1^{3}}]$}\hss}\cg%
\cx{\hskip 16 true pt\flip{$20_p{\times}[{1^{3}}]$}\hss}\cg%
\cx{\hskip 16 true pt\flip{$30_p{\times}[{1^{3}}]$}\hss}\cg%
\cx{\hskip 16 true pt\flip{$64_p{\times}[{1^{3}}]$}\hss}\cg%
\cx{\hskip 16 true pt\flip{$81_p{\times}[{1^{3}}]$}\hss}\cg%
\cx{\hskip 16 true pt\flip{$15_q{\times}[{1^{3}}]$}\hss}\cg%
\cx{\hskip 16 true pt\flip{$24_p{\times}[{1^{3}}]$}\hss}\cg%
\cx{\hskip 16 true pt\flip{$60_p{\times}[{1^{3}}]$}\hss}\cg%
\cx{\hskip 16 true pt\flip{$1_p^{*}{\times}[{1^{3}}]$}\hss}\cg%
\cx{\hskip 16 true pt\flip{$6_p^{*}{\times}[{1^{3}}]$}\hss}\cg%
\cx{\hskip 16 true pt\flip{$15_p^{*}{\times}[{1^{3}}]$}\hss}\cg%
\cx{\hskip 16 true pt\flip{$20_p^{*}{\times}[{1^{3}}]$}\hss}\cg%
\cx{\hskip 16 true pt\flip{$30_p^{*}{\times}[{1^{3}}]$}\hss}\cg%
\cx{\hskip 16 true pt\flip{$64_p^{*}{\times}[{1^{3}}]$}\hss}\cg%
\cx{\hskip 16 true pt\flip{$81_p^{*}{\times}[{1^{3}}]$}\hss}\cg%
\cx{\hskip 16 true pt\flip{$15_q^{*}{\times}[{1^{3}}]$}\hss}\cg%
\cx{\hskip 16 true pt\flip{$24_p^{*}{\times}[{1^{3}}]$}\hss}\cg%
\cx{\hskip 16 true pt\flip{$60_p^{*}{\times}[{1^{3}}]$}\hss}\cg%
\cx{\hskip 16 true pt\flip{$20_s{\times}[{1^{3}}]$}\hss}\cg%
\cx{\hskip 16 true pt\flip{$90_s{\times}[{1^{3}}]$}\hss}\cg%
\cx{\hskip 16 true pt\flip{$80_s{\times}[{1^{3}}]$}\hss}\cg%
\cx{\hskip 16 true pt\flip{$60_s{\times}[{1^{3}}]$}\hss}\cg%
\cx{\hskip 16 true pt\flip{$10_s{\times}[{1^{3}}]$}\hss}\cg%
\eol}}\rg%
%
%
\rx{\vss\hfull{%
\rlx{\hss{$4096_{x}$}}\cg%
\e{0}%
\e{0}%
\e{1}%
\e{1}%
\e{1}%
\e{2}%
\e{1}%
\e{0}%
\e{1}%
\e{1}%
\e{0}%
\e{0}%
\e{0}%
\e{0}%
\e{0}%
\e{0}%
\e{0}%
\e{0}%
\e{0}%
\e{0}%
\e{0}%
\e{1}%
\e{0}%
\e{0}%
\e{0}%
\eol}\vss}\rg%
%
%
\rx{\vss\hfull{%
\rlx{\hss{$4200_{x}$}}\cg%
\e{0}%
\e{0}%
\e{0}%
\e{0}%
\e{1}%
\e{1}%
\e{1}%
\e{1}%
\e{0}%
\e{1}%
\e{0}%
\e{0}%
\e{0}%
\e{0}%
\e{0}%
\e{0}%
\e{0}%
\e{0}%
\e{0}%
\e{0}%
\e{0}%
\e{1}%
\e{1}%
\e{1}%
\e{0}%
\eol}\vss}\rg%
%
%
\rx{\vss\hfull{%
\rlx{\hss{$2240_{x}$}}\cg%
\e{0}%
\e{0}%
\e{0}%
\e{0}%
\e{0}%
\e{1}%
\e{0}%
\e{0}%
\e{0}%
\e{1}%
\e{0}%
\e{0}%
\e{0}%
\e{0}%
\e{0}%
\e{0}%
\e{0}%
\e{0}%
\e{0}%
\e{0}%
\e{0}%
\e{0}%
\e{1}%
\e{0}%
\e{0}%
\eol}\vss}\rg%
%
%
\rx{\vss\hfull{%
\rlx{\hss{$2835_{x}$}}\cg%
\e{0}%
\e{0}%
\e{0}%
\e{0}%
\e{0}%
\e{0}%
\e{1}%
\e{0}%
\e{0}%
\e{1}%
\e{0}%
\e{0}%
\e{0}%
\e{0}%
\e{0}%
\e{0}%
\e{0}%
\e{0}%
\e{0}%
\e{1}%
\e{0}%
\e{0}%
\e{1}%
\e{1}%
\e{1}%
\eol}\vss}\rg%
%
%
\rx{\vss\hfull{%
\rlx{\hss{$6075_{x}$}}\cg%
\e{0}%
\e{0}%
\e{0}%
\e{1}%
\e{1}%
\e{2}%
\e{2}%
\e{1}%
\e{1}%
\e{2}%
\e{0}%
\e{0}%
\e{0}%
\e{0}%
\e{0}%
\e{0}%
\e{1}%
\e{0}%
\e{0}%
\e{0}%
\e{0}%
\e{1}%
\e{1}%
\e{1}%
\e{0}%
\eol}\vss}\rg%
%
%
\rx{\vss\hfull{%
\rlx{\hss{$3200_{x}$}}\cg%
\e{0}%
\e{0}%
\e{0}%
\e{0}%
\e{0}%
\e{1}%
\e{1}%
\e{0}%
\e{1}%
\e{0}%
\e{0}%
\e{0}%
\e{0}%
\e{0}%
\e{0}%
\e{0}%
\e{1}%
\e{0}%
\e{0}%
\e{0}%
\e{1}%
\e{1}%
\e{1}%
\e{0}%
\e{0}%
\eol}\vss}\rg%
%
%
\rx{\vss\hfull{%
\rlx{\hss{$70_{y}$}}\cg%
\e{0}%
\e{0}%
\e{1}%
\e{0}%
\e{0}%
\e{0}%
\e{0}%
\e{0}%
\e{0}%
\e{0}%
\e{0}%
\e{0}%
\e{0}%
\e{0}%
\e{0}%
\e{0}%
\e{0}%
\e{0}%
\e{0}%
\e{0}%
\e{0}%
\e{0}%
\e{0}%
\e{0}%
\e{0}%
\eol}\vss}\rg%
%
%
\rx{\vss\hfull{%
\rlx{\hss{$1134_{y}$}}\cg%
\e{0}%
\e{0}%
\e{0}%
\e{0}%
\e{0}%
\e{0}%
\e{0}%
\e{0}%
\e{1}%
\e{1}%
\e{0}%
\e{0}%
\e{1}%
\e{0}%
\e{0}%
\e{0}%
\e{0}%
\e{0}%
\e{0}%
\e{0}%
\e{0}%
\e{1}%
\e{0}%
\e{0}%
\e{0}%
\eol}\vss}\rg%
%
%
\rx{\vss\hfull{%
\rlx{\hss{$1680_{y}$}}\cg%
\e{0}%
\e{0}%
\e{1}%
\e{0}%
\e{1}%
\e{1}%
\e{1}%
\e{0}%
\e{0}%
\e{0}%
\e{0}%
\e{0}%
\e{0}%
\e{0}%
\e{0}%
\e{0}%
\e{0}%
\e{0}%
\e{0}%
\e{0}%
\e{1}%
\e{1}%
\e{0}%
\e{0}%
\e{0}%
\eol}\vss}\rg%
%
%
\rx{\vss\hfull{%
\rlx{\hss{$168_{y}$}}\cg%
\e{0}%
\e{0}%
\e{0}%
\e{0}%
\e{0}%
\e{0}%
\e{0}%
\e{0}%
\e{0}%
\e{0}%
\e{0}%
\e{0}%
\e{0}%
\e{0}%
\e{0}%
\e{0}%
\e{0}%
\e{0}%
\e{1}%
\e{0}%
\e{0}%
\e{0}%
\e{0}%
\e{0}%
\e{0}%
\eol}\vss}\rg%
%
%
\rx{\vss\hfull{%
\rlx{\hss{$420_{y}$}}\cg%
\e{0}%
\e{0}%
\e{0}%
\e{0}%
\e{0}%
\e{0}%
\e{0}%
\e{0}%
\e{0}%
\e{0}%
\e{0}%
\e{0}%
\e{0}%
\e{0}%
\e{0}%
\e{0}%
\e{0}%
\e{0}%
\e{0}%
\e{1}%
\e{0}%
\e{0}%
\e{0}%
\e{0}%
\e{0}%
\eol}\vss}\rg%
%
%
\rx{\vss\hfull{%
\rlx{\hss{$3150_{y}$}}\cg%
\e{0}%
\e{0}%
\e{0}%
\e{0}%
\e{0}%
\e{0}%
\e{0}%
\e{0}%
\e{1}%
\e{1}%
\e{0}%
\e{0}%
\e{0}%
\e{0}%
\e{1}%
\e{0}%
\e{1}%
\e{0}%
\e{0}%
\e{1}%
\e{0}%
\e{1}%
\e{1}%
\e{1}%
\e{1}%
\eol}\vss}\rg%
%
%
\rx{\vss\hfull{%
\rlx{\hss{$4200_{y}$}}\cg%
\e{0}%
\e{0}%
\e{0}%
\e{0}%
\e{0}%
\e{0}%
\e{1}%
\e{0}%
\e{0}%
\e{1}%
\e{0}%
\e{0}%
\e{0}%
\e{0}%
\e{0}%
\e{1}%
\e{1}%
\e{0}%
\e{1}%
\e{1}%
\e{0}%
\e{1}%
\e{1}%
\e{2}%
\e{0}%
\eol}\vss}\rg%
%
%
\rx{\vss\hfull{%
\rlx{\hss{$2688_{y}$}}\cg%
\e{0}%
\e{0}%
\e{0}%
\e{0}%
\e{0}%
\e{0}%
\e{1}%
\e{0}%
\e{0}%
\e{1}%
\e{0}%
\e{0}%
\e{0}%
\e{0}%
\e{0}%
\e{0}%
\e{1}%
\e{0}%
\e{0}%
\e{1}%
\e{0}%
\e{1}%
\e{0}%
\e{1}%
\e{0}%
\eol}\vss}\rg%
%
%
\rx{\vss\hfull{%
\rlx{\hss{$2100_{y}$}}\cg%
\e{0}%
\e{0}%
\e{0}%
\e{1}%
\e{0}%
\e{1}%
\e{1}%
\e{0}%
\e{1}%
\e{0}%
\e{0}%
\e{0}%
\e{0}%
\e{0}%
\e{0}%
\e{0}%
\e{1}%
\e{0}%
\e{0}%
\e{0}%
\e{0}%
\e{1}%
\e{0}%
\e{0}%
\e{0}%
\eol}\vss}\rg%
%
%
\rx{\vss\hfull{%
\rlx{\hss{$1400_{y}$}}\cg%
\e{0}%
\e{0}%
\e{1}%
\e{0}%
\e{1}%
\e{0}%
\e{1}%
\e{0}%
\e{0}%
\e{0}%
\e{0}%
\e{0}%
\e{0}%
\e{0}%
\e{0}%
\e{0}%
\e{0}%
\e{0}%
\e{1}%
\e{0}%
\e{0}%
\e{1}%
\e{0}%
\e{0}%
\e{0}%
\eol}\vss}\rg%
%
%
\rx{\vss\hfull{%
\rlx{\hss{$4536_{y}$}}\cg%
\e{0}%
\e{0}%
\e{1}%
\e{0}%
\e{1}%
\e{1}%
\e{2}%
\e{0}%
\e{0}%
\e{1}%
\e{0}%
\e{0}%
\e{0}%
\e{0}%
\e{0}%
\e{0}%
\e{1}%
\e{0}%
\e{1}%
\e{1}%
\e{0}%
\e{2}%
\e{1}%
\e{0}%
\e{0}%
\eol}\vss}\rg%
%
%
\rx{\vss\hfull{%
\rlx{\hss{$5670_{y}$}}\cg%
\e{0}%
\e{0}%
\e{0}%
\e{0}%
\e{1}%
\e{2}%
\e{2}%
\e{0}%
\e{0}%
\e{1}%
\e{0}%
\e{0}%
\e{0}%
\e{0}%
\e{0}%
\e{1}%
\e{1}%
\e{0}%
\e{0}%
\e{0}%
\e{1}%
\e{2}%
\e{2}%
\e{1}%
\e{0}%
\eol}\vss}\rg%
%
%
\rx{\vss\hfull{%
\rlx{\hss{$4480_{y}$}}\cg%
\e{0}%
\e{0}%
\e{0}%
\e{0}%
\e{0}%
\e{1}%
\e{1}%
\e{0}%
\e{0}%
\e{1}%
\e{0}%
\e{0}%
\e{0}%
\e{0}%
\e{0}%
\e{1}%
\e{1}%
\e{0}%
\e{0}%
\e{1}%
\e{0}%
\e{1}%
\e{2}%
\e{1}%
\e{0}%
\eol}\vss}\rg%
%
%
\rx{\vss\hfull{%
\rlx{\hss{$8_{z}$}}\cg%
\e{0}%
\e{0}%
\e{0}%
\e{0}%
\e{0}%
\e{0}%
\e{0}%
\e{0}%
\e{0}%
\e{0}%
\e{0}%
\e{0}%
\e{0}%
\e{0}%
\e{0}%
\e{0}%
\e{0}%
\e{0}%
\e{0}%
\e{0}%
\e{0}%
\e{0}%
\e{0}%
\e{0}%
\e{0}%
\eol}\vss}\rg%
%
%
\rx{\vss\hfull{%
\rlx{\hss{$56_{z}$}}\cg%
\e{0}%
\e{1}%
\e{0}%
\e{0}%
\e{0}%
\e{0}%
\e{0}%
\e{0}%
\e{0}%
\e{0}%
\e{0}%
\e{0}%
\e{0}%
\e{0}%
\e{0}%
\e{0}%
\e{0}%
\e{0}%
\e{0}%
\e{0}%
\e{0}%
\e{0}%
\e{0}%
\e{0}%
\e{0}%
\eol}\vss}\rg%
%
%
\rx{\vss\hfull{%
\rlx{\hss{$160_{z}$}}\cg%
\e{0}%
\e{1}%
\e{0}%
\e{0}%
\e{0}%
\e{0}%
\e{0}%
\e{0}%
\e{0}%
\e{0}%
\e{0}%
\e{0}%
\e{0}%
\e{0}%
\e{0}%
\e{0}%
\e{0}%
\e{0}%
\e{0}%
\e{0}%
\e{0}%
\e{0}%
\e{0}%
\e{0}%
\e{0}%
\eol}\vss}\rg%
%
%
\rx{\vss\hfull{%
\rlx{\hss{$112_{z}$}}\cg%
\e{1}%
\e{0}%
\e{0}%
\e{0}%
\e{0}%
\e{0}%
\e{0}%
\e{0}%
\e{0}%
\e{0}%
\e{0}%
\e{0}%
\e{0}%
\e{0}%
\e{0}%
\e{0}%
\e{0}%
\e{0}%
\e{0}%
\e{0}%
\e{0}%
\e{0}%
\e{0}%
\e{0}%
\e{0}%
\eol}\vss}\rg%
%
%
\rx{\vss\hfull{%
\rlx{\hss{$840_{z}$}}\cg%
\e{0}%
\e{0}%
\e{0}%
\e{0}%
\e{0}%
\e{1}%
\e{0}%
\e{0}%
\e{0}%
\e{0}%
\e{0}%
\e{0}%
\e{0}%
\e{0}%
\e{0}%
\e{0}%
\e{0}%
\e{0}%
\e{0}%
\e{0}%
\e{1}%
\e{0}%
\e{0}%
\e{0}%
\e{0}%
\eol}\vss}\rg%
%
%
\rx{\vss\hfull{%
\rlx{\hss{$1296_{z}$}}\cg%
\e{0}%
\e{1}%
\e{1}%
\e{1}%
\e{1}%
\e{1}%
\e{0}%
\e{0}%
\e{0}%
\e{0}%
\e{0}%
\e{0}%
\e{0}%
\e{0}%
\e{0}%
\e{0}%
\e{0}%
\e{0}%
\e{0}%
\e{0}%
\e{0}%
\e{0}%
\e{0}%
\e{0}%
\e{0}%
\eol}\vss}\rg%
%
%
\rx{\vss\hfull{%
\rlx{\hss{$1400_{z}$}}\cg%
\e{0}%
\e{1}%
\e{0}%
\e{1}%
\e{0}%
\e{1}%
\e{0}%
\e{0}%
\e{0}%
\e{0}%
\e{0}%
\e{0}%
\e{0}%
\e{0}%
\e{0}%
\e{0}%
\e{0}%
\e{0}%
\e{0}%
\e{0}%
\e{0}%
\e{0}%
\e{0}%
\e{0}%
\e{0}%
\eol}\vss}\rg%
%
%
\rx{\vss\hfull{%
\rlx{\hss{$1008_{z}$}}\cg%
\e{1}%
\e{1}%
\e{1}%
\e{1}%
\e{1}%
\e{0}%
\e{0}%
\e{0}%
\e{0}%
\e{0}%
\e{0}%
\e{0}%
\e{0}%
\e{0}%
\e{0}%
\e{0}%
\e{0}%
\e{0}%
\e{0}%
\e{0}%
\e{0}%
\e{0}%
\e{0}%
\e{0}%
\e{0}%
\eol}\vss}\rg%
%
%
\rx{\vss\hfull{%
\rlx{\hss{$560_{z}$}}\cg%
\e{0}%
\e{0}%
\e{0}%
\e{1}%
\e{0}%
\e{0}%
\e{0}%
\e{0}%
\e{0}%
\e{0}%
\e{0}%
\e{0}%
\e{0}%
\e{0}%
\e{0}%
\e{0}%
\e{0}%
\e{0}%
\e{0}%
\e{0}%
\e{0}%
\e{0}%
\e{0}%
\e{0}%
\e{0}%
\eol}\vss}\rg%
%
%
\rx{\vss\hfull{%
\rlx{\hss{$1400_{zz}$}}\cg%
\e{0}%
\e{0}%
\e{0}%
\e{0}%
\e{0}%
\e{0}%
\e{0}%
\e{0}%
\e{0}%
\e{1}%
\e{0}%
\e{0}%
\e{0}%
\e{0}%
\e{0}%
\e{0}%
\e{0}%
\e{0}%
\e{0}%
\e{0}%
\e{0}%
\e{0}%
\e{0}%
\e{1}%
\e{0}%
\eol}\vss}\rg%
%
%
\rx{\vss\hfull{%
\rlx{\hss{$4200_{z}$}}\cg%
\e{0}%
\e{0}%
\e{0}%
\e{0}%
\e{1}%
\e{0}%
\e{1}%
\e{1}%
\e{0}%
\e{2}%
\e{0}%
\e{0}%
\e{0}%
\e{0}%
\e{0}%
\e{0}%
\e{0}%
\e{0}%
\e{1}%
\e{1}%
\e{0}%
\e{1}%
\e{1}%
\e{1}%
\e{1}%
\eol}\vss}\rg%
%
%
\rx{\vss\hfull{%
\rlx{\hss{$400_{z}$}}\cg%
\e{0}%
\e{0}%
\e{0}%
\e{0}%
\e{0}%
\e{0}%
\e{0}%
\e{1}%
\e{0}%
\e{0}%
\e{0}%
\e{0}%
\e{0}%
\e{0}%
\e{0}%
\e{0}%
\e{0}%
\e{0}%
\e{0}%
\e{0}%
\e{0}%
\e{0}%
\e{0}%
\e{0}%
\e{0}%
\eol}\vss}\rg%
%
%
\rx{\vss\hfull{%
\rlx{\hss{$3240_{z}$}}\cg%
\e{0}%
\e{0}%
\e{0}%
\e{1}%
\e{1}%
\e{1}%
\e{1}%
\e{1}%
\e{0}%
\e{1}%
\e{0}%
\e{0}%
\e{0}%
\e{0}%
\e{0}%
\e{0}%
\e{0}%
\e{0}%
\e{0}%
\e{0}%
\e{0}%
\e{0}%
\e{0}%
\e{0}%
\e{0}%
\eol}\vss}\rg%
%
%
\rx{\vss\hfull{%
\rlx{\hss{$4536_{z}$}}\cg%
\e{0}%
\e{0}%
\e{0}%
\e{0}%
\e{0}%
\e{1}%
\e{1}%
\e{0}%
\e{1}%
\e{1}%
\e{0}%
\e{0}%
\e{0}%
\e{0}%
\e{0}%
\e{0}%
\e{1}%
\e{0}%
\e{0}%
\e{0}%
\e{0}%
\e{1}%
\e{1}%
\e{1}%
\e{0}%
\eol}\vss}\rg%
%
%
\rx{\vss\hfull{%
\rlx{\hss{$2400_{z}$}}\cg%
\e{0}%
\e{0}%
\e{1}%
\e{1}%
\e{1}%
\e{1}%
\e{1}%
\e{0}%
\e{0}%
\e{1}%
\e{0}%
\e{0}%
\e{0}%
\e{0}%
\e{0}%
\e{0}%
\e{0}%
\e{0}%
\e{0}%
\e{0}%
\e{0}%
\e{1}%
\e{0}%
\e{0}%
\e{0}%
\eol}\vss}\rg%
\eop
\eject
\tablecont%
%
%
%
%
%
%
\rowpts=18 true pt%
\colpts=18 true pt%
\rowlabpts=40 true pt%
\collabpts=60 true pt%
\clx{\vss\hfull{%
\rlx{\hss{$ $}}\cg%
\cx{\hskip 16 true pt\flip{$1_p{\times}[{1^{3}}]$}\hss}\cg%
\cx{\hskip 16 true pt\flip{$6_p{\times}[{1^{3}}]$}\hss}\cg%
\cx{\hskip 16 true pt\flip{$15_p{\times}[{1^{3}}]$}\hss}\cg%
\cx{\hskip 16 true pt\flip{$20_p{\times}[{1^{3}}]$}\hss}\cg%
\cx{\hskip 16 true pt\flip{$30_p{\times}[{1^{3}}]$}\hss}\cg%
\cx{\hskip 16 true pt\flip{$64_p{\times}[{1^{3}}]$}\hss}\cg%
\cx{\hskip 16 true pt\flip{$81_p{\times}[{1^{3}}]$}\hss}\cg%
\cx{\hskip 16 true pt\flip{$15_q{\times}[{1^{3}}]$}\hss}\cg%
\cx{\hskip 16 true pt\flip{$24_p{\times}[{1^{3}}]$}\hss}\cg%
\cx{\hskip 16 true pt\flip{$60_p{\times}[{1^{3}}]$}\hss}\cg%
\cx{\hskip 16 true pt\flip{$1_p^{*}{\times}[{1^{3}}]$}\hss}\cg%
\cx{\hskip 16 true pt\flip{$6_p^{*}{\times}[{1^{3}}]$}\hss}\cg%
\cx{\hskip 16 true pt\flip{$15_p^{*}{\times}[{1^{3}}]$}\hss}\cg%
\cx{\hskip 16 true pt\flip{$20_p^{*}{\times}[{1^{3}}]$}\hss}\cg%
\cx{\hskip 16 true pt\flip{$30_p^{*}{\times}[{1^{3}}]$}\hss}\cg%
\cx{\hskip 16 true pt\flip{$64_p^{*}{\times}[{1^{3}}]$}\hss}\cg%
\cx{\hskip 16 true pt\flip{$81_p^{*}{\times}[{1^{3}}]$}\hss}\cg%
\cx{\hskip 16 true pt\flip{$15_q^{*}{\times}[{1^{3}}]$}\hss}\cg%
\cx{\hskip 16 true pt\flip{$24_p^{*}{\times}[{1^{3}}]$}\hss}\cg%
\cx{\hskip 16 true pt\flip{$60_p^{*}{\times}[{1^{3}}]$}\hss}\cg%
\cx{\hskip 16 true pt\flip{$20_s{\times}[{1^{3}}]$}\hss}\cg%
\cx{\hskip 16 true pt\flip{$90_s{\times}[{1^{3}}]$}\hss}\cg%
\cx{\hskip 16 true pt\flip{$80_s{\times}[{1^{3}}]$}\hss}\cg%
\cx{\hskip 16 true pt\flip{$60_s{\times}[{1^{3}}]$}\hss}\cg%
\cx{\hskip 16 true pt\flip{$10_s{\times}[{1^{3}}]$}\hss}\cg%
\eol}}\rg%
%
%
\rx{\vss\hfull{%
\rlx{\hss{$3360_{z}$}}\cg%
\e{0}%
\e{0}%
\e{0}%
\e{0}%
\e{1}%
\e{1}%
\e{1}%
\e{1}%
\e{0}%
\e{1}%
\e{0}%
\e{0}%
\e{0}%
\e{0}%
\e{0}%
\e{0}%
\e{0}%
\e{0}%
\e{0}%
\e{0}%
\e{0}%
\e{0}%
\e{1}%
\e{1}%
\e{0}%
\eol}\vss}\rg%
%
%
\rx{\vss\hfull{%
\rlx{\hss{$2800_{z}$}}\cg%
\e{0}%
\e{1}%
\e{0}%
\e{1}%
\e{1}%
\e{2}%
\e{1}%
\e{0}%
\e{0}%
\e{0}%
\e{0}%
\e{0}%
\e{0}%
\e{0}%
\e{0}%
\e{0}%
\e{0}%
\e{0}%
\e{0}%
\e{0}%
\e{0}%
\e{0}%
\e{1}%
\e{0}%
\e{0}%
\eol}\vss}\rg%
%
%
\rx{\vss\hfull{%
\rlx{\hss{$4096_{z}$}}\cg%
\e{0}%
\e{0}%
\e{1}%
\e{1}%
\e{1}%
\e{2}%
\e{1}%
\e{0}%
\e{1}%
\e{1}%
\e{0}%
\e{0}%
\e{0}%
\e{0}%
\e{0}%
\e{0}%
\e{0}%
\e{0}%
\e{0}%
\e{0}%
\e{0}%
\e{1}%
\e{0}%
\e{0}%
\e{0}%
\eol}\vss}\rg%
%
%
\rx{\vss\hfull{%
\rlx{\hss{$5600_{z}$}}\cg%
\e{0}%
\e{0}%
\e{1}%
\e{1}%
\e{1}%
\e{2}%
\e{2}%
\e{0}%
\e{1}%
\e{1}%
\e{0}%
\e{0}%
\e{0}%
\e{0}%
\e{0}%
\e{0}%
\e{1}%
\e{0}%
\e{0}%
\e{0}%
\e{1}%
\e{2}%
\e{1}%
\e{0}%
\e{0}%
\eol}\vss}\rg%
%
%
\rx{\vss\hfull{%
\rlx{\hss{$448_{z}$}}\cg%
\e{0}%
\e{0}%
\e{0}%
\e{0}%
\e{0}%
\e{0}%
\e{0}%
\e{0}%
\e{1}%
\e{0}%
\e{0}%
\e{0}%
\e{0}%
\e{0}%
\e{0}%
\e{0}%
\e{0}%
\e{0}%
\e{0}%
\e{0}%
\e{0}%
\e{0}%
\e{0}%
\e{0}%
\e{0}%
\eol}\vss}\rg%
%
%
\rx{\vss\hfull{%
\rlx{\hss{$448_{w}$}}\cg%
\e{0}%
\e{0}%
\e{0}%
\e{0}%
\e{0}%
\e{1}%
\e{0}%
\e{0}%
\e{0}%
\e{0}%
\e{0}%
\e{0}%
\e{0}%
\e{0}%
\e{0}%
\e{0}%
\e{0}%
\e{0}%
\e{0}%
\e{0}%
\e{1}%
\e{0}%
\e{0}%
\e{0}%
\e{0}%
\eol}\vss}\rg%
%
%
\rx{\vss\hfull{%
\rlx{\hss{$1344_{w}$}}\cg%
\e{0}%
\e{0}%
\e{0}%
\e{0}%
\e{0}%
\e{0}%
\e{0}%
\e{0}%
\e{0}%
\e{0}%
\e{0}%
\e{0}%
\e{0}%
\e{0}%
\e{0}%
\e{1}%
\e{0}%
\e{0}%
\e{0}%
\e{0}%
\e{0}%
\e{0}%
\e{1}%
\e{1}%
\e{0}%
\eol}\vss}\rg%
%
%
\rx{\vss\hfull{%
\rlx{\hss{$5600_{w}$}}\cg%
\e{0}%
\e{0}%
\e{0}%
\e{0}%
\e{0}%
\e{2}%
\e{2}%
\e{1}%
\e{0}%
\e{1}%
\e{0}%
\e{0}%
\e{0}%
\e{0}%
\e{0}%
\e{1}%
\e{1}%
\e{0}%
\e{0}%
\e{0}%
\e{1}%
\e{2}%
\e{2}%
\e{1}%
\e{0}%
\eol}\vss}\rg%
%
%
\rx{\vss\hfull{%
\rlx{\hss{$2016_{w}$}}\cg%
\e{0}%
\e{0}%
\e{0}%
\e{0}%
\e{0}%
\e{0}%
\e{0}%
\e{0}%
\e{0}%
\e{0}%
\e{0}%
\e{0}%
\e{0}%
\e{0}%
\e{0}%
\e{0}%
\e{1}%
\e{1}%
\e{0}%
\e{1}%
\e{0}%
\e{0}%
\e{1}%
\e{1}%
\e{1}%
\eol}\vss}\rg%
%
%
\rx{\vss\hfull{%
\rlx{\hss{$7168_{w}$}}\cg%
\e{0}%
\e{0}%
\e{0}%
\e{0}%
\e{0}%
\e{1}%
\e{2}%
\e{0}%
\e{0}%
\e{2}%
\e{0}%
\e{0}%
\e{0}%
\e{0}%
\e{0}%
\e{1}%
\e{2}%
\e{0}%
\e{0}%
\e{2}%
\e{0}%
\e{2}%
\e{2}%
\e{2}%
\e{0}%
\eol}\vss}\rg%
\tableclose%
%
%
%
%
%
%
\tableopen{Induce/restrict matrix for $W({E_{7}}{A_{1}})\,\subset\,W(E_{8})$}%
%
%
%
%
%
%
\rowpts=18 true pt%
\colpts=18 true pt%
\rowlabpts=40 true pt%
\collabpts=70 true pt%
\clx{\vss\hfull{%
\rlx{\hss{$ $}}\cg%
\cx{\hskip 16 true pt\flip{$1_a{\times}[{2}]$}\hss}\cg%
\cx{\hskip 16 true pt\flip{$7_a{\times}[{2}]$}\hss}\cg%
\cx{\hskip 16 true pt\flip{$27_a{\times}[{2}]$}\hss}\cg%
\cx{\hskip 16 true pt\flip{$21_a{\times}[{2}]$}\hss}\cg%
\cx{\hskip 16 true pt\flip{$35_a{\times}[{2}]$}\hss}\cg%
\cx{\hskip 16 true pt\flip{$105_a{\times}[{2}]$}\hss}\cg%
\cx{\hskip 16 true pt\flip{$189_a{\times}[{2}]$}\hss}\cg%
\cx{\hskip 16 true pt\flip{$21_b{\times}[{2}]$}\hss}\cg%
\cx{\hskip 16 true pt\flip{$35_b{\times}[{2}]$}\hss}\cg%
\cx{\hskip 16 true pt\flip{$189_b{\times}[{2}]$}\hss}\cg%
\cx{\hskip 16 true pt\flip{$189_c{\times}[{2}]$}\hss}\cg%
\cx{\hskip 16 true pt\flip{$15_a{\times}[{2}]$}\hss}\cg%
\cx{\hskip 16 true pt\flip{$105_b{\times}[{2}]$}\hss}\cg%
\cx{\hskip 16 true pt\flip{$105_c{\times}[{2}]$}\hss}\cg%
\cx{\hskip 16 true pt\flip{$315_a{\times}[{2}]$}\hss}\cg%
\cx{\hskip 16 true pt\flip{$405_a{\times}[{2}]$}\hss}\cg%
\cx{\hskip 16 true pt\flip{$168_a{\times}[{2}]$}\hss}\cg%
\cx{\hskip 16 true pt\flip{$56_a{\times}[{2}]$}\hss}\cg%
\cx{\hskip 16 true pt\flip{$120_a{\times}[{2}]$}\hss}\cg%
\cx{\hskip 16 true pt\flip{$210_a{\times}[{2}]$}\hss}\cg%
\cx{\hskip 16 true pt\flip{$280_a{\times}[{2}]$}\hss}\cg%
\cx{\hskip 16 true pt\flip{$336_a{\times}[{2}]$}\hss}\cg%
\cx{\hskip 16 true pt\flip{$216_a{\times}[{2}]$}\hss}\cg%
\cx{\hskip 16 true pt\flip{$512_a{\times}[{2}]$}\hss}\cg%
\eol}}\rg%
%
%
\rx{\vss\hfull{%
\rlx{\hss{$1_{x}$}}\cg%
\e{1}%
\e{0}%
\e{0}%
\e{0}%
\e{0}%
\e{0}%
\e{0}%
\e{0}%
\e{0}%
\e{0}%
\e{0}%
\e{0}%
\e{0}%
\e{0}%
\e{0}%
\e{0}%
\e{0}%
\e{0}%
\e{0}%
\e{0}%
\e{0}%
\e{0}%
\e{0}%
\e{0}%
\eol}\vss}\rg%
%
%
\rx{\vss\hfull{%
\rlx{\hss{$28_{x}$}}\cg%
\e{0}%
\e{0}%
\e{0}%
\e{1}%
\e{0}%
\e{0}%
\e{0}%
\e{0}%
\e{0}%
\e{0}%
\e{0}%
\e{0}%
\e{0}%
\e{0}%
\e{0}%
\e{0}%
\e{0}%
\e{0}%
\e{0}%
\e{0}%
\e{0}%
\e{0}%
\e{0}%
\e{0}%
\eol}\vss}\rg%
%
%
\rx{\vss\hfull{%
\rlx{\hss{$35_{x}$}}\cg%
\e{1}%
\e{0}%
\e{1}%
\e{0}%
\e{0}%
\e{0}%
\e{0}%
\e{0}%
\e{0}%
\e{0}%
\e{0}%
\e{0}%
\e{0}%
\e{0}%
\e{0}%
\e{0}%
\e{0}%
\e{0}%
\e{0}%
\e{0}%
\e{0}%
\e{0}%
\e{0}%
\e{0}%
\eol}\vss}\rg%
%
%
\rx{\vss\hfull{%
\rlx{\hss{$84_{x}$}}\cg%
\e{1}%
\e{0}%
\e{1}%
\e{0}%
\e{0}%
\e{0}%
\e{0}%
\e{0}%
\e{1}%
\e{0}%
\e{0}%
\e{0}%
\e{0}%
\e{0}%
\e{0}%
\e{0}%
\e{0}%
\e{0}%
\e{0}%
\e{0}%
\e{0}%
\e{0}%
\e{0}%
\e{0}%
\eol}\vss}\rg%
%
%
\rx{\vss\hfull{%
\rlx{\hss{$50_{x}$}}\cg%
\e{0}%
\e{0}%
\e{0}%
\e{0}%
\e{0}%
\e{0}%
\e{0}%
\e{0}%
\e{1}%
\e{0}%
\e{0}%
\e{0}%
\e{0}%
\e{0}%
\e{0}%
\e{0}%
\e{0}%
\e{0}%
\e{0}%
\e{0}%
\e{0}%
\e{0}%
\e{0}%
\e{0}%
\eol}\vss}\rg%
%
%
\rx{\vss\hfull{%
\rlx{\hss{$350_{x}$}}\cg%
\e{0}%
\e{0}%
\e{0}%
\e{1}%
\e{0}%
\e{0}%
\e{1}%
\e{0}%
\e{0}%
\e{0}%
\e{0}%
\e{0}%
\e{0}%
\e{0}%
\e{0}%
\e{0}%
\e{0}%
\e{0}%
\e{0}%
\e{0}%
\e{0}%
\e{0}%
\e{0}%
\e{0}%
\eol}\vss}\rg%
%
%
\rx{\vss\hfull{%
\rlx{\hss{$300_{x}$}}\cg%
\e{0}%
\e{0}%
\e{1}%
\e{0}%
\e{0}%
\e{0}%
\e{0}%
\e{0}%
\e{0}%
\e{0}%
\e{0}%
\e{0}%
\e{0}%
\e{0}%
\e{0}%
\e{0}%
\e{1}%
\e{0}%
\e{0}%
\e{0}%
\e{0}%
\e{0}%
\e{0}%
\e{0}%
\eol}\vss}\rg%
%
%
\rx{\vss\hfull{%
\rlx{\hss{$567_{x}$}}\cg%
\e{0}%
\e{0}%
\e{1}%
\e{1}%
\e{0}%
\e{0}%
\e{0}%
\e{0}%
\e{0}%
\e{0}%
\e{0}%
\e{0}%
\e{0}%
\e{0}%
\e{0}%
\e{0}%
\e{0}%
\e{0}%
\e{1}%
\e{1}%
\e{0}%
\e{0}%
\e{0}%
\e{0}%
\eol}\vss}\rg%
%
%
\rx{\vss\hfull{%
\rlx{\hss{$210_{x}$}}\cg%
\e{0}%
\e{0}%
\e{1}%
\e{0}%
\e{0}%
\e{0}%
\e{0}%
\e{0}%
\e{0}%
\e{0}%
\e{0}%
\e{0}%
\e{0}%
\e{0}%
\e{0}%
\e{0}%
\e{0}%
\e{0}%
\e{1}%
\e{0}%
\e{0}%
\e{0}%
\e{0}%
\e{0}%
\eol}\vss}\rg%
%
%
\rx{\vss\hfull{%
\rlx{\hss{$840_{x}$}}\cg%
\e{0}%
\e{0}%
\e{0}%
\e{0}%
\e{0}%
\e{0}%
\e{0}%
\e{0}%
\e{0}%
\e{0}%
\e{0}%
\e{0}%
\e{0}%
\e{0}%
\e{0}%
\e{0}%
\e{1}%
\e{0}%
\e{0}%
\e{0}%
\e{0}%
\e{0}%
\e{0}%
\e{0}%
\eol}\vss}\rg%
%
%
\rx{\vss\hfull{%
\rlx{\hss{$700_{x}$}}\cg%
\e{0}%
\e{0}%
\e{1}%
\e{0}%
\e{0}%
\e{0}%
\e{0}%
\e{0}%
\e{1}%
\e{0}%
\e{0}%
\e{0}%
\e{1}%
\e{0}%
\e{0}%
\e{0}%
\e{1}%
\e{0}%
\e{1}%
\e{0}%
\e{0}%
\e{0}%
\e{0}%
\e{0}%
\eol}\vss}\rg%
%
%
\rx{\vss\hfull{%
\rlx{\hss{$175_{x}$}}\cg%
\e{0}%
\e{0}%
\e{0}%
\e{0}%
\e{0}%
\e{0}%
\e{0}%
\e{0}%
\e{0}%
\e{0}%
\e{0}%
\e{0}%
\e{1}%
\e{0}%
\e{0}%
\e{0}%
\e{0}%
\e{0}%
\e{0}%
\e{0}%
\e{0}%
\e{0}%
\e{0}%
\e{0}%
\eol}\vss}\rg%
%
%
\rx{\vss\hfull{%
\rlx{\hss{$1400_{x}$}}\cg%
\e{0}%
\e{0}%
\e{0}%
\e{0}%
\e{0}%
\e{0}%
\e{0}%
\e{0}%
\e{0}%
\e{0}%
\e{0}%
\e{0}%
\e{1}%
\e{0}%
\e{0}%
\e{1}%
\e{0}%
\e{0}%
\e{1}%
\e{1}%
\e{0}%
\e{0}%
\e{0}%
\e{0}%
\eol}\vss}\rg%
%
%
\rx{\vss\hfull{%
\rlx{\hss{$1050_{x}$}}\cg%
\e{0}%
\e{0}%
\e{0}%
\e{0}%
\e{0}%
\e{0}%
\e{0}%
\e{0}%
\e{1}%
\e{0}%
\e{0}%
\e{0}%
\e{1}%
\e{0}%
\e{0}%
\e{0}%
\e{0}%
\e{0}%
\e{0}%
\e{1}%
\e{0}%
\e{0}%
\e{0}%
\e{0}%
\eol}\vss}\rg%
\eop
\eject
\tablecont%
%
%
%
%
%
%
\rowpts=18 true pt%
\colpts=18 true pt%
\rowlabpts=40 true pt%
\collabpts=70 true pt%
\clx{\vss\hfull{%
\rlx{\hss{$ $}}\cg%
\cx{\hskip 16 true pt\flip{$1_a{\times}[{2}]$}\hss}\cg%
\cx{\hskip 16 true pt\flip{$7_a{\times}[{2}]$}\hss}\cg%
\cx{\hskip 16 true pt\flip{$27_a{\times}[{2}]$}\hss}\cg%
\cx{\hskip 16 true pt\flip{$21_a{\times}[{2}]$}\hss}\cg%
\cx{\hskip 16 true pt\flip{$35_a{\times}[{2}]$}\hss}\cg%
\cx{\hskip 16 true pt\flip{$105_a{\times}[{2}]$}\hss}\cg%
\cx{\hskip 16 true pt\flip{$189_a{\times}[{2}]$}\hss}\cg%
\cx{\hskip 16 true pt\flip{$21_b{\times}[{2}]$}\hss}\cg%
\cx{\hskip 16 true pt\flip{$35_b{\times}[{2}]$}\hss}\cg%
\cx{\hskip 16 true pt\flip{$189_b{\times}[{2}]$}\hss}\cg%
\cx{\hskip 16 true pt\flip{$189_c{\times}[{2}]$}\hss}\cg%
\cx{\hskip 16 true pt\flip{$15_a{\times}[{2}]$}\hss}\cg%
\cx{\hskip 16 true pt\flip{$105_b{\times}[{2}]$}\hss}\cg%
\cx{\hskip 16 true pt\flip{$105_c{\times}[{2}]$}\hss}\cg%
\cx{\hskip 16 true pt\flip{$315_a{\times}[{2}]$}\hss}\cg%
\cx{\hskip 16 true pt\flip{$405_a{\times}[{2}]$}\hss}\cg%
\cx{\hskip 16 true pt\flip{$168_a{\times}[{2}]$}\hss}\cg%
\cx{\hskip 16 true pt\flip{$56_a{\times}[{2}]$}\hss}\cg%
\cx{\hskip 16 true pt\flip{$120_a{\times}[{2}]$}\hss}\cg%
\cx{\hskip 16 true pt\flip{$210_a{\times}[{2}]$}\hss}\cg%
\cx{\hskip 16 true pt\flip{$280_a{\times}[{2}]$}\hss}\cg%
\cx{\hskip 16 true pt\flip{$336_a{\times}[{2}]$}\hss}\cg%
\cx{\hskip 16 true pt\flip{$216_a{\times}[{2}]$}\hss}\cg%
\cx{\hskip 16 true pt\flip{$512_a{\times}[{2}]$}\hss}\cg%
\eol}}\rg%
%
%
\rx{\vss\hfull{%
\rlx{\hss{$1575_{x}$}}\cg%
\e{0}%
\e{0}%
\e{0}%
\e{1}%
\e{0}%
\e{0}%
\e{1}%
\e{0}%
\e{0}%
\e{0}%
\e{0}%
\e{0}%
\e{0}%
\e{0}%
\e{0}%
\e{1}%
\e{0}%
\e{0}%
\e{1}%
\e{1}%
\e{0}%
\e{0}%
\e{0}%
\e{0}%
\eol}\vss}\rg%
%
%
\rx{\vss\hfull{%
\rlx{\hss{$1344_{x}$}}\cg%
\e{0}%
\e{0}%
\e{1}%
\e{0}%
\e{0}%
\e{0}%
\e{0}%
\e{0}%
\e{1}%
\e{0}%
\e{0}%
\e{0}%
\e{0}%
\e{0}%
\e{0}%
\e{0}%
\e{1}%
\e{0}%
\e{1}%
\e{1}%
\e{0}%
\e{0}%
\e{0}%
\e{0}%
\eol}\vss}\rg%
%
%
\rx{\vss\hfull{%
\rlx{\hss{$2100_{x}$}}\cg%
\e{0}%
\e{0}%
\e{0}%
\e{0}%
\e{0}%
\e{0}%
\e{1}%
\e{0}%
\e{0}%
\e{0}%
\e{0}%
\e{0}%
\e{0}%
\e{0}%
\e{0}%
\e{0}%
\e{0}%
\e{0}%
\e{0}%
\e{1}%
\e{0}%
\e{1}%
\e{0}%
\e{0}%
\eol}\vss}\rg%
%
%
\rx{\vss\hfull{%
\rlx{\hss{$2268_{x}$}}\cg%
\e{0}%
\e{0}%
\e{0}%
\e{0}%
\e{0}%
\e{0}%
\e{0}%
\e{0}%
\e{0}%
\e{0}%
\e{0}%
\e{0}%
\e{0}%
\e{0}%
\e{0}%
\e{1}%
\e{1}%
\e{0}%
\e{1}%
\e{1}%
\e{0}%
\e{0}%
\e{0}%
\e{0}%
\eol}\vss}\rg%
%
%
\rx{\vss\hfull{%
\rlx{\hss{$525_{x}$}}\cg%
\e{0}%
\e{0}%
\e{0}%
\e{0}%
\e{0}%
\e{0}%
\e{0}%
\e{0}%
\e{0}%
\e{0}%
\e{0}%
\e{0}%
\e{0}%
\e{1}%
\e{0}%
\e{0}%
\e{0}%
\e{0}%
\e{0}%
\e{1}%
\e{0}%
\e{0}%
\e{0}%
\e{0}%
\eol}\vss}\rg%
%
%
\rx{\vss\hfull{%
\rlx{\hss{$700_{xx}$}}\cg%
\e{0}%
\e{0}%
\e{0}%
\e{0}%
\e{0}%
\e{0}%
\e{0}%
\e{0}%
\e{0}%
\e{0}%
\e{0}%
\e{0}%
\e{0}%
\e{1}%
\e{0}%
\e{0}%
\e{0}%
\e{0}%
\e{0}%
\e{0}%
\e{0}%
\e{0}%
\e{0}%
\e{0}%
\eol}\vss}\rg%
%
%
\rx{\vss\hfull{%
\rlx{\hss{$972_{x}$}}\cg%
\e{0}%
\e{0}%
\e{0}%
\e{0}%
\e{0}%
\e{0}%
\e{0}%
\e{0}%
\e{1}%
\e{0}%
\e{0}%
\e{0}%
\e{0}%
\e{0}%
\e{0}%
\e{0}%
\e{1}%
\e{0}%
\e{0}%
\e{0}%
\e{0}%
\e{0}%
\e{0}%
\e{0}%
\eol}\vss}\rg%
%
%
\rx{\vss\hfull{%
\rlx{\hss{$4096_{x}$}}\cg%
\e{0}%
\e{0}%
\e{0}%
\e{0}%
\e{0}%
\e{0}%
\e{1}%
\e{0}%
\e{0}%
\e{0}%
\e{0}%
\e{0}%
\e{0}%
\e{0}%
\e{0}%
\e{1}%
\e{1}%
\e{0}%
\e{1}%
\e{1}%
\e{0}%
\e{0}%
\e{0}%
\e{1}%
\eol}\vss}\rg%
%
%
\rx{\vss\hfull{%
\rlx{\hss{$4200_{x}$}}\cg%
\e{0}%
\e{0}%
\e{0}%
\e{0}%
\e{0}%
\e{0}%
\e{0}%
\e{0}%
\e{0}%
\e{0}%
\e{0}%
\e{0}%
\e{1}%
\e{0}%
\e{0}%
\e{1}%
\e{1}%
\e{0}%
\e{0}%
\e{1}%
\e{0}%
\e{0}%
\e{0}%
\e{1}%
\eol}\vss}\rg%
%
%
\rx{\vss\hfull{%
\rlx{\hss{$2240_{x}$}}\cg%
\e{0}%
\e{0}%
\e{0}%
\e{0}%
\e{0}%
\e{0}%
\e{0}%
\e{0}%
\e{0}%
\e{0}%
\e{0}%
\e{0}%
\e{1}%
\e{0}%
\e{0}%
\e{1}%
\e{1}%
\e{0}%
\e{1}%
\e{0}%
\e{0}%
\e{0}%
\e{0}%
\e{0}%
\eol}\vss}\rg%
%
%
\rx{\vss\hfull{%
\rlx{\hss{$2835_{x}$}}\cg%
\e{0}%
\e{0}%
\e{0}%
\e{0}%
\e{0}%
\e{0}%
\e{0}%
\e{0}%
\e{0}%
\e{0}%
\e{0}%
\e{0}%
\e{1}%
\e{0}%
\e{0}%
\e{1}%
\e{0}%
\e{0}%
\e{0}%
\e{0}%
\e{0}%
\e{0}%
\e{0}%
\e{1}%
\eol}\vss}\rg%
%
%
\rx{\vss\hfull{%
\rlx{\hss{$6075_{x}$}}\cg%
\e{0}%
\e{0}%
\e{0}%
\e{0}%
\e{0}%
\e{0}%
\e{1}%
\e{0}%
\e{0}%
\e{0}%
\e{0}%
\e{0}%
\e{0}%
\e{1}%
\e{0}%
\e{2}%
\e{0}%
\e{0}%
\e{0}%
\e{1}%
\e{0}%
\e{1}%
\e{0}%
\e{1}%
\eol}\vss}\rg%
%
%
\rx{\vss\hfull{%
\rlx{\hss{$3200_{x}$}}\cg%
\e{0}%
\e{0}%
\e{0}%
\e{0}%
\e{0}%
\e{0}%
\e{0}%
\e{0}%
\e{0}%
\e{0}%
\e{0}%
\e{0}%
\e{0}%
\e{0}%
\e{0}%
\e{0}%
\e{1}%
\e{0}%
\e{0}%
\e{0}%
\e{0}%
\e{0}%
\e{1}%
\e{1}%
\eol}\vss}\rg%
%
%
\rx{\vss\hfull{%
\rlx{\hss{$70_{y}$}}\cg%
\e{0}%
\e{0}%
\e{0}%
\e{0}%
\e{1}%
\e{0}%
\e{0}%
\e{0}%
\e{0}%
\e{0}%
\e{0}%
\e{0}%
\e{0}%
\e{0}%
\e{0}%
\e{0}%
\e{0}%
\e{0}%
\e{0}%
\e{0}%
\e{0}%
\e{0}%
\e{0}%
\e{0}%
\eol}\vss}\rg%
%
%
\rx{\vss\hfull{%
\rlx{\hss{$1134_{y}$}}\cg%
\e{0}%
\e{0}%
\e{0}%
\e{0}%
\e{0}%
\e{0}%
\e{1}%
\e{0}%
\e{0}%
\e{0}%
\e{0}%
\e{0}%
\e{0}%
\e{0}%
\e{0}%
\e{0}%
\e{0}%
\e{0}%
\e{0}%
\e{0}%
\e{0}%
\e{0}%
\e{0}%
\e{0}%
\eol}\vss}\rg%
%
%
\rx{\vss\hfull{%
\rlx{\hss{$1680_{y}$}}\cg%
\e{0}%
\e{0}%
\e{0}%
\e{0}%
\e{1}%
\e{0}%
\e{1}%
\e{0}%
\e{0}%
\e{0}%
\e{0}%
\e{0}%
\e{0}%
\e{0}%
\e{0}%
\e{0}%
\e{0}%
\e{0}%
\e{0}%
\e{0}%
\e{1}%
\e{1}%
\e{0}%
\e{0}%
\eol}\vss}\rg%
%
%
\rx{\vss\hfull{%
\rlx{\hss{$168_{y}$}}\cg%
\e{0}%
\e{0}%
\e{0}%
\e{0}%
\e{0}%
\e{0}%
\e{0}%
\e{0}%
\e{0}%
\e{0}%
\e{0}%
\e{0}%
\e{0}%
\e{0}%
\e{0}%
\e{0}%
\e{0}%
\e{0}%
\e{0}%
\e{0}%
\e{0}%
\e{0}%
\e{0}%
\e{0}%
\eol}\vss}\rg%
%
%
\rx{\vss\hfull{%
\rlx{\hss{$420_{y}$}}\cg%
\e{0}%
\e{0}%
\e{0}%
\e{0}%
\e{0}%
\e{0}%
\e{0}%
\e{0}%
\e{0}%
\e{0}%
\e{0}%
\e{0}%
\e{0}%
\e{0}%
\e{0}%
\e{0}%
\e{0}%
\e{0}%
\e{0}%
\e{0}%
\e{0}%
\e{0}%
\e{0}%
\e{0}%
\eol}\vss}\rg%
%
%
\rx{\vss\hfull{%
\rlx{\hss{$3150_{y}$}}\cg%
\e{0}%
\e{0}%
\e{0}%
\e{0}%
\e{0}%
\e{0}%
\e{0}%
\e{0}%
\e{0}%
\e{0}%
\e{0}%
\e{0}%
\e{0}%
\e{0}%
\e{0}%
\e{1}%
\e{0}%
\e{0}%
\e{0}%
\e{0}%
\e{0}%
\e{0}%
\e{0}%
\e{1}%
\eol}\vss}\rg%
%
%
\rx{\vss\hfull{%
\rlx{\hss{$4200_{y}$}}\cg%
\e{0}%
\e{0}%
\e{0}%
\e{0}%
\e{0}%
\e{0}%
\e{0}%
\e{0}%
\e{0}%
\e{0}%
\e{0}%
\e{0}%
\e{0}%
\e{0}%
\e{0}%
\e{0}%
\e{0}%
\e{0}%
\e{0}%
\e{0}%
\e{0}%
\e{0}%
\e{1}%
\e{1}%
\eol}\vss}\rg%
%
%
\rx{\vss\hfull{%
\rlx{\hss{$2688_{y}$}}\cg%
\e{0}%
\e{0}%
\e{0}%
\e{0}%
\e{0}%
\e{0}%
\e{0}%
\e{0}%
\e{0}%
\e{0}%
\e{0}%
\e{0}%
\e{0}%
\e{0}%
\e{0}%
\e{0}%
\e{0}%
\e{0}%
\e{0}%
\e{0}%
\e{0}%
\e{1}%
\e{1}%
\e{1}%
\eol}\vss}\rg%
%
%
\rx{\vss\hfull{%
\rlx{\hss{$2100_{y}$}}\cg%
\e{0}%
\e{0}%
\e{0}%
\e{0}%
\e{0}%
\e{0}%
\e{0}%
\e{0}%
\e{0}%
\e{0}%
\e{1}%
\e{0}%
\e{0}%
\e{1}%
\e{0}%
\e{0}%
\e{0}%
\e{0}%
\e{0}%
\e{0}%
\e{0}%
\e{1}%
\e{0}%
\e{0}%
\eol}\vss}\rg%
%
%
\rx{\vss\hfull{%
\rlx{\hss{$1400_{y}$}}\cg%
\e{0}%
\e{0}%
\e{0}%
\e{0}%
\e{0}%
\e{0}%
\e{0}%
\e{0}%
\e{0}%
\e{0}%
\e{0}%
\e{0}%
\e{0}%
\e{0}%
\e{0}%
\e{0}%
\e{0}%
\e{0}%
\e{0}%
\e{0}%
\e{1}%
\e{0}%
\e{0}%
\e{0}%
\eol}\vss}\rg%
%
%
\rx{\vss\hfull{%
\rlx{\hss{$4536_{y}$}}\cg%
\e{0}%
\e{0}%
\e{0}%
\e{0}%
\e{0}%
\e{0}%
\e{0}%
\e{0}%
\e{0}%
\e{0}%
\e{0}%
\e{0}%
\e{0}%
\e{0}%
\e{1}%
\e{1}%
\e{0}%
\e{0}%
\e{0}%
\e{0}%
\e{1}%
\e{1}%
\e{0}%
\e{1}%
\eol}\vss}\rg%
%
%
\rx{\vss\hfull{%
\rlx{\hss{$5670_{y}$}}\cg%
\e{0}%
\e{0}%
\e{0}%
\e{0}%
\e{0}%
\e{0}%
\e{1}%
\e{0}%
\e{0}%
\e{0}%
\e{0}%
\e{0}%
\e{0}%
\e{0}%
\e{1}%
\e{1}%
\e{0}%
\e{0}%
\e{0}%
\e{0}%
\e{1}%
\e{1}%
\e{0}%
\e{1}%
\eol}\vss}\rg%
%
%
\rx{\vss\hfull{%
\rlx{\hss{$4480_{y}$}}\cg%
\e{0}%
\e{0}%
\e{0}%
\e{0}%
\e{0}%
\e{0}%
\e{0}%
\e{0}%
\e{0}%
\e{0}%
\e{0}%
\e{0}%
\e{0}%
\e{0}%
\e{1}%
\e{1}%
\e{0}%
\e{0}%
\e{0}%
\e{0}%
\e{0}%
\e{0}%
\e{0}%
\e{1}%
\eol}\vss}\rg%
%
%
\rx{\vss\hfull{%
\rlx{\hss{$8_{z}$}}\cg%
\e{0}%
\e{0}%
\e{0}%
\e{0}%
\e{0}%
\e{0}%
\e{0}%
\e{0}%
\e{0}%
\e{0}%
\e{0}%
\e{0}%
\e{0}%
\e{0}%
\e{0}%
\e{0}%
\e{0}%
\e{0}%
\e{0}%
\e{0}%
\e{0}%
\e{0}%
\e{0}%
\e{0}%
\eol}\vss}\rg%
%
%
\rx{\vss\hfull{%
\rlx{\hss{$56_{z}$}}\cg%
\e{0}%
\e{0}%
\e{0}%
\e{0}%
\e{0}%
\e{0}%
\e{0}%
\e{0}%
\e{0}%
\e{0}%
\e{0}%
\e{0}%
\e{0}%
\e{0}%
\e{0}%
\e{0}%
\e{0}%
\e{0}%
\e{0}%
\e{0}%
\e{0}%
\e{0}%
\e{0}%
\e{0}%
\eol}\vss}\rg%
%
%
\rx{\vss\hfull{%
\rlx{\hss{$160_{z}$}}\cg%
\e{0}%
\e{0}%
\e{0}%
\e{0}%
\e{0}%
\e{0}%
\e{0}%
\e{0}%
\e{0}%
\e{0}%
\e{0}%
\e{0}%
\e{0}%
\e{0}%
\e{0}%
\e{0}%
\e{0}%
\e{0}%
\e{0}%
\e{0}%
\e{0}%
\e{0}%
\e{0}%
\e{0}%
\eol}\vss}\rg%
%
%
\rx{\vss\hfull{%
\rlx{\hss{$112_{z}$}}\cg%
\e{0}%
\e{0}%
\e{0}%
\e{0}%
\e{0}%
\e{0}%
\e{0}%
\e{0}%
\e{0}%
\e{0}%
\e{0}%
\e{0}%
\e{0}%
\e{0}%
\e{0}%
\e{0}%
\e{0}%
\e{0}%
\e{0}%
\e{0}%
\e{0}%
\e{0}%
\e{0}%
\e{0}%
\eol}\vss}\rg%
%
%
\rx{\vss\hfull{%
\rlx{\hss{$840_{z}$}}\cg%
\e{0}%
\e{0}%
\e{0}%
\e{0}%
\e{0}%
\e{0}%
\e{0}%
\e{0}%
\e{0}%
\e{0}%
\e{0}%
\e{0}%
\e{0}%
\e{0}%
\e{0}%
\e{0}%
\e{0}%
\e{0}%
\e{0}%
\e{0}%
\e{0}%
\e{0}%
\e{0}%
\e{0}%
\eol}\vss}\rg%
\eop
\eject
\tablecont%
%
%
%
%
%
%
\rowpts=18 true pt%
\colpts=18 true pt%
\rowlabpts=40 true pt%
\collabpts=70 true pt%
\clx{\vss\hfull{%
\rlx{\hss{$ $}}\cg%
\cx{\hskip 16 true pt\flip{$1_a{\times}[{2}]$}\hss}\cg%
\cx{\hskip 16 true pt\flip{$7_a{\times}[{2}]$}\hss}\cg%
\cx{\hskip 16 true pt\flip{$27_a{\times}[{2}]$}\hss}\cg%
\cx{\hskip 16 true pt\flip{$21_a{\times}[{2}]$}\hss}\cg%
\cx{\hskip 16 true pt\flip{$35_a{\times}[{2}]$}\hss}\cg%
\cx{\hskip 16 true pt\flip{$105_a{\times}[{2}]$}\hss}\cg%
\cx{\hskip 16 true pt\flip{$189_a{\times}[{2}]$}\hss}\cg%
\cx{\hskip 16 true pt\flip{$21_b{\times}[{2}]$}\hss}\cg%
\cx{\hskip 16 true pt\flip{$35_b{\times}[{2}]$}\hss}\cg%
\cx{\hskip 16 true pt\flip{$189_b{\times}[{2}]$}\hss}\cg%
\cx{\hskip 16 true pt\flip{$189_c{\times}[{2}]$}\hss}\cg%
\cx{\hskip 16 true pt\flip{$15_a{\times}[{2}]$}\hss}\cg%
\cx{\hskip 16 true pt\flip{$105_b{\times}[{2}]$}\hss}\cg%
\cx{\hskip 16 true pt\flip{$105_c{\times}[{2}]$}\hss}\cg%
\cx{\hskip 16 true pt\flip{$315_a{\times}[{2}]$}\hss}\cg%
\cx{\hskip 16 true pt\flip{$405_a{\times}[{2}]$}\hss}\cg%
\cx{\hskip 16 true pt\flip{$168_a{\times}[{2}]$}\hss}\cg%
\cx{\hskip 16 true pt\flip{$56_a{\times}[{2}]$}\hss}\cg%
\cx{\hskip 16 true pt\flip{$120_a{\times}[{2}]$}\hss}\cg%
\cx{\hskip 16 true pt\flip{$210_a{\times}[{2}]$}\hss}\cg%
\cx{\hskip 16 true pt\flip{$280_a{\times}[{2}]$}\hss}\cg%
\cx{\hskip 16 true pt\flip{$336_a{\times}[{2}]$}\hss}\cg%
\cx{\hskip 16 true pt\flip{$216_a{\times}[{2}]$}\hss}\cg%
\cx{\hskip 16 true pt\flip{$512_a{\times}[{2}]$}\hss}\cg%
\eol}}\rg%
%
%
\rx{\vss\hfull{%
\rlx{\hss{$1296_{z}$}}\cg%
\e{0}%
\e{0}%
\e{0}%
\e{0}%
\e{0}%
\e{0}%
\e{0}%
\e{0}%
\e{0}%
\e{0}%
\e{0}%
\e{0}%
\e{0}%
\e{0}%
\e{0}%
\e{0}%
\e{0}%
\e{0}%
\e{0}%
\e{0}%
\e{0}%
\e{0}%
\e{0}%
\e{0}%
\eol}\vss}\rg%
%
%
\rx{\vss\hfull{%
\rlx{\hss{$1400_{z}$}}\cg%
\e{0}%
\e{0}%
\e{0}%
\e{0}%
\e{0}%
\e{0}%
\e{0}%
\e{0}%
\e{0}%
\e{0}%
\e{0}%
\e{0}%
\e{0}%
\e{0}%
\e{0}%
\e{0}%
\e{0}%
\e{0}%
\e{0}%
\e{0}%
\e{0}%
\e{0}%
\e{0}%
\e{0}%
\eol}\vss}\rg%
%
%
\rx{\vss\hfull{%
\rlx{\hss{$1008_{z}$}}\cg%
\e{0}%
\e{0}%
\e{0}%
\e{0}%
\e{0}%
\e{0}%
\e{0}%
\e{0}%
\e{0}%
\e{0}%
\e{0}%
\e{0}%
\e{0}%
\e{0}%
\e{0}%
\e{0}%
\e{0}%
\e{0}%
\e{0}%
\e{0}%
\e{0}%
\e{0}%
\e{0}%
\e{0}%
\eol}\vss}\rg%
%
%
\rx{\vss\hfull{%
\rlx{\hss{$560_{z}$}}\cg%
\e{0}%
\e{0}%
\e{0}%
\e{0}%
\e{0}%
\e{0}%
\e{0}%
\e{0}%
\e{0}%
\e{0}%
\e{0}%
\e{0}%
\e{0}%
\e{0}%
\e{0}%
\e{0}%
\e{0}%
\e{0}%
\e{0}%
\e{0}%
\e{0}%
\e{0}%
\e{0}%
\e{0}%
\eol}\vss}\rg%
%
%
\rx{\vss\hfull{%
\rlx{\hss{$1400_{zz}$}}\cg%
\e{0}%
\e{0}%
\e{0}%
\e{0}%
\e{0}%
\e{0}%
\e{0}%
\e{0}%
\e{0}%
\e{0}%
\e{0}%
\e{0}%
\e{0}%
\e{0}%
\e{0}%
\e{0}%
\e{0}%
\e{0}%
\e{0}%
\e{0}%
\e{0}%
\e{0}%
\e{0}%
\e{0}%
\eol}\vss}\rg%
%
%
\rx{\vss\hfull{%
\rlx{\hss{$4200_{z}$}}\cg%
\e{0}%
\e{0}%
\e{0}%
\e{0}%
\e{0}%
\e{0}%
\e{0}%
\e{0}%
\e{0}%
\e{0}%
\e{0}%
\e{0}%
\e{0}%
\e{0}%
\e{0}%
\e{0}%
\e{0}%
\e{0}%
\e{0}%
\e{0}%
\e{0}%
\e{0}%
\e{0}%
\e{0}%
\eol}\vss}\rg%
%
%
\rx{\vss\hfull{%
\rlx{\hss{$400_{z}$}}\cg%
\e{0}%
\e{0}%
\e{0}%
\e{0}%
\e{0}%
\e{0}%
\e{0}%
\e{0}%
\e{0}%
\e{0}%
\e{0}%
\e{0}%
\e{0}%
\e{0}%
\e{0}%
\e{0}%
\e{0}%
\e{0}%
\e{0}%
\e{0}%
\e{0}%
\e{0}%
\e{0}%
\e{0}%
\eol}\vss}\rg%
%
%
\rx{\vss\hfull{%
\rlx{\hss{$3240_{z}$}}\cg%
\e{0}%
\e{0}%
\e{0}%
\e{0}%
\e{0}%
\e{0}%
\e{0}%
\e{0}%
\e{0}%
\e{0}%
\e{0}%
\e{0}%
\e{0}%
\e{0}%
\e{0}%
\e{0}%
\e{0}%
\e{0}%
\e{0}%
\e{0}%
\e{0}%
\e{0}%
\e{0}%
\e{0}%
\eol}\vss}\rg%
%
%
\rx{\vss\hfull{%
\rlx{\hss{$4536_{z}$}}\cg%
\e{0}%
\e{0}%
\e{0}%
\e{0}%
\e{0}%
\e{0}%
\e{0}%
\e{0}%
\e{0}%
\e{0}%
\e{0}%
\e{0}%
\e{0}%
\e{0}%
\e{0}%
\e{0}%
\e{0}%
\e{0}%
\e{0}%
\e{0}%
\e{0}%
\e{0}%
\e{0}%
\e{0}%
\eol}\vss}\rg%
%
%
\rx{\vss\hfull{%
\rlx{\hss{$2400_{z}$}}\cg%
\e{0}%
\e{0}%
\e{0}%
\e{0}%
\e{0}%
\e{0}%
\e{0}%
\e{0}%
\e{0}%
\e{0}%
\e{0}%
\e{0}%
\e{0}%
\e{0}%
\e{0}%
\e{0}%
\e{0}%
\e{0}%
\e{0}%
\e{0}%
\e{0}%
\e{0}%
\e{0}%
\e{0}%
\eol}\vss}\rg%
%
%
\rx{\vss\hfull{%
\rlx{\hss{$3360_{z}$}}\cg%
\e{0}%
\e{0}%
\e{0}%
\e{0}%
\e{0}%
\e{0}%
\e{0}%
\e{0}%
\e{0}%
\e{0}%
\e{0}%
\e{0}%
\e{0}%
\e{0}%
\e{0}%
\e{0}%
\e{0}%
\e{0}%
\e{0}%
\e{0}%
\e{0}%
\e{0}%
\e{0}%
\e{0}%
\eol}\vss}\rg%
%
%
\rx{\vss\hfull{%
\rlx{\hss{$2800_{z}$}}\cg%
\e{0}%
\e{0}%
\e{0}%
\e{0}%
\e{0}%
\e{0}%
\e{0}%
\e{0}%
\e{0}%
\e{0}%
\e{0}%
\e{0}%
\e{0}%
\e{0}%
\e{0}%
\e{0}%
\e{0}%
\e{0}%
\e{0}%
\e{0}%
\e{0}%
\e{0}%
\e{0}%
\e{0}%
\eol}\vss}\rg%
%
%
\rx{\vss\hfull{%
\rlx{\hss{$4096_{z}$}}\cg%
\e{0}%
\e{0}%
\e{0}%
\e{0}%
\e{0}%
\e{0}%
\e{0}%
\e{0}%
\e{0}%
\e{0}%
\e{0}%
\e{0}%
\e{0}%
\e{0}%
\e{0}%
\e{0}%
\e{0}%
\e{0}%
\e{0}%
\e{0}%
\e{0}%
\e{0}%
\e{0}%
\e{0}%
\eol}\vss}\rg%
%
%
\rx{\vss\hfull{%
\rlx{\hss{$5600_{z}$}}\cg%
\e{0}%
\e{0}%
\e{0}%
\e{0}%
\e{0}%
\e{0}%
\e{0}%
\e{0}%
\e{0}%
\e{0}%
\e{0}%
\e{0}%
\e{0}%
\e{0}%
\e{0}%
\e{0}%
\e{0}%
\e{0}%
\e{0}%
\e{0}%
\e{0}%
\e{0}%
\e{0}%
\e{0}%
\eol}\vss}\rg%
%
%
\rx{\vss\hfull{%
\rlx{\hss{$448_{z}$}}\cg%
\e{0}%
\e{0}%
\e{0}%
\e{0}%
\e{0}%
\e{0}%
\e{0}%
\e{0}%
\e{0}%
\e{0}%
\e{0}%
\e{0}%
\e{0}%
\e{0}%
\e{0}%
\e{0}%
\e{0}%
\e{0}%
\e{0}%
\e{0}%
\e{0}%
\e{0}%
\e{0}%
\e{0}%
\eol}\vss}\rg%
%
%
\rx{\vss\hfull{%
\rlx{\hss{$448_{w}$}}\cg%
\e{0}%
\e{0}%
\e{0}%
\e{0}%
\e{0}%
\e{0}%
\e{0}%
\e{0}%
\e{0}%
\e{0}%
\e{0}%
\e{0}%
\e{0}%
\e{0}%
\e{0}%
\e{0}%
\e{0}%
\e{0}%
\e{0}%
\e{0}%
\e{0}%
\e{0}%
\e{0}%
\e{0}%
\eol}\vss}\rg%
%
%
\rx{\vss\hfull{%
\rlx{\hss{$1344_{w}$}}\cg%
\e{0}%
\e{0}%
\e{0}%
\e{0}%
\e{0}%
\e{0}%
\e{0}%
\e{0}%
\e{0}%
\e{0}%
\e{0}%
\e{0}%
\e{0}%
\e{0}%
\e{0}%
\e{0}%
\e{0}%
\e{0}%
\e{0}%
\e{0}%
\e{0}%
\e{0}%
\e{0}%
\e{0}%
\eol}\vss}\rg%
%
%
\rx{\vss\hfull{%
\rlx{\hss{$5600_{w}$}}\cg%
\e{0}%
\e{0}%
\e{0}%
\e{0}%
\e{0}%
\e{0}%
\e{0}%
\e{0}%
\e{0}%
\e{0}%
\e{0}%
\e{0}%
\e{0}%
\e{0}%
\e{0}%
\e{0}%
\e{0}%
\e{0}%
\e{0}%
\e{0}%
\e{0}%
\e{0}%
\e{0}%
\e{0}%
\eol}\vss}\rg%
%
%
\rx{\vss\hfull{%
\rlx{\hss{$2016_{w}$}}\cg%
\e{0}%
\e{0}%
\e{0}%
\e{0}%
\e{0}%
\e{0}%
\e{0}%
\e{0}%
\e{0}%
\e{0}%
\e{0}%
\e{0}%
\e{0}%
\e{0}%
\e{0}%
\e{0}%
\e{0}%
\e{0}%
\e{0}%
\e{0}%
\e{0}%
\e{0}%
\e{0}%
\e{0}%
\eol}\vss}\rg%
%
%
\rx{\vss\hfull{%
\rlx{\hss{$7168_{w}$}}\cg%
\e{0}%
\e{0}%
\e{0}%
\e{0}%
\e{0}%
\e{0}%
\e{0}%
\e{0}%
\e{0}%
\e{0}%
\e{0}%
\e{0}%
\e{0}%
\e{0}%
\e{0}%
\e{0}%
\e{0}%
\e{0}%
\e{0}%
\e{0}%
\e{0}%
\e{0}%
\e{0}%
\e{0}%
\eol}\vss}\rg%
%
%
%
%
%
%
\rowpts=18 true pt%
\colpts=18 true pt%
\rowlabpts=40 true pt%
\collabpts=70 true pt%
\clx{\vss\hfull{%
\rlx{\hss{$ $}}\cg%
\cx{\hskip 16 true pt\flip{$378_a{\times}[{2}]$}\hss}\cg%
\cx{\hskip 16 true pt\flip{$84_a{\times}[{2}]$}\hss}\cg%
\cx{\hskip 16 true pt\flip{$420_a{\times}[{2}]$}\hss}\cg%
\cx{\hskip 16 true pt\flip{$280_b{\times}[{2}]$}\hss}\cg%
\cx{\hskip 16 true pt\flip{$210_b{\times}[{2}]$}\hss}\cg%
\cx{\hskip 16 true pt\flip{$70_a{\times}[{2}]$}\hss}\cg%
\cx{\hskip 16 true pt\flip{$1_a^{*}{\times}[{2}]$}\hss}\cg%
\cx{\hskip 16 true pt\flip{$7_a^{*}{\times}[{2}]$}\hss}\cg%
\cx{\hskip 16 true pt\flip{$27_a^{*}{\times}[{2}]$}\hss}\cg%
\cx{\hskip 16 true pt\flip{$21_a^{*}{\times}[{2}]$}\hss}\cg%
\cx{\hskip 16 true pt\flip{$35_a^{*}{\times}[{2}]$}\hss}\cg%
\cx{\hskip 16 true pt\flip{$105_a^{*}{\times}[{2}]$}\hss}\cg%
\cx{\hskip 16 true pt\flip{$189_a^{*}{\times}[{2}]$}\hss}\cg%
\cx{\hskip 16 true pt\flip{$21_b^{*}{\times}[{2}]$}\hss}\cg%
\cx{\hskip 16 true pt\flip{$35_b^{*}{\times}[{2}]$}\hss}\cg%
\cx{\hskip 16 true pt\flip{$189_b^{*}{\times}[{2}]$}\hss}\cg%
\cx{\hskip 16 true pt\flip{$189_c^{*}{\times}[{2}]$}\hss}\cg%
\cx{\hskip 16 true pt\flip{$15_a^{*}{\times}[{2}]$}\hss}\cg%
\cx{\hskip 16 true pt\flip{$105_b^{*}{\times}[{2}]$}\hss}\cg%
\cx{\hskip 16 true pt\flip{$105_c^{*}{\times}[{2}]$}\hss}\cg%
\cx{\hskip 16 true pt\flip{$315_a^{*}{\times}[{2}]$}\hss}\cg%
\cx{\hskip 16 true pt\flip{$405_a^{*}{\times}[{2}]$}\hss}\cg%
\cx{\hskip 16 true pt\flip{$168_a^{*}{\times}[{2}]$}\hss}\cg%
\cx{\hskip 16 true pt\flip{$56_a^{*}{\times}[{2}]$}\hss}\cg%
\eol}}\rg%
%
%
\rx{\vss\hfull{%
\rlx{\hss{$1_{x}$}}\cg%
\e{0}%
\e{0}%
\e{0}%
\e{0}%
\e{0}%
\e{0}%
\e{0}%
\e{0}%
\e{0}%
\e{0}%
\e{0}%
\e{0}%
\e{0}%
\e{0}%
\e{0}%
\e{0}%
\e{0}%
\e{0}%
\e{0}%
\e{0}%
\e{0}%
\e{0}%
\e{0}%
\e{0}%
\eol}\vss}\rg%
%
%
\rx{\vss\hfull{%
\rlx{\hss{$28_{x}$}}\cg%
\e{0}%
\e{0}%
\e{0}%
\e{0}%
\e{0}%
\e{0}%
\e{0}%
\e{0}%
\e{0}%
\e{0}%
\e{0}%
\e{0}%
\e{0}%
\e{0}%
\e{0}%
\e{0}%
\e{0}%
\e{0}%
\e{0}%
\e{0}%
\e{0}%
\e{0}%
\e{0}%
\e{0}%
\eol}\vss}\rg%
%
%
\rx{\vss\hfull{%
\rlx{\hss{$35_{x}$}}\cg%
\e{0}%
\e{0}%
\e{0}%
\e{0}%
\e{0}%
\e{0}%
\e{0}%
\e{0}%
\e{0}%
\e{0}%
\e{0}%
\e{0}%
\e{0}%
\e{0}%
\e{0}%
\e{0}%
\e{0}%
\e{0}%
\e{0}%
\e{0}%
\e{0}%
\e{0}%
\e{0}%
\e{0}%
\eol}\vss}\rg%
%
%
\rx{\vss\hfull{%
\rlx{\hss{$84_{x}$}}\cg%
\e{0}%
\e{0}%
\e{0}%
\e{0}%
\e{0}%
\e{0}%
\e{0}%
\e{0}%
\e{0}%
\e{0}%
\e{0}%
\e{0}%
\e{0}%
\e{0}%
\e{0}%
\e{0}%
\e{0}%
\e{0}%
\e{0}%
\e{0}%
\e{0}%
\e{0}%
\e{0}%
\e{0}%
\eol}\vss}\rg%
%
%
\rx{\vss\hfull{%
\rlx{\hss{$50_{x}$}}\cg%
\e{0}%
\e{0}%
\e{0}%
\e{0}%
\e{0}%
\e{0}%
\e{0}%
\e{0}%
\e{0}%
\e{0}%
\e{0}%
\e{0}%
\e{0}%
\e{0}%
\e{0}%
\e{0}%
\e{0}%
\e{0}%
\e{0}%
\e{0}%
\e{0}%
\e{0}%
\e{0}%
\e{0}%
\eol}\vss}\rg%
%
%
\rx{\vss\hfull{%
\rlx{\hss{$350_{x}$}}\cg%
\e{0}%
\e{0}%
\e{0}%
\e{0}%
\e{0}%
\e{0}%
\e{0}%
\e{0}%
\e{0}%
\e{0}%
\e{0}%
\e{0}%
\e{0}%
\e{0}%
\e{0}%
\e{0}%
\e{0}%
\e{0}%
\e{0}%
\e{0}%
\e{0}%
\e{0}%
\e{0}%
\e{0}%
\eol}\vss}\rg%
%
%
\rx{\vss\hfull{%
\rlx{\hss{$300_{x}$}}\cg%
\e{0}%
\e{0}%
\e{0}%
\e{0}%
\e{0}%
\e{0}%
\e{0}%
\e{0}%
\e{0}%
\e{0}%
\e{0}%
\e{0}%
\e{0}%
\e{0}%
\e{0}%
\e{0}%
\e{0}%
\e{0}%
\e{0}%
\e{0}%
\e{0}%
\e{0}%
\e{0}%
\e{0}%
\eol}\vss}\rg%
\eop
\eject
\tablecont%
%
%
%
%
%
%
\rowpts=18 true pt%
\colpts=18 true pt%
\rowlabpts=40 true pt%
\collabpts=70 true pt%
\clx{\vss\hfull{%
\rlx{\hss{$ $}}\cg%
\cx{\hskip 16 true pt\flip{$378_a{\times}[{2}]$}\hss}\cg%
\cx{\hskip 16 true pt\flip{$84_a{\times}[{2}]$}\hss}\cg%
\cx{\hskip 16 true pt\flip{$420_a{\times}[{2}]$}\hss}\cg%
\cx{\hskip 16 true pt\flip{$280_b{\times}[{2}]$}\hss}\cg%
\cx{\hskip 16 true pt\flip{$210_b{\times}[{2}]$}\hss}\cg%
\cx{\hskip 16 true pt\flip{$70_a{\times}[{2}]$}\hss}\cg%
\cx{\hskip 16 true pt\flip{$1_a^{*}{\times}[{2}]$}\hss}\cg%
\cx{\hskip 16 true pt\flip{$7_a^{*}{\times}[{2}]$}\hss}\cg%
\cx{\hskip 16 true pt\flip{$27_a^{*}{\times}[{2}]$}\hss}\cg%
\cx{\hskip 16 true pt\flip{$21_a^{*}{\times}[{2}]$}\hss}\cg%
\cx{\hskip 16 true pt\flip{$35_a^{*}{\times}[{2}]$}\hss}\cg%
\cx{\hskip 16 true pt\flip{$105_a^{*}{\times}[{2}]$}\hss}\cg%
\cx{\hskip 16 true pt\flip{$189_a^{*}{\times}[{2}]$}\hss}\cg%
\cx{\hskip 16 true pt\flip{$21_b^{*}{\times}[{2}]$}\hss}\cg%
\cx{\hskip 16 true pt\flip{$35_b^{*}{\times}[{2}]$}\hss}\cg%
\cx{\hskip 16 true pt\flip{$189_b^{*}{\times}[{2}]$}\hss}\cg%
\cx{\hskip 16 true pt\flip{$189_c^{*}{\times}[{2}]$}\hss}\cg%
\cx{\hskip 16 true pt\flip{$15_a^{*}{\times}[{2}]$}\hss}\cg%
\cx{\hskip 16 true pt\flip{$105_b^{*}{\times}[{2}]$}\hss}\cg%
\cx{\hskip 16 true pt\flip{$105_c^{*}{\times}[{2}]$}\hss}\cg%
\cx{\hskip 16 true pt\flip{$315_a^{*}{\times}[{2}]$}\hss}\cg%
\cx{\hskip 16 true pt\flip{$405_a^{*}{\times}[{2}]$}\hss}\cg%
\cx{\hskip 16 true pt\flip{$168_a^{*}{\times}[{2}]$}\hss}\cg%
\cx{\hskip 16 true pt\flip{$56_a^{*}{\times}[{2}]$}\hss}\cg%
\eol}}\rg%
%
%
\rx{\vss\hfull{%
\rlx{\hss{$567_{x}$}}\cg%
\e{0}%
\e{0}%
\e{0}%
\e{0}%
\e{0}%
\e{0}%
\e{0}%
\e{0}%
\e{0}%
\e{0}%
\e{0}%
\e{0}%
\e{0}%
\e{0}%
\e{0}%
\e{0}%
\e{0}%
\e{0}%
\e{0}%
\e{0}%
\e{0}%
\e{0}%
\e{0}%
\e{0}%
\eol}\vss}\rg%
%
%
\rx{\vss\hfull{%
\rlx{\hss{$210_{x}$}}\cg%
\e{0}%
\e{0}%
\e{0}%
\e{0}%
\e{0}%
\e{0}%
\e{0}%
\e{0}%
\e{0}%
\e{0}%
\e{0}%
\e{0}%
\e{0}%
\e{0}%
\e{0}%
\e{0}%
\e{0}%
\e{0}%
\e{0}%
\e{0}%
\e{0}%
\e{0}%
\e{0}%
\e{0}%
\eol}\vss}\rg%
%
%
\rx{\vss\hfull{%
\rlx{\hss{$840_{x}$}}\cg%
\e{0}%
\e{1}%
\e{0}%
\e{0}%
\e{1}%
\e{0}%
\e{0}%
\e{0}%
\e{0}%
\e{0}%
\e{0}%
\e{0}%
\e{0}%
\e{0}%
\e{0}%
\e{0}%
\e{0}%
\e{0}%
\e{0}%
\e{0}%
\e{0}%
\e{0}%
\e{0}%
\e{0}%
\eol}\vss}\rg%
%
%
\rx{\vss\hfull{%
\rlx{\hss{$700_{x}$}}\cg%
\e{0}%
\e{0}%
\e{0}%
\e{0}%
\e{0}%
\e{0}%
\e{0}%
\e{0}%
\e{0}%
\e{0}%
\e{0}%
\e{0}%
\e{0}%
\e{0}%
\e{0}%
\e{0}%
\e{0}%
\e{0}%
\e{0}%
\e{0}%
\e{0}%
\e{0}%
\e{0}%
\e{0}%
\eol}\vss}\rg%
%
%
\rx{\vss\hfull{%
\rlx{\hss{$175_{x}$}}\cg%
\e{0}%
\e{0}%
\e{0}%
\e{0}%
\e{0}%
\e{0}%
\e{0}%
\e{0}%
\e{0}%
\e{0}%
\e{0}%
\e{0}%
\e{0}%
\e{0}%
\e{0}%
\e{0}%
\e{0}%
\e{0}%
\e{0}%
\e{0}%
\e{0}%
\e{0}%
\e{0}%
\e{0}%
\eol}\vss}\rg%
%
%
\rx{\vss\hfull{%
\rlx{\hss{$1400_{x}$}}\cg%
\e{0}%
\e{0}%
\e{0}%
\e{0}%
\e{0}%
\e{0}%
\e{0}%
\e{0}%
\e{0}%
\e{0}%
\e{0}%
\e{0}%
\e{0}%
\e{0}%
\e{0}%
\e{0}%
\e{0}%
\e{0}%
\e{0}%
\e{0}%
\e{0}%
\e{0}%
\e{0}%
\e{0}%
\eol}\vss}\rg%
%
%
\rx{\vss\hfull{%
\rlx{\hss{$1050_{x}$}}\cg%
\e{0}%
\e{0}%
\e{0}%
\e{1}%
\e{0}%
\e{0}%
\e{0}%
\e{0}%
\e{0}%
\e{0}%
\e{0}%
\e{0}%
\e{0}%
\e{0}%
\e{0}%
\e{0}%
\e{0}%
\e{0}%
\e{0}%
\e{0}%
\e{0}%
\e{0}%
\e{0}%
\e{0}%
\eol}\vss}\rg%
%
%
\rx{\vss\hfull{%
\rlx{\hss{$1575_{x}$}}\cg%
\e{0}%
\e{0}%
\e{0}%
\e{0}%
\e{0}%
\e{0}%
\e{0}%
\e{0}%
\e{0}%
\e{0}%
\e{0}%
\e{0}%
\e{0}%
\e{0}%
\e{0}%
\e{0}%
\e{0}%
\e{0}%
\e{0}%
\e{0}%
\e{0}%
\e{0}%
\e{0}%
\e{0}%
\eol}\vss}\rg%
%
%
\rx{\vss\hfull{%
\rlx{\hss{$1344_{x}$}}\cg%
\e{0}%
\e{0}%
\e{0}%
\e{1}%
\e{0}%
\e{0}%
\e{0}%
\e{0}%
\e{0}%
\e{0}%
\e{0}%
\e{0}%
\e{0}%
\e{0}%
\e{0}%
\e{0}%
\e{0}%
\e{0}%
\e{0}%
\e{0}%
\e{0}%
\e{0}%
\e{0}%
\e{0}%
\eol}\vss}\rg%
%
%
\rx{\vss\hfull{%
\rlx{\hss{$2100_{x}$}}\cg%
\e{0}%
\e{0}%
\e{1}%
\e{0}%
\e{0}%
\e{0}%
\e{0}%
\e{0}%
\e{0}%
\e{0}%
\e{0}%
\e{0}%
\e{0}%
\e{0}%
\e{0}%
\e{0}%
\e{0}%
\e{0}%
\e{0}%
\e{0}%
\e{0}%
\e{0}%
\e{0}%
\e{0}%
\eol}\vss}\rg%
%
%
\rx{\vss\hfull{%
\rlx{\hss{$2268_{x}$}}\cg%
\e{0}%
\e{0}%
\e{1}%
\e{0}%
\e{0}%
\e{0}%
\e{0}%
\e{0}%
\e{0}%
\e{0}%
\e{0}%
\e{0}%
\e{0}%
\e{0}%
\e{0}%
\e{0}%
\e{0}%
\e{0}%
\e{0}%
\e{0}%
\e{0}%
\e{0}%
\e{0}%
\e{0}%
\eol}\vss}\rg%
%
%
\rx{\vss\hfull{%
\rlx{\hss{$525_{x}$}}\cg%
\e{0}%
\e{0}%
\e{0}%
\e{0}%
\e{0}%
\e{0}%
\e{0}%
\e{0}%
\e{0}%
\e{0}%
\e{0}%
\e{0}%
\e{0}%
\e{0}%
\e{0}%
\e{0}%
\e{0}%
\e{0}%
\e{0}%
\e{0}%
\e{0}%
\e{0}%
\e{0}%
\e{0}%
\eol}\vss}\rg%
%
%
\rx{\vss\hfull{%
\rlx{\hss{$700_{xx}$}}\cg%
\e{0}%
\e{0}%
\e{0}%
\e{1}%
\e{0}%
\e{0}%
\e{0}%
\e{0}%
\e{0}%
\e{0}%
\e{0}%
\e{0}%
\e{0}%
\e{0}%
\e{0}%
\e{0}%
\e{0}%
\e{0}%
\e{0}%
\e{0}%
\e{0}%
\e{0}%
\e{0}%
\e{0}%
\eol}\vss}\rg%
%
%
\rx{\vss\hfull{%
\rlx{\hss{$972_{x}$}}\cg%
\e{0}%
\e{1}%
\e{0}%
\e{1}%
\e{0}%
\e{0}%
\e{0}%
\e{0}%
\e{0}%
\e{0}%
\e{0}%
\e{0}%
\e{0}%
\e{0}%
\e{0}%
\e{0}%
\e{0}%
\e{0}%
\e{0}%
\e{0}%
\e{0}%
\e{0}%
\e{0}%
\e{0}%
\eol}\vss}\rg%
%
%
\rx{\vss\hfull{%
\rlx{\hss{$4096_{x}$}}\cg%
\e{0}%
\e{0}%
\e{1}%
\e{1}%
\e{0}%
\e{0}%
\e{0}%
\e{0}%
\e{0}%
\e{0}%
\e{0}%
\e{0}%
\e{0}%
\e{0}%
\e{0}%
\e{0}%
\e{0}%
\e{0}%
\e{0}%
\e{0}%
\e{0}%
\e{0}%
\e{0}%
\e{0}%
\eol}\vss}\rg%
%
%
\rx{\vss\hfull{%
\rlx{\hss{$4200_{x}$}}\cg%
\e{0}%
\e{0}%
\e{1}%
\e{1}%
\e{1}%
\e{0}%
\e{0}%
\e{0}%
\e{0}%
\e{0}%
\e{0}%
\e{0}%
\e{0}%
\e{0}%
\e{0}%
\e{0}%
\e{0}%
\e{0}%
\e{0}%
\e{0}%
\e{0}%
\e{0}%
\e{0}%
\e{0}%
\eol}\vss}\rg%
%
%
\rx{\vss\hfull{%
\rlx{\hss{$2240_{x}$}}\cg%
\e{0}%
\e{0}%
\e{0}%
\e{1}%
\e{1}%
\e{0}%
\e{0}%
\e{0}%
\e{0}%
\e{0}%
\e{0}%
\e{0}%
\e{0}%
\e{0}%
\e{0}%
\e{0}%
\e{0}%
\e{0}%
\e{0}%
\e{0}%
\e{0}%
\e{0}%
\e{0}%
\e{0}%
\eol}\vss}\rg%
%
%
\rx{\vss\hfull{%
\rlx{\hss{$2835_{x}$}}\cg%
\e{0}%
\e{0}%
\e{0}%
\e{1}%
\e{1}%
\e{0}%
\e{0}%
\e{0}%
\e{0}%
\e{0}%
\e{0}%
\e{0}%
\e{0}%
\e{0}%
\e{0}%
\e{0}%
\e{0}%
\e{0}%
\e{0}%
\e{0}%
\e{0}%
\e{0}%
\e{0}%
\e{0}%
\eol}\vss}\rg%
%
%
\rx{\vss\hfull{%
\rlx{\hss{$6075_{x}$}}\cg%
\e{1}%
\e{0}%
\e{1}%
\e{1}%
\e{0}%
\e{0}%
\e{0}%
\e{0}%
\e{0}%
\e{0}%
\e{0}%
\e{0}%
\e{0}%
\e{0}%
\e{0}%
\e{0}%
\e{0}%
\e{0}%
\e{0}%
\e{0}%
\e{0}%
\e{0}%
\e{0}%
\e{0}%
\eol}\vss}\rg%
%
%
\rx{\vss\hfull{%
\rlx{\hss{$3200_{x}$}}\cg%
\e{0}%
\e{1}%
\e{1}%
\e{1}%
\e{0}%
\e{0}%
\e{0}%
\e{0}%
\e{0}%
\e{0}%
\e{0}%
\e{0}%
\e{0}%
\e{0}%
\e{0}%
\e{0}%
\e{0}%
\e{0}%
\e{0}%
\e{0}%
\e{0}%
\e{0}%
\e{0}%
\e{0}%
\eol}\vss}\rg%
%
%
\rx{\vss\hfull{%
\rlx{\hss{$70_{y}$}}\cg%
\e{0}%
\e{0}%
\e{0}%
\e{0}%
\e{0}%
\e{0}%
\e{0}%
\e{0}%
\e{0}%
\e{0}%
\e{0}%
\e{0}%
\e{0}%
\e{0}%
\e{0}%
\e{0}%
\e{0}%
\e{0}%
\e{0}%
\e{0}%
\e{0}%
\e{0}%
\e{0}%
\e{0}%
\eol}\vss}\rg%
%
%
\rx{\vss\hfull{%
\rlx{\hss{$1134_{y}$}}\cg%
\e{1}%
\e{0}%
\e{0}%
\e{0}%
\e{0}%
\e{0}%
\e{0}%
\e{0}%
\e{0}%
\e{0}%
\e{0}%
\e{0}%
\e{0}%
\e{0}%
\e{0}%
\e{0}%
\e{0}%
\e{0}%
\e{0}%
\e{0}%
\e{0}%
\e{0}%
\e{0}%
\e{0}%
\eol}\vss}\rg%
%
%
\rx{\vss\hfull{%
\rlx{\hss{$1680_{y}$}}\cg%
\e{0}%
\e{0}%
\e{0}%
\e{0}%
\e{0}%
\e{0}%
\e{0}%
\e{0}%
\e{0}%
\e{0}%
\e{0}%
\e{0}%
\e{0}%
\e{0}%
\e{0}%
\e{0}%
\e{0}%
\e{0}%
\e{0}%
\e{0}%
\e{0}%
\e{0}%
\e{0}%
\e{0}%
\eol}\vss}\rg%
%
%
\rx{\vss\hfull{%
\rlx{\hss{$168_{y}$}}\cg%
\e{0}%
\e{1}%
\e{0}%
\e{0}%
\e{0}%
\e{0}%
\e{0}%
\e{0}%
\e{0}%
\e{0}%
\e{0}%
\e{0}%
\e{0}%
\e{0}%
\e{0}%
\e{0}%
\e{0}%
\e{0}%
\e{0}%
\e{0}%
\e{0}%
\e{0}%
\e{0}%
\e{0}%
\eol}\vss}\rg%
%
%
\rx{\vss\hfull{%
\rlx{\hss{$420_{y}$}}\cg%
\e{0}%
\e{0}%
\e{0}%
\e{0}%
\e{1}%
\e{0}%
\e{0}%
\e{0}%
\e{0}%
\e{0}%
\e{0}%
\e{0}%
\e{0}%
\e{0}%
\e{0}%
\e{0}%
\e{0}%
\e{0}%
\e{0}%
\e{0}%
\e{0}%
\e{0}%
\e{0}%
\e{0}%
\eol}\vss}\rg%
%
%
\rx{\vss\hfull{%
\rlx{\hss{$3150_{y}$}}\cg%
\e{1}%
\e{0}%
\e{0}%
\e{0}%
\e{1}%
\e{1}%
\e{0}%
\e{0}%
\e{0}%
\e{0}%
\e{0}%
\e{0}%
\e{0}%
\e{0}%
\e{0}%
\e{0}%
\e{0}%
\e{0}%
\e{0}%
\e{0}%
\e{0}%
\e{0}%
\e{0}%
\e{0}%
\eol}\vss}\rg%
%
%
\rx{\vss\hfull{%
\rlx{\hss{$4200_{y}$}}\cg%
\e{1}%
\e{1}%
\e{1}%
\e{1}%
\e{1}%
\e{0}%
\e{0}%
\e{0}%
\e{0}%
\e{0}%
\e{0}%
\e{0}%
\e{0}%
\e{0}%
\e{0}%
\e{0}%
\e{0}%
\e{0}%
\e{0}%
\e{0}%
\e{0}%
\e{0}%
\e{0}%
\e{0}%
\eol}\vss}\rg%
%
%
\rx{\vss\hfull{%
\rlx{\hss{$2688_{y}$}}\cg%
\e{0}%
\e{0}%
\e{0}%
\e{1}%
\e{0}%
\e{0}%
\e{0}%
\e{0}%
\e{0}%
\e{0}%
\e{0}%
\e{0}%
\e{0}%
\e{0}%
\e{0}%
\e{0}%
\e{0}%
\e{0}%
\e{0}%
\e{0}%
\e{0}%
\e{0}%
\e{0}%
\e{0}%
\eol}\vss}\rg%
%
%
\rx{\vss\hfull{%
\rlx{\hss{$2100_{y}$}}\cg%
\e{0}%
\e{0}%
\e{1}%
\e{0}%
\e{0}%
\e{0}%
\e{0}%
\e{0}%
\e{0}%
\e{0}%
\e{0}%
\e{0}%
\e{0}%
\e{0}%
\e{0}%
\e{0}%
\e{0}%
\e{0}%
\e{0}%
\e{0}%
\e{0}%
\e{0}%
\e{0}%
\e{0}%
\eol}\vss}\rg%
%
%
\rx{\vss\hfull{%
\rlx{\hss{$1400_{y}$}}\cg%
\e{0}%
\e{0}%
\e{1}%
\e{0}%
\e{0}%
\e{0}%
\e{0}%
\e{0}%
\e{0}%
\e{0}%
\e{0}%
\e{0}%
\e{0}%
\e{0}%
\e{0}%
\e{0}%
\e{0}%
\e{0}%
\e{0}%
\e{0}%
\e{0}%
\e{0}%
\e{0}%
\e{0}%
\eol}\vss}\rg%
%
%
\rx{\vss\hfull{%
\rlx{\hss{$4536_{y}$}}\cg%
\e{0}%
\e{0}%
\e{1}%
\e{0}%
\e{0}%
\e{0}%
\e{0}%
\e{0}%
\e{0}%
\e{0}%
\e{0}%
\e{0}%
\e{0}%
\e{0}%
\e{0}%
\e{0}%
\e{0}%
\e{0}%
\e{0}%
\e{0}%
\e{0}%
\e{0}%
\e{0}%
\e{0}%
\eol}\vss}\rg%
%
%
\rx{\vss\hfull{%
\rlx{\hss{$5670_{y}$}}\cg%
\e{1}%
\e{0}%
\e{1}%
\e{0}%
\e{0}%
\e{0}%
\e{0}%
\e{0}%
\e{0}%
\e{0}%
\e{0}%
\e{0}%
\e{0}%
\e{0}%
\e{0}%
\e{0}%
\e{0}%
\e{0}%
\e{0}%
\e{0}%
\e{0}%
\e{0}%
\e{0}%
\e{0}%
\eol}\vss}\rg%
%
%
\rx{\vss\hfull{%
\rlx{\hss{$4480_{y}$}}\cg%
\e{1}%
\e{0}%
\e{1}%
\e{0}%
\e{1}%
\e{0}%
\e{0}%
\e{0}%
\e{0}%
\e{0}%
\e{0}%
\e{0}%
\e{0}%
\e{0}%
\e{0}%
\e{0}%
\e{0}%
\e{0}%
\e{0}%
\e{0}%
\e{0}%
\e{0}%
\e{0}%
\e{0}%
\eol}\vss}\rg%
\eop
\eject
\tablecont%
%
%
%
%
%
%
\rowpts=18 true pt%
\colpts=18 true pt%
\rowlabpts=40 true pt%
\collabpts=70 true pt%
\clx{\vss\hfull{%
\rlx{\hss{$ $}}\cg%
\cx{\hskip 16 true pt\flip{$378_a{\times}[{2}]$}\hss}\cg%
\cx{\hskip 16 true pt\flip{$84_a{\times}[{2}]$}\hss}\cg%
\cx{\hskip 16 true pt\flip{$420_a{\times}[{2}]$}\hss}\cg%
\cx{\hskip 16 true pt\flip{$280_b{\times}[{2}]$}\hss}\cg%
\cx{\hskip 16 true pt\flip{$210_b{\times}[{2}]$}\hss}\cg%
\cx{\hskip 16 true pt\flip{$70_a{\times}[{2}]$}\hss}\cg%
\cx{\hskip 16 true pt\flip{$1_a^{*}{\times}[{2}]$}\hss}\cg%
\cx{\hskip 16 true pt\flip{$7_a^{*}{\times}[{2}]$}\hss}\cg%
\cx{\hskip 16 true pt\flip{$27_a^{*}{\times}[{2}]$}\hss}\cg%
\cx{\hskip 16 true pt\flip{$21_a^{*}{\times}[{2}]$}\hss}\cg%
\cx{\hskip 16 true pt\flip{$35_a^{*}{\times}[{2}]$}\hss}\cg%
\cx{\hskip 16 true pt\flip{$105_a^{*}{\times}[{2}]$}\hss}\cg%
\cx{\hskip 16 true pt\flip{$189_a^{*}{\times}[{2}]$}\hss}\cg%
\cx{\hskip 16 true pt\flip{$21_b^{*}{\times}[{2}]$}\hss}\cg%
\cx{\hskip 16 true pt\flip{$35_b^{*}{\times}[{2}]$}\hss}\cg%
\cx{\hskip 16 true pt\flip{$189_b^{*}{\times}[{2}]$}\hss}\cg%
\cx{\hskip 16 true pt\flip{$189_c^{*}{\times}[{2}]$}\hss}\cg%
\cx{\hskip 16 true pt\flip{$15_a^{*}{\times}[{2}]$}\hss}\cg%
\cx{\hskip 16 true pt\flip{$105_b^{*}{\times}[{2}]$}\hss}\cg%
\cx{\hskip 16 true pt\flip{$105_c^{*}{\times}[{2}]$}\hss}\cg%
\cx{\hskip 16 true pt\flip{$315_a^{*}{\times}[{2}]$}\hss}\cg%
\cx{\hskip 16 true pt\flip{$405_a^{*}{\times}[{2}]$}\hss}\cg%
\cx{\hskip 16 true pt\flip{$168_a^{*}{\times}[{2}]$}\hss}\cg%
\cx{\hskip 16 true pt\flip{$56_a^{*}{\times}[{2}]$}\hss}\cg%
\eol}}\rg%
%
%
\rx{\vss\hfull{%
\rlx{\hss{$8_{z}$}}\cg%
\e{0}%
\e{0}%
\e{0}%
\e{0}%
\e{0}%
\e{0}%
\e{0}%
\e{1}%
\e{0}%
\e{0}%
\e{0}%
\e{0}%
\e{0}%
\e{0}%
\e{0}%
\e{0}%
\e{0}%
\e{0}%
\e{0}%
\e{0}%
\e{0}%
\e{0}%
\e{0}%
\e{0}%
\eol}\vss}\rg%
%
%
\rx{\vss\hfull{%
\rlx{\hss{$56_{z}$}}\cg%
\e{0}%
\e{0}%
\e{0}%
\e{0}%
\e{0}%
\e{0}%
\e{0}%
\e{0}%
\e{0}%
\e{0}%
\e{1}%
\e{0}%
\e{0}%
\e{0}%
\e{0}%
\e{0}%
\e{0}%
\e{0}%
\e{0}%
\e{0}%
\e{0}%
\e{0}%
\e{0}%
\e{0}%
\eol}\vss}\rg%
%
%
\rx{\vss\hfull{%
\rlx{\hss{$160_{z}$}}\cg%
\e{0}%
\e{0}%
\e{0}%
\e{0}%
\e{0}%
\e{0}%
\e{0}%
\e{1}%
\e{0}%
\e{0}%
\e{0}%
\e{1}%
\e{0}%
\e{0}%
\e{0}%
\e{0}%
\e{0}%
\e{0}%
\e{0}%
\e{0}%
\e{0}%
\e{0}%
\e{0}%
\e{0}%
\eol}\vss}\rg%
%
%
\rx{\vss\hfull{%
\rlx{\hss{$112_{z}$}}\cg%
\e{0}%
\e{0}%
\e{0}%
\e{0}%
\e{0}%
\e{0}%
\e{0}%
\e{1}%
\e{0}%
\e{0}%
\e{0}%
\e{0}%
\e{0}%
\e{1}%
\e{0}%
\e{0}%
\e{0}%
\e{0}%
\e{0}%
\e{0}%
\e{0}%
\e{0}%
\e{0}%
\e{1}%
\eol}\vss}\rg%
%
%
\rx{\vss\hfull{%
\rlx{\hss{$840_{z}$}}\cg%
\e{0}%
\e{0}%
\e{0}%
\e{0}%
\e{0}%
\e{0}%
\e{0}%
\e{0}%
\e{0}%
\e{0}%
\e{0}%
\e{1}%
\e{0}%
\e{0}%
\e{0}%
\e{0}%
\e{0}%
\e{0}%
\e{0}%
\e{0}%
\e{0}%
\e{0}%
\e{0}%
\e{0}%
\eol}\vss}\rg%
%
%
\rx{\vss\hfull{%
\rlx{\hss{$1296_{z}$}}\cg%
\e{0}%
\e{0}%
\e{0}%
\e{0}%
\e{0}%
\e{0}%
\e{0}%
\e{0}%
\e{0}%
\e{0}%
\e{1}%
\e{1}%
\e{0}%
\e{0}%
\e{0}%
\e{0}%
\e{0}%
\e{0}%
\e{0}%
\e{0}%
\e{0}%
\e{0}%
\e{0}%
\e{0}%
\eol}\vss}\rg%
%
%
\rx{\vss\hfull{%
\rlx{\hss{$1400_{z}$}}\cg%
\e{0}%
\e{0}%
\e{0}%
\e{0}%
\e{0}%
\e{0}%
\e{0}%
\e{0}%
\e{0}%
\e{0}%
\e{0}%
\e{1}%
\e{0}%
\e{1}%
\e{0}%
\e{1}%
\e{1}%
\e{0}%
\e{0}%
\e{0}%
\e{1}%
\e{0}%
\e{0}%
\e{1}%
\eol}\vss}\rg%
%
%
\rx{\vss\hfull{%
\rlx{\hss{$1008_{z}$}}\cg%
\e{0}%
\e{0}%
\e{0}%
\e{0}%
\e{0}%
\e{0}%
\e{0}%
\e{0}%
\e{0}%
\e{0}%
\e{0}%
\e{1}%
\e{0}%
\e{0}%
\e{0}%
\e{0}%
\e{1}%
\e{0}%
\e{0}%
\e{0}%
\e{0}%
\e{0}%
\e{0}%
\e{1}%
\eol}\vss}\rg%
%
%
\rx{\vss\hfull{%
\rlx{\hss{$560_{z}$}}\cg%
\e{0}%
\e{0}%
\e{0}%
\e{0}%
\e{0}%
\e{0}%
\e{0}%
\e{1}%
\e{0}%
\e{0}%
\e{0}%
\e{1}%
\e{0}%
\e{1}%
\e{0}%
\e{1}%
\e{0}%
\e{0}%
\e{0}%
\e{0}%
\e{0}%
\e{0}%
\e{0}%
\e{1}%
\eol}\vss}\rg%
%
%
\rx{\vss\hfull{%
\rlx{\hss{$1400_{zz}$}}\cg%
\e{0}%
\e{0}%
\e{0}%
\e{0}%
\e{0}%
\e{0}%
\e{0}%
\e{0}%
\e{0}%
\e{0}%
\e{0}%
\e{0}%
\e{0}%
\e{0}%
\e{0}%
\e{1}%
\e{0}%
\e{1}%
\e{0}%
\e{0}%
\e{1}%
\e{0}%
\e{0}%
\e{0}%
\eol}\vss}\rg%
%
%
\rx{\vss\hfull{%
\rlx{\hss{$4200_{z}$}}\cg%
\e{0}%
\e{0}%
\e{0}%
\e{0}%
\e{0}%
\e{0}%
\e{0}%
\e{0}%
\e{0}%
\e{0}%
\e{0}%
\e{0}%
\e{0}%
\e{0}%
\e{0}%
\e{0}%
\e{0}%
\e{0}%
\e{0}%
\e{0}%
\e{1}%
\e{0}%
\e{0}%
\e{0}%
\eol}\vss}\rg%
%
%
\rx{\vss\hfull{%
\rlx{\hss{$400_{z}$}}\cg%
\e{0}%
\e{0}%
\e{0}%
\e{0}%
\e{0}%
\e{0}%
\e{0}%
\e{0}%
\e{0}%
\e{0}%
\e{0}%
\e{0}%
\e{0}%
\e{0}%
\e{0}%
\e{1}%
\e{0}%
\e{1}%
\e{0}%
\e{0}%
\e{0}%
\e{0}%
\e{0}%
\e{1}%
\eol}\vss}\rg%
%
%
\rx{\vss\hfull{%
\rlx{\hss{$3240_{z}$}}\cg%
\e{0}%
\e{0}%
\e{0}%
\e{0}%
\e{0}%
\e{0}%
\e{0}%
\e{0}%
\e{0}%
\e{0}%
\e{0}%
\e{1}%
\e{0}%
\e{0}%
\e{0}%
\e{2}%
\e{1}%
\e{0}%
\e{0}%
\e{0}%
\e{1}%
\e{0}%
\e{0}%
\e{1}%
\eol}\vss}\rg%
%
%
\rx{\vss\hfull{%
\rlx{\hss{$4536_{z}$}}\cg%
\e{0}%
\e{0}%
\e{0}%
\e{0}%
\e{0}%
\e{0}%
\e{0}%
\e{0}%
\e{0}%
\e{0}%
\e{0}%
\e{0}%
\e{0}%
\e{0}%
\e{0}%
\e{1}%
\e{1}%
\e{0}%
\e{0}%
\e{0}%
\e{1}%
\e{0}%
\e{0}%
\e{0}%
\eol}\vss}\rg%
%
%
\rx{\vss\hfull{%
\rlx{\hss{$2400_{z}$}}\cg%
\e{0}%
\e{0}%
\e{0}%
\e{0}%
\e{0}%
\e{0}%
\e{0}%
\e{0}%
\e{0}%
\e{0}%
\e{1}%
\e{0}%
\e{1}%
\e{0}%
\e{0}%
\e{0}%
\e{0}%
\e{0}%
\e{0}%
\e{0}%
\e{0}%
\e{0}%
\e{0}%
\e{0}%
\eol}\vss}\rg%
%
%
\rx{\vss\hfull{%
\rlx{\hss{$3360_{z}$}}\cg%
\e{0}%
\e{0}%
\e{0}%
\e{0}%
\e{0}%
\e{0}%
\e{0}%
\e{0}%
\e{0}%
\e{0}%
\e{0}%
\e{0}%
\e{0}%
\e{0}%
\e{0}%
\e{1}%
\e{0}%
\e{0}%
\e{0}%
\e{0}%
\e{1}%
\e{0}%
\e{0}%
\e{0}%
\eol}\vss}\rg%
%
%
\rx{\vss\hfull{%
\rlx{\hss{$2800_{z}$}}\cg%
\e{0}%
\e{0}%
\e{0}%
\e{0}%
\e{0}%
\e{0}%
\e{0}%
\e{0}%
\e{0}%
\e{0}%
\e{0}%
\e{0}%
\e{0}%
\e{0}%
\e{0}%
\e{0}%
\e{1}%
\e{0}%
\e{0}%
\e{0}%
\e{1}%
\e{0}%
\e{0}%
\e{0}%
\eol}\vss}\rg%
%
%
\rx{\vss\hfull{%
\rlx{\hss{$4096_{z}$}}\cg%
\e{0}%
\e{0}%
\e{0}%
\e{0}%
\e{0}%
\e{0}%
\e{0}%
\e{0}%
\e{0}%
\e{0}%
\e{0}%
\e{1}%
\e{0}%
\e{0}%
\e{0}%
\e{1}%
\e{1}%
\e{0}%
\e{0}%
\e{0}%
\e{1}%
\e{0}%
\e{0}%
\e{0}%
\eol}\vss}\rg%
%
%
\rx{\vss\hfull{%
\rlx{\hss{$5600_{z}$}}\cg%
\e{0}%
\e{0}%
\e{0}%
\e{0}%
\e{0}%
\e{0}%
\e{0}%
\e{0}%
\e{0}%
\e{0}%
\e{0}%
\e{0}%
\e{0}%
\e{0}%
\e{0}%
\e{0}%
\e{1}%
\e{0}%
\e{0}%
\e{1}%
\e{1}%
\e{1}%
\e{0}%
\e{0}%
\eol}\vss}\rg%
%
%
\rx{\vss\hfull{%
\rlx{\hss{$448_{z}$}}\cg%
\e{0}%
\e{0}%
\e{0}%
\e{0}%
\e{0}%
\e{0}%
\e{0}%
\e{0}%
\e{0}%
\e{0}%
\e{0}%
\e{0}%
\e{0}%
\e{1}%
\e{0}%
\e{1}%
\e{0}%
\e{0}%
\e{0}%
\e{0}%
\e{0}%
\e{0}%
\e{0}%
\e{0}%
\eol}\vss}\rg%
%
%
\rx{\vss\hfull{%
\rlx{\hss{$448_{w}$}}\cg%
\e{0}%
\e{0}%
\e{0}%
\e{0}%
\e{0}%
\e{0}%
\e{0}%
\e{0}%
\e{0}%
\e{0}%
\e{1}%
\e{0}%
\e{1}%
\e{0}%
\e{0}%
\e{0}%
\e{0}%
\e{0}%
\e{0}%
\e{0}%
\e{0}%
\e{0}%
\e{0}%
\e{0}%
\eol}\vss}\rg%
%
%
\rx{\vss\hfull{%
\rlx{\hss{$1344_{w}$}}\cg%
\e{0}%
\e{0}%
\e{0}%
\e{0}%
\e{0}%
\e{0}%
\e{0}%
\e{0}%
\e{0}%
\e{0}%
\e{0}%
\e{0}%
\e{0}%
\e{0}%
\e{0}%
\e{0}%
\e{0}%
\e{0}%
\e{0}%
\e{0}%
\e{0}%
\e{0}%
\e{0}%
\e{0}%
\eol}\vss}\rg%
%
%
\rx{\vss\hfull{%
\rlx{\hss{$5600_{w}$}}\cg%
\e{0}%
\e{0}%
\e{0}%
\e{0}%
\e{0}%
\e{0}%
\e{0}%
\e{0}%
\e{0}%
\e{0}%
\e{0}%
\e{0}%
\e{1}%
\e{0}%
\e{0}%
\e{0}%
\e{0}%
\e{0}%
\e{0}%
\e{0}%
\e{0}%
\e{1}%
\e{0}%
\e{0}%
\eol}\vss}\rg%
%
%
\rx{\vss\hfull{%
\rlx{\hss{$2016_{w}$}}\cg%
\e{0}%
\e{0}%
\e{0}%
\e{0}%
\e{0}%
\e{0}%
\e{0}%
\e{0}%
\e{0}%
\e{0}%
\e{0}%
\e{0}%
\e{0}%
\e{0}%
\e{0}%
\e{0}%
\e{0}%
\e{0}%
\e{0}%
\e{0}%
\e{0}%
\e{0}%
\e{0}%
\e{0}%
\eol}\vss}\rg%
%
%
\rx{\vss\hfull{%
\rlx{\hss{$7168_{w}$}}\cg%
\e{0}%
\e{0}%
\e{0}%
\e{0}%
\e{0}%
\e{0}%
\e{0}%
\e{0}%
\e{0}%
\e{0}%
\e{0}%
\e{0}%
\e{0}%
\e{0}%
\e{0}%
\e{0}%
\e{0}%
\e{0}%
\e{0}%
\e{0}%
\e{1}%
\e{1}%
\e{0}%
\e{0}%
\eol}\vss}\rg%
\eop
\eject
\tablecont%
%
%
%
%
%
%
\rowpts=18 true pt%
\colpts=18 true pt%
\rowlabpts=40 true pt%
\collabpts=70 true pt%
\clx{\vss\hfull{%
\rlx{\hss{$ $}}\cg%
\cx{\hskip 16 true pt\flip{$120_a^{*}{\times}[{2}]$}\hss}\cg%
\cx{\hskip 16 true pt\flip{$210_a^{*}{\times}[{2}]$}\hss}\cg%
\cx{\hskip 16 true pt\flip{$280_a^{*}{\times}[{2}]$}\hss}\cg%
\cx{\hskip 16 true pt\flip{$336_a^{*}{\times}[{2}]$}\hss}\cg%
\cx{\hskip 16 true pt\flip{$216_a^{*}{\times}[{2}]$}\hss}\cg%
\cx{\hskip 16 true pt\flip{$512_a^{*}{\times}[{2}]$}\hss}\cg%
\cx{\hskip 16 true pt\flip{$378_a^{*}{\times}[{2}]$}\hss}\cg%
\cx{\hskip 16 true pt\flip{$84_a^{*}{\times}[{2}]$}\hss}\cg%
\cx{\hskip 16 true pt\flip{$420_a^{*}{\times}[{2}]$}\hss}\cg%
\cx{\hskip 16 true pt\flip{$280_b^{*}{\times}[{2}]$}\hss}\cg%
\cx{\hskip 16 true pt\flip{$210_b^{*}{\times}[{2}]$}\hss}\cg%
\cx{\hskip 16 true pt\flip{$70_a^{*}{\times}[{2}]$}\hss}\cg%
\cx{\hskip 16 true pt\flip{$1_a{\times}[{1^{2}}]$}\hss}\cg%
\cx{\hskip 16 true pt\flip{$7_a{\times}[{1^{2}}]$}\hss}\cg%
\cx{\hskip 16 true pt\flip{$27_a{\times}[{1^{2}}]$}\hss}\cg%
\cx{\hskip 16 true pt\flip{$21_a{\times}[{1^{2}}]$}\hss}\cg%
\cx{\hskip 16 true pt\flip{$35_a{\times}[{1^{2}}]$}\hss}\cg%
\cx{\hskip 16 true pt\flip{$105_a{\times}[{1^{2}}]$}\hss}\cg%
\cx{\hskip 16 true pt\flip{$189_a{\times}[{1^{2}}]$}\hss}\cg%
\cx{\hskip 16 true pt\flip{$21_b{\times}[{1^{2}}]$}\hss}\cg%
\cx{\hskip 16 true pt\flip{$35_b{\times}[{1^{2}}]$}\hss}\cg%
\cx{\hskip 16 true pt\flip{$189_b{\times}[{1^{2}}]$}\hss}\cg%
\cx{\hskip 16 true pt\flip{$189_c{\times}[{1^{2}}]$}\hss}\cg%
\cx{\hskip 16 true pt\flip{$15_a{\times}[{1^{2}}]$}\hss}\cg%
\eol}}\rg%
%
%
\rx{\vss\hfull{%
\rlx{\hss{$1_{x}$}}\cg%
\e{0}%
\e{0}%
\e{0}%
\e{0}%
\e{0}%
\e{0}%
\e{0}%
\e{0}%
\e{0}%
\e{0}%
\e{0}%
\e{0}%
\e{0}%
\e{0}%
\e{0}%
\e{0}%
\e{0}%
\e{0}%
\e{0}%
\e{0}%
\e{0}%
\e{0}%
\e{0}%
\e{0}%
\eol}\vss}\rg%
%
%
\rx{\vss\hfull{%
\rlx{\hss{$28_{x}$}}\cg%
\e{0}%
\e{0}%
\e{0}%
\e{0}%
\e{0}%
\e{0}%
\e{0}%
\e{0}%
\e{0}%
\e{0}%
\e{0}%
\e{0}%
\e{0}%
\e{0}%
\e{0}%
\e{0}%
\e{0}%
\e{0}%
\e{0}%
\e{0}%
\e{0}%
\e{0}%
\e{0}%
\e{0}%
\eol}\vss}\rg%
%
%
\rx{\vss\hfull{%
\rlx{\hss{$35_{x}$}}\cg%
\e{0}%
\e{0}%
\e{0}%
\e{0}%
\e{0}%
\e{0}%
\e{0}%
\e{0}%
\e{0}%
\e{0}%
\e{0}%
\e{0}%
\e{0}%
\e{0}%
\e{0}%
\e{0}%
\e{0}%
\e{0}%
\e{0}%
\e{0}%
\e{0}%
\e{0}%
\e{0}%
\e{0}%
\eol}\vss}\rg%
%
%
\rx{\vss\hfull{%
\rlx{\hss{$84_{x}$}}\cg%
\e{0}%
\e{0}%
\e{0}%
\e{0}%
\e{0}%
\e{0}%
\e{0}%
\e{0}%
\e{0}%
\e{0}%
\e{0}%
\e{0}%
\e{0}%
\e{0}%
\e{0}%
\e{0}%
\e{0}%
\e{0}%
\e{0}%
\e{0}%
\e{0}%
\e{0}%
\e{0}%
\e{0}%
\eol}\vss}\rg%
%
%
\rx{\vss\hfull{%
\rlx{\hss{$50_{x}$}}\cg%
\e{0}%
\e{0}%
\e{0}%
\e{0}%
\e{0}%
\e{0}%
\e{0}%
\e{0}%
\e{0}%
\e{0}%
\e{0}%
\e{0}%
\e{0}%
\e{0}%
\e{0}%
\e{0}%
\e{0}%
\e{0}%
\e{0}%
\e{0}%
\e{0}%
\e{0}%
\e{0}%
\e{0}%
\eol}\vss}\rg%
%
%
\rx{\vss\hfull{%
\rlx{\hss{$350_{x}$}}\cg%
\e{0}%
\e{0}%
\e{0}%
\e{0}%
\e{0}%
\e{0}%
\e{0}%
\e{0}%
\e{0}%
\e{0}%
\e{0}%
\e{0}%
\e{0}%
\e{0}%
\e{0}%
\e{0}%
\e{0}%
\e{0}%
\e{0}%
\e{0}%
\e{0}%
\e{0}%
\e{0}%
\e{0}%
\eol}\vss}\rg%
%
%
\rx{\vss\hfull{%
\rlx{\hss{$300_{x}$}}\cg%
\e{0}%
\e{0}%
\e{0}%
\e{0}%
\e{0}%
\e{0}%
\e{0}%
\e{0}%
\e{0}%
\e{0}%
\e{0}%
\e{0}%
\e{0}%
\e{0}%
\e{0}%
\e{0}%
\e{0}%
\e{0}%
\e{0}%
\e{0}%
\e{0}%
\e{0}%
\e{0}%
\e{0}%
\eol}\vss}\rg%
%
%
\rx{\vss\hfull{%
\rlx{\hss{$567_{x}$}}\cg%
\e{0}%
\e{0}%
\e{0}%
\e{0}%
\e{0}%
\e{0}%
\e{0}%
\e{0}%
\e{0}%
\e{0}%
\e{0}%
\e{0}%
\e{0}%
\e{0}%
\e{0}%
\e{0}%
\e{0}%
\e{0}%
\e{0}%
\e{0}%
\e{0}%
\e{0}%
\e{0}%
\e{0}%
\eol}\vss}\rg%
%
%
\rx{\vss\hfull{%
\rlx{\hss{$210_{x}$}}\cg%
\e{0}%
\e{0}%
\e{0}%
\e{0}%
\e{0}%
\e{0}%
\e{0}%
\e{0}%
\e{0}%
\e{0}%
\e{0}%
\e{0}%
\e{0}%
\e{0}%
\e{0}%
\e{0}%
\e{0}%
\e{0}%
\e{0}%
\e{0}%
\e{0}%
\e{0}%
\e{0}%
\e{0}%
\eol}\vss}\rg%
%
%
\rx{\vss\hfull{%
\rlx{\hss{$840_{x}$}}\cg%
\e{0}%
\e{0}%
\e{0}%
\e{0}%
\e{0}%
\e{0}%
\e{0}%
\e{0}%
\e{0}%
\e{0}%
\e{0}%
\e{0}%
\e{0}%
\e{0}%
\e{0}%
\e{0}%
\e{0}%
\e{0}%
\e{0}%
\e{0}%
\e{0}%
\e{0}%
\e{0}%
\e{0}%
\eol}\vss}\rg%
%
%
\rx{\vss\hfull{%
\rlx{\hss{$700_{x}$}}\cg%
\e{0}%
\e{0}%
\e{0}%
\e{0}%
\e{0}%
\e{0}%
\e{0}%
\e{0}%
\e{0}%
\e{0}%
\e{0}%
\e{0}%
\e{0}%
\e{0}%
\e{0}%
\e{0}%
\e{0}%
\e{0}%
\e{0}%
\e{0}%
\e{0}%
\e{0}%
\e{0}%
\e{0}%
\eol}\vss}\rg%
%
%
\rx{\vss\hfull{%
\rlx{\hss{$175_{x}$}}\cg%
\e{0}%
\e{0}%
\e{0}%
\e{0}%
\e{0}%
\e{0}%
\e{0}%
\e{0}%
\e{0}%
\e{0}%
\e{0}%
\e{0}%
\e{0}%
\e{0}%
\e{0}%
\e{0}%
\e{0}%
\e{0}%
\e{0}%
\e{0}%
\e{0}%
\e{0}%
\e{0}%
\e{0}%
\eol}\vss}\rg%
%
%
\rx{\vss\hfull{%
\rlx{\hss{$1400_{x}$}}\cg%
\e{0}%
\e{0}%
\e{0}%
\e{0}%
\e{0}%
\e{0}%
\e{0}%
\e{0}%
\e{0}%
\e{0}%
\e{0}%
\e{0}%
\e{0}%
\e{0}%
\e{0}%
\e{0}%
\e{0}%
\e{0}%
\e{0}%
\e{0}%
\e{0}%
\e{0}%
\e{0}%
\e{0}%
\eol}\vss}\rg%
%
%
\rx{\vss\hfull{%
\rlx{\hss{$1050_{x}$}}\cg%
\e{0}%
\e{0}%
\e{0}%
\e{0}%
\e{0}%
\e{0}%
\e{0}%
\e{0}%
\e{0}%
\e{0}%
\e{0}%
\e{0}%
\e{0}%
\e{0}%
\e{0}%
\e{0}%
\e{0}%
\e{0}%
\e{0}%
\e{0}%
\e{0}%
\e{0}%
\e{0}%
\e{0}%
\eol}\vss}\rg%
%
%
\rx{\vss\hfull{%
\rlx{\hss{$1575_{x}$}}\cg%
\e{0}%
\e{0}%
\e{0}%
\e{0}%
\e{0}%
\e{0}%
\e{0}%
\e{0}%
\e{0}%
\e{0}%
\e{0}%
\e{0}%
\e{0}%
\e{0}%
\e{0}%
\e{0}%
\e{0}%
\e{0}%
\e{0}%
\e{0}%
\e{0}%
\e{0}%
\e{0}%
\e{0}%
\eol}\vss}\rg%
%
%
\rx{\vss\hfull{%
\rlx{\hss{$1344_{x}$}}\cg%
\e{0}%
\e{0}%
\e{0}%
\e{0}%
\e{0}%
\e{0}%
\e{0}%
\e{0}%
\e{0}%
\e{0}%
\e{0}%
\e{0}%
\e{0}%
\e{0}%
\e{0}%
\e{0}%
\e{0}%
\e{0}%
\e{0}%
\e{0}%
\e{0}%
\e{0}%
\e{0}%
\e{0}%
\eol}\vss}\rg%
%
%
\rx{\vss\hfull{%
\rlx{\hss{$2100_{x}$}}\cg%
\e{0}%
\e{0}%
\e{0}%
\e{0}%
\e{0}%
\e{0}%
\e{0}%
\e{0}%
\e{0}%
\e{0}%
\e{0}%
\e{0}%
\e{0}%
\e{0}%
\e{0}%
\e{0}%
\e{0}%
\e{0}%
\e{0}%
\e{0}%
\e{0}%
\e{0}%
\e{0}%
\e{0}%
\eol}\vss}\rg%
%
%
\rx{\vss\hfull{%
\rlx{\hss{$2268_{x}$}}\cg%
\e{0}%
\e{0}%
\e{0}%
\e{0}%
\e{0}%
\e{0}%
\e{0}%
\e{0}%
\e{0}%
\e{0}%
\e{0}%
\e{0}%
\e{0}%
\e{0}%
\e{0}%
\e{0}%
\e{0}%
\e{0}%
\e{0}%
\e{0}%
\e{0}%
\e{0}%
\e{0}%
\e{0}%
\eol}\vss}\rg%
%
%
\rx{\vss\hfull{%
\rlx{\hss{$525_{x}$}}\cg%
\e{0}%
\e{0}%
\e{0}%
\e{0}%
\e{0}%
\e{0}%
\e{0}%
\e{0}%
\e{0}%
\e{0}%
\e{0}%
\e{0}%
\e{0}%
\e{0}%
\e{0}%
\e{0}%
\e{0}%
\e{0}%
\e{0}%
\e{0}%
\e{0}%
\e{0}%
\e{0}%
\e{0}%
\eol}\vss}\rg%
%
%
\rx{\vss\hfull{%
\rlx{\hss{$700_{xx}$}}\cg%
\e{0}%
\e{0}%
\e{0}%
\e{0}%
\e{0}%
\e{0}%
\e{0}%
\e{0}%
\e{0}%
\e{0}%
\e{0}%
\e{0}%
\e{0}%
\e{0}%
\e{0}%
\e{0}%
\e{0}%
\e{0}%
\e{0}%
\e{0}%
\e{0}%
\e{0}%
\e{0}%
\e{0}%
\eol}\vss}\rg%
%
%
\rx{\vss\hfull{%
\rlx{\hss{$972_{x}$}}\cg%
\e{0}%
\e{0}%
\e{0}%
\e{0}%
\e{0}%
\e{0}%
\e{0}%
\e{0}%
\e{0}%
\e{0}%
\e{0}%
\e{0}%
\e{0}%
\e{0}%
\e{0}%
\e{0}%
\e{0}%
\e{0}%
\e{0}%
\e{0}%
\e{0}%
\e{0}%
\e{0}%
\e{0}%
\eol}\vss}\rg%
%
%
\rx{\vss\hfull{%
\rlx{\hss{$4096_{x}$}}\cg%
\e{0}%
\e{0}%
\e{0}%
\e{0}%
\e{0}%
\e{0}%
\e{0}%
\e{0}%
\e{0}%
\e{0}%
\e{0}%
\e{0}%
\e{0}%
\e{0}%
\e{0}%
\e{0}%
\e{0}%
\e{0}%
\e{0}%
\e{0}%
\e{0}%
\e{0}%
\e{0}%
\e{0}%
\eol}\vss}\rg%
%
%
\rx{\vss\hfull{%
\rlx{\hss{$4200_{x}$}}\cg%
\e{0}%
\e{0}%
\e{0}%
\e{0}%
\e{0}%
\e{0}%
\e{0}%
\e{0}%
\e{0}%
\e{0}%
\e{0}%
\e{0}%
\e{0}%
\e{0}%
\e{0}%
\e{0}%
\e{0}%
\e{0}%
\e{0}%
\e{0}%
\e{0}%
\e{0}%
\e{0}%
\e{0}%
\eol}\vss}\rg%
%
%
\rx{\vss\hfull{%
\rlx{\hss{$2240_{x}$}}\cg%
\e{0}%
\e{0}%
\e{0}%
\e{0}%
\e{0}%
\e{0}%
\e{0}%
\e{0}%
\e{0}%
\e{0}%
\e{0}%
\e{0}%
\e{0}%
\e{0}%
\e{0}%
\e{0}%
\e{0}%
\e{0}%
\e{0}%
\e{0}%
\e{0}%
\e{0}%
\e{0}%
\e{0}%
\eol}\vss}\rg%
%
%
\rx{\vss\hfull{%
\rlx{\hss{$2835_{x}$}}\cg%
\e{0}%
\e{0}%
\e{0}%
\e{0}%
\e{0}%
\e{0}%
\e{0}%
\e{0}%
\e{0}%
\e{0}%
\e{0}%
\e{0}%
\e{0}%
\e{0}%
\e{0}%
\e{0}%
\e{0}%
\e{0}%
\e{0}%
\e{0}%
\e{0}%
\e{0}%
\e{0}%
\e{0}%
\eol}\vss}\rg%
%
%
\rx{\vss\hfull{%
\rlx{\hss{$6075_{x}$}}\cg%
\e{0}%
\e{0}%
\e{0}%
\e{0}%
\e{0}%
\e{0}%
\e{0}%
\e{0}%
\e{0}%
\e{0}%
\e{0}%
\e{0}%
\e{0}%
\e{0}%
\e{0}%
\e{0}%
\e{0}%
\e{0}%
\e{0}%
\e{0}%
\e{0}%
\e{0}%
\e{0}%
\e{0}%
\eol}\vss}\rg%
%
%
\rx{\vss\hfull{%
\rlx{\hss{$3200_{x}$}}\cg%
\e{0}%
\e{0}%
\e{0}%
\e{0}%
\e{0}%
\e{0}%
\e{0}%
\e{0}%
\e{0}%
\e{0}%
\e{0}%
\e{0}%
\e{0}%
\e{0}%
\e{0}%
\e{0}%
\e{0}%
\e{0}%
\e{0}%
\e{0}%
\e{0}%
\e{0}%
\e{0}%
\e{0}%
\eol}\vss}\rg%
%
%
\rx{\vss\hfull{%
\rlx{\hss{$70_{y}$}}\cg%
\e{0}%
\e{0}%
\e{0}%
\e{0}%
\e{0}%
\e{0}%
\e{0}%
\e{0}%
\e{0}%
\e{0}%
\e{0}%
\e{0}%
\e{0}%
\e{0}%
\e{0}%
\e{0}%
\e{0}%
\e{0}%
\e{0}%
\e{0}%
\e{0}%
\e{0}%
\e{0}%
\e{0}%
\eol}\vss}\rg%
%
%
\rx{\vss\hfull{%
\rlx{\hss{$1134_{y}$}}\cg%
\e{0}%
\e{0}%
\e{0}%
\e{0}%
\e{0}%
\e{0}%
\e{0}%
\e{0}%
\e{0}%
\e{0}%
\e{0}%
\e{0}%
\e{0}%
\e{0}%
\e{0}%
\e{0}%
\e{0}%
\e{0}%
\e{0}%
\e{0}%
\e{0}%
\e{0}%
\e{0}%
\e{0}%
\eol}\vss}\rg%
%
%
\rx{\vss\hfull{%
\rlx{\hss{$1680_{y}$}}\cg%
\e{0}%
\e{0}%
\e{0}%
\e{0}%
\e{0}%
\e{0}%
\e{0}%
\e{0}%
\e{0}%
\e{0}%
\e{0}%
\e{0}%
\e{0}%
\e{0}%
\e{0}%
\e{0}%
\e{0}%
\e{0}%
\e{0}%
\e{0}%
\e{0}%
\e{0}%
\e{0}%
\e{0}%
\eol}\vss}\rg%
%
%
\rx{\vss\hfull{%
\rlx{\hss{$168_{y}$}}\cg%
\e{0}%
\e{0}%
\e{0}%
\e{0}%
\e{0}%
\e{0}%
\e{0}%
\e{0}%
\e{0}%
\e{0}%
\e{0}%
\e{0}%
\e{0}%
\e{0}%
\e{0}%
\e{0}%
\e{0}%
\e{0}%
\e{0}%
\e{0}%
\e{0}%
\e{0}%
\e{0}%
\e{0}%
\eol}\vss}\rg%
%
%
\rx{\vss\hfull{%
\rlx{\hss{$420_{y}$}}\cg%
\e{0}%
\e{0}%
\e{0}%
\e{0}%
\e{0}%
\e{0}%
\e{0}%
\e{0}%
\e{0}%
\e{0}%
\e{0}%
\e{0}%
\e{0}%
\e{0}%
\e{0}%
\e{0}%
\e{0}%
\e{0}%
\e{0}%
\e{0}%
\e{0}%
\e{0}%
\e{0}%
\e{0}%
\eol}\vss}\rg%
%
%
\rx{\vss\hfull{%
\rlx{\hss{$3150_{y}$}}\cg%
\e{0}%
\e{0}%
\e{0}%
\e{0}%
\e{0}%
\e{0}%
\e{0}%
\e{0}%
\e{0}%
\e{0}%
\e{0}%
\e{0}%
\e{0}%
\e{0}%
\e{0}%
\e{0}%
\e{0}%
\e{0}%
\e{0}%
\e{0}%
\e{0}%
\e{0}%
\e{0}%
\e{0}%
\eol}\vss}\rg%
%
%
\rx{\vss\hfull{%
\rlx{\hss{$4200_{y}$}}\cg%
\e{0}%
\e{0}%
\e{0}%
\e{0}%
\e{0}%
\e{0}%
\e{0}%
\e{0}%
\e{0}%
\e{0}%
\e{0}%
\e{0}%
\e{0}%
\e{0}%
\e{0}%
\e{0}%
\e{0}%
\e{0}%
\e{0}%
\e{0}%
\e{0}%
\e{0}%
\e{0}%
\e{0}%
\eol}\vss}\rg%
\eop
\eject
\tablecont%
%
%
%
%
%
%
\rowpts=18 true pt%
\colpts=18 true pt%
\rowlabpts=40 true pt%
\collabpts=70 true pt%
\clx{\vss\hfull{%
\rlx{\hss{$ $}}\cg%
\cx{\hskip 16 true pt\flip{$120_a^{*}{\times}[{2}]$}\hss}\cg%
\cx{\hskip 16 true pt\flip{$210_a^{*}{\times}[{2}]$}\hss}\cg%
\cx{\hskip 16 true pt\flip{$280_a^{*}{\times}[{2}]$}\hss}\cg%
\cx{\hskip 16 true pt\flip{$336_a^{*}{\times}[{2}]$}\hss}\cg%
\cx{\hskip 16 true pt\flip{$216_a^{*}{\times}[{2}]$}\hss}\cg%
\cx{\hskip 16 true pt\flip{$512_a^{*}{\times}[{2}]$}\hss}\cg%
\cx{\hskip 16 true pt\flip{$378_a^{*}{\times}[{2}]$}\hss}\cg%
\cx{\hskip 16 true pt\flip{$84_a^{*}{\times}[{2}]$}\hss}\cg%
\cx{\hskip 16 true pt\flip{$420_a^{*}{\times}[{2}]$}\hss}\cg%
\cx{\hskip 16 true pt\flip{$280_b^{*}{\times}[{2}]$}\hss}\cg%
\cx{\hskip 16 true pt\flip{$210_b^{*}{\times}[{2}]$}\hss}\cg%
\cx{\hskip 16 true pt\flip{$70_a^{*}{\times}[{2}]$}\hss}\cg%
\cx{\hskip 16 true pt\flip{$1_a{\times}[{1^{2}}]$}\hss}\cg%
\cx{\hskip 16 true pt\flip{$7_a{\times}[{1^{2}}]$}\hss}\cg%
\cx{\hskip 16 true pt\flip{$27_a{\times}[{1^{2}}]$}\hss}\cg%
\cx{\hskip 16 true pt\flip{$21_a{\times}[{1^{2}}]$}\hss}\cg%
\cx{\hskip 16 true pt\flip{$35_a{\times}[{1^{2}}]$}\hss}\cg%
\cx{\hskip 16 true pt\flip{$105_a{\times}[{1^{2}}]$}\hss}\cg%
\cx{\hskip 16 true pt\flip{$189_a{\times}[{1^{2}}]$}\hss}\cg%
\cx{\hskip 16 true pt\flip{$21_b{\times}[{1^{2}}]$}\hss}\cg%
\cx{\hskip 16 true pt\flip{$35_b{\times}[{1^{2}}]$}\hss}\cg%
\cx{\hskip 16 true pt\flip{$189_b{\times}[{1^{2}}]$}\hss}\cg%
\cx{\hskip 16 true pt\flip{$189_c{\times}[{1^{2}}]$}\hss}\cg%
\cx{\hskip 16 true pt\flip{$15_a{\times}[{1^{2}}]$}\hss}\cg%
\eol}}\rg%
%
%
\rx{\vss\hfull{%
\rlx{\hss{$2688_{y}$}}\cg%
\e{0}%
\e{0}%
\e{0}%
\e{0}%
\e{0}%
\e{0}%
\e{0}%
\e{0}%
\e{0}%
\e{0}%
\e{0}%
\e{0}%
\e{0}%
\e{0}%
\e{0}%
\e{0}%
\e{0}%
\e{0}%
\e{0}%
\e{0}%
\e{0}%
\e{0}%
\e{0}%
\e{0}%
\eol}\vss}\rg%
%
%
\rx{\vss\hfull{%
\rlx{\hss{$2100_{y}$}}\cg%
\e{0}%
\e{0}%
\e{0}%
\e{0}%
\e{0}%
\e{0}%
\e{0}%
\e{0}%
\e{0}%
\e{0}%
\e{0}%
\e{0}%
\e{0}%
\e{0}%
\e{0}%
\e{0}%
\e{0}%
\e{0}%
\e{0}%
\e{0}%
\e{0}%
\e{0}%
\e{0}%
\e{0}%
\eol}\vss}\rg%
%
%
\rx{\vss\hfull{%
\rlx{\hss{$1400_{y}$}}\cg%
\e{0}%
\e{0}%
\e{0}%
\e{0}%
\e{0}%
\e{0}%
\e{0}%
\e{0}%
\e{0}%
\e{0}%
\e{0}%
\e{0}%
\e{0}%
\e{0}%
\e{0}%
\e{0}%
\e{0}%
\e{0}%
\e{0}%
\e{0}%
\e{0}%
\e{0}%
\e{0}%
\e{0}%
\eol}\vss}\rg%
%
%
\rx{\vss\hfull{%
\rlx{\hss{$4536_{y}$}}\cg%
\e{0}%
\e{0}%
\e{0}%
\e{0}%
\e{0}%
\e{0}%
\e{0}%
\e{0}%
\e{0}%
\e{0}%
\e{0}%
\e{0}%
\e{0}%
\e{0}%
\e{0}%
\e{0}%
\e{0}%
\e{0}%
\e{0}%
\e{0}%
\e{0}%
\e{0}%
\e{0}%
\e{0}%
\eol}\vss}\rg%
%
%
\rx{\vss\hfull{%
\rlx{\hss{$5670_{y}$}}\cg%
\e{0}%
\e{0}%
\e{0}%
\e{0}%
\e{0}%
\e{0}%
\e{0}%
\e{0}%
\e{0}%
\e{0}%
\e{0}%
\e{0}%
\e{0}%
\e{0}%
\e{0}%
\e{0}%
\e{0}%
\e{0}%
\e{0}%
\e{0}%
\e{0}%
\e{0}%
\e{0}%
\e{0}%
\eol}\vss}\rg%
%
%
\rx{\vss\hfull{%
\rlx{\hss{$4480_{y}$}}\cg%
\e{0}%
\e{0}%
\e{0}%
\e{0}%
\e{0}%
\e{0}%
\e{0}%
\e{0}%
\e{0}%
\e{0}%
\e{0}%
\e{0}%
\e{0}%
\e{0}%
\e{0}%
\e{0}%
\e{0}%
\e{0}%
\e{0}%
\e{0}%
\e{0}%
\e{0}%
\e{0}%
\e{0}%
\eol}\vss}\rg%
%
%
\rx{\vss\hfull{%
\rlx{\hss{$8_{z}$}}\cg%
\e{0}%
\e{0}%
\e{0}%
\e{0}%
\e{0}%
\e{0}%
\e{0}%
\e{0}%
\e{0}%
\e{0}%
\e{0}%
\e{0}%
\e{1}%
\e{0}%
\e{0}%
\e{0}%
\e{0}%
\e{0}%
\e{0}%
\e{0}%
\e{0}%
\e{0}%
\e{0}%
\e{0}%
\eol}\vss}\rg%
%
%
\rx{\vss\hfull{%
\rlx{\hss{$56_{z}$}}\cg%
\e{0}%
\e{0}%
\e{0}%
\e{0}%
\e{0}%
\e{0}%
\e{0}%
\e{0}%
\e{0}%
\e{0}%
\e{0}%
\e{0}%
\e{0}%
\e{0}%
\e{0}%
\e{1}%
\e{0}%
\e{0}%
\e{0}%
\e{0}%
\e{0}%
\e{0}%
\e{0}%
\e{0}%
\eol}\vss}\rg%
%
%
\rx{\vss\hfull{%
\rlx{\hss{$160_{z}$}}\cg%
\e{0}%
\e{0}%
\e{0}%
\e{0}%
\e{0}%
\e{0}%
\e{0}%
\e{0}%
\e{0}%
\e{0}%
\e{0}%
\e{0}%
\e{0}%
\e{0}%
\e{1}%
\e{1}%
\e{0}%
\e{0}%
\e{0}%
\e{0}%
\e{0}%
\e{0}%
\e{0}%
\e{0}%
\eol}\vss}\rg%
%
%
\rx{\vss\hfull{%
\rlx{\hss{$112_{z}$}}\cg%
\e{0}%
\e{0}%
\e{0}%
\e{0}%
\e{0}%
\e{0}%
\e{0}%
\e{0}%
\e{0}%
\e{0}%
\e{0}%
\e{0}%
\e{1}%
\e{0}%
\e{1}%
\e{0}%
\e{0}%
\e{0}%
\e{0}%
\e{0}%
\e{0}%
\e{0}%
\e{0}%
\e{0}%
\eol}\vss}\rg%
%
%
\rx{\vss\hfull{%
\rlx{\hss{$840_{z}$}}\cg%
\e{0}%
\e{0}%
\e{0}%
\e{0}%
\e{0}%
\e{0}%
\e{1}%
\e{0}%
\e{0}%
\e{0}%
\e{0}%
\e{0}%
\e{0}%
\e{0}%
\e{0}%
\e{0}%
\e{0}%
\e{0}%
\e{1}%
\e{0}%
\e{0}%
\e{0}%
\e{0}%
\e{0}%
\eol}\vss}\rg%
%
%
\rx{\vss\hfull{%
\rlx{\hss{$1296_{z}$}}\cg%
\e{0}%
\e{0}%
\e{1}%
\e{1}%
\e{0}%
\e{0}%
\e{0}%
\e{0}%
\e{0}%
\e{0}%
\e{0}%
\e{0}%
\e{0}%
\e{0}%
\e{0}%
\e{1}%
\e{0}%
\e{0}%
\e{1}%
\e{0}%
\e{0}%
\e{0}%
\e{0}%
\e{0}%
\eol}\vss}\rg%
%
%
\rx{\vss\hfull{%
\rlx{\hss{$1400_{z}$}}\cg%
\e{0}%
\e{0}%
\e{0}%
\e{0}%
\e{0}%
\e{0}%
\e{0}%
\e{0}%
\e{0}%
\e{0}%
\e{0}%
\e{0}%
\e{0}%
\e{0}%
\e{1}%
\e{0}%
\e{0}%
\e{0}%
\e{0}%
\e{0}%
\e{0}%
\e{0}%
\e{0}%
\e{0}%
\eol}\vss}\rg%
%
%
\rx{\vss\hfull{%
\rlx{\hss{$1008_{z}$}}\cg%
\e{0}%
\e{0}%
\e{1}%
\e{0}%
\e{0}%
\e{0}%
\e{0}%
\e{0}%
\e{0}%
\e{0}%
\e{0}%
\e{0}%
\e{0}%
\e{0}%
\e{1}%
\e{1}%
\e{0}%
\e{0}%
\e{0}%
\e{0}%
\e{0}%
\e{0}%
\e{0}%
\e{0}%
\eol}\vss}\rg%
%
%
\rx{\vss\hfull{%
\rlx{\hss{$560_{z}$}}\cg%
\e{0}%
\e{0}%
\e{0}%
\e{0}%
\e{0}%
\e{0}%
\e{0}%
\e{0}%
\e{0}%
\e{0}%
\e{0}%
\e{0}%
\e{0}%
\e{0}%
\e{1}%
\e{0}%
\e{0}%
\e{0}%
\e{0}%
\e{0}%
\e{1}%
\e{0}%
\e{0}%
\e{0}%
\eol}\vss}\rg%
%
%
\rx{\vss\hfull{%
\rlx{\hss{$1400_{zz}$}}\cg%
\e{0}%
\e{0}%
\e{0}%
\e{0}%
\e{1}%
\e{0}%
\e{0}%
\e{0}%
\e{0}%
\e{0}%
\e{0}%
\e{1}%
\e{0}%
\e{0}%
\e{0}%
\e{0}%
\e{0}%
\e{0}%
\e{0}%
\e{0}%
\e{0}%
\e{0}%
\e{0}%
\e{0}%
\eol}\vss}\rg%
%
%
\rx{\vss\hfull{%
\rlx{\hss{$4200_{z}$}}\cg%
\e{0}%
\e{0}%
\e{1}%
\e{0}%
\e{1}%
\e{1}%
\e{1}%
\e{1}%
\e{1}%
\e{0}%
\e{0}%
\e{0}%
\e{0}%
\e{0}%
\e{0}%
\e{0}%
\e{0}%
\e{0}%
\e{0}%
\e{0}%
\e{0}%
\e{0}%
\e{0}%
\e{0}%
\eol}\vss}\rg%
%
%
\rx{\vss\hfull{%
\rlx{\hss{$400_{z}$}}\cg%
\e{0}%
\e{0}%
\e{0}%
\e{0}%
\e{0}%
\e{0}%
\e{0}%
\e{0}%
\e{0}%
\e{0}%
\e{0}%
\e{0}%
\e{0}%
\e{0}%
\e{0}%
\e{0}%
\e{0}%
\e{0}%
\e{0}%
\e{0}%
\e{1}%
\e{0}%
\e{0}%
\e{0}%
\eol}\vss}\rg%
%
%
\rx{\vss\hfull{%
\rlx{\hss{$3240_{z}$}}\cg%
\e{0}%
\e{0}%
\e{1}%
\e{0}%
\e{1}%
\e{0}%
\e{1}%
\e{0}%
\e{0}%
\e{0}%
\e{0}%
\e{0}%
\e{0}%
\e{0}%
\e{0}%
\e{0}%
\e{0}%
\e{0}%
\e{0}%
\e{0}%
\e{1}%
\e{0}%
\e{0}%
\e{0}%
\eol}\vss}\rg%
%
%
\rx{\vss\hfull{%
\rlx{\hss{$4536_{z}$}}\cg%
\e{0}%
\e{0}%
\e{0}%
\e{0}%
\e{1}%
\e{1}%
\e{2}%
\e{0}%
\e{0}%
\e{0}%
\e{1}%
\e{1}%
\e{0}%
\e{0}%
\e{0}%
\e{0}%
\e{0}%
\e{0}%
\e{0}%
\e{0}%
\e{0}%
\e{0}%
\e{0}%
\e{0}%
\eol}\vss}\rg%
%
%
\rx{\vss\hfull{%
\rlx{\hss{$2400_{z}$}}\cg%
\e{0}%
\e{0}%
\e{1}%
\e{1}%
\e{0}%
\e{0}%
\e{0}%
\e{0}%
\e{1}%
\e{0}%
\e{0}%
\e{0}%
\e{0}%
\e{0}%
\e{0}%
\e{0}%
\e{0}%
\e{0}%
\e{1}%
\e{0}%
\e{0}%
\e{0}%
\e{0}%
\e{0}%
\eol}\vss}\rg%
%
%
\rx{\vss\hfull{%
\rlx{\hss{$3360_{z}$}}\cg%
\e{0}%
\e{0}%
\e{1}%
\e{1}%
\e{1}%
\e{1}%
\e{0}%
\e{0}%
\e{0}%
\e{0}%
\e{0}%
\e{0}%
\e{0}%
\e{0}%
\e{0}%
\e{0}%
\e{0}%
\e{0}%
\e{0}%
\e{0}%
\e{0}%
\e{0}%
\e{0}%
\e{0}%
\eol}\vss}\rg%
%
%
\rx{\vss\hfull{%
\rlx{\hss{$2800_{z}$}}\cg%
\e{0}%
\e{0}%
\e{1}%
\e{1}%
\e{0}%
\e{0}%
\e{0}%
\e{0}%
\e{1}%
\e{0}%
\e{0}%
\e{0}%
\e{0}%
\e{0}%
\e{0}%
\e{0}%
\e{0}%
\e{0}%
\e{0}%
\e{0}%
\e{0}%
\e{0}%
\e{0}%
\e{0}%
\eol}\vss}\rg%
%
%
\rx{\vss\hfull{%
\rlx{\hss{$4096_{z}$}}\cg%
\e{0}%
\e{0}%
\e{1}%
\e{1}%
\e{0}%
\e{1}%
\e{1}%
\e{0}%
\e{0}%
\e{0}%
\e{0}%
\e{0}%
\e{0}%
\e{0}%
\e{0}%
\e{0}%
\e{0}%
\e{0}%
\e{1}%
\e{0}%
\e{0}%
\e{0}%
\e{0}%
\e{0}%
\eol}\vss}\rg%
%
%
\rx{\vss\hfull{%
\rlx{\hss{$5600_{z}$}}\cg%
\e{0}%
\e{0}%
\e{1}%
\e{1}%
\e{0}%
\e{1}%
\e{1}%
\e{0}%
\e{1}%
\e{0}%
\e{0}%
\e{0}%
\e{0}%
\e{0}%
\e{0}%
\e{0}%
\e{0}%
\e{0}%
\e{1}%
\e{0}%
\e{0}%
\e{0}%
\e{0}%
\e{0}%
\eol}\vss}\rg%
%
%
\rx{\vss\hfull{%
\rlx{\hss{$448_{z}$}}\cg%
\e{0}%
\e{0}%
\e{0}%
\e{0}%
\e{0}%
\e{0}%
\e{0}%
\e{0}%
\e{0}%
\e{0}%
\e{0}%
\e{1}%
\e{0}%
\e{0}%
\e{0}%
\e{0}%
\e{0}%
\e{0}%
\e{0}%
\e{0}%
\e{0}%
\e{0}%
\e{0}%
\e{0}%
\eol}\vss}\rg%
%
%
\rx{\vss\hfull{%
\rlx{\hss{$448_{w}$}}\cg%
\e{0}%
\e{0}%
\e{0}%
\e{0}%
\e{0}%
\e{0}%
\e{0}%
\e{0}%
\e{0}%
\e{0}%
\e{0}%
\e{0}%
\e{0}%
\e{0}%
\e{0}%
\e{0}%
\e{1}%
\e{0}%
\e{1}%
\e{0}%
\e{0}%
\e{0}%
\e{0}%
\e{0}%
\eol}\vss}\rg%
%
%
\rx{\vss\hfull{%
\rlx{\hss{$1344_{w}$}}\cg%
\e{0}%
\e{0}%
\e{0}%
\e{0}%
\e{0}%
\e{0}%
\e{1}%
\e{1}%
\e{0}%
\e{0}%
\e{1}%
\e{0}%
\e{0}%
\e{0}%
\e{0}%
\e{0}%
\e{0}%
\e{0}%
\e{0}%
\e{0}%
\e{0}%
\e{0}%
\e{0}%
\e{0}%
\eol}\vss}\rg%
%
%
\rx{\vss\hfull{%
\rlx{\hss{$5600_{w}$}}\cg%
\e{0}%
\e{0}%
\e{1}%
\e{1}%
\e{0}%
\e{1}%
\e{1}%
\e{0}%
\e{1}%
\e{1}%
\e{0}%
\e{0}%
\e{0}%
\e{0}%
\e{0}%
\e{0}%
\e{0}%
\e{0}%
\e{1}%
\e{0}%
\e{0}%
\e{0}%
\e{0}%
\e{0}%
\eol}\vss}\rg%
%
%
\rx{\vss\hfull{%
\rlx{\hss{$2016_{w}$}}\cg%
\e{0}%
\e{0}%
\e{0}%
\e{0}%
\e{1}%
\e{1}%
\e{0}%
\e{0}%
\e{0}%
\e{0}%
\e{1}%
\e{1}%
\e{0}%
\e{0}%
\e{0}%
\e{0}%
\e{0}%
\e{0}%
\e{0}%
\e{0}%
\e{0}%
\e{0}%
\e{0}%
\e{0}%
\eol}\vss}\rg%
%
%
\rx{\vss\hfull{%
\rlx{\hss{$7168_{w}$}}\cg%
\e{0}%
\e{0}%
\e{0}%
\e{1}%
\e{1}%
\e{2}%
\e{1}%
\e{0}%
\e{1}%
\e{1}%
\e{1}%
\e{0}%
\e{0}%
\e{0}%
\e{0}%
\e{0}%
\e{0}%
\e{0}%
\e{0}%
\e{0}%
\e{0}%
\e{0}%
\e{0}%
\e{0}%
\eol}\vss}\rg%
\eop
\eject
\tablecont%
%
%
%
%
%
%
\rowpts=18 true pt%
\colpts=18 true pt%
\rowlabpts=40 true pt%
\collabpts=70 true pt%
\clx{\vss\hfull{%
\rlx{\hss{$ $}}\cg%
\cx{\hskip 16 true pt\flip{$105_b{\times}[{1^{2}}]$}\hss}\cg%
\cx{\hskip 16 true pt\flip{$105_c{\times}[{1^{2}}]$}\hss}\cg%
\cx{\hskip 16 true pt\flip{$315_a{\times}[{1^{2}}]$}\hss}\cg%
\cx{\hskip 16 true pt\flip{$405_a{\times}[{1^{2}}]$}\hss}\cg%
\cx{\hskip 16 true pt\flip{$168_a{\times}[{1^{2}}]$}\hss}\cg%
\cx{\hskip 16 true pt\flip{$56_a{\times}[{1^{2}}]$}\hss}\cg%
\cx{\hskip 16 true pt\flip{$120_a{\times}[{1^{2}}]$}\hss}\cg%
\cx{\hskip 16 true pt\flip{$210_a{\times}[{1^{2}}]$}\hss}\cg%
\cx{\hskip 16 true pt\flip{$280_a{\times}[{1^{2}}]$}\hss}\cg%
\cx{\hskip 16 true pt\flip{$336_a{\times}[{1^{2}}]$}\hss}\cg%
\cx{\hskip 16 true pt\flip{$216_a{\times}[{1^{2}}]$}\hss}\cg%
\cx{\hskip 16 true pt\flip{$512_a{\times}[{1^{2}}]$}\hss}\cg%
\cx{\hskip 16 true pt\flip{$378_a{\times}[{1^{2}}]$}\hss}\cg%
\cx{\hskip 16 true pt\flip{$84_a{\times}[{1^{2}}]$}\hss}\cg%
\cx{\hskip 16 true pt\flip{$420_a{\times}[{1^{2}}]$}\hss}\cg%
\cx{\hskip 16 true pt\flip{$280_b{\times}[{1^{2}}]$}\hss}\cg%
\cx{\hskip 16 true pt\flip{$210_b{\times}[{1^{2}}]$}\hss}\cg%
\cx{\hskip 16 true pt\flip{$70_a{\times}[{1^{2}}]$}\hss}\cg%
\cx{\hskip 16 true pt\flip{$1_a^{*}{\times}[{1^{2}}]$}\hss}\cg%
\cx{\hskip 16 true pt\flip{$7_a^{*}{\times}[{1^{2}}]$}\hss}\cg%
\cx{\hskip 16 true pt\flip{$27_a^{*}{\times}[{1^{2}}]$}\hss}\cg%
\cx{\hskip 16 true pt\flip{$21_a^{*}{\times}[{1^{2}}]$}\hss}\cg%
\cx{\hskip 16 true pt\flip{$35_a^{*}{\times}[{1^{2}}]$}\hss}\cg%
\cx{\hskip 16 true pt\flip{$105_a^{*}{\times}[{1^{2}}]$}\hss}\cg%
\eol}}\rg%
%
%
\rx{\vss\hfull{%
\rlx{\hss{$1_{x}$}}\cg%
\e{0}%
\e{0}%
\e{0}%
\e{0}%
\e{0}%
\e{0}%
\e{0}%
\e{0}%
\e{0}%
\e{0}%
\e{0}%
\e{0}%
\e{0}%
\e{0}%
\e{0}%
\e{0}%
\e{0}%
\e{0}%
\e{0}%
\e{0}%
\e{0}%
\e{0}%
\e{0}%
\e{0}%
\eol}\vss}\rg%
%
%
\rx{\vss\hfull{%
\rlx{\hss{$28_{x}$}}\cg%
\e{0}%
\e{0}%
\e{0}%
\e{0}%
\e{0}%
\e{0}%
\e{0}%
\e{0}%
\e{0}%
\e{0}%
\e{0}%
\e{0}%
\e{0}%
\e{0}%
\e{0}%
\e{0}%
\e{0}%
\e{0}%
\e{0}%
\e{1}%
\e{0}%
\e{0}%
\e{0}%
\e{0}%
\eol}\vss}\rg%
%
%
\rx{\vss\hfull{%
\rlx{\hss{$35_{x}$}}\cg%
\e{0}%
\e{0}%
\e{0}%
\e{0}%
\e{0}%
\e{0}%
\e{0}%
\e{0}%
\e{0}%
\e{0}%
\e{0}%
\e{0}%
\e{0}%
\e{0}%
\e{0}%
\e{0}%
\e{0}%
\e{0}%
\e{0}%
\e{1}%
\e{0}%
\e{0}%
\e{0}%
\e{0}%
\eol}\vss}\rg%
%
%
\rx{\vss\hfull{%
\rlx{\hss{$84_{x}$}}\cg%
\e{0}%
\e{0}%
\e{0}%
\e{0}%
\e{0}%
\e{0}%
\e{0}%
\e{0}%
\e{0}%
\e{0}%
\e{0}%
\e{0}%
\e{0}%
\e{0}%
\e{0}%
\e{0}%
\e{0}%
\e{0}%
\e{0}%
\e{0}%
\e{0}%
\e{0}%
\e{0}%
\e{0}%
\eol}\vss}\rg%
%
%
\rx{\vss\hfull{%
\rlx{\hss{$50_{x}$}}\cg%
\e{0}%
\e{0}%
\e{0}%
\e{0}%
\e{0}%
\e{0}%
\e{0}%
\e{0}%
\e{0}%
\e{0}%
\e{0}%
\e{0}%
\e{0}%
\e{0}%
\e{0}%
\e{0}%
\e{0}%
\e{0}%
\e{0}%
\e{0}%
\e{0}%
\e{0}%
\e{0}%
\e{0}%
\eol}\vss}\rg%
%
%
\rx{\vss\hfull{%
\rlx{\hss{$350_{x}$}}\cg%
\e{0}%
\e{0}%
\e{0}%
\e{0}%
\e{0}%
\e{0}%
\e{0}%
\e{0}%
\e{0}%
\e{0}%
\e{0}%
\e{0}%
\e{0}%
\e{0}%
\e{0}%
\e{0}%
\e{0}%
\e{0}%
\e{0}%
\e{0}%
\e{0}%
\e{0}%
\e{1}%
\e{1}%
\eol}\vss}\rg%
%
%
\rx{\vss\hfull{%
\rlx{\hss{$300_{x}$}}\cg%
\e{0}%
\e{0}%
\e{0}%
\e{0}%
\e{0}%
\e{0}%
\e{0}%
\e{0}%
\e{0}%
\e{0}%
\e{0}%
\e{0}%
\e{0}%
\e{0}%
\e{0}%
\e{0}%
\e{0}%
\e{0}%
\e{0}%
\e{0}%
\e{0}%
\e{0}%
\e{0}%
\e{1}%
\eol}\vss}\rg%
%
%
\rx{\vss\hfull{%
\rlx{\hss{$567_{x}$}}\cg%
\e{0}%
\e{0}%
\e{0}%
\e{0}%
\e{0}%
\e{0}%
\e{0}%
\e{0}%
\e{0}%
\e{0}%
\e{0}%
\e{0}%
\e{0}%
\e{0}%
\e{0}%
\e{0}%
\e{0}%
\e{0}%
\e{0}%
\e{1}%
\e{0}%
\e{0}%
\e{0}%
\e{1}%
\eol}\vss}\rg%
%
%
\rx{\vss\hfull{%
\rlx{\hss{$210_{x}$}}\cg%
\e{0}%
\e{0}%
\e{0}%
\e{0}%
\e{0}%
\e{0}%
\e{0}%
\e{0}%
\e{0}%
\e{0}%
\e{0}%
\e{0}%
\e{0}%
\e{0}%
\e{0}%
\e{0}%
\e{0}%
\e{0}%
\e{0}%
\e{1}%
\e{0}%
\e{0}%
\e{0}%
\e{0}%
\eol}\vss}\rg%
%
%
\rx{\vss\hfull{%
\rlx{\hss{$840_{x}$}}\cg%
\e{0}%
\e{0}%
\e{0}%
\e{0}%
\e{0}%
\e{0}%
\e{0}%
\e{0}%
\e{0}%
\e{0}%
\e{0}%
\e{0}%
\e{0}%
\e{0}%
\e{0}%
\e{0}%
\e{0}%
\e{0}%
\e{0}%
\e{0}%
\e{0}%
\e{0}%
\e{0}%
\e{0}%
\eol}\vss}\rg%
%
%
\rx{\vss\hfull{%
\rlx{\hss{$700_{x}$}}\cg%
\e{0}%
\e{0}%
\e{0}%
\e{0}%
\e{0}%
\e{0}%
\e{0}%
\e{0}%
\e{0}%
\e{0}%
\e{0}%
\e{0}%
\e{0}%
\e{0}%
\e{0}%
\e{0}%
\e{0}%
\e{0}%
\e{0}%
\e{0}%
\e{0}%
\e{0}%
\e{0}%
\e{0}%
\eol}\vss}\rg%
%
%
\rx{\vss\hfull{%
\rlx{\hss{$175_{x}$}}\cg%
\e{0}%
\e{0}%
\e{0}%
\e{0}%
\e{0}%
\e{0}%
\e{0}%
\e{0}%
\e{0}%
\e{0}%
\e{0}%
\e{0}%
\e{0}%
\e{0}%
\e{0}%
\e{0}%
\e{0}%
\e{0}%
\e{0}%
\e{0}%
\e{0}%
\e{0}%
\e{0}%
\e{0}%
\eol}\vss}\rg%
%
%
\rx{\vss\hfull{%
\rlx{\hss{$1400_{x}$}}\cg%
\e{0}%
\e{0}%
\e{0}%
\e{0}%
\e{0}%
\e{0}%
\e{0}%
\e{0}%
\e{0}%
\e{0}%
\e{0}%
\e{0}%
\e{0}%
\e{0}%
\e{0}%
\e{0}%
\e{0}%
\e{0}%
\e{0}%
\e{0}%
\e{0}%
\e{0}%
\e{0}%
\e{0}%
\eol}\vss}\rg%
%
%
\rx{\vss\hfull{%
\rlx{\hss{$1050_{x}$}}\cg%
\e{0}%
\e{0}%
\e{0}%
\e{0}%
\e{0}%
\e{0}%
\e{0}%
\e{0}%
\e{0}%
\e{0}%
\e{0}%
\e{0}%
\e{0}%
\e{0}%
\e{0}%
\e{0}%
\e{0}%
\e{0}%
\e{0}%
\e{0}%
\e{0}%
\e{0}%
\e{0}%
\e{0}%
\eol}\vss}\rg%
%
%
\rx{\vss\hfull{%
\rlx{\hss{$1575_{x}$}}\cg%
\e{0}%
\e{0}%
\e{0}%
\e{0}%
\e{0}%
\e{0}%
\e{0}%
\e{0}%
\e{0}%
\e{0}%
\e{0}%
\e{0}%
\e{0}%
\e{0}%
\e{0}%
\e{0}%
\e{0}%
\e{0}%
\e{0}%
\e{0}%
\e{0}%
\e{0}%
\e{0}%
\e{1}%
\eol}\vss}\rg%
%
%
\rx{\vss\hfull{%
\rlx{\hss{$1344_{x}$}}\cg%
\e{0}%
\e{0}%
\e{0}%
\e{0}%
\e{0}%
\e{0}%
\e{0}%
\e{0}%
\e{0}%
\e{0}%
\e{0}%
\e{0}%
\e{0}%
\e{0}%
\e{0}%
\e{0}%
\e{0}%
\e{0}%
\e{0}%
\e{0}%
\e{0}%
\e{0}%
\e{0}%
\e{1}%
\eol}\vss}\rg%
%
%
\rx{\vss\hfull{%
\rlx{\hss{$2100_{x}$}}\cg%
\e{0}%
\e{0}%
\e{0}%
\e{0}%
\e{0}%
\e{0}%
\e{0}%
\e{0}%
\e{0}%
\e{0}%
\e{0}%
\e{0}%
\e{0}%
\e{0}%
\e{0}%
\e{0}%
\e{0}%
\e{0}%
\e{0}%
\e{0}%
\e{0}%
\e{0}%
\e{1}%
\e{1}%
\eol}\vss}\rg%
%
%
\rx{\vss\hfull{%
\rlx{\hss{$2268_{x}$}}\cg%
\e{0}%
\e{0}%
\e{0}%
\e{0}%
\e{0}%
\e{0}%
\e{0}%
\e{0}%
\e{0}%
\e{0}%
\e{0}%
\e{0}%
\e{0}%
\e{0}%
\e{0}%
\e{0}%
\e{0}%
\e{0}%
\e{0}%
\e{0}%
\e{0}%
\e{0}%
\e{0}%
\e{1}%
\eol}\vss}\rg%
%
%
\rx{\vss\hfull{%
\rlx{\hss{$525_{x}$}}\cg%
\e{0}%
\e{0}%
\e{0}%
\e{0}%
\e{0}%
\e{0}%
\e{0}%
\e{0}%
\e{0}%
\e{0}%
\e{0}%
\e{0}%
\e{0}%
\e{0}%
\e{0}%
\e{0}%
\e{0}%
\e{0}%
\e{0}%
\e{0}%
\e{0}%
\e{0}%
\e{0}%
\e{0}%
\eol}\vss}\rg%
%
%
\rx{\vss\hfull{%
\rlx{\hss{$700_{xx}$}}\cg%
\e{0}%
\e{0}%
\e{0}%
\e{0}%
\e{0}%
\e{0}%
\e{0}%
\e{0}%
\e{0}%
\e{0}%
\e{0}%
\e{0}%
\e{0}%
\e{0}%
\e{0}%
\e{0}%
\e{0}%
\e{0}%
\e{0}%
\e{0}%
\e{0}%
\e{0}%
\e{0}%
\e{0}%
\eol}\vss}\rg%
%
%
\rx{\vss\hfull{%
\rlx{\hss{$972_{x}$}}\cg%
\e{0}%
\e{0}%
\e{0}%
\e{0}%
\e{0}%
\e{0}%
\e{0}%
\e{0}%
\e{0}%
\e{0}%
\e{0}%
\e{0}%
\e{0}%
\e{0}%
\e{0}%
\e{0}%
\e{0}%
\e{0}%
\e{0}%
\e{0}%
\e{0}%
\e{0}%
\e{0}%
\e{0}%
\eol}\vss}\rg%
%
%
\rx{\vss\hfull{%
\rlx{\hss{$4096_{x}$}}\cg%
\e{0}%
\e{0}%
\e{0}%
\e{0}%
\e{0}%
\e{0}%
\e{0}%
\e{0}%
\e{0}%
\e{0}%
\e{0}%
\e{0}%
\e{0}%
\e{0}%
\e{0}%
\e{0}%
\e{0}%
\e{0}%
\e{0}%
\e{0}%
\e{0}%
\e{0}%
\e{0}%
\e{1}%
\eol}\vss}\rg%
%
%
\rx{\vss\hfull{%
\rlx{\hss{$4200_{x}$}}\cg%
\e{0}%
\e{0}%
\e{0}%
\e{0}%
\e{0}%
\e{0}%
\e{0}%
\e{0}%
\e{0}%
\e{0}%
\e{0}%
\e{0}%
\e{0}%
\e{0}%
\e{0}%
\e{0}%
\e{0}%
\e{0}%
\e{0}%
\e{0}%
\e{0}%
\e{0}%
\e{0}%
\e{0}%
\eol}\vss}\rg%
%
%
\rx{\vss\hfull{%
\rlx{\hss{$2240_{x}$}}\cg%
\e{0}%
\e{0}%
\e{0}%
\e{0}%
\e{0}%
\e{0}%
\e{0}%
\e{0}%
\e{0}%
\e{0}%
\e{0}%
\e{0}%
\e{0}%
\e{0}%
\e{0}%
\e{0}%
\e{0}%
\e{0}%
\e{0}%
\e{0}%
\e{0}%
\e{0}%
\e{0}%
\e{0}%
\eol}\vss}\rg%
%
%
\rx{\vss\hfull{%
\rlx{\hss{$2835_{x}$}}\cg%
\e{0}%
\e{0}%
\e{0}%
\e{0}%
\e{0}%
\e{0}%
\e{0}%
\e{0}%
\e{0}%
\e{0}%
\e{0}%
\e{0}%
\e{0}%
\e{0}%
\e{0}%
\e{0}%
\e{0}%
\e{0}%
\e{0}%
\e{0}%
\e{0}%
\e{0}%
\e{0}%
\e{0}%
\eol}\vss}\rg%
%
%
\rx{\vss\hfull{%
\rlx{\hss{$6075_{x}$}}\cg%
\e{0}%
\e{0}%
\e{0}%
\e{0}%
\e{0}%
\e{0}%
\e{0}%
\e{0}%
\e{0}%
\e{0}%
\e{0}%
\e{0}%
\e{0}%
\e{0}%
\e{0}%
\e{0}%
\e{0}%
\e{0}%
\e{0}%
\e{0}%
\e{0}%
\e{0}%
\e{0}%
\e{0}%
\eol}\vss}\rg%
%
%
\rx{\vss\hfull{%
\rlx{\hss{$3200_{x}$}}\cg%
\e{0}%
\e{0}%
\e{0}%
\e{0}%
\e{0}%
\e{0}%
\e{0}%
\e{0}%
\e{0}%
\e{0}%
\e{0}%
\e{0}%
\e{0}%
\e{0}%
\e{0}%
\e{0}%
\e{0}%
\e{0}%
\e{0}%
\e{0}%
\e{0}%
\e{0}%
\e{0}%
\e{0}%
\eol}\vss}\rg%
%
%
\rx{\vss\hfull{%
\rlx{\hss{$70_{y}$}}\cg%
\e{0}%
\e{0}%
\e{0}%
\e{0}%
\e{0}%
\e{0}%
\e{0}%
\e{0}%
\e{0}%
\e{0}%
\e{0}%
\e{0}%
\e{0}%
\e{0}%
\e{0}%
\e{0}%
\e{0}%
\e{0}%
\e{0}%
\e{0}%
\e{0}%
\e{0}%
\e{1}%
\e{0}%
\eol}\vss}\rg%
%
%
\rx{\vss\hfull{%
\rlx{\hss{$1134_{y}$}}\cg%
\e{0}%
\e{0}%
\e{0}%
\e{0}%
\e{0}%
\e{0}%
\e{0}%
\e{0}%
\e{0}%
\e{0}%
\e{0}%
\e{0}%
\e{0}%
\e{0}%
\e{0}%
\e{0}%
\e{0}%
\e{0}%
\e{0}%
\e{0}%
\e{0}%
\e{0}%
\e{0}%
\e{0}%
\eol}\vss}\rg%
%
%
\rx{\vss\hfull{%
\rlx{\hss{$1680_{y}$}}\cg%
\e{0}%
\e{0}%
\e{0}%
\e{0}%
\e{0}%
\e{0}%
\e{0}%
\e{0}%
\e{0}%
\e{0}%
\e{0}%
\e{0}%
\e{0}%
\e{0}%
\e{0}%
\e{0}%
\e{0}%
\e{0}%
\e{0}%
\e{0}%
\e{0}%
\e{0}%
\e{1}%
\e{0}%
\eol}\vss}\rg%
%
%
\rx{\vss\hfull{%
\rlx{\hss{$168_{y}$}}\cg%
\e{0}%
\e{0}%
\e{0}%
\e{0}%
\e{0}%
\e{0}%
\e{0}%
\e{0}%
\e{0}%
\e{0}%
\e{0}%
\e{0}%
\e{0}%
\e{0}%
\e{0}%
\e{0}%
\e{0}%
\e{0}%
\e{0}%
\e{0}%
\e{0}%
\e{0}%
\e{0}%
\e{0}%
\eol}\vss}\rg%
%
%
\rx{\vss\hfull{%
\rlx{\hss{$420_{y}$}}\cg%
\e{0}%
\e{0}%
\e{0}%
\e{0}%
\e{0}%
\e{0}%
\e{0}%
\e{0}%
\e{0}%
\e{0}%
\e{0}%
\e{0}%
\e{0}%
\e{0}%
\e{0}%
\e{0}%
\e{0}%
\e{0}%
\e{0}%
\e{0}%
\e{0}%
\e{0}%
\e{0}%
\e{0}%
\eol}\vss}\rg%
%
%
\rx{\vss\hfull{%
\rlx{\hss{$3150_{y}$}}\cg%
\e{0}%
\e{0}%
\e{0}%
\e{0}%
\e{0}%
\e{0}%
\e{0}%
\e{0}%
\e{0}%
\e{0}%
\e{0}%
\e{0}%
\e{0}%
\e{0}%
\e{0}%
\e{0}%
\e{0}%
\e{0}%
\e{0}%
\e{0}%
\e{0}%
\e{0}%
\e{0}%
\e{0}%
\eol}\vss}\rg%
%
%
\rx{\vss\hfull{%
\rlx{\hss{$4200_{y}$}}\cg%
\e{0}%
\e{0}%
\e{0}%
\e{0}%
\e{0}%
\e{0}%
\e{0}%
\e{0}%
\e{0}%
\e{0}%
\e{0}%
\e{0}%
\e{0}%
\e{0}%
\e{0}%
\e{0}%
\e{0}%
\e{0}%
\e{0}%
\e{0}%
\e{0}%
\e{0}%
\e{0}%
\e{0}%
\eol}\vss}\rg%
\eop
\eject
\tablecont%
%
%
%
%
%
%
\rowpts=18 true pt%
\colpts=18 true pt%
\rowlabpts=40 true pt%
\collabpts=70 true pt%
\clx{\vss\hfull{%
\rlx{\hss{$ $}}\cg%
\cx{\hskip 16 true pt\flip{$105_b{\times}[{1^{2}}]$}\hss}\cg%
\cx{\hskip 16 true pt\flip{$105_c{\times}[{1^{2}}]$}\hss}\cg%
\cx{\hskip 16 true pt\flip{$315_a{\times}[{1^{2}}]$}\hss}\cg%
\cx{\hskip 16 true pt\flip{$405_a{\times}[{1^{2}}]$}\hss}\cg%
\cx{\hskip 16 true pt\flip{$168_a{\times}[{1^{2}}]$}\hss}\cg%
\cx{\hskip 16 true pt\flip{$56_a{\times}[{1^{2}}]$}\hss}\cg%
\cx{\hskip 16 true pt\flip{$120_a{\times}[{1^{2}}]$}\hss}\cg%
\cx{\hskip 16 true pt\flip{$210_a{\times}[{1^{2}}]$}\hss}\cg%
\cx{\hskip 16 true pt\flip{$280_a{\times}[{1^{2}}]$}\hss}\cg%
\cx{\hskip 16 true pt\flip{$336_a{\times}[{1^{2}}]$}\hss}\cg%
\cx{\hskip 16 true pt\flip{$216_a{\times}[{1^{2}}]$}\hss}\cg%
\cx{\hskip 16 true pt\flip{$512_a{\times}[{1^{2}}]$}\hss}\cg%
\cx{\hskip 16 true pt\flip{$378_a{\times}[{1^{2}}]$}\hss}\cg%
\cx{\hskip 16 true pt\flip{$84_a{\times}[{1^{2}}]$}\hss}\cg%
\cx{\hskip 16 true pt\flip{$420_a{\times}[{1^{2}}]$}\hss}\cg%
\cx{\hskip 16 true pt\flip{$280_b{\times}[{1^{2}}]$}\hss}\cg%
\cx{\hskip 16 true pt\flip{$210_b{\times}[{1^{2}}]$}\hss}\cg%
\cx{\hskip 16 true pt\flip{$70_a{\times}[{1^{2}}]$}\hss}\cg%
\cx{\hskip 16 true pt\flip{$1_a^{*}{\times}[{1^{2}}]$}\hss}\cg%
\cx{\hskip 16 true pt\flip{$7_a^{*}{\times}[{1^{2}}]$}\hss}\cg%
\cx{\hskip 16 true pt\flip{$27_a^{*}{\times}[{1^{2}}]$}\hss}\cg%
\cx{\hskip 16 true pt\flip{$21_a^{*}{\times}[{1^{2}}]$}\hss}\cg%
\cx{\hskip 16 true pt\flip{$35_a^{*}{\times}[{1^{2}}]$}\hss}\cg%
\cx{\hskip 16 true pt\flip{$105_a^{*}{\times}[{1^{2}}]$}\hss}\cg%
\eol}}\rg%
%
%
\rx{\vss\hfull{%
\rlx{\hss{$2688_{y}$}}\cg%
\e{0}%
\e{0}%
\e{0}%
\e{0}%
\e{0}%
\e{0}%
\e{0}%
\e{0}%
\e{0}%
\e{0}%
\e{0}%
\e{0}%
\e{0}%
\e{0}%
\e{0}%
\e{0}%
\e{0}%
\e{0}%
\e{0}%
\e{0}%
\e{0}%
\e{0}%
\e{0}%
\e{0}%
\eol}\vss}\rg%
%
%
\rx{\vss\hfull{%
\rlx{\hss{$2100_{y}$}}\cg%
\e{0}%
\e{0}%
\e{0}%
\e{0}%
\e{0}%
\e{0}%
\e{0}%
\e{0}%
\e{0}%
\e{0}%
\e{0}%
\e{0}%
\e{0}%
\e{0}%
\e{0}%
\e{0}%
\e{0}%
\e{0}%
\e{0}%
\e{0}%
\e{0}%
\e{0}%
\e{0}%
\e{0}%
\eol}\vss}\rg%
%
%
\rx{\vss\hfull{%
\rlx{\hss{$1400_{y}$}}\cg%
\e{0}%
\e{0}%
\e{0}%
\e{0}%
\e{0}%
\e{0}%
\e{0}%
\e{0}%
\e{0}%
\e{0}%
\e{0}%
\e{0}%
\e{0}%
\e{0}%
\e{0}%
\e{0}%
\e{0}%
\e{0}%
\e{0}%
\e{0}%
\e{0}%
\e{0}%
\e{0}%
\e{0}%
\eol}\vss}\rg%
%
%
\rx{\vss\hfull{%
\rlx{\hss{$4536_{y}$}}\cg%
\e{0}%
\e{0}%
\e{0}%
\e{0}%
\e{0}%
\e{0}%
\e{0}%
\e{0}%
\e{0}%
\e{0}%
\e{0}%
\e{0}%
\e{0}%
\e{0}%
\e{0}%
\e{0}%
\e{0}%
\e{0}%
\e{0}%
\e{0}%
\e{0}%
\e{0}%
\e{0}%
\e{0}%
\eol}\vss}\rg%
%
%
\rx{\vss\hfull{%
\rlx{\hss{$5670_{y}$}}\cg%
\e{0}%
\e{0}%
\e{0}%
\e{0}%
\e{0}%
\e{0}%
\e{0}%
\e{0}%
\e{0}%
\e{0}%
\e{0}%
\e{0}%
\e{0}%
\e{0}%
\e{0}%
\e{0}%
\e{0}%
\e{0}%
\e{0}%
\e{0}%
\e{0}%
\e{0}%
\e{0}%
\e{0}%
\eol}\vss}\rg%
%
%
\rx{\vss\hfull{%
\rlx{\hss{$4480_{y}$}}\cg%
\e{0}%
\e{0}%
\e{0}%
\e{0}%
\e{0}%
\e{0}%
\e{0}%
\e{0}%
\e{0}%
\e{0}%
\e{0}%
\e{0}%
\e{0}%
\e{0}%
\e{0}%
\e{0}%
\e{0}%
\e{0}%
\e{0}%
\e{0}%
\e{0}%
\e{0}%
\e{0}%
\e{0}%
\eol}\vss}\rg%
%
%
\rx{\vss\hfull{%
\rlx{\hss{$8_{z}$}}\cg%
\e{0}%
\e{0}%
\e{0}%
\e{0}%
\e{0}%
\e{0}%
\e{0}%
\e{0}%
\e{0}%
\e{0}%
\e{0}%
\e{0}%
\e{0}%
\e{0}%
\e{0}%
\e{0}%
\e{0}%
\e{0}%
\e{0}%
\e{0}%
\e{0}%
\e{0}%
\e{0}%
\e{0}%
\eol}\vss}\rg%
%
%
\rx{\vss\hfull{%
\rlx{\hss{$56_{z}$}}\cg%
\e{0}%
\e{0}%
\e{0}%
\e{0}%
\e{0}%
\e{0}%
\e{0}%
\e{0}%
\e{0}%
\e{0}%
\e{0}%
\e{0}%
\e{0}%
\e{0}%
\e{0}%
\e{0}%
\e{0}%
\e{0}%
\e{0}%
\e{0}%
\e{0}%
\e{0}%
\e{0}%
\e{0}%
\eol}\vss}\rg%
%
%
\rx{\vss\hfull{%
\rlx{\hss{$160_{z}$}}\cg%
\e{0}%
\e{0}%
\e{0}%
\e{0}%
\e{0}%
\e{0}%
\e{0}%
\e{0}%
\e{0}%
\e{0}%
\e{0}%
\e{0}%
\e{0}%
\e{0}%
\e{0}%
\e{0}%
\e{0}%
\e{0}%
\e{0}%
\e{0}%
\e{0}%
\e{0}%
\e{0}%
\e{0}%
\eol}\vss}\rg%
%
%
\rx{\vss\hfull{%
\rlx{\hss{$112_{z}$}}\cg%
\e{0}%
\e{0}%
\e{0}%
\e{0}%
\e{0}%
\e{0}%
\e{0}%
\e{0}%
\e{0}%
\e{0}%
\e{0}%
\e{0}%
\e{0}%
\e{0}%
\e{0}%
\e{0}%
\e{0}%
\e{0}%
\e{0}%
\e{0}%
\e{0}%
\e{0}%
\e{0}%
\e{0}%
\eol}\vss}\rg%
%
%
\rx{\vss\hfull{%
\rlx{\hss{$840_{z}$}}\cg%
\e{0}%
\e{0}%
\e{0}%
\e{0}%
\e{1}%
\e{0}%
\e{0}%
\e{0}%
\e{0}%
\e{0}%
\e{0}%
\e{0}%
\e{0}%
\e{0}%
\e{0}%
\e{0}%
\e{0}%
\e{0}%
\e{0}%
\e{0}%
\e{0}%
\e{0}%
\e{0}%
\e{0}%
\eol}\vss}\rg%
%
%
\rx{\vss\hfull{%
\rlx{\hss{$1296_{z}$}}\cg%
\e{0}%
\e{0}%
\e{0}%
\e{0}%
\e{0}%
\e{0}%
\e{1}%
\e{1}%
\e{0}%
\e{0}%
\e{0}%
\e{0}%
\e{0}%
\e{0}%
\e{0}%
\e{0}%
\e{0}%
\e{0}%
\e{0}%
\e{0}%
\e{0}%
\e{0}%
\e{0}%
\e{0}%
\eol}\vss}\rg%
%
%
\rx{\vss\hfull{%
\rlx{\hss{$1400_{z}$}}\cg%
\e{0}%
\e{0}%
\e{0}%
\e{0}%
\e{1}%
\e{0}%
\e{1}%
\e{1}%
\e{0}%
\e{0}%
\e{0}%
\e{0}%
\e{0}%
\e{0}%
\e{0}%
\e{0}%
\e{0}%
\e{0}%
\e{0}%
\e{0}%
\e{0}%
\e{0}%
\e{0}%
\e{0}%
\eol}\vss}\rg%
%
%
\rx{\vss\hfull{%
\rlx{\hss{$1008_{z}$}}\cg%
\e{0}%
\e{0}%
\e{0}%
\e{0}%
\e{0}%
\e{0}%
\e{1}%
\e{1}%
\e{0}%
\e{0}%
\e{0}%
\e{0}%
\e{0}%
\e{0}%
\e{0}%
\e{0}%
\e{0}%
\e{0}%
\e{0}%
\e{0}%
\e{0}%
\e{0}%
\e{0}%
\e{0}%
\eol}\vss}\rg%
%
%
\rx{\vss\hfull{%
\rlx{\hss{$560_{z}$}}\cg%
\e{0}%
\e{0}%
\e{0}%
\e{0}%
\e{0}%
\e{0}%
\e{1}%
\e{0}%
\e{0}%
\e{0}%
\e{0}%
\e{0}%
\e{0}%
\e{0}%
\e{0}%
\e{0}%
\e{0}%
\e{0}%
\e{0}%
\e{0}%
\e{0}%
\e{0}%
\e{0}%
\e{0}%
\eol}\vss}\rg%
%
%
\rx{\vss\hfull{%
\rlx{\hss{$1400_{zz}$}}\cg%
\e{1}%
\e{0}%
\e{0}%
\e{0}%
\e{0}%
\e{0}%
\e{0}%
\e{0}%
\e{0}%
\e{0}%
\e{0}%
\e{0}%
\e{0}%
\e{0}%
\e{0}%
\e{1}%
\e{1}%
\e{0}%
\e{0}%
\e{0}%
\e{0}%
\e{0}%
\e{0}%
\e{0}%
\eol}\vss}\rg%
%
%
\rx{\vss\hfull{%
\rlx{\hss{$4200_{z}$}}\cg%
\e{1}%
\e{1}%
\e{0}%
\e{1}%
\e{0}%
\e{0}%
\e{0}%
\e{0}%
\e{0}%
\e{0}%
\e{0}%
\e{1}%
\e{1}%
\e{0}%
\e{0}%
\e{1}%
\e{1}%
\e{0}%
\e{0}%
\e{0}%
\e{0}%
\e{0}%
\e{0}%
\e{0}%
\eol}\vss}\rg%
%
%
\rx{\vss\hfull{%
\rlx{\hss{$400_{z}$}}\cg%
\e{1}%
\e{0}%
\e{0}%
\e{0}%
\e{0}%
\e{0}%
\e{0}%
\e{0}%
\e{0}%
\e{0}%
\e{0}%
\e{0}%
\e{0}%
\e{0}%
\e{0}%
\e{0}%
\e{0}%
\e{0}%
\e{0}%
\e{0}%
\e{0}%
\e{0}%
\e{0}%
\e{0}%
\eol}\vss}\rg%
%
%
\rx{\vss\hfull{%
\rlx{\hss{$3240_{z}$}}\cg%
\e{1}%
\e{0}%
\e{0}%
\e{1}%
\e{1}%
\e{0}%
\e{1}%
\e{1}%
\e{0}%
\e{0}%
\e{0}%
\e{0}%
\e{0}%
\e{0}%
\e{0}%
\e{1}%
\e{0}%
\e{0}%
\e{0}%
\e{0}%
\e{0}%
\e{0}%
\e{0}%
\e{0}%
\eol}\vss}\rg%
%
%
\rx{\vss\hfull{%
\rlx{\hss{$4536_{z}$}}\cg%
\e{0}%
\e{0}%
\e{0}%
\e{1}%
\e{1}%
\e{0}%
\e{0}%
\e{0}%
\e{0}%
\e{0}%
\e{0}%
\e{1}%
\e{0}%
\e{1}%
\e{1}%
\e{1}%
\e{1}%
\e{0}%
\e{0}%
\e{0}%
\e{0}%
\e{0}%
\e{0}%
\e{0}%
\eol}\vss}\rg%
%
%
\rx{\vss\hfull{%
\rlx{\hss{$2400_{z}$}}\cg%
\e{0}%
\e{0}%
\e{0}%
\e{1}%
\e{0}%
\e{0}%
\e{0}%
\e{1}%
\e{0}%
\e{1}%
\e{0}%
\e{0}%
\e{0}%
\e{0}%
\e{0}%
\e{0}%
\e{0}%
\e{0}%
\e{0}%
\e{0}%
\e{0}%
\e{0}%
\e{0}%
\e{0}%
\eol}\vss}\rg%
%
%
\rx{\vss\hfull{%
\rlx{\hss{$3360_{z}$}}\cg%
\e{1}%
\e{0}%
\e{0}%
\e{1}%
\e{0}%
\e{0}%
\e{0}%
\e{1}%
\e{0}%
\e{0}%
\e{0}%
\e{1}%
\e{0}%
\e{0}%
\e{0}%
\e{1}%
\e{0}%
\e{0}%
\e{0}%
\e{0}%
\e{0}%
\e{0}%
\e{0}%
\e{0}%
\eol}\vss}\rg%
%
%
\rx{\vss\hfull{%
\rlx{\hss{$2800_{z}$}}\cg%
\e{0}%
\e{1}%
\e{0}%
\e{1}%
\e{0}%
\e{0}%
\e{1}%
\e{1}%
\e{0}%
\e{0}%
\e{0}%
\e{0}%
\e{0}%
\e{0}%
\e{1}%
\e{0}%
\e{0}%
\e{0}%
\e{0}%
\e{0}%
\e{0}%
\e{0}%
\e{0}%
\e{0}%
\eol}\vss}\rg%
%
%
\rx{\vss\hfull{%
\rlx{\hss{$4096_{z}$}}\cg%
\e{0}%
\e{0}%
\e{0}%
\e{1}%
\e{1}%
\e{0}%
\e{1}%
\e{1}%
\e{0}%
\e{0}%
\e{0}%
\e{0}%
\e{0}%
\e{0}%
\e{1}%
\e{1}%
\e{0}%
\e{0}%
\e{0}%
\e{0}%
\e{0}%
\e{0}%
\e{0}%
\e{0}%
\eol}\vss}\rg%
%
%
\rx{\vss\hfull{%
\rlx{\hss{$5600_{z}$}}\cg%
\e{0}%
\e{0}%
\e{0}%
\e{1}%
\e{1}%
\e{0}%
\e{0}%
\e{1}%
\e{0}%
\e{1}%
\e{0}%
\e{1}%
\e{0}%
\e{0}%
\e{2}%
\e{0}%
\e{0}%
\e{0}%
\e{0}%
\e{0}%
\e{0}%
\e{0}%
\e{0}%
\e{0}%
\eol}\vss}\rg%
%
%
\rx{\vss\hfull{%
\rlx{\hss{$448_{z}$}}\cg%
\e{0}%
\e{0}%
\e{0}%
\e{0}%
\e{1}%
\e{0}%
\e{0}%
\e{0}%
\e{0}%
\e{0}%
\e{0}%
\e{0}%
\e{0}%
\e{0}%
\e{0}%
\e{0}%
\e{0}%
\e{0}%
\e{0}%
\e{0}%
\e{0}%
\e{0}%
\e{0}%
\e{0}%
\eol}\vss}\rg%
%
%
\rx{\vss\hfull{%
\rlx{\hss{$448_{w}$}}\cg%
\e{0}%
\e{0}%
\e{0}%
\e{0}%
\e{0}%
\e{0}%
\e{0}%
\e{0}%
\e{0}%
\e{0}%
\e{0}%
\e{0}%
\e{0}%
\e{0}%
\e{0}%
\e{0}%
\e{0}%
\e{0}%
\e{0}%
\e{0}%
\e{0}%
\e{0}%
\e{0}%
\e{0}%
\eol}\vss}\rg%
%
%
\rx{\vss\hfull{%
\rlx{\hss{$1344_{w}$}}\cg%
\e{0}%
\e{0}%
\e{0}%
\e{0}%
\e{0}%
\e{0}%
\e{0}%
\e{0}%
\e{0}%
\e{0}%
\e{0}%
\e{0}%
\e{1}%
\e{1}%
\e{0}%
\e{0}%
\e{1}%
\e{0}%
\e{0}%
\e{0}%
\e{0}%
\e{0}%
\e{0}%
\e{0}%
\eol}\vss}\rg%
%
%
\rx{\vss\hfull{%
\rlx{\hss{$5600_{w}$}}\cg%
\e{0}%
\e{0}%
\e{0}%
\e{1}%
\e{0}%
\e{0}%
\e{0}%
\e{0}%
\e{1}%
\e{1}%
\e{0}%
\e{1}%
\e{1}%
\e{0}%
\e{1}%
\e{1}%
\e{0}%
\e{0}%
\e{0}%
\e{0}%
\e{0}%
\e{0}%
\e{0}%
\e{0}%
\eol}\vss}\rg%
%
%
\rx{\vss\hfull{%
\rlx{\hss{$2016_{w}$}}\cg%
\e{0}%
\e{0}%
\e{0}%
\e{0}%
\e{0}%
\e{0}%
\e{0}%
\e{0}%
\e{0}%
\e{0}%
\e{1}%
\e{1}%
\e{0}%
\e{0}%
\e{0}%
\e{0}%
\e{1}%
\e{1}%
\e{0}%
\e{0}%
\e{0}%
\e{0}%
\e{0}%
\e{0}%
\eol}\vss}\rg%
%
%
\rx{\vss\hfull{%
\rlx{\hss{$7168_{w}$}}\cg%
\e{0}%
\e{0}%
\e{1}%
\e{1}%
\e{0}%
\e{0}%
\e{0}%
\e{0}%
\e{0}%
\e{1}%
\e{1}%
\e{2}%
\e{1}%
\e{0}%
\e{1}%
\e{1}%
\e{1}%
\e{0}%
\e{0}%
\e{0}%
\e{0}%
\e{0}%
\e{0}%
\e{0}%
\eol}\vss}\rg%
\eop
\eject
\tablecont%
%
%
%
%
%
%
\rowpts=18 true pt%
\colpts=18 true pt%
\rowlabpts=40 true pt%
\collabpts=70 true pt%
\clx{\vss\hfull{%
\rlx{\hss{$ $}}\cg%
\cx{\hskip 16 true pt\flip{$189_a^{*}{\times}[{1^{2}}]$}\hss}\cg%
\cx{\hskip 16 true pt\flip{$21_b^{*}{\times}[{1^{2}}]$}\hss}\cg%
\cx{\hskip 16 true pt\flip{$35_b^{*}{\times}[{1^{2}}]$}\hss}\cg%
\cx{\hskip 16 true pt\flip{$189_b^{*}{\times}[{1^{2}}]$}\hss}\cg%
\cx{\hskip 16 true pt\flip{$189_c^{*}{\times}[{1^{2}}]$}\hss}\cg%
\cx{\hskip 16 true pt\flip{$15_a^{*}{\times}[{1^{2}}]$}\hss}\cg%
\cx{\hskip 16 true pt\flip{$105_b^{*}{\times}[{1^{2}}]$}\hss}\cg%
\cx{\hskip 16 true pt\flip{$105_c^{*}{\times}[{1^{2}}]$}\hss}\cg%
\cx{\hskip 16 true pt\flip{$315_a^{*}{\times}[{1^{2}}]$}\hss}\cg%
\cx{\hskip 16 true pt\flip{$405_a^{*}{\times}[{1^{2}}]$}\hss}\cg%
\cx{\hskip 16 true pt\flip{$168_a^{*}{\times}[{1^{2}}]$}\hss}\cg%
\cx{\hskip 16 true pt\flip{$56_a^{*}{\times}[{1^{2}}]$}\hss}\cg%
\cx{\hskip 16 true pt\flip{$120_a^{*}{\times}[{1^{2}}]$}\hss}\cg%
\cx{\hskip 16 true pt\flip{$210_a^{*}{\times}[{1^{2}}]$}\hss}\cg%
\cx{\hskip 16 true pt\flip{$280_a^{*}{\times}[{1^{2}}]$}\hss}\cg%
\cx{\hskip 16 true pt\flip{$336_a^{*}{\times}[{1^{2}}]$}\hss}\cg%
\cx{\hskip 16 true pt\flip{$216_a^{*}{\times}[{1^{2}}]$}\hss}\cg%
\cx{\hskip 16 true pt\flip{$512_a^{*}{\times}[{1^{2}}]$}\hss}\cg%
\cx{\hskip 16 true pt\flip{$378_a^{*}{\times}[{1^{2}}]$}\hss}\cg%
\cx{\hskip 16 true pt\flip{$84_a^{*}{\times}[{1^{2}}]$}\hss}\cg%
\cx{\hskip 16 true pt\flip{$420_a^{*}{\times}[{1^{2}}]$}\hss}\cg%
\cx{\hskip 16 true pt\flip{$280_b^{*}{\times}[{1^{2}}]$}\hss}\cg%
\cx{\hskip 16 true pt\flip{$210_b^{*}{\times}[{1^{2}}]$}\hss}\cg%
\cx{\hskip 16 true pt\flip{$70_a^{*}{\times}[{1^{2}}]$}\hss}\cg%
\eol}}\rg%
%
%
\rx{\vss\hfull{%
\rlx{\hss{$1_{x}$}}\cg%
\e{0}%
\e{0}%
\e{0}%
\e{0}%
\e{0}%
\e{0}%
\e{0}%
\e{0}%
\e{0}%
\e{0}%
\e{0}%
\e{0}%
\e{0}%
\e{0}%
\e{0}%
\e{0}%
\e{0}%
\e{0}%
\e{0}%
\e{0}%
\e{0}%
\e{0}%
\e{0}%
\e{0}%
\eol}\vss}\rg%
%
%
\rx{\vss\hfull{%
\rlx{\hss{$28_{x}$}}\cg%
\e{0}%
\e{0}%
\e{0}%
\e{0}%
\e{0}%
\e{0}%
\e{0}%
\e{0}%
\e{0}%
\e{0}%
\e{0}%
\e{0}%
\e{0}%
\e{0}%
\e{0}%
\e{0}%
\e{0}%
\e{0}%
\e{0}%
\e{0}%
\e{0}%
\e{0}%
\e{0}%
\e{0}%
\eol}\vss}\rg%
%
%
\rx{\vss\hfull{%
\rlx{\hss{$35_{x}$}}\cg%
\e{0}%
\e{0}%
\e{0}%
\e{0}%
\e{0}%
\e{0}%
\e{0}%
\e{0}%
\e{0}%
\e{0}%
\e{0}%
\e{0}%
\e{0}%
\e{0}%
\e{0}%
\e{0}%
\e{0}%
\e{0}%
\e{0}%
\e{0}%
\e{0}%
\e{0}%
\e{0}%
\e{0}%
\eol}\vss}\rg%
%
%
\rx{\vss\hfull{%
\rlx{\hss{$84_{x}$}}\cg%
\e{0}%
\e{1}%
\e{0}%
\e{0}%
\e{0}%
\e{0}%
\e{0}%
\e{0}%
\e{0}%
\e{0}%
\e{0}%
\e{0}%
\e{0}%
\e{0}%
\e{0}%
\e{0}%
\e{0}%
\e{0}%
\e{0}%
\e{0}%
\e{0}%
\e{0}%
\e{0}%
\e{0}%
\eol}\vss}\rg%
%
%
\rx{\vss\hfull{%
\rlx{\hss{$50_{x}$}}\cg%
\e{0}%
\e{0}%
\e{0}%
\e{0}%
\e{0}%
\e{1}%
\e{0}%
\e{0}%
\e{0}%
\e{0}%
\e{0}%
\e{0}%
\e{0}%
\e{0}%
\e{0}%
\e{0}%
\e{0}%
\e{0}%
\e{0}%
\e{0}%
\e{0}%
\e{0}%
\e{0}%
\e{0}%
\eol}\vss}\rg%
%
%
\rx{\vss\hfull{%
\rlx{\hss{$350_{x}$}}\cg%
\e{0}%
\e{0}%
\e{0}%
\e{0}%
\e{0}%
\e{0}%
\e{0}%
\e{0}%
\e{0}%
\e{0}%
\e{0}%
\e{0}%
\e{0}%
\e{0}%
\e{0}%
\e{0}%
\e{0}%
\e{0}%
\e{0}%
\e{0}%
\e{0}%
\e{0}%
\e{0}%
\e{0}%
\eol}\vss}\rg%
%
%
\rx{\vss\hfull{%
\rlx{\hss{$300_{x}$}}\cg%
\e{0}%
\e{0}%
\e{0}%
\e{0}%
\e{0}%
\e{0}%
\e{0}%
\e{0}%
\e{0}%
\e{0}%
\e{0}%
\e{0}%
\e{0}%
\e{0}%
\e{0}%
\e{0}%
\e{0}%
\e{0}%
\e{0}%
\e{0}%
\e{0}%
\e{0}%
\e{0}%
\e{0}%
\eol}\vss}\rg%
%
%
\rx{\vss\hfull{%
\rlx{\hss{$567_{x}$}}\cg%
\e{0}%
\e{1}%
\e{0}%
\e{0}%
\e{0}%
\e{0}%
\e{0}%
\e{0}%
\e{0}%
\e{0}%
\e{0}%
\e{1}%
\e{0}%
\e{0}%
\e{0}%
\e{0}%
\e{0}%
\e{0}%
\e{0}%
\e{0}%
\e{0}%
\e{0}%
\e{0}%
\e{0}%
\eol}\vss}\rg%
%
%
\rx{\vss\hfull{%
\rlx{\hss{$210_{x}$}}\cg%
\e{0}%
\e{0}%
\e{0}%
\e{0}%
\e{0}%
\e{0}%
\e{0}%
\e{0}%
\e{0}%
\e{0}%
\e{0}%
\e{1}%
\e{0}%
\e{0}%
\e{0}%
\e{0}%
\e{0}%
\e{0}%
\e{0}%
\e{0}%
\e{0}%
\e{0}%
\e{0}%
\e{0}%
\eol}\vss}\rg%
%
%
\rx{\vss\hfull{%
\rlx{\hss{$840_{x}$}}\cg%
\e{0}%
\e{0}%
\e{0}%
\e{0}%
\e{0}%
\e{0}%
\e{0}%
\e{0}%
\e{0}%
\e{0}%
\e{0}%
\e{0}%
\e{0}%
\e{0}%
\e{0}%
\e{0}%
\e{0}%
\e{0}%
\e{1}%
\e{0}%
\e{0}%
\e{0}%
\e{0}%
\e{0}%
\eol}\vss}\rg%
%
%
\rx{\vss\hfull{%
\rlx{\hss{$700_{x}$}}\cg%
\e{0}%
\e{0}%
\e{0}%
\e{1}%
\e{0}%
\e{0}%
\e{0}%
\e{0}%
\e{0}%
\e{0}%
\e{0}%
\e{1}%
\e{0}%
\e{0}%
\e{0}%
\e{0}%
\e{0}%
\e{0}%
\e{0}%
\e{0}%
\e{0}%
\e{0}%
\e{0}%
\e{0}%
\eol}\vss}\rg%
%
%
\rx{\vss\hfull{%
\rlx{\hss{$175_{x}$}}\cg%
\e{0}%
\e{0}%
\e{0}%
\e{0}%
\e{0}%
\e{0}%
\e{0}%
\e{0}%
\e{0}%
\e{0}%
\e{0}%
\e{0}%
\e{0}%
\e{0}%
\e{0}%
\e{0}%
\e{0}%
\e{0}%
\e{0}%
\e{0}%
\e{0}%
\e{0}%
\e{0}%
\e{1}%
\eol}\vss}\rg%
%
%
\rx{\vss\hfull{%
\rlx{\hss{$1400_{x}$}}\cg%
\e{0}%
\e{0}%
\e{0}%
\e{1}%
\e{0}%
\e{0}%
\e{0}%
\e{0}%
\e{1}%
\e{0}%
\e{0}%
\e{1}%
\e{0}%
\e{0}%
\e{0}%
\e{0}%
\e{0}%
\e{0}%
\e{0}%
\e{0}%
\e{0}%
\e{0}%
\e{0}%
\e{0}%
\eol}\vss}\rg%
%
%
\rx{\vss\hfull{%
\rlx{\hss{$1050_{x}$}}\cg%
\e{0}%
\e{0}%
\e{0}%
\e{1}%
\e{0}%
\e{1}%
\e{0}%
\e{0}%
\e{0}%
\e{0}%
\e{0}%
\e{0}%
\e{0}%
\e{0}%
\e{0}%
\e{0}%
\e{1}%
\e{0}%
\e{0}%
\e{0}%
\e{0}%
\e{0}%
\e{0}%
\e{0}%
\eol}\vss}\rg%
%
%
\rx{\vss\hfull{%
\rlx{\hss{$1575_{x}$}}\cg%
\e{0}%
\e{0}%
\e{0}%
\e{1}%
\e{0}%
\e{0}%
\e{0}%
\e{0}%
\e{0}%
\e{0}%
\e{0}%
\e{1}%
\e{0}%
\e{0}%
\e{1}%
\e{0}%
\e{0}%
\e{0}%
\e{0}%
\e{0}%
\e{0}%
\e{0}%
\e{0}%
\e{0}%
\eol}\vss}\rg%
%
%
\rx{\vss\hfull{%
\rlx{\hss{$1344_{x}$}}\cg%
\e{0}%
\e{1}%
\e{0}%
\e{1}%
\e{1}%
\e{0}%
\e{0}%
\e{0}%
\e{0}%
\e{0}%
\e{0}%
\e{0}%
\e{0}%
\e{0}%
\e{0}%
\e{0}%
\e{0}%
\e{0}%
\e{0}%
\e{0}%
\e{0}%
\e{0}%
\e{0}%
\e{0}%
\eol}\vss}\rg%
%
%
\rx{\vss\hfull{%
\rlx{\hss{$2100_{x}$}}\cg%
\e{0}%
\e{0}%
\e{0}%
\e{0}%
\e{1}%
\e{0}%
\e{0}%
\e{0}%
\e{0}%
\e{0}%
\e{0}%
\e{0}%
\e{0}%
\e{0}%
\e{1}%
\e{1}%
\e{0}%
\e{0}%
\e{0}%
\e{0}%
\e{0}%
\e{0}%
\e{0}%
\e{0}%
\eol}\vss}\rg%
%
%
\rx{\vss\hfull{%
\rlx{\hss{$2268_{x}$}}\cg%
\e{0}%
\e{0}%
\e{0}%
\e{0}%
\e{1}%
\e{0}%
\e{0}%
\e{0}%
\e{1}%
\e{0}%
\e{0}%
\e{1}%
\e{0}%
\e{0}%
\e{1}%
\e{0}%
\e{0}%
\e{0}%
\e{0}%
\e{0}%
\e{0}%
\e{0}%
\e{0}%
\e{0}%
\eol}\vss}\rg%
%
%
\rx{\vss\hfull{%
\rlx{\hss{$525_{x}$}}\cg%
\e{0}%
\e{1}%
\e{0}%
\e{0}%
\e{1}%
\e{0}%
\e{0}%
\e{0}%
\e{0}%
\e{0}%
\e{0}%
\e{0}%
\e{0}%
\e{0}%
\e{0}%
\e{0}%
\e{0}%
\e{0}%
\e{0}%
\e{0}%
\e{0}%
\e{0}%
\e{0}%
\e{0}%
\eol}\vss}\rg%
%
%
\rx{\vss\hfull{%
\rlx{\hss{$700_{xx}$}}\cg%
\e{0}%
\e{0}%
\e{0}%
\e{0}%
\e{0}%
\e{1}%
\e{0}%
\e{0}%
\e{0}%
\e{0}%
\e{0}%
\e{0}%
\e{0}%
\e{0}%
\e{0}%
\e{0}%
\e{1}%
\e{0}%
\e{0}%
\e{1}%
\e{0}%
\e{0}%
\e{0}%
\e{0}%
\eol}\vss}\rg%
%
%
\rx{\vss\hfull{%
\rlx{\hss{$972_{x}$}}\cg%
\e{0}%
\e{0}%
\e{0}%
\e{0}%
\e{1}%
\e{0}%
\e{0}%
\e{0}%
\e{0}%
\e{0}%
\e{0}%
\e{0}%
\e{0}%
\e{0}%
\e{0}%
\e{0}%
\e{1}%
\e{0}%
\e{0}%
\e{0}%
\e{0}%
\e{0}%
\e{0}%
\e{0}%
\eol}\vss}\rg%
%
%
\rx{\vss\hfull{%
\rlx{\hss{$4096_{x}$}}\cg%
\e{0}%
\e{0}%
\e{0}%
\e{1}%
\e{1}%
\e{0}%
\e{0}%
\e{0}%
\e{1}%
\e{0}%
\e{0}%
\e{0}%
\e{0}%
\e{0}%
\e{1}%
\e{1}%
\e{0}%
\e{0}%
\e{1}%
\e{0}%
\e{0}%
\e{0}%
\e{0}%
\e{0}%
\eol}\vss}\rg%
%
%
\rx{\vss\hfull{%
\rlx{\hss{$4200_{x}$}}\cg%
\e{0}%
\e{0}%
\e{0}%
\e{1}%
\e{0}%
\e{0}%
\e{0}%
\e{0}%
\e{1}%
\e{0}%
\e{0}%
\e{0}%
\e{0}%
\e{0}%
\e{1}%
\e{0}%
\e{1}%
\e{1}%
\e{1}%
\e{0}%
\e{0}%
\e{0}%
\e{0}%
\e{0}%
\eol}\vss}\rg%
%
%
\rx{\vss\hfull{%
\rlx{\hss{$2240_{x}$}}\cg%
\e{0}%
\e{0}%
\e{0}%
\e{1}%
\e{0}%
\e{0}%
\e{0}%
\e{0}%
\e{1}%
\e{0}%
\e{0}%
\e{0}%
\e{0}%
\e{0}%
\e{0}%
\e{0}%
\e{0}%
\e{0}%
\e{1}%
\e{0}%
\e{0}%
\e{0}%
\e{0}%
\e{1}%
\eol}\vss}\rg%
%
%
\rx{\vss\hfull{%
\rlx{\hss{$2835_{x}$}}\cg%
\e{0}%
\e{0}%
\e{0}%
\e{0}%
\e{0}%
\e{0}%
\e{0}%
\e{0}%
\e{1}%
\e{0}%
\e{0}%
\e{0}%
\e{0}%
\e{0}%
\e{0}%
\e{0}%
\e{1}%
\e{1}%
\e{0}%
\e{0}%
\e{0}%
\e{0}%
\e{1}%
\e{1}%
\eol}\vss}\rg%
%
%
\rx{\vss\hfull{%
\rlx{\hss{$6075_{x}$}}\cg%
\e{0}%
\e{0}%
\e{0}%
\e{1}%
\e{1}%
\e{0}%
\e{0}%
\e{0}%
\e{1}%
\e{0}%
\e{0}%
\e{0}%
\e{0}%
\e{0}%
\e{1}%
\e{1}%
\e{1}%
\e{1}%
\e{1}%
\e{0}%
\e{1}%
\e{0}%
\e{0}%
\e{0}%
\eol}\vss}\rg%
%
%
\rx{\vss\hfull{%
\rlx{\hss{$3200_{x}$}}\cg%
\e{0}%
\e{0}%
\e{0}%
\e{0}%
\e{1}%
\e{0}%
\e{0}%
\e{1}%
\e{0}%
\e{0}%
\e{0}%
\e{0}%
\e{0}%
\e{0}%
\e{0}%
\e{1}%
\e{0}%
\e{1}%
\e{1}%
\e{0}%
\e{0}%
\e{0}%
\e{0}%
\e{0}%
\eol}\vss}\rg%
%
%
\rx{\vss\hfull{%
\rlx{\hss{$70_{y}$}}\cg%
\e{0}%
\e{0}%
\e{0}%
\e{0}%
\e{0}%
\e{0}%
\e{0}%
\e{0}%
\e{0}%
\e{0}%
\e{0}%
\e{0}%
\e{0}%
\e{0}%
\e{0}%
\e{0}%
\e{0}%
\e{0}%
\e{0}%
\e{0}%
\e{0}%
\e{0}%
\e{0}%
\e{0}%
\eol}\vss}\rg%
%
%
\rx{\vss\hfull{%
\rlx{\hss{$1134_{y}$}}\cg%
\e{1}%
\e{0}%
\e{0}%
\e{0}%
\e{0}%
\e{0}%
\e{0}%
\e{0}%
\e{0}%
\e{0}%
\e{0}%
\e{0}%
\e{0}%
\e{0}%
\e{0}%
\e{0}%
\e{0}%
\e{0}%
\e{1}%
\e{0}%
\e{0}%
\e{0}%
\e{0}%
\e{0}%
\eol}\vss}\rg%
%
%
\rx{\vss\hfull{%
\rlx{\hss{$1680_{y}$}}\cg%
\e{1}%
\e{0}%
\e{0}%
\e{0}%
\e{0}%
\e{0}%
\e{0}%
\e{0}%
\e{0}%
\e{0}%
\e{0}%
\e{0}%
\e{0}%
\e{0}%
\e{1}%
\e{1}%
\e{0}%
\e{0}%
\e{0}%
\e{0}%
\e{0}%
\e{0}%
\e{0}%
\e{0}%
\eol}\vss}\rg%
%
%
\rx{\vss\hfull{%
\rlx{\hss{$168_{y}$}}\cg%
\e{0}%
\e{0}%
\e{0}%
\e{0}%
\e{0}%
\e{0}%
\e{0}%
\e{0}%
\e{0}%
\e{0}%
\e{0}%
\e{0}%
\e{0}%
\e{0}%
\e{0}%
\e{0}%
\e{0}%
\e{0}%
\e{0}%
\e{1}%
\e{0}%
\e{0}%
\e{0}%
\e{0}%
\eol}\vss}\rg%
%
%
\rx{\vss\hfull{%
\rlx{\hss{$420_{y}$}}\cg%
\e{0}%
\e{0}%
\e{0}%
\e{0}%
\e{0}%
\e{0}%
\e{0}%
\e{0}%
\e{0}%
\e{0}%
\e{0}%
\e{0}%
\e{0}%
\e{0}%
\e{0}%
\e{0}%
\e{0}%
\e{0}%
\e{0}%
\e{0}%
\e{0}%
\e{0}%
\e{1}%
\e{0}%
\eol}\vss}\rg%
%
%
\rx{\vss\hfull{%
\rlx{\hss{$3150_{y}$}}\cg%
\e{0}%
\e{0}%
\e{0}%
\e{0}%
\e{0}%
\e{0}%
\e{0}%
\e{0}%
\e{0}%
\e{1}%
\e{0}%
\e{0}%
\e{0}%
\e{0}%
\e{0}%
\e{0}%
\e{0}%
\e{1}%
\e{1}%
\e{0}%
\e{0}%
\e{0}%
\e{1}%
\e{1}%
\eol}\vss}\rg%
%
%
\rx{\vss\hfull{%
\rlx{\hss{$4200_{y}$}}\cg%
\e{0}%
\e{0}%
\e{0}%
\e{0}%
\e{0}%
\e{0}%
\e{0}%
\e{0}%
\e{0}%
\e{0}%
\e{0}%
\e{0}%
\e{0}%
\e{0}%
\e{0}%
\e{0}%
\e{1}%
\e{1}%
\e{1}%
\e{1}%
\e{1}%
\e{1}%
\e{1}%
\e{0}%
\eol}\vss}\rg%
\eop
\eject
\tablecont%
%
%
%
%
%
%
\rowpts=18 true pt%
\colpts=18 true pt%
\rowlabpts=40 true pt%
\collabpts=70 true pt%
\clx{\vss\hfull{%
\rlx{\hss{$ $}}\cg%
\cx{\hskip 16 true pt\flip{$189_a^{*}{\times}[{1^{2}}]$}\hss}\cg%
\cx{\hskip 16 true pt\flip{$21_b^{*}{\times}[{1^{2}}]$}\hss}\cg%
\cx{\hskip 16 true pt\flip{$35_b^{*}{\times}[{1^{2}}]$}\hss}\cg%
\cx{\hskip 16 true pt\flip{$189_b^{*}{\times}[{1^{2}}]$}\hss}\cg%
\cx{\hskip 16 true pt\flip{$189_c^{*}{\times}[{1^{2}}]$}\hss}\cg%
\cx{\hskip 16 true pt\flip{$15_a^{*}{\times}[{1^{2}}]$}\hss}\cg%
\cx{\hskip 16 true pt\flip{$105_b^{*}{\times}[{1^{2}}]$}\hss}\cg%
\cx{\hskip 16 true pt\flip{$105_c^{*}{\times}[{1^{2}}]$}\hss}\cg%
\cx{\hskip 16 true pt\flip{$315_a^{*}{\times}[{1^{2}}]$}\hss}\cg%
\cx{\hskip 16 true pt\flip{$405_a^{*}{\times}[{1^{2}}]$}\hss}\cg%
\cx{\hskip 16 true pt\flip{$168_a^{*}{\times}[{1^{2}}]$}\hss}\cg%
\cx{\hskip 16 true pt\flip{$56_a^{*}{\times}[{1^{2}}]$}\hss}\cg%
\cx{\hskip 16 true pt\flip{$120_a^{*}{\times}[{1^{2}}]$}\hss}\cg%
\cx{\hskip 16 true pt\flip{$210_a^{*}{\times}[{1^{2}}]$}\hss}\cg%
\cx{\hskip 16 true pt\flip{$280_a^{*}{\times}[{1^{2}}]$}\hss}\cg%
\cx{\hskip 16 true pt\flip{$336_a^{*}{\times}[{1^{2}}]$}\hss}\cg%
\cx{\hskip 16 true pt\flip{$216_a^{*}{\times}[{1^{2}}]$}\hss}\cg%
\cx{\hskip 16 true pt\flip{$512_a^{*}{\times}[{1^{2}}]$}\hss}\cg%
\cx{\hskip 16 true pt\flip{$378_a^{*}{\times}[{1^{2}}]$}\hss}\cg%
\cx{\hskip 16 true pt\flip{$84_a^{*}{\times}[{1^{2}}]$}\hss}\cg%
\cx{\hskip 16 true pt\flip{$420_a^{*}{\times}[{1^{2}}]$}\hss}\cg%
\cx{\hskip 16 true pt\flip{$280_b^{*}{\times}[{1^{2}}]$}\hss}\cg%
\cx{\hskip 16 true pt\flip{$210_b^{*}{\times}[{1^{2}}]$}\hss}\cg%
\cx{\hskip 16 true pt\flip{$70_a^{*}{\times}[{1^{2}}]$}\hss}\cg%
\eol}}\rg%
%
%
\rx{\vss\hfull{%
\rlx{\hss{$2688_{y}$}}\cg%
\e{0}%
\e{0}%
\e{0}%
\e{0}%
\e{0}%
\e{0}%
\e{0}%
\e{0}%
\e{0}%
\e{0}%
\e{0}%
\e{0}%
\e{0}%
\e{0}%
\e{0}%
\e{1}%
\e{1}%
\e{1}%
\e{0}%
\e{0}%
\e{0}%
\e{1}%
\e{0}%
\e{0}%
\eol}\vss}\rg%
%
%
\rx{\vss\hfull{%
\rlx{\hss{$2100_{y}$}}\cg%
\e{0}%
\e{0}%
\e{0}%
\e{0}%
\e{1}%
\e{0}%
\e{0}%
\e{1}%
\e{0}%
\e{0}%
\e{0}%
\e{0}%
\e{0}%
\e{0}%
\e{0}%
\e{1}%
\e{0}%
\e{0}%
\e{0}%
\e{0}%
\e{1}%
\e{0}%
\e{0}%
\e{0}%
\eol}\vss}\rg%
%
%
\rx{\vss\hfull{%
\rlx{\hss{$1400_{y}$}}\cg%
\e{0}%
\e{0}%
\e{0}%
\e{0}%
\e{0}%
\e{0}%
\e{0}%
\e{0}%
\e{0}%
\e{0}%
\e{0}%
\e{0}%
\e{0}%
\e{0}%
\e{1}%
\e{0}%
\e{0}%
\e{0}%
\e{0}%
\e{0}%
\e{1}%
\e{0}%
\e{0}%
\e{0}%
\eol}\vss}\rg%
%
%
\rx{\vss\hfull{%
\rlx{\hss{$4536_{y}$}}\cg%
\e{0}%
\e{0}%
\e{0}%
\e{0}%
\e{0}%
\e{0}%
\e{0}%
\e{0}%
\e{1}%
\e{1}%
\e{0}%
\e{0}%
\e{0}%
\e{0}%
\e{1}%
\e{1}%
\e{0}%
\e{1}%
\e{0}%
\e{0}%
\e{1}%
\e{0}%
\e{0}%
\e{0}%
\eol}\vss}\rg%
%
%
\rx{\vss\hfull{%
\rlx{\hss{$5670_{y}$}}\cg%
\e{1}%
\e{0}%
\e{0}%
\e{0}%
\e{0}%
\e{0}%
\e{0}%
\e{0}%
\e{1}%
\e{1}%
\e{0}%
\e{0}%
\e{0}%
\e{0}%
\e{1}%
\e{1}%
\e{0}%
\e{1}%
\e{1}%
\e{0}%
\e{1}%
\e{0}%
\e{0}%
\e{0}%
\eol}\vss}\rg%
%
%
\rx{\vss\hfull{%
\rlx{\hss{$4480_{y}$}}\cg%
\e{0}%
\e{0}%
\e{0}%
\e{0}%
\e{0}%
\e{0}%
\e{0}%
\e{0}%
\e{1}%
\e{1}%
\e{0}%
\e{0}%
\e{0}%
\e{0}%
\e{0}%
\e{0}%
\e{0}%
\e{1}%
\e{1}%
\e{0}%
\e{1}%
\e{0}%
\e{1}%
\e{0}%
\eol}\vss}\rg%
%
%
\rx{\vss\hfull{%
\rlx{\hss{$8_{z}$}}\cg%
\e{0}%
\e{0}%
\e{0}%
\e{0}%
\e{0}%
\e{0}%
\e{0}%
\e{0}%
\e{0}%
\e{0}%
\e{0}%
\e{0}%
\e{0}%
\e{0}%
\e{0}%
\e{0}%
\e{0}%
\e{0}%
\e{0}%
\e{0}%
\e{0}%
\e{0}%
\e{0}%
\e{0}%
\eol}\vss}\rg%
%
%
\rx{\vss\hfull{%
\rlx{\hss{$56_{z}$}}\cg%
\e{0}%
\e{0}%
\e{0}%
\e{0}%
\e{0}%
\e{0}%
\e{0}%
\e{0}%
\e{0}%
\e{0}%
\e{0}%
\e{0}%
\e{0}%
\e{0}%
\e{0}%
\e{0}%
\e{0}%
\e{0}%
\e{0}%
\e{0}%
\e{0}%
\e{0}%
\e{0}%
\e{0}%
\eol}\vss}\rg%
%
%
\rx{\vss\hfull{%
\rlx{\hss{$160_{z}$}}\cg%
\e{0}%
\e{0}%
\e{0}%
\e{0}%
\e{0}%
\e{0}%
\e{0}%
\e{0}%
\e{0}%
\e{0}%
\e{0}%
\e{0}%
\e{0}%
\e{0}%
\e{0}%
\e{0}%
\e{0}%
\e{0}%
\e{0}%
\e{0}%
\e{0}%
\e{0}%
\e{0}%
\e{0}%
\eol}\vss}\rg%
%
%
\rx{\vss\hfull{%
\rlx{\hss{$112_{z}$}}\cg%
\e{0}%
\e{0}%
\e{0}%
\e{0}%
\e{0}%
\e{0}%
\e{0}%
\e{0}%
\e{0}%
\e{0}%
\e{0}%
\e{0}%
\e{0}%
\e{0}%
\e{0}%
\e{0}%
\e{0}%
\e{0}%
\e{0}%
\e{0}%
\e{0}%
\e{0}%
\e{0}%
\e{0}%
\eol}\vss}\rg%
%
%
\rx{\vss\hfull{%
\rlx{\hss{$840_{z}$}}\cg%
\e{0}%
\e{0}%
\e{0}%
\e{0}%
\e{0}%
\e{0}%
\e{0}%
\e{0}%
\e{0}%
\e{0}%
\e{0}%
\e{0}%
\e{0}%
\e{0}%
\e{0}%
\e{0}%
\e{0}%
\e{0}%
\e{0}%
\e{0}%
\e{0}%
\e{0}%
\e{0}%
\e{0}%
\eol}\vss}\rg%
%
%
\rx{\vss\hfull{%
\rlx{\hss{$1296_{z}$}}\cg%
\e{0}%
\e{0}%
\e{0}%
\e{0}%
\e{0}%
\e{0}%
\e{0}%
\e{0}%
\e{0}%
\e{0}%
\e{0}%
\e{0}%
\e{0}%
\e{0}%
\e{0}%
\e{0}%
\e{0}%
\e{0}%
\e{0}%
\e{0}%
\e{0}%
\e{0}%
\e{0}%
\e{0}%
\eol}\vss}\rg%
%
%
\rx{\vss\hfull{%
\rlx{\hss{$1400_{z}$}}\cg%
\e{0}%
\e{0}%
\e{0}%
\e{0}%
\e{0}%
\e{0}%
\e{0}%
\e{0}%
\e{0}%
\e{0}%
\e{0}%
\e{0}%
\e{0}%
\e{0}%
\e{0}%
\e{0}%
\e{0}%
\e{0}%
\e{0}%
\e{0}%
\e{0}%
\e{0}%
\e{0}%
\e{0}%
\eol}\vss}\rg%
%
%
\rx{\vss\hfull{%
\rlx{\hss{$1008_{z}$}}\cg%
\e{0}%
\e{0}%
\e{0}%
\e{0}%
\e{0}%
\e{0}%
\e{0}%
\e{0}%
\e{0}%
\e{0}%
\e{0}%
\e{0}%
\e{0}%
\e{0}%
\e{0}%
\e{0}%
\e{0}%
\e{0}%
\e{0}%
\e{0}%
\e{0}%
\e{0}%
\e{0}%
\e{0}%
\eol}\vss}\rg%
%
%
\rx{\vss\hfull{%
\rlx{\hss{$560_{z}$}}\cg%
\e{0}%
\e{0}%
\e{0}%
\e{0}%
\e{0}%
\e{0}%
\e{0}%
\e{0}%
\e{0}%
\e{0}%
\e{0}%
\e{0}%
\e{0}%
\e{0}%
\e{0}%
\e{0}%
\e{0}%
\e{0}%
\e{0}%
\e{0}%
\e{0}%
\e{0}%
\e{0}%
\e{0}%
\eol}\vss}\rg%
%
%
\rx{\vss\hfull{%
\rlx{\hss{$1400_{zz}$}}\cg%
\e{0}%
\e{0}%
\e{0}%
\e{0}%
\e{0}%
\e{0}%
\e{0}%
\e{0}%
\e{0}%
\e{0}%
\e{0}%
\e{0}%
\e{0}%
\e{0}%
\e{0}%
\e{0}%
\e{0}%
\e{0}%
\e{0}%
\e{0}%
\e{0}%
\e{0}%
\e{0}%
\e{0}%
\eol}\vss}\rg%
%
%
\rx{\vss\hfull{%
\rlx{\hss{$4200_{z}$}}\cg%
\e{0}%
\e{0}%
\e{0}%
\e{0}%
\e{0}%
\e{0}%
\e{0}%
\e{0}%
\e{0}%
\e{0}%
\e{0}%
\e{0}%
\e{0}%
\e{0}%
\e{0}%
\e{0}%
\e{0}%
\e{0}%
\e{0}%
\e{0}%
\e{0}%
\e{0}%
\e{0}%
\e{0}%
\eol}\vss}\rg%
%
%
\rx{\vss\hfull{%
\rlx{\hss{$400_{z}$}}\cg%
\e{0}%
\e{0}%
\e{0}%
\e{0}%
\e{0}%
\e{0}%
\e{0}%
\e{0}%
\e{0}%
\e{0}%
\e{0}%
\e{0}%
\e{0}%
\e{0}%
\e{0}%
\e{0}%
\e{0}%
\e{0}%
\e{0}%
\e{0}%
\e{0}%
\e{0}%
\e{0}%
\e{0}%
\eol}\vss}\rg%
%
%
\rx{\vss\hfull{%
\rlx{\hss{$3240_{z}$}}\cg%
\e{0}%
\e{0}%
\e{0}%
\e{0}%
\e{0}%
\e{0}%
\e{0}%
\e{0}%
\e{0}%
\e{0}%
\e{0}%
\e{0}%
\e{0}%
\e{0}%
\e{0}%
\e{0}%
\e{0}%
\e{0}%
\e{0}%
\e{0}%
\e{0}%
\e{0}%
\e{0}%
\e{0}%
\eol}\vss}\rg%
%
%
\rx{\vss\hfull{%
\rlx{\hss{$4536_{z}$}}\cg%
\e{0}%
\e{0}%
\e{0}%
\e{0}%
\e{0}%
\e{0}%
\e{0}%
\e{0}%
\e{0}%
\e{0}%
\e{0}%
\e{0}%
\e{0}%
\e{0}%
\e{0}%
\e{0}%
\e{0}%
\e{0}%
\e{0}%
\e{0}%
\e{0}%
\e{0}%
\e{0}%
\e{0}%
\eol}\vss}\rg%
%
%
\rx{\vss\hfull{%
\rlx{\hss{$2400_{z}$}}\cg%
\e{0}%
\e{0}%
\e{0}%
\e{0}%
\e{0}%
\e{0}%
\e{0}%
\e{0}%
\e{0}%
\e{0}%
\e{0}%
\e{0}%
\e{0}%
\e{0}%
\e{0}%
\e{0}%
\e{0}%
\e{0}%
\e{0}%
\e{0}%
\e{0}%
\e{0}%
\e{0}%
\e{0}%
\eol}\vss}\rg%
%
%
\rx{\vss\hfull{%
\rlx{\hss{$3360_{z}$}}\cg%
\e{0}%
\e{0}%
\e{0}%
\e{0}%
\e{0}%
\e{0}%
\e{0}%
\e{0}%
\e{0}%
\e{0}%
\e{0}%
\e{0}%
\e{0}%
\e{0}%
\e{0}%
\e{0}%
\e{0}%
\e{0}%
\e{0}%
\e{0}%
\e{0}%
\e{0}%
\e{0}%
\e{0}%
\eol}\vss}\rg%
%
%
\rx{\vss\hfull{%
\rlx{\hss{$2800_{z}$}}\cg%
\e{0}%
\e{0}%
\e{0}%
\e{0}%
\e{0}%
\e{0}%
\e{0}%
\e{0}%
\e{0}%
\e{0}%
\e{0}%
\e{0}%
\e{0}%
\e{0}%
\e{0}%
\e{0}%
\e{0}%
\e{0}%
\e{0}%
\e{0}%
\e{0}%
\e{0}%
\e{0}%
\e{0}%
\eol}\vss}\rg%
%
%
\rx{\vss\hfull{%
\rlx{\hss{$4096_{z}$}}\cg%
\e{0}%
\e{0}%
\e{0}%
\e{0}%
\e{0}%
\e{0}%
\e{0}%
\e{0}%
\e{0}%
\e{0}%
\e{0}%
\e{0}%
\e{0}%
\e{0}%
\e{0}%
\e{0}%
\e{0}%
\e{0}%
\e{0}%
\e{0}%
\e{0}%
\e{0}%
\e{0}%
\e{0}%
\eol}\vss}\rg%
%
%
\rx{\vss\hfull{%
\rlx{\hss{$5600_{z}$}}\cg%
\e{0}%
\e{0}%
\e{0}%
\e{0}%
\e{0}%
\e{0}%
\e{0}%
\e{0}%
\e{0}%
\e{0}%
\e{0}%
\e{0}%
\e{0}%
\e{0}%
\e{0}%
\e{0}%
\e{0}%
\e{0}%
\e{0}%
\e{0}%
\e{0}%
\e{0}%
\e{0}%
\e{0}%
\eol}\vss}\rg%
%
%
\rx{\vss\hfull{%
\rlx{\hss{$448_{z}$}}\cg%
\e{0}%
\e{0}%
\e{0}%
\e{0}%
\e{0}%
\e{0}%
\e{0}%
\e{0}%
\e{0}%
\e{0}%
\e{0}%
\e{0}%
\e{0}%
\e{0}%
\e{0}%
\e{0}%
\e{0}%
\e{0}%
\e{0}%
\e{0}%
\e{0}%
\e{0}%
\e{0}%
\e{0}%
\eol}\vss}\rg%
%
%
\rx{\vss\hfull{%
\rlx{\hss{$448_{w}$}}\cg%
\e{0}%
\e{0}%
\e{0}%
\e{0}%
\e{0}%
\e{0}%
\e{0}%
\e{0}%
\e{0}%
\e{0}%
\e{0}%
\e{0}%
\e{0}%
\e{0}%
\e{0}%
\e{0}%
\e{0}%
\e{0}%
\e{0}%
\e{0}%
\e{0}%
\e{0}%
\e{0}%
\e{0}%
\eol}\vss}\rg%
%
%
\rx{\vss\hfull{%
\rlx{\hss{$1344_{w}$}}\cg%
\e{0}%
\e{0}%
\e{0}%
\e{0}%
\e{0}%
\e{0}%
\e{0}%
\e{0}%
\e{0}%
\e{0}%
\e{0}%
\e{0}%
\e{0}%
\e{0}%
\e{0}%
\e{0}%
\e{0}%
\e{0}%
\e{0}%
\e{0}%
\e{0}%
\e{0}%
\e{0}%
\e{0}%
\eol}\vss}\rg%
%
%
\rx{\vss\hfull{%
\rlx{\hss{$5600_{w}$}}\cg%
\e{0}%
\e{0}%
\e{0}%
\e{0}%
\e{0}%
\e{0}%
\e{0}%
\e{0}%
\e{0}%
\e{0}%
\e{0}%
\e{0}%
\e{0}%
\e{0}%
\e{0}%
\e{0}%
\e{0}%
\e{0}%
\e{0}%
\e{0}%
\e{0}%
\e{0}%
\e{0}%
\e{0}%
\eol}\vss}\rg%
%
%
\rx{\vss\hfull{%
\rlx{\hss{$2016_{w}$}}\cg%
\e{0}%
\e{0}%
\e{0}%
\e{0}%
\e{0}%
\e{0}%
\e{0}%
\e{0}%
\e{0}%
\e{0}%
\e{0}%
\e{0}%
\e{0}%
\e{0}%
\e{0}%
\e{0}%
\e{0}%
\e{0}%
\e{0}%
\e{0}%
\e{0}%
\e{0}%
\e{0}%
\e{0}%
\eol}\vss}\rg%
%
%
\rx{\vss\hfull{%
\rlx{\hss{$7168_{w}$}}\cg%
\e{0}%
\e{0}%
\e{0}%
\e{0}%
\e{0}%
\e{0}%
\e{0}%
\e{0}%
\e{0}%
\e{0}%
\e{0}%
\e{0}%
\e{0}%
\e{0}%
\e{0}%
\e{0}%
\e{0}%
\e{0}%
\e{0}%
\e{0}%
\e{0}%
\e{0}%
\e{0}%
\e{0}%
\eol}\vss}\rg%
\tableclose%
%
%
%
%
%
%
\eop
\eject
\tableopen{Induce/restrict matrix for $W(D_{7})\,\subset\,W(E_{8})$}%
%
%
%
%
%
%
\rowpts=18 true pt%
\colpts=18 true pt%
\rowlabpts=40 true pt%
\collabpts=65 true pt%
\clx{\vss\hfull{%
\rlx{\hss{$ $}}\cg%
\cx{\hskip 16 true pt\flip{$[{7}:-]$}\hss}\cg%
\cx{\hskip 16 true pt\flip{$[{6}{1}:-]$}\hss}\cg%
\cx{\hskip 16 true pt\flip{$[{5}{2}:-]$}\hss}\cg%
\cx{\hskip 16 true pt\flip{$[{5}{1^{2}}:-]$}\hss}\cg%
\cx{\hskip 16 true pt\flip{$[{4}{3}:-]$}\hss}\cg%
\cx{\hskip 16 true pt\flip{$[{4}{2}{1}:-]$}\hss}\cg%
\cx{\hskip 16 true pt\flip{$[{4}{1^{3}}:-]$}\hss}\cg%
\cx{\hskip 16 true pt\flip{$[{3^{2}}{1}:-]$}\hss}\cg%
\cx{\hskip 16 true pt\flip{$[{3}{2^{2}}:-]$}\hss}\cg%
\cx{\hskip 16 true pt\flip{$[{3}{2}{1^{2}}:-]$}\hss}\cg%
\cx{\hskip 16 true pt\flip{$[{3}{1^{4}}:-]$}\hss}\cg%
\cx{\hskip 16 true pt\flip{$[{2^{3}}{1}:-]$}\hss}\cg%
\cx{\hskip 16 true pt\flip{$[{2^{2}}{1^{3}}:-]$}\hss}\cg%
\cx{\hskip 16 true pt\flip{$[{2}{1^{5}}:-]$}\hss}\cg%
\cx{\hskip 16 true pt\flip{$[{1^{7}}:-]$}\hss}\cg%
\cx{\hskip 16 true pt\flip{$[{6}:{1}]$}\hss}\cg%
\cx{\hskip 16 true pt\flip{$[{5}{1}:{1}]$}\hss}\cg%
\cx{\hskip 16 true pt\flip{$[{4}{2}:{1}]$}\hss}\cg%
\cx{\hskip 16 true pt\flip{$[{4}{1^{2}}:{1}]$}\hss}\cg%
\eol}}\rg%
%
%
\rx{\vss\hfull{%
\rlx{\hss{$1_{x}$}}\cg%
\e{1}%
\e{0}%
\e{0}%
\e{0}%
\e{0}%
\e{0}%
\e{0}%
\e{0}%
\e{0}%
\e{0}%
\e{0}%
\e{0}%
\e{0}%
\e{0}%
\e{0}%
\e{0}%
\e{0}%
\e{0}%
\e{0}%
\eol}\vss}\rg%
%
%
\rx{\vss\hfull{%
\rlx{\hss{$28_{x}$}}\cg%
\e{0}%
\e{0}%
\e{0}%
\e{0}%
\e{0}%
\e{0}%
\e{0}%
\e{0}%
\e{0}%
\e{0}%
\e{0}%
\e{0}%
\e{0}%
\e{0}%
\e{0}%
\e{1}%
\e{0}%
\e{0}%
\e{0}%
\eol}\vss}\rg%
%
%
\rx{\vss\hfull{%
\rlx{\hss{$35_{x}$}}\cg%
\e{1}%
\e{1}%
\e{0}%
\e{0}%
\e{0}%
\e{0}%
\e{0}%
\e{0}%
\e{0}%
\e{0}%
\e{0}%
\e{0}%
\e{0}%
\e{0}%
\e{0}%
\e{1}%
\e{0}%
\e{0}%
\e{0}%
\eol}\vss}\rg%
%
%
\rx{\vss\hfull{%
\rlx{\hss{$84_{x}$}}\cg%
\e{1}%
\e{1}%
\e{1}%
\e{0}%
\e{0}%
\e{0}%
\e{0}%
\e{0}%
\e{0}%
\e{0}%
\e{0}%
\e{0}%
\e{0}%
\e{0}%
\e{0}%
\e{1}%
\e{0}%
\e{0}%
\e{0}%
\eol}\vss}\rg%
%
%
\rx{\vss\hfull{%
\rlx{\hss{$50_{x}$}}\cg%
\e{1}%
\e{0}%
\e{0}%
\e{0}%
\e{1}%
\e{0}%
\e{0}%
\e{0}%
\e{0}%
\e{0}%
\e{0}%
\e{0}%
\e{0}%
\e{0}%
\e{0}%
\e{0}%
\e{0}%
\e{0}%
\e{0}%
\eol}\vss}\rg%
%
%
\rx{\vss\hfull{%
\rlx{\hss{$350_{x}$}}\cg%
\e{0}%
\e{0}%
\e{0}%
\e{0}%
\e{0}%
\e{0}%
\e{0}%
\e{0}%
\e{0}%
\e{0}%
\e{0}%
\e{0}%
\e{0}%
\e{0}%
\e{0}%
\e{0}%
\e{1}%
\e{0}%
\e{0}%
\eol}\vss}\rg%
%
%
\rx{\vss\hfull{%
\rlx{\hss{$300_{x}$}}\cg%
\e{0}%
\e{1}%
\e{1}%
\e{0}%
\e{0}%
\e{0}%
\e{0}%
\e{0}%
\e{0}%
\e{0}%
\e{0}%
\e{0}%
\e{0}%
\e{0}%
\e{0}%
\e{0}%
\e{1}%
\e{0}%
\e{0}%
\eol}\vss}\rg%
%
%
\rx{\vss\hfull{%
\rlx{\hss{$567_{x}$}}\cg%
\e{0}%
\e{1}%
\e{0}%
\e{1}%
\e{0}%
\e{0}%
\e{0}%
\e{0}%
\e{0}%
\e{0}%
\e{0}%
\e{0}%
\e{0}%
\e{0}%
\e{0}%
\e{2}%
\e{2}%
\e{0}%
\e{0}%
\eol}\vss}\rg%
%
%
\rx{\vss\hfull{%
\rlx{\hss{$210_{x}$}}\cg%
\e{1}%
\e{1}%
\e{0}%
\e{0}%
\e{0}%
\e{0}%
\e{0}%
\e{0}%
\e{0}%
\e{0}%
\e{0}%
\e{0}%
\e{0}%
\e{0}%
\e{0}%
\e{1}%
\e{1}%
\e{0}%
\e{0}%
\eol}\vss}\rg%
%
%
\rx{\vss\hfull{%
\rlx{\hss{$840_{x}$}}\cg%
\e{0}%
\e{0}%
\e{1}%
\e{0}%
\e{1}%
\e{0}%
\e{0}%
\e{0}%
\e{0}%
\e{0}%
\e{0}%
\e{0}%
\e{0}%
\e{0}%
\e{0}%
\e{0}%
\e{0}%
\e{1}%
\e{0}%
\eol}\vss}\rg%
%
%
\rx{\vss\hfull{%
\rlx{\hss{$700_{x}$}}\cg%
\e{1}%
\e{1}%
\e{1}%
\e{0}%
\e{1}%
\e{0}%
\e{0}%
\e{0}%
\e{0}%
\e{0}%
\e{0}%
\e{0}%
\e{0}%
\e{0}%
\e{0}%
\e{1}%
\e{1}%
\e{1}%
\e{0}%
\eol}\vss}\rg%
%
%
\rx{\vss\hfull{%
\rlx{\hss{$175_{x}$}}\cg%
\e{0}%
\e{0}%
\e{0}%
\e{0}%
\e{0}%
\e{0}%
\e{0}%
\e{0}%
\e{0}%
\e{0}%
\e{0}%
\e{0}%
\e{0}%
\e{0}%
\e{0}%
\e{0}%
\e{0}%
\e{0}%
\e{0}%
\eol}\vss}\rg%
%
%
\rx{\vss\hfull{%
\rlx{\hss{$1400_{x}$}}\cg%
\e{0}%
\e{0}%
\e{0}%
\e{0}%
\e{0}%
\e{0}%
\e{0}%
\e{0}%
\e{0}%
\e{0}%
\e{0}%
\e{0}%
\e{0}%
\e{0}%
\e{0}%
\e{1}%
\e{1}%
\e{1}%
\e{1}%
\eol}\vss}\rg%
%
%
\rx{\vss\hfull{%
\rlx{\hss{$1050_{x}$}}\cg%
\e{0}%
\e{0}%
\e{0}%
\e{0}%
\e{1}%
\e{1}%
\e{0}%
\e{1}%
\e{0}%
\e{0}%
\e{0}%
\e{0}%
\e{0}%
\e{0}%
\e{0}%
\e{1}%
\e{0}%
\e{1}%
\e{0}%
\eol}\vss}\rg%
%
%
\rx{\vss\hfull{%
\rlx{\hss{$1575_{x}$}}\cg%
\e{0}%
\e{0}%
\e{0}%
\e{0}%
\e{0}%
\e{0}%
\e{0}%
\e{0}%
\e{0}%
\e{0}%
\e{0}%
\e{0}%
\e{0}%
\e{0}%
\e{0}%
\e{1}%
\e{2}%
\e{1}%
\e{1}%
\eol}\vss}\rg%
%
%
\rx{\vss\hfull{%
\rlx{\hss{$1344_{x}$}}\cg%
\e{0}%
\e{1}%
\e{2}%
\e{1}%
\e{0}%
\e{1}%
\e{0}%
\e{0}%
\e{0}%
\e{0}%
\e{0}%
\e{0}%
\e{0}%
\e{0}%
\e{0}%
\e{1}%
\e{2}%
\e{1}%
\e{0}%
\eol}\vss}\rg%
%
%
\rx{\vss\hfull{%
\rlx{\hss{$2100_{x}$}}\cg%
\e{0}%
\e{0}%
\e{0}%
\e{1}%
\e{0}%
\e{0}%
\e{1}%
\e{0}%
\e{0}%
\e{0}%
\e{0}%
\e{0}%
\e{0}%
\e{0}%
\e{0}%
\e{0}%
\e{1}%
\e{0}%
\e{2}%
\eol}\vss}\rg%
%
%
\rx{\vss\hfull{%
\rlx{\hss{$2268_{x}$}}\cg%
\e{0}%
\e{1}%
\e{0}%
\e{1}%
\e{0}%
\e{0}%
\e{0}%
\e{0}%
\e{0}%
\e{0}%
\e{0}%
\e{0}%
\e{0}%
\e{0}%
\e{0}%
\e{0}%
\e{2}%
\e{1}%
\e{2}%
\eol}\vss}\rg%
%
%
\rx{\vss\hfull{%
\rlx{\hss{$525_{x}$}}\cg%
\e{0}%
\e{0}%
\e{0}%
\e{1}%
\e{0}%
\e{0}%
\e{1}%
\e{0}%
\e{0}%
\e{0}%
\e{0}%
\e{0}%
\e{0}%
\e{0}%
\e{0}%
\e{0}%
\e{1}%
\e{0}%
\e{0}%
\eol}\vss}\rg%
%
%
\rx{\vss\hfull{%
\rlx{\hss{$700_{xx}$}}\cg%
\e{0}%
\e{0}%
\e{0}%
\e{0}%
\e{1}%
\e{0}%
\e{0}%
\e{1}%
\e{0}%
\e{1}%
\e{0}%
\e{0}%
\e{0}%
\e{0}%
\e{0}%
\e{0}%
\e{0}%
\e{0}%
\e{0}%
\eol}\vss}\rg%
%
%
\rx{\vss\hfull{%
\rlx{\hss{$972_{x}$}}\cg%
\e{0}%
\e{1}%
\e{1}%
\e{0}%
\e{1}%
\e{1}%
\e{0}%
\e{0}%
\e{1}%
\e{0}%
\e{0}%
\e{0}%
\e{0}%
\e{0}%
\e{0}%
\e{0}%
\e{0}%
\e{1}%
\e{0}%
\eol}\vss}\rg%
%
%
\rx{\vss\hfull{%
\rlx{\hss{$4096_{x}$}}\cg%
\e{0}%
\e{0}%
\e{1}%
\e{1}%
\e{0}%
\e{1}%
\e{0}%
\e{0}%
\e{0}%
\e{0}%
\e{0}%
\e{0}%
\e{0}%
\e{0}%
\e{0}%
\e{0}%
\e{2}%
\e{2}%
\e{2}%
\eol}\vss}\rg%
%
%
\rx{\vss\hfull{%
\rlx{\hss{$4200_{x}$}}\cg%
\e{0}%
\e{0}%
\e{0}%
\e{0}%
\e{1}%
\e{1}%
\e{0}%
\e{1}%
\e{0}%
\e{0}%
\e{0}%
\e{0}%
\e{0}%
\e{0}%
\e{0}%
\e{0}%
\e{1}%
\e{2}%
\e{1}%
\eol}\vss}\rg%
%
%
\rx{\vss\hfull{%
\rlx{\hss{$2240_{x}$}}\cg%
\e{0}%
\e{0}%
\e{1}%
\e{0}%
\e{1}%
\e{0}%
\e{0}%
\e{0}%
\e{0}%
\e{0}%
\e{0}%
\e{0}%
\e{0}%
\e{0}%
\e{0}%
\e{0}%
\e{1}%
\e{2}%
\e{0}%
\eol}\vss}\rg%
%
%
\rx{\vss\hfull{%
\rlx{\hss{$2835_{x}$}}\cg%
\e{0}%
\e{0}%
\e{0}%
\e{0}%
\e{0}%
\e{0}%
\e{0}%
\e{1}%
\e{1}%
\e{0}%
\e{0}%
\e{0}%
\e{0}%
\e{0}%
\e{0}%
\e{0}%
\e{0}%
\e{1}%
\e{0}%
\eol}\vss}\rg%
%
%
\rx{\vss\hfull{%
\rlx{\hss{$6075_{x}$}}\cg%
\e{0}%
\e{0}%
\e{0}%
\e{0}%
\e{0}%
\e{1}%
\e{1}%
\e{0}%
\e{0}%
\e{1}%
\e{0}%
\e{0}%
\e{0}%
\e{0}%
\e{0}%
\e{0}%
\e{1}%
\e{2}%
\e{2}%
\eol}\vss}\rg%
%
%
\rx{\vss\hfull{%
\rlx{\hss{$3200_{x}$}}\cg%
\e{0}%
\e{0}%
\e{1}%
\e{1}%
\e{0}%
\e{2}%
\e{0}%
\e{0}%
\e{1}%
\e{0}%
\e{0}%
\e{0}%
\e{0}%
\e{0}%
\e{0}%
\e{0}%
\e{0}%
\e{1}%
\e{1}%
\eol}\vss}\rg%
%
%
\rx{\vss\hfull{%
\rlx{\hss{$70_{y}$}}\cg%
\e{0}%
\e{0}%
\e{0}%
\e{0}%
\e{0}%
\e{0}%
\e{0}%
\e{0}%
\e{0}%
\e{0}%
\e{0}%
\e{0}%
\e{0}%
\e{0}%
\e{0}%
\e{0}%
\e{0}%
\e{0}%
\e{0}%
\eol}\vss}\rg%
%
%
\rx{\vss\hfull{%
\rlx{\hss{$1134_{y}$}}\cg%
\e{0}%
\e{0}%
\e{0}%
\e{0}%
\e{0}%
\e{0}%
\e{0}%
\e{0}%
\e{0}%
\e{0}%
\e{0}%
\e{0}%
\e{0}%
\e{0}%
\e{0}%
\e{0}%
\e{0}%
\e{1}%
\e{0}%
\eol}\vss}\rg%
%
%
\rx{\vss\hfull{%
\rlx{\hss{$1680_{y}$}}\cg%
\e{0}%
\e{0}%
\e{0}%
\e{0}%
\e{0}%
\e{0}%
\e{0}%
\e{0}%
\e{0}%
\e{0}%
\e{0}%
\e{0}%
\e{0}%
\e{0}%
\e{0}%
\e{0}%
\e{0}%
\e{0}%
\e{1}%
\eol}\vss}\rg%
%
%
\rx{\vss\hfull{%
\rlx{\hss{$168_{y}$}}\cg%
\e{0}%
\e{0}%
\e{0}%
\e{0}%
\e{1}%
\e{0}%
\e{0}%
\e{0}%
\e{0}%
\e{0}%
\e{0}%
\e{1}%
\e{0}%
\e{0}%
\e{0}%
\e{0}%
\e{0}%
\e{0}%
\e{0}%
\eol}\vss}\rg%
%
%
\rx{\vss\hfull{%
\rlx{\hss{$420_{y}$}}\cg%
\e{0}%
\e{0}%
\e{0}%
\e{0}%
\e{0}%
\e{0}%
\e{0}%
\e{0}%
\e{0}%
\e{0}%
\e{0}%
\e{0}%
\e{0}%
\e{0}%
\e{0}%
\e{0}%
\e{0}%
\e{0}%
\e{0}%
\eol}\vss}\rg%
%
%
\rx{\vss\hfull{%
\rlx{\hss{$3150_{y}$}}\cg%
\e{0}%
\e{0}%
\e{0}%
\e{0}%
\e{0}%
\e{0}%
\e{0}%
\e{0}%
\e{0}%
\e{0}%
\e{0}%
\e{0}%
\e{0}%
\e{0}%
\e{0}%
\e{0}%
\e{0}%
\e{1}%
\e{0}%
\eol}\vss}\rg%
%
%
\rx{\vss\hfull{%
\rlx{\hss{$4200_{y}$}}\cg%
\e{0}%
\e{0}%
\e{0}%
\e{0}%
\e{1}%
\e{1}%
\e{0}%
\e{1}%
\e{1}%
\e{1}%
\e{0}%
\e{1}%
\e{0}%
\e{0}%
\e{0}%
\e{0}%
\e{0}%
\e{1}%
\e{0}%
\eol}\vss}\rg%
\eop
\eject
\tablecont%
%
%
%
%
%
%
\rowpts=18 true pt%
\colpts=18 true pt%
\rowlabpts=40 true pt%
\collabpts=65 true pt%
\clx{\vss\hfull{%
\rlx{\hss{$ $}}\cg%
\cx{\hskip 16 true pt\flip{$[{7}:-]$}\hss}\cg%
\cx{\hskip 16 true pt\flip{$[{6}{1}:-]$}\hss}\cg%
\cx{\hskip 16 true pt\flip{$[{5}{2}:-]$}\hss}\cg%
\cx{\hskip 16 true pt\flip{$[{5}{1^{2}}:-]$}\hss}\cg%
\cx{\hskip 16 true pt\flip{$[{4}{3}:-]$}\hss}\cg%
\cx{\hskip 16 true pt\flip{$[{4}{2}{1}:-]$}\hss}\cg%
\cx{\hskip 16 true pt\flip{$[{4}{1^{3}}:-]$}\hss}\cg%
\cx{\hskip 16 true pt\flip{$[{3^{2}}{1}:-]$}\hss}\cg%
\cx{\hskip 16 true pt\flip{$[{3}{2^{2}}:-]$}\hss}\cg%
\cx{\hskip 16 true pt\flip{$[{3}{2}{1^{2}}:-]$}\hss}\cg%
\cx{\hskip 16 true pt\flip{$[{3}{1^{4}}:-]$}\hss}\cg%
\cx{\hskip 16 true pt\flip{$[{2^{3}}{1}:-]$}\hss}\cg%
\cx{\hskip 16 true pt\flip{$[{2^{2}}{1^{3}}:-]$}\hss}\cg%
\cx{\hskip 16 true pt\flip{$[{2}{1^{5}}:-]$}\hss}\cg%
\cx{\hskip 16 true pt\flip{$[{1^{7}}:-]$}\hss}\cg%
\cx{\hskip 16 true pt\flip{$[{6}:{1}]$}\hss}\cg%
\cx{\hskip 16 true pt\flip{$[{5}{1}:{1}]$}\hss}\cg%
\cx{\hskip 16 true pt\flip{$[{4}{2}:{1}]$}\hss}\cg%
\cx{\hskip 16 true pt\flip{$[{4}{1^{2}}:{1}]$}\hss}\cg%
\eol}}\rg%
%
%
\rx{\vss\hfull{%
\rlx{\hss{$2688_{y}$}}\cg%
\e{0}%
\e{0}%
\e{0}%
\e{0}%
\e{0}%
\e{1}%
\e{0}%
\e{1}%
\e{1}%
\e{1}%
\e{0}%
\e{0}%
\e{0}%
\e{0}%
\e{0}%
\e{0}%
\e{0}%
\e{0}%
\e{0}%
\eol}\vss}\rg%
%
%
\rx{\vss\hfull{%
\rlx{\hss{$2100_{y}$}}\cg%
\e{0}%
\e{0}%
\e{0}%
\e{1}%
\e{0}%
\e{0}%
\e{2}%
\e{0}%
\e{0}%
\e{0}%
\e{1}%
\e{0}%
\e{0}%
\e{0}%
\e{0}%
\e{0}%
\e{0}%
\e{0}%
\e{1}%
\eol}\vss}\rg%
%
%
\rx{\vss\hfull{%
\rlx{\hss{$1400_{y}$}}\cg%
\e{0}%
\e{0}%
\e{0}%
\e{0}%
\e{0}%
\e{0}%
\e{0}%
\e{0}%
\e{0}%
\e{0}%
\e{0}%
\e{0}%
\e{0}%
\e{0}%
\e{0}%
\e{0}%
\e{0}%
\e{0}%
\e{1}%
\eol}\vss}\rg%
%
%
\rx{\vss\hfull{%
\rlx{\hss{$4536_{y}$}}\cg%
\e{0}%
\e{0}%
\e{0}%
\e{0}%
\e{0}%
\e{0}%
\e{0}%
\e{0}%
\e{0}%
\e{0}%
\e{0}%
\e{0}%
\e{0}%
\e{0}%
\e{0}%
\e{0}%
\e{0}%
\e{0}%
\e{2}%
\eol}\vss}\rg%
%
%
\rx{\vss\hfull{%
\rlx{\hss{$5670_{y}$}}\cg%
\e{0}%
\e{0}%
\e{0}%
\e{0}%
\e{0}%
\e{0}%
\e{0}%
\e{0}%
\e{0}%
\e{0}%
\e{0}%
\e{0}%
\e{0}%
\e{0}%
\e{0}%
\e{0}%
\e{0}%
\e{1}%
\e{2}%
\eol}\vss}\rg%
%
%
\rx{\vss\hfull{%
\rlx{\hss{$4480_{y}$}}\cg%
\e{0}%
\e{0}%
\e{0}%
\e{0}%
\e{0}%
\e{0}%
\e{0}%
\e{0}%
\e{0}%
\e{0}%
\e{0}%
\e{0}%
\e{0}%
\e{0}%
\e{0}%
\e{0}%
\e{0}%
\e{1}%
\e{1}%
\eol}\vss}\rg%
%
%
\rx{\vss\hfull{%
\rlx{\hss{$8_{z}$}}\cg%
\e{1}%
\e{0}%
\e{0}%
\e{0}%
\e{0}%
\e{0}%
\e{0}%
\e{0}%
\e{0}%
\e{0}%
\e{0}%
\e{0}%
\e{0}%
\e{0}%
\e{0}%
\e{1}%
\e{0}%
\e{0}%
\e{0}%
\eol}\vss}\rg%
%
%
\rx{\vss\hfull{%
\rlx{\hss{$56_{z}$}}\cg%
\e{0}%
\e{0}%
\e{0}%
\e{0}%
\e{0}%
\e{0}%
\e{0}%
\e{0}%
\e{0}%
\e{0}%
\e{0}%
\e{0}%
\e{0}%
\e{0}%
\e{0}%
\e{0}%
\e{0}%
\e{0}%
\e{0}%
\eol}\vss}\rg%
%
%
\rx{\vss\hfull{%
\rlx{\hss{$160_{z}$}}\cg%
\e{0}%
\e{1}%
\e{0}%
\e{0}%
\e{0}%
\e{0}%
\e{0}%
\e{0}%
\e{0}%
\e{0}%
\e{0}%
\e{0}%
\e{0}%
\e{0}%
\e{0}%
\e{1}%
\e{1}%
\e{0}%
\e{0}%
\eol}\vss}\rg%
%
%
\rx{\vss\hfull{%
\rlx{\hss{$112_{z}$}}\cg%
\e{1}%
\e{1}%
\e{0}%
\e{0}%
\e{0}%
\e{0}%
\e{0}%
\e{0}%
\e{0}%
\e{0}%
\e{0}%
\e{0}%
\e{0}%
\e{0}%
\e{0}%
\e{2}%
\e{1}%
\e{0}%
\e{0}%
\eol}\vss}\rg%
%
%
\rx{\vss\hfull{%
\rlx{\hss{$840_{z}$}}\cg%
\e{0}%
\e{0}%
\e{1}%
\e{0}%
\e{0}%
\e{0}%
\e{0}%
\e{0}%
\e{0}%
\e{0}%
\e{0}%
\e{0}%
\e{0}%
\e{0}%
\e{0}%
\e{0}%
\e{1}%
\e{1}%
\e{0}%
\eol}\vss}\rg%
%
%
\rx{\vss\hfull{%
\rlx{\hss{$1296_{z}$}}\cg%
\e{0}%
\e{0}%
\e{0}%
\e{1}%
\e{0}%
\e{0}%
\e{0}%
\e{0}%
\e{0}%
\e{0}%
\e{0}%
\e{0}%
\e{0}%
\e{0}%
\e{0}%
\e{0}%
\e{1}%
\e{0}%
\e{1}%
\eol}\vss}\rg%
%
%
\rx{\vss\hfull{%
\rlx{\hss{$1400_{z}$}}\cg%
\e{0}%
\e{1}%
\e{1}%
\e{1}%
\e{0}%
\e{0}%
\e{0}%
\e{0}%
\e{0}%
\e{0}%
\e{0}%
\e{0}%
\e{0}%
\e{0}%
\e{0}%
\e{1}%
\e{3}%
\e{1}%
\e{1}%
\eol}\vss}\rg%
%
%
\rx{\vss\hfull{%
\rlx{\hss{$1008_{z}$}}\cg%
\e{0}%
\e{1}%
\e{0}%
\e{1}%
\e{0}%
\e{0}%
\e{0}%
\e{0}%
\e{0}%
\e{0}%
\e{0}%
\e{0}%
\e{0}%
\e{0}%
\e{0}%
\e{1}%
\e{2}%
\e{0}%
\e{1}%
\eol}\vss}\rg%
%
%
\rx{\vss\hfull{%
\rlx{\hss{$560_{z}$}}\cg%
\e{1}%
\e{1}%
\e{1}%
\e{0}%
\e{0}%
\e{0}%
\e{0}%
\e{0}%
\e{0}%
\e{0}%
\e{0}%
\e{0}%
\e{0}%
\e{0}%
\e{0}%
\e{2}%
\e{2}%
\e{1}%
\e{0}%
\eol}\vss}\rg%
%
%
\rx{\vss\hfull{%
\rlx{\hss{$1400_{zz}$}}\cg%
\e{0}%
\e{0}%
\e{0}%
\e{0}%
\e{1}%
\e{0}%
\e{0}%
\e{1}%
\e{0}%
\e{0}%
\e{0}%
\e{0}%
\e{0}%
\e{0}%
\e{0}%
\e{0}%
\e{0}%
\e{1}%
\e{0}%
\eol}\vss}\rg%
%
%
\rx{\vss\hfull{%
\rlx{\hss{$4200_{z}$}}\cg%
\e{0}%
\e{0}%
\e{0}%
\e{0}%
\e{1}%
\e{0}%
\e{0}%
\e{1}%
\e{0}%
\e{1}%
\e{0}%
\e{0}%
\e{0}%
\e{0}%
\e{0}%
\e{0}%
\e{0}%
\e{1}%
\e{0}%
\eol}\vss}\rg%
%
%
\rx{\vss\hfull{%
\rlx{\hss{$400_{z}$}}\cg%
\e{1}%
\e{0}%
\e{0}%
\e{0}%
\e{1}%
\e{0}%
\e{0}%
\e{0}%
\e{0}%
\e{0}%
\e{0}%
\e{0}%
\e{0}%
\e{0}%
\e{0}%
\e{1}%
\e{0}%
\e{1}%
\e{0}%
\eol}\vss}\rg%
%
%
\rx{\vss\hfull{%
\rlx{\hss{$3240_{z}$}}\cg%
\e{0}%
\e{1}%
\e{1}%
\e{0}%
\e{1}%
\e{1}%
\e{0}%
\e{0}%
\e{0}%
\e{0}%
\e{0}%
\e{0}%
\e{0}%
\e{0}%
\e{0}%
\e{1}%
\e{2}%
\e{3}%
\e{1}%
\eol}\vss}\rg%
%
%
\rx{\vss\hfull{%
\rlx{\hss{$4536_{z}$}}\cg%
\e{0}%
\e{0}%
\e{1}%
\e{0}%
\e{1}%
\e{1}%
\e{0}%
\e{0}%
\e{1}%
\e{0}%
\e{0}%
\e{0}%
\e{0}%
\e{0}%
\e{0}%
\e{0}%
\e{1}%
\e{3}%
\e{1}%
\eol}\vss}\rg%
%
%
\rx{\vss\hfull{%
\rlx{\hss{$2400_{z}$}}\cg%
\e{0}%
\e{0}%
\e{0}%
\e{0}%
\e{0}%
\e{0}%
\e{1}%
\e{0}%
\e{0}%
\e{0}%
\e{0}%
\e{0}%
\e{0}%
\e{0}%
\e{0}%
\e{0}%
\e{0}%
\e{0}%
\e{1}%
\eol}\vss}\rg%
%
%
\rx{\vss\hfull{%
\rlx{\hss{$3360_{z}$}}\cg%
\e{0}%
\e{0}%
\e{0}%
\e{0}%
\e{0}%
\e{1}%
\e{0}%
\e{1}%
\e{0}%
\e{0}%
\e{0}%
\e{0}%
\e{0}%
\e{0}%
\e{0}%
\e{0}%
\e{0}%
\e{1}%
\e{1}%
\eol}\vss}\rg%
%
%
\rx{\vss\hfull{%
\rlx{\hss{$2800_{z}$}}\cg%
\e{0}%
\e{0}%
\e{0}%
\e{1}%
\e{0}%
\e{0}%
\e{1}%
\e{0}%
\e{0}%
\e{0}%
\e{0}%
\e{0}%
\e{0}%
\e{0}%
\e{0}%
\e{0}%
\e{1}%
\e{0}%
\e{2}%
\eol}\vss}\rg%
%
%
\rx{\vss\hfull{%
\rlx{\hss{$4096_{z}$}}\cg%
\e{0}%
\e{0}%
\e{1}%
\e{1}%
\e{0}%
\e{1}%
\e{0}%
\e{0}%
\e{0}%
\e{0}%
\e{0}%
\e{0}%
\e{0}%
\e{0}%
\e{0}%
\e{0}%
\e{2}%
\e{2}%
\e{2}%
\eol}\vss}\rg%
%
%
\rx{\vss\hfull{%
\rlx{\hss{$5600_{z}$}}\cg%
\e{0}%
\e{0}%
\e{0}%
\e{1}%
\e{0}%
\e{1}%
\e{1}%
\e{0}%
\e{0}%
\e{0}%
\e{0}%
\e{0}%
\e{0}%
\e{0}%
\e{0}%
\e{0}%
\e{1}%
\e{1}%
\e{3}%
\eol}\vss}\rg%
%
%
\rx{\vss\hfull{%
\rlx{\hss{$448_{z}$}}\cg%
\e{0}%
\e{0}%
\e{1}%
\e{0}%
\e{0}%
\e{0}%
\e{0}%
\e{0}%
\e{0}%
\e{0}%
\e{0}%
\e{0}%
\e{0}%
\e{0}%
\e{0}%
\e{0}%
\e{1}%
\e{1}%
\e{0}%
\eol}\vss}\rg%
%
%
\rx{\vss\hfull{%
\rlx{\hss{$448_{w}$}}\cg%
\e{0}%
\e{0}%
\e{0}%
\e{0}%
\e{0}%
\e{0}%
\e{0}%
\e{0}%
\e{0}%
\e{0}%
\e{0}%
\e{0}%
\e{0}%
\e{0}%
\e{0}%
\e{0}%
\e{0}%
\e{0}%
\e{0}%
\eol}\vss}\rg%
%
%
\rx{\vss\hfull{%
\rlx{\hss{$1344_{w}$}}\cg%
\e{0}%
\e{0}%
\e{0}%
\e{0}%
\e{1}%
\e{0}%
\e{0}%
\e{0}%
\e{0}%
\e{0}%
\e{0}%
\e{1}%
\e{0}%
\e{0}%
\e{0}%
\e{0}%
\e{0}%
\e{1}%
\e{0}%
\eol}\vss}\rg%
%
%
\rx{\vss\hfull{%
\rlx{\hss{$5600_{w}$}}\cg%
\e{0}%
\e{0}%
\e{0}%
\e{0}%
\e{0}%
\e{1}%
\e{0}%
\e{0}%
\e{0}%
\e{1}%
\e{0}%
\e{0}%
\e{0}%
\e{0}%
\e{0}%
\e{0}%
\e{0}%
\e{1}%
\e{1}%
\eol}\vss}\rg%
%
%
\rx{\vss\hfull{%
\rlx{\hss{$2016_{w}$}}\cg%
\e{0}%
\e{0}%
\e{0}%
\e{0}%
\e{0}%
\e{0}%
\e{0}%
\e{1}%
\e{1}%
\e{0}%
\e{0}%
\e{0}%
\e{0}%
\e{0}%
\e{0}%
\e{0}%
\e{0}%
\e{0}%
\e{0}%
\eol}\vss}\rg%
%
%
\rx{\vss\hfull{%
\rlx{\hss{$7168_{w}$}}\cg%
\e{0}%
\e{0}%
\e{0}%
\e{0}%
\e{0}%
\e{1}%
\e{0}%
\e{1}%
\e{1}%
\e{1}%
\e{0}%
\e{0}%
\e{0}%
\e{0}%
\e{0}%
\e{0}%
\e{0}%
\e{1}%
\e{1}%
\eol}\vss}\rg%
\eop
\eject
\tablecont%
%
%
%
%
%
%
\rowpts=18 true pt%
\colpts=18 true pt%
\rowlabpts=40 true pt%
\collabpts=65 true pt%
\clx{\vss\hfull{%
\rlx{\hss{$ $}}\cg%
\cx{\hskip 16 true pt\flip{$[{3^{2}}:{1}]$}\hss}\cg%
\cx{\hskip 16 true pt\flip{$[{3}{2}{1}:{1}]$}\hss}\cg%
\cx{\hskip 16 true pt\flip{$[{3}{1^{3}}:{1}]$}\hss}\cg%
\cx{\hskip 16 true pt\flip{$[{2^{3}}:{1}]$}\hss}\cg%
\cx{\hskip 16 true pt\flip{$[{2^{2}}{1^{2}}:{1}]$}\hss}\cg%
\cx{\hskip 16 true pt\flip{$[{2}{1^{4}}:{1}]$}\hss}\cg%
\cx{\hskip 16 true pt\flip{$[{1^{6}}:{1}]$}\hss}\cg%
\cx{\hskip 16 true pt\flip{$[{5}:{2}]$}\hss}\cg%
\cx{\hskip 16 true pt\flip{$[{5}:{1^{2}}]$}\hss}\cg%
\cx{\hskip 16 true pt\flip{$[{4}{1}:{2}]$}\hss}\cg%
\cx{\hskip 16 true pt\flip{$[{4}{1}:{1^{2}}]$}\hss}\cg%
\cx{\hskip 16 true pt\flip{$[{3}{2}:{2}]$}\hss}\cg%
\cx{\hskip 16 true pt\flip{$[{3}{2}:{1^{2}}]$}\hss}\cg%
\cx{\hskip 16 true pt\flip{$[{3}{1^{2}}:{2}]$}\hss}\cg%
\cx{\hskip 16 true pt\flip{$[{3}{1^{2}}:{1^{2}}]$}\hss}\cg%
\cx{\hskip 16 true pt\flip{$[{2^{2}}{1}:{2}]$}\hss}\cg%
\cx{\hskip 16 true pt\flip{$[{2^{2}}{1}:{1^{2}}]$}\hss}\cg%
\cx{\hskip 16 true pt\flip{$[{2}{1^{3}}:{2}]$}\hss}\cg%
\eol}}\rg%
%
%
\rx{\vss\hfull{%
\rlx{\hss{$1_{x}$}}\cg%
\e{0}%
\e{0}%
\e{0}%
\e{0}%
\e{0}%
\e{0}%
\e{0}%
\e{0}%
\e{0}%
\e{0}%
\e{0}%
\e{0}%
\e{0}%
\e{0}%
\e{0}%
\e{0}%
\e{0}%
\e{0}%
\eol}\vss}\rg%
%
%
\rx{\vss\hfull{%
\rlx{\hss{$28_{x}$}}\cg%
\e{0}%
\e{0}%
\e{0}%
\e{0}%
\e{0}%
\e{0}%
\e{0}%
\e{0}%
\e{1}%
\e{0}%
\e{0}%
\e{0}%
\e{0}%
\e{0}%
\e{0}%
\e{0}%
\e{0}%
\e{0}%
\eol}\vss}\rg%
%
%
\rx{\vss\hfull{%
\rlx{\hss{$35_{x}$}}\cg%
\e{0}%
\e{0}%
\e{0}%
\e{0}%
\e{0}%
\e{0}%
\e{0}%
\e{1}%
\e{0}%
\e{0}%
\e{0}%
\e{0}%
\e{0}%
\e{0}%
\e{0}%
\e{0}%
\e{0}%
\e{0}%
\eol}\vss}\rg%
%
%
\rx{\vss\hfull{%
\rlx{\hss{$84_{x}$}}\cg%
\e{0}%
\e{0}%
\e{0}%
\e{0}%
\e{0}%
\e{0}%
\e{0}%
\e{1}%
\e{0}%
\e{0}%
\e{0}%
\e{0}%
\e{0}%
\e{0}%
\e{0}%
\e{0}%
\e{0}%
\e{0}%
\eol}\vss}\rg%
%
%
\rx{\vss\hfull{%
\rlx{\hss{$50_{x}$}}\cg%
\e{0}%
\e{0}%
\e{0}%
\e{0}%
\e{0}%
\e{0}%
\e{0}%
\e{0}%
\e{0}%
\e{0}%
\e{0}%
\e{0}%
\e{0}%
\e{0}%
\e{0}%
\e{0}%
\e{0}%
\e{0}%
\eol}\vss}\rg%
%
%
\rx{\vss\hfull{%
\rlx{\hss{$350_{x}$}}\cg%
\e{0}%
\e{0}%
\e{0}%
\e{0}%
\e{0}%
\e{0}%
\e{0}%
\e{0}%
\e{1}%
\e{0}%
\e{1}%
\e{0}%
\e{0}%
\e{0}%
\e{0}%
\e{0}%
\e{0}%
\e{0}%
\eol}\vss}\rg%
%
%
\rx{\vss\hfull{%
\rlx{\hss{$300_{x}$}}\cg%
\e{0}%
\e{0}%
\e{0}%
\e{0}%
\e{0}%
\e{0}%
\e{0}%
\e{1}%
\e{0}%
\e{1}%
\e{0}%
\e{0}%
\e{0}%
\e{0}%
\e{0}%
\e{0}%
\e{0}%
\e{0}%
\eol}\vss}\rg%
%
%
\rx{\vss\hfull{%
\rlx{\hss{$567_{x}$}}\cg%
\e{0}%
\e{0}%
\e{0}%
\e{0}%
\e{0}%
\e{0}%
\e{0}%
\e{2}%
\e{2}%
\e{1}%
\e{1}%
\e{0}%
\e{0}%
\e{0}%
\e{0}%
\e{0}%
\e{0}%
\e{0}%
\eol}\vss}\rg%
%
%
\rx{\vss\hfull{%
\rlx{\hss{$210_{x}$}}\cg%
\e{0}%
\e{0}%
\e{0}%
\e{0}%
\e{0}%
\e{0}%
\e{0}%
\e{1}%
\e{1}%
\e{1}%
\e{0}%
\e{0}%
\e{0}%
\e{0}%
\e{0}%
\e{0}%
\e{0}%
\e{0}%
\eol}\vss}\rg%
%
%
\rx{\vss\hfull{%
\rlx{\hss{$840_{x}$}}\cg%
\e{0}%
\e{0}%
\e{0}%
\e{1}%
\e{0}%
\e{0}%
\e{0}%
\e{0}%
\e{0}%
\e{1}%
\e{0}%
\e{1}%
\e{0}%
\e{0}%
\e{0}%
\e{1}%
\e{0}%
\e{0}%
\eol}\vss}\rg%
%
%
\rx{\vss\hfull{%
\rlx{\hss{$700_{x}$}}\cg%
\e{0}%
\e{0}%
\e{0}%
\e{0}%
\e{0}%
\e{0}%
\e{0}%
\e{2}%
\e{0}%
\e{2}%
\e{0}%
\e{1}%
\e{0}%
\e{0}%
\e{0}%
\e{0}%
\e{0}%
\e{0}%
\eol}\vss}\rg%
%
%
\rx{\vss\hfull{%
\rlx{\hss{$175_{x}$}}\cg%
\e{1}%
\e{0}%
\e{0}%
\e{0}%
\e{0}%
\e{0}%
\e{0}%
\e{0}%
\e{0}%
\e{0}%
\e{0}%
\e{1}%
\e{0}%
\e{0}%
\e{0}%
\e{0}%
\e{0}%
\e{0}%
\eol}\vss}\rg%
%
%
\rx{\vss\hfull{%
\rlx{\hss{$1400_{x}$}}\cg%
\e{1}%
\e{0}%
\e{0}%
\e{0}%
\e{0}%
\e{0}%
\e{0}%
\e{1}%
\e{1}%
\e{2}%
\e{1}%
\e{1}%
\e{1}%
\e{1}%
\e{0}%
\e{0}%
\e{0}%
\e{0}%
\eol}\vss}\rg%
%
%
\rx{\vss\hfull{%
\rlx{\hss{$1050_{x}$}}\cg%
\e{1}%
\e{0}%
\e{0}%
\e{0}%
\e{0}%
\e{0}%
\e{0}%
\e{1}%
\e{0}%
\e{1}%
\e{0}%
\e{1}%
\e{1}%
\e{0}%
\e{0}%
\e{0}%
\e{0}%
\e{0}%
\eol}\vss}\rg%
%
%
\rx{\vss\hfull{%
\rlx{\hss{$1575_{x}$}}\cg%
\e{0}%
\e{0}%
\e{0}%
\e{0}%
\e{0}%
\e{0}%
\e{0}%
\e{1}%
\e{2}%
\e{2}%
\e{2}%
\e{0}%
\e{1}%
\e{1}%
\e{0}%
\e{0}%
\e{0}%
\e{0}%
\eol}\vss}\rg%
%
%
\rx{\vss\hfull{%
\rlx{\hss{$1344_{x}$}}\cg%
\e{0}%
\e{0}%
\e{0}%
\e{0}%
\e{0}%
\e{0}%
\e{0}%
\e{2}%
\e{1}%
\e{2}%
\e{1}%
\e{1}%
\e{0}%
\e{0}%
\e{0}%
\e{0}%
\e{0}%
\e{0}%
\eol}\vss}\rg%
%
%
\rx{\vss\hfull{%
\rlx{\hss{$2100_{x}$}}\cg%
\e{0}%
\e{0}%
\e{0}%
\e{0}%
\e{0}%
\e{0}%
\e{0}%
\e{0}%
\e{1}%
\e{1}%
\e{2}%
\e{0}%
\e{0}%
\e{1}%
\e{1}%
\e{0}%
\e{0}%
\e{0}%
\eol}\vss}\rg%
%
%
\rx{\vss\hfull{%
\rlx{\hss{$2268_{x}$}}\cg%
\e{0}%
\e{0}%
\e{0}%
\e{0}%
\e{0}%
\e{0}%
\e{0}%
\e{1}%
\e{1}%
\e{3}%
\e{2}%
\e{1}%
\e{0}%
\e{1}%
\e{1}%
\e{0}%
\e{0}%
\e{0}%
\eol}\vss}\rg%
%
%
\rx{\vss\hfull{%
\rlx{\hss{$525_{x}$}}\cg%
\e{0}%
\e{0}%
\e{0}%
\e{0}%
\e{0}%
\e{0}%
\e{0}%
\e{0}%
\e{1}%
\e{0}%
\e{1}%
\e{0}%
\e{0}%
\e{0}%
\e{0}%
\e{0}%
\e{0}%
\e{0}%
\eol}\vss}\rg%
%
%
\rx{\vss\hfull{%
\rlx{\hss{$700_{xx}$}}\cg%
\e{1}%
\e{0}%
\e{0}%
\e{0}%
\e{0}%
\e{0}%
\e{0}%
\e{0}%
\e{0}%
\e{0}%
\e{0}%
\e{0}%
\e{1}%
\e{0}%
\e{0}%
\e{0}%
\e{0}%
\e{0}%
\eol}\vss}\rg%
%
%
\rx{\vss\hfull{%
\rlx{\hss{$972_{x}$}}\cg%
\e{0}%
\e{0}%
\e{0}%
\e{0}%
\e{0}%
\e{0}%
\e{0}%
\e{0}%
\e{0}%
\e{1}%
\e{0}%
\e{1}%
\e{0}%
\e{0}%
\e{0}%
\e{0}%
\e{0}%
\e{0}%
\eol}\vss}\rg%
%
%
\rx{\vss\hfull{%
\rlx{\hss{$4096_{x}$}}\cg%
\e{0}%
\e{1}%
\e{0}%
\e{0}%
\e{0}%
\e{0}%
\e{0}%
\e{1}%
\e{1}%
\e{3}%
\e{3}%
\e{2}%
\e{1}%
\e{2}%
\e{1}%
\e{1}%
\e{0}%
\e{0}%
\eol}\vss}\rg%
%
%
\rx{\vss\hfull{%
\rlx{\hss{$4200_{x}$}}\cg%
\e{1}%
\e{2}%
\e{0}%
\e{0}%
\e{0}%
\e{0}%
\e{0}%
\e{1}%
\e{0}%
\e{3}%
\e{1}%
\e{3}%
\e{2}%
\e{2}%
\e{1}%
\e{1}%
\e{1}%
\e{0}%
\eol}\vss}\rg%
%
%
\rx{\vss\hfull{%
\rlx{\hss{$2240_{x}$}}\cg%
\e{1}%
\e{1}%
\e{0}%
\e{0}%
\e{0}%
\e{0}%
\e{0}%
\e{1}%
\e{0}%
\e{2}%
\e{1}%
\e{3}%
\e{1}%
\e{1}%
\e{0}%
\e{1}%
\e{0}%
\e{0}%
\eol}\vss}\rg%
%
%
\rx{\vss\hfull{%
\rlx{\hss{$2835_{x}$}}\cg%
\e{2}%
\e{2}%
\e{0}%
\e{0}%
\e{0}%
\e{0}%
\e{0}%
\e{0}%
\e{0}%
\e{1}%
\e{0}%
\e{3}%
\e{2}%
\e{1}%
\e{1}%
\e{1}%
\e{1}%
\e{0}%
\eol}\vss}\rg%
%
%
\rx{\vss\hfull{%
\rlx{\hss{$6075_{x}$}}\cg%
\e{1}%
\e{2}%
\e{1}%
\e{0}%
\e{0}%
\e{0}%
\e{0}%
\e{0}%
\e{1}%
\e{2}%
\e{3}%
\e{2}%
\e{3}%
\e{3}%
\e{2}%
\e{1}%
\e{1}%
\e{1}%
\eol}\vss}\rg%
%
%
\rx{\vss\hfull{%
\rlx{\hss{$3200_{x}$}}\cg%
\e{0}%
\e{1}%
\e{0}%
\e{1}%
\e{0}%
\e{0}%
\e{0}%
\e{0}%
\e{0}%
\e{1}%
\e{1}%
\e{2}%
\e{0}%
\e{1}%
\e{1}%
\e{2}%
\e{0}%
\e{0}%
\eol}\vss}\rg%
%
%
\rx{\vss\hfull{%
\rlx{\hss{$70_{y}$}}\cg%
\e{0}%
\e{0}%
\e{0}%
\e{0}%
\e{0}%
\e{0}%
\e{0}%
\e{0}%
\e{0}%
\e{0}%
\e{0}%
\e{0}%
\e{0}%
\e{0}%
\e{0}%
\e{0}%
\e{0}%
\e{0}%
\eol}\vss}\rg%
%
%
\rx{\vss\hfull{%
\rlx{\hss{$1134_{y}$}}\cg%
\e{0}%
\e{0}%
\e{0}%
\e{0}%
\e{1}%
\e{0}%
\e{0}%
\e{0}%
\e{0}%
\e{0}%
\e{1}%
\e{0}%
\e{1}%
\e{0}%
\e{0}%
\e{1}%
\e{0}%
\e{1}%
\eol}\vss}\rg%
%
%
\rx{\vss\hfull{%
\rlx{\hss{$1680_{y}$}}\cg%
\e{0}%
\e{0}%
\e{1}%
\e{0}%
\e{0}%
\e{0}%
\e{0}%
\e{0}%
\e{0}%
\e{0}%
\e{1}%
\e{0}%
\e{0}%
\e{1}%
\e{1}%
\e{0}%
\e{0}%
\e{1}%
\eol}\vss}\rg%
%
%
\rx{\vss\hfull{%
\rlx{\hss{$168_{y}$}}\cg%
\e{0}%
\e{0}%
\e{0}%
\e{0}%
\e{0}%
\e{0}%
\e{0}%
\e{0}%
\e{0}%
\e{0}%
\e{0}%
\e{0}%
\e{0}%
\e{0}%
\e{0}%
\e{0}%
\e{0}%
\e{0}%
\eol}\vss}\rg%
%
%
\rx{\vss\hfull{%
\rlx{\hss{$420_{y}$}}\cg%
\e{1}%
\e{0}%
\e{0}%
\e{1}%
\e{0}%
\e{0}%
\e{0}%
\e{0}%
\e{0}%
\e{0}%
\e{0}%
\e{1}%
\e{0}%
\e{0}%
\e{0}%
\e{0}%
\e{1}%
\e{0}%
\eol}\vss}\rg%
%
%
\rx{\vss\hfull{%
\rlx{\hss{$3150_{y}$}}\cg%
\e{1}%
\e{2}%
\e{0}%
\e{1}%
\e{1}%
\e{0}%
\e{0}%
\e{0}%
\e{0}%
\e{0}%
\e{1}%
\e{2}%
\e{2}%
\e{1}%
\e{1}%
\e{2}%
\e{2}%
\e{1}%
\eol}\vss}\rg%
%
%
\rx{\vss\hfull{%
\rlx{\hss{$4200_{y}$}}\cg%
\e{1}%
\e{2}%
\e{0}%
\e{1}%
\e{1}%
\e{0}%
\e{0}%
\e{0}%
\e{0}%
\e{1}%
\e{0}%
\e{2}%
\e{2}%
\e{1}%
\e{1}%
\e{2}%
\e{2}%
\e{0}%
\eol}\vss}\rg%
\eop
\eject
\tablecont%
%
%
%
%
%
%
\rowpts=18 true pt%
\colpts=18 true pt%
\rowlabpts=40 true pt%
\collabpts=65 true pt%
\clx{\vss\hfull{%
\rlx{\hss{$ $}}\cg%
\cx{\hskip 16 true pt\flip{$[{3^{2}}:{1}]$}\hss}\cg%
\cx{\hskip 16 true pt\flip{$[{3}{2}{1}:{1}]$}\hss}\cg%
\cx{\hskip 16 true pt\flip{$[{3}{1^{3}}:{1}]$}\hss}\cg%
\cx{\hskip 16 true pt\flip{$[{2^{3}}:{1}]$}\hss}\cg%
\cx{\hskip 16 true pt\flip{$[{2^{2}}{1^{2}}:{1}]$}\hss}\cg%
\cx{\hskip 16 true pt\flip{$[{2}{1^{4}}:{1}]$}\hss}\cg%
\cx{\hskip 16 true pt\flip{$[{1^{6}}:{1}]$}\hss}\cg%
\cx{\hskip 16 true pt\flip{$[{5}:{2}]$}\hss}\cg%
\cx{\hskip 16 true pt\flip{$[{5}:{1^{2}}]$}\hss}\cg%
\cx{\hskip 16 true pt\flip{$[{4}{1}:{2}]$}\hss}\cg%
\cx{\hskip 16 true pt\flip{$[{4}{1}:{1^{2}}]$}\hss}\cg%
\cx{\hskip 16 true pt\flip{$[{3}{2}:{2}]$}\hss}\cg%
\cx{\hskip 16 true pt\flip{$[{3}{2}:{1^{2}}]$}\hss}\cg%
\cx{\hskip 16 true pt\flip{$[{3}{1^{2}}:{2}]$}\hss}\cg%
\cx{\hskip 16 true pt\flip{$[{3}{1^{2}}:{1^{2}}]$}\hss}\cg%
\cx{\hskip 16 true pt\flip{$[{2^{2}}{1}:{2}]$}\hss}\cg%
\cx{\hskip 16 true pt\flip{$[{2^{2}}{1}:{1^{2}}]$}\hss}\cg%
\cx{\hskip 16 true pt\flip{$[{2}{1^{3}}:{2}]$}\hss}\cg%
\eol}}\rg%
%
%
\rx{\vss\hfull{%
\rlx{\hss{$2688_{y}$}}\cg%
\e{0}%
\e{2}%
\e{0}%
\e{0}%
\e{0}%
\e{0}%
\e{0}%
\e{0}%
\e{0}%
\e{0}%
\e{0}%
\e{1}%
\e{1}%
\e{1}%
\e{1}%
\e{1}%
\e{1}%
\e{0}%
\eol}\vss}\rg%
%
%
\rx{\vss\hfull{%
\rlx{\hss{$2100_{y}$}}\cg%
\e{0}%
\e{0}%
\e{1}%
\e{0}%
\e{0}%
\e{0}%
\e{0}%
\e{0}%
\e{0}%
\e{0}%
\e{1}%
\e{0}%
\e{0}%
\e{1}%
\e{1}%
\e{0}%
\e{0}%
\e{1}%
\eol}\vss}\rg%
%
%
\rx{\vss\hfull{%
\rlx{\hss{$1400_{y}$}}\cg%
\e{0}%
\e{0}%
\e{1}%
\e{0}%
\e{0}%
\e{0}%
\e{0}%
\e{0}%
\e{0}%
\e{1}%
\e{0}%
\e{0}%
\e{0}%
\e{1}%
\e{1}%
\e{0}%
\e{0}%
\e{0}%
\eol}\vss}\rg%
%
%
\rx{\vss\hfull{%
\rlx{\hss{$4536_{y}$}}\cg%
\e{0}%
\e{2}%
\e{2}%
\e{0}%
\e{0}%
\e{0}%
\e{0}%
\e{0}%
\e{0}%
\e{1}%
\e{1}%
\e{1}%
\e{1}%
\e{3}%
\e{3}%
\e{1}%
\e{1}%
\e{1}%
\eol}\vss}\rg%
%
%
\rx{\vss\hfull{%
\rlx{\hss{$5670_{y}$}}\cg%
\e{0}%
\e{2}%
\e{2}%
\e{0}%
\e{1}%
\e{0}%
\e{0}%
\e{0}%
\e{0}%
\e{1}%
\e{2}%
\e{1}%
\e{2}%
\e{3}%
\e{3}%
\e{2}%
\e{1}%
\e{2}%
\eol}\vss}\rg%
%
%
\rx{\vss\hfull{%
\rlx{\hss{$4480_{y}$}}\cg%
\e{1}%
\e{2}%
\e{1}%
\e{1}%
\e{1}%
\e{0}%
\e{0}%
\e{0}%
\e{0}%
\e{1}%
\e{1}%
\e{2}%
\e{2}%
\e{2}%
\e{2}%
\e{2}%
\e{2}%
\e{1}%
\eol}\vss}\rg%
%
%
\rx{\vss\hfull{%
\rlx{\hss{$8_{z}$}}\cg%
\e{0}%
\e{0}%
\e{0}%
\e{0}%
\e{0}%
\e{0}%
\e{0}%
\e{0}%
\e{0}%
\e{0}%
\e{0}%
\e{0}%
\e{0}%
\e{0}%
\e{0}%
\e{0}%
\e{0}%
\e{0}%
\eol}\vss}\rg%
%
%
\rx{\vss\hfull{%
\rlx{\hss{$56_{z}$}}\cg%
\e{0}%
\e{0}%
\e{0}%
\e{0}%
\e{0}%
\e{0}%
\e{0}%
\e{0}%
\e{1}%
\e{0}%
\e{0}%
\e{0}%
\e{0}%
\e{0}%
\e{0}%
\e{0}%
\e{0}%
\e{0}%
\eol}\vss}\rg%
%
%
\rx{\vss\hfull{%
\rlx{\hss{$160_{z}$}}\cg%
\e{0}%
\e{0}%
\e{0}%
\e{0}%
\e{0}%
\e{0}%
\e{0}%
\e{1}%
\e{1}%
\e{0}%
\e{0}%
\e{0}%
\e{0}%
\e{0}%
\e{0}%
\e{0}%
\e{0}%
\e{0}%
\eol}\vss}\rg%
%
%
\rx{\vss\hfull{%
\rlx{\hss{$112_{z}$}}\cg%
\e{0}%
\e{0}%
\e{0}%
\e{0}%
\e{0}%
\e{0}%
\e{0}%
\e{1}%
\e{0}%
\e{0}%
\e{0}%
\e{0}%
\e{0}%
\e{0}%
\e{0}%
\e{0}%
\e{0}%
\e{0}%
\eol}\vss}\rg%
%
%
\rx{\vss\hfull{%
\rlx{\hss{$840_{z}$}}\cg%
\e{0}%
\e{0}%
\e{0}%
\e{0}%
\e{0}%
\e{0}%
\e{0}%
\e{0}%
\e{0}%
\e{1}%
\e{1}%
\e{0}%
\e{0}%
\e{0}%
\e{0}%
\e{1}%
\e{0}%
\e{0}%
\eol}\vss}\rg%
%
%
\rx{\vss\hfull{%
\rlx{\hss{$1296_{z}$}}\cg%
\e{0}%
\e{0}%
\e{0}%
\e{0}%
\e{0}%
\e{0}%
\e{0}%
\e{1}%
\e{2}%
\e{1}%
\e{2}%
\e{0}%
\e{0}%
\e{1}%
\e{0}%
\e{0}%
\e{0}%
\e{0}%
\eol}\vss}\rg%
%
%
\rx{\vss\hfull{%
\rlx{\hss{$1400_{z}$}}\cg%
\e{0}%
\e{0}%
\e{0}%
\e{0}%
\e{0}%
\e{0}%
\e{0}%
\e{2}%
\e{1}%
\e{2}%
\e{1}%
\e{1}%
\e{0}%
\e{0}%
\e{0}%
\e{0}%
\e{0}%
\e{0}%
\eol}\vss}\rg%
%
%
\rx{\vss\hfull{%
\rlx{\hss{$1008_{z}$}}\cg%
\e{0}%
\e{0}%
\e{0}%
\e{0}%
\e{0}%
\e{0}%
\e{0}%
\e{1}%
\e{2}%
\e{2}%
\e{1}%
\e{0}%
\e{0}%
\e{0}%
\e{0}%
\e{0}%
\e{0}%
\e{0}%
\eol}\vss}\rg%
%
%
\rx{\vss\hfull{%
\rlx{\hss{$560_{z}$}}\cg%
\e{0}%
\e{0}%
\e{0}%
\e{0}%
\e{0}%
\e{0}%
\e{0}%
\e{2}%
\e{1}%
\e{1}%
\e{0}%
\e{0}%
\e{0}%
\e{0}%
\e{0}%
\e{0}%
\e{0}%
\e{0}%
\eol}\vss}\rg%
%
%
\rx{\vss\hfull{%
\rlx{\hss{$1400_{zz}$}}\cg%
\e{2}%
\e{1}%
\e{0}%
\e{0}%
\e{0}%
\e{0}%
\e{0}%
\e{1}%
\e{0}%
\e{1}%
\e{0}%
\e{2}%
\e{1}%
\e{0}%
\e{0}%
\e{0}%
\e{0}%
\e{0}%
\eol}\vss}\rg%
%
%
\rx{\vss\hfull{%
\rlx{\hss{$4200_{z}$}}\cg%
\e{2}%
\e{2}%
\e{1}%
\e{0}%
\e{1}%
\e{0}%
\e{0}%
\e{0}%
\e{0}%
\e{2}%
\e{1}%
\e{2}%
\e{3}%
\e{2}%
\e{1}%
\e{1}%
\e{1}%
\e{0}%
\eol}\vss}\rg%
%
%
\rx{\vss\hfull{%
\rlx{\hss{$400_{z}$}}\cg%
\e{1}%
\e{0}%
\e{0}%
\e{0}%
\e{0}%
\e{0}%
\e{0}%
\e{1}%
\e{0}%
\e{1}%
\e{0}%
\e{0}%
\e{0}%
\e{0}%
\e{0}%
\e{0}%
\e{0}%
\e{0}%
\eol}\vss}\rg%
%
%
\rx{\vss\hfull{%
\rlx{\hss{$3240_{z}$}}\cg%
\e{1}%
\e{1}%
\e{0}%
\e{0}%
\e{0}%
\e{0}%
\e{0}%
\e{2}%
\e{1}%
\e{4}%
\e{2}%
\e{2}%
\e{1}%
\e{1}%
\e{0}%
\e{0}%
\e{0}%
\e{0}%
\eol}\vss}\rg%
%
%
\rx{\vss\hfull{%
\rlx{\hss{$4536_{z}$}}\cg%
\e{1}%
\e{2}%
\e{0}%
\e{1}%
\e{0}%
\e{0}%
\e{0}%
\e{0}%
\e{0}%
\e{2}%
\e{1}%
\e{4}%
\e{2}%
\e{1}%
\e{1}%
\e{2}%
\e{1}%
\e{0}%
\eol}\vss}\rg%
%
%
\rx{\vss\hfull{%
\rlx{\hss{$2400_{z}$}}\cg%
\e{0}%
\e{0}%
\e{1}%
\e{0}%
\e{0}%
\e{0}%
\e{0}%
\e{0}%
\e{1}%
\e{1}%
\e{2}%
\e{0}%
\e{1}%
\e{2}%
\e{1}%
\e{0}%
\e{0}%
\e{1}%
\eol}\vss}\rg%
%
%
\rx{\vss\hfull{%
\rlx{\hss{$3360_{z}$}}\cg%
\e{1}%
\e{2}%
\e{0}%
\e{0}%
\e{0}%
\e{0}%
\e{0}%
\e{1}%
\e{1}%
\e{2}%
\e{1}%
\e{2}%
\e{2}%
\e{2}%
\e{1}%
\e{0}%
\e{0}%
\e{0}%
\eol}\vss}\rg%
%
%
\rx{\vss\hfull{%
\rlx{\hss{$2800_{z}$}}\cg%
\e{0}%
\e{0}%
\e{1}%
\e{0}%
\e{0}%
\e{0}%
\e{0}%
\e{0}%
\e{1}%
\e{2}%
\e{2}%
\e{1}%
\e{1}%
\e{2}%
\e{1}%
\e{0}%
\e{0}%
\e{0}%
\eol}\vss}\rg%
%
%
\rx{\vss\hfull{%
\rlx{\hss{$4096_{z}$}}\cg%
\e{0}%
\e{1}%
\e{0}%
\e{0}%
\e{0}%
\e{0}%
\e{0}%
\e{1}%
\e{1}%
\e{3}%
\e{3}%
\e{2}%
\e{1}%
\e{2}%
\e{1}%
\e{1}%
\e{0}%
\e{0}%
\eol}\vss}\rg%
%
%
\rx{\vss\hfull{%
\rlx{\hss{$5600_{z}$}}\cg%
\e{0}%
\e{1}%
\e{1}%
\e{0}%
\e{0}%
\e{0}%
\e{0}%
\e{0}%
\e{0}%
\e{2}%
\e{3}%
\e{2}%
\e{1}%
\e{3}%
\e{3}%
\e{2}%
\e{1}%
\e{1}%
\eol}\vss}\rg%
%
%
\rx{\vss\hfull{%
\rlx{\hss{$448_{z}$}}\cg%
\e{0}%
\e{0}%
\e{0}%
\e{0}%
\e{0}%
\e{0}%
\e{0}%
\e{1}%
\e{0}%
\e{0}%
\e{0}%
\e{1}%
\e{0}%
\e{0}%
\e{0}%
\e{0}%
\e{0}%
\e{0}%
\eol}\vss}\rg%
%
%
\rx{\vss\hfull{%
\rlx{\hss{$448_{w}$}}\cg%
\e{0}%
\e{0}%
\e{0}%
\e{0}%
\e{0}%
\e{0}%
\e{0}%
\e{0}%
\e{0}%
\e{0}%
\e{1}%
\e{0}%
\e{0}%
\e{0}%
\e{0}%
\e{0}%
\e{0}%
\e{1}%
\eol}\vss}\rg%
%
%
\rx{\vss\hfull{%
\rlx{\hss{$1344_{w}$}}\cg%
\e{1}%
\e{0}%
\e{0}%
\e{1}%
\e{1}%
\e{0}%
\e{0}%
\e{0}%
\e{0}%
\e{0}%
\e{0}%
\e{1}%
\e{1}%
\e{0}%
\e{0}%
\e{1}%
\e{1}%
\e{0}%
\eol}\vss}\rg%
%
%
\rx{\vss\hfull{%
\rlx{\hss{$5600_{w}$}}\cg%
\e{0}%
\e{2}%
\e{1}%
\e{0}%
\e{1}%
\e{0}%
\e{0}%
\e{0}%
\e{0}%
\e{1}%
\e{2}%
\e{1}%
\e{2}%
\e{3}%
\e{3}%
\e{2}%
\e{1}%
\e{2}%
\eol}\vss}\rg%
%
%
\rx{\vss\hfull{%
\rlx{\hss{$2016_{w}$}}\cg%
\e{1}%
\e{2}%
\e{0}%
\e{1}%
\e{0}%
\e{0}%
\e{0}%
\e{0}%
\e{0}%
\e{0}%
\e{0}%
\e{2}%
\e{1}%
\e{0}%
\e{0}%
\e{1}%
\e{2}%
\e{0}%
\eol}\vss}\rg%
%
%
\rx{\vss\hfull{%
\rlx{\hss{$7168_{w}$}}\cg%
\e{1}%
\e{4}%
\e{1}%
\e{1}%
\e{1}%
\e{0}%
\e{0}%
\e{0}%
\e{0}%
\e{1}%
\e{1}%
\e{3}%
\e{3}%
\e{3}%
\e{3}%
\e{3}%
\e{3}%
\e{1}%
\eol}\vss}\rg%
\eop
\eject
\tablecont%
%
%
%
%
%
%
\rowpts=18 true pt%
\colpts=18 true pt%
\rowlabpts=40 true pt%
\collabpts=65 true pt%
\clx{\vss\hfull{%
\rlx{\hss{$ $}}\cg%
\cx{\hskip 16 true pt\flip{$[{2}{1^{3}}:{1^{2}}]$}\hss}\cg%
\cx{\hskip 16 true pt\flip{$[{1^{5}}:{2}]$}\hss}\cg%
\cx{\hskip 16 true pt\flip{$[{1^{5}}:{1^{2}}]$}\hss}\cg%
\cx{\hskip 16 true pt\flip{$[{4}:{3}]$}\hss}\cg%
\cx{\hskip 16 true pt\flip{$[{4}:{2}{1}]$}\hss}\cg%
\cx{\hskip 16 true pt\flip{$[{4}:{1^{3}}]$}\hss}\cg%
\cx{\hskip 16 true pt\flip{$[{3}{1}:{3}]$}\hss}\cg%
\cx{\hskip 16 true pt\flip{$[{3}{1}:{2}{1}]$}\hss}\cg%
\cx{\hskip 16 true pt\flip{$[{3}{1}:{1^{3}}]$}\hss}\cg%
\cx{\hskip 16 true pt\flip{$[{2^{2}}:{3}]$}\hss}\cg%
\cx{\hskip 16 true pt\flip{$[{2^{2}}:{2}{1}]$}\hss}\cg%
\cx{\hskip 16 true pt\flip{$[{2^{2}}:{1^{3}}]$}\hss}\cg%
\cx{\hskip 16 true pt\flip{$[{2}{1^{2}}:{3}]$}\hss}\cg%
\cx{\hskip 16 true pt\flip{$[{2}{1^{2}}:{2}{1}]$}\hss}\cg%
\cx{\hskip 16 true pt\flip{$[{2}{1^{2}}:{1^{3}}]$}\hss}\cg%
\cx{\hskip 16 true pt\flip{$[{1^{4}}:{3}]$}\hss}\cg%
\cx{\hskip 16 true pt\flip{$[{1^{4}}:{2}{1}]$}\hss}\cg%
\cx{\hskip 16 true pt\flip{$[{1^{4}}:{1^{3}}]$}\hss}\cg%
\eol}}\rg%
%
%
\rx{\vss\hfull{%
\rlx{\hss{$1_{x}$}}\cg%
\e{0}%
\e{0}%
\e{0}%
\e{0}%
\e{0}%
\e{0}%
\e{0}%
\e{0}%
\e{0}%
\e{0}%
\e{0}%
\e{0}%
\e{0}%
\e{0}%
\e{0}%
\e{0}%
\e{0}%
\e{0}%
\eol}\vss}\rg%
%
%
\rx{\vss\hfull{%
\rlx{\hss{$28_{x}$}}\cg%
\e{0}%
\e{0}%
\e{0}%
\e{0}%
\e{0}%
\e{0}%
\e{0}%
\e{0}%
\e{0}%
\e{0}%
\e{0}%
\e{0}%
\e{0}%
\e{0}%
\e{0}%
\e{0}%
\e{0}%
\e{0}%
\eol}\vss}\rg%
%
%
\rx{\vss\hfull{%
\rlx{\hss{$35_{x}$}}\cg%
\e{0}%
\e{0}%
\e{0}%
\e{0}%
\e{0}%
\e{0}%
\e{0}%
\e{0}%
\e{0}%
\e{0}%
\e{0}%
\e{0}%
\e{0}%
\e{0}%
\e{0}%
\e{0}%
\e{0}%
\e{0}%
\eol}\vss}\rg%
%
%
\rx{\vss\hfull{%
\rlx{\hss{$84_{x}$}}\cg%
\e{0}%
\e{0}%
\e{0}%
\e{1}%
\e{0}%
\e{0}%
\e{0}%
\e{0}%
\e{0}%
\e{0}%
\e{0}%
\e{0}%
\e{0}%
\e{0}%
\e{0}%
\e{0}%
\e{0}%
\e{0}%
\eol}\vss}\rg%
%
%
\rx{\vss\hfull{%
\rlx{\hss{$50_{x}$}}\cg%
\e{0}%
\e{0}%
\e{0}%
\e{1}%
\e{0}%
\e{0}%
\e{0}%
\e{0}%
\e{0}%
\e{0}%
\e{0}%
\e{0}%
\e{0}%
\e{0}%
\e{0}%
\e{0}%
\e{0}%
\e{0}%
\eol}\vss}\rg%
%
%
\rx{\vss\hfull{%
\rlx{\hss{$350_{x}$}}\cg%
\e{0}%
\e{0}%
\e{0}%
\e{0}%
\e{1}%
\e{1}%
\e{0}%
\e{0}%
\e{0}%
\e{0}%
\e{0}%
\e{0}%
\e{1}%
\e{0}%
\e{0}%
\e{0}%
\e{0}%
\e{0}%
\eol}\vss}\rg%
%
%
\rx{\vss\hfull{%
\rlx{\hss{$300_{x}$}}\cg%
\e{0}%
\e{0}%
\e{0}%
\e{0}%
\e{1}%
\e{0}%
\e{0}%
\e{0}%
\e{0}%
\e{1}%
\e{0}%
\e{0}%
\e{0}%
\e{0}%
\e{0}%
\e{0}%
\e{0}%
\e{0}%
\eol}\vss}\rg%
%
%
\rx{\vss\hfull{%
\rlx{\hss{$567_{x}$}}\cg%
\e{0}%
\e{0}%
\e{0}%
\e{1}%
\e{1}%
\e{0}%
\e{1}%
\e{0}%
\e{0}%
\e{0}%
\e{0}%
\e{0}%
\e{0}%
\e{0}%
\e{0}%
\e{0}%
\e{0}%
\e{0}%
\eol}\vss}\rg%
%
%
\rx{\vss\hfull{%
\rlx{\hss{$210_{x}$}}\cg%
\e{0}%
\e{0}%
\e{0}%
\e{1}%
\e{0}%
\e{0}%
\e{0}%
\e{0}%
\e{0}%
\e{0}%
\e{0}%
\e{0}%
\e{0}%
\e{0}%
\e{0}%
\e{0}%
\e{0}%
\e{0}%
\eol}\vss}\rg%
%
%
\rx{\vss\hfull{%
\rlx{\hss{$840_{x}$}}\cg%
\e{0}%
\e{0}%
\e{0}%
\e{0}%
\e{0}%
\e{0}%
\e{0}%
\e{1}%
\e{0}%
\e{1}%
\e{1}%
\e{0}%
\e{0}%
\e{0}%
\e{0}%
\e{0}%
\e{0}%
\e{0}%
\eol}\vss}\rg%
%
%
\rx{\vss\hfull{%
\rlx{\hss{$700_{x}$}}\cg%
\e{0}%
\e{0}%
\e{0}%
\e{2}%
\e{1}%
\e{0}%
\e{1}%
\e{0}%
\e{0}%
\e{0}%
\e{0}%
\e{0}%
\e{0}%
\e{0}%
\e{0}%
\e{0}%
\e{0}%
\e{0}%
\eol}\vss}\rg%
%
%
\rx{\vss\hfull{%
\rlx{\hss{$175_{x}$}}\cg%
\e{0}%
\e{0}%
\e{0}%
\e{1}%
\e{0}%
\e{0}%
\e{0}%
\e{0}%
\e{0}%
\e{0}%
\e{0}%
\e{0}%
\e{0}%
\e{0}%
\e{0}%
\e{0}%
\e{0}%
\e{0}%
\eol}\vss}\rg%
%
%
\rx{\vss\hfull{%
\rlx{\hss{$1400_{x}$}}\cg%
\e{0}%
\e{0}%
\e{0}%
\e{2}%
\e{1}%
\e{0}%
\e{2}%
\e{1}%
\e{0}%
\e{0}%
\e{0}%
\e{0}%
\e{0}%
\e{0}%
\e{0}%
\e{0}%
\e{0}%
\e{0}%
\eol}\vss}\rg%
%
%
\rx{\vss\hfull{%
\rlx{\hss{$1050_{x}$}}\cg%
\e{0}%
\e{0}%
\e{0}%
\e{2}%
\e{1}%
\e{0}%
\e{2}%
\e{1}%
\e{0}%
\e{0}%
\e{0}%
\e{0}%
\e{0}%
\e{0}%
\e{0}%
\e{0}%
\e{0}%
\e{0}%
\eol}\vss}\rg%
%
%
\rx{\vss\hfull{%
\rlx{\hss{$1575_{x}$}}\cg%
\e{0}%
\e{0}%
\e{0}%
\e{1}%
\e{2}%
\e{1}%
\e{2}%
\e{1}%
\e{0}%
\e{0}%
\e{0}%
\e{0}%
\e{1}%
\e{0}%
\e{0}%
\e{0}%
\e{0}%
\e{0}%
\eol}\vss}\rg%
%
%
\rx{\vss\hfull{%
\rlx{\hss{$1344_{x}$}}\cg%
\e{0}%
\e{0}%
\e{0}%
\e{2}%
\e{2}%
\e{0}%
\e{2}%
\e{1}%
\e{0}%
\e{1}%
\e{0}%
\e{0}%
\e{0}%
\e{0}%
\e{0}%
\e{0}%
\e{0}%
\e{0}%
\eol}\vss}\rg%
%
%
\rx{\vss\hfull{%
\rlx{\hss{$2100_{x}$}}\cg%
\e{0}%
\e{0}%
\e{0}%
\e{0}%
\e{2}%
\e{2}%
\e{1}%
\e{2}%
\e{1}%
\e{1}%
\e{0}%
\e{0}%
\e{2}%
\e{1}%
\e{0}%
\e{1}%
\e{0}%
\e{0}%
\eol}\vss}\rg%
%
%
\rx{\vss\hfull{%
\rlx{\hss{$2268_{x}$}}\cg%
\e{0}%
\e{0}%
\e{0}%
\e{1}%
\e{2}%
\e{1}%
\e{2}%
\e{2}%
\e{0}%
\e{1}%
\e{1}%
\e{0}%
\e{1}%
\e{0}%
\e{0}%
\e{0}%
\e{0}%
\e{0}%
\eol}\vss}\rg%
%
%
\rx{\vss\hfull{%
\rlx{\hss{$525_{x}$}}\cg%
\e{0}%
\e{0}%
\e{0}%
\e{1}%
\e{0}%
\e{0}%
\e{1}%
\e{1}%
\e{0}%
\e{0}%
\e{0}%
\e{0}%
\e{0}%
\e{0}%
\e{0}%
\e{0}%
\e{0}%
\e{0}%
\eol}\vss}\rg%
%
%
\rx{\vss\hfull{%
\rlx{\hss{$700_{xx}$}}\cg%
\e{0}%
\e{0}%
\e{0}%
\e{1}%
\e{0}%
\e{0}%
\e{1}%
\e{1}%
\e{0}%
\e{0}%
\e{1}%
\e{0}%
\e{0}%
\e{0}%
\e{0}%
\e{0}%
\e{0}%
\e{0}%
\eol}\vss}\rg%
%
%
\rx{\vss\hfull{%
\rlx{\hss{$972_{x}$}}\cg%
\e{0}%
\e{0}%
\e{0}%
\e{1}%
\e{1}%
\e{0}%
\e{1}%
\e{1}%
\e{0}%
\e{1}%
\e{1}%
\e{0}%
\e{0}%
\e{0}%
\e{0}%
\e{0}%
\e{0}%
\e{0}%
\eol}\vss}\rg%
%
%
\rx{\vss\hfull{%
\rlx{\hss{$4096_{x}$}}\cg%
\e{0}%
\e{0}%
\e{0}%
\e{1}%
\e{3}%
\e{1}%
\e{3}%
\e{4}%
\e{1}%
\e{2}%
\e{1}%
\e{0}%
\e{2}%
\e{1}%
\e{0}%
\e{0}%
\e{0}%
\e{0}%
\eol}\vss}\rg%
%
%
\rx{\vss\hfull{%
\rlx{\hss{$4200_{x}$}}\cg%
\e{0}%
\e{0}%
\e{0}%
\e{1}%
\e{2}%
\e{0}%
\e{3}%
\e{4}%
\e{1}%
\e{2}%
\e{2}%
\e{0}%
\e{1}%
\e{1}%
\e{0}%
\e{0}%
\e{0}%
\e{0}%
\eol}\vss}\rg%
%
%
\rx{\vss\hfull{%
\rlx{\hss{$2240_{x}$}}\cg%
\e{0}%
\e{0}%
\e{0}%
\e{2}%
\e{1}%
\e{0}%
\e{2}%
\e{2}%
\e{0}%
\e{1}%
\e{1}%
\e{0}%
\e{0}%
\e{0}%
\e{0}%
\e{0}%
\e{0}%
\e{0}%
\eol}\vss}\rg%
%
%
\rx{\vss\hfull{%
\rlx{\hss{$2835_{x}$}}\cg%
\e{0}%
\e{0}%
\e{0}%
\e{1}%
\e{1}%
\e{0}%
\e{2}%
\e{2}%
\e{0}%
\e{1}%
\e{2}%
\e{1}%
\e{0}%
\e{1}%
\e{0}%
\e{0}%
\e{0}%
\e{0}%
\eol}\vss}\rg%
%
%
\rx{\vss\hfull{%
\rlx{\hss{$6075_{x}$}}\cg%
\e{0}%
\e{0}%
\e{0}%
\e{1}%
\e{2}%
\e{1}%
\e{4}%
\e{6}%
\e{2}%
\e{1}%
\e{2}%
\e{1}%
\e{3}%
\e{3}%
\e{0}%
\e{0}%
\e{0}%
\e{0}%
\eol}\vss}\rg%
%
%
\rx{\vss\hfull{%
\rlx{\hss{$3200_{x}$}}\cg%
\e{0}%
\e{0}%
\e{0}%
\e{0}%
\e{1}%
\e{0}%
\e{1}%
\e{3}%
\e{1}%
\e{2}%
\e{2}%
\e{0}%
\e{1}%
\e{2}%
\e{1}%
\e{0}%
\e{0}%
\e{0}%
\eol}\vss}\rg%
%
%
\rx{\vss\hfull{%
\rlx{\hss{$70_{y}$}}\cg%
\e{0}%
\e{0}%
\e{0}%
\e{0}%
\e{0}%
\e{1}%
\e{0}%
\e{0}%
\e{0}%
\e{0}%
\e{0}%
\e{0}%
\e{0}%
\e{0}%
\e{0}%
\e{1}%
\e{0}%
\e{0}%
\eol}\vss}\rg%
%
%
\rx{\vss\hfull{%
\rlx{\hss{$1134_{y}$}}\cg%
\e{0}%
\e{0}%
\e{0}%
\e{0}%
\e{0}%
\e{0}%
\e{0}%
\e{1}%
\e{1}%
\e{0}%
\e{0}%
\e{0}%
\e{1}%
\e{1}%
\e{0}%
\e{0}%
\e{0}%
\e{0}%
\eol}\vss}\rg%
%
%
\rx{\vss\hfull{%
\rlx{\hss{$1680_{y}$}}\cg%
\e{0}%
\e{0}%
\e{0}%
\e{0}%
\e{1}%
\e{2}%
\e{0}%
\e{1}%
\e{2}%
\e{0}%
\e{0}%
\e{0}%
\e{2}%
\e{1}%
\e{0}%
\e{2}%
\e{1}%
\e{0}%
\eol}\vss}\rg%
%
%
\rx{\vss\hfull{%
\rlx{\hss{$168_{y}$}}\cg%
\e{0}%
\e{0}%
\e{0}%
\e{0}%
\e{0}%
\e{0}%
\e{0}%
\e{0}%
\e{0}%
\e{0}%
\e{1}%
\e{0}%
\e{0}%
\e{0}%
\e{0}%
\e{0}%
\e{0}%
\e{0}%
\eol}\vss}\rg%
%
%
\rx{\vss\hfull{%
\rlx{\hss{$420_{y}$}}\cg%
\e{0}%
\e{0}%
\e{0}%
\e{0}%
\e{0}%
\e{0}%
\e{0}%
\e{0}%
\e{0}%
\e{0}%
\e{1}%
\e{0}%
\e{0}%
\e{0}%
\e{0}%
\e{0}%
\e{0}%
\e{0}%
\eol}\vss}\rg%
%
%
\rx{\vss\hfull{%
\rlx{\hss{$3150_{y}$}}\cg%
\e{0}%
\e{0}%
\e{0}%
\e{0}%
\e{0}%
\e{0}%
\e{1}%
\e{2}%
\e{0}%
\e{1}%
\e{2}%
\e{1}%
\e{0}%
\e{2}%
\e{1}%
\e{0}%
\e{0}%
\e{0}%
\eol}\vss}\rg%
%
%
\rx{\vss\hfull{%
\rlx{\hss{$4200_{y}$}}\cg%
\e{1}%
\e{0}%
\e{0}%
\e{0}%
\e{0}%
\e{0}%
\e{1}%
\e{3}%
\e{1}%
\e{1}%
\e{4}%
\e{1}%
\e{1}%
\e{3}%
\e{1}%
\e{0}%
\e{0}%
\e{0}%
\eol}\vss}\rg%
\eop
\eject
\tablecont%
%
%
%
%
%
%
\rowpts=18 true pt%
\colpts=18 true pt%
\rowlabpts=40 true pt%
\collabpts=65 true pt%
\clx{\vss\hfull{%
\rlx{\hss{$ $}}\cg%
\cx{\hskip 16 true pt\flip{$[{2}{1^{3}}:{1^{2}}]$}\hss}\cg%
\cx{\hskip 16 true pt\flip{$[{1^{5}}:{2}]$}\hss}\cg%
\cx{\hskip 16 true pt\flip{$[{1^{5}}:{1^{2}}]$}\hss}\cg%
\cx{\hskip 16 true pt\flip{$[{4}:{3}]$}\hss}\cg%
\cx{\hskip 16 true pt\flip{$[{4}:{2}{1}]$}\hss}\cg%
\cx{\hskip 16 true pt\flip{$[{4}:{1^{3}}]$}\hss}\cg%
\cx{\hskip 16 true pt\flip{$[{3}{1}:{3}]$}\hss}\cg%
\cx{\hskip 16 true pt\flip{$[{3}{1}:{2}{1}]$}\hss}\cg%
\cx{\hskip 16 true pt\flip{$[{3}{1}:{1^{3}}]$}\hss}\cg%
\cx{\hskip 16 true pt\flip{$[{2^{2}}:{3}]$}\hss}\cg%
\cx{\hskip 16 true pt\flip{$[{2^{2}}:{2}{1}]$}\hss}\cg%
\cx{\hskip 16 true pt\flip{$[{2^{2}}:{1^{3}}]$}\hss}\cg%
\cx{\hskip 16 true pt\flip{$[{2}{1^{2}}:{3}]$}\hss}\cg%
\cx{\hskip 16 true pt\flip{$[{2}{1^{2}}:{2}{1}]$}\hss}\cg%
\cx{\hskip 16 true pt\flip{$[{2}{1^{2}}:{1^{3}}]$}\hss}\cg%
\cx{\hskip 16 true pt\flip{$[{1^{4}}:{3}]$}\hss}\cg%
\cx{\hskip 16 true pt\flip{$[{1^{4}}:{2}{1}]$}\hss}\cg%
\cx{\hskip 16 true pt\flip{$[{1^{4}}:{1^{3}}]$}\hss}\cg%
\eol}}\rg%
%
%
\rx{\vss\hfull{%
\rlx{\hss{$2688_{y}$}}\cg%
\e{0}%
\e{0}%
\e{0}%
\e{0}%
\e{1}%
\e{0}%
\e{1}%
\e{2}%
\e{1}%
\e{1}%
\e{1}%
\e{1}%
\e{1}%
\e{2}%
\e{1}%
\e{0}%
\e{1}%
\e{0}%
\eol}\vss}\rg%
%
%
\rx{\vss\hfull{%
\rlx{\hss{$2100_{y}$}}\cg%
\e{0}%
\e{0}%
\e{0}%
\e{0}%
\e{0}%
\e{1}%
\e{1}%
\e{2}%
\e{1}%
\e{0}%
\e{1}%
\e{0}%
\e{1}%
\e{2}%
\e{1}%
\e{1}%
\e{0}%
\e{0}%
\eol}\vss}\rg%
%
%
\rx{\vss\hfull{%
\rlx{\hss{$1400_{y}$}}\cg%
\e{1}%
\e{0}%
\e{0}%
\e{0}%
\e{0}%
\e{1}%
\e{0}%
\e{1}%
\e{1}%
\e{0}%
\e{1}%
\e{0}%
\e{1}%
\e{1}%
\e{0}%
\e{1}%
\e{0}%
\e{0}%
\eol}\vss}\rg%
%
%
\rx{\vss\hfull{%
\rlx{\hss{$4536_{y}$}}\cg%
\e{1}%
\e{0}%
\e{0}%
\e{0}%
\e{1}%
\e{1}%
\e{1}%
\e{3}%
\e{2}%
\e{1}%
\e{2}%
\e{1}%
\e{2}%
\e{3}%
\e{1}%
\e{1}%
\e{1}%
\e{0}%
\eol}\vss}\rg%
%
%
\rx{\vss\hfull{%
\rlx{\hss{$5670_{y}$}}\cg%
\e{1}%
\e{0}%
\e{0}%
\e{0}%
\e{1}%
\e{1}%
\e{1}%
\e{4}%
\e{3}%
\e{1}%
\e{2}%
\e{1}%
\e{3}%
\e{4}%
\e{1}%
\e{1}%
\e{1}%
\e{0}%
\eol}\vss}\rg%
%
%
\rx{\vss\hfull{%
\rlx{\hss{$4480_{y}$}}\cg%
\e{1}%
\e{0}%
\e{0}%
\e{0}%
\e{0}%
\e{0}%
\e{1}%
\e{3}%
\e{1}%
\e{1}%
\e{3}%
\e{1}%
\e{1}%
\e{3}%
\e{1}%
\e{0}%
\e{0}%
\e{0}%
\eol}\vss}\rg%
%
%
\rx{\vss\hfull{%
\rlx{\hss{$8_{z}$}}\cg%
\e{0}%
\e{0}%
\e{0}%
\e{0}%
\e{0}%
\e{0}%
\e{0}%
\e{0}%
\e{0}%
\e{0}%
\e{0}%
\e{0}%
\e{0}%
\e{0}%
\e{0}%
\e{0}%
\e{0}%
\e{0}%
\eol}\vss}\rg%
%
%
\rx{\vss\hfull{%
\rlx{\hss{$56_{z}$}}\cg%
\e{0}%
\e{0}%
\e{0}%
\e{0}%
\e{0}%
\e{1}%
\e{0}%
\e{0}%
\e{0}%
\e{0}%
\e{0}%
\e{0}%
\e{0}%
\e{0}%
\e{0}%
\e{0}%
\e{0}%
\e{0}%
\eol}\vss}\rg%
%
%
\rx{\vss\hfull{%
\rlx{\hss{$160_{z}$}}\cg%
\e{0}%
\e{0}%
\e{0}%
\e{0}%
\e{1}%
\e{0}%
\e{0}%
\e{0}%
\e{0}%
\e{0}%
\e{0}%
\e{0}%
\e{0}%
\e{0}%
\e{0}%
\e{0}%
\e{0}%
\e{0}%
\eol}\vss}\rg%
%
%
\rx{\vss\hfull{%
\rlx{\hss{$112_{z}$}}\cg%
\e{0}%
\e{0}%
\e{0}%
\e{1}%
\e{0}%
\e{0}%
\e{0}%
\e{0}%
\e{0}%
\e{0}%
\e{0}%
\e{0}%
\e{0}%
\e{0}%
\e{0}%
\e{0}%
\e{0}%
\e{0}%
\eol}\vss}\rg%
%
%
\rx{\vss\hfull{%
\rlx{\hss{$840_{z}$}}\cg%
\e{0}%
\e{0}%
\e{0}%
\e{0}%
\e{1}%
\e{0}%
\e{0}%
\e{1}%
\e{0}%
\e{1}%
\e{0}%
\e{0}%
\e{1}%
\e{0}%
\e{0}%
\e{0}%
\e{0}%
\e{0}%
\eol}\vss}\rg%
%
%
\rx{\vss\hfull{%
\rlx{\hss{$1296_{z}$}}\cg%
\e{0}%
\e{0}%
\e{0}%
\e{0}%
\e{2}%
\e{2}%
\e{1}%
\e{1}%
\e{1}%
\e{0}%
\e{0}%
\e{0}%
\e{1}%
\e{0}%
\e{0}%
\e{0}%
\e{0}%
\e{0}%
\eol}\vss}\rg%
%
%
\rx{\vss\hfull{%
\rlx{\hss{$1400_{z}$}}\cg%
\e{0}%
\e{0}%
\e{0}%
\e{2}%
\e{2}%
\e{0}%
\e{2}%
\e{1}%
\e{0}%
\e{1}%
\e{0}%
\e{0}%
\e{0}%
\e{0}%
\e{0}%
\e{0}%
\e{0}%
\e{0}%
\eol}\vss}\rg%
%
%
\rx{\vss\hfull{%
\rlx{\hss{$1008_{z}$}}\cg%
\e{0}%
\e{0}%
\e{0}%
\e{1}%
\e{2}%
\e{1}%
\e{1}%
\e{1}%
\e{0}%
\e{0}%
\e{0}%
\e{0}%
\e{0}%
\e{0}%
\e{0}%
\e{0}%
\e{0}%
\e{0}%
\eol}\vss}\rg%
%
%
\rx{\vss\hfull{%
\rlx{\hss{$560_{z}$}}\cg%
\e{0}%
\e{0}%
\e{0}%
\e{2}%
\e{1}%
\e{0}%
\e{1}%
\e{0}%
\e{0}%
\e{0}%
\e{0}%
\e{0}%
\e{0}%
\e{0}%
\e{0}%
\e{0}%
\e{0}%
\e{0}%
\eol}\vss}\rg%
%
%
\rx{\vss\hfull{%
\rlx{\hss{$1400_{zz}$}}\cg%
\e{0}%
\e{0}%
\e{0}%
\e{2}%
\e{0}%
\e{0}%
\e{2}%
\e{1}%
\e{0}%
\e{1}%
\e{1}%
\e{0}%
\e{0}%
\e{0}%
\e{0}%
\e{0}%
\e{0}%
\e{0}%
\eol}\vss}\rg%
%
%
\rx{\vss\hfull{%
\rlx{\hss{$4200_{z}$}}\cg%
\e{0}%
\e{0}%
\e{0}%
\e{1}%
\e{1}%
\e{0}%
\e{2}%
\e{4}%
\e{1}%
\e{0}%
\e{3}%
\e{1}%
\e{1}%
\e{2}%
\e{0}%
\e{0}%
\e{0}%
\e{0}%
\eol}\vss}\rg%
%
%
\rx{\vss\hfull{%
\rlx{\hss{$400_{z}$}}\cg%
\e{0}%
\e{0}%
\e{0}%
\e{2}%
\e{0}%
\e{0}%
\e{1}%
\e{0}%
\e{0}%
\e{0}%
\e{0}%
\e{0}%
\e{0}%
\e{0}%
\e{0}%
\e{0}%
\e{0}%
\e{0}%
\eol}\vss}\rg%
%
%
\rx{\vss\hfull{%
\rlx{\hss{$3240_{z}$}}\cg%
\e{0}%
\e{0}%
\e{0}%
\e{3}%
\e{3}%
\e{0}%
\e{4}%
\e{3}%
\e{0}%
\e{1}%
\e{1}%
\e{0}%
\e{1}%
\e{0}%
\e{0}%
\e{0}%
\e{0}%
\e{0}%
\eol}\vss}\rg%
%
%
\rx{\vss\hfull{%
\rlx{\hss{$4536_{z}$}}\cg%
\e{0}%
\e{0}%
\e{0}%
\e{1}%
\e{1}%
\e{0}%
\e{3}%
\e{4}%
\e{0}%
\e{3}%
\e{3}%
\e{0}%
\e{1}%
\e{2}%
\e{0}%
\e{0}%
\e{0}%
\e{0}%
\eol}\vss}\rg%
%
%
\rx{\vss\hfull{%
\rlx{\hss{$2400_{z}$}}\cg%
\e{0}%
\e{0}%
\e{0}%
\e{0}%
\e{1}%
\e{2}%
\e{1}%
\e{2}%
\e{2}%
\e{0}%
\e{0}%
\e{1}%
\e{2}%
\e{1}%
\e{0}%
\e{1}%
\e{0}%
\e{0}%
\eol}\vss}\rg%
%
%
\rx{\vss\hfull{%
\rlx{\hss{$3360_{z}$}}\cg%
\e{0}%
\e{0}%
\e{0}%
\e{1}%
\e{2}%
\e{0}%
\e{3}%
\e{3}%
\e{1}%
\e{1}%
\e{1}%
\e{1}%
\e{1}%
\e{1}%
\e{0}%
\e{0}%
\e{0}%
\e{0}%
\eol}\vss}\rg%
%
%
\rx{\vss\hfull{%
\rlx{\hss{$2800_{z}$}}\cg%
\e{0}%
\e{0}%
\e{0}%
\e{1}%
\e{1}%
\e{2}%
\e{2}%
\e{3}%
\e{1}%
\e{0}%
\e{1}%
\e{0}%
\e{1}%
\e{1}%
\e{0}%
\e{0}%
\e{0}%
\e{0}%
\eol}\vss}\rg%
%
%
\rx{\vss\hfull{%
\rlx{\hss{$4096_{z}$}}\cg%
\e{0}%
\e{0}%
\e{0}%
\e{1}%
\e{3}%
\e{1}%
\e{3}%
\e{4}%
\e{1}%
\e{2}%
\e{1}%
\e{0}%
\e{2}%
\e{1}%
\e{0}%
\e{0}%
\e{0}%
\e{0}%
\eol}\vss}\rg%
%
%
\rx{\vss\hfull{%
\rlx{\hss{$5600_{z}$}}\cg%
\e{0}%
\e{0}%
\e{0}%
\e{0}%
\e{2}%
\e{1}%
\e{2}%
\e{5}%
\e{2}%
\e{2}%
\e{2}%
\e{0}%
\e{3}%
\e{3}%
\e{1}%
\e{1}%
\e{0}%
\e{0}%
\eol}\vss}\rg%
%
%
\rx{\vss\hfull{%
\rlx{\hss{$448_{z}$}}\cg%
\e{0}%
\e{0}%
\e{0}%
\e{1}%
\e{0}%
\e{0}%
\e{1}%
\e{0}%
\e{0}%
\e{1}%
\e{0}%
\e{0}%
\e{0}%
\e{0}%
\e{0}%
\e{0}%
\e{0}%
\e{0}%
\eol}\vss}\rg%
%
%
\rx{\vss\hfull{%
\rlx{\hss{$448_{w}$}}\cg%
\e{0}%
\e{0}%
\e{0}%
\e{0}%
\e{0}%
\e{1}%
\e{0}%
\e{0}%
\e{1}%
\e{0}%
\e{0}%
\e{0}%
\e{1}%
\e{0}%
\e{0}%
\e{1}%
\e{0}%
\e{0}%
\eol}\vss}\rg%
%
%
\rx{\vss\hfull{%
\rlx{\hss{$1344_{w}$}}\cg%
\e{0}%
\e{0}%
\e{0}%
\e{0}%
\e{0}%
\e{0}%
\e{0}%
\e{1}%
\e{0}%
\e{0}%
\e{2}%
\e{0}%
\e{0}%
\e{1}%
\e{0}%
\e{0}%
\e{0}%
\e{0}%
\eol}\vss}\rg%
%
%
\rx{\vss\hfull{%
\rlx{\hss{$5600_{w}$}}\cg%
\e{1}%
\e{0}%
\e{0}%
\e{0}%
\e{1}%
\e{1}%
\e{1}%
\e{4}%
\e{3}%
\e{1}%
\e{2}%
\e{1}%
\e{3}%
\e{4}%
\e{1}%
\e{1}%
\e{1}%
\e{0}%
\eol}\vss}\rg%
%
%
\rx{\vss\hfull{%
\rlx{\hss{$2016_{w}$}}\cg%
\e{0}%
\e{0}%
\e{0}%
\e{0}%
\e{0}%
\e{0}%
\e{1}%
\e{1}%
\e{0}%
\e{1}%
\e{2}%
\e{1}%
\e{0}%
\e{1}%
\e{1}%
\e{0}%
\e{0}%
\e{0}%
\eol}\vss}\rg%
%
%
\rx{\vss\hfull{%
\rlx{\hss{$7168_{w}$}}\cg%
\e{1}%
\e{0}%
\e{0}%
\e{0}%
\e{1}%
\e{0}%
\e{2}%
\e{5}%
\e{2}%
\e{2}%
\e{4}%
\e{2}%
\e{2}%
\e{5}%
\e{2}%
\e{0}%
\e{1}%
\e{0}%
\eol}\vss}\rg%
\tableclose%
%
%
%
%
%
%
\eop
\eject
\tableopen{Induce/restrict matrix for $W(A_{7})\,\subset\,W(E_{8})$}%
%
%
%
%
%
%
\rowpts=18 true pt%
\colpts=18 true pt%
\rowlabpts=40 true pt%
\collabpts=50 true pt%
\clx{\vss\hfull{%
\rlx{\hss{$ $}}\cg%
\cx{\hskip 16 true pt\flip{$[{8}]$}\hss}\cg%
\cx{\hskip 16 true pt\flip{$[{7}{1}]$}\hss}\cg%
\cx{\hskip 16 true pt\flip{$[{6}{2}]$}\hss}\cg%
\cx{\hskip 16 true pt\flip{$[{6}{1^{2}}]$}\hss}\cg%
\cx{\hskip 16 true pt\flip{$[{5}{3}]$}\hss}\cg%
\cx{\hskip 16 true pt\flip{$[{5}{2}{1}]$}\hss}\cg%
\cx{\hskip 16 true pt\flip{$[{5}{1^{3}}]$}\hss}\cg%
\cx{\hskip 16 true pt\flip{$[{4^{2}}]$}\hss}\cg%
\cx{\hskip 16 true pt\flip{$[{4}{3}{1}]$}\hss}\cg%
\cx{\hskip 16 true pt\flip{$[{4}{2^{2}}]$}\hss}\cg%
\cx{\hskip 16 true pt\flip{$[{4}{2}{1^{2}}]$}\hss}\cg%
\cx{\hskip 16 true pt\flip{$[{4}{1^{4}}]$}\hss}\cg%
\cx{\hskip 16 true pt\flip{$[{3^{2}}{2}]$}\hss}\cg%
\cx{\hskip 16 true pt\flip{$[{3^{2}}{1^{2}}]$}\hss}\cg%
\cx{\hskip 16 true pt\flip{$[{3}{2^{2}}{1}]$}\hss}\cg%
\cx{\hskip 16 true pt\flip{$[{3}{2}{1^{3}}]$}\hss}\cg%
\cx{\hskip 16 true pt\flip{$[{3}{1^{5}}]$}\hss}\cg%
\cx{\hskip 16 true pt\flip{$[{2^{4}}]$}\hss}\cg%
\cx{\hskip 16 true pt\flip{$[{2^{3}}{1^{2}}]$}\hss}\cg%
\cx{\hskip 16 true pt\flip{$[{2^{2}}{1^{4}}]$}\hss}\cg%
\cx{\hskip 16 true pt\flip{$[{2}{1^{6}}]$}\hss}\cg%
\cx{\hskip 16 true pt\flip{$[{1^{8}}]$}\hss}\cg%
\eol}}\rg%
%
%
\rx{\vss\hfull{%
\rlx{\hss{$1_{x}$}}\cg%
\e{1}%
\e{0}%
\e{0}%
\e{0}%
\e{0}%
\e{0}%
\e{0}%
\e{0}%
\e{0}%
\e{0}%
\e{0}%
\e{0}%
\e{0}%
\e{0}%
\e{0}%
\e{0}%
\e{0}%
\e{0}%
\e{0}%
\e{0}%
\e{0}%
\e{0}%
\eol}\vss}\rg%
%
%
\rx{\vss\hfull{%
\rlx{\hss{$28_{x}$}}\cg%
\e{0}%
\e{1}%
\e{0}%
\e{1}%
\e{0}%
\e{0}%
\e{0}%
\e{0}%
\e{0}%
\e{0}%
\e{0}%
\e{0}%
\e{0}%
\e{0}%
\e{0}%
\e{0}%
\e{0}%
\e{0}%
\e{0}%
\e{0}%
\e{0}%
\e{0}%
\eol}\vss}\rg%
%
%
\rx{\vss\hfull{%
\rlx{\hss{$35_{x}$}}\cg%
\e{1}%
\e{2}%
\e{1}%
\e{0}%
\e{0}%
\e{0}%
\e{0}%
\e{0}%
\e{0}%
\e{0}%
\e{0}%
\e{0}%
\e{0}%
\e{0}%
\e{0}%
\e{0}%
\e{0}%
\e{0}%
\e{0}%
\e{0}%
\e{0}%
\e{0}%
\eol}\vss}\rg%
%
%
\rx{\vss\hfull{%
\rlx{\hss{$84_{x}$}}\cg%
\e{2}%
\e{2}%
\e{2}%
\e{0}%
\e{1}%
\e{0}%
\e{0}%
\e{0}%
\e{0}%
\e{0}%
\e{0}%
\e{0}%
\e{0}%
\e{0}%
\e{0}%
\e{0}%
\e{0}%
\e{0}%
\e{0}%
\e{0}%
\e{0}%
\e{0}%
\eol}\vss}\rg%
%
%
\rx{\vss\hfull{%
\rlx{\hss{$50_{x}$}}\cg%
\e{1}%
\e{1}%
\e{0}%
\e{0}%
\e{1}%
\e{0}%
\e{0}%
\e{1}%
\e{0}%
\e{0}%
\e{0}%
\e{0}%
\e{0}%
\e{0}%
\e{0}%
\e{0}%
\e{0}%
\e{0}%
\e{0}%
\e{0}%
\e{0}%
\e{0}%
\eol}\vss}\rg%
%
%
\rx{\vss\hfull{%
\rlx{\hss{$350_{x}$}}\cg%
\e{0}%
\e{0}%
\e{1}%
\e{2}%
\e{0}%
\e{2}%
\e{2}%
\e{0}%
\e{0}%
\e{0}%
\e{1}%
\e{0}%
\e{0}%
\e{0}%
\e{0}%
\e{0}%
\e{0}%
\e{0}%
\e{0}%
\e{0}%
\e{0}%
\e{0}%
\eol}\vss}\rg%
%
%
\rx{\vss\hfull{%
\rlx{\hss{$300_{x}$}}\cg%
\e{0}%
\e{1}%
\e{3}%
\e{1}%
\e{1}%
\e{2}%
\e{0}%
\e{0}%
\e{0}%
\e{1}%
\e{0}%
\e{0}%
\e{0}%
\e{0}%
\e{0}%
\e{0}%
\e{0}%
\e{0}%
\e{0}%
\e{0}%
\e{0}%
\e{0}%
\eol}\vss}\rg%
%
%
\rx{\vss\hfull{%
\rlx{\hss{$567_{x}$}}\cg%
\e{1}%
\e{4}%
\e{4}%
\e{5}%
\e{2}%
\e{3}%
\e{1}%
\e{0}%
\e{1}%
\e{0}%
\e{0}%
\e{0}%
\e{0}%
\e{0}%
\e{0}%
\e{0}%
\e{0}%
\e{0}%
\e{0}%
\e{0}%
\e{0}%
\e{0}%
\eol}\vss}\rg%
%
%
\rx{\vss\hfull{%
\rlx{\hss{$210_{x}$}}\cg%
\e{1}%
\e{3}%
\e{2}%
\e{2}%
\e{1}%
\e{1}%
\e{0}%
\e{1}%
\e{0}%
\e{0}%
\e{0}%
\e{0}%
\e{0}%
\e{0}%
\e{0}%
\e{0}%
\e{0}%
\e{0}%
\e{0}%
\e{0}%
\e{0}%
\e{0}%
\eol}\vss}\rg%
%
%
\rx{\vss\hfull{%
\rlx{\hss{$840_{x}$}}\cg%
\e{0}%
\e{0}%
\e{1}%
\e{0}%
\e{3}%
\e{2}%
\e{0}%
\e{1}%
\e{2}%
\e{3}%
\e{1}%
\e{0}%
\e{1}%
\e{0}%
\e{2}%
\e{0}%
\e{0}%
\e{1}%
\e{0}%
\e{0}%
\e{0}%
\e{0}%
\eol}\vss}\rg%
%
%
\rx{\vss\hfull{%
\rlx{\hss{$700_{x}$}}\cg%
\e{2}%
\e{4}%
\e{5}%
\e{2}%
\e{4}%
\e{3}%
\e{0}%
\e{2}%
\e{2}%
\e{1}%
\e{0}%
\e{0}%
\e{0}%
\e{0}%
\e{0}%
\e{0}%
\e{0}%
\e{0}%
\e{0}%
\e{0}%
\e{0}%
\e{0}%
\eol}\vss}\rg%
%
%
\rx{\vss\hfull{%
\rlx{\hss{$175_{x}$}}\cg%
\e{1}%
\e{0}%
\e{1}%
\e{0}%
\e{1}%
\e{0}%
\e{0}%
\e{1}%
\e{1}%
\e{0}%
\e{0}%
\e{0}%
\e{1}%
\e{0}%
\e{0}%
\e{0}%
\e{0}%
\e{0}%
\e{0}%
\e{0}%
\e{0}%
\e{0}%
\eol}\vss}\rg%
%
%
\rx{\vss\hfull{%
\rlx{\hss{$1400_{x}$}}\cg%
\e{1}%
\e{3}%
\e{4}%
\e{4}%
\e{4}%
\e{5}%
\e{2}%
\e{2}%
\e{5}%
\e{1}%
\e{2}%
\e{0}%
\e{1}%
\e{1}%
\e{0}%
\e{0}%
\e{0}%
\e{0}%
\e{0}%
\e{0}%
\e{0}%
\e{0}%
\eol}\vss}\rg%
%
%
\rx{\vss\hfull{%
\rlx{\hss{$1050_{x}$}}\cg%
\e{1}%
\e{3}%
\e{3}%
\e{2}%
\e{4}%
\e{3}%
\e{0}%
\e{2}%
\e{5}%
\e{1}%
\e{1}%
\e{0}%
\e{1}%
\e{1}%
\e{0}%
\e{0}%
\e{0}%
\e{0}%
\e{0}%
\e{0}%
\e{0}%
\e{0}%
\eol}\vss}\rg%
%
%
\rx{\vss\hfull{%
\rlx{\hss{$1575_{x}$}}\cg%
\e{0}%
\e{3}%
\e{4}%
\e{6}%
\e{3}%
\e{7}%
\e{4}%
\e{1}%
\e{4}%
\e{1}%
\e{3}%
\e{0}%
\e{0}%
\e{1}%
\e{0}%
\e{0}%
\e{0}%
\e{0}%
\e{0}%
\e{0}%
\e{0}%
\e{0}%
\eol}\vss}\rg%
%
%
\rx{\vss\hfull{%
\rlx{\hss{$1344_{x}$}}\cg%
\e{1}%
\e{4}%
\e{7}%
\e{4}%
\e{5}%
\e{7}%
\e{1}%
\e{1}%
\e{3}%
\e{2}%
\e{1}%
\e{0}%
\e{1}%
\e{0}%
\e{0}%
\e{0}%
\e{0}%
\e{0}%
\e{0}%
\e{0}%
\e{0}%
\e{0}%
\eol}\vss}\rg%
%
%
\rx{\vss\hfull{%
\rlx{\hss{$2100_{x}$}}\cg%
\e{0}%
\e{0}%
\e{2}%
\e{4}%
\e{1}%
\e{7}%
\e{7}%
\e{0}%
\e{3}%
\e{3}%
\e{6}%
\e{3}%
\e{1}%
\e{1}%
\e{1}%
\e{1}%
\e{0}%
\e{0}%
\e{0}%
\e{0}%
\e{0}%
\e{0}%
\eol}\vss}\rg%
%
%
\rx{\vss\hfull{%
\rlx{\hss{$2268_{x}$}}\cg%
\e{0}%
\e{2}%
\e{5}%
\e{6}%
\e{4}%
\e{9}%
\e{5}%
\e{2}%
\e{5}%
\e{4}%
\e{4}%
\e{1}%
\e{1}%
\e{1}%
\e{1}%
\e{0}%
\e{0}%
\e{0}%
\e{0}%
\e{0}%
\e{0}%
\e{0}%
\eol}\vss}\rg%
%
%
\rx{\vss\hfull{%
\rlx{\hss{$525_{x}$}}\cg%
\e{0}%
\e{1}%
\e{1}%
\e{2}%
\e{2}%
\e{2}%
\e{2}%
\e{0}%
\e{1}%
\e{0}%
\e{1}%
\e{0}%
\e{1}%
\e{0}%
\e{0}%
\e{0}%
\e{0}%
\e{0}%
\e{0}%
\e{0}%
\e{0}%
\e{0}%
\eol}\vss}\rg%
%
%
\rx{\vss\hfull{%
\rlx{\hss{$700_{xx}$}}\cg%
\e{0}%
\e{1}%
\e{0}%
\e{1}%
\e{2}%
\e{1}%
\e{0}%
\e{2}%
\e{3}%
\e{0}%
\e{1}%
\e{0}%
\e{1}%
\e{2}%
\e{1}%
\e{0}%
\e{0}%
\e{0}%
\e{0}%
\e{0}%
\e{0}%
\e{0}%
\eol}\vss}\rg%
%
%
\rx{\vss\hfull{%
\rlx{\hss{$972_{x}$}}\cg%
\e{0}%
\e{1}%
\e{3}%
\e{1}%
\e{3}%
\e{3}%
\e{0}%
\e{2}%
\e{3}%
\e{3}%
\e{1}%
\e{0}%
\e{1}%
\e{0}%
\e{1}%
\e{0}%
\e{0}%
\e{0}%
\e{0}%
\e{0}%
\e{0}%
\e{0}%
\eol}\vss}\rg%
%
%
\rx{\vss\hfull{%
\rlx{\hss{$4096_{x}$}}\cg%
\e{0}%
\e{2}%
\e{6}%
\e{6}%
\e{6}%
\e{15}%
\e{6}%
\e{2}%
\e{9}%
\e{7}%
\e{9}%
\e{2}%
\e{3}%
\e{3}%
\e{3}%
\e{1}%
\e{0}%
\e{0}%
\e{0}%
\e{0}%
\e{0}%
\e{0}%
\eol}\vss}\rg%
%
%
\rx{\vss\hfull{%
\rlx{\hss{$4200_{x}$}}\cg%
\e{0}%
\e{2}%
\e{5}%
\e{4}%
\e{7}%
\e{11}%
\e{3}%
\e{3}%
\e{12}%
\e{8}%
\e{8}%
\e{1}%
\e{5}%
\e{5}%
\e{4}%
\e{2}%
\e{0}%
\e{1}%
\e{0}%
\e{0}%
\e{0}%
\e{0}%
\eol}\vss}\rg%
%
%
\rx{\vss\hfull{%
\rlx{\hss{$2240_{x}$}}\cg%
\e{1}%
\e{2}%
\e{5}%
\e{2}%
\e{7}%
\e{7}%
\e{1}%
\e{3}%
\e{7}%
\e{4}%
\e{3}%
\e{0}%
\e{3}%
\e{2}%
\e{2}%
\e{0}%
\e{0}%
\e{0}%
\e{0}%
\e{0}%
\e{0}%
\e{0}%
\eol}\vss}\rg%
%
%
\rx{\vss\hfull{%
\rlx{\hss{$2835_{x}$}}\cg%
\e{0}%
\e{1}%
\e{3}%
\e{1}%
\e{4}%
\e{5}%
\e{1}%
\e{3}%
\e{9}%
\e{5}%
\e{4}%
\e{0}%
\e{6}%
\e{5}%
\e{4}%
\e{2}%
\e{0}%
\e{0}%
\e{1}%
\e{0}%
\e{0}%
\e{0}%
\eol}\vss}\rg%
%
%
\rx{\vss\hfull{%
\rlx{\hss{$6075_{x}$}}\cg%
\e{0}%
\e{1}%
\e{4}%
\e{6}%
\e{6}%
\e{14}%
\e{8}%
\e{3}%
\e{14}%
\e{8}%
\e{16}%
\e{4}%
\e{6}%
\e{8}%
\e{6}%
\e{5}%
\e{0}%
\e{0}%
\e{1}%
\e{0}%
\e{0}%
\e{0}%
\eol}\vss}\rg%
%
%
\rx{\vss\hfull{%
\rlx{\hss{$3200_{x}$}}\cg%
\e{0}%
\e{0}%
\e{2}%
\e{1}%
\e{3}%
\e{8}%
\e{2}%
\e{1}%
\e{6}%
\e{8}%
\e{7}%
\e{2}%
\e{4}%
\e{2}%
\e{6}%
\e{2}%
\e{1}%
\e{1}%
\e{1}%
\e{0}%
\e{0}%
\e{0}%
\eol}\vss}\rg%
%
%
\rx{\vss\hfull{%
\rlx{\hss{$70_{y}$}}\cg%
\e{0}%
\e{0}%
\e{0}%
\e{0}%
\e{0}%
\e{0}%
\e{1}%
\e{0}%
\e{0}%
\e{0}%
\e{0}%
\e{1}%
\e{0}%
\e{0}%
\e{0}%
\e{0}%
\e{0}%
\e{0}%
\e{0}%
\e{0}%
\e{0}%
\e{0}%
\eol}\vss}\rg%
%
%
\rx{\vss\hfull{%
\rlx{\hss{$1134_{y}$}}\cg%
\e{0}%
\e{0}%
\e{0}%
\e{0}%
\e{1}%
\e{2}%
\e{1}%
\e{0}%
\e{2}%
\e{1}%
\e{4}%
\e{1}%
\e{0}%
\e{1}%
\e{2}%
\e{2}%
\e{0}%
\e{0}%
\e{1}%
\e{0}%
\e{0}%
\e{0}%
\eol}\vss}\rg%
%
%
\rx{\vss\hfull{%
\rlx{\hss{$1680_{y}$}}\cg%
\e{0}%
\e{0}%
\e{0}%
\e{1}%
\e{0}%
\e{3}%
\e{5}%
\e{0}%
\e{1}%
\e{2}%
\e{6}%
\e{5}%
\e{0}%
\e{2}%
\e{1}%
\e{3}%
\e{1}%
\e{0}%
\e{0}%
\e{0}%
\e{0}%
\e{0}%
\eol}\vss}\rg%
%
%
\rx{\vss\hfull{%
\rlx{\hss{$168_{y}$}}\cg%
\e{0}%
\e{0}%
\e{0}%
\e{0}%
\e{0}%
\e{0}%
\e{0}%
\e{1}%
\e{1}%
\e{0}%
\e{0}%
\e{0}%
\e{0}%
\e{0}%
\e{1}%
\e{0}%
\e{0}%
\e{1}%
\e{0}%
\e{0}%
\e{0}%
\e{0}%
\eol}\vss}\rg%
%
%
\rx{\vss\hfull{%
\rlx{\hss{$420_{y}$}}\cg%
\e{0}%
\e{0}%
\e{0}%
\e{0}%
\e{1}%
\e{0}%
\e{0}%
\e{1}%
\e{1}%
\e{1}%
\e{0}%
\e{0}%
\e{2}%
\e{1}%
\e{1}%
\e{0}%
\e{0}%
\e{1}%
\e{1}%
\e{0}%
\e{0}%
\e{0}%
\eol}\vss}\rg%
%
%
\rx{\vss\hfull{%
\rlx{\hss{$3150_{y}$}}\cg%
\e{0}%
\e{0}%
\e{1}%
\e{0}%
\e{3}%
\e{4}%
\e{1}%
\e{1}%
\e{7}%
\e{5}%
\e{6}%
\e{1}%
\e{6}%
\e{5}%
\e{7}%
\e{4}%
\e{0}%
\e{1}%
\e{3}%
\e{1}%
\e{0}%
\e{0}%
\eol}\vss}\rg%
%
%
\rx{\vss\hfull{%
\rlx{\hss{$4200_{y}$}}\cg%
\e{0}%
\e{0}%
\e{1}%
\e{1}%
\e{3}%
\e{5}%
\e{1}%
\e{3}%
\e{10}%
\e{7}%
\e{8}%
\e{1}%
\e{6}%
\e{7}%
\e{10}%
\e{5}%
\e{1}%
\e{3}%
\e{3}%
\e{1}%
\e{0}%
\e{0}%
\eol}\vss}\rg%
\eop
\eject
\tablecont%
%
%
%
%
%
%
\rowpts=18 true pt%
\colpts=18 true pt%
\rowlabpts=40 true pt%
\collabpts=50 true pt%
\clx{\vss\hfull{%
\rlx{\hss{$ $}}\cg%
\cx{\hskip 16 true pt\flip{$[{8}]$}\hss}\cg%
\cx{\hskip 16 true pt\flip{$[{7}{1}]$}\hss}\cg%
\cx{\hskip 16 true pt\flip{$[{6}{2}]$}\hss}\cg%
\cx{\hskip 16 true pt\flip{$[{6}{1^{2}}]$}\hss}\cg%
\cx{\hskip 16 true pt\flip{$[{5}{3}]$}\hss}\cg%
\cx{\hskip 16 true pt\flip{$[{5}{2}{1}]$}\hss}\cg%
\cx{\hskip 16 true pt\flip{$[{5}{1^{3}}]$}\hss}\cg%
\cx{\hskip 16 true pt\flip{$[{4^{2}}]$}\hss}\cg%
\cx{\hskip 16 true pt\flip{$[{4}{3}{1}]$}\hss}\cg%
\cx{\hskip 16 true pt\flip{$[{4}{2^{2}}]$}\hss}\cg%
\cx{\hskip 16 true pt\flip{$[{4}{2}{1^{2}}]$}\hss}\cg%
\cx{\hskip 16 true pt\flip{$[{4}{1^{4}}]$}\hss}\cg%
\cx{\hskip 16 true pt\flip{$[{3^{2}}{2}]$}\hss}\cg%
\cx{\hskip 16 true pt\flip{$[{3^{2}}{1^{2}}]$}\hss}\cg%
\cx{\hskip 16 true pt\flip{$[{3}{2^{2}}{1}]$}\hss}\cg%
\cx{\hskip 16 true pt\flip{$[{3}{2}{1^{3}}]$}\hss}\cg%
\cx{\hskip 16 true pt\flip{$[{3}{1^{5}}]$}\hss}\cg%
\cx{\hskip 16 true pt\flip{$[{2^{4}}]$}\hss}\cg%
\cx{\hskip 16 true pt\flip{$[{2^{3}}{1^{2}}]$}\hss}\cg%
\cx{\hskip 16 true pt\flip{$[{2^{2}}{1^{4}}]$}\hss}\cg%
\cx{\hskip 16 true pt\flip{$[{2}{1^{6}}]$}\hss}\cg%
\cx{\hskip 16 true pt\flip{$[{1^{8}}]$}\hss}\cg%
\eol}}\rg%
%
%
\rx{\vss\hfull{%
\rlx{\hss{$2688_{y}$}}\cg%
\e{0}%
\e{0}%
\e{1}%
\e{1}%
\e{1}%
\e{4}%
\e{1}%
\e{0}%
\e{5}%
\e{5}%
\e{6}%
\e{1}%
\e{4}%
\e{5}%
\e{5}%
\e{4}%
\e{1}%
\e{0}%
\e{1}%
\e{1}%
\e{0}%
\e{0}%
\eol}\vss}\rg%
%
%
\rx{\vss\hfull{%
\rlx{\hss{$2100_{y}$}}\cg%
\e{0}%
\e{0}%
\e{0}%
\e{1}%
\e{1}%
\e{3}%
\e{5}%
\e{0}%
\e{3}%
\e{2}%
\e{6}%
\e{5}%
\e{2}%
\e{2}%
\e{3}%
\e{3}%
\e{1}%
\e{0}%
\e{1}%
\e{0}%
\e{0}%
\e{0}%
\eol}\vss}\rg%
%
%
\rx{\vss\hfull{%
\rlx{\hss{$1400_{y}$}}\cg%
\e{0}%
\e{0}%
\e{0}%
\e{1}%
\e{0}%
\e{2}%
\e{3}%
\e{1}%
\e{2}%
\e{2}%
\e{4}%
\e{3}%
\e{0}%
\e{2}%
\e{2}%
\e{2}%
\e{1}%
\e{1}%
\e{0}%
\e{0}%
\e{0}%
\e{0}%
\eol}\vss}\rg%
%
%
\rx{\vss\hfull{%
\rlx{\hss{$4536_{y}$}}\cg%
\e{0}%
\e{0}%
\e{1}%
\e{2}%
\e{1}%
\e{7}%
\e{6}%
\e{1}%
\e{7}%
\e{7}%
\e{12}%
\e{6}%
\e{4}%
\e{7}%
\e{7}%
\e{7}%
\e{2}%
\e{1}%
\e{1}%
\e{1}%
\e{0}%
\e{0}%
\eol}\vss}\rg%
%
%
\rx{\vss\hfull{%
\rlx{\hss{$5670_{y}$}}\cg%
\e{0}%
\e{0}%
\e{1}%
\e{2}%
\e{2}%
\e{9}%
\e{7}%
\e{1}%
\e{9}%
\e{8}%
\e{16}%
\e{7}%
\e{4}%
\e{8}%
\e{9}%
\e{9}%
\e{2}%
\e{1}%
\e{2}%
\e{1}%
\e{0}%
\e{0}%
\eol}\vss}\rg%
%
%
\rx{\vss\hfull{%
\rlx{\hss{$4480_{y}$}}\cg%
\e{0}%
\e{0}%
\e{1}%
\e{1}%
\e{3}%
\e{6}%
\e{3}%
\e{2}%
\e{9}%
\e{7}%
\e{10}%
\e{3}%
\e{6}%
\e{7}%
\e{9}%
\e{6}%
\e{1}%
\e{2}%
\e{3}%
\e{1}%
\e{0}%
\e{0}%
\eol}\vss}\rg%
%
%
\rx{\vss\hfull{%
\rlx{\hss{$8_{z}$}}\cg%
\e{1}%
\e{1}%
\e{0}%
\e{0}%
\e{0}%
\e{0}%
\e{0}%
\e{0}%
\e{0}%
\e{0}%
\e{0}%
\e{0}%
\e{0}%
\e{0}%
\e{0}%
\e{0}%
\e{0}%
\e{0}%
\e{0}%
\e{0}%
\e{0}%
\e{0}%
\eol}\vss}\rg%
%
%
\rx{\vss\hfull{%
\rlx{\hss{$56_{z}$}}\cg%
\e{0}%
\e{0}%
\e{0}%
\e{1}%
\e{0}%
\e{0}%
\e{1}%
\e{0}%
\e{0}%
\e{0}%
\e{0}%
\e{0}%
\e{0}%
\e{0}%
\e{0}%
\e{0}%
\e{0}%
\e{0}%
\e{0}%
\e{0}%
\e{0}%
\e{0}%
\eol}\vss}\rg%
%
%
\rx{\vss\hfull{%
\rlx{\hss{$160_{z}$}}\cg%
\e{0}%
\e{2}%
\e{2}%
\e{2}%
\e{0}%
\e{1}%
\e{0}%
\e{0}%
\e{0}%
\e{0}%
\e{0}%
\e{0}%
\e{0}%
\e{0}%
\e{0}%
\e{0}%
\e{0}%
\e{0}%
\e{0}%
\e{0}%
\e{0}%
\e{0}%
\eol}\vss}\rg%
%
%
\rx{\vss\hfull{%
\rlx{\hss{$112_{z}$}}\cg%
\e{2}%
\e{3}%
\e{2}%
\e{1}%
\e{1}%
\e{0}%
\e{0}%
\e{0}%
\e{0}%
\e{0}%
\e{0}%
\e{0}%
\e{0}%
\e{0}%
\e{0}%
\e{0}%
\e{0}%
\e{0}%
\e{0}%
\e{0}%
\e{0}%
\e{0}%
\eol}\vss}\rg%
%
%
\rx{\vss\hfull{%
\rlx{\hss{$840_{z}$}}\cg%
\e{0}%
\e{0}%
\e{2}%
\e{1}%
\e{2}%
\e{4}%
\e{1}%
\e{0}%
\e{1}%
\e{2}%
\e{2}%
\e{0}%
\e{0}%
\e{0}%
\e{1}%
\e{0}%
\e{0}%
\e{0}%
\e{0}%
\e{0}%
\e{0}%
\e{0}%
\eol}\vss}\rg%
%
%
\rx{\vss\hfull{%
\rlx{\hss{$1296_{z}$}}\cg%
\e{0}%
\e{1}%
\e{2}%
\e{5}%
\e{1}%
\e{6}%
\e{5}%
\e{0}%
\e{2}%
\e{1}%
\e{3}%
\e{1}%
\e{0}%
\e{1}%
\e{0}%
\e{0}%
\e{0}%
\e{0}%
\e{0}%
\e{0}%
\e{0}%
\e{0}%
\eol}\vss}\rg%
%
%
\rx{\vss\hfull{%
\rlx{\hss{$1400_{z}$}}\cg%
\e{1}%
\e{4}%
\e{7}%
\e{5}%
\e{5}%
\e{7}%
\e{2}%
\e{1}%
\e{3}%
\e{2}%
\e{1}%
\e{0}%
\e{1}%
\e{0}%
\e{0}%
\e{0}%
\e{0}%
\e{0}%
\e{0}%
\e{0}%
\e{0}%
\e{0}%
\eol}\vss}\rg%
%
%
\rx{\vss\hfull{%
\rlx{\hss{$1008_{z}$}}\cg%
\e{0}%
\e{3}%
\e{4}%
\e{6}%
\e{2}%
\e{5}%
\e{3}%
\e{1}%
\e{2}%
\e{1}%
\e{1}%
\e{0}%
\e{0}%
\e{0}%
\e{0}%
\e{0}%
\e{0}%
\e{0}%
\e{0}%
\e{0}%
\e{0}%
\e{0}%
\eol}\vss}\rg%
%
%
\rx{\vss\hfull{%
\rlx{\hss{$560_{z}$}}\cg%
\e{2}%
\e{5}%
\e{5}%
\e{3}%
\e{3}%
\e{3}%
\e{0}%
\e{1}%
\e{1}%
\e{0}%
\e{0}%
\e{0}%
\e{0}%
\e{0}%
\e{0}%
\e{0}%
\e{0}%
\e{0}%
\e{0}%
\e{0}%
\e{0}%
\e{0}%
\eol}\vss}\rg%
%
%
\rx{\vss\hfull{%
\rlx{\hss{$1400_{zz}$}}\cg%
\e{1}%
\e{2}%
\e{3}%
\e{1}%
\e{5}%
\e{3}%
\e{0}%
\e{3}%
\e{6}%
\e{2}%
\e{1}%
\e{0}%
\e{3}%
\e{2}%
\e{1}%
\e{0}%
\e{0}%
\e{0}%
\e{0}%
\e{0}%
\e{0}%
\e{0}%
\eol}\vss}\rg%
%
%
\rx{\vss\hfull{%
\rlx{\hss{$4200_{z}$}}\cg%
\e{0}%
\e{1}%
\e{2}%
\e{3}%
\e{5}%
\e{7}%
\e{3}%
\e{4}%
\e{12}%
\e{5}%
\e{9}%
\e{1}%
\e{5}%
\e{8}%
\e{6}%
\e{4}%
\e{0}%
\e{1}%
\e{1}%
\e{0}%
\e{0}%
\e{0}%
\eol}\vss}\rg%
%
%
\rx{\vss\hfull{%
\rlx{\hss{$400_{z}$}}\cg%
\e{2}%
\e{3}%
\e{2}%
\e{1}%
\e{3}%
\e{1}%
\e{0}%
\e{2}%
\e{2}%
\e{0}%
\e{0}%
\e{0}%
\e{0}%
\e{0}%
\e{0}%
\e{0}%
\e{0}%
\e{0}%
\e{0}%
\e{0}%
\e{0}%
\e{0}%
\eol}\vss}\rg%
%
%
\rx{\vss\hfull{%
\rlx{\hss{$3240_{z}$}}\cg%
\e{1}%
\e{5}%
\e{9}%
\e{7}%
\e{9}%
\e{12}%
\e{3}%
\e{4}%
\e{10}%
\e{5}%
\e{5}%
\e{0}%
\e{2}%
\e{2}%
\e{1}%
\e{0}%
\e{0}%
\e{0}%
\e{0}%
\e{0}%
\e{0}%
\e{0}%
\eol}\vss}\rg%
%
%
\rx{\vss\hfull{%
\rlx{\hss{$4536_{z}$}}\cg%
\e{0}%
\e{1}%
\e{5}%
\e{2}%
\e{8}%
\e{11}%
\e{2}%
\e{4}%
\e{12}%
\e{10}%
\e{8}%
\e{1}%
\e{7}%
\e{4}%
\e{7}%
\e{2}%
\e{0}%
\e{1}%
\e{1}%
\e{0}%
\e{0}%
\e{0}%
\eol}\vss}\rg%
%
%
\rx{\vss\hfull{%
\rlx{\hss{$2400_{z}$}}\cg%
\e{0}%
\e{0}%
\e{1}%
\e{3}%
\e{1}%
\e{5}%
\e{7}%
\e{0}%
\e{4}%
\e{2}%
\e{8}%
\e{4}%
\e{1}%
\e{3}%
\e{1}%
\e{3}%
\e{0}%
\e{0}%
\e{0}%
\e{0}%
\e{0}%
\e{0}%
\eol}\vss}\rg%
%
%
\rx{\vss\hfull{%
\rlx{\hss{$3360_{z}$}}\cg%
\e{0}%
\e{2}%
\e{4}%
\e{4}%
\e{4}%
\e{9}%
\e{3}%
\e{2}%
\e{10}%
\e{5}%
\e{7}%
\e{1}%
\e{4}%
\e{5}%
\e{2}%
\e{2}%
\e{0}%
\e{0}%
\e{0}%
\e{0}%
\e{0}%
\e{0}%
\eol}\vss}\rg%
%
%
\rx{\vss\hfull{%
\rlx{\hss{$2800_{z}$}}\cg%
\e{0}%
\e{1}%
\e{2}%
\e{5}%
\e{3}%
\e{8}%
\e{7}%
\e{2}%
\e{6}%
\e{3}%
\e{7}%
\e{3}%
\e{2}%
\e{3}%
\e{2}%
\e{1}%
\e{0}%
\e{0}%
\e{0}%
\e{0}%
\e{0}%
\e{0}%
\eol}\vss}\rg%
%
%
\rx{\vss\hfull{%
\rlx{\hss{$4096_{z}$}}\cg%
\e{0}%
\e{2}%
\e{6}%
\e{6}%
\e{6}%
\e{15}%
\e{6}%
\e{2}%
\e{9}%
\e{7}%
\e{9}%
\e{2}%
\e{3}%
\e{3}%
\e{3}%
\e{1}%
\e{0}%
\e{0}%
\e{0}%
\e{0}%
\e{0}%
\e{0}%
\eol}\vss}\rg%
%
%
\rx{\vss\hfull{%
\rlx{\hss{$5600_{z}$}}\cg%
\e{0}%
\e{0}%
\e{3}%
\e{4}%
\e{4}%
\e{13}%
\e{9}%
\e{1}%
\e{10}%
\e{10}%
\e{15}%
\e{6}%
\e{5}%
\e{5}%
\e{7}%
\e{5}%
\e{1}%
\e{1}%
\e{1}%
\e{0}%
\e{0}%
\e{0}%
\eol}\vss}\rg%
%
%
\rx{\vss\hfull{%
\rlx{\hss{$448_{z}$}}\cg%
\e{1}%
\e{1}%
\e{3}%
\e{0}%
\e{3}%
\e{2}%
\e{0}%
\e{0}%
\e{1}%
\e{1}%
\e{0}%
\e{0}%
\e{1}%
\e{0}%
\e{0}%
\e{0}%
\e{0}%
\e{0}%
\e{0}%
\e{0}%
\e{0}%
\e{0}%
\eol}\vss}\rg%
%
%
\rx{\vss\hfull{%
\rlx{\hss{$448_{w}$}}\cg%
\e{0}%
\e{0}%
\e{0}%
\e{0}%
\e{0}%
\e{1}%
\e{2}%
\e{0}%
\e{0}%
\e{0}%
\e{2}%
\e{2}%
\e{0}%
\e{0}%
\e{0}%
\e{1}%
\e{0}%
\e{0}%
\e{0}%
\e{0}%
\e{0}%
\e{0}%
\eol}\vss}\rg%
%
%
\rx{\vss\hfull{%
\rlx{\hss{$1344_{w}$}}\cg%
\e{0}%
\e{0}%
\e{0}%
\e{0}%
\e{2}%
\e{1}%
\e{0}%
\e{2}%
\e{4}%
\e{2}%
\e{2}%
\e{0}%
\e{2}%
\e{2}%
\e{4}%
\e{1}%
\e{0}%
\e{2}%
\e{2}%
\e{0}%
\e{0}%
\e{0}%
\eol}\vss}\rg%
%
%
\rx{\vss\hfull{%
\rlx{\hss{$5600_{w}$}}\cg%
\e{0}%
\e{0}%
\e{1}%
\e{2}%
\e{2}%
\e{9}%
\e{6}%
\e{1}%
\e{9}%
\e{8}%
\e{16}%
\e{6}%
\e{4}%
\e{8}%
\e{9}%
\e{9}%
\e{2}%
\e{1}%
\e{2}%
\e{1}%
\e{0}%
\e{0}%
\eol}\vss}\rg%
%
%
\rx{\vss\hfull{%
\rlx{\hss{$2016_{w}$}}\cg%
\e{0}%
\e{0}%
\e{1}%
\e{0}%
\e{2}%
\e{2}%
\e{0}%
\e{1}%
\e{5}%
\e{4}%
\e{2}%
\e{0}%
\e{6}%
\e{4}%
\e{5}%
\e{2}%
\e{0}%
\e{1}%
\e{2}%
\e{1}%
\e{0}%
\e{0}%
\eol}\vss}\rg%
%
%
\rx{\vss\hfull{%
\rlx{\hss{$7168_{w}$}}\cg%
\e{0}%
\e{0}%
\e{2}%
\e{2}%
\e{4}%
\e{10}%
\e{4}%
\e{2}%
\e{14}%
\e{12}%
\e{16}%
\e{4}%
\e{10}%
\e{12}%
\e{14}%
\e{10}%
\e{2}%
\e{2}%
\e{4}%
\e{2}%
\e{0}%
\e{0}%
\eol}\vss}\rg%
\tableclose%
%
%
%
%
%
%
\eop
\eject
\tableopen{Induce/restrict matrix for $W({A_{6}}{A_{1}})\,\subset\,W(E_{8})$}%
%
%
%
%
%
%
\rowpts=18 true pt%
\colpts=18 true pt%
\rowlabpts=40 true pt%
\collabpts=70 true pt%
\clx{\vss\hfull{%
\rlx{\hss{$ $}}\cg%
\cx{\hskip 16 true pt\flip{$[{7}]{\times}[{2}]$}\hss}\cg%
\cx{\hskip 16 true pt\flip{$[{6}{1}]{\times}[{2}]$}\hss}\cg%
\cx{\hskip 16 true pt\flip{$[{5}{2}]{\times}[{2}]$}\hss}\cg%
\cx{\hskip 16 true pt\flip{$[{5}{1^{2}}]{\times}[{2}]$}\hss}\cg%
\cx{\hskip 16 true pt\flip{$[{4}{3}]{\times}[{2}]$}\hss}\cg%
\cx{\hskip 16 true pt\flip{$[{4}{2}{1}]{\times}[{2}]$}\hss}\cg%
\cx{\hskip 16 true pt\flip{$[{4}{1^{3}}]{\times}[{2}]$}\hss}\cg%
\cx{\hskip 16 true pt\flip{$[{3^{2}}{1}]{\times}[{2}]$}\hss}\cg%
\cx{\hskip 16 true pt\flip{$[{3}{2^{2}}]{\times}[{2}]$}\hss}\cg%
\cx{\hskip 16 true pt\flip{$[{3}{2}{1^{2}}]{\times}[{2}]$}\hss}\cg%
\cx{\hskip 16 true pt\flip{$[{3}{1^{4}}]{\times}[{2}]$}\hss}\cg%
\cx{\hskip 16 true pt\flip{$[{2^{3}}{1}]{\times}[{2}]$}\hss}\cg%
\cx{\hskip 16 true pt\flip{$[{2^{2}}{1^{3}}]{\times}[{2}]$}\hss}\cg%
\cx{\hskip 16 true pt\flip{$[{2}{1^{5}}]{\times}[{2}]$}\hss}\cg%
\cx{\hskip 16 true pt\flip{$[{1^{7}}]{\times}[{2}]$}\hss}\cg%
\eol}}\rg%
%
%
\rx{\vss\hfull{%
\rlx{\hss{$1_{x}$}}\cg%
\e{1}%
\e{0}%
\e{0}%
\e{0}%
\e{0}%
\e{0}%
\e{0}%
\e{0}%
\e{0}%
\e{0}%
\e{0}%
\e{0}%
\e{0}%
\e{0}%
\e{0}%
\eol}\vss}\rg%
%
%
\rx{\vss\hfull{%
\rlx{\hss{$28_{x}$}}\cg%
\e{0}%
\e{1}%
\e{0}%
\e{1}%
\e{0}%
\e{0}%
\e{0}%
\e{0}%
\e{0}%
\e{0}%
\e{0}%
\e{0}%
\e{0}%
\e{0}%
\e{0}%
\eol}\vss}\rg%
%
%
\rx{\vss\hfull{%
\rlx{\hss{$35_{x}$}}\cg%
\e{2}%
\e{2}%
\e{1}%
\e{0}%
\e{0}%
\e{0}%
\e{0}%
\e{0}%
\e{0}%
\e{0}%
\e{0}%
\e{0}%
\e{0}%
\e{0}%
\e{0}%
\eol}\vss}\rg%
%
%
\rx{\vss\hfull{%
\rlx{\hss{$84_{x}$}}\cg%
\e{3}%
\e{3}%
\e{2}%
\e{0}%
\e{1}%
\e{0}%
\e{0}%
\e{0}%
\e{0}%
\e{0}%
\e{0}%
\e{0}%
\e{0}%
\e{0}%
\e{0}%
\eol}\vss}\rg%
%
%
\rx{\vss\hfull{%
\rlx{\hss{$50_{x}$}}\cg%
\e{1}%
\e{1}%
\e{1}%
\e{0}%
\e{1}%
\e{0}%
\e{0}%
\e{0}%
\e{0}%
\e{0}%
\e{0}%
\e{0}%
\e{0}%
\e{0}%
\e{0}%
\eol}\vss}\rg%
%
%
\rx{\vss\hfull{%
\rlx{\hss{$350_{x}$}}\cg%
\e{0}%
\e{1}%
\e{1}%
\e{3}%
\e{0}%
\e{2}%
\e{2}%
\e{0}%
\e{0}%
\e{1}%
\e{0}%
\e{0}%
\e{0}%
\e{0}%
\e{0}%
\eol}\vss}\rg%
%
%
\rx{\vss\hfull{%
\rlx{\hss{$300_{x}$}}\cg%
\e{1}%
\e{3}%
\e{4}%
\e{1}%
\e{1}%
\e{2}%
\e{0}%
\e{0}%
\e{1}%
\e{0}%
\e{0}%
\e{0}%
\e{0}%
\e{0}%
\e{0}%
\eol}\vss}\rg%
%
%
\rx{\vss\hfull{%
\rlx{\hss{$567_{x}$}}\cg%
\e{2}%
\e{7}%
\e{5}%
\e{6}%
\e{2}%
\e{3}%
\e{1}%
\e{1}%
\e{0}%
\e{0}%
\e{0}%
\e{0}%
\e{0}%
\e{0}%
\e{0}%
\eol}\vss}\rg%
%
%
\rx{\vss\hfull{%
\rlx{\hss{$210_{x}$}}\cg%
\e{2}%
\e{4}%
\e{3}%
\e{2}%
\e{1}%
\e{1}%
\e{0}%
\e{0}%
\e{0}%
\e{0}%
\e{0}%
\e{0}%
\e{0}%
\e{0}%
\e{0}%
\eol}\vss}\rg%
%
%
\rx{\vss\hfull{%
\rlx{\hss{$840_{x}$}}\cg%
\e{0}%
\e{1}%
\e{4}%
\e{1}%
\e{4}%
\e{4}%
\e{0}%
\e{2}%
\e{4}%
\e{1}%
\e{0}%
\e{2}%
\e{0}%
\e{0}%
\e{0}%
\eol}\vss}\rg%
%
%
\rx{\vss\hfull{%
\rlx{\hss{$700_{x}$}}\cg%
\e{4}%
\e{7}%
\e{8}%
\e{3}%
\e{5}%
\e{4}%
\e{0}%
\e{1}%
\e{1}%
\e{0}%
\e{0}%
\e{0}%
\e{0}%
\e{0}%
\e{0}%
\eol}\vss}\rg%
%
%
\rx{\vss\hfull{%
\rlx{\hss{$175_{x}$}}\cg%
\e{1}%
\e{1}%
\e{1}%
\e{0}%
\e{2}%
\e{1}%
\e{0}%
\e{1}%
\e{0}%
\e{0}%
\e{0}%
\e{0}%
\e{0}%
\e{0}%
\e{0}%
\eol}\vss}\rg%
%
%
\rx{\vss\hfull{%
\rlx{\hss{$1400_{x}$}}\cg%
\e{2}%
\e{6}%
\e{7}%
\e{7}%
\e{6}%
\e{8}%
\e{3}%
\e{4}%
\e{1}%
\e{2}%
\e{0}%
\e{0}%
\e{0}%
\e{0}%
\e{0}%
\eol}\vss}\rg%
%
%
\rx{\vss\hfull{%
\rlx{\hss{$1050_{x}$}}\cg%
\e{2}%
\e{5}%
\e{6}%
\e{4}%
\e{6}%
\e{6}%
\e{1}%
\e{4}%
\e{1}%
\e{1}%
\e{0}%
\e{0}%
\e{0}%
\e{0}%
\e{0}%
\eol}\vss}\rg%
%
%
\rx{\vss\hfull{%
\rlx{\hss{$1575_{x}$}}\cg%
\e{1}%
\e{6}%
\e{7}%
\e{10}%
\e{4}%
\e{9}%
\e{5}%
\e{3}%
\e{1}%
\e{3}%
\e{0}%
\e{0}%
\e{0}%
\e{0}%
\e{0}%
\eol}\vss}\rg%
%
%
\rx{\vss\hfull{%
\rlx{\hss{$1344_{x}$}}\cg%
\e{3}%
\e{9}%
\e{11}%
\e{7}%
\e{6}%
\e{8}%
\e{1}%
\e{3}%
\e{2}%
\e{1}%
\e{0}%
\e{0}%
\e{0}%
\e{0}%
\e{0}%
\eol}\vss}\rg%
%
%
\rx{\vss\hfull{%
\rlx{\hss{$2100_{x}$}}\cg%
\e{0}%
\e{2}%
\e{4}%
\e{8}%
\e{2}%
\e{10}%
\e{9}%
\e{3}%
\e{3}%
\e{6}%
\e{3}%
\e{1}%
\e{1}%
\e{0}%
\e{0}%
\eol}\vss}\rg%
%
%
\rx{\vss\hfull{%
\rlx{\hss{$2268_{x}$}}\cg%
\e{1}%
\e{6}%
\e{10}%
\e{10}%
\e{6}%
\e{13}%
\e{6}%
\e{4}%
\e{4}%
\e{4}%
\e{1}%
\e{1}%
\e{0}%
\e{0}%
\e{0}%
\eol}\vss}\rg%
%
%
\rx{\vss\hfull{%
\rlx{\hss{$525_{x}$}}\cg%
\e{0}%
\e{2}%
\e{2}%
\e{4}%
\e{2}%
\e{2}%
\e{2}%
\e{2}%
\e{0}%
\e{1}%
\e{0}%
\e{0}%
\e{0}%
\e{0}%
\e{0}%
\eol}\vss}\rg%
%
%
\rx{\vss\hfull{%
\rlx{\hss{$700_{xx}$}}\cg%
\e{0}%
\e{1}%
\e{2}%
\e{2}%
\e{3}%
\e{3}%
\e{1}%
\e{3}%
\e{1}%
\e{2}%
\e{0}%
\e{0}%
\e{0}%
\e{0}%
\e{0}%
\eol}\vss}\rg%
%
%
\rx{\vss\hfull{%
\rlx{\hss{$972_{x}$}}\cg%
\e{1}%
\e{3}%
\e{6}%
\e{2}%
\e{5}%
\e{6}%
\e{0}%
\e{2}%
\e{3}%
\e{1}%
\e{0}%
\e{1}%
\e{0}%
\e{0}%
\e{0}%
\eol}\vss}\rg%
%
%
\rx{\vss\hfull{%
\rlx{\hss{$4096_{x}$}}\cg%
\e{1}%
\e{7}%
\e{14}%
\e{14}%
\e{9}%
\e{22}%
\e{9}%
\e{9}%
\e{8}%
\e{10}%
\e{2}%
\e{2}%
\e{1}%
\e{0}%
\e{0}%
\eol}\vss}\rg%
%
%
\rx{\vss\hfull{%
\rlx{\hss{$4200_{x}$}}\cg%
\e{1}%
\e{6}%
\e{13}%
\e{10}%
\e{12}%
\e{21}%
\e{7}%
\e{12}%
\e{10}%
\e{10}%
\e{2}%
\e{3}%
\e{1}%
\e{0}%
\e{0}%
\eol}\vss}\rg%
%
%
\rx{\vss\hfull{%
\rlx{\hss{$2240_{x}$}}\cg%
\e{2}%
\e{6}%
\e{11}%
\e{6}%
\e{10}%
\e{12}%
\e{2}%
\e{7}%
\e{5}%
\e{4}%
\e{0}%
\e{1}%
\e{0}%
\e{0}%
\e{0}%
\eol}\vss}\rg%
%
%
\rx{\vss\hfull{%
\rlx{\hss{$2835_{x}$}}\cg%
\e{1}%
\e{3}%
\e{7}%
\e{4}%
\e{9}%
\e{13}%
\e{3}%
\e{10}%
\e{7}%
\e{8}%
\e{1}%
\e{2}%
\e{1}%
\e{0}%
\e{0}%
\eol}\vss}\rg%
%
%
\rx{\vss\hfull{%
\rlx{\hss{$6075_{x}$}}\cg%
\e{0}%
\e{5}%
\e{11}%
\e{15}%
\e{11}%
\e{27}%
\e{16}%
\e{15}%
\e{10}%
\e{20}%
\e{5}%
\e{4}%
\e{4}%
\e{0}%
\e{0}%
\eol}\vss}\rg%
%
%
\rx{\vss\hfull{%
\rlx{\hss{$3200_{x}$}}\cg%
\e{0}%
\e{2}%
\e{7}%
\e{5}%
\e{6}%
\e{15}%
\e{4}%
\e{7}%
\e{10}%
\e{9}%
\e{2}%
\e{5}%
\e{2}%
\e{1}%
\e{0}%
\eol}\vss}\rg%
%
%
\rx{\vss\hfull{%
\rlx{\hss{$70_{y}$}}\cg%
\e{0}%
\e{0}%
\e{0}%
\e{0}%
\e{0}%
\e{0}%
\e{1}%
\e{0}%
\e{0}%
\e{0}%
\e{1}%
\e{0}%
\e{0}%
\e{0}%
\e{0}%
\eol}\vss}\rg%
%
%
\rx{\vss\hfull{%
\rlx{\hss{$1134_{y}$}}\cg%
\e{0}%
\e{0}%
\e{1}%
\e{2}%
\e{1}%
\e{4}%
\e{3}%
\e{2}%
\e{1}%
\e{5}%
\e{1}%
\e{2}%
\e{2}%
\e{0}%
\e{0}%
\eol}\vss}\rg%
%
%
\rx{\vss\hfull{%
\rlx{\hss{$1680_{y}$}}\cg%
\e{0}%
\e{0}%
\e{1}%
\e{3}%
\e{0}%
\e{5}%
\e{8}%
\e{1}%
\e{2}%
\e{7}%
\e{6}%
\e{1}%
\e{2}%
\e{1}%
\e{0}%
\eol}\vss}\rg%
%
%
\rx{\vss\hfull{%
\rlx{\hss{$168_{y}$}}\cg%
\e{0}%
\e{0}%
\e{0}%
\e{0}%
\e{1}%
\e{1}%
\e{0}%
\e{0}%
\e{1}%
\e{0}%
\e{0}%
\e{1}%
\e{0}%
\e{0}%
\e{0}%
\eol}\vss}\rg%
%
%
\rx{\vss\hfull{%
\rlx{\hss{$420_{y}$}}\cg%
\e{0}%
\e{0}%
\e{1}%
\e{0}%
\e{2}%
\e{1}%
\e{0}%
\e{2}%
\e{2}%
\e{1}%
\e{0}%
\e{1}%
\e{0}%
\e{0}%
\e{0}%
\eol}\vss}\rg%
%
%
\rx{\vss\hfull{%
\rlx{\hss{$3150_{y}$}}\cg%
\e{0}%
\e{1}%
\e{4}%
\e{3}%
\e{6}%
\e{11}%
\e{4}%
\e{10}%
\e{8}%
\e{11}%
\e{2}%
\e{5}%
\e{4}%
\e{0}%
\e{0}%
\eol}\vss}\rg%
%
%
\rx{\vss\hfull{%
\rlx{\hss{$4200_{y}$}}\cg%
\e{0}%
\e{1}%
\e{5}%
\e{4}%
\e{8}%
\e{16}%
\e{5}%
\e{11}%
\e{12}%
\e{14}%
\e{3}%
\e{8}%
\e{4}%
\e{1}%
\e{0}%
\eol}\vss}\rg%
\eop
\eject
\tablecont%
%
%
%
%
%
%
\rowpts=18 true pt%
\colpts=18 true pt%
\rowlabpts=40 true pt%
\collabpts=70 true pt%
\clx{\vss\hfull{%
\rlx{\hss{$ $}}\cg%
\cx{\hskip 16 true pt\flip{$[{7}]{\times}[{2}]$}\hss}\cg%
\cx{\hskip 16 true pt\flip{$[{6}{1}]{\times}[{2}]$}\hss}\cg%
\cx{\hskip 16 true pt\flip{$[{5}{2}]{\times}[{2}]$}\hss}\cg%
\cx{\hskip 16 true pt\flip{$[{5}{1^{2}}]{\times}[{2}]$}\hss}\cg%
\cx{\hskip 16 true pt\flip{$[{4}{3}]{\times}[{2}]$}\hss}\cg%
\cx{\hskip 16 true pt\flip{$[{4}{2}{1}]{\times}[{2}]$}\hss}\cg%
\cx{\hskip 16 true pt\flip{$[{4}{1^{3}}]{\times}[{2}]$}\hss}\cg%
\cx{\hskip 16 true pt\flip{$[{3^{2}}{1}]{\times}[{2}]$}\hss}\cg%
\cx{\hskip 16 true pt\flip{$[{3}{2^{2}}]{\times}[{2}]$}\hss}\cg%
\cx{\hskip 16 true pt\flip{$[{3}{2}{1^{2}}]{\times}[{2}]$}\hss}\cg%
\cx{\hskip 16 true pt\flip{$[{3}{1^{4}}]{\times}[{2}]$}\hss}\cg%
\cx{\hskip 16 true pt\flip{$[{2^{3}}{1}]{\times}[{2}]$}\hss}\cg%
\cx{\hskip 16 true pt\flip{$[{2^{2}}{1^{3}}]{\times}[{2}]$}\hss}\cg%
\cx{\hskip 16 true pt\flip{$[{2}{1^{5}}]{\times}[{2}]$}\hss}\cg%
\cx{\hskip 16 true pt\flip{$[{1^{7}}]{\times}[{2}]$}\hss}\cg%
\eol}}\rg%
%
%
\rx{\vss\hfull{%
\rlx{\hss{$2688_{y}$}}\cg%
\e{0}%
\e{1}%
\e{3}%
\e{3}%
\e{3}%
\e{10}%
\e{4}%
\e{7}%
\e{7}%
\e{10}%
\e{3}%
\e{3}%
\e{3}%
\e{1}%
\e{0}%
\eol}\vss}\rg%
%
%
\rx{\vss\hfull{%
\rlx{\hss{$2100_{y}$}}\cg%
\e{0}%
\e{0}%
\e{1}%
\e{4}%
\e{2}%
\e{6}%
\e{8}%
\e{4}%
\e{3}%
\e{8}%
\e{5}%
\e{2}%
\e{3}%
\e{1}%
\e{0}%
\eol}\vss}\rg%
%
%
\rx{\vss\hfull{%
\rlx{\hss{$1400_{y}$}}\cg%
\e{0}%
\e{0}%
\e{1}%
\e{2}%
\e{1}%
\e{5}%
\e{5}%
\e{1}%
\e{3}%
\e{5}%
\e{4}%
\e{2}%
\e{1}%
\e{1}%
\e{0}%
\eol}\vss}\rg%
%
%
\rx{\vss\hfull{%
\rlx{\hss{$4536_{y}$}}\cg%
\e{0}%
\e{1}%
\e{4}%
\e{6}%
\e{4}%
\e{16}%
\e{12}%
\e{8}%
\e{10}%
\e{17}%
\e{9}%
\e{5}%
\e{5}%
\e{2}%
\e{0}%
\eol}\vss}\rg%
%
%
\rx{\vss\hfull{%
\rlx{\hss{$5670_{y}$}}\cg%
\e{0}%
\e{1}%
\e{5}%
\e{8}%
\e{5}%
\e{20}%
\e{15}%
\e{10}%
\e{11}%
\e{22}%
\e{10}%
\e{7}%
\e{7}%
\e{2}%
\e{0}%
\eol}\vss}\rg%
%
%
\rx{\vss\hfull{%
\rlx{\hss{$4480_{y}$}}\cg%
\e{0}%
\e{1}%
\e{5}%
\e{5}%
\e{7}%
\e{16}%
\e{8}%
\e{11}%
\e{11}%
\e{16}%
\e{5}%
\e{7}%
\e{5}%
\e{1}%
\e{0}%
\eol}\vss}\rg%
%
%
\rx{\vss\hfull{%
\rlx{\hss{$8_{z}$}}\cg%
\e{1}%
\e{1}%
\e{0}%
\e{0}%
\e{0}%
\e{0}%
\e{0}%
\e{0}%
\e{0}%
\e{0}%
\e{0}%
\e{0}%
\e{0}%
\e{0}%
\e{0}%
\eol}\vss}\rg%
%
%
\rx{\vss\hfull{%
\rlx{\hss{$56_{z}$}}\cg%
\e{0}%
\e{0}%
\e{0}%
\e{1}%
\e{0}%
\e{0}%
\e{1}%
\e{0}%
\e{0}%
\e{0}%
\e{0}%
\e{0}%
\e{0}%
\e{0}%
\e{0}%
\eol}\vss}\rg%
%
%
\rx{\vss\hfull{%
\rlx{\hss{$160_{z}$}}\cg%
\e{1}%
\e{3}%
\e{2}%
\e{2}%
\e{0}%
\e{1}%
\e{0}%
\e{0}%
\e{0}%
\e{0}%
\e{0}%
\e{0}%
\e{0}%
\e{0}%
\e{0}%
\eol}\vss}\rg%
%
%
\rx{\vss\hfull{%
\rlx{\hss{$112_{z}$}}\cg%
\e{3}%
\e{4}%
\e{2}%
\e{1}%
\e{1}%
\e{0}%
\e{0}%
\e{0}%
\e{0}%
\e{0}%
\e{0}%
\e{0}%
\e{0}%
\e{0}%
\e{0}%
\eol}\vss}\rg%
%
%
\rx{\vss\hfull{%
\rlx{\hss{$840_{z}$}}\cg%
\e{0}%
\e{2}%
\e{4}%
\e{3}%
\e{2}%
\e{5}%
\e{1}%
\e{1}%
\e{2}%
\e{2}%
\e{0}%
\e{1}%
\e{0}%
\e{0}%
\e{0}%
\eol}\vss}\rg%
%
%
\rx{\vss\hfull{%
\rlx{\hss{$1296_{z}$}}\cg%
\e{0}%
\e{3}%
\e{4}%
\e{8}%
\e{1}%
\e{7}%
\e{6}%
\e{2}%
\e{1}%
\e{3}%
\e{1}%
\e{0}%
\e{0}%
\e{0}%
\e{0}%
\eol}\vss}\rg%
%
%
\rx{\vss\hfull{%
\rlx{\hss{$1400_{z}$}}\cg%
\e{3}%
\e{9}%
\e{11}%
\e{8}%
\e{6}%
\e{8}%
\e{2}%
\e{3}%
\e{2}%
\e{1}%
\e{0}%
\e{0}%
\e{0}%
\e{0}%
\e{0}%
\eol}\vss}\rg%
%
%
\rx{\vss\hfull{%
\rlx{\hss{$1008_{z}$}}\cg%
\e{1}%
\e{6}%
\e{6}%
\e{8}%
\e{3}%
\e{6}%
\e{3}%
\e{1}%
\e{1}%
\e{1}%
\e{0}%
\e{0}%
\e{0}%
\e{0}%
\e{0}%
\eol}\vss}\rg%
%
%
\rx{\vss\hfull{%
\rlx{\hss{$560_{z}$}}\cg%
\e{4}%
\e{8}%
\e{7}%
\e{4}%
\e{3}%
\e{3}%
\e{0}%
\e{1}%
\e{0}%
\e{0}%
\e{0}%
\e{0}%
\e{0}%
\e{0}%
\e{0}%
\eol}\vss}\rg%
%
%
\rx{\vss\hfull{%
\rlx{\hss{$1400_{zz}$}}\cg%
\e{2}%
\e{4}%
\e{7}%
\e{3}%
\e{8}%
\e{7}%
\e{1}%
\e{6}%
\e{3}%
\e{2}%
\e{0}%
\e{0}%
\e{0}%
\e{0}%
\e{0}%
\eol}\vss}\rg%
%
%
\rx{\vss\hfull{%
\rlx{\hss{$4200_{z}$}}\cg%
\e{0}%
\e{3}%
\e{7}%
\e{8}%
\e{10}%
\e{18}%
\e{8}%
\e{12}%
\e{8}%
\e{14}%
\e{3}%
\e{4}%
\e{2}%
\e{0}%
\e{0}%
\eol}\vss}\rg%
%
%
\rx{\vss\hfull{%
\rlx{\hss{$400_{z}$}}\cg%
\e{3}%
\e{4}%
\e{4}%
\e{2}%
\e{4}%
\e{2}%
\e{0}%
\e{1}%
\e{0}%
\e{0}%
\e{0}%
\e{0}%
\e{0}%
\e{0}%
\e{0}%
\eol}\vss}\rg%
%
%
\rx{\vss\hfull{%
\rlx{\hss{$3240_{z}$}}\cg%
\e{3}%
\e{12}%
\e{17}%
\e{13}%
\e{13}%
\e{19}%
\e{5}%
\e{8}%
\e{5}%
\e{5}%
\e{0}%
\e{1}%
\e{0}%
\e{0}%
\e{0}%
\eol}\vss}\rg%
%
%
\rx{\vss\hfull{%
\rlx{\hss{$4536_{z}$}}\cg%
\e{1}%
\e{5}%
\e{14}%
\e{8}%
\e{14}%
\e{22}%
\e{5}%
\e{13}%
\e{13}%
\e{11}%
\e{1}%
\e{5}%
\e{2}%
\e{0}%
\e{0}%
\eol}\vss}\rg%
%
%
\rx{\vss\hfull{%
\rlx{\hss{$2400_{z}$}}\cg%
\e{0}%
\e{1}%
\e{2}%
\e{7}%
\e{2}%
\e{9}%
\e{11}%
\e{4}%
\e{2}%
\e{9}%
\e{5}%
\e{1}%
\e{2}%
\e{0}%
\e{0}%
\eol}\vss}\rg%
%
%
\rx{\vss\hfull{%
\rlx{\hss{$3360_{z}$}}\cg%
\e{1}%
\e{5}%
\e{9}%
\e{9}%
\e{8}%
\e{17}%
\e{7}%
\e{10}%
\e{6}%
\e{9}%
\e{2}%
\e{1}%
\e{1}%
\e{0}%
\e{0}%
\eol}\vss}\rg%
%
%
\rx{\vss\hfull{%
\rlx{\hss{$2800_{z}$}}\cg%
\e{0}%
\e{3}%
\e{6}%
\e{10}%
\e{5}%
\e{13}%
\e{10}%
\e{6}%
\e{4}%
\e{8}%
\e{3}%
\e{1}%
\e{1}%
\e{0}%
\e{0}%
\eol}\vss}\rg%
%
%
\rx{\vss\hfull{%
\rlx{\hss{$4096_{z}$}}\cg%
\e{1}%
\e{7}%
\e{14}%
\e{14}%
\e{9}%
\e{22}%
\e{9}%
\e{9}%
\e{8}%
\e{10}%
\e{2}%
\e{2}%
\e{1}%
\e{0}%
\e{0}%
\eol}\vss}\rg%
%
%
\rx{\vss\hfull{%
\rlx{\hss{$5600_{z}$}}\cg%
\e{0}%
\e{3}%
\e{9}%
\e{12}%
\e{8}%
\e{24}%
\e{15}%
\e{11}%
\e{12}%
\e{18}%
\e{7}%
\e{6}%
\e{4}%
\e{1}%
\e{0}%
\eol}\vss}\rg%
%
%
\rx{\vss\hfull{%
\rlx{\hss{$448_{z}$}}\cg%
\e{2}%
\e{3}%
\e{5}%
\e{1}%
\e{3}%
\e{2}%
\e{0}%
\e{2}%
\e{1}%
\e{0}%
\e{0}%
\e{0}%
\e{0}%
\e{0}%
\e{0}%
\eol}\vss}\rg%
%
%
\rx{\vss\hfull{%
\rlx{\hss{$448_{w}$}}\cg%
\e{0}%
\e{0}%
\e{0}%
\e{1}%
\e{0}%
\e{1}%
\e{3}%
\e{0}%
\e{0}%
\e{2}%
\e{2}%
\e{0}%
\e{1}%
\e{0}%
\e{0}%
\eol}\vss}\rg%
%
%
\rx{\vss\hfull{%
\rlx{\hss{$1344_{w}$}}\cg%
\e{0}%
\e{0}%
\e{2}%
\e{1}%
\e{4}%
\e{5}%
\e{1}%
\e{4}%
\e{4}%
\e{4}%
\e{0}%
\e{4}%
\e{1}%
\e{0}%
\e{0}%
\eol}\vss}\rg%
%
%
\rx{\vss\hfull{%
\rlx{\hss{$5600_{w}$}}\cg%
\e{0}%
\e{1}%
\e{5}%
\e{8}%
\e{5}%
\e{20}%
\e{14}%
\e{10}%
\e{11}%
\e{22}%
\e{9}%
\e{7}%
\e{7}%
\e{2}%
\e{0}%
\eol}\vss}\rg%
%
%
\rx{\vss\hfull{%
\rlx{\hss{$2016_{w}$}}\cg%
\e{0}%
\e{1}%
\e{3}%
\e{1}%
\e{5}%
\e{7}%
\e{1}%
\e{8}%
\e{7}%
\e{6}%
\e{1}%
\e{3}%
\e{2}%
\e{0}%
\e{0}%
\eol}\vss}\rg%
%
%
\rx{\vss\hfull{%
\rlx{\hss{$7168_{w}$}}\cg%
\e{0}%
\e{2}%
\e{8}%
\e{8}%
\e{10}%
\e{26}%
\e{12}%
\e{18}%
\e{18}%
\e{26}%
\e{8}%
\e{10}%
\e{8}%
\e{2}%
\e{0}%
\eol}\vss}\rg%
\eop
\eject
\tablecont%
%
%
%
%
%
%
\rowpts=18 true pt%
\colpts=18 true pt%
\rowlabpts=40 true pt%
\collabpts=70 true pt%
\clx{\vss\hfull{%
\rlx{\hss{$ $}}\cg%
\cx{\hskip 16 true pt\flip{$[{7}]{\times}[{1^{2}}]$}\hss}\cg%
\cx{\hskip 16 true pt\flip{$[{6}{1}]{\times}[{1^{2}}]$}\hss}\cg%
\cx{\hskip 16 true pt\flip{$[{5}{2}]{\times}[{1^{2}}]$}\hss}\cg%
\cx{\hskip 16 true pt\flip{$[{5}{1^{2}}]{\times}[{1^{2}}]$}\hss}\cg%
\cx{\hskip 16 true pt\flip{$[{4}{3}]{\times}[{1^{2}}]$}\hss}\cg%
\cx{\hskip 16 true pt\flip{$[{4}{2}{1}]{\times}[{1^{2}}]$}\hss}\cg%
\cx{\hskip 16 true pt\flip{$[{4}{1^{3}}]{\times}[{1^{2}}]$}\hss}\cg%
\cx{\hskip 16 true pt\flip{$[{3^{2}}{1}]{\times}[{1^{2}}]$}\hss}\cg%
\cx{\hskip 16 true pt\flip{$[{3}{2^{2}}]{\times}[{1^{2}}]$}\hss}\cg%
\cx{\hskip 16 true pt\flip{$[{3}{2}{1^{2}}]{\times}[{1^{2}}]$}\hss}\cg%
\cx{\hskip 16 true pt\flip{$[{3}{1^{4}}]{\times}[{1^{2}}]$}\hss}\cg%
\cx{\hskip 16 true pt\flip{$[{2^{3}}{1}]{\times}[{1^{2}}]$}\hss}\cg%
\cx{\hskip 16 true pt\flip{$[{2^{2}}{1^{3}}]{\times}[{1^{2}}]$}\hss}\cg%
\cx{\hskip 16 true pt\flip{$[{2}{1^{5}}]{\times}[{1^{2}}]$}\hss}\cg%
\cx{\hskip 16 true pt\flip{$[{1^{7}}]{\times}[{1^{2}}]$}\hss}\cg%
\eol}}\rg%
%
%
\rx{\vss\hfull{%
\rlx{\hss{$1_{x}$}}\cg%
\e{0}%
\e{0}%
\e{0}%
\e{0}%
\e{0}%
\e{0}%
\e{0}%
\e{0}%
\e{0}%
\e{0}%
\e{0}%
\e{0}%
\e{0}%
\e{0}%
\e{0}%
\eol}\vss}\rg%
%
%
\rx{\vss\hfull{%
\rlx{\hss{$28_{x}$}}\cg%
\e{1}%
\e{1}%
\e{0}%
\e{0}%
\e{0}%
\e{0}%
\e{0}%
\e{0}%
\e{0}%
\e{0}%
\e{0}%
\e{0}%
\e{0}%
\e{0}%
\e{0}%
\eol}\vss}\rg%
%
%
\rx{\vss\hfull{%
\rlx{\hss{$35_{x}$}}\cg%
\e{1}%
\e{1}%
\e{0}%
\e{0}%
\e{0}%
\e{0}%
\e{0}%
\e{0}%
\e{0}%
\e{0}%
\e{0}%
\e{0}%
\e{0}%
\e{0}%
\e{0}%
\eol}\vss}\rg%
%
%
\rx{\vss\hfull{%
\rlx{\hss{$84_{x}$}}\cg%
\e{1}%
\e{1}%
\e{1}%
\e{0}%
\e{0}%
\e{0}%
\e{0}%
\e{0}%
\e{0}%
\e{0}%
\e{0}%
\e{0}%
\e{0}%
\e{0}%
\e{0}%
\eol}\vss}\rg%
%
%
\rx{\vss\hfull{%
\rlx{\hss{$50_{x}$}}\cg%
\e{1}%
\e{0}%
\e{0}%
\e{0}%
\e{1}%
\e{0}%
\e{0}%
\e{0}%
\e{0}%
\e{0}%
\e{0}%
\e{0}%
\e{0}%
\e{0}%
\e{0}%
\eol}\vss}\rg%
%
%
\rx{\vss\hfull{%
\rlx{\hss{$350_{x}$}}\cg%
\e{0}%
\e{2}%
\e{2}%
\e{3}%
\e{0}%
\e{1}%
\e{1}%
\e{0}%
\e{0}%
\e{0}%
\e{0}%
\e{0}%
\e{0}%
\e{0}%
\e{0}%
\eol}\vss}\rg%
%
%
\rx{\vss\hfull{%
\rlx{\hss{$300_{x}$}}\cg%
\e{0}%
\e{2}%
\e{2}%
\e{2}%
\e{0}%
\e{1}%
\e{0}%
\e{0}%
\e{0}%
\e{0}%
\e{0}%
\e{0}%
\e{0}%
\e{0}%
\e{0}%
\eol}\vss}\rg%
%
%
\rx{\vss\hfull{%
\rlx{\hss{$567_{x}$}}\cg%
\e{3}%
\e{6}%
\e{4}%
\e{3}%
\e{1}%
\e{1}%
\e{0}%
\e{0}%
\e{0}%
\e{0}%
\e{0}%
\e{0}%
\e{0}%
\e{0}%
\e{0}%
\eol}\vss}\rg%
%
%
\rx{\vss\hfull{%
\rlx{\hss{$210_{x}$}}\cg%
\e{2}%
\e{3}%
\e{1}%
\e{1}%
\e{1}%
\e{0}%
\e{0}%
\e{0}%
\e{0}%
\e{0}%
\e{0}%
\e{0}%
\e{0}%
\e{0}%
\e{0}%
\eol}\vss}\rg%
%
%
\rx{\vss\hfull{%
\rlx{\hss{$840_{x}$}}\cg%
\e{0}%
\e{0}%
\e{2}%
\e{1}%
\e{2}%
\e{4}%
\e{1}%
\e{1}%
\e{2}%
\e{2}%
\e{0}%
\e{1}%
\e{0}%
\e{0}%
\e{0}%
\eol}\vss}\rg%
%
%
\rx{\vss\hfull{%
\rlx{\hss{$700_{x}$}}\cg%
\e{2}%
\e{4}%
\e{4}%
\e{2}%
\e{3}%
\e{2}%
\e{0}%
\e{1}%
\e{0}%
\e{0}%
\e{0}%
\e{0}%
\e{0}%
\e{0}%
\e{0}%
\eol}\vss}\rg%
%
%
\rx{\vss\hfull{%
\rlx{\hss{$175_{x}$}}\cg%
\e{0}%
\e{0}%
\e{1}%
\e{0}%
\e{1}%
\e{0}%
\e{0}%
\e{1}%
\e{1}%
\e{0}%
\e{0}%
\e{0}%
\e{0}%
\e{0}%
\e{0}%
\eol}\vss}\rg%
%
%
\rx{\vss\hfull{%
\rlx{\hss{$1400_{x}$}}\cg%
\e{2}%
\e{5}%
\e{6}%
\e{4}%
\e{5}%
\e{5}%
\e{1}%
\e{3}%
\e{1}%
\e{1}%
\e{0}%
\e{0}%
\e{0}%
\e{0}%
\e{0}%
\eol}\vss}\rg%
%
%
\rx{\vss\hfull{%
\rlx{\hss{$1050_{x}$}}\cg%
\e{2}%
\e{3}%
\e{4}%
\e{1}%
\e{5}%
\e{4}%
\e{0}%
\e{3}%
\e{1}%
\e{1}%
\e{0}%
\e{0}%
\e{0}%
\e{0}%
\e{0}%
\eol}\vss}\rg%
%
%
\rx{\vss\hfull{%
\rlx{\hss{$1575_{x}$}}\cg%
\e{2}%
\e{7}%
\e{7}%
\e{7}%
\e{4}%
\e{6}%
\e{2}%
\e{2}%
\e{0}%
\e{1}%
\e{0}%
\e{0}%
\e{0}%
\e{0}%
\e{0}%
\eol}\vss}\rg%
%
%
\rx{\vss\hfull{%
\rlx{\hss{$1344_{x}$}}\cg%
\e{2}%
\e{6}%
\e{8}%
\e{5}%
\e{3}%
\e{5}%
\e{1}%
\e{1}%
\e{1}%
\e{0}%
\e{0}%
\e{0}%
\e{0}%
\e{0}%
\e{0}%
\eol}\vss}\rg%
%
%
\rx{\vss\hfull{%
\rlx{\hss{$2100_{x}$}}\cg%
\e{0}%
\e{4}%
\e{6}%
\e{10}%
\e{2}%
\e{9}%
\e{7}%
\e{2}%
\e{2}%
\e{3}%
\e{1}%
\e{0}%
\e{0}%
\e{0}%
\e{0}%
\eol}\vss}\rg%
%
%
\rx{\vss\hfull{%
\rlx{\hss{$2268_{x}$}}\cg%
\e{1}%
\e{7}%
\e{8}%
\e{10}%
\e{5}%
\e{9}%
\e{4}%
\e{3}%
\e{2}%
\e{2}%
\e{0}%
\e{0}%
\e{0}%
\e{0}%
\e{0}%
\eol}\vss}\rg%
%
%
\rx{\vss\hfull{%
\rlx{\hss{$525_{x}$}}\cg%
\e{1}%
\e{2}%
\e{3}%
\e{2}%
\e{1}%
\e{2}%
\e{1}%
\e{0}%
\e{1}%
\e{0}%
\e{0}%
\e{0}%
\e{0}%
\e{0}%
\e{0}%
\eol}\vss}\rg%
%
%
\rx{\vss\hfull{%
\rlx{\hss{$700_{xx}$}}\cg%
\e{1}%
\e{1}%
\e{1}%
\e{0}%
\e{4}%
\e{2}%
\e{0}%
\e{3}%
\e{1}%
\e{2}%
\e{0}%
\e{1}%
\e{0}%
\e{0}%
\e{0}%
\eol}\vss}\rg%
%
%
\rx{\vss\hfull{%
\rlx{\hss{$972_{x}$}}\cg%
\e{0}%
\e{2}%
\e{3}%
\e{2}%
\e{3}%
\e{4}%
\e{1}%
\e{2}%
\e{2}%
\e{1}%
\e{0}%
\e{0}%
\e{0}%
\e{0}%
\e{0}%
\eol}\vss}\rg%
%
%
\rx{\vss\hfull{%
\rlx{\hss{$4096_{x}$}}\cg%
\e{1}%
\e{7}%
\e{13}%
\e{13}%
\e{8}%
\e{18}%
\e{8}%
\e{6}%
\e{5}%
\e{6}%
\e{1}%
\e{1}%
\e{0}%
\e{0}%
\e{0}%
\eol}\vss}\rg%
%
%
\rx{\vss\hfull{%
\rlx{\hss{$4200_{x}$}}\cg%
\e{1}%
\e{5}%
\e{10}%
\e{8}%
\e{10}%
\e{18}%
\e{5}%
\e{10}%
\e{7}%
\e{9}%
\e{1}%
\e{2}%
\e{1}%
\e{0}%
\e{0}%
\eol}\vss}\rg%
%
%
\rx{\vss\hfull{%
\rlx{\hss{$2240_{x}$}}\cg%
\e{1}%
\e{3}%
\e{8}%
\e{4}%
\e{7}%
\e{9}%
\e{2}%
\e{5}%
\e{4}%
\e{3}%
\e{0}%
\e{1}%
\e{0}%
\e{0}%
\e{0}%
\eol}\vss}\rg%
%
%
\rx{\vss\hfull{%
\rlx{\hss{$2835_{x}$}}\cg%
\e{0}%
\e{2}%
\e{5}%
\e{3}%
\e{7}%
\e{10}%
\e{2}%
\e{10}%
\e{8}%
\e{7}%
\e{1}%
\e{3}%
\e{2}%
\e{0}%
\e{0}%
\eol}\vss}\rg%
%
%
\rx{\vss\hfull{%
\rlx{\hss{$6075_{x}$}}\cg%
\e{1}%
\e{6}%
\e{13}%
\e{13}%
\e{12}%
\e{25}%
\e{12}%
\e{13}%
\e{10}%
\e{15}%
\e{4}%
\e{3}%
\e{2}%
\e{0}%
\e{0}%
\eol}\vss}\rg%
%
%
\rx{\vss\hfull{%
\rlx{\hss{$3200_{x}$}}\cg%
\e{0}%
\e{1}%
\e{6}%
\e{6}%
\e{4}%
\e{14}%
\e{7}%
\e{5}%
\e{8}%
\e{8}%
\e{3}%
\e{3}%
\e{1}%
\e{0}%
\e{0}%
\eol}\vss}\rg%
%
%
\rx{\vss\hfull{%
\rlx{\hss{$70_{y}$}}\cg%
\e{0}%
\e{0}%
\e{0}%
\e{1}%
\e{0}%
\e{0}%
\e{1}%
\e{0}%
\e{0}%
\e{0}%
\e{0}%
\e{0}%
\e{0}%
\e{0}%
\e{0}%
\eol}\vss}\rg%
%
%
\rx{\vss\hfull{%
\rlx{\hss{$1134_{y}$}}\cg%
\e{0}%
\e{0}%
\e{2}%
\e{1}%
\e{2}%
\e{5}%
\e{3}%
\e{1}%
\e{2}%
\e{4}%
\e{2}%
\e{1}%
\e{1}%
\e{0}%
\e{0}%
\eol}\vss}\rg%
%
%
\rx{\vss\hfull{%
\rlx{\hss{$1680_{y}$}}\cg%
\e{0}%
\e{1}%
\e{2}%
\e{6}%
\e{1}%
\e{7}%
\e{8}%
\e{2}%
\e{1}%
\e{5}%
\e{3}%
\e{0}%
\e{1}%
\e{0}%
\e{0}%
\eol}\vss}\rg%
%
%
\rx{\vss\hfull{%
\rlx{\hss{$168_{y}$}}\cg%
\e{0}%
\e{0}%
\e{0}%
\e{0}%
\e{1}%
\e{0}%
\e{0}%
\e{1}%
\e{0}%
\e{1}%
\e{0}%
\e{1}%
\e{0}%
\e{0}%
\e{0}%
\eol}\vss}\rg%
%
%
\rx{\vss\hfull{%
\rlx{\hss{$420_{y}$}}\cg%
\e{0}%
\e{0}%
\e{0}%
\e{0}%
\e{1}%
\e{1}%
\e{0}%
\e{2}%
\e{2}%
\e{1}%
\e{0}%
\e{2}%
\e{1}%
\e{0}%
\e{0}%
\eol}\vss}\rg%
%
%
\rx{\vss\hfull{%
\rlx{\hss{$3150_{y}$}}\cg%
\e{0}%
\e{0}%
\e{4}%
\e{2}%
\e{5}%
\e{11}%
\e{4}%
\e{8}%
\e{10}%
\e{11}%
\e{3}%
\e{6}%
\e{4}%
\e{1}%
\e{0}%
\eol}\vss}\rg%
%
%
\rx{\vss\hfull{%
\rlx{\hss{$4200_{y}$}}\cg%
\e{0}%
\e{1}%
\e{4}%
\e{3}%
\e{8}%
\e{14}%
\e{5}%
\e{12}%
\e{11}%
\e{16}%
\e{4}%
\e{8}%
\e{5}%
\e{1}%
\e{0}%
\eol}\vss}\rg%
\eop
\eject
\tablecont%
%
%
%
%
%
%
\rowpts=18 true pt%
\colpts=18 true pt%
\rowlabpts=40 true pt%
\collabpts=70 true pt%
\clx{\vss\hfull{%
\rlx{\hss{$ $}}\cg%
\cx{\hskip 16 true pt\flip{$[{7}]{\times}[{1^{2}}]$}\hss}\cg%
\cx{\hskip 16 true pt\flip{$[{6}{1}]{\times}[{1^{2}}]$}\hss}\cg%
\cx{\hskip 16 true pt\flip{$[{5}{2}]{\times}[{1^{2}}]$}\hss}\cg%
\cx{\hskip 16 true pt\flip{$[{5}{1^{2}}]{\times}[{1^{2}}]$}\hss}\cg%
\cx{\hskip 16 true pt\flip{$[{4}{3}]{\times}[{1^{2}}]$}\hss}\cg%
\cx{\hskip 16 true pt\flip{$[{4}{2}{1}]{\times}[{1^{2}}]$}\hss}\cg%
\cx{\hskip 16 true pt\flip{$[{4}{1^{3}}]{\times}[{1^{2}}]$}\hss}\cg%
\cx{\hskip 16 true pt\flip{$[{3^{2}}{1}]{\times}[{1^{2}}]$}\hss}\cg%
\cx{\hskip 16 true pt\flip{$[{3}{2^{2}}]{\times}[{1^{2}}]$}\hss}\cg%
\cx{\hskip 16 true pt\flip{$[{3}{2}{1^{2}}]{\times}[{1^{2}}]$}\hss}\cg%
\cx{\hskip 16 true pt\flip{$[{3}{1^{4}}]{\times}[{1^{2}}]$}\hss}\cg%
\cx{\hskip 16 true pt\flip{$[{2^{3}}{1}]{\times}[{1^{2}}]$}\hss}\cg%
\cx{\hskip 16 true pt\flip{$[{2^{2}}{1^{3}}]{\times}[{1^{2}}]$}\hss}\cg%
\cx{\hskip 16 true pt\flip{$[{2}{1^{5}}]{\times}[{1^{2}}]$}\hss}\cg%
\cx{\hskip 16 true pt\flip{$[{1^{7}}]{\times}[{1^{2}}]$}\hss}\cg%
\eol}}\rg%
%
%
\rx{\vss\hfull{%
\rlx{\hss{$2688_{y}$}}\cg%
\e{0}%
\e{1}%
\e{3}%
\e{3}%
\e{3}%
\e{10}%
\e{4}%
\e{7}%
\e{7}%
\e{10}%
\e{3}%
\e{3}%
\e{3}%
\e{1}%
\e{0}%
\eol}\vss}\rg%
%
%
\rx{\vss\hfull{%
\rlx{\hss{$2100_{y}$}}\cg%
\e{0}%
\e{1}%
\e{3}%
\e{5}%
\e{2}%
\e{8}%
\e{8}%
\e{3}%
\e{4}%
\e{6}%
\e{4}%
\e{2}%
\e{1}%
\e{0}%
\e{0}%
\eol}\vss}\rg%
%
%
\rx{\vss\hfull{%
\rlx{\hss{$1400_{y}$}}\cg%
\e{0}%
\e{1}%
\e{1}%
\e{4}%
\e{2}%
\e{5}%
\e{5}%
\e{3}%
\e{1}%
\e{5}%
\e{2}%
\e{1}%
\e{1}%
\e{0}%
\e{0}%
\eol}\vss}\rg%
%
%
\rx{\vss\hfull{%
\rlx{\hss{$4536_{y}$}}\cg%
\e{0}%
\e{2}%
\e{5}%
\e{9}%
\e{5}%
\e{17}%
\e{12}%
\e{10}%
\e{8}%
\e{16}%
\e{6}%
\e{4}%
\e{4}%
\e{1}%
\e{0}%
\eol}\vss}\rg%
%
%
\rx{\vss\hfull{%
\rlx{\hss{$5670_{y}$}}\cg%
\e{0}%
\e{2}%
\e{7}%
\e{10}%
\e{7}%
\e{22}%
\e{15}%
\e{11}%
\e{10}%
\e{20}%
\e{8}%
\e{5}%
\e{5}%
\e{1}%
\e{0}%
\eol}\vss}\rg%
%
%
\rx{\vss\hfull{%
\rlx{\hss{$4480_{y}$}}\cg%
\e{0}%
\e{1}%
\e{5}%
\e{5}%
\e{7}%
\e{16}%
\e{8}%
\e{11}%
\e{11}%
\e{16}%
\e{5}%
\e{7}%
\e{5}%
\e{1}%
\e{0}%
\eol}\vss}\rg%
%
%
\rx{\vss\hfull{%
\rlx{\hss{$8_{z}$}}\cg%
\e{1}%
\e{0}%
\e{0}%
\e{0}%
\e{0}%
\e{0}%
\e{0}%
\e{0}%
\e{0}%
\e{0}%
\e{0}%
\e{0}%
\e{0}%
\e{0}%
\e{0}%
\eol}\vss}\rg%
%
%
\rx{\vss\hfull{%
\rlx{\hss{$56_{z}$}}\cg%
\e{0}%
\e{1}%
\e{0}%
\e{1}%
\e{0}%
\e{0}%
\e{0}%
\e{0}%
\e{0}%
\e{0}%
\e{0}%
\e{0}%
\e{0}%
\e{0}%
\e{0}%
\eol}\vss}\rg%
%
%
\rx{\vss\hfull{%
\rlx{\hss{$160_{z}$}}\cg%
\e{1}%
\e{3}%
\e{1}%
\e{1}%
\e{0}%
\e{0}%
\e{0}%
\e{0}%
\e{0}%
\e{0}%
\e{0}%
\e{0}%
\e{0}%
\e{0}%
\e{0}%
\eol}\vss}\rg%
%
%
\rx{\vss\hfull{%
\rlx{\hss{$112_{z}$}}\cg%
\e{2}%
\e{2}%
\e{1}%
\e{0}%
\e{0}%
\e{0}%
\e{0}%
\e{0}%
\e{0}%
\e{0}%
\e{0}%
\e{0}%
\e{0}%
\e{0}%
\e{0}%
\eol}\vss}\rg%
%
%
\rx{\vss\hfull{%
\rlx{\hss{$840_{z}$}}\cg%
\e{0}%
\e{1}%
\e{4}%
\e{3}%
\e{1}%
\e{4}%
\e{2}%
\e{0}%
\e{1}%
\e{1}%
\e{0}%
\e{0}%
\e{0}%
\e{0}%
\e{0}%
\eol}\vss}\rg%
%
%
\rx{\vss\hfull{%
\rlx{\hss{$1296_{z}$}}\cg%
\e{1}%
\e{5}%
\e{5}%
\e{8}%
\e{2}%
\e{5}%
\e{3}%
\e{1}%
\e{0}%
\e{1}%
\e{0}%
\e{0}%
\e{0}%
\e{0}%
\e{0}%
\eol}\vss}\rg%
%
%
\rx{\vss\hfull{%
\rlx{\hss{$1400_{z}$}}\cg%
\e{2}%
\e{7}%
\e{8}%
\e{6}%
\e{3}%
\e{5}%
\e{1}%
\e{1}%
\e{1}%
\e{0}%
\e{0}%
\e{0}%
\e{0}%
\e{0}%
\e{0}%
\eol}\vss}\rg%
%
%
\rx{\vss\hfull{%
\rlx{\hss{$1008_{z}$}}\cg%
\e{2}%
\e{7}%
\e{5}%
\e{6}%
\e{2}%
\e{3}%
\e{1}%
\e{1}%
\e{0}%
\e{0}%
\e{0}%
\e{0}%
\e{0}%
\e{0}%
\e{0}%
\eol}\vss}\rg%
%
%
\rx{\vss\hfull{%
\rlx{\hss{$560_{z}$}}\cg%
\e{3}%
\e{5}%
\e{4}%
\e{2}%
\e{2}%
\e{1}%
\e{0}%
\e{0}%
\e{0}%
\e{0}%
\e{0}%
\e{0}%
\e{0}%
\e{0}%
\e{0}%
\eol}\vss}\rg%
%
%
\rx{\vss\hfull{%
\rlx{\hss{$1400_{zz}$}}\cg%
\e{1}%
\e{2}%
\e{4}%
\e{1}%
\e{6}%
\e{5}%
\e{0}%
\e{5}%
\e{3}%
\e{2}%
\e{0}%
\e{1}%
\e{0}%
\e{0}%
\e{0}%
\eol}\vss}\rg%
%
%
\rx{\vss\hfull{%
\rlx{\hss{$4200_{z}$}}\cg%
\e{1}%
\e{3}%
\e{7}%
\e{5}%
\e{11}%
\e{15}%
\e{5}%
\e{13}%
\e{8}%
\e{13}%
\e{2}%
\e{4}%
\e{3}%
\e{0}%
\e{0}%
\eol}\vss}\rg%
%
%
\rx{\vss\hfull{%
\rlx{\hss{$400_{z}$}}\cg%
\e{2}%
\e{2}%
\e{2}%
\e{0}%
\e{3}%
\e{1}%
\e{0}%
\e{1}%
\e{0}%
\e{0}%
\e{0}%
\e{0}%
\e{0}%
\e{0}%
\e{0}%
\eol}\vss}\rg%
%
%
\rx{\vss\hfull{%
\rlx{\hss{$3240_{z}$}}\cg%
\e{3}%
\e{9}%
\e{13}%
\e{9}%
\e{10}%
\e{13}%
\e{3}%
\e{6}%
\e{3}%
\e{3}%
\e{0}%
\e{0}%
\e{0}%
\e{0}%
\e{0}%
\eol}\vss}\rg%
%
%
\rx{\vss\hfull{%
\rlx{\hss{$4536_{z}$}}\cg%
\e{0}%
\e{3}%
\e{10}%
\e{7}%
\e{10}%
\e{19}%
\e{6}%
\e{10}%
\e{11}%
\e{10}%
\e{2}%
\e{4}%
\e{1}%
\e{0}%
\e{0}%
\eol}\vss}\rg%
%
%
\rx{\vss\hfull{%
\rlx{\hss{$2400_{z}$}}\cg%
\e{0}%
\e{3}%
\e{5}%
\e{8}%
\e{3}%
\e{10}%
\e{8}%
\e{4}%
\e{2}%
\e{6}%
\e{2}%
\e{0}%
\e{1}%
\e{0}%
\e{0}%
\eol}\vss}\rg%
%
%
\rx{\vss\hfull{%
\rlx{\hss{$3360_{z}$}}\cg%
\e{1}%
\e{5}%
\e{8}%
\e{7}%
\e{8}%
\e{14}%
\e{4}%
\e{9}%
\e{5}%
\e{7}%
\e{1}%
\e{1}%
\e{1}%
\e{0}%
\e{0}%
\eol}\vss}\rg%
%
%
\rx{\vss\hfull{%
\rlx{\hss{$2800_{z}$}}\cg%
\e{1}%
\e{5}%
\e{7}%
\e{10}%
\e{6}%
\e{11}%
\e{7}%
\e{5}%
\e{3}%
\e{5}%
\e{1}%
\e{1}%
\e{0}%
\e{0}%
\e{0}%
\eol}\vss}\rg%
%
%
\rx{\vss\hfull{%
\rlx{\hss{$4096_{z}$}}\cg%
\e{1}%
\e{7}%
\e{13}%
\e{13}%
\e{8}%
\e{18}%
\e{8}%
\e{6}%
\e{5}%
\e{6}%
\e{1}%
\e{1}%
\e{0}%
\e{0}%
\e{0}%
\eol}\vss}\rg%
%
%
\rx{\vss\hfull{%
\rlx{\hss{$5600_{z}$}}\cg%
\e{0}%
\e{4}%
\e{11}%
\e{14}%
\e{7}%
\e{24}%
\e{15}%
\e{9}%
\e{10}%
\e{14}%
\e{5}%
\e{3}%
\e{2}%
\e{0}%
\e{0}%
\eol}\vss}\rg%
%
%
\rx{\vss\hfull{%
\rlx{\hss{$448_{z}$}}\cg%
\e{0}%
\e{1}%
\e{3}%
\e{1}%
\e{1}%
\e{2}%
\e{0}%
\e{0}%
\e{1}%
\e{0}%
\e{0}%
\e{0}%
\e{0}%
\e{0}%
\e{0}%
\eol}\vss}\rg%
%
%
\rx{\vss\hfull{%
\rlx{\hss{$448_{w}$}}\cg%
\e{0}%
\e{0}%
\e{1}%
\e{2}%
\e{0}%
\e{2}%
\e{3}%
\e{0}%
\e{0}%
\e{1}%
\e{1}%
\e{0}%
\e{0}%
\e{0}%
\e{0}%
\eol}\vss}\rg%
%
%
\rx{\vss\hfull{%
\rlx{\hss{$1344_{w}$}}\cg%
\e{0}%
\e{0}%
\e{1}%
\e{0}%
\e{4}%
\e{4}%
\e{1}%
\e{4}%
\e{4}%
\e{5}%
\e{1}%
\e{4}%
\e{2}%
\e{0}%
\e{0}%
\eol}\vss}\rg%
%
%
\rx{\vss\hfull{%
\rlx{\hss{$5600_{w}$}}\cg%
\e{0}%
\e{2}%
\e{7}%
\e{9}%
\e{7}%
\e{22}%
\e{14}%
\e{11}%
\e{10}%
\e{20}%
\e{8}%
\e{5}%
\e{5}%
\e{1}%
\e{0}%
\eol}\vss}\rg%
%
%
\rx{\vss\hfull{%
\rlx{\hss{$2016_{w}$}}\cg%
\e{0}%
\e{0}%
\e{2}%
\e{1}%
\e{3}%
\e{6}%
\e{1}%
\e{7}%
\e{8}%
\e{7}%
\e{1}%
\e{5}%
\e{3}%
\e{1}%
\e{0}%
\eol}\vss}\rg%
%
%
\rx{\vss\hfull{%
\rlx{\hss{$7168_{w}$}}\cg%
\e{0}%
\e{2}%
\e{8}%
\e{8}%
\e{10}%
\e{26}%
\e{12}%
\e{18}%
\e{18}%
\e{26}%
\e{8}%
\e{10}%
\e{8}%
\e{2}%
\e{0}%
\eol}\vss}\rg%
\tableclose%
%
%
%
%
%
%
\eop
\eject
\tableopen{Induce/restrict matrix for $W({A_{4}}{A_{2}}{A_{1}})\,\subset\,W(E_{8})$}%
%
%
%
%
%
%
\rowpts=18 true pt%
\colpts=18 true pt%
\rowlabpts=40 true pt%
\collabpts=90 true pt%
\clx{\vss\hfull{%
\rlx{\hss{$ $}}\cg%
\cx{\hskip 16 true pt\flip{$[{5}]{\times}[{3}]{\times}[{2}]$}\hss}\cg%
\cx{\hskip 16 true pt\flip{$[{4}{1}]{\times}[{3}]{\times}[{2}]$}\hss}\cg%
\cx{\hskip 16 true pt\flip{$[{3}{2}]{\times}[{3}]{\times}[{2}]$}\hss}\cg%
\cx{\hskip 16 true pt\flip{$[{3}{1^{2}}]{\times}[{3}]{\times}[{2}]$}\hss}\cg%
\cx{\hskip 16 true pt\flip{$[{2^{2}}{1}]{\times}[{3}]{\times}[{2}]$}\hss}\cg%
\cx{\hskip 16 true pt\flip{$[{2}{1^{3}}]{\times}[{3}]{\times}[{2}]$}\hss}\cg%
\cx{\hskip 16 true pt\flip{$[{1^{5}}]{\times}[{3}]{\times}[{2}]$}\hss}\cg%
\cx{\hskip 16 true pt\flip{$[{5}]{\times}[{2}{1}]{\times}[{2}]$}\hss}\cg%
\cx{\hskip 16 true pt\flip{$[{4}{1}]{\times}[{2}{1}]{\times}[{2}]$}\hss}\cg%
\cx{\hskip 16 true pt\flip{$[{3}{2}]{\times}[{2}{1}]{\times}[{2}]$}\hss}\cg%
\cx{\hskip 16 true pt\flip{$[{3}{1^{2}}]{\times}[{2}{1}]{\times}[{2}]$}\hss}\cg%
\cx{\hskip 16 true pt\flip{$[{2^{2}}{1}]{\times}[{2}{1}]{\times}[{2}]$}\hss}\cg%
\cx{\hskip 16 true pt\flip{$[{2}{1^{3}}]{\times}[{2}{1}]{\times}[{2}]$}\hss}\cg%
\cx{\hskip 16 true pt\flip{$[{1^{5}}]{\times}[{2}{1}]{\times}[{2}]$}\hss}\cg%
\cx{\hskip 16 true pt\flip{$[{5}]{\times}[{1^{3}}]{\times}[{2}]$}\hss}\cg%
\cx{\hskip 16 true pt\flip{$[{4}{1}]{\times}[{1^{3}}]{\times}[{2}]$}\hss}\cg%
\cx{\hskip 16 true pt\flip{$[{3}{2}]{\times}[{1^{3}}]{\times}[{2}]$}\hss}\cg%
\cx{\hskip 16 true pt\flip{$[{3}{1^{2}}]{\times}[{1^{3}}]{\times}[{2}]$}\hss}\cg%
\cx{\hskip 16 true pt\flip{$[{2^{2}}{1}]{\times}[{1^{3}}]{\times}[{2}]$}\hss}\cg%
\cx{\hskip 16 true pt\flip{$[{2}{1^{3}}]{\times}[{1^{3}}]{\times}[{2}]$}\hss}\cg%
\cx{\hskip 16 true pt\flip{$[{1^{5}}]{\times}[{1^{3}}]{\times}[{2}]$}\hss}\cg%
\eol}}\rg%
%
%
\rx{\vss\hfull{%
\rlx{\hss{$1_{x}$}}\cg%
\e{1}%
\e{0}%
\e{0}%
\e{0}%
\e{0}%
\e{0}%
\e{0}%
\e{0}%
\e{0}%
\e{0}%
\e{0}%
\e{0}%
\e{0}%
\e{0}%
\e{0}%
\e{0}%
\e{0}%
\e{0}%
\e{0}%
\e{0}%
\e{0}%
\eol}\vss}\rg%
%
%
\rx{\vss\hfull{%
\rlx{\hss{$28_{x}$}}\cg%
\e{0}%
\e{1}%
\e{0}%
\e{1}%
\e{0}%
\e{0}%
\e{0}%
\e{1}%
\e{1}%
\e{0}%
\e{0}%
\e{0}%
\e{0}%
\e{0}%
\e{1}%
\e{0}%
\e{0}%
\e{0}%
\e{0}%
\e{0}%
\e{0}%
\eol}\vss}\rg%
%
%
\rx{\vss\hfull{%
\rlx{\hss{$35_{x}$}}\cg%
\e{3}%
\e{2}%
\e{1}%
\e{0}%
\e{0}%
\e{0}%
\e{0}%
\e{2}%
\e{1}%
\e{0}%
\e{0}%
\e{0}%
\e{0}%
\e{0}%
\e{0}%
\e{0}%
\e{0}%
\e{0}%
\e{0}%
\e{0}%
\e{0}%
\eol}\vss}\rg%
%
%
\rx{\vss\hfull{%
\rlx{\hss{$84_{x}$}}\cg%
\e{5}%
\e{4}%
\e{2}%
\e{0}%
\e{0}%
\e{0}%
\e{0}%
\e{3}%
\e{2}%
\e{1}%
\e{0}%
\e{0}%
\e{0}%
\e{0}%
\e{0}%
\e{0}%
\e{0}%
\e{0}%
\e{0}%
\e{0}%
\e{0}%
\eol}\vss}\rg%
%
%
\rx{\vss\hfull{%
\rlx{\hss{$50_{x}$}}\cg%
\e{2}%
\e{2}%
\e{1}%
\e{0}%
\e{0}%
\e{0}%
\e{0}%
\e{1}%
\e{1}%
\e{1}%
\e{0}%
\e{0}%
\e{0}%
\e{0}%
\e{0}%
\e{0}%
\e{0}%
\e{0}%
\e{0}%
\e{0}%
\e{0}%
\eol}\vss}\rg%
%
%
\rx{\vss\hfull{%
\rlx{\hss{$350_{x}$}}\cg%
\e{0}%
\e{2}%
\e{1}%
\e{4}%
\e{2}%
\e{2}%
\e{0}%
\e{2}%
\e{5}%
\e{2}%
\e{4}%
\e{1}%
\e{1}%
\e{0}%
\e{2}%
\e{3}%
\e{1}%
\e{1}%
\e{0}%
\e{0}%
\e{0}%
\eol}\vss}\rg%
%
%
\rx{\vss\hfull{%
\rlx{\hss{$300_{x}$}}\cg%
\e{4}%
\e{5}%
\e{5}%
\e{1}%
\e{2}%
\e{0}%
\e{0}%
\e{4}%
\e{6}%
\e{3}%
\e{2}%
\e{1}%
\e{0}%
\e{0}%
\e{0}%
\e{1}%
\e{0}%
\e{1}%
\e{0}%
\e{0}%
\e{0}%
\eol}\vss}\rg%
%
%
\rx{\vss\hfull{%
\rlx{\hss{$567_{x}$}}\cg%
\e{5}%
\e{11}%
\e{6}%
\e{7}%
\e{2}%
\e{1}%
\e{0}%
\e{9}%
\e{12}%
\e{5}%
\e{4}%
\e{1}%
\e{0}%
\e{0}%
\e{4}%
\e{3}%
\e{1}%
\e{0}%
\e{0}%
\e{0}%
\e{0}%
\eol}\vss}\rg%
%
%
\rx{\vss\hfull{%
\rlx{\hss{$210_{x}$}}\cg%
\e{4}%
\e{6}%
\e{3}%
\e{2}%
\e{1}%
\e{0}%
\e{0}%
\e{5}%
\e{5}%
\e{2}%
\e{1}%
\e{0}%
\e{0}%
\e{0}%
\e{1}%
\e{1}%
\e{0}%
\e{0}%
\e{0}%
\e{0}%
\e{0}%
\eol}\vss}\rg%
%
%
\rx{\vss\hfull{%
\rlx{\hss{$840_{x}$}}\cg%
\e{3}%
\e{6}%
\e{8}%
\e{2}%
\e{4}%
\e{0}%
\e{0}%
\e{2}%
\e{8}%
\e{10}%
\e{6}%
\e{7}%
\e{1}%
\e{0}%
\e{0}%
\e{1}%
\e{1}%
\e{3}%
\e{2}%
\e{2}%
\e{0}%
\eol}\vss}\rg%
%
%
\rx{\vss\hfull{%
\rlx{\hss{$700_{x}$}}\cg%
\e{10}%
\e{14}%
\e{10}%
\e{4}%
\e{3}%
\e{0}%
\e{0}%
\e{9}%
\e{13}%
\e{9}%
\e{4}%
\e{2}%
\e{0}%
\e{0}%
\e{1}%
\e{2}%
\e{1}%
\e{1}%
\e{0}%
\e{0}%
\e{0}%
\eol}\vss}\rg%
%
%
\rx{\vss\hfull{%
\rlx{\hss{$175_{x}$}}\cg%
\e{2}%
\e{3}%
\e{2}%
\e{1}%
\e{0}%
\e{0}%
\e{0}%
\e{1}%
\e{2}%
\e{3}%
\e{1}%
\e{1}%
\e{0}%
\e{0}%
\e{0}%
\e{0}%
\e{1}%
\e{0}%
\e{0}%
\e{0}%
\e{0}%
\eol}\vss}\rg%
%
%
\rx{\vss\hfull{%
\rlx{\hss{$1400_{x}$}}\cg%
\e{6}%
\e{15}%
\e{11}%
\e{12}%
\e{5}%
\e{3}%
\e{0}%
\e{9}%
\e{19}%
\e{15}%
\e{12}%
\e{6}%
\e{2}%
\e{0}%
\e{4}%
\e{6}%
\e{5}%
\e{2}%
\e{1}%
\e{0}%
\e{0}%
\eol}\vss}\rg%
%
%
\rx{\vss\hfull{%
\rlx{\hss{$1050_{x}$}}\cg%
\e{6}%
\e{13}%
\e{10}%
\e{8}%
\e{3}%
\e{1}%
\e{0}%
\e{7}%
\e{14}%
\e{13}%
\e{8}%
\e{5}%
\e{1}%
\e{0}%
\e{2}%
\e{3}%
\e{4}%
\e{1}%
\e{1}%
\e{0}%
\e{0}%
\eol}\vss}\rg%
%
%
\rx{\vss\hfull{%
\rlx{\hss{$1575_{x}$}}\cg%
\e{4}%
\e{14}%
\e{10}%
\e{15}%
\e{7}%
\e{5}%
\e{0}%
\e{10}%
\e{22}%
\e{14}%
\e{15}%
\e{6}%
\e{3}%
\e{0}%
\e{6}%
\e{9}%
\e{5}%
\e{3}%
\e{1}%
\e{0}%
\e{0}%
\eol}\vss}\rg%
%
%
\rx{\vss\hfull{%
\rlx{\hss{$1344_{x}$}}\cg%
\e{10}%
\e{19}%
\e{15}%
\e{10}%
\e{6}%
\e{1}%
\e{0}%
\e{13}%
\e{22}%
\e{15}%
\e{10}%
\e{5}%
\e{1}%
\e{0}%
\e{3}%
\e{5}%
\e{3}%
\e{2}%
\e{1}%
\e{0}%
\e{0}%
\eol}\vss}\rg%
%
%
\rx{\vss\hfull{%
\rlx{\hss{$2100_{x}$}}\cg%
\e{1}%
\e{7}%
\e{8}%
\e{14}%
\e{9}%
\e{10}%
\e{2}%
\e{5}%
\e{19}%
\e{14}%
\e{22}%
\e{11}%
\e{9}%
\e{1}%
\e{5}%
\e{12}%
\e{7}%
\e{9}%
\e{3}%
\e{2}%
\e{0}%
\eol}\vss}\rg%
%
%
\rx{\vss\hfull{%
\rlx{\hss{$2268_{x}$}}\cg%
\e{6}%
\e{17}%
\e{16}%
\e{16}%
\e{11}%
\e{6}%
\e{1}%
\e{11}%
\e{28}%
\e{21}%
\e{21}%
\e{11}%
\e{5}%
\e{0}%
\e{5}%
\e{11}%
\e{6}%
\e{7}%
\e{2}%
\e{1}%
\e{0}%
\eol}\vss}\rg%
%
%
\rx{\vss\hfull{%
\rlx{\hss{$525_{x}$}}\cg%
\e{1}%
\e{5}%
\e{3}%
\e{6}%
\e{1}%
\e{2}%
\e{0}%
\e{3}%
\e{7}%
\e{5}%
\e{5}%
\e{2}%
\e{1}%
\e{0}%
\e{3}%
\e{3}%
\e{2}%
\e{0}%
\e{1}%
\e{0}%
\e{0}%
\eol}\vss}\rg%
%
%
\rx{\vss\hfull{%
\rlx{\hss{$700_{xx}$}}\cg%
\e{1}%
\e{5}%
\e{4}%
\e{5}%
\e{2}%
\e{1}%
\e{0}%
\e{2}%
\e{6}%
\e{8}%
\e{6}%
\e{4}%
\e{2}%
\e{0}%
\e{1}%
\e{2}%
\e{3}%
\e{1}%
\e{2}%
\e{0}%
\e{0}%
\eol}\vss}\rg%
%
%
\rx{\vss\hfull{%
\rlx{\hss{$972_{x}$}}\cg%
\e{5}%
\e{10}%
\e{10}%
\e{4}%
\e{5}%
\e{0}%
\e{0}%
\e{5}%
\e{12}%
\e{12}%
\e{7}%
\e{6}%
\e{1}%
\e{0}%
\e{0}%
\e{2}%
\e{2}%
\e{3}%
\e{1}%
\e{1}%
\e{0}%
\eol}\vss}\rg%
%
%
\rx{\vss\hfull{%
\rlx{\hss{$4096_{x}$}}\cg%
\e{7}%
\e{24}%
\e{25}%
\e{26}%
\e{19}%
\e{10}%
\e{1}%
\e{15}%
\e{42}%
\e{37}%
\e{38}%
\e{23}%
\e{12}%
\e{1}%
\e{6}%
\e{17}%
\e{13}%
\e{13}%
\e{7}%
\e{3}%
\e{0}%
\eol}\vss}\rg%
%
%
\rx{\vss\hfull{%
\rlx{\hss{$4200_{x}$}}\cg%
\e{8}%
\e{24}%
\e{27}%
\e{23}%
\e{17}%
\e{8}%
\e{1}%
\e{12}%
\e{38}%
\e{41}%
\e{37}%
\e{27}%
\e{12}%
\e{1}%
\e{4}%
\e{13}%
\e{14}%
\e{13}%
\e{9}%
\e{4}%
\e{0}%
\eol}\vss}\rg%
%
%
\rx{\vss\hfull{%
\rlx{\hss{$2240_{x}$}}\cg%
\e{9}%
\e{20}%
\e{19}%
\e{13}%
\e{9}%
\e{2}%
\e{0}%
\e{10}%
\e{25}%
\e{26}%
\e{18}%
\e{13}%
\e{4}%
\e{0}%
\e{2}%
\e{6}%
\e{7}%
\e{5}%
\e{4}%
\e{1}%
\e{0}%
\eol}\vss}\rg%
%
%
\rx{\vss\hfull{%
\rlx{\hss{$2835_{x}$}}\cg%
\e{5}%
\e{15}%
\e{17}%
\e{14}%
\e{10}%
\e{4}%
\e{0}%
\e{6}%
\e{21}%
\e{29}%
\e{24}%
\e{20}%
\e{9}%
\e{1}%
\e{1}%
\e{6}%
\e{11}%
\e{8}%
\e{8}%
\e{3}%
\e{0}%
\eol}\vss}\rg%
%
%
\rx{\vss\hfull{%
\rlx{\hss{$6075_{x}$}}\cg%
\e{4}%
\e{23}%
\e{25}%
\e{36}%
\e{23}%
\e{19}%
\e{2}%
\e{12}%
\e{46}%
\e{49}%
\e{57}%
\e{38}%
\e{26}%
\e{3}%
\e{8}%
\e{23}%
\e{24}%
\e{21}%
\e{15}%
\e{7}%
\e{1}%
\eol}\vss}\rg%
%
%
\rx{\vss\hfull{%
\rlx{\hss{$3200_{x}$}}\cg%
\e{3}%
\e{12}%
\e{17}%
\e{13}%
\e{15}%
\e{5}%
\e{1}%
\e{6}%
\e{23}%
\e{28}%
\e{27}%
\e{24}%
\e{12}%
\e{2}%
\e{1}%
\e{8}%
\e{9}%
\e{13}%
\e{9}%
\e{7}%
\e{1}%
\eol}\vss}\rg%
\eop
\eject
\tablecont%
%
%
%
%
%
%
\rowpts=18 true pt%
\colpts=18 true pt%
\rowlabpts=40 true pt%
\collabpts=90 true pt%
\clx{\vss\hfull{%
\rlx{\hss{$ $}}\cg%
\cx{\hskip 16 true pt\flip{$[{5}]{\times}[{3}]{\times}[{2}]$}\hss}\cg%
\cx{\hskip 16 true pt\flip{$[{4}{1}]{\times}[{3}]{\times}[{2}]$}\hss}\cg%
\cx{\hskip 16 true pt\flip{$[{3}{2}]{\times}[{3}]{\times}[{2}]$}\hss}\cg%
\cx{\hskip 16 true pt\flip{$[{3}{1^{2}}]{\times}[{3}]{\times}[{2}]$}\hss}\cg%
\cx{\hskip 16 true pt\flip{$[{2^{2}}{1}]{\times}[{3}]{\times}[{2}]$}\hss}\cg%
\cx{\hskip 16 true pt\flip{$[{2}{1^{3}}]{\times}[{3}]{\times}[{2}]$}\hss}\cg%
\cx{\hskip 16 true pt\flip{$[{1^{5}}]{\times}[{3}]{\times}[{2}]$}\hss}\cg%
\cx{\hskip 16 true pt\flip{$[{5}]{\times}[{2}{1}]{\times}[{2}]$}\hss}\cg%
\cx{\hskip 16 true pt\flip{$[{4}{1}]{\times}[{2}{1}]{\times}[{2}]$}\hss}\cg%
\cx{\hskip 16 true pt\flip{$[{3}{2}]{\times}[{2}{1}]{\times}[{2}]$}\hss}\cg%
\cx{\hskip 16 true pt\flip{$[{3}{1^{2}}]{\times}[{2}{1}]{\times}[{2}]$}\hss}\cg%
\cx{\hskip 16 true pt\flip{$[{2^{2}}{1}]{\times}[{2}{1}]{\times}[{2}]$}\hss}\cg%
\cx{\hskip 16 true pt\flip{$[{2}{1^{3}}]{\times}[{2}{1}]{\times}[{2}]$}\hss}\cg%
\cx{\hskip 16 true pt\flip{$[{1^{5}}]{\times}[{2}{1}]{\times}[{2}]$}\hss}\cg%
\cx{\hskip 16 true pt\flip{$[{5}]{\times}[{1^{3}}]{\times}[{2}]$}\hss}\cg%
\cx{\hskip 16 true pt\flip{$[{4}{1}]{\times}[{1^{3}}]{\times}[{2}]$}\hss}\cg%
\cx{\hskip 16 true pt\flip{$[{3}{2}]{\times}[{1^{3}}]{\times}[{2}]$}\hss}\cg%
\cx{\hskip 16 true pt\flip{$[{3}{1^{2}}]{\times}[{1^{3}}]{\times}[{2}]$}\hss}\cg%
\cx{\hskip 16 true pt\flip{$[{2^{2}}{1}]{\times}[{1^{3}}]{\times}[{2}]$}\hss}\cg%
\cx{\hskip 16 true pt\flip{$[{2}{1^{3}}]{\times}[{1^{3}}]{\times}[{2}]$}\hss}\cg%
\cx{\hskip 16 true pt\flip{$[{1^{5}}]{\times}[{1^{3}}]{\times}[{2}]$}\hss}\cg%
\eol}}\rg%
%
%
\rx{\vss\hfull{%
\rlx{\hss{$70_{y}$}}\cg%
\e{0}%
\e{0}%
\e{0}%
\e{0}%
\e{0}%
\e{1}%
\e{1}%
\e{0}%
\e{0}%
\e{0}%
\e{1}%
\e{0}%
\e{1}%
\e{0}%
\e{0}%
\e{1}%
\e{0}%
\e{1}%
\e{0}%
\e{0}%
\e{0}%
\eol}\vss}\rg%
%
%
\rx{\vss\hfull{%
\rlx{\hss{$1134_{y}$}}\cg%
\e{0}%
\e{2}%
\e{2}%
\e{6}%
\e{4}%
\e{4}%
\e{0}%
\e{1}%
\e{6}%
\e{7}%
\e{10}%
\e{8}%
\e{7}%
\e{1}%
\e{1}%
\e{4}%
\e{5}%
\e{4}%
\e{4}%
\e{3}%
\e{1}%
\eol}\vss}\rg%
%
%
\rx{\vss\hfull{%
\rlx{\hss{$1680_{y}$}}\cg%
\e{0}%
\e{1}%
\e{2}%
\e{7}%
\e{6}%
\e{10}%
\e{4}%
\e{1}%
\e{8}%
\e{7}%
\e{17}%
\e{9}%
\e{13}%
\e{3}%
\e{2}%
\e{9}%
\e{6}%
\e{11}%
\e{4}%
\e{4}%
\e{0}%
\eol}\vss}\rg%
%
%
\rx{\vss\hfull{%
\rlx{\hss{$168_{y}$}}\cg%
\e{0}%
\e{1}%
\e{1}%
\e{0}%
\e{1}%
\e{0}%
\e{0}%
\e{0}%
\e{1}%
\e{2}%
\e{1}%
\e{2}%
\e{0}%
\e{0}%
\e{0}%
\e{0}%
\e{0}%
\e{1}%
\e{0}%
\e{1}%
\e{0}%
\eol}\vss}\rg%
%
%
\rx{\vss\hfull{%
\rlx{\hss{$420_{y}$}}\cg%
\e{1}%
\e{2}%
\e{3}%
\e{1}%
\e{1}%
\e{0}%
\e{0}%
\e{0}%
\e{2}%
\e{5}%
\e{3}%
\e{4}%
\e{1}%
\e{0}%
\e{0}%
\e{0}%
\e{1}%
\e{1}%
\e{2}%
\e{1}%
\e{0}%
\eol}\vss}\rg%
%
%
\rx{\vss\hfull{%
\rlx{\hss{$3150_{y}$}}\cg%
\e{2}%
\e{9}%
\e{13}%
\e{13}%
\e{10}%
\e{6}%
\e{0}%
\e{3}%
\e{16}%
\e{26}%
\e{26}%
\e{25}%
\e{15}%
\e{2}%
\e{1}%
\e{6}%
\e{12}%
\e{10}%
\e{13}%
\e{7}%
\e{2}%
\eol}\vss}\rg%
%
%
\rx{\vss\hfull{%
\rlx{\hss{$4200_{y}$}}\cg%
\e{2}%
\e{12}%
\e{17}%
\e{16}%
\e{16}%
\e{7}%
\e{1}%
\e{4}%
\e{22}%
\e{35}%
\e{34}%
\e{34}%
\e{19}%
\e{3}%
\e{1}%
\e{8}%
\e{14}%
\e{16}%
\e{15}%
\e{11}%
\e{2}%
\eol}\vss}\rg%
%
%
\rx{\vss\hfull{%
\rlx{\hss{$2688_{y}$}}\cg%
\e{1}%
\e{6}%
\e{10}%
\e{11}%
\e{10}%
\e{6}%
\e{1}%
\e{3}%
\e{14}%
\e{20}%
\e{23}%
\e{20}%
\e{14}%
\e{3}%
\e{1}%
\e{6}%
\e{10}%
\e{11}%
\e{10}%
\e{6}%
\e{1}%
\eol}\vss}\rg%
%
%
\rx{\vss\hfull{%
\rlx{\hss{$2100_{y}$}}\cg%
\e{0}%
\e{3}%
\e{4}%
\e{10}%
\e{6}%
\e{10}%
\e{3}%
\e{1}%
\e{10}%
\e{12}%
\e{20}%
\e{13}%
\e{14}%
\e{3}%
\e{3}%
\e{9}%
\e{8}%
\e{10}%
\e{7}%
\e{5}%
\e{1}%
\eol}\vss}\rg%
%
%
\rx{\vss\hfull{%
\rlx{\hss{$1400_{y}$}}\cg%
\e{0}%
\e{2}%
\e{3}%
\e{5}%
\e{6}%
\e{6}%
\e{3}%
\e{1}%
\e{7}%
\e{8}%
\e{13}%
\e{9}%
\e{9}%
\e{2}%
\e{1}%
\e{6}%
\e{4}%
\e{9}%
\e{3}%
\e{4}%
\e{0}%
\eol}\vss}\rg%
%
%
\rx{\vss\hfull{%
\rlx{\hss{$4536_{y}$}}\cg%
\e{1}%
\e{8}%
\e{13}%
\e{18}%
\e{17}%
\e{16}%
\e{5}%
\e{4}%
\e{23}%
\e{29}%
\e{41}%
\e{31}%
\e{27}%
\e{6}%
\e{3}%
\e{15}%
\e{16}%
\e{23}%
\e{14}%
\e{11}%
\e{1}%
\eol}\vss}\rg%
%
%
\rx{\vss\hfull{%
\rlx{\hss{$5670_{y}$}}\cg%
\e{1}%
\e{10}%
\e{15}%
\e{24}%
\e{21}%
\e{20}%
\e{5}%
\e{5}%
\e{29}%
\e{36}%
\e{51}%
\e{39}%
\e{34}%
\e{7}%
\e{4}%
\e{19}%
\e{21}%
\e{27}%
\e{18}%
\e{14}%
\e{2}%
\eol}\vss}\rg%
%
%
\rx{\vss\hfull{%
\rlx{\hss{$4480_{y}$}}\cg%
\e{2}%
\e{11}%
\e{16}%
\e{18}%
\e{16}%
\e{11}%
\e{2}%
\e{4}%
\e{23}%
\e{34}%
\e{38}%
\e{34}%
\e{23}%
\e{4}%
\e{2}%
\e{11}%
\e{16}%
\e{18}%
\e{16}%
\e{11}%
\e{2}%
\eol}\vss}\rg%
%
%
\rx{\vss\hfull{%
\rlx{\hss{$8_{z}$}}\cg%
\e{1}%
\e{1}%
\e{0}%
\e{0}%
\e{0}%
\e{0}%
\e{0}%
\e{1}%
\e{0}%
\e{0}%
\e{0}%
\e{0}%
\e{0}%
\e{0}%
\e{0}%
\e{0}%
\e{0}%
\e{0}%
\e{0}%
\e{0}%
\e{0}%
\eol}\vss}\rg%
%
%
\rx{\vss\hfull{%
\rlx{\hss{$56_{z}$}}\cg%
\e{0}%
\e{0}%
\e{0}%
\e{1}%
\e{0}%
\e{1}%
\e{0}%
\e{0}%
\e{1}%
\e{0}%
\e{1}%
\e{0}%
\e{0}%
\e{0}%
\e{1}%
\e{1}%
\e{0}%
\e{0}%
\e{0}%
\e{0}%
\e{0}%
\eol}\vss}\rg%
%
%
\rx{\vss\hfull{%
\rlx{\hss{$160_{z}$}}\cg%
\e{2}%
\e{4}%
\e{2}%
\e{2}%
\e{1}%
\e{0}%
\e{0}%
\e{4}%
\e{4}%
\e{1}%
\e{1}%
\e{0}%
\e{0}%
\e{0}%
\e{1}%
\e{1}%
\e{0}%
\e{0}%
\e{0}%
\e{0}%
\e{0}%
\eol}\vss}\rg%
%
%
\rx{\vss\hfull{%
\rlx{\hss{$112_{z}$}}\cg%
\e{5}%
\e{5}%
\e{2}%
\e{1}%
\e{0}%
\e{0}%
\e{0}%
\e{4}%
\e{3}%
\e{1}%
\e{0}%
\e{0}%
\e{0}%
\e{0}%
\e{1}%
\e{0}%
\e{0}%
\e{0}%
\e{0}%
\e{0}%
\e{0}%
\eol}\vss}\rg%
%
%
\rx{\vss\hfull{%
\rlx{\hss{$840_{z}$}}\cg%
\e{2}%
\e{6}%
\e{6}%
\e{5}%
\e{5}%
\e{1}%
\e{0}%
\e{4}%
\e{10}%
\e{8}%
\e{7}%
\e{5}%
\e{2}%
\e{0}%
\e{1}%
\e{3}%
\e{2}%
\e{3}%
\e{1}%
\e{1}%
\e{0}%
\eol}\vss}\rg%
%
%
\rx{\vss\hfull{%
\rlx{\hss{$1296_{z}$}}\cg%
\e{1}%
\e{7}%
\e{6}%
\e{12}%
\e{6}%
\e{6}%
\e{1}%
\e{6}%
\e{16}%
\e{9}%
\e{14}%
\e{5}%
\e{4}%
\e{0}%
\e{5}%
\e{9}%
\e{4}%
\e{4}%
\e{1}%
\e{0}%
\e{0}%
\eol}\vss}\rg%
%
%
\rx{\vss\hfull{%
\rlx{\hss{$1400_{z}$}}\cg%
\e{10}%
\e{19}%
\e{15}%
\e{11}%
\e{6}%
\e{2}%
\e{0}%
\e{13}%
\e{23}%
\e{15}%
\e{11}%
\e{5}%
\e{1}%
\e{0}%
\e{4}%
\e{6}%
\e{3}%
\e{2}%
\e{1}%
\e{0}%
\e{0}%
\eol}\vss}\rg%
%
%
\rx{\vss\hfull{%
\rlx{\hss{$1008_{z}$}}\cg%
\e{4}%
\e{12}%
\e{8}%
\e{10}%
\e{5}%
\e{3}%
\e{0}%
\e{9}%
\e{17}%
\e{9}%
\e{9}%
\e{3}%
\e{1}%
\e{0}%
\e{5}%
\e{6}%
\e{2}%
\e{2}%
\e{0}%
\e{0}%
\e{0}%
\eol}\vss}\rg%
%
%
\rx{\vss\hfull{%
\rlx{\hss{$560_{z}$}}\cg%
\e{9}%
\e{13}%
\e{8}%
\e{5}%
\e{2}%
\e{0}%
\e{0}%
\e{10}%
\e{12}%
\e{6}%
\e{3}%
\e{1}%
\e{0}%
\e{0}%
\e{2}%
\e{2}%
\e{1}%
\e{0}%
\e{0}%
\e{0}%
\e{0}%
\eol}\vss}\rg%
%
%
\rx{\vss\hfull{%
\rlx{\hss{$1400_{zz}$}}\cg%
\e{7}%
\e{14}%
\e{13}%
\e{8}%
\e{4}%
\e{1}%
\e{0}%
\e{6}%
\e{15}%
\e{18}%
\e{11}%
\e{8}%
\e{2}%
\e{0}%
\e{1}%
\e{3}%
\e{5}%
\e{2}%
\e{3}%
\e{0}%
\e{0}%
\eol}\vss}\rg%
%
%
\rx{\vss\hfull{%
\rlx{\hss{$4200_{z}$}}\cg%
\e{3}%
\e{17}%
\e{18}%
\e{23}%
\e{15}%
\e{10}%
\e{1}%
\e{7}%
\e{29}%
\e{37}%
\e{37}%
\e{29}%
\e{17}%
\e{2}%
\e{4}%
\e{12}%
\e{17}%
\e{14}%
\e{11}%
\e{6}%
\e{0}%
\eol}\vss}\rg%
%
%
\rx{\vss\hfull{%
\rlx{\hss{$400_{z}$}}\cg%
\e{6}%
\e{9}%
\e{5}%
\e{3}%
\e{1}%
\e{0}%
\e{0}%
\e{5}%
\e{7}%
\e{6}%
\e{2}%
\e{1}%
\e{0}%
\e{0}%
\e{1}%
\e{1}%
\e{1}%
\e{0}%
\e{0}%
\e{0}%
\e{0}%
\eol}\vss}\rg%
%
%
\rx{\vss\hfull{%
\rlx{\hss{$3240_{z}$}}\cg%
\e{13}%
\e{32}%
\e{27}%
\e{23}%
\e{14}%
\e{5}%
\e{0}%
\e{19}%
\e{42}%
\e{35}%
\e{27}%
\e{16}%
\e{5}%
\e{0}%
\e{6}%
\e{12}%
\e{10}%
\e{7}%
\e{3}%
\e{1}%
\e{0}%
\eol}\vss}\rg%
%
%
\rx{\vss\hfull{%
\rlx{\hss{$4536_{z}$}}\cg%
\e{9}%
\e{25}%
\e{30}%
\e{21}%
\e{19}%
\e{6}%
\e{0}%
\e{11}%
\e{38}%
\e{46}%
\e{38}%
\e{32}%
\e{13}%
\e{1}%
\e{2}%
\e{11}%
\e{14}%
\e{14}%
\e{12}%
\e{6}%
\e{1}%
\eol}\vss}\rg%
%
%
\rx{\vss\hfull{%
\rlx{\hss{$2400_{z}$}}\cg%
\e{0}%
\e{5}%
\e{5}%
\e{15}%
\e{8}%
\e{13}%
\e{3}%
\e{3}%
\e{16}%
\e{14}%
\e{25}%
\e{13}%
\e{14}%
\e{2}%
\e{5}%
\e{13}%
\e{10}%
\e{11}%
\e{5}%
\e{3}%
\e{0}%
\eol}\vss}\rg%
\eop
\eject
\tablecont%
%
%
%
%
%
%
\rowpts=18 true pt%
\colpts=18 true pt%
\rowlabpts=40 true pt%
\collabpts=90 true pt%
\clx{\vss\hfull{%
\rlx{\hss{$ $}}\cg%
\cx{\hskip 16 true pt\flip{$[{5}]{\times}[{3}]{\times}[{2}]$}\hss}\cg%
\cx{\hskip 16 true pt\flip{$[{4}{1}]{\times}[{3}]{\times}[{2}]$}\hss}\cg%
\cx{\hskip 16 true pt\flip{$[{3}{2}]{\times}[{3}]{\times}[{2}]$}\hss}\cg%
\cx{\hskip 16 true pt\flip{$[{3}{1^{2}}]{\times}[{3}]{\times}[{2}]$}\hss}\cg%
\cx{\hskip 16 true pt\flip{$[{2^{2}}{1}]{\times}[{3}]{\times}[{2}]$}\hss}\cg%
\cx{\hskip 16 true pt\flip{$[{2}{1^{3}}]{\times}[{3}]{\times}[{2}]$}\hss}\cg%
\cx{\hskip 16 true pt\flip{$[{1^{5}}]{\times}[{3}]{\times}[{2}]$}\hss}\cg%
\cx{\hskip 16 true pt\flip{$[{5}]{\times}[{2}{1}]{\times}[{2}]$}\hss}\cg%
\cx{\hskip 16 true pt\flip{$[{4}{1}]{\times}[{2}{1}]{\times}[{2}]$}\hss}\cg%
\cx{\hskip 16 true pt\flip{$[{3}{2}]{\times}[{2}{1}]{\times}[{2}]$}\hss}\cg%
\cx{\hskip 16 true pt\flip{$[{3}{1^{2}}]{\times}[{2}{1}]{\times}[{2}]$}\hss}\cg%
\cx{\hskip 16 true pt\flip{$[{2^{2}}{1}]{\times}[{2}{1}]{\times}[{2}]$}\hss}\cg%
\cx{\hskip 16 true pt\flip{$[{2}{1^{3}}]{\times}[{2}{1}]{\times}[{2}]$}\hss}\cg%
\cx{\hskip 16 true pt\flip{$[{1^{5}}]{\times}[{2}{1}]{\times}[{2}]$}\hss}\cg%
\cx{\hskip 16 true pt\flip{$[{5}]{\times}[{1^{3}}]{\times}[{2}]$}\hss}\cg%
\cx{\hskip 16 true pt\flip{$[{4}{1}]{\times}[{1^{3}}]{\times}[{2}]$}\hss}\cg%
\cx{\hskip 16 true pt\flip{$[{3}{2}]{\times}[{1^{3}}]{\times}[{2}]$}\hss}\cg%
\cx{\hskip 16 true pt\flip{$[{3}{1^{2}}]{\times}[{1^{3}}]{\times}[{2}]$}\hss}\cg%
\cx{\hskip 16 true pt\flip{$[{2^{2}}{1}]{\times}[{1^{3}}]{\times}[{2}]$}\hss}\cg%
\cx{\hskip 16 true pt\flip{$[{2}{1^{3}}]{\times}[{1^{3}}]{\times}[{2}]$}\hss}\cg%
\cx{\hskip 16 true pt\flip{$[{1^{5}}]{\times}[{1^{3}}]{\times}[{2}]$}\hss}\cg%
\eol}}\rg%
%
%
\rx{\vss\hfull{%
\rlx{\hss{$3360_{z}$}}\cg%
\e{5}%
\e{18}%
\e{19}%
\e{21}%
\e{13}%
\e{8}%
\e{1}%
\e{10}%
\e{30}%
\e{31}%
\e{31}%
\e{20}%
\e{11}%
\e{1}%
\e{4}%
\e{12}%
\e{13}%
\e{10}%
\e{7}%
\e{2}%
\e{0}%
\eol}\vss}\rg%
%
%
\rx{\vss\hfull{%
\rlx{\hss{$2800_{z}$}}\cg%
\e{2}%
\e{12}%
\e{12}%
\e{19}%
\e{11}%
\e{11}%
\e{2}%
\e{7}%
\e{25}%
\e{22}%
\e{28}%
\e{15}%
\e{11}%
\e{1}%
\e{6}%
\e{14}%
\e{10}%
\e{10}%
\e{5}%
\e{2}%
\e{0}%
\eol}\vss}\rg%
%
%
\rx{\vss\hfull{%
\rlx{\hss{$4096_{z}$}}\cg%
\e{7}%
\e{24}%
\e{25}%
\e{26}%
\e{19}%
\e{10}%
\e{1}%
\e{15}%
\e{42}%
\e{37}%
\e{38}%
\e{23}%
\e{12}%
\e{1}%
\e{6}%
\e{17}%
\e{13}%
\e{13}%
\e{7}%
\e{3}%
\e{0}%
\eol}\vss}\rg%
%
%
\rx{\vss\hfull{%
\rlx{\hss{$5600_{z}$}}\cg%
\e{3}%
\e{17}%
\e{22}%
\e{28}%
\e{23}%
\e{18}%
\e{4}%
\e{9}%
\e{39}%
\e{42}%
\e{52}%
\e{37}%
\e{26}%
\e{4}%
\e{6}%
\e{21}%
\e{19}%
\e{24}%
\e{14}%
\e{10}%
\e{1}%
\eol}\vss}\rg%
%
%
\rx{\vss\hfull{%
\rlx{\hss{$448_{z}$}}\cg%
\e{6}%
\e{7}%
\e{7}%
\e{2}%
\e{1}%
\e{0}%
\e{0}%
\e{4}%
\e{7}%
\e{6}%
\e{3}%
\e{2}%
\e{0}%
\e{0}%
\e{0}%
\e{1}%
\e{1}%
\e{0}%
\e{1}%
\e{0}%
\e{0}%
\eol}\vss}\rg%
%
%
\rx{\vss\hfull{%
\rlx{\hss{$448_{w}$}}\cg%
\e{0}%
\e{0}%
\e{0}%
\e{2}%
\e{1}%
\e{4}%
\e{1}%
\e{0}%
\e{2}%
\e{1}%
\e{5}%
\e{2}%
\e{4}%
\e{1}%
\e{1}%
\e{3}%
\e{2}%
\e{3}%
\e{1}%
\e{1}%
\e{0}%
\eol}\vss}\rg%
%
%
\rx{\vss\hfull{%
\rlx{\hss{$1344_{w}$}}\cg%
\e{1}%
\e{5}%
\e{6}%
\e{5}%
\e{5}%
\e{1}%
\e{0}%
\e{1}%
\e{7}%
\e{13}%
\e{10}%
\e{12}%
\e{5}%
\e{0}%
\e{0}%
\e{2}%
\e{4}%
\e{4}%
\e{5}%
\e{4}%
\e{1}%
\eol}\vss}\rg%
%
%
\rx{\vss\hfull{%
\rlx{\hss{$5600_{w}$}}\cg%
\e{1}%
\e{10}%
\e{15}%
\e{24}%
\e{21}%
\e{19}%
\e{4}%
\e{5}%
\e{29}%
\e{36}%
\e{50}%
\e{39}%
\e{33}%
\e{7}%
\e{4}%
\e{18}%
\e{21}%
\e{26}%
\e{18}%
\e{14}%
\e{2}%
\eol}\vss}\rg%
%
%
\rx{\vss\hfull{%
\rlx{\hss{$2016_{w}$}}\cg%
\e{2}%
\e{7}%
\e{11}%
\e{7}%
\e{6}%
\e{2}%
\e{0}%
\e{2}%
\e{10}%
\e{19}%
\e{16}%
\e{17}%
\e{8}%
\e{1}%
\e{0}%
\e{2}%
\e{7}%
\e{6}%
\e{9}%
\e{4}%
\e{1}%
\eol}\vss}\rg%
%
%
\rx{\vss\hfull{%
\rlx{\hss{$7168_{w}$}}\cg%
\e{3}%
\e{17}%
\e{26}%
\e{29}%
\e{26}%
\e{17}%
\e{3}%
\e{7}%
\e{37}%
\e{54}%
\e{61}%
\e{54}%
\e{37}%
\e{7}%
\e{3}%
\e{17}%
\e{26}%
\e{29}%
\e{26}%
\e{17}%
\e{3}%
\eol}\vss}\rg%
%
%
%
%
%
%
\rowpts=18 true pt%
\colpts=18 true pt%
\rowlabpts=40 true pt%
\collabpts=90 true pt%
\clx{\vss\hfull{%
\rlx{\hss{$ $}}\cg%
\cx{\hskip 16 true pt\flip{$[{5}]{\times}[{3}]{\times}[{1^{2}}]$}\hss}\cg%
\cx{\hskip 16 true pt\flip{$[{4}{1}]{\times}[{3}]{\times}[{1^{2}}]$}\hss}\cg%
\cx{\hskip 16 true pt\flip{$[{3}{2}]{\times}[{3}]{\times}[{1^{2}}]$}\hss}\cg%
\cx{\hskip 16 true pt\flip{$[{3}{1^{2}}]{\times}[{3}]{\times}[{1^{2}}]$}\hss}\cg%
\cx{\hskip 16 true pt\flip{$[{2^{2}}{1}]{\times}[{3}]{\times}[{1^{2}}]$}\hss}\cg%
\cx{\hskip 16 true pt\flip{$[{2}{1^{3}}]{\times}[{3}]{\times}[{1^{2}}]$}\hss}\cg%
\cx{\hskip 16 true pt\flip{$[{1^{5}}]{\times}[{3}]{\times}[{1^{2}}]$}\hss}\cg%
\cx{\hskip 16 true pt\flip{$[{5}]{\times}[{2}{1}]{\times}[{1^{2}}]$}\hss}\cg%
\cx{\hskip 16 true pt\flip{$[{4}{1}]{\times}[{2}{1}]{\times}[{1^{2}}]$}\hss}\cg%
\cx{\hskip 16 true pt\flip{$[{3}{2}]{\times}[{2}{1}]{\times}[{1^{2}}]$}\hss}\cg%
\cx{\hskip 16 true pt\flip{$[{3}{1^{2}}]{\times}[{2}{1}]{\times}[{1^{2}}]$}\hss}\cg%
\cx{\hskip 16 true pt\flip{$[{2^{2}}{1}]{\times}[{2}{1}]{\times}[{1^{2}}]$}\hss}\cg%
\cx{\hskip 16 true pt\flip{$[{2}{1^{3}}]{\times}[{2}{1}]{\times}[{1^{2}}]$}\hss}\cg%
\cx{\hskip 16 true pt\flip{$[{1^{5}}]{\times}[{2}{1}]{\times}[{1^{2}}]$}\hss}\cg%
\cx{\hskip 16 true pt\flip{$[{5}]{\times}[{1^{3}}]{\times}[{1^{2}}]$}\hss}\cg%
\cx{\hskip 16 true pt\flip{$[{4}{1}]{\times}[{1^{3}}]{\times}[{1^{2}}]$}\hss}\cg%
\cx{\hskip 16 true pt\flip{$[{3}{2}]{\times}[{1^{3}}]{\times}[{1^{2}}]$}\hss}\cg%
\cx{\hskip 16 true pt\flip{$[{3}{1^{2}}]{\times}[{1^{3}}]{\times}[{1^{2}}]$}\hss}\cg%
\cx{\hskip 16 true pt\flip{$[{2^{2}}{1}]{\times}[{1^{3}}]{\times}[{1^{2}}]$}\hss}\cg%
\cx{\hskip 16 true pt\flip{$[{2}{1^{3}}]{\times}[{1^{3}}]{\times}[{1^{2}}]$}\hss}\cg%
\cx{\hskip 16 true pt\flip{$[{1^{5}}]{\times}[{1^{3}}]{\times}[{1^{2}}]$}\hss}\cg%
\eol}}\rg%
%
%
\rx{\vss\hfull{%
\rlx{\hss{$1_{x}$}}\cg%
\e{0}%
\e{0}%
\e{0}%
\e{0}%
\e{0}%
\e{0}%
\e{0}%
\e{0}%
\e{0}%
\e{0}%
\e{0}%
\e{0}%
\e{0}%
\e{0}%
\e{0}%
\e{0}%
\e{0}%
\e{0}%
\e{0}%
\e{0}%
\e{0}%
\eol}\vss}\rg%
%
%
\rx{\vss\hfull{%
\rlx{\hss{$28_{x}$}}\cg%
\e{1}%
\e{1}%
\e{0}%
\e{0}%
\e{0}%
\e{0}%
\e{0}%
\e{1}%
\e{0}%
\e{0}%
\e{0}%
\e{0}%
\e{0}%
\e{0}%
\e{0}%
\e{0}%
\e{0}%
\e{0}%
\e{0}%
\e{0}%
\e{0}%
\eol}\vss}\rg%
%
%
\rx{\vss\hfull{%
\rlx{\hss{$35_{x}$}}\cg%
\e{1}%
\e{1}%
\e{0}%
\e{0}%
\e{0}%
\e{0}%
\e{0}%
\e{1}%
\e{0}%
\e{0}%
\e{0}%
\e{0}%
\e{0}%
\e{0}%
\e{0}%
\e{0}%
\e{0}%
\e{0}%
\e{0}%
\e{0}%
\e{0}%
\eol}\vss}\rg%
%
%
\rx{\vss\hfull{%
\rlx{\hss{$84_{x}$}}\cg%
\e{2}%
\e{1}%
\e{1}%
\e{0}%
\e{0}%
\e{0}%
\e{0}%
\e{1}%
\e{1}%
\e{0}%
\e{0}%
\e{0}%
\e{0}%
\e{0}%
\e{0}%
\e{0}%
\e{0}%
\e{0}%
\e{0}%
\e{0}%
\e{0}%
\eol}\vss}\rg%
%
%
\rx{\vss\hfull{%
\rlx{\hss{$50_{x}$}}\cg%
\e{1}%
\e{1}%
\e{0}%
\e{0}%
\e{0}%
\e{0}%
\e{0}%
\e{0}%
\e{0}%
\e{1}%
\e{0}%
\e{0}%
\e{0}%
\e{0}%
\e{0}%
\e{0}%
\e{0}%
\e{0}%
\e{0}%
\e{0}%
\e{0}%
\eol}\vss}\rg%
%
%
\rx{\vss\hfull{%
\rlx{\hss{$350_{x}$}}\cg%
\e{1}%
\e{3}%
\e{2}%
\e{3}%
\e{1}%
\e{1}%
\e{0}%
\e{3}%
\e{5}%
\e{1}%
\e{2}%
\e{0}%
\e{0}%
\e{0}%
\e{2}%
\e{2}%
\e{0}%
\e{0}%
\e{0}%
\e{0}%
\e{0}%
\eol}\vss}\rg%
%
%
\rx{\vss\hfull{%
\rlx{\hss{$300_{x}$}}\cg%
\e{1}%
\e{3}%
\e{2}%
\e{2}%
\e{1}%
\e{0}%
\e{0}%
\e{3}%
\e{4}%
\e{1}%
\e{1}%
\e{0}%
\e{0}%
\e{0}%
\e{1}%
\e{1}%
\e{0}%
\e{0}%
\e{0}%
\e{0}%
\e{0}%
\eol}\vss}\rg%
%
%
\rx{\vss\hfull{%
\rlx{\hss{$567_{x}$}}\cg%
\e{6}%
\e{8}%
\e{4}%
\e{3}%
\e{1}%
\e{0}%
\e{0}%
\e{7}%
\e{7}%
\e{2}%
\e{1}%
\e{0}%
\e{0}%
\e{0}%
\e{2}%
\e{1}%
\e{0}%
\e{0}%
\e{0}%
\e{0}%
\e{0}%
\eol}\vss}\rg%
%
%
\rx{\vss\hfull{%
\rlx{\hss{$210_{x}$}}\cg%
\e{3}%
\e{4}%
\e{1}%
\e{1}%
\e{0}%
\e{0}%
\e{0}%
\e{3}%
\e{2}%
\e{1}%
\e{0}%
\e{0}%
\e{0}%
\e{0}%
\e{1}%
\e{0}%
\e{0}%
\e{0}%
\e{0}%
\e{0}%
\e{0}%
\eol}\vss}\rg%
%
%
\rx{\vss\hfull{%
\rlx{\hss{$840_{x}$}}\cg%
\e{1}%
\e{3}%
\e{4}%
\e{3}%
\e{4}%
\e{1}%
\e{0}%
\e{1}%
\e{6}%
\e{7}%
\e{6}%
\e{5}%
\e{2}%
\e{0}%
\e{0}%
\e{2}%
\e{2}%
\e{3}%
\e{1}%
\e{1}%
\e{0}%
\eol}\vss}\rg%
%
%
\rx{\vss\hfull{%
\rlx{\hss{$700_{x}$}}\cg%
\e{5}%
\e{8}%
\e{5}%
\e{3}%
\e{1}%
\e{0}%
\e{0}%
\e{5}%
\e{7}%
\e{5}%
\e{2}%
\e{1}%
\e{0}%
\e{0}%
\e{1}%
\e{1}%
\e{1}%
\e{0}%
\e{0}%
\e{0}%
\e{0}%
\eol}\vss}\rg%
%
%
\rx{\vss\hfull{%
\rlx{\hss{$175_{x}$}}\cg%
\e{1}%
\e{1}%
\e{2}%
\e{0}%
\e{0}%
\e{0}%
\e{0}%
\e{0}%
\e{1}%
\e{2}%
\e{1}%
\e{1}%
\e{0}%
\e{0}%
\e{0}%
\e{0}%
\e{0}%
\e{0}%
\e{1}%
\e{0}%
\e{0}%
\eol}\vss}\rg%
%
%
\rx{\vss\hfull{%
\rlx{\hss{$1400_{x}$}}\cg%
\e{6}%
\e{12}%
\e{9}%
\e{7}%
\e{3}%
\e{1}%
\e{0}%
\e{7}%
\e{13}%
\e{11}%
\e{7}%
\e{4}%
\e{1}%
\e{0}%
\e{2}%
\e{3}%
\e{3}%
\e{1}%
\e{1}%
\e{0}%
\e{0}%
\eol}\vss}\rg%
%
%
\rx{\vss\hfull{%
\rlx{\hss{$1050_{x}$}}\cg%
\e{5}%
\e{9}%
\e{7}%
\e{4}%
\e{2}%
\e{0}%
\e{0}%
\e{4}%
\e{8}%
\e{10}%
\e{5}%
\e{4}%
\e{1}%
\e{0}%
\e{0}%
\e{1}%
\e{3}%
\e{1}%
\e{1}%
\e{0}%
\e{0}%
\eol}\vss}\rg%
\eop
\eject
\tablecont%
%
%
%
%
%
%
\rowpts=18 true pt%
\colpts=18 true pt%
\rowlabpts=40 true pt%
\collabpts=90 true pt%
\clx{\vss\hfull{%
\rlx{\hss{$ $}}\cg%
\cx{\hskip 16 true pt\flip{$[{5}]{\times}[{3}]{\times}[{1^{2}}]$}\hss}\cg%
\cx{\hskip 16 true pt\flip{$[{4}{1}]{\times}[{3}]{\times}[{1^{2}}]$}\hss}\cg%
\cx{\hskip 16 true pt\flip{$[{3}{2}]{\times}[{3}]{\times}[{1^{2}}]$}\hss}\cg%
\cx{\hskip 16 true pt\flip{$[{3}{1^{2}}]{\times}[{3}]{\times}[{1^{2}}]$}\hss}\cg%
\cx{\hskip 16 true pt\flip{$[{2^{2}}{1}]{\times}[{3}]{\times}[{1^{2}}]$}\hss}\cg%
\cx{\hskip 16 true pt\flip{$[{2}{1^{3}}]{\times}[{3}]{\times}[{1^{2}}]$}\hss}\cg%
\cx{\hskip 16 true pt\flip{$[{1^{5}}]{\times}[{3}]{\times}[{1^{2}}]$}\hss}\cg%
\cx{\hskip 16 true pt\flip{$[{5}]{\times}[{2}{1}]{\times}[{1^{2}}]$}\hss}\cg%
\cx{\hskip 16 true pt\flip{$[{4}{1}]{\times}[{2}{1}]{\times}[{1^{2}}]$}\hss}\cg%
\cx{\hskip 16 true pt\flip{$[{3}{2}]{\times}[{2}{1}]{\times}[{1^{2}}]$}\hss}\cg%
\cx{\hskip 16 true pt\flip{$[{3}{1^{2}}]{\times}[{2}{1}]{\times}[{1^{2}}]$}\hss}\cg%
\cx{\hskip 16 true pt\flip{$[{2^{2}}{1}]{\times}[{2}{1}]{\times}[{1^{2}}]$}\hss}\cg%
\cx{\hskip 16 true pt\flip{$[{2}{1^{3}}]{\times}[{2}{1}]{\times}[{1^{2}}]$}\hss}\cg%
\cx{\hskip 16 true pt\flip{$[{1^{5}}]{\times}[{2}{1}]{\times}[{1^{2}}]$}\hss}\cg%
\cx{\hskip 16 true pt\flip{$[{5}]{\times}[{1^{3}}]{\times}[{1^{2}}]$}\hss}\cg%
\cx{\hskip 16 true pt\flip{$[{4}{1}]{\times}[{1^{3}}]{\times}[{1^{2}}]$}\hss}\cg%
\cx{\hskip 16 true pt\flip{$[{3}{2}]{\times}[{1^{3}}]{\times}[{1^{2}}]$}\hss}\cg%
\cx{\hskip 16 true pt\flip{$[{3}{1^{2}}]{\times}[{1^{3}}]{\times}[{1^{2}}]$}\hss}\cg%
\cx{\hskip 16 true pt\flip{$[{2^{2}}{1}]{\times}[{1^{3}}]{\times}[{1^{2}}]$}\hss}\cg%
\cx{\hskip 16 true pt\flip{$[{2}{1^{3}}]{\times}[{1^{3}}]{\times}[{1^{2}}]$}\hss}\cg%
\cx{\hskip 16 true pt\flip{$[{1^{5}}]{\times}[{1^{3}}]{\times}[{1^{2}}]$}\hss}\cg%
\eol}}\rg%
%
%
\rx{\vss\hfull{%
\rlx{\hss{$1575_{x}$}}\cg%
\e{6}%
\e{14}%
\e{9}%
\e{10}%
\e{4}%
\e{2}%
\e{0}%
\e{10}%
\e{17}%
\e{10}%
\e{8}%
\e{3}%
\e{1}%
\e{0}%
\e{4}%
\e{5}%
\e{3}%
\e{1}%
\e{0}%
\e{0}%
\e{0}%
\eol}\vss}\rg%
%
%
\rx{\vss\hfull{%
\rlx{\hss{$1344_{x}$}}\cg%
\e{7}%
\e{12}%
\e{10}%
\e{6}%
\e{4}%
\e{1}%
\e{0}%
\e{9}%
\e{15}%
\e{8}%
\e{6}%
\e{2}%
\e{0}%
\e{0}%
\e{2}%
\e{4}%
\e{1}%
\e{1}%
\e{0}%
\e{0}%
\e{0}%
\eol}\vss}\rg%
%
%
\rx{\vss\hfull{%
\rlx{\hss{$2100_{x}$}}\cg%
\e{2}%
\e{10}%
\e{9}%
\e{14}%
\e{8}%
\e{7}%
\e{1}%
\e{8}%
\e{21}%
\e{12}%
\e{17}%
\e{6}%
\e{4}%
\e{0}%
\e{6}%
\e{11}%
\e{4}%
\e{5}%
\e{1}%
\e{0}%
\e{0}%
\eol}\vss}\rg%
%
%
\rx{\vss\hfull{%
\rlx{\hss{$2268_{x}$}}\cg%
\e{5}%
\e{16}%
\e{12}%
\e{14}%
\e{7}%
\e{4}%
\e{0}%
\e{11}%
\e{23}%
\e{15}%
\e{14}%
\e{6}%
\e{2}%
\e{0}%
\e{6}%
\e{8}%
\e{4}%
\e{3}%
\e{1}%
\e{0}%
\e{0}%
\eol}\vss}\rg%
%
%
\rx{\vss\hfull{%
\rlx{\hss{$525_{x}$}}\cg%
\e{3}%
\e{4}%
\e{4}%
\e{2}%
\e{2}%
\e{1}%
\e{0}%
\e{3}%
\e{6}%
\e{3}%
\e{3}%
\e{1}%
\e{0}%
\e{0}%
\e{1}%
\e{2}%
\e{0}%
\e{1}%
\e{0}%
\e{0}%
\e{0}%
\eol}\vss}\rg%
%
%
\rx{\vss\hfull{%
\rlx{\hss{$700_{xx}$}}\cg%
\e{2}%
\e{5}%
\e{3}%
\e{3}%
\e{1}%
\e{0}%
\e{0}%
\e{1}%
\e{3}%
\e{8}%
\e{4}%
\e{5}%
\e{2}%
\e{0}%
\e{0}%
\e{0}%
\e{3}%
\e{1}%
\e{2}%
\e{1}%
\e{0}%
\eol}\vss}\rg%
%
%
\rx{\vss\hfull{%
\rlx{\hss{$972_{x}$}}\cg%
\e{2}%
\e{6}%
\e{6}%
\e{4}%
\e{3}%
\e{1}%
\e{0}%
\e{3}%
\e{8}%
\e{8}%
\e{6}%
\e{4}%
\e{1}%
\e{0}%
\e{1}%
\e{2}%
\e{2}%
\e{2}%
\e{1}%
\e{0}%
\e{0}%
\eol}\vss}\rg%
%
%
\rx{\vss\hfull{%
\rlx{\hss{$4096_{x}$}}\cg%
\e{7}%
\e{22}%
\e{21}%
\e{22}%
\e{15}%
\e{8}%
\e{1}%
\e{14}%
\e{37}%
\e{29}%
\e{29}%
\e{15}%
\e{7}%
\e{0}%
\e{6}%
\e{15}%
\e{9}%
\e{9}%
\e{3}%
\e{1}%
\e{0}%
\eol}\vss}\rg%
%
%
\rx{\vss\hfull{%
\rlx{\hss{$4200_{x}$}}\cg%
\e{6}%
\e{20}%
\e{21}%
\e{20}%
\e{14}%
\e{6}%
\e{0}%
\e{10}%
\e{31}%
\e{34}%
\e{30}%
\e{22}%
\e{10}%
\e{1}%
\e{3}%
\e{10}%
\e{13}%
\e{10}%
\e{7}%
\e{3}%
\e{0}%
\eol}\vss}\rg%
%
%
\rx{\vss\hfull{%
\rlx{\hss{$2240_{x}$}}\cg%
\e{6}%
\e{13}%
\e{14}%
\e{9}%
\e{7}%
\e{2}%
\e{0}%
\e{6}%
\e{18}%
\e{19}%
\e{14}%
\e{10}%
\e{3}%
\e{0}%
\e{1}%
\e{5}%
\e{5}%
\e{4}%
\e{3}%
\e{1}%
\e{0}%
\eol}\vss}\rg%
%
%
\rx{\vss\hfull{%
\rlx{\hss{$2835_{x}$}}\cg%
\e{3}%
\e{11}%
\e{15}%
\e{11}%
\e{8}%
\e{3}%
\e{0}%
\e{4}%
\e{16}%
\e{25}%
\e{21}%
\e{20}%
\e{9}%
\e{1}%
\e{1}%
\e{4}%
\e{9}%
\e{7}%
\e{10}%
\e{4}%
\e{1}%
\eol}\vss}\rg%
%
%
\rx{\vss\hfull{%
\rlx{\hss{$6075_{x}$}}\cg%
\e{7}%
\e{25}%
\e{27}%
\e{30}%
\e{21}%
\e{14}%
\e{2}%
\e{13}%
\e{44}%
\e{46}%
\e{48}%
\e{32}%
\e{19}%
\e{2}%
\e{6}%
\e{19}%
\e{19}%
\e{18}%
\e{11}%
\e{5}%
\e{0}%
\eol}\vss}\rg%
%
%
\rx{\vss\hfull{%
\rlx{\hss{$3200_{x}$}}\cg%
\e{2}%
\e{9}%
\e{14}%
\e{13}%
\e{14}%
\e{8}%
\e{2}%
\e{5}%
\e{22}%
\e{23}%
\e{27}%
\e{19}%
\e{11}%
\e{1}%
\e{2}%
\e{11}%
\e{8}%
\e{13}%
\e{6}%
\e{4}%
\e{0}%
\eol}\vss}\rg%
%
%
\rx{\vss\hfull{%
\rlx{\hss{$70_{y}$}}\cg%
\e{0}%
\e{0}%
\e{0}%
\e{1}%
\e{0}%
\e{1}%
\e{0}%
\e{0}%
\e{1}%
\e{0}%
\e{1}%
\e{0}%
\e{0}%
\e{0}%
\e{1}%
\e{1}%
\e{0}%
\e{0}%
\e{0}%
\e{0}%
\e{0}%
\eol}\vss}\rg%
%
%
\rx{\vss\hfull{%
\rlx{\hss{$1134_{y}$}}\cg%
\e{1}%
\e{3}%
\e{4}%
\e{4}%
\e{5}%
\e{4}%
\e{1}%
\e{1}%
\e{7}%
\e{8}%
\e{10}%
\e{7}%
\e{6}%
\e{1}%
\e{0}%
\e{4}%
\e{4}%
\e{6}%
\e{2}%
\e{2}%
\e{0}%
\eol}\vss}\rg%
%
%
\rx{\vss\hfull{%
\rlx{\hss{$1680_{y}$}}\cg%
\e{0}%
\e{4}%
\e{4}%
\e{11}%
\e{6}%
\e{9}%
\e{2}%
\e{3}%
\e{13}%
\e{9}%
\e{17}%
\e{7}%
\e{8}%
\e{1}%
\e{4}%
\e{10}%
\e{6}%
\e{7}%
\e{2}%
\e{1}%
\e{0}%
\eol}\vss}\rg%
%
%
\rx{\vss\hfull{%
\rlx{\hss{$168_{y}$}}\cg%
\e{0}%
\e{1}%
\e{0}%
\e{1}%
\e{0}%
\e{0}%
\e{0}%
\e{0}%
\e{0}%
\e{2}%
\e{1}%
\e{2}%
\e{1}%
\e{0}%
\e{0}%
\e{0}%
\e{1}%
\e{0}%
\e{1}%
\e{1}%
\e{0}%
\eol}\vss}\rg%
%
%
\rx{\vss\hfull{%
\rlx{\hss{$420_{y}$}}\cg%
\e{0}%
\e{1}%
\e{2}%
\e{1}%
\e{1}%
\e{0}%
\e{0}%
\e{0}%
\e{1}%
\e{4}%
\e{3}%
\e{5}%
\e{2}%
\e{0}%
\e{0}%
\e{0}%
\e{1}%
\e{1}%
\e{3}%
\e{2}%
\e{1}%
\eol}\vss}\rg%
%
%
\rx{\vss\hfull{%
\rlx{\hss{$3150_{y}$}}\cg%
\e{2}%
\e{7}%
\e{13}%
\e{10}%
\e{12}%
\e{6}%
\e{1}%
\e{2}%
\e{15}%
\e{25}%
\e{26}%
\e{26}%
\e{16}%
\e{3}%
\e{0}%
\e{6}%
\e{10}%
\e{13}%
\e{13}%
\e{9}%
\e{2}%
\eol}\vss}\rg%
%
%
\rx{\vss\hfull{%
\rlx{\hss{$4200_{y}$}}\cg%
\e{2}%
\e{11}%
\e{15}%
\e{16}%
\e{14}%
\e{8}%
\e{1}%
\e{3}%
\e{19}%
\e{34}%
\e{34}%
\e{35}%
\e{22}%
\e{4}%
\e{1}%
\e{7}%
\e{16}%
\e{16}%
\e{17}%
\e{12}%
\e{2}%
\eol}\vss}\rg%
%
%
\rx{\vss\hfull{%
\rlx{\hss{$2688_{y}$}}\cg%
\e{1}%
\e{6}%
\e{10}%
\e{11}%
\e{10}%
\e{6}%
\e{1}%
\e{3}%
\e{14}%
\e{20}%
\e{23}%
\e{20}%
\e{14}%
\e{3}%
\e{1}%
\e{6}%
\e{10}%
\e{11}%
\e{10}%
\e{6}%
\e{1}%
\eol}\vss}\rg%
%
%
\rx{\vss\hfull{%
\rlx{\hss{$2100_{y}$}}\cg%
\e{1}%
\e{5}%
\e{7}%
\e{10}%
\e{8}%
\e{9}%
\e{3}%
\e{3}%
\e{14}%
\e{13}%
\e{20}%
\e{12}%
\e{10}%
\e{1}%
\e{3}%
\e{10}%
\e{6}%
\e{10}%
\e{4}%
\e{3}%
\e{0}%
\eol}\vss}\rg%
%
%
\rx{\vss\hfull{%
\rlx{\hss{$1400_{y}$}}\cg%
\e{0}%
\e{4}%
\e{3}%
\e{9}%
\e{4}%
\e{6}%
\e{1}%
\e{2}%
\e{9}%
\e{9}%
\e{13}%
\e{8}%
\e{7}%
\e{1}%
\e{3}%
\e{6}%
\e{6}%
\e{5}%
\e{3}%
\e{2}%
\e{0}%
\eol}\vss}\rg%
%
%
\rx{\vss\hfull{%
\rlx{\hss{$4536_{y}$}}\cg%
\e{1}%
\e{11}%
\e{14}%
\e{23}%
\e{16}%
\e{15}%
\e{3}%
\e{6}%
\e{27}%
\e{31}%
\e{41}%
\e{29}%
\e{23}%
\e{4}%
\e{5}%
\e{16}%
\e{17}%
\e{18}%
\e{13}%
\e{8}%
\e{1}%
\eol}\vss}\rg%
%
%
\rx{\vss\hfull{%
\rlx{\hss{$5670_{y}$}}\cg%
\e{2}%
\e{14}%
\e{18}%
\e{27}%
\e{21}%
\e{19}%
\e{4}%
\e{7}%
\e{34}%
\e{39}%
\e{51}%
\e{36}%
\e{29}%
\e{5}%
\e{5}%
\e{20}%
\e{21}%
\e{24}%
\e{15}%
\e{10}%
\e{1}%
\eol}\vss}\rg%
%
%
\rx{\vss\hfull{%
\rlx{\hss{$4480_{y}$}}\cg%
\e{2}%
\e{11}%
\e{16}%
\e{18}%
\e{16}%
\e{11}%
\e{2}%
\e{4}%
\e{23}%
\e{34}%
\e{38}%
\e{34}%
\e{23}%
\e{4}%
\e{2}%
\e{11}%
\e{16}%
\e{18}%
\e{16}%
\e{11}%
\e{2}%
\eol}\vss}\rg%
%
%
\rx{\vss\hfull{%
\rlx{\hss{$8_{z}$}}\cg%
\e{1}%
\e{0}%
\e{0}%
\e{0}%
\e{0}%
\e{0}%
\e{0}%
\e{0}%
\e{0}%
\e{0}%
\e{0}%
\e{0}%
\e{0}%
\e{0}%
\e{0}%
\e{0}%
\e{0}%
\e{0}%
\e{0}%
\e{0}%
\e{0}%
\eol}\vss}\rg%
%
%
\rx{\vss\hfull{%
\rlx{\hss{$56_{z}$}}\cg%
\e{0}%
\e{1}%
\e{0}%
\e{1}%
\e{0}%
\e{0}%
\e{0}%
\e{1}%
\e{1}%
\e{0}%
\e{0}%
\e{0}%
\e{0}%
\e{0}%
\e{1}%
\e{0}%
\e{0}%
\e{0}%
\e{0}%
\e{0}%
\e{0}%
\eol}\vss}\rg%
%
%
\rx{\vss\hfull{%
\rlx{\hss{$160_{z}$}}\cg%
\e{2}%
\e{3}%
\e{1}%
\e{1}%
\e{0}%
\e{0}%
\e{0}%
\e{3}%
\e{2}%
\e{0}%
\e{0}%
\e{0}%
\e{0}%
\e{0}%
\e{1}%
\e{0}%
\e{0}%
\e{0}%
\e{0}%
\e{0}%
\e{0}%
\eol}\vss}\rg%
%
%
\rx{\vss\hfull{%
\rlx{\hss{$112_{z}$}}\cg%
\e{3}%
\e{2}%
\e{1}%
\e{0}%
\e{0}%
\e{0}%
\e{0}%
\e{2}%
\e{1}%
\e{0}%
\e{0}%
\e{0}%
\e{0}%
\e{0}%
\e{0}%
\e{0}%
\e{0}%
\e{0}%
\e{0}%
\e{0}%
\e{0}%
\eol}\vss}\rg%
%
%
\rx{\vss\hfull{%
\rlx{\hss{$840_{z}$}}\cg%
\e{2}%
\e{4}%
\e{5}%
\e{4}%
\e{4}%
\e{2}%
\e{0}%
\e{3}%
\e{9}%
\e{5}%
\e{6}%
\e{2}%
\e{1}%
\e{0}%
\e{1}%
\e{4}%
\e{1}%
\e{2}%
\e{0}%
\e{0}%
\e{0}%
\eol}\vss}\rg%
\eop
\eject
\tablecont%
%
%
%
%
%
%
\rowpts=18 true pt%
\colpts=18 true pt%
\rowlabpts=40 true pt%
\collabpts=90 true pt%
\clx{\vss\hfull{%
\rlx{\hss{$ $}}\cg%
\cx{\hskip 16 true pt\flip{$[{5}]{\times}[{3}]{\times}[{1^{2}}]$}\hss}\cg%
\cx{\hskip 16 true pt\flip{$[{4}{1}]{\times}[{3}]{\times}[{1^{2}}]$}\hss}\cg%
\cx{\hskip 16 true pt\flip{$[{3}{2}]{\times}[{3}]{\times}[{1^{2}}]$}\hss}\cg%
\cx{\hskip 16 true pt\flip{$[{3}{1^{2}}]{\times}[{3}]{\times}[{1^{2}}]$}\hss}\cg%
\cx{\hskip 16 true pt\flip{$[{2^{2}}{1}]{\times}[{3}]{\times}[{1^{2}}]$}\hss}\cg%
\cx{\hskip 16 true pt\flip{$[{2}{1^{3}}]{\times}[{3}]{\times}[{1^{2}}]$}\hss}\cg%
\cx{\hskip 16 true pt\flip{$[{1^{5}}]{\times}[{3}]{\times}[{1^{2}}]$}\hss}\cg%
\cx{\hskip 16 true pt\flip{$[{5}]{\times}[{2}{1}]{\times}[{1^{2}}]$}\hss}\cg%
\cx{\hskip 16 true pt\flip{$[{4}{1}]{\times}[{2}{1}]{\times}[{1^{2}}]$}\hss}\cg%
\cx{\hskip 16 true pt\flip{$[{3}{2}]{\times}[{2}{1}]{\times}[{1^{2}}]$}\hss}\cg%
\cx{\hskip 16 true pt\flip{$[{3}{1^{2}}]{\times}[{2}{1}]{\times}[{1^{2}}]$}\hss}\cg%
\cx{\hskip 16 true pt\flip{$[{2^{2}}{1}]{\times}[{2}{1}]{\times}[{1^{2}}]$}\hss}\cg%
\cx{\hskip 16 true pt\flip{$[{2}{1^{3}}]{\times}[{2}{1}]{\times}[{1^{2}}]$}\hss}\cg%
\cx{\hskip 16 true pt\flip{$[{1^{5}}]{\times}[{2}{1}]{\times}[{1^{2}}]$}\hss}\cg%
\cx{\hskip 16 true pt\flip{$[{5}]{\times}[{1^{3}}]{\times}[{1^{2}}]$}\hss}\cg%
\cx{\hskip 16 true pt\flip{$[{4}{1}]{\times}[{1^{3}}]{\times}[{1^{2}}]$}\hss}\cg%
\cx{\hskip 16 true pt\flip{$[{3}{2}]{\times}[{1^{3}}]{\times}[{1^{2}}]$}\hss}\cg%
\cx{\hskip 16 true pt\flip{$[{3}{1^{2}}]{\times}[{1^{3}}]{\times}[{1^{2}}]$}\hss}\cg%
\cx{\hskip 16 true pt\flip{$[{2^{2}}{1}]{\times}[{1^{3}}]{\times}[{1^{2}}]$}\hss}\cg%
\cx{\hskip 16 true pt\flip{$[{2}{1^{3}}]{\times}[{1^{3}}]{\times}[{1^{2}}]$}\hss}\cg%
\cx{\hskip 16 true pt\flip{$[{1^{5}}]{\times}[{1^{3}}]{\times}[{1^{2}}]$}\hss}\cg%
\eol}}\rg%
%
%
\rx{\vss\hfull{%
\rlx{\hss{$1296_{z}$}}\cg%
\e{3}%
\e{10}%
\e{6}%
\e{10}%
\e{4}%
\e{3}%
\e{0}%
\e{8}%
\e{15}%
\e{7}%
\e{8}%
\e{2}%
\e{1}%
\e{0}%
\e{5}%
\e{6}%
\e{2}%
\e{1}%
\e{0}%
\e{0}%
\e{0}%
\eol}\vss}\rg%
%
%
\rx{\vss\hfull{%
\rlx{\hss{$1400_{z}$}}\cg%
\e{7}%
\e{13}%
\e{10}%
\e{7}%
\e{4}%
\e{1}%
\e{0}%
\e{10}%
\e{16}%
\e{8}%
\e{6}%
\e{2}%
\e{0}%
\e{0}%
\e{3}%
\e{4}%
\e{1}%
\e{1}%
\e{0}%
\e{0}%
\e{0}%
\eol}\vss}\rg%
%
%
\rx{\vss\hfull{%
\rlx{\hss{$1008_{z}$}}\cg%
\e{5}%
\e{11}%
\e{6}%
\e{7}%
\e{2}%
\e{1}%
\e{0}%
\e{9}%
\e{12}%
\e{5}%
\e{4}%
\e{1}%
\e{0}%
\e{0}%
\e{4}%
\e{3}%
\e{1}%
\e{0}%
\e{0}%
\e{0}%
\e{0}%
\eol}\vss}\rg%
%
%
\rx{\vss\hfull{%
\rlx{\hss{$560_{z}$}}\cg%
\e{6}%
\e{8}%
\e{4}%
\e{2}%
\e{1}%
\e{0}%
\e{0}%
\e{6}%
\e{6}%
\e{3}%
\e{1}%
\e{0}%
\e{0}%
\e{0}%
\e{1}%
\e{1}%
\e{0}%
\e{0}%
\e{0}%
\e{0}%
\e{0}%
\eol}\vss}\rg%
%
%
\rx{\vss\hfull{%
\rlx{\hss{$1400_{zz}$}}\cg%
\e{4}%
\e{9}%
\e{9}%
\e{5}%
\e{3}%
\e{0}%
\e{0}%
\e{3}%
\e{9}%
\e{14}%
\e{8}%
\e{8}%
\e{2}%
\e{0}%
\e{0}%
\e{1}%
\e{4}%
\e{2}%
\e{3}%
\e{1}%
\e{0}%
\eol}\vss}\rg%
%
%
\rx{\vss\hfull{%
\rlx{\hss{$4200_{z}$}}\cg%
\e{5}%
\e{17}%
\e{18}%
\e{19}%
\e{12}%
\e{7}%
\e{0}%
\e{6}%
\e{24}%
\e{36}%
\e{32}%
\e{28}%
\e{16}%
\e{2}%
\e{2}%
\e{8}%
\e{16}%
\e{11}%
\e{13}%
\e{6}%
\e{1}%
\eol}\vss}\rg%
%
%
\rx{\vss\hfull{%
\rlx{\hss{$400_{z}$}}\cg%
\e{4}%
\e{5}%
\e{3}%
\e{1}%
\e{0}%
\e{0}%
\e{0}%
\e{2}%
\e{3}%
\e{4}%
\e{1}%
\e{1}%
\e{0}%
\e{0}%
\e{0}%
\e{0}%
\e{1}%
\e{0}%
\e{0}%
\e{0}%
\e{0}%
\eol}\vss}\rg%
%
%
\rx{\vss\hfull{%
\rlx{\hss{$3240_{z}$}}\cg%
\e{11}%
\e{24}%
\e{20}%
\e{16}%
\e{9}%
\e{3}%
\e{0}%
\e{14}%
\e{30}%
\e{25}%
\e{18}%
\e{10}%
\e{3}%
\e{0}%
\e{4}%
\e{8}%
\e{7}%
\e{4}%
\e{2}%
\e{0}%
\e{0}%
\eol}\vss}\rg%
%
%
\rx{\vss\hfull{%
\rlx{\hss{$4536_{z}$}}\cg%
\e{5}%
\e{18}%
\e{23}%
\e{18}%
\e{17}%
\e{7}%
\e{1}%
\e{8}%
\e{31}%
\e{37}%
\e{34}%
\e{27}%
\e{12}%
\e{1}%
\e{2}%
\e{11}%
\e{12}%
\e{14}%
\e{9}%
\e{5}%
\e{0}%
\eol}\vss}\rg%
%
%
\rx{\vss\hfull{%
\rlx{\hss{$2400_{z}$}}\cg%
\e{2}%
\e{9}%
\e{9}%
\e{15}%
\e{8}%
\e{9}%
\e{1}%
\e{6}%
\e{20}%
\e{15}%
\e{21}%
\e{10}%
\e{8}%
\e{1}%
\e{5}%
\e{11}%
\e{8}%
\e{7}%
\e{3}%
\e{1}%
\e{0}%
\eol}\vss}\rg%
%
%
\rx{\vss\hfull{%
\rlx{\hss{$3360_{z}$}}\cg%
\e{5}%
\e{17}%
\e{17}%
\e{17}%
\e{10}%
\e{5}%
\e{0}%
\e{9}%
\e{25}%
\e{27}%
\e{24}%
\e{17}%
\e{8}%
\e{1}%
\e{3}%
\e{8}%
\e{11}%
\e{7}%
\e{6}%
\e{2}%
\e{0}%
\eol}\vss}\rg%
%
%
\rx{\vss\hfull{%
\rlx{\hss{$2800_{z}$}}\cg%
\e{4}%
\e{15}%
\e{12}%
\e{17}%
\e{9}%
\e{7}%
\e{1}%
\e{9}%
\e{24}%
\e{20}%
\e{21}%
\e{11}%
\e{6}%
\e{0}%
\e{6}%
\e{11}%
\e{7}%
\e{6}%
\e{3}%
\e{1}%
\e{0}%
\eol}\vss}\rg%
%
%
\rx{\vss\hfull{%
\rlx{\hss{$4096_{z}$}}\cg%
\e{7}%
\e{22}%
\e{21}%
\e{22}%
\e{15}%
\e{8}%
\e{1}%
\e{14}%
\e{37}%
\e{29}%
\e{29}%
\e{15}%
\e{7}%
\e{0}%
\e{6}%
\e{15}%
\e{9}%
\e{9}%
\e{3}%
\e{1}%
\e{0}%
\eol}\vss}\rg%
%
%
\rx{\vss\hfull{%
\rlx{\hss{$5600_{z}$}}\cg%
\e{4}%
\e{18}%
\e{23}%
\e{28}%
\e{22}%
\e{17}%
\e{3}%
\e{11}%
\e{42}%
\e{38}%
\e{48}%
\e{29}%
\e{19}%
\e{2}%
\e{7}%
\e{22}%
\e{16}%
\e{20}%
\e{9}%
\e{5}%
\e{0}%
\eol}\vss}\rg%
%
%
\rx{\vss\hfull{%
\rlx{\hss{$448_{z}$}}\cg%
\e{2}%
\e{3}%
\e{4}%
\e{1}%
\e{2}%
\e{0}%
\e{0}%
\e{2}%
\e{5}%
\e{3}%
\e{2}%
\e{1}%
\e{0}%
\e{0}%
\e{0}%
\e{1}%
\e{0}%
\e{1}%
\e{0}%
\e{0}%
\e{0}%
\eol}\vss}\rg%
%
%
\rx{\vss\hfull{%
\rlx{\hss{$448_{w}$}}\cg%
\e{0}%
\e{1}%
\e{1}%
\e{3}%
\e{2}%
\e{3}%
\e{1}%
\e{1}%
\e{4}%
\e{2}%
\e{5}%
\e{1}%
\e{2}%
\e{0}%
\e{1}%
\e{4}%
\e{1}%
\e{2}%
\e{0}%
\e{0}%
\e{0}%
\eol}\vss}\rg%
%
%
\rx{\vss\hfull{%
\rlx{\hss{$1344_{w}$}}\cg%
\e{1}%
\e{4}%
\e{5}%
\e{4}%
\e{4}%
\e{2}%
\e{0}%
\e{0}%
\e{5}%
\e{12}%
\e{10}%
\e{13}%
\e{7}%
\e{1}%
\e{0}%
\e{1}%
\e{5}%
\e{5}%
\e{6}%
\e{5}%
\e{1}%
\eol}\vss}\rg%
%
%
\rx{\vss\hfull{%
\rlx{\hss{$5600_{w}$}}\cg%
\e{2}%
\e{14}%
\e{18}%
\e{26}%
\e{21}%
\e{18}%
\e{4}%
\e{7}%
\e{33}%
\e{39}%
\e{50}%
\e{36}%
\e{29}%
\e{5}%
\e{4}%
\e{19}%
\e{21}%
\e{24}%
\e{15}%
\e{10}%
\e{1}%
\eol}\vss}\rg%
%
%
\rx{\vss\hfull{%
\rlx{\hss{$2016_{w}$}}\cg%
\e{1}%
\e{4}%
\e{9}%
\e{6}%
\e{7}%
\e{2}%
\e{0}%
\e{1}%
\e{8}%
\e{17}%
\e{16}%
\e{19}%
\e{10}%
\e{2}%
\e{0}%
\e{2}%
\e{6}%
\e{7}%
\e{11}%
\e{7}%
\e{2}%
\eol}\vss}\rg%
%
%
\rx{\vss\hfull{%
\rlx{\hss{$7168_{w}$}}\cg%
\e{3}%
\e{17}%
\e{26}%
\e{29}%
\e{26}%
\e{17}%
\e{3}%
\e{7}%
\e{37}%
\e{54}%
\e{61}%
\e{54}%
\e{37}%
\e{7}%
\e{3}%
\e{17}%
\e{26}%
\e{29}%
\e{26}%
\e{17}%
\e{3}%
\eol}\vss}\rg%
\tableclose%
%
%
%
%
%
%
\tableopen{Induce/restrict matrix for $W({A_{4}}{A_{3}})\,\subset\,W(E_{8})$}%
%
%
%
%
%
%
\rowpts=18 true pt%
\colpts=18 true pt%
\rowlabpts=40 true pt%
\collabpts=70 true pt%
\clx{\vss\hfull{%
\rlx{\hss{$ $}}\cg%
\cx{\hskip 16 true pt\flip{$[{5}]{\times}[{4}]$}\hss}\cg%
\cx{\hskip 16 true pt\flip{$[{4}{1}]{\times}[{4}]$}\hss}\cg%
\cx{\hskip 16 true pt\flip{$[{3}{2}]{\times}[{4}]$}\hss}\cg%
\cx{\hskip 16 true pt\flip{$[{3}{1^{2}}]{\times}[{4}]$}\hss}\cg%
\cx{\hskip 16 true pt\flip{$[{2^{2}}{1}]{\times}[{4}]$}\hss}\cg%
\cx{\hskip 16 true pt\flip{$[{2}{1^{3}}]{\times}[{4}]$}\hss}\cg%
\cx{\hskip 16 true pt\flip{$[{1^{5}}]{\times}[{4}]$}\hss}\cg%
\cx{\hskip 16 true pt\flip{$[{5}]{\times}[{3}{1}]$}\hss}\cg%
\cx{\hskip 16 true pt\flip{$[{4}{1}]{\times}[{3}{1}]$}\hss}\cg%
\cx{\hskip 16 true pt\flip{$[{3}{2}]{\times}[{3}{1}]$}\hss}\cg%
\cx{\hskip 16 true pt\flip{$[{3}{1^{2}}]{\times}[{3}{1}]$}\hss}\cg%
\cx{\hskip 16 true pt\flip{$[{2^{2}}{1}]{\times}[{3}{1}]$}\hss}\cg%
\cx{\hskip 16 true pt\flip{$[{2}{1^{3}}]{\times}[{3}{1}]$}\hss}\cg%
\cx{\hskip 16 true pt\flip{$[{1^{5}}]{\times}[{3}{1}]$}\hss}\cg%
\cx{\hskip 16 true pt\flip{$[{5}]{\times}[{2^{2}}]$}\hss}\cg%
\cx{\hskip 16 true pt\flip{$[{4}{1}]{\times}[{2^{2}}]$}\hss}\cg%
\cx{\hskip 16 true pt\flip{$[{3}{2}]{\times}[{2^{2}}]$}\hss}\cg%
\cx{\hskip 16 true pt\flip{$[{3}{1^{2}}]{\times}[{2^{2}}]$}\hss}\cg%
\eol}}\rg%
%
%
\rx{\vss\hfull{%
\rlx{\hss{$1_{x}$}}\cg%
\e{1}%
\e{0}%
\e{0}%
\e{0}%
\e{0}%
\e{0}%
\e{0}%
\e{0}%
\e{0}%
\e{0}%
\e{0}%
\e{0}%
\e{0}%
\e{0}%
\e{0}%
\e{0}%
\e{0}%
\e{0}%
\eol}\vss}\rg%
%
%
\rx{\vss\hfull{%
\rlx{\hss{$28_{x}$}}\cg%
\e{0}%
\e{1}%
\e{0}%
\e{1}%
\e{0}%
\e{0}%
\e{0}%
\e{1}%
\e{1}%
\e{0}%
\e{0}%
\e{0}%
\e{0}%
\e{0}%
\e{0}%
\e{0}%
\e{0}%
\e{0}%
\eol}\vss}\rg%
%
%
\rx{\vss\hfull{%
\rlx{\hss{$35_{x}$}}\cg%
\e{2}%
\e{2}%
\e{1}%
\e{0}%
\e{0}%
\e{0}%
\e{0}%
\e{2}%
\e{1}%
\e{0}%
\e{0}%
\e{0}%
\e{0}%
\e{0}%
\e{1}%
\e{0}%
\e{0}%
\e{0}%
\eol}\vss}\rg%
%
%
\rx{\vss\hfull{%
\rlx{\hss{$84_{x}$}}\cg%
\e{4}%
\e{3}%
\e{2}%
\e{0}%
\e{0}%
\e{0}%
\e{0}%
\e{3}%
\e{2}%
\e{1}%
\e{0}%
\e{0}%
\e{0}%
\e{0}%
\e{1}%
\e{1}%
\e{0}%
\e{0}%
\eol}\vss}\rg%
%
%
\rx{\vss\hfull{%
\rlx{\hss{$50_{x}$}}\cg%
\e{2}%
\e{2}%
\e{0}%
\e{0}%
\e{0}%
\e{0}%
\e{0}%
\e{1}%
\e{1}%
\e{1}%
\e{0}%
\e{0}%
\e{0}%
\e{0}%
\e{0}%
\e{0}%
\e{1}%
\e{0}%
\eol}\vss}\rg%
%
%
\rx{\vss\hfull{%
\rlx{\hss{$350_{x}$}}\cg%
\e{0}%
\e{1}%
\e{1}%
\e{3}%
\e{2}%
\e{2}%
\e{0}%
\e{1}%
\e{4}%
\e{2}%
\e{4}%
\e{1}%
\e{1}%
\e{0}%
\e{1}%
\e{2}%
\e{0}%
\e{1}%
\eol}\vss}\rg%
%
%
\rx{\vss\hfull{%
\rlx{\hss{$300_{x}$}}\cg%
\e{2}%
\e{3}%
\e{4}%
\e{1}%
\e{2}%
\e{0}%
\e{0}%
\e{3}%
\e{5}%
\e{3}%
\e{2}%
\e{1}%
\e{0}%
\e{0}%
\e{3}%
\e{3}%
\e{1}%
\e{0}%
\eol}\vss}\rg%
\eop
\eject
\tablecont%
%
%
%
%
%
%
\rowpts=18 true pt%
\colpts=18 true pt%
\rowlabpts=40 true pt%
\collabpts=70 true pt%
\clx{\vss\hfull{%
\rlx{\hss{$ $}}\cg%
\cx{\hskip 16 true pt\flip{$[{5}]{\times}[{4}]$}\hss}\cg%
\cx{\hskip 16 true pt\flip{$[{4}{1}]{\times}[{4}]$}\hss}\cg%
\cx{\hskip 16 true pt\flip{$[{3}{2}]{\times}[{4}]$}\hss}\cg%
\cx{\hskip 16 true pt\flip{$[{3}{1^{2}}]{\times}[{4}]$}\hss}\cg%
\cx{\hskip 16 true pt\flip{$[{2^{2}}{1}]{\times}[{4}]$}\hss}\cg%
\cx{\hskip 16 true pt\flip{$[{2}{1^{3}}]{\times}[{4}]$}\hss}\cg%
\cx{\hskip 16 true pt\flip{$[{1^{5}}]{\times}[{4}]$}\hss}\cg%
\cx{\hskip 16 true pt\flip{$[{5}]{\times}[{3}{1}]$}\hss}\cg%
\cx{\hskip 16 true pt\flip{$[{4}{1}]{\times}[{3}{1}]$}\hss}\cg%
\cx{\hskip 16 true pt\flip{$[{3}{2}]{\times}[{3}{1}]$}\hss}\cg%
\cx{\hskip 16 true pt\flip{$[{3}{1^{2}}]{\times}[{3}{1}]$}\hss}\cg%
\cx{\hskip 16 true pt\flip{$[{2^{2}}{1}]{\times}[{3}{1}]$}\hss}\cg%
\cx{\hskip 16 true pt\flip{$[{2}{1^{3}}]{\times}[{3}{1}]$}\hss}\cg%
\cx{\hskip 16 true pt\flip{$[{1^{5}}]{\times}[{3}{1}]$}\hss}\cg%
\cx{\hskip 16 true pt\flip{$[{5}]{\times}[{2^{2}}]$}\hss}\cg%
\cx{\hskip 16 true pt\flip{$[{4}{1}]{\times}[{2^{2}}]$}\hss}\cg%
\cx{\hskip 16 true pt\flip{$[{3}{2}]{\times}[{2^{2}}]$}\hss}\cg%
\cx{\hskip 16 true pt\flip{$[{3}{1^{2}}]{\times}[{2^{2}}]$}\hss}\cg%
\eol}}\rg%
%
%
\rx{\vss\hfull{%
\rlx{\hss{$567_{x}$}}\cg%
\e{3}%
\e{8}%
\e{5}%
\e{6}%
\e{2}%
\e{1}%
\e{0}%
\e{8}%
\e{11}%
\e{5}%
\e{4}%
\e{1}%
\e{0}%
\e{0}%
\e{3}%
\e{4}%
\e{1}%
\e{1}%
\eol}\vss}\rg%
%
%
\rx{\vss\hfull{%
\rlx{\hss{$210_{x}$}}\cg%
\e{3}%
\e{5}%
\e{2}%
\e{2}%
\e{1}%
\e{0}%
\e{0}%
\e{4}%
\e{5}%
\e{2}%
\e{1}%
\e{0}%
\e{0}%
\e{0}%
\e{2}%
\e{1}%
\e{1}%
\e{0}%
\eol}\vss}\rg%
%
%
\rx{\vss\hfull{%
\rlx{\hss{$840_{x}$}}\cg%
\e{2}%
\e{3}%
\e{4}%
\e{1}%
\e{2}%
\e{0}%
\e{0}%
\e{2}%
\e{6}%
\e{8}%
\e{4}%
\e{6}%
\e{1}%
\e{0}%
\e{1}%
\e{5}%
\e{6}%
\e{3}%
\eol}\vss}\rg%
%
%
\rx{\vss\hfull{%
\rlx{\hss{$700_{x}$}}\cg%
\e{7}%
\e{10}%
\e{7}%
\e{3}%
\e{2}%
\e{0}%
\e{0}%
\e{8}%
\e{12}%
\e{8}%
\e{4}%
\e{2}%
\e{0}%
\e{0}%
\e{4}%
\e{5}%
\e{4}%
\e{1}%
\eol}\vss}\rg%
%
%
\rx{\vss\hfull{%
\rlx{\hss{$175_{x}$}}\cg%
\e{2}%
\e{2}%
\e{1}%
\e{0}%
\e{0}%
\e{0}%
\e{0}%
\e{1}%
\e{2}%
\e{3}%
\e{1}%
\e{0}%
\e{0}%
\e{0}%
\e{0}%
\e{1}%
\e{1}%
\e{1}%
\eol}\vss}\rg%
%
%
\rx{\vss\hfull{%
\rlx{\hss{$1400_{x}$}}\cg%
\e{4}%
\e{10}%
\e{7}%
\e{8}%
\e{3}%
\e{2}%
\e{0}%
\e{8}%
\e{17}%
\e{13}%
\e{11}%
\e{5}%
\e{2}%
\e{0}%
\e{3}%
\e{7}%
\e{6}%
\e{5}%
\eol}\vss}\rg%
%
%
\rx{\vss\hfull{%
\rlx{\hss{$1050_{x}$}}\cg%
\e{4}%
\e{9}%
\e{6}%
\e{5}%
\e{1}%
\e{0}%
\e{0}%
\e{7}%
\e{13}%
\e{11}%
\e{7}%
\e{4}%
\e{1}%
\e{0}%
\e{2}%
\e{5}%
\e{6}%
\e{4}%
\eol}\vss}\rg%
%
%
\rx{\vss\hfull{%
\rlx{\hss{$1575_{x}$}}\cg%
\e{2}%
\e{9}%
\e{7}%
\e{11}%
\e{5}%
\e{4}%
\e{0}%
\e{8}%
\e{19}%
\e{12}%
\e{14}%
\e{6}%
\e{3}%
\e{0}%
\e{4}%
\e{8}%
\e{5}%
\e{5}%
\eol}\vss}\rg%
%
%
\rx{\vss\hfull{%
\rlx{\hss{$1344_{x}$}}\cg%
\e{6}%
\e{12}%
\e{11}%
\e{7}%
\e{5}%
\e{1}%
\e{0}%
\e{11}%
\e{19}%
\e{14}%
\e{9}%
\e{5}%
\e{1}%
\e{0}%
\e{6}%
\e{10}%
\e{5}%
\e{4}%
\eol}\vss}\rg%
%
%
\rx{\vss\hfull{%
\rlx{\hss{$2100_{x}$}}\cg%
\e{0}%
\e{3}%
\e{5}%
\e{9}%
\e{7}%
\e{8}%
\e{2}%
\e{3}%
\e{14}%
\e{12}%
\e{19}%
\e{10}%
\e{9}%
\e{1}%
\e{3}%
\e{9}%
\e{5}%
\e{8}%
\eol}\vss}\rg%
%
%
\rx{\vss\hfull{%
\rlx{\hss{$2268_{x}$}}\cg%
\e{3}%
\e{10}%
\e{10}%
\e{11}%
\e{8}%
\e{5}%
\e{1}%
\e{8}%
\e{23}%
\e{18}%
\e{19}%
\e{10}%
\e{5}%
\e{0}%
\e{6}%
\e{12}%
\e{9}%
\e{7}%
\eol}\vss}\rg%
%
%
\rx{\vss\hfull{%
\rlx{\hss{$525_{x}$}}\cg%
\e{1}%
\e{3}%
\e{2}%
\e{4}%
\e{1}%
\e{2}%
\e{0}%
\e{3}%
\e{6}%
\e{5}%
\e{4}%
\e{2}%
\e{1}%
\e{0}%
\e{0}%
\e{3}%
\e{1}%
\e{3}%
\eol}\vss}\rg%
%
%
\rx{\vss\hfull{%
\rlx{\hss{$700_{xx}$}}\cg%
\e{1}%
\e{4}%
\e{1}%
\e{3}%
\e{0}%
\e{0}%
\e{0}%
\e{2}%
\e{6}%
\e{6}%
\e{5}%
\e{3}%
\e{1}%
\e{0}%
\e{0}%
\e{1}%
\e{5}%
\e{3}%
\eol}\vss}\rg%
%
%
\rx{\vss\hfull{%
\rlx{\hss{$972_{x}$}}\cg%
\e{3}%
\e{6}%
\e{6}%
\e{2}%
\e{3}%
\e{0}%
\e{0}%
\e{4}%
\e{10}%
\e{10}%
\e{6}%
\e{5}%
\e{1}%
\e{0}%
\e{3}%
\e{6}%
\e{6}%
\e{3}%
\eol}\vss}\rg%
%
%
\rx{\vss\hfull{%
\rlx{\hss{$4096_{x}$}}\cg%
\e{3}%
\e{13}%
\e{15}%
\e{16}%
\e{13}%
\e{7}%
\e{1}%
\e{11}%
\e{33}%
\e{31}%
\e{32}%
\e{21}%
\e{11}%
\e{1}%
\e{8}%
\e{20}%
\e{16}%
\e{16}%
\eol}\vss}\rg%
%
%
\rx{\vss\hfull{%
\rlx{\hss{$4200_{x}$}}\cg%
\e{4}%
\e{13}%
\e{15}%
\e{13}%
\e{9}%
\e{4}%
\e{0}%
\e{10}%
\e{31}%
\e{33}%
\e{30}%
\e{22}%
\e{10}%
\e{1}%
\e{6}%
\e{18}%
\e{20}%
\e{17}%
\eol}\vss}\rg%
%
%
\rx{\vss\hfull{%
\rlx{\hss{$2240_{x}$}}\cg%
\e{6}%
\e{12}%
\e{11}%
\e{7}%
\e{5}%
\e{1}%
\e{0}%
\e{9}%
\e{21}%
\e{22}%
\e{15}%
\e{11}%
\e{3}%
\e{0}%
\e{4}%
\e{12}%
\e{12}%
\e{9}%
\eol}\vss}\rg%
%
%
\rx{\vss\hfull{%
\rlx{\hss{$2835_{x}$}}\cg%
\e{3}%
\e{8}%
\e{9}%
\e{6}%
\e{4}%
\e{1}%
\e{0}%
\e{5}%
\e{18}%
\e{23}%
\e{19}%
\e{14}%
\e{6}%
\e{0}%
\e{3}%
\e{10}%
\e{14}%
\e{13}%
\eol}\vss}\rg%
%
%
\rx{\vss\hfull{%
\rlx{\hss{$6075_{x}$}}\cg%
\e{2}%
\e{12}%
\e{13}%
\e{20}%
\e{13}%
\e{11}%
\e{1}%
\e{9}%
\e{36}%
\e{39}%
\e{46}%
\e{31}%
\e{22}%
\e{3}%
\e{5}%
\e{21}%
\e{22}%
\e{27}%
\eol}\vss}\rg%
%
%
\rx{\vss\hfull{%
\rlx{\hss{$3200_{x}$}}\cg%
\e{1}%
\e{5}%
\e{9}%
\e{6}%
\e{9}%
\e{3}%
\e{1}%
\e{4}%
\e{16}%
\e{22}%
\e{20}%
\e{20}%
\e{10}%
\e{2}%
\e{4}%
\e{14}%
\e{14}%
\e{14}%
\eol}\vss}\rg%
%
%
\rx{\vss\hfull{%
\rlx{\hss{$70_{y}$}}\cg%
\e{0}%
\e{0}%
\e{0}%
\e{0}%
\e{0}%
\e{1}%
\e{1}%
\e{0}%
\e{0}%
\e{0}%
\e{1}%
\e{0}%
\e{1}%
\e{0}%
\e{0}%
\e{0}%
\e{0}%
\e{0}%
\eol}\vss}\rg%
%
%
\rx{\vss\hfull{%
\rlx{\hss{$1134_{y}$}}\cg%
\e{0}%
\e{1}%
\e{1}%
\e{3}%
\e{2}%
\e{2}%
\e{0}%
\e{1}%
\e{4}%
\e{5}%
\e{7}%
\e{7}%
\e{6}%
\e{1}%
\e{0}%
\e{3}%
\e{3}%
\e{6}%
\eol}\vss}\rg%
%
%
\rx{\vss\hfull{%
\rlx{\hss{$1680_{y}$}}\cg%
\e{0}%
\e{0}%
\e{1}%
\e{4}%
\e{4}%
\e{7}%
\e{3}%
\e{0}%
\e{5}%
\e{5}%
\e{14}%
\e{8}%
\e{12}%
\e{3}%
\e{1}%
\e{4}%
\e{3}%
\e{6}%
\eol}\vss}\rg%
%
%
\rx{\vss\hfull{%
\rlx{\hss{$168_{y}$}}\cg%
\e{0}%
\e{1}%
\e{0}%
\e{0}%
\e{0}%
\e{0}%
\e{0}%
\e{0}%
\e{1}%
\e{1}%
\e{1}%
\e{1}%
\e{0}%
\e{0}%
\e{0}%
\e{0}%
\e{2}%
\e{0}%
\eol}\vss}\rg%
%
%
\rx{\vss\hfull{%
\rlx{\hss{$420_{y}$}}\cg%
\e{1}%
\e{1}%
\e{1}%
\e{0}%
\e{0}%
\e{0}%
\e{0}%
\e{0}%
\e{2}%
\e{4}%
\e{2}%
\e{2}%
\e{0}%
\e{0}%
\e{0}%
\e{1}%
\e{3}%
\e{2}%
\eol}\vss}\rg%
%
%
\rx{\vss\hfull{%
\rlx{\hss{$3150_{y}$}}\cg%
\e{1}%
\e{4}%
\e{6}%
\e{5}%
\e{4}%
\e{2}%
\e{0}%
\e{3}%
\e{12}%
\e{20}%
\e{18}%
\e{18}%
\e{10}%
\e{1}%
\e{1}%
\e{9}%
\e{13}%
\e{16}%
\eol}\vss}\rg%
%
%
\rx{\vss\hfull{%
\rlx{\hss{$4200_{y}$}}\cg%
\e{1}%
\e{6}%
\e{7}%
\e{7}%
\e{6}%
\e{2}%
\e{0}%
\e{3}%
\e{17}%
\e{25}%
\e{25}%
\e{24}%
\e{13}%
\e{2}%
\e{2}%
\e{11}%
\e{20}%
\e{18}%
\eol}\vss}\rg%
%
%
\rx{\vss\hfull{%
\rlx{\hss{$2688_{y}$}}\cg%
\e{0}%
\e{2}%
\e{5}%
\e{5}%
\e{5}%
\e{2}%
\e{0}%
\e{2}%
\e{10}%
\e{15}%
\e{17}%
\e{15}%
\e{10}%
\e{2}%
\e{2}%
\e{8}%
\e{10}%
\e{12}%
\eol}\vss}\rg%
%
%
\rx{\vss\hfull{%
\rlx{\hss{$2100_{y}$}}\cg%
\e{0}%
\e{1}%
\e{2}%
\e{5}%
\e{3}%
\e{7}%
\e{3}%
\e{1}%
\e{7}%
\e{9}%
\e{15}%
\e{11}%
\e{12}%
\e{3}%
\e{0}%
\e{5}%
\e{5}%
\e{10}%
\eol}\vss}\rg%
%
%
\rx{\vss\hfull{%
\rlx{\hss{$1400_{y}$}}\cg%
\e{0}%
\e{1}%
\e{1}%
\e{3}%
\e{3}%
\e{4}%
\e{2}%
\e{0}%
\e{5}%
\e{5}%
\e{11}%
\e{7}%
\e{8}%
\e{2}%
\e{1}%
\e{3}%
\e{5}%
\e{4}%
\eol}\vss}\rg%
%
%
\rx{\vss\hfull{%
\rlx{\hss{$4536_{y}$}}\cg%
\e{0}%
\e{3}%
\e{6}%
\e{9}%
\e{9}%
\e{9}%
\e{3}%
\e{2}%
\e{16}%
\e{21}%
\e{32}%
\e{24}%
\e{22}%
\e{5}%
\e{3}%
\e{12}%
\e{15}%
\e{18}%
\eol}\vss}\rg%
%
%
\rx{\vss\hfull{%
\rlx{\hss{$5670_{y}$}}\cg%
\e{0}%
\e{4}%
\e{7}%
\e{12}%
\e{11}%
\e{11}%
\e{3}%
\e{3}%
\e{20}%
\e{26}%
\e{39}%
\e{31}%
\e{28}%
\e{6}%
\e{3}%
\e{15}%
\e{18}%
\e{24}%
\eol}\vss}\rg%
%
%
\rx{\vss\hfull{%
\rlx{\hss{$4480_{y}$}}\cg%
\e{1}%
\e{5}%
\e{7}%
\e{8}%
\e{7}%
\e{5}%
\e{1}%
\e{3}%
\e{17}%
\e{25}%
\e{28}%
\e{25}%
\e{17}%
\e{3}%
\e{2}%
\e{12}%
\e{18}%
\e{20}%
\eol}\vss}\rg%
\eop
\eject
\tablecont%
%
%
%
%
%
%
\rowpts=18 true pt%
\colpts=18 true pt%
\rowlabpts=40 true pt%
\collabpts=70 true pt%
\clx{\vss\hfull{%
\rlx{\hss{$ $}}\cg%
\cx{\hskip 16 true pt\flip{$[{5}]{\times}[{4}]$}\hss}\cg%
\cx{\hskip 16 true pt\flip{$[{4}{1}]{\times}[{4}]$}\hss}\cg%
\cx{\hskip 16 true pt\flip{$[{3}{2}]{\times}[{4}]$}\hss}\cg%
\cx{\hskip 16 true pt\flip{$[{3}{1^{2}}]{\times}[{4}]$}\hss}\cg%
\cx{\hskip 16 true pt\flip{$[{2^{2}}{1}]{\times}[{4}]$}\hss}\cg%
\cx{\hskip 16 true pt\flip{$[{2}{1^{3}}]{\times}[{4}]$}\hss}\cg%
\cx{\hskip 16 true pt\flip{$[{1^{5}}]{\times}[{4}]$}\hss}\cg%
\cx{\hskip 16 true pt\flip{$[{5}]{\times}[{3}{1}]$}\hss}\cg%
\cx{\hskip 16 true pt\flip{$[{4}{1}]{\times}[{3}{1}]$}\hss}\cg%
\cx{\hskip 16 true pt\flip{$[{3}{2}]{\times}[{3}{1}]$}\hss}\cg%
\cx{\hskip 16 true pt\flip{$[{3}{1^{2}}]{\times}[{3}{1}]$}\hss}\cg%
\cx{\hskip 16 true pt\flip{$[{2^{2}}{1}]{\times}[{3}{1}]$}\hss}\cg%
\cx{\hskip 16 true pt\flip{$[{2}{1^{3}}]{\times}[{3}{1}]$}\hss}\cg%
\cx{\hskip 16 true pt\flip{$[{1^{5}}]{\times}[{3}{1}]$}\hss}\cg%
\cx{\hskip 16 true pt\flip{$[{5}]{\times}[{2^{2}}]$}\hss}\cg%
\cx{\hskip 16 true pt\flip{$[{4}{1}]{\times}[{2^{2}}]$}\hss}\cg%
\cx{\hskip 16 true pt\flip{$[{3}{2}]{\times}[{2^{2}}]$}\hss}\cg%
\cx{\hskip 16 true pt\flip{$[{3}{1^{2}}]{\times}[{2^{2}}]$}\hss}\cg%
\eol}}\rg%
%
%
\rx{\vss\hfull{%
\rlx{\hss{$8_{z}$}}\cg%
\e{1}%
\e{1}%
\e{0}%
\e{0}%
\e{0}%
\e{0}%
\e{0}%
\e{1}%
\e{0}%
\e{0}%
\e{0}%
\e{0}%
\e{0}%
\e{0}%
\e{0}%
\e{0}%
\e{0}%
\e{0}%
\eol}\vss}\rg%
%
%
\rx{\vss\hfull{%
\rlx{\hss{$56_{z}$}}\cg%
\e{0}%
\e{0}%
\e{0}%
\e{1}%
\e{0}%
\e{1}%
\e{0}%
\e{0}%
\e{1}%
\e{0}%
\e{1}%
\e{0}%
\e{0}%
\e{0}%
\e{0}%
\e{0}%
\e{0}%
\e{0}%
\eol}\vss}\rg%
%
%
\rx{\vss\hfull{%
\rlx{\hss{$160_{z}$}}\cg%
\e{1}%
\e{3}%
\e{2}%
\e{2}%
\e{1}%
\e{0}%
\e{0}%
\e{3}%
\e{4}%
\e{1}%
\e{1}%
\e{0}%
\e{0}%
\e{0}%
\e{2}%
\e{1}%
\e{0}%
\e{0}%
\eol}\vss}\rg%
%
%
\rx{\vss\hfull{%
\rlx{\hss{$112_{z}$}}\cg%
\e{4}%
\e{4}%
\e{2}%
\e{1}%
\e{0}%
\e{0}%
\e{0}%
\e{4}%
\e{3}%
\e{1}%
\e{0}%
\e{0}%
\e{0}%
\e{0}%
\e{1}%
\e{1}%
\e{0}%
\e{0}%
\eol}\vss}\rg%
%
%
\rx{\vss\hfull{%
\rlx{\hss{$840_{z}$}}\cg%
\e{1}%
\e{3}%
\e{4}%
\e{3}%
\e{4}%
\e{1}%
\e{0}%
\e{3}%
\e{7}%
\e{7}%
\e{6}%
\e{5}%
\e{2}%
\e{0}%
\e{2}%
\e{6}%
\e{3}%
\e{3}%
\eol}\vss}\rg%
%
%
\rx{\vss\hfull{%
\rlx{\hss{$1296_{z}$}}\cg%
\e{0}%
\e{4}%
\e{4}%
\e{9}%
\e{5}%
\e{5}%
\e{1}%
\e{4}%
\e{13}%
\e{8}%
\e{13}%
\e{5}%
\e{4}%
\e{0}%
\e{3}%
\e{6}%
\e{3}%
\e{4}%
\eol}\vss}\rg%
%
%
\rx{\vss\hfull{%
\rlx{\hss{$1400_{z}$}}\cg%
\e{6}%
\e{12}%
\e{11}%
\e{8}%
\e{5}%
\e{2}%
\e{0}%
\e{11}%
\e{20}%
\e{14}%
\e{10}%
\e{5}%
\e{1}%
\e{0}%
\e{6}%
\e{10}%
\e{5}%
\e{4}%
\eol}\vss}\rg%
%
%
\rx{\vss\hfull{%
\rlx{\hss{$1008_{z}$}}\cg%
\e{2}%
\e{8}%
\e{6}%
\e{8}%
\e{4}%
\e{3}%
\e{0}%
\e{7}%
\e{15}%
\e{8}%
\e{9}%
\e{3}%
\e{1}%
\e{0}%
\e{4}%
\e{6}%
\e{3}%
\e{2}%
\eol}\vss}\rg%
%
%
\rx{\vss\hfull{%
\rlx{\hss{$560_{z}$}}\cg%
\e{6}%
\e{10}%
\e{6}%
\e{4}%
\e{2}%
\e{0}%
\e{0}%
\e{9}%
\e{11}%
\e{6}%
\e{3}%
\e{1}%
\e{0}%
\e{0}%
\e{4}%
\e{4}%
\e{2}%
\e{1}%
\eol}\vss}\rg%
%
%
\rx{\vss\hfull{%
\rlx{\hss{$1400_{zz}$}}\cg%
\e{5}%
\e{9}%
\e{7}%
\e{4}%
\e{1}%
\e{0}%
\e{0}%
\e{6}%
\e{14}%
\e{15}%
\e{9}%
\e{6}%
\e{1}%
\e{0}%
\e{2}%
\e{6}%
\e{9}%
\e{6}%
\eol}\vss}\rg%
%
%
\rx{\vss\hfull{%
\rlx{\hss{$4200_{z}$}}\cg%
\e{2}%
\e{10}%
\e{8}%
\e{12}%
\e{6}%
\e{4}%
\e{0}%
\e{6}%
\e{24}%
\e{28}%
\e{30}%
\e{21}%
\e{13}%
\e{1}%
\e{2}%
\e{12}%
\e{19}%
\e{18}%
\eol}\vss}\rg%
%
%
\rx{\vss\hfull{%
\rlx{\hss{$400_{z}$}}\cg%
\e{5}%
\e{7}%
\e{3}%
\e{2}%
\e{0}%
\e{0}%
\e{0}%
\e{5}%
\e{7}%
\e{5}%
\e{2}%
\e{1}%
\e{0}%
\e{0}%
\e{1}%
\e{2}%
\e{3}%
\e{1}%
\eol}\vss}\rg%
%
%
\rx{\vss\hfull{%
\rlx{\hss{$3240_{z}$}}\cg%
\e{8}%
\e{20}%
\e{17}%
\e{15}%
\e{9}%
\e{3}%
\e{0}%
\e{16}%
\e{36}%
\e{30}%
\e{24}%
\e{14}%
\e{5}%
\e{0}%
\e{8}%
\e{18}%
\e{15}%
\e{11}%
\eol}\vss}\rg%
%
%
\rx{\vss\hfull{%
\rlx{\hss{$4536_{z}$}}\cg%
\e{5}%
\e{13}%
\e{16}%
\e{10}%
\e{10}%
\e{3}%
\e{0}%
\e{9}%
\e{30}%
\e{37}%
\e{29}%
\e{26}%
\e{10}%
\e{1}%
\e{6}%
\e{20}%
\e{23}%
\e{20}%
\eol}\vss}\rg%
%
%
\rx{\vss\hfull{%
\rlx{\hss{$2400_{z}$}}\cg%
\e{0}%
\e{2}%
\e{3}%
\e{9}%
\e{5}%
\e{9}%
\e{2}%
\e{2}%
\e{12}%
\e{11}%
\e{21}%
\e{11}%
\e{13}%
\e{2}%
\e{1}%
\e{7}%
\e{5}%
\e{10}%
\eol}\vss}\rg%
%
%
\rx{\vss\hfull{%
\rlx{\hss{$3360_{z}$}}\cg%
\e{2}%
\e{10}%
\e{11}%
\e{12}%
\e{7}%
\e{4}%
\e{0}%
\e{8}%
\e{25}%
\e{25}%
\e{26}%
\e{16}%
\e{9}%
\e{1}%
\e{5}%
\e{13}%
\e{14}%
\e{14}%
\eol}\vss}\rg%
%
%
\rx{\vss\hfull{%
\rlx{\hss{$2800_{z}$}}\cg%
\e{1}%
\e{7}%
\e{6}%
\e{12}%
\e{7}%
\e{8}%
\e{2}%
\e{5}%
\e{20}%
\e{18}%
\e{24}%
\e{13}%
\e{10}%
\e{1}%
\e{3}%
\e{10}%
\e{10}%
\e{11}%
\eol}\vss}\rg%
%
%
\rx{\vss\hfull{%
\rlx{\hss{$4096_{z}$}}\cg%
\e{3}%
\e{13}%
\e{15}%
\e{16}%
\e{13}%
\e{7}%
\e{1}%
\e{11}%
\e{33}%
\e{31}%
\e{32}%
\e{21}%
\e{11}%
\e{1}%
\e{8}%
\e{20}%
\e{16}%
\e{16}%
\eol}\vss}\rg%
%
%
\rx{\vss\hfull{%
\rlx{\hss{$5600_{z}$}}\cg%
\e{1}%
\e{7}%
\e{12}%
\e{15}%
\e{14}%
\e{12}%
\e{3}%
\e{6}%
\e{28}%
\e{33}%
\e{41}%
\e{31}%
\e{23}%
\e{4}%
\e{5}%
\e{21}%
\e{19}%
\e{24}%
\eol}\vss}\rg%
%
%
\rx{\vss\hfull{%
\rlx{\hss{$448_{z}$}}\cg%
\e{4}%
\e{4}%
\e{5}%
\e{1}%
\e{1}%
\e{0}%
\e{0}%
\e{4}%
\e{6}%
\e{6}%
\e{2}%
\e{2}%
\e{0}%
\e{0}%
\e{2}%
\e{4}%
\e{2}%
\e{2}%
\eol}\vss}\rg%
%
%
\rx{\vss\hfull{%
\rlx{\hss{$448_{w}$}}\cg%
\e{0}%
\e{0}%
\e{0}%
\e{1}%
\e{1}%
\e{3}%
\e{1}%
\e{0}%
\e{1}%
\e{1}%
\e{4}%
\e{2}%
\e{4}%
\e{1}%
\e{0}%
\e{1}%
\e{0}%
\e{2}%
\eol}\vss}\rg%
%
%
\rx{\vss\hfull{%
\rlx{\hss{$1344_{w}$}}\cg%
\e{1}%
\e{3}%
\e{2}%
\e{2}%
\e{1}%
\e{0}%
\e{0}%
\e{1}%
\e{6}%
\e{9}%
\e{7}%
\e{8}%
\e{3}%
\e{0}%
\e{0}%
\e{3}%
\e{8}%
\e{6}%
\eol}\vss}\rg%
%
%
\rx{\vss\hfull{%
\rlx{\hss{$5600_{w}$}}\cg%
\e{0}%
\e{4}%
\e{7}%
\e{12}%
\e{11}%
\e{10}%
\e{2}%
\e{3}%
\e{20}%
\e{26}%
\e{38}%
\e{31}%
\e{27}%
\e{6}%
\e{3}%
\e{15}%
\e{18}%
\e{24}%
\eol}\vss}\rg%
%
%
\rx{\vss\hfull{%
\rlx{\hss{$2016_{w}$}}\cg%
\e{1}%
\e{3}%
\e{5}%
\e{2}%
\e{2}%
\e{0}%
\e{0}%
\e{2}%
\e{8}%
\e{15}%
\e{11}%
\e{11}%
\e{4}%
\e{0}%
\e{1}%
\e{6}%
\e{10}%
\e{10}%
\eol}\vss}\rg%
%
%
\rx{\vss\hfull{%
\rlx{\hss{$7168_{w}$}}\cg%
\e{1}%
\e{7}%
\e{12}%
\e{13}%
\e{12}%
\e{7}%
\e{1}%
\e{5}%
\e{27}%
\e{40}%
\e{45}%
\e{40}%
\e{27}%
\e{5}%
\e{4}%
\e{20}%
\e{28}%
\e{32}%
\eol}\vss}\rg%
\eop
\eject
\tablecont%
%
%
%
%
%
%
\rowpts=18 true pt%
\colpts=18 true pt%
\rowlabpts=40 true pt%
\collabpts=70 true pt%
\clx{\vss\hfull{%
\rlx{\hss{$ $}}\cg%
\cx{\hskip 16 true pt\flip{$[{2^{2}}{1}]{\times}[{2^{2}}]$}\hss}\cg%
\cx{\hskip 16 true pt\flip{$[{2}{1^{3}}]{\times}[{2^{2}}]$}\hss}\cg%
\cx{\hskip 16 true pt\flip{$[{1^{5}}]{\times}[{2^{2}}]$}\hss}\cg%
\cx{\hskip 16 true pt\flip{$[{5}]{\times}[{2}{1^{2}}]$}\hss}\cg%
\cx{\hskip 16 true pt\flip{$[{4}{1}]{\times}[{2}{1^{2}}]$}\hss}\cg%
\cx{\hskip 16 true pt\flip{$[{3}{2}]{\times}[{2}{1^{2}}]$}\hss}\cg%
\cx{\hskip 16 true pt\flip{$[{3}{1^{2}}]{\times}[{2}{1^{2}}]$}\hss}\cg%
\cx{\hskip 16 true pt\flip{$[{2^{2}}{1}]{\times}[{2}{1^{2}}]$}\hss}\cg%
\cx{\hskip 16 true pt\flip{$[{2}{1^{3}}]{\times}[{2}{1^{2}}]$}\hss}\cg%
\cx{\hskip 16 true pt\flip{$[{1^{5}}]{\times}[{2}{1^{2}}]$}\hss}\cg%
\cx{\hskip 16 true pt\flip{$[{5}]{\times}[{1^{4}}]$}\hss}\cg%
\cx{\hskip 16 true pt\flip{$[{4}{1}]{\times}[{1^{4}}]$}\hss}\cg%
\cx{\hskip 16 true pt\flip{$[{3}{2}]{\times}[{1^{4}}]$}\hss}\cg%
\cx{\hskip 16 true pt\flip{$[{3}{1^{2}}]{\times}[{1^{4}}]$}\hss}\cg%
\cx{\hskip 16 true pt\flip{$[{2^{2}}{1}]{\times}[{1^{4}}]$}\hss}\cg%
\cx{\hskip 16 true pt\flip{$[{2}{1^{3}}]{\times}[{1^{4}}]$}\hss}\cg%
\cx{\hskip 16 true pt\flip{$[{1^{5}}]{\times}[{1^{4}}]$}\hss}\cg%
\eol}}\rg%
%
%
\rx{\vss\hfull{%
\rlx{\hss{$1_{x}$}}\cg%
\e{0}%
\e{0}%
\e{0}%
\e{0}%
\e{0}%
\e{0}%
\e{0}%
\e{0}%
\e{0}%
\e{0}%
\e{0}%
\e{0}%
\e{0}%
\e{0}%
\e{0}%
\e{0}%
\e{0}%
\eol}\vss}\rg%
%
%
\rx{\vss\hfull{%
\rlx{\hss{$28_{x}$}}\cg%
\e{0}%
\e{0}%
\e{0}%
\e{1}%
\e{0}%
\e{0}%
\e{0}%
\e{0}%
\e{0}%
\e{0}%
\e{0}%
\e{0}%
\e{0}%
\e{0}%
\e{0}%
\e{0}%
\e{0}%
\eol}\vss}\rg%
%
%
\rx{\vss\hfull{%
\rlx{\hss{$35_{x}$}}\cg%
\e{0}%
\e{0}%
\e{0}%
\e{0}%
\e{0}%
\e{0}%
\e{0}%
\e{0}%
\e{0}%
\e{0}%
\e{0}%
\e{0}%
\e{0}%
\e{0}%
\e{0}%
\e{0}%
\e{0}%
\eol}\vss}\rg%
%
%
\rx{\vss\hfull{%
\rlx{\hss{$84_{x}$}}\cg%
\e{0}%
\e{0}%
\e{0}%
\e{0}%
\e{0}%
\e{0}%
\e{0}%
\e{0}%
\e{0}%
\e{0}%
\e{0}%
\e{0}%
\e{0}%
\e{0}%
\e{0}%
\e{0}%
\e{0}%
\eol}\vss}\rg%
%
%
\rx{\vss\hfull{%
\rlx{\hss{$50_{x}$}}\cg%
\e{0}%
\e{0}%
\e{0}%
\e{0}%
\e{0}%
\e{0}%
\e{0}%
\e{0}%
\e{0}%
\e{0}%
\e{0}%
\e{0}%
\e{0}%
\e{0}%
\e{0}%
\e{0}%
\e{0}%
\eol}\vss}\rg%
%
%
\rx{\vss\hfull{%
\rlx{\hss{$350_{x}$}}\cg%
\e{0}%
\e{0}%
\e{0}%
\e{3}%
\e{4}%
\e{1}%
\e{1}%
\e{0}%
\e{0}%
\e{0}%
\e{1}%
\e{1}%
\e{0}%
\e{0}%
\e{0}%
\e{0}%
\e{0}%
\eol}\vss}\rg%
%
%
\rx{\vss\hfull{%
\rlx{\hss{$300_{x}$}}\cg%
\e{0}%
\e{0}%
\e{0}%
\e{1}%
\e{2}%
\e{0}%
\e{1}%
\e{0}%
\e{0}%
\e{0}%
\e{0}%
\e{0}%
\e{0}%
\e{0}%
\e{0}%
\e{0}%
\e{0}%
\eol}\vss}\rg%
%
%
\rx{\vss\hfull{%
\rlx{\hss{$567_{x}$}}\cg%
\e{0}%
\e{0}%
\e{0}%
\e{5}%
\e{4}%
\e{1}%
\e{0}%
\e{0}%
\e{0}%
\e{0}%
\e{1}%
\e{0}%
\e{0}%
\e{0}%
\e{0}%
\e{0}%
\e{0}%
\eol}\vss}\rg%
%
%
\rx{\vss\hfull{%
\rlx{\hss{$210_{x}$}}\cg%
\e{0}%
\e{0}%
\e{0}%
\e{2}%
\e{1}%
\e{0}%
\e{0}%
\e{0}%
\e{0}%
\e{0}%
\e{0}%
\e{0}%
\e{0}%
\e{0}%
\e{0}%
\e{0}%
\e{0}%
\eol}\vss}\rg%
%
%
\rx{\vss\hfull{%
\rlx{\hss{$840_{x}$}}\cg%
\e{3}%
\e{0}%
\e{0}%
\e{0}%
\e{3}%
\e{3}%
\e{5}%
\e{3}%
\e{2}%
\e{0}%
\e{0}%
\e{0}%
\e{0}%
\e{1}%
\e{0}%
\e{1}%
\e{0}%
\eol}\vss}\rg%
%
%
\rx{\vss\hfull{%
\rlx{\hss{$700_{x}$}}\cg%
\e{1}%
\e{0}%
\e{0}%
\e{2}%
\e{3}%
\e{2}%
\e{1}%
\e{0}%
\e{0}%
\e{0}%
\e{0}%
\e{0}%
\e{0}%
\e{0}%
\e{0}%
\e{0}%
\e{0}%
\eol}\vss}\rg%
%
%
\rx{\vss\hfull{%
\rlx{\hss{$175_{x}$}}\cg%
\e{1}%
\e{0}%
\e{0}%
\e{0}%
\e{0}%
\e{1}%
\e{0}%
\e{1}%
\e{0}%
\e{0}%
\e{0}%
\e{0}%
\e{0}%
\e{0}%
\e{0}%
\e{0}%
\e{0}%
\eol}\vss}\rg%
%
%
\rx{\vss\hfull{%
\rlx{\hss{$1400_{x}$}}\cg%
\e{3}%
\e{1}%
\e{0}%
\e{5}%
\e{8}%
\e{7}%
\e{3}%
\e{2}%
\e{0}%
\e{0}%
\e{1}%
\e{1}%
\e{1}%
\e{0}%
\e{0}%
\e{0}%
\e{0}%
\eol}\vss}\rg%
%
%
\rx{\vss\hfull{%
\rlx{\hss{$1050_{x}$}}\cg%
\e{3}%
\e{1}%
\e{0}%
\e{2}%
\e{4}%
\e{6}%
\e{2}%
\e{2}%
\e{0}%
\e{0}%
\e{0}%
\e{0}%
\e{1}%
\e{0}%
\e{0}%
\e{0}%
\e{0}%
\eol}\vss}\rg%
%
%
\rx{\vss\hfull{%
\rlx{\hss{$1575_{x}$}}\cg%
\e{2}%
\e{1}%
\e{0}%
\e{8}%
\e{12}%
\e{7}%
\e{4}%
\e{1}%
\e{0}%
\e{0}%
\e{2}%
\e{2}%
\e{1}%
\e{0}%
\e{0}%
\e{0}%
\e{0}%
\eol}\vss}\rg%
%
%
\rx{\vss\hfull{%
\rlx{\hss{$1344_{x}$}}\cg%
\e{1}%
\e{0}%
\e{0}%
\e{5}%
\e{8}%
\e{4}%
\e{3}%
\e{1}%
\e{0}%
\e{0}%
\e{0}%
\e{1}%
\e{0}%
\e{0}%
\e{0}%
\e{0}%
\e{0}%
\eol}\vss}\rg%
%
%
\rx{\vss\hfull{%
\rlx{\hss{$2100_{x}$}}\cg%
\e{3}%
\e{2}%
\e{0}%
\e{7}%
\e{17}%
\e{9}%
\e{12}%
\e{4}%
\e{2}%
\e{0}%
\e{4}%
\e{6}%
\e{2}%
\e{2}%
\e{0}%
\e{0}%
\e{0}%
\eol}\vss}\rg%
%
%
\rx{\vss\hfull{%
\rlx{\hss{$2268_{x}$}}\cg%
\e{4}%
\e{1}%
\e{0}%
\e{8}%
\e{16}%
\e{9}%
\e{9}%
\e{3}%
\e{1}%
\e{0}%
\e{3}%
\e{3}%
\e{1}%
\e{1}%
\e{0}%
\e{0}%
\e{0}%
\eol}\vss}\rg%
%
%
\rx{\vss\hfull{%
\rlx{\hss{$525_{x}$}}\cg%
\e{0}%
\e{0}%
\e{0}%
\e{3}%
\e{4}%
\e{2}%
\e{1}%
\e{1}%
\e{0}%
\e{0}%
\e{1}%
\e{1}%
\e{0}%
\e{0}%
\e{0}%
\e{0}%
\e{0}%
\eol}\vss}\rg%
%
%
\rx{\vss\hfull{%
\rlx{\hss{$700_{xx}$}}\cg%
\e{3}%
\e{2}%
\e{0}%
\e{1}%
\e{2}%
\e{5}%
\e{2}%
\e{3}%
\e{1}%
\e{0}%
\e{0}%
\e{0}%
\e{1}%
\e{0}%
\e{1}%
\e{0}%
\e{0}%
\eol}\vss}\rg%
%
%
\rx{\vss\hfull{%
\rlx{\hss{$972_{x}$}}\cg%
\e{3}%
\e{0}%
\e{0}%
\e{1}%
\e{4}%
\e{4}%
\e{4}%
\e{2}%
\e{1}%
\e{0}%
\e{0}%
\e{0}%
\e{0}%
\e{1}%
\e{0}%
\e{0}%
\e{0}%
\eol}\vss}\rg%
%
%
\rx{\vss\hfull{%
\rlx{\hss{$4096_{x}$}}\cg%
\e{8}%
\e{4}%
\e{0}%
\e{10}%
\e{26}%
\e{19}%
\e{19}%
\e{9}%
\e{4}%
\e{0}%
\e{2}%
\e{6}%
\e{3}%
\e{3}%
\e{1}%
\e{0}%
\e{0}%
\eol}\vss}\rg%
%
%
\rx{\vss\hfull{%
\rlx{\hss{$4200_{x}$}}\cg%
\e{13}%
\e{6}%
\e{1}%
\e{6}%
\e{20}%
\e{22}%
\e{20}%
\e{14}%
\e{6}%
\e{0}%
\e{1}%
\e{3}%
\e{5}%
\e{3}%
\e{2}%
\e{1}%
\e{0}%
\eol}\vss}\rg%
%
%
\rx{\vss\hfull{%
\rlx{\hss{$2240_{x}$}}\cg%
\e{6}%
\e{2}%
\e{0}%
\e{3}%
\e{10}%
\e{11}%
\e{8}%
\e{6}%
\e{2}%
\e{0}%
\e{0}%
\e{1}%
\e{1}%
\e{1}%
\e{1}%
\e{0}%
\e{0}%
\eol}\vss}\rg%
%
%
\rx{\vss\hfull{%
\rlx{\hss{$2835_{x}$}}\cg%
\e{12}%
\e{6}%
\e{1}%
\e{2}%
\e{9}%
\e{17}%
\e{13}%
\e{14}%
\e{6}%
\e{1}%
\e{0}%
\e{1}%
\e{3}%
\e{2}%
\e{4}%
\e{1}%
\e{0}%
\eol}\vss}\rg%
%
%
\rx{\vss\hfull{%
\rlx{\hss{$6075_{x}$}}\cg%
\e{17}%
\e{12}%
\e{1}%
\e{11}%
\e{33}%
\e{34}%
\e{32}%
\e{22}%
\e{11}%
\e{1}%
\e{3}%
\e{9}%
\e{9}%
\e{7}%
\e{4}%
\e{1}%
\e{0}%
\eol}\vss}\rg%
%
%
\rx{\vss\hfull{%
\rlx{\hss{$3200_{x}$}}\cg%
\e{10}%
\e{4}%
\e{0}%
\e{3}%
\e{15}%
\e{15}%
\e{20}%
\e{13}%
\e{9}%
\e{1}%
\e{0}%
\e{4}%
\e{2}%
\e{6}%
\e{2}%
\e{2}%
\e{0}%
\eol}\vss}\rg%
%
%
\rx{\vss\hfull{%
\rlx{\hss{$70_{y}$}}\cg%
\e{0}%
\e{0}%
\e{0}%
\e{0}%
\e{1}%
\e{0}%
\e{1}%
\e{0}%
\e{0}%
\e{0}%
\e{1}%
\e{1}%
\e{0}%
\e{0}%
\e{0}%
\e{0}%
\e{0}%
\eol}\vss}\rg%
%
%
\rx{\vss\hfull{%
\rlx{\hss{$1134_{y}$}}\cg%
\e{3}%
\e{3}%
\e{0}%
\e{1}%
\e{6}%
\e{7}%
\e{7}%
\e{5}%
\e{4}%
\e{1}%
\e{0}%
\e{2}%
\e{2}%
\e{3}%
\e{1}%
\e{1}%
\e{0}%
\eol}\vss}\rg%
%
%
\rx{\vss\hfull{%
\rlx{\hss{$1680_{y}$}}\cg%
\e{3}%
\e{4}%
\e{1}%
\e{3}%
\e{12}%
\e{8}%
\e{14}%
\e{5}%
\e{5}%
\e{0}%
\e{3}%
\e{7}%
\e{4}%
\e{4}%
\e{1}%
\e{0}%
\e{0}%
\eol}\vss}\rg%
%
%
\rx{\vss\hfull{%
\rlx{\hss{$168_{y}$}}\cg%
\e{2}%
\e{0}%
\e{0}%
\e{0}%
\e{0}%
\e{1}%
\e{1}%
\e{1}%
\e{1}%
\e{0}%
\e{0}%
\e{0}%
\e{0}%
\e{0}%
\e{0}%
\e{1}%
\e{0}%
\eol}\vss}\rg%
%
%
\rx{\vss\hfull{%
\rlx{\hss{$420_{y}$}}\cg%
\e{3}%
\e{1}%
\e{0}%
\e{0}%
\e{0}%
\e{2}%
\e{2}%
\e{4}%
\e{2}%
\e{0}%
\e{0}%
\e{0}%
\e{0}%
\e{0}%
\e{1}%
\e{1}%
\e{1}%
\eol}\vss}\rg%
%
%
\rx{\vss\hfull{%
\rlx{\hss{$3150_{y}$}}\cg%
\e{13}%
\e{9}%
\e{1}%
\e{1}%
\e{10}%
\e{18}%
\e{18}%
\e{20}%
\e{12}%
\e{3}%
\e{0}%
\e{2}%
\e{4}%
\e{5}%
\e{6}%
\e{4}%
\e{1}%
\eol}\vss}\rg%
%
%
\rx{\vss\hfull{%
\rlx{\hss{$4200_{y}$}}\cg%
\e{20}%
\e{11}%
\e{2}%
\e{2}%
\e{13}%
\e{24}%
\e{25}%
\e{25}%
\e{17}%
\e{3}%
\e{0}%
\e{2}%
\e{6}%
\e{7}%
\e{7}%
\e{6}%
\e{1}%
\eol}\vss}\rg%
\eop
\eject
\tablecont%
%
%
%
%
%
%
\rowpts=18 true pt%
\colpts=18 true pt%
\rowlabpts=40 true pt%
\collabpts=70 true pt%
\clx{\vss\hfull{%
\rlx{\hss{$ $}}\cg%
\cx{\hskip 16 true pt\flip{$[{2^{2}}{1}]{\times}[{2^{2}}]$}\hss}\cg%
\cx{\hskip 16 true pt\flip{$[{2}{1^{3}}]{\times}[{2^{2}}]$}\hss}\cg%
\cx{\hskip 16 true pt\flip{$[{1^{5}}]{\times}[{2^{2}}]$}\hss}\cg%
\cx{\hskip 16 true pt\flip{$[{5}]{\times}[{2}{1^{2}}]$}\hss}\cg%
\cx{\hskip 16 true pt\flip{$[{4}{1}]{\times}[{2}{1^{2}}]$}\hss}\cg%
\cx{\hskip 16 true pt\flip{$[{3}{2}]{\times}[{2}{1^{2}}]$}\hss}\cg%
\cx{\hskip 16 true pt\flip{$[{3}{1^{2}}]{\times}[{2}{1^{2}}]$}\hss}\cg%
\cx{\hskip 16 true pt\flip{$[{2^{2}}{1}]{\times}[{2}{1^{2}}]$}\hss}\cg%
\cx{\hskip 16 true pt\flip{$[{2}{1^{3}}]{\times}[{2}{1^{2}}]$}\hss}\cg%
\cx{\hskip 16 true pt\flip{$[{1^{5}}]{\times}[{2}{1^{2}}]$}\hss}\cg%
\cx{\hskip 16 true pt\flip{$[{5}]{\times}[{1^{4}}]$}\hss}\cg%
\cx{\hskip 16 true pt\flip{$[{4}{1}]{\times}[{1^{4}}]$}\hss}\cg%
\cx{\hskip 16 true pt\flip{$[{3}{2}]{\times}[{1^{4}}]$}\hss}\cg%
\cx{\hskip 16 true pt\flip{$[{3}{1^{2}}]{\times}[{1^{4}}]$}\hss}\cg%
\cx{\hskip 16 true pt\flip{$[{2^{2}}{1}]{\times}[{1^{4}}]$}\hss}\cg%
\cx{\hskip 16 true pt\flip{$[{2}{1^{3}}]{\times}[{1^{4}}]$}\hss}\cg%
\cx{\hskip 16 true pt\flip{$[{1^{5}}]{\times}[{1^{4}}]$}\hss}\cg%
\eol}}\rg%
%
%
\rx{\vss\hfull{%
\rlx{\hss{$2688_{y}$}}\cg%
\e{10}%
\e{8}%
\e{2}%
\e{2}%
\e{10}%
\e{15}%
\e{17}%
\e{15}%
\e{10}%
\e{2}%
\e{0}%
\e{2}%
\e{5}%
\e{5}%
\e{5}%
\e{2}%
\e{0}%
\eol}\vss}\rg%
%
%
\rx{\vss\hfull{%
\rlx{\hss{$2100_{y}$}}\cg%
\e{5}%
\e{5}%
\e{0}%
\e{3}%
\e{12}%
\e{11}%
\e{15}%
\e{9}%
\e{7}%
\e{1}%
\e{3}%
\e{7}%
\e{3}%
\e{5}%
\e{2}%
\e{1}%
\e{0}%
\eol}\vss}\rg%
%
%
\rx{\vss\hfull{%
\rlx{\hss{$1400_{y}$}}\cg%
\e{5}%
\e{3}%
\e{1}%
\e{2}%
\e{8}%
\e{7}%
\e{11}%
\e{5}%
\e{5}%
\e{0}%
\e{2}%
\e{4}%
\e{3}%
\e{3}%
\e{1}%
\e{1}%
\e{0}%
\eol}\vss}\rg%
%
%
\rx{\vss\hfull{%
\rlx{\hss{$4536_{y}$}}\cg%
\e{15}%
\e{12}%
\e{3}%
\e{5}%
\e{22}%
\e{24}%
\e{32}%
\e{21}%
\e{16}%
\e{2}%
\e{3}%
\e{9}%
\e{9}%
\e{9}%
\e{6}%
\e{3}%
\e{0}%
\eol}\vss}\rg%
%
%
\rx{\vss\hfull{%
\rlx{\hss{$5670_{y}$}}\cg%
\e{18}%
\e{15}%
\e{3}%
\e{6}%
\e{28}%
\e{31}%
\e{39}%
\e{26}%
\e{20}%
\e{3}%
\e{3}%
\e{11}%
\e{11}%
\e{12}%
\e{7}%
\e{4}%
\e{0}%
\eol}\vss}\rg%
%
%
\rx{\vss\hfull{%
\rlx{\hss{$4480_{y}$}}\cg%
\e{18}%
\e{12}%
\e{2}%
\e{3}%
\e{17}%
\e{25}%
\e{28}%
\e{25}%
\e{17}%
\e{3}%
\e{1}%
\e{5}%
\e{7}%
\e{8}%
\e{7}%
\e{5}%
\e{1}%
\eol}\vss}\rg%
%
%
\rx{\vss\hfull{%
\rlx{\hss{$8_{z}$}}\cg%
\e{0}%
\e{0}%
\e{0}%
\e{0}%
\e{0}%
\e{0}%
\e{0}%
\e{0}%
\e{0}%
\e{0}%
\e{0}%
\e{0}%
\e{0}%
\e{0}%
\e{0}%
\e{0}%
\e{0}%
\eol}\vss}\rg%
%
%
\rx{\vss\hfull{%
\rlx{\hss{$56_{z}$}}\cg%
\e{0}%
\e{0}%
\e{0}%
\e{1}%
\e{1}%
\e{0}%
\e{0}%
\e{0}%
\e{0}%
\e{0}%
\e{1}%
\e{0}%
\e{0}%
\e{0}%
\e{0}%
\e{0}%
\e{0}%
\eol}\vss}\rg%
%
%
\rx{\vss\hfull{%
\rlx{\hss{$160_{z}$}}\cg%
\e{0}%
\e{0}%
\e{0}%
\e{2}%
\e{1}%
\e{0}%
\e{0}%
\e{0}%
\e{0}%
\e{0}%
\e{0}%
\e{0}%
\e{0}%
\e{0}%
\e{0}%
\e{0}%
\e{0}%
\eol}\vss}\rg%
%
%
\rx{\vss\hfull{%
\rlx{\hss{$112_{z}$}}\cg%
\e{0}%
\e{0}%
\e{0}%
\e{1}%
\e{0}%
\e{0}%
\e{0}%
\e{0}%
\e{0}%
\e{0}%
\e{0}%
\e{0}%
\e{0}%
\e{0}%
\e{0}%
\e{0}%
\e{0}%
\eol}\vss}\rg%
%
%
\rx{\vss\hfull{%
\rlx{\hss{$840_{z}$}}\cg%
\e{1}%
\e{0}%
\e{0}%
\e{2}%
\e{6}%
\e{3}%
\e{4}%
\e{1}%
\e{1}%
\e{0}%
\e{0}%
\e{1}%
\e{0}%
\e{1}%
\e{0}%
\e{0}%
\e{0}%
\eol}\vss}\rg%
%
%
\rx{\vss\hfull{%
\rlx{\hss{$1296_{z}$}}\cg%
\e{1}%
\e{1}%
\e{0}%
\e{7}%
\e{12}%
\e{5}%
\e{5}%
\e{1}%
\e{0}%
\e{0}%
\e{3}%
\e{3}%
\e{1}%
\e{0}%
\e{0}%
\e{0}%
\e{0}%
\eol}\vss}\rg%
%
%
\rx{\vss\hfull{%
\rlx{\hss{$1400_{z}$}}\cg%
\e{1}%
\e{0}%
\e{0}%
\e{6}%
\e{9}%
\e{4}%
\e{3}%
\e{1}%
\e{0}%
\e{0}%
\e{1}%
\e{1}%
\e{0}%
\e{0}%
\e{0}%
\e{0}%
\e{0}%
\eol}\vss}\rg%
%
%
\rx{\vss\hfull{%
\rlx{\hss{$1008_{z}$}}\cg%
\e{1}%
\e{0}%
\e{0}%
\e{7}%
\e{8}%
\e{3}%
\e{2}%
\e{0}%
\e{0}%
\e{0}%
\e{2}%
\e{1}%
\e{0}%
\e{0}%
\e{0}%
\e{0}%
\e{0}%
\eol}\vss}\rg%
%
%
\rx{\vss\hfull{%
\rlx{\hss{$560_{z}$}}\cg%
\e{0}%
\e{0}%
\e{0}%
\e{3}%
\e{3}%
\e{1}%
\e{0}%
\e{0}%
\e{0}%
\e{0}%
\e{0}%
\e{0}%
\e{0}%
\e{0}%
\e{0}%
\e{0}%
\e{0}%
\eol}\vss}\rg%
%
%
\rx{\vss\hfull{%
\rlx{\hss{$1400_{zz}$}}\cg%
\e{5}%
\e{2}%
\e{0}%
\e{1}%
\e{4}%
\e{8}%
\e{4}%
\e{5}%
\e{1}%
\e{0}%
\e{0}%
\e{0}%
\e{1}%
\e{0}%
\e{1}%
\e{0}%
\e{0}%
\eol}\vss}\rg%
%
%
\rx{\vss\hfull{%
\rlx{\hss{$4200_{z}$}}\cg%
\e{17}%
\e{10}%
\e{2}%
\e{5}%
\e{17}%
\e{26}%
\e{21}%
\e{19}%
\e{10}%
\e{1}%
\e{1}%
\e{3}%
\e{7}%
\e{4}%
\e{5}%
\e{2}%
\e{0}%
\eol}\vss}\rg%
%
%
\rx{\vss\hfull{%
\rlx{\hss{$400_{z}$}}\cg%
\e{1}%
\e{0}%
\e{0}%
\e{1}%
\e{1}%
\e{2}%
\e{0}%
\e{0}%
\e{0}%
\e{0}%
\e{0}%
\e{0}%
\e{0}%
\e{0}%
\e{0}%
\e{0}%
\e{0}%
\eol}\vss}\rg%
%
%
\rx{\vss\hfull{%
\rlx{\hss{$3240_{z}$}}\cg%
\e{7}%
\e{2}%
\e{0}%
\e{9}%
\e{18}%
\e{15}%
\e{10}%
\e{5}%
\e{1}%
\e{0}%
\e{1}%
\e{2}%
\e{2}%
\e{1}%
\e{0}%
\e{0}%
\e{0}%
\eol}\vss}\rg%
%
%
\rx{\vss\hfull{%
\rlx{\hss{$4536_{z}$}}\cg%
\e{15}%
\e{6}%
\e{0}%
\e{4}%
\e{19}%
\e{23}%
\e{23}%
\e{18}%
\e{9}%
\e{1}%
\e{0}%
\e{3}%
\e{3}%
\e{5}%
\e{3}%
\e{2}%
\e{0}%
\eol}\vss}\rg%
%
%
\rx{\vss\hfull{%
\rlx{\hss{$2400_{z}$}}\cg%
\e{5}%
\e{5}%
\e{1}%
\e{6}%
\e{17}%
\e{13}%
\e{15}%
\e{7}%
\e{4}%
\e{0}%
\e{4}%
\e{7}%
\e{5}%
\e{3}%
\e{1}%
\e{0}%
\e{0}%
\eol}\vss}\rg%
%
%
\rx{\vss\hfull{%
\rlx{\hss{$3360_{z}$}}\cg%
\e{10}%
\e{6}%
\e{1}%
\e{6}%
\e{17}%
\e{19}%
\e{15}%
\e{11}%
\e{4}%
\e{0}%
\e{1}%
\e{3}%
\e{5}%
\e{2}%
\e{2}%
\e{0}%
\e{0}%
\eol}\vss}\rg%
%
%
\rx{\vss\hfull{%
\rlx{\hss{$2800_{z}$}}\cg%
\e{6}%
\e{4}%
\e{0}%
\e{8}%
\e{19}%
\e{14}%
\e{14}%
\e{7}%
\e{3}%
\e{0}%
\e{4}%
\e{6}%
\e{3}%
\e{2}%
\e{1}%
\e{0}%
\e{0}%
\eol}\vss}\rg%
%
%
\rx{\vss\hfull{%
\rlx{\hss{$4096_{z}$}}\cg%
\e{8}%
\e{4}%
\e{0}%
\e{10}%
\e{26}%
\e{19}%
\e{19}%
\e{9}%
\e{4}%
\e{0}%
\e{2}%
\e{6}%
\e{3}%
\e{3}%
\e{1}%
\e{0}%
\e{0}%
\eol}\vss}\rg%
%
%
\rx{\vss\hfull{%
\rlx{\hss{$5600_{z}$}}\cg%
\e{15}%
\e{9}%
\e{1}%
\e{9}%
\e{32}%
\e{28}%
\e{35}%
\e{20}%
\e{13}%
\e{1}%
\e{4}%
\e{11}%
\e{7}%
\e{9}%
\e{3}%
\e{2}%
\e{0}%
\eol}\vss}\rg%
%
%
\rx{\vss\hfull{%
\rlx{\hss{$448_{z}$}}\cg%
\e{0}%
\e{0}%
\e{0}%
\e{0}%
\e{2}%
\e{1}%
\e{1}%
\e{1}%
\e{0}%
\e{0}%
\e{0}%
\e{0}%
\e{0}%
\e{0}%
\e{0}%
\e{0}%
\e{0}%
\eol}\vss}\rg%
%
%
\rx{\vss\hfull{%
\rlx{\hss{$448_{w}$}}\cg%
\e{0}%
\e{1}%
\e{0}%
\e{1}%
\e{4}%
\e{2}%
\e{4}%
\e{1}%
\e{1}%
\e{0}%
\e{1}%
\e{3}%
\e{1}%
\e{1}%
\e{0}%
\e{0}%
\e{0}%
\eol}\vss}\rg%
%
%
\rx{\vss\hfull{%
\rlx{\hss{$1344_{w}$}}\cg%
\e{8}%
\e{3}%
\e{0}%
\e{0}%
\e{3}%
\e{8}%
\e{7}%
\e{9}%
\e{6}%
\e{1}%
\e{0}%
\e{0}%
\e{1}%
\e{2}%
\e{2}%
\e{3}%
\e{1}%
\eol}\vss}\rg%
%
%
\rx{\vss\hfull{%
\rlx{\hss{$5600_{w}$}}\cg%
\e{18}%
\e{15}%
\e{3}%
\e{6}%
\e{27}%
\e{31}%
\e{38}%
\e{26}%
\e{20}%
\e{3}%
\e{2}%
\e{10}%
\e{11}%
\e{12}%
\e{7}%
\e{4}%
\e{0}%
\eol}\vss}\rg%
%
%
\rx{\vss\hfull{%
\rlx{\hss{$2016_{w}$}}\cg%
\e{10}%
\e{6}%
\e{1}%
\e{0}%
\e{4}%
\e{11}%
\e{11}%
\e{15}%
\e{8}%
\e{2}%
\e{0}%
\e{0}%
\e{2}%
\e{2}%
\e{5}%
\e{3}%
\e{1}%
\eol}\vss}\rg%
%
%
\rx{\vss\hfull{%
\rlx{\hss{$7168_{w}$}}\cg%
\e{28}%
\e{20}%
\e{4}%
\e{5}%
\e{27}%
\e{40}%
\e{45}%
\e{40}%
\e{27}%
\e{5}%
\e{1}%
\e{7}%
\e{12}%
\e{13}%
\e{12}%
\e{7}%
\e{1}%
\eol}\vss}\rg%
\tableclose%
%
%
%
%
%
%
\eop
\eject
\tableopen{Induce/restrict matrix for $W({D_{5}}{A_{2}})\,\subset\,W(E_{8})$}%
%
%
%
%
%
%
\rowpts=18 true pt%
\colpts=18 true pt%
\rowlabpts=40 true pt%
\collabpts=85 true pt%
\clx{\vss\hfull{%
\rlx{\hss{$ $}}\cg%
\cx{\hskip 16 true pt\flip{$[{5}:-]{\times}[{3}]$}\hss}\cg%
\cx{\hskip 16 true pt\flip{$[{4}{1}:-]{\times}[{3}]$}\hss}\cg%
\cx{\hskip 16 true pt\flip{$[{3}{2}:-]{\times}[{3}]$}\hss}\cg%
\cx{\hskip 16 true pt\flip{$[{3}{1^{2}}:-]{\times}[{3}]$}\hss}\cg%
\cx{\hskip 16 true pt\flip{$[{2^{2}}{1}:-]{\times}[{3}]$}\hss}\cg%
\cx{\hskip 16 true pt\flip{$[{2}{1^{3}}:-]{\times}[{3}]$}\hss}\cg%
\cx{\hskip 16 true pt\flip{$[{1^{5}}:-]{\times}[{3}]$}\hss}\cg%
\cx{\hskip 16 true pt\flip{$[{4}:{1}]{\times}[{3}]$}\hss}\cg%
\cx{\hskip 16 true pt\flip{$[{3}{1}:{1}]{\times}[{3}]$}\hss}\cg%
\cx{\hskip 16 true pt\flip{$[{2^{2}}:{1}]{\times}[{3}]$}\hss}\cg%
\cx{\hskip 16 true pt\flip{$[{2}{1^{2}}:{1}]{\times}[{3}]$}\hss}\cg%
\cx{\hskip 16 true pt\flip{$[{1^{4}}:{1}]{\times}[{3}]$}\hss}\cg%
\cx{\hskip 16 true pt\flip{$[{3}:{2}]{\times}[{3}]$}\hss}\cg%
\cx{\hskip 16 true pt\flip{$[{3}:{1^{2}}]{\times}[{3}]$}\hss}\cg%
\cx{\hskip 16 true pt\flip{$[{2}{1}:{2}]{\times}[{3}]$}\hss}\cg%
\cx{\hskip 16 true pt\flip{$[{2}{1}:{1^{2}}]{\times}[{3}]$}\hss}\cg%
\cx{\hskip 16 true pt\flip{$[{1^{3}}:{2}]{\times}[{3}]$}\hss}\cg%
\cx{\hskip 16 true pt\flip{$[{1^{3}}:{1^{2}}]{\times}[{3}]$}\hss}\cg%
\eol}}\rg%
%
%
\rx{\vss\hfull{%
\rlx{\hss{$1_{x}$}}\cg%
\e{1}%
\e{0}%
\e{0}%
\e{0}%
\e{0}%
\e{0}%
\e{0}%
\e{0}%
\e{0}%
\e{0}%
\e{0}%
\e{0}%
\e{0}%
\e{0}%
\e{0}%
\e{0}%
\e{0}%
\e{0}%
\eol}\vss}\rg%
%
%
\rx{\vss\hfull{%
\rlx{\hss{$28_{x}$}}\cg%
\e{0}%
\e{0}%
\e{0}%
\e{0}%
\e{0}%
\e{0}%
\e{0}%
\e{1}%
\e{0}%
\e{0}%
\e{0}%
\e{0}%
\e{0}%
\e{1}%
\e{0}%
\e{0}%
\e{0}%
\e{0}%
\eol}\vss}\rg%
%
%
\rx{\vss\hfull{%
\rlx{\hss{$35_{x}$}}\cg%
\e{2}%
\e{1}%
\e{0}%
\e{0}%
\e{0}%
\e{0}%
\e{0}%
\e{1}%
\e{0}%
\e{0}%
\e{0}%
\e{0}%
\e{1}%
\e{0}%
\e{0}%
\e{0}%
\e{0}%
\e{0}%
\eol}\vss}\rg%
%
%
\rx{\vss\hfull{%
\rlx{\hss{$84_{x}$}}\cg%
\e{3}%
\e{1}%
\e{1}%
\e{0}%
\e{0}%
\e{0}%
\e{0}%
\e{2}%
\e{0}%
\e{0}%
\e{0}%
\e{0}%
\e{2}%
\e{0}%
\e{0}%
\e{0}%
\e{0}%
\e{0}%
\eol}\vss}\rg%
%
%
\rx{\vss\hfull{%
\rlx{\hss{$50_{x}$}}\cg%
\e{1}%
\e{1}%
\e{0}%
\e{0}%
\e{0}%
\e{0}%
\e{0}%
\e{1}%
\e{0}%
\e{0}%
\e{0}%
\e{0}%
\e{1}%
\e{0}%
\e{0}%
\e{0}%
\e{0}%
\e{0}%
\eol}\vss}\rg%
%
%
\rx{\vss\hfull{%
\rlx{\hss{$350_{x}$}}\cg%
\e{0}%
\e{0}%
\e{0}%
\e{0}%
\e{0}%
\e{0}%
\e{0}%
\e{1}%
\e{1}%
\e{0}%
\e{1}%
\e{0}%
\e{0}%
\e{2}%
\e{1}%
\e{1}%
\e{1}%
\e{0}%
\eol}\vss}\rg%
%
%
\rx{\vss\hfull{%
\rlx{\hss{$300_{x}$}}\cg%
\e{2}%
\e{2}%
\e{1}%
\e{0}%
\e{0}%
\e{0}%
\e{0}%
\e{1}%
\e{1}%
\e{1}%
\e{0}%
\e{0}%
\e{2}%
\e{0}%
\e{2}%
\e{0}%
\e{0}%
\e{0}%
\eol}\vss}\rg%
%
%
\rx{\vss\hfull{%
\rlx{\hss{$567_{x}$}}\cg%
\e{2}%
\e{2}%
\e{0}%
\e{1}%
\e{0}%
\e{0}%
\e{0}%
\e{5}%
\e{3}%
\e{0}%
\e{0}%
\e{0}%
\e{4}%
\e{3}%
\e{2}%
\e{1}%
\e{0}%
\e{0}%
\eol}\vss}\rg%
%
%
\rx{\vss\hfull{%
\rlx{\hss{$210_{x}$}}\cg%
\e{2}%
\e{2}%
\e{0}%
\e{0}%
\e{0}%
\e{0}%
\e{0}%
\e{3}%
\e{1}%
\e{0}%
\e{0}%
\e{0}%
\e{2}%
\e{1}%
\e{1}%
\e{0}%
\e{0}%
\e{0}%
\eol}\vss}\rg%
%
%
\rx{\vss\hfull{%
\rlx{\hss{$840_{x}$}}\cg%
\e{1}%
\e{2}%
\e{2}%
\e{0}%
\e{1}%
\e{0}%
\e{0}%
\e{1}%
\e{1}%
\e{3}%
\e{0}%
\e{0}%
\e{2}%
\e{0}%
\e{3}%
\e{1}%
\e{0}%
\e{0}%
\eol}\vss}\rg%
%
%
\rx{\vss\hfull{%
\rlx{\hss{$700_{x}$}}\cg%
\e{4}%
\e{4}%
\e{2}%
\e{0}%
\e{0}%
\e{0}%
\e{0}%
\e{5}%
\e{3}%
\e{1}%
\e{0}%
\e{0}%
\e{6}%
\e{1}%
\e{3}%
\e{0}%
\e{0}%
\e{0}%
\eol}\vss}\rg%
%
%
\rx{\vss\hfull{%
\rlx{\hss{$175_{x}$}}\cg%
\e{0}%
\e{0}%
\e{1}%
\e{0}%
\e{0}%
\e{0}%
\e{0}%
\e{1}%
\e{1}%
\e{0}%
\e{0}%
\e{0}%
\e{2}%
\e{0}%
\e{0}%
\e{0}%
\e{0}%
\e{0}%
\eol}\vss}\rg%
%
%
\rx{\vss\hfull{%
\rlx{\hss{$1400_{x}$}}\cg%
\e{1}%
\e{2}%
\e{1}%
\e{1}%
\e{0}%
\e{0}%
\e{0}%
\e{5}%
\e{6}%
\e{1}%
\e{1}%
\e{0}%
\e{6}%
\e{4}%
\e{4}%
\e{2}%
\e{1}%
\e{0}%
\eol}\vss}\rg%
%
%
\rx{\vss\hfull{%
\rlx{\hss{$1050_{x}$}}\cg%
\e{1}%
\e{2}%
\e{2}%
\e{1}%
\e{0}%
\e{0}%
\e{0}%
\e{4}%
\e{4}%
\e{1}%
\e{0}%
\e{0}%
\e{6}%
\e{3}%
\e{3}%
\e{1}%
\e{0}%
\e{0}%
\eol}\vss}\rg%
%
%
\rx{\vss\hfull{%
\rlx{\hss{$1575_{x}$}}\cg%
\e{1}%
\e{2}%
\e{0}%
\e{1}%
\e{0}%
\e{0}%
\e{0}%
\e{5}%
\e{6}%
\e{1}%
\e{2}%
\e{0}%
\e{4}%
\e{6}%
\e{5}%
\e{3}%
\e{2}%
\e{0}%
\eol}\vss}\rg%
%
%
\rx{\vss\hfull{%
\rlx{\hss{$1344_{x}$}}\cg%
\e{3}%
\e{4}%
\e{3}%
\e{1}%
\e{1}%
\e{0}%
\e{0}%
\e{6}%
\e{5}%
\e{2}%
\e{0}%
\e{0}%
\e{8}%
\e{3}%
\e{5}%
\e{2}%
\e{0}%
\e{0}%
\eol}\vss}\rg%
%
%
\rx{\vss\hfull{%
\rlx{\hss{$2100_{x}$}}\cg%
\e{0}%
\e{1}%
\e{0}%
\e{2}%
\e{0}%
\e{1}%
\e{0}%
\e{1}%
\e{4}%
\e{1}%
\e{4}%
\e{1}%
\e{2}%
\e{4}%
\e{5}%
\e{5}%
\e{4}%
\e{2}%
\eol}\vss}\rg%
%
%
\rx{\vss\hfull{%
\rlx{\hss{$2268_{x}$}}\cg%
\e{1}%
\e{4}%
\e{1}%
\e{2}%
\e{0}%
\e{0}%
\e{0}%
\e{4}%
\e{7}%
\e{2}%
\e{3}%
\e{0}%
\e{6}%
\e{4}%
\e{8}%
\e{4}%
\e{2}%
\e{1}%
\eol}\vss}\rg%
%
%
\rx{\vss\hfull{%
\rlx{\hss{$525_{x}$}}\cg%
\e{0}%
\e{0}%
\e{0}%
\e{1}%
\e{0}%
\e{1}%
\e{0}%
\e{2}%
\e{2}%
\e{0}%
\e{0}%
\e{0}%
\e{2}%
\e{2}%
\e{1}%
\e{2}%
\e{0}%
\e{0}%
\eol}\vss}\rg%
%
%
\rx{\vss\hfull{%
\rlx{\hss{$700_{xx}$}}\cg%
\e{0}%
\e{1}%
\e{0}%
\e{1}%
\e{0}%
\e{0}%
\e{0}%
\e{1}%
\e{2}%
\e{0}%
\e{0}%
\e{0}%
\e{2}%
\e{2}%
\e{2}%
\e{1}%
\e{0}%
\e{0}%
\eol}\vss}\rg%
%
%
\rx{\vss\hfull{%
\rlx{\hss{$972_{x}$}}\cg%
\e{1}%
\e{3}%
\e{3}%
\e{0}%
\e{1}%
\e{0}%
\e{0}%
\e{2}%
\e{2}%
\e{2}%
\e{0}%
\e{0}%
\e{4}%
\e{1}%
\e{4}%
\e{1}%
\e{0}%
\e{0}%
\eol}\vss}\rg%
%
%
\rx{\vss\hfull{%
\rlx{\hss{$4096_{x}$}}\cg%
\e{1}%
\e{4}%
\e{3}%
\e{3}%
\e{2}%
\e{0}%
\e{0}%
\e{5}%
\e{10}%
\e{5}%
\e{5}%
\e{0}%
\e{8}%
\e{7}%
\e{12}%
\e{8}%
\e{3}%
\e{2}%
\eol}\vss}\rg%
%
%
\rx{\vss\hfull{%
\rlx{\hss{$4200_{x}$}}\cg%
\e{1}%
\e{4}%
\e{4}%
\e{3}%
\e{1}%
\e{0}%
\e{0}%
\e{4}%
\e{10}%
\e{6}%
\e{4}%
\e{0}%
\e{9}%
\e{5}%
\e{12}%
\e{7}%
\e{2}%
\e{1}%
\eol}\vss}\rg%
%
%
\rx{\vss\hfull{%
\rlx{\hss{$2240_{x}$}}\cg%
\e{2}%
\e{3}%
\e{4}%
\e{1}%
\e{1}%
\e{0}%
\e{0}%
\e{5}%
\e{7}%
\e{4}%
\e{1}%
\e{0}%
\e{8}%
\e{3}%
\e{7}%
\e{3}%
\e{0}%
\e{0}%
\eol}\vss}\rg%
%
%
\rx{\vss\hfull{%
\rlx{\hss{$2835_{x}$}}\cg%
\e{0}%
\e{1}%
\e{4}%
\e{1}%
\e{1}%
\e{0}%
\e{0}%
\e{2}%
\e{7}%
\e{4}%
\e{2}%
\e{0}%
\e{6}%
\e{3}%
\e{7}%
\e{4}%
\e{1}%
\e{0}%
\eol}\vss}\rg%
%
%
\rx{\vss\hfull{%
\rlx{\hss{$6075_{x}$}}\cg%
\e{0}%
\e{2}%
\e{2}%
\e{4}%
\e{2}%
\e{2}%
\e{0}%
\e{4}%
\e{12}%
\e{5}%
\e{8}%
\e{1}%
\e{7}%
\e{10}%
\e{13}%
\e{13}%
\e{6}%
\e{3}%
\eol}\vss}\rg%
%
%
\rx{\vss\hfull{%
\rlx{\hss{$3200_{x}$}}\cg%
\e{0}%
\e{2}%
\e{4}%
\e{2}%
\e{4}%
\e{0}%
\e{0}%
\e{1}%
\e{4}%
\e{5}%
\e{3}%
\e{0}%
\e{4}%
\e{2}%
\e{8}%
\e{6}%
\e{1}%
\e{3}%
\eol}\vss}\rg%
\eop
\eject
\tablecont%
%
%
%
%
%
%
\rowpts=18 true pt%
\colpts=18 true pt%
\rowlabpts=40 true pt%
\collabpts=85 true pt%
\clx{\vss\hfull{%
\rlx{\hss{$ $}}\cg%
\cx{\hskip 16 true pt\flip{$[{5}:-]{\times}[{3}]$}\hss}\cg%
\cx{\hskip 16 true pt\flip{$[{4}{1}:-]{\times}[{3}]$}\hss}\cg%
\cx{\hskip 16 true pt\flip{$[{3}{2}:-]{\times}[{3}]$}\hss}\cg%
\cx{\hskip 16 true pt\flip{$[{3}{1^{2}}:-]{\times}[{3}]$}\hss}\cg%
\cx{\hskip 16 true pt\flip{$[{2^{2}}{1}:-]{\times}[{3}]$}\hss}\cg%
\cx{\hskip 16 true pt\flip{$[{2}{1^{3}}:-]{\times}[{3}]$}\hss}\cg%
\cx{\hskip 16 true pt\flip{$[{1^{5}}:-]{\times}[{3}]$}\hss}\cg%
\cx{\hskip 16 true pt\flip{$[{4}:{1}]{\times}[{3}]$}\hss}\cg%
\cx{\hskip 16 true pt\flip{$[{3}{1}:{1}]{\times}[{3}]$}\hss}\cg%
\cx{\hskip 16 true pt\flip{$[{2^{2}}:{1}]{\times}[{3}]$}\hss}\cg%
\cx{\hskip 16 true pt\flip{$[{2}{1^{2}}:{1}]{\times}[{3}]$}\hss}\cg%
\cx{\hskip 16 true pt\flip{$[{1^{4}}:{1}]{\times}[{3}]$}\hss}\cg%
\cx{\hskip 16 true pt\flip{$[{3}:{2}]{\times}[{3}]$}\hss}\cg%
\cx{\hskip 16 true pt\flip{$[{3}:{1^{2}}]{\times}[{3}]$}\hss}\cg%
\cx{\hskip 16 true pt\flip{$[{2}{1}:{2}]{\times}[{3}]$}\hss}\cg%
\cx{\hskip 16 true pt\flip{$[{2}{1}:{1^{2}}]{\times}[{3}]$}\hss}\cg%
\cx{\hskip 16 true pt\flip{$[{1^{3}}:{2}]{\times}[{3}]$}\hss}\cg%
\cx{\hskip 16 true pt\flip{$[{1^{3}}:{1^{2}}]{\times}[{3}]$}\hss}\cg%
\eol}}\rg%
%
%
\rx{\vss\hfull{%
\rlx{\hss{$70_{y}$}}\cg%
\e{0}%
\e{0}%
\e{0}%
\e{0}%
\e{0}%
\e{0}%
\e{0}%
\e{0}%
\e{0}%
\e{0}%
\e{0}%
\e{1}%
\e{0}%
\e{0}%
\e{0}%
\e{0}%
\e{1}%
\e{0}%
\eol}\vss}\rg%
%
%
\rx{\vss\hfull{%
\rlx{\hss{$1134_{y}$}}\cg%
\e{0}%
\e{0}%
\e{0}%
\e{0}%
\e{1}%
\e{1}%
\e{0}%
\e{1}%
\e{1}%
\e{1}%
\e{2}%
\e{0}%
\e{0}%
\e{2}%
\e{1}%
\e{3}%
\e{1}%
\e{1}%
\eol}\vss}\rg%
%
%
\rx{\vss\hfull{%
\rlx{\hss{$1680_{y}$}}\cg%
\e{0}%
\e{0}%
\e{0}%
\e{1}%
\e{0}%
\e{1}%
\e{0}%
\e{0}%
\e{1}%
\e{0}%
\e{4}%
\e{3}%
\e{0}%
\e{2}%
\e{2}%
\e{3}%
\e{5}%
\e{3}%
\eol}\vss}\rg%
%
%
\rx{\vss\hfull{%
\rlx{\hss{$168_{y}$}}\cg%
\e{0}%
\e{1}%
\e{0}%
\e{0}%
\e{0}%
\e{0}%
\e{0}%
\e{0}%
\e{0}%
\e{0}%
\e{0}%
\e{0}%
\e{0}%
\e{0}%
\e{1}%
\e{0}%
\e{0}%
\e{0}%
\eol}\vss}\rg%
%
%
\rx{\vss\hfull{%
\rlx{\hss{$420_{y}$}}\cg%
\e{0}%
\e{0}%
\e{1}%
\e{0}%
\e{0}%
\e{0}%
\e{0}%
\e{0}%
\e{1}%
\e{1}%
\e{0}%
\e{0}%
\e{1}%
\e{0}%
\e{1}%
\e{0}%
\e{0}%
\e{0}%
\eol}\vss}\rg%
%
%
\rx{\vss\hfull{%
\rlx{\hss{$3150_{y}$}}\cg%
\e{0}%
\e{0}%
\e{2}%
\e{1}%
\e{2}%
\e{1}%
\e{0}%
\e{1}%
\e{5}%
\e{5}%
\e{3}%
\e{0}%
\e{3}%
\e{2}%
\e{5}%
\e{6}%
\e{1}%
\e{1}%
\eol}\vss}\rg%
%
%
\rx{\vss\hfull{%
\rlx{\hss{$4200_{y}$}}\cg%
\e{0}%
\e{2}%
\e{3}%
\e{2}%
\e{3}%
\e{0}%
\e{0}%
\e{1}%
\e{5}%
\e{5}%
\e{4}%
\e{0}%
\e{3}%
\e{3}%
\e{9}%
\e{7}%
\e{2}%
\e{2}%
\eol}\vss}\rg%
%
%
\rx{\vss\hfull{%
\rlx{\hss{$2688_{y}$}}\cg%
\e{0}%
\e{0}%
\e{2}%
\e{2}%
\e{2}%
\e{0}%
\e{0}%
\e{0}%
\e{3}%
\e{3}%
\e{3}%
\e{0}%
\e{2}%
\e{2}%
\e{5}%
\e{5}%
\e{2}%
\e{2}%
\eol}\vss}\rg%
%
%
\rx{\vss\hfull{%
\rlx{\hss{$2100_{y}$}}\cg%
\e{0}%
\e{0}%
\e{0}%
\e{2}%
\e{0}%
\e{3}%
\e{1}%
\e{0}%
\e{2}%
\e{0}%
\e{3}%
\e{2}%
\e{1}%
\e{2}%
\e{3}%
\e{5}%
\e{3}%
\e{3}%
\eol}\vss}\rg%
%
%
\rx{\vss\hfull{%
\rlx{\hss{$1400_{y}$}}\cg%
\e{0}%
\e{1}%
\e{0}%
\e{1}%
\e{0}%
\e{0}%
\e{0}%
\e{0}%
\e{1}%
\e{0}%
\e{3}%
\e{2}%
\e{0}%
\e{1}%
\e{3}%
\e{2}%
\e{3}%
\e{2}%
\eol}\vss}\rg%
%
%
\rx{\vss\hfull{%
\rlx{\hss{$4536_{y}$}}\cg%
\e{0}%
\e{1}%
\e{1}%
\e{3}%
\e{1}%
\e{1}%
\e{0}%
\e{0}%
\e{5}%
\e{3}%
\e{8}%
\e{3}%
\e{2}%
\e{3}%
\e{8}%
\e{8}%
\e{6}%
\e{5}%
\eol}\vss}\rg%
%
%
\rx{\vss\hfull{%
\rlx{\hss{$5670_{y}$}}\cg%
\e{0}%
\e{1}%
\e{1}%
\e{3}%
\e{2}%
\e{2}%
\e{0}%
\e{1}%
\e{6}%
\e{4}%
\e{10}%
\e{3}%
\e{2}%
\e{5}%
\e{9}%
\e{11}%
\e{7}%
\e{6}%
\eol}\vss}\rg%
%
%
\rx{\vss\hfull{%
\rlx{\hss{$4480_{y}$}}\cg%
\e{0}%
\e{1}%
\e{2}%
\e{2}%
\e{2}%
\e{1}%
\e{0}%
\e{1}%
\e{6}%
\e{5}%
\e{6}%
\e{1}%
\e{3}%
\e{3}%
\e{8}%
\e{8}%
\e{3}%
\e{3}%
\eol}\vss}\rg%
%
%
\rx{\vss\hfull{%
\rlx{\hss{$8_{z}$}}\cg%
\e{1}%
\e{0}%
\e{0}%
\e{0}%
\e{0}%
\e{0}%
\e{0}%
\e{1}%
\e{0}%
\e{0}%
\e{0}%
\e{0}%
\e{0}%
\e{0}%
\e{0}%
\e{0}%
\e{0}%
\e{0}%
\eol}\vss}\rg%
%
%
\rx{\vss\hfull{%
\rlx{\hss{$56_{z}$}}\cg%
\e{0}%
\e{0}%
\e{0}%
\e{0}%
\e{0}%
\e{0}%
\e{0}%
\e{0}%
\e{0}%
\e{0}%
\e{0}%
\e{0}%
\e{0}%
\e{1}%
\e{0}%
\e{0}%
\e{1}%
\e{0}%
\eol}\vss}\rg%
%
%
\rx{\vss\hfull{%
\rlx{\hss{$160_{z}$}}\cg%
\e{1}%
\e{1}%
\e{0}%
\e{0}%
\e{0}%
\e{0}%
\e{0}%
\e{2}%
\e{1}%
\e{0}%
\e{0}%
\e{0}%
\e{1}%
\e{1}%
\e{1}%
\e{0}%
\e{0}%
\e{0}%
\eol}\vss}\rg%
%
%
\rx{\vss\hfull{%
\rlx{\hss{$112_{z}$}}\cg%
\e{3}%
\e{1}%
\e{0}%
\e{0}%
\e{0}%
\e{0}%
\e{0}%
\e{3}%
\e{1}%
\e{0}%
\e{0}%
\e{0}%
\e{2}%
\e{0}%
\e{0}%
\e{0}%
\e{0}%
\e{0}%
\eol}\vss}\rg%
%
%
\rx{\vss\hfull{%
\rlx{\hss{$840_{z}$}}\cg%
\e{1}%
\e{1}%
\e{1}%
\e{0}%
\e{1}%
\e{0}%
\e{0}%
\e{2}%
\e{2}%
\e{2}%
\e{1}%
\e{0}%
\e{1}%
\e{1}%
\e{3}%
\e{2}%
\e{0}%
\e{0}%
\eol}\vss}\rg%
%
%
\rx{\vss\hfull{%
\rlx{\hss{$1296_{z}$}}\cg%
\e{0}%
\e{1}%
\e{0}%
\e{2}%
\e{0}%
\e{0}%
\e{0}%
\e{2}%
\e{3}%
\e{0}%
\e{2}%
\e{0}%
\e{2}%
\e{5}%
\e{4}%
\e{3}%
\e{3}%
\e{1}%
\eol}\vss}\rg%
%
%
\rx{\vss\hfull{%
\rlx{\hss{$1400_{z}$}}\cg%
\e{3}%
\e{4}%
\e{2}%
\e{1}%
\e{0}%
\e{0}%
\e{0}%
\e{6}%
\e{6}%
\e{2}%
\e{1}%
\e{0}%
\e{8}%
\e{3}%
\e{5}%
\e{2}%
\e{0}%
\e{0}%
\eol}\vss}\rg%
%
%
\rx{\vss\hfull{%
\rlx{\hss{$1008_{z}$}}\cg%
\e{1}%
\e{3}%
\e{0}%
\e{1}%
\e{0}%
\e{0}%
\e{0}%
\e{4}%
\e{4}%
\e{0}%
\e{1}%
\e{0}%
\e{4}%
\e{4}%
\e{4}%
\e{2}%
\e{1}%
\e{0}%
\eol}\vss}\rg%
%
%
\rx{\vss\hfull{%
\rlx{\hss{$560_{z}$}}\cg%
\e{4}%
\e{3}%
\e{1}%
\e{0}%
\e{0}%
\e{0}%
\e{0}%
\e{6}%
\e{3}%
\e{1}%
\e{0}%
\e{0}%
\e{5}%
\e{2}%
\e{2}%
\e{0}%
\e{0}%
\e{0}%
\eol}\vss}\rg%
%
%
\rx{\vss\hfull{%
\rlx{\hss{$1400_{zz}$}}\cg%
\e{1}%
\e{2}%
\e{3}%
\e{1}%
\e{0}%
\e{0}%
\e{0}%
\e{3}%
\e{5}%
\e{2}%
\e{0}%
\e{0}%
\e{7}%
\e{2}%
\e{4}%
\e{1}%
\e{0}%
\e{0}%
\eol}\vss}\rg%
%
%
\rx{\vss\hfull{%
\rlx{\hss{$4200_{z}$}}\cg%
\e{0}%
\e{2}%
\e{2}%
\e{3}%
\e{1}%
\e{1}%
\e{0}%
\e{3}%
\e{8}%
\e{3}%
\e{4}%
\e{0}%
\e{5}%
\e{6}%
\e{10}%
\e{8}%
\e{3}%
\e{1}%
\eol}\vss}\rg%
%
%
\rx{\vss\hfull{%
\rlx{\hss{$400_{z}$}}\cg%
\e{2}%
\e{2}%
\e{1}%
\e{0}%
\e{0}%
\e{0}%
\e{0}%
\e{4}%
\e{2}%
\e{0}%
\e{0}%
\e{0}%
\e{4}%
\e{1}%
\e{1}%
\e{0}%
\e{0}%
\e{0}%
\eol}\vss}\rg%
%
%
\rx{\vss\hfull{%
\rlx{\hss{$3240_{z}$}}\cg%
\e{3}%
\e{6}%
\e{4}%
\e{2}%
\e{1}%
\e{0}%
\e{0}%
\e{9}%
\e{11}%
\e{4}%
\e{2}%
\e{0}%
\e{12}%
\e{7}%
\e{11}%
\e{5}%
\e{1}%
\e{0}%
\eol}\vss}\rg%
%
%
\rx{\vss\hfull{%
\rlx{\hss{$4536_{z}$}}\cg%
\e{1}%
\e{4}%
\e{6}%
\e{1}%
\e{3}%
\e{0}%
\e{0}%
\e{4}%
\e{10}%
\e{9}%
\e{4}%
\e{0}%
\e{9}%
\e{4}%
\e{12}%
\e{7}%
\e{1}%
\e{1}%
\eol}\vss}\rg%
%
%
\rx{\vss\hfull{%
\rlx{\hss{$2400_{z}$}}\cg%
\e{0}%
\e{0}%
\e{0}%
\e{2}%
\e{0}%
\e{2}%
\e{0}%
\e{1}%
\e{3}%
\e{0}%
\e{4}%
\e{2}%
\e{1}%
\e{5}%
\e{4}%
\e{6}%
\e{6}%
\e{2}%
\eol}\vss}\rg%
%
%
\rx{\vss\hfull{%
\rlx{\hss{$3360_{z}$}}\cg%
\e{0}%
\e{2}%
\e{3}%
\e{3}%
\e{1}%
\e{0}%
\e{0}%
\e{3}%
\e{8}%
\e{3}%
\e{3}%
\e{0}%
\e{7}%
\e{6}%
\e{9}%
\e{6}%
\e{3}%
\e{1}%
\eol}\vss}\rg%
%
%
\rx{\vss\hfull{%
\rlx{\hss{$2800_{z}$}}\cg%
\e{0}%
\e{2}%
\e{0}%
\e{3}%
\e{0}%
\e{1}%
\e{0}%
\e{2}%
\e{6}%
\e{1}%
\e{4}%
\e{1}%
\e{4}%
\e{6}%
\e{7}%
\e{6}%
\e{4}%
\e{2}%
\eol}\vss}\rg%
%
%
\rx{\vss\hfull{%
\rlx{\hss{$4096_{z}$}}\cg%
\e{1}%
\e{4}%
\e{3}%
\e{3}%
\e{2}%
\e{0}%
\e{0}%
\e{5}%
\e{10}%
\e{5}%
\e{5}%
\e{0}%
\e{8}%
\e{7}%
\e{12}%
\e{8}%
\e{3}%
\e{2}%
\eol}\vss}\rg%
%
%
\rx{\vss\hfull{%
\rlx{\hss{$5600_{z}$}}\cg%
\e{0}%
\e{2}%
\e{2}%
\e{4}%
\e{2}%
\e{2}%
\e{0}%
\e{2}%
\e{9}%
\e{5}%
\e{9}%
\e{2}%
\e{5}%
\e{5}%
\e{12}%
\e{12}%
\e{5}%
\e{5}%
\eol}\vss}\rg%
%
%
\rx{\vss\hfull{%
\rlx{\hss{$448_{z}$}}\cg%
\e{2}%
\e{1}%
\e{2}%
\e{0}%
\e{0}%
\e{0}%
\e{0}%
\e{2}%
\e{2}%
\e{2}%
\e{0}%
\e{0}%
\e{4}%
\e{0}%
\e{1}%
\e{0}%
\e{0}%
\e{0}%
\eol}\vss}\rg%
\eop
\eject
\tablecont%
%
%
%
%
%
%
\rowpts=18 true pt%
\colpts=18 true pt%
\rowlabpts=40 true pt%
\collabpts=85 true pt%
\clx{\vss\hfull{%
\rlx{\hss{$ $}}\cg%
\cx{\hskip 16 true pt\flip{$[{5}:-]{\times}[{3}]$}\hss}\cg%
\cx{\hskip 16 true pt\flip{$[{4}{1}:-]{\times}[{3}]$}\hss}\cg%
\cx{\hskip 16 true pt\flip{$[{3}{2}:-]{\times}[{3}]$}\hss}\cg%
\cx{\hskip 16 true pt\flip{$[{3}{1^{2}}:-]{\times}[{3}]$}\hss}\cg%
\cx{\hskip 16 true pt\flip{$[{2^{2}}{1}:-]{\times}[{3}]$}\hss}\cg%
\cx{\hskip 16 true pt\flip{$[{2}{1^{3}}:-]{\times}[{3}]$}\hss}\cg%
\cx{\hskip 16 true pt\flip{$[{1^{5}}:-]{\times}[{3}]$}\hss}\cg%
\cx{\hskip 16 true pt\flip{$[{4}:{1}]{\times}[{3}]$}\hss}\cg%
\cx{\hskip 16 true pt\flip{$[{3}{1}:{1}]{\times}[{3}]$}\hss}\cg%
\cx{\hskip 16 true pt\flip{$[{2^{2}}:{1}]{\times}[{3}]$}\hss}\cg%
\cx{\hskip 16 true pt\flip{$[{2}{1^{2}}:{1}]{\times}[{3}]$}\hss}\cg%
\cx{\hskip 16 true pt\flip{$[{1^{4}}:{1}]{\times}[{3}]$}\hss}\cg%
\cx{\hskip 16 true pt\flip{$[{3}:{2}]{\times}[{3}]$}\hss}\cg%
\cx{\hskip 16 true pt\flip{$[{3}:{1^{2}}]{\times}[{3}]$}\hss}\cg%
\cx{\hskip 16 true pt\flip{$[{2}{1}:{2}]{\times}[{3}]$}\hss}\cg%
\cx{\hskip 16 true pt\flip{$[{2}{1}:{1^{2}}]{\times}[{3}]$}\hss}\cg%
\cx{\hskip 16 true pt\flip{$[{1^{3}}:{2}]{\times}[{3}]$}\hss}\cg%
\cx{\hskip 16 true pt\flip{$[{1^{3}}:{1^{2}}]{\times}[{3}]$}\hss}\cg%
\eol}}\rg%
%
%
\rx{\vss\hfull{%
\rlx{\hss{$448_{w}$}}\cg%
\e{0}%
\e{0}%
\e{0}%
\e{0}%
\e{0}%
\e{1}%
\e{0}%
\e{0}%
\e{0}%
\e{0}%
\e{1}%
\e{1}%
\e{0}%
\e{1}%
\e{0}%
\e{1}%
\e{2}%
\e{1}%
\eol}\vss}\rg%
%
%
\rx{\vss\hfull{%
\rlx{\hss{$1344_{w}$}}\cg%
\e{0}%
\e{1}%
\e{1}%
\e{0}%
\e{1}%
\e{0}%
\e{0}%
\e{1}%
\e{2}%
\e{2}%
\e{1}%
\e{0}%
\e{1}%
\e{1}%
\e{3}%
\e{2}%
\e{0}%
\e{0}%
\eol}\vss}\rg%
%
%
\rx{\vss\hfull{%
\rlx{\hss{$5600_{w}$}}\cg%
\e{0}%
\e{1}%
\e{1}%
\e{3}%
\e{3}%
\e{2}%
\e{0}%
\e{1}%
\e{6}%
\e{4}%
\e{9}%
\e{2}%
\e{2}%
\e{5}%
\e{9}%
\e{11}%
\e{7}%
\e{6}%
\eol}\vss}\rg%
%
%
\rx{\vss\hfull{%
\rlx{\hss{$2016_{w}$}}\cg%
\e{0}%
\e{0}%
\e{3}%
\e{1}%
\e{1}%
\e{0}%
\e{0}%
\e{0}%
\e{3}%
\e{4}%
\e{1}%
\e{0}%
\e{3}%
\e{1}%
\e{4}%
\e{3}%
\e{0}%
\e{0}%
\eol}\vss}\rg%
%
%
\rx{\vss\hfull{%
\rlx{\hss{$7168_{w}$}}\cg%
\e{0}%
\e{1}%
\e{4}%
\e{4}%
\e{4}%
\e{1}%
\e{0}%
\e{1}%
\e{9}%
\e{8}%
\e{9}%
\e{1}%
\e{5}%
\e{5}%
\e{13}%
\e{13}%
\e{5}%
\e{5}%
\eol}\vss}\rg%
%
%
%
%
%
%
\rowpts=18 true pt%
\colpts=18 true pt%
\rowlabpts=40 true pt%
\collabpts=85 true pt%
\clx{\vss\hfull{%
\rlx{\hss{$ $}}\cg%
\cx{\hskip 16 true pt\flip{$[{5}:-]{\times}[{2}{1}]$}\hss}\cg%
\cx{\hskip 16 true pt\flip{$[{4}{1}:-]{\times}[{2}{1}]$}\hss}\cg%
\cx{\hskip 16 true pt\flip{$[{3}{2}:-]{\times}[{2}{1}]$}\hss}\cg%
\cx{\hskip 16 true pt\flip{$[{3}{1^{2}}:-]{\times}[{2}{1}]$}\hss}\cg%
\cx{\hskip 16 true pt\flip{$[{2^{2}}{1}:-]{\times}[{2}{1}]$}\hss}\cg%
\cx{\hskip 16 true pt\flip{$[{2}{1^{3}}:-]{\times}[{2}{1}]$}\hss}\cg%
\cx{\hskip 16 true pt\flip{$[{1^{5}}:-]{\times}[{2}{1}]$}\hss}\cg%
\cx{\hskip 16 true pt\flip{$[{4}:{1}]{\times}[{2}{1}]$}\hss}\cg%
\cx{\hskip 16 true pt\flip{$[{3}{1}:{1}]{\times}[{2}{1}]$}\hss}\cg%
\cx{\hskip 16 true pt\flip{$[{2^{2}}:{1}]{\times}[{2}{1}]$}\hss}\cg%
\cx{\hskip 16 true pt\flip{$[{2}{1^{2}}:{1}]{\times}[{2}{1}]$}\hss}\cg%
\cx{\hskip 16 true pt\flip{$[{1^{4}}:{1}]{\times}[{2}{1}]$}\hss}\cg%
\cx{\hskip 16 true pt\flip{$[{3}:{2}]{\times}[{2}{1}]$}\hss}\cg%
\cx{\hskip 16 true pt\flip{$[{3}:{1^{2}}]{\times}[{2}{1}]$}\hss}\cg%
\cx{\hskip 16 true pt\flip{$[{2}{1}:{2}]{\times}[{2}{1}]$}\hss}\cg%
\cx{\hskip 16 true pt\flip{$[{2}{1}:{1^{2}}]{\times}[{2}{1}]$}\hss}\cg%
\cx{\hskip 16 true pt\flip{$[{1^{3}}:{2}]{\times}[{2}{1}]$}\hss}\cg%
\cx{\hskip 16 true pt\flip{$[{1^{3}}:{1^{2}}]{\times}[{2}{1}]$}\hss}\cg%
\eol}}\rg%
%
%
\rx{\vss\hfull{%
\rlx{\hss{$1_{x}$}}\cg%
\e{0}%
\e{0}%
\e{0}%
\e{0}%
\e{0}%
\e{0}%
\e{0}%
\e{0}%
\e{0}%
\e{0}%
\e{0}%
\e{0}%
\e{0}%
\e{0}%
\e{0}%
\e{0}%
\e{0}%
\e{0}%
\eol}\vss}\rg%
%
%
\rx{\vss\hfull{%
\rlx{\hss{$28_{x}$}}\cg%
\e{1}%
\e{0}%
\e{0}%
\e{0}%
\e{0}%
\e{0}%
\e{0}%
\e{1}%
\e{0}%
\e{0}%
\e{0}%
\e{0}%
\e{0}%
\e{0}%
\e{0}%
\e{0}%
\e{0}%
\e{0}%
\eol}\vss}\rg%
%
%
\rx{\vss\hfull{%
\rlx{\hss{$35_{x}$}}\cg%
\e{2}%
\e{0}%
\e{0}%
\e{0}%
\e{0}%
\e{0}%
\e{0}%
\e{1}%
\e{0}%
\e{0}%
\e{0}%
\e{0}%
\e{0}%
\e{0}%
\e{0}%
\e{0}%
\e{0}%
\e{0}%
\eol}\vss}\rg%
%
%
\rx{\vss\hfull{%
\rlx{\hss{$84_{x}$}}\cg%
\e{2}%
\e{1}%
\e{0}%
\e{0}%
\e{0}%
\e{0}%
\e{0}%
\e{1}%
\e{0}%
\e{0}%
\e{0}%
\e{0}%
\e{1}%
\e{0}%
\e{0}%
\e{0}%
\e{0}%
\e{0}%
\eol}\vss}\rg%
%
%
\rx{\vss\hfull{%
\rlx{\hss{$50_{x}$}}\cg%
\e{0}%
\e{0}%
\e{1}%
\e{0}%
\e{0}%
\e{0}%
\e{0}%
\e{0}%
\e{0}%
\e{0}%
\e{0}%
\e{0}%
\e{1}%
\e{0}%
\e{0}%
\e{0}%
\e{0}%
\e{0}%
\eol}\vss}\rg%
%
%
\rx{\vss\hfull{%
\rlx{\hss{$350_{x}$}}\cg%
\e{1}%
\e{1}%
\e{0}%
\e{0}%
\e{0}%
\e{0}%
\e{0}%
\e{3}%
\e{1}%
\e{0}%
\e{0}%
\e{0}%
\e{1}%
\e{3}%
\e{1}%
\e{0}%
\e{1}%
\e{0}%
\eol}\vss}\rg%
%
%
\rx{\vss\hfull{%
\rlx{\hss{$300_{x}$}}\cg%
\e{2}%
\e{2}%
\e{0}%
\e{0}%
\e{0}%
\e{0}%
\e{0}%
\e{3}%
\e{1}%
\e{0}%
\e{0}%
\e{0}%
\e{2}%
\e{1}%
\e{1}%
\e{0}%
\e{0}%
\e{0}%
\eol}\vss}\rg%
%
%
\rx{\vss\hfull{%
\rlx{\hss{$567_{x}$}}\cg%
\e{5}%
\e{3}%
\e{0}%
\e{0}%
\e{0}%
\e{0}%
\e{0}%
\e{7}%
\e{2}%
\e{0}%
\e{0}%
\e{0}%
\e{4}%
\e{2}%
\e{1}%
\e{0}%
\e{0}%
\e{0}%
\eol}\vss}\rg%
%
%
\rx{\vss\hfull{%
\rlx{\hss{$210_{x}$}}\cg%
\e{3}%
\e{1}%
\e{0}%
\e{0}%
\e{0}%
\e{0}%
\e{0}%
\e{3}%
\e{1}%
\e{0}%
\e{0}%
\e{0}%
\e{2}%
\e{0}%
\e{0}%
\e{0}%
\e{0}%
\e{0}%
\eol}\vss}\rg%
%
%
\rx{\vss\hfull{%
\rlx{\hss{$840_{x}$}}\cg%
\e{0}%
\e{2}%
\e{2}%
\e{0}%
\e{1}%
\e{0}%
\e{0}%
\e{1}%
\e{3}%
\e{3}%
\e{1}%
\e{0}%
\e{2}%
\e{1}%
\e{5}%
\e{2}%
\e{0}%
\e{0}%
\eol}\vss}\rg%
%
%
\rx{\vss\hfull{%
\rlx{\hss{$700_{x}$}}\cg%
\e{3}%
\e{3}%
\e{2}%
\e{0}%
\e{0}%
\e{0}%
\e{0}%
\e{5}%
\e{3}%
\e{1}%
\e{0}%
\e{0}%
\e{6}%
\e{1}%
\e{2}%
\e{0}%
\e{0}%
\e{0}%
\eol}\vss}\rg%
%
%
\rx{\vss\hfull{%
\rlx{\hss{$175_{x}$}}\cg%
\e{0}%
\e{0}%
\e{1}%
\e{0}%
\e{0}%
\e{0}%
\e{0}%
\e{0}%
\e{1}%
\e{1}%
\e{0}%
\e{0}%
\e{1}%
\e{0}%
\e{1}%
\e{0}%
\e{0}%
\e{0}%
\eol}\vss}\rg%
%
%
\rx{\vss\hfull{%
\rlx{\hss{$1400_{x}$}}\cg%
\e{2}%
\e{3}%
\e{2}%
\e{1}%
\e{0}%
\e{0}%
\e{0}%
\e{6}%
\e{7}%
\e{2}%
\e{1}%
\e{0}%
\e{8}%
\e{3}%
\e{5}%
\e{2}%
\e{0}%
\e{0}%
\eol}\vss}\rg%
%
%
\rx{\vss\hfull{%
\rlx{\hss{$1050_{x}$}}\cg%
\e{1}%
\e{2}%
\e{4}%
\e{1}%
\e{1}%
\e{0}%
\e{0}%
\e{3}%
\e{4}%
\e{2}%
\e{0}%
\e{0}%
\e{7}%
\e{2}%
\e{4}%
\e{2}%
\e{0}%
\e{0}%
\eol}\vss}\rg%
%
%
\rx{\vss\hfull{%
\rlx{\hss{$1575_{x}$}}\cg%
\e{3}%
\e{4}%
\e{1}%
\e{1}%
\e{0}%
\e{0}%
\e{0}%
\e{9}%
\e{7}%
\e{1}%
\e{1}%
\e{0}%
\e{8}%
\e{6}%
\e{5}%
\e{2}%
\e{1}%
\e{0}%
\eol}\vss}\rg%
%
%
\rx{\vss\hfull{%
\rlx{\hss{$1344_{x}$}}\cg%
\e{5}%
\e{6}%
\e{2}%
\e{1}%
\e{0}%
\e{0}%
\e{0}%
\e{8}%
\e{5}%
\e{1}%
\e{0}%
\e{0}%
\e{9}%
\e{4}%
\e{5}%
\e{1}%
\e{0}%
\e{0}%
\eol}\vss}\rg%
%
%
\rx{\vss\hfull{%
\rlx{\hss{$2100_{x}$}}\cg%
\e{1}%
\e{4}%
\e{0}%
\e{3}%
\e{0}%
\e{0}%
\e{0}%
\e{6}%
\e{7}%
\e{0}%
\e{3}%
\e{0}%
\e{6}%
\e{9}%
\e{8}%
\e{5}%
\e{4}%
\e{1}%
\eol}\vss}\rg%
%
%
\rx{\vss\hfull{%
\rlx{\hss{$2268_{x}$}}\cg%
\e{3}%
\e{6}%
\e{1}%
\e{2}%
\e{0}%
\e{0}%
\e{0}%
\e{9}%
\e{10}%
\e{2}%
\e{2}%
\e{0}%
\e{10}%
\e{7}%
\e{9}%
\e{4}%
\e{1}%
\e{0}%
\eol}\vss}\rg%
%
%
\rx{\vss\hfull{%
\rlx{\hss{$525_{x}$}}\cg%
\e{1}%
\e{2}%
\e{0}%
\e{1}%
\e{0}%
\e{0}%
\e{0}%
\e{2}%
\e{2}%
\e{0}%
\e{0}%
\e{0}%
\e{3}%
\e{2}%
\e{2}%
\e{1}%
\e{0}%
\e{0}%
\eol}\vss}\rg%
%
%
\rx{\vss\hfull{%
\rlx{\hss{$700_{xx}$}}\cg%
\e{0}%
\e{0}%
\e{3}%
\e{1}%
\e{1}%
\e{1}%
\e{0}%
\e{0}%
\e{2}%
\e{2}%
\e{0}%
\e{0}%
\e{3}%
\e{1}%
\e{3}%
\e{3}%
\e{0}%
\e{0}%
\eol}\vss}\rg%
%
%
\rx{\vss\hfull{%
\rlx{\hss{$972_{x}$}}\cg%
\e{1}%
\e{3}%
\e{3}%
\e{1}%
\e{1}%
\e{0}%
\e{0}%
\e{2}%
\e{3}%
\e{2}%
\e{0}%
\e{0}%
\e{5}%
\e{2}%
\e{5}%
\e{2}%
\e{0}%
\e{0}%
\eol}\vss}\rg%
\eop
\eject
\tablecont%
%
%
%
%
%
%
\rowpts=18 true pt%
\colpts=18 true pt%
\rowlabpts=40 true pt%
\collabpts=85 true pt%
\clx{\vss\hfull{%
\rlx{\hss{$ $}}\cg%
\cx{\hskip 16 true pt\flip{$[{5}:-]{\times}[{2}{1}]$}\hss}\cg%
\cx{\hskip 16 true pt\flip{$[{4}{1}:-]{\times}[{2}{1}]$}\hss}\cg%
\cx{\hskip 16 true pt\flip{$[{3}{2}:-]{\times}[{2}{1}]$}\hss}\cg%
\cx{\hskip 16 true pt\flip{$[{3}{1^{2}}:-]{\times}[{2}{1}]$}\hss}\cg%
\cx{\hskip 16 true pt\flip{$[{2^{2}}{1}:-]{\times}[{2}{1}]$}\hss}\cg%
\cx{\hskip 16 true pt\flip{$[{2}{1^{3}}:-]{\times}[{2}{1}]$}\hss}\cg%
\cx{\hskip 16 true pt\flip{$[{1^{5}}:-]{\times}[{2}{1}]$}\hss}\cg%
\cx{\hskip 16 true pt\flip{$[{4}:{1}]{\times}[{2}{1}]$}\hss}\cg%
\cx{\hskip 16 true pt\flip{$[{3}{1}:{1}]{\times}[{2}{1}]$}\hss}\cg%
\cx{\hskip 16 true pt\flip{$[{2^{2}}:{1}]{\times}[{2}{1}]$}\hss}\cg%
\cx{\hskip 16 true pt\flip{$[{2}{1^{2}}:{1}]{\times}[{2}{1}]$}\hss}\cg%
\cx{\hskip 16 true pt\flip{$[{1^{4}}:{1}]{\times}[{2}{1}]$}\hss}\cg%
\cx{\hskip 16 true pt\flip{$[{3}:{2}]{\times}[{2}{1}]$}\hss}\cg%
\cx{\hskip 16 true pt\flip{$[{3}:{1^{2}}]{\times}[{2}{1}]$}\hss}\cg%
\cx{\hskip 16 true pt\flip{$[{2}{1}:{2}]{\times}[{2}{1}]$}\hss}\cg%
\cx{\hskip 16 true pt\flip{$[{2}{1}:{1^{2}}]{\times}[{2}{1}]$}\hss}\cg%
\cx{\hskip 16 true pt\flip{$[{1^{3}}:{2}]{\times}[{2}{1}]$}\hss}\cg%
\cx{\hskip 16 true pt\flip{$[{1^{3}}:{1^{2}}]{\times}[{2}{1}]$}\hss}\cg%
\eol}}\rg%
%
%
\rx{\vss\hfull{%
\rlx{\hss{$4096_{x}$}}\cg%
\e{3}%
\e{8}%
\e{4}%
\e{4}%
\e{1}%
\e{0}%
\e{0}%
\e{11}%
\e{16}%
\e{5}%
\e{5}%
\e{0}%
\e{15}%
\e{12}%
\e{17}%
\e{9}%
\e{4}%
\e{1}%
\eol}\vss}\rg%
%
%
\rx{\vss\hfull{%
\rlx{\hss{$4200_{x}$}}\cg%
\e{1}%
\e{5}%
\e{7}%
\e{4}%
\e{3}%
\e{0}%
\e{0}%
\e{7}%
\e{16}%
\e{9}%
\e{6}%
\e{0}%
\e{14}%
\e{9}%
\e{18}%
\e{11}%
\e{3}%
\e{2}%
\eol}\vss}\rg%
%
%
\rx{\vss\hfull{%
\rlx{\hss{$2240_{x}$}}\cg%
\e{1}%
\e{4}%
\e{5}%
\e{1}%
\e{1}%
\e{0}%
\e{0}%
\e{5}%
\e{10}%
\e{6}%
\e{2}%
\e{0}%
\e{10}%
\e{4}%
\e{10}%
\e{4}%
\e{1}%
\e{0}%
\eol}\vss}\rg%
%
%
\rx{\vss\hfull{%
\rlx{\hss{$2835_{x}$}}\cg%
\e{0}%
\e{1}%
\e{6}%
\e{3}%
\e{3}%
\e{0}%
\e{0}%
\e{2}%
\e{10}%
\e{9}%
\e{5}%
\e{0}%
\e{8}%
\e{4}%
\e{12}%
\e{9}%
\e{2}%
\e{2}%
\eol}\vss}\rg%
%
%
\rx{\vss\hfull{%
\rlx{\hss{$6075_{x}$}}\cg%
\e{1}%
\e{6}%
\e{6}%
\e{7}%
\e{3}%
\e{2}%
\e{0}%
\e{8}%
\e{21}%
\e{9}%
\e{11}%
\e{1}%
\e{16}%
\e{15}%
\e{24}%
\e{19}%
\e{8}%
\e{4}%
\eol}\vss}\rg%
%
%
\rx{\vss\hfull{%
\rlx{\hss{$3200_{x}$}}\cg%
\e{1}%
\e{5}%
\e{4}%
\e{4}%
\e{3}%
\e{0}%
\e{0}%
\e{3}%
\e{9}%
\e{6}%
\e{6}%
\e{0}%
\e{7}%
\e{6}%
\e{15}%
\e{10}%
\e{4}%
\e{3}%
\eol}\vss}\rg%
%
%
\rx{\vss\hfull{%
\rlx{\hss{$70_{y}$}}\cg%
\e{0}%
\e{0}%
\e{0}%
\e{0}%
\e{0}%
\e{0}%
\e{0}%
\e{0}%
\e{0}%
\e{0}%
\e{0}%
\e{0}%
\e{0}%
\e{1}%
\e{0}%
\e{0}%
\e{1}%
\e{0}%
\eol}\vss}\rg%
%
%
\rx{\vss\hfull{%
\rlx{\hss{$1134_{y}$}}\cg%
\e{0}%
\e{1}%
\e{1}%
\e{0}%
\e{1}%
\e{1}%
\e{0}%
\e{1}%
\e{3}%
\e{2}%
\e{3}%
\e{1}%
\e{1}%
\e{3}%
\e{4}%
\e{4}%
\e{3}%
\e{1}%
\eol}\vss}\rg%
%
%
\rx{\vss\hfull{%
\rlx{\hss{$1680_{y}$}}\cg%
\e{0}%
\e{1}%
\e{0}%
\e{2}%
\e{0}%
\e{1}%
\e{0}%
\e{2}%
\e{4}%
\e{0}%
\e{4}%
\e{2}%
\e{2}%
\e{7}%
\e{5}%
\e{5}%
\e{7}%
\e{2}%
\eol}\vss}\rg%
%
%
\rx{\vss\hfull{%
\rlx{\hss{$168_{y}$}}\cg%
\e{0}%
\e{0}%
\e{1}%
\e{0}%
\e{1}%
\e{0}%
\e{0}%
\e{0}%
\e{0}%
\e{1}%
\e{0}%
\e{0}%
\e{0}%
\e{0}%
\e{1}%
\e{1}%
\e{0}%
\e{0}%
\eol}\vss}\rg%
%
%
\rx{\vss\hfull{%
\rlx{\hss{$420_{y}$}}\cg%
\e{0}%
\e{0}%
\e{1}%
\e{0}%
\e{1}%
\e{0}%
\e{0}%
\e{0}%
\e{1}%
\e{3}%
\e{1}%
\e{0}%
\e{0}%
\e{0}%
\e{2}%
\e{2}%
\e{0}%
\e{0}%
\eol}\vss}\rg%
%
%
\rx{\vss\hfull{%
\rlx{\hss{$3150_{y}$}}\cg%
\e{0}%
\e{1}%
\e{4}%
\e{2}%
\e{4}%
\e{1}%
\e{0}%
\e{1}%
\e{9}%
\e{10}%
\e{9}%
\e{1}%
\e{4}%
\e{4}%
\e{12}%
\e{12}%
\e{4}%
\e{4}%
\eol}\vss}\rg%
%
%
\rx{\vss\hfull{%
\rlx{\hss{$4200_{y}$}}\cg%
\e{0}%
\e{2}%
\e{7}%
\e{4}%
\e{7}%
\e{2}%
\e{0}%
\e{1}%
\e{10}%
\e{12}%
\e{10}%
\e{1}%
\e{6}%
\e{5}%
\e{17}%
\e{17}%
\e{5}%
\e{6}%
\eol}\vss}\rg%
%
%
\rx{\vss\hfull{%
\rlx{\hss{$2688_{y}$}}\cg%
\e{0}%
\e{1}%
\e{4}%
\e{4}%
\e{4}%
\e{1}%
\e{0}%
\e{1}%
\e{6}%
\e{5}%
\e{6}%
\e{1}%
\e{5}%
\e{5}%
\e{10}%
\e{10}%
\e{5}%
\e{5}%
\eol}\vss}\rg%
%
%
\rx{\vss\hfull{%
\rlx{\hss{$2100_{y}$}}\cg%
\e{0}%
\e{2}%
\e{0}%
\e{4}%
\e{0}%
\e{2}%
\e{0}%
\e{1}%
\e{5}%
\e{1}%
\e{5}%
\e{1}%
\e{3}%
\e{5}%
\e{8}%
\e{8}%
\e{5}%
\e{3}%
\eol}\vss}\rg%
%
%
\rx{\vss\hfull{%
\rlx{\hss{$1400_{y}$}}\cg%
\e{0}%
\e{1}%
\e{0}%
\e{2}%
\e{0}%
\e{1}%
\e{0}%
\e{1}%
\e{4}%
\e{1}%
\e{4}%
\e{1}%
\e{2}%
\e{3}%
\e{5}%
\e{5}%
\e{3}%
\e{2}%
\eol}\vss}\rg%
%
%
\rx{\vss\hfull{%
\rlx{\hss{$4536_{y}$}}\cg%
\e{0}%
\e{2}%
\e{2}%
\e{6}%
\e{2}%
\e{2}%
\e{0}%
\e{3}%
\e{13}%
\e{6}%
\e{13}%
\e{3}%
\e{7}%
\e{9}%
\e{16}%
\e{16}%
\e{9}%
\e{7}%
\eol}\vss}\rg%
%
%
\rx{\vss\hfull{%
\rlx{\hss{$5670_{y}$}}\cg%
\e{0}%
\e{3}%
\e{3}%
\e{6}%
\e{3}%
\e{3}%
\e{0}%
\e{4}%
\e{16}%
\e{8}%
\e{16}%
\e{4}%
\e{8}%
\e{12}%
\e{20}%
\e{20}%
\e{12}%
\e{8}%
\eol}\vss}\rg%
%
%
\rx{\vss\hfull{%
\rlx{\hss{$4480_{y}$}}\cg%
\e{0}%
\e{2}%
\e{4}%
\e{4}%
\e{4}%
\e{2}%
\e{0}%
\e{2}%
\e{13}%
\e{11}%
\e{13}%
\e{2}%
\e{6}%
\e{6}%
\e{17}%
\e{17}%
\e{6}%
\e{6}%
\eol}\vss}\rg%
%
%
\rx{\vss\hfull{%
\rlx{\hss{$8_{z}$}}\cg%
\e{1}%
\e{0}%
\e{0}%
\e{0}%
\e{0}%
\e{0}%
\e{0}%
\e{0}%
\e{0}%
\e{0}%
\e{0}%
\e{0}%
\e{0}%
\e{0}%
\e{0}%
\e{0}%
\e{0}%
\e{0}%
\eol}\vss}\rg%
%
%
\rx{\vss\hfull{%
\rlx{\hss{$56_{z}$}}\cg%
\e{0}%
\e{0}%
\e{0}%
\e{0}%
\e{0}%
\e{0}%
\e{0}%
\e{1}%
\e{0}%
\e{0}%
\e{0}%
\e{0}%
\e{0}%
\e{1}%
\e{0}%
\e{0}%
\e{0}%
\e{0}%
\eol}\vss}\rg%
%
%
\rx{\vss\hfull{%
\rlx{\hss{$160_{z}$}}\cg%
\e{3}%
\e{1}%
\e{0}%
\e{0}%
\e{0}%
\e{0}%
\e{0}%
\e{3}%
\e{0}%
\e{0}%
\e{0}%
\e{0}%
\e{1}%
\e{1}%
\e{0}%
\e{0}%
\e{0}%
\e{0}%
\eol}\vss}\rg%
%
%
\rx{\vss\hfull{%
\rlx{\hss{$112_{z}$}}\cg%
\e{3}%
\e{1}%
\e{0}%
\e{0}%
\e{0}%
\e{0}%
\e{0}%
\e{2}%
\e{0}%
\e{0}%
\e{0}%
\e{0}%
\e{1}%
\e{0}%
\e{0}%
\e{0}%
\e{0}%
\e{0}%
\eol}\vss}\rg%
%
%
\rx{\vss\hfull{%
\rlx{\hss{$840_{z}$}}\cg%
\e{1}%
\e{3}%
\e{1}%
\e{0}%
\e{0}%
\e{0}%
\e{0}%
\e{3}%
\e{3}%
\e{1}%
\e{1}%
\e{0}%
\e{3}%
\e{3}%
\e{4}%
\e{1}%
\e{1}%
\e{0}%
\eol}\vss}\rg%
%
%
\rx{\vss\hfull{%
\rlx{\hss{$1296_{z}$}}\cg%
\e{2}%
\e{3}%
\e{0}%
\e{1}%
\e{0}%
\e{0}%
\e{0}%
\e{7}%
\e{5}%
\e{0}%
\e{1}%
\e{0}%
\e{5}%
\e{7}%
\e{4}%
\e{2}%
\e{2}%
\e{0}%
\eol}\vss}\rg%
%
%
\rx{\vss\hfull{%
\rlx{\hss{$1400_{z}$}}\cg%
\e{5}%
\e{6}%
\e{1}%
\e{1}%
\e{0}%
\e{0}%
\e{0}%
\e{9}%
\e{6}%
\e{1}%
\e{0}%
\e{0}%
\e{9}%
\e{4}%
\e{5}%
\e{1}%
\e{0}%
\e{0}%
\eol}\vss}\rg%
%
%
\rx{\vss\hfull{%
\rlx{\hss{$1008_{z}$}}\cg%
\e{4}%
\e{4}%
\e{0}%
\e{1}%
\e{0}%
\e{0}%
\e{0}%
\e{8}%
\e{4}%
\e{0}%
\e{0}%
\e{0}%
\e{6}%
\e{4}%
\e{3}%
\e{1}%
\e{0}%
\e{0}%
\eol}\vss}\rg%
%
%
\rx{\vss\hfull{%
\rlx{\hss{$560_{z}$}}\cg%
\e{5}%
\e{3}%
\e{1}%
\e{0}%
\e{0}%
\e{0}%
\e{0}%
\e{6}%
\e{2}%
\e{0}%
\e{0}%
\e{0}%
\e{5}%
\e{1}%
\e{1}%
\e{0}%
\e{0}%
\e{0}%
\eol}\vss}\rg%
%
%
\rx{\vss\hfull{%
\rlx{\hss{$1400_{zz}$}}\cg%
\e{0}%
\e{1}%
\e{5}%
\e{1}%
\e{1}%
\e{0}%
\e{0}%
\e{2}%
\e{6}%
\e{5}%
\e{1}%
\e{0}%
\e{7}%
\e{2}%
\e{6}%
\e{3}%
\e{0}%
\e{0}%
\eol}\vss}\rg%
\eop
\eject
\tablecont%
%
%
%
%
%
%
\rowpts=18 true pt%
\colpts=18 true pt%
\rowlabpts=40 true pt%
\collabpts=85 true pt%
\clx{\vss\hfull{%
\rlx{\hss{$ $}}\cg%
\cx{\hskip 16 true pt\flip{$[{5}:-]{\times}[{2}{1}]$}\hss}\cg%
\cx{\hskip 16 true pt\flip{$[{4}{1}:-]{\times}[{2}{1}]$}\hss}\cg%
\cx{\hskip 16 true pt\flip{$[{3}{2}:-]{\times}[{2}{1}]$}\hss}\cg%
\cx{\hskip 16 true pt\flip{$[{3}{1^{2}}:-]{\times}[{2}{1}]$}\hss}\cg%
\cx{\hskip 16 true pt\flip{$[{2^{2}}{1}:-]{\times}[{2}{1}]$}\hss}\cg%
\cx{\hskip 16 true pt\flip{$[{2}{1^{3}}:-]{\times}[{2}{1}]$}\hss}\cg%
\cx{\hskip 16 true pt\flip{$[{1^{5}}:-]{\times}[{2}{1}]$}\hss}\cg%
\cx{\hskip 16 true pt\flip{$[{4}:{1}]{\times}[{2}{1}]$}\hss}\cg%
\cx{\hskip 16 true pt\flip{$[{3}{1}:{1}]{\times}[{2}{1}]$}\hss}\cg%
\cx{\hskip 16 true pt\flip{$[{2^{2}}:{1}]{\times}[{2}{1}]$}\hss}\cg%
\cx{\hskip 16 true pt\flip{$[{2}{1^{2}}:{1}]{\times}[{2}{1}]$}\hss}\cg%
\cx{\hskip 16 true pt\flip{$[{1^{4}}:{1}]{\times}[{2}{1}]$}\hss}\cg%
\cx{\hskip 16 true pt\flip{$[{3}:{2}]{\times}[{2}{1}]$}\hss}\cg%
\cx{\hskip 16 true pt\flip{$[{3}:{1^{2}}]{\times}[{2}{1}]$}\hss}\cg%
\cx{\hskip 16 true pt\flip{$[{2}{1}:{2}]{\times}[{2}{1}]$}\hss}\cg%
\cx{\hskip 16 true pt\flip{$[{2}{1}:{1^{2}}]{\times}[{2}{1}]$}\hss}\cg%
\cx{\hskip 16 true pt\flip{$[{1^{3}}:{2}]{\times}[{2}{1}]$}\hss}\cg%
\cx{\hskip 16 true pt\flip{$[{1^{3}}:{1^{2}}]{\times}[{2}{1}]$}\hss}\cg%
\eol}}\rg%
%
%
\rx{\vss\hfull{%
\rlx{\hss{$4200_{z}$}}\cg%
\e{0}%
\e{2}%
\e{7}%
\e{4}%
\e{4}%
\e{2}%
\e{0}%
\e{3}%
\e{14}%
\e{10}%
\e{8}%
\e{1}%
\e{10}%
\e{7}%
\e{17}%
\e{15}%
\e{4}%
\e{3}%
\eol}\vss}\rg%
%
%
\rx{\vss\hfull{%
\rlx{\hss{$400_{z}$}}\cg%
\e{1}%
\e{1}%
\e{2}%
\e{0}%
\e{0}%
\e{0}%
\e{0}%
\e{2}%
\e{2}%
\e{1}%
\e{0}%
\e{0}%
\e{4}%
\e{0}%
\e{1}%
\e{0}%
\e{0}%
\e{0}%
\eol}\vss}\rg%
%
%
\rx{\vss\hfull{%
\rlx{\hss{$3240_{z}$}}\cg%
\e{4}%
\e{8}%
\e{6}%
\e{2}%
\e{1}%
\e{0}%
\e{0}%
\e{12}%
\e{14}%
\e{5}%
\e{2}%
\e{0}%
\e{17}%
\e{8}%
\e{13}%
\e{5}%
\e{1}%
\e{0}%
\eol}\vss}\rg%
%
%
\rx{\vss\hfull{%
\rlx{\hss{$4536_{z}$}}\cg%
\e{1}%
\e{6}%
\e{8}%
\e{4}%
\e{4}%
\e{0}%
\e{0}%
\e{5}%
\e{16}%
\e{12}%
\e{7}%
\e{0}%
\e{13}%
\e{8}%
\e{21}%
\e{13}%
\e{3}%
\e{2}%
\eol}\vss}\rg%
%
%
\rx{\vss\hfull{%
\rlx{\hss{$2400_{z}$}}\cg%
\e{0}%
\e{2}%
\e{0}%
\e{3}%
\e{0}%
\e{1}%
\e{0}%
\e{4}%
\e{8}%
\e{1}%
\e{5}%
\e{1}%
\e{5}%
\e{9}%
\e{8}%
\e{7}%
\e{6}%
\e{2}%
\eol}\vss}\rg%
%
%
\rx{\vss\hfull{%
\rlx{\hss{$3360_{z}$}}\cg%
\e{1}%
\e{3}%
\e{5}%
\e{4}%
\e{2}%
\e{0}%
\e{0}%
\e{6}%
\e{13}%
\e{6}%
\e{5}%
\e{0}%
\e{12}%
\e{8}%
\e{13}%
\e{9}%
\e{3}%
\e{2}%
\eol}\vss}\rg%
%
%
\rx{\vss\hfull{%
\rlx{\hss{$2800_{z}$}}\cg%
\e{1}%
\e{4}%
\e{1}%
\e{4}%
\e{0}%
\e{1}%
\e{0}%
\e{6}%
\e{11}%
\e{2}%
\e{4}%
\e{0}%
\e{9}%
\e{8}%
\e{11}%
\e{8}%
\e{3}%
\e{1}%
\eol}\vss}\rg%
%
%
\rx{\vss\hfull{%
\rlx{\hss{$4096_{z}$}}\cg%
\e{3}%
\e{8}%
\e{4}%
\e{4}%
\e{1}%
\e{0}%
\e{0}%
\e{11}%
\e{16}%
\e{5}%
\e{5}%
\e{0}%
\e{15}%
\e{12}%
\e{17}%
\e{9}%
\e{4}%
\e{1}%
\eol}\vss}\rg%
%
%
\rx{\vss\hfull{%
\rlx{\hss{$5600_{z}$}}\cg%
\e{1}%
\e{7}%
\e{3}%
\e{7}%
\e{2}%
\e{1}%
\e{0}%
\e{7}%
\e{18}%
\e{7}%
\e{12}%
\e{1}%
\e{12}%
\e{14}%
\e{23}%
\e{17}%
\e{9}%
\e{5}%
\eol}\vss}\rg%
%
%
\rx{\vss\hfull{%
\rlx{\hss{$448_{z}$}}\cg%
\e{1}%
\e{2}%
\e{1}%
\e{0}%
\e{0}%
\e{0}%
\e{0}%
\e{2}%
\e{2}%
\e{1}%
\e{0}%
\e{0}%
\e{3}%
\e{1}%
\e{2}%
\e{0}%
\e{0}%
\e{0}%
\eol}\vss}\rg%
%
%
\rx{\vss\hfull{%
\rlx{\hss{$448_{w}$}}\cg%
\e{0}%
\e{0}%
\e{0}%
\e{0}%
\e{0}%
\e{0}%
\e{0}%
\e{1}%
\e{1}%
\e{0}%
\e{1}%
\e{1}%
\e{0}%
\e{3}%
\e{1}%
\e{1}%
\e{3}%
\e{0}%
\eol}\vss}\rg%
%
%
\rx{\vss\hfull{%
\rlx{\hss{$1344_{w}$}}\cg%
\e{0}%
\e{1}%
\e{3}%
\e{0}%
\e{3}%
\e{1}%
\e{0}%
\e{0}%
\e{3}%
\e{6}%
\e{3}%
\e{0}%
\e{1}%
\e{1}%
\e{6}%
\e{6}%
\e{1}%
\e{1}%
\eol}\vss}\rg%
%
%
\rx{\vss\hfull{%
\rlx{\hss{$5600_{w}$}}\cg%
\e{0}%
\e{3}%
\e{4}%
\e{6}%
\e{4}%
\e{3}%
\e{0}%
\e{4}%
\e{15}%
\e{8}%
\e{15}%
\e{4}%
\e{8}%
\e{12}%
\e{20}%
\e{20}%
\e{12}%
\e{8}%
\eol}\vss}\rg%
%
%
\rx{\vss\hfull{%
\rlx{\hss{$2016_{w}$}}\cg%
\e{0}%
\e{0}%
\e{4}%
\e{2}%
\e{4}%
\e{0}%
\e{0}%
\e{0}%
\e{5}%
\e{8}%
\e{5}%
\e{0}%
\e{3}%
\e{2}%
\e{8}%
\e{8}%
\e{2}%
\e{3}%
\eol}\vss}\rg%
%
%
\rx{\vss\hfull{%
\rlx{\hss{$7168_{w}$}}\cg%
\e{0}%
\e{3}%
\e{8}%
\e{8}%
\e{8}%
\e{3}%
\e{0}%
\e{3}%
\e{19}%
\e{16}%
\e{19}%
\e{3}%
\e{11}%
\e{11}%
\e{27}%
\e{27}%
\e{11}%
\e{11}%
\eol}\vss}\rg%
%
%
%
%
%
%
\rowpts=18 true pt%
\colpts=18 true pt%
\rowlabpts=40 true pt%
\collabpts=85 true pt%
\clx{\vss\hfull{%
\rlx{\hss{$ $}}\cg%
\cx{\hskip 16 true pt\flip{$[{5}:-]{\times}[{1^{3}}]$}\hss}\cg%
\cx{\hskip 16 true pt\flip{$[{4}{1}:-]{\times}[{1^{3}}]$}\hss}\cg%
\cx{\hskip 16 true pt\flip{$[{3}{2}:-]{\times}[{1^{3}}]$}\hss}\cg%
\cx{\hskip 16 true pt\flip{$[{3}{1^{2}}:-]{\times}[{1^{3}}]$}\hss}\cg%
\cx{\hskip 16 true pt\flip{$[{2^{2}}{1}:-]{\times}[{1^{3}}]$}\hss}\cg%
\cx{\hskip 16 true pt\flip{$[{2}{1^{3}}:-]{\times}[{1^{3}}]$}\hss}\cg%
\cx{\hskip 16 true pt\flip{$[{1^{5}}:-]{\times}[{1^{3}}]$}\hss}\cg%
\cx{\hskip 16 true pt\flip{$[{4}:{1}]{\times}[{1^{3}}]$}\hss}\cg%
\cx{\hskip 16 true pt\flip{$[{3}{1}:{1}]{\times}[{1^{3}}]$}\hss}\cg%
\cx{\hskip 16 true pt\flip{$[{2^{2}}:{1}]{\times}[{1^{3}}]$}\hss}\cg%
\cx{\hskip 16 true pt\flip{$[{2}{1^{2}}:{1}]{\times}[{1^{3}}]$}\hss}\cg%
\cx{\hskip 16 true pt\flip{$[{1^{4}}:{1}]{\times}[{1^{3}}]$}\hss}\cg%
\cx{\hskip 16 true pt\flip{$[{3}:{2}]{\times}[{1^{3}}]$}\hss}\cg%
\cx{\hskip 16 true pt\flip{$[{3}:{1^{2}}]{\times}[{1^{3}}]$}\hss}\cg%
\cx{\hskip 16 true pt\flip{$[{2}{1}:{2}]{\times}[{1^{3}}]$}\hss}\cg%
\cx{\hskip 16 true pt\flip{$[{2}{1}:{1^{2}}]{\times}[{1^{3}}]$}\hss}\cg%
\cx{\hskip 16 true pt\flip{$[{1^{3}}:{2}]{\times}[{1^{3}}]$}\hss}\cg%
\cx{\hskip 16 true pt\flip{$[{1^{3}}:{1^{2}}]{\times}[{1^{3}}]$}\hss}\cg%
\eol}}\rg%
%
%
\rx{\vss\hfull{%
\rlx{\hss{$1_{x}$}}\cg%
\e{0}%
\e{0}%
\e{0}%
\e{0}%
\e{0}%
\e{0}%
\e{0}%
\e{0}%
\e{0}%
\e{0}%
\e{0}%
\e{0}%
\e{0}%
\e{0}%
\e{0}%
\e{0}%
\e{0}%
\e{0}%
\eol}\vss}\rg%
%
%
\rx{\vss\hfull{%
\rlx{\hss{$28_{x}$}}\cg%
\e{1}%
\e{0}%
\e{0}%
\e{0}%
\e{0}%
\e{0}%
\e{0}%
\e{0}%
\e{0}%
\e{0}%
\e{0}%
\e{0}%
\e{0}%
\e{0}%
\e{0}%
\e{0}%
\e{0}%
\e{0}%
\eol}\vss}\rg%
%
%
\rx{\vss\hfull{%
\rlx{\hss{$35_{x}$}}\cg%
\e{0}%
\e{0}%
\e{0}%
\e{0}%
\e{0}%
\e{0}%
\e{0}%
\e{0}%
\e{0}%
\e{0}%
\e{0}%
\e{0}%
\e{0}%
\e{0}%
\e{0}%
\e{0}%
\e{0}%
\e{0}%
\eol}\vss}\rg%
%
%
\rx{\vss\hfull{%
\rlx{\hss{$84_{x}$}}\cg%
\e{0}%
\e{0}%
\e{0}%
\e{0}%
\e{0}%
\e{0}%
\e{0}%
\e{0}%
\e{0}%
\e{0}%
\e{0}%
\e{0}%
\e{0}%
\e{0}%
\e{0}%
\e{0}%
\e{0}%
\e{0}%
\eol}\vss}\rg%
%
%
\rx{\vss\hfull{%
\rlx{\hss{$50_{x}$}}\cg%
\e{0}%
\e{0}%
\e{0}%
\e{0}%
\e{0}%
\e{0}%
\e{0}%
\e{0}%
\e{0}%
\e{0}%
\e{0}%
\e{0}%
\e{0}%
\e{0}%
\e{0}%
\e{0}%
\e{0}%
\e{0}%
\eol}\vss}\rg%
%
%
\rx{\vss\hfull{%
\rlx{\hss{$350_{x}$}}\cg%
\e{1}%
\e{1}%
\e{0}%
\e{0}%
\e{0}%
\e{0}%
\e{0}%
\e{2}%
\e{0}%
\e{0}%
\e{0}%
\e{0}%
\e{1}%
\e{1}%
\e{0}%
\e{0}%
\e{0}%
\e{0}%
\eol}\vss}\rg%
%
%
\rx{\vss\hfull{%
\rlx{\hss{$300_{x}$}}\cg%
\e{0}%
\e{0}%
\e{0}%
\e{0}%
\e{0}%
\e{0}%
\e{0}%
\e{1}%
\e{0}%
\e{0}%
\e{0}%
\e{0}%
\e{0}%
\e{1}%
\e{0}%
\e{0}%
\e{0}%
\e{0}%
\eol}\vss}\rg%
\eop
\eject
\tablecont%
%
%
%
%
%
%
\rowpts=18 true pt%
\colpts=18 true pt%
\rowlabpts=40 true pt%
\collabpts=85 true pt%
\clx{\vss\hfull{%
\rlx{\hss{$ $}}\cg%
\cx{\hskip 16 true pt\flip{$[{5}:-]{\times}[{1^{3}}]$}\hss}\cg%
\cx{\hskip 16 true pt\flip{$[{4}{1}:-]{\times}[{1^{3}}]$}\hss}\cg%
\cx{\hskip 16 true pt\flip{$[{3}{2}:-]{\times}[{1^{3}}]$}\hss}\cg%
\cx{\hskip 16 true pt\flip{$[{3}{1^{2}}:-]{\times}[{1^{3}}]$}\hss}\cg%
\cx{\hskip 16 true pt\flip{$[{2^{2}}{1}:-]{\times}[{1^{3}}]$}\hss}\cg%
\cx{\hskip 16 true pt\flip{$[{2}{1^{3}}:-]{\times}[{1^{3}}]$}\hss}\cg%
\cx{\hskip 16 true pt\flip{$[{1^{5}}:-]{\times}[{1^{3}}]$}\hss}\cg%
\cx{\hskip 16 true pt\flip{$[{4}:{1}]{\times}[{1^{3}}]$}\hss}\cg%
\cx{\hskip 16 true pt\flip{$[{3}{1}:{1}]{\times}[{1^{3}}]$}\hss}\cg%
\cx{\hskip 16 true pt\flip{$[{2^{2}}:{1}]{\times}[{1^{3}}]$}\hss}\cg%
\cx{\hskip 16 true pt\flip{$[{2}{1^{2}}:{1}]{\times}[{1^{3}}]$}\hss}\cg%
\cx{\hskip 16 true pt\flip{$[{1^{4}}:{1}]{\times}[{1^{3}}]$}\hss}\cg%
\cx{\hskip 16 true pt\flip{$[{3}:{2}]{\times}[{1^{3}}]$}\hss}\cg%
\cx{\hskip 16 true pt\flip{$[{3}:{1^{2}}]{\times}[{1^{3}}]$}\hss}\cg%
\cx{\hskip 16 true pt\flip{$[{2}{1}:{2}]{\times}[{1^{3}}]$}\hss}\cg%
\cx{\hskip 16 true pt\flip{$[{2}{1}:{1^{2}}]{\times}[{1^{3}}]$}\hss}\cg%
\cx{\hskip 16 true pt\flip{$[{1^{3}}:{2}]{\times}[{1^{3}}]$}\hss}\cg%
\cx{\hskip 16 true pt\flip{$[{1^{3}}:{1^{2}}]{\times}[{1^{3}}]$}\hss}\cg%
\eol}}\rg%
%
%
\rx{\vss\hfull{%
\rlx{\hss{$567_{x}$}}\cg%
\e{3}%
\e{1}%
\e{0}%
\e{0}%
\e{0}%
\e{0}%
\e{0}%
\e{2}%
\e{0}%
\e{0}%
\e{0}%
\e{0}%
\e{1}%
\e{0}%
\e{0}%
\e{0}%
\e{0}%
\e{0}%
\eol}\vss}\rg%
%
%
\rx{\vss\hfull{%
\rlx{\hss{$210_{x}$}}\cg%
\e{1}%
\e{0}%
\e{0}%
\e{0}%
\e{0}%
\e{0}%
\e{0}%
\e{1}%
\e{0}%
\e{0}%
\e{0}%
\e{0}%
\e{0}%
\e{0}%
\e{0}%
\e{0}%
\e{0}%
\e{0}%
\eol}\vss}\rg%
%
%
\rx{\vss\hfull{%
\rlx{\hss{$840_{x}$}}\cg%
\e{0}%
\e{0}%
\e{0}%
\e{0}%
\e{0}%
\e{0}%
\e{0}%
\e{0}%
\e{1}%
\e{0}%
\e{1}%
\e{0}%
\e{0}%
\e{1}%
\e{1}%
\e{1}%
\e{1}%
\e{0}%
\eol}\vss}\rg%
%
%
\rx{\vss\hfull{%
\rlx{\hss{$700_{x}$}}\cg%
\e{0}%
\e{0}%
\e{0}%
\e{0}%
\e{0}%
\e{0}%
\e{0}%
\e{1}%
\e{1}%
\e{0}%
\e{0}%
\e{0}%
\e{1}%
\e{0}%
\e{0}%
\e{0}%
\e{0}%
\e{0}%
\eol}\vss}\rg%
%
%
\rx{\vss\hfull{%
\rlx{\hss{$175_{x}$}}\cg%
\e{0}%
\e{0}%
\e{0}%
\e{0}%
\e{0}%
\e{0}%
\e{0}%
\e{0}%
\e{0}%
\e{1}%
\e{0}%
\e{0}%
\e{0}%
\e{0}%
\e{0}%
\e{0}%
\e{0}%
\e{0}%
\eol}\vss}\rg%
%
%
\rx{\vss\hfull{%
\rlx{\hss{$1400_{x}$}}\cg%
\e{1}%
\e{1}%
\e{1}%
\e{0}%
\e{0}%
\e{0}%
\e{0}%
\e{2}%
\e{2}%
\e{1}%
\e{0}%
\e{0}%
\e{3}%
\e{0}%
\e{1}%
\e{0}%
\e{0}%
\e{0}%
\eol}\vss}\rg%
%
%
\rx{\vss\hfull{%
\rlx{\hss{$1050_{x}$}}\cg%
\e{0}%
\e{0}%
\e{2}%
\e{0}%
\e{0}%
\e{0}%
\e{0}%
\e{0}%
\e{1}%
\e{1}%
\e{0}%
\e{0}%
\e{2}%
\e{0}%
\e{1}%
\e{0}%
\e{0}%
\e{0}%
\eol}\vss}\rg%
%
%
\rx{\vss\hfull{%
\rlx{\hss{$1575_{x}$}}\cg%
\e{2}%
\e{2}%
\e{1}%
\e{0}%
\e{0}%
\e{0}%
\e{0}%
\e{4}%
\e{2}%
\e{0}%
\e{0}%
\e{0}%
\e{4}%
\e{1}%
\e{1}%
\e{0}%
\e{0}%
\e{0}%
\eol}\vss}\rg%
%
%
\rx{\vss\hfull{%
\rlx{\hss{$1344_{x}$}}\cg%
\e{1}%
\e{2}%
\e{0}%
\e{0}%
\e{0}%
\e{0}%
\e{0}%
\e{2}%
\e{1}%
\e{0}%
\e{0}%
\e{0}%
\e{2}%
\e{1}%
\e{1}%
\e{0}%
\e{0}%
\e{0}%
\eol}\vss}\rg%
%
%
\rx{\vss\hfull{%
\rlx{\hss{$2100_{x}$}}\cg%
\e{2}%
\e{3}%
\e{0}%
\e{1}%
\e{0}%
\e{0}%
\e{0}%
\e{5}%
\e{3}%
\e{0}%
\e{0}%
\e{0}%
\e{4}%
\e{5}%
\e{3}%
\e{1}%
\e{1}%
\e{0}%
\eol}\vss}\rg%
%
%
\rx{\vss\hfull{%
\rlx{\hss{$2268_{x}$}}\cg%
\e{2}%
\e{2}%
\e{0}%
\e{1}%
\e{0}%
\e{0}%
\e{0}%
\e{5}%
\e{3}%
\e{0}%
\e{0}%
\e{0}%
\e{4}%
\e{3}%
\e{2}%
\e{1}%
\e{0}%
\e{0}%
\eol}\vss}\rg%
%
%
\rx{\vss\hfull{%
\rlx{\hss{$525_{x}$}}\cg%
\e{2}%
\e{2}%
\e{0}%
\e{0}%
\e{0}%
\e{0}%
\e{0}%
\e{1}%
\e{0}%
\e{0}%
\e{0}%
\e{0}%
\e{1}%
\e{0}%
\e{1}%
\e{0}%
\e{0}%
\e{0}%
\eol}\vss}\rg%
%
%
\rx{\vss\hfull{%
\rlx{\hss{$700_{xx}$}}\cg%
\e{0}%
\e{0}%
\e{2}%
\e{0}%
\e{1}%
\e{0}%
\e{0}%
\e{0}%
\e{0}%
\e{1}%
\e{0}%
\e{0}%
\e{1}%
\e{0}%
\e{1}%
\e{1}%
\e{0}%
\e{0}%
\eol}\vss}\rg%
%
%
\rx{\vss\hfull{%
\rlx{\hss{$972_{x}$}}\cg%
\e{0}%
\e{0}%
\e{0}%
\e{1}%
\e{0}%
\e{0}%
\e{0}%
\e{0}%
\e{1}%
\e{0}%
\e{0}%
\e{0}%
\e{1}%
\e{1}%
\e{1}%
\e{1}%
\e{0}%
\e{0}%
\eol}\vss}\rg%
%
%
\rx{\vss\hfull{%
\rlx{\hss{$4096_{x}$}}\cg%
\e{1}%
\e{4}%
\e{1}%
\e{1}%
\e{0}%
\e{0}%
\e{0}%
\e{5}%
\e{6}%
\e{1}%
\e{1}%
\e{0}%
\e{6}%
\e{5}%
\e{6}%
\e{2}%
\e{1}%
\e{0}%
\eol}\vss}\rg%
%
%
\rx{\vss\hfull{%
\rlx{\hss{$4200_{x}$}}\cg%
\e{0}%
\e{1}%
\e{3}%
\e{1}%
\e{1}%
\e{0}%
\e{0}%
\e{2}%
\e{6}%
\e{3}%
\e{2}%
\e{0}%
\e{5}%
\e{3}%
\e{6}%
\e{4}%
\e{1}%
\e{0}%
\eol}\vss}\rg%
%
%
\rx{\vss\hfull{%
\rlx{\hss{$2240_{x}$}}\cg%
\e{0}%
\e{1}%
\e{1}%
\e{0}%
\e{0}%
\e{0}%
\e{0}%
\e{1}%
\e{3}%
\e{2}%
\e{1}%
\e{0}%
\e{2}%
\e{1}%
\e{3}%
\e{1}%
\e{0}%
\e{0}%
\eol}\vss}\rg%
%
%
\rx{\vss\hfull{%
\rlx{\hss{$2835_{x}$}}\cg%
\e{0}%
\e{0}%
\e{2}%
\e{1}%
\e{2}%
\e{0}%
\e{0}%
\e{0}%
\e{3}%
\e{5}%
\e{2}%
\e{0}%
\e{2}%
\e{1}%
\e{4}%
\e{4}%
\e{0}%
\e{1}%
\eol}\vss}\rg%
%
%
\rx{\vss\hfull{%
\rlx{\hss{$6075_{x}$}}\cg%
\e{1}%
\e{4}%
\e{4}%
\e{3}%
\e{1}%
\e{0}%
\e{0}%
\e{4}%
\e{9}%
\e{4}%
\e{3}%
\e{0}%
\e{9}%
\e{5}%
\e{11}%
\e{6}%
\e{2}%
\e{1}%
\eol}\vss}\rg%
%
%
\rx{\vss\hfull{%
\rlx{\hss{$3200_{x}$}}\cg%
\e{0}%
\e{2}%
\e{0}%
\e{2}%
\e{0}%
\e{0}%
\e{0}%
\e{1}%
\e{4}%
\e{1}%
\e{3}%
\e{0}%
\e{2}%
\e{4}%
\e{6}%
\e{4}%
\e{3}%
\e{1}%
\eol}\vss}\rg%
%
%
\rx{\vss\hfull{%
\rlx{\hss{$70_{y}$}}\cg%
\e{0}%
\e{0}%
\e{0}%
\e{0}%
\e{0}%
\e{0}%
\e{0}%
\e{1}%
\e{0}%
\e{0}%
\e{0}%
\e{0}%
\e{0}%
\e{1}%
\e{0}%
\e{0}%
\e{0}%
\e{0}%
\eol}\vss}\rg%
%
%
\rx{\vss\hfull{%
\rlx{\hss{$1134_{y}$}}\cg%
\e{0}%
\e{1}%
\e{1}%
\e{0}%
\e{0}%
\e{0}%
\e{0}%
\e{0}%
\e{2}%
\e{1}%
\e{1}%
\e{1}%
\e{1}%
\e{1}%
\e{3}%
\e{1}%
\e{2}%
\e{0}%
\eol}\vss}\rg%
%
%
\rx{\vss\hfull{%
\rlx{\hss{$1680_{y}$}}\cg%
\e{0}%
\e{1}%
\e{0}%
\e{1}%
\e{0}%
\e{0}%
\e{0}%
\e{3}%
\e{4}%
\e{0}%
\e{1}%
\e{0}%
\e{3}%
\e{5}%
\e{3}%
\e{2}%
\e{2}%
\e{0}%
\eol}\vss}\rg%
%
%
\rx{\vss\hfull{%
\rlx{\hss{$168_{y}$}}\cg%
\e{0}%
\e{0}%
\e{0}%
\e{0}%
\e{0}%
\e{1}%
\e{0}%
\e{0}%
\e{0}%
\e{0}%
\e{0}%
\e{0}%
\e{0}%
\e{0}%
\e{0}%
\e{1}%
\e{0}%
\e{0}%
\eol}\vss}\rg%
%
%
\rx{\vss\hfull{%
\rlx{\hss{$420_{y}$}}\cg%
\e{0}%
\e{0}%
\e{0}%
\e{0}%
\e{1}%
\e{0}%
\e{0}%
\e{0}%
\e{0}%
\e{1}%
\e{1}%
\e{0}%
\e{0}%
\e{0}%
\e{0}%
\e{1}%
\e{0}%
\e{1}%
\eol}\vss}\rg%
%
%
\rx{\vss\hfull{%
\rlx{\hss{$3150_{y}$}}\cg%
\e{0}%
\e{1}%
\e{2}%
\e{1}%
\e{2}%
\e{0}%
\e{0}%
\e{0}%
\e{3}%
\e{5}%
\e{5}%
\e{1}%
\e{1}%
\e{1}%
\e{6}%
\e{5}%
\e{2}%
\e{3}%
\eol}\vss}\rg%
%
%
\rx{\vss\hfull{%
\rlx{\hss{$4200_{y}$}}\cg%
\e{0}%
\e{0}%
\e{3}%
\e{2}%
\e{3}%
\e{2}%
\e{0}%
\e{0}%
\e{4}%
\e{5}%
\e{5}%
\e{1}%
\e{2}%
\e{2}%
\e{7}%
\e{9}%
\e{3}%
\e{3}%
\eol}\vss}\rg%
%
%
\rx{\vss\hfull{%
\rlx{\hss{$2688_{y}$}}\cg%
\e{0}%
\e{0}%
\e{2}%
\e{2}%
\e{2}%
\e{0}%
\e{0}%
\e{0}%
\e{3}%
\e{3}%
\e{3}%
\e{0}%
\e{2}%
\e{2}%
\e{5}%
\e{5}%
\e{2}%
\e{2}%
\eol}\vss}\rg%
%
%
\rx{\vss\hfull{%
\rlx{\hss{$2100_{y}$}}\cg%
\e{1}%
\e{3}%
\e{0}%
\e{2}%
\e{0}%
\e{0}%
\e{0}%
\e{2}%
\e{3}%
\e{0}%
\e{2}%
\e{0}%
\e{3}%
\e{3}%
\e{5}%
\e{3}%
\e{2}%
\e{1}%
\eol}\vss}\rg%
%
%
\rx{\vss\hfull{%
\rlx{\hss{$1400_{y}$}}\cg%
\e{0}%
\e{0}%
\e{0}%
\e{1}%
\e{0}%
\e{1}%
\e{0}%
\e{2}%
\e{3}%
\e{0}%
\e{1}%
\e{0}%
\e{2}%
\e{3}%
\e{2}%
\e{3}%
\e{1}%
\e{0}%
\eol}\vss}\rg%
%
%
\rx{\vss\hfull{%
\rlx{\hss{$4536_{y}$}}\cg%
\e{0}%
\e{1}%
\e{1}%
\e{3}%
\e{1}%
\e{1}%
\e{0}%
\e{3}%
\e{8}%
\e{3}%
\e{5}%
\e{0}%
\e{5}%
\e{6}%
\e{8}%
\e{8}%
\e{3}%
\e{2}%
\eol}\vss}\rg%
%
%
\rx{\vss\hfull{%
\rlx{\hss{$5670_{y}$}}\cg%
\e{0}%
\e{2}%
\e{2}%
\e{3}%
\e{1}%
\e{1}%
\e{0}%
\e{3}%
\e{10}%
\e{4}%
\e{6}%
\e{1}%
\e{6}%
\e{7}%
\e{11}%
\e{9}%
\e{5}%
\e{2}%
\eol}\vss}\rg%
%
%
\rx{\vss\hfull{%
\rlx{\hss{$4480_{y}$}}\cg%
\e{0}%
\e{1}%
\e{2}%
\e{2}%
\e{2}%
\e{1}%
\e{0}%
\e{1}%
\e{6}%
\e{5}%
\e{6}%
\e{1}%
\e{3}%
\e{3}%
\e{8}%
\e{8}%
\e{3}%
\e{3}%
\eol}\vss}\rg%
\eop
\eject
\tablecont%
%
%
%
%
%
%
\rowpts=18 true pt%
\colpts=18 true pt%
\rowlabpts=40 true pt%
\collabpts=85 true pt%
\clx{\vss\hfull{%
\rlx{\hss{$ $}}\cg%
\cx{\hskip 16 true pt\flip{$[{5}:-]{\times}[{1^{3}}]$}\hss}\cg%
\cx{\hskip 16 true pt\flip{$[{4}{1}:-]{\times}[{1^{3}}]$}\hss}\cg%
\cx{\hskip 16 true pt\flip{$[{3}{2}:-]{\times}[{1^{3}}]$}\hss}\cg%
\cx{\hskip 16 true pt\flip{$[{3}{1^{2}}:-]{\times}[{1^{3}}]$}\hss}\cg%
\cx{\hskip 16 true pt\flip{$[{2^{2}}{1}:-]{\times}[{1^{3}}]$}\hss}\cg%
\cx{\hskip 16 true pt\flip{$[{2}{1^{3}}:-]{\times}[{1^{3}}]$}\hss}\cg%
\cx{\hskip 16 true pt\flip{$[{1^{5}}:-]{\times}[{1^{3}}]$}\hss}\cg%
\cx{\hskip 16 true pt\flip{$[{4}:{1}]{\times}[{1^{3}}]$}\hss}\cg%
\cx{\hskip 16 true pt\flip{$[{3}{1}:{1}]{\times}[{1^{3}}]$}\hss}\cg%
\cx{\hskip 16 true pt\flip{$[{2^{2}}:{1}]{\times}[{1^{3}}]$}\hss}\cg%
\cx{\hskip 16 true pt\flip{$[{2}{1^{2}}:{1}]{\times}[{1^{3}}]$}\hss}\cg%
\cx{\hskip 16 true pt\flip{$[{1^{4}}:{1}]{\times}[{1^{3}}]$}\hss}\cg%
\cx{\hskip 16 true pt\flip{$[{3}:{2}]{\times}[{1^{3}}]$}\hss}\cg%
\cx{\hskip 16 true pt\flip{$[{3}:{1^{2}}]{\times}[{1^{3}}]$}\hss}\cg%
\cx{\hskip 16 true pt\flip{$[{2}{1}:{2}]{\times}[{1^{3}}]$}\hss}\cg%
\cx{\hskip 16 true pt\flip{$[{2}{1}:{1^{2}}]{\times}[{1^{3}}]$}\hss}\cg%
\cx{\hskip 16 true pt\flip{$[{1^{3}}:{2}]{\times}[{1^{3}}]$}\hss}\cg%
\cx{\hskip 16 true pt\flip{$[{1^{3}}:{1^{2}}]{\times}[{1^{3}}]$}\hss}\cg%
\eol}}\rg%
%
%
\rx{\vss\hfull{%
\rlx{\hss{$8_{z}$}}\cg%
\e{0}%
\e{0}%
\e{0}%
\e{0}%
\e{0}%
\e{0}%
\e{0}%
\e{0}%
\e{0}%
\e{0}%
\e{0}%
\e{0}%
\e{0}%
\e{0}%
\e{0}%
\e{0}%
\e{0}%
\e{0}%
\eol}\vss}\rg%
%
%
\rx{\vss\hfull{%
\rlx{\hss{$56_{z}$}}\cg%
\e{1}%
\e{0}%
\e{0}%
\e{0}%
\e{0}%
\e{0}%
\e{0}%
\e{1}%
\e{0}%
\e{0}%
\e{0}%
\e{0}%
\e{0}%
\e{0}%
\e{0}%
\e{0}%
\e{0}%
\e{0}%
\eol}\vss}\rg%
%
%
\rx{\vss\hfull{%
\rlx{\hss{$160_{z}$}}\cg%
\e{1}%
\e{0}%
\e{0}%
\e{0}%
\e{0}%
\e{0}%
\e{0}%
\e{1}%
\e{0}%
\e{0}%
\e{0}%
\e{0}%
\e{0}%
\e{0}%
\e{0}%
\e{0}%
\e{0}%
\e{0}%
\eol}\vss}\rg%
%
%
\rx{\vss\hfull{%
\rlx{\hss{$112_{z}$}}\cg%
\e{1}%
\e{0}%
\e{0}%
\e{0}%
\e{0}%
\e{0}%
\e{0}%
\e{0}%
\e{0}%
\e{0}%
\e{0}%
\e{0}%
\e{0}%
\e{0}%
\e{0}%
\e{0}%
\e{0}%
\e{0}%
\eol}\vss}\rg%
%
%
\rx{\vss\hfull{%
\rlx{\hss{$840_{z}$}}\cg%
\e{0}%
\e{1}%
\e{0}%
\e{0}%
\e{0}%
\e{0}%
\e{0}%
\e{1}%
\e{1}%
\e{0}%
\e{0}%
\e{0}%
\e{1}%
\e{2}%
\e{1}%
\e{0}%
\e{1}%
\e{0}%
\eol}\vss}\rg%
%
%
\rx{\vss\hfull{%
\rlx{\hss{$1296_{z}$}}\cg%
\e{2}%
\e{2}%
\e{0}%
\e{0}%
\e{0}%
\e{0}%
\e{0}%
\e{5}%
\e{2}%
\e{0}%
\e{0}%
\e{0}%
\e{3}%
\e{2}%
\e{1}%
\e{0}%
\e{0}%
\e{0}%
\eol}\vss}\rg%
%
%
\rx{\vss\hfull{%
\rlx{\hss{$1400_{z}$}}\cg%
\e{2}%
\e{2}%
\e{0}%
\e{0}%
\e{0}%
\e{0}%
\e{0}%
\e{3}%
\e{1}%
\e{0}%
\e{0}%
\e{0}%
\e{2}%
\e{1}%
\e{1}%
\e{0}%
\e{0}%
\e{0}%
\eol}\vss}\rg%
%
%
\rx{\vss\hfull{%
\rlx{\hss{$1008_{z}$}}\cg%
\e{3}%
\e{1}%
\e{0}%
\e{0}%
\e{0}%
\e{0}%
\e{0}%
\e{4}%
\e{1}%
\e{0}%
\e{0}%
\e{0}%
\e{2}%
\e{1}%
\e{0}%
\e{0}%
\e{0}%
\e{0}%
\eol}\vss}\rg%
%
%
\rx{\vss\hfull{%
\rlx{\hss{$560_{z}$}}\cg%
\e{1}%
\e{1}%
\e{0}%
\e{0}%
\e{0}%
\e{0}%
\e{0}%
\e{1}%
\e{0}%
\e{0}%
\e{0}%
\e{0}%
\e{1}%
\e{0}%
\e{0}%
\e{0}%
\e{0}%
\e{0}%
\eol}\vss}\rg%
%
%
\rx{\vss\hfull{%
\rlx{\hss{$1400_{zz}$}}\cg%
\e{0}%
\e{0}%
\e{2}%
\e{0}%
\e{1}%
\e{0}%
\e{0}%
\e{0}%
\e{1}%
\e{2}%
\e{0}%
\e{0}%
\e{1}%
\e{0}%
\e{2}%
\e{1}%
\e{0}%
\e{0}%
\eol}\vss}\rg%
%
%
\rx{\vss\hfull{%
\rlx{\hss{$4200_{z}$}}\cg%
\e{0}%
\e{0}%
\e{4}%
\e{1}%
\e{2}%
\e{1}%
\e{0}%
\e{1}%
\e{6}%
\e{6}%
\e{3}%
\e{0}%
\e{5}%
\e{2}%
\e{6}%
\e{6}%
\e{1}%
\e{1}%
\eol}\vss}\rg%
%
%
\rx{\vss\hfull{%
\rlx{\hss{$400_{z}$}}\cg%
\e{0}%
\e{0}%
\e{1}%
\e{0}%
\e{0}%
\e{0}%
\e{0}%
\e{0}%
\e{0}%
\e{0}%
\e{0}%
\e{0}%
\e{1}%
\e{0}%
\e{0}%
\e{0}%
\e{0}%
\e{0}%
\eol}\vss}\rg%
%
%
\rx{\vss\hfull{%
\rlx{\hss{$3240_{z}$}}\cg%
\e{1}%
\e{2}%
\e{2}%
\e{1}%
\e{0}%
\e{0}%
\e{0}%
\e{3}%
\e{4}%
\e{1}%
\e{0}%
\e{0}%
\e{6}%
\e{2}%
\e{3}%
\e{1}%
\e{0}%
\e{0}%
\eol}\vss}\rg%
%
%
\rx{\vss\hfull{%
\rlx{\hss{$4536_{z}$}}\cg%
\e{0}%
\e{2}%
\e{2}%
\e{2}%
\e{1}%
\e{0}%
\e{0}%
\e{1}%
\e{5}%
\e{3}%
\e{3}%
\e{0}%
\e{3}%
\e{3}%
\e{8}%
\e{5}%
\e{2}%
\e{1}%
\eol}\vss}\rg%
%
%
\rx{\vss\hfull{%
\rlx{\hss{$2400_{z}$}}\cg%
\e{1}%
\e{2}%
\e{1}%
\e{1}%
\e{0}%
\e{0}%
\e{0}%
\e{4}%
\e{5}%
\e{1}%
\e{1}%
\e{0}%
\e{5}%
\e{4}%
\e{4}%
\e{2}%
\e{1}%
\e{0}%
\eol}\vss}\rg%
%
%
\rx{\vss\hfull{%
\rlx{\hss{$3360_{z}$}}\cg%
\e{0}%
\e{1}%
\e{3}%
\e{1}%
\e{1}%
\e{0}%
\e{0}%
\e{2}%
\e{5}%
\e{3}%
\e{1}%
\e{0}%
\e{5}%
\e{2}%
\e{5}%
\e{3}%
\e{0}%
\e{0}%
\eol}\vss}\rg%
%
%
\rx{\vss\hfull{%
\rlx{\hss{$2800_{z}$}}\cg%
\e{2}%
\e{3}%
\e{0}%
\e{1}%
\e{0}%
\e{0}%
\e{0}%
\e{5}%
\e{5}%
\e{1}%
\e{1}%
\e{0}%
\e{5}%
\e{3}%
\e{4}%
\e{2}%
\e{0}%
\e{0}%
\eol}\vss}\rg%
%
%
\rx{\vss\hfull{%
\rlx{\hss{$4096_{z}$}}\cg%
\e{1}%
\e{4}%
\e{1}%
\e{1}%
\e{0}%
\e{0}%
\e{0}%
\e{5}%
\e{6}%
\e{1}%
\e{1}%
\e{0}%
\e{6}%
\e{5}%
\e{6}%
\e{2}%
\e{1}%
\e{0}%
\eol}\vss}\rg%
%
%
\rx{\vss\hfull{%
\rlx{\hss{$5600_{z}$}}\cg%
\e{1}%
\e{4}%
\e{1}%
\e{3}%
\e{0}%
\e{0}%
\e{0}%
\e{5}%
\e{9}%
\e{2}%
\e{4}%
\e{0}%
\e{7}%
\e{8}%
\e{10}%
\e{6}%
\e{4}%
\e{1}%
\eol}\vss}\rg%
%
%
\rx{\vss\hfull{%
\rlx{\hss{$448_{z}$}}\cg%
\e{0}%
\e{1}%
\e{0}%
\e{0}%
\e{0}%
\e{0}%
\e{0}%
\e{0}%
\e{0}%
\e{0}%
\e{0}%
\e{0}%
\e{0}%
\e{0}%
\e{1}%
\e{0}%
\e{0}%
\e{0}%
\eol}\vss}\rg%
%
%
\rx{\vss\hfull{%
\rlx{\hss{$448_{w}$}}\cg%
\e{0}%
\e{1}%
\e{0}%
\e{0}%
\e{0}%
\e{0}%
\e{0}%
\e{1}%
\e{1}%
\e{0}%
\e{0}%
\e{0}%
\e{1}%
\e{2}%
\e{1}%
\e{0}%
\e{1}%
\e{0}%
\eol}\vss}\rg%
%
%
\rx{\vss\hfull{%
\rlx{\hss{$1344_{w}$}}\cg%
\e{0}%
\e{0}%
\e{1}%
\e{0}%
\e{1}%
\e{1}%
\e{0}%
\e{0}%
\e{1}%
\e{2}%
\e{2}%
\e{1}%
\e{0}%
\e{0}%
\e{2}%
\e{3}%
\e{1}%
\e{1}%
\eol}\vss}\rg%
%
%
\rx{\vss\hfull{%
\rlx{\hss{$5600_{w}$}}\cg%
\e{0}%
\e{2}%
\e{3}%
\e{3}%
\e{1}%
\e{1}%
\e{0}%
\e{2}%
\e{9}%
\e{4}%
\e{6}%
\e{1}%
\e{6}%
\e{7}%
\e{11}%
\e{9}%
\e{5}%
\e{2}%
\eol}\vss}\rg%
%
%
\rx{\vss\hfull{%
\rlx{\hss{$2016_{w}$}}\cg%
\e{0}%
\e{0}%
\e{1}%
\e{1}%
\e{3}%
\e{0}%
\e{0}%
\e{0}%
\e{1}%
\e{4}%
\e{3}%
\e{0}%
\e{0}%
\e{0}%
\e{3}%
\e{4}%
\e{1}%
\e{3}%
\eol}\vss}\rg%
%
%
\rx{\vss\hfull{%
\rlx{\hss{$7168_{w}$}}\cg%
\e{0}%
\e{1}%
\e{4}%
\e{4}%
\e{4}%
\e{1}%
\e{0}%
\e{1}%
\e{9}%
\e{8}%
\e{9}%
\e{1}%
\e{5}%
\e{5}%
\e{13}%
\e{13}%
\e{5}%
\e{5}%
\eol}\vss}\rg%
\tableclose%
%
%
%
%
%
%
\eop
\eject
\tableopen{Induce/restrict matrix for $W({E_{6}}{A_{1}})\,\subset\,W(E_{8})$}%
%
%
%
%
%
%
\rowpts=18 true pt%
\colpts=18 true pt%
\rowlabpts=40 true pt%
\collabpts=60 true pt%
\clx{\vss\hfull{%
\rlx{\hss{$ $}}\cg%
\cx{\hskip 16 true pt\flip{$1_p{\times}[{2}]$}\hss}\cg%
\cx{\hskip 16 true pt\flip{$6_p{\times}[{2}]$}\hss}\cg%
\cx{\hskip 16 true pt\flip{$15_p{\times}[{2}]$}\hss}\cg%
\cx{\hskip 16 true pt\flip{$20_p{\times}[{2}]$}\hss}\cg%
\cx{\hskip 16 true pt\flip{$30_p{\times}[{2}]$}\hss}\cg%
\cx{\hskip 16 true pt\flip{$64_p{\times}[{2}]$}\hss}\cg%
\cx{\hskip 16 true pt\flip{$81_p{\times}[{2}]$}\hss}\cg%
\cx{\hskip 16 true pt\flip{$15_q{\times}[{2}]$}\hss}\cg%
\cx{\hskip 16 true pt\flip{$24_p{\times}[{2}]$}\hss}\cg%
\cx{\hskip 16 true pt\flip{$60_p{\times}[{2}]$}\hss}\cg%
\cx{\hskip 16 true pt\flip{$1_p^{*}{\times}[{2}]$}\hss}\cg%
\cx{\hskip 16 true pt\flip{$6_p^{*}{\times}[{2}]$}\hss}\cg%
\cx{\hskip 16 true pt\flip{$15_p^{*}{\times}[{2}]$}\hss}\cg%
\cx{\hskip 16 true pt\flip{$20_p^{*}{\times}[{2}]$}\hss}\cg%
\cx{\hskip 16 true pt\flip{$30_p^{*}{\times}[{2}]$}\hss}\cg%
\cx{\hskip 16 true pt\flip{$64_p^{*}{\times}[{2}]$}\hss}\cg%
\cx{\hskip 16 true pt\flip{$81_p^{*}{\times}[{2}]$}\hss}\cg%
\cx{\hskip 16 true pt\flip{$15_q^{*}{\times}[{2}]$}\hss}\cg%
\cx{\hskip 16 true pt\flip{$24_p^{*}{\times}[{2}]$}\hss}\cg%
\cx{\hskip 16 true pt\flip{$60_p^{*}{\times}[{2}]$}\hss}\cg%
\cx{\hskip 16 true pt\flip{$20_s{\times}[{2}]$}\hss}\cg%
\cx{\hskip 16 true pt\flip{$90_s{\times}[{2}]$}\hss}\cg%
\cx{\hskip 16 true pt\flip{$80_s{\times}[{2}]$}\hss}\cg%
\cx{\hskip 16 true pt\flip{$60_s{\times}[{2}]$}\hss}\cg%
\cx{\hskip 16 true pt\flip{$10_s{\times}[{2}]$}\hss}\cg%
\eol}}\rg%
%
%
\rx{\vss\hfull{%
\rlx{\hss{$1_{x}$}}\cg%
\e{1}%
\e{0}%
\e{0}%
\e{0}%
\e{0}%
\e{0}%
\e{0}%
\e{0}%
\e{0}%
\e{0}%
\e{0}%
\e{0}%
\e{0}%
\e{0}%
\e{0}%
\e{0}%
\e{0}%
\e{0}%
\e{0}%
\e{0}%
\e{0}%
\e{0}%
\e{0}%
\e{0}%
\e{0}%
\eol}\vss}\rg%
%
%
\rx{\vss\hfull{%
\rlx{\hss{$28_{x}$}}\cg%
\e{0}%
\e{1}%
\e{1}%
\e{0}%
\e{0}%
\e{0}%
\e{0}%
\e{0}%
\e{0}%
\e{0}%
\e{0}%
\e{0}%
\e{0}%
\e{0}%
\e{0}%
\e{0}%
\e{0}%
\e{0}%
\e{0}%
\e{0}%
\e{0}%
\e{0}%
\e{0}%
\e{0}%
\e{0}%
\eol}\vss}\rg%
%
%
\rx{\vss\hfull{%
\rlx{\hss{$35_{x}$}}\cg%
\e{2}%
\e{1}%
\e{0}%
\e{1}%
\e{0}%
\e{0}%
\e{0}%
\e{0}%
\e{0}%
\e{0}%
\e{0}%
\e{0}%
\e{0}%
\e{0}%
\e{0}%
\e{0}%
\e{0}%
\e{0}%
\e{0}%
\e{0}%
\e{0}%
\e{0}%
\e{0}%
\e{0}%
\e{0}%
\eol}\vss}\rg%
%
%
\rx{\vss\hfull{%
\rlx{\hss{$84_{x}$}}\cg%
\e{2}%
\e{1}%
\e{0}%
\e{2}%
\e{0}%
\e{0}%
\e{0}%
\e{1}%
\e{0}%
\e{0}%
\e{0}%
\e{0}%
\e{0}%
\e{0}%
\e{0}%
\e{0}%
\e{0}%
\e{0}%
\e{0}%
\e{0}%
\e{0}%
\e{0}%
\e{0}%
\e{0}%
\e{0}%
\eol}\vss}\rg%
%
%
\rx{\vss\hfull{%
\rlx{\hss{$50_{x}$}}\cg%
\e{0}%
\e{0}%
\e{0}%
\e{1}%
\e{0}%
\e{0}%
\e{0}%
\e{1}%
\e{0}%
\e{0}%
\e{0}%
\e{0}%
\e{0}%
\e{0}%
\e{0}%
\e{0}%
\e{0}%
\e{0}%
\e{0}%
\e{0}%
\e{0}%
\e{0}%
\e{0}%
\e{0}%
\e{0}%
\eol}\vss}\rg%
%
%
\rx{\vss\hfull{%
\rlx{\hss{$350_{x}$}}\cg%
\e{0}%
\e{1}%
\e{2}%
\e{0}%
\e{0}%
\e{1}%
\e{0}%
\e{0}%
\e{0}%
\e{0}%
\e{0}%
\e{0}%
\e{0}%
\e{0}%
\e{0}%
\e{0}%
\e{0}%
\e{0}%
\e{0}%
\e{0}%
\e{1}%
\e{1}%
\e{0}%
\e{0}%
\e{0}%
\eol}\vss}\rg%
%
%
\rx{\vss\hfull{%
\rlx{\hss{$300_{x}$}}\cg%
\e{1}%
\e{1}%
\e{0}%
\e{2}%
\e{0}%
\e{1}%
\e{0}%
\e{0}%
\e{1}%
\e{1}%
\e{0}%
\e{0}%
\e{0}%
\e{0}%
\e{0}%
\e{0}%
\e{0}%
\e{0}%
\e{0}%
\e{0}%
\e{0}%
\e{0}%
\e{0}%
\e{0}%
\e{0}%
\eol}\vss}\rg%
%
%
\rx{\vss\hfull{%
\rlx{\hss{$567_{x}$}}\cg%
\e{1}%
\e{3}%
\e{2}%
\e{3}%
\e{2}%
\e{2}%
\e{1}%
\e{0}%
\e{0}%
\e{0}%
\e{0}%
\e{0}%
\e{0}%
\e{0}%
\e{0}%
\e{0}%
\e{0}%
\e{0}%
\e{0}%
\e{0}%
\e{0}%
\e{0}%
\e{0}%
\e{0}%
\e{0}%
\eol}\vss}\rg%
%
%
\rx{\vss\hfull{%
\rlx{\hss{$210_{x}$}}\cg%
\e{1}%
\e{2}%
\e{0}%
\e{2}%
\e{1}%
\e{1}%
\e{0}%
\e{0}%
\e{0}%
\e{0}%
\e{0}%
\e{0}%
\e{0}%
\e{0}%
\e{0}%
\e{0}%
\e{0}%
\e{0}%
\e{0}%
\e{0}%
\e{0}%
\e{0}%
\e{0}%
\e{0}%
\e{0}%
\eol}\vss}\rg%
%
%
\rx{\vss\hfull{%
\rlx{\hss{$840_{x}$}}\cg%
\e{0}%
\e{0}%
\e{0}%
\e{1}%
\e{0}%
\e{1}%
\e{0}%
\e{0}%
\e{2}%
\e{2}%
\e{0}%
\e{0}%
\e{0}%
\e{0}%
\e{0}%
\e{0}%
\e{0}%
\e{0}%
\e{0}%
\e{0}%
\e{0}%
\e{0}%
\e{1}%
\e{2}%
\e{1}%
\eol}\vss}\rg%
%
%
\rx{\vss\hfull{%
\rlx{\hss{$700_{x}$}}\cg%
\e{1}%
\e{2}%
\e{0}%
\e{4}%
\e{2}%
\e{2}%
\e{0}%
\e{2}%
\e{1}%
\e{2}%
\e{0}%
\e{0}%
\e{0}%
\e{0}%
\e{0}%
\e{0}%
\e{0}%
\e{0}%
\e{0}%
\e{0}%
\e{0}%
\e{0}%
\e{0}%
\e{0}%
\e{0}%
\eol}\vss}\rg%
%
%
\rx{\vss\hfull{%
\rlx{\hss{$175_{x}$}}\cg%
\e{0}%
\e{0}%
\e{0}%
\e{0}%
\e{1}%
\e{0}%
\e{0}%
\e{1}%
\e{0}%
\e{1}%
\e{0}%
\e{0}%
\e{0}%
\e{0}%
\e{0}%
\e{0}%
\e{0}%
\e{0}%
\e{0}%
\e{0}%
\e{0}%
\e{0}%
\e{0}%
\e{0}%
\e{0}%
\eol}\vss}\rg%
%
%
\rx{\vss\hfull{%
\rlx{\hss{$1400_{x}$}}\cg%
\e{0}%
\e{1}%
\e{1}%
\e{2}%
\e{4}%
\e{3}%
\e{2}%
\e{1}%
\e{0}%
\e{2}%
\e{0}%
\e{0}%
\e{0}%
\e{0}%
\e{0}%
\e{0}%
\e{0}%
\e{0}%
\e{0}%
\e{0}%
\e{0}%
\e{1}%
\e{1}%
\e{0}%
\e{0}%
\eol}\vss}\rg%
%
%
\rx{\vss\hfull{%
\rlx{\hss{$1050_{x}$}}\cg%
\e{0}%
\e{0}%
\e{1}%
\e{2}%
\e{2}%
\e{2}%
\e{2}%
\e{3}%
\e{0}%
\e{2}%
\e{0}%
\e{0}%
\e{0}%
\e{0}%
\e{0}%
\e{0}%
\e{0}%
\e{0}%
\e{0}%
\e{0}%
\e{0}%
\e{0}%
\e{0}%
\e{1}%
\e{0}%
\eol}\vss}\rg%
%
%
\rx{\vss\hfull{%
\rlx{\hss{$1575_{x}$}}\cg%
\e{0}%
\e{2}%
\e{3}%
\e{2}%
\e{3}%
\e{4}%
\e{2}%
\e{0}%
\e{0}%
\e{1}%
\e{0}%
\e{0}%
\e{0}%
\e{0}%
\e{0}%
\e{0}%
\e{0}%
\e{0}%
\e{0}%
\e{0}%
\e{1}%
\e{2}%
\e{1}%
\e{0}%
\e{0}%
\eol}\vss}\rg%
%
%
\rx{\vss\hfull{%
\rlx{\hss{$1344_{x}$}}\cg%
\e{1}%
\e{2}%
\e{1}%
\e{5}%
\e{2}%
\e{4}%
\e{2}%
\e{2}%
\e{1}%
\e{2}%
\e{0}%
\e{0}%
\e{0}%
\e{0}%
\e{0}%
\e{0}%
\e{0}%
\e{0}%
\e{0}%
\e{0}%
\e{0}%
\e{0}%
\e{0}%
\e{1}%
\e{0}%
\eol}\vss}\rg%
%
%
\rx{\vss\hfull{%
\rlx{\hss{$2100_{x}$}}\cg%
\e{0}%
\e{0}%
\e{2}%
\e{1}%
\e{1}%
\e{3}%
\e{3}%
\e{0}%
\e{1}%
\e{0}%
\e{0}%
\e{0}%
\e{0}%
\e{0}%
\e{0}%
\e{1}%
\e{2}%
\e{0}%
\e{0}%
\e{0}%
\e{2}%
\e{3}%
\e{1}%
\e{0}%
\e{0}%
\eol}\vss}\rg%
%
%
\rx{\vss\hfull{%
\rlx{\hss{$2268_{x}$}}\cg%
\e{0}%
\e{1}%
\e{1}%
\e{3}%
\e{3}%
\e{5}%
\e{3}%
\e{0}%
\e{2}%
\e{2}%
\e{0}%
\e{0}%
\e{0}%
\e{0}%
\e{0}%
\e{0}%
\e{1}%
\e{0}%
\e{0}%
\e{0}%
\e{0}%
\e{2}%
\e{2}%
\e{0}%
\e{0}%
\eol}\vss}\rg%
%
%
\rx{\vss\hfull{%
\rlx{\hss{$525_{x}$}}\cg%
\e{0}%
\e{0}%
\e{1}%
\e{1}%
\e{1}%
\e{1}%
\e{2}%
\e{0}%
\e{0}%
\e{0}%
\e{0}%
\e{0}%
\e{0}%
\e{0}%
\e{0}%
\e{0}%
\e{0}%
\e{0}%
\e{1}%
\e{0}%
\e{0}%
\e{0}%
\e{0}%
\e{0}%
\e{0}%
\eol}\vss}\rg%
%
%
\rx{\vss\hfull{%
\rlx{\hss{$700_{xx}$}}\cg%
\e{0}%
\e{0}%
\e{0}%
\e{0}%
\e{0}%
\e{1}%
\e{2}%
\e{1}%
\e{0}%
\e{1}%
\e{0}%
\e{0}%
\e{0}%
\e{0}%
\e{0}%
\e{0}%
\e{0}%
\e{0}%
\e{1}%
\e{0}%
\e{0}%
\e{0}%
\e{0}%
\e{1}%
\e{0}%
\eol}\vss}\rg%
%
%
\rx{\vss\hfull{%
\rlx{\hss{$972_{x}$}}\cg%
\e{0}%
\e{0}%
\e{0}%
\e{2}%
\e{0}%
\e{2}%
\e{1}%
\e{2}%
\e{2}%
\e{2}%
\e{0}%
\e{0}%
\e{0}%
\e{0}%
\e{0}%
\e{0}%
\e{0}%
\e{0}%
\e{0}%
\e{0}%
\e{0}%
\e{0}%
\e{0}%
\e{2}%
\e{0}%
\eol}\vss}\rg%
%
%
\rx{\vss\hfull{%
\rlx{\hss{$4096_{x}$}}\cg%
\e{0}%
\e{1}%
\e{2}%
\e{3}%
\e{3}%
\e{7}%
\e{5}%
\e{1}%
\e{2}%
\e{4}%
\e{0}%
\e{0}%
\e{0}%
\e{0}%
\e{0}%
\e{0}%
\e{2}%
\e{0}%
\e{0}%
\e{1}%
\e{1}%
\e{4}%
\e{3}%
\e{2}%
\e{0}%
\eol}\vss}\rg%
%
%
\rx{\vss\hfull{%
\rlx{\hss{$4200_{x}$}}\cg%
\e{0}%
\e{0}%
\e{1}%
\e{2}%
\e{3}%
\e{5}%
\e{5}%
\e{2}%
\e{2}%
\e{6}%
\e{0}%
\e{0}%
\e{0}%
\e{0}%
\e{0}%
\e{0}%
\e{2}%
\e{0}%
\e{0}%
\e{1}%
\e{0}%
\e{3}%
\e{4}%
\e{3}%
\e{1}%
\eol}\vss}\rg%
%
%
\rx{\vss\hfull{%
\rlx{\hss{$2240_{x}$}}\cg%
\e{0}%
\e{1}%
\e{0}%
\e{2}%
\e{3}%
\e{4}%
\e{2}%
\e{2}%
\e{1}%
\e{5}%
\e{0}%
\e{0}%
\e{0}%
\e{0}%
\e{0}%
\e{0}%
\e{0}%
\e{0}%
\e{0}%
\e{0}%
\e{0}%
\e{1}%
\e{2}%
\e{2}%
\e{1}%
\eol}\vss}\rg%
%
%
\rx{\vss\hfull{%
\rlx{\hss{$2835_{x}$}}\cg%
\e{0}%
\e{0}%
\e{0}%
\e{0}%
\e{2}%
\e{2}%
\e{3}%
\e{2}%
\e{0}%
\e{5}%
\e{0}%
\e{0}%
\e{0}%
\e{0}%
\e{0}%
\e{0}%
\e{1}%
\e{0}%
\e{0}%
\e{1}%
\e{0}%
\e{2}%
\e{3}%
\e{3}%
\e{1}%
\eol}\vss}\rg%
%
%
\rx{\vss\hfull{%
\rlx{\hss{$6075_{x}$}}\cg%
\e{0}%
\e{0}%
\e{2}%
\e{1}%
\e{3}%
\e{6}%
\e{8}%
\e{1}%
\e{1}%
\e{4}%
\e{0}%
\e{0}%
\e{0}%
\e{0}%
\e{0}%
\e{2}%
\e{3}%
\e{0}%
\e{2}%
\e{2}%
\e{2}%
\e{7}%
\e{5}%
\e{3}%
\e{0}%
\eol}\vss}\rg%
%
%
\rx{\vss\hfull{%
\rlx{\hss{$3200_{x}$}}\cg%
\e{0}%
\e{0}%
\e{0}%
\e{1}%
\e{0}%
\e{3}%
\e{3}%
\e{1}%
\e{3}%
\e{3}%
\e{0}%
\e{0}%
\e{0}%
\e{0}%
\e{0}%
\e{0}%
\e{3}%
\e{1}%
\e{0}%
\e{2}%
\e{0}%
\e{2}%
\e{2}%
\e{4}%
\e{0}%
\eol}\vss}\rg%
%
%
\rx{\vss\hfull{%
\rlx{\hss{$70_{y}$}}\cg%
\e{0}%
\e{0}%
\e{0}%
\e{0}%
\e{0}%
\e{0}%
\e{0}%
\e{0}%
\e{0}%
\e{0}%
\e{0}%
\e{0}%
\e{1}%
\e{0}%
\e{0}%
\e{0}%
\e{0}%
\e{0}%
\e{0}%
\e{0}%
\e{1}%
\e{0}%
\e{0}%
\e{0}%
\e{0}%
\eol}\vss}\rg%
%
%
\rx{\vss\hfull{%
\rlx{\hss{$1134_{y}$}}\cg%
\e{0}%
\e{0}%
\e{1}%
\e{0}%
\e{0}%
\e{1}%
\e{0}%
\e{0}%
\e{0}%
\e{0}%
\e{0}%
\e{0}%
\e{0}%
\e{0}%
\e{0}%
\e{1}%
\e{0}%
\e{0}%
\e{1}%
\e{1}%
\e{1}%
\e{2}%
\e{1}%
\e{1}%
\e{0}%
\eol}\vss}\rg%
%
%
\rx{\vss\hfull{%
\rlx{\hss{$1680_{y}$}}\cg%
\e{0}%
\e{0}%
\e{1}%
\e{0}%
\e{0}%
\e{1}%
\e{1}%
\e{0}%
\e{0}%
\e{0}%
\e{0}%
\e{0}%
\e{2}%
\e{0}%
\e{1}%
\e{2}%
\e{2}%
\e{0}%
\e{0}%
\e{0}%
\e{3}%
\e{3}%
\e{0}%
\e{0}%
\e{0}%
\eol}\vss}\rg%
%
%
\rx{\vss\hfull{%
\rlx{\hss{$168_{y}$}}\cg%
\e{0}%
\e{0}%
\e{0}%
\e{0}%
\e{0}%
\e{0}%
\e{0}%
\e{0}%
\e{1}%
\e{0}%
\e{0}%
\e{0}%
\e{0}%
\e{0}%
\e{0}%
\e{0}%
\e{0}%
\e{0}%
\e{0}%
\e{0}%
\e{0}%
\e{0}%
\e{0}%
\e{1}%
\e{0}%
\eol}\vss}\rg%
%
%
\rx{\vss\hfull{%
\rlx{\hss{$420_{y}$}}\cg%
\e{0}%
\e{0}%
\e{0}%
\e{0}%
\e{0}%
\e{0}%
\e{0}%
\e{0}%
\e{0}%
\e{1}%
\e{0}%
\e{0}%
\e{0}%
\e{0}%
\e{0}%
\e{0}%
\e{0}%
\e{0}%
\e{0}%
\e{0}%
\e{0}%
\e{0}%
\e{1}%
\e{1}%
\e{1}%
\eol}\vss}\rg%
%
%
\rx{\vss\hfull{%
\rlx{\hss{$3150_{y}$}}\cg%
\e{0}%
\e{0}%
\e{0}%
\e{0}%
\e{1}%
\e{1}%
\e{2}%
\e{0}%
\e{0}%
\e{3}%
\e{0}%
\e{0}%
\e{0}%
\e{0}%
\e{0}%
\e{1}%
\e{1}%
\e{0}%
\e{1}%
\e{3}%
\e{0}%
\e{3}%
\e{4}%
\e{3}%
\e{2}%
\eol}\vss}\rg%
%
%
\rx{\vss\hfull{%
\rlx{\hss{$4200_{y}$}}\cg%
\e{0}%
\e{0}%
\e{0}%
\e{0}%
\e{0}%
\e{2}%
\e{3}%
\e{1}%
\e{2}%
\e{3}%
\e{0}%
\e{0}%
\e{0}%
\e{0}%
\e{0}%
\e{1}%
\e{3}%
\e{1}%
\e{1}%
\e{3}%
\e{0}%
\e{3}%
\e{4}%
\e{6}%
\e{1}%
\eol}\vss}\rg%
\eop
\eject
\tablecont%
%
%
%
%
%
%
\rowpts=18 true pt%
\colpts=18 true pt%
\rowlabpts=40 true pt%
\collabpts=60 true pt%
\clx{\vss\hfull{%
\rlx{\hss{$ $}}\cg%
\cx{\hskip 16 true pt\flip{$1_p{\times}[{2}]$}\hss}\cg%
\cx{\hskip 16 true pt\flip{$6_p{\times}[{2}]$}\hss}\cg%
\cx{\hskip 16 true pt\flip{$15_p{\times}[{2}]$}\hss}\cg%
\cx{\hskip 16 true pt\flip{$20_p{\times}[{2}]$}\hss}\cg%
\cx{\hskip 16 true pt\flip{$30_p{\times}[{2}]$}\hss}\cg%
\cx{\hskip 16 true pt\flip{$64_p{\times}[{2}]$}\hss}\cg%
\cx{\hskip 16 true pt\flip{$81_p{\times}[{2}]$}\hss}\cg%
\cx{\hskip 16 true pt\flip{$15_q{\times}[{2}]$}\hss}\cg%
\cx{\hskip 16 true pt\flip{$24_p{\times}[{2}]$}\hss}\cg%
\cx{\hskip 16 true pt\flip{$60_p{\times}[{2}]$}\hss}\cg%
\cx{\hskip 16 true pt\flip{$1_p^{*}{\times}[{2}]$}\hss}\cg%
\cx{\hskip 16 true pt\flip{$6_p^{*}{\times}[{2}]$}\hss}\cg%
\cx{\hskip 16 true pt\flip{$15_p^{*}{\times}[{2}]$}\hss}\cg%
\cx{\hskip 16 true pt\flip{$20_p^{*}{\times}[{2}]$}\hss}\cg%
\cx{\hskip 16 true pt\flip{$30_p^{*}{\times}[{2}]$}\hss}\cg%
\cx{\hskip 16 true pt\flip{$64_p^{*}{\times}[{2}]$}\hss}\cg%
\cx{\hskip 16 true pt\flip{$81_p^{*}{\times}[{2}]$}\hss}\cg%
\cx{\hskip 16 true pt\flip{$15_q^{*}{\times}[{2}]$}\hss}\cg%
\cx{\hskip 16 true pt\flip{$24_p^{*}{\times}[{2}]$}\hss}\cg%
\cx{\hskip 16 true pt\flip{$60_p^{*}{\times}[{2}]$}\hss}\cg%
\cx{\hskip 16 true pt\flip{$20_s{\times}[{2}]$}\hss}\cg%
\cx{\hskip 16 true pt\flip{$90_s{\times}[{2}]$}\hss}\cg%
\cx{\hskip 16 true pt\flip{$80_s{\times}[{2}]$}\hss}\cg%
\cx{\hskip 16 true pt\flip{$60_s{\times}[{2}]$}\hss}\cg%
\cx{\hskip 16 true pt\flip{$10_s{\times}[{2}]$}\hss}\cg%
\eol}}\rg%
%
%
\rx{\vss\hfull{%
\rlx{\hss{$2688_{y}$}}\cg%
\e{0}%
\e{0}%
\e{0}%
\e{0}%
\e{0}%
\e{1}%
\e{3}%
\e{1}%
\e{0}%
\e{2}%
\e{0}%
\e{0}%
\e{0}%
\e{0}%
\e{0}%
\e{1}%
\e{3}%
\e{1}%
\e{0}%
\e{2}%
\e{1}%
\e{2}%
\e{1}%
\e{3}%
\e{0}%
\eol}\vss}\rg%
%
%
\rx{\vss\hfull{%
\rlx{\hss{$2100_{y}$}}\cg%
\e{0}%
\e{0}%
\e{0}%
\e{0}%
\e{0}%
\e{1}%
\e{3}%
\e{0}%
\e{1}%
\e{0}%
\e{0}%
\e{0}%
\e{0}%
\e{1}%
\e{0}%
\e{2}%
\e{3}%
\e{0}%
\e{2}%
\e{0}%
\e{1}%
\e{2}%
\e{1}%
\e{0}%
\e{0}%
\eol}\vss}\rg%
%
%
\rx{\vss\hfull{%
\rlx{\hss{$1400_{y}$}}\cg%
\e{0}%
\e{0}%
\e{0}%
\e{0}%
\e{0}%
\e{1}%
\e{1}%
\e{0}%
\e{1}%
\e{0}%
\e{0}%
\e{0}%
\e{1}%
\e{0}%
\e{1}%
\e{1}%
\e{2}%
\e{0}%
\e{0}%
\e{0}%
\e{0}%
\e{2}%
\e{1}%
\e{0}%
\e{0}%
\eol}\vss}\rg%
%
%
\rx{\vss\hfull{%
\rlx{\hss{$4536_{y}$}}\cg%
\e{0}%
\e{0}%
\e{0}%
\e{0}%
\e{1}%
\e{2}%
\e{4}%
\e{0}%
\e{1}%
\e{2}%
\e{0}%
\e{0}%
\e{1}%
\e{0}%
\e{2}%
\e{3}%
\e{5}%
\e{0}%
\e{0}%
\e{2}%
\e{1}%
\e{5}%
\e{4}%
\e{1}%
\e{0}%
\eol}\vss}\rg%
%
%
\rx{\vss\hfull{%
\rlx{\hss{$5670_{y}$}}\cg%
\e{0}%
\e{0}%
\e{1}%
\e{0}%
\e{1}%
\e{3}%
\e{4}%
\e{0}%
\e{1}%
\e{2}%
\e{0}%
\e{0}%
\e{1}%
\e{0}%
\e{2}%
\e{4}%
\e{5}%
\e{0}%
\e{1}%
\e{3}%
\e{2}%
\e{7}%
\e{5}%
\e{2}%
\e{0}%
\eol}\vss}\rg%
%
%
\rx{\vss\hfull{%
\rlx{\hss{$4480_{y}$}}\cg%
\e{0}%
\e{0}%
\e{0}%
\e{0}%
\e{1}%
\e{2}%
\e{3}%
\e{0}%
\e{1}%
\e{3}%
\e{0}%
\e{0}%
\e{0}%
\e{0}%
\e{1}%
\e{2}%
\e{3}%
\e{0}%
\e{1}%
\e{3}%
\e{0}%
\e{4}%
\e{6}%
\e{3}%
\e{1}%
\eol}\vss}\rg%
%
%
\rx{\vss\hfull{%
\rlx{\hss{$8_{z}$}}\cg%
\e{1}%
\e{1}%
\e{0}%
\e{0}%
\e{0}%
\e{0}%
\e{0}%
\e{0}%
\e{0}%
\e{0}%
\e{0}%
\e{0}%
\e{0}%
\e{0}%
\e{0}%
\e{0}%
\e{0}%
\e{0}%
\e{0}%
\e{0}%
\e{0}%
\e{0}%
\e{0}%
\e{0}%
\e{0}%
\eol}\vss}\rg%
%
%
\rx{\vss\hfull{%
\rlx{\hss{$56_{z}$}}\cg%
\e{0}%
\e{0}%
\e{1}%
\e{0}%
\e{0}%
\e{0}%
\e{0}%
\e{0}%
\e{0}%
\e{0}%
\e{0}%
\e{0}%
\e{0}%
\e{0}%
\e{0}%
\e{0}%
\e{0}%
\e{0}%
\e{0}%
\e{0}%
\e{1}%
\e{0}%
\e{0}%
\e{0}%
\e{0}%
\eol}\vss}\rg%
%
%
\rx{\vss\hfull{%
\rlx{\hss{$160_{z}$}}\cg%
\e{1}%
\e{2}%
\e{1}%
\e{1}%
\e{0}%
\e{1}%
\e{0}%
\e{0}%
\e{0}%
\e{0}%
\e{0}%
\e{0}%
\e{0}%
\e{0}%
\e{0}%
\e{0}%
\e{0}%
\e{0}%
\e{0}%
\e{0}%
\e{0}%
\e{0}%
\e{0}%
\e{0}%
\e{0}%
\eol}\vss}\rg%
%
%
\rx{\vss\hfull{%
\rlx{\hss{$112_{z}$}}\cg%
\e{2}%
\e{2}%
\e{0}%
\e{2}%
\e{1}%
\e{0}%
\e{0}%
\e{0}%
\e{0}%
\e{0}%
\e{0}%
\e{0}%
\e{0}%
\e{0}%
\e{0}%
\e{0}%
\e{0}%
\e{0}%
\e{0}%
\e{0}%
\e{0}%
\e{0}%
\e{0}%
\e{0}%
\e{0}%
\eol}\vss}\rg%
%
%
\rx{\vss\hfull{%
\rlx{\hss{$840_{z}$}}\cg%
\e{0}%
\e{1}%
\e{1}%
\e{1}%
\e{0}%
\e{2}%
\e{0}%
\e{0}%
\e{1}%
\e{1}%
\e{0}%
\e{0}%
\e{0}%
\e{0}%
\e{0}%
\e{0}%
\e{0}%
\e{0}%
\e{0}%
\e{0}%
\e{0}%
\e{1}%
\e{1}%
\e{1}%
\e{0}%
\eol}\vss}\rg%
%
%
\rx{\vss\hfull{%
\rlx{\hss{$1296_{z}$}}\cg%
\e{0}%
\e{1}%
\e{3}%
\e{1}%
\e{1}%
\e{3}%
\e{2}%
\e{0}%
\e{0}%
\e{0}%
\e{0}%
\e{0}%
\e{0}%
\e{0}%
\e{0}%
\e{0}%
\e{1}%
\e{0}%
\e{0}%
\e{0}%
\e{2}%
\e{2}%
\e{0}%
\e{0}%
\e{0}%
\eol}\vss}\rg%
%
%
\rx{\vss\hfull{%
\rlx{\hss{$1400_{z}$}}\cg%
\e{1}%
\e{2}%
\e{1}%
\e{5}%
\e{3}%
\e{4}%
\e{2}%
\e{1}%
\e{1}%
\e{2}%
\e{0}%
\e{0}%
\e{0}%
\e{0}%
\e{0}%
\e{0}%
\e{0}%
\e{0}%
\e{0}%
\e{0}%
\e{0}%
\e{0}%
\e{1}%
\e{0}%
\e{0}%
\eol}\vss}\rg%
%
%
\rx{\vss\hfull{%
\rlx{\hss{$1008_{z}$}}\cg%
\e{0}%
\e{2}%
\e{2}%
\e{3}%
\e{2}%
\e{3}%
\e{2}%
\e{0}%
\e{1}%
\e{0}%
\e{0}%
\e{0}%
\e{0}%
\e{0}%
\e{0}%
\e{0}%
\e{0}%
\e{0}%
\e{0}%
\e{0}%
\e{0}%
\e{1}%
\e{0}%
\e{0}%
\e{0}%
\eol}\vss}\rg%
%
%
\rx{\vss\hfull{%
\rlx{\hss{$560_{z}$}}\cg%
\e{2}%
\e{3}%
\e{1}%
\e{4}%
\e{2}%
\e{2}%
\e{0}%
\e{1}%
\e{0}%
\e{1}%
\e{0}%
\e{0}%
\e{0}%
\e{0}%
\e{0}%
\e{0}%
\e{0}%
\e{0}%
\e{0}%
\e{0}%
\e{0}%
\e{0}%
\e{0}%
\e{0}%
\e{0}%
\eol}\vss}\rg%
%
%
\rx{\vss\hfull{%
\rlx{\hss{$1400_{zz}$}}\cg%
\e{0}%
\e{0}%
\e{0}%
\e{1}%
\e{2}%
\e{2}%
\e{2}%
\e{3}%
\e{0}%
\e{4}%
\e{0}%
\e{0}%
\e{0}%
\e{0}%
\e{0}%
\e{0}%
\e{0}%
\e{0}%
\e{0}%
\e{0}%
\e{0}%
\e{0}%
\e{1}%
\e{1}%
\e{1}%
\eol}\vss}\rg%
%
%
\rx{\vss\hfull{%
\rlx{\hss{$4200_{z}$}}\cg%
\e{0}%
\e{0}%
\e{1}%
\e{0}%
\e{2}%
\e{3}%
\e{5}%
\e{1}%
\e{1}%
\e{4}%
\e{0}%
\e{0}%
\e{0}%
\e{0}%
\e{0}%
\e{1}%
\e{2}%
\e{0}%
\e{2}%
\e{1}%
\e{0}%
\e{4}%
\e{4}%
\e{4}%
\e{0}%
\eol}\vss}\rg%
%
%
\rx{\vss\hfull{%
\rlx{\hss{$400_{z}$}}\cg%
\e{0}%
\e{1}%
\e{0}%
\e{2}%
\e{2}%
\e{1}%
\e{0}%
\e{2}%
\e{0}%
\e{1}%
\e{0}%
\e{0}%
\e{0}%
\e{0}%
\e{0}%
\e{0}%
\e{0}%
\e{0}%
\e{0}%
\e{0}%
\e{0}%
\e{0}%
\e{0}%
\e{0}%
\e{0}%
\eol}\vss}\rg%
%
%
\rx{\vss\hfull{%
\rlx{\hss{$3240_{z}$}}\cg%
\e{0}%
\e{2}%
\e{2}%
\e{5}%
\e{5}%
\e{7}%
\e{4}%
\e{3}%
\e{2}%
\e{5}%
\e{0}%
\e{0}%
\e{0}%
\e{0}%
\e{0}%
\e{0}%
\e{0}%
\e{0}%
\e{0}%
\e{0}%
\e{0}%
\e{2}%
\e{2}%
\e{2}%
\e{0}%
\eol}\vss}\rg%
%
%
\rx{\vss\hfull{%
\rlx{\hss{$4536_{z}$}}\cg%
\e{0}%
\e{0}%
\e{0}%
\e{2}%
\e{2}%
\e{5}%
\e{4}%
\e{2}%
\e{3}%
\e{7}%
\e{0}%
\e{0}%
\e{0}%
\e{0}%
\e{0}%
\e{0}%
\e{1}%
\e{0}%
\e{0}%
\e{2}%
\e{0}%
\e{3}%
\e{5}%
\e{5}%
\e{2}%
\eol}\vss}\rg%
%
%
\rx{\vss\hfull{%
\rlx{\hss{$2400_{z}$}}\cg%
\e{0}%
\e{0}%
\e{2}%
\e{0}%
\e{1}%
\e{2}%
\e{3}%
\e{0}%
\e{0}%
\e{0}%
\e{0}%
\e{0}%
\e{1}%
\e{0}%
\e{0}%
\e{2}%
\e{2}%
\e{0}%
\e{1}%
\e{0}%
\e{3}%
\e{4}%
\e{1}%
\e{0}%
\e{0}%
\eol}\vss}\rg%
%
%
\rx{\vss\hfull{%
\rlx{\hss{$3360_{z}$}}\cg%
\e{0}%
\e{0}%
\e{1}%
\e{1}%
\e{3}%
\e{4}%
\e{5}%
\e{2}%
\e{0}%
\e{4}%
\e{0}%
\e{0}%
\e{0}%
\e{0}%
\e{0}%
\e{0}%
\e{2}%
\e{0}%
\e{0}%
\e{1}%
\e{1}%
\e{3}%
\e{2}%
\e{2}%
\e{0}%
\eol}\vss}\rg%
%
%
\rx{\vss\hfull{%
\rlx{\hss{$2800_{z}$}}\cg%
\e{0}%
\e{0}%
\e{1}%
\e{1}%
\e{2}%
\e{4}%
\e{5}%
\e{0}%
\e{1}%
\e{1}%
\e{0}%
\e{0}%
\e{0}%
\e{0}%
\e{0}%
\e{1}%
\e{2}%
\e{0}%
\e{1}%
\e{0}%
\e{1}%
\e{3}%
\e{2}%
\e{0}%
\e{0}%
\eol}\vss}\rg%
%
%
\rx{\vss\hfull{%
\rlx{\hss{$4096_{z}$}}\cg%
\e{0}%
\e{1}%
\e{2}%
\e{3}%
\e{3}%
\e{7}%
\e{5}%
\e{1}%
\e{2}%
\e{4}%
\e{0}%
\e{0}%
\e{0}%
\e{0}%
\e{0}%
\e{0}%
\e{2}%
\e{0}%
\e{0}%
\e{1}%
\e{1}%
\e{4}%
\e{3}%
\e{2}%
\e{0}%
\eol}\vss}\rg%
%
%
\rx{\vss\hfull{%
\rlx{\hss{$5600_{z}$}}\cg%
\e{0}%
\e{0}%
\e{1}%
\e{1}%
\e{2}%
\e{5}%
\e{6}%
\e{0}%
\e{3}%
\e{3}%
\e{0}%
\e{0}%
\e{0}%
\e{0}%
\e{1}%
\e{2}%
\e{5}%
\e{0}%
\e{1}%
\e{2}%
\e{1}%
\e{6}%
\e{5}%
\e{2}%
\e{0}%
\eol}\vss}\rg%
%
%
\rx{\vss\hfull{%
\rlx{\hss{$448_{z}$}}\cg%
\e{1}%
\e{0}%
\e{0}%
\e{2}%
\e{1}%
\e{1}%
\e{0}%
\e{1}%
\e{0}%
\e{2}%
\e{0}%
\e{0}%
\e{0}%
\e{0}%
\e{0}%
\e{0}%
\e{0}%
\e{0}%
\e{0}%
\e{0}%
\e{0}%
\e{0}%
\e{0}%
\e{0}%
\e{1}%
\eol}\vss}\rg%
%
%
\rx{\vss\hfull{%
\rlx{\hss{$448_{w}$}}\cg%
\e{0}%
\e{0}%
\e{1}%
\e{0}%
\e{0}%
\e{0}%
\e{0}%
\e{0}%
\e{0}%
\e{0}%
\e{0}%
\e{0}%
\e{1}%
\e{0}%
\e{0}%
\e{1}%
\e{0}%
\e{0}%
\e{0}%
\e{0}%
\e{2}%
\e{1}%
\e{0}%
\e{0}%
\e{0}%
\eol}\vss}\rg%
%
%
\rx{\vss\hfull{%
\rlx{\hss{$1344_{w}$}}\cg%
\e{0}%
\e{0}%
\e{0}%
\e{0}%
\e{0}%
\e{1}%
\e{0}%
\e{0}%
\e{1}%
\e{1}%
\e{0}%
\e{0}%
\e{0}%
\e{0}%
\e{0}%
\e{0}%
\e{0}%
\e{0}%
\e{1}%
\e{1}%
\e{0}%
\e{1}%
\e{2}%
\e{3}%
\e{1}%
\eol}\vss}\rg%
%
%
\rx{\vss\hfull{%
\rlx{\hss{$5600_{w}$}}\cg%
\e{0}%
\e{0}%
\e{1}%
\e{0}%
\e{1}%
\e{3}%
\e{4}%
\e{0}%
\e{1}%
\e{2}%
\e{0}%
\e{0}%
\e{1}%
\e{0}%
\e{1}%
\e{4}%
\e{5}%
\e{1}%
\e{1}%
\e{3}%
\e{2}%
\e{7}%
\e{4}%
\e{3}%
\e{0}%
\eol}\vss}\rg%
%
%
\rx{\vss\hfull{%
\rlx{\hss{$2016_{w}$}}\cg%
\e{0}%
\e{0}%
\e{0}%
\e{0}%
\e{0}%
\e{0}%
\e{2}%
\e{1}%
\e{0}%
\e{3}%
\e{0}%
\e{0}%
\e{0}%
\e{0}%
\e{0}%
\e{0}%
\e{1}%
\e{0}%
\e{0}%
\e{2}%
\e{0}%
\e{1}%
\e{2}%
\e{3}%
\e{2}%
\eol}\vss}\rg%
%
%
\rx{\vss\hfull{%
\rlx{\hss{$7168_{w}$}}\cg%
\e{0}%
\e{0}%
\e{0}%
\e{0}%
\e{1}%
\e{3}%
\e{6}%
\e{1}%
\e{1}%
\e{5}%
\e{0}%
\e{0}%
\e{0}%
\e{0}%
\e{1}%
\e{3}%
\e{6}%
\e{1}%
\e{1}%
\e{5}%
\e{1}%
\e{6}%
\e{7}%
\e{6}%
\e{1}%
\eol}\vss}\rg%
\eop
\eject
\tablecont%
%
%
%
%
%
%
\rowpts=18 true pt%
\colpts=18 true pt%
\rowlabpts=40 true pt%
\collabpts=60 true pt%
\clx{\vss\hfull{%
\rlx{\hss{$ $}}\cg%
\cx{\hskip 16 true pt\flip{$1_p{\times}[{1^{2}}]$}\hss}\cg%
\cx{\hskip 16 true pt\flip{$6_p{\times}[{1^{2}}]$}\hss}\cg%
\cx{\hskip 16 true pt\flip{$15_p{\times}[{1^{2}}]$}\hss}\cg%
\cx{\hskip 16 true pt\flip{$20_p{\times}[{1^{2}}]$}\hss}\cg%
\cx{\hskip 16 true pt\flip{$30_p{\times}[{1^{2}}]$}\hss}\cg%
\cx{\hskip 16 true pt\flip{$64_p{\times}[{1^{2}}]$}\hss}\cg%
\cx{\hskip 16 true pt\flip{$81_p{\times}[{1^{2}}]$}\hss}\cg%
\cx{\hskip 16 true pt\flip{$15_q{\times}[{1^{2}}]$}\hss}\cg%
\cx{\hskip 16 true pt\flip{$24_p{\times}[{1^{2}}]$}\hss}\cg%
\cx{\hskip 16 true pt\flip{$60_p{\times}[{1^{2}}]$}\hss}\cg%
\cx{\hskip 16 true pt\flip{$1_p^{*}{\times}[{1^{2}}]$}\hss}\cg%
\cx{\hskip 16 true pt\flip{$6_p^{*}{\times}[{1^{2}}]$}\hss}\cg%
\cx{\hskip 16 true pt\flip{$15_p^{*}{\times}[{1^{2}}]$}\hss}\cg%
\cx{\hskip 16 true pt\flip{$20_p^{*}{\times}[{1^{2}}]$}\hss}\cg%
\cx{\hskip 16 true pt\flip{$30_p^{*}{\times}[{1^{2}}]$}\hss}\cg%
\cx{\hskip 16 true pt\flip{$64_p^{*}{\times}[{1^{2}}]$}\hss}\cg%
\cx{\hskip 16 true pt\flip{$81_p^{*}{\times}[{1^{2}}]$}\hss}\cg%
\cx{\hskip 16 true pt\flip{$15_q^{*}{\times}[{1^{2}}]$}\hss}\cg%
\cx{\hskip 16 true pt\flip{$24_p^{*}{\times}[{1^{2}}]$}\hss}\cg%
\cx{\hskip 16 true pt\flip{$60_p^{*}{\times}[{1^{2}}]$}\hss}\cg%
\cx{\hskip 16 true pt\flip{$20_s{\times}[{1^{2}}]$}\hss}\cg%
\cx{\hskip 16 true pt\flip{$90_s{\times}[{1^{2}}]$}\hss}\cg%
\cx{\hskip 16 true pt\flip{$80_s{\times}[{1^{2}}]$}\hss}\cg%
\cx{\hskip 16 true pt\flip{$60_s{\times}[{1^{2}}]$}\hss}\cg%
\cx{\hskip 16 true pt\flip{$10_s{\times}[{1^{2}}]$}\hss}\cg%
\eol}}\rg%
%
%
\rx{\vss\hfull{%
\rlx{\hss{$1_{x}$}}\cg%
\e{0}%
\e{0}%
\e{0}%
\e{0}%
\e{0}%
\e{0}%
\e{0}%
\e{0}%
\e{0}%
\e{0}%
\e{0}%
\e{0}%
\e{0}%
\e{0}%
\e{0}%
\e{0}%
\e{0}%
\e{0}%
\e{0}%
\e{0}%
\e{0}%
\e{0}%
\e{0}%
\e{0}%
\e{0}%
\eol}\vss}\rg%
%
%
\rx{\vss\hfull{%
\rlx{\hss{$28_{x}$}}\cg%
\e{1}%
\e{1}%
\e{0}%
\e{0}%
\e{0}%
\e{0}%
\e{0}%
\e{0}%
\e{0}%
\e{0}%
\e{0}%
\e{0}%
\e{0}%
\e{0}%
\e{0}%
\e{0}%
\e{0}%
\e{0}%
\e{0}%
\e{0}%
\e{0}%
\e{0}%
\e{0}%
\e{0}%
\e{0}%
\eol}\vss}\rg%
%
%
\rx{\vss\hfull{%
\rlx{\hss{$35_{x}$}}\cg%
\e{1}%
\e{1}%
\e{0}%
\e{0}%
\e{0}%
\e{0}%
\e{0}%
\e{0}%
\e{0}%
\e{0}%
\e{0}%
\e{0}%
\e{0}%
\e{0}%
\e{0}%
\e{0}%
\e{0}%
\e{0}%
\e{0}%
\e{0}%
\e{0}%
\e{0}%
\e{0}%
\e{0}%
\e{0}%
\eol}\vss}\rg%
%
%
\rx{\vss\hfull{%
\rlx{\hss{$84_{x}$}}\cg%
\e{1}%
\e{0}%
\e{0}%
\e{1}%
\e{0}%
\e{0}%
\e{0}%
\e{0}%
\e{0}%
\e{0}%
\e{0}%
\e{0}%
\e{0}%
\e{0}%
\e{0}%
\e{0}%
\e{0}%
\e{0}%
\e{0}%
\e{0}%
\e{0}%
\e{0}%
\e{0}%
\e{0}%
\e{0}%
\eol}\vss}\rg%
%
%
\rx{\vss\hfull{%
\rlx{\hss{$50_{x}$}}\cg%
\e{0}%
\e{0}%
\e{0}%
\e{0}%
\e{0}%
\e{0}%
\e{0}%
\e{1}%
\e{0}%
\e{0}%
\e{0}%
\e{0}%
\e{0}%
\e{0}%
\e{0}%
\e{0}%
\e{0}%
\e{0}%
\e{0}%
\e{0}%
\e{0}%
\e{0}%
\e{0}%
\e{0}%
\e{0}%
\eol}\vss}\rg%
%
%
\rx{\vss\hfull{%
\rlx{\hss{$350_{x}$}}\cg%
\e{0}%
\e{1}%
\e{2}%
\e{1}%
\e{0}%
\e{1}%
\e{0}%
\e{0}%
\e{0}%
\e{0}%
\e{0}%
\e{0}%
\e{0}%
\e{0}%
\e{0}%
\e{0}%
\e{0}%
\e{0}%
\e{0}%
\e{0}%
\e{1}%
\e{0}%
\e{0}%
\e{0}%
\e{0}%
\eol}\vss}\rg%
%
%
\rx{\vss\hfull{%
\rlx{\hss{$300_{x}$}}\cg%
\e{0}%
\e{1}%
\e{1}%
\e{1}%
\e{0}%
\e{1}%
\e{0}%
\e{0}%
\e{0}%
\e{0}%
\e{0}%
\e{0}%
\e{0}%
\e{0}%
\e{0}%
\e{0}%
\e{0}%
\e{0}%
\e{0}%
\e{0}%
\e{0}%
\e{0}%
\e{0}%
\e{0}%
\e{0}%
\eol}\vss}\rg%
%
%
\rx{\vss\hfull{%
\rlx{\hss{$567_{x}$}}\cg%
\e{2}%
\e{3}%
\e{1}%
\e{3}%
\e{1}%
\e{1}%
\e{0}%
\e{0}%
\e{0}%
\e{0}%
\e{0}%
\e{0}%
\e{0}%
\e{0}%
\e{0}%
\e{0}%
\e{0}%
\e{0}%
\e{0}%
\e{0}%
\e{0}%
\e{0}%
\e{0}%
\e{0}%
\e{0}%
\eol}\vss}\rg%
%
%
\rx{\vss\hfull{%
\rlx{\hss{$210_{x}$}}\cg%
\e{1}%
\e{2}%
\e{0}%
\e{1}%
\e{1}%
\e{0}%
\e{0}%
\e{0}%
\e{0}%
\e{0}%
\e{0}%
\e{0}%
\e{0}%
\e{0}%
\e{0}%
\e{0}%
\e{0}%
\e{0}%
\e{0}%
\e{0}%
\e{0}%
\e{0}%
\e{0}%
\e{0}%
\e{0}%
\eol}\vss}\rg%
%
%
\rx{\vss\hfull{%
\rlx{\hss{$840_{x}$}}\cg%
\e{0}%
\e{0}%
\e{0}%
\e{0}%
\e{0}%
\e{1}%
\e{0}%
\e{0}%
\e{1}%
\e{1}%
\e{0}%
\e{0}%
\e{0}%
\e{0}%
\e{0}%
\e{0}%
\e{0}%
\e{0}%
\e{0}%
\e{0}%
\e{0}%
\e{1}%
\e{1}%
\e{1}%
\e{0}%
\eol}\vss}\rg%
%
%
\rx{\vss\hfull{%
\rlx{\hss{$700_{x}$}}\cg%
\e{0}%
\e{1}%
\e{0}%
\e{2}%
\e{2}%
\e{1}%
\e{0}%
\e{1}%
\e{0}%
\e{1}%
\e{0}%
\e{0}%
\e{0}%
\e{0}%
\e{0}%
\e{0}%
\e{0}%
\e{0}%
\e{0}%
\e{0}%
\e{0}%
\e{0}%
\e{0}%
\e{0}%
\e{0}%
\eol}\vss}\rg%
%
%
\rx{\vss\hfull{%
\rlx{\hss{$175_{x}$}}\cg%
\e{0}%
\e{0}%
\e{0}%
\e{0}%
\e{0}%
\e{0}%
\e{0}%
\e{0}%
\e{0}%
\e{1}%
\e{0}%
\e{0}%
\e{0}%
\e{0}%
\e{0}%
\e{0}%
\e{0}%
\e{0}%
\e{0}%
\e{0}%
\e{0}%
\e{0}%
\e{0}%
\e{0}%
\e{1}%
\eol}\vss}\rg%
%
%
\rx{\vss\hfull{%
\rlx{\hss{$1400_{x}$}}\cg%
\e{0}%
\e{1}%
\e{0}%
\e{2}%
\e{3}%
\e{2}%
\e{1}%
\e{1}%
\e{0}%
\e{2}%
\e{0}%
\e{0}%
\e{0}%
\e{0}%
\e{0}%
\e{0}%
\e{0}%
\e{0}%
\e{0}%
\e{0}%
\e{0}%
\e{0}%
\e{1}%
\e{0}%
\e{0}%
\eol}\vss}\rg%
%
%
\rx{\vss\hfull{%
\rlx{\hss{$1050_{x}$}}\cg%
\e{0}%
\e{0}%
\e{0}%
\e{1}%
\e{1}%
\e{1}%
\e{1}%
\e{3}%
\e{0}%
\e{2}%
\e{0}%
\e{0}%
\e{0}%
\e{0}%
\e{0}%
\e{0}%
\e{0}%
\e{0}%
\e{0}%
\e{0}%
\e{0}%
\e{0}%
\e{0}%
\e{1}%
\e{0}%
\eol}\vss}\rg%
%
%
\rx{\vss\hfull{%
\rlx{\hss{$1575_{x}$}}\cg%
\e{0}%
\e{2}%
\e{2}%
\e{3}%
\e{3}%
\e{3}%
\e{1}%
\e{1}%
\e{0}%
\e{1}%
\e{0}%
\e{0}%
\e{0}%
\e{0}%
\e{0}%
\e{0}%
\e{0}%
\e{0}%
\e{0}%
\e{0}%
\e{0}%
\e{1}%
\e{0}%
\e{0}%
\e{0}%
\eol}\vss}\rg%
%
%
\rx{\vss\hfull{%
\rlx{\hss{$1344_{x}$}}\cg%
\e{1}%
\e{1}%
\e{1}%
\e{4}%
\e{1}%
\e{3}%
\e{1}%
\e{1}%
\e{1}%
\e{1}%
\e{0}%
\e{0}%
\e{0}%
\e{0}%
\e{0}%
\e{0}%
\e{0}%
\e{0}%
\e{0}%
\e{0}%
\e{0}%
\e{0}%
\e{0}%
\e{0}%
\e{0}%
\eol}\vss}\rg%
%
%
\rx{\vss\hfull{%
\rlx{\hss{$2100_{x}$}}\cg%
\e{0}%
\e{1}%
\e{3}%
\e{2}%
\e{1}%
\e{4}%
\e{3}%
\e{0}%
\e{1}%
\e{0}%
\e{0}%
\e{0}%
\e{0}%
\e{0}%
\e{0}%
\e{0}%
\e{1}%
\e{0}%
\e{0}%
\e{0}%
\e{2}%
\e{2}%
\e{0}%
\e{0}%
\e{0}%
\eol}\vss}\rg%
%
%
\rx{\vss\hfull{%
\rlx{\hss{$2268_{x}$}}\cg%
\e{0}%
\e{2}%
\e{2}%
\e{3}%
\e{3}%
\e{4}%
\e{3}%
\e{0}%
\e{1}%
\e{1}%
\e{0}%
\e{0}%
\e{0}%
\e{0}%
\e{0}%
\e{0}%
\e{0}%
\e{0}%
\e{0}%
\e{0}%
\e{0}%
\e{1}%
\e{1}%
\e{0}%
\e{0}%
\eol}\vss}\rg%
%
%
\rx{\vss\hfull{%
\rlx{\hss{$525_{x}$}}\cg%
\e{1}%
\e{0}%
\e{0}%
\e{2}%
\e{0}%
\e{1}%
\e{1}%
\e{0}%
\e{1}%
\e{0}%
\e{0}%
\e{0}%
\e{0}%
\e{0}%
\e{0}%
\e{0}%
\e{0}%
\e{0}%
\e{0}%
\e{0}%
\e{0}%
\e{0}%
\e{0}%
\e{0}%
\e{0}%
\eol}\vss}\rg%
%
%
\rx{\vss\hfull{%
\rlx{\hss{$700_{xx}$}}\cg%
\e{0}%
\e{0}%
\e{0}%
\e{0}%
\e{0}%
\e{0}%
\e{1}%
\e{2}%
\e{0}%
\e{1}%
\e{0}%
\e{0}%
\e{0}%
\e{0}%
\e{0}%
\e{0}%
\e{0}%
\e{0}%
\e{1}%
\e{0}%
\e{0}%
\e{0}%
\e{0}%
\e{2}%
\e{0}%
\eol}\vss}\rg%
%
%
\rx{\vss\hfull{%
\rlx{\hss{$972_{x}$}}\cg%
\e{0}%
\e{0}%
\e{0}%
\e{1}%
\e{0}%
\e{1}%
\e{2}%
\e{1}%
\e{1}%
\e{1}%
\e{0}%
\e{0}%
\e{0}%
\e{0}%
\e{0}%
\e{0}%
\e{0}%
\e{0}%
\e{0}%
\e{0}%
\e{0}%
\e{0}%
\e{0}%
\e{1}%
\e{0}%
\eol}\vss}\rg%
%
%
\rx{\vss\hfull{%
\rlx{\hss{$4096_{x}$}}\cg%
\e{0}%
\e{1}%
\e{2}%
\e{3}%
\e{3}%
\e{7}%
\e{4}%
\e{1}%
\e{2}%
\e{3}%
\e{0}%
\e{0}%
\e{0}%
\e{0}%
\e{0}%
\e{0}%
\e{1}%
\e{0}%
\e{0}%
\e{0}%
\e{1}%
\e{3}%
\e{2}%
\e{1}%
\e{0}%
\eol}\vss}\rg%
%
%
\rx{\vss\hfull{%
\rlx{\hss{$4200_{x}$}}\cg%
\e{0}%
\e{0}%
\e{1}%
\e{1}%
\e{3}%
\e{4}%
\e{4}%
\e{2}%
\e{1}%
\e{5}%
\e{0}%
\e{0}%
\e{0}%
\e{0}%
\e{0}%
\e{0}%
\e{1}%
\e{0}%
\e{0}%
\e{1}%
\e{0}%
\e{3}%
\e{3}%
\e{3}%
\e{0}%
\eol}\vss}\rg%
%
%
\rx{\vss\hfull{%
\rlx{\hss{$2240_{x}$}}\cg%
\e{0}%
\e{0}%
\e{0}%
\e{1}%
\e{2}%
\e{3}%
\e{1}%
\e{1}%
\e{1}%
\e{4}%
\e{0}%
\e{0}%
\e{0}%
\e{0}%
\e{0}%
\e{0}%
\e{0}%
\e{0}%
\e{0}%
\e{0}%
\e{0}%
\e{1}%
\e{2}%
\e{1}%
\e{1}%
\eol}\vss}\rg%
%
%
\rx{\vss\hfull{%
\rlx{\hss{$2835_{x}$}}\cg%
\e{0}%
\e{0}%
\e{0}%
\e{0}%
\e{1}%
\e{1}%
\e{3}%
\e{1}%
\e{0}%
\e{4}%
\e{0}%
\e{0}%
\e{0}%
\e{0}%
\e{0}%
\e{0}%
\e{1}%
\e{0}%
\e{0}%
\e{2}%
\e{0}%
\e{1}%
\e{3}%
\e{3}%
\e{2}%
\eol}\vss}\rg%
%
%
\rx{\vss\hfull{%
\rlx{\hss{$6075_{x}$}}\cg%
\e{0}%
\e{0}%
\e{1}%
\e{2}%
\e{3}%
\e{6}%
\e{7}%
\e{2}%
\e{2}%
\e{5}%
\e{0}%
\e{0}%
\e{0}%
\e{0}%
\e{0}%
\e{1}%
\e{3}%
\e{0}%
\e{1}%
\e{1}%
\e{1}%
\e{5}%
\e{4}%
\e{3}%
\e{0}%
\eol}\vss}\rg%
%
%
\rx{\vss\hfull{%
\rlx{\hss{$3200_{x}$}}\cg%
\e{0}%
\e{0}%
\e{0}%
\e{1}%
\e{0}%
\e{3}%
\e{3}%
\e{0}%
\e{3}%
\e{2}%
\e{0}%
\e{0}%
\e{0}%
\e{0}%
\e{0}%
\e{0}%
\e{3}%
\e{0}%
\e{0}%
\e{1}%
\e{1}%
\e{3}%
\e{2}%
\e{2}%
\e{0}%
\eol}\vss}\rg%
%
%
\rx{\vss\hfull{%
\rlx{\hss{$70_{y}$}}\cg%
\e{0}%
\e{0}%
\e{1}%
\e{0}%
\e{0}%
\e{0}%
\e{0}%
\e{0}%
\e{0}%
\e{0}%
\e{0}%
\e{0}%
\e{0}%
\e{0}%
\e{0}%
\e{0}%
\e{0}%
\e{0}%
\e{0}%
\e{0}%
\e{1}%
\e{0}%
\e{0}%
\e{0}%
\e{0}%
\eol}\vss}\rg%
%
%
\rx{\vss\hfull{%
\rlx{\hss{$1134_{y}$}}\cg%
\e{0}%
\e{0}%
\e{0}%
\e{0}%
\e{0}%
\e{1}%
\e{0}%
\e{0}%
\e{1}%
\e{1}%
\e{0}%
\e{0}%
\e{1}%
\e{0}%
\e{0}%
\e{1}%
\e{0}%
\e{0}%
\e{0}%
\e{0}%
\e{1}%
\e{2}%
\e{1}%
\e{1}%
\e{0}%
\eol}\vss}\rg%
%
%
\rx{\vss\hfull{%
\rlx{\hss{$1680_{y}$}}\cg%
\e{0}%
\e{0}%
\e{2}%
\e{0}%
\e{1}%
\e{2}%
\e{2}%
\e{0}%
\e{0}%
\e{0}%
\e{0}%
\e{0}%
\e{1}%
\e{0}%
\e{0}%
\e{1}%
\e{1}%
\e{0}%
\e{0}%
\e{0}%
\e{3}%
\e{3}%
\e{0}%
\e{0}%
\e{0}%
\eol}\vss}\rg%
%
%
\rx{\vss\hfull{%
\rlx{\hss{$168_{y}$}}\cg%
\e{0}%
\e{0}%
\e{0}%
\e{0}%
\e{0}%
\e{0}%
\e{0}%
\e{0}%
\e{0}%
\e{0}%
\e{0}%
\e{0}%
\e{0}%
\e{0}%
\e{0}%
\e{0}%
\e{0}%
\e{0}%
\e{1}%
\e{0}%
\e{0}%
\e{0}%
\e{0}%
\e{1}%
\e{0}%
\eol}\vss}\rg%
%
%
\rx{\vss\hfull{%
\rlx{\hss{$420_{y}$}}\cg%
\e{0}%
\e{0}%
\e{0}%
\e{0}%
\e{0}%
\e{0}%
\e{0}%
\e{0}%
\e{0}%
\e{0}%
\e{0}%
\e{0}%
\e{0}%
\e{0}%
\e{0}%
\e{0}%
\e{0}%
\e{0}%
\e{0}%
\e{1}%
\e{0}%
\e{0}%
\e{1}%
\e{1}%
\e{1}%
\eol}\vss}\rg%
%
%
\rx{\vss\hfull{%
\rlx{\hss{$3150_{y}$}}\cg%
\e{0}%
\e{0}%
\e{0}%
\e{0}%
\e{0}%
\e{1}%
\e{1}%
\e{0}%
\e{1}%
\e{3}%
\e{0}%
\e{0}%
\e{0}%
\e{0}%
\e{1}%
\e{1}%
\e{2}%
\e{0}%
\e{0}%
\e{3}%
\e{0}%
\e{3}%
\e{4}%
\e{3}%
\e{2}%
\eol}\vss}\rg%
%
%
\rx{\vss\hfull{%
\rlx{\hss{$4200_{y}$}}\cg%
\e{0}%
\e{0}%
\e{0}%
\e{0}%
\e{0}%
\e{1}%
\e{3}%
\e{1}%
\e{1}%
\e{3}%
\e{0}%
\e{0}%
\e{0}%
\e{0}%
\e{0}%
\e{2}%
\e{3}%
\e{1}%
\e{2}%
\e{3}%
\e{0}%
\e{3}%
\e{4}%
\e{6}%
\e{1}%
\eol}\vss}\rg%
\eop
\eject
\tablecont%
%
%
%
%
%
%
\rowpts=18 true pt%
\colpts=18 true pt%
\rowlabpts=40 true pt%
\collabpts=60 true pt%
\clx{\vss\hfull{%
\rlx{\hss{$ $}}\cg%
\cx{\hskip 16 true pt\flip{$1_p{\times}[{1^{2}}]$}\hss}\cg%
\cx{\hskip 16 true pt\flip{$6_p{\times}[{1^{2}}]$}\hss}\cg%
\cx{\hskip 16 true pt\flip{$15_p{\times}[{1^{2}}]$}\hss}\cg%
\cx{\hskip 16 true pt\flip{$20_p{\times}[{1^{2}}]$}\hss}\cg%
\cx{\hskip 16 true pt\flip{$30_p{\times}[{1^{2}}]$}\hss}\cg%
\cx{\hskip 16 true pt\flip{$64_p{\times}[{1^{2}}]$}\hss}\cg%
\cx{\hskip 16 true pt\flip{$81_p{\times}[{1^{2}}]$}\hss}\cg%
\cx{\hskip 16 true pt\flip{$15_q{\times}[{1^{2}}]$}\hss}\cg%
\cx{\hskip 16 true pt\flip{$24_p{\times}[{1^{2}}]$}\hss}\cg%
\cx{\hskip 16 true pt\flip{$60_p{\times}[{1^{2}}]$}\hss}\cg%
\cx{\hskip 16 true pt\flip{$1_p^{*}{\times}[{1^{2}}]$}\hss}\cg%
\cx{\hskip 16 true pt\flip{$6_p^{*}{\times}[{1^{2}}]$}\hss}\cg%
\cx{\hskip 16 true pt\flip{$15_p^{*}{\times}[{1^{2}}]$}\hss}\cg%
\cx{\hskip 16 true pt\flip{$20_p^{*}{\times}[{1^{2}}]$}\hss}\cg%
\cx{\hskip 16 true pt\flip{$30_p^{*}{\times}[{1^{2}}]$}\hss}\cg%
\cx{\hskip 16 true pt\flip{$64_p^{*}{\times}[{1^{2}}]$}\hss}\cg%
\cx{\hskip 16 true pt\flip{$81_p^{*}{\times}[{1^{2}}]$}\hss}\cg%
\cx{\hskip 16 true pt\flip{$15_q^{*}{\times}[{1^{2}}]$}\hss}\cg%
\cx{\hskip 16 true pt\flip{$24_p^{*}{\times}[{1^{2}}]$}\hss}\cg%
\cx{\hskip 16 true pt\flip{$60_p^{*}{\times}[{1^{2}}]$}\hss}\cg%
\cx{\hskip 16 true pt\flip{$20_s{\times}[{1^{2}}]$}\hss}\cg%
\cx{\hskip 16 true pt\flip{$90_s{\times}[{1^{2}}]$}\hss}\cg%
\cx{\hskip 16 true pt\flip{$80_s{\times}[{1^{2}}]$}\hss}\cg%
\cx{\hskip 16 true pt\flip{$60_s{\times}[{1^{2}}]$}\hss}\cg%
\cx{\hskip 16 true pt\flip{$10_s{\times}[{1^{2}}]$}\hss}\cg%
\eol}}\rg%
%
%
\rx{\vss\hfull{%
\rlx{\hss{$2688_{y}$}}\cg%
\e{0}%
\e{0}%
\e{0}%
\e{0}%
\e{0}%
\e{1}%
\e{3}%
\e{1}%
\e{0}%
\e{2}%
\e{0}%
\e{0}%
\e{0}%
\e{0}%
\e{0}%
\e{1}%
\e{3}%
\e{1}%
\e{0}%
\e{2}%
\e{1}%
\e{2}%
\e{1}%
\e{3}%
\e{0}%
\eol}\vss}\rg%
%
%
\rx{\vss\hfull{%
\rlx{\hss{$2100_{y}$}}\cg%
\e{0}%
\e{0}%
\e{0}%
\e{1}%
\e{0}%
\e{2}%
\e{3}%
\e{0}%
\e{2}%
\e{0}%
\e{0}%
\e{0}%
\e{0}%
\e{0}%
\e{0}%
\e{1}%
\e{3}%
\e{0}%
\e{1}%
\e{0}%
\e{1}%
\e{2}%
\e{1}%
\e{0}%
\e{0}%
\eol}\vss}\rg%
%
%
\rx{\vss\hfull{%
\rlx{\hss{$1400_{y}$}}\cg%
\e{0}%
\e{0}%
\e{1}%
\e{0}%
\e{1}%
\e{1}%
\e{2}%
\e{0}%
\e{0}%
\e{0}%
\e{0}%
\e{0}%
\e{0}%
\e{0}%
\e{0}%
\e{1}%
\e{1}%
\e{0}%
\e{1}%
\e{0}%
\e{0}%
\e{2}%
\e{1}%
\e{0}%
\e{0}%
\eol}\vss}\rg%
%
%
\rx{\vss\hfull{%
\rlx{\hss{$4536_{y}$}}\cg%
\e{0}%
\e{0}%
\e{1}%
\e{0}%
\e{2}%
\e{3}%
\e{5}%
\e{0}%
\e{0}%
\e{2}%
\e{0}%
\e{0}%
\e{0}%
\e{0}%
\e{1}%
\e{2}%
\e{4}%
\e{0}%
\e{1}%
\e{2}%
\e{1}%
\e{5}%
\e{4}%
\e{1}%
\e{0}%
\eol}\vss}\rg%
%
%
\rx{\vss\hfull{%
\rlx{\hss{$5670_{y}$}}\cg%
\e{0}%
\e{0}%
\e{1}%
\e{0}%
\e{2}%
\e{4}%
\e{5}%
\e{0}%
\e{1}%
\e{3}%
\e{0}%
\e{0}%
\e{1}%
\e{0}%
\e{1}%
\e{3}%
\e{4}%
\e{0}%
\e{1}%
\e{2}%
\e{2}%
\e{7}%
\e{5}%
\e{2}%
\e{0}%
\eol}\vss}\rg%
%
%
\rx{\vss\hfull{%
\rlx{\hss{$4480_{y}$}}\cg%
\e{0}%
\e{0}%
\e{0}%
\e{0}%
\e{1}%
\e{2}%
\e{3}%
\e{0}%
\e{1}%
\e{3}%
\e{0}%
\e{0}%
\e{0}%
\e{0}%
\e{1}%
\e{2}%
\e{3}%
\e{0}%
\e{1}%
\e{3}%
\e{0}%
\e{4}%
\e{6}%
\e{3}%
\e{1}%
\eol}\vss}\rg%
%
%
\rx{\vss\hfull{%
\rlx{\hss{$8_{z}$}}\cg%
\e{1}%
\e{0}%
\e{0}%
\e{0}%
\e{0}%
\e{0}%
\e{0}%
\e{0}%
\e{0}%
\e{0}%
\e{0}%
\e{0}%
\e{0}%
\e{0}%
\e{0}%
\e{0}%
\e{0}%
\e{0}%
\e{0}%
\e{0}%
\e{0}%
\e{0}%
\e{0}%
\e{0}%
\e{0}%
\eol}\vss}\rg%
%
%
\rx{\vss\hfull{%
\rlx{\hss{$56_{z}$}}\cg%
\e{0}%
\e{1}%
\e{1}%
\e{0}%
\e{0}%
\e{0}%
\e{0}%
\e{0}%
\e{0}%
\e{0}%
\e{0}%
\e{0}%
\e{0}%
\e{0}%
\e{0}%
\e{0}%
\e{0}%
\e{0}%
\e{0}%
\e{0}%
\e{0}%
\e{0}%
\e{0}%
\e{0}%
\e{0}%
\eol}\vss}\rg%
%
%
\rx{\vss\hfull{%
\rlx{\hss{$160_{z}$}}\cg%
\e{1}%
\e{2}%
\e{1}%
\e{1}%
\e{0}%
\e{0}%
\e{0}%
\e{0}%
\e{0}%
\e{0}%
\e{0}%
\e{0}%
\e{0}%
\e{0}%
\e{0}%
\e{0}%
\e{0}%
\e{0}%
\e{0}%
\e{0}%
\e{0}%
\e{0}%
\e{0}%
\e{0}%
\e{0}%
\eol}\vss}\rg%
%
%
\rx{\vss\hfull{%
\rlx{\hss{$112_{z}$}}\cg%
\e{2}%
\e{1}%
\e{0}%
\e{1}%
\e{0}%
\e{0}%
\e{0}%
\e{0}%
\e{0}%
\e{0}%
\e{0}%
\e{0}%
\e{0}%
\e{0}%
\e{0}%
\e{0}%
\e{0}%
\e{0}%
\e{0}%
\e{0}%
\e{0}%
\e{0}%
\e{0}%
\e{0}%
\e{0}%
\eol}\vss}\rg%
%
%
\rx{\vss\hfull{%
\rlx{\hss{$840_{z}$}}\cg%
\e{0}%
\e{0}%
\e{1}%
\e{1}%
\e{0}%
\e{2}%
\e{0}%
\e{0}%
\e{1}%
\e{1}%
\e{0}%
\e{0}%
\e{0}%
\e{0}%
\e{0}%
\e{0}%
\e{0}%
\e{0}%
\e{0}%
\e{0}%
\e{1}%
\e{1}%
\e{0}%
\e{0}%
\e{0}%
\eol}\vss}\rg%
%
%
\rx{\vss\hfull{%
\rlx{\hss{$1296_{z}$}}\cg%
\e{0}%
\e{2}%
\e{3}%
\e{2}%
\e{2}%
\e{3}%
\e{1}%
\e{0}%
\e{0}%
\e{0}%
\e{0}%
\e{0}%
\e{0}%
\e{0}%
\e{0}%
\e{0}%
\e{0}%
\e{0}%
\e{0}%
\e{0}%
\e{1}%
\e{1}%
\e{0}%
\e{0}%
\e{0}%
\eol}\vss}\rg%
%
%
\rx{\vss\hfull{%
\rlx{\hss{$1400_{z}$}}\cg%
\e{1}%
\e{2}%
\e{1}%
\e{4}%
\e{2}%
\e{3}%
\e{1}%
\e{0}%
\e{1}%
\e{1}%
\e{0}%
\e{0}%
\e{0}%
\e{0}%
\e{0}%
\e{0}%
\e{0}%
\e{0}%
\e{0}%
\e{0}%
\e{0}%
\e{0}%
\e{0}%
\e{0}%
\e{0}%
\eol}\vss}\rg%
%
%
\rx{\vss\hfull{%
\rlx{\hss{$1008_{z}$}}\cg%
\e{1}%
\e{3}%
\e{2}%
\e{3}%
\e{2}%
\e{2}%
\e{1}%
\e{0}%
\e{0}%
\e{0}%
\e{0}%
\e{0}%
\e{0}%
\e{0}%
\e{0}%
\e{0}%
\e{0}%
\e{0}%
\e{0}%
\e{0}%
\e{0}%
\e{0}%
\e{0}%
\e{0}%
\e{0}%
\eol}\vss}\rg%
%
%
\rx{\vss\hfull{%
\rlx{\hss{$560_{z}$}}\cg%
\e{1}%
\e{2}%
\e{0}%
\e{3}%
\e{1}%
\e{1}%
\e{0}%
\e{1}%
\e{0}%
\e{0}%
\e{0}%
\e{0}%
\e{0}%
\e{0}%
\e{0}%
\e{0}%
\e{0}%
\e{0}%
\e{0}%
\e{0}%
\e{0}%
\e{0}%
\e{0}%
\e{0}%
\e{0}%
\eol}\vss}\rg%
%
%
\rx{\vss\hfull{%
\rlx{\hss{$1400_{zz}$}}\cg%
\e{0}%
\e{0}%
\e{0}%
\e{0}%
\e{1}%
\e{1}%
\e{1}%
\e{2}%
\e{0}%
\e{3}%
\e{0}%
\e{0}%
\e{0}%
\e{0}%
\e{0}%
\e{0}%
\e{0}%
\e{0}%
\e{0}%
\e{0}%
\e{0}%
\e{0}%
\e{1}%
\e{2}%
\e{1}%
\eol}\vss}\rg%
%
%
\rx{\vss\hfull{%
\rlx{\hss{$4200_{z}$}}\cg%
\e{0}%
\e{0}%
\e{0}%
\e{0}%
\e{2}%
\e{2}%
\e{4}%
\e{2}%
\e{0}%
\e{5}%
\e{0}%
\e{0}%
\e{0}%
\e{0}%
\e{0}%
\e{1}%
\e{1}%
\e{0}%
\e{2}%
\e{2}%
\e{0}%
\e{3}%
\e{4}%
\e{4}%
\e{1}%
\eol}\vss}\rg%
%
%
\rx{\vss\hfull{%
\rlx{\hss{$400_{z}$}}\cg%
\e{0}%
\e{0}%
\e{0}%
\e{1}%
\e{1}%
\e{0}%
\e{0}%
\e{2}%
\e{0}%
\e{1}%
\e{0}%
\e{0}%
\e{0}%
\e{0}%
\e{0}%
\e{0}%
\e{0}%
\e{0}%
\e{0}%
\e{0}%
\e{0}%
\e{0}%
\e{0}%
\e{0}%
\e{0}%
\eol}\vss}\rg%
%
%
\rx{\vss\hfull{%
\rlx{\hss{$3240_{z}$}}\cg%
\e{0}%
\e{1}%
\e{1}%
\e{4}%
\e{4}%
\e{5}%
\e{3}%
\e{3}%
\e{1}%
\e{4}%
\e{0}%
\e{0}%
\e{0}%
\e{0}%
\e{0}%
\e{0}%
\e{0}%
\e{0}%
\e{0}%
\e{0}%
\e{0}%
\e{1}%
\e{1}%
\e{1}%
\e{0}%
\eol}\vss}\rg%
%
%
\rx{\vss\hfull{%
\rlx{\hss{$4536_{z}$}}\cg%
\e{0}%
\e{0}%
\e{0}%
\e{1}%
\e{1}%
\e{4}%
\e{4}%
\e{1}%
\e{3}%
\e{5}%
\e{0}%
\e{0}%
\e{0}%
\e{0}%
\e{0}%
\e{0}%
\e{2}%
\e{0}%
\e{0}%
\e{1}%
\e{0}%
\e{3}%
\e{4}%
\e{4}%
\e{1}%
\eol}\vss}\rg%
%
%
\rx{\vss\hfull{%
\rlx{\hss{$2400_{z}$}}\cg%
\e{0}%
\e{0}%
\e{2}%
\e{1}%
\e{2}%
\e{3}%
\e{3}%
\e{0}%
\e{0}%
\e{1}%
\e{0}%
\e{0}%
\e{0}%
\e{0}%
\e{0}%
\e{1}%
\e{1}%
\e{0}%
\e{0}%
\e{0}%
\e{2}%
\e{3}%
\e{1}%
\e{0}%
\e{0}%
\eol}\vss}\rg%
%
%
\rx{\vss\hfull{%
\rlx{\hss{$3360_{z}$}}\cg%
\e{0}%
\e{0}%
\e{1}%
\e{1}%
\e{3}%
\e{3}%
\e{4}%
\e{2}%
\e{0}%
\e{4}%
\e{0}%
\e{0}%
\e{0}%
\e{0}%
\e{0}%
\e{0}%
\e{1}%
\e{0}%
\e{0}%
\e{1}%
\e{0}%
\e{2}%
\e{2}%
\e{2}%
\e{0}%
\eol}\vss}\rg%
%
%
\rx{\vss\hfull{%
\rlx{\hss{$2800_{z}$}}\cg%
\e{0}%
\e{1}%
\e{1}%
\e{2}%
\e{3}%
\e{4}%
\e{4}%
\e{0}%
\e{1}%
\e{1}%
\e{0}%
\e{0}%
\e{0}%
\e{0}%
\e{0}%
\e{0}%
\e{1}%
\e{0}%
\e{1}%
\e{0}%
\e{0}%
\e{2}%
\e{2}%
\e{0}%
\e{0}%
\eol}\vss}\rg%
%
%
\rx{\vss\hfull{%
\rlx{\hss{$4096_{z}$}}\cg%
\e{0}%
\e{1}%
\e{2}%
\e{3}%
\e{3}%
\e{7}%
\e{4}%
\e{1}%
\e{2}%
\e{3}%
\e{0}%
\e{0}%
\e{0}%
\e{0}%
\e{0}%
\e{0}%
\e{1}%
\e{0}%
\e{0}%
\e{0}%
\e{1}%
\e{3}%
\e{2}%
\e{1}%
\e{0}%
\eol}\vss}\rg%
%
%
\rx{\vss\hfull{%
\rlx{\hss{$5600_{z}$}}\cg%
\e{0}%
\e{0}%
\e{2}%
\e{2}%
\e{2}%
\e{6}%
\e{6}%
\e{0}%
\e{3}%
\e{3}%
\e{0}%
\e{0}%
\e{0}%
\e{0}%
\e{0}%
\e{1}%
\e{4}%
\e{0}%
\e{0}%
\e{1}%
\e{2}%
\e{6}%
\e{4}%
\e{1}%
\e{0}%
\eol}\vss}\rg%
%
%
\rx{\vss\hfull{%
\rlx{\hss{$448_{z}$}}\cg%
\e{0}%
\e{0}%
\e{0}%
\e{1}%
\e{0}%
\e{1}%
\e{0}%
\e{0}%
\e{1}%
\e{1}%
\e{0}%
\e{0}%
\e{0}%
\e{0}%
\e{0}%
\e{0}%
\e{0}%
\e{0}%
\e{0}%
\e{0}%
\e{0}%
\e{0}%
\e{0}%
\e{0}%
\e{0}%
\eol}\vss}\rg%
%
%
\rx{\vss\hfull{%
\rlx{\hss{$448_{w}$}}\cg%
\e{0}%
\e{0}%
\e{1}%
\e{0}%
\e{0}%
\e{1}%
\e{0}%
\e{0}%
\e{0}%
\e{0}%
\e{0}%
\e{0}%
\e{1}%
\e{0}%
\e{0}%
\e{0}%
\e{0}%
\e{0}%
\e{0}%
\e{0}%
\e{2}%
\e{1}%
\e{0}%
\e{0}%
\e{0}%
\eol}\vss}\rg%
%
%
\rx{\vss\hfull{%
\rlx{\hss{$1344_{w}$}}\cg%
\e{0}%
\e{0}%
\e{0}%
\e{0}%
\e{0}%
\e{0}%
\e{0}%
\e{0}%
\e{1}%
\e{1}%
\e{0}%
\e{0}%
\e{0}%
\e{0}%
\e{0}%
\e{1}%
\e{0}%
\e{0}%
\e{1}%
\e{1}%
\e{0}%
\e{1}%
\e{2}%
\e{3}%
\e{1}%
\eol}\vss}\rg%
%
%
\rx{\vss\hfull{%
\rlx{\hss{$5600_{w}$}}\cg%
\e{0}%
\e{0}%
\e{1}%
\e{0}%
\e{1}%
\e{4}%
\e{5}%
\e{1}%
\e{1}%
\e{3}%
\e{0}%
\e{0}%
\e{1}%
\e{0}%
\e{1}%
\e{3}%
\e{4}%
\e{0}%
\e{1}%
\e{2}%
\e{2}%
\e{7}%
\e{4}%
\e{3}%
\e{0}%
\eol}\vss}\rg%
%
%
\rx{\vss\hfull{%
\rlx{\hss{$2016_{w}$}}\cg%
\e{0}%
\e{0}%
\e{0}%
\e{0}%
\e{0}%
\e{0}%
\e{1}%
\e{0}%
\e{0}%
\e{2}%
\e{0}%
\e{0}%
\e{0}%
\e{0}%
\e{0}%
\e{0}%
\e{2}%
\e{1}%
\e{0}%
\e{3}%
\e{0}%
\e{1}%
\e{2}%
\e{3}%
\e{2}%
\eol}\vss}\rg%
%
%
\rx{\vss\hfull{%
\rlx{\hss{$7168_{w}$}}\cg%
\e{0}%
\e{0}%
\e{0}%
\e{0}%
\e{1}%
\e{3}%
\e{6}%
\e{1}%
\e{1}%
\e{5}%
\e{0}%
\e{0}%
\e{0}%
\e{0}%
\e{1}%
\e{3}%
\e{6}%
\e{1}%
\e{1}%
\e{5}%
\e{1}%
\e{6}%
\e{7}%
\e{6}%
\e{1}%
\eol}\vss}\rg%
\tableclose%
%
%
%
%
%
%
\eop
\eject
\tableopen{Induce/restrict matrix for $W(E_{7})\,\subset\,W(E_{8})$}%
%
%
%
%
%
%
\rowpts=18 true pt%
\colpts=18 true pt%
\rowlabpts=40 true pt%
\collabpts=45 true pt%
\clx{\vss\hfull{%
\rlx{\hss{$ $}}\cg%
\cx{\hskip 16 true pt\flip{$1_a$}\hss}\cg%
\cx{\hskip 16 true pt\flip{$7_a$}\hss}\cg%
\cx{\hskip 16 true pt\flip{$27_a$}\hss}\cg%
\cx{\hskip 16 true pt\flip{$21_a$}\hss}\cg%
\cx{\hskip 16 true pt\flip{$35_a$}\hss}\cg%
\cx{\hskip 16 true pt\flip{$105_a$}\hss}\cg%
\cx{\hskip 16 true pt\flip{$189_a$}\hss}\cg%
\cx{\hskip 16 true pt\flip{$21_b$}\hss}\cg%
\cx{\hskip 16 true pt\flip{$35_b$}\hss}\cg%
\cx{\hskip 16 true pt\flip{$189_b$}\hss}\cg%
\cx{\hskip 16 true pt\flip{$189_c$}\hss}\cg%
\cx{\hskip 16 true pt\flip{$15_a$}\hss}\cg%
\cx{\hskip 16 true pt\flip{$105_b$}\hss}\cg%
\cx{\hskip 16 true pt\flip{$105_c$}\hss}\cg%
\cx{\hskip 16 true pt\flip{$315_a$}\hss}\cg%
\cx{\hskip 16 true pt\flip{$405_a$}\hss}\cg%
\cx{\hskip 16 true pt\flip{$168_a$}\hss}\cg%
\cx{\hskip 16 true pt\flip{$56_a$}\hss}\cg%
\cx{\hskip 16 true pt\flip{$120_a$}\hss}\cg%
\cx{\hskip 16 true pt\flip{$210_a$}\hss}\cg%
\eol}}\rg%
%
%
\rx{\vss\hfull{%
\rlx{\hss{$1_{x}$}}\cg%
\e{1}%
\e{0}%
\e{0}%
\e{0}%
\e{0}%
\e{0}%
\e{0}%
\e{0}%
\e{0}%
\e{0}%
\e{0}%
\e{0}%
\e{0}%
\e{0}%
\e{0}%
\e{0}%
\e{0}%
\e{0}%
\e{0}%
\e{0}%
\eol}\vss}\rg%
%
%
\rx{\vss\hfull{%
\rlx{\hss{$28_{x}$}}\cg%
\e{0}%
\e{0}%
\e{0}%
\e{1}%
\e{0}%
\e{0}%
\e{0}%
\e{0}%
\e{0}%
\e{0}%
\e{0}%
\e{0}%
\e{0}%
\e{0}%
\e{0}%
\e{0}%
\e{0}%
\e{0}%
\e{0}%
\e{0}%
\eol}\vss}\rg%
%
%
\rx{\vss\hfull{%
\rlx{\hss{$35_{x}$}}\cg%
\e{1}%
\e{0}%
\e{1}%
\e{0}%
\e{0}%
\e{0}%
\e{0}%
\e{0}%
\e{0}%
\e{0}%
\e{0}%
\e{0}%
\e{0}%
\e{0}%
\e{0}%
\e{0}%
\e{0}%
\e{0}%
\e{0}%
\e{0}%
\eol}\vss}\rg%
%
%
\rx{\vss\hfull{%
\rlx{\hss{$84_{x}$}}\cg%
\e{1}%
\e{0}%
\e{1}%
\e{0}%
\e{0}%
\e{0}%
\e{0}%
\e{0}%
\e{1}%
\e{0}%
\e{0}%
\e{0}%
\e{0}%
\e{0}%
\e{0}%
\e{0}%
\e{0}%
\e{0}%
\e{0}%
\e{0}%
\eol}\vss}\rg%
%
%
\rx{\vss\hfull{%
\rlx{\hss{$50_{x}$}}\cg%
\e{0}%
\e{0}%
\e{0}%
\e{0}%
\e{0}%
\e{0}%
\e{0}%
\e{0}%
\e{1}%
\e{0}%
\e{0}%
\e{0}%
\e{0}%
\e{0}%
\e{0}%
\e{0}%
\e{0}%
\e{0}%
\e{0}%
\e{0}%
\eol}\vss}\rg%
%
%
\rx{\vss\hfull{%
\rlx{\hss{$350_{x}$}}\cg%
\e{0}%
\e{0}%
\e{0}%
\e{1}%
\e{0}%
\e{0}%
\e{1}%
\e{0}%
\e{0}%
\e{0}%
\e{0}%
\e{0}%
\e{0}%
\e{0}%
\e{0}%
\e{0}%
\e{0}%
\e{0}%
\e{0}%
\e{0}%
\eol}\vss}\rg%
%
%
\rx{\vss\hfull{%
\rlx{\hss{$300_{x}$}}\cg%
\e{0}%
\e{0}%
\e{1}%
\e{0}%
\e{0}%
\e{0}%
\e{0}%
\e{0}%
\e{0}%
\e{0}%
\e{0}%
\e{0}%
\e{0}%
\e{0}%
\e{0}%
\e{0}%
\e{1}%
\e{0}%
\e{0}%
\e{0}%
\eol}\vss}\rg%
%
%
\rx{\vss\hfull{%
\rlx{\hss{$567_{x}$}}\cg%
\e{0}%
\e{0}%
\e{1}%
\e{1}%
\e{0}%
\e{0}%
\e{0}%
\e{0}%
\e{0}%
\e{0}%
\e{0}%
\e{0}%
\e{0}%
\e{0}%
\e{0}%
\e{0}%
\e{0}%
\e{0}%
\e{1}%
\e{1}%
\eol}\vss}\rg%
%
%
\rx{\vss\hfull{%
\rlx{\hss{$210_{x}$}}\cg%
\e{0}%
\e{0}%
\e{1}%
\e{0}%
\e{0}%
\e{0}%
\e{0}%
\e{0}%
\e{0}%
\e{0}%
\e{0}%
\e{0}%
\e{0}%
\e{0}%
\e{0}%
\e{0}%
\e{0}%
\e{0}%
\e{1}%
\e{0}%
\eol}\vss}\rg%
%
%
\rx{\vss\hfull{%
\rlx{\hss{$840_{x}$}}\cg%
\e{0}%
\e{0}%
\e{0}%
\e{0}%
\e{0}%
\e{0}%
\e{0}%
\e{0}%
\e{0}%
\e{0}%
\e{0}%
\e{0}%
\e{0}%
\e{0}%
\e{0}%
\e{0}%
\e{1}%
\e{0}%
\e{0}%
\e{0}%
\eol}\vss}\rg%
%
%
\rx{\vss\hfull{%
\rlx{\hss{$700_{x}$}}\cg%
\e{0}%
\e{0}%
\e{1}%
\e{0}%
\e{0}%
\e{0}%
\e{0}%
\e{0}%
\e{1}%
\e{0}%
\e{0}%
\e{0}%
\e{1}%
\e{0}%
\e{0}%
\e{0}%
\e{1}%
\e{0}%
\e{1}%
\e{0}%
\eol}\vss}\rg%
%
%
\rx{\vss\hfull{%
\rlx{\hss{$175_{x}$}}\cg%
\e{0}%
\e{0}%
\e{0}%
\e{0}%
\e{0}%
\e{0}%
\e{0}%
\e{0}%
\e{0}%
\e{0}%
\e{0}%
\e{0}%
\e{1}%
\e{0}%
\e{0}%
\e{0}%
\e{0}%
\e{0}%
\e{0}%
\e{0}%
\eol}\vss}\rg%
%
%
\rx{\vss\hfull{%
\rlx{\hss{$1400_{x}$}}\cg%
\e{0}%
\e{0}%
\e{0}%
\e{0}%
\e{0}%
\e{0}%
\e{0}%
\e{0}%
\e{0}%
\e{0}%
\e{0}%
\e{0}%
\e{1}%
\e{0}%
\e{0}%
\e{1}%
\e{0}%
\e{0}%
\e{1}%
\e{1}%
\eol}\vss}\rg%
%
%
\rx{\vss\hfull{%
\rlx{\hss{$1050_{x}$}}\cg%
\e{0}%
\e{0}%
\e{0}%
\e{0}%
\e{0}%
\e{0}%
\e{0}%
\e{0}%
\e{1}%
\e{0}%
\e{0}%
\e{0}%
\e{1}%
\e{0}%
\e{0}%
\e{0}%
\e{0}%
\e{0}%
\e{0}%
\e{1}%
\eol}\vss}\rg%
%
%
\rx{\vss\hfull{%
\rlx{\hss{$1575_{x}$}}\cg%
\e{0}%
\e{0}%
\e{0}%
\e{1}%
\e{0}%
\e{0}%
\e{1}%
\e{0}%
\e{0}%
\e{0}%
\e{0}%
\e{0}%
\e{0}%
\e{0}%
\e{0}%
\e{1}%
\e{0}%
\e{0}%
\e{1}%
\e{1}%
\eol}\vss}\rg%
%
%
\rx{\vss\hfull{%
\rlx{\hss{$1344_{x}$}}\cg%
\e{0}%
\e{0}%
\e{1}%
\e{0}%
\e{0}%
\e{0}%
\e{0}%
\e{0}%
\e{1}%
\e{0}%
\e{0}%
\e{0}%
\e{0}%
\e{0}%
\e{0}%
\e{0}%
\e{1}%
\e{0}%
\e{1}%
\e{1}%
\eol}\vss}\rg%
%
%
\rx{\vss\hfull{%
\rlx{\hss{$2100_{x}$}}\cg%
\e{0}%
\e{0}%
\e{0}%
\e{0}%
\e{0}%
\e{0}%
\e{1}%
\e{0}%
\e{0}%
\e{0}%
\e{0}%
\e{0}%
\e{0}%
\e{0}%
\e{0}%
\e{0}%
\e{0}%
\e{0}%
\e{0}%
\e{1}%
\eol}\vss}\rg%
%
%
\rx{\vss\hfull{%
\rlx{\hss{$2268_{x}$}}\cg%
\e{0}%
\e{0}%
\e{0}%
\e{0}%
\e{0}%
\e{0}%
\e{0}%
\e{0}%
\e{0}%
\e{0}%
\e{0}%
\e{0}%
\e{0}%
\e{0}%
\e{0}%
\e{1}%
\e{1}%
\e{0}%
\e{1}%
\e{1}%
\eol}\vss}\rg%
%
%
\rx{\vss\hfull{%
\rlx{\hss{$525_{x}$}}\cg%
\e{0}%
\e{0}%
\e{0}%
\e{0}%
\e{0}%
\e{0}%
\e{0}%
\e{0}%
\e{0}%
\e{0}%
\e{0}%
\e{0}%
\e{0}%
\e{1}%
\e{0}%
\e{0}%
\e{0}%
\e{0}%
\e{0}%
\e{1}%
\eol}\vss}\rg%
%
%
\rx{\vss\hfull{%
\rlx{\hss{$700_{xx}$}}\cg%
\e{0}%
\e{0}%
\e{0}%
\e{0}%
\e{0}%
\e{0}%
\e{0}%
\e{0}%
\e{0}%
\e{0}%
\e{0}%
\e{0}%
\e{0}%
\e{1}%
\e{0}%
\e{0}%
\e{0}%
\e{0}%
\e{0}%
\e{0}%
\eol}\vss}\rg%
%
%
\rx{\vss\hfull{%
\rlx{\hss{$972_{x}$}}\cg%
\e{0}%
\e{0}%
\e{0}%
\e{0}%
\e{0}%
\e{0}%
\e{0}%
\e{0}%
\e{1}%
\e{0}%
\e{0}%
\e{0}%
\e{0}%
\e{0}%
\e{0}%
\e{0}%
\e{1}%
\e{0}%
\e{0}%
\e{0}%
\eol}\vss}\rg%
%
%
\rx{\vss\hfull{%
\rlx{\hss{$4096_{x}$}}\cg%
\e{0}%
\e{0}%
\e{0}%
\e{0}%
\e{0}%
\e{0}%
\e{1}%
\e{0}%
\e{0}%
\e{0}%
\e{0}%
\e{0}%
\e{0}%
\e{0}%
\e{0}%
\e{1}%
\e{1}%
\e{0}%
\e{1}%
\e{1}%
\eol}\vss}\rg%
%
%
\rx{\vss\hfull{%
\rlx{\hss{$4200_{x}$}}\cg%
\e{0}%
\e{0}%
\e{0}%
\e{0}%
\e{0}%
\e{0}%
\e{0}%
\e{0}%
\e{0}%
\e{0}%
\e{0}%
\e{0}%
\e{1}%
\e{0}%
\e{0}%
\e{1}%
\e{1}%
\e{0}%
\e{0}%
\e{1}%
\eol}\vss}\rg%
%
%
\rx{\vss\hfull{%
\rlx{\hss{$2240_{x}$}}\cg%
\e{0}%
\e{0}%
\e{0}%
\e{0}%
\e{0}%
\e{0}%
\e{0}%
\e{0}%
\e{0}%
\e{0}%
\e{0}%
\e{0}%
\e{1}%
\e{0}%
\e{0}%
\e{1}%
\e{1}%
\e{0}%
\e{1}%
\e{0}%
\eol}\vss}\rg%
%
%
\rx{\vss\hfull{%
\rlx{\hss{$2835_{x}$}}\cg%
\e{0}%
\e{0}%
\e{0}%
\e{0}%
\e{0}%
\e{0}%
\e{0}%
\e{0}%
\e{0}%
\e{0}%
\e{0}%
\e{0}%
\e{1}%
\e{0}%
\e{0}%
\e{1}%
\e{0}%
\e{0}%
\e{0}%
\e{0}%
\eol}\vss}\rg%
%
%
\rx{\vss\hfull{%
\rlx{\hss{$6075_{x}$}}\cg%
\e{0}%
\e{0}%
\e{0}%
\e{0}%
\e{0}%
\e{0}%
\e{1}%
\e{0}%
\e{0}%
\e{0}%
\e{0}%
\e{0}%
\e{0}%
\e{1}%
\e{0}%
\e{2}%
\e{0}%
\e{0}%
\e{0}%
\e{1}%
\eol}\vss}\rg%
%
%
\rx{\vss\hfull{%
\rlx{\hss{$3200_{x}$}}\cg%
\e{0}%
\e{0}%
\e{0}%
\e{0}%
\e{0}%
\e{0}%
\e{0}%
\e{0}%
\e{0}%
\e{0}%
\e{0}%
\e{0}%
\e{0}%
\e{0}%
\e{0}%
\e{0}%
\e{1}%
\e{0}%
\e{0}%
\e{0}%
\eol}\vss}\rg%
%
%
\rx{\vss\hfull{%
\rlx{\hss{$70_{y}$}}\cg%
\e{0}%
\e{0}%
\e{0}%
\e{0}%
\e{1}%
\e{0}%
\e{0}%
\e{0}%
\e{0}%
\e{0}%
\e{0}%
\e{0}%
\e{0}%
\e{0}%
\e{0}%
\e{0}%
\e{0}%
\e{0}%
\e{0}%
\e{0}%
\eol}\vss}\rg%
%
%
\rx{\vss\hfull{%
\rlx{\hss{$1134_{y}$}}\cg%
\e{0}%
\e{0}%
\e{0}%
\e{0}%
\e{0}%
\e{0}%
\e{1}%
\e{0}%
\e{0}%
\e{0}%
\e{0}%
\e{0}%
\e{0}%
\e{0}%
\e{0}%
\e{0}%
\e{0}%
\e{0}%
\e{0}%
\e{0}%
\eol}\vss}\rg%
%
%
\rx{\vss\hfull{%
\rlx{\hss{$1680_{y}$}}\cg%
\e{0}%
\e{0}%
\e{0}%
\e{0}%
\e{1}%
\e{0}%
\e{1}%
\e{0}%
\e{0}%
\e{0}%
\e{0}%
\e{0}%
\e{0}%
\e{0}%
\e{0}%
\e{0}%
\e{0}%
\e{0}%
\e{0}%
\e{0}%
\eol}\vss}\rg%
%
%
\rx{\vss\hfull{%
\rlx{\hss{$168_{y}$}}\cg%
\e{0}%
\e{0}%
\e{0}%
\e{0}%
\e{0}%
\e{0}%
\e{0}%
\e{0}%
\e{0}%
\e{0}%
\e{0}%
\e{0}%
\e{0}%
\e{0}%
\e{0}%
\e{0}%
\e{0}%
\e{0}%
\e{0}%
\e{0}%
\eol}\vss}\rg%
%
%
\rx{\vss\hfull{%
\rlx{\hss{$420_{y}$}}\cg%
\e{0}%
\e{0}%
\e{0}%
\e{0}%
\e{0}%
\e{0}%
\e{0}%
\e{0}%
\e{0}%
\e{0}%
\e{0}%
\e{0}%
\e{0}%
\e{0}%
\e{0}%
\e{0}%
\e{0}%
\e{0}%
\e{0}%
\e{0}%
\eol}\vss}\rg%
%
%
\rx{\vss\hfull{%
\rlx{\hss{$3150_{y}$}}\cg%
\e{0}%
\e{0}%
\e{0}%
\e{0}%
\e{0}%
\e{0}%
\e{0}%
\e{0}%
\e{0}%
\e{0}%
\e{0}%
\e{0}%
\e{0}%
\e{0}%
\e{0}%
\e{1}%
\e{0}%
\e{0}%
\e{0}%
\e{0}%
\eol}\vss}\rg%
%
%
\rx{\vss\hfull{%
\rlx{\hss{$4200_{y}$}}\cg%
\e{0}%
\e{0}%
\e{0}%
\e{0}%
\e{0}%
\e{0}%
\e{0}%
\e{0}%
\e{0}%
\e{0}%
\e{0}%
\e{0}%
\e{0}%
\e{0}%
\e{0}%
\e{0}%
\e{0}%
\e{0}%
\e{0}%
\e{0}%
\eol}\vss}\rg%
\eop
\eject
\tablecont%
%
%
%
%
%
%
\rowpts=18 true pt%
\colpts=18 true pt%
\rowlabpts=40 true pt%
\collabpts=45 true pt%
\clx{\vss\hfull{%
\rlx{\hss{$ $}}\cg%
\cx{\hskip 16 true pt\flip{$1_a$}\hss}\cg%
\cx{\hskip 16 true pt\flip{$7_a$}\hss}\cg%
\cx{\hskip 16 true pt\flip{$27_a$}\hss}\cg%
\cx{\hskip 16 true pt\flip{$21_a$}\hss}\cg%
\cx{\hskip 16 true pt\flip{$35_a$}\hss}\cg%
\cx{\hskip 16 true pt\flip{$105_a$}\hss}\cg%
\cx{\hskip 16 true pt\flip{$189_a$}\hss}\cg%
\cx{\hskip 16 true pt\flip{$21_b$}\hss}\cg%
\cx{\hskip 16 true pt\flip{$35_b$}\hss}\cg%
\cx{\hskip 16 true pt\flip{$189_b$}\hss}\cg%
\cx{\hskip 16 true pt\flip{$189_c$}\hss}\cg%
\cx{\hskip 16 true pt\flip{$15_a$}\hss}\cg%
\cx{\hskip 16 true pt\flip{$105_b$}\hss}\cg%
\cx{\hskip 16 true pt\flip{$105_c$}\hss}\cg%
\cx{\hskip 16 true pt\flip{$315_a$}\hss}\cg%
\cx{\hskip 16 true pt\flip{$405_a$}\hss}\cg%
\cx{\hskip 16 true pt\flip{$168_a$}\hss}\cg%
\cx{\hskip 16 true pt\flip{$56_a$}\hss}\cg%
\cx{\hskip 16 true pt\flip{$120_a$}\hss}\cg%
\cx{\hskip 16 true pt\flip{$210_a$}\hss}\cg%
\eol}}\rg%
%
%
\rx{\vss\hfull{%
\rlx{\hss{$2688_{y}$}}\cg%
\e{0}%
\e{0}%
\e{0}%
\e{0}%
\e{0}%
\e{0}%
\e{0}%
\e{0}%
\e{0}%
\e{0}%
\e{0}%
\e{0}%
\e{0}%
\e{0}%
\e{0}%
\e{0}%
\e{0}%
\e{0}%
\e{0}%
\e{0}%
\eol}\vss}\rg%
%
%
\rx{\vss\hfull{%
\rlx{\hss{$2100_{y}$}}\cg%
\e{0}%
\e{0}%
\e{0}%
\e{0}%
\e{0}%
\e{0}%
\e{0}%
\e{0}%
\e{0}%
\e{0}%
\e{1}%
\e{0}%
\e{0}%
\e{1}%
\e{0}%
\e{0}%
\e{0}%
\e{0}%
\e{0}%
\e{0}%
\eol}\vss}\rg%
%
%
\rx{\vss\hfull{%
\rlx{\hss{$1400_{y}$}}\cg%
\e{0}%
\e{0}%
\e{0}%
\e{0}%
\e{0}%
\e{0}%
\e{0}%
\e{0}%
\e{0}%
\e{0}%
\e{0}%
\e{0}%
\e{0}%
\e{0}%
\e{0}%
\e{0}%
\e{0}%
\e{0}%
\e{0}%
\e{0}%
\eol}\vss}\rg%
%
%
\rx{\vss\hfull{%
\rlx{\hss{$4536_{y}$}}\cg%
\e{0}%
\e{0}%
\e{0}%
\e{0}%
\e{0}%
\e{0}%
\e{0}%
\e{0}%
\e{0}%
\e{0}%
\e{0}%
\e{0}%
\e{0}%
\e{0}%
\e{1}%
\e{1}%
\e{0}%
\e{0}%
\e{0}%
\e{0}%
\eol}\vss}\rg%
%
%
\rx{\vss\hfull{%
\rlx{\hss{$5670_{y}$}}\cg%
\e{0}%
\e{0}%
\e{0}%
\e{0}%
\e{0}%
\e{0}%
\e{1}%
\e{0}%
\e{0}%
\e{0}%
\e{0}%
\e{0}%
\e{0}%
\e{0}%
\e{1}%
\e{1}%
\e{0}%
\e{0}%
\e{0}%
\e{0}%
\eol}\vss}\rg%
%
%
\rx{\vss\hfull{%
\rlx{\hss{$4480_{y}$}}\cg%
\e{0}%
\e{0}%
\e{0}%
\e{0}%
\e{0}%
\e{0}%
\e{0}%
\e{0}%
\e{0}%
\e{0}%
\e{0}%
\e{0}%
\e{0}%
\e{0}%
\e{1}%
\e{1}%
\e{0}%
\e{0}%
\e{0}%
\e{0}%
\eol}\vss}\rg%
%
%
\rx{\vss\hfull{%
\rlx{\hss{$8_{z}$}}\cg%
\e{1}%
\e{0}%
\e{0}%
\e{0}%
\e{0}%
\e{0}%
\e{0}%
\e{0}%
\e{0}%
\e{0}%
\e{0}%
\e{0}%
\e{0}%
\e{0}%
\e{0}%
\e{0}%
\e{0}%
\e{0}%
\e{0}%
\e{0}%
\eol}\vss}\rg%
%
%
\rx{\vss\hfull{%
\rlx{\hss{$56_{z}$}}\cg%
\e{0}%
\e{0}%
\e{0}%
\e{1}%
\e{0}%
\e{0}%
\e{0}%
\e{0}%
\e{0}%
\e{0}%
\e{0}%
\e{0}%
\e{0}%
\e{0}%
\e{0}%
\e{0}%
\e{0}%
\e{0}%
\e{0}%
\e{0}%
\eol}\vss}\rg%
%
%
\rx{\vss\hfull{%
\rlx{\hss{$160_{z}$}}\cg%
\e{0}%
\e{0}%
\e{1}%
\e{1}%
\e{0}%
\e{0}%
\e{0}%
\e{0}%
\e{0}%
\e{0}%
\e{0}%
\e{0}%
\e{0}%
\e{0}%
\e{0}%
\e{0}%
\e{0}%
\e{0}%
\e{0}%
\e{0}%
\eol}\vss}\rg%
%
%
\rx{\vss\hfull{%
\rlx{\hss{$112_{z}$}}\cg%
\e{1}%
\e{0}%
\e{1}%
\e{0}%
\e{0}%
\e{0}%
\e{0}%
\e{0}%
\e{0}%
\e{0}%
\e{0}%
\e{0}%
\e{0}%
\e{0}%
\e{0}%
\e{0}%
\e{0}%
\e{0}%
\e{0}%
\e{0}%
\eol}\vss}\rg%
%
%
\rx{\vss\hfull{%
\rlx{\hss{$840_{z}$}}\cg%
\e{0}%
\e{0}%
\e{0}%
\e{0}%
\e{0}%
\e{0}%
\e{1}%
\e{0}%
\e{0}%
\e{0}%
\e{0}%
\e{0}%
\e{0}%
\e{0}%
\e{0}%
\e{0}%
\e{1}%
\e{0}%
\e{0}%
\e{0}%
\eol}\vss}\rg%
%
%
\rx{\vss\hfull{%
\rlx{\hss{$1296_{z}$}}\cg%
\e{0}%
\e{0}%
\e{0}%
\e{1}%
\e{0}%
\e{0}%
\e{1}%
\e{0}%
\e{0}%
\e{0}%
\e{0}%
\e{0}%
\e{0}%
\e{0}%
\e{0}%
\e{0}%
\e{0}%
\e{0}%
\e{1}%
\e{1}%
\eol}\vss}\rg%
%
%
\rx{\vss\hfull{%
\rlx{\hss{$1400_{z}$}}\cg%
\e{0}%
\e{0}%
\e{1}%
\e{0}%
\e{0}%
\e{0}%
\e{0}%
\e{0}%
\e{0}%
\e{0}%
\e{0}%
\e{0}%
\e{0}%
\e{0}%
\e{0}%
\e{0}%
\e{1}%
\e{0}%
\e{1}%
\e{1}%
\eol}\vss}\rg%
%
%
\rx{\vss\hfull{%
\rlx{\hss{$1008_{z}$}}\cg%
\e{0}%
\e{0}%
\e{1}%
\e{1}%
\e{0}%
\e{0}%
\e{0}%
\e{0}%
\e{0}%
\e{0}%
\e{0}%
\e{0}%
\e{0}%
\e{0}%
\e{0}%
\e{0}%
\e{0}%
\e{0}%
\e{1}%
\e{1}%
\eol}\vss}\rg%
%
%
\rx{\vss\hfull{%
\rlx{\hss{$560_{z}$}}\cg%
\e{0}%
\e{0}%
\e{1}%
\e{0}%
\e{0}%
\e{0}%
\e{0}%
\e{0}%
\e{1}%
\e{0}%
\e{0}%
\e{0}%
\e{0}%
\e{0}%
\e{0}%
\e{0}%
\e{0}%
\e{0}%
\e{1}%
\e{0}%
\eol}\vss}\rg%
%
%
\rx{\vss\hfull{%
\rlx{\hss{$1400_{zz}$}}\cg%
\e{0}%
\e{0}%
\e{0}%
\e{0}%
\e{0}%
\e{0}%
\e{0}%
\e{0}%
\e{0}%
\e{0}%
\e{0}%
\e{0}%
\e{1}%
\e{0}%
\e{0}%
\e{0}%
\e{0}%
\e{0}%
\e{0}%
\e{0}%
\eol}\vss}\rg%
%
%
\rx{\vss\hfull{%
\rlx{\hss{$4200_{z}$}}\cg%
\e{0}%
\e{0}%
\e{0}%
\e{0}%
\e{0}%
\e{0}%
\e{0}%
\e{0}%
\e{0}%
\e{0}%
\e{0}%
\e{0}%
\e{1}%
\e{1}%
\e{0}%
\e{1}%
\e{0}%
\e{0}%
\e{0}%
\e{0}%
\eol}\vss}\rg%
%
%
\rx{\vss\hfull{%
\rlx{\hss{$400_{z}$}}\cg%
\e{0}%
\e{0}%
\e{0}%
\e{0}%
\e{0}%
\e{0}%
\e{0}%
\e{0}%
\e{1}%
\e{0}%
\e{0}%
\e{0}%
\e{1}%
\e{0}%
\e{0}%
\e{0}%
\e{0}%
\e{0}%
\e{0}%
\e{0}%
\eol}\vss}\rg%
%
%
\rx{\vss\hfull{%
\rlx{\hss{$3240_{z}$}}\cg%
\e{0}%
\e{0}%
\e{0}%
\e{0}%
\e{0}%
\e{0}%
\e{0}%
\e{0}%
\e{1}%
\e{0}%
\e{0}%
\e{0}%
\e{1}%
\e{0}%
\e{0}%
\e{1}%
\e{1}%
\e{0}%
\e{1}%
\e{1}%
\eol}\vss}\rg%
%
%
\rx{\vss\hfull{%
\rlx{\hss{$4536_{z}$}}\cg%
\e{0}%
\e{0}%
\e{0}%
\e{0}%
\e{0}%
\e{0}%
\e{0}%
\e{0}%
\e{0}%
\e{0}%
\e{0}%
\e{0}%
\e{0}%
\e{0}%
\e{0}%
\e{1}%
\e{1}%
\e{0}%
\e{0}%
\e{0}%
\eol}\vss}\rg%
%
%
\rx{\vss\hfull{%
\rlx{\hss{$2400_{z}$}}\cg%
\e{0}%
\e{0}%
\e{0}%
\e{0}%
\e{0}%
\e{0}%
\e{1}%
\e{0}%
\e{0}%
\e{0}%
\e{0}%
\e{0}%
\e{0}%
\e{0}%
\e{0}%
\e{1}%
\e{0}%
\e{0}%
\e{0}%
\e{1}%
\eol}\vss}\rg%
%
%
\rx{\vss\hfull{%
\rlx{\hss{$3360_{z}$}}\cg%
\e{0}%
\e{0}%
\e{0}%
\e{0}%
\e{0}%
\e{0}%
\e{0}%
\e{0}%
\e{0}%
\e{0}%
\e{0}%
\e{0}%
\e{1}%
\e{0}%
\e{0}%
\e{1}%
\e{0}%
\e{0}%
\e{0}%
\e{1}%
\eol}\vss}\rg%
%
%
\rx{\vss\hfull{%
\rlx{\hss{$2800_{z}$}}\cg%
\e{0}%
\e{0}%
\e{0}%
\e{0}%
\e{0}%
\e{0}%
\e{0}%
\e{0}%
\e{0}%
\e{0}%
\e{0}%
\e{0}%
\e{0}%
\e{1}%
\e{0}%
\e{1}%
\e{0}%
\e{0}%
\e{1}%
\e{1}%
\eol}\vss}\rg%
%
%
\rx{\vss\hfull{%
\rlx{\hss{$4096_{z}$}}\cg%
\e{0}%
\e{0}%
\e{0}%
\e{0}%
\e{0}%
\e{0}%
\e{1}%
\e{0}%
\e{0}%
\e{0}%
\e{0}%
\e{0}%
\e{0}%
\e{0}%
\e{0}%
\e{1}%
\e{1}%
\e{0}%
\e{1}%
\e{1}%
\eol}\vss}\rg%
%
%
\rx{\vss\hfull{%
\rlx{\hss{$5600_{z}$}}\cg%
\e{0}%
\e{0}%
\e{0}%
\e{0}%
\e{0}%
\e{0}%
\e{1}%
\e{0}%
\e{0}%
\e{0}%
\e{0}%
\e{0}%
\e{0}%
\e{0}%
\e{0}%
\e{1}%
\e{1}%
\e{0}%
\e{0}%
\e{1}%
\eol}\vss}\rg%
%
%
\rx{\vss\hfull{%
\rlx{\hss{$448_{z}$}}\cg%
\e{0}%
\e{0}%
\e{0}%
\e{0}%
\e{0}%
\e{0}%
\e{0}%
\e{0}%
\e{0}%
\e{0}%
\e{0}%
\e{0}%
\e{0}%
\e{0}%
\e{0}%
\e{0}%
\e{1}%
\e{0}%
\e{0}%
\e{0}%
\eol}\vss}\rg%
%
%
\rx{\vss\hfull{%
\rlx{\hss{$448_{w}$}}\cg%
\e{0}%
\e{0}%
\e{0}%
\e{0}%
\e{1}%
\e{0}%
\e{1}%
\e{0}%
\e{0}%
\e{0}%
\e{0}%
\e{0}%
\e{0}%
\e{0}%
\e{0}%
\e{0}%
\e{0}%
\e{0}%
\e{0}%
\e{0}%
\eol}\vss}\rg%
%
%
\rx{\vss\hfull{%
\rlx{\hss{$1344_{w}$}}\cg%
\e{0}%
\e{0}%
\e{0}%
\e{0}%
\e{0}%
\e{0}%
\e{0}%
\e{0}%
\e{0}%
\e{0}%
\e{0}%
\e{0}%
\e{0}%
\e{0}%
\e{0}%
\e{0}%
\e{0}%
\e{0}%
\e{0}%
\e{0}%
\eol}\vss}\rg%
%
%
\rx{\vss\hfull{%
\rlx{\hss{$5600_{w}$}}\cg%
\e{0}%
\e{0}%
\e{0}%
\e{0}%
\e{0}%
\e{0}%
\e{1}%
\e{0}%
\e{0}%
\e{0}%
\e{0}%
\e{0}%
\e{0}%
\e{0}%
\e{0}%
\e{1}%
\e{0}%
\e{0}%
\e{0}%
\e{0}%
\eol}\vss}\rg%
%
%
\rx{\vss\hfull{%
\rlx{\hss{$2016_{w}$}}\cg%
\e{0}%
\e{0}%
\e{0}%
\e{0}%
\e{0}%
\e{0}%
\e{0}%
\e{0}%
\e{0}%
\e{0}%
\e{0}%
\e{0}%
\e{0}%
\e{0}%
\e{0}%
\e{0}%
\e{0}%
\e{0}%
\e{0}%
\e{0}%
\eol}\vss}\rg%
%
%
\rx{\vss\hfull{%
\rlx{\hss{$7168_{w}$}}\cg%
\e{0}%
\e{0}%
\e{0}%
\e{0}%
\e{0}%
\e{0}%
\e{0}%
\e{0}%
\e{0}%
\e{0}%
\e{0}%
\e{0}%
\e{0}%
\e{0}%
\e{1}%
\e{1}%
\e{0}%
\e{0}%
\e{0}%
\e{0}%
\eol}\vss}\rg%
\eop
\eject
\tablecont%
%
%
%
%
%
%
\rowpts=18 true pt%
\colpts=18 true pt%
\rowlabpts=40 true pt%
\collabpts=45 true pt%
\clx{\vss\hfull{%
\rlx{\hss{$ $}}\cg%
\cx{\hskip 16 true pt\flip{$280_a$}\hss}\cg%
\cx{\hskip 16 true pt\flip{$336_a$}\hss}\cg%
\cx{\hskip 16 true pt\flip{$216_a$}\hss}\cg%
\cx{\hskip 16 true pt\flip{$512_a$}\hss}\cg%
\cx{\hskip 16 true pt\flip{$378_a$}\hss}\cg%
\cx{\hskip 16 true pt\flip{$84_a$}\hss}\cg%
\cx{\hskip 16 true pt\flip{$420_a$}\hss}\cg%
\cx{\hskip 16 true pt\flip{$280_b$}\hss}\cg%
\cx{\hskip 16 true pt\flip{$210_b$}\hss}\cg%
\cx{\hskip 16 true pt\flip{$70_a$}\hss}\cg%
\cx{\hskip 16 true pt\flip{$1_a^{*}$}\hss}\cg%
\cx{\hskip 16 true pt\flip{$7_a^{*}$}\hss}\cg%
\cx{\hskip 16 true pt\flip{$27_a^{*}$}\hss}\cg%
\cx{\hskip 16 true pt\flip{$21_a^{*}$}\hss}\cg%
\cx{\hskip 16 true pt\flip{$35_a^{*}$}\hss}\cg%
\cx{\hskip 16 true pt\flip{$105_a^{*}$}\hss}\cg%
\cx{\hskip 16 true pt\flip{$189_a^{*}$}\hss}\cg%
\cx{\hskip 16 true pt\flip{$21_b^{*}$}\hss}\cg%
\cx{\hskip 16 true pt\flip{$35_b^{*}$}\hss}\cg%
\cx{\hskip 16 true pt\flip{$189_b^{*}$}\hss}\cg%
\eol}}\rg%
%
%
\rx{\vss\hfull{%
\rlx{\hss{$1_{x}$}}\cg%
\e{0}%
\e{0}%
\e{0}%
\e{0}%
\e{0}%
\e{0}%
\e{0}%
\e{0}%
\e{0}%
\e{0}%
\e{0}%
\e{0}%
\e{0}%
\e{0}%
\e{0}%
\e{0}%
\e{0}%
\e{0}%
\e{0}%
\e{0}%
\eol}\vss}\rg%
%
%
\rx{\vss\hfull{%
\rlx{\hss{$28_{x}$}}\cg%
\e{0}%
\e{0}%
\e{0}%
\e{0}%
\e{0}%
\e{0}%
\e{0}%
\e{0}%
\e{0}%
\e{0}%
\e{0}%
\e{1}%
\e{0}%
\e{0}%
\e{0}%
\e{0}%
\e{0}%
\e{0}%
\e{0}%
\e{0}%
\eol}\vss}\rg%
%
%
\rx{\vss\hfull{%
\rlx{\hss{$35_{x}$}}\cg%
\e{0}%
\e{0}%
\e{0}%
\e{0}%
\e{0}%
\e{0}%
\e{0}%
\e{0}%
\e{0}%
\e{0}%
\e{0}%
\e{1}%
\e{0}%
\e{0}%
\e{0}%
\e{0}%
\e{0}%
\e{0}%
\e{0}%
\e{0}%
\eol}\vss}\rg%
%
%
\rx{\vss\hfull{%
\rlx{\hss{$84_{x}$}}\cg%
\e{0}%
\e{0}%
\e{0}%
\e{0}%
\e{0}%
\e{0}%
\e{0}%
\e{0}%
\e{0}%
\e{0}%
\e{0}%
\e{0}%
\e{0}%
\e{0}%
\e{0}%
\e{0}%
\e{0}%
\e{1}%
\e{0}%
\e{0}%
\eol}\vss}\rg%
%
%
\rx{\vss\hfull{%
\rlx{\hss{$50_{x}$}}\cg%
\e{0}%
\e{0}%
\e{0}%
\e{0}%
\e{0}%
\e{0}%
\e{0}%
\e{0}%
\e{0}%
\e{0}%
\e{0}%
\e{0}%
\e{0}%
\e{0}%
\e{0}%
\e{0}%
\e{0}%
\e{0}%
\e{0}%
\e{0}%
\eol}\vss}\rg%
%
%
\rx{\vss\hfull{%
\rlx{\hss{$350_{x}$}}\cg%
\e{0}%
\e{0}%
\e{0}%
\e{0}%
\e{0}%
\e{0}%
\e{0}%
\e{0}%
\e{0}%
\e{0}%
\e{0}%
\e{0}%
\e{0}%
\e{0}%
\e{1}%
\e{1}%
\e{0}%
\e{0}%
\e{0}%
\e{0}%
\eol}\vss}\rg%
%
%
\rx{\vss\hfull{%
\rlx{\hss{$300_{x}$}}\cg%
\e{0}%
\e{0}%
\e{0}%
\e{0}%
\e{0}%
\e{0}%
\e{0}%
\e{0}%
\e{0}%
\e{0}%
\e{0}%
\e{0}%
\e{0}%
\e{0}%
\e{0}%
\e{1}%
\e{0}%
\e{0}%
\e{0}%
\e{0}%
\eol}\vss}\rg%
%
%
\rx{\vss\hfull{%
\rlx{\hss{$567_{x}$}}\cg%
\e{0}%
\e{0}%
\e{0}%
\e{0}%
\e{0}%
\e{0}%
\e{0}%
\e{0}%
\e{0}%
\e{0}%
\e{0}%
\e{1}%
\e{0}%
\e{0}%
\e{0}%
\e{1}%
\e{0}%
\e{1}%
\e{0}%
\e{0}%
\eol}\vss}\rg%
%
%
\rx{\vss\hfull{%
\rlx{\hss{$210_{x}$}}\cg%
\e{0}%
\e{0}%
\e{0}%
\e{0}%
\e{0}%
\e{0}%
\e{0}%
\e{0}%
\e{0}%
\e{0}%
\e{0}%
\e{1}%
\e{0}%
\e{0}%
\e{0}%
\e{0}%
\e{0}%
\e{0}%
\e{0}%
\e{0}%
\eol}\vss}\rg%
%
%
\rx{\vss\hfull{%
\rlx{\hss{$840_{x}$}}\cg%
\e{0}%
\e{0}%
\e{0}%
\e{0}%
\e{0}%
\e{1}%
\e{0}%
\e{0}%
\e{1}%
\e{0}%
\e{0}%
\e{0}%
\e{0}%
\e{0}%
\e{0}%
\e{0}%
\e{0}%
\e{0}%
\e{0}%
\e{0}%
\eol}\vss}\rg%
%
%
\rx{\vss\hfull{%
\rlx{\hss{$700_{x}$}}\cg%
\e{0}%
\e{0}%
\e{0}%
\e{0}%
\e{0}%
\e{0}%
\e{0}%
\e{0}%
\e{0}%
\e{0}%
\e{0}%
\e{0}%
\e{0}%
\e{0}%
\e{0}%
\e{0}%
\e{0}%
\e{0}%
\e{0}%
\e{1}%
\eol}\vss}\rg%
%
%
\rx{\vss\hfull{%
\rlx{\hss{$175_{x}$}}\cg%
\e{0}%
\e{0}%
\e{0}%
\e{0}%
\e{0}%
\e{0}%
\e{0}%
\e{0}%
\e{0}%
\e{0}%
\e{0}%
\e{0}%
\e{0}%
\e{0}%
\e{0}%
\e{0}%
\e{0}%
\e{0}%
\e{0}%
\e{0}%
\eol}\vss}\rg%
%
%
\rx{\vss\hfull{%
\rlx{\hss{$1400_{x}$}}\cg%
\e{0}%
\e{0}%
\e{0}%
\e{0}%
\e{0}%
\e{0}%
\e{0}%
\e{0}%
\e{0}%
\e{0}%
\e{0}%
\e{0}%
\e{0}%
\e{0}%
\e{0}%
\e{0}%
\e{0}%
\e{0}%
\e{0}%
\e{1}%
\eol}\vss}\rg%
%
%
\rx{\vss\hfull{%
\rlx{\hss{$1050_{x}$}}\cg%
\e{0}%
\e{0}%
\e{0}%
\e{0}%
\e{0}%
\e{0}%
\e{0}%
\e{1}%
\e{0}%
\e{0}%
\e{0}%
\e{0}%
\e{0}%
\e{0}%
\e{0}%
\e{0}%
\e{0}%
\e{0}%
\e{0}%
\e{1}%
\eol}\vss}\rg%
%
%
\rx{\vss\hfull{%
\rlx{\hss{$1575_{x}$}}\cg%
\e{0}%
\e{0}%
\e{0}%
\e{0}%
\e{0}%
\e{0}%
\e{0}%
\e{0}%
\e{0}%
\e{0}%
\e{0}%
\e{0}%
\e{0}%
\e{0}%
\e{0}%
\e{1}%
\e{0}%
\e{0}%
\e{0}%
\e{1}%
\eol}\vss}\rg%
%
%
\rx{\vss\hfull{%
\rlx{\hss{$1344_{x}$}}\cg%
\e{0}%
\e{0}%
\e{0}%
\e{0}%
\e{0}%
\e{0}%
\e{0}%
\e{1}%
\e{0}%
\e{0}%
\e{0}%
\e{0}%
\e{0}%
\e{0}%
\e{0}%
\e{1}%
\e{0}%
\e{1}%
\e{0}%
\e{1}%
\eol}\vss}\rg%
%
%
\rx{\vss\hfull{%
\rlx{\hss{$2100_{x}$}}\cg%
\e{0}%
\e{1}%
\e{0}%
\e{0}%
\e{0}%
\e{0}%
\e{1}%
\e{0}%
\e{0}%
\e{0}%
\e{0}%
\e{0}%
\e{0}%
\e{0}%
\e{1}%
\e{1}%
\e{0}%
\e{0}%
\e{0}%
\e{0}%
\eol}\vss}\rg%
%
%
\rx{\vss\hfull{%
\rlx{\hss{$2268_{x}$}}\cg%
\e{0}%
\e{0}%
\e{0}%
\e{0}%
\e{0}%
\e{0}%
\e{1}%
\e{0}%
\e{0}%
\e{0}%
\e{0}%
\e{0}%
\e{0}%
\e{0}%
\e{0}%
\e{1}%
\e{0}%
\e{0}%
\e{0}%
\e{0}%
\eol}\vss}\rg%
%
%
\rx{\vss\hfull{%
\rlx{\hss{$525_{x}$}}\cg%
\e{0}%
\e{0}%
\e{0}%
\e{0}%
\e{0}%
\e{0}%
\e{0}%
\e{0}%
\e{0}%
\e{0}%
\e{0}%
\e{0}%
\e{0}%
\e{0}%
\e{0}%
\e{0}%
\e{0}%
\e{1}%
\e{0}%
\e{0}%
\eol}\vss}\rg%
%
%
\rx{\vss\hfull{%
\rlx{\hss{$700_{xx}$}}\cg%
\e{0}%
\e{0}%
\e{0}%
\e{0}%
\e{0}%
\e{0}%
\e{0}%
\e{1}%
\e{0}%
\e{0}%
\e{0}%
\e{0}%
\e{0}%
\e{0}%
\e{0}%
\e{0}%
\e{0}%
\e{0}%
\e{0}%
\e{0}%
\eol}\vss}\rg%
%
%
\rx{\vss\hfull{%
\rlx{\hss{$972_{x}$}}\cg%
\e{0}%
\e{0}%
\e{0}%
\e{0}%
\e{0}%
\e{1}%
\e{0}%
\e{1}%
\e{0}%
\e{0}%
\e{0}%
\e{0}%
\e{0}%
\e{0}%
\e{0}%
\e{0}%
\e{0}%
\e{0}%
\e{0}%
\e{0}%
\eol}\vss}\rg%
%
%
\rx{\vss\hfull{%
\rlx{\hss{$4096_{x}$}}\cg%
\e{0}%
\e{0}%
\e{0}%
\e{1}%
\e{0}%
\e{0}%
\e{1}%
\e{1}%
\e{0}%
\e{0}%
\e{0}%
\e{0}%
\e{0}%
\e{0}%
\e{0}%
\e{1}%
\e{0}%
\e{0}%
\e{0}%
\e{1}%
\eol}\vss}\rg%
%
%
\rx{\vss\hfull{%
\rlx{\hss{$4200_{x}$}}\cg%
\e{0}%
\e{0}%
\e{0}%
\e{1}%
\e{0}%
\e{0}%
\e{1}%
\e{1}%
\e{1}%
\e{0}%
\e{0}%
\e{0}%
\e{0}%
\e{0}%
\e{0}%
\e{0}%
\e{0}%
\e{0}%
\e{0}%
\e{1}%
\eol}\vss}\rg%
%
%
\rx{\vss\hfull{%
\rlx{\hss{$2240_{x}$}}\cg%
\e{0}%
\e{0}%
\e{0}%
\e{0}%
\e{0}%
\e{0}%
\e{0}%
\e{1}%
\e{1}%
\e{0}%
\e{0}%
\e{0}%
\e{0}%
\e{0}%
\e{0}%
\e{0}%
\e{0}%
\e{0}%
\e{0}%
\e{1}%
\eol}\vss}\rg%
%
%
\rx{\vss\hfull{%
\rlx{\hss{$2835_{x}$}}\cg%
\e{0}%
\e{0}%
\e{0}%
\e{1}%
\e{0}%
\e{0}%
\e{0}%
\e{1}%
\e{1}%
\e{0}%
\e{0}%
\e{0}%
\e{0}%
\e{0}%
\e{0}%
\e{0}%
\e{0}%
\e{0}%
\e{0}%
\e{0}%
\eol}\vss}\rg%
%
%
\rx{\vss\hfull{%
\rlx{\hss{$6075_{x}$}}\cg%
\e{0}%
\e{1}%
\e{0}%
\e{1}%
\e{1}%
\e{0}%
\e{1}%
\e{1}%
\e{0}%
\e{0}%
\e{0}%
\e{0}%
\e{0}%
\e{0}%
\e{0}%
\e{0}%
\e{0}%
\e{0}%
\e{0}%
\e{1}%
\eol}\vss}\rg%
%
%
\rx{\vss\hfull{%
\rlx{\hss{$3200_{x}$}}\cg%
\e{0}%
\e{0}%
\e{1}%
\e{1}%
\e{0}%
\e{1}%
\e{1}%
\e{1}%
\e{0}%
\e{0}%
\e{0}%
\e{0}%
\e{0}%
\e{0}%
\e{0}%
\e{0}%
\e{0}%
\e{0}%
\e{0}%
\e{0}%
\eol}\vss}\rg%
%
%
\rx{\vss\hfull{%
\rlx{\hss{$70_{y}$}}\cg%
\e{0}%
\e{0}%
\e{0}%
\e{0}%
\e{0}%
\e{0}%
\e{0}%
\e{0}%
\e{0}%
\e{0}%
\e{0}%
\e{0}%
\e{0}%
\e{0}%
\e{1}%
\e{0}%
\e{0}%
\e{0}%
\e{0}%
\e{0}%
\eol}\vss}\rg%
%
%
\rx{\vss\hfull{%
\rlx{\hss{$1134_{y}$}}\cg%
\e{0}%
\e{0}%
\e{0}%
\e{0}%
\e{1}%
\e{0}%
\e{0}%
\e{0}%
\e{0}%
\e{0}%
\e{0}%
\e{0}%
\e{0}%
\e{0}%
\e{0}%
\e{0}%
\e{1}%
\e{0}%
\e{0}%
\e{0}%
\eol}\vss}\rg%
%
%
\rx{\vss\hfull{%
\rlx{\hss{$1680_{y}$}}\cg%
\e{1}%
\e{1}%
\e{0}%
\e{0}%
\e{0}%
\e{0}%
\e{0}%
\e{0}%
\e{0}%
\e{0}%
\e{0}%
\e{0}%
\e{0}%
\e{0}%
\e{1}%
\e{0}%
\e{1}%
\e{0}%
\e{0}%
\e{0}%
\eol}\vss}\rg%
%
%
\rx{\vss\hfull{%
\rlx{\hss{$168_{y}$}}\cg%
\e{0}%
\e{0}%
\e{0}%
\e{0}%
\e{0}%
\e{1}%
\e{0}%
\e{0}%
\e{0}%
\e{0}%
\e{0}%
\e{0}%
\e{0}%
\e{0}%
\e{0}%
\e{0}%
\e{0}%
\e{0}%
\e{0}%
\e{0}%
\eol}\vss}\rg%
%
%
\rx{\vss\hfull{%
\rlx{\hss{$420_{y}$}}\cg%
\e{0}%
\e{0}%
\e{0}%
\e{0}%
\e{0}%
\e{0}%
\e{0}%
\e{0}%
\e{1}%
\e{0}%
\e{0}%
\e{0}%
\e{0}%
\e{0}%
\e{0}%
\e{0}%
\e{0}%
\e{0}%
\e{0}%
\e{0}%
\eol}\vss}\rg%
%
%
\rx{\vss\hfull{%
\rlx{\hss{$3150_{y}$}}\cg%
\e{0}%
\e{0}%
\e{0}%
\e{1}%
\e{1}%
\e{0}%
\e{0}%
\e{0}%
\e{1}%
\e{1}%
\e{0}%
\e{0}%
\e{0}%
\e{0}%
\e{0}%
\e{0}%
\e{0}%
\e{0}%
\e{0}%
\e{0}%
\eol}\vss}\rg%
%
%
\rx{\vss\hfull{%
\rlx{\hss{$4200_{y}$}}\cg%
\e{0}%
\e{0}%
\e{1}%
\e{1}%
\e{1}%
\e{1}%
\e{1}%
\e{1}%
\e{1}%
\e{0}%
\e{0}%
\e{0}%
\e{0}%
\e{0}%
\e{0}%
\e{0}%
\e{0}%
\e{0}%
\e{0}%
\e{0}%
\eol}\vss}\rg%
\eop
\eject
\tablecont%
%
%
%
%
%
%
\rowpts=18 true pt%
\colpts=18 true pt%
\rowlabpts=40 true pt%
\collabpts=45 true pt%
\clx{\vss\hfull{%
\rlx{\hss{$ $}}\cg%
\cx{\hskip 16 true pt\flip{$280_a$}\hss}\cg%
\cx{\hskip 16 true pt\flip{$336_a$}\hss}\cg%
\cx{\hskip 16 true pt\flip{$216_a$}\hss}\cg%
\cx{\hskip 16 true pt\flip{$512_a$}\hss}\cg%
\cx{\hskip 16 true pt\flip{$378_a$}\hss}\cg%
\cx{\hskip 16 true pt\flip{$84_a$}\hss}\cg%
\cx{\hskip 16 true pt\flip{$420_a$}\hss}\cg%
\cx{\hskip 16 true pt\flip{$280_b$}\hss}\cg%
\cx{\hskip 16 true pt\flip{$210_b$}\hss}\cg%
\cx{\hskip 16 true pt\flip{$70_a$}\hss}\cg%
\cx{\hskip 16 true pt\flip{$1_a^{*}$}\hss}\cg%
\cx{\hskip 16 true pt\flip{$7_a^{*}$}\hss}\cg%
\cx{\hskip 16 true pt\flip{$27_a^{*}$}\hss}\cg%
\cx{\hskip 16 true pt\flip{$21_a^{*}$}\hss}\cg%
\cx{\hskip 16 true pt\flip{$35_a^{*}$}\hss}\cg%
\cx{\hskip 16 true pt\flip{$105_a^{*}$}\hss}\cg%
\cx{\hskip 16 true pt\flip{$189_a^{*}$}\hss}\cg%
\cx{\hskip 16 true pt\flip{$21_b^{*}$}\hss}\cg%
\cx{\hskip 16 true pt\flip{$35_b^{*}$}\hss}\cg%
\cx{\hskip 16 true pt\flip{$189_b^{*}$}\hss}\cg%
\eol}}\rg%
%
%
\rx{\vss\hfull{%
\rlx{\hss{$2688_{y}$}}\cg%
\e{0}%
\e{1}%
\e{1}%
\e{1}%
\e{0}%
\e{0}%
\e{0}%
\e{1}%
\e{0}%
\e{0}%
\e{0}%
\e{0}%
\e{0}%
\e{0}%
\e{0}%
\e{0}%
\e{0}%
\e{0}%
\e{0}%
\e{0}%
\eol}\vss}\rg%
%
%
\rx{\vss\hfull{%
\rlx{\hss{$2100_{y}$}}\cg%
\e{0}%
\e{1}%
\e{0}%
\e{0}%
\e{0}%
\e{0}%
\e{1}%
\e{0}%
\e{0}%
\e{0}%
\e{0}%
\e{0}%
\e{0}%
\e{0}%
\e{0}%
\e{0}%
\e{0}%
\e{0}%
\e{0}%
\e{0}%
\eol}\vss}\rg%
%
%
\rx{\vss\hfull{%
\rlx{\hss{$1400_{y}$}}\cg%
\e{1}%
\e{0}%
\e{0}%
\e{0}%
\e{0}%
\e{0}%
\e{1}%
\e{0}%
\e{0}%
\e{0}%
\e{0}%
\e{0}%
\e{0}%
\e{0}%
\e{0}%
\e{0}%
\e{0}%
\e{0}%
\e{0}%
\e{0}%
\eol}\vss}\rg%
%
%
\rx{\vss\hfull{%
\rlx{\hss{$4536_{y}$}}\cg%
\e{1}%
\e{1}%
\e{0}%
\e{1}%
\e{0}%
\e{0}%
\e{1}%
\e{0}%
\e{0}%
\e{0}%
\e{0}%
\e{0}%
\e{0}%
\e{0}%
\e{0}%
\e{0}%
\e{0}%
\e{0}%
\e{0}%
\e{0}%
\eol}\vss}\rg%
%
%
\rx{\vss\hfull{%
\rlx{\hss{$5670_{y}$}}\cg%
\e{1}%
\e{1}%
\e{0}%
\e{1}%
\e{1}%
\e{0}%
\e{1}%
\e{0}%
\e{0}%
\e{0}%
\e{0}%
\e{0}%
\e{0}%
\e{0}%
\e{0}%
\e{0}%
\e{1}%
\e{0}%
\e{0}%
\e{0}%
\eol}\vss}\rg%
%
%
\rx{\vss\hfull{%
\rlx{\hss{$4480_{y}$}}\cg%
\e{0}%
\e{0}%
\e{0}%
\e{1}%
\e{1}%
\e{0}%
\e{1}%
\e{0}%
\e{1}%
\e{0}%
\e{0}%
\e{0}%
\e{0}%
\e{0}%
\e{0}%
\e{0}%
\e{0}%
\e{0}%
\e{0}%
\e{0}%
\eol}\vss}\rg%
%
%
\rx{\vss\hfull{%
\rlx{\hss{$8_{z}$}}\cg%
\e{0}%
\e{0}%
\e{0}%
\e{0}%
\e{0}%
\e{0}%
\e{0}%
\e{0}%
\e{0}%
\e{0}%
\e{0}%
\e{1}%
\e{0}%
\e{0}%
\e{0}%
\e{0}%
\e{0}%
\e{0}%
\e{0}%
\e{0}%
\eol}\vss}\rg%
%
%
\rx{\vss\hfull{%
\rlx{\hss{$56_{z}$}}\cg%
\e{0}%
\e{0}%
\e{0}%
\e{0}%
\e{0}%
\e{0}%
\e{0}%
\e{0}%
\e{0}%
\e{0}%
\e{0}%
\e{0}%
\e{0}%
\e{0}%
\e{1}%
\e{0}%
\e{0}%
\e{0}%
\e{0}%
\e{0}%
\eol}\vss}\rg%
%
%
\rx{\vss\hfull{%
\rlx{\hss{$160_{z}$}}\cg%
\e{0}%
\e{0}%
\e{0}%
\e{0}%
\e{0}%
\e{0}%
\e{0}%
\e{0}%
\e{0}%
\e{0}%
\e{0}%
\e{1}%
\e{0}%
\e{0}%
\e{0}%
\e{1}%
\e{0}%
\e{0}%
\e{0}%
\e{0}%
\eol}\vss}\rg%
%
%
\rx{\vss\hfull{%
\rlx{\hss{$112_{z}$}}\cg%
\e{0}%
\e{0}%
\e{0}%
\e{0}%
\e{0}%
\e{0}%
\e{0}%
\e{0}%
\e{0}%
\e{0}%
\e{0}%
\e{1}%
\e{0}%
\e{0}%
\e{0}%
\e{0}%
\e{0}%
\e{1}%
\e{0}%
\e{0}%
\eol}\vss}\rg%
%
%
\rx{\vss\hfull{%
\rlx{\hss{$840_{z}$}}\cg%
\e{0}%
\e{0}%
\e{0}%
\e{0}%
\e{0}%
\e{0}%
\e{0}%
\e{0}%
\e{0}%
\e{0}%
\e{0}%
\e{0}%
\e{0}%
\e{0}%
\e{0}%
\e{1}%
\e{0}%
\e{0}%
\e{0}%
\e{0}%
\eol}\vss}\rg%
%
%
\rx{\vss\hfull{%
\rlx{\hss{$1296_{z}$}}\cg%
\e{0}%
\e{0}%
\e{0}%
\e{0}%
\e{0}%
\e{0}%
\e{0}%
\e{0}%
\e{0}%
\e{0}%
\e{0}%
\e{0}%
\e{0}%
\e{0}%
\e{1}%
\e{1}%
\e{0}%
\e{0}%
\e{0}%
\e{0}%
\eol}\vss}\rg%
%
%
\rx{\vss\hfull{%
\rlx{\hss{$1400_{z}$}}\cg%
\e{0}%
\e{0}%
\e{0}%
\e{0}%
\e{0}%
\e{0}%
\e{0}%
\e{0}%
\e{0}%
\e{0}%
\e{0}%
\e{0}%
\e{0}%
\e{0}%
\e{0}%
\e{1}%
\e{0}%
\e{1}%
\e{0}%
\e{1}%
\eol}\vss}\rg%
%
%
\rx{\vss\hfull{%
\rlx{\hss{$1008_{z}$}}\cg%
\e{0}%
\e{0}%
\e{0}%
\e{0}%
\e{0}%
\e{0}%
\e{0}%
\e{0}%
\e{0}%
\e{0}%
\e{0}%
\e{0}%
\e{0}%
\e{0}%
\e{0}%
\e{1}%
\e{0}%
\e{0}%
\e{0}%
\e{0}%
\eol}\vss}\rg%
%
%
\rx{\vss\hfull{%
\rlx{\hss{$560_{z}$}}\cg%
\e{0}%
\e{0}%
\e{0}%
\e{0}%
\e{0}%
\e{0}%
\e{0}%
\e{0}%
\e{0}%
\e{0}%
\e{0}%
\e{1}%
\e{0}%
\e{0}%
\e{0}%
\e{1}%
\e{0}%
\e{1}%
\e{0}%
\e{1}%
\eol}\vss}\rg%
%
%
\rx{\vss\hfull{%
\rlx{\hss{$1400_{zz}$}}\cg%
\e{0}%
\e{0}%
\e{0}%
\e{0}%
\e{0}%
\e{0}%
\e{0}%
\e{1}%
\e{1}%
\e{0}%
\e{0}%
\e{0}%
\e{0}%
\e{0}%
\e{0}%
\e{0}%
\e{0}%
\e{0}%
\e{0}%
\e{1}%
\eol}\vss}\rg%
%
%
\rx{\vss\hfull{%
\rlx{\hss{$4200_{z}$}}\cg%
\e{0}%
\e{0}%
\e{0}%
\e{1}%
\e{1}%
\e{0}%
\e{0}%
\e{1}%
\e{1}%
\e{0}%
\e{0}%
\e{0}%
\e{0}%
\e{0}%
\e{0}%
\e{0}%
\e{0}%
\e{0}%
\e{0}%
\e{0}%
\eol}\vss}\rg%
%
%
\rx{\vss\hfull{%
\rlx{\hss{$400_{z}$}}\cg%
\e{0}%
\e{0}%
\e{0}%
\e{0}%
\e{0}%
\e{0}%
\e{0}%
\e{0}%
\e{0}%
\e{0}%
\e{0}%
\e{0}%
\e{0}%
\e{0}%
\e{0}%
\e{0}%
\e{0}%
\e{0}%
\e{0}%
\e{1}%
\eol}\vss}\rg%
%
%
\rx{\vss\hfull{%
\rlx{\hss{$3240_{z}$}}\cg%
\e{0}%
\e{0}%
\e{0}%
\e{0}%
\e{0}%
\e{0}%
\e{0}%
\e{1}%
\e{0}%
\e{0}%
\e{0}%
\e{0}%
\e{0}%
\e{0}%
\e{0}%
\e{1}%
\e{0}%
\e{0}%
\e{0}%
\e{2}%
\eol}\vss}\rg%
%
%
\rx{\vss\hfull{%
\rlx{\hss{$4536_{z}$}}\cg%
\e{0}%
\e{0}%
\e{0}%
\e{1}%
\e{0}%
\e{1}%
\e{1}%
\e{1}%
\e{1}%
\e{0}%
\e{0}%
\e{0}%
\e{0}%
\e{0}%
\e{0}%
\e{0}%
\e{0}%
\e{0}%
\e{0}%
\e{1}%
\eol}\vss}\rg%
%
%
\rx{\vss\hfull{%
\rlx{\hss{$2400_{z}$}}\cg%
\e{0}%
\e{1}%
\e{0}%
\e{0}%
\e{0}%
\e{0}%
\e{0}%
\e{0}%
\e{0}%
\e{0}%
\e{0}%
\e{0}%
\e{0}%
\e{0}%
\e{1}%
\e{0}%
\e{1}%
\e{0}%
\e{0}%
\e{0}%
\eol}\vss}\rg%
%
%
\rx{\vss\hfull{%
\rlx{\hss{$3360_{z}$}}\cg%
\e{0}%
\e{0}%
\e{0}%
\e{1}%
\e{0}%
\e{0}%
\e{0}%
\e{1}%
\e{0}%
\e{0}%
\e{0}%
\e{0}%
\e{0}%
\e{0}%
\e{0}%
\e{0}%
\e{0}%
\e{0}%
\e{0}%
\e{1}%
\eol}\vss}\rg%
%
%
\rx{\vss\hfull{%
\rlx{\hss{$2800_{z}$}}\cg%
\e{0}%
\e{0}%
\e{0}%
\e{0}%
\e{0}%
\e{0}%
\e{1}%
\e{0}%
\e{0}%
\e{0}%
\e{0}%
\e{0}%
\e{0}%
\e{0}%
\e{0}%
\e{0}%
\e{0}%
\e{0}%
\e{0}%
\e{0}%
\eol}\vss}\rg%
%
%
\rx{\vss\hfull{%
\rlx{\hss{$4096_{z}$}}\cg%
\e{0}%
\e{0}%
\e{0}%
\e{0}%
\e{0}%
\e{0}%
\e{1}%
\e{1}%
\e{0}%
\e{0}%
\e{0}%
\e{0}%
\e{0}%
\e{0}%
\e{0}%
\e{1}%
\e{0}%
\e{0}%
\e{0}%
\e{1}%
\eol}\vss}\rg%
%
%
\rx{\vss\hfull{%
\rlx{\hss{$5600_{z}$}}\cg%
\e{0}%
\e{1}%
\e{0}%
\e{1}%
\e{0}%
\e{0}%
\e{2}%
\e{0}%
\e{0}%
\e{0}%
\e{0}%
\e{0}%
\e{0}%
\e{0}%
\e{0}%
\e{0}%
\e{0}%
\e{0}%
\e{0}%
\e{0}%
\eol}\vss}\rg%
%
%
\rx{\vss\hfull{%
\rlx{\hss{$448_{z}$}}\cg%
\e{0}%
\e{0}%
\e{0}%
\e{0}%
\e{0}%
\e{0}%
\e{0}%
\e{0}%
\e{0}%
\e{0}%
\e{0}%
\e{0}%
\e{0}%
\e{0}%
\e{0}%
\e{0}%
\e{0}%
\e{1}%
\e{0}%
\e{1}%
\eol}\vss}\rg%
%
%
\rx{\vss\hfull{%
\rlx{\hss{$448_{w}$}}\cg%
\e{0}%
\e{0}%
\e{0}%
\e{0}%
\e{0}%
\e{0}%
\e{0}%
\e{0}%
\e{0}%
\e{0}%
\e{0}%
\e{0}%
\e{0}%
\e{0}%
\e{1}%
\e{0}%
\e{1}%
\e{0}%
\e{0}%
\e{0}%
\eol}\vss}\rg%
%
%
\rx{\vss\hfull{%
\rlx{\hss{$1344_{w}$}}\cg%
\e{0}%
\e{0}%
\e{0}%
\e{0}%
\e{1}%
\e{1}%
\e{0}%
\e{0}%
\e{1}%
\e{0}%
\e{0}%
\e{0}%
\e{0}%
\e{0}%
\e{0}%
\e{0}%
\e{0}%
\e{0}%
\e{0}%
\e{0}%
\eol}\vss}\rg%
%
%
\rx{\vss\hfull{%
\rlx{\hss{$5600_{w}$}}\cg%
\e{1}%
\e{1}%
\e{0}%
\e{1}%
\e{1}%
\e{0}%
\e{1}%
\e{1}%
\e{0}%
\e{0}%
\e{0}%
\e{0}%
\e{0}%
\e{0}%
\e{0}%
\e{0}%
\e{1}%
\e{0}%
\e{0}%
\e{0}%
\eol}\vss}\rg%
%
%
\rx{\vss\hfull{%
\rlx{\hss{$2016_{w}$}}\cg%
\e{0}%
\e{0}%
\e{1}%
\e{1}%
\e{0}%
\e{0}%
\e{0}%
\e{0}%
\e{1}%
\e{1}%
\e{0}%
\e{0}%
\e{0}%
\e{0}%
\e{0}%
\e{0}%
\e{0}%
\e{0}%
\e{0}%
\e{0}%
\eol}\vss}\rg%
%
%
\rx{\vss\hfull{%
\rlx{\hss{$7168_{w}$}}\cg%
\e{0}%
\e{1}%
\e{1}%
\e{2}%
\e{1}%
\e{0}%
\e{1}%
\e{1}%
\e{1}%
\e{0}%
\e{0}%
\e{0}%
\e{0}%
\e{0}%
\e{0}%
\e{0}%
\e{0}%
\e{0}%
\e{0}%
\e{0}%
\eol}\vss}\rg%
\eop
\eject
\tablecont%
%
%
%
%
%
%
\rowpts=18 true pt%
\colpts=18 true pt%
\rowlabpts=40 true pt%
\collabpts=45 true pt%
\clx{\vss\hfull{%
\rlx{\hss{$ $}}\cg%
\cx{\hskip 16 true pt\flip{$189_c^{*}$}\hss}\cg%
\cx{\hskip 16 true pt\flip{$15_a^{*}$}\hss}\cg%
\cx{\hskip 16 true pt\flip{$105_b^{*}$}\hss}\cg%
\cx{\hskip 16 true pt\flip{$105_c^{*}$}\hss}\cg%
\cx{\hskip 16 true pt\flip{$315_a^{*}$}\hss}\cg%
\cx{\hskip 16 true pt\flip{$405_a^{*}$}\hss}\cg%
\cx{\hskip 16 true pt\flip{$168_a^{*}$}\hss}\cg%
\cx{\hskip 16 true pt\flip{$56_a^{*}$}\hss}\cg%
\cx{\hskip 16 true pt\flip{$120_a^{*}$}\hss}\cg%
\cx{\hskip 16 true pt\flip{$210_a^{*}$}\hss}\cg%
\cx{\hskip 16 true pt\flip{$280_a^{*}$}\hss}\cg%
\cx{\hskip 16 true pt\flip{$336_a^{*}$}\hss}\cg%
\cx{\hskip 16 true pt\flip{$216_a^{*}$}\hss}\cg%
\cx{\hskip 16 true pt\flip{$512_a^{*}$}\hss}\cg%
\cx{\hskip 16 true pt\flip{$378_a^{*}$}\hss}\cg%
\cx{\hskip 16 true pt\flip{$84_a^{*}$}\hss}\cg%
\cx{\hskip 16 true pt\flip{$420_a^{*}$}\hss}\cg%
\cx{\hskip 16 true pt\flip{$280_b^{*}$}\hss}\cg%
\cx{\hskip 16 true pt\flip{$210_b^{*}$}\hss}\cg%
\cx{\hskip 16 true pt\flip{$70_a^{*}$}\hss}\cg%
\eol}}\rg%
%
%
\rx{\vss\hfull{%
\rlx{\hss{$1_{x}$}}\cg%
\e{0}%
\e{0}%
\e{0}%
\e{0}%
\e{0}%
\e{0}%
\e{0}%
\e{0}%
\e{0}%
\e{0}%
\e{0}%
\e{0}%
\e{0}%
\e{0}%
\e{0}%
\e{0}%
\e{0}%
\e{0}%
\e{0}%
\e{0}%
\eol}\vss}\rg%
%
%
\rx{\vss\hfull{%
\rlx{\hss{$28_{x}$}}\cg%
\e{0}%
\e{0}%
\e{0}%
\e{0}%
\e{0}%
\e{0}%
\e{0}%
\e{0}%
\e{0}%
\e{0}%
\e{0}%
\e{0}%
\e{0}%
\e{0}%
\e{0}%
\e{0}%
\e{0}%
\e{0}%
\e{0}%
\e{0}%
\eol}\vss}\rg%
%
%
\rx{\vss\hfull{%
\rlx{\hss{$35_{x}$}}\cg%
\e{0}%
\e{0}%
\e{0}%
\e{0}%
\e{0}%
\e{0}%
\e{0}%
\e{0}%
\e{0}%
\e{0}%
\e{0}%
\e{0}%
\e{0}%
\e{0}%
\e{0}%
\e{0}%
\e{0}%
\e{0}%
\e{0}%
\e{0}%
\eol}\vss}\rg%
%
%
\rx{\vss\hfull{%
\rlx{\hss{$84_{x}$}}\cg%
\e{0}%
\e{0}%
\e{0}%
\e{0}%
\e{0}%
\e{0}%
\e{0}%
\e{0}%
\e{0}%
\e{0}%
\e{0}%
\e{0}%
\e{0}%
\e{0}%
\e{0}%
\e{0}%
\e{0}%
\e{0}%
\e{0}%
\e{0}%
\eol}\vss}\rg%
%
%
\rx{\vss\hfull{%
\rlx{\hss{$50_{x}$}}\cg%
\e{0}%
\e{1}%
\e{0}%
\e{0}%
\e{0}%
\e{0}%
\e{0}%
\e{0}%
\e{0}%
\e{0}%
\e{0}%
\e{0}%
\e{0}%
\e{0}%
\e{0}%
\e{0}%
\e{0}%
\e{0}%
\e{0}%
\e{0}%
\eol}\vss}\rg%
%
%
\rx{\vss\hfull{%
\rlx{\hss{$350_{x}$}}\cg%
\e{0}%
\e{0}%
\e{0}%
\e{0}%
\e{0}%
\e{0}%
\e{0}%
\e{0}%
\e{0}%
\e{0}%
\e{0}%
\e{0}%
\e{0}%
\e{0}%
\e{0}%
\e{0}%
\e{0}%
\e{0}%
\e{0}%
\e{0}%
\eol}\vss}\rg%
%
%
\rx{\vss\hfull{%
\rlx{\hss{$300_{x}$}}\cg%
\e{0}%
\e{0}%
\e{0}%
\e{0}%
\e{0}%
\e{0}%
\e{0}%
\e{0}%
\e{0}%
\e{0}%
\e{0}%
\e{0}%
\e{0}%
\e{0}%
\e{0}%
\e{0}%
\e{0}%
\e{0}%
\e{0}%
\e{0}%
\eol}\vss}\rg%
%
%
\rx{\vss\hfull{%
\rlx{\hss{$567_{x}$}}\cg%
\e{0}%
\e{0}%
\e{0}%
\e{0}%
\e{0}%
\e{0}%
\e{0}%
\e{1}%
\e{0}%
\e{0}%
\e{0}%
\e{0}%
\e{0}%
\e{0}%
\e{0}%
\e{0}%
\e{0}%
\e{0}%
\e{0}%
\e{0}%
\eol}\vss}\rg%
%
%
\rx{\vss\hfull{%
\rlx{\hss{$210_{x}$}}\cg%
\e{0}%
\e{0}%
\e{0}%
\e{0}%
\e{0}%
\e{0}%
\e{0}%
\e{1}%
\e{0}%
\e{0}%
\e{0}%
\e{0}%
\e{0}%
\e{0}%
\e{0}%
\e{0}%
\e{0}%
\e{0}%
\e{0}%
\e{0}%
\eol}\vss}\rg%
%
%
\rx{\vss\hfull{%
\rlx{\hss{$840_{x}$}}\cg%
\e{0}%
\e{0}%
\e{0}%
\e{0}%
\e{0}%
\e{0}%
\e{0}%
\e{0}%
\e{0}%
\e{0}%
\e{0}%
\e{0}%
\e{0}%
\e{0}%
\e{1}%
\e{0}%
\e{0}%
\e{0}%
\e{0}%
\e{0}%
\eol}\vss}\rg%
%
%
\rx{\vss\hfull{%
\rlx{\hss{$700_{x}$}}\cg%
\e{0}%
\e{0}%
\e{0}%
\e{0}%
\e{0}%
\e{0}%
\e{0}%
\e{1}%
\e{0}%
\e{0}%
\e{0}%
\e{0}%
\e{0}%
\e{0}%
\e{0}%
\e{0}%
\e{0}%
\e{0}%
\e{0}%
\e{0}%
\eol}\vss}\rg%
%
%
\rx{\vss\hfull{%
\rlx{\hss{$175_{x}$}}\cg%
\e{0}%
\e{0}%
\e{0}%
\e{0}%
\e{0}%
\e{0}%
\e{0}%
\e{0}%
\e{0}%
\e{0}%
\e{0}%
\e{0}%
\e{0}%
\e{0}%
\e{0}%
\e{0}%
\e{0}%
\e{0}%
\e{0}%
\e{1}%
\eol}\vss}\rg%
%
%
\rx{\vss\hfull{%
\rlx{\hss{$1400_{x}$}}\cg%
\e{0}%
\e{0}%
\e{0}%
\e{0}%
\e{1}%
\e{0}%
\e{0}%
\e{1}%
\e{0}%
\e{0}%
\e{0}%
\e{0}%
\e{0}%
\e{0}%
\e{0}%
\e{0}%
\e{0}%
\e{0}%
\e{0}%
\e{0}%
\eol}\vss}\rg%
%
%
\rx{\vss\hfull{%
\rlx{\hss{$1050_{x}$}}\cg%
\e{0}%
\e{1}%
\e{0}%
\e{0}%
\e{0}%
\e{0}%
\e{0}%
\e{0}%
\e{0}%
\e{0}%
\e{0}%
\e{0}%
\e{1}%
\e{0}%
\e{0}%
\e{0}%
\e{0}%
\e{0}%
\e{0}%
\e{0}%
\eol}\vss}\rg%
%
%
\rx{\vss\hfull{%
\rlx{\hss{$1575_{x}$}}\cg%
\e{0}%
\e{0}%
\e{0}%
\e{0}%
\e{0}%
\e{0}%
\e{0}%
\e{1}%
\e{0}%
\e{0}%
\e{1}%
\e{0}%
\e{0}%
\e{0}%
\e{0}%
\e{0}%
\e{0}%
\e{0}%
\e{0}%
\e{0}%
\eol}\vss}\rg%
%
%
\rx{\vss\hfull{%
\rlx{\hss{$1344_{x}$}}\cg%
\e{1}%
\e{0}%
\e{0}%
\e{0}%
\e{0}%
\e{0}%
\e{0}%
\e{0}%
\e{0}%
\e{0}%
\e{0}%
\e{0}%
\e{0}%
\e{0}%
\e{0}%
\e{0}%
\e{0}%
\e{0}%
\e{0}%
\e{0}%
\eol}\vss}\rg%
%
%
\rx{\vss\hfull{%
\rlx{\hss{$2100_{x}$}}\cg%
\e{1}%
\e{0}%
\e{0}%
\e{0}%
\e{0}%
\e{0}%
\e{0}%
\e{0}%
\e{0}%
\e{0}%
\e{1}%
\e{1}%
\e{0}%
\e{0}%
\e{0}%
\e{0}%
\e{0}%
\e{0}%
\e{0}%
\e{0}%
\eol}\vss}\rg%
%
%
\rx{\vss\hfull{%
\rlx{\hss{$2268_{x}$}}\cg%
\e{1}%
\e{0}%
\e{0}%
\e{0}%
\e{1}%
\e{0}%
\e{0}%
\e{1}%
\e{0}%
\e{0}%
\e{1}%
\e{0}%
\e{0}%
\e{0}%
\e{0}%
\e{0}%
\e{0}%
\e{0}%
\e{0}%
\e{0}%
\eol}\vss}\rg%
%
%
\rx{\vss\hfull{%
\rlx{\hss{$525_{x}$}}\cg%
\e{1}%
\e{0}%
\e{0}%
\e{0}%
\e{0}%
\e{0}%
\e{0}%
\e{0}%
\e{0}%
\e{0}%
\e{0}%
\e{0}%
\e{0}%
\e{0}%
\e{0}%
\e{0}%
\e{0}%
\e{0}%
\e{0}%
\e{0}%
\eol}\vss}\rg%
%
%
\rx{\vss\hfull{%
\rlx{\hss{$700_{xx}$}}\cg%
\e{0}%
\e{1}%
\e{0}%
\e{0}%
\e{0}%
\e{0}%
\e{0}%
\e{0}%
\e{0}%
\e{0}%
\e{0}%
\e{0}%
\e{1}%
\e{0}%
\e{0}%
\e{1}%
\e{0}%
\e{0}%
\e{0}%
\e{0}%
\eol}\vss}\rg%
%
%
\rx{\vss\hfull{%
\rlx{\hss{$972_{x}$}}\cg%
\e{1}%
\e{0}%
\e{0}%
\e{0}%
\e{0}%
\e{0}%
\e{0}%
\e{0}%
\e{0}%
\e{0}%
\e{0}%
\e{0}%
\e{1}%
\e{0}%
\e{0}%
\e{0}%
\e{0}%
\e{0}%
\e{0}%
\e{0}%
\eol}\vss}\rg%
%
%
\rx{\vss\hfull{%
\rlx{\hss{$4096_{x}$}}\cg%
\e{1}%
\e{0}%
\e{0}%
\e{0}%
\e{1}%
\e{0}%
\e{0}%
\e{0}%
\e{0}%
\e{0}%
\e{1}%
\e{1}%
\e{0}%
\e{0}%
\e{1}%
\e{0}%
\e{0}%
\e{0}%
\e{0}%
\e{0}%
\eol}\vss}\rg%
%
%
\rx{\vss\hfull{%
\rlx{\hss{$4200_{x}$}}\cg%
\e{0}%
\e{0}%
\e{0}%
\e{0}%
\e{1}%
\e{0}%
\e{0}%
\e{0}%
\e{0}%
\e{0}%
\e{1}%
\e{0}%
\e{1}%
\e{1}%
\e{1}%
\e{0}%
\e{0}%
\e{0}%
\e{0}%
\e{0}%
\eol}\vss}\rg%
%
%
\rx{\vss\hfull{%
\rlx{\hss{$2240_{x}$}}\cg%
\e{0}%
\e{0}%
\e{0}%
\e{0}%
\e{1}%
\e{0}%
\e{0}%
\e{0}%
\e{0}%
\e{0}%
\e{0}%
\e{0}%
\e{0}%
\e{0}%
\e{1}%
\e{0}%
\e{0}%
\e{0}%
\e{0}%
\e{1}%
\eol}\vss}\rg%
%
%
\rx{\vss\hfull{%
\rlx{\hss{$2835_{x}$}}\cg%
\e{0}%
\e{0}%
\e{0}%
\e{0}%
\e{1}%
\e{0}%
\e{0}%
\e{0}%
\e{0}%
\e{0}%
\e{0}%
\e{0}%
\e{1}%
\e{1}%
\e{0}%
\e{0}%
\e{0}%
\e{0}%
\e{1}%
\e{1}%
\eol}\vss}\rg%
%
%
\rx{\vss\hfull{%
\rlx{\hss{$6075_{x}$}}\cg%
\e{1}%
\e{0}%
\e{0}%
\e{0}%
\e{1}%
\e{0}%
\e{0}%
\e{0}%
\e{0}%
\e{0}%
\e{1}%
\e{1}%
\e{1}%
\e{1}%
\e{1}%
\e{0}%
\e{1}%
\e{0}%
\e{0}%
\e{0}%
\eol}\vss}\rg%
%
%
\rx{\vss\hfull{%
\rlx{\hss{$3200_{x}$}}\cg%
\e{1}%
\e{0}%
\e{0}%
\e{1}%
\e{0}%
\e{0}%
\e{0}%
\e{0}%
\e{0}%
\e{0}%
\e{0}%
\e{1}%
\e{0}%
\e{1}%
\e{1}%
\e{0}%
\e{0}%
\e{0}%
\e{0}%
\e{0}%
\eol}\vss}\rg%
%
%
\rx{\vss\hfull{%
\rlx{\hss{$70_{y}$}}\cg%
\e{0}%
\e{0}%
\e{0}%
\e{0}%
\e{0}%
\e{0}%
\e{0}%
\e{0}%
\e{0}%
\e{0}%
\e{0}%
\e{0}%
\e{0}%
\e{0}%
\e{0}%
\e{0}%
\e{0}%
\e{0}%
\e{0}%
\e{0}%
\eol}\vss}\rg%
%
%
\rx{\vss\hfull{%
\rlx{\hss{$1134_{y}$}}\cg%
\e{0}%
\e{0}%
\e{0}%
\e{0}%
\e{0}%
\e{0}%
\e{0}%
\e{0}%
\e{0}%
\e{0}%
\e{0}%
\e{0}%
\e{0}%
\e{0}%
\e{1}%
\e{0}%
\e{0}%
\e{0}%
\e{0}%
\e{0}%
\eol}\vss}\rg%
%
%
\rx{\vss\hfull{%
\rlx{\hss{$1680_{y}$}}\cg%
\e{0}%
\e{0}%
\e{0}%
\e{0}%
\e{0}%
\e{0}%
\e{0}%
\e{0}%
\e{0}%
\e{0}%
\e{1}%
\e{1}%
\e{0}%
\e{0}%
\e{0}%
\e{0}%
\e{0}%
\e{0}%
\e{0}%
\e{0}%
\eol}\vss}\rg%
%
%
\rx{\vss\hfull{%
\rlx{\hss{$168_{y}$}}\cg%
\e{0}%
\e{0}%
\e{0}%
\e{0}%
\e{0}%
\e{0}%
\e{0}%
\e{0}%
\e{0}%
\e{0}%
\e{0}%
\e{0}%
\e{0}%
\e{0}%
\e{0}%
\e{1}%
\e{0}%
\e{0}%
\e{0}%
\e{0}%
\eol}\vss}\rg%
%
%
\rx{\vss\hfull{%
\rlx{\hss{$420_{y}$}}\cg%
\e{0}%
\e{0}%
\e{0}%
\e{0}%
\e{0}%
\e{0}%
\e{0}%
\e{0}%
\e{0}%
\e{0}%
\e{0}%
\e{0}%
\e{0}%
\e{0}%
\e{0}%
\e{0}%
\e{0}%
\e{0}%
\e{1}%
\e{0}%
\eol}\vss}\rg%
%
%
\rx{\vss\hfull{%
\rlx{\hss{$3150_{y}$}}\cg%
\e{0}%
\e{0}%
\e{0}%
\e{0}%
\e{0}%
\e{1}%
\e{0}%
\e{0}%
\e{0}%
\e{0}%
\e{0}%
\e{0}%
\e{0}%
\e{1}%
\e{1}%
\e{0}%
\e{0}%
\e{0}%
\e{1}%
\e{1}%
\eol}\vss}\rg%
%
%
\rx{\vss\hfull{%
\rlx{\hss{$4200_{y}$}}\cg%
\e{0}%
\e{0}%
\e{0}%
\e{0}%
\e{0}%
\e{0}%
\e{0}%
\e{0}%
\e{0}%
\e{0}%
\e{0}%
\e{0}%
\e{1}%
\e{1}%
\e{1}%
\e{1}%
\e{1}%
\e{1}%
\e{1}%
\e{0}%
\eol}\vss}\rg%
\eop
\eject
\tablecont%
%
%
%
%
%
%
\rowpts=18 true pt%
\colpts=18 true pt%
\rowlabpts=40 true pt%
\collabpts=45 true pt%
\clx{\vss\hfull{%
\rlx{\hss{$ $}}\cg%
\cx{\hskip 16 true pt\flip{$189_c^{*}$}\hss}\cg%
\cx{\hskip 16 true pt\flip{$15_a^{*}$}\hss}\cg%
\cx{\hskip 16 true pt\flip{$105_b^{*}$}\hss}\cg%
\cx{\hskip 16 true pt\flip{$105_c^{*}$}\hss}\cg%
\cx{\hskip 16 true pt\flip{$315_a^{*}$}\hss}\cg%
\cx{\hskip 16 true pt\flip{$405_a^{*}$}\hss}\cg%
\cx{\hskip 16 true pt\flip{$168_a^{*}$}\hss}\cg%
\cx{\hskip 16 true pt\flip{$56_a^{*}$}\hss}\cg%
\cx{\hskip 16 true pt\flip{$120_a^{*}$}\hss}\cg%
\cx{\hskip 16 true pt\flip{$210_a^{*}$}\hss}\cg%
\cx{\hskip 16 true pt\flip{$280_a^{*}$}\hss}\cg%
\cx{\hskip 16 true pt\flip{$336_a^{*}$}\hss}\cg%
\cx{\hskip 16 true pt\flip{$216_a^{*}$}\hss}\cg%
\cx{\hskip 16 true pt\flip{$512_a^{*}$}\hss}\cg%
\cx{\hskip 16 true pt\flip{$378_a^{*}$}\hss}\cg%
\cx{\hskip 16 true pt\flip{$84_a^{*}$}\hss}\cg%
\cx{\hskip 16 true pt\flip{$420_a^{*}$}\hss}\cg%
\cx{\hskip 16 true pt\flip{$280_b^{*}$}\hss}\cg%
\cx{\hskip 16 true pt\flip{$210_b^{*}$}\hss}\cg%
\cx{\hskip 16 true pt\flip{$70_a^{*}$}\hss}\cg%
\eol}}\rg%
%
%
\rx{\vss\hfull{%
\rlx{\hss{$2688_{y}$}}\cg%
\e{0}%
\e{0}%
\e{0}%
\e{0}%
\e{0}%
\e{0}%
\e{0}%
\e{0}%
\e{0}%
\e{0}%
\e{0}%
\e{1}%
\e{1}%
\e{1}%
\e{0}%
\e{0}%
\e{0}%
\e{1}%
\e{0}%
\e{0}%
\eol}\vss}\rg%
%
%
\rx{\vss\hfull{%
\rlx{\hss{$2100_{y}$}}\cg%
\e{1}%
\e{0}%
\e{0}%
\e{1}%
\e{0}%
\e{0}%
\e{0}%
\e{0}%
\e{0}%
\e{0}%
\e{0}%
\e{1}%
\e{0}%
\e{0}%
\e{0}%
\e{0}%
\e{1}%
\e{0}%
\e{0}%
\e{0}%
\eol}\vss}\rg%
%
%
\rx{\vss\hfull{%
\rlx{\hss{$1400_{y}$}}\cg%
\e{0}%
\e{0}%
\e{0}%
\e{0}%
\e{0}%
\e{0}%
\e{0}%
\e{0}%
\e{0}%
\e{0}%
\e{1}%
\e{0}%
\e{0}%
\e{0}%
\e{0}%
\e{0}%
\e{1}%
\e{0}%
\e{0}%
\e{0}%
\eol}\vss}\rg%
%
%
\rx{\vss\hfull{%
\rlx{\hss{$4536_{y}$}}\cg%
\e{0}%
\e{0}%
\e{0}%
\e{0}%
\e{1}%
\e{1}%
\e{0}%
\e{0}%
\e{0}%
\e{0}%
\e{1}%
\e{1}%
\e{0}%
\e{1}%
\e{0}%
\e{0}%
\e{1}%
\e{0}%
\e{0}%
\e{0}%
\eol}\vss}\rg%
%
%
\rx{\vss\hfull{%
\rlx{\hss{$5670_{y}$}}\cg%
\e{0}%
\e{0}%
\e{0}%
\e{0}%
\e{1}%
\e{1}%
\e{0}%
\e{0}%
\e{0}%
\e{0}%
\e{1}%
\e{1}%
\e{0}%
\e{1}%
\e{1}%
\e{0}%
\e{1}%
\e{0}%
\e{0}%
\e{0}%
\eol}\vss}\rg%
%
%
\rx{\vss\hfull{%
\rlx{\hss{$4480_{y}$}}\cg%
\e{0}%
\e{0}%
\e{0}%
\e{0}%
\e{1}%
\e{1}%
\e{0}%
\e{0}%
\e{0}%
\e{0}%
\e{0}%
\e{0}%
\e{0}%
\e{1}%
\e{1}%
\e{0}%
\e{1}%
\e{0}%
\e{1}%
\e{0}%
\eol}\vss}\rg%
%
%
\rx{\vss\hfull{%
\rlx{\hss{$8_{z}$}}\cg%
\e{0}%
\e{0}%
\e{0}%
\e{0}%
\e{0}%
\e{0}%
\e{0}%
\e{0}%
\e{0}%
\e{0}%
\e{0}%
\e{0}%
\e{0}%
\e{0}%
\e{0}%
\e{0}%
\e{0}%
\e{0}%
\e{0}%
\e{0}%
\eol}\vss}\rg%
%
%
\rx{\vss\hfull{%
\rlx{\hss{$56_{z}$}}\cg%
\e{0}%
\e{0}%
\e{0}%
\e{0}%
\e{0}%
\e{0}%
\e{0}%
\e{0}%
\e{0}%
\e{0}%
\e{0}%
\e{0}%
\e{0}%
\e{0}%
\e{0}%
\e{0}%
\e{0}%
\e{0}%
\e{0}%
\e{0}%
\eol}\vss}\rg%
%
%
\rx{\vss\hfull{%
\rlx{\hss{$160_{z}$}}\cg%
\e{0}%
\e{0}%
\e{0}%
\e{0}%
\e{0}%
\e{0}%
\e{0}%
\e{0}%
\e{0}%
\e{0}%
\e{0}%
\e{0}%
\e{0}%
\e{0}%
\e{0}%
\e{0}%
\e{0}%
\e{0}%
\e{0}%
\e{0}%
\eol}\vss}\rg%
%
%
\rx{\vss\hfull{%
\rlx{\hss{$112_{z}$}}\cg%
\e{0}%
\e{0}%
\e{0}%
\e{0}%
\e{0}%
\e{0}%
\e{0}%
\e{1}%
\e{0}%
\e{0}%
\e{0}%
\e{0}%
\e{0}%
\e{0}%
\e{0}%
\e{0}%
\e{0}%
\e{0}%
\e{0}%
\e{0}%
\eol}\vss}\rg%
%
%
\rx{\vss\hfull{%
\rlx{\hss{$840_{z}$}}\cg%
\e{0}%
\e{0}%
\e{0}%
\e{0}%
\e{0}%
\e{0}%
\e{0}%
\e{0}%
\e{0}%
\e{0}%
\e{0}%
\e{0}%
\e{0}%
\e{0}%
\e{1}%
\e{0}%
\e{0}%
\e{0}%
\e{0}%
\e{0}%
\eol}\vss}\rg%
%
%
\rx{\vss\hfull{%
\rlx{\hss{$1296_{z}$}}\cg%
\e{0}%
\e{0}%
\e{0}%
\e{0}%
\e{0}%
\e{0}%
\e{0}%
\e{0}%
\e{0}%
\e{0}%
\e{1}%
\e{1}%
\e{0}%
\e{0}%
\e{0}%
\e{0}%
\e{0}%
\e{0}%
\e{0}%
\e{0}%
\eol}\vss}\rg%
%
%
\rx{\vss\hfull{%
\rlx{\hss{$1400_{z}$}}\cg%
\e{1}%
\e{0}%
\e{0}%
\e{0}%
\e{1}%
\e{0}%
\e{0}%
\e{1}%
\e{0}%
\e{0}%
\e{0}%
\e{0}%
\e{0}%
\e{0}%
\e{0}%
\e{0}%
\e{0}%
\e{0}%
\e{0}%
\e{0}%
\eol}\vss}\rg%
%
%
\rx{\vss\hfull{%
\rlx{\hss{$1008_{z}$}}\cg%
\e{1}%
\e{0}%
\e{0}%
\e{0}%
\e{0}%
\e{0}%
\e{0}%
\e{1}%
\e{0}%
\e{0}%
\e{1}%
\e{0}%
\e{0}%
\e{0}%
\e{0}%
\e{0}%
\e{0}%
\e{0}%
\e{0}%
\e{0}%
\eol}\vss}\rg%
%
%
\rx{\vss\hfull{%
\rlx{\hss{$560_{z}$}}\cg%
\e{0}%
\e{0}%
\e{0}%
\e{0}%
\e{0}%
\e{0}%
\e{0}%
\e{1}%
\e{0}%
\e{0}%
\e{0}%
\e{0}%
\e{0}%
\e{0}%
\e{0}%
\e{0}%
\e{0}%
\e{0}%
\e{0}%
\e{0}%
\eol}\vss}\rg%
%
%
\rx{\vss\hfull{%
\rlx{\hss{$1400_{zz}$}}\cg%
\e{0}%
\e{1}%
\e{0}%
\e{0}%
\e{1}%
\e{0}%
\e{0}%
\e{0}%
\e{0}%
\e{0}%
\e{0}%
\e{0}%
\e{1}%
\e{0}%
\e{0}%
\e{0}%
\e{0}%
\e{0}%
\e{0}%
\e{1}%
\eol}\vss}\rg%
%
%
\rx{\vss\hfull{%
\rlx{\hss{$4200_{z}$}}\cg%
\e{0}%
\e{0}%
\e{0}%
\e{0}%
\e{1}%
\e{0}%
\e{0}%
\e{0}%
\e{0}%
\e{0}%
\e{1}%
\e{0}%
\e{1}%
\e{1}%
\e{1}%
\e{1}%
\e{1}%
\e{0}%
\e{0}%
\e{0}%
\eol}\vss}\rg%
%
%
\rx{\vss\hfull{%
\rlx{\hss{$400_{z}$}}\cg%
\e{0}%
\e{1}%
\e{0}%
\e{0}%
\e{0}%
\e{0}%
\e{0}%
\e{1}%
\e{0}%
\e{0}%
\e{0}%
\e{0}%
\e{0}%
\e{0}%
\e{0}%
\e{0}%
\e{0}%
\e{0}%
\e{0}%
\e{0}%
\eol}\vss}\rg%
%
%
\rx{\vss\hfull{%
\rlx{\hss{$3240_{z}$}}\cg%
\e{1}%
\e{0}%
\e{0}%
\e{0}%
\e{1}%
\e{0}%
\e{0}%
\e{1}%
\e{0}%
\e{0}%
\e{1}%
\e{0}%
\e{1}%
\e{0}%
\e{1}%
\e{0}%
\e{0}%
\e{0}%
\e{0}%
\e{0}%
\eol}\vss}\rg%
%
%
\rx{\vss\hfull{%
\rlx{\hss{$4536_{z}$}}\cg%
\e{1}%
\e{0}%
\e{0}%
\e{0}%
\e{1}%
\e{0}%
\e{0}%
\e{0}%
\e{0}%
\e{0}%
\e{0}%
\e{0}%
\e{1}%
\e{1}%
\e{2}%
\e{0}%
\e{0}%
\e{0}%
\e{1}%
\e{1}%
\eol}\vss}\rg%
%
%
\rx{\vss\hfull{%
\rlx{\hss{$2400_{z}$}}\cg%
\e{0}%
\e{0}%
\e{0}%
\e{0}%
\e{0}%
\e{0}%
\e{0}%
\e{0}%
\e{0}%
\e{0}%
\e{1}%
\e{1}%
\e{0}%
\e{0}%
\e{0}%
\e{0}%
\e{1}%
\e{0}%
\e{0}%
\e{0}%
\eol}\vss}\rg%
%
%
\rx{\vss\hfull{%
\rlx{\hss{$3360_{z}$}}\cg%
\e{0}%
\e{0}%
\e{0}%
\e{0}%
\e{1}%
\e{0}%
\e{0}%
\e{0}%
\e{0}%
\e{0}%
\e{1}%
\e{1}%
\e{1}%
\e{1}%
\e{0}%
\e{0}%
\e{0}%
\e{0}%
\e{0}%
\e{0}%
\eol}\vss}\rg%
%
%
\rx{\vss\hfull{%
\rlx{\hss{$2800_{z}$}}\cg%
\e{1}%
\e{0}%
\e{0}%
\e{0}%
\e{1}%
\e{0}%
\e{0}%
\e{0}%
\e{0}%
\e{0}%
\e{1}%
\e{1}%
\e{0}%
\e{0}%
\e{0}%
\e{0}%
\e{1}%
\e{0}%
\e{0}%
\e{0}%
\eol}\vss}\rg%
%
%
\rx{\vss\hfull{%
\rlx{\hss{$4096_{z}$}}\cg%
\e{1}%
\e{0}%
\e{0}%
\e{0}%
\e{1}%
\e{0}%
\e{0}%
\e{0}%
\e{0}%
\e{0}%
\e{1}%
\e{1}%
\e{0}%
\e{1}%
\e{1}%
\e{0}%
\e{0}%
\e{0}%
\e{0}%
\e{0}%
\eol}\vss}\rg%
%
%
\rx{\vss\hfull{%
\rlx{\hss{$5600_{z}$}}\cg%
\e{1}%
\e{0}%
\e{0}%
\e{1}%
\e{1}%
\e{1}%
\e{0}%
\e{0}%
\e{0}%
\e{0}%
\e{1}%
\e{1}%
\e{0}%
\e{1}%
\e{1}%
\e{0}%
\e{1}%
\e{0}%
\e{0}%
\e{0}%
\eol}\vss}\rg%
%
%
\rx{\vss\hfull{%
\rlx{\hss{$448_{z}$}}\cg%
\e{0}%
\e{0}%
\e{0}%
\e{0}%
\e{0}%
\e{0}%
\e{0}%
\e{0}%
\e{0}%
\e{0}%
\e{0}%
\e{0}%
\e{0}%
\e{0}%
\e{0}%
\e{0}%
\e{0}%
\e{0}%
\e{0}%
\e{1}%
\eol}\vss}\rg%
%
%
\rx{\vss\hfull{%
\rlx{\hss{$448_{w}$}}\cg%
\e{0}%
\e{0}%
\e{0}%
\e{0}%
\e{0}%
\e{0}%
\e{0}%
\e{0}%
\e{0}%
\e{0}%
\e{0}%
\e{0}%
\e{0}%
\e{0}%
\e{0}%
\e{0}%
\e{0}%
\e{0}%
\e{0}%
\e{0}%
\eol}\vss}\rg%
%
%
\rx{\vss\hfull{%
\rlx{\hss{$1344_{w}$}}\cg%
\e{0}%
\e{0}%
\e{0}%
\e{0}%
\e{0}%
\e{0}%
\e{0}%
\e{0}%
\e{0}%
\e{0}%
\e{0}%
\e{0}%
\e{0}%
\e{0}%
\e{1}%
\e{1}%
\e{0}%
\e{0}%
\e{1}%
\e{0}%
\eol}\vss}\rg%
%
%
\rx{\vss\hfull{%
\rlx{\hss{$5600_{w}$}}\cg%
\e{0}%
\e{0}%
\e{0}%
\e{0}%
\e{0}%
\e{1}%
\e{0}%
\e{0}%
\e{0}%
\e{0}%
\e{1}%
\e{1}%
\e{0}%
\e{1}%
\e{1}%
\e{0}%
\e{1}%
\e{1}%
\e{0}%
\e{0}%
\eol}\vss}\rg%
%
%
\rx{\vss\hfull{%
\rlx{\hss{$2016_{w}$}}\cg%
\e{0}%
\e{0}%
\e{0}%
\e{0}%
\e{0}%
\e{0}%
\e{0}%
\e{0}%
\e{0}%
\e{0}%
\e{0}%
\e{0}%
\e{1}%
\e{1}%
\e{0}%
\e{0}%
\e{0}%
\e{0}%
\e{1}%
\e{1}%
\eol}\vss}\rg%
%
%
\rx{\vss\hfull{%
\rlx{\hss{$7168_{w}$}}\cg%
\e{0}%
\e{0}%
\e{0}%
\e{0}%
\e{1}%
\e{1}%
\e{0}%
\e{0}%
\e{0}%
\e{0}%
\e{0}%
\e{1}%
\e{1}%
\e{2}%
\e{1}%
\e{0}%
\e{1}%
\e{1}%
\e{1}%
\e{0}%
\eol}\vss}\rg%
\tableclose%
%
%
%
%
%
%
\eop
\eject
\tableopen{Induce/restrict matrix for $W({C_{3}}{A_{1}})\,\subset\,W(F_{4})$}%
%
%
%
%
%
%
\rowpts=18 true pt%
\colpts=18 true pt%
\rowlabpts=40 true pt%
\collabpts=80 true pt%
\clx{\vss\hfull{%
\rlx{\hss{$ $}}\cg%
\cx{\hskip 16 true pt\flip{$[{3}:-]{\times}[{2}]$}\hss}\cg%
\cx{\hskip 16 true pt\flip{$[{2}{1}:-]{\times}[{2}]$}\hss}\cg%
\cx{\hskip 16 true pt\flip{$[{1^{3}}:-]{\times}[{2}]$}\hss}\cg%
\cx{\hskip 16 true pt\flip{$[{2}:{1}]{\times}[{2}]$}\hss}\cg%
\cx{\hskip 16 true pt\flip{$[{1^{2}}:{1}]{\times}[{2}]$}\hss}\cg%
\cx{\hskip 16 true pt\flip{$[{1}:{2}]{\times}[{2}]$}\hss}\cg%
\cx{\hskip 16 true pt\flip{$[{1}:{1^{2}}]{\times}[{2}]$}\hss}\cg%
\cx{\hskip 16 true pt\flip{$[-:{3}]{\times}[{2}]$}\hss}\cg%
\cx{\hskip 16 true pt\flip{$[-:{2}{1}]{\times}[{2}]$}\hss}\cg%
\cx{\hskip 16 true pt\flip{$[-:{1^{3}}]{\times}[{2}]$}\hss}\cg%
\cx{\hskip 16 true pt\flip{$[{3}:-]{\times}[{1^{2}}]$}\hss}\cg%
\cx{\hskip 16 true pt\flip{$[{2}{1}:-]{\times}[{1^{2}}]$}\hss}\cg%
\cx{\hskip 16 true pt\flip{$[{1^{3}}:-]{\times}[{1^{2}}]$}\hss}\cg%
\cx{\hskip 16 true pt\flip{$[{2}:{1}]{\times}[{1^{2}}]$}\hss}\cg%
\cx{\hskip 16 true pt\flip{$[{1^{2}}:{1}]{\times}[{1^{2}}]$}\hss}\cg%
\cx{\hskip 16 true pt\flip{$[{1}:{2}]{\times}[{1^{2}}]$}\hss}\cg%
\cx{\hskip 16 true pt\flip{$[{1}:{1^{2}}]{\times}[{1^{2}}]$}\hss}\cg%
\cx{\hskip 16 true pt\flip{$[-:{3}]{\times}[{1^{2}}]$}\hss}\cg%
\cx{\hskip 16 true pt\flip{$[-:{2}{1}]{\times}[{1^{2}}]$}\hss}\cg%
\cx{\hskip 16 true pt\flip{$[-:{1^{3}}]{\times}[{1^{2}}]$}\hss}\cg%
\eol}}\rg%
%
%
\rx{\vss\hfull{%
\rlx{\hss{$\chi_{1,1}$}}\cg%
\e{1}%
\e{0}%
\e{0}%
\e{0}%
\e{0}%
\e{0}%
\e{0}%
\e{0}%
\e{0}%
\e{0}%
\e{0}%
\e{0}%
\e{0}%
\e{0}%
\e{0}%
\e{0}%
\e{0}%
\e{0}%
\e{0}%
\e{0}%
\eol}\vss}\rg%
%
%
\rx{\vss\hfull{%
\rlx{\hss{$\chi_{1,2}$}}\cg%
\e{0}%
\e{0}%
\e{1}%
\e{0}%
\e{0}%
\e{0}%
\e{0}%
\e{0}%
\e{0}%
\e{0}%
\e{0}%
\e{0}%
\e{0}%
\e{0}%
\e{0}%
\e{0}%
\e{0}%
\e{0}%
\e{0}%
\e{0}%
\eol}\vss}\rg%
%
%
\rx{\vss\hfull{%
\rlx{\hss{$\chi_{1,3}$}}\cg%
\e{0}%
\e{0}%
\e{0}%
\e{0}%
\e{0}%
\e{0}%
\e{0}%
\e{0}%
\e{0}%
\e{0}%
\e{0}%
\e{0}%
\e{0}%
\e{0}%
\e{0}%
\e{0}%
\e{0}%
\e{1}%
\e{0}%
\e{0}%
\eol}\vss}\rg%
%
%
\rx{\vss\hfull{%
\rlx{\hss{$\chi_{1,4}$}}\cg%
\e{0}%
\e{0}%
\e{0}%
\e{0}%
\e{0}%
\e{0}%
\e{0}%
\e{0}%
\e{0}%
\e{0}%
\e{0}%
\e{0}%
\e{0}%
\e{0}%
\e{0}%
\e{0}%
\e{0}%
\e{0}%
\e{0}%
\e{1}%
\eol}\vss}\rg%
%
%
\rx{\vss\hfull{%
\rlx{\hss{$\chi_{2,1}$}}\cg%
\e{0}%
\e{1}%
\e{0}%
\e{0}%
\e{0}%
\e{0}%
\e{0}%
\e{0}%
\e{0}%
\e{0}%
\e{0}%
\e{0}%
\e{0}%
\e{0}%
\e{0}%
\e{0}%
\e{0}%
\e{0}%
\e{0}%
\e{0}%
\eol}\vss}\rg%
%
%
\rx{\vss\hfull{%
\rlx{\hss{$\chi_{2,2}$}}\cg%
\e{0}%
\e{0}%
\e{0}%
\e{0}%
\e{0}%
\e{0}%
\e{0}%
\e{0}%
\e{0}%
\e{0}%
\e{0}%
\e{0}%
\e{0}%
\e{0}%
\e{0}%
\e{0}%
\e{0}%
\e{0}%
\e{1}%
\e{0}%
\eol}\vss}\rg%
%
%
\rx{\vss\hfull{%
\rlx{\hss{$\chi_{2,3}$}}\cg%
\e{1}%
\e{0}%
\e{0}%
\e{0}%
\e{0}%
\e{0}%
\e{0}%
\e{0}%
\e{0}%
\e{0}%
\e{0}%
\e{0}%
\e{0}%
\e{0}%
\e{0}%
\e{0}%
\e{0}%
\e{1}%
\e{0}%
\e{0}%
\eol}\vss}\rg%
%
%
\rx{\vss\hfull{%
\rlx{\hss{$\chi_{2,4}$}}\cg%
\e{0}%
\e{0}%
\e{1}%
\e{0}%
\e{0}%
\e{0}%
\e{0}%
\e{0}%
\e{0}%
\e{0}%
\e{0}%
\e{0}%
\e{0}%
\e{0}%
\e{0}%
\e{0}%
\e{0}%
\e{0}%
\e{0}%
\e{1}%
\eol}\vss}\rg%
%
%
\rx{\vss\hfull{%
\rlx{\hss{$\chi_{4,1}$}}\cg%
\e{0}%
\e{1}%
\e{0}%
\e{0}%
\e{0}%
\e{0}%
\e{0}%
\e{0}%
\e{0}%
\e{0}%
\e{0}%
\e{0}%
\e{0}%
\e{0}%
\e{0}%
\e{0}%
\e{0}%
\e{0}%
\e{1}%
\e{0}%
\eol}\vss}\rg%
%
%
\rx{\vss\hfull{%
\rlx{\hss{$\chi_{9,1}$}}\cg%
\e{1}%
\e{1}%
\e{0}%
\e{0}%
\e{0}%
\e{1}%
\e{0}%
\e{0}%
\e{0}%
\e{0}%
\e{0}%
\e{0}%
\e{0}%
\e{1}%
\e{0}%
\e{0}%
\e{0}%
\e{0}%
\e{0}%
\e{0}%
\eol}\vss}\rg%
%
%
\rx{\vss\hfull{%
\rlx{\hss{$\chi_{9,2}$}}\cg%
\e{0}%
\e{1}%
\e{1}%
\e{0}%
\e{0}%
\e{0}%
\e{1}%
\e{0}%
\e{0}%
\e{0}%
\e{0}%
\e{0}%
\e{0}%
\e{0}%
\e{1}%
\e{0}%
\e{0}%
\e{0}%
\e{0}%
\e{0}%
\eol}\vss}\rg%
%
%
\rx{\vss\hfull{%
\rlx{\hss{$\chi_{9,3}$}}\cg%
\e{0}%
\e{0}%
\e{0}%
\e{0}%
\e{0}%
\e{1}%
\e{0}%
\e{0}%
\e{0}%
\e{0}%
\e{0}%
\e{0}%
\e{0}%
\e{1}%
\e{0}%
\e{0}%
\e{0}%
\e{1}%
\e{1}%
\e{0}%
\eol}\vss}\rg%
%
%
\rx{\vss\hfull{%
\rlx{\hss{$\chi_{9,4}$}}\cg%
\e{0}%
\e{0}%
\e{0}%
\e{0}%
\e{0}%
\e{0}%
\e{1}%
\e{0}%
\e{0}%
\e{0}%
\e{0}%
\e{0}%
\e{0}%
\e{0}%
\e{1}%
\e{0}%
\e{0}%
\e{0}%
\e{1}%
\e{1}%
\eol}\vss}\rg%
%
%
\rx{\vss\hfull{%
\rlx{\hss{$\chi_{6,1}$}}\cg%
\e{0}%
\e{0}%
\e{0}%
\e{0}%
\e{0}%
\e{1}%
\e{0}%
\e{0}%
\e{0}%
\e{0}%
\e{0}%
\e{0}%
\e{0}%
\e{0}%
\e{1}%
\e{0}%
\e{0}%
\e{0}%
\e{0}%
\e{0}%
\eol}\vss}\rg%
%
%
\rx{\vss\hfull{%
\rlx{\hss{$\chi_{6,2}$}}\cg%
\e{0}%
\e{0}%
\e{0}%
\e{0}%
\e{0}%
\e{0}%
\e{1}%
\e{0}%
\e{0}%
\e{0}%
\e{0}%
\e{0}%
\e{0}%
\e{1}%
\e{0}%
\e{0}%
\e{0}%
\e{0}%
\e{0}%
\e{0}%
\eol}\vss}\rg%
%
%
\rx{\vss\hfull{%
\rlx{\hss{$\chi_{12,1}$}}\cg%
\e{0}%
\e{0}%
\e{0}%
\e{0}%
\e{0}%
\e{1}%
\e{1}%
\e{0}%
\e{0}%
\e{0}%
\e{0}%
\e{0}%
\e{0}%
\e{1}%
\e{1}%
\e{0}%
\e{0}%
\e{0}%
\e{0}%
\e{0}%
\eol}\vss}\rg%
%
%
\rx{\vss\hfull{%
\rlx{\hss{$\chi_{4,2}$}}\cg%
\e{0}%
\e{0}%
\e{0}%
\e{1}%
\e{0}%
\e{0}%
\e{0}%
\e{0}%
\e{0}%
\e{0}%
\e{1}%
\e{0}%
\e{0}%
\e{0}%
\e{0}%
\e{0}%
\e{0}%
\e{0}%
\e{0}%
\e{0}%
\eol}\vss}\rg%
%
%
\rx{\vss\hfull{%
\rlx{\hss{$\chi_{4,3}$}}\cg%
\e{0}%
\e{0}%
\e{0}%
\e{0}%
\e{1}%
\e{0}%
\e{0}%
\e{0}%
\e{0}%
\e{0}%
\e{0}%
\e{0}%
\e{1}%
\e{0}%
\e{0}%
\e{0}%
\e{0}%
\e{0}%
\e{0}%
\e{0}%
\eol}\vss}\rg%
%
%
\rx{\vss\hfull{%
\rlx{\hss{$\chi_{4,4}$}}\cg%
\e{0}%
\e{0}%
\e{0}%
\e{0}%
\e{0}%
\e{0}%
\e{0}%
\e{1}%
\e{0}%
\e{0}%
\e{0}%
\e{0}%
\e{0}%
\e{0}%
\e{0}%
\e{1}%
\e{0}%
\e{0}%
\e{0}%
\e{0}%
\eol}\vss}\rg%
%
%
\rx{\vss\hfull{%
\rlx{\hss{$\chi_{4,5}$}}\cg%
\e{0}%
\e{0}%
\e{0}%
\e{0}%
\e{0}%
\e{0}%
\e{0}%
\e{0}%
\e{0}%
\e{1}%
\e{0}%
\e{0}%
\e{0}%
\e{0}%
\e{0}%
\e{0}%
\e{1}%
\e{0}%
\e{0}%
\e{0}%
\eol}\vss}\rg%
%
%
\rx{\vss\hfull{%
\rlx{\hss{$\chi_{8,1}$}}\cg%
\e{0}%
\e{0}%
\e{0}%
\e{1}%
\e{1}%
\e{0}%
\e{0}%
\e{0}%
\e{0}%
\e{0}%
\e{0}%
\e{1}%
\e{0}%
\e{0}%
\e{0}%
\e{0}%
\e{0}%
\e{0}%
\e{0}%
\e{0}%
\eol}\vss}\rg%
%
%
\rx{\vss\hfull{%
\rlx{\hss{$\chi_{8,2}$}}\cg%
\e{0}%
\e{0}%
\e{0}%
\e{0}%
\e{0}%
\e{0}%
\e{0}%
\e{0}%
\e{1}%
\e{0}%
\e{0}%
\e{0}%
\e{0}%
\e{0}%
\e{0}%
\e{1}%
\e{1}%
\e{0}%
\e{0}%
\e{0}%
\eol}\vss}\rg%
%
%
\rx{\vss\hfull{%
\rlx{\hss{$\chi_{8,3}$}}\cg%
\e{0}%
\e{0}%
\e{0}%
\e{1}%
\e{0}%
\e{0}%
\e{0}%
\e{1}%
\e{0}%
\e{0}%
\e{1}%
\e{0}%
\e{0}%
\e{0}%
\e{0}%
\e{1}%
\e{0}%
\e{0}%
\e{0}%
\e{0}%
\eol}\vss}\rg%
%
%
\rx{\vss\hfull{%
\rlx{\hss{$\chi_{8,4}$}}\cg%
\e{0}%
\e{0}%
\e{0}%
\e{0}%
\e{1}%
\e{0}%
\e{0}%
\e{0}%
\e{0}%
\e{1}%
\e{0}%
\e{0}%
\e{1}%
\e{0}%
\e{0}%
\e{0}%
\e{1}%
\e{0}%
\e{0}%
\e{0}%
\eol}\vss}\rg%
%
%
\rx{\vss\hfull{%
\rlx{\hss{$\chi_{16,1}$}}\cg%
\e{0}%
\e{0}%
\e{0}%
\e{1}%
\e{1}%
\e{0}%
\e{0}%
\e{0}%
\e{1}%
\e{0}%
\e{0}%
\e{1}%
\e{0}%
\e{0}%
\e{0}%
\e{1}%
\e{1}%
\e{0}%
\e{0}%
\e{0}%
\eol}\vss}\rg%
\tableclose%
%
%
%
%
%
%
\eop
\eject
\tableopen{Induce/restrict matrix for $W({A_{2}}{A_{2}^{*}})\,\subset\,W(F_{4})$}%
%
%
%
%
%
%
\rowpts=18 true pt%
\colpts=18 true pt%
\rowlabpts=40 true pt%
\collabpts=60 true pt%
\clx{\vss\hfull{%
\rlx{\hss{$ $}}\cg%
\cx{\hskip 16 true pt\flip{$[{3}]{\times}[{3}]$}\hss}\cg%
\cx{\hskip 16 true pt\flip{$[{2}{1}]{\times}[{3}]$}\hss}\cg%
\cx{\hskip 16 true pt\flip{$[{1^{3}}]{\times}[{3}]$}\hss}\cg%
\cx{\hskip 16 true pt\flip{$[{3}]{\times}[{2}{1}]$}\hss}\cg%
\cx{\hskip 16 true pt\flip{$[{2}{1}]{\times}[{2}{1}]$}\hss}\cg%
\cx{\hskip 16 true pt\flip{$[{1^{3}}]{\times}[{2}{1}]$}\hss}\cg%
\cx{\hskip 16 true pt\flip{$[{3}]{\times}[{1^{3}}]$}\hss}\cg%
\cx{\hskip 16 true pt\flip{$[{2}{1}]{\times}[{1^{3}}]$}\hss}\cg%
\cx{\hskip 16 true pt\flip{$[{1^{3}}]{\times}[{1^{3}}]$}\hss}\cg%
\eol}}\rg%
%
%
\rx{\vss\hfull{%
\rlx{\hss{$\chi_{1,1}$}}\cg%
\e{1}%
\e{0}%
\e{0}%
\e{0}%
\e{0}%
\e{0}%
\e{0}%
\e{0}%
\e{0}%
\eol}\vss}\rg%
%
%
\rx{\vss\hfull{%
\rlx{\hss{$\chi_{1,2}$}}\cg%
\e{0}%
\e{0}%
\e{0}%
\e{0}%
\e{0}%
\e{0}%
\e{1}%
\e{0}%
\e{0}%
\eol}\vss}\rg%
%
%
\rx{\vss\hfull{%
\rlx{\hss{$\chi_{1,3}$}}\cg%
\e{0}%
\e{0}%
\e{1}%
\e{0}%
\e{0}%
\e{0}%
\e{0}%
\e{0}%
\e{0}%
\eol}\vss}\rg%
%
%
\rx{\vss\hfull{%
\rlx{\hss{$\chi_{1,4}$}}\cg%
\e{0}%
\e{0}%
\e{0}%
\e{0}%
\e{0}%
\e{0}%
\e{0}%
\e{0}%
\e{1}%
\eol}\vss}\rg%
%
%
\rx{\vss\hfull{%
\rlx{\hss{$\chi_{2,1}$}}\cg%
\e{0}%
\e{0}%
\e{0}%
\e{1}%
\e{0}%
\e{0}%
\e{0}%
\e{0}%
\e{0}%
\eol}\vss}\rg%
%
%
\rx{\vss\hfull{%
\rlx{\hss{$\chi_{2,2}$}}\cg%
\e{0}%
\e{0}%
\e{0}%
\e{0}%
\e{0}%
\e{1}%
\e{0}%
\e{0}%
\e{0}%
\eol}\vss}\rg%
%
%
\rx{\vss\hfull{%
\rlx{\hss{$\chi_{2,3}$}}\cg%
\e{0}%
\e{1}%
\e{0}%
\e{0}%
\e{0}%
\e{0}%
\e{0}%
\e{0}%
\e{0}%
\eol}\vss}\rg%
%
%
\rx{\vss\hfull{%
\rlx{\hss{$\chi_{2,4}$}}\cg%
\e{0}%
\e{0}%
\e{0}%
\e{0}%
\e{0}%
\e{0}%
\e{0}%
\e{1}%
\e{0}%
\eol}\vss}\rg%
%
%
\rx{\vss\hfull{%
\rlx{\hss{$\chi_{4,1}$}}\cg%
\e{0}%
\e{0}%
\e{0}%
\e{0}%
\e{1}%
\e{0}%
\e{0}%
\e{0}%
\e{0}%
\eol}\vss}\rg%
%
%
\rx{\vss\hfull{%
\rlx{\hss{$\chi_{9,1}$}}\cg%
\e{1}%
\e{1}%
\e{0}%
\e{1}%
\e{1}%
\e{0}%
\e{0}%
\e{0}%
\e{0}%
\eol}\vss}\rg%
%
%
\rx{\vss\hfull{%
\rlx{\hss{$\chi_{9,2}$}}\cg%
\e{0}%
\e{0}%
\e{0}%
\e{1}%
\e{1}%
\e{0}%
\e{1}%
\e{1}%
\e{0}%
\eol}\vss}\rg%
%
%
\rx{\vss\hfull{%
\rlx{\hss{$\chi_{9,3}$}}\cg%
\e{0}%
\e{1}%
\e{1}%
\e{0}%
\e{1}%
\e{1}%
\e{0}%
\e{0}%
\e{0}%
\eol}\vss}\rg%
%
%
\rx{\vss\hfull{%
\rlx{\hss{$\chi_{9,4}$}}\cg%
\e{0}%
\e{0}%
\e{0}%
\e{0}%
\e{1}%
\e{1}%
\e{0}%
\e{1}%
\e{1}%
\eol}\vss}\rg%
%
%
\rx{\vss\hfull{%
\rlx{\hss{$\chi_{6,1}$}}\cg%
\e{1}%
\e{0}%
\e{0}%
\e{0}%
\e{1}%
\e{0}%
\e{0}%
\e{0}%
\e{1}%
\eol}\vss}\rg%
%
%
\rx{\vss\hfull{%
\rlx{\hss{$\chi_{6,2}$}}\cg%
\e{0}%
\e{0}%
\e{1}%
\e{0}%
\e{1}%
\e{0}%
\e{1}%
\e{0}%
\e{0}%
\eol}\vss}\rg%
%
%
\rx{\vss\hfull{%
\rlx{\hss{$\chi_{12,1}$}}\cg%
\e{0}%
\e{1}%
\e{0}%
\e{1}%
\e{1}%
\e{1}%
\e{0}%
\e{1}%
\e{0}%
\eol}\vss}\rg%
%
%
\rx{\vss\hfull{%
\rlx{\hss{$\chi_{4,2}$}}\cg%
\e{0}%
\e{1}%
\e{0}%
\e{1}%
\e{0}%
\e{0}%
\e{0}%
\e{0}%
\e{0}%
\eol}\vss}\rg%
%
%
\rx{\vss\hfull{%
\rlx{\hss{$\chi_{4,3}$}}\cg%
\e{0}%
\e{0}%
\e{0}%
\e{1}%
\e{0}%
\e{0}%
\e{0}%
\e{1}%
\e{0}%
\eol}\vss}\rg%
%
%
\rx{\vss\hfull{%
\rlx{\hss{$\chi_{4,4}$}}\cg%
\e{0}%
\e{1}%
\e{0}%
\e{0}%
\e{0}%
\e{1}%
\e{0}%
\e{0}%
\e{0}%
\eol}\vss}\rg%
%
%
\rx{\vss\hfull{%
\rlx{\hss{$\chi_{4,5}$}}\cg%
\e{0}%
\e{0}%
\e{0}%
\e{0}%
\e{0}%
\e{1}%
\e{0}%
\e{1}%
\e{0}%
\eol}\vss}\rg%
%
%
\rx{\vss\hfull{%
\rlx{\hss{$\chi_{8,1}$}}\cg%
\e{1}%
\e{0}%
\e{0}%
\e{1}%
\e{1}%
\e{0}%
\e{1}%
\e{0}%
\e{0}%
\eol}\vss}\rg%
%
%
\rx{\vss\hfull{%
\rlx{\hss{$\chi_{8,2}$}}\cg%
\e{0}%
\e{0}%
\e{1}%
\e{0}%
\e{1}%
\e{1}%
\e{0}%
\e{0}%
\e{1}%
\eol}\vss}\rg%
%
%
\rx{\vss\hfull{%
\rlx{\hss{$\chi_{8,3}$}}\cg%
\e{1}%
\e{1}%
\e{1}%
\e{0}%
\e{1}%
\e{0}%
\e{0}%
\e{0}%
\e{0}%
\eol}\vss}\rg%
%
%
\rx{\vss\hfull{%
\rlx{\hss{$\chi_{8,4}$}}\cg%
\e{0}%
\e{0}%
\e{0}%
\e{0}%
\e{1}%
\e{0}%
\e{1}%
\e{1}%
\e{1}%
\eol}\vss}\rg%
%
%
\rx{\vss\hfull{%
\rlx{\hss{$\chi_{16,1}$}}\cg%
\e{0}%
\e{1}%
\e{0}%
\e{1}%
\e{2}%
\e{1}%
\e{0}%
\e{1}%
\e{0}%
\eol}\vss}\rg%
\tableclose%
%
%
%
%
%
%
\eop
\eject
\tableopen{Induce/restrict matrix for $W({A_{3}}{A_{1}})\,\subset\,W(F_{4})$}%
%
%
%
%
%
%
\rowpts=18 true pt%
\colpts=18 true pt%
\rowlabpts=40 true pt%
\collabpts=60 true pt%
\clx{\vss\hfull{%
\rlx{\hss{$ $}}\cg%
\cx{\hskip 16 true pt\flip{$[{4}]{\times}[{2}]$}\hss}\cg%
\cx{\hskip 16 true pt\flip{$[{3}{1}]{\times}[{2}]$}\hss}\cg%
\cx{\hskip 16 true pt\flip{$[{2^{2}}]{\times}[{2}]$}\hss}\cg%
\cx{\hskip 16 true pt\flip{$[{2}{1^{2}}]{\times}[{2}]$}\hss}\cg%
\cx{\hskip 16 true pt\flip{$[{1^{4}}]{\times}[{2}]$}\hss}\cg%
\cx{\hskip 16 true pt\flip{$[{4}]{\times}[{1^{2}}]$}\hss}\cg%
\cx{\hskip 16 true pt\flip{$[{3}{1}]{\times}[{1^{2}}]$}\hss}\cg%
\cx{\hskip 16 true pt\flip{$[{2^{2}}]{\times}[{1^{2}}]$}\hss}\cg%
\cx{\hskip 16 true pt\flip{$[{2}{1^{2}}]{\times}[{1^{2}}]$}\hss}\cg%
\cx{\hskip 16 true pt\flip{$[{1^{4}}]{\times}[{1^{2}}]$}\hss}\cg%
\eol}}\rg%
%
%
\rx{\vss\hfull{%
\rlx{\hss{$\chi_{1,1}$}}\cg%
\e{1}%
\e{0}%
\e{0}%
\e{0}%
\e{0}%
\e{0}%
\e{0}%
\e{0}%
\e{0}%
\e{0}%
\eol}\vss}\rg%
%
%
\rx{\vss\hfull{%
\rlx{\hss{$\chi_{1,2}$}}\cg%
\e{0}%
\e{0}%
\e{0}%
\e{0}%
\e{0}%
\e{1}%
\e{0}%
\e{0}%
\e{0}%
\e{0}%
\eol}\vss}\rg%
%
%
\rx{\vss\hfull{%
\rlx{\hss{$\chi_{1,3}$}}\cg%
\e{0}%
\e{0}%
\e{0}%
\e{0}%
\e{1}%
\e{0}%
\e{0}%
\e{0}%
\e{0}%
\e{0}%
\eol}\vss}\rg%
%
%
\rx{\vss\hfull{%
\rlx{\hss{$\chi_{1,4}$}}\cg%
\e{0}%
\e{0}%
\e{0}%
\e{0}%
\e{0}%
\e{0}%
\e{0}%
\e{0}%
\e{0}%
\e{1}%
\eol}\vss}\rg%
%
%
\rx{\vss\hfull{%
\rlx{\hss{$\chi_{2,1}$}}\cg%
\e{1}%
\e{0}%
\e{0}%
\e{0}%
\e{0}%
\e{1}%
\e{0}%
\e{0}%
\e{0}%
\e{0}%
\eol}\vss}\rg%
%
%
\rx{\vss\hfull{%
\rlx{\hss{$\chi_{2,2}$}}\cg%
\e{0}%
\e{0}%
\e{0}%
\e{0}%
\e{1}%
\e{0}%
\e{0}%
\e{0}%
\e{0}%
\e{1}%
\eol}\vss}\rg%
%
%
\rx{\vss\hfull{%
\rlx{\hss{$\chi_{2,3}$}}\cg%
\e{0}%
\e{0}%
\e{1}%
\e{0}%
\e{0}%
\e{0}%
\e{0}%
\e{0}%
\e{0}%
\e{0}%
\eol}\vss}\rg%
%
%
\rx{\vss\hfull{%
\rlx{\hss{$\chi_{2,4}$}}\cg%
\e{0}%
\e{0}%
\e{0}%
\e{0}%
\e{0}%
\e{0}%
\e{0}%
\e{1}%
\e{0}%
\e{0}%
\eol}\vss}\rg%
%
%
\rx{\vss\hfull{%
\rlx{\hss{$\chi_{4,1}$}}\cg%
\e{0}%
\e{0}%
\e{1}%
\e{0}%
\e{0}%
\e{0}%
\e{0}%
\e{1}%
\e{0}%
\e{0}%
\eol}\vss}\rg%
%
%
\rx{\vss\hfull{%
\rlx{\hss{$\chi_{9,1}$}}\cg%
\e{1}%
\e{1}%
\e{1}%
\e{0}%
\e{0}%
\e{0}%
\e{1}%
\e{0}%
\e{0}%
\e{0}%
\eol}\vss}\rg%
%
%
\rx{\vss\hfull{%
\rlx{\hss{$\chi_{9,2}$}}\cg%
\e{0}%
\e{1}%
\e{0}%
\e{0}%
\e{0}%
\e{1}%
\e{1}%
\e{1}%
\e{0}%
\e{0}%
\eol}\vss}\rg%
%
%
\rx{\vss\hfull{%
\rlx{\hss{$\chi_{9,3}$}}\cg%
\e{0}%
\e{0}%
\e{1}%
\e{1}%
\e{1}%
\e{0}%
\e{0}%
\e{0}%
\e{1}%
\e{0}%
\eol}\vss}\rg%
%
%
\rx{\vss\hfull{%
\rlx{\hss{$\chi_{9,4}$}}\cg%
\e{0}%
\e{0}%
\e{0}%
\e{1}%
\e{0}%
\e{0}%
\e{0}%
\e{1}%
\e{1}%
\e{1}%
\eol}\vss}\rg%
%
%
\rx{\vss\hfull{%
\rlx{\hss{$\chi_{6,1}$}}\cg%
\e{0}%
\e{1}%
\e{0}%
\e{0}%
\e{0}%
\e{0}%
\e{0}%
\e{0}%
\e{1}%
\e{0}%
\eol}\vss}\rg%
%
%
\rx{\vss\hfull{%
\rlx{\hss{$\chi_{6,2}$}}\cg%
\e{0}%
\e{0}%
\e{0}%
\e{1}%
\e{0}%
\e{0}%
\e{1}%
\e{0}%
\e{0}%
\e{0}%
\eol}\vss}\rg%
%
%
\rx{\vss\hfull{%
\rlx{\hss{$\chi_{12,1}$}}\cg%
\e{0}%
\e{1}%
\e{0}%
\e{1}%
\e{0}%
\e{0}%
\e{1}%
\e{0}%
\e{1}%
\e{0}%
\eol}\vss}\rg%
%
%
\rx{\vss\hfull{%
\rlx{\hss{$\chi_{4,2}$}}\cg%
\e{0}%
\e{1}%
\e{0}%
\e{0}%
\e{0}%
\e{1}%
\e{0}%
\e{0}%
\e{0}%
\e{0}%
\eol}\vss}\rg%
%
%
\rx{\vss\hfull{%
\rlx{\hss{$\chi_{4,3}$}}\cg%
\e{1}%
\e{0}%
\e{0}%
\e{0}%
\e{0}%
\e{0}%
\e{1}%
\e{0}%
\e{0}%
\e{0}%
\eol}\vss}\rg%
%
%
\rx{\vss\hfull{%
\rlx{\hss{$\chi_{4,4}$}}\cg%
\e{0}%
\e{0}%
\e{0}%
\e{1}%
\e{0}%
\e{0}%
\e{0}%
\e{0}%
\e{0}%
\e{1}%
\eol}\vss}\rg%
%
%
\rx{\vss\hfull{%
\rlx{\hss{$\chi_{4,5}$}}\cg%
\e{0}%
\e{0}%
\e{0}%
\e{0}%
\e{1}%
\e{0}%
\e{0}%
\e{0}%
\e{1}%
\e{0}%
\eol}\vss}\rg%
%
%
\rx{\vss\hfull{%
\rlx{\hss{$\chi_{8,1}$}}\cg%
\e{1}%
\e{1}%
\e{0}%
\e{0}%
\e{0}%
\e{1}%
\e{1}%
\e{0}%
\e{0}%
\e{0}%
\eol}\vss}\rg%
%
%
\rx{\vss\hfull{%
\rlx{\hss{$\chi_{8,2}$}}\cg%
\e{0}%
\e{0}%
\e{0}%
\e{1}%
\e{1}%
\e{0}%
\e{0}%
\e{0}%
\e{1}%
\e{1}%
\eol}\vss}\rg%
%
%
\rx{\vss\hfull{%
\rlx{\hss{$\chi_{8,3}$}}\cg%
\e{0}%
\e{1}%
\e{0}%
\e{1}%
\e{0}%
\e{0}%
\e{0}%
\e{1}%
\e{0}%
\e{0}%
\eol}\vss}\rg%
%
%
\rx{\vss\hfull{%
\rlx{\hss{$\chi_{8,4}$}}\cg%
\e{0}%
\e{0}%
\e{1}%
\e{0}%
\e{0}%
\e{0}%
\e{1}%
\e{0}%
\e{1}%
\e{0}%
\eol}\vss}\rg%
%
%
\rx{\vss\hfull{%
\rlx{\hss{$\chi_{16,1}$}}\cg%
\e{0}%
\e{1}%
\e{1}%
\e{1}%
\e{0}%
\e{0}%
\e{1}%
\e{1}%
\e{1}%
\e{0}%
\eol}\vss}\rg%
\tableclose%
%
%
%
%
%
%
\eop
\eject
\tableopen{Induce/restrict matrix for $W(B_{4})\,\subset\,W(F_{4})$}%
%
%
%
%
%
%
\rowpts=18 true pt%
\colpts=18 true pt%
\rowlabpts=40 true pt%
\collabpts=60 true pt%
\clx{\vss\hfull{%
\rlx{\hss{$ $}}\cg%
\cx{\hskip 16 true pt\flip{$[{4}:-]$}\hss}\cg%
\cx{\hskip 16 true pt\flip{$[{3}{1}:-]$}\hss}\cg%
\cx{\hskip 16 true pt\flip{$[{2^{2}}:-]$}\hss}\cg%
\cx{\hskip 16 true pt\flip{$[{2}{1^{2}}:-]$}\hss}\cg%
\cx{\hskip 16 true pt\flip{$[{1^{4}}:-]$}\hss}\cg%
\cx{\hskip 16 true pt\flip{$[{3}:{1}]$}\hss}\cg%
\cx{\hskip 16 true pt\flip{$[{2}{1}:{1}]$}\hss}\cg%
\cx{\hskip 16 true pt\flip{$[{1^{3}}:{1}]$}\hss}\cg%
\cx{\hskip 16 true pt\flip{$[{2}:{2}]$}\hss}\cg%
\cx{\hskip 16 true pt\flip{$[{2}:{1^{2}}]$}\hss}\cg%
\cx{\hskip 16 true pt\flip{$[{1^{2}}:{2}]$}\hss}\cg%
\cx{\hskip 16 true pt\flip{$[{1^{2}}:{1^{2}}]$}\hss}\cg%
\cx{\hskip 16 true pt\flip{$[{1}:{3}]$}\hss}\cg%
\cx{\hskip 16 true pt\flip{$[{1}:{2}{1}]$}\hss}\cg%
\cx{\hskip 16 true pt\flip{$[{1}:{1^{3}}]$}\hss}\cg%
\cx{\hskip 16 true pt\flip{$[-:{4}]$}\hss}\cg%
\cx{\hskip 16 true pt\flip{$[-:{3}{1}]$}\hss}\cg%
\cx{\hskip 16 true pt\flip{$[-:{2^{2}}]$}\hss}\cg%
\cx{\hskip 16 true pt\flip{$[-:{2}{1^{2}}]$}\hss}\cg%
\cx{\hskip 16 true pt\flip{$[-:{1^{4}}]$}\hss}\cg%
\eol}}\rg%
%
%
\rx{\vss\hfull{%
\rlx{\hss{$\chi_{1,1}$}}\cg%
\e{1}%
\e{0}%
\e{0}%
\e{0}%
\e{0}%
\e{0}%
\e{0}%
\e{0}%
\e{0}%
\e{0}%
\e{0}%
\e{0}%
\e{0}%
\e{0}%
\e{0}%
\e{0}%
\e{0}%
\e{0}%
\e{0}%
\e{0}%
\eol}\vss}\rg%
%
%
\rx{\vss\hfull{%
\rlx{\hss{$\chi_{1,2}$}}\cg%
\e{0}%
\e{0}%
\e{0}%
\e{0}%
\e{0}%
\e{0}%
\e{0}%
\e{0}%
\e{0}%
\e{0}%
\e{0}%
\e{0}%
\e{0}%
\e{0}%
\e{0}%
\e{1}%
\e{0}%
\e{0}%
\e{0}%
\e{0}%
\eol}\vss}\rg%
%
%
\rx{\vss\hfull{%
\rlx{\hss{$\chi_{1,3}$}}\cg%
\e{0}%
\e{0}%
\e{0}%
\e{0}%
\e{1}%
\e{0}%
\e{0}%
\e{0}%
\e{0}%
\e{0}%
\e{0}%
\e{0}%
\e{0}%
\e{0}%
\e{0}%
\e{0}%
\e{0}%
\e{0}%
\e{0}%
\e{0}%
\eol}\vss}\rg%
%
%
\rx{\vss\hfull{%
\rlx{\hss{$\chi_{1,4}$}}\cg%
\e{0}%
\e{0}%
\e{0}%
\e{0}%
\e{0}%
\e{0}%
\e{0}%
\e{0}%
\e{0}%
\e{0}%
\e{0}%
\e{0}%
\e{0}%
\e{0}%
\e{0}%
\e{0}%
\e{0}%
\e{0}%
\e{0}%
\e{1}%
\eol}\vss}\rg%
%
%
\rx{\vss\hfull{%
\rlx{\hss{$\chi_{2,1}$}}\cg%
\e{1}%
\e{0}%
\e{0}%
\e{0}%
\e{0}%
\e{0}%
\e{0}%
\e{0}%
\e{0}%
\e{0}%
\e{0}%
\e{0}%
\e{0}%
\e{0}%
\e{0}%
\e{1}%
\e{0}%
\e{0}%
\e{0}%
\e{0}%
\eol}\vss}\rg%
%
%
\rx{\vss\hfull{%
\rlx{\hss{$\chi_{2,2}$}}\cg%
\e{0}%
\e{0}%
\e{0}%
\e{0}%
\e{1}%
\e{0}%
\e{0}%
\e{0}%
\e{0}%
\e{0}%
\e{0}%
\e{0}%
\e{0}%
\e{0}%
\e{0}%
\e{0}%
\e{0}%
\e{0}%
\e{0}%
\e{1}%
\eol}\vss}\rg%
%
%
\rx{\vss\hfull{%
\rlx{\hss{$\chi_{2,3}$}}\cg%
\e{0}%
\e{0}%
\e{1}%
\e{0}%
\e{0}%
\e{0}%
\e{0}%
\e{0}%
\e{0}%
\e{0}%
\e{0}%
\e{0}%
\e{0}%
\e{0}%
\e{0}%
\e{0}%
\e{0}%
\e{0}%
\e{0}%
\e{0}%
\eol}\vss}\rg%
%
%
\rx{\vss\hfull{%
\rlx{\hss{$\chi_{2,4}$}}\cg%
\e{0}%
\e{0}%
\e{0}%
\e{0}%
\e{0}%
\e{0}%
\e{0}%
\e{0}%
\e{0}%
\e{0}%
\e{0}%
\e{0}%
\e{0}%
\e{0}%
\e{0}%
\e{0}%
\e{0}%
\e{1}%
\e{0}%
\e{0}%
\eol}\vss}\rg%
%
%
\rx{\vss\hfull{%
\rlx{\hss{$\chi_{4,1}$}}\cg%
\e{0}%
\e{0}%
\e{1}%
\e{0}%
\e{0}%
\e{0}%
\e{0}%
\e{0}%
\e{0}%
\e{0}%
\e{0}%
\e{0}%
\e{0}%
\e{0}%
\e{0}%
\e{0}%
\e{0}%
\e{1}%
\e{0}%
\e{0}%
\eol}\vss}\rg%
%
%
\rx{\vss\hfull{%
\rlx{\hss{$\chi_{9,1}$}}\cg%
\e{0}%
\e{1}%
\e{0}%
\e{0}%
\e{0}%
\e{0}%
\e{0}%
\e{0}%
\e{1}%
\e{0}%
\e{0}%
\e{0}%
\e{0}%
\e{0}%
\e{0}%
\e{0}%
\e{0}%
\e{0}%
\e{0}%
\e{0}%
\eol}\vss}\rg%
%
%
\rx{\vss\hfull{%
\rlx{\hss{$\chi_{9,2}$}}\cg%
\e{0}%
\e{0}%
\e{0}%
\e{0}%
\e{0}%
\e{0}%
\e{0}%
\e{0}%
\e{1}%
\e{0}%
\e{0}%
\e{0}%
\e{0}%
\e{0}%
\e{0}%
\e{0}%
\e{1}%
\e{0}%
\e{0}%
\e{0}%
\eol}\vss}\rg%
%
%
\rx{\vss\hfull{%
\rlx{\hss{$\chi_{9,3}$}}\cg%
\e{0}%
\e{0}%
\e{0}%
\e{1}%
\e{0}%
\e{0}%
\e{0}%
\e{0}%
\e{0}%
\e{0}%
\e{0}%
\e{1}%
\e{0}%
\e{0}%
\e{0}%
\e{0}%
\e{0}%
\e{0}%
\e{0}%
\e{0}%
\eol}\vss}\rg%
%
%
\rx{\vss\hfull{%
\rlx{\hss{$\chi_{9,4}$}}\cg%
\e{0}%
\e{0}%
\e{0}%
\e{0}%
\e{0}%
\e{0}%
\e{0}%
\e{0}%
\e{0}%
\e{0}%
\e{0}%
\e{1}%
\e{0}%
\e{0}%
\e{0}%
\e{0}%
\e{0}%
\e{0}%
\e{1}%
\e{0}%
\eol}\vss}\rg%
%
%
\rx{\vss\hfull{%
\rlx{\hss{$\chi_{6,1}$}}\cg%
\e{0}%
\e{0}%
\e{0}%
\e{0}%
\e{0}%
\e{0}%
\e{0}%
\e{0}%
\e{0}%
\e{0}%
\e{1}%
\e{0}%
\e{0}%
\e{0}%
\e{0}%
\e{0}%
\e{0}%
\e{0}%
\e{0}%
\e{0}%
\eol}\vss}\rg%
%
%
\rx{\vss\hfull{%
\rlx{\hss{$\chi_{6,2}$}}\cg%
\e{0}%
\e{0}%
\e{0}%
\e{0}%
\e{0}%
\e{0}%
\e{0}%
\e{0}%
\e{0}%
\e{1}%
\e{0}%
\e{0}%
\e{0}%
\e{0}%
\e{0}%
\e{0}%
\e{0}%
\e{0}%
\e{0}%
\e{0}%
\eol}\vss}\rg%
%
%
\rx{\vss\hfull{%
\rlx{\hss{$\chi_{12,1}$}}\cg%
\e{0}%
\e{0}%
\e{0}%
\e{0}%
\e{0}%
\e{0}%
\e{0}%
\e{0}%
\e{0}%
\e{1}%
\e{1}%
\e{0}%
\e{0}%
\e{0}%
\e{0}%
\e{0}%
\e{0}%
\e{0}%
\e{0}%
\e{0}%
\eol}\vss}\rg%
%
%
\rx{\vss\hfull{%
\rlx{\hss{$\chi_{4,2}$}}\cg%
\e{0}%
\e{0}%
\e{0}%
\e{0}%
\e{0}%
\e{1}%
\e{0}%
\e{0}%
\e{0}%
\e{0}%
\e{0}%
\e{0}%
\e{0}%
\e{0}%
\e{0}%
\e{0}%
\e{0}%
\e{0}%
\e{0}%
\e{0}%
\eol}\vss}\rg%
%
%
\rx{\vss\hfull{%
\rlx{\hss{$\chi_{4,3}$}}\cg%
\e{0}%
\e{0}%
\e{0}%
\e{0}%
\e{0}%
\e{0}%
\e{0}%
\e{0}%
\e{0}%
\e{0}%
\e{0}%
\e{0}%
\e{1}%
\e{0}%
\e{0}%
\e{0}%
\e{0}%
\e{0}%
\e{0}%
\e{0}%
\eol}\vss}\rg%
%
%
\rx{\vss\hfull{%
\rlx{\hss{$\chi_{4,4}$}}\cg%
\e{0}%
\e{0}%
\e{0}%
\e{0}%
\e{0}%
\e{0}%
\e{0}%
\e{1}%
\e{0}%
\e{0}%
\e{0}%
\e{0}%
\e{0}%
\e{0}%
\e{0}%
\e{0}%
\e{0}%
\e{0}%
\e{0}%
\e{0}%
\eol}\vss}\rg%
%
%
\rx{\vss\hfull{%
\rlx{\hss{$\chi_{4,5}$}}\cg%
\e{0}%
\e{0}%
\e{0}%
\e{0}%
\e{0}%
\e{0}%
\e{0}%
\e{0}%
\e{0}%
\e{0}%
\e{0}%
\e{0}%
\e{0}%
\e{0}%
\e{1}%
\e{0}%
\e{0}%
\e{0}%
\e{0}%
\e{0}%
\eol}\vss}\rg%
%
%
\rx{\vss\hfull{%
\rlx{\hss{$\chi_{8,1}$}}\cg%
\e{0}%
\e{0}%
\e{0}%
\e{0}%
\e{0}%
\e{1}%
\e{0}%
\e{0}%
\e{0}%
\e{0}%
\e{0}%
\e{0}%
\e{1}%
\e{0}%
\e{0}%
\e{0}%
\e{0}%
\e{0}%
\e{0}%
\e{0}%
\eol}\vss}\rg%
%
%
\rx{\vss\hfull{%
\rlx{\hss{$\chi_{8,2}$}}\cg%
\e{0}%
\e{0}%
\e{0}%
\e{0}%
\e{0}%
\e{0}%
\e{0}%
\e{1}%
\e{0}%
\e{0}%
\e{0}%
\e{0}%
\e{0}%
\e{0}%
\e{1}%
\e{0}%
\e{0}%
\e{0}%
\e{0}%
\e{0}%
\eol}\vss}\rg%
%
%
\rx{\vss\hfull{%
\rlx{\hss{$\chi_{8,3}$}}\cg%
\e{0}%
\e{0}%
\e{0}%
\e{0}%
\e{0}%
\e{0}%
\e{1}%
\e{0}%
\e{0}%
\e{0}%
\e{0}%
\e{0}%
\e{0}%
\e{0}%
\e{0}%
\e{0}%
\e{0}%
\e{0}%
\e{0}%
\e{0}%
\eol}\vss}\rg%
%
%
\rx{\vss\hfull{%
\rlx{\hss{$\chi_{8,4}$}}\cg%
\e{0}%
\e{0}%
\e{0}%
\e{0}%
\e{0}%
\e{0}%
\e{0}%
\e{0}%
\e{0}%
\e{0}%
\e{0}%
\e{0}%
\e{0}%
\e{1}%
\e{0}%
\e{0}%
\e{0}%
\e{0}%
\e{0}%
\e{0}%
\eol}\vss}\rg%
%
%
\rx{\vss\hfull{%
\rlx{\hss{$\chi_{16,1}$}}\cg%
\e{0}%
\e{0}%
\e{0}%
\e{0}%
\e{0}%
\e{0}%
\e{1}%
\e{0}%
\e{0}%
\e{0}%
\e{0}%
\e{0}%
\e{0}%
\e{1}%
\e{0}%
\e{0}%
\e{0}%
\e{0}%
\e{0}%
\e{0}%
\eol}\vss}\rg%
\tableclose%
%
%
%
%
%
%
\eop
\eject
\tableopen{Induce/restrict matrix for $W(B_{3})\,\subset\,W(F_{4})$}%
%
%
%
%
%
%
\rowpts=18 true pt%
\colpts=18 true pt%
\rowlabpts=40 true pt%
\collabpts=65 true pt%
\clx{\vss\hfull{%
\rlx{\hss{$ $}}\cg%
\cx{\hskip 16 true pt\flip{$[{3}:-]$}\hss}\cg%
\cx{\hskip 16 true pt\flip{$[{2}{1}:-]$}\hss}\cg%
\cx{\hskip 16 true pt\flip{$[{1^{3}}:-]$}\hss}\cg%
\cx{\hskip 16 true pt\flip{$[{2}:{1}]$}\hss}\cg%
\cx{\hskip 16 true pt\flip{$[{1^{2}}:{1}]$}\hss}\cg%
\cx{\hskip 16 true pt\flip{$[{1}:{2}]$}\hss}\cg%
\cx{\hskip 16 true pt\flip{$[{1}:{1^{2}}]$}\hss}\cg%
\cx{\hskip 16 true pt\flip{$[-:{3}]$}\hss}\cg%
\cx{\hskip 16 true pt\flip{$[-:{2}{1}]$}\hss}\cg%
\cx{\hskip 16 true pt\flip{$[-:{1^{3}}]$}\hss}\cg%
\eol}}\rg%
%
%
\rx{\vss\hfull{%
\rlx{\hss{$\chi_{1,1}$}}\cg%
\e{1}%
\e{0}%
\e{0}%
\e{0}%
\e{0}%
\e{0}%
\e{0}%
\e{0}%
\e{0}%
\e{0}%
\eol}\vss}\rg%
%
%
\rx{\vss\hfull{%
\rlx{\hss{$\chi_{1,2}$}}\cg%
\e{0}%
\e{0}%
\e{0}%
\e{0}%
\e{0}%
\e{0}%
\e{0}%
\e{1}%
\e{0}%
\e{0}%
\eol}\vss}\rg%
%
%
\rx{\vss\hfull{%
\rlx{\hss{$\chi_{1,3}$}}\cg%
\e{0}%
\e{0}%
\e{1}%
\e{0}%
\e{0}%
\e{0}%
\e{0}%
\e{0}%
\e{0}%
\e{0}%
\eol}\vss}\rg%
%
%
\rx{\vss\hfull{%
\rlx{\hss{$\chi_{1,4}$}}\cg%
\e{0}%
\e{0}%
\e{0}%
\e{0}%
\e{0}%
\e{0}%
\e{0}%
\e{0}%
\e{0}%
\e{1}%
\eol}\vss}\rg%
%
%
\rx{\vss\hfull{%
\rlx{\hss{$\chi_{2,1}$}}\cg%
\e{1}%
\e{0}%
\e{0}%
\e{0}%
\e{0}%
\e{0}%
\e{0}%
\e{1}%
\e{0}%
\e{0}%
\eol}\vss}\rg%
%
%
\rx{\vss\hfull{%
\rlx{\hss{$\chi_{2,2}$}}\cg%
\e{0}%
\e{0}%
\e{1}%
\e{0}%
\e{0}%
\e{0}%
\e{0}%
\e{0}%
\e{0}%
\e{1}%
\eol}\vss}\rg%
%
%
\rx{\vss\hfull{%
\rlx{\hss{$\chi_{2,3}$}}\cg%
\e{0}%
\e{1}%
\e{0}%
\e{0}%
\e{0}%
\e{0}%
\e{0}%
\e{0}%
\e{0}%
\e{0}%
\eol}\vss}\rg%
%
%
\rx{\vss\hfull{%
\rlx{\hss{$\chi_{2,4}$}}\cg%
\e{0}%
\e{0}%
\e{0}%
\e{0}%
\e{0}%
\e{0}%
\e{0}%
\e{0}%
\e{1}%
\e{0}%
\eol}\vss}\rg%
%
%
\rx{\vss\hfull{%
\rlx{\hss{$\chi_{4,1}$}}\cg%
\e{0}%
\e{1}%
\e{0}%
\e{0}%
\e{0}%
\e{0}%
\e{0}%
\e{0}%
\e{1}%
\e{0}%
\eol}\vss}\rg%
%
%
\rx{\vss\hfull{%
\rlx{\hss{$\chi_{9,1}$}}\cg%
\e{1}%
\e{1}%
\e{0}%
\e{1}%
\e{0}%
\e{1}%
\e{0}%
\e{0}%
\e{0}%
\e{0}%
\eol}\vss}\rg%
%
%
\rx{\vss\hfull{%
\rlx{\hss{$\chi_{9,2}$}}\cg%
\e{0}%
\e{0}%
\e{0}%
\e{1}%
\e{0}%
\e{1}%
\e{0}%
\e{1}%
\e{1}%
\e{0}%
\eol}\vss}\rg%
%
%
\rx{\vss\hfull{%
\rlx{\hss{$\chi_{9,3}$}}\cg%
\e{0}%
\e{1}%
\e{1}%
\e{0}%
\e{1}%
\e{0}%
\e{1}%
\e{0}%
\e{0}%
\e{0}%
\eol}\vss}\rg%
%
%
\rx{\vss\hfull{%
\rlx{\hss{$\chi_{9,4}$}}\cg%
\e{0}%
\e{0}%
\e{0}%
\e{0}%
\e{1}%
\e{0}%
\e{1}%
\e{0}%
\e{1}%
\e{1}%
\eol}\vss}\rg%
%
%
\rx{\vss\hfull{%
\rlx{\hss{$\chi_{6,1}$}}\cg%
\e{0}%
\e{0}%
\e{0}%
\e{0}%
\e{1}%
\e{1}%
\e{0}%
\e{0}%
\e{0}%
\e{0}%
\eol}\vss}\rg%
%
%
\rx{\vss\hfull{%
\rlx{\hss{$\chi_{6,2}$}}\cg%
\e{0}%
\e{0}%
\e{0}%
\e{1}%
\e{0}%
\e{0}%
\e{1}%
\e{0}%
\e{0}%
\e{0}%
\eol}\vss}\rg%
%
%
\rx{\vss\hfull{%
\rlx{\hss{$\chi_{12,1}$}}\cg%
\e{0}%
\e{0}%
\e{0}%
\e{1}%
\e{1}%
\e{1}%
\e{1}%
\e{0}%
\e{0}%
\e{0}%
\eol}\vss}\rg%
%
%
\rx{\vss\hfull{%
\rlx{\hss{$\chi_{4,2}$}}\cg%
\e{1}%
\e{0}%
\e{0}%
\e{1}%
\e{0}%
\e{0}%
\e{0}%
\e{0}%
\e{0}%
\e{0}%
\eol}\vss}\rg%
%
%
\rx{\vss\hfull{%
\rlx{\hss{$\chi_{4,3}$}}\cg%
\e{0}%
\e{0}%
\e{0}%
\e{0}%
\e{0}%
\e{1}%
\e{0}%
\e{1}%
\e{0}%
\e{0}%
\eol}\vss}\rg%
%
%
\rx{\vss\hfull{%
\rlx{\hss{$\chi_{4,4}$}}\cg%
\e{0}%
\e{0}%
\e{1}%
\e{0}%
\e{1}%
\e{0}%
\e{0}%
\e{0}%
\e{0}%
\e{0}%
\eol}\vss}\rg%
%
%
\rx{\vss\hfull{%
\rlx{\hss{$\chi_{4,5}$}}\cg%
\e{0}%
\e{0}%
\e{0}%
\e{0}%
\e{0}%
\e{0}%
\e{1}%
\e{0}%
\e{0}%
\e{1}%
\eol}\vss}\rg%
%
%
\rx{\vss\hfull{%
\rlx{\hss{$\chi_{8,1}$}}\cg%
\e{1}%
\e{0}%
\e{0}%
\e{1}%
\e{0}%
\e{1}%
\e{0}%
\e{1}%
\e{0}%
\e{0}%
\eol}\vss}\rg%
%
%
\rx{\vss\hfull{%
\rlx{\hss{$\chi_{8,2}$}}\cg%
\e{0}%
\e{0}%
\e{1}%
\e{0}%
\e{1}%
\e{0}%
\e{1}%
\e{0}%
\e{0}%
\e{1}%
\eol}\vss}\rg%
%
%
\rx{\vss\hfull{%
\rlx{\hss{$\chi_{8,3}$}}\cg%
\e{0}%
\e{1}%
\e{0}%
\e{1}%
\e{1}%
\e{0}%
\e{0}%
\e{0}%
\e{0}%
\e{0}%
\eol}\vss}\rg%
%
%
\rx{\vss\hfull{%
\rlx{\hss{$\chi_{8,4}$}}\cg%
\e{0}%
\e{0}%
\e{0}%
\e{0}%
\e{0}%
\e{1}%
\e{1}%
\e{0}%
\e{1}%
\e{0}%
\eol}\vss}\rg%
%
%
\rx{\vss\hfull{%
\rlx{\hss{$\chi_{16,1}$}}\cg%
\e{0}%
\e{1}%
\e{0}%
\e{1}%
\e{1}%
\e{1}%
\e{1}%
\e{0}%
\e{1}%
\e{0}%
\eol}\vss}\rg%
\tableclose%
%
%
%
%
%
%
\eop
\eject
\tableopen{Induce/restrict matrix for $W({A_{2}}{A_{1}^{*}})\,\subset\,W(F_{4})$}%
%
%
%
%
%
%
\rowpts=18 true pt%
\colpts=18 true pt%
\rowlabpts=40 true pt%
\collabpts=60 true pt%
\clx{\vss\hfull{%
\rlx{\hss{$ $}}\cg%
\cx{\hskip 16 true pt\flip{$[{3}]{\times}[{2}]$}\hss}\cg%
\cx{\hskip 16 true pt\flip{$[{2}{1}]{\times}[{2}]$}\hss}\cg%
\cx{\hskip 16 true pt\flip{$[{1^{3}}]{\times}[{2}]$}\hss}\cg%
\cx{\hskip 16 true pt\flip{$[{3}]{\times}[{1^{2}}]$}\hss}\cg%
\cx{\hskip 16 true pt\flip{$[{2}{1}]{\times}[{1^{2}}]$}\hss}\cg%
\cx{\hskip 16 true pt\flip{$[{1^{3}}]{\times}[{1^{2}}]$}\hss}\cg%
\eol}}\rg%
%
%
\rx{\vss\hfull{%
\rlx{\hss{$\chi_{1,1}$}}\cg%
\e{1}%
\e{0}%
\e{0}%
\e{0}%
\e{0}%
\e{0}%
\eol}\vss}\rg%
%
%
\rx{\vss\hfull{%
\rlx{\hss{$\chi_{1,2}$}}\cg%
\e{0}%
\e{0}%
\e{0}%
\e{1}%
\e{0}%
\e{0}%
\eol}\vss}\rg%
%
%
\rx{\vss\hfull{%
\rlx{\hss{$\chi_{1,3}$}}\cg%
\e{0}%
\e{0}%
\e{1}%
\e{0}%
\e{0}%
\e{0}%
\eol}\vss}\rg%
%
%
\rx{\vss\hfull{%
\rlx{\hss{$\chi_{1,4}$}}\cg%
\e{0}%
\e{0}%
\e{0}%
\e{0}%
\e{0}%
\e{1}%
\eol}\vss}\rg%
%
%
\rx{\vss\hfull{%
\rlx{\hss{$\chi_{2,1}$}}\cg%
\e{1}%
\e{0}%
\e{0}%
\e{1}%
\e{0}%
\e{0}%
\eol}\vss}\rg%
%
%
\rx{\vss\hfull{%
\rlx{\hss{$\chi_{2,2}$}}\cg%
\e{0}%
\e{0}%
\e{1}%
\e{0}%
\e{0}%
\e{1}%
\eol}\vss}\rg%
%
%
\rx{\vss\hfull{%
\rlx{\hss{$\chi_{2,3}$}}\cg%
\e{0}%
\e{1}%
\e{0}%
\e{0}%
\e{0}%
\e{0}%
\eol}\vss}\rg%
%
%
\rx{\vss\hfull{%
\rlx{\hss{$\chi_{2,4}$}}\cg%
\e{0}%
\e{0}%
\e{0}%
\e{0}%
\e{1}%
\e{0}%
\eol}\vss}\rg%
%
%
\rx{\vss\hfull{%
\rlx{\hss{$\chi_{4,1}$}}\cg%
\e{0}%
\e{1}%
\e{0}%
\e{0}%
\e{1}%
\e{0}%
\eol}\vss}\rg%
%
%
\rx{\vss\hfull{%
\rlx{\hss{$\chi_{9,1}$}}\cg%
\e{2}%
\e{2}%
\e{0}%
\e{1}%
\e{1}%
\e{0}%
\eol}\vss}\rg%
%
%
\rx{\vss\hfull{%
\rlx{\hss{$\chi_{9,2}$}}\cg%
\e{1}%
\e{1}%
\e{0}%
\e{2}%
\e{2}%
\e{0}%
\eol}\vss}\rg%
%
%
\rx{\vss\hfull{%
\rlx{\hss{$\chi_{9,3}$}}\cg%
\e{0}%
\e{2}%
\e{2}%
\e{0}%
\e{1}%
\e{1}%
\eol}\vss}\rg%
%
%
\rx{\vss\hfull{%
\rlx{\hss{$\chi_{9,4}$}}\cg%
\e{0}%
\e{1}%
\e{1}%
\e{0}%
\e{2}%
\e{2}%
\eol}\vss}\rg%
%
%
\rx{\vss\hfull{%
\rlx{\hss{$\chi_{6,1}$}}\cg%
\e{1}%
\e{1}%
\e{0}%
\e{0}%
\e{1}%
\e{1}%
\eol}\vss}\rg%
%
%
\rx{\vss\hfull{%
\rlx{\hss{$\chi_{6,2}$}}\cg%
\e{0}%
\e{1}%
\e{1}%
\e{1}%
\e{1}%
\e{0}%
\eol}\vss}\rg%
%
%
\rx{\vss\hfull{%
\rlx{\hss{$\chi_{12,1}$}}\cg%
\e{1}%
\e{2}%
\e{1}%
\e{1}%
\e{2}%
\e{1}%
\eol}\vss}\rg%
%
%
\rx{\vss\hfull{%
\rlx{\hss{$\chi_{4,2}$}}\cg%
\e{1}%
\e{1}%
\e{0}%
\e{1}%
\e{0}%
\e{0}%
\eol}\vss}\rg%
%
%
\rx{\vss\hfull{%
\rlx{\hss{$\chi_{4,3}$}}\cg%
\e{1}%
\e{0}%
\e{0}%
\e{1}%
\e{1}%
\e{0}%
\eol}\vss}\rg%
%
%
\rx{\vss\hfull{%
\rlx{\hss{$\chi_{4,4}$}}\cg%
\e{0}%
\e{1}%
\e{1}%
\e{0}%
\e{0}%
\e{1}%
\eol}\vss}\rg%
%
%
\rx{\vss\hfull{%
\rlx{\hss{$\chi_{4,5}$}}\cg%
\e{0}%
\e{0}%
\e{1}%
\e{0}%
\e{1}%
\e{1}%
\eol}\vss}\rg%
%
%
\rx{\vss\hfull{%
\rlx{\hss{$\chi_{8,1}$}}\cg%
\e{2}%
\e{1}%
\e{0}%
\e{2}%
\e{1}%
\e{0}%
\eol}\vss}\rg%
%
%
\rx{\vss\hfull{%
\rlx{\hss{$\chi_{8,2}$}}\cg%
\e{0}%
\e{1}%
\e{2}%
\e{0}%
\e{1}%
\e{2}%
\eol}\vss}\rg%
%
%
\rx{\vss\hfull{%
\rlx{\hss{$\chi_{8,3}$}}\cg%
\e{1}%
\e{2}%
\e{1}%
\e{0}%
\e{1}%
\e{0}%
\eol}\vss}\rg%
%
%
\rx{\vss\hfull{%
\rlx{\hss{$\chi_{8,4}$}}\cg%
\e{0}%
\e{1}%
\e{0}%
\e{1}%
\e{2}%
\e{1}%
\eol}\vss}\rg%
%
%
\rx{\vss\hfull{%
\rlx{\hss{$\chi_{16,1}$}}\cg%
\e{1}%
\e{3}%
\e{1}%
\e{1}%
\e{3}%
\e{1}%
\eol}\vss}\rg%
\tableclose%
%
%
%
%
%
%
\eop
\eject
\tableopen{Induce/restrict matrix for $W({A_{1}}{A_{1}^{*}})\,\subset\,W(G_{2})$}%
%
%
%
%
%
%
\rowpts=18 true pt%
\colpts=18 true pt%
\rowlabpts=40 true pt%
\collabpts=60 true pt%
\clx{\vss\hfull{%
\rlx{\hss{$ $}}\cg%
\cx{\hskip 16 true pt\flip{$[{2}]{\times}[{2}]$}\hss}\cg%
\cx{\hskip 16 true pt\flip{$[{1^{2}}]{\times}[{2}]$}\hss}\cg%
\cx{\hskip 16 true pt\flip{$[{2}]{\times}[{1^{2}}]$}\hss}\cg%
\cx{\hskip 16 true pt\flip{$[{1^{2}}]{\times}[{1^{2}}]$}\hss}\cg%
\eol}}\rg%
%
%
\rx{\vss\hfull{%
\rlx{\hss{$\chi_{1,1}$}}\cg%
\e{1}%
\e{0}%
\e{0}%
\e{0}%
\eol}\vss}\rg%
%
%
\rx{\vss\hfull{%
\rlx{\hss{$\chi_{1,2}$}}\cg%
\e{0}%
\e{0}%
\e{0}%
\e{1}%
\eol}\vss}\rg%
%
%
\rx{\vss\hfull{%
\rlx{\hss{$\chi_{1,3}$}}\cg%
\e{0}%
\e{1}%
\e{0}%
\e{0}%
\eol}\vss}\rg%
%
%
\rx{\vss\hfull{%
\rlx{\hss{$\chi_{1,4}$}}\cg%
\e{0}%
\e{0}%
\e{1}%
\e{0}%
\eol}\vss}\rg%
%
%
\rx{\vss\hfull{%
\rlx{\hss{$\chi_{2,1}$}}\cg%
\e{0}%
\e{1}%
\e{1}%
\e{0}%
\eol}\vss}\rg%
%
%
\rx{\vss\hfull{%
\rlx{\hss{$\chi_{2,2}$}}\cg%
\e{1}%
\e{0}%
\e{0}%
\e{1}%
\eol}\vss}\rg%
\tableclose%
%
%
%
%
%
%
\tableopen{Induce/restrict matrix for $W(A_{2})\,\subset\,W(G_{2})$}%
%
%
%
%
%
%
\rowpts=18 true pt%
\colpts=18 true pt%
\rowlabpts=40 true pt%
\collabpts=45 true pt%
\clx{\vss\hfull{%
\rlx{\hss{$ $}}\cg%
\cx{\hskip 16 true pt\flip{$[{3}]$}\hss}\cg%
\cx{\hskip 16 true pt\flip{$[{2}{1}]$}\hss}\cg%
\cx{\hskip 16 true pt\flip{$[{1^{3}}]$}\hss}\cg%
\eol}}\rg%
%
%
\rx{\vss\hfull{%
\rlx{\hss{$\chi_{1,1}$}}\cg%
\e{1}%
\e{0}%
\e{0}%
\eol}\vss}\rg%
%
%
\rx{\vss\hfull{%
\rlx{\hss{$\chi_{1,2}$}}\cg%
\e{0}%
\e{0}%
\e{1}%
\eol}\vss}\rg%
%
%
\rx{\vss\hfull{%
\rlx{\hss{$\chi_{1,3}$}}\cg%
\e{0}%
\e{0}%
\e{1}%
\eol}\vss}\rg%
%
%
\rx{\vss\hfull{%
\rlx{\hss{$\chi_{1,4}$}}\cg%
\e{1}%
\e{0}%
\e{0}%
\eol}\vss}\rg%
%
%
\rx{\vss\hfull{%
\rlx{\hss{$\chi_{2,1}$}}\cg%
\e{0}%
\e{1}%
\e{0}%
\eol}\vss}\rg%
%
%
\rx{\vss\hfull{%
\rlx{\hss{$\chi_{2,2}$}}\cg%
\e{0}%
\e{1}%
\e{0}%
\eol}\vss}\rg%
\tableclose%
%
%
%
%
%
%
\tableopen{Induce/restrict matrix for $W(A_{1})\,\subset\,W(G_{2})$}%
%
%
%
%
%
%
\rowpts=18 true pt%
\colpts=18 true pt%
\rowlabpts=40 true pt%
\collabpts=45 true pt%
\clx{\vss\hfull{%
\rlx{\hss{$ $}}\cg%
\cx{\hskip 16 true pt\flip{$[{2}]$}\hss}\cg%
\cx{\hskip 16 true pt\flip{$[{1^{2}}]$}\hss}\cg%
\eol}}\rg%
%
%
\rx{\vss\hfull{%
\rlx{\hss{$\chi_{1,1}$}}\cg%
\e{1}%
\e{0}%
\eol}\vss}\rg%
%
%
\rx{\vss\hfull{%
\rlx{\hss{$\chi_{1,2}$}}\cg%
\e{0}%
\e{1}%
\eol}\vss}\rg%
%
%
\rx{\vss\hfull{%
\rlx{\hss{$\chi_{1,3}$}}\cg%
\e{0}%
\e{1}%
\eol}\vss}\rg%
%
%
\rx{\vss\hfull{%
\rlx{\hss{$\chi_{1,4}$}}\cg%
\e{1}%
\e{0}%
\eol}\vss}\rg%
%
%
\rx{\vss\hfull{%
\rlx{\hss{$\chi_{2,1}$}}\cg%
\e{1}%
\e{1}%
\eol}\vss}\rg%
%
%
\rx{\vss\hfull{%
\rlx{\hss{$\chi_{2,2}$}}\cg%
\e{1}%
\e{1}%
\eol}\vss}\rg%
\tableclose%
\vfil
\eject
\sectionhead{Intertwining Numbers}

\bigskip
Let $\Gamma_0$ and $\Gamma_1$ be maximal subdiagrams of
$\widetilde{\Gamma}$ or $\Gamma$, and let
$\varphi_0$, $\varphi_1$ be irreducible characters
of $W(\Gamma_0)$ and $W(\Gamma_1)$, respectively.
The inner product, or intertwining number,
$$(Ind^{W}_{W_0}\varphi_0,Ind^{W}_{W_1}\varphi_1)_{W}$$
can be computed by forming the inner product of the
columns of the induce/restrict matrices
corresponding to $\varphi_0$
and $\varphi_1$.
The intertwining numbers will not
be reproduced in their entirety here.
Instead, we consider only the simplest cases, in
which $\varphi_0$ and $\varphi_1$
are the unit, sign or reflection characters
of $W_0$ and $W_1$, respectively.
The reflection character of
$W$
will be denoted
$\rho$,
and the reflection character of
$W_{i}$
will be denoted
$\rho_{i}$,
$i = 0$, $1$.

\smallskip
The intertwining number
$(Ind^{W}_{W_0}1,Ind^{W}_{W_1}1)_{\down{2 true pt}{W}}$
is equal to
the number of $(W_0,W_1)$-double
cosets in $W$, as is well-known.

\smallskip
The tables below are symmetric, as can
be seen from the identities
$Ind^{W}_{W_0}\varepsilon = \varepsilon\cdot Ind^{W}_{W_0}1$
and
$Ind^{W}_{W_i}\rho_{i} =
(\rho - \delta_{i}\cdot 1)\cdot
Ind^{W}_{W_i}1$,
where $\delta_{i} = |\Gamma| - |\Gamma_i|$.
Only the entries
on or above the diagonal are displayed.

\smallskip
The component $A_2^{*}$ of
${A_1}{A_2^{*}}\,\subset\,F_4$
corresponds to short roots.
The subgroup
$W(A_1^{*})$
of $W(G_2)$
contains the reflection with respect
to a short root.
\begingroup
\def\e#1{\cx{\hss{$#1$}\hss}\cg}%
\def\eol{\hskip 0pt plus 20 true pt\break}%
\vskip 5 true pt plus 10 true pt
\goodbreak
\tableopen{Intertwining numbers $(Ind^{W}_{W_0}1,Ind^{W}_{W_1}1)_{\down{2 true pt}{W}}$ for $W = W(E_{6})$}%
\rowpts=19 true pt%
\rowlabpts=120 true pt%
\collabpts=50 true pt%
\clx{\vss\hfull{%
\rlx{\hss{$ $}}\cg%
\hbox to 20 true pt{\hskip 12 true pt\twist{${A_{5}}{A_{1}}$}\hss}\cg%
\hbox to 20 true pt{\hskip 12 true pt\twist{${A_{2}}{A_{2}}{A_{2}}$}\hss}\cg%
\hbox to 1.5\colpts{\hfil}\rlx{\hss{$ $}}\cg%
\hbox to 20 true pt{\hskip 12 true pt\twist{$D_{5}$}\hss}\cg%
\hbox to 20 true pt{\hskip 12 true pt\twist{$A_{5}$}\hss}\cg%
\hbox to 20 true pt{\hskip 12 true pt\twist{${A_{4}}{A_{1}}$}\hss}\cg%
\hbox to 20 true pt{\hskip 12 true pt\twist{${A_{2}}{A_{2}}{A_{1}}$}\hss}\cg%
\eol}}\rg%
\rx{\vss\hfull{%
\rlx{\hss{${A_{5}}{A_{1}}$}}\cg%
\hbox to 20 true pt{\hss{$3$}}\cg%
\hbox to 20 true pt{\hss{$4$}}\cg%
\hbox to 1.5\colpts{\hfil}\rlx{\hss{$D_{5}$}}\cg%
\hbox to 20 true pt{\hss{$3$}}\cg%
\hbox to 20 true pt{\hss{$3$}}\cg%
\hbox to 20 true pt{\hss{$4$}}\cg%
\hbox to 20 true pt{\hss{$5$}}\cg%
\eol}\vss}\rg%
\rx{\vss\hfull{%
\rlx{\hss{${A_{2}}{A_{2}}{A_{2}}$}}\cg%
\hbox to 20 true pt{\hfil}\cg%
\hbox to 20 true pt{\hss{$13$}}\cg%
\hbox to 1.5\colpts{\hfil}\rlx{\hss{$A_{5}$}}\cg%
\hbox to 20 true pt{\hfil}\cg%
\hbox to 20 true pt{\hss{$5$}}\cg%
\hbox to 20 true pt{\hss{$6$}}\cg%
\hbox to 20 true pt{\hss{$9$}}\cg%
\eol}\vss}\rg%
\rx{\vss\hfull{%
\rlx{\hss{$ $}}\cg%
\hbox to 20 true pt{\hfil}\cg%
\hbox to 20 true pt{\hfil}\cg%
\hbox to 1.5\colpts{\hfil}\rlx{\hss{${A_{4}}{A_{1}}$}}\cg%
\hbox to 20 true pt{\hfil}\cg%
\hbox to 20 true pt{\hfil}\cg%
\hbox to 20 true pt{\hss{$10$}}\cg%
\hbox to 20 true pt{\hss{$17$}}\cg%
\eol}\vss}\rg%
\rx{\vss\hfull{%
\rlx{\hss{$ $}}\cg%
\hbox to 20 true pt{\hfil}\cg%
\hbox to 20 true pt{\hfil}\cg%
\hbox to 1.5\colpts{\hfil}\rlx{\hss{${A_{2}}{A_{2}}{A_{1}}$}}\cg%
\hbox to 20 true pt{\hfil}\cg%
\hbox to 20 true pt{\hfil}\cg%
\hbox to 20 true pt{\hfil}\cg%
\hbox to 20 true pt{\hss{$37$}}\cg%
\eol}\vss}\rg%
\tableclose%

\vskip 5 true pt plus 10 true pt
\goodbreak
\tableopen{Intertwining numbers $(Ind^{W}_{W_0}1,Ind^{W}_{W_1}1)_{\down{2 true pt}{W}}$ for $W = W(E_{7})$}%
\rowpts=19 true pt%
\rowlabpts=85 true pt%
\collabpts=50 true pt%
\clx{\vss\hfull{%
\rlx{\hss{$ $}}\cg%
\hbox to 20 true pt{\hskip 12 true pt\twist{${D_{6}}{A_{1}}$}\hss}\cg%
\hbox to 20 true pt{\hskip 12 true pt\twist{$A_{7}$}\hss}\cg%
\hbox to 20 true pt{\hskip 12 true pt\twist{${A_{5}}{A_{2}}$}\hss}\cg%
\hbox to 20 true pt{\hskip 12 true pt\twist{${A_{3}}{A_{3}}{A_{1}}$}\hss}\cg%
\hbox to 1.5\colpts{\hfil}\rlx{\hss{$ $}}\cg%
\hbox to 20 true pt{\hskip 12 true pt\twist{$D_{6}$}\hss}\cg%
\hbox to 20 true pt{\hskip 12 true pt\twist{$A_{6}$}\hss}\cg%
\hbox to 20 true pt{\hskip 12 true pt\twist{${A_{5}}{A_{1}}$}\hss}\cg%
\hbox to 21 true pt{\hskip 12 true pt\twist{${A_{3}}{A_{2}}{A_{1}}$}\hss}\cg%
\hbox to 20 true pt{\hskip 12 true pt\twist{${A_{4}}{A_{2}}$}\hss}\cg%
\hbox to 20 true pt{\hskip 12 true pt\twist{${D_{5}}{A_{1}}$}\hss}\cg%
\hbox to 20 true pt{\hskip 12 true pt\twist{$E_{6}$}\hss}\cg%
\eol}}\rg%
\rx{\vss\hfull{%
\rlx{\hss{${D_{6}}{A_{1}}$}}\cg%
\hbox to 20 true pt{\hss{$3$}}\cg%
\hbox to 20 true pt{\hss{$2$}}\cg%
\hbox to 20 true pt{\hss{$3$}}\cg%
\hbox to 20 true pt{\hss{$5$}}\cg%
\hbox to 1.5\colpts{\hfil}\rlx{\hss{$D_{6}$}}\cg%
\hbox to 20 true pt{\hss{$5$}}\cg%
\hbox to 20 true pt{\hss{$5$}}\cg%
\hbox to 20 true pt{\hss{$8$}}\cg%
\hbox to 21 true pt{\hss{$11$}}\cg%
\hbox to 20 true pt{\hss{$8$}}\cg%
\hbox to 20 true pt{\hss{$6$}}\cg%
\hbox to 20 true pt{\hss{$3$}}\cg%
\eol}\vss}\rg%
\rx{\vss\hfull{%
\rlx{\hss{$A_{7}$}}\cg%
\hbox to 20 true pt{\hfil}\cg%
\hbox to 20 true pt{\hss{$4$}}\cg%
\hbox to 20 true pt{\hss{$3$}}\cg%
\hbox to 20 true pt{\hss{$5$}}\cg%
\hbox to 1.5\colpts{\hfil}\rlx{\hss{$A_{6}$}}\cg%
\hbox to 20 true pt{\hfil}\cg%
\hbox to 20 true pt{\hss{$10$}}\cg%
\hbox to 20 true pt{\hss{$13$}}\cg%
\hbox to 21 true pt{\hss{$25$}}\cg%
\hbox to 20 true pt{\hss{$18$}}\cg%
\hbox to 20 true pt{\hss{$9$}}\cg%
\hbox to 20 true pt{\hss{$4$}}\cg%
\eol}\vss}\rg%
\rx{\vss\hfull{%
\rlx{\hss{${A_{5}}{A_{2}}$}}\cg%
\hbox to 20 true pt{\hfil}\cg%
\hbox to 20 true pt{\hfil}\cg%
\hbox to 20 true pt{\hss{$9$}}\cg%
\hbox to 20 true pt{\hss{$13$}}\cg%
\hbox to 1.5\colpts{\hfil}\rlx{\hss{${A_{5}}{A_{1}}$}}\cg%
\hbox to 20 true pt{\hfil}\cg%
\hbox to 20 true pt{\hfil}\cg%
\hbox to 20 true pt{\hss{$27$}}\cg%
\hbox to 21 true pt{\hss{$54$}}\cg%
\hbox to 20 true pt{\hss{$33$}}\cg%
\hbox to 20 true pt{\hss{$16$}}\cg%
\hbox to 20 true pt{\hss{$5$}}\cg%
\eol}\vss}\rg%
\rx{\vss\hfull{%
\rlx{\hss{${A_{3}}{A_{3}}{A_{1}}$}}\cg%
\hbox to 20 true pt{\hfil}\cg%
\hbox to 20 true pt{\hfil}\cg%
\hbox to 20 true pt{\hfil}\cg%
\hbox to 20 true pt{\hss{$29$}}\cg%
\hbox to 1.5\colpts{\hfil}\rlx{\hss{${A_{3}}{A_{2}}{A_{1}}$}}\cg%
\hbox to 20 true pt{\hfil}\cg%
\hbox to 20 true pt{\hfil}\cg%
\hbox to 20 true pt{\hfil}\cg%
\hbox to 21 true pt{\hss{$149$}}\cg%
\hbox to 20 true pt{\hss{$82$}}\cg%
\hbox to 20 true pt{\hss{$30$}}\cg%
\hbox to 20 true pt{\hss{$7$}}\cg%
\eol}\vss}\rg%
\rx{\vss\hfull{%
\rlx{\hss{$ $}}\cg%
\hbox to 20 true pt{\hfil}\cg%
\hbox to 20 true pt{\hfil}\cg%
\hbox to 20 true pt{\hfil}\cg%
\hbox to 20 true pt{\hfil}\cg%
\hbox to 1.5\colpts{\hfil}\rlx{\hss{${A_{4}}{A_{2}}$}}\cg%
\hbox to 20 true pt{\hfil}\cg%
\hbox to 20 true pt{\hfil}\cg%
\hbox to 20 true pt{\hfil}\cg%
\hbox to 21 true pt{\hfil}\cg%
\hbox to 20 true pt{\hss{$50$}}\cg%
\hbox to 20 true pt{\hss{$20$}}\cg%
\hbox to 20 true pt{\hss{$6$}}\cg%
\eol}\vss}\rg%
\rx{\vss\hfull{%
\rlx{\hss{$ $}}\cg%
\hbox to 20 true pt{\hfil}\cg%
\hbox to 20 true pt{\hfil}\cg%
\hbox to 20 true pt{\hfil}\cg%
\hbox to 20 true pt{\hfil}\cg%
\hbox to 1.5\colpts{\hfil}\rlx{\hss{${D_{5}}{A_{1}}$}}\cg%
\hbox to 20 true pt{\hfil}\cg%
\hbox to 20 true pt{\hfil}\cg%
\hbox to 20 true pt{\hfil}\cg%
\hbox to 21 true pt{\hfil}\cg%
\hbox to 20 true pt{\hfil}\cg%
\hbox to 20 true pt{\hss{$13$}}\cg%
\hbox to 20 true pt{\hss{$5$}}\cg%
\eol}\vss}\rg%
\rx{\vss\hfull{%
\rlx{\hss{$ $}}\cg%
\hbox to 20 true pt{\hfil}\cg%
\hbox to 20 true pt{\hfil}\cg%
\hbox to 20 true pt{\hfil}\cg%
\hbox to 20 true pt{\hfil}\cg%
\hbox to 1.5\colpts{\hfil}\rlx{\hss{$E_{6}$}}\cg%
\hbox to 20 true pt{\hfil}\cg%
\hbox to 20 true pt{\hfil}\cg%
\hbox to 20 true pt{\hfil}\cg%
\hbox to 21 true pt{\hfil}\cg%
\hbox to 20 true pt{\hfil}\cg%
\hbox to 20 true pt{\hfil}\cg%
\hbox to 20 true pt{\hss{$4$}}\cg%
\eol}\vss}\rg%
\tableclose%

\vskip 5 true pt plus 10 true pt
\goodbreak
\tableopen{Intertwining numbers $(Ind^{W}_{W_0}1,Ind^{W}_{W_1}1)_{\down{2 true pt}{W}}$ for $W = W(E_{8})$}%
\rowpts=19 true pt%
\rowlabpts=50 true pt%
\collabpts=50 true pt%
\clx{\vss\hfull{%
\rlx{\hss{$ $}}\cg%
\hbox to 20 true pt{\hskip 12 true pt\twist{$D_{8}$}\hss}\cg%
\hbox to 20 true pt{\hskip 12 true pt\twist{$A_{8}$}\hss}\cg%
\hbox to 20 true pt{\hskip 12 true pt\twist{${A_{7}}{A_{1}}$}\hss}\cg%
\hbox to 21 true pt{\hskip 12 true pt\twist{${A_{5}}{A_{2}}{A_{1}}$}\hss}\cg%
\hbox to 20 true pt{\hskip 12 true pt\twist{${A_{4}}{A_{4}}$}\hss}\cg%
\hbox to 20 true pt{\hskip 12 true pt\twist{${D_{5}}{A_{3}}$}\hss}\cg%
\hbox to 20 true pt{\hskip 12 true pt\twist{${E_{6}}{A_{2}}$}\hss}\cg%
\hbox to 20 true pt{\hskip 12 true pt\twist{${E_{7}}{A_{1}}$}\hss}\cg%
\hbox to 1.5\colpts{\hfil}\rlx{\hss{$ $}}\cg%
\hbox to 20 true pt{\hskip 12 true pt\twist{$D_{7}$}\hss}\cg%
\hbox to 20 true pt{\hskip 12 true pt\twist{$A_{7}$}\hss}\cg%
\hbox to 21 true pt{\hskip 12 true pt\twist{${A_{6}}{A_{1}}$}\hss}\cg%
\hbox to 28 true pt{\hskip 16 true pt\twist{${A_{4}}{A_{2}}{A_{1}}$}\hss}\cg%
\hbox to 21 true pt{\hskip 12 true pt\twist{${A_{4}}{A_{3}}$}\hss}\cg%
\hbox to 21 true pt{\hskip 12 true pt\twist{${D_{5}}{A_{2}}$}\hss}\cg%
\hbox to 20 true pt{\hskip 12 true pt\twist{${E_{6}}{A_{1}}$}\hss}\cg%
\hbox to 20 true pt{\hskip 12 true pt\twist{$E_{7}$}\hss}\cg%
\eol}}\rg%
\rx{\vss\hfull{%
\rlx{\hss{$D_{8}$}}\cg%
\hbox to 20 true pt{\hss{$3$}}\cg%
\hbox to 20 true pt{\hss{$2$}}\cg%
\hbox to 20 true pt{\hss{$4$}}\cg%
\hbox to 21 true pt{\hss{$5$}}\cg%
\hbox to 20 true pt{\hss{$4$}}\cg%
\hbox to 20 true pt{\hss{$4$}}\cg%
\hbox to 20 true pt{\hss{$2$}}\cg%
\hbox to 20 true pt{\hss{$2$}}\cg%
\hbox to 1.5\colpts{\hfil}\rlx{\hss{$D_{7}$}}\cg%
\hbox to 20 true pt{\hss{$10$}}\cg%
\hbox to 20 true pt{\hss{$15$}}\cg%
\hbox to 21 true pt{\hss{$24$}}\cg%
\hbox to 28 true pt{\hss{$46$}}\cg%
\hbox to 21 true pt{\hss{$35$}}\cg%
\hbox to 21 true pt{\hss{$23$}}\cg%
\hbox to 20 true pt{\hss{$12$}}\cg%
\hbox to 20 true pt{\hss{$5$}}\cg%
\eol}\vss}\rg%
\rx{\vss\hfull{%
\rlx{\hss{$A_{8}$}}\cg%
\hbox to 20 true pt{\hfil}\cg%
\hbox to 20 true pt{\hss{$7$}}\cg%
\hbox to 20 true pt{\hss{$7$}}\cg%
\hbox to 21 true pt{\hss{$16$}}\cg%
\hbox to 20 true pt{\hss{$13$}}\cg%
\hbox to 20 true pt{\hss{$8$}}\cg%
\hbox to 20 true pt{\hss{$5$}}\cg%
\hbox to 20 true pt{\hss{$2$}}\cg%
\hbox to 1.5\colpts{\hfil}\rlx{\hss{$A_{7}$}}\cg%
\hbox to 20 true pt{\hfil}\cg%
\hbox to 20 true pt{\hss{$35$}}\cg%
\hbox to 21 true pt{\hss{$63$}}\cg%
\hbox to 28 true pt{\hss{$155$}}\cg%
\hbox to 21 true pt{\hss{$109$}}\cg%
\hbox to 21 true pt{\hss{$57$}}\cg%
\hbox to 20 true pt{\hss{$23$}}\cg%
\hbox to 20 true pt{\hss{$7$}}\cg%
\eol}\vss}\rg%
\rx{\vss\hfull{%
\rlx{\hss{${A_{7}}{A_{1}}$}}\cg%
\hbox to 20 true pt{\hfil}\cg%
\hbox to 20 true pt{\hfil}\cg%
\hbox to 20 true pt{\hss{$17$}}\cg%
\hbox to 21 true pt{\hss{$36$}}\cg%
\hbox to 20 true pt{\hss{$26$}}\cg%
\hbox to 20 true pt{\hss{$17$}}\cg%
\hbox to 20 true pt{\hss{$8$}}\cg%
\hbox to 20 true pt{\hss{$4$}}\cg%
\hbox to 1.5\colpts{\hfil}\rlx{\hss{${A_{6}}{A_{1}}$}}\cg%
\hbox to 20 true pt{\hfil}\cg%
\hbox to 20 true pt{\hfil}\cg%
\hbox to 21 true pt{\hss{$134$}}\cg%
\hbox to 28 true pt{\hss{$387$}}\cg%
\hbox to 21 true pt{\hss{$252$}}\cg%
\hbox to 21 true pt{\hss{$123$}}\cg%
\hbox to 20 true pt{\hss{$44$}}\cg%
\hbox to 20 true pt{\hss{$10$}}\cg%
\eol}\vss}\rg%
\rx{\vss\hfull{%
\rlx{\hss{${A_{5}}{A_{2}}{A_{1}}$}}\cg%
\hbox to 20 true pt{\hfil}\cg%
\hbox to 20 true pt{\hfil}\cg%
\hbox to 20 true pt{\hfil}\cg%
\hbox to 21 true pt{\hss{$125$}}\cg%
\hbox to 20 true pt{\hss{$86$}}\cg%
\hbox to 20 true pt{\hss{$47$}}\cg%
\hbox to 20 true pt{\hss{$20$}}\cg%
\hbox to 20 true pt{\hss{$6$}}\cg%
\hbox to 1.5\colpts{\hfil}\rlx{\hss{${A_{4}}{A_{2}}{A_{1}}$}}\cg%
\hbox to 20 true pt{\hfil}\cg%
\hbox to 20 true pt{\hfil}\cg%
\hbox to 21 true pt{\hfil}\cg%
\hbox to 28 true pt{\hss{$1437$}}\cg%
\hbox to 21 true pt{\hss{$870$}}\cg%
\hbox to 21 true pt{\hss{$359$}}\cg%
\hbox to 20 true pt{\hss{$97$}}\cg%
\hbox to 20 true pt{\hss{$15$}}\cg%
\eol}\vss}\rg%
\rx{\vss\hfull{%
\rlx{\hss{${A_{4}}{A_{4}}$}}\cg%
\hbox to 20 true pt{\hfil}\cg%
\hbox to 20 true pt{\hfil}\cg%
\hbox to 20 true pt{\hfil}\cg%
\hbox to 21 true pt{\hfil}\cg%
\hbox to 20 true pt{\hss{$67$}}\cg%
\hbox to 20 true pt{\hss{$35$}}\cg%
\hbox to 20 true pt{\hss{$14$}}\cg%
\hbox to 20 true pt{\hss{$4$}}\cg%
\hbox to 1.5\colpts{\hfil}\rlx{\hss{${A_{4}}{A_{3}}$}}\cg%
\hbox to 20 true pt{\hfil}\cg%
\hbox to 20 true pt{\hfil}\cg%
\hbox to 21 true pt{\hfil}\cg%
\hbox to 28 true pt{\hfil}\cg%
\hbox to 21 true pt{\hss{$547$}}\cg%
\hbox to 21 true pt{\hss{$232$}}\cg%
\hbox to 20 true pt{\hss{$66$}}\cg%
\hbox to 20 true pt{\hss{$12$}}\cg%
\eol}\vss}\rg%
\rx{\vss\hfull{%
\rlx{\hss{${D_{5}}{A_{3}}$}}\cg%
\hbox to 20 true pt{\hfil}\cg%
\hbox to 20 true pt{\hfil}\cg%
\hbox to 20 true pt{\hfil}\cg%
\hbox to 21 true pt{\hfil}\cg%
\hbox to 20 true pt{\hfil}\cg%
\hbox to 20 true pt{\hss{$24$}}\cg%
\hbox to 20 true pt{\hss{$10$}}\cg%
\hbox to 20 true pt{\hss{$4$}}\cg%
\hbox to 1.5\colpts{\hfil}\rlx{\hss{${D_{5}}{A_{2}}$}}\cg%
\hbox to 20 true pt{\hfil}\cg%
\hbox to 20 true pt{\hfil}\cg%
\hbox to 21 true pt{\hfil}\cg%
\hbox to 28 true pt{\hfil}\cg%
\hbox to 21 true pt{\hfil}\cg%
\hbox to 21 true pt{\hss{$122$}}\cg%
\hbox to 20 true pt{\hss{$45$}}\cg%
\hbox to 20 true pt{\hss{$10$}}\cg%
\eol}\vss}\rg%
\rx{\vss\hfull{%
\rlx{\hss{${E_{6}}{A_{2}}$}}\cg%
\hbox to 20 true pt{\hfil}\cg%
\hbox to 20 true pt{\hfil}\cg%
\hbox to 20 true pt{\hfil}\cg%
\hbox to 21 true pt{\hfil}\cg%
\hbox to 20 true pt{\hfil}\cg%
\hbox to 20 true pt{\hfil}\cg%
\hbox to 20 true pt{\hss{$8$}}\cg%
\hbox to 20 true pt{\hss{$3$}}\cg%
\hbox to 1.5\colpts{\hfil}\rlx{\hss{${E_{6}}{A_{1}}$}}\cg%
\hbox to 20 true pt{\hfil}\cg%
\hbox to 20 true pt{\hfil}\cg%
\hbox to 21 true pt{\hfil}\cg%
\hbox to 28 true pt{\hfil}\cg%
\hbox to 21 true pt{\hfil}\cg%
\hbox to 21 true pt{\hfil}\cg%
\hbox to 20 true pt{\hss{$26$}}\cg%
\hbox to 20 true pt{\hss{$8$}}\cg%
\eol}\vss}\rg%
\rx{\vss\hfull{%
\rlx{\hss{${E_{7}}{A_{1}}$}}\cg%
\hbox to 20 true pt{\hfil}\cg%
\hbox to 20 true pt{\hfil}\cg%
\hbox to 20 true pt{\hfil}\cg%
\hbox to 21 true pt{\hfil}\cg%
\hbox to 20 true pt{\hfil}\cg%
\hbox to 20 true pt{\hfil}\cg%
\hbox to 20 true pt{\hfil}\cg%
\hbox to 20 true pt{\hss{$3$}}\cg%
\hbox to 1.5\colpts{\hfil}\rlx{\hss{$E_{7}$}}\cg%
\hbox to 20 true pt{\hfil}\cg%
\hbox to 20 true pt{\hfil}\cg%
\hbox to 21 true pt{\hfil}\cg%
\hbox to 28 true pt{\hfil}\cg%
\hbox to 21 true pt{\hfil}\cg%
\hbox to 21 true pt{\hfil}\cg%
\hbox to 20 true pt{\hfil}\cg%
\hbox to 20 true pt{\hss{$5$}}\cg%
\eol}\vss}\rg%
\tableclose%

\vskip 5 true pt plus 10 true pt
\goodbreak
\tableopen{Intertwining numbers $(Ind^{W}_{W_0}1,Ind^{W}_{W_1}1)_{\down{2 true pt}{W}}$ for $W = W(F_{4})$}%
\rowpts=19 true pt%
\rowlabpts=105 true pt%
\collabpts=45 true pt%
\clx{\vss\hfull{%
\rlx{\hss{$ $}}\cg%
\hbox to 20 true pt{\hskip 12 true pt\twist{${C_{3}}{A_{1}}$}\hss}\cg%
\hbox to 20 true pt{\hskip 12 true pt\twist{${A_{2}}{A_{2}^{*}}$}\hss}\cg%
\hbox to 20 true pt{\hskip 12 true pt\twist{${A_{3}}{A_{1}}$}\hss}\cg%
\hbox to 20 true pt{\hskip 12 true pt\twist{$B_{4}$}\hss}\cg%
\hbox to 1.5\colpts{\hfil}\rlx{\hss{$ $}}\cg%
\hbox to 20 true pt{\hskip 12 true pt\twist{$C_{3}$}\hss}\cg%
\hbox to 20 true pt{\hskip 12 true pt\twist{${A_{1}}{A_{2}^{*}}$}\hss}\cg%
\hbox to 20 true pt{\hskip 12 true pt\twist{${A_{2}}{A_{1}^{*}}$}\hss}\cg%
\hbox to 20 true pt{\hskip 12 true pt\twist{$B_{3}$}\hss}\cg%
\eol}}\rg%
\rx{\vss\hfull{%
\rlx{\hss{${C_{3}}{A_{1}}$}}\cg%
\hbox to 20 true pt{\hss{$3$}}\cg%
\hbox to 20 true pt{\hss{$2$}}\cg%
\hbox to 20 true pt{\hss{$2$}}\cg%
\hbox to 20 true pt{\hss{$1$}}\cg%
\hbox to 1.5\colpts{\hfil}\rlx{\hss{$C_{3}$}}\cg%
\hbox to 20 true pt{\hss{$5$}}\cg%
\hbox to 20 true pt{\hss{$7$}}\cg%
\hbox to 20 true pt{\hss{$5$}}\cg%
\hbox to 20 true pt{\hss{$3$}}\cg%
\eol}\vss}\rg%
\rx{\vss\hfull{%
\rlx{\hss{${A_{2}}{A_{2}^{*}}$}}\cg%
\hbox to 20 true pt{\hfil}\cg%
\hbox to 20 true pt{\hss{$5$}}\cg%
\hbox to 20 true pt{\hss{$3$}}\cg%
\hbox to 20 true pt{\hss{$1$}}\cg%
\hbox to 1.5\colpts{\hfil}\rlx{\hss{${A_{1}}{A_{2}^{*}}$}}\cg%
\hbox to 20 true pt{\hfil}\cg%
\hbox to 20 true pt{\hss{$17$}}\cg%
\hbox to 20 true pt{\hss{$13$}}\cg%
\hbox to 20 true pt{\hss{$5$}}\cg%
\eol}\vss}\rg%
\rx{\vss\hfull{%
\rlx{\hss{${A_{3}}{A_{1}}$}}\cg%
\hbox to 20 true pt{\hfil}\cg%
\hbox to 20 true pt{\hfil}\cg%
\hbox to 20 true pt{\hss{$5$}}\cg%
\hbox to 20 true pt{\hss{$2$}}\cg%
\hbox to 1.5\colpts{\hfil}\rlx{\hss{${A_{2}}{A_{1}^{*}}$}}\cg%
\hbox to 20 true pt{\hfil}\cg%
\hbox to 20 true pt{\hfil}\cg%
\hbox to 20 true pt{\hss{$17$}}\cg%
\hbox to 20 true pt{\hss{$7$}}\cg%
\eol}\vss}\rg%
\rx{\vss\hfull{%
\rlx{\hss{$B_{4}$}}\cg%
\hbox to 20 true pt{\hfil}\cg%
\hbox to 20 true pt{\hfil}\cg%
\hbox to 20 true pt{\hfil}\cg%
\hbox to 20 true pt{\hss{$2$}}\cg%
\hbox to 1.5\colpts{\hfil}\rlx{\hss{$B_{3}$}}\cg%
\hbox to 20 true pt{\hfil}\cg%
\hbox to 20 true pt{\hfil}\cg%
\hbox to 20 true pt{\hfil}\cg%
\hbox to 20 true pt{\hss{$5$}}\cg%
\eol}\vss}\rg%
\tableclose%

\vskip 5 true pt plus 10 true pt
\goodbreak
\tableopen{Intertwining numbers $(Ind^{W}_{W_0}1,Ind^{W}_{W_1}1)_{\down{2 true pt}{W}}$ for $W = W(G_{2})$}%
\rowpts=19 true pt%
\rowlabpts=122 true pt%
\collabpts=40 true pt%
\clx{\vss\hfull{%
\rlx{\hss{$ $}}\cg%
\hbox to 20 true pt{\hskip 12 true pt\twist{${A_{1}}{A_{1}^{*}}$}\hss}\cg%
\hbox to 20 true pt{\hskip 12 true pt\twist{$A_{2}$}\hss}\cg%
\hbox to 1.5\colpts{\hfil}\rlx{\hss{$ $}}\cg%
\hbox to 20 true pt{\hskip 12 true pt\twist{$A_{1}$}\hss}\cg%
\hbox to 20 true pt{\hskip 12 true pt\twist{$A_{1}^{*}$}\hss}\cg%
\eol}}\rg%
\rx{\vss\hfull{%
\rlx{\hss{${A_{1}}{A_{1}^{*}}$}}\cg%
\hbox to 20 true pt{\hss{$2$}}\cg%
\hbox to 20 true pt{\hss{$1$}}\cg%
\hbox to 1.5\colpts{\hfil}\rlx{\hss{$A_{1}$}}\cg%
\hbox to 20 true pt{\hss{$4$}}\cg%
\hbox to 20 true pt{\hss{$3$}}\cg%
\eol}\vss}\rg%
\rx{\vss\hfull{%
\rlx{\hss{$A_{2}$}}\cg%
\hbox to 20 true pt{\hfil}\cg%
\hbox to 20 true pt{\hss{$2$}}\cg%
\hbox to 1.5\colpts{\hfil}\rlx{\hss{$A_{1}^{*}$}}\cg%
\hbox to 20 true pt{\hfil}\cg%
\hbox to 20 true pt{\hss{$4$}}\cg%
\eol}\vss}\rg%
\tableclose%

\vskip 5 true pt plus 10 true pt
\goodbreak
\tableopen{Intertwining numbers $(Ind^{W}_{W_0}1,Ind^{W}_{W_1}\varepsilon)_{\down{2 true pt}{W}}$ for $W = W(E_{6})$}%
\rowpts=19 true pt%
\rowlabpts=120 true pt%
\collabpts=50 true pt%
\clx{\vss\hfull{%
\rlx{\hss{$ $}}\cg%
\hbox to 20 true pt{\hskip 12 true pt\twist{${A_{5}}{A_{1}}$}\hss}\cg%
\hbox to 20 true pt{\hskip 12 true pt\twist{${A_{2}}{A_{2}}{A_{2}}$}\hss}\cg%
\hbox to 1.5\colpts{\hfil}\rlx{\hss{$ $}}\cg%
\hbox to 20 true pt{\hskip 12 true pt\twist{$D_{5}$}\hss}\cg%
\hbox to 20 true pt{\hskip 12 true pt\twist{$A_{5}$}\hss}\cg%
\hbox to 20 true pt{\hskip 12 true pt\twist{${A_{4}}{A_{1}}$}\hss}\cg%
\hbox to 20 true pt{\hskip 12 true pt\twist{${A_{2}}{A_{2}}{A_{1}}$}\hss}\cg%
\eol}}\rg%
\rx{\vss\hfull{%
\rlx{\hss{${A_{5}}{A_{1}}$}}\cg%
\hbox to 20 true pt{\hss{$0$}}\cg%
\hbox to 20 true pt{\hss{$0$}}\cg%
\hbox to 1.5\colpts{\hfil}\rlx{\hss{$D_{5}$}}\cg%
\hbox to 20 true pt{\hss{$0$}}\cg%
\hbox to 20 true pt{\hss{$0$}}\cg%
\hbox to 20 true pt{\hss{$0$}}\cg%
\hbox to 20 true pt{\hss{$0$}}\cg%
\eol}\vss}\rg%
\rx{\vss\hfull{%
\rlx{\hss{${A_{2}}{A_{2}}{A_{2}}$}}\cg%
\hbox to 20 true pt{\hfil}\cg%
\hbox to 20 true pt{\hss{$1$}}\cg%
\hbox to 1.5\colpts{\hfil}\rlx{\hss{$A_{5}$}}\cg%
\hbox to 20 true pt{\hfil}\cg%
\hbox to 20 true pt{\hss{$0$}}\cg%
\hbox to 20 true pt{\hss{$0$}}\cg%
\hbox to 20 true pt{\hss{$0$}}\cg%
\eol}\vss}\rg%
\rx{\vss\hfull{%
\rlx{\hss{$ $}}\cg%
\hbox to 20 true pt{\hfil}\cg%
\hbox to 20 true pt{\hfil}\cg%
\hbox to 1.5\colpts{\hfil}\rlx{\hss{${A_{4}}{A_{1}}$}}\cg%
\hbox to 20 true pt{\hfil}\cg%
\hbox to 20 true pt{\hfil}\cg%
\hbox to 20 true pt{\hss{$0$}}\cg%
\hbox to 20 true pt{\hss{$0$}}\cg%
\eol}\vss}\rg%
\rx{\vss\hfull{%
\rlx{\hss{$ $}}\cg%
\hbox to 20 true pt{\hfil}\cg%
\hbox to 20 true pt{\hfil}\cg%
\hbox to 1.5\colpts{\hfil}\rlx{\hss{${A_{2}}{A_{2}}{A_{1}}$}}\cg%
\hbox to 20 true pt{\hfil}\cg%
\hbox to 20 true pt{\hfil}\cg%
\hbox to 20 true pt{\hfil}\cg%
\hbox to 20 true pt{\hss{$3$}}\cg%
\eol}\vss}\rg%
\tableclose%

\vskip 5 true pt plus 10 true pt
\goodbreak
\tableopen{Intertwining numbers $(Ind^{W}_{W_0}1,Ind^{W}_{W_1}\varepsilon)_{\down{2 true pt}{W}}$ for $W = W(E_{7})$}%
\rowpts=19 true pt%
\rowlabpts=85 true pt%
\collabpts=50 true pt%
\clx{\vss\hfull{%
\rlx{\hss{$ $}}\cg%
\hbox to 20 true pt{\hskip 12 true pt\twist{${D_{6}}{A_{1}}$}\hss}\cg%
\hbox to 20 true pt{\hskip 12 true pt\twist{$A_{7}$}\hss}\cg%
\hbox to 20 true pt{\hskip 12 true pt\twist{${A_{5}}{A_{2}}$}\hss}\cg%
\hbox to 20 true pt{\hskip 12 true pt\twist{${A_{3}}{A_{3}}{A_{1}}$}\hss}\cg%
\hbox to 1.5\colpts{\hfil}\rlx{\hss{$ $}}\cg%
\hbox to 20 true pt{\hskip 12 true pt\twist{$D_{6}$}\hss}\cg%
\hbox to 20 true pt{\hskip 12 true pt\twist{$A_{6}$}\hss}\cg%
\hbox to 20 true pt{\hskip 12 true pt\twist{${A_{5}}{A_{1}}$}\hss}\cg%
\hbox to 20 true pt{\hskip 12 true pt\twist{${A_{3}}{A_{2}}{A_{1}}$}\hss}\cg%
\hbox to 20 true pt{\hskip 12 true pt\twist{${A_{4}}{A_{2}}$}\hss}\cg%
\hbox to 20 true pt{\hskip 12 true pt\twist{${D_{5}}{A_{1}}$}\hss}\cg%
\hbox to 20 true pt{\hskip 12 true pt\twist{$E_{6}$}\hss}\cg%
\eol}}\rg%
\rx{\vss\hfull{%
\rlx{\hss{${D_{6}}{A_{1}}$}}\cg%
\hbox to 20 true pt{\hss{$0$}}\cg%
\hbox to 20 true pt{\hss{$0$}}\cg%
\hbox to 20 true pt{\hss{$0$}}\cg%
\hbox to 20 true pt{\hss{$0$}}\cg%
\hbox to 1.5\colpts{\hfil}\rlx{\hss{$D_{6}$}}\cg%
\hbox to 20 true pt{\hss{$0$}}\cg%
\hbox to 20 true pt{\hss{$0$}}\cg%
\hbox to 20 true pt{\hss{$0$}}\cg%
\hbox to 20 true pt{\hss{$0$}}\cg%
\hbox to 20 true pt{\hss{$0$}}\cg%
\hbox to 20 true pt{\hss{$0$}}\cg%
\hbox to 20 true pt{\hss{$0$}}\cg%
\eol}\vss}\rg%
\rx{\vss\hfull{%
\rlx{\hss{$A_{7}$}}\cg%
\hbox to 20 true pt{\hfil}\cg%
\hbox to 20 true pt{\hss{$0$}}\cg%
\hbox to 20 true pt{\hss{$0$}}\cg%
\hbox to 20 true pt{\hss{$0$}}\cg%
\hbox to 1.5\colpts{\hfil}\rlx{\hss{$A_{6}$}}\cg%
\hbox to 20 true pt{\hfil}\cg%
\hbox to 20 true pt{\hss{$0$}}\cg%
\hbox to 20 true pt{\hss{$0$}}\cg%
\hbox to 20 true pt{\hss{$0$}}\cg%
\hbox to 20 true pt{\hss{$0$}}\cg%
\hbox to 20 true pt{\hss{$0$}}\cg%
\hbox to 20 true pt{\hss{$0$}}\cg%
\eol}\vss}\rg%
\rx{\vss\hfull{%
\rlx{\hss{${A_{5}}{A_{2}}$}}\cg%
\hbox to 20 true pt{\hfil}\cg%
\hbox to 20 true pt{\hfil}\cg%
\hbox to 20 true pt{\hss{$0$}}\cg%
\hbox to 20 true pt{\hss{$0$}}\cg%
\hbox to 1.5\colpts{\hfil}\rlx{\hss{${A_{5}}{A_{1}}$}}\cg%
\hbox to 20 true pt{\hfil}\cg%
\hbox to 20 true pt{\hfil}\cg%
\hbox to 20 true pt{\hss{$0$}}\cg%
\hbox to 20 true pt{\hss{$0$}}\cg%
\hbox to 20 true pt{\hss{$0$}}\cg%
\hbox to 20 true pt{\hss{$0$}}\cg%
\hbox to 20 true pt{\hss{$0$}}\cg%
\eol}\vss}\rg%
\rx{\vss\hfull{%
\rlx{\hss{${A_{3}}{A_{3}}{A_{1}}$}}\cg%
\hbox to 20 true pt{\hfil}\cg%
\hbox to 20 true pt{\hfil}\cg%
\hbox to 20 true pt{\hfil}\cg%
\hbox to 20 true pt{\hss{$0$}}\cg%
\hbox to 1.5\colpts{\hfil}\rlx{\hss{${A_{3}}{A_{2}}{A_{1}}$}}\cg%
\hbox to 20 true pt{\hfil}\cg%
\hbox to 20 true pt{\hfil}\cg%
\hbox to 20 true pt{\hfil}\cg%
\hbox to 20 true pt{\hss{$6$}}\cg%
\hbox to 20 true pt{\hss{$1$}}\cg%
\hbox to 20 true pt{\hss{$0$}}\cg%
\hbox to 20 true pt{\hss{$0$}}\cg%
\eol}\vss}\rg%
\rx{\vss\hfull{%
\rlx{\hss{$ $}}\cg%
\hbox to 20 true pt{\hfil}\cg%
\hbox to 20 true pt{\hfil}\cg%
\hbox to 20 true pt{\hfil}\cg%
\hbox to 20 true pt{\hfil}\cg%
\hbox to 1.5\colpts{\hfil}\rlx{\hss{${A_{4}}{A_{2}}$}}\cg%
\hbox to 20 true pt{\hfil}\cg%
\hbox to 20 true pt{\hfil}\cg%
\hbox to 20 true pt{\hfil}\cg%
\hbox to 20 true pt{\hfil}\cg%
\hbox to 20 true pt{\hss{$0$}}\cg%
\hbox to 20 true pt{\hss{$0$}}\cg%
\hbox to 20 true pt{\hss{$0$}}\cg%
\eol}\vss}\rg%
\rx{\vss\hfull{%
\rlx{\hss{$ $}}\cg%
\hbox to 20 true pt{\hfil}\cg%
\hbox to 20 true pt{\hfil}\cg%
\hbox to 20 true pt{\hfil}\cg%
\hbox to 20 true pt{\hfil}\cg%
\hbox to 1.5\colpts{\hfil}\rlx{\hss{${D_{5}}{A_{1}}$}}\cg%
\hbox to 20 true pt{\hfil}\cg%
\hbox to 20 true pt{\hfil}\cg%
\hbox to 20 true pt{\hfil}\cg%
\hbox to 20 true pt{\hfil}\cg%
\hbox to 20 true pt{\hfil}\cg%
\hbox to 20 true pt{\hss{$0$}}\cg%
\hbox to 20 true pt{\hss{$0$}}\cg%
\eol}\vss}\rg%
\rx{\vss\hfull{%
\rlx{\hss{$ $}}\cg%
\hbox to 20 true pt{\hfil}\cg%
\hbox to 20 true pt{\hfil}\cg%
\hbox to 20 true pt{\hfil}\cg%
\hbox to 20 true pt{\hfil}\cg%
\hbox to 1.5\colpts{\hfil}\rlx{\hss{$E_{6}$}}\cg%
\hbox to 20 true pt{\hfil}\cg%
\hbox to 20 true pt{\hfil}\cg%
\hbox to 20 true pt{\hfil}\cg%
\hbox to 20 true pt{\hfil}\cg%
\hbox to 20 true pt{\hfil}\cg%
\hbox to 20 true pt{\hfil}\cg%
\hbox to 20 true pt{\hss{$0$}}\cg%
\eol}\vss}\rg%
\tableclose%

\vskip 5 true pt plus 10 true pt
\goodbreak
\tableopen{Intertwining numbers $(Ind^{W}_{W_0}1,Ind^{W}_{W_1}\varepsilon)_{\down{2 true pt}{W}}$ for $W = W(E_{8})$}%
\rowpts=19 true pt%
\rowlabpts=50 true pt%
\collabpts=50 true pt%
\clx{\vss\hfull{%
\rlx{\hss{$ $}}\cg%
\hbox to 20 true pt{\hskip 12 true pt\twist{$D_{8}$}\hss}\cg%
\hbox to 20 true pt{\hskip 12 true pt\twist{$A_{8}$}\hss}\cg%
\hbox to 20 true pt{\hskip 12 true pt\twist{${A_{7}}{A_{1}}$}\hss}\cg%
\hbox to 20 true pt{\hskip 12 true pt\twist{${A_{5}}{A_{2}}{A_{1}}$}\hss}\cg%
\hbox to 20 true pt{\hskip 12 true pt\twist{${A_{4}}{A_{4}}$}\hss}\cg%
\hbox to 20 true pt{\hskip 12 true pt\twist{${D_{5}}{A_{3}}$}\hss}\cg%
\hbox to 20 true pt{\hskip 12 true pt\twist{${E_{6}}{A_{2}}$}\hss}\cg%
\hbox to 20 true pt{\hskip 12 true pt\twist{${E_{7}}{A_{1}}$}\hss}\cg%
\hbox to 1.5\colpts{\hfil}\rlx{\hss{$ $}}\cg%
\hbox to 20 true pt{\hskip 12 true pt\twist{$D_{7}$}\hss}\cg%
\hbox to 20 true pt{\hskip 12 true pt\twist{$A_{7}$}\hss}\cg%
\hbox to 20 true pt{\hskip 12 true pt\twist{${A_{6}}{A_{1}}$}\hss}\cg%
\hbox to 20 true pt{\hskip 12 true pt\twist{${A_{4}}{A_{2}}{A_{1}}$}\hss}\cg%
\hbox to 20 true pt{\hskip 12 true pt\twist{${A_{4}}{A_{3}}$}\hss}\cg%
\hbox to 20 true pt{\hskip 12 true pt\twist{${D_{5}}{A_{2}}$}\hss}\cg%
\hbox to 20 true pt{\hskip 12 true pt\twist{${E_{6}}{A_{1}}$}\hss}\cg%
\hbox to 20 true pt{\hskip 12 true pt\twist{$E_{7}$}\hss}\cg%
\eol}}\rg%
\rx{\vss\hfull{%
\rlx{\hss{$D_{8}$}}\cg%
\hbox to 20 true pt{\hss{$0$}}\cg%
\hbox to 20 true pt{\hss{$0$}}\cg%
\hbox to 20 true pt{\hss{$0$}}\cg%
\hbox to 20 true pt{\hss{$0$}}\cg%
\hbox to 20 true pt{\hss{$0$}}\cg%
\hbox to 20 true pt{\hss{$0$}}\cg%
\hbox to 20 true pt{\hss{$0$}}\cg%
\hbox to 20 true pt{\hss{$0$}}\cg%
\hbox to 1.5\colpts{\hfil}\rlx{\hss{$D_{7}$}}\cg%
\hbox to 20 true pt{\hss{$0$}}\cg%
\hbox to 20 true pt{\hss{$0$}}\cg%
\hbox to 20 true pt{\hss{$0$}}\cg%
\hbox to 20 true pt{\hss{$0$}}\cg%
\hbox to 20 true pt{\hss{$0$}}\cg%
\hbox to 20 true pt{\hss{$0$}}\cg%
\hbox to 20 true pt{\hss{$0$}}\cg%
\hbox to 20 true pt{\hss{$0$}}\cg%
\eol}\vss}\rg%
\rx{\vss\hfull{%
\rlx{\hss{$A_{8}$}}\cg%
\hbox to 20 true pt{\hfil}\cg%
\hbox to 20 true pt{\hss{$0$}}\cg%
\hbox to 20 true pt{\hss{$0$}}\cg%
\hbox to 20 true pt{\hss{$0$}}\cg%
\hbox to 20 true pt{\hss{$0$}}\cg%
\hbox to 20 true pt{\hss{$0$}}\cg%
\hbox to 20 true pt{\hss{$0$}}\cg%
\hbox to 20 true pt{\hss{$0$}}\cg%
\hbox to 1.5\colpts{\hfil}\rlx{\hss{$A_{7}$}}\cg%
\hbox to 20 true pt{\hfil}\cg%
\hbox to 20 true pt{\hss{$0$}}\cg%
\hbox to 20 true pt{\hss{$0$}}\cg%
\hbox to 20 true pt{\hss{$0$}}\cg%
\hbox to 20 true pt{\hss{$0$}}\cg%
\hbox to 20 true pt{\hss{$0$}}\cg%
\hbox to 20 true pt{\hss{$0$}}\cg%
\hbox to 20 true pt{\hss{$0$}}\cg%
\eol}\vss}\rg%
\rx{\vss\hfull{%
\rlx{\hss{${A_{7}}{A_{1}}$}}\cg%
\hbox to 20 true pt{\hfil}\cg%
\hbox to 20 true pt{\hfil}\cg%
\hbox to 20 true pt{\hss{$0$}}\cg%
\hbox to 20 true pt{\hss{$0$}}\cg%
\hbox to 20 true pt{\hss{$0$}}\cg%
\hbox to 20 true pt{\hss{$0$}}\cg%
\hbox to 20 true pt{\hss{$0$}}\cg%
\hbox to 20 true pt{\hss{$0$}}\cg%
\hbox to 1.5\colpts{\hfil}\rlx{\hss{${A_{6}}{A_{1}}$}}\cg%
\hbox to 20 true pt{\hfil}\cg%
\hbox to 20 true pt{\hfil}\cg%
\hbox to 20 true pt{\hss{$0$}}\cg%
\hbox to 20 true pt{\hss{$1$}}\cg%
\hbox to 20 true pt{\hss{$0$}}\cg%
\hbox to 20 true pt{\hss{$0$}}\cg%
\hbox to 20 true pt{\hss{$0$}}\cg%
\hbox to 20 true pt{\hss{$0$}}\cg%
\eol}\vss}\rg%
\rx{\vss\hfull{%
\rlx{\hss{${A_{5}}{A_{2}}{A_{1}}$}}\cg%
\hbox to 20 true pt{\hfil}\cg%
\hbox to 20 true pt{\hfil}\cg%
\hbox to 20 true pt{\hfil}\cg%
\hbox to 20 true pt{\hss{$1$}}\cg%
\hbox to 20 true pt{\hss{$0$}}\cg%
\hbox to 20 true pt{\hss{$0$}}\cg%
\hbox to 20 true pt{\hss{$0$}}\cg%
\hbox to 20 true pt{\hss{$0$}}\cg%
\hbox to 1.5\colpts{\hfil}\rlx{\hss{${A_{4}}{A_{2}}{A_{1}}$}}\cg%
\hbox to 20 true pt{\hfil}\cg%
\hbox to 20 true pt{\hfil}\cg%
\hbox to 20 true pt{\hfil}\cg%
\hbox to 20 true pt{\hss{$47$}}\cg%
\hbox to 20 true pt{\hss{$18$}}\cg%
\hbox to 20 true pt{\hss{$0$}}\cg%
\hbox to 20 true pt{\hss{$0$}}\cg%
\hbox to 20 true pt{\hss{$0$}}\cg%
\eol}\vss}\rg%
\rx{\vss\hfull{%
\rlx{\hss{${A_{4}}{A_{4}}$}}\cg%
\hbox to 20 true pt{\hfil}\cg%
\hbox to 20 true pt{\hfil}\cg%
\hbox to 20 true pt{\hfil}\cg%
\hbox to 20 true pt{\hfil}\cg%
\hbox to 20 true pt{\hss{$1$}}\cg%
\hbox to 20 true pt{\hss{$0$}}\cg%
\hbox to 20 true pt{\hss{$0$}}\cg%
\hbox to 20 true pt{\hss{$0$}}\cg%
\hbox to 1.5\colpts{\hfil}\rlx{\hss{${A_{4}}{A_{3}}$}}\cg%
\hbox to 20 true pt{\hfil}\cg%
\hbox to 20 true pt{\hfil}\cg%
\hbox to 20 true pt{\hfil}\cg%
\hbox to 20 true pt{\hfil}\cg%
\hbox to 20 true pt{\hss{$7$}}\cg%
\hbox to 20 true pt{\hss{$0$}}\cg%
\hbox to 20 true pt{\hss{$0$}}\cg%
\hbox to 20 true pt{\hss{$0$}}\cg%
\eol}\vss}\rg%
\rx{\vss\hfull{%
\rlx{\hss{${D_{5}}{A_{3}}$}}\cg%
\hbox to 20 true pt{\hfil}\cg%
\hbox to 20 true pt{\hfil}\cg%
\hbox to 20 true pt{\hfil}\cg%
\hbox to 20 true pt{\hfil}\cg%
\hbox to 20 true pt{\hfil}\cg%
\hbox to 20 true pt{\hss{$0$}}\cg%
\hbox to 20 true pt{\hss{$0$}}\cg%
\hbox to 20 true pt{\hss{$0$}}\cg%
\hbox to 1.5\colpts{\hfil}\rlx{\hss{${D_{5}}{A_{2}}$}}\cg%
\hbox to 20 true pt{\hfil}\cg%
\hbox to 20 true pt{\hfil}\cg%
\hbox to 20 true pt{\hfil}\cg%
\hbox to 20 true pt{\hfil}\cg%
\hbox to 20 true pt{\hfil}\cg%
\hbox to 20 true pt{\hss{$0$}}\cg%
\hbox to 20 true pt{\hss{$0$}}\cg%
\hbox to 20 true pt{\hss{$0$}}\cg%
\eol}\vss}\rg%
\rx{\vss\hfull{%
\rlx{\hss{${E_{6}}{A_{2}}$}}\cg%
\hbox to 20 true pt{\hfil}\cg%
\hbox to 20 true pt{\hfil}\cg%
\hbox to 20 true pt{\hfil}\cg%
\hbox to 20 true pt{\hfil}\cg%
\hbox to 20 true pt{\hfil}\cg%
\hbox to 20 true pt{\hfil}\cg%
\hbox to 20 true pt{\hss{$0$}}\cg%
\hbox to 20 true pt{\hss{$0$}}\cg%
\hbox to 1.5\colpts{\hfil}\rlx{\hss{${E_{6}}{A_{1}}$}}\cg%
\hbox to 20 true pt{\hfil}\cg%
\hbox to 20 true pt{\hfil}\cg%
\hbox to 20 true pt{\hfil}\cg%
\hbox to 20 true pt{\hfil}\cg%
\hbox to 20 true pt{\hfil}\cg%
\hbox to 20 true pt{\hfil}\cg%
\hbox to 20 true pt{\hss{$0$}}\cg%
\hbox to 20 true pt{\hss{$0$}}\cg%
\eol}\vss}\rg%
\rx{\vss\hfull{%
\rlx{\hss{${E_{7}}{A_{1}}$}}\cg%
\hbox to 20 true pt{\hfil}\cg%
\hbox to 20 true pt{\hfil}\cg%
\hbox to 20 true pt{\hfil}\cg%
\hbox to 20 true pt{\hfil}\cg%
\hbox to 20 true pt{\hfil}\cg%
\hbox to 20 true pt{\hfil}\cg%
\hbox to 20 true pt{\hfil}\cg%
\hbox to 20 true pt{\hss{$0$}}\cg%
\hbox to 1.5\colpts{\hfil}\rlx{\hss{$E_{7}$}}\cg%
\hbox to 20 true pt{\hfil}\cg%
\hbox to 20 true pt{\hfil}\cg%
\hbox to 20 true pt{\hfil}\cg%
\hbox to 20 true pt{\hfil}\cg%
\hbox to 20 true pt{\hfil}\cg%
\hbox to 20 true pt{\hfil}\cg%
\hbox to 20 true pt{\hfil}\cg%
\hbox to 20 true pt{\hss{$0$}}\cg%
\eol}\vss}\rg%
\tableclose%

\vskip 5 true pt plus 10 true pt
\goodbreak
\tableopen{Intertwining numbers $(Ind^{W}_{W_0}1,Ind^{W}_{W_1}\varepsilon)_{\down{2 true pt}{W}}$ for $W = W(F_{4})$}%
\rowpts=19 true pt%
\rowlabpts=105 true pt%
\collabpts=45 true pt%
\clx{\vss\hfull{%
\rlx{\hss{$ $}}\cg%
\hbox to 20 true pt{\hskip 12 true pt\twist{${C_{3}}{A_{1}}$}\hss}\cg%
\hbox to 20 true pt{\hskip 12 true pt\twist{${A_{2}}{A_{2}^{*}}$}\hss}\cg%
\hbox to 20 true pt{\hskip 12 true pt\twist{${A_{3}}{A_{1}}$}\hss}\cg%
\hbox to 20 true pt{\hskip 12 true pt\twist{$B_{4}$}\hss}\cg%
\hbox to 1.5\colpts{\hfil}\rlx{\hss{$ $}}\cg%
\hbox to 20 true pt{\hskip 12 true pt\twist{$C_{3}$}\hss}\cg%
\hbox to 20 true pt{\hskip 12 true pt\twist{${A_{1}}{A_{2}^{*}}$}\hss}\cg%
\hbox to 20 true pt{\hskip 12 true pt\twist{${A_{2}}{A_{1}^{*}}$}\hss}\cg%
\hbox to 20 true pt{\hskip 12 true pt\twist{$B_{3}$}\hss}\cg%
\eol}}\rg%
\rx{\vss\hfull{%
\rlx{\hss{${C_{3}}{A_{1}}$}}\cg%
\hbox to 20 true pt{\hss{$0$}}\cg%
\hbox to 20 true pt{\hss{$0$}}\cg%
\hbox to 20 true pt{\hss{$0$}}\cg%
\hbox to 20 true pt{\hss{$0$}}\cg%
\hbox to 1.5\colpts{\hfil}\rlx{\hss{$C_{3}$}}\cg%
\hbox to 20 true pt{\hss{$0$}}\cg%
\hbox to 20 true pt{\hss{$0$}}\cg%
\hbox to 20 true pt{\hss{$0$}}\cg%
\hbox to 20 true pt{\hss{$0$}}\cg%
\eol}\vss}\rg%
\rx{\vss\hfull{%
\rlx{\hss{${A_{2}}{A_{2}^{*}}$}}\cg%
\hbox to 20 true pt{\hfil}\cg%
\hbox to 20 true pt{\hss{$1$}}\cg%
\hbox to 20 true pt{\hss{$0$}}\cg%
\hbox to 20 true pt{\hss{$0$}}\cg%
\hbox to 1.5\colpts{\hfil}\rlx{\hss{${A_{1}}{A_{2}^{*}}$}}\cg%
\hbox to 20 true pt{\hfil}\cg%
\hbox to 20 true pt{\hss{$3$}}\cg%
\hbox to 20 true pt{\hss{$5$}}\cg%
\hbox to 20 true pt{\hss{$0$}}\cg%
\eol}\vss}\rg%
\rx{\vss\hfull{%
\rlx{\hss{${A_{3}}{A_{1}}$}}\cg%
\hbox to 20 true pt{\hfil}\cg%
\hbox to 20 true pt{\hfil}\cg%
\hbox to 20 true pt{\hss{$0$}}\cg%
\hbox to 20 true pt{\hss{$0$}}\cg%
\hbox to 1.5\colpts{\hfil}\rlx{\hss{${A_{2}}{A_{1}^{*}}$}}\cg%
\hbox to 20 true pt{\hfil}\cg%
\hbox to 20 true pt{\hfil}\cg%
\hbox to 20 true pt{\hss{$3$}}\cg%
\hbox to 20 true pt{\hss{$0$}}\cg%
\eol}\vss}\rg%
\rx{\vss\hfull{%
\rlx{\hss{$B_{4}$}}\cg%
\hbox to 20 true pt{\hfil}\cg%
\hbox to 20 true pt{\hfil}\cg%
\hbox to 20 true pt{\hfil}\cg%
\hbox to 20 true pt{\hss{$0$}}\cg%
\hbox to 1.5\colpts{\hfil}\rlx{\hss{$B_{3}$}}\cg%
\hbox to 20 true pt{\hfil}\cg%
\hbox to 20 true pt{\hfil}\cg%
\hbox to 20 true pt{\hfil}\cg%
\hbox to 20 true pt{\hss{$0$}}\cg%
\eol}\vss}\rg%
\tableclose%

\vskip 5 true pt plus 10 true pt
\goodbreak
\tableopen{Intertwining numbers $(Ind^{W}_{W_0}1,Ind^{W}_{W_1}\varepsilon)_{\down{2 true pt}{W}}$ for $W = W(G_{2})$}%
\rowpts=19 true pt%
\rowlabpts=122 true pt%
\collabpts=40 true pt%
\clx{\vss\hfull{%
\rlx{\hss{$ $}}\cg%
\hbox to 20 true pt{\hskip 12 true pt\twist{${A_{1}}{A_{1}^{*}}$}\hss}\cg%
\hbox to 20 true pt{\hskip 12 true pt\twist{$A_{2}$}\hss}\cg%
\hbox to 1.5\colpts{\hfil}\rlx{\hss{$ $}}\cg%
\hbox to 20 true pt{\hskip 12 true pt\twist{$A_{1}$}\hss}\cg%
\hbox to 20 true pt{\hskip 12 true pt\twist{$A_{1}^{*}$}\hss}\cg%
\eol}}\rg%
\rx{\vss\hfull{%
\rlx{\hss{${A_{1}}{A_{1}^{*}}$}}\cg%
\hbox to 20 true pt{\hss{$1$}}\cg%
\hbox to 20 true pt{\hss{$0$}}\cg%
\hbox to 1.5\colpts{\hfil}\rlx{\hss{$A_{1}$}}\cg%
\hbox to 20 true pt{\hss{$2$}}\cg%
\hbox to 20 true pt{\hss{$3$}}\cg%
\eol}\vss}\rg%
\rx{\vss\hfull{%
\rlx{\hss{$A_{2}$}}\cg%
\hbox to 20 true pt{\hfil}\cg%
\hbox to 20 true pt{\hss{$0$}}\cg%
\hbox to 1.5\colpts{\hfil}\rlx{\hss{$A_{1}^{*}$}}\cg%
\hbox to 20 true pt{\hfil}\cg%
\hbox to 20 true pt{\hss{$2$}}\cg%
\eol}\vss}\rg%
\tableclose%

\vskip 5 true pt plus 10 true pt
\goodbreak
\tableopen{Intertwining numbers $(Ind^{W}_{W_0}1,Ind^{W}_{W_1}\rho_{1})_{\down{2 true pt}{W}}$ for $W = W(E_{6})$}%
\rowpts=19 true pt%
\rowlabpts=120 true pt%
\collabpts=50 true pt%
\clx{\vss\hfull{%
\rlx{\hss{$ $}}\cg%
\hbox to 20 true pt{\hskip 12 true pt\twist{${A_{5}}{A_{1}}$}\hss}\cg%
\hbox to 20 true pt{\hskip 12 true pt\twist{${A_{2}}{A_{2}}{A_{2}}$}\hss}\cg%
\hbox to 1.5\colpts{\hfil}\rlx{\hss{$ $}}\cg%
\hbox to 20 true pt{\hskip 12 true pt\twist{$D_{5}$}\hss}\cg%
\hbox to 20 true pt{\hskip 12 true pt\twist{$A_{5}$}\hss}\cg%
\hbox to 20 true pt{\hskip 12 true pt\twist{${A_{4}}{A_{1}}$}\hss}\cg%
\hbox to 20 true pt{\hskip 12 true pt\twist{${A_{2}}{A_{2}}{A_{1}}$}\hss}\cg%
\eol}}\rg%
\rx{\vss\hfull{%
\rlx{\hss{${A_{5}}{A_{1}}$}}\cg%
\hbox to 20 true pt{\hss{$1$}}\cg%
\hbox to 20 true pt{\hss{$6$}}\cg%
\hbox to 1.5\colpts{\hfil}\rlx{\hss{$D_{5}$}}\cg%
\hbox to 20 true pt{\hss{$2$}}\cg%
\hbox to 20 true pt{\hss{$3$}}\cg%
\hbox to 20 true pt{\hss{$5$}}\cg%
\hbox to 20 true pt{\hss{$8$}}\cg%
\eol}\vss}\rg%
\rx{\vss\hfull{%
\rlx{\hss{${A_{2}}{A_{2}}{A_{2}}$}}\cg%
\hbox to 20 true pt{\hfil}\cg%
\hbox to 20 true pt{\hss{$18$}}\cg%
\hbox to 1.5\colpts{\hfil}\rlx{\hss{$A_{5}$}}\cg%
\hbox to 20 true pt{\hfil}\cg%
\hbox to 20 true pt{\hss{$4$}}\cg%
\hbox to 20 true pt{\hss{$9$}}\cg%
\hbox to 20 true pt{\hss{$17$}}\cg%
\eol}\vss}\rg%
\rx{\vss\hfull{%
\rlx{\hss{$ $}}\cg%
\hbox to 20 true pt{\hfil}\cg%
\hbox to 20 true pt{\hfil}\cg%
\hbox to 1.5\colpts{\hfil}\rlx{\hss{${A_{4}}{A_{1}}$}}\cg%
\hbox to 20 true pt{\hfil}\cg%
\hbox to 20 true pt{\hfil}\cg%
\hbox to 20 true pt{\hss{$18$}}\cg%
\hbox to 20 true pt{\hss{$39$}}\cg%
\eol}\vss}\rg%
\rx{\vss\hfull{%
\rlx{\hss{$ $}}\cg%
\hbox to 20 true pt{\hfil}\cg%
\hbox to 20 true pt{\hfil}\cg%
\hbox to 1.5\colpts{\hfil}\rlx{\hss{${A_{2}}{A_{2}}{A_{1}}$}}\cg%
\hbox to 20 true pt{\hfil}\cg%
\hbox to 20 true pt{\hfil}\cg%
\hbox to 20 true pt{\hfil}\cg%
\hbox to 20 true pt{\hss{$97$}}\cg%
\eol}\vss}\rg%
\tableclose%

\vskip 5 true pt plus 10 true pt
\goodbreak
\tableopen{Intertwining numbers $(Ind^{W}_{W_0}1,Ind^{W}_{W_1}\rho_{1})_{\down{2 true pt}{W}}$ for $W = W(E_{7})$}%
\rowpts=19 true pt%
\rowlabpts=85 true pt%
\collabpts=50 true pt%
\clx{\vss\hfull{%
\rlx{\hss{$ $}}\cg%
\hbox to 20 true pt{\hskip 12 true pt\twist{${D_{6}}{A_{1}}$}\hss}\cg%
\hbox to 20 true pt{\hskip 12 true pt\twist{$A_{7}$}\hss}\cg%
\hbox to 20 true pt{\hskip 12 true pt\twist{${A_{5}}{A_{2}}$}\hss}\cg%
\hbox to 20 true pt{\hskip 12 true pt\twist{${A_{3}}{A_{3}}{A_{1}}$}\hss}\cg%
\hbox to 1.5\colpts{\hfil}\rlx{\hss{$ $}}\cg%
\hbox to 20 true pt{\hskip 12 true pt\twist{$D_{6}$}\hss}\cg%
\hbox to 20 true pt{\hskip 12 true pt\twist{$A_{6}$}\hss}\cg%
\hbox to 20 true pt{\hskip 12 true pt\twist{${A_{5}}{A_{1}}$}\hss}\cg%
\hbox to 21 true pt{\hskip 12 true pt\twist{${A_{3}}{A_{2}}{A_{1}}$}\hss}\cg%
\hbox to 21 true pt{\hskip 12 true pt\twist{${A_{4}}{A_{2}}$}\hss}\cg%
\hbox to 20 true pt{\hskip 12 true pt\twist{${D_{5}}{A_{1}}$}\hss}\cg%
\hbox to 20 true pt{\hskip 12 true pt\twist{$E_{6}$}\hss}\cg%
\eol}}\rg%
\rx{\vss\hfull{%
\rlx{\hss{${D_{6}}{A_{1}}$}}\cg%
\hbox to 20 true pt{\hss{$0$}}\cg%
\hbox to 20 true pt{\hss{$1$}}\cg%
\hbox to 20 true pt{\hss{$4$}}\cg%
\hbox to 20 true pt{\hss{$5$}}\cg%
\hbox to 1.5\colpts{\hfil}\rlx{\hss{$D_{6}$}}\cg%
\hbox to 20 true pt{\hss{$3$}}\cg%
\hbox to 20 true pt{\hss{$6$}}\cg%
\hbox to 20 true pt{\hss{$11$}}\cg%
\hbox to 21 true pt{\hss{$21$}}\cg%
\hbox to 21 true pt{\hss{$14$}}\cg%
\hbox to 20 true pt{\hss{$8$}}\cg%
\hbox to 20 true pt{\hss{$3$}}\cg%
\eol}\vss}\rg%
\rx{\vss\hfull{%
\rlx{\hss{$A_{7}$}}\cg%
\hbox to 20 true pt{\hfil}\cg%
\hbox to 20 true pt{\hss{$0$}}\cg%
\hbox to 20 true pt{\hss{$4$}}\cg%
\hbox to 20 true pt{\hss{$7$}}\cg%
\hbox to 1.5\colpts{\hfil}\rlx{\hss{$A_{6}$}}\cg%
\hbox to 20 true pt{\hfil}\cg%
\hbox to 20 true pt{\hss{$12$}}\cg%
\hbox to 20 true pt{\hss{$26$}}\cg%
\hbox to 21 true pt{\hss{$62$}}\cg%
\hbox to 21 true pt{\hss{$36$}}\cg%
\hbox to 20 true pt{\hss{$16$}}\cg%
\hbox to 20 true pt{\hss{$4$}}\cg%
\eol}\vss}\rg%
\rx{\vss\hfull{%
\rlx{\hss{${A_{5}}{A_{2}}$}}\cg%
\hbox to 20 true pt{\hfil}\cg%
\hbox to 20 true pt{\hfil}\cg%
\hbox to 20 true pt{\hss{$12$}}\cg%
\hbox to 20 true pt{\hss{$31$}}\cg%
\hbox to 1.5\colpts{\hfil}\rlx{\hss{${A_{5}}{A_{1}}$}}\cg%
\hbox to 20 true pt{\hfil}\cg%
\hbox to 20 true pt{\hfil}\cg%
\hbox to 20 true pt{\hss{$57$}}\cg%
\hbox to 21 true pt{\hss{$156$}}\cg%
\hbox to 21 true pt{\hss{$86$}}\cg%
\hbox to 20 true pt{\hss{$33$}}\cg%
\hbox to 20 true pt{\hss{$7$}}\cg%
\eol}\vss}\rg%
\rx{\vss\hfull{%
\rlx{\hss{${A_{3}}{A_{3}}{A_{1}}$}}\cg%
\hbox to 20 true pt{\hfil}\cg%
\hbox to 20 true pt{\hfil}\cg%
\hbox to 20 true pt{\hfil}\cg%
\hbox to 20 true pt{\hss{$70$}}\cg%
\hbox to 1.5\colpts{\hfil}\rlx{\hss{${A_{3}}{A_{2}}{A_{1}}$}}\cg%
\hbox to 20 true pt{\hfil}\cg%
\hbox to 20 true pt{\hfil}\cg%
\hbox to 20 true pt{\hfil}\cg%
\hbox to 21 true pt{\hss{$517$}}\cg%
\hbox to 21 true pt{\hss{$258$}}\cg%
\hbox to 20 true pt{\hss{$77$}}\cg%
\hbox to 20 true pt{\hss{$12$}}\cg%
\eol}\vss}\rg%
\rx{\vss\hfull{%
\rlx{\hss{$ $}}\cg%
\hbox to 20 true pt{\hfil}\cg%
\hbox to 20 true pt{\hfil}\cg%
\hbox to 20 true pt{\hfil}\cg%
\hbox to 20 true pt{\hfil}\cg%
\hbox to 1.5\colpts{\hfil}\rlx{\hss{${A_{4}}{A_{2}}$}}\cg%
\hbox to 20 true pt{\hfil}\cg%
\hbox to 20 true pt{\hfil}\cg%
\hbox to 20 true pt{\hfil}\cg%
\hbox to 21 true pt{\hfil}\cg%
\hbox to 21 true pt{\hss{$134$}}\cg%
\hbox to 20 true pt{\hss{$46$}}\cg%
\hbox to 20 true pt{\hss{$8$}}\cg%
\eol}\vss}\rg%
\rx{\vss\hfull{%
\rlx{\hss{$ $}}\cg%
\hbox to 20 true pt{\hfil}\cg%
\hbox to 20 true pt{\hfil}\cg%
\hbox to 20 true pt{\hfil}\cg%
\hbox to 20 true pt{\hfil}\cg%
\hbox to 1.5\colpts{\hfil}\rlx{\hss{${D_{5}}{A_{1}}$}}\cg%
\hbox to 20 true pt{\hfil}\cg%
\hbox to 20 true pt{\hfil}\cg%
\hbox to 20 true pt{\hfil}\cg%
\hbox to 21 true pt{\hfil}\cg%
\hbox to 21 true pt{\hfil}\cg%
\hbox to 20 true pt{\hss{$21$}}\cg%
\hbox to 20 true pt{\hss{$6$}}\cg%
\eol}\vss}\rg%
\rx{\vss\hfull{%
\rlx{\hss{$ $}}\cg%
\hbox to 20 true pt{\hfil}\cg%
\hbox to 20 true pt{\hfil}\cg%
\hbox to 20 true pt{\hfil}\cg%
\hbox to 20 true pt{\hfil}\cg%
\hbox to 1.5\colpts{\hfil}\rlx{\hss{$E_{6}$}}\cg%
\hbox to 20 true pt{\hfil}\cg%
\hbox to 20 true pt{\hfil}\cg%
\hbox to 20 true pt{\hfil}\cg%
\hbox to 21 true pt{\hfil}\cg%
\hbox to 21 true pt{\hfil}\cg%
\hbox to 20 true pt{\hfil}\cg%
\hbox to 20 true pt{\hss{$2$}}\cg%
\eol}\vss}\rg%
\tableclose%

\vskip 5 true pt plus 10 true pt
\goodbreak
\tableopen{Intertwining numbers $(Ind^{W}_{W_0}1,Ind^{W}_{W_1}\rho_{1})_{\down{2 true pt}{W}}$ for $W = W(E_{8})$}%
\rowpts=19 true pt%
\rowlabpts=50 true pt%
\collabpts=50 true pt%
\clx{\vss\hfull{%
\rlx{\hss{$ $}}\cg%
\hbox to 20 true pt{\hskip 12 true pt\twist{$D_{8}$}\hss}\cg%
\hbox to 20 true pt{\hskip 12 true pt\twist{$A_{8}$}\hss}\cg%
\hbox to 20 true pt{\hskip 12 true pt\twist{${A_{7}}{A_{1}}$}\hss}\cg%
\hbox to 21 true pt{\hskip 12 true pt\twist{${A_{5}}{A_{2}}{A_{1}}$}\hss}\cg%
\hbox to 21 true pt{\hskip 12 true pt\twist{${A_{4}}{A_{4}}$}\hss}\cg%
\hbox to 21 true pt{\hskip 12 true pt\twist{${D_{5}}{A_{3}}$}\hss}\cg%
\hbox to 20 true pt{\hskip 12 true pt\twist{${E_{6}}{A_{2}}$}\hss}\cg%
\hbox to 20 true pt{\hskip 12 true pt\twist{${E_{7}}{A_{1}}$}\hss}\cg%
\hbox to 1.5\colpts{\hfil}\rlx{\hss{$ $}}\cg%
\hbox to 20 true pt{\hskip 12 true pt\twist{$D_{7}$}\hss}\cg%
\hbox to 20 true pt{\hskip 12 true pt\twist{$A_{7}$}\hss}\cg%
\hbox to 21 true pt{\hskip 12 true pt\twist{${A_{6}}{A_{1}}$}\hss}\cg%
\hbox to 28 true pt{\hskip 16 true pt\twist{${A_{4}}{A_{2}}{A_{1}}$}\hss}\cg%
\hbox to 28 true pt{\hskip 16 true pt\twist{${A_{4}}{A_{3}}$}\hss}\cg%
\hbox to 28 true pt{\hskip 16 true pt\twist{${D_{5}}{A_{2}}$}\hss}\cg%
\hbox to 21 true pt{\hskip 12 true pt\twist{${E_{6}}{A_{1}}$}\hss}\cg%
\hbox to 20 true pt{\hskip 12 true pt\twist{$E_{7}$}\hss}\cg%
\eol}}\rg%
\rx{\vss\hfull{%
\rlx{\hss{$D_{8}$}}\cg%
\hbox to 20 true pt{\hss{$0$}}\cg%
\hbox to 20 true pt{\hss{$2$}}\cg%
\hbox to 20 true pt{\hss{$3$}}\cg%
\hbox to 21 true pt{\hss{$7$}}\cg%
\hbox to 21 true pt{\hss{$6$}}\cg%
\hbox to 21 true pt{\hss{$3$}}\cg%
\hbox to 20 true pt{\hss{$2$}}\cg%
\hbox to 20 true pt{\hss{$0$}}\cg%
\hbox to 1.5\colpts{\hfil}\rlx{\hss{$D_{7}$}}\cg%
\hbox to 20 true pt{\hss{$10$}}\cg%
\hbox to 20 true pt{\hss{$24$}}\cg%
\hbox to 21 true pt{\hss{$49$}}\cg%
\hbox to 28 true pt{\hss{$122$}}\cg%
\hbox to 28 true pt{\hss{$81$}}\cg%
\hbox to 28 true pt{\hss{$46$}}\cg%
\hbox to 21 true pt{\hss{$20$}}\cg%
\hbox to 20 true pt{\hss{$5$}}\cg%
\eol}\vss}\rg%
\rx{\vss\hfull{%
\rlx{\hss{$A_{8}$}}\cg%
\hbox to 20 true pt{\hfil}\cg%
\hbox to 20 true pt{\hss{$4$}}\cg%
\hbox to 20 true pt{\hss{$13$}}\cg%
\hbox to 21 true pt{\hss{$38$}}\cg%
\hbox to 21 true pt{\hss{$28$}}\cg%
\hbox to 21 true pt{\hss{$16$}}\cg%
\hbox to 20 true pt{\hss{$6$}}\cg%
\hbox to 20 true pt{\hss{$2$}}\cg%
\hbox to 1.5\colpts{\hfil}\rlx{\hss{$A_{7}$}}\cg%
\hbox to 20 true pt{\hfil}\cg%
\hbox to 20 true pt{\hss{$69$}}\cg%
\hbox to 21 true pt{\hss{$163$}}\cg%
\hbox to 28 true pt{\hss{$517$}}\cg%
\hbox to 28 true pt{\hss{$323$}}\cg%
\hbox to 28 true pt{\hss{$150$}}\cg%
\hbox to 21 true pt{\hss{$48$}}\cg%
\hbox to 20 true pt{\hss{$8$}}\cg%
\eol}\vss}\rg%
\rx{\vss\hfull{%
\rlx{\hss{${A_{7}}{A_{1}}$}}\cg%
\hbox to 20 true pt{\hfil}\cg%
\hbox to 20 true pt{\hfil}\cg%
\hbox to 20 true pt{\hss{$28$}}\cg%
\hbox to 21 true pt{\hss{$105$}}\cg%
\hbox to 21 true pt{\hss{$74$}}\cg%
\hbox to 21 true pt{\hss{$38$}}\cg%
\hbox to 20 true pt{\hss{$15$}}\cg%
\hbox to 20 true pt{\hss{$3$}}\cg%
\hbox to 1.5\colpts{\hfil}\rlx{\hss{${A_{6}}{A_{1}}$}}\cg%
\hbox to 20 true pt{\hfil}\cg%
\hbox to 20 true pt{\hfil}\cg%
\hbox to 21 true pt{\hss{$412$}}\cg%
\hbox to 28 true pt{\hss{$1481$}}\cg%
\hbox to 28 true pt{\hss{$890$}}\cg%
\hbox to 28 true pt{\hss{$382$}}\cg%
\hbox to 21 true pt{\hss{$107$}}\cg%
\hbox to 20 true pt{\hss{$15$}}\cg%
\eol}\vss}\rg%
\rx{\vss\hfull{%
\rlx{\hss{${A_{5}}{A_{2}}{A_{1}}$}}\cg%
\hbox to 20 true pt{\hfil}\cg%
\hbox to 20 true pt{\hfil}\cg%
\hbox to 20 true pt{\hfil}\cg%
\hbox to 21 true pt{\hss{$468$}}\cg%
\hbox to 21 true pt{\hss{$318$}}\cg%
\hbox to 21 true pt{\hss{$148$}}\cg%
\hbox to 20 true pt{\hss{$45$}}\cg%
\hbox to 20 true pt{\hss{$7$}}\cg%
\hbox to 1.5\colpts{\hfil}\rlx{\hss{${A_{4}}{A_{2}}{A_{1}}$}}\cg%
\hbox to 20 true pt{\hfil}\cg%
\hbox to 20 true pt{\hfil}\cg%
\hbox to 21 true pt{\hfil}\cg%
\hbox to 28 true pt{\hss{$6493$}}\cg%
\hbox to 28 true pt{\hss{$3706$}}\cg%
\hbox to 28 true pt{\hss{$1353$}}\cg%
\hbox to 21 true pt{\hss{$292$}}\cg%
\hbox to 20 true pt{\hss{$29$}}\cg%
\eol}\vss}\rg%
\rx{\vss\hfull{%
\rlx{\hss{${A_{4}}{A_{4}}$}}\cg%
\hbox to 20 true pt{\hfil}\cg%
\hbox to 20 true pt{\hfil}\cg%
\hbox to 20 true pt{\hfil}\cg%
\hbox to 21 true pt{\hfil}\cg%
\hbox to 21 true pt{\hss{$216$}}\cg%
\hbox to 21 true pt{\hss{$104$}}\cg%
\hbox to 20 true pt{\hss{$32$}}\cg%
\hbox to 20 true pt{\hss{$6$}}\cg%
\hbox to 1.5\colpts{\hfil}\rlx{\hss{${A_{4}}{A_{3}}$}}\cg%
\hbox to 20 true pt{\hfil}\cg%
\hbox to 20 true pt{\hfil}\cg%
\hbox to 21 true pt{\hfil}\cg%
\hbox to 28 true pt{\hfil}\cg%
\hbox to 28 true pt{\hss{$2155$}}\cg%
\hbox to 28 true pt{\hss{$812$}}\cg%
\hbox to 21 true pt{\hss{$185$}}\cg%
\hbox to 20 true pt{\hss{$20$}}\cg%
\eol}\vss}\rg%
\rx{\vss\hfull{%
\rlx{\hss{${D_{5}}{A_{3}}$}}\cg%
\hbox to 20 true pt{\hfil}\cg%
\hbox to 20 true pt{\hfil}\cg%
\hbox to 20 true pt{\hfil}\cg%
\hbox to 21 true pt{\hfil}\cg%
\hbox to 21 true pt{\hfil}\cg%
\hbox to 21 true pt{\hss{$50$}}\cg%
\hbox to 20 true pt{\hss{$20$}}\cg%
\hbox to 20 true pt{\hss{$3$}}\cg%
\hbox to 1.5\colpts{\hfil}\rlx{\hss{${D_{5}}{A_{2}}$}}\cg%
\hbox to 20 true pt{\hfil}\cg%
\hbox to 20 true pt{\hfil}\cg%
\hbox to 21 true pt{\hfil}\cg%
\hbox to 28 true pt{\hfil}\cg%
\hbox to 28 true pt{\hfil}\cg%
\hbox to 28 true pt{\hss{$362$}}\cg%
\hbox to 21 true pt{\hss{$107$}}\cg%
\hbox to 20 true pt{\hss{$15$}}\cg%
\eol}\vss}\rg%
\rx{\vss\hfull{%
\rlx{\hss{${E_{6}}{A_{2}}$}}\cg%
\hbox to 20 true pt{\hfil}\cg%
\hbox to 20 true pt{\hfil}\cg%
\hbox to 20 true pt{\hfil}\cg%
\hbox to 21 true pt{\hfil}\cg%
\hbox to 21 true pt{\hfil}\cg%
\hbox to 21 true pt{\hfil}\cg%
\hbox to 20 true pt{\hss{$8$}}\cg%
\hbox to 20 true pt{\hss{$3$}}\cg%
\hbox to 1.5\colpts{\hfil}\rlx{\hss{${E_{6}}{A_{1}}$}}\cg%
\hbox to 20 true pt{\hfil}\cg%
\hbox to 20 true pt{\hfil}\cg%
\hbox to 21 true pt{\hfil}\cg%
\hbox to 28 true pt{\hfil}\cg%
\hbox to 28 true pt{\hfil}\cg%
\hbox to 28 true pt{\hfil}\cg%
\hbox to 21 true pt{\hss{$46$}}\cg%
\hbox to 20 true pt{\hss{$10$}}\cg%
\eol}\vss}\rg%
\rx{\vss\hfull{%
\rlx{\hss{${E_{7}}{A_{1}}$}}\cg%
\hbox to 20 true pt{\hfil}\cg%
\hbox to 20 true pt{\hfil}\cg%
\hbox to 20 true pt{\hfil}\cg%
\hbox to 21 true pt{\hfil}\cg%
\hbox to 21 true pt{\hfil}\cg%
\hbox to 21 true pt{\hfil}\cg%
\hbox to 20 true pt{\hfil}\cg%
\hbox to 20 true pt{\hss{$0$}}\cg%
\hbox to 1.5\colpts{\hfil}\rlx{\hss{$E_{7}$}}\cg%
\hbox to 20 true pt{\hfil}\cg%
\hbox to 20 true pt{\hfil}\cg%
\hbox to 21 true pt{\hfil}\cg%
\hbox to 28 true pt{\hfil}\cg%
\hbox to 28 true pt{\hfil}\cg%
\hbox to 28 true pt{\hfil}\cg%
\hbox to 21 true pt{\hfil}\cg%
\hbox to 20 true pt{\hss{$3$}}\cg%
\eol}\vss}\rg%
\tableclose%

\vskip 5 true pt plus 10 true pt
\goodbreak
\tableopen{Intertwining numbers $(Ind^{W}_{W_0}1,Ind^{W}_{W_1}\rho_{1})_{\down{2 true pt}{W}}$ for $W = W(F_{4})$}%
\rowpts=19 true pt%
\rowlabpts=105 true pt%
\collabpts=45 true pt%
\clx{\vss\hfull{%
\rlx{\hss{$ $}}\cg%
\hbox to 20 true pt{\hskip 12 true pt\twist{${C_{3}}{A_{1}}$}\hss}\cg%
\hbox to 20 true pt{\hskip 12 true pt\twist{${A_{2}}{A_{2}^{*}}$}\hss}\cg%
\hbox to 20 true pt{\hskip 12 true pt\twist{${A_{3}}{A_{1}}$}\hss}\cg%
\hbox to 20 true pt{\hskip 12 true pt\twist{$B_{4}$}\hss}\cg%
\hbox to 1.5\colpts{\hfil}\rlx{\hss{$ $}}\cg%
\hbox to 20 true pt{\hskip 12 true pt\twist{$C_{3}$}\hss}\cg%
\hbox to 20 true pt{\hskip 12 true pt\twist{${A_{1}}{A_{2}^{*}}$}\hss}\cg%
\hbox to 20 true pt{\hskip 12 true pt\twist{${A_{2}}{A_{1}^{*}}$}\hss}\cg%
\hbox to 20 true pt{\hskip 12 true pt\twist{$B_{3}$}\hss}\cg%
\eol}}\rg%
\rx{\vss\hfull{%
\rlx{\hss{${C_{3}}{A_{1}}$}}\cg%
\hbox to 20 true pt{\hss{$0$}}\cg%
\hbox to 20 true pt{\hss{$3$}}\cg%
\hbox to 20 true pt{\hss{$1$}}\cg%
\hbox to 20 true pt{\hss{$0$}}\cg%
\hbox to 1.5\colpts{\hfil}\rlx{\hss{$C_{3}$}}\cg%
\hbox to 20 true pt{\hss{$3$}}\cg%
\hbox to 20 true pt{\hss{$9$}}\cg%
\hbox to 20 true pt{\hss{$8$}}\cg%
\hbox to 20 true pt{\hss{$3$}}\cg%
\eol}\vss}\rg%
\rx{\vss\hfull{%
\rlx{\hss{${A_{2}}{A_{2}^{*}}$}}\cg%
\hbox to 20 true pt{\hfil}\cg%
\hbox to 20 true pt{\hss{$4$}}\cg%
\hbox to 20 true pt{\hss{$5$}}\cg%
\hbox to 20 true pt{\hss{$1$}}\cg%
\hbox to 1.5\colpts{\hfil}\rlx{\hss{${A_{1}}{A_{2}^{*}}$}}\cg%
\hbox to 20 true pt{\hfil}\cg%
\hbox to 20 true pt{\hss{$29$}}\cg%
\hbox to 20 true pt{\hss{$27$}}\cg%
\hbox to 20 true pt{\hss{$8$}}\cg%
\eol}\vss}\rg%
\rx{\vss\hfull{%
\rlx{\hss{${A_{3}}{A_{1}}$}}\cg%
\hbox to 20 true pt{\hfil}\cg%
\hbox to 20 true pt{\hfil}\cg%
\hbox to 20 true pt{\hss{$4$}}\cg%
\hbox to 20 true pt{\hss{$1$}}\cg%
\hbox to 1.5\colpts{\hfil}\rlx{\hss{${A_{2}}{A_{1}^{*}}$}}\cg%
\hbox to 20 true pt{\hfil}\cg%
\hbox to 20 true pt{\hfil}\cg%
\hbox to 20 true pt{\hss{$29$}}\cg%
\hbox to 20 true pt{\hss{$9$}}\cg%
\eol}\vss}\rg%
\rx{\vss\hfull{%
\rlx{\hss{$B_{4}$}}\cg%
\hbox to 20 true pt{\hfil}\cg%
\hbox to 20 true pt{\hfil}\cg%
\hbox to 20 true pt{\hfil}\cg%
\hbox to 20 true pt{\hss{$0$}}\cg%
\hbox to 1.5\colpts{\hfil}\rlx{\hss{$B_{3}$}}\cg%
\hbox to 20 true pt{\hfil}\cg%
\hbox to 20 true pt{\hfil}\cg%
\hbox to 20 true pt{\hfil}\cg%
\hbox to 20 true pt{\hss{$3$}}\cg%
\eol}\vss}\rg%
\tableclose%

\vskip 5 true pt plus 10 true pt
\goodbreak
\tableopen{Intertwining numbers $(Ind^{W}_{W_0}1,Ind^{W}_{W_1}\rho_{1})_{\down{2 true pt}{W}}$ for $W = W(G_{2})$}%
\rowpts=19 true pt%
\rowlabpts=122 true pt%
\collabpts=40 true pt%
\clx{\vss\hfull{%
\rlx{\hss{$ $}}\cg%
\hbox to 20 true pt{\hskip 12 true pt\twist{${A_{1}}{A_{1}^{*}}$}\hss}\cg%
\hbox to 20 true pt{\hskip 12 true pt\twist{$A_{2}$}\hss}\cg%
\hbox to 1.5\colpts{\hfil}\rlx{\hss{$ $}}\cg%
\hbox to 20 true pt{\hskip 12 true pt\twist{$A_{1}$}\hss}\cg%
\hbox to 20 true pt{\hskip 12 true pt\twist{$A_{1}^{*}$}\hss}\cg%
\eol}}\rg%
\rx{\vss\hfull{%
\rlx{\hss{${A_{1}}{A_{1}^{*}}$}}\cg%
\hbox to 20 true pt{\hss{$0$}}\cg%
\hbox to 20 true pt{\hss{$1$}}\cg%
\hbox to 1.5\colpts{\hfil}\rlx{\hss{$A_{1}$}}\cg%
\hbox to 20 true pt{\hss{$2$}}\cg%
\hbox to 20 true pt{\hss{$3$}}\cg%
\eol}\vss}\rg%
\rx{\vss\hfull{%
\rlx{\hss{$A_{2}$}}\cg%
\hbox to 20 true pt{\hfil}\cg%
\hbox to 20 true pt{\hss{$0$}}\cg%
\hbox to 1.5\colpts{\hfil}\rlx{\hss{$A_{1}^{*}$}}\cg%
\hbox to 20 true pt{\hfil}\cg%
\hbox to 20 true pt{\hss{$2$}}\cg%
\eol}\vss}\rg%
\tableclose%

\vskip 5 true pt plus 10 true pt
\goodbreak
\tableopen{Intertwining numbers $(Ind^{W}_{W_0}\varepsilon,Ind^{W}_{W_1}\rho_{1})_{\down{2 true pt}{W}}$ for $W = W(E_{6})$}%
\rowpts=19 true pt%
\rowlabpts=120 true pt%
\collabpts=50 true pt%
\clx{\vss\hfull{%
\rlx{\hss{$ $}}\cg%
\hbox to 20 true pt{\hskip 12 true pt\twist{${A_{5}}{A_{1}}$}\hss}\cg%
\hbox to 20 true pt{\hskip 12 true pt\twist{${A_{2}}{A_{2}}{A_{2}}$}\hss}\cg%
\hbox to 1.5\colpts{\hfil}\rlx{\hss{$ $}}\cg%
\hbox to 20 true pt{\hskip 12 true pt\twist{$D_{5}$}\hss}\cg%
\hbox to 20 true pt{\hskip 12 true pt\twist{$A_{5}$}\hss}\cg%
\hbox to 20 true pt{\hskip 12 true pt\twist{${A_{4}}{A_{1}}$}\hss}\cg%
\hbox to 20 true pt{\hskip 12 true pt\twist{${A_{2}}{A_{2}}{A_{1}}$}\hss}\cg%
\eol}}\rg%
\rx{\vss\hfull{%
\rlx{\hss{${A_{5}}{A_{1}}$}}\cg%
\hbox to 20 true pt{\hss{$0$}}\cg%
\hbox to 20 true pt{\hss{$0$}}\cg%
\hbox to 1.5\colpts{\hfil}\rlx{\hss{$D_{5}$}}\cg%
\hbox to 20 true pt{\hss{$0$}}\cg%
\hbox to 20 true pt{\hss{$0$}}\cg%
\hbox to 20 true pt{\hss{$0$}}\cg%
\hbox to 20 true pt{\hss{$0$}}\cg%
\eol}\vss}\rg%
\rx{\vss\hfull{%
\rlx{\hss{${A_{2}}{A_{2}}{A_{2}}$}}\cg%
\hbox to 20 true pt{\hfil}\cg%
\hbox to 20 true pt{\hss{$0$}}\cg%
\hbox to 1.5\colpts{\hfil}\rlx{\hss{$A_{5}$}}\cg%
\hbox to 20 true pt{\hfil}\cg%
\hbox to 20 true pt{\hss{$0$}}\cg%
\hbox to 20 true pt{\hss{$0$}}\cg%
\hbox to 20 true pt{\hss{$0$}}\cg%
\eol}\vss}\rg%
\rx{\vss\hfull{%
\rlx{\hss{$ $}}\cg%
\hbox to 20 true pt{\hfil}\cg%
\hbox to 20 true pt{\hfil}\cg%
\hbox to 1.5\colpts{\hfil}\rlx{\hss{${A_{4}}{A_{1}}$}}\cg%
\hbox to 20 true pt{\hfil}\cg%
\hbox to 20 true pt{\hfil}\cg%
\hbox to 20 true pt{\hss{$0$}}\cg%
\hbox to 20 true pt{\hss{$2$}}\cg%
\eol}\vss}\rg%
\rx{\vss\hfull{%
\rlx{\hss{$ $}}\cg%
\hbox to 20 true pt{\hfil}\cg%
\hbox to 20 true pt{\hfil}\cg%
\hbox to 1.5\colpts{\hfil}\rlx{\hss{${A_{2}}{A_{2}}{A_{1}}$}}\cg%
\hbox to 20 true pt{\hfil}\cg%
\hbox to 20 true pt{\hfil}\cg%
\hbox to 20 true pt{\hfil}\cg%
\hbox to 20 true pt{\hss{$19$}}\cg%
\eol}\vss}\rg%
\tableclose%

\vskip 5 true pt plus 10 true pt
\goodbreak
\tableopen{Intertwining numbers $(Ind^{W}_{W_0}\varepsilon,Ind^{W}_{W_1}\rho_{1})_{\down{2 true pt}{W}}$ for $W = W(E_{7})$}%
\rowpts=19 true pt%
\rowlabpts=85 true pt%
\collabpts=50 true pt%
\clx{\vss\hfull{%
\rlx{\hss{$ $}}\cg%
\hbox to 20 true pt{\hskip 12 true pt\twist{${D_{6}}{A_{1}}$}\hss}\cg%
\hbox to 20 true pt{\hskip 12 true pt\twist{$A_{7}$}\hss}\cg%
\hbox to 20 true pt{\hskip 12 true pt\twist{${A_{5}}{A_{2}}$}\hss}\cg%
\hbox to 20 true pt{\hskip 12 true pt\twist{${A_{3}}{A_{3}}{A_{1}}$}\hss}\cg%
\hbox to 1.5\colpts{\hfil}\rlx{\hss{$ $}}\cg%
\hbox to 20 true pt{\hskip 12 true pt\twist{$D_{6}$}\hss}\cg%
\hbox to 20 true pt{\hskip 12 true pt\twist{$A_{6}$}\hss}\cg%
\hbox to 20 true pt{\hskip 12 true pt\twist{${A_{5}}{A_{1}}$}\hss}\cg%
\hbox to 20 true pt{\hskip 12 true pt\twist{${A_{3}}{A_{2}}{A_{1}}$}\hss}\cg%
\hbox to 20 true pt{\hskip 12 true pt\twist{${A_{4}}{A_{2}}$}\hss}\cg%
\hbox to 20 true pt{\hskip 12 true pt\twist{${D_{5}}{A_{1}}$}\hss}\cg%
\hbox to 20 true pt{\hskip 12 true pt\twist{$E_{6}$}\hss}\cg%
\eol}}\rg%
\rx{\vss\hfull{%
\rlx{\hss{${D_{6}}{A_{1}}$}}\cg%
\hbox to 20 true pt{\hss{$0$}}\cg%
\hbox to 20 true pt{\hss{$0$}}\cg%
\hbox to 20 true pt{\hss{$0$}}\cg%
\hbox to 20 true pt{\hss{$0$}}\cg%
\hbox to 1.5\colpts{\hfil}\rlx{\hss{$D_{6}$}}\cg%
\hbox to 20 true pt{\hss{$0$}}\cg%
\hbox to 20 true pt{\hss{$0$}}\cg%
\hbox to 20 true pt{\hss{$0$}}\cg%
\hbox to 20 true pt{\hss{$0$}}\cg%
\hbox to 20 true pt{\hss{$0$}}\cg%
\hbox to 20 true pt{\hss{$0$}}\cg%
\hbox to 20 true pt{\hss{$0$}}\cg%
\eol}\vss}\rg%
\rx{\vss\hfull{%
\rlx{\hss{$A_{7}$}}\cg%
\hbox to 20 true pt{\hfil}\cg%
\hbox to 20 true pt{\hss{$0$}}\cg%
\hbox to 20 true pt{\hss{$0$}}\cg%
\hbox to 20 true pt{\hss{$0$}}\cg%
\hbox to 1.5\colpts{\hfil}\rlx{\hss{$A_{6}$}}\cg%
\hbox to 20 true pt{\hfil}\cg%
\hbox to 20 true pt{\hss{$0$}}\cg%
\hbox to 20 true pt{\hss{$0$}}\cg%
\hbox to 20 true pt{\hss{$0$}}\cg%
\hbox to 20 true pt{\hss{$0$}}\cg%
\hbox to 20 true pt{\hss{$0$}}\cg%
\hbox to 20 true pt{\hss{$0$}}\cg%
\eol}\vss}\rg%
\rx{\vss\hfull{%
\rlx{\hss{${A_{5}}{A_{2}}$}}\cg%
\hbox to 20 true pt{\hfil}\cg%
\hbox to 20 true pt{\hfil}\cg%
\hbox to 20 true pt{\hss{$0$}}\cg%
\hbox to 20 true pt{\hss{$0$}}\cg%
\hbox to 1.5\colpts{\hfil}\rlx{\hss{${A_{5}}{A_{1}}$}}\cg%
\hbox to 20 true pt{\hfil}\cg%
\hbox to 20 true pt{\hfil}\cg%
\hbox to 20 true pt{\hss{$0$}}\cg%
\hbox to 20 true pt{\hss{$3$}}\cg%
\hbox to 20 true pt{\hss{$0$}}\cg%
\hbox to 20 true pt{\hss{$0$}}\cg%
\hbox to 20 true pt{\hss{$0$}}\cg%
\eol}\vss}\rg%
\rx{\vss\hfull{%
\rlx{\hss{${A_{3}}{A_{3}}{A_{1}}$}}\cg%
\hbox to 20 true pt{\hfil}\cg%
\hbox to 20 true pt{\hfil}\cg%
\hbox to 20 true pt{\hfil}\cg%
\hbox to 20 true pt{\hss{$3$}}\cg%
\hbox to 1.5\colpts{\hfil}\rlx{\hss{${A_{3}}{A_{2}}{A_{1}}$}}\cg%
\hbox to 20 true pt{\hfil}\cg%
\hbox to 20 true pt{\hfil}\cg%
\hbox to 20 true pt{\hfil}\cg%
\hbox to 20 true pt{\hss{$60$}}\cg%
\hbox to 20 true pt{\hss{$14$}}\cg%
\hbox to 20 true pt{\hss{$0$}}\cg%
\hbox to 20 true pt{\hss{$0$}}\cg%
\eol}\vss}\rg%
\rx{\vss\hfull{%
\rlx{\hss{$ $}}\cg%
\hbox to 20 true pt{\hfil}\cg%
\hbox to 20 true pt{\hfil}\cg%
\hbox to 20 true pt{\hfil}\cg%
\hbox to 20 true pt{\hfil}\cg%
\hbox to 1.5\colpts{\hfil}\rlx{\hss{${A_{4}}{A_{2}}$}}\cg%
\hbox to 20 true pt{\hfil}\cg%
\hbox to 20 true pt{\hfil}\cg%
\hbox to 20 true pt{\hfil}\cg%
\hbox to 20 true pt{\hfil}\cg%
\hbox to 20 true pt{\hss{$2$}}\cg%
\hbox to 20 true pt{\hss{$0$}}\cg%
\hbox to 20 true pt{\hss{$0$}}\cg%
\eol}\vss}\rg%
\rx{\vss\hfull{%
\rlx{\hss{$ $}}\cg%
\hbox to 20 true pt{\hfil}\cg%
\hbox to 20 true pt{\hfil}\cg%
\hbox to 20 true pt{\hfil}\cg%
\hbox to 20 true pt{\hfil}\cg%
\hbox to 1.5\colpts{\hfil}\rlx{\hss{${D_{5}}{A_{1}}$}}\cg%
\hbox to 20 true pt{\hfil}\cg%
\hbox to 20 true pt{\hfil}\cg%
\hbox to 20 true pt{\hfil}\cg%
\hbox to 20 true pt{\hfil}\cg%
\hbox to 20 true pt{\hfil}\cg%
\hbox to 20 true pt{\hss{$0$}}\cg%
\hbox to 20 true pt{\hss{$0$}}\cg%
\eol}\vss}\rg%
\rx{\vss\hfull{%
\rlx{\hss{$ $}}\cg%
\hbox to 20 true pt{\hfil}\cg%
\hbox to 20 true pt{\hfil}\cg%
\hbox to 20 true pt{\hfil}\cg%
\hbox to 20 true pt{\hfil}\cg%
\hbox to 1.5\colpts{\hfil}\rlx{\hss{$E_{6}$}}\cg%
\hbox to 20 true pt{\hfil}\cg%
\hbox to 20 true pt{\hfil}\cg%
\hbox to 20 true pt{\hfil}\cg%
\hbox to 20 true pt{\hfil}\cg%
\hbox to 20 true pt{\hfil}\cg%
\hbox to 20 true pt{\hfil}\cg%
\hbox to 20 true pt{\hss{$0$}}\cg%
\eol}\vss}\rg%
\tableclose%

\vskip 5 true pt plus 10 true pt
\goodbreak
\tableopen{Intertwining numbers $(Ind^{W}_{W_0}\varepsilon,Ind^{W}_{W_1}\rho_{1})_{\down{2 true pt}{W}}$ for $W = W(E_{8})$}%
\rowpts=19 true pt%
\rowlabpts=50 true pt%
\collabpts=50 true pt%
\clx{\vss\hfull{%
\rlx{\hss{$ $}}\cg%
\hbox to 20 true pt{\hskip 12 true pt\twist{$D_{8}$}\hss}\cg%
\hbox to 20 true pt{\hskip 12 true pt\twist{$A_{8}$}\hss}\cg%
\hbox to 20 true pt{\hskip 12 true pt\twist{${A_{7}}{A_{1}}$}\hss}\cg%
\hbox to 20 true pt{\hskip 12 true pt\twist{${A_{5}}{A_{2}}{A_{1}}$}\hss}\cg%
\hbox to 20 true pt{\hskip 12 true pt\twist{${A_{4}}{A_{4}}$}\hss}\cg%
\hbox to 20 true pt{\hskip 12 true pt\twist{${D_{5}}{A_{3}}$}\hss}\cg%
\hbox to 20 true pt{\hskip 12 true pt\twist{${E_{6}}{A_{2}}$}\hss}\cg%
\hbox to 20 true pt{\hskip 12 true pt\twist{${E_{7}}{A_{1}}$}\hss}\cg%
\hbox to 1.5\colpts{\hfil}\rlx{\hss{$ $}}\cg%
\hbox to 20 true pt{\hskip 12 true pt\twist{$D_{7}$}\hss}\cg%
\hbox to 20 true pt{\hskip 12 true pt\twist{$A_{7}$}\hss}\cg%
\hbox to 20 true pt{\hskip 12 true pt\twist{${A_{6}}{A_{1}}$}\hss}\cg%
\hbox to 21 true pt{\hskip 12 true pt\twist{${A_{4}}{A_{2}}{A_{1}}$}\hss}\cg%
\hbox to 21 true pt{\hskip 12 true pt\twist{${A_{4}}{A_{3}}$}\hss}\cg%
\hbox to 20 true pt{\hskip 12 true pt\twist{${D_{5}}{A_{2}}$}\hss}\cg%
\hbox to 20 true pt{\hskip 12 true pt\twist{${E_{6}}{A_{1}}$}\hss}\cg%
\hbox to 20 true pt{\hskip 12 true pt\twist{$E_{7}$}\hss}\cg%
\eol}}\rg%
\rx{\vss\hfull{%
\rlx{\hss{$D_{8}$}}\cg%
\hbox to 20 true pt{\hss{$0$}}\cg%
\hbox to 20 true pt{\hss{$0$}}\cg%
\hbox to 20 true pt{\hss{$0$}}\cg%
\hbox to 20 true pt{\hss{$0$}}\cg%
\hbox to 20 true pt{\hss{$0$}}\cg%
\hbox to 20 true pt{\hss{$0$}}\cg%
\hbox to 20 true pt{\hss{$0$}}\cg%
\hbox to 20 true pt{\hss{$0$}}\cg%
\hbox to 1.5\colpts{\hfil}\rlx{\hss{$D_{7}$}}\cg%
\hbox to 20 true pt{\hss{$0$}}\cg%
\hbox to 20 true pt{\hss{$0$}}\cg%
\hbox to 20 true pt{\hss{$0$}}\cg%
\hbox to 21 true pt{\hss{$0$}}\cg%
\hbox to 21 true pt{\hss{$0$}}\cg%
\hbox to 20 true pt{\hss{$0$}}\cg%
\hbox to 20 true pt{\hss{$0$}}\cg%
\hbox to 20 true pt{\hss{$0$}}\cg%
\eol}\vss}\rg%
\rx{\vss\hfull{%
\rlx{\hss{$A_{8}$}}\cg%
\hbox to 20 true pt{\hfil}\cg%
\hbox to 20 true pt{\hss{$0$}}\cg%
\hbox to 20 true pt{\hss{$0$}}\cg%
\hbox to 20 true pt{\hss{$0$}}\cg%
\hbox to 20 true pt{\hss{$0$}}\cg%
\hbox to 20 true pt{\hss{$0$}}\cg%
\hbox to 20 true pt{\hss{$0$}}\cg%
\hbox to 20 true pt{\hss{$0$}}\cg%
\hbox to 1.5\colpts{\hfil}\rlx{\hss{$A_{7}$}}\cg%
\hbox to 20 true pt{\hfil}\cg%
\hbox to 20 true pt{\hss{$0$}}\cg%
\hbox to 20 true pt{\hss{$0$}}\cg%
\hbox to 21 true pt{\hss{$2$}}\cg%
\hbox to 21 true pt{\hss{$0$}}\cg%
\hbox to 20 true pt{\hss{$0$}}\cg%
\hbox to 20 true pt{\hss{$0$}}\cg%
\hbox to 20 true pt{\hss{$0$}}\cg%
\eol}\vss}\rg%
\rx{\vss\hfull{%
\rlx{\hss{${A_{7}}{A_{1}}$}}\cg%
\hbox to 20 true pt{\hfil}\cg%
\hbox to 20 true pt{\hfil}\cg%
\hbox to 20 true pt{\hss{$0$}}\cg%
\hbox to 20 true pt{\hss{$0$}}\cg%
\hbox to 20 true pt{\hss{$0$}}\cg%
\hbox to 20 true pt{\hss{$0$}}\cg%
\hbox to 20 true pt{\hss{$0$}}\cg%
\hbox to 20 true pt{\hss{$0$}}\cg%
\hbox to 1.5\colpts{\hfil}\rlx{\hss{${A_{6}}{A_{1}}$}}\cg%
\hbox to 20 true pt{\hfil}\cg%
\hbox to 20 true pt{\hfil}\cg%
\hbox to 20 true pt{\hss{$0$}}\cg%
\hbox to 21 true pt{\hss{$25$}}\cg%
\hbox to 21 true pt{\hss{$6$}}\cg%
\hbox to 20 true pt{\hss{$0$}}\cg%
\hbox to 20 true pt{\hss{$0$}}\cg%
\hbox to 20 true pt{\hss{$0$}}\cg%
\eol}\vss}\rg%
\rx{\vss\hfull{%
\rlx{\hss{${A_{5}}{A_{2}}{A_{1}}$}}\cg%
\hbox to 20 true pt{\hfil}\cg%
\hbox to 20 true pt{\hfil}\cg%
\hbox to 20 true pt{\hfil}\cg%
\hbox to 20 true pt{\hss{$2$}}\cg%
\hbox to 20 true pt{\hss{$2$}}\cg%
\hbox to 20 true pt{\hss{$0$}}\cg%
\hbox to 20 true pt{\hss{$0$}}\cg%
\hbox to 20 true pt{\hss{$0$}}\cg%
\hbox to 1.5\colpts{\hfil}\rlx{\hss{${A_{4}}{A_{2}}{A_{1}}$}}\cg%
\hbox to 20 true pt{\hfil}\cg%
\hbox to 20 true pt{\hfil}\cg%
\hbox to 20 true pt{\hfil}\cg%
\hbox to 21 true pt{\hss{$587$}}\cg%
\hbox to 21 true pt{\hss{$234$}}\cg%
\hbox to 20 true pt{\hss{$17$}}\cg%
\hbox to 20 true pt{\hss{$0$}}\cg%
\hbox to 20 true pt{\hss{$0$}}\cg%
\eol}\vss}\rg%
\rx{\vss\hfull{%
\rlx{\hss{${A_{4}}{A_{4}}$}}\cg%
\hbox to 20 true pt{\hfil}\cg%
\hbox to 20 true pt{\hfil}\cg%
\hbox to 20 true pt{\hfil}\cg%
\hbox to 20 true pt{\hfil}\cg%
\hbox to 20 true pt{\hss{$0$}}\cg%
\hbox to 20 true pt{\hss{$0$}}\cg%
\hbox to 20 true pt{\hss{$0$}}\cg%
\hbox to 20 true pt{\hss{$0$}}\cg%
\hbox to 1.5\colpts{\hfil}\rlx{\hss{${A_{4}}{A_{3}}$}}\cg%
\hbox to 20 true pt{\hfil}\cg%
\hbox to 20 true pt{\hfil}\cg%
\hbox to 20 true pt{\hfil}\cg%
\hbox to 21 true pt{\hfil}\cg%
\hbox to 21 true pt{\hss{$89$}}\cg%
\hbox to 20 true pt{\hss{$5$}}\cg%
\hbox to 20 true pt{\hss{$0$}}\cg%
\hbox to 20 true pt{\hss{$0$}}\cg%
\eol}\vss}\rg%
\rx{\vss\hfull{%
\rlx{\hss{${D_{5}}{A_{3}}$}}\cg%
\hbox to 20 true pt{\hfil}\cg%
\hbox to 20 true pt{\hfil}\cg%
\hbox to 20 true pt{\hfil}\cg%
\hbox to 20 true pt{\hfil}\cg%
\hbox to 20 true pt{\hfil}\cg%
\hbox to 20 true pt{\hss{$0$}}\cg%
\hbox to 20 true pt{\hss{$0$}}\cg%
\hbox to 20 true pt{\hss{$0$}}\cg%
\hbox to 1.5\colpts{\hfil}\rlx{\hss{${D_{5}}{A_{2}}$}}\cg%
\hbox to 20 true pt{\hfil}\cg%
\hbox to 20 true pt{\hfil}\cg%
\hbox to 20 true pt{\hfil}\cg%
\hbox to 21 true pt{\hfil}\cg%
\hbox to 21 true pt{\hfil}\cg%
\hbox to 20 true pt{\hss{$0$}}\cg%
\hbox to 20 true pt{\hss{$0$}}\cg%
\hbox to 20 true pt{\hss{$0$}}\cg%
\eol}\vss}\rg%
\rx{\vss\hfull{%
\rlx{\hss{${E_{6}}{A_{2}}$}}\cg%
\hbox to 20 true pt{\hfil}\cg%
\hbox to 20 true pt{\hfil}\cg%
\hbox to 20 true pt{\hfil}\cg%
\hbox to 20 true pt{\hfil}\cg%
\hbox to 20 true pt{\hfil}\cg%
\hbox to 20 true pt{\hfil}\cg%
\hbox to 20 true pt{\hss{$0$}}\cg%
\hbox to 20 true pt{\hss{$0$}}\cg%
\hbox to 1.5\colpts{\hfil}\rlx{\hss{${E_{6}}{A_{1}}$}}\cg%
\hbox to 20 true pt{\hfil}\cg%
\hbox to 20 true pt{\hfil}\cg%
\hbox to 20 true pt{\hfil}\cg%
\hbox to 21 true pt{\hfil}\cg%
\hbox to 21 true pt{\hfil}\cg%
\hbox to 20 true pt{\hfil}\cg%
\hbox to 20 true pt{\hss{$0$}}\cg%
\hbox to 20 true pt{\hss{$0$}}\cg%
\eol}\vss}\rg%
\rx{\vss\hfull{%
\rlx{\hss{${E_{7}}{A_{1}}$}}\cg%
\hbox to 20 true pt{\hfil}\cg%
\hbox to 20 true pt{\hfil}\cg%
\hbox to 20 true pt{\hfil}\cg%
\hbox to 20 true pt{\hfil}\cg%
\hbox to 20 true pt{\hfil}\cg%
\hbox to 20 true pt{\hfil}\cg%
\hbox to 20 true pt{\hfil}\cg%
\hbox to 20 true pt{\hss{$0$}}\cg%
\hbox to 1.5\colpts{\hfil}\rlx{\hss{$E_{7}$}}\cg%
\hbox to 20 true pt{\hfil}\cg%
\hbox to 20 true pt{\hfil}\cg%
\hbox to 20 true pt{\hfil}\cg%
\hbox to 21 true pt{\hfil}\cg%
\hbox to 21 true pt{\hfil}\cg%
\hbox to 20 true pt{\hfil}\cg%
\hbox to 20 true pt{\hfil}\cg%
\hbox to 20 true pt{\hss{$0$}}\cg%
\eol}\vss}\rg%
\tableclose%

\vskip 5 true pt plus 10 true pt
\goodbreak
\tableopen{Intertwining numbers $(Ind^{W}_{W_0}\varepsilon,Ind^{W}_{W_1}\rho_{1})_{\down{2 true pt}{W}}$ for $W = W(F_{4})$}%
\rowpts=19 true pt%
\rowlabpts=105 true pt%
\collabpts=45 true pt%
\clx{\vss\hfull{%
\rlx{\hss{$ $}}\cg%
\hbox to 20 true pt{\hskip 12 true pt\twist{${C_{3}}{A_{1}}$}\hss}\cg%
\hbox to 20 true pt{\hskip 12 true pt\twist{${A_{2}}{A_{2}^{*}}$}\hss}\cg%
\hbox to 20 true pt{\hskip 12 true pt\twist{${A_{3}}{A_{1}}$}\hss}\cg%
\hbox to 20 true pt{\hskip 12 true pt\twist{$B_{4}$}\hss}\cg%
\hbox to 1.5\colpts{\hfil}\rlx{\hss{$ $}}\cg%
\hbox to 20 true pt{\hskip 12 true pt\twist{$C_{3}$}\hss}\cg%
\hbox to 20 true pt{\hskip 12 true pt\twist{${A_{1}}{A_{2}^{*}}$}\hss}\cg%
\hbox to 20 true pt{\hskip 12 true pt\twist{${A_{2}}{A_{1}^{*}}$}\hss}\cg%
\hbox to 20 true pt{\hskip 12 true pt\twist{$B_{3}$}\hss}\cg%
\eol}}\rg%
\rx{\vss\hfull{%
\rlx{\hss{${C_{3}}{A_{1}}$}}\cg%
\hbox to 20 true pt{\hss{$0$}}\cg%
\hbox to 20 true pt{\hss{$0$}}\cg%
\hbox to 20 true pt{\hss{$0$}}\cg%
\hbox to 20 true pt{\hss{$0$}}\cg%
\hbox to 1.5\colpts{\hfil}\rlx{\hss{$C_{3}$}}\cg%
\hbox to 20 true pt{\hss{$0$}}\cg%
\hbox to 20 true pt{\hss{$2$}}\cg%
\hbox to 20 true pt{\hss{$3$}}\cg%
\hbox to 20 true pt{\hss{$0$}}\cg%
\eol}\vss}\rg%
\rx{\vss\hfull{%
\rlx{\hss{${A_{2}}{A_{2}^{*}}$}}\cg%
\hbox to 20 true pt{\hfil}\cg%
\hbox to 20 true pt{\hss{$0$}}\cg%
\hbox to 20 true pt{\hss{$1$}}\cg%
\hbox to 20 true pt{\hss{$0$}}\cg%
\hbox to 1.5\colpts{\hfil}\rlx{\hss{${A_{1}}{A_{2}^{*}}$}}\cg%
\hbox to 20 true pt{\hfil}\cg%
\hbox to 20 true pt{\hss{$17$}}\cg%
\hbox to 20 true pt{\hss{$19$}}\cg%
\hbox to 20 true pt{\hss{$3$}}\cg%
\eol}\vss}\rg%
\rx{\vss\hfull{%
\rlx{\hss{${A_{3}}{A_{1}}$}}\cg%
\hbox to 20 true pt{\hfil}\cg%
\hbox to 20 true pt{\hfil}\cg%
\hbox to 20 true pt{\hss{$0$}}\cg%
\hbox to 20 true pt{\hss{$0$}}\cg%
\hbox to 1.5\colpts{\hfil}\rlx{\hss{${A_{2}}{A_{1}^{*}}$}}\cg%
\hbox to 20 true pt{\hfil}\cg%
\hbox to 20 true pt{\hfil}\cg%
\hbox to 20 true pt{\hss{$17$}}\cg%
\hbox to 20 true pt{\hss{$2$}}\cg%
\eol}\vss}\rg%
\rx{\vss\hfull{%
\rlx{\hss{$B_{4}$}}\cg%
\hbox to 20 true pt{\hfil}\cg%
\hbox to 20 true pt{\hfil}\cg%
\hbox to 20 true pt{\hfil}\cg%
\hbox to 20 true pt{\hss{$0$}}\cg%
\hbox to 1.5\colpts{\hfil}\rlx{\hss{$B_{3}$}}\cg%
\hbox to 20 true pt{\hfil}\cg%
\hbox to 20 true pt{\hfil}\cg%
\hbox to 20 true pt{\hfil}\cg%
\hbox to 20 true pt{\hss{$0$}}\cg%
\eol}\vss}\rg%
\tableclose%

\vskip 5 true pt plus 10 true pt
\goodbreak
\tableopen{Intertwining numbers $(Ind^{W}_{W_0}\varepsilon,Ind^{W}_{W_1}\rho_{1})_{\down{2 true pt}{W}}$ for $W = W(G_{2})$}%
\rowpts=19 true pt%
\rowlabpts=122 true pt%
\collabpts=40 true pt%
\clx{\vss\hfull{%
\rlx{\hss{$ $}}\cg%
\hbox to 20 true pt{\hskip 12 true pt\twist{${A_{1}}{A_{1}^{*}}$}\hss}\cg%
\hbox to 20 true pt{\hskip 12 true pt\twist{$A_{2}$}\hss}\cg%
\hbox to 1.5\colpts{\hfil}\rlx{\hss{$ $}}\cg%
\hbox to 20 true pt{\hskip 12 true pt\twist{$A_{1}$}\hss}\cg%
\hbox to 20 true pt{\hskip 12 true pt\twist{$A_{1}^{*}$}\hss}\cg%
\eol}}\rg%
\rx{\vss\hfull{%
\rlx{\hss{${A_{1}}{A_{1}^{*}}$}}\cg%
\hbox to 20 true pt{\hss{$0$}}\cg%
\hbox to 20 true pt{\hss{$1$}}\cg%
\hbox to 1.5\colpts{\hfil}\rlx{\hss{$A_{1}$}}\cg%
\hbox to 20 true pt{\hss{$4$}}\cg%
\hbox to 20 true pt{\hss{$3$}}\cg%
\eol}\vss}\rg%
\rx{\vss\hfull{%
\rlx{\hss{$A_{2}$}}\cg%
\hbox to 20 true pt{\hfil}\cg%
\hbox to 20 true pt{\hss{$0$}}\cg%
\hbox to 1.5\colpts{\hfil}\rlx{\hss{$A_{1}^{*}$}}\cg%
\hbox to 20 true pt{\hfil}\cg%
\hbox to 20 true pt{\hss{$4$}}\cg%
\eol}\vss}\rg%
\tableclose%

\vskip 5 true pt plus 10 true pt
\goodbreak
\tableopen{Intertwining numbers $(Ind^{W}_{W_0}\rho_{0},Ind^{W}_{W_1}\rho_{1})_{\down{2 true pt}{W}}$ for $W = W(E_{6})$}%
\rowpts=19 true pt%
\rowlabpts=120 true pt%
\collabpts=50 true pt%
\clx{\vss\hfull{%
\rlx{\hss{$ $}}\cg%
\hbox to 20 true pt{\hskip 12 true pt\twist{${A_{5}}{A_{1}}$}\hss}\cg%
\hbox to 21 true pt{\hskip 12 true pt\twist{${A_{2}}{A_{2}}{A_{2}}$}\hss}\cg%
\hbox to 1.5\colpts{\hfil}\rlx{\hss{$ $}}\cg%
\hbox to 20 true pt{\hskip 12 true pt\twist{$D_{5}$}\hss}\cg%
\hbox to 20 true pt{\hskip 12 true pt\twist{$A_{5}$}\hss}\cg%
\hbox to 20 true pt{\hskip 12 true pt\twist{${A_{4}}{A_{1}}$}\hss}\cg%
\hbox to 21 true pt{\hskip 12 true pt\twist{${A_{2}}{A_{2}}{A_{1}}$}\hss}\cg%
\eol}}\rg%
\rx{\vss\hfull{%
\rlx{\hss{${A_{5}}{A_{1}}$}}\cg%
\hbox to 20 true pt{\hss{$11$}}\cg%
\hbox to 21 true pt{\hss{$24$}}\cg%
\hbox to 1.5\colpts{\hfil}\rlx{\hss{$D_{5}$}}\cg%
\hbox to 20 true pt{\hss{$5$}}\cg%
\hbox to 20 true pt{\hss{$7$}}\cg%
\hbox to 20 true pt{\hss{$14$}}\cg%
\hbox to 21 true pt{\hss{$27$}}\cg%
\eol}\vss}\rg%
\rx{\vss\hfull{%
\rlx{\hss{${A_{2}}{A_{2}}{A_{2}}$}}\cg%
\hbox to 20 true pt{\hfil}\cg%
\hbox to 21 true pt{\hss{$108$}}\cg%
\hbox to 1.5\colpts{\hfil}\rlx{\hss{$A_{5}$}}\cg%
\hbox to 20 true pt{\hfil}\cg%
\hbox to 20 true pt{\hss{$13$}}\cg%
\hbox to 20 true pt{\hss{$26$}}\cg%
\hbox to 21 true pt{\hss{$57$}}\cg%
\eol}\vss}\rg%
\rx{\vss\hfull{%
\rlx{\hss{$ $}}\cg%
\hbox to 20 true pt{\hfil}\cg%
\hbox to 21 true pt{\hfil}\cg%
\hbox to 1.5\colpts{\hfil}\rlx{\hss{${A_{4}}{A_{1}}$}}\cg%
\hbox to 20 true pt{\hfil}\cg%
\hbox to 20 true pt{\hfil}\cg%
\hbox to 20 true pt{\hss{$60$}}\cg%
\hbox to 21 true pt{\hss{$143$}}\cg%
\eol}\vss}\rg%
\rx{\vss\hfull{%
\rlx{\hss{$ $}}\cg%
\hbox to 20 true pt{\hfil}\cg%
\hbox to 21 true pt{\hfil}\cg%
\hbox to 1.5\colpts{\hfil}\rlx{\hss{${A_{2}}{A_{2}}{A_{1}}$}}\cg%
\hbox to 20 true pt{\hfil}\cg%
\hbox to 20 true pt{\hfil}\cg%
\hbox to 20 true pt{\hfil}\cg%
\hbox to 21 true pt{\hss{$389$}}\cg%
\eol}\vss}\rg%
\tableclose%

\vskip 5 true pt plus 10 true pt
\goodbreak
\tableopen{Intertwining numbers $(Ind^{W}_{W_0}\rho_{0},Ind^{W}_{W_1}\rho_{1})_{\down{2 true pt}{W}}$ for $W = W(E_{7})$}%
\rowpts=19 true pt%
\rowlabpts=85 true pt%
\collabpts=50 true pt%
\clx{\vss\hfull{%
\rlx{\hss{$ $}}\cg%
\hbox to 20 true pt{\hskip 12 true pt\twist{${D_{6}}{A_{1}}$}\hss}\cg%
\hbox to 20 true pt{\hskip 12 true pt\twist{$A_{7}$}\hss}\cg%
\hbox to 20 true pt{\hskip 12 true pt\twist{${A_{5}}{A_{2}}$}\hss}\cg%
\hbox to 21 true pt{\hskip 12 true pt\twist{${A_{3}}{A_{3}}{A_{1}}$}\hss}\cg%
\hbox to 1.5\colpts{\hfil}\rlx{\hss{$ $}}\cg%
\hbox to 20 true pt{\hskip 12 true pt\twist{$D_{6}$}\hss}\cg%
\hbox to 20 true pt{\hskip 12 true pt\twist{$A_{6}$}\hss}\cg%
\hbox to 21 true pt{\hskip 12 true pt\twist{${A_{5}}{A_{1}}$}\hss}\cg%
\hbox to 28 true pt{\hskip 16 true pt\twist{${A_{3}}{A_{2}}{A_{1}}$}\hss}\cg%
\hbox to 28 true pt{\hskip 16 true pt\twist{${A_{4}}{A_{2}}$}\hss}\cg%
\hbox to 21 true pt{\hskip 12 true pt\twist{${D_{5}}{A_{1}}$}\hss}\cg%
\hbox to 20 true pt{\hskip 12 true pt\twist{$E_{6}$}\hss}\cg%
\eol}}\rg%
\rx{\vss\hfull{%
\rlx{\hss{${D_{6}}{A_{1}}$}}\cg%
\hbox to 20 true pt{\hss{$11$}}\cg%
\hbox to 20 true pt{\hss{$5$}}\cg%
\hbox to 20 true pt{\hss{$14$}}\cg%
\hbox to 21 true pt{\hss{$31$}}\cg%
\hbox to 1.5\colpts{\hfil}\rlx{\hss{$D_{6}$}}\cg%
\hbox to 20 true pt{\hss{$9$}}\cg%
\hbox to 20 true pt{\hss{$15$}}\cg%
\hbox to 21 true pt{\hss{$32$}}\cg%
\hbox to 28 true pt{\hss{$74$}}\cg%
\hbox to 28 true pt{\hss{$43$}}\cg%
\hbox to 21 true pt{\hss{$22$}}\cg%
\hbox to 20 true pt{\hss{$6$}}\cg%
\eol}\vss}\rg%
\rx{\vss\hfull{%
\rlx{\hss{$A_{7}$}}\cg%
\hbox to 20 true pt{\hfil}\cg%
\hbox to 20 true pt{\hss{$6$}}\cg%
\hbox to 20 true pt{\hss{$13$}}\cg%
\hbox to 21 true pt{\hss{$27$}}\cg%
\hbox to 1.5\colpts{\hfil}\rlx{\hss{$A_{6}$}}\cg%
\hbox to 20 true pt{\hfil}\cg%
\hbox to 20 true pt{\hss{$36$}}\cg%
\hbox to 21 true pt{\hss{$82$}}\cg%
\hbox to 28 true pt{\hss{$230$}}\cg%
\hbox to 28 true pt{\hss{$124$}}\cg%
\hbox to 21 true pt{\hss{$47$}}\cg%
\hbox to 20 true pt{\hss{$10$}}\cg%
\eol}\vss}\rg%
\rx{\vss\hfull{%
\rlx{\hss{${A_{5}}{A_{2}}$}}\cg%
\hbox to 20 true pt{\hfil}\cg%
\hbox to 20 true pt{\hfil}\cg%
\hbox to 20 true pt{\hss{$59$}}\cg%
\hbox to 21 true pt{\hss{$131$}}\cg%
\hbox to 1.5\colpts{\hfil}\rlx{\hss{${A_{5}}{A_{1}}$}}\cg%
\hbox to 20 true pt{\hfil}\cg%
\hbox to 20 true pt{\hfil}\cg%
\hbox to 21 true pt{\hss{$206$}}\cg%
\hbox to 28 true pt{\hss{$638$}}\cg%
\hbox to 28 true pt{\hss{$324$}}\cg%
\hbox to 21 true pt{\hss{$110$}}\cg%
\hbox to 20 true pt{\hss{$20$}}\cg%
\eol}\vss}\rg%
\rx{\vss\hfull{%
\rlx{\hss{${A_{3}}{A_{3}}{A_{1}}$}}\cg%
\hbox to 20 true pt{\hfil}\cg%
\hbox to 20 true pt{\hfil}\cg%
\hbox to 20 true pt{\hfil}\cg%
\hbox to 21 true pt{\hss{$361$}}\cg%
\hbox to 1.5\colpts{\hfil}\rlx{\hss{${A_{3}}{A_{2}}{A_{1}}$}}\cg%
\hbox to 20 true pt{\hfil}\cg%
\hbox to 20 true pt{\hfil}\cg%
\hbox to 21 true pt{\hfil}\cg%
\hbox to 28 true pt{\hss{$2346$}}\cg%
\hbox to 28 true pt{\hss{$1098$}}\cg%
\hbox to 21 true pt{\hss{$302$}}\cg%
\hbox to 20 true pt{\hss{$41$}}\cg%
\eol}\vss}\rg%
\rx{\vss\hfull{%
\rlx{\hss{$ $}}\cg%
\hbox to 20 true pt{\hfil}\cg%
\hbox to 20 true pt{\hfil}\cg%
\hbox to 20 true pt{\hfil}\cg%
\hbox to 21 true pt{\hfil}\cg%
\hbox to 1.5\colpts{\hfil}\rlx{\hss{${A_{4}}{A_{2}}$}}\cg%
\hbox to 20 true pt{\hfil}\cg%
\hbox to 20 true pt{\hfil}\cg%
\hbox to 21 true pt{\hfil}\cg%
\hbox to 28 true pt{\hfil}\cg%
\hbox to 28 true pt{\hss{$538$}}\cg%
\hbox to 21 true pt{\hss{$162$}}\cg%
\hbox to 20 true pt{\hss{$26$}}\cg%
\eol}\vss}\rg%
\rx{\vss\hfull{%
\rlx{\hss{$ $}}\cg%
\hbox to 20 true pt{\hfil}\cg%
\hbox to 20 true pt{\hfil}\cg%
\hbox to 20 true pt{\hfil}\cg%
\hbox to 21 true pt{\hfil}\cg%
\hbox to 1.5\colpts{\hfil}\rlx{\hss{${D_{5}}{A_{1}}$}}\cg%
\hbox to 20 true pt{\hfil}\cg%
\hbox to 20 true pt{\hfil}\cg%
\hbox to 21 true pt{\hfil}\cg%
\hbox to 28 true pt{\hfil}\cg%
\hbox to 28 true pt{\hfil}\cg%
\hbox to 21 true pt{\hss{$69$}}\cg%
\hbox to 20 true pt{\hss{$15$}}\cg%
\eol}\vss}\rg%
\rx{\vss\hfull{%
\rlx{\hss{$ $}}\cg%
\hbox to 20 true pt{\hfil}\cg%
\hbox to 20 true pt{\hfil}\cg%
\hbox to 20 true pt{\hfil}\cg%
\hbox to 21 true pt{\hfil}\cg%
\hbox to 1.5\colpts{\hfil}\rlx{\hss{$E_{6}$}}\cg%
\hbox to 20 true pt{\hfil}\cg%
\hbox to 20 true pt{\hfil}\cg%
\hbox to 21 true pt{\hfil}\cg%
\hbox to 28 true pt{\hfil}\cg%
\hbox to 28 true pt{\hfil}\cg%
\hbox to 21 true pt{\hfil}\cg%
\hbox to 20 true pt{\hss{$6$}}\cg%
\eol}\vss}\rg%
\tableclose%

\vskip 5 true pt plus 10 true pt
\goodbreak
\tableopen{Intertwining numbers $(Ind^{W}_{W_0}\rho_{0},Ind^{W}_{W_1}\rho_{1})_{\down{2 true pt}{W}}$ for $W = W(E_{8})$}%
\rowpts=19 true pt%
\rowlabpts=50 true pt%
\collabpts=50 true pt%
\clx{\vss\hfull{%
\rlx{\hss{$ $}}\cg%
\hbox to 20 true pt{\hskip 12 true pt\twist{$D_{8}$}\hss}\cg%
\hbox to 20 true pt{\hskip 12 true pt\twist{$A_{8}$}\hss}\cg%
\hbox to 21 true pt{\hskip 12 true pt\twist{${A_{7}}{A_{1}}$}\hss}\cg%
\hbox to 28 true pt{\hskip 16 true pt\twist{${A_{5}}{A_{2}}{A_{1}}$}\hss}\cg%
\hbox to 28 true pt{\hskip 16 true pt\twist{${A_{4}}{A_{4}}$}\hss}\cg%
\hbox to 21 true pt{\hskip 12 true pt\twist{${D_{5}}{A_{3}}$}\hss}\cg%
\hbox to 21 true pt{\hskip 12 true pt\twist{${E_{6}}{A_{2}}$}\hss}\cg%
\hbox to 20 true pt{\hskip 12 true pt\twist{${E_{7}}{A_{1}}$}\hss}\cg%
\hbox to 1.5\colpts{\hfil}\rlx{\hss{$ $}}\cg%
\hbox to 20 true pt{\hskip 12 true pt\twist{$D_{7}$}\hss}\cg%
\hbox to 21 true pt{\hskip 12 true pt\twist{$A_{7}$}\hss}\cg%
\hbox to 28 true pt{\hskip 16 true pt\twist{${A_{6}}{A_{1}}$}\hss}\cg%
\hbox to 35 true pt{\hskip 21 true pt\twist{${A_{4}}{A_{2}}{A_{1}}$}\hss}\cg%
\hbox to 35 true pt{\hskip 21 true pt\twist{${A_{4}}{A_{3}}$}\hss}\cg%
\hbox to 28 true pt{\hskip 16 true pt\twist{${D_{5}}{A_{2}}$}\hss}\cg%
\hbox to 28 true pt{\hskip 16 true pt\twist{${E_{6}}{A_{1}}$}\hss}\cg%
\hbox to 21 true pt{\hskip 12 true pt\twist{$E_{7}$}\hss}\cg%
\eol}}\rg%
\rx{\vss\hfull{%
\rlx{\hss{$D_{8}$}}\cg%
\hbox to 20 true pt{\hss{$4$}}\cg%
\hbox to 20 true pt{\hss{$5$}}\cg%
\hbox to 21 true pt{\hss{$12$}}\cg%
\hbox to 28 true pt{\hss{$30$}}\cg%
\hbox to 28 true pt{\hss{$20$}}\cg%
\hbox to 21 true pt{\hss{$15$}}\cg%
\hbox to 21 true pt{\hss{$7$}}\cg%
\hbox to 20 true pt{\hss{$5$}}\cg%
\hbox to 1.5\colpts{\hfil}\rlx{\hss{$D_{7}$}}\cg%
\hbox to 20 true pt{\hss{$27$}}\cg%
\hbox to 21 true pt{\hss{$67$}}\cg%
\hbox to 28 true pt{\hss{$158$}}\cg%
\hbox to 35 true pt{\hss{$475$}}\cg%
\hbox to 35 true pt{\hss{$293$}}\cg%
\hbox to 28 true pt{\hss{$152$}}\cg%
\hbox to 28 true pt{\hss{$57$}}\cg%
\hbox to 21 true pt{\hss{$11$}}\cg%
\eol}\vss}\rg%
\rx{\vss\hfull{%
\rlx{\hss{$A_{8}$}}\cg%
\hbox to 20 true pt{\hfil}\cg%
\hbox to 20 true pt{\hss{$20$}}\cg%
\hbox to 21 true pt{\hss{$44$}}\cg%
\hbox to 28 true pt{\hss{$158$}}\cg%
\hbox to 28 true pt{\hss{$108$}}\cg%
\hbox to 21 true pt{\hss{$57$}}\cg%
\hbox to 21 true pt{\hss{$23$}}\cg%
\hbox to 20 true pt{\hss{$6$}}\cg%
\hbox to 1.5\colpts{\hfil}\rlx{\hss{$A_{7}$}}\cg%
\hbox to 20 true pt{\hfil}\cg%
\hbox to 21 true pt{\hss{$231$}}\cg%
\hbox to 28 true pt{\hss{$605$}}\cg%
\hbox to 35 true pt{\hss{$2261$}}\cg%
\hbox to 35 true pt{\hss{$1333$}}\cg%
\hbox to 28 true pt{\hss{$561$}}\cg%
\hbox to 28 true pt{\hss{$153$}}\cg%
\hbox to 21 true pt{\hss{$21$}}\cg%
\eol}\vss}\rg%
\rx{\vss\hfull{%
\rlx{\hss{${A_{7}}{A_{1}}$}}\cg%
\hbox to 20 true pt{\hfil}\cg%
\hbox to 20 true pt{\hfil}\cg%
\hbox to 21 true pt{\hss{$123$}}\cg%
\hbox to 28 true pt{\hss{$483$}}\cg%
\hbox to 28 true pt{\hss{$320$}}\cg%
\hbox to 21 true pt{\hss{$163$}}\cg%
\hbox to 21 true pt{\hss{$58$}}\cg%
\hbox to 20 true pt{\hss{$15$}}\cg%
\hbox to 1.5\colpts{\hfil}\rlx{\hss{${A_{6}}{A_{1}}$}}\cg%
\hbox to 20 true pt{\hfil}\cg%
\hbox to 21 true pt{\hfil}\cg%
\hbox to 28 true pt{\hss{$1702$}}\cg%
\hbox to 35 true pt{\hss{$7041$}}\cg%
\hbox to 35 true pt{\hss{$4030$}}\cg%
\hbox to 28 true pt{\hss{$1573$}}\cg%
\hbox to 28 true pt{\hss{$381$}}\cg%
\hbox to 21 true pt{\hss{$44$}}\cg%
\eol}\vss}\rg%
\rx{\vss\hfull{%
\rlx{\hss{${A_{5}}{A_{2}}{A_{1}}$}}\cg%
\hbox to 20 true pt{\hfil}\cg%
\hbox to 20 true pt{\hfil}\cg%
\hbox to 21 true pt{\hfil}\cg%
\hbox to 28 true pt{\hss{$2482$}}\cg%
\hbox to 28 true pt{\hss{$1596$}}\cg%
\hbox to 21 true pt{\hss{$711$}}\cg%
\hbox to 21 true pt{\hss{$210$}}\cg%
\hbox to 20 true pt{\hss{$37$}}\cg%
\hbox to 1.5\colpts{\hfil}\rlx{\hss{${A_{4}}{A_{2}}{A_{1}}$}}\cg%
\hbox to 20 true pt{\hfil}\cg%
\hbox to 21 true pt{\hfil}\cg%
\hbox to 28 true pt{\hfil}\cg%
\hbox to 35 true pt{\hss{$34447$}}\cg%
\hbox to 35 true pt{\hss{$18948$}}\cg%
\hbox to 28 true pt{\hss{$6404$}}\cg%
\hbox to 28 true pt{\hss{$1221$}}\cg%
\hbox to 21 true pt{\hss{$102$}}\cg%
\eol}\vss}\rg%
\rx{\vss\hfull{%
\rlx{\hss{${A_{4}}{A_{4}}$}}\cg%
\hbox to 20 true pt{\hfil}\cg%
\hbox to 20 true pt{\hfil}\cg%
\hbox to 21 true pt{\hfil}\cg%
\hbox to 28 true pt{\hfil}\cg%
\hbox to 28 true pt{\hss{$1056$}}\cg%
\hbox to 21 true pt{\hss{$466$}}\cg%
\hbox to 21 true pt{\hss{$134$}}\cg%
\hbox to 20 true pt{\hss{$22$}}\cg%
\hbox to 1.5\colpts{\hfil}\rlx{\hss{${A_{4}}{A_{3}}$}}\cg%
\hbox to 20 true pt{\hfil}\cg%
\hbox to 21 true pt{\hfil}\cg%
\hbox to 28 true pt{\hfil}\cg%
\hbox to 35 true pt{\hfil}\cg%
\hbox to 35 true pt{\hss{$10559$}}\cg%
\hbox to 28 true pt{\hss{$3668$}}\cg%
\hbox to 28 true pt{\hss{$731$}}\cg%
\hbox to 21 true pt{\hss{$66$}}\cg%
\eol}\vss}\rg%
\rx{\vss\hfull{%
\rlx{\hss{${D_{5}}{A_{3}}$}}\cg%
\hbox to 20 true pt{\hfil}\cg%
\hbox to 20 true pt{\hfil}\cg%
\hbox to 21 true pt{\hfil}\cg%
\hbox to 28 true pt{\hfil}\cg%
\hbox to 28 true pt{\hfil}\cg%
\hbox to 21 true pt{\hss{$239$}}\cg%
\hbox to 21 true pt{\hss{$78$}}\cg%
\hbox to 20 true pt{\hss{$20$}}\cg%
\hbox to 1.5\colpts{\hfil}\rlx{\hss{${D_{5}}{A_{2}}$}}\cg%
\hbox to 20 true pt{\hfil}\cg%
\hbox to 21 true pt{\hfil}\cg%
\hbox to 28 true pt{\hfil}\cg%
\hbox to 35 true pt{\hfil}\cg%
\hbox to 35 true pt{\hfil}\cg%
\hbox to 28 true pt{\hss{$1483$}}\cg%
\hbox to 28 true pt{\hss{$377$}}\cg%
\hbox to 21 true pt{\hss{$45$}}\cg%
\eol}\vss}\rg%
\rx{\vss\hfull{%
\rlx{\hss{${E_{6}}{A_{2}}$}}\cg%
\hbox to 20 true pt{\hfil}\cg%
\hbox to 20 true pt{\hfil}\cg%
\hbox to 21 true pt{\hfil}\cg%
\hbox to 28 true pt{\hfil}\cg%
\hbox to 28 true pt{\hfil}\cg%
\hbox to 21 true pt{\hfil}\cg%
\hbox to 21 true pt{\hss{$43$}}\cg%
\hbox to 20 true pt{\hss{$12$}}\cg%
\hbox to 1.5\colpts{\hfil}\rlx{\hss{${E_{6}}{A_{1}}$}}\cg%
\hbox to 20 true pt{\hfil}\cg%
\hbox to 21 true pt{\hfil}\cg%
\hbox to 28 true pt{\hfil}\cg%
\hbox to 35 true pt{\hfil}\cg%
\hbox to 35 true pt{\hfil}\cg%
\hbox to 28 true pt{\hfil}\cg%
\hbox to 28 true pt{\hss{$141$}}\cg%
\hbox to 21 true pt{\hss{$25$}}\cg%
\eol}\vss}\rg%
\rx{\vss\hfull{%
\rlx{\hss{${E_{7}}{A_{1}}$}}\cg%
\hbox to 20 true pt{\hfil}\cg%
\hbox to 20 true pt{\hfil}\cg%
\hbox to 21 true pt{\hfil}\cg%
\hbox to 28 true pt{\hfil}\cg%
\hbox to 28 true pt{\hfil}\cg%
\hbox to 21 true pt{\hfil}\cg%
\hbox to 21 true pt{\hfil}\cg%
\hbox to 20 true pt{\hss{$10$}}\cg%
\hbox to 1.5\colpts{\hfil}\rlx{\hss{$E_{7}$}}\cg%
\hbox to 20 true pt{\hfil}\cg%
\hbox to 21 true pt{\hfil}\cg%
\hbox to 28 true pt{\hfil}\cg%
\hbox to 35 true pt{\hfil}\cg%
\hbox to 35 true pt{\hfil}\cg%
\hbox to 28 true pt{\hfil}\cg%
\hbox to 28 true pt{\hfil}\cg%
\hbox to 21 true pt{\hss{$8$}}\cg%
\eol}\vss}\rg%
\tableclose%

\vskip 5 true pt plus 10 true pt
\goodbreak
\tableopen{Intertwining numbers $(Ind^{W}_{W_0}\rho_{0},Ind^{W}_{W_1}\rho_{1})_{\down{2 true pt}{W}}$ for $W = W(F_{4})$}%
\rowpts=19 true pt%
\rowlabpts=105 true pt%
\collabpts=45 true pt%
\clx{\vss\hfull{%
\rlx{\hss{$ $}}\cg%
\hbox to 20 true pt{\hskip 12 true pt\twist{${C_{3}}{A_{1}}$}\hss}\cg%
\hbox to 20 true pt{\hskip 12 true pt\twist{${A_{2}}{A_{2}^{*}}$}\hss}\cg%
\hbox to 20 true pt{\hskip 12 true pt\twist{${A_{3}}{A_{1}}$}\hss}\cg%
\hbox to 20 true pt{\hskip 12 true pt\twist{$B_{4}$}\hss}\cg%
\hbox to 1.5\colpts{\hfil}\rlx{\hss{$ $}}\cg%
\hbox to 20 true pt{\hskip 12 true pt\twist{$C_{3}$}\hss}\cg%
\hbox to 20 true pt{\hskip 12 true pt\twist{${A_{1}}{A_{2}^{*}}$}\hss}\cg%
\hbox to 20 true pt{\hskip 12 true pt\twist{${A_{2}}{A_{1}^{*}}$}\hss}\cg%
\hbox to 20 true pt{\hskip 12 true pt\twist{$B_{3}$}\hss}\cg%
\eol}}\rg%
\rx{\vss\hfull{%
\rlx{\hss{${C_{3}}{A_{1}}$}}\cg%
\hbox to 20 true pt{\hss{$10$}}\cg%
\hbox to 20 true pt{\hss{$9$}}\cg%
\hbox to 20 true pt{\hss{$9$}}\cg%
\hbox to 20 true pt{\hss{$3$}}\cg%
\hbox to 1.5\colpts{\hfil}\rlx{\hss{$C_{3}$}}\cg%
\hbox to 20 true pt{\hss{$8$}}\cg%
\hbox to 20 true pt{\hss{$22$}}\cg%
\hbox to 20 true pt{\hss{$21$}}\cg%
\hbox to 20 true pt{\hss{$7$}}\cg%
\eol}\vss}\rg%
\rx{\vss\hfull{%
\rlx{\hss{${A_{2}}{A_{2}^{*}}$}}\cg%
\hbox to 20 true pt{\hfil}\cg%
\hbox to 20 true pt{\hss{$24$}}\cg%
\hbox to 20 true pt{\hss{$16$}}\cg%
\hbox to 20 true pt{\hss{$3$}}\cg%
\hbox to 1.5\colpts{\hfil}\rlx{\hss{${A_{1}}{A_{2}^{*}}$}}\cg%
\hbox to 20 true pt{\hfil}\cg%
\hbox to 20 true pt{\hss{$81$}}\cg%
\hbox to 20 true pt{\hss{$77$}}\cg%
\hbox to 20 true pt{\hss{$21$}}\cg%
\eol}\vss}\rg%
\rx{\vss\hfull{%
\rlx{\hss{${A_{3}}{A_{1}}$}}\cg%
\hbox to 20 true pt{\hfil}\cg%
\hbox to 20 true pt{\hfil}\cg%
\hbox to 20 true pt{\hss{$19$}}\cg%
\hbox to 20 true pt{\hss{$4$}}\cg%
\hbox to 1.5\colpts{\hfil}\rlx{\hss{${A_{2}}{A_{1}^{*}}$}}\cg%
\hbox to 20 true pt{\hfil}\cg%
\hbox to 20 true pt{\hfil}\cg%
\hbox to 20 true pt{\hss{$81$}}\cg%
\hbox to 20 true pt{\hss{$22$}}\cg%
\eol}\vss}\rg%
\rx{\vss\hfull{%
\rlx{\hss{$B_{4}$}}\cg%
\hbox to 20 true pt{\hfil}\cg%
\hbox to 20 true pt{\hfil}\cg%
\hbox to 20 true pt{\hfil}\cg%
\hbox to 20 true pt{\hss{$2$}}\cg%
\hbox to 1.5\colpts{\hfil}\rlx{\hss{$B_{3}$}}\cg%
\hbox to 20 true pt{\hfil}\cg%
\hbox to 20 true pt{\hfil}\cg%
\hbox to 20 true pt{\hfil}\cg%
\hbox to 20 true pt{\hss{$8$}}\cg%
\eol}\vss}\rg%
\tableclose%

\vskip 5 true pt plus 10 true pt
\goodbreak
\tableopen{Intertwining numbers $(Ind^{W}_{W_0}\rho_{0},Ind^{W}_{W_1}\rho_{1})_{\down{2 true pt}{W}}$ for $W = W(G_{2})$}%
\rowpts=19 true pt%
\rowlabpts=122 true pt%
\collabpts=40 true pt%
\clx{\vss\hfull{%
\rlx{\hss{$ $}}\cg%
\hbox to 20 true pt{\hskip 12 true pt\twist{${A_{1}}{A_{1}^{*}}$}\hss}\cg%
\hbox to 20 true pt{\hskip 12 true pt\twist{$A_{2}$}\hss}\cg%
\hbox to 1.5\colpts{\hfil}\rlx{\hss{$ $}}\cg%
\hbox to 20 true pt{\hskip 12 true pt\twist{$A_{1}$}\hss}\cg%
\hbox to 20 true pt{\hskip 12 true pt\twist{$A_{1}^{*}$}\hss}\cg%
\eol}}\rg%
\rx{\vss\hfull{%
\rlx{\hss{${A_{1}}{A_{1}^{*}}$}}\cg%
\hbox to 20 true pt{\hss{$6$}}\cg%
\hbox to 20 true pt{\hss{$2$}}\cg%
\hbox to 1.5\colpts{\hfil}\rlx{\hss{$A_{1}$}}\cg%
\hbox to 20 true pt{\hss{$4$}}\cg%
\hbox to 20 true pt{\hss{$3$}}\cg%
\eol}\vss}\rg%
\rx{\vss\hfull{%
\rlx{\hss{$A_{2}$}}\cg%
\hbox to 20 true pt{\hfil}\cg%
\hbox to 20 true pt{\hss{$2$}}\cg%
\hbox to 1.5\colpts{\hfil}\rlx{\hss{$A_{1}^{*}$}}\cg%
\hbox to 20 true pt{\hfil}\cg%
\hbox to 20 true pt{\hss{$4$}}\cg%
\eol}\vss}\rg%
\tableclose%
\endgroup
\vfil
\eject
%
%

\vfil
\eject

\otherhead{References}

\bigskip
\frenchspacing

\item{[A]} D.~Alvis, ``Duality and generic degress'',
{\sl Proc.~Amer.~Math.~Soc.} {\bf 91} (1984), 532-536.

\item{[AL]} D.~Alvis and G.~Lusztig, ``On Springer's
correspondence for simple groups of type
$E_n$, ($n = 6$, $7$, $8$)'',
{\sl Math.~Proc.~Camb.~Phil.~Soc.}
{\bf 92} (1982), 65-72.

\item{[BL]} W.~M.~Beynon and G.~Lusztig, ``Some numerical
results on the characters of exceptional Weyl groups'',
{\sl Math.~Proc.~Camb.~Phil.~Soc.} {\bf84} (1978), 417--426.

\item{[B]} N.~Bourbaki,
{\sl Groupes et alg\`ebres de Lie},
Chap. 4, 5, 6,
Paris, Hermann, 1968.

\item{[C]} R.~Carter,
``Conjugacy classes in the Weyl group'',
in {\sl Seminar on Algebraic Groups and Related
Finite Groups}, Lecture Notes in Mathematics, vol. 131,
Springer-Verlag, Berlin, 1970.

\item{[CR]} C.~W.~Curtis and I.~Reiner,
{\sl Methods of Representation Theory,
with Applications to Finite Groups
and Orders, vol. I},
Wiley, New York, 1981.

\item{[F1]} J.~S.~Frame, ``The classes and
representations of the groups of 27 lines
and 28 bitangents'',
{\sl Ann.~Mat.~Pura.~Appl. \bf{4}} (1951), 83--119.

\item{[F2]}
J.~S.~Frame,
``The characters of the Weyl group
$E_8$'',
in
{\sl Computational problems in Abstract Algebra},
ed. J.~Leech,
Pergamon Press,
Oxford, 1970.

\item{[JK]} G.~James and A.~Kerber,
{\sl The Representation Theory of the
Symmetric Group},
vol. 16 of
{\sl Encyclopedia of Mathematics and Its
Applications},
Addison-Wesley,
London, 1981.

\item{[K]} T.~Kondo, ``The characters of the Weyl group of
type $F_4$'',
{\sl J.~Fac.~Sci. Univ.~Tokyo}
Sect. I {\bf 11}
(1965),
145--153.

\item{[S]} N.~Spaltenstein,
``Appendix'', {\sl Math.~Proc.~Camb.~Phil.~Soc.},
{\bf 92} (1982), 65-72.

\item{[Y]} A.~Young,
``On quantitative substitutional analysis (fifth paper)'',
{\sl Proc.~London~Math.~Soc.}
(2) {\bf 31} (1930), 273-288.

\vfil
\eject
%
%
\vfil
\eject
\immediate\write\toctxt{TheEnd}
\immediate\closeout\toctxt
\immediate\closeout\tocnum
\pageno=-1
\leftline{\sectionfont Table of Contents}

\bigskip
\baselineskip=\tocskip

\newread\toctxt
\openin\toctxt toctxt
\newread\tocnum
\openin\tocnum tocnum
\newwrite\toc
\immediate\openout\toc toc.tex
\newif\ifmore
\moretrue
\def\eof{TheEnd }
\loop{%
  \read\toctxt to \txt%
  \ifx\txt\eof%
    \global\morefalse
  \else{%
    \read\tocnum to \num%
    \immediate\write\toc{%
      \line{\txt\lotsofdots\num}%
      }%
    }%
  \fi
  }%
\ifmore{\relax}
\repeat

\closein\toctxt
\closein\tocnum
\immediate\closeout\toc

\hbox to\hsize {{\fam \bffam \tenbf 1.\kern .5em Introduction} \leaders \hbox to 0.7em{\hss .\hss }\hfil 1 }
\hbox to\hsize {{\fam \bffam \tenbf 2.\kern .5em Conjugacy Classes} \leaders \hbox to 0.7em{\hss .\hss }\hfil 3 }
\hbox to\hsize {{\fam \bffam \tenbf 3.\kern .5em Class Inclusions} \leaders \hbox to 0.7em{\hss .\hss }\hfil 7 }
\hbox to\hsize {{\fam \bffam \tenbf \hskip 2em\relax Table 1.\kern .5em }{\fam \slfam \tensl Class inclusions for $W({A_{5}}{A_{1}})\mskip \thinmuskip \subset \mskip \thinmuskip W(E_{6})$} \leaders \hbox to 0.7em{\hss .\hss }\hfil 8 }
\hbox to\hsize {{\fam \bffam \tenbf \hskip 2em\relax Table 2.\kern .5em }{\fam \slfam \tensl Class inclusions for $W({A_{2}}{A_{2}}{A_{2}})\mskip \thinmuskip \subset \mskip \thinmuskip W(E_{6})$} \leaders \hbox to 0.7em{\hss .\hss }\hfil 8 }
\hbox to\hsize {{\fam \bffam \tenbf \hskip 2em\relax Table 3.\kern .5em }{\fam \slfam \tensl Class inclusions for $W(D_{5})\mskip \thinmuskip \subset \mskip \thinmuskip W(E_{6})$} \leaders \hbox to 0.7em{\hss .\hss }\hfil 9 }
\hbox to\hsize {{\fam \bffam \tenbf \hskip 2em\relax Table 4.\kern .5em }{\fam \slfam \tensl Class inclusions for $W({D_{6}}{A_{1}})\mskip \thinmuskip \subset \mskip \thinmuskip W(E_{7})$} \leaders \hbox to 0.7em{\hss .\hss }\hfil 10 }
\hbox to\hsize {{\fam \bffam \tenbf \hskip 2em\relax Table 5.\kern .5em }{\fam \slfam \tensl Class inclusions for $W(A_{7})\mskip \thinmuskip \subset \mskip \thinmuskip W(E_{7})$} \leaders \hbox to 0.7em{\hss .\hss }\hfil 11 }
\hbox to\hsize {{\fam \bffam \tenbf \hskip 2em\relax Table 6.\kern .5em }{\fam \slfam \tensl Class inclusions for $W({A_{5}}{A_{2}})\mskip \thinmuskip \subset \mskip \thinmuskip W(E_{7})$} \leaders \hbox to 0.7em{\hss .\hss }\hfil 11 }
\hbox to\hsize {{\fam \bffam \tenbf \hskip 2em\relax Table 7.\kern .5em }{\fam \slfam \tensl Class inclusions for $W({A_{3}}{A_{3}}{A_{1}})\mskip \thinmuskip \subset \mskip \thinmuskip W(E_{7})$} \leaders \hbox to 0.7em{\hss .\hss }\hfil 12 }
\hbox to\hsize {{\fam \bffam \tenbf \hskip 2em\relax Table 8.\kern .5em }{\fam \slfam \tensl Class inclusions for $W(E_{6})\mskip \thinmuskip \subset \mskip \thinmuskip W(E_{7})$} \leaders \hbox to 0.7em{\hss .\hss }\hfil 12 }
\hbox to\hsize {{\fam \bffam \tenbf \hskip 2em\relax Table 9.\kern .5em }{\fam \slfam \tensl Class inclusions for $W(D_{8})\mskip \thinmuskip \subset \mskip \thinmuskip W(E_{8})$} \leaders \hbox to 0.7em{\hss .\hss }\hfil 14 }
\hbox to\hsize {{\fam \bffam \tenbf \hskip 2em\relax Table 10.\kern .5em }{\fam \slfam \tensl Class inclusions for $W(A_{8})\mskip \thinmuskip \subset \mskip \thinmuskip W(E_{8})$} \leaders \hbox to 0.7em{\hss .\hss }\hfil 15 }
\hbox to\hsize {{\fam \bffam \tenbf \hskip 2em\relax Table 11.\kern .5em }{\fam \slfam \tensl Class inclusions for $W({A_{7}}{A_{1}})\mskip \thinmuskip \subset \mskip \thinmuskip W(E_{8})$} \leaders \hbox to 0.7em{\hss .\hss }\hfil 15 }
\hbox to\hsize {{\fam \bffam \tenbf \hskip 2em\relax Table 12.\kern .5em }{\fam \slfam \tensl Class inclusions for $W({A_{5}}{A_{2}}{A_{1}})\mskip \thinmuskip \subset \mskip \thinmuskip W(E_{8})$} \leaders \hbox to 0.7em{\hss .\hss }\hfil 16 }
\hbox to\hsize {{\fam \bffam \tenbf \hskip 2em\relax Table 13.\kern .5em }{\fam \slfam \tensl Class inclusions for $W({A_{4}}{A_{4}})\mskip \thinmuskip \subset \mskip \thinmuskip W(E_{8})$} \leaders \hbox to 0.7em{\hss .\hss }\hfil 17 }
\hbox to\hsize {{\fam \bffam \tenbf \hskip 2em\relax Table 14.\kern .5em }{\fam \slfam \tensl Class inclusions for $W({D_{5}}{A_{3}})\mskip \thinmuskip \subset \mskip \thinmuskip W(E_{8})$} \leaders \hbox to 0.7em{\hss .\hss }\hfil 18 }
\hbox to\hsize {{\fam \bffam \tenbf \hskip 2em\relax Table 15.\kern .5em }{\fam \slfam \tensl Class inclusions for $W({E_{6}}{A_{2}})\mskip \thinmuskip \subset \mskip \thinmuskip W(E_{8})$} \leaders \hbox to 0.7em{\hss .\hss }\hfil 19 }
\hbox to\hsize {{\fam \bffam \tenbf \hskip 2em\relax Table 16.\kern .5em }{\fam \slfam \tensl Class inclusions for $W({E_{7}}{A_{1}})\mskip \thinmuskip \subset \mskip \thinmuskip W(E_{8})$} \leaders \hbox to 0.7em{\hss .\hss }\hfil 20 }
\hbox to\hsize {{\fam \bffam \tenbf \hskip 2em\relax Table 17.\kern .5em }{\fam \slfam \tensl Class inclusions for $W({C_{3}}{A_{1}})\mskip \thinmuskip \subset \mskip \thinmuskip W(F_{4})$} \leaders \hbox to 0.7em{\hss .\hss }\hfil 22 }
\hbox to\hsize {{\fam \bffam \tenbf \hskip 2em\relax Table 18.\kern .5em }{\fam \slfam \tensl Class inclusions for $W({A_{2}}{A_{2}^{*}})\mskip \thinmuskip \subset \mskip \thinmuskip W(F_{4})$} \leaders \hbox to 0.7em{\hss .\hss }\hfil 22 }
\hbox to\hsize {{\fam \bffam \tenbf \hskip 2em\relax Table 19.\kern .5em }{\fam \slfam \tensl Class inclusions for $W({A_{3}}{A_{1}})\mskip \thinmuskip \subset \mskip \thinmuskip W(F_{4})$} \leaders \hbox to 0.7em{\hss .\hss }\hfil 23 }
\hbox to\hsize {{\fam \bffam \tenbf \hskip 2em\relax Table 20.\kern .5em }{\fam \slfam \tensl Class inclusions for $W(B_{4})\mskip \thinmuskip \subset \mskip \thinmuskip W(F_{4})$} \leaders \hbox to 0.7em{\hss .\hss }\hfil 23 }
\hbox to\hsize {{\fam \bffam \tenbf \hskip 2em\relax Table 21.\kern .5em }{\fam \slfam \tensl Class inclusions for $W({A_{1}}{A_{1}^{*}})\mskip \thinmuskip \subset \mskip \thinmuskip W(G_{2})$} \leaders \hbox to 0.7em{\hss .\hss }\hfil 23 }
\hbox to\hsize {{\fam \bffam \tenbf \hskip 2em\relax Table 22.\kern .5em }{\fam \slfam \tensl Class inclusions for $W(A_{2})\mskip \thinmuskip \subset \mskip \thinmuskip W(G_{2})$} \leaders \hbox to 0.7em{\hss .\hss }\hfil 23 }
\hbox to\hsize {{\fam \bffam \tenbf 4.\kern .5em Irreducible Characters} \leaders \hbox to 0.7em{\hss .\hss }\hfil 24 }
\hbox to\hsize {{\fam \bffam \tenbf 5.\kern .5em Induce/Restrict Matrices} \leaders \hbox to 0.7em{\hss .\hss }\hfil 27 }
\hbox to\hsize {{\fam \bffam \tenbf \hskip 2em\relax Table 23.\kern .5em }{\fam \slfam \tensl Induce/restrict matrix for $W({A_{5}}{A_{1}})\mskip \thinmuskip \subset \mskip \thinmuskip W(E_{6})$} \leaders \hbox to 0.7em{\hss .\hss }\hfil 29 }
\hbox to\hsize {{\fam \bffam \tenbf \hskip 2em\relax Table 24.\kern .5em }{\fam \slfam \tensl Induce/restrict matrix for $W({A_{2}}{A_{2}}{A_{2}})\mskip \thinmuskip \subset \mskip \thinmuskip W(E_{6})$} \leaders \hbox to 0.7em{\hss .\hss }\hfil 29 }
\hbox to\hsize {{\fam \bffam \tenbf \hskip 2em\relax Table 25.\kern .5em }{\fam \slfam \tensl Induce/restrict matrix for $W(D_{5})\mskip \thinmuskip \subset \mskip \thinmuskip W(E_{6})$} \leaders \hbox to 0.7em{\hss .\hss }\hfil 31 }
\hbox to\hsize {{\fam \bffam \tenbf \hskip 2em\relax Table 26.\kern .5em }{\fam \slfam \tensl Induce/restrict matrix for $W(A_{5})\mskip \thinmuskip \subset \mskip \thinmuskip W(E_{6})$} \leaders \hbox to 0.7em{\hss .\hss }\hfil 31 }
\hbox to\hsize {{\fam \bffam \tenbf \hskip 2em\relax Table 27.\kern .5em }{\fam \slfam \tensl Induce/restrict matrix for $W({A_{4}}{A_{1}})\mskip \thinmuskip \subset \mskip \thinmuskip W(E_{6})$} \leaders \hbox to 0.7em{\hss .\hss }\hfil 32 }
\hbox to\hsize {{\fam \bffam \tenbf \hskip 2em\relax Table 28.\kern .5em }{\fam \slfam \tensl Induce/restrict matrix for $W({A_{2}}{A_{2}}{A_{1}})\mskip \thinmuskip \subset \mskip \thinmuskip W(E_{6})$} \leaders \hbox to 0.7em{\hss .\hss }\hfil 32 }
\hbox to\hsize {{\fam \bffam \tenbf \hskip 2em\relax Table 29.\kern .5em }{\fam \slfam \tensl Induce/restrict matrix for $W({D_{6}}{A_{1}})\mskip \thinmuskip \subset \mskip \thinmuskip W(E_{7})$} \leaders \hbox to 0.7em{\hss .\hss }\hfil 33 }
\hbox to\hsize {{\fam \bffam \tenbf \hskip 2em\relax Table 30.\kern .5em }{\fam \slfam \tensl Induce/restrict matrix for $W(A_{7})\mskip \thinmuskip \subset \mskip \thinmuskip W(E_{7})$} \leaders \hbox to 0.7em{\hss .\hss }\hfil 36 }
\hbox to\hsize {{\fam \bffam \tenbf \hskip 2em\relax Table 31.\kern .5em }{\fam \slfam \tensl Induce/restrict matrix for $W({A_{5}}{A_{2}})\mskip \thinmuskip \subset \mskip \thinmuskip W(E_{7})$} \leaders \hbox to 0.7em{\hss .\hss }\hfil 38 }
\hbox to\hsize {{\fam \bffam \tenbf \hskip 2em\relax Table 32.\kern .5em }{\fam \slfam \tensl Induce/restrict matrix for $W({A_{3}}{A_{3}}{A_{1}})\mskip \thinmuskip \subset \mskip \thinmuskip W(E_{7})$} \leaders \hbox to 0.7em{\hss .\hss }\hfil 40 }
\hbox to\hsize {{\fam \bffam \tenbf \hskip 2em\relax Table 33.\kern .5em }{\fam \slfam \tensl Induce/restrict matrix for $W(D_{6})\mskip \thinmuskip \subset \mskip \thinmuskip W(E_{7})$} \leaders \hbox to 0.7em{\hss .\hss }\hfil 42 }
\hbox to\hsize {{\fam \bffam \tenbf \hskip 2em\relax Table 34.\kern .5em }{\fam \slfam \tensl Induce/restrict matrix for $W(A_{6})\mskip \thinmuskip \subset \mskip \thinmuskip W(E_{7})$} \leaders \hbox to 0.7em{\hss .\hss }\hfil 44 }
\hbox to\hsize {{\fam \bffam \tenbf \hskip 2em\relax Table 35.\kern .5em }{\fam \slfam \tensl Induce/restrict matrix for $W({A_{5}}{A_{1}})\mskip \thinmuskip \subset \mskip \thinmuskip W(E_{7})$} \leaders \hbox to 0.7em{\hss .\hss }\hfil 45 }
\hbox to\hsize {{\fam \bffam \tenbf \hskip 2em\relax Table 36.\kern .5em }{\fam \slfam \tensl Induce/restrict matrix for $W({A_{3}}{A_{2}}{A_{1}})\mskip \thinmuskip \subset \mskip \thinmuskip W(E_{7})$} \leaders \hbox to 0.7em{\hss .\hss }\hfil 46 }
\hbox to\hsize {{\fam \bffam \tenbf \hskip 2em\relax Table 37.\kern .5em }{\fam \slfam \tensl Induce/restrict matrix for $W({A_{4}}{A_{2}})\mskip \thinmuskip \subset \mskip \thinmuskip W(E_{7})$} \leaders \hbox to 0.7em{\hss .\hss }\hfil 48 }
\hbox to\hsize {{\fam \bffam \tenbf \hskip 2em\relax Table 38.\kern .5em }{\fam \slfam \tensl Induce/restrict matrix for $W({D_{5}}{A_{1}})\mskip \thinmuskip \subset \mskip \thinmuskip W(E_{7})$} \leaders \hbox to 0.7em{\hss .\hss }\hfil 49 }
\hbox to\hsize {{\fam \bffam \tenbf \hskip 2em\relax Table 39.\kern .5em }{\fam \slfam \tensl Induce/restrict matrix for $W(E_{6})\mskip \thinmuskip \subset \mskip \thinmuskip W(E_{7})$} \leaders \hbox to 0.7em{\hss .\hss }\hfil 51 }
\hbox to\hsize {{\fam \bffam \tenbf \hskip 2em\relax Table 40.\kern .5em }{\fam \slfam \tensl Induce/restrict matrix for $W(D_{8})\mskip \thinmuskip \subset \mskip \thinmuskip W(E_{8})$} \leaders \hbox to 0.7em{\hss .\hss }\hfil 52 }
\hbox to\hsize {{\fam \bffam \tenbf \hskip 2em\relax Table 41.\kern .5em }{\fam \slfam \tensl Induce/restrict matrix for $W(A_{8})\mskip \thinmuskip \subset \mskip \thinmuskip W(E_{8})$} \leaders \hbox to 0.7em{\hss .\hss }\hfil 60 }
\hbox to\hsize {{\fam \bffam \tenbf \hskip 2em\relax Table 42.\kern .5em }{\fam \slfam \tensl Induce/restrict matrix for $W({A_{7}}{A_{1}})\mskip \thinmuskip \subset \mskip \thinmuskip W(E_{8})$} \leaders \hbox to 0.7em{\hss .\hss }\hfil 64 }
\hbox to\hsize {{\fam \bffam \tenbf \hskip 2em\relax Table 43.\kern .5em }{\fam \slfam \tensl Induce/restrict matrix for $W({A_{5}}{A_{2}}{A_{1}})\mskip \thinmuskip \subset \mskip \thinmuskip W(E_{8})$} \leaders \hbox to 0.7em{\hss .\hss }\hfil 68 }
\hbox to\hsize {{\fam \bffam \tenbf \hskip 2em\relax Table 44.\kern .5em }{\fam \slfam \tensl Induce/restrict matrix for $W({A_{4}}{A_{4}})\mskip \thinmuskip \subset \mskip \thinmuskip W(E_{8})$} \leaders \hbox to 0.7em{\hss .\hss }\hfil 75 }
\hbox to\hsize {{\fam \bffam \tenbf \hskip 2em\relax Table 45.\kern .5em }{\fam \slfam \tensl Induce/restrict matrix for $W({D_{5}}{A_{3}})\mskip \thinmuskip \subset \mskip \thinmuskip W(E_{8})$} \leaders \hbox to 0.7em{\hss .\hss }\hfil 79 }
\hbox to\hsize {{\fam \bffam \tenbf \hskip 2em\relax Table 46.\kern .5em }{\fam \slfam \tensl Induce/restrict matrix for $W({E_{6}}{A_{2}})\mskip \thinmuskip \subset \mskip \thinmuskip W(E_{8})$} \leaders \hbox to 0.7em{\hss .\hss }\hfil 88 }
\hbox to\hsize {{\fam \bffam \tenbf \hskip 2em\relax Table 47.\kern .5em }{\fam \slfam \tensl Induce/restrict matrix for $W({E_{7}}{A_{1}})\mskip \thinmuskip \subset \mskip \thinmuskip W(E_{8})$} \leaders \hbox to 0.7em{\hss .\hss }\hfil 94 }
\hbox to\hsize {{\fam \bffam \tenbf \hskip 2em\relax Table 48.\kern .5em }{\fam \slfam \tensl Induce/restrict matrix for $W(D_{7})\mskip \thinmuskip \subset \mskip \thinmuskip W(E_{8})$} \leaders \hbox to 0.7em{\hss .\hss }\hfil 105 }
\hbox to\hsize {{\fam \bffam \tenbf \hskip 2em\relax Table 49.\kern .5em }{\fam \slfam \tensl Induce/restrict matrix for $W(A_{7})\mskip \thinmuskip \subset \mskip \thinmuskip W(E_{8})$} \leaders \hbox to 0.7em{\hss .\hss }\hfil 111 }
\hbox to\hsize {{\fam \bffam \tenbf \hskip 2em\relax Table 50.\kern .5em }{\fam \slfam \tensl Induce/restrict matrix for $W({A_{6}}{A_{1}})\mskip \thinmuskip \subset \mskip \thinmuskip W(E_{8})$} \leaders \hbox to 0.7em{\hss .\hss }\hfil 113 }
\hbox to\hsize {{\fam \bffam \tenbf \hskip 2em\relax Table 51.\kern .5em }{\fam \slfam \tensl Induce/restrict matrix for $W({A_{4}}{A_{2}}{A_{1}})\mskip \thinmuskip \subset \mskip \thinmuskip W(E_{8})$} \leaders \hbox to 0.7em{\hss .\hss }\hfil 117 }
\hbox to\hsize {{\fam \bffam \tenbf \hskip 2em\relax Table 52.\kern .5em }{\fam \slfam \tensl Induce/restrict matrix for $W({A_{4}}{A_{3}})\mskip \thinmuskip \subset \mskip \thinmuskip W(E_{8})$} \leaders \hbox to 0.7em{\hss .\hss }\hfil 121 }
\hbox to\hsize {{\fam \bffam \tenbf \hskip 2em\relax Table 53.\kern .5em }{\fam \slfam \tensl Induce/restrict matrix for $W({D_{5}}{A_{2}})\mskip \thinmuskip \subset \mskip \thinmuskip W(E_{8})$} \leaders \hbox to 0.7em{\hss .\hss }\hfil 126 }
\hbox to\hsize {{\fam \bffam \tenbf \hskip 2em\relax Table 54.\kern .5em }{\fam \slfam \tensl Induce/restrict matrix for $W({E_{6}}{A_{1}})\mskip \thinmuskip \subset \mskip \thinmuskip W(E_{8})$} \leaders \hbox to 0.7em{\hss .\hss }\hfil 133 }
\hbox to\hsize {{\fam \bffam \tenbf \hskip 2em\relax Table 55.\kern .5em }{\fam \slfam \tensl Induce/restrict matrix for $W(E_{7})\mskip \thinmuskip \subset \mskip \thinmuskip W(E_{8})$} \leaders \hbox to 0.7em{\hss .\hss }\hfil 137 }
\hbox to\hsize {{\fam \bffam \tenbf \hskip 2em\relax Table 56.\kern .5em }{\fam \slfam \tensl Induce/restrict matrix for $W({C_{3}}{A_{1}})\mskip \thinmuskip \subset \mskip \thinmuskip W(F_{4})$} \leaders \hbox to 0.7em{\hss .\hss }\hfil 143 }
\hbox to\hsize {{\fam \bffam \tenbf \hskip 2em\relax Table 57.\kern .5em }{\fam \slfam \tensl Induce/restrict matrix for $W({A_{2}}{A_{2}^{*}})\mskip \thinmuskip \subset \mskip \thinmuskip W(F_{4})$} \leaders \hbox to 0.7em{\hss .\hss }\hfil 144 }
\hbox to\hsize {{\fam \bffam \tenbf \hskip 2em\relax Table 58.\kern .5em }{\fam \slfam \tensl Induce/restrict matrix for $W({A_{3}}{A_{1}})\mskip \thinmuskip \subset \mskip \thinmuskip W(F_{4})$} \leaders \hbox to 0.7em{\hss .\hss }\hfil 145 }
\hbox to\hsize {{\fam \bffam \tenbf \hskip 2em\relax Table 59.\kern .5em }{\fam \slfam \tensl Induce/restrict matrix for $W(B_{4})\mskip \thinmuskip \subset \mskip \thinmuskip W(F_{4})$} \leaders \hbox to 0.7em{\hss .\hss }\hfil 146 }
\hbox to\hsize {{\fam \bffam \tenbf \hskip 2em\relax Table 60.\kern .5em }{\fam \slfam \tensl Induce/restrict matrix for $W(B_{3})\mskip \thinmuskip \subset \mskip \thinmuskip W(F_{4})$} \leaders \hbox to 0.7em{\hss .\hss }\hfil 147 }
\hbox to\hsize {{\fam \bffam \tenbf \hskip 2em\relax Table 61.\kern .5em }{\fam \slfam \tensl Induce/restrict matrix for $W({A_{2}}{A_{1}^{*}})\mskip \thinmuskip \subset \mskip \thinmuskip W(F_{4})$} \leaders \hbox to 0.7em{\hss .\hss }\hfil 148 }
\hbox to\hsize {{\fam \bffam \tenbf \hskip 2em\relax Table 62.\kern .5em }{\fam \slfam \tensl Induce/restrict matrix for $W({A_{1}}{A_{1}^{*}})\mskip \thinmuskip \subset \mskip \thinmuskip W(G_{2})$} \leaders \hbox to 0.7em{\hss .\hss }\hfil 149 }
\hbox to\hsize {{\fam \bffam \tenbf \hskip 2em\relax Table 63.\kern .5em }{\fam \slfam \tensl Induce/restrict matrix for $W(A_{2})\mskip \thinmuskip \subset \mskip \thinmuskip W(G_{2})$} \leaders \hbox to 0.7em{\hss .\hss }\hfil 149 }
\hbox to\hsize {{\fam \bffam \tenbf \hskip 2em\relax Table 64.\kern .5em }{\fam \slfam \tensl Induce/restrict matrix for $W(A_{1})\mskip \thinmuskip \subset \mskip \thinmuskip W(G_{2})$} \leaders \hbox to 0.7em{\hss .\hss }\hfil 149 }
\hbox to\hsize {{\fam \bffam \tenbf 6.\kern .5em Intertwining Numbers} \leaders \hbox to 0.7em{\hss .\hss }\hfil 150 }
\hbox to\hsize {{\fam \bffam \tenbf \hskip 2em\relax Table 65.\kern .5em }{\fam \slfam \tensl Intertwining numbers $(Ind^{W}_{W_0}1,Ind^{W}_{W_1}1)_{\lower 2 true pt\hbox {$\scriptstyle {W}$}}$ for $W = W(E_{6})$} \leaders \hbox to 0.7em{\hss .\hss }\hfil 150 }
\hbox to\hsize {{\fam \bffam \tenbf \hskip 2em\relax Table 66.\kern .5em }{\fam \slfam \tensl Intertwining numbers $(Ind^{W}_{W_0}1,Ind^{W}_{W_1}1)_{\lower 2 true pt\hbox {$\scriptstyle {W}$}}$ for $W = W(E_{7})$} \leaders \hbox to 0.7em{\hss .\hss }\hfil 151 }
\hbox to\hsize {{\fam \bffam \tenbf \hskip 2em\relax Table 67.\kern .5em }{\fam \slfam \tensl Intertwining numbers $(Ind^{W}_{W_0}1,Ind^{W}_{W_1}1)_{\lower 2 true pt\hbox {$\scriptstyle {W}$}}$ for $W = W(E_{8})$} \leaders \hbox to 0.7em{\hss .\hss }\hfil 151 }
\hbox to\hsize {{\fam \bffam \tenbf \hskip 2em\relax Table 68.\kern .5em }{\fam \slfam \tensl Intertwining numbers $(Ind^{W}_{W_0}1,Ind^{W}_{W_1}1)_{\lower 2 true pt\hbox {$\scriptstyle {W}$}}$ for $W = W(F_{4})$} \leaders \hbox to 0.7em{\hss .\hss }\hfil 151 }
\hbox to\hsize {{\fam \bffam \tenbf \hskip 2em\relax Table 69.\kern .5em }{\fam \slfam \tensl Intertwining numbers $(Ind^{W}_{W_0}1,Ind^{W}_{W_1}1)_{\lower 2 true pt\hbox {$\scriptstyle {W}$}}$ for $W = W(G_{2})$} \leaders \hbox to 0.7em{\hss .\hss }\hfil 152 }
\hbox to\hsize {{\fam \bffam \tenbf \hskip 2em\relax Table 70.\kern .5em }{\fam \slfam \tensl Intertwining numbers $(Ind^{W}_{W_0}1,Ind^{W}_{W_1}\varepsilon )_{\lower 2 true pt\hbox {$\scriptstyle {W}$}}$ for $W = W(E_{6})$} \leaders \hbox to 0.7em{\hss .\hss }\hfil 152 }
\hbox to\hsize {{\fam \bffam \tenbf \hskip 2em\relax Table 71.\kern .5em }{\fam \slfam \tensl Intertwining numbers $(Ind^{W}_{W_0}1,Ind^{W}_{W_1}\varepsilon )_{\lower 2 true pt\hbox {$\scriptstyle {W}$}}$ for $W = W(E_{7})$} \leaders \hbox to 0.7em{\hss .\hss }\hfil 152 }
\hbox to\hsize {{\fam \bffam \tenbf \hskip 2em\relax Table 72.\kern .5em }{\fam \slfam \tensl Intertwining numbers $(Ind^{W}_{W_0}1,Ind^{W}_{W_1}\varepsilon )_{\lower 2 true pt\hbox {$\scriptstyle {W}$}}$ for $W = W(E_{8})$} \leaders \hbox to 0.7em{\hss .\hss }\hfil 153 }
\hbox to\hsize {{\fam \bffam \tenbf \hskip 2em\relax Table 73.\kern .5em }{\fam \slfam \tensl Intertwining numbers $(Ind^{W}_{W_0}1,Ind^{W}_{W_1}\varepsilon )_{\lower 2 true pt\hbox {$\scriptstyle {W}$}}$ for $W = W(F_{4})$} \leaders \hbox to 0.7em{\hss .\hss }\hfil 153 }
\hbox to\hsize {{\fam \bffam \tenbf \hskip 2em\relax Table 74.\kern .5em }{\fam \slfam \tensl Intertwining numbers $(Ind^{W}_{W_0}1,Ind^{W}_{W_1}\varepsilon )_{\lower 2 true pt\hbox {$\scriptstyle {W}$}}$ for $W = W(G_{2})$} \leaders \hbox to 0.7em{\hss .\hss }\hfil 153 }
\hbox to\hsize {{\fam \bffam \tenbf \hskip 2em\relax Table 75.\kern .5em }{\fam \slfam \tensl Intertwining numbers $(Ind^{W}_{W_0}1,Ind^{W}_{W_1}\rho _{1})_{\lower 2 true pt\hbox {$\scriptstyle {W}$}}$ for $W = W(E_{6})$} \leaders \hbox to 0.7em{\hss .\hss }\hfil 153 }
\hbox to\hsize {{\fam \bffam \tenbf \hskip 2em\relax Table 76.\kern .5em }{\fam \slfam \tensl Intertwining numbers $(Ind^{W}_{W_0}1,Ind^{W}_{W_1}\rho _{1})_{\lower 2 true pt\hbox {$\scriptstyle {W}$}}$ for $W = W(E_{7})$} \leaders \hbox to 0.7em{\hss .\hss }\hfil 154 }
\hbox to\hsize {{\fam \bffam \tenbf \hskip 2em\relax Table 77.\kern .5em }{\fam \slfam \tensl Intertwining numbers $(Ind^{W}_{W_0}1,Ind^{W}_{W_1}\rho _{1})_{\lower 2 true pt\hbox {$\scriptstyle {W}$}}$ for $W = W(E_{8})$} \leaders \hbox to 0.7em{\hss .\hss }\hfil 154 }
\hbox to\hsize {{\fam \bffam \tenbf \hskip 2em\relax Table 78.\kern .5em }{\fam \slfam \tensl Intertwining numbers $(Ind^{W}_{W_0}1,Ind^{W}_{W_1}\rho _{1})_{\lower 2 true pt\hbox {$\scriptstyle {W}$}}$ for $W = W(F_{4})$} \leaders \hbox to 0.7em{\hss .\hss }\hfil 154 }
\hbox to\hsize {{\fam \bffam \tenbf \hskip 2em\relax Table 79.\kern .5em }{\fam \slfam \tensl Intertwining numbers $(Ind^{W}_{W_0}1,Ind^{W}_{W_1}\rho _{1})_{\lower 2 true pt\hbox {$\scriptstyle {W}$}}$ for $W = W(G_{2})$} \leaders \hbox to 0.7em{\hss .\hss }\hfil 155 }
\hbox to\hsize {{\fam \bffam \tenbf \hskip 2em\relax Table 80.\kern .5em }{\fam \slfam \tensl Intertwining numbers $(Ind^{W}_{W_0}\varepsilon ,Ind^{W}_{W_1}\rho _{1})_{\lower 2 true pt\hbox {$\scriptstyle {W}$}}$ for $W = W(E_{6})$} \leaders \hbox to 0.7em{\hss .\hss }\hfil 155 }
\hbox to\hsize {{\fam \bffam \tenbf \hskip 2em\relax Table 81.\kern .5em }{\fam \slfam \tensl Intertwining numbers $(Ind^{W}_{W_0}\varepsilon ,Ind^{W}_{W_1}\rho _{1})_{\lower 2 true pt\hbox {$\scriptstyle {W}$}}$ for $W = W(E_{7})$} \leaders \hbox to 0.7em{\hss .\hss }\hfil 155 }
\hbox to\hsize {{\fam \bffam \tenbf \hskip 2em\relax Table 82.\kern .5em }{\fam \slfam \tensl Intertwining numbers $(Ind^{W}_{W_0}\varepsilon ,Ind^{W}_{W_1}\rho _{1})_{\lower 2 true pt\hbox {$\scriptstyle {W}$}}$ for $W = W(E_{8})$} \leaders \hbox to 0.7em{\hss .\hss }\hfil 156 }
\hbox to\hsize {{\fam \bffam \tenbf \hskip 2em\relax Table 83.\kern .5em }{\fam \slfam \tensl Intertwining numbers $(Ind^{W}_{W_0}\varepsilon ,Ind^{W}_{W_1}\rho _{1})_{\lower 2 true pt\hbox {$\scriptstyle {W}$}}$ for $W = W(F_{4})$} \leaders \hbox to 0.7em{\hss .\hss }\hfil 156 }
\hbox to\hsize {{\fam \bffam \tenbf \hskip 2em\relax Table 84.\kern .5em }{\fam \slfam \tensl Intertwining numbers $(Ind^{W}_{W_0}\varepsilon ,Ind^{W}_{W_1}\rho _{1})_{\lower 2 true pt\hbox {$\scriptstyle {W}$}}$ for $W = W(G_{2})$} \leaders \hbox to 0.7em{\hss .\hss }\hfil 156 }
\hbox to\hsize {{\fam \bffam \tenbf \hskip 2em\relax Table 85.\kern .5em }{\fam \slfam \tensl Intertwining numbers $(Ind^{W}_{W_0}\rho _{0},Ind^{W}_{W_1}\rho _{1})_{\lower 2 true pt\hbox {$\scriptstyle {W}$}}$ for $W = W(E_{6})$} \leaders \hbox to 0.7em{\hss .\hss }\hfil 156 }
\hbox to\hsize {{\fam \bffam \tenbf \hskip 2em\relax Table 86.\kern .5em }{\fam \slfam \tensl Intertwining numbers $(Ind^{W}_{W_0}\rho _{0},Ind^{W}_{W_1}\rho _{1})_{\lower 2 true pt\hbox {$\scriptstyle {W}$}}$ for $W = W(E_{7})$} \leaders \hbox to 0.7em{\hss .\hss }\hfil 157 }
\hbox to\hsize {{\fam \bffam \tenbf \hskip 2em\relax Table 87.\kern .5em }{\fam \slfam \tensl Intertwining numbers $(Ind^{W}_{W_0}\rho _{0},Ind^{W}_{W_1}\rho _{1})_{\lower 2 true pt\hbox {$\scriptstyle {W}$}}$ for $W = W(E_{8})$} \leaders \hbox to 0.7em{\hss .\hss }\hfil 157 }
\hbox to\hsize {{\fam \bffam \tenbf \hskip 2em\relax Table 88.\kern .5em }{\fam \slfam \tensl Intertwining numbers $(Ind^{W}_{W_0}\rho _{0},Ind^{W}_{W_1}\rho _{1})_{\lower 2 true pt\hbox {$\scriptstyle {W}$}}$ for $W = W(F_{4})$} \leaders \hbox to 0.7em{\hss .\hss }\hfil 157 }
\hbox to\hsize {{\fam \bffam \tenbf \hskip 2em\relax Table 89.\kern .5em }{\fam \slfam \tensl Intertwining numbers $(Ind^{W}_{W_0}\rho _{0},Ind^{W}_{W_1}\rho _{1})_{\lower 2 true pt\hbox {$\scriptstyle {W}$}}$ for $W = W(G_{2})$} \leaders \hbox to 0.7em{\hss .\hss }\hfil 157 }
\hbox to\hsize {{\fam \bffam \tenbf References} \leaders \hbox to 0.7em{\hss .\hss }\hfil 158 }

\vfil
\eject
\end